%% file: EqK2.tex
\documentclass{amsbook}

\usepackage{sty/preamble}

\title{\vspace{5\baselineskip}\scalebox{2.7}{Equivariant Algebraic}\protect\\ \scalebox{2.7}{$K$-Theories}}
\author{\vspace{2.5\baselineskip}\scalebox{2}{Donald Yau}}

\hypersetup{pdfauthor=\authors}

\date{\today}

\begin{document}

\newcommand{\pagecol}{white}
\pagecolor{\pagecol} 

\frontmatter


\maketitle

\cleardoublepage
\thispagestyle{empty}
\vspace*{13.5pc}
\begin{center}
The author dedicates this work to Jacqueline.  
\end{center}
\cleardoublepage

\pdfbookmark{\contentsname}{Contents}
\tableofcontents

\chapter*{Preface}
\input{chap/preface.tex}

\chapter*{Summary of Main Results}
\label{ch:intro}
\input{chap/intro}

\chapter*{Glossary of Notation}
\label{ch:notation}
\input{chap/notation.tex}

\mainmatter
\part{Multifunctorial Equivariant $K$-Theory}
\label{part:kgo_shi}

\chapter{Multifunctorial Equivariant $K$-Theory via $J$-Theory}
\label{ch:kgo}
\input{chap/kgo} 
\input{chap/kgo_ii} 
\input{chap/kgo_iii} 
\input{chap/kgo_iv} 
\input{chap/kgo_v} 
\input{chap/kgo_vi} 
\input{chap/kgo_vii} 

\chapter{$\GGG$-Categories}
\label{ch:ggcatg}
\input{chap/ggcatg}

\chapter{Equivariant $K$-Theory via $H$-Theory}
\label{ch:hgo}
\input{chap/hgo}

\part{Shimakawa Equivariant $K$-Theory}
\label{part:shim_k}

\chapter{Shimakawa $H$-Theory and $J$-Theory} 
\label{ch:shih}
\input{chap/shih}

\chapter{Shimakawa $K$-Theory}
\label{ch:shimakawa_K}
\input{chap/shimakawa_K}

\chapter{Topological Equivalence of Shimakawa and Multifunctorial $K$-Theories}
\label{ch:shim_top}
\input{chap/shim_top}

\part{Categorical Equivalence of Shimakawa and Multifunctorial $K$-Theories}
\label{part:kgo_shi_comp}

\chapter{Comparison of $H$-Theories}
\label{ch:h_comparison}
\input{chap/h_comparison}

\chapter{Shimakawa Strong $H$-Theory and Twisted Products}
\label{ch:sgoprod}
\input{chap/sgoprod}

\chapter{Strong $H$-Theory and Twisted Products}
\label{ch:hgoprod}
\input{chap/hgoprod}

\chapter{Strong $H$-Theory Comparison Weak $G$-Equivalence}
\label{ch:compgen}
\input{chap/compgen}

\chapter{Special Objects and Weak $G$-Equivalences}
\label{ch:special}
\input{chap/special}

\part{Equivalences of Shimakawa, GMMO, and Schwede Global $K$-Theories}
\label{part:shi_gmmo_gl}

\chapter{Equivalence of Shimakawa and GMMO $K$-Theories}
\label{ch:kgmmo}
\input{chap/kgmmo}

\input{chap/gmmo_shi}

\chapter{Schwede Global $K$-Theory}
\label{ch:kgl}
\input{chap/kgl}

\chapter{Equivalence of GMMO and Global $K$-Theories}
\label{ch:kgl_gmmo}
\input{chap/kgl_gmmo}

\appendix

\chapter{Categories and Operads}
\label{ch:mon_cat}
\input{chap/mon_cat.tex}


\backmatter
\bibliographystyle{sty/amsalpha3}
\bibliography{references}
\printindex
\end{document}

%% file: chap/preface.tex
A cornerstone of algebraic $K$-theory is the equivalence between the $K$-theory machines of May, Segal, and Elmendorf and Mandell.  Equivariant algebraic $K$-theory enriches the theory with group actions, making it more powerful and complex.  There are a number of equivariant $K$-theory machines that turn equivariant categorical data into equivariant spectra, the main objects of study in equivariant stable homotopy theory.  

This work proves that the following four equivariant $K$-theory machines are appropriately equivalent: Shimakawa equivariant $K$-theory \cite{shimakawa89,shimakawa91}; the author's enriched multifunctorial equivariant $K$-theory \cite{yau-eqk}; the equivariant $K$-theory of Guillou, May, Merling, and Osorno \cite{gmmo23}; and Schwede global equivariant $K$-theory \cite{schwede_global}.  \cref{part:kgo_shi,part:shim_k} prove the topological equivalence between Shimakawa and multifunctorial equivariant $K$-theories.  \cref{part:kgo_shi_comp} proves that their categorical parts are equivalent.  \cref{part:shi_gmmo_gl} proves that the equivariant $K$-theory of Guillou, May, Merling, and Osorno is equivalent to Shimakawa $K$-theory and Schwede global $K$-theory for each finite group.

%% file: chap/intro.tex
\subsection*{Objective}
For each finite group $G$, each of the following four equivariant $K$-theory machines turns equivariant categorical data into equivariant spectra, the main objects of study in equivariant stable homotopy theory:
\begin{enumerate}
\item\label{machine_shi} Shimakawa equivariant $K$-theory \cite{shimakawa89,shimakawa91};
\item\label{machine_yau} the author's enriched multifunctorial equivariant $K$-theory \cite{yau-eqk};
\item\label{machine_gmmo} the equivariant $K$-theory of Guillou, May, Merling, and Osorno \cite{gmmo23}, called GMMO $K$-theory; and
\item\label{machine_sch} Schwede global equivariant $K$-theory \cite{schwede_global}.
\end{enumerate}
This work proves that these four equivariant $K$-theory machines are appropriately equivalent.  Without assuming any prior knowledge of algebraic $K$-theory, the main text thoroughly explains these four equivariant $K$-theory machines and their comparisons.

\organization
The rest of this introductory chapter summarizes the main results of this work and provides references to the main text.  \cref{sec:intro_tables} is a nontechnical overview of the four equivariant $K$-theory machines and their comparisons.  \cref{sec:intro_shim} summarizes the equivalence between Shimakawa and multifunctorial equivariant $K$-theories.  \cref{sec:intro_shi_gmmo} summarizes the equivalence between Shimakawa and GMMO $K$-theories.  \cref{sec:intro_gmmo_gl} summarizes the equivalence between GMMO and Schwede global $K$-theories.  \cref{sec:intro_summ} summarizes each part and each chapter in the main text. 

\subsection*{Reading Suggestion}
This work is mostly self-contained and has a substantial amount of detailed explanation.  The references in \cref{sec:intro_tables} can be used to jump to a specific equivariant $K$-theory machine and comparison result of interest.  The extensive cross references throughout this work can be used to trace backward for necessary details.  After this introductory chapter, one way to skim this work is to read only the chapter introductions.  In the main text, each chapter introduction summarizes the objectives of that chapter, its connection with other chapters, and the content of each section.  In a similar manner, the beginning of each section has a summary and a list of key references for that section.  Relevant concepts of category theory are reviewed in \cref{ch:mon_cat}.

\section{An Overview of Equivariant Algebraic $K$-Theory}
\label{sec:intro_tables}

This section provides a bird's-eye view of (non)equivariant algebraic $K$-theory \pcref{comparison_schematic,eqmachines} and the comparisons of the four equivariant $K$-theory machines \pcref{table.eqmach,table.foureqmach}, along with references to subsequent sections and the main text.

\subsection*{Nonequivariant Algebraic $K$-Theory}
Algebraic $K$-theory machines---including the operadic machine of May \cite{may}, the Segal machine \cite{segal}, and the multifunctorial Elmendorf-Mandell machine \cite{elmendorf-mandell}---are functors that send small permutative categories to spectra, the main objects of study in stable homotopy theory.  The work \cite{may-thomason} of May and Thomason proves the topological equivalence between the May machine and the Segal machine.  The work \cite{elmendorf-mandell} of Elmendorf and Mandell proves the equivalence between their machine and Segal's.  Thus, several ways to construct spectra from either categorical or topological data are equivalent.

Each algebraic $K$-theory machine has its own relative advantages.  The Segal machine \cite{segal} is the easiest one to define.  Its categorical part, which sends small permutative categories to special $\Fsk$-categories, is reminiscent of lower $K$-groups.  However, the Segal machine does not generally preserve multiplicative structures.  With a slightly more elaborate construction than the Segal machine, the Elmendorf-Mandell machine \cite{elmendorf-mandell} is an enriched multifunctor.  It preserves all algebraic structures parametrized by operads, including $E_n$-algebras for $1\leq n \leq \infty$ and their modules.  The book \cite{cerberusIII} has a detailed exposition of the Segal machine, the Elmendorf-Mandell machine, and their equivalence.

Relative to the Segal machine and the Elmendorf-Mandell machine, May's operadic machine \cite{may} has a more categorical flavor and uses categorical bar constructions extensively.  It is a homotopical functor, so it sends levelwise weak homotopy equivalences to equivalences of spectra.  Earlier work \cite{may-einfinity,may-multiplicative} on multiplicative properties of the May machine requires further tuning, as discussed in \cite{may-construction,may-good,may-precisely}.  With the later adjustment, the May machine sends bipermutative categories to $\Einf$-ring spaces and then to $\Einf$-ring spectra.  It should be noted that $\Einf$-ring spaces and $\Einf$-ring spectra are not defined as algebras over some $\Einf$-operads.

\subsection*{Equivariant Algebraic $K$-Theory}
Equivariant algebraic $K$-theory enriches the theory with group actions, making it more powerful and complex.  The nonequivariant $K$-theory machines of May, Segal, and Elmendorf and Mandell generalize to the equivariant $K$-theory machines of, respectively, Guillou and May \cite{gm17}, Shimakawa \cite{shimakawa89,shimakawa91}, and the author \cite{yau-eqk}.  Starting from topological data, the work \cite{mmo} of May, Merling, and Osorno establishes the equivalence between the Guillou-May machine and the Shimakawa machine.  The equivariant $K$-theory of Guillou, May, Merling, and Osorno \cite{gmmo23}, which we call GMMO $K$-theory, preserves algebras over nonsymmetric operads, a property that Guillou-May and Shimakawa $K$-theories lack.  Additionally, Schwede's global equivariant $K$-theory machine \cite{schwede_global} keeps track of group actions for all finite groups.  

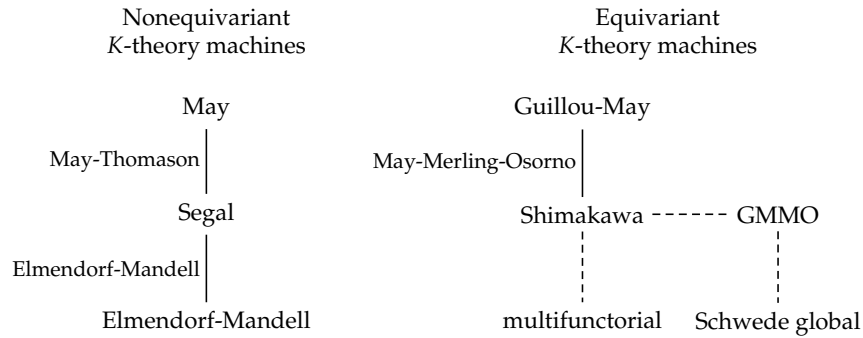
\begin{figure}[H] 
\centering
\begin{tikzpicture}
\def\u{-1} \def\v{-1.4} \def\w{-1.45} \def\a{-1} \def\h{2.6} \def\t{.2}
\draw[0cell=.9]
(0,0) node (a01) [align=center] {\text{Nonequivariant}\\ \text{$K$-theory machines}}
(a01)++(0,\u) node (a11) {\text{May}}
(a11)++(0,\v) node (a21) {\text{Segal}}
(a21)++(0,\v) node (a31) {\text{Elmendorf-Mandell}}
(a01)++(6,0) node (a02) [align=center] {\text{Equivariant}\\ \text{$K$-theory machines}}
(a02)++(\a,\u) node (a12) {\text{Guillou-May}}
(a12)++(0,\v) node (a22) {\text{Shimakawa}}
(a22)++(0,\t) node (a22') {\phantom{X}}
(a22)++(0,-\t) node (a22'') {\phantom{X}}
(a22)++(0,\v) node (a32) {\text{multifunctorial}}
(a22)++(\h,0) node (b1) {\text{GMMO}}
(b1)++(0,\w) node (b2) {\text{Schwede global}}
;
\draw[1cell=.9]
(a11) edge[-,semithick] node[swap] {\text{May-Thomason}} (a21)
(a21) edge[-,semithick] node[swap] {\text{Elmendorf-Mandell}} (a31)
(a12) edge[-,semithick] node[swap] {\text{May-Merling-Osorno}} (a22)
(a22) edge[-,densely dashed,semithick] (a32)
(a22) edge[-,densely dashed,semithick] (b1)
(b1) edge[-,densely dashed,semithick] (b2)
;
\end{tikzpicture}
\caption{Equivalences of (non)equivariant $K$-theory machines.}
\label{comparison_schematic}
\end{figure}

\Cref{comparison_schematic} summarizes the (non)equivariant machines and their equivalences mentioned so far.
\noindent
The three solid lines represent the equivalences established in \cite{elmendorf-mandell,may-thomason,mmo}.  The three dashed lines represent the equivalences of equivariant $K$-theory machines established in this work, as summarized in \crefrange{sec:intro_shim}{sec:intro_gmmo_gl}.  Given the May-Merling-Osorno equivalence between Guillou-May and Shimakawa $K$-theories, this work considers only the latter.

\subsection*{Summary of Equivariant Machines}
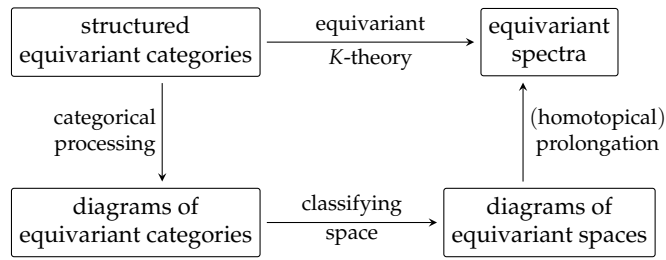
\begin{figure}[H] 
\centering
\begin{tikzpicture}[x=55mm, y=24mm,
block/.style ={rectangle, draw=black,
align=center, rounded corners=1pt, outer sep=1mm}]
\def\t{1em}
\draw[0cell=.9]
(0,0) node[block] (a11) {\text{structured}\\ \text{equivariant categories}}
(1,0) node[block] (a12) {\text{equivariant}\\ \text{spectra}}
(0,-1) node[block] (a21) {\text{diagrams of}\\ \text{equivariant categories}}
(1,-1) node[block] (a22) {\text{diagrams of}\\ \text{equivariant spaces}}
;
\draw[1cell=.9]
(a11) edge node {\text{equivariant}} node[swap] {\text{$K$-theory}} (a12)
(a11) edge[transform canvas={xshift=\t}] node[swap,align=center] {\text{categorical}\\ \text{processing}} (a21)
(a21) edge node{\text{classifying}} node[swap] {\text{space}} (a22)
(a22) edge[transform canvas={xshift=-\t}] node[swap,align=center] {\text{$\mathrm{(homotopical)}$}\\ \text{prolongation}} (a12);
\end{tikzpicture}
\caption{Steps of equivariant $K$-theory machines.}
\label{eqmachines}
\end{figure}

Each equivariant $K$-theory machine considered in this work sends equivariant categorical input data to equivariant spectral output data in three main steps, as displayed in \cref{eqmachines}.
\noindent
The first main step of each machine turns structured equivariant categories into diagrams of equivariant categories indexed by some indexing categories or equivariant categories.  The first step varies substantially among different equivariant $K$-theory machines.  It is a major source of difficulty in their comparison.  After applying the classifying space functor levelwise, the last step turns diagrams of equivariant spaces into equivariant spectra by prolongation or homotopical prolongation. 

This work proves that the four equivariant $K$-theory machines in \cref{table.eqmach} are appropriately equivalent for each finite group $G$.
\begin{figure}[H] 
\centering
\resizebox{\linewidth}{!}{%
{\renewcommand{\arraystretch}{1.3}%
{\setlength{\tabcolsep}{1ex}
\begin{tabular}{c|c|c|c|c|c|c}
\multirow{2}{*}{\makecell{equivariant \\ $K$-theories}} 
& \multirow{2}{*}{\makecell{main \\ definitions}} 
& \multirow{2}{*}{\makecell{explicit \\ description}} 
& \multirow{2}{*}{\makecell{equivariant categorical \\ input data}} 
& \multirow{2}{*}{\makecell{indexing \\ categories}} 
& \multirow{2}{*}{\makecell{equivariant spectral \\ output data}} 
& \multirow{2}{*}{\makecell{notable \\ properties}} 
\\ 
&&&&&& \\ \hline
\multirow{2}{*}{\makecell{Shimakawa \\ \cref{kshosg_intro}, \cref{khshosg_intro}}} 
& \multirow{2}{*}{\makecell{\cref{ksho_kshosg} \\ \cref{khsho_khshosg}}}
& \multirow{2}{*}{\makecell{\ref{expl:ksho} \\ \ref{expl:khsho}}}
& \multirow{2}{*}{\makecell{$\Op$-pseudoalgebras \eqref{def:pseudoalgebra} \\ 1-connected $\Gcat$-operad $\Op$ \cref{i_connected}}} 
& \multirow{2}{*}{\makecell{$\FG$ \eqref{def:FG}}}
& \multirow{2}{*}{\makecell{orthogonal \\ $G$-spectra \eqref{def:gsp_module}}}
& \multirow{2}{*}{\makecell{based on \\ Segal $K$-theory}} 
\\
&&&&&& \\ \hline 
\multirow{2}{*}{\makecell{multifunctorial \\ \cref{kgosg_intro}, \cref{khgosg_intro}}} 
& \multirow{2}{*}{\makecell{\cref{Kgo_functors} \\ \cref{Khgo_functors}}}
& \multirow{2}{*}{\makecell{\cref{kgoav} \\ \cref{khgoav}}}
& \multirow{2}{*}{\makecell{$\Op$-pseudoalgebras \\ $\Tinf$-operad $\Op$ \eqref{as:OpA}}} 
& \multirow{2}{*}{\makecell{$\Gsk$ \eqref{def:Gsk} \\ $\GG$ \eqref{def:GG}}} 
& \multirow{2}{*}{\makecell{orthogonal \\ $G$-spectra}}
& \multirow{2}{*}{\makecell{enriched \\ multifunctor}} 
\\
&&&&&& \\
\end{tabular}}}}
\caption{The four equivariant $K$-theory machines proved to be equivalent in this work.}
\label{table.eqmach}
\end{figure}
\noindent
Each of Shimakawa $K$-theory \cite{shimakawa89,shimakawa91}, multifunctorial equivariant $K$-theory \cite{yau-eqk}, and GMMO $K$-theory \cite{gmmo23} sends operadic (pseudo)algebras to orthogonal $G$-spectra.  Shimakawa $K$-theory is based on nonequivariant Segal $K$-theory.  The indexing $G$-category $\FG$ of pointed finite $G$-sets replaces the indexing category $\Fsk$ of pointed finite sets used in Segal $K$-theory.  The author's multifunctorial equivariant $K$-theory is an enriched multifunctor.  It preserves equivariant algebras parametrized by equivariant operads, including equivariant $\Einf$-algebras.  The equivariant $K$-theory of Guillou, May, Merling, and Osorno satisfies a multiplicative equivariant Barratt-Priddy-Quillen Theorem.  Schwede global equivariant $K$-theory \cite{schwede_global} sends parsummable categories to topological symmetric spectra.  The word global means that the machine keeps track of $G$-actions for all finite groups $G$ in a compatible manner.

\subsection*{Summary of Comparisons}

\Cref{table.foureqmach} summarizes the comparisons of the four equivariant $K$-theory machines established in this work.  
\begin{figure}[H] 
\centering
\resizebox{\linewidth}{!}{%
{\renewcommand{\arraystretch}{1.3}%
{\setlength{\tabcolsep}{1ex}
\begin{tabular}{c|c|c|c|c|c|c}
\multirow{2}{*}{\makecell{equivariant \\ $K$-theories}} 
& \multirow{2}{*}{\makecell{comparison \\ theorems}}  
& \multirow{2}{*}{\makecell{explicit \\ description}}  
& \multirow{2}{*}{\makecell{comparison \\ diagrams}} 
& \multirow{2}{*}{\makecell{categorical \\ comparisons}} 
& \multirow{2}{*}{\makecell{properties of \\ comparisons}}
& \multirow{2}{*}{\makecell{method of \\ proofs}}  
\\
&&&&&& \\ \hline
\multirow{5}{*}{\makecell{Shimakawa and \\ multifunctorial}}
& \multirow{5}{*}{\makecell{topological \\ \ref{kfgkgg_compare}, \ref{kfgbkgg_compare}, \ref{kggkfgb_compare} \\ \phantom{X} \\ categorical \\ \ref{HgoSgo}, \ref{thm:PistAequivalence}, \ref{thm:pistweakgeq}}}
& \multirow{5}{*}{\makecell{\ref{expl:retn}, \ref{expl:kfgkgg_compare} \\ \phantom{X} \\ \phantom{X} \\ \cref{PistA}}}
& \multirow{5}{*}{\makecell{\cref{KgoKsh_comparison} \\ \cref{Pistsg_prod_diag}}}
& \multirow{5}{*}{\makecell{$\Pistsg_{\A}$ \cref{PistsgA} \\ $\Op$-pseudoalgebras $\A$ \\ $\Uinf$-operad $\Op$ \eqref{as:OpA'}}} 
& \multirow{5}{*}{\makecell{categorical weak \\ $G$-equivalences for \\ $\Oph = \Catg(\EG,\Op)$ \\ $\Ah = \Catg(\EG,\A)$ \\ \eqref{catgego}, \eqref{def:cat_weakg}}} 
& \multirow{5}{*}{\makecell{embed into \\ $G$-thickening \\ $\Catg(\EG,-)$ \\ \cref{pistsgahnbe_inc}}} 
\\
&&&&&& \\
&&&&&& \\ 
&&&&&& \\  
&&&&&& \\ \hline 
\multirow{3}{*}{\makecell{Shimakawa and \\ GMMO}}
& \multirow{3}{*}{\makecell{\ref{thm:gmmo_shi}}}
& \multirow{3}{*}{\makecell{\cref{strvx_mal} \\ \cref{cgsam_def}, \cref{cgsax_comp}}}
& \multirow{3}{*}{\makecell{\cref{kgmmo_shi_intro} \\ \cref{gmmo_shi_zig}}}
& \multirow{3}{*}{\makecell{$\cgs_{\A}$ \cref{cgsa} \\ $\Op$-algebras $\A$ \\ chaotic $\Einfg$-operad $\Op$}} 
& \multirow{3}{*}{\makecell{same as above}} & \multirow{3}{*}{\makecell{same as above \\ \cref{zdsg_catweakg}, \cref{cgs_ah_wgeq}}}
\\
\end{tabular}}}}
\caption{Comparisons of the four equivariant $K$-theory machines.}
\label{table.foureqmach}
\end{figure}
\begin{description}
\item[Shimakawa, multifunctorial, and GMMO]
In the top two rows in \cref{table.foureqmach}, for general input data $(\Op,\A)$ with $\A$ an $\Op$-(pseudo)algebra, the comparison morphisms $\Pistsg_\A$ and $\cgs_\A$ compare Shimakawa $K$-theory with, respectively, multifunctorial equivariant $K$-theory and GMMO $K$-theory.  Each component of each of $\Pistsg_\A$ and $\cgs_\A$ is a pointed $G$-functor and a nonequivariant equivalence of categories.  See \cref{Pist_to_prod,thm:PistAequivalence,thm:zbsg_eq,thm:zdsg_eq,cgsam_def}.

In order for the comparison morphisms $\Pistsg$ and $\cgs$ to be componentwise categorical weak $G$-equivalences \pcref{def:cat_weakg}, the input pair $(\Op,\A)$ has to be replaced by its $G$-thickening $(\Oph,\Ah)$ consisting of the $\Gcat$-operad \pcref{catgego}
\[\Oph = \Catg(\EG,\Op)\] 
and the $\Oph$-(pseudo)algebra
\[\Ah = \Catg(\EG,\A).\]  
Here, $\Catg(-,-)$ is the $G$-category of functors and natural transformations with the conjugation $G$-action \pcref{def:Catg}.  The category $\EG$ is the translation category of the group $G$ \pcref{def:translation_cat}, and $\A$ is an $\Op$-(pseudo)algebra.  For example, the $G$-thickening of the Barratt-Eccles operad $\BE$ is the $G$-Barratt-Eccles operad 
\[\GBE = \Catg(\EG,\BE),\] 
whose algebras and pseudoalgebras are, respectively, genuine permutative $G$-categories and genuine symmetric monoidal $G$-categories \pcref{def:BE,def:GBE,def:GBE_pseudoalg}.  The proofs of the comparison categorical weak $G$-equivalences $\Pistsg_{\Ah}$ and $\cgs_{\Ah}$ are structurally similar.  Each proof involves embedding into the $G$-thickening construction $\Catg(\EG,-)$ and a 2-out-of-3 argument.  See \cref{pistsgahnbe_inc,zdsg_catweakg,cgs_ah_wgeq}.
\item[GMMO and Schwede Global]
In the bottom row in \cref{table.foureqmach}, the comparison $\cglg$ compares Schwede global $K$-theory with GMMO $K$-theory.  The comparison $\cglg$ has stronger properties than the comparisons $\Pistsg$ and $\cgs$ in the following sense.  For each input parsummable category $\C$, each component of the comparison $\cglgc$ is an adjoint equivalence in the pointed $G$-equivariant sense \pcref{def:adjointGeq}.  The proof constructs an explicit right $G$-adjoint inverse, along with invertible $G$-equivariant unit and counit that satisfy the triangle identities for a $G$-adjoint equivalence.  See \cref{cglgc_equiv}.
\end{description}
\crefrange{sec:intro_shim}{sec:intro_gmmo_gl} have more elaborate summaries of these equivariant $K$-theory machines and their comparisons.

\section{Shimakawa and Multifunctorial Equivariant $K$-Theories}
\label{sec:intro_shim}

This section summarizes \crefrange{part:kgo_shi}{part:kgo_shi_comp} on the comparison between Shimakawa equivariant $K$-theory \cite{shimakawa89,shimakawa91} and the author's multifunctorial equivariant $K$-theory \cite{yau-eqk}.  The first half of this section briefly summarizes these equivariant $K$-theories.  The second half of this section summarizes their topological and categorical comparisons.  The key comparison diagram is \cref{KgoKsh_comparison}.

\subsection*{Shimakawa Equivariant $K$-Theory}
Shimakawa $K$-theory is based on nonequivariant Segal $K$-theory \cite{segal,segal_eht}.  The latter machine sends small permutative categories to special $\Fsk$-spaces and then to connective spectra.  To compare Shimakawa $K$-theory with multifunctorial equivariant $K$-theory \cite{yau-eqk} and GMMO $K$-theory \cite{gmmo23}, this work extends Shimakawa's original construction so that it accepts operadic pseudoalgebras as input data.  This introductory chapter focuses on the strong variants.

For a compact Lie group $G$ and a 1-connected $\Gcat$-operad $\Op$ \pcref{def:Catg,as:Op_iconn,def:voperad}, \emph{Shimakawa strong $K$-theory} is the composite functor \cref{ksho_kshosg} 
\begin{equation}\label{kshosg_intro}
\begin{tikzpicture}[vcenter]
\def\h{2} \def\u{.6}
\draw[0cell=.9]
(0,0) node (a1) {\AlgpspsO}
(a1)++(1.1*\h,0) node (a2) {\FGCatg}
(a2)++(\h,0) node (a3) {\FGTopg}
(a3)++(\h,0) node (a4) {\Gspec.} 
;
\draw[1cell=.85]
(a1) edge node {\Sgosg} (a2)
(a2) edge node {\clast} (a3)
(a3) edge node {\Kfg} (a4)
(a1) [rounded corners=2pt, shorten <=-.2ex] |- ($(a2)+(0,\u)$) -- node {\Kshosg} ($(a3)+(0,\u)$) -| (a4)
;
\end{tikzpicture}
\end{equation}
\begin{description}
\item[Categories]
The domain of the functor $\Kshosg$ is the 2-category $\AlgpspsO$ of $\Op$-pseudoalgebras, $\Op$-pseudomorphisms, and $\Op$-transformations \pcref{oalgps_twocat}.  The codomain is the category $\Gspec$ of orthogonal $G$-spectra and $G$-morphisms \pcref{def:gsp_morphism}.  In the intermediate categories, the indexing $G$-category $\FG$ has pointed finite $G$-sets as objects and pointed functions as morphisms \pcref{def:FG}.  The pointed $G$-category $\Catgst$ \cref{Catgst_iicat} has small pointed $G$-categories as objects and pointed functors as morphisms.  The group $G$ acts by conjugation on the morphisms of $\FG$ and $\Catgst$.   The 2-category $\FGCatg$ has pointed $G$-functors $\FG \to \Catgst$, called \emph{$\FGG$-categories}, as objects, $G$-natural transformations as 1-cells, and $G$-modifications as 2-cells \pcref{def:fgcatg}.  The category $\FGTopg$ is defined analogously using pointed $G$-spaces \pcref{def:ggtopg}.
\item[Functors]
The 2-functor $\Sgosg$ in \cref{kshosg_intro}, called \emph{Shimakawa strong $H$-theory} \pcref{ch:shih}, sends $\Op$-pseudoalgebras to $\FGG$-categories, $\Op$-pseudomorphisms to $G$-natural transformations, and $\Op$-transformations to $G$-modifications.  The functor $\clast$ \cref{clast_ggcatg} is given levelwise by the classifying space functor $\cla$ \cref{classifying_space}.  The prolongation functor $\Kfg$ in \cref{kshosg_intro} produces orthogonal $G$-spectra from $\FGG$-spaces.
\item[Notational Convention]
In the functors $\Kshosg$ and $\Sgosg$, the subscript \texttt{Sh} refers to Shimakawa.  The tilde decoration refers to the strong variant.  A functor that builds equivariant spectra from either categorical or topological data, such as $\Kshosg$ and $\Kfg$, is called \emph{$K$-theory} with corresponding notation.  A functor that sends equivariant categorical data to equivariant categories indexed by an equivariant category, such as $\Sgosg$, is called \emph{$H$-theory} with corresponding notation.
\end{description}
While Shimakawa strong $K$-theory can be defined for compact Lie groups $G$, some further properties only hold for finite groups $G$.

There are two notable characteristics of Shimakawa strong $K$-theory $\Kshosg$.  First, Shimakawa strong $H$-theory $\Sgosg$ \cref{def:sgosg} sends pseudo structures in its domain $\AlgpspsO$ to strict structures in its codomain $\FGCatg$.  However, the 2-functor $\Sgosg$ does \emph{not} involve any strictification functor that would first strictify the pseudo structures in its domain.  Instead, the associativity constraint of an $\Op$-pseudoalgebra is directly incorporated into the associativity axiom \cref{nsys_associativity} in the construction of the 2-functor $\Sgosg$.

Second, the prolongation functor $\Kfg$ is not homotopical.  In other words, it does not generally preserve componentwise weak $G$-equivalences, unless one restricts to \emph{proper} $\FGG$-spaces \pcref{thm:Kfg_inv} and a finite group $G$.  One remedy is to precompose the prolongation functor $\Kfg$ with the bar construction $\Bc$ \pcref{bar_proper}.  This idea leads to the homotopical variant of Shimakawa strong $K$-theory, which is discussed next.

\subsection*{Homotopical Shimakawa Machine}
The \emph{homotopical Shimakawa strong $K$-theory} is the composite functor \cref{khsho_khshosg}
\begin{equation}\label{khshosg_intro}
\begin{tikzpicture}[vcenter]
\def\h{2} \def\u{.6}
\draw[0cell=.9]
(0,0) node (a1) {\phantom{\AlgpspsO}}
(a1)++(0,0) node (a1') {\AlgpspsO}
(a1)++(1.1*\h,0) node (a2) {\FGCatg}
(a2)++(\h,0) node (a3) {\FGTopg}
(a3)++(\h,0) node (a4) {\FGTopg}
(a4)++(\h,0) node (a5) {\Gspec} 
;
\draw[1cell=.85]
(a1) edge node {\Sgosg} (a2)
(a2) edge node {\clast} (a3)
(a3) edge node {\Bc} (a4)
(a4) edge node {\Kfg} (a5)
(a1) [rounded corners=2pt, shorten <=-.2ex] |- ($(a2)+(0,\u)$) -- node {\Khshosg} ($(a4)+(0,\u)$) -| (a5)
;
\end{tikzpicture}
\end{equation}
that includes the bar functor $\Bc$ \cref{bar_functor_FG} before the prolongation functor $\Kfg$.  The letter $h$ in the superscript of $\Khshosg$ is an abbreviation of \emph{homotopical}.  For a finite group $G$, via the retraction $\retn \cn \Bc \to 1_{\FGTopg}$ \cref{retn_barFG_id}, the functors $\Kshosg$ \cref{kshosg_intro} and $\Khshosg$ \cref{khshosg_intro} are naturally componentwise weakly $G$-equivalent \pcref{thm:ksho_khsho}.  Shimakawa's original construction in \cite{shimakawa89} is naturally componentwise weakly $G$-equivalent to the homotopical Shimakawa strong $K$-theory $\Khshbesg$ for the Barratt-Eccles operad $\BE$ \pcref{def:BE}.  \cref{ch:shimakawa_K} constructs Shimakawa's functors $\Kshosg$ and $\Khshosg$.

At the object level, Shimakawa's functors $\Kshosg$ and $\Khshosg$ send $\Op$-pseudoalgebras to orthogonal $G$-spectra \pcref{def:pseudoalgebra,def:gsp_module}.  Similar to nonequivariant Segal $K$-theory, the functors $\Sgosg$, $\Kshosg$, and $\Khshosg$ do \emph{not} generally preserve multiplicative structures such as monoids and their modules.  From the point of view of multiplicative equivariant stable homotopy theory, it is desirable to have a multiplicative alternative to Shimakawa's $K$-theory functors.

\subsection*{Multifunctorial Equivariant $K$-Theory}

In the work \cite{yau-eqk}, the author constructs an equivariant $K$-theory machine that preserves equivariant multiplicative structures, including equivariant $\Einf$-algebras.  A \emph{$\Tinf$-operad} is a $G$-categorically enriched pseudo-commutative operad $\Op$ \pcref{def:pseudocom_operad} such that $\Op(1)$ is a terminal $G$-category \pcref{as:OpA}.  For a compact Lie group $G$ and a $\Tinf$-operad $\Op$, the equivariant $K$-theory constructed in \cite[\namecref{EqK:thm:Kgo_multi} \ref*{EqK:thm:Kgo_multi}]{yau-eqk} is the composite \emph{enriched multifunctor}
\begin{equation}\label{Kgosg_intro}
\begin{tikzpicture}[vcenter]
\def\h{2.3} \def\u{.6}
\draw[0cell=.9]
(0,0) node (a1) {\MultpspsO}
(a1)++(\h,0) node (a2) {\GGCat}
(a2)++(\h,0) node (a3) {\GGTop}
(a3)++(.85*\h,0) node (a4) {\GSp.} 
;
\draw[1cell=.85]
(a1) edge node {\Jgosg} (a2)
(a2) edge node {\clast} (a3)
(a3) edge node {\Kg} (a4)
(a1) [rounded corners=2pt, shorten <=-.2ex] |- ($(a2)+(0,\u)$) -- node {\Kgosg} ($(a3)+(0,\u)$) -| (a4)
;
\end{tikzpicture}
\end{equation}
\begin{description}
\item[Multicategories]
The domain is the $\Gcat$-enriched multicategory $\MultpspsO$ of $\Op$-pseudoalgebras.  The codomain is the $\Gtop$-enriched multicategory $\GSp$ of orthogonal $G$-spectra.  The indexing category $\Gsk$ has finite tuples of pointed finite sets as objects.  Its morphisms are generated by injections and permutations of entries, along with pointed functions in each entry \pcref{def:Gsk}.  The $\Gcat$-multicategory $\GGCat$ has pointed functors $\Gsk \to \Gcatst$, called \emph{$\Gskg$-categories}, as objects.  The $\Gtop$-multicategory $\GGTop$ has pointed functors $\Gsk \to \Gtopst$, called \emph{$\Gskg$-spaces}, as objects.  
\item[Multifunctors]
At the object level, the $\Gcat$-multifunctor $\Jgosg$ in \cref{Kgosg_intro} sends $\Op$-pseudoalgebras to $\Gskg$-categories.  The $\Gtop$-multifunctor $\clast$ is given by postcomposing with the classifying space functor $\cla$.  The prolongation $\Gtop$-multifunctor $\Kg$ in \cref{Kgosg_intro} sends $\Gskg$-spaces to orthogonal $G$-spectra.  The composite $\Gtop$-multifunctor $\Kgosg$ sends $\Op$-pseudoalgebras to orthogonal $G$-spectra.
\item[Notational Convention]
A (multi)functor that sends equivariant categorical data to equivariant categories indexed by a nonequivariant category, such as $\Jgosg$, is called \emph{$J$-theory} with corresponding notation.
\end{description}

The enriched multifunctor $\Kgosg$ preserves the $\Gcat$-multicategory structure of its domain $\MultpspsO$ and the $\Gtop$-multicategory structure of its codomain $\GSp$.  Thus, $\Kgosg$ sends equivariant $\Einf$-algebras of $\Op$-pseudoalgebras to equivariant $\Einf$-algebras of orthogonal $G$-spectra.  The main purpose of \crefrange{part:kgo_shi}{part:kgo_shi_comp} is to prove that Shimakawa's equivariant $K$-theory functors $\Kshosg$ \cref{kshosg_intro} and $\Khshosg$ \cref{khshosg_intro} are equivalent to $\Kgosg$ \cref{Kgosg_intro} in an appropriate sense.

\subsection*{Underlying Functor of Multifunctorial Equivariant $K$-Theory}
Since Shimakawa's equivariant $K$-theory functors $\Kshosg$ \cref{kshosg_intro} and $\Khshosg$ \cref{khshosg_intro} are not multifunctors, the comparison with the author's equivariant $K$-theory enriched multifunctor $\Kgosg$ \cref{Kgosg_intro} occurs at the functor level.  Thus, it suffices to consider the underlying functor 
\begin{equation}\label{kgosg_intro}
\begin{tikzpicture}[vcenter]
\def\h{2.2} \def\u{.6}
\draw[0cell=.9]
(0,0) node (a1) {\AlgpspsO}
(a1)++(\h,0) node (a2) {\GGCatii}
(a2)++(\h,0) node (a3) {\GGTopii}
(a3)++(.85*\h,0) node (a4) {\Gspec} 
;
\draw[1cell=.85]
(a1) edge node {\Jgosg} (a2)
(a2) edge node {\clast} (a3)
(a3) edge node {\Kg} (a4)
(a1) [rounded corners=2pt, shorten <=-.2ex] |- ($(a2)+(0,\u)$) -- node {\Kgosg} ($(a3)+(0,\u)$) -| (a4)
;
\end{tikzpicture}
\end{equation}
of the $\Gtop$-multifunctor $\Kgosg$.  Each of the four categories in \cref{kgosg_intro} has the same objects as the corresponding enriched multicategory in \cref{Kgosg_intro}.  Each of the four functors in \cref{kgosg_intro} has the same object assignment as the corresponding enriched multifunctor in \cref{Kgosg_intro}.  At the morphism level, each of the four functors in \cref{kgosg_intro} is the restriction of the corresponding enriched multifunctor to either
\begin{itemize}
\item 1-ary $G$-fixed 1-cells and 1-ary $G$-fixed 2-cells in the $\Gcat$-enrichment or
\item 1-ary $G$-fixed subspaces in the $\Gtop$-enrichment.
\end{itemize}  
\cref{ch:kgo} discusses the functor $\Kgosg$ \cref{kgosg_intro} in detail.

\subsection*{An Isomorphic Variant of the Functor $\Kgosg$}
To compare Shimakawa's functors with the functor $\Kgosg$ \cref{kgosg_intro}, it is necessary to reconcile the difference between
\begin{itemize}
\item the indexing $G$-category $\FG$, with pointed finite $G$-sets as object, used in Shimakawa's functors $\Kshosg$ and $\Khshosg$, and
\item the indexing category $\Gsk$, with finite tuples of pointed finite sets as objects, used in $\Kgosg$.  
\end{itemize}
Mediating between the indexing $G$-category $\FG$ and the indexing category $\Gsk$ is the indexing $G$-category $\GG$ \pcref{def:GG}.  Its objects are finite tuples of pointed finite $G$-sets.  There are pointed functors
\begin{equation}\label{ifg_intro}
\begin{tikzpicture}[vcenter]
\def\h{1.8}
\draw[0cell]
(0,0) node (a1) {\FG}
(a1)++(\h,0) node (a2) {\GG}
(a2)++(.9*\h,0) node (a3) {\Gsk}
;
\draw[1cell=.9]
(a1) edge[right hook->, transform canvas={yshift=-.5ex}] node[inner sep=2pt] {\ifg} (a2)
(a2) edge[bend right=25, transform canvas={yshift=0ex}] node[swap] {\sma} (a1)
(a3) edge[left hook->] node[swap] {\ig} (a2)
;
\end{tikzpicture}
\end{equation}
relating the categories $\FG$, $\GG$, and $\Gsk$.  The full subcategory inclusion pointed $G$-functor $\ifg$ regards each pointed finite $G$-set as a length-1 object in $\GG$.  It admits a retraction given by the smash functor $\sma$, so 
\[\sma\ifg = 1_{\FG}.\]  
The full subcategory inclusion pointed functor $\ig$ equips each pointed finite set in each object in $\Gsk$ with the trivial $G$-action.  

The functor $\Kgosg$ \cref{kgosg_intro} is naturally isomorphic to the composite functor 
\begin{equation}\label{khgosg_intro}
\begin{tikzpicture}[vcenter]
\def\h{2.2} \def\u{.6}
\draw[0cell=.9]
(0,0) node (a1) {\AlgpspsO}
(a1)++(\h,0) node (a2) {\GGCatg}
(a2)++(\h,0) node (a3) {\GGTopg}
(a3)++(.9*\h,0) node (a4) {\Gspec} 
;
\draw[1cell=.85]
(a1) edge node {\Hgosg} (a2)
(a2) edge node {\clast} (a3)
(a3) edge node {\Kgg} (a4)
(a1) [rounded corners=2pt, shorten <=-.2ex] |- ($(a2)+(0,\u)$) -- node {\Khgosg} ($(a3)+(0,\u)$) -| (a4)
;
\end{tikzpicture}
\end{equation}
that involves the indexing $G$-category $\GG$ \pcref{KgoKhgo}.  \cref{ch:ggcatg} discusses the indexing $G$-category $\GG$ and the 2-category $\GGCatg$.  \cref{ch:hgo} constructs the 2-functor $\Hgosg$, called \emph{strong $H$-theory}, the prolongation functor $\Kgg$, and the equivariant $K$-theory functor $\Khgosg$.  It also proves that the functor $\Khgosg$ \cref{khgosg_intro} is naturally isomorphic to the functor $\Kgosg$ \cref{kgosg_intro}.

\subsection*{Comparison Diagram: Shimakawa and Multifunctorial}
The diagram \cref{KgoKsh_comparison} summarizes the topological and categorical comparisons between Shimakawa and multifunctorial equivariant $K$-theories.
\begin{itemize}
\item The top row in \cref{KgoKsh_comparison} displays Shimakawa's equivariant $K$-theory functors $\Khshosg$ \cref{khshosg_intro}, which involves the bar functor $\Bc$, and $\Kshosg$ \cref{kshosg_intro}, which does not involve the bar functor.
\item The bottom row in \cref{KgoKsh_comparison} displays the functor $\Kgosg$ \cref{kgosg_intro}, which is the underlying functor of the author's equivariant $K$-theory enriched multifunctor \cref{Kgosg_intro}.
\item The middle row in \cref{KgoKsh_comparison} displays the functor $\Khgosg$ \cref{khgosg_intro}, which is naturally isomorphic to the functor $\Kgosg$ in the bottom row.
\end{itemize}
\begin{equation}\label{KgoKsh_comparison}

\end{equation}
\begin{description}
\item[Bottom] In the bottom half of the diagram \cref{KgoKsh_comparison}, each functor $\igst$ is the pullback along the inclusion functor $\ig \cn \Gsk \to \GG$ \cref{ifg_intro}.  The functor $\igst$ admits a left adjoint inverse $\Lg$ \pcref{thm:ggtop_ggtopg_iieq}.  The left and middle regions strictly commute \pcref{jgohgoigst}.  The right region commutes up to natural isomorphisms \pcref{KgKgg}:
\[\begin{split}
\GGTopii & \fto{\Kg \iso \Kgg\Lg} \Gspec \andspace\\
\GGTopg & \fto{\Kgg \iso \Kg\igst} \Gspec.
\end{split}\]
The second natural isomorphism follows from the first one and the counit of the adjoint equivalence $(\Lg,\igst)$.
\item[Top] In the top half of the diagram \cref{KgoKsh_comparison}, the functors $\smast$ and $\ifgst$ are the pullbacks along the functors $\sma$ and $\ifg$ \cref{ifg_intro}.  They exhibit the category $\FGCatg$ as a retract of the category $\GGCatg$ in the sense that
\[\ifgst\smast = 1_{\FGCatg}.\]  
The functors $\ifgl$ and $\smal$ are the left adjoints of, respectively, the topological analogues of the pullback functors $\ifgst$ and $\smast$.  They exhibit the category $\FGTopg$ as a retract of the category $\GGTopg$ in the sense that there is a natural isomorphism \pcref{ex:smai_adj}
\[\smal\ifgl \iso 1_{\FGTopg}.\]
\end{description}
The comparison of Shimakawa's functors $\Kshosg$ and $\Khshosg$ with the functor $\Khgosg$ consists of two parts.
\begin{description}
\item[Topological] The topological comparison computes Shimakawa's functors $\Kfg$ and $\Kfg\Bc$ in terms of the functor $\Kgg$ and vice versa.  This comparison is displayed in the top right regions in \cref{KgoKsh_comparison}.  It is summarized in the upcoming subsection \nameref{subsec:top_comp}.
\item[Categorical] The categorical comparison computes Shimakawa strong $H$-theory $\Sgosg$ in terms of strong $H$-theory $\Hgosg$ and vice versa.  This comparison is displayed in the top left region in \cref{KgoKsh_comparison}.  It is summarized in the subsections after the next one.  
\end{description}

\subsection*{Topological Comparison}
\label{subsec:top_comp}

This subsection summarizes the topological comparisons in \cref{ch:shim_top}.  The main topological comparison results are \cref{kfgkgg_compare,kfgbkgg_compare,kggkfgb_compare}.  Suppose $G$ is a compact Lie group.  Some of the results require $G$ to be finite.
\begin{description}
\item[Comparing $\Kfg$ and $\Kgg$]
In the upper right part of the comparison diagram \cref{KgoKsh_comparison}, the functors $\Kfg$ and $\Kgg$ factor through each other up to natural isomorphisms \pcref{kfgkgg_compare}:
\[\begin{split}
\FGTopg & \fto{\Kfg \iso \Kgg\ifgl} \Gspec \andspace\\
\GGTopg & \fto{\Kgg \iso \Kfg\smal\!} \Gspec.
\end{split}\]
Thus, any orthogonal $G$-spectrum produced by Shimakawa's prolongation functor $\Kfg$ is also produced by $\Kgg$ and vice versa.
\item[Computing $\Kfg\Bc$ using $\Kgg$]
Taking the bar functor $\Bc$ \cref{bar_functor_FG} into account, the retraction $\retn \cn \Bc \to 1_{\FGTopg}$ \cref{retn_barFG_id} induces a natural transformation \pcref{kfgbkgg_compare}
\begin{equation}\label{kbgg_intro}
\begin{tikzpicture}[vcenter]
\def\h{2.2}
\draw[0cell]
(0,0) node (a1) {\Kfg\Bc}
(a1)++(\h,0) node (a2) {\Kfg} 
(a2)++(.8*\h,0) node (a3) {\phantom{\Kgg\ifgl}} 
(a3)++(0,-.04) node (a3') {\Kgg\ifgl}
;
\draw[1cell=.9]
(a1) edge node {\Kfg\retn} (a2)
(a2) edge node {\iso} (a3)
;
\draw[1cell=.9]
(a1) [rounded corners=2pt] |- ($(a2)+(-1,.6)$) -- node {\kbgg} ($(a2)+(1,.6)$) -| (a3')
 ;
\end{tikzpicture}
\end{equation}
that computes Shimakawa's homotopical prolongation functor $\Kfg\Bc$ in terms of $\Kgg$ in the following sense.  For a finite group $G$ and a proper $\FGG$-space $X$ \pcref{def:proper_fgg}, the $G$-morphism of orthogonal $G$-spectra 
\[\Kfg\Bc X \fto{\kbgg_X} \Kgg\ifgl X\]
is componentwise a weak $G$-equivalence between $G$-spaces \pcref{def:weakG_top}.  For example, each $\FGG$-space $X$ in the image of the levelwise classifying space functor $\clast$ is proper \pcref{reast_proper}, so $\kbgg_X$ is componentwise a weak $G$-equivalence for such $X$.
\item[Computing $\Kgg$ using $\Kfg\Bc$]
Conversely, there is a natural transformation \pcref{kggkfgb_compare}
\begin{equation}\label{kgbk_intro}
\begin{tikzpicture}[vcenter]
\def\h{2.7} \def\u{.65}
\draw[0cell]
(0,0) node (a1) {\Kfg\Bc\smal}
(a1)++(\h,0) node (a2) {\Kfg\smal} 
(a2)++(.7*\h,0) node (a3) {\phantom{\Kgg}}
(a3)++(0,.04) node (a3') {\Kgg}
;
\draw[1cell=.9]
(a1) edge node {\Kfg\retn_{\smal}} (a2)
(a2) edge node {\iso} (a3)
;
\draw[1cell=.9]
(a1) [rounded corners=2pt] |- ($(a2)+(-.3*\h,\u)$) -- node {\kgbk} ($(a2)+(0,\u)$) -| (a3')
 ;
\end{tikzpicture}
\end{equation}
that computes $\Kgg$ in terms of Shimakawa's homotopical prolongation functor $\Kfg\Bc$ in the following sense.  For a finite group $G$ and a proper $\GGG$-space $Y$ \pcref{def:proper_ggg}, the $G$-morphism of orthogonal $G$-spectra 
\[\Kfg\Bc\!\smal\! Y \fto{\kgbk_Y} \Kgg Y\]
is componentwise a weak $G$-equivalence between $G$-spaces.  For example, each $\GGG$-space $Y$ in the image of the levelwise classifying space functor $\clast$ is proper (\cref{reast_properGG} \cref{reast_properGG_ii}), so $\kgbk_Y$ is componentwise a weak $G$-equivalence for such $Y$.
\end{description}  
In summary, Shimakawa's homotopical prolongation functor $\Kfg\Bc$ and the functor $\Kgg$ compute each other in the following sense.  For a finite group $G$ and proper objects, the functors $\Kfg\Bc$ and $\Kgg\ifgl$ are naturally componentwise weakly $G$-equivalent, and so are the functors $\Kfg\Bc\smal$ and $\Kgg$.

\subsection*{Categorical Comparison: Computing $\Sgosg$ using $\Hgosg$}
\label{subsec:cat_comp_i}
This subsection and the next one summarize the categorical comparison between Shimakawa $K$-theory and $\Khgosg$ that is discussed in detail in \cref{part:kgo_shi_comp}.  This subsection summarizes the first categorical comparison result.   

Suppose that $\Op$ is a $\Tinf$-operad \pcref{as:OpA} for a group $G$. 
The first categorical comparison is the equality of 2-functors \pcref{HgoSgo}
\[\AlgpspsO \fto{\Sgosg = \ifgst\Hgosg} \FGCatg.\]
In other words, Shimakawa strong $H$-theory $\Sgosg$ is the restriction of strong $H$-theory $\Hgosg$ to length-1 objects in $\GG$.  Thus, each $\FGG$-category produced by Shimakawa strong $H$-theory $\Sgosg$ can also be produced by strong $H$-theory $\Hgosg$.  This categorical comparison is discussed in \cref{sec:factor_hj}.

\subsection*{Categorical Comparison: Computing $\Hgosg$ using $\Sgosg$}
\label{subsec:cat_comp_ii}

The second categorical comparison between Shimakawa $K$-theory and $\Khgosg$, called the \emph{strong $H$-theory comparison}, is the 2-natural transformation \cref{Pistarsg_twonat}
\begin{equation}\label{pistsg_intro}
\begin{tikzpicture}[vcenter]
\def\h{2.3} \def\v{.7} \def\s{15}
\draw[0cell]
(0,0) node (a1) {\AlgpspsO}
(a1)++(\h,\v) node (a2) {\FGCatg}
(a1)++(\h,-\v) node (a3) {\GGCatg}
;
\draw[1cell=.9]
(a1) edge[bend left=\s] node {\Sgosg} (a2)
(a2) edge[bend left=\s, shorten <=-.5ex] node {\smast} (a3)
(a1) edge[bend right=\s] node[swap] {\Hgosg} (a3)
;
\draw[2cell]
node[between=a2 and a3 at .55, shift={(-.3*\h,0)}, rotate=-120, 2labelw={below,\Pistsg,1pt}] {\Rightarrow}
;
\end{tikzpicture}
\end{equation}
that expresses $\Hgosg$ in terms of Shimakawa strong $H$-theory $\Sgosg$.  All of \cref{part:kgo_shi_comp} after \cref{sec:factor_hj} studies the comparison 2-natural transformation $\Pistsg$ and related constructions.  The main results about $\Pistsg$ are \cref{thm:PistAequivalence,thm:pistweakgeq}, which are summarized next.  
\begin{description}
\item[Comparison nonequivariant equivalences]
A \emph{$\Uinf$-operad} is a 1-connected $\Gcat$-operad that is levelwise a nonempty translation category \pcref{as:OpA'}.  \cref{thm:PistAequivalence} proves that, for each $\Uinf$-operad $\Op$, $\Op$-pseudoalgebra $\A$ \pcref{def:pseudoalgebra}, and object $\nbe \in \GG$ \cref{GG_objects}, the $\nbe$-component pointed $G$-functor
\begin{equation}\label{pistsganbe_intro}
(\smast\Sgosg\A)\angordnbe \fto[\sim]{\Pistsg_{\A,\angordnbe}} (\Hgosg\A)\angordnbe
\end{equation}
is a nonequivariant equivalence of categories.  The proof of \cref{thm:PistAequivalence} occupies \cref{sec:pistar,sec:pistar_proof,sec:hcomp_prod,sec:sgo_tprod,sec:sgo_prod_eq,sec:hgo_tprod,sec:hgo_unit_counit,sec:hgo_prod_eq}.  It involves a factorization of pointed $G$-functors \cref{Pistsg_prod_diag}
\[\zbsg = \zdsg \comp \Pistsg_{\A,\angordnbe}\]
such that $\zbsg$ and $\zdsg$ are nonequivariant equivalences \pcref{thm:zbsg_eq,thm:zdsg_eq}.  Applying the classifying space functor $\cla$ yields a pointed $G$-morphism between $G$-spaces
\[\cla(\smast\Sgosg\A)\angordnbe \fto{\cla\Pistsg_{\A,\angordnbe}} \cla(\Hgosg\A)\angordnbe\]
that is a nonequivariant homotopy equivalence.  Specializing to the Barratt-Eccles operad $\BE$ and the trivial group $G$, the equivalence \cref{pistsganbe_intro} recovers the categorical equivalence between nonequivariant Segal $K$-theory and Elmendorf-Mandell $K$-theory proved in \cite[4.6]{elmendorf-mandell}.  
\item[Pseudo equivariant adjoints]
For a nontrivial group $G$, the comparison pointed $G$-functor $\Pistsg_{\A,\angordnbe}$ in \cref{pistsganbe_intro} is \emph{not} known to be a $G$-equivalence of $G$-categories.  There are no known $G$-functors going backward such that the two composites are $G$-naturally isomorphic to the respective identity functors.  The essential difficulty is that the left adjoint inverses of the pointed $G$-functors $\zbsg$ and $\zdsg$, denoted by $\zbsgad$ \cref{zbsgad} and $\zdsgad$ \cref{zdsgad}, are only \emph{pseudo} $G$-equivariant \pcref{def:pseudoG,thm:zbsgad_pseudo,thm:zdsgad_pseudo} but not strictly $G$-equivariant \pcref{zbsgad_constraint,zdsgad_constraint}.  

In fact, for a nontrivial group $G$, the left adjoint inverses $\zbsgad$ and $\zdsgad$ are not $G$-equivariant even for the Barratt-Eccles operad $\BE$, on which $G$ acts trivially \pcref{ex:al_g_s}.  Thus, the $G$-morphism $\cla\Pistsg_{\A,\angordnbe}$ is not expected to be a weak $G$-equivalence between $G$-spaces in general.  Something more subtle is true, as we discuss next.  
\item[$G$-thickening]
Instead of a general $\Uinf$-operad $\Op$ and an $\Op$-pseudoalgebra $\A$, we consider 
\begin{itemize}
\item the $G$-category $\Catg(-,-)$ of functors and natural transformations with the conjugation $G$-action \pcref{def:Catg},
\item the translation category $\EG$ of a group $G$ \pcref{def:translation_cat}, 
\item the $\Uinf$-operad \pcref{catgego}
\[\Oph =  \Catg(\EG,\Op),\] 
and
\item the $\Oph$-pseudoalgebra 
\[\Ah = \Catg(\EG,\A).\]
\end{itemize}
We call $\Catg(\EG,-)$ the $G$-thickening construction \pcref{def:gcat_inc}.  For example, the $G$-thickening of the Barratt-Eccles operad $\BE$ is the $G$-Barratt-Eccles operad 
\[\GBE = \Catg(\EG,\BE),\] 
whose pseudoalgebras are genuine symmetric monoidal $G$-categories \pcref{def:GBE,def:GBE_pseudoalg}.  

The group $G$ acts freely and transitively on the translation category $\EG$.   Moreover, $\EG$ is nonequivariantly equivalent to the terminal category.  For a small $G$-category $\C$, the $G$-thickening $\Catg(\EG,\C)$ gives the \emph{homotopy points} of $\C$.  The $G$-fixed point category $\Catg(\EG,\C)^G$ is then the \emph{homotopy fixed point category} of $\C$.  By replacing the terminal category with the translation category $\EG$, the $G$-thickening construction $\Catg(\EG,-)$ improves the equivariant properties of small $G$-categories and $G$-functors.
\item[Comparison categorical weak $G$-equivalences]
\cref{thm:pistweakgeq} proves that the strong $H$-theory comparison pointed $G$-functor for the input pair $(\Oph,\Ah)$ 
\begin{equation}\label{pistsgahnbe_intro}
(\smast\Sgohsg\Ah)\angordnbe \fto{\Pistsg_{\Ah,\angordnbe}} (\Hgohsg\Ah)\angordnbe
\end{equation}
is a \emph{categorical weak $G$-equivalence} \pcref{def:cat_weakg} for each object $\nbe \in \GG$.  This means that, for each subgroup $H \subseteq G$, the $H$-fixed subfunctor 
\[\big((\smast\Sgohsg\Ah)\angordnbe \big)^H 
\fto[\sim]{(\Pistsg_{\Ah,\nbe})^H} 
\big((\Hgohsg\Ah)\angordnbe \big)^H\]
is an equivalence of categories.  The proof of \cref{thm:pistweakgeq} involves a 2-out-of-3 argument, the nonequivariant equivalence pointed $G$-functor in \cref{pistsganbe_intro}, and the concept of an adjoint $G$-equivalence \pcref{def:adjointGeq}.  The introduction of \cref{ch:compgen} has an outline of this proof.  Since the classifying space functor $\cla$ commutes with taking $H$-fixed points for each subgroup $H \subseteq G$, applying $\cla$ to the pointed $G$-functor in \cref{pistsgahnbe_intro} yields a weak $G$-equivalence
\[\cla(\smast\Sgohsg\Ah)\angordnbe \fto{\cla\Pistsg_{\Ah,\angordnbe}} \cla(\Hgohsg\Ah)\angordnbe\]
between $G$-spaces.
\end{description}
In summary, for an $\Oph$-pseudoalgebra of the form $\Ah$, each component of the strong $H$-theory comparison $\Pistsg_{\Ah}$ \cref{pistsgahnbe_intro} is a categorical weak $G$-equivalence.

\subsection*{Specialness}
Special $\FGG$-spaces play a crucial role in Shimakawa equivariant $K$-theory \cite[Theorem B]{shimakawa89}; see \cref{expl:sp_gggcat}.  This work generalizes special $\FGG$-spaces to \emph{special $\GGG$-categories}.  A $\GGG$-category is special if each of its Segal functors is a categorical weak $G$-equivalence \pcref{def:cat_weakg,def:sp_gggcat}.

\cref{thm:h_special,cor:h_special_shi} prove that the domain and the codomain of the strong $H$-theory comparison for $(\Oph,\Ah)$, 
\begin{equation}\label{pistsgah_intro}
\smast\Sgohsg\Ah \fto{\Pistsg_{\Ah}} \Hgohsg\Ah \inspace \GGCatg,
\end{equation}
are special $\GGG$-categories.  The proof of the specialness of the $\GGG$-category $\Hgohsg\Ah$ in \cref{thm:h_special} involves a 2-out-of-3 argument similar to the one used in the proof of \cref{thm:pistweakgeq}.  See the commutative diagrams \cref{pistsgahnbe_inc,zdsg_catweakg}.  Moreover, the diagram \cref{zdsg_catweakg} is also used in the proof of \cref{thm:gmmo_shi} \cref{thm:gmmo_shi_iii} to prove the equivalence between Shimakawa and GMMO $K$-theories.  See \cref{sec:intro_shi_gmmo,rk:zdsg_gmmo} for further discussion.

\section{Shimakawa and GMMO $K$-Theories}
\label{sec:intro_shi_gmmo}

This section summarizes \cref{ch:kgmmo} on the equivalence between the homotopical Shimakawa strong $K$-theory and GMMO $K$-theory.  The first subsection briefly summarizes GMMO $K$-theory.  The key comparison diagram is \cref{kgmmo_shi_intro}.

\subsection*{GMMO $K$-Theory}
For a finite group $G$ and a chaotic $\Einfg$-operad $\Op$ enriched in $\Gcat$ \pcref{def:chaotic_einf}, the equivariant $K$-theory functor $\Kgmmo$ of Guillou, May, Merling, and Osorno \cite{gmmo23} is the following composite functor \cref{kgmmo_diag}.
\begin{equation}\label{intro_kgmmo}
\begin{tikzpicture}[vcenter]
\def\v{1.4} \def\h{2.8}
\draw[0cell]
(0,0) node (a1) {\AlgstpsO}
(a1)++(0,-\v) node (a2) {\DGCatg}
(a2)++(0,-\v) node (a3) {\FGCatgps}
(a3)++(\h,0) node (a4) {\FGCatg}
(a4)++(\h,0) node (a5) {\FGTopg}
(a5)++(0,\v) node (a6) {\FGTopg}
(a6)++(0,\v) node (a7) {\Gspec}
;
\draw[1cell=.9]
(a1) edge node {\Kgmmo} (a7)
(a1) edge node[swap] {\Rg} (a2)
(a2) edge node[swap] {\gzest} (a3)
(a3) edge node {\str} (a4)
(a4) edge node {\clast} (a5)
(a5) edge node[swap] {\Bc} (a6)
(a6) edge node[swap] {\Kfg} (a7)
;
\end{tikzpicture}
\end{equation}
The last three functors---the levelwise classifying space functor $\clast$, the bar functor $\Bc$, and the prolongation functor $\Kfg$---also appear in the homotopical Shimakawa strong $K$-theory $\Khshosg$ \cref{khshosg_intro}.  The first three functors in $\Kgmmo$ are discussed in detail in \crefrange{sec:dgcatg}{sec:kgmmo} and briefly summarized next.
\begin{description}
\item[$\Rg$] The domain 2-category $\AlgstpsO$ consists of $\Op$-algebras, $\Op$-pseudomorphisms, and $\Op$-transformations \pcref{def:algstpsO}.  The $\Gcatst$-category $\DG$ has pointed finite $G$-sets as objects and hom pointed $G$-categories constructed from $\Op$ \cref{dg_mn}.  The objects of the 2-category $\DGCatg$ are pointed $\Gcatst$-functors $\DG \to \Gcatst$, called $\DGG$-categories \cref{dgcatg_obj}.  Its 1-cells are $\Pig$-strict $\Gcatst$-pseudotransformations \pcref{def:pig_strict}. Its 2-cells  are $\Gcatst$-modifications \pcref{def:dgcatg_iicell}.  For an $\Op$-algebra $\A$, the $\DGG$-category $\Rg\A$ sends a pointed finite $G$-set $\mal$ to the $\mal$-twisted product \cref{rga_mal}
\[(\Rg\A)\mal = \proAmal.\]
The 1-cell assignment \cref{rga_objx} and 2-cell assignment \cref{rgaf} of $\Rg\A$ use the $\Op$-algebra structure on $\A$.
\item[$\gzest$] The 2-category $\FGCatgps$ has pointed $\Gcatst$-pseudofunctors $\FG \to \Gcatst$, called pseudo $\FGG$-categories, as objects; $\Gcatst$-pseudotransformations as 1-cells; and $\Gcatst$-modifications as 2-cells \pcref{def:psfggcat,def:psfggcat_icell,def:psfggcat_iicell}.  The 2-functor $\gzest$ in $\Kgmmo$ is the pullback along a $\Gcatst$-pseudofunctor $\gze \cn \FG \to \DG$ \cref{gze}.  The construction of $\gze$ uses the assumption that $\Op$ is a chaotic $\Einfg$-operad.
\item[$\str$] The strictification 2-functor $\str$ in $\Kgmmo$ is the left 2-adjoint of a 2-adjunction \cref{str_iifunctor}
\begin{equation}\label{StIn_intro}
\begin{tikzpicture}
\draw[0cell]
(0,0) node (a1) {\FGCatgps}
(a1)++(2.7,0) node (a2) {\FGCatg}
;
\draw[1cell=.9]
(a1) edge[transform canvas={yshift=.5ex}] node {\str} (a2)
(a2) edge[transform canvas={yshift=-.4ex}] node {\Incj} (a1)
;
\end{tikzpicture}
\end{equation}
whose right 2-adjoint is the inclusion 2-functor $\Incj$ \cref{fgcatg_psfgcatg}.
\end{description}
GMMO $K$-theory extends to a nonsymmetric multifunctor \pcref{rk:gmmo_mult} and satisfies a multiplicative equivariant version of the Barratt-Priddy-Quillen theorem.  The last property is used in the work \cite{gm-presheaves} to express equivariant spectra as spectral Mackey functors.  This work does not use those properties of GMMO $K$-theory.

\subsection*{Comparison Diagram: Shimakawa and GMMO}
The diagram \cref{kgmmo_shi_intro} displays the homotopical Shimakawa strong $K$-theory $\Khshosg$ \cref{khshosg_intro} along the top boundary, where $\algi$ is the inclusion functor, and GMMO $K$-theory $\Kgmmo$ \cref{intro_kgmmo} along the left-bottom-right boundary.
\begin{equation}\label{kgmmo_shi_intro}
\begin{tikzpicture}[vcenter]
\def\v{1.4} \def\h{2.8} \def\u{-1.2} \def\w{1}
\draw[0cell]
(0,0) node (a1) {\AlgstpsO}
(a1)++(\h,0) node (b) {\AlgpspsO}
(a1)++(0,-\v) node (a2) {\DGCatg}
(a2)++(0,.2) node (a2') {\phantom{\DGCatg}}
(a2)++(0,-\v) node (a3) {\FGCatgps}
(a3)++(\h,0) node (a4) {\FGCatg}
(a4)++(\h,0) node (a5) {\FGTopg}
(a5)++(0,\v) node (a6) {\FGTopg}
(a6)++(0,\v) node (a7) {\Gspec}
;
\draw[1cell=.9]
(a1) edge[right hook->] node {\algi} (b)
(b) edge node {\Khshosg} (a7)
(b) edge node {\Sgosg} (a4)
(a1) edge node[swap] {\Rg} (a2)
(a2) edge node[swap] {\gzest} (a3)
(a3) edge[transform canvas={yshift=-.3ex}] node[swap] {\str} (a4)
(a4) edge[transform canvas={yshift=.6ex}, left hook->] node[swap] {\Incj} (a3)
(a4) edge node {\clast} (a5)
(a5) edge node {\Bc} (a6)
(a6) edge node {\Kfg} (a7)
(a4) edge[bend right=10,shorten >=.2ex] node[swap] {\gxist} (a2')
(a1) [rounded corners=2pt] -| ($(a2)+(\u,0)$) -- node[pos=0,swap] {\Kgmmo} ($(a3)+(\u,0)$) |- ($(a5)+(0,-.6)$) -| node[pos=0,swap] {\phantom{x}} ($(a5)+(\w,0)$) |- (a7)
;
\draw[2cell=.9]
node[between=a1 and a2 at .55, shift={(.4*\h,0)}, rotate=180, 2label={below,\cgs}] {\Rightarrow}
;
\end{tikzpicture}
\end{equation}
Both equivariant $K$-theories $\Khshosg$ and $\Kgmmo$ end with the composite functor 
\[\FGCatg \fto{\Kfg\Bc\clast} \Gspec.\]
Thus, it suffices to compare their passages from $\Op$-algebras to $\FGG$-categories via the functors
\begin{equation}\label{hiszr_intro}
\begin{tikzpicture}[vcenter]
\draw[0cell]
(0,0) node (a1) {\AlgstpsO}
(a1)++(3,0) node (a2) {\FGCatg.}
;
\draw[1cell=.9]
(a1) edge[transform canvas={yshift=.6ex}] node {\Sgosg\algi} (a2)
(a1) edge[transform canvas={yshift=-.3ex}] node[swap] {\str\gzest\Rg} (a2)
;
\end{tikzpicture}
\end{equation}
However, it is not at all obvious how the two functors in \cref{hiszr_intro} are related.  One main difficulty is that the morphisms of the category $\FGCatg$---namely, $G$-natural transformations---are too rigid for an explicit comparison between $\Khshosg$ and $\Kgmmo$.  Instead, the comparison happens in the 2-category $\DGCatg$, whose 1-cells are $\Pig$-strict $\Gcatst$-pseudotransformations.  The construction of the comparison $\cgs$ \pcref{def:cgs} makes full use of the flexibility offered by enriched pseudotransformations; see \cref{cgsax}.

\subsection*{Comparison Categorical Weak $G$-Equivalences}
The main comparison between $\Khshosg$ and $\Kgmmo$ is the 2-natural transformation \pcref{sec:cgs} 
\begin{equation}\label{cgs_iinat_intro}
\begin{tikzpicture}[vcenter]
\def\t{25}
\draw[0cell]
(0,0) node (a1) {\phantom{A}}
(a1)++(2,0) node (a2) {\phantom{A}}
(a1)++(-.5,0) node (a1') {\AlgstpsO}
(a2)++(.7,0) node (a2') {\DGCatg}
;
\draw[1cell=.9]
(a1) edge[bend left=\t] node {\gxist\Sgosg\algi} (a2)
(a1) edge[bend right=\t] node[swap] {\Rg} (a2)
;
\draw[2cell=.9]
node[between=a1 and a2 at .42, rotate=-90, 2label={above,\cgs}] {\Rightarrow}
;
\end{tikzpicture}
\end{equation}
that appears in the left region in \cref{kgmmo_shi_intro}.  It compares their first steps, namely, Shimakawa strong $H$-theory $\Sgosg$ and $\Rg$.  The 2-functor $\gxist$ is the pullback along the $\Gcatst$-functor $\gxi \cn \DG \to \FG$ that comes with the $\Gcatst$-category $\DG$ \cref{dgcoop}.  For a chaotic $\Einfg$-operad of the form 
\[\Oph = \Catg(\EG,\Op)\] 
with $\EG$ the translation category of $G$ and an $\Oph$-algebra 
\[\Ah = \Catg(\EG,\A)\] 
with $\A$ an $\Op$-algebra, the comparison morphism \cref{cgs_ah_wgeq} 
\[\gxist\Sgohsg\algi\Ah \fto{\cgs_{\Ah}} \Rg\Ah \inspace \DGCatg\] 
is componentwise a categorical weak $G$-equivalence \pcref{def:cat_weakg}.  This means that, for each pointed finite $G$-set $\mal$, the $\mal$-component pointed $G$-functor
\begin{equation}\label{intro_cgsah_mal}
(\gxist\Sgohsg\algi\Ah)\mal \fto{\cgs_{\Ah,\mal}} (\Rg\Ah)\mal
\end{equation}
is a categorical weak $G$-equivalence.  In other words, for each subgroup $H \subseteq G$, the $H$-fixed subfunctor
\[\big((\gxist\Sgohsg\algi\Ah)\mal\big)^H 
\fto[\sim]{(\cgs_{\Ah,\mal})^H} 
\big((\Rg\Ah)\mal\big)^H\] 
is an equivalence of categories.  This nontrivial fact uses \cref{zdsg_catweakg}, which involves results from \crefrange{ch:hgoprod}{ch:special}.

\subsection*{Comparison Weak $G$-Equivalences}
Using the counit of the 2-adjunction $(\str,\Incj)$ and the comparison $\cgs$, \cref{thm:gmmo_shi} establishes a natural zigzag 
\[\Khshosg\algi \fot[]{} \bdot \fto{} \Kgmmo\]
between the homotopical Shimakawa strong $K$-theory $\Khshosg$ and GMMO $K$-theory $\Kgmmo$.  Applying this zigzag to the input pair $(\Oph,\Ah)$ yields componentwise weak $G$-equivalences between orthogonal $G$-spectra
\[\Khshohsg\algi\Ah \fot[]{} \bdot \fto{} \Kgmmoh\Ah.\]
Thus,  for $\Oph$-algebras of the form $\Ah$, the homotopical Shimakawa strong $K$-theory and GMMO $K$-theory produce weakly $G$-equivalent orthogonal $G$-spectra.  

For example, consider the Barratt-Eccles $\Gcat$-operad $\BE$ and the $G$-Barratt-Eccles operad $\GBE = \Catg(\EG,\BE)$ \pcref{def:BE,def:GBE}.  Their algebras are, respectively, naive and genuine permutative $G$-categories.  By \cref{thm:gmmo_shi}, for each genuine permutative $G$-category of the form $\Ah$ with $\A$ a naive permutative $G$-category, there is a zigzag between orthogonal $G$-spectra 
\[\Khshgbesg\algi\Ah \fot[]{} \bdot \fto{} \Kgmmogbe\Ah\]
consisting of componentwise weak $G$-equivalences \pcref{ex:gmmo_shi_gbe}.

\section{GMMO and Schwede Global $K$-Theories}
\label{sec:intro_gmmo_gl}

This section summarizes \cref{ch:kgl_gmmo} on the equivalence between GMMO $K$-theory \cref{intro_kgmmo} and Schwede global $K$-theory \cref{intro_ksc} for each finite group $G$.  The first half of this section briefly summarizes Schwede global $K$-theory, which is reviewed in detail in \cref{ch:kgl}.  The key comparison diagram is \cref{cglg_intro}.

\subsection*{Schwede Global $K$-Theory}
Each of the equivariant $K$-theories---Shimakawa, multifunctorial, and GMMO---summarized so far is defined for an arbitrary but fixed finite group $G$.  In contrast, Schwede global $K$-theory $\Ksc$ \cite{schwede_global}, defined as the following composite functor, keeps track of $G$-actions for all finite groups $G$.
\begin{equation}\label{intro_ksc}
\begin{tikzpicture}[vcenter]
\def\h{2.5} \def\v{1.4}
\draw[0cell]
(0,0) node (a1) {\Parcat}
(a1)++(0,-\v) node (a2) {\FMCat}
(a2)++(\h,0) node (a3) {\FICat}
(a3)++(\h,0) node (a4) {\phantom{\FITop}}
(a4)++(0,-.035) node (a4') {\FITop}
(a4)++(0,\v) node (a5) {\Sptop}
;
\draw[1cell=.9]
(a1) edge node {\Ksc} (a5)
(a1) edge node[swap] {\Jsc} (a2) 
(a2) edge node {\uprst} (a3)
(a3) edge node {\clast} (a4)
(a4') edge node[swap] {\Kfi} (a5)
;
\end{tikzpicture}
\end{equation}
The domain category $\Parcat$ of $\Ksc$ consists of parsummable categories and parsummable functors \pcref{def:mcat,def:bxtimes,def:parcat}.  A parsummable category is a commutative monoid in the symmetric monoidal category of tame $\schm$-categories, where $\schm$ is the translation category of the monoid of injections 
\[\{0,1,2,\ldots\} = \ome \to \ome.\]  
In the intermediate categories, $\Fsk$ is the category of pointed finite sets and pointed morphisms, and $\bdI$ is the category of finite sets and injections.  The intermediate diagram categories---$\FMCat$, $\FICat$, and $\FITop$---are defined in \cref{def:fmcat,def:ficat}.  The codomain category $\Sptop$ of global $K$-theory consists of symmetric spectra based on topological spaces, with structure morphisms indexed by $\bdI$ \pcref{def:sptop}. 

For a finite group $G$, each symmetric spectrum $X$ yields a $G$-symmetric spectrum $X_G$, called the underlying $G$-spectrum of $X$, with the trivial $G$-action.  Thus, the global $K$-theory $\Ksc\C$ of a parsummable category $\C$ yields a $G$-symmetric spectrum $(\Ksc\C)_G$ for each finite group $G$.  Note that a $G$-symmetric spectrum with the trivial $G$-action, such as $(\Ksc\C)_G$, can still be $G$-stably equivalent to a $G$-symmetric spectrum with a nontrivial $G$-action.  

The comparison of GMMO and global $K$-theories consists of two main steps.  First, we replace global $K$-theory for each finite group $G$ with what we call \emph{Schwede $K$-theory} \cref{kscgb_functor}.  Then we prove that Schwede $K$-theory is equivalent to GMMO $K$-theory.  The rest of this section briefly summarizes these two steps.

\subsection*{Schwede $K$-Theory}
For a finite group $G$, we define Schwede $K$-theory \cref{kscgb_functor} as the composite functor
\begin{equation}\label{kscgb_fun_intro}
\begin{tikzpicture}[vcenter]
\def\h{2.3} \def\g{1.6} \def\u{.7} \def\d{.05}
\draw[0cell]
(0,0) node (a1) {\phantom{\Parcat}}
(a1)++(0,.02) node (a1') {\Parcat}
(a1)++(\h,0) node (a2) {\FGCat}
(a2)++(\h,0) node (a3) {\FGTop}
(a3)++(\h,0) node (a4) {\Sptopg}
;
\draw[1cell=.9]
(a1) edge node {\Jscg} (a2)
(a2) edge node {\clast} (a3)
(a3) edge node {\Kfgsi} (a4)
(a1') [rounded corners=2pt] |- ($(a2)+(0,\u)$) -- node {\Kscgb} ($(a3)+(0,\u)$) -| (a4)
;
\end{tikzpicture}
\end{equation}
that sends parsummable categories to $G$-symmetric spectra with generally nontrivial $G$-actions.  The first functor $\Jscg$ \cref{Jscg}, called Schwede $J$-theory, sends parsummable categories to $\Fskg$-categories.  It is analogous to the first functor $\Jsc$ in global $K$-theory $\Ksc$ \cref{intro_ksc}, but it takes the finite group $G$ into account.  By Schwede's \cref{sch4.15} and \cref{Isc_geq}, for each parsummable category $\C$, there are natural $G$-stable equivalences between $G$-symmetric spectra
\begin{equation}\label{IscC_Gequiv_intro}
(\Ksc\C)_G \fto{\ascgc} \gacomegs \fot{\bscgc} \kcg \fto{\Isc^\C} \Kscgb\C
\end{equation}
that connect the underlying $G$-spectrum $(\Ksc\C)_G$ of the global $K$-theory of $\C$ and the Schwede $K$-theory $\Kscgb\C$ of $\C$.  Note that the group $G$ acts trivially on $(\Ksc\C)_G$ but nontrivially on $\Kscgb\C$ in general.  The intermediate $G$-symmetric spectra $\kcg$ and $\gacomegs$ are defined in \cref{def:sch4.15,expl:sch3.13}.

\subsection*{$\Einfg$-Categories from Parsummable Categories}
In order to compare GMMO and global $K$-theories, it is necessary to reconcile their different input data.  For this purpose, we consider the categorical injection $G$-operad $\Iopg$ \pcref{def:global_einf}.  It is a chaotic $\Einfg$-operad \pcref{iopg_einf}, so its algebras are $\Einfg$-categories.  There is a functor \cref{Ig_functor}
\[\Parcat \fto{\Ig} \Algiopg\]
that sends parsummable categories, the input data for global $K$-theory, to $\Iopg$-algebras, the input data for GMMO $K$-theory.  \cref{sec:iopg} discusses the chaotic $\Einfg$-operad $\Iopg$ and the functor $\Ig$ in detail.

\subsection*{Comparison Diagram: GMMO and Global}
Recall that global $K$-theory at a finite group $G$ is naturally $G$-stably equivalent to Schwede $K$-theory \cref{IscC_Gequiv_intro}.  To compare global and GMMO $K$-theories, it suffices to compare Schwede and GMMO $K$-theories as follows.
\begin{equation}\label{cglg_intro}

\end{equation}
The top-right composite $\Kscgb$ is Schwede $K$-theory \cref{kscgb_fun_intro}.  The left-bottom-right composite $\Kiopgb$ is GMMO $K$-theory for the $\Einfg$-operad $\Iopg$, restricted to $\Algiopg$ in the domain and prolonged to the category $\Sptopg$ of $G$-symmetric spectra in the codomain.  The functor $\gxist$ is the pullback along the $\Gcatst$-functor \cref{gxist} 
\[\DG \fto{\gxi} \FG.\]  
Each instance of $\Lg$ is an equivalence of categories \pcref{thm:fgcat_fgcatg_iieq,fgtop_fgtopg_eq}, and $\retn \cn \Bc \to 1$ is the retraction for the bar functor \cref{retn_barFG_id}.  

The category $\FGCatg$ is not convenient to compare Schwede and GMMO $K$-theories because its morphisms, which are $G$-natural transformations, are too rigid.  Instead, the comparison happens in the category $\DGCatg$, whose morphisms are $\Pig$-strict $\Gcatst$-pseudotransformations \pcref{def:dgcatg_iicat}.  This aspect of the comparison is conceptually similar to the comparison between Shimakawa and GMMO $K$-theories \cref{kgmmo_shi_intro}.

\subsection*{Comparison $G$-Equivalences}
The main comparison between Schwede and GMMO $K$-theories is the natural transformation
\cref{cglg_nat} 
\begin{equation}\label{cglg_comp_intro}
\begin{tikzpicture}[vcenter]
\def\t{25}
\draw[0cell]
(0,0) node (a1) {\phantom{A}}
(a1)++(2,0) node (a2) {\phantom{A}}
(a1)++(-.38,0) node (a1') {\Parcat}
(a2)++(.7,-.02) node (a2') {\DGCatg}
;
\draw[1cell=.9]
(a1) edge[bend left=\t] node {\gxist\Lg\Jscg} (a2)
(a1) edge[bend right=\t] node[swap] {\Rg\Ig} (a2)
;
\draw[2cell=.9]
node[between=a1 and a2 at .42, rotate=-90, 2label={above,\cglg}] {\Rightarrow}
;
\end{tikzpicture}
\end{equation}
that compares their first steps, namely, Schwede $J$-theory $\Jscg$ \cref{Jscg} and the functor $\Rg$ \cref{rg_twofunctor}.  \cref{cglgc_equiv} proves that, for each finite group $G$ and parsummable category $\C$, the comparison morphism
\[\gxist\Lg\Jscg\C \fto{\cglgc} \Rg\Ig\C \inspace \DGCatg\]
is componentwise an adjoint pointed $G$-equivalence between pointed $G$-categories \pcref{def:adjointGeq}.  This means that, for each pointed finite $G$-set $\nbeta$, the $\nbeta$-component pointed $G$-functor
\begin{equation}\label{cglgcnbeta_intro}
(\gxist\Lg\Jscg\C)\nbeta \fto{\cglgcnbeta} (\Rg\Ig\C)\nbeta
\end{equation}
admits a pointed $G$-functor going backward, along with unit and counit pointed $G$-natural isomorphisms \pcref{def:cgglc,def:cglgu,def:cglgv} that satisfy the triangle identities for an adjoint $G$-equivalence.  \cref{gl_gmmo_geq} explains how \cref{IscC_Gequiv_intro,cglg_intro,cglgcnbeta_intro} yield a zigzag of natural $G$-stable equivalences
\begin{equation}\label{kglgmmo_zigzag}
(\Ksc\C)_G \to \Cdots \to \Kscgb\C \leftarrow \Cdots \leftarrow \Kiopgb\Ig\C
\end{equation}
between the underlying $G$-spectrum $(\Ksc\C)_G$ of the global $K$-theory of $\C$ and the GMMO $K$-theory $\Kiopgb\Ig\C$ of the $\Einfg$-category $\Ig\C$.  \cref{expl:gl_gmmo_geq} explicitly describes the components of the zigzag \cref{kglgmmo_zigzag}.

\subsection*{Lenz's Comparison of Shimakawa and Global $K$-Theories}
Lenz's work \cite[Cor.\! 4.1.44]{lenz-global} compares Shimakawa $K$-theory and global $K$-theory.  The context and method of Lenz's comparison are substantially different from this work.  Lenz's comparison starts with small permutative categories instead of parsummable categories.  It involves $G$-global algebraic $K$-theory, quasi-categories, and Schwede's constructions $(-)^{\mathrm{sat}}$ and $\Phi$ in \cite[7.23 and 11.1]{schwede_global}.  Our direct comparison of GMMO and global $K$-theories \pcref{gl_gmmo_geq} is \emph{not} a combination of Lenz's comparison and our \cref{thm:gmmo_shi}, which directly compares Shimakawa and GMMO $K$-theories.  Since this work does not use Lenz's result in any way, it will not be discussed any further.

\section{Chapter Summaries}
\label{sec:intro_summ}

This section briefly summarizes each chapter.  In the main text, each chapter and each section start with a summary.

\prefacepartNumName{part:kgo_shi}
This part (\crefrange{ch:kgo}{ch:hgo}) discusses the underlying functor \cref{kgosg_intro} 
\begin{equation}\label{kgo_chsum}
\begin{tikzpicture}[vcenter]
\def\h{2.2} \def\u{.6}
\draw[0cell=.9]
(0,0) node (a1) {\AlgpspsO}
(a1)++(\h,0) node (a2) {\GGCatii}
(a2)++(\h,0) node (a3) {\GGTopii}
(a3)++(.85*\h,0) node (a4) {\Gspec} 
;
\draw[1cell=.85]
(a1) edge node {\Jgosg} (a2)
(a2) edge node {\clast} (a3)
(a3) edge node {\Kg} (a4)
(a1) [rounded corners=2pt, shorten <=-.2ex] |- ($(a2)+(0,\u)$) -- node {\Kgosg} ($(a3)+(0,\u)$) -| (a4)
;
\end{tikzpicture}
\end{equation}
of the author's enriched multifunctorial equivariant $K$-theory for a $\Tinf$-operad $\Op$ and proves that it is naturally isomorphic to the functor \cref{khgosg_intro}
\begin{equation}\label{khgosg_chsum}
\begin{tikzpicture}[vcenter]
\def\h{2.2} \def\u{.6}
\draw[0cell=.9]
(0,0) node (a1) {\AlgpspsO}
(a1)++(\h,0) node (a2) {\GGCatg}
(a2)++(\h,0) node (a3) {\GGTopg}
(a3)++(.9*\h,0) node (a4) {\Gspec.}
;
\draw[1cell=.85]
(a1) edge node {\Hgosg} (a2)
(a2) edge node {\clast} (a3)
(a3) edge node {\Kgg} (a4)
(a1) [rounded corners=2pt, shorten <=-.2ex] |- ($(a2)+(0,\u)$) -- node {\Khgosg} ($(a3)+(0,\u)$) -| (a4)
;
\end{tikzpicture}
\end{equation}

\prefacechapNumName{ch:kgo} 
To facilitate the comparison with Shimakawa $K$-theory, this chapter reviews the underlying functor $\Kgosg$ \cref{kgosg_intro} of our enriched multifunctorial equivariant $K$-theory.  This material is adapted from \cite{yau-eqk} in a self-contained manner, having all the necessary definitions in full.

\prefacechapNumName{ch:ggcatg}
This chapter constructs the 2-category $\GGCatg$ and the 2-equivalence
\begin{equation}\label{Ligst_chsum}
\begin{tikzpicture}[vcenter]
\draw[0cell]
(0,0) node (a1) {\GGCatii}
(a1)++(2.5,0) node (a2) {\GGCatg}
;
\draw[1cell=.9]
(a1) edge[transform canvas={yshift=.5ex}] node {\Lg} (a2)
(a2) edge[transform canvas={yshift=-.4ex}] node {\igst} (a1)
;
\end{tikzpicture}
\end{equation}
induced by the full subcategory inclusion $\ig \cn \Gsk \to \GG$.  The 2-categories $\GGCatii$ and $\GGCatg$ are intermediate stops for, respectively, the functors $\Kgosg$ and $\Khgosg$.

\prefacechapNumName{ch:hgo}
This chapter constructs the functor $\Khgosg$ \cref{khgosg_intro}, which involves the category $\GGCatg$, and proves that it is naturally isomorphic to $\Kgosg$.  The main feature of $\Khgosg$ is that it can be effectively compared with Shimakawa $K$-theory, which involves the category $\FGCatg$.

\prefacepartNumName{part:shim_k}
This part (\crefrange{ch:shih}{ch:shim_top}) constructs the operadic extension of Shimakawa equivariant $K$-theory and topologically compares that with the last step $\Kgg$ of the functor $\Khgosg$.

\prefacechapNumName{ch:shih}
This chapter constructs Shimakawa strong $H$-theory 2-functor
\[\AlgpspsO \fto{\Sgosg} \FGCatg\]
for a group $G$ and a 1-connected $\Gcat$-operad $\Op$.  For the Barratt-Eccles operad $\BE$, $\Sgosgbe$ recovers Shimakawa's original construction in \cite[pp.\ 251--252]{shimakawa89}.

\prefacechapNumName{ch:shimakawa_K}
This chapter constructs Shimakawa strong $K$-theory \cref{kshosg_intro} 
\begin{equation}\label{kshosg_chsum}
\begin{tikzpicture}[vcenter]
\def\h{2} \def\u{.6}
\draw[0cell=.9]
(0,0) node (a1) {\AlgpspsO}
(a1)++(1.1*\h,0) node (a2) {\FGCatg}
(a2)++(\h,0) node (a3) {\FGTopg}
(a3)++(\h,0) node (a4) {\Gspec} 
;
\draw[1cell=.85]
(a1) edge node {\Sgosg} (a2)
(a2) edge node {\clast} (a3)
(a3) edge node {\Kfg} (a4)
(a1) [rounded corners=2pt, shorten <=-.2ex] |- ($(a2)+(0,\u)$) -- node {\Kshosg} ($(a3)+(0,\u)$) -| (a4)
;
\end{tikzpicture}
\end{equation}
and its homotopical variant \cref{khshosg_intro}
\begin{equation}\label{khshosg_chsum}
\begin{tikzpicture}[vcenter]
\def\h{2} \def\u{.6}
\draw[0cell=.9]
(0,0) node (a1) {\phantom{\AlgpspsO}}
(a1)++(0,0) node (a1') {\AlgpspsO}
(a1)++(1.1*\h,0) node (a2) {\FGCatg}
(a2)++(\h,0) node (a3) {\FGTopg}
(a3)++(\h,0) node (a4) {\FGTopg}
(a4)++(\h,0) node (a5) {\Gspec.} 
;
\draw[1cell=.85]
(a1) edge node {\Sgosg} (a2)
(a2) edge node {\clast} (a3)
(a3) edge node {\Bc} (a4)
(a4) edge node {\Kfg} (a5)
(a1) [rounded corners=2pt, shorten <=-.2ex] |- ($(a2)+(0,\u)$) -- node {\Khshosg} ($(a4)+(0,\u)$) -| (a5)
;
\end{tikzpicture}
\end{equation}
For a finite group $G$, these two functors produce naturally componentwise weakly $G$-equivalent orthogonal $G$-spectra \pcref{thm:ksho_khsho}.  For the Barratt-Eccles operad $\BE$, $\Khshbesg$ computes Shimakawa's original equivariant $K$-theory functor in \cite{shimakawa89} up to a natural componentwise weak $G$-equivalence \pcref{expl:khsho_shi,}

\prefacechapNumName{ch:shim_top}
This chapter compares the topological parts of Shimakawa $K$-theory and multifunctorial equivariant $K$-theory, namely, the functors
\begin{equation}\label{topeq_chsum}
\begin{tikzpicture}[vcenter]
\def\h{2}
\draw[0cell]
(0,0) node (a1) {\FGTopg}
(a1)++(2.5,0) node (a2) {\Gspec}
(a2)++(1.3,0) node (and) {\text{and}}
(and)++(1.6,0) node (b1) {\GGTopg}
(b1)++(2.3,0) node (b2) {\Gspec.}
;
\draw[1cell=.9]
(a1) edge[transform canvas={yshift=.5ex}] node {\Kfg} (a2)
(a1) edge[transform canvas={yshift=-.4ex}] node[swap] {\Kfg\Bc} (a2)
(b1) edge node {\Kgg} (b2)
;
\end{tikzpicture}
\end{equation}
The functors $\Kfg$ and $\Kgg$ factor through each other up to natural isomorphisms \pcref{kfgkgg_compare}.  Thus, any orthogonal $G$-spectrum produced by one of these two functors is also produced by the other functor.  For a finite group $G$ and proper objects, the functors $\Kfg\Bc$ and $\Kgg$ factor through each other up to natural componentwise weak $G$-equivalences \pcref{kfgbkgg_compare,kggkfgb_compare}.

\prefacepartNumName{part:kgo_shi_comp}
This part (\crefrange{ch:h_comparison}{ch:special}) compares the categorical parts of Shimakawa $K$-theory and multifunctorial equivariant $K$-theory, namely, Shimakawa strong $H$-theory \pcref{ch:shih} 
\[\AlgpspsO \fto{\Sgosg} \FGCatg\]
and strong $H$-theory \pcref{ch:hgo}
\[\AlgpspsO \fto{\Hgosg} \GGCatg.\]

\prefacechapNumName{ch:h_comparison}
This chapter first proves that Shimakawa strong $H$-theory $\Sgosg$ factors through strong $H$-theory as $\ifgst\Hgosg$, where $\ifgst$ is the pullback along the inclusion functor $\ifg \cn \FG \to \GG$.  Thus, for each $\Op$-pseudoalgebra $\A$, the $\FGG$-category $\Sgosg\A$ is the restriction to length-1 objects of the $\GGG$-category $\Hgosg\A$.  The remainder of \cref{part:kgo_shi_comp} studies the reverse comparison, computing $\Hgosg$ in terms of $\Sgosg$.  The rest of this chapter constructs the strong $H$-theory comparison 2-natural transformation \cref{pistsg_intro}
\[\smast\Sgosg \fto{\Pistsg} \Hgosg\]
and a commutative diagram of pointed $G$-functors \cref{pistsganbe_intro}
\begin{equation}\label{pist_chsum}
\begin{tikzpicture}[vcenter]
\def\h{3.5}
\draw[0cell]
(0,0) node (a11) {(\smast\Sgosg\A)\nbe}
(a11)++(\h,0) node (a12) {(\Hgosg\A)\nbe}
(a11)++(\h/2,-1.1) node (a2) {\phantom{\proAnbe}}
(a2)++(0,.15) node (a2') {\proAnbe}
;
\draw[1cell=.9]
(a11) edge node {\Pistsg_{\A,\angordnbe}} (a12)
(a11) edge node[swap] {\zbsg} (a2)
(a12) edge node {\zdsg} (a2)
;
\end{tikzpicture}
\end{equation}
for a group $G$, an $\Op$-pseudoalgebra $\A$, and an object $\angordnbe \in \GG$.  The pointed $G$-category $\proAnbe$ is the $\nbe$-twisted product of the pointed $G$-category $\A$ \pcref{def:proCnbe}, taking into account the $G$-actions on the entries of $\nbe$.  \cref{sec:twprod_fixed} computes the $G$-fixed point subcategory of any $\nbe$-twisted product.  The main result about the comparison $\Pistsg$ is \cref{thm:pistweakgeq}, whose proof involves most of \crefrange{ch:sgoprod}{ch:compgen}.

\prefacechapNumName{ch:sgoprod}
This chapter proves that, for a group $G$ and a $\Uinf$-operad $\Op$ \pcref{as:OpA'}, the $G$-functor 
\[(\smast\Sgosg\A)\nbe \fto{\zbsg} \proAnbe\]
constructed in \cref{ch:h_comparison} is the right adjoint of an adjoint equivalence \pcref{thm:zbsg_eq}.  While the functor $\zbsg$ is $G$-equivariant, its left adjoint inverse $\zbsgad$ is only \emph{pseudo} $G$-equivariant \pcref{thm:zbsgad_pseudo}.  This means that $\zbsgad$ commutes with the $G$-actions of its domain and codomain up to coherent natural isomorphisms.  For a nontrivial group $G$, the left adjoint $\zbsgad$ is not $G$-equivariant even for the Barratt-Eccles operad $\BE$, on which $G$ acts trivially \pcref{ex:al_g_s}.

\prefacechapNumName{ch:hgoprod}
This chapter is the analogue of \cref{ch:sgoprod} for the pointed $G$-functor 
\[(\Hgosg\A)\nbe \fto{\zdsg} \proAnbe\]
constructed in \cref{ch:h_comparison}.  \cref{thm:zdsg_eq} proves that $\zdsg$ is the right adjoint of an adjoint equivalence.  Its left adjoint inverse $\zdsgad$ is only pseudo $G$-equivariant \pcref{thm:zdsgad_pseudo}.  The fact that $\zbsg$ and $\zdsg$ are both equivalences implies that the strong $H$-theory comparison $\Pistsg_{\A,\angordnbe}$ \cref{pistsganbe_intro} is also an equivalence of categories in the nonequivariant sense \pcref{thm:PistAequivalence}.

\prefacechapNumName{ch:compgen} 
This chapter proves the main result about the strong $H$-theory comparison \pcref{thm:pistweakgeq}: The pointed $G$-functor \cref{pistsgahnbe_intro} 
\[(\smast\Sgohsg\Ah)\angordnbe \fto{\Pistsg_{\Ah,\angordnbe}} (\Hgohsg\Ah)\angordnbe\] 
is a categorical weak $G$-equivalence for the $\Uinf$-operad $\Oph = \Catg(\EG,\Op)$ with $\EG$ the translation category of $G$, an $\Oph$-pseudoalgebra of the form $\Ah = \Catg(\EG,\A)$, and an object $\nbe \in \GG$.  Therefore, up to componentwise categorical weak $G$-equivalences, the strong $H$-theory $\Hgohsg\Ah$ of $\Ah$ is computed by the Shimakawa strong $H$-theory $\smast\Sgohsg\Ah$ of $\Ah$.  The basic strategy of the proof of \cref{thm:pistweakgeq} is a 2-out-of-3 argument outlined in the introduction of \cref{ch:compgen}.

\prefacechapNumName{ch:special}
This chapter discusses special objects and weak $G$-equivalences in the categories $\GGCatg$ and $\GGCatii$.  There are two main observations.  First, the domain and codomain of the strong $H$-theory comparison $\Pistsg_{\Ah}$ \cref{pistsgah_intro} are special $\GGG$-categories, meaning that their Segal functors are categorical weak $G$-equivalences \pcref{thm:h_special,cor:h_special_shi}.  Thus, $\Pistsg_{\Ah}$ is a weak $G$-equivalence between special objects in $\GGCatg$.  Second, for a finite group $G$, each of the 2-equivalences $(\Lg,\igst)$ in \cref{ch:ggcatg} preserves and reflects special objects and weak $G$-equivalences \pcref{ggcat_sp_LX,gggcat_sp_iX,Lg_weq,igst_weq}.

\prefacepartNumName{part:shi_gmmo_gl}
This part (\crefrange{ch:kgmmo}{ch:kgl_gmmo}) proves the equivalences between Shimakawa $K$-theory, GMMO $K$-theory, and Schwede global $K$-theory for each finite group $G$.

\prefacechapNumName{ch:kgmmo}
For a finite group $G$ and a chaotic $\Einfg$-operad $\Op$ in $\Gcat$, this chapter constructs a natural zigzag 
\[\Khshosg\algi \fot[]{} \bdot \fto{} \Kgmmo\]
between the homotopical Shimakawa strong $K$-theory $\Khshosg$ \cref{khshosg_intro} and GMMO $K$-theory $\Kgmmo$ \cref{intro_kgmmo}.  For a chaotic $\Einfg$-operad of the form $\Oph = \Catg(\EG,\Op)$ with $\EG$ the translation category of $G$ and an $\Oph$-algebra $\Ah = \Catg(\EG,\A)$ with $\A$ an $\Op$-algebra, the zigzag comparison between $\Khshohsg\algi\Ah$ and $\Kgmmoh\Ah$ consists of componentwise weak $G$-equivalences of orthogonal $G$-spectra.  See \cref{thm:gmmo_shi}.  Applying this result to the $G$-Barratt-Eccles operad $\GBE$, for a naive permutative $G$-category $\A$ and the genuine permutative $G$-category $\Ah$, Shimakawa and GMMO $K$-theories produce weakly $G$-equivalent orthogonal $G$-spectra $\Khshgbesg\algi\Ah$ and $\Kgmmogbe\Ah$.  \crefrange{sec:dgcatg}{sec:kgmmo} review GMMO $K$-theory in detail.  \cref{sec:cgs,sec:gmmo_shi_weq} compare $\Khshosg$ and $\Kgmmo$.

\prefacechapNumName{ch:kgl} 
To prepare for the comparison with GMMO $K$-theory in \cref{ch:kgl_gmmo}, this chapter reviews Schwede global $K$-theory \cref{intro_ksc}
\[\Parcat \fto{\Ksc} \Sptop.\]
Denote by $\schm$ the translation category of the monoid of injections 
\[\{0,1,2,\ldots\} = \ome \to \ome.\]
The input data for global $K$-theory, called parsummable categories, are commutative monoids in the symmetric monoidal category of tame $\schm$-categories.  Global $K$-theory produces topological symmetric spectra indexed by finite sets and injections.  For each parsummable category $\C$ and finite group $G$, the symmetric spectrum $\Ksc\C$ yields a $G$-symmetric spectrum $(\Ksc\C)_G$, called the underlying $G$-spectrum, with the trivial $G$-action.

\prefacechapNumName{ch:kgl_gmmo} 
This chapter proves that, for each finite group $G$ and parsummable category $\C$, there is a natural zigzag of $G$-stable equivalences
\[(\Ksc\C)_G \to \Cdots \to \Kscgb\C \leftarrow \Cdots \leftarrow \Kiopgb\Ig\C\]
connecting the underlying $G$-spectrum $(\Ksc\C)_G$ of the global $K$-theory of $\C$ and the GMMO $K$-theory $\Kiopgb\Ig\C$ of the $\Einfg$-category $\Ig\C$ associated to $\C$.  See \cref{gl_gmmo_geq}.  The $\Einfg$-category $\Ig\C$ and GMMO $K$-theory $\Kiopgb$ are defined using the categorical injection $G$-operad $\Iopg$ \cref{iopgn}, which is a chaotic $\Einfg$-operad.  \cref{gl_gmmo_geq} is independent of \cref{thm:gmmo_shi}, which compares Shimakawa $K$-theory with GMMO $K$-theory.

\prefacechapNumName{ch:mon_cat}
To make this work as self-contained as possible, this appendix reviews basic concepts of category theory, monoidal categories, enriched categories, 2-categories, and enriched operads.

%% file: chap/notation.tex
\def\s{.85} \def\v{-.08ex} \def\hor{\kern .1em}
\newcommand{\blob}{\>}
\newcommand{\entry}[3]{\scalebox{\s}{#1} \hor : \scalebox{\s}{#3}, \scalebox{\s}{\pageref{#2}} \blob\\[\v]}
\newcommand{\newchHeader}[1]{\blob\\ \textbf{\Cref{#1}}. \nameref{#1}\blob\\}

\begin{tabbing}
\phantom{\textbf{Notation}} \= \kill

\newchHeader{ch:kgo} 
\entry{$(\Op,\ga,\opu$)}{def:intrinsic_pairing}{operad}
\entry{$\intr$, $\intr_{j,k}$}{intr_jk}{intrinsic pairing}
\entry{$\Delta^j$, $\Delta$}{intr_jk}{diagonal}
\entry{$\As$}{ex:as_intrinsic}{associative operad}
\entry{$\phi\ang{k_1,\ldots,k_n}$}{as_gamma}{block permutation}
\entry{$\phi_1 \oplus \cdots \oplus \phi_n$}{as_gamma}{block sum}
\entry{$\ufs{n}$}{ufsn}{unpointed finite set $\{1,2,\ldots,n\}$}
\entry{$\txsum_{j \in \ufs{n}}$}{ufsn}{sum $\txsum_{j=1}^n$}
\entry{$\twist_{j,k}$}{eq:transpose_perm}{transpose permutation}
\entry{$\lambda_{j,k}$}{lex_bijection}{lexicographic ordering}
\entry{$\Cat$}{def:GCat}{2-category of small categories}
\entry{$G$}{def:GCat}{a group, also regarded as a (2-)category with one object and morphisms $G$}
\entry{$\Gcat$}{not:Gcat}{2-category of small $G$-categories}
\entry{$\boldone$}{not:bone}{terminal category}
\entry{$\Catg(-,-)$}{def:Catg}{$G$-category of functors and natural transformations}
\entry{$(-)^G$}{catg_gfixed}{$G$-fixed subcategory}
\entry{$\pcom$, $\pcom_{j,k}$}{pseudocom_isos}{pseudo-commutativity isomorphism}
\entry{$\twist$}{pseudocom_isos}{swap permutation}
\entry{$\tn S$}{def:translation_cat}{translation category of $S$}
\entry{$[b,a] , ! \cn a \to b$}{not:b_to_a}{unique morphism from $a$ to $b$}
\entry{$\BE$}{def:BE}{Barratt-Eccles operad}
\entry{$\GBE$}{def:GBE}{$G$-Barratt-Eccles operad}
\entry{$[-,-]$}{def:GBE}{$G$-set of functions}
\entry{$(\A,\gaA,\phiA)$}{def:pseudoalgebra}{$\Op$-pseudoalgebra}
\entry{$\gaA_n$}{gaAn}{$n$-th $\Op$-action $G$-functor}
\entry{$\zero$, $\gaA_0(*)$}{pseudoalg_zero}{basepoint}
\entry{$\phiA$}{phiA}{associativity constraint}
\entry{$\Oph$}{ng_i}{$\Catg(\EG,\Op)$ for a $\Gcat$-operad $\Op$}
\entry{$\Ah$}{ng_ii}{$\Catg(\EG,\A)$ for an $\Op$-pseudoalgebra $\A$}
\entry{$(f,\actf)$}{def:laxmorphism}{lax $\Op$-morphism with action constraint $\actf$}
\entry{$\AlglaxO$, $\AlgpspsO$, $\AlgstO$}{oalgps_twocat}{2-categories of $\Op$-pseudoalgebras}
\entry{$(\A,\otimes,\zero,\alpha,\lambda,\rho,\xi)$}{def:naive_smGcat}{naive symmetric monoidal $G$-category}
\entry{$(f,f^2,f^0)$}{def:smGfunctor}{symmetric monoidal $G$-functor}
\entry{$\smgcat$, $\smgcatsg$, $\smgcatst$}{def:smGcat_twocat}{2-categories of naive symmetric monoidal $G$-categories}
\entry{$(S, \bp)$}{def:Fsk}{pointed set $S$ with basepoint $\bp$}
\entry{$\ordn$}{ordn}{pointed finite set $\{0 < 1 < \cdots < n\}$}
\entry{$\Fsk$}{not:Fsk}{category of $\ordn$ for $n \geq 0$ and pointed functions}
\entry{$(\sma, \ord{1}, \xi)$}{def:Fsk_permutative}{permutative structure on $\Fsk$}
\entry{$\sma$}{m-sma-n}{smash product}
\entry{$\vstar$}{not:Fsk_smash_zero}{initial-terminal basepoint of $\Fsk^{(0)}$}
\entry{$\ang{}$}{not:Fsk_smash_zero}{empty tuple}
\entry{$\Fsk^{(q)}$}{not:Fskq}{$q$-fold smash power of $\Fsk$}
\entry{$\ang{\ord{n}}$, $\ang{\ord{n}_i}_{i \in \ufs{q}}$}{angordn}{$q$-tuple $(\ord{n}_1, \ldots, \ord{n}_q)$}
\entry{$h_*$}{not:reind_func}{reindexing functor for an injection $h \cn \ufs{q} \to \ufsr$}
\entry{$\ord{n}_{\emptyset}$, $\psi_{\emptyset}$}{ordn_empty}{$\ord{1}$ and $1_{\ord{1}}$}
\entry{$\Gsk$}{def:Gsk}{category with finite tuples of pointed finite sets as objects}
\entry{$\vstar$}{Gsk_objects}{initial-terminal basepoint of $\Gsk$}
\entry{$\angordm$}{Gsk_morphisms}{object in $\Gsk$}
\entry{$(f, \ang{\psi})$}{fangpsi}{morphism in $\Gsk$}
\entry{$\ifg$}{ifg}{length-1 inclusion $\Fsk \to \Gsk$}
\entry{$(\oplus, \ang{}, \xi)$}{def:Gsk_permutative}{permutative structure on $\Gsk$}
\entry{$\tau_{p,q}$}{Gsk_braiding}{block permutation interchanging $p$ and $q$ elements}
\entry{$\sma$}{def:smashFskGsk}{smash functor $\Gsk \to \Fsk$}
\entry{$\Gcatst$}{Gcatst}{2-category of small pointed $G$-categories}
\entry{$\bonep$}{gcatst_unit}{discrete trivial $G$-category with two objects}
\entry{$\Catgst$}{catgst_cd}{internal hom for $\Gcatst$}
\entry{$\GGCatii$}{def:GGCat}{2-category of $\Gskg$-categories}
\entry{$\ang{s} \compk t$}{compk}{replace the $k$-th entry of $\ang{s}$ by $t$}
\entry{$\coprod_{j \in \ufs{p}}\, S_j$}{partition}{partition of a set $S$}
\entry{$(a,\glu)$}{nsystem}{$\angordn$-system}
\entry{$\ang{s}$, $\ang{s_j}_{j \in \ufs{q}}$}{marker}{marker $\ang{s_j \subseteq \ufs{n}_j}_{j \in \ufs{q}}$}
\entry{$a_{\ang{s}}$}{a_angs}{$\angs$-component object}
\entry{$\glu_{x; \ang{s} \csp k, \ang{s_{k,i}}_{i \in \ufs{r}}}$}{gluing-morphism}{gluing morphism}
\entry{$(\zero,1_\zero)$}{not:basesystem}{base $\angordn$-system}
\entry{$\theta_{\angs}$}{theta_angs}{$\angs$-component morphism}
\entry{$\Aangordn$}{not:Aangordn}{category of $\angordn$-systems in $\A$}
\entry{$\Aangordnsg$}{not:Aangordn}{category of strong $\angordn$-systems in $\A$}
\entry{$g \cdot (a,\glu)$, $(ga, g\glu)$}{nsystem_gaction}{$g$-action on $(a,\glu)$}
\entry{$g\theta$}{gtheta}{$g$-action on $\theta$}
\entry{$\ftil$}{ftil_functor}{pointed $G$-functor $\Aangordm \to \Afangordm$}
\entry{$\left(\atil, \glutil\right)$}{ftil_aglu}{$f_*\angordm$-system $\ftil(a,\glu)$}
\entry{$\ftil_*\ang{s}$}{ftil_angs}{marker $\ang{s_{f(i)} \subseteq \ufs{m}_i}_{i \in \ufs{p}}$}
\entry{$\thatil$}{thatil_component}{morphism of $f_*\angordm$-systems $\ftil(\theta)$}
\entry{$\psitil$}{psitil_functor}{pointed $G$-functor $\Afangordm \to \Aangordn$}
\entry{$(a^{\psitil}, \glu^{\psitil})$}{angs_ordtun}{$\angordn$-system $\psitil(a,\glu)$}
\entry{$\psiinv\ang{s}$}{apsitil_angs}{marker $\ang{\psiinv_j s_j \subseteq \ufs{m}_{\finv(j)}}_{j \in \ufs{q}}$}
\entry{$\tha^{\psitil}$}{thapsitil_component}{morphism of $\angordn$-systems $\psitil(\theta)$}
\entry{$\Aupom$, $\Afpsi$}{AF}{pointed $G$-functor $\Aangordm \to \Aangordn$}
\entry{$\Aupomsg$}{AF_sg}{strong variant of $\Aupom$}
\entry{$\Jgo\A$, $\Adash$}{A_ptfunctor}{$J$-theory of an $\Op$-pseudoalgebra $\A$}
\entry{$\Jgosg\A$, $\Asgdash$}{A_ptfunctor}{strong $J$-theory of an $\Op$-pseudoalgebra $\A$}
\entry{$\Jgo f$}{Jgo_f}{$J$-theory of a lax $\Op$-morphism $f$}
\entry{$\Jgosg f$}{Jgosg_f}{strong $J$-theory of an $\Op$-pseudomorphism $f$}
\entry{$\Jgo\theta$}{Jgotheta}{$J$-theory of an $\Op$-transformation $\theta$}
\entry{$\Jgosg\theta$}{Jgosgtheta}{strong $J$-theory of an $\Op$-transformation $\theta$}
\entry{$\Jgo$}{thm:Jgo_twofunctor}{$J$-theory 2-functor $\AlglaxO \to \GGCatii$}
\entry{$\Jgosg$}{thm:Jgo_twofunctor}{strong $J$-theory 2-functor $\AlgpspsO \to \GGCatii$}
\entry{$\Top$}{def:Gtop}{category of spaces and morphisms}
\entry{$\Gtop$}{not:Gtop}{category of $G$-spaces and $G$-morphisms}
\entry{$\Topg(-,-)$}{not:Topg}{$G$-space of morphisms}
\entry{$(-)^G$}{def:weakG_top}{$G$-fixed point space}
\entry{$\eqg$}{not:eqg}{weak $G$-equivalence}
\entry{$\Gtopst$}{Gtopst_smc}{category of pointed $G$-spaces and pointed $G$-morphisms}
\entry{$\stplus$}{not:stplus}{smash unit}
\entry{$\Topgst(-,-)$}{not:Topgst}{pointed $G$-space of pointed morphisms}
\entry{$\Topgst$}{topgst_gtopst_enr}{pointed $G$-category of pointed $G$-spaces and pointed morphisms}
\entry{$\GGTopii$}{def:ggtop_smc}{category of $\Gskg$-spaces}
\entry{$\FGTop$}{not:FGTop}{category of $\Fskg$-spaces}
\entry{$\cla$, $\Ner$, $\Rea$}{classifying_space}{classifying space, nerve, and realization}
\entry{$\sset$}{classifying_space}{category of simplicial sets}
\entry{$\clast$}{clast}{functor induced by $\cla$}
\entry{$\univ$}{not:univ}{complete $G$-universe}
\entry{$\IU$}{not:IU}{$\Gtop$-category of $G$-inner product spaces isomorphic to indexing $G$-spaces in $\univ$}
\entry{$(\oplus,0,\xi)$}{def:iu_spaces_i}{symmetric monoidal structure on $\IU$}
\entry{$\IU(-,-)$}{def:iu_spaces_ii}{$G$-space of linear isometric isomorphisms}
\entry{$\IUsk$}{def:iu_spaces_iv}{small skeleton of $\IU$ of indexing $G$-spaces in $\univ$}
\entry{$\upphi_V$}{VV'}{$G$-linear isometric isomorphism $V \fiso V'$ with $V' \in \IUsk$}
\entry{$(\pSet, \sma, \stplus)$}{not:pSet}{symmetric monoidal closed category of pointed sets and pointed functions}
\entry{$X_V$}{expl:iu_space}{value of an $\IU$-space $X$ at $V \in \IU$}
\entry{$X_f$}{iu_space_xf}{value of an $\IU$-space $X$ at a linear isometric isomorphism $f \in \IU$}
\entry{$\IUT$}{IUT}{category of $\IU$-spaces and $\IU$-morphisms}
\entry{$\GIUT$}{GIUT}{category of $\IU$-spaces and $G$-equivariant $\IU$-morphisms}
\entry{$S^V$}{def:g_sphere}{$V$-sphere for $V \in \IU$}
\entry{$\mu_{V,W}$}{gsphere_multiplication}{pointed $G$-homeomorphism $S^V \sma S^W \fiso S^{V \oplus W}$}
\entry{$(X,\umu)$}{def:gsp_module}{orthogonal $G$-spectrum with sphere action $\umu$}
\entry{$\GSp$}{not:GSp}{$\Gtop$-category of orthogonal $G$-spectra and morphisms}
\entry{$\Gspec$}{not:Gspec}{$\Top$-category of orthogonal $G$-spectra and $G$-morphisms}
\entry{$\Kg X$}{Kg_object}{orthogonal $G$-spectrum associated to a $\Gskg$-space $X$}
\entry{$(S^V)^{\sma\angordn}$}{SVsman}{pointed $G$-space $\Topgst(\sma\angordn, S^V)$}
\entry{$\assm_{\angordn}$}{assembly_sph}{pointed $G$-morphism $(S^V)^{\sma\angordn} \sma S^W \to (S^{V \oplus W})^{\sma\angordn}$}
\entry{$\Kg$}{Kg_functor}{prolongation functor $\GGTopii \to \Gspec$}
\entry{$\Kgo$}{Kgo_functors}{equivariant $K$-theory $\Kg\clast\Jgo \cn \AlglaxO \to \Gspec$}
\entry{$\Kgosg$}{Kgo_functors}{equivariant $K$-theory $\Kg\clast\Jgosg \cn \AlgpspsO \to \Gspec$}

\newchHeader{ch:ggcatg} 
\entry{$\Aut(S)$}{not:AutS}{group of self-bijections of $S$}
\entry{$\nbeta$}{ordn_be}{pointed finite $G$-set $(\ordn,\be)$}
\entry{$\FG$}{def:FG}{pointed $G$-category of pointed finite $G$-sets and pointed functions}
\entry{$(\sma, \ord{1}, \xi)$}{def:FG_permutative}{naive permutative $G$-category structure on $\FG$}
\entry{$\FGsma{q}$}{def:FG_smashpower}{$q$-fold smash power of $\FG$}
\entry{$\vstar$}{def:FG_smashpower}{initial-terminal basepoint of $\FG^{(0)}$}
\entry{$\ang{}$}{def:FG_smashpower}{empty tuple}
\entry{$\ang{\ord{n}^\be}$, $\ang{\ord{n}_i^{\be_i}}_{i \in \ufs{q}}$}{ordnbeta}{$q$-tuple $\big( \ord{n}_1^{\be_1}, \ldots, \ord{n}_q^{\be_q}\big)$ of pointed finite $G$-sets}
\entry{$h_*$}{def:injectionsFG}{reindexing $G$-functor $\FGsma{q} \to \FGsma{r}$ for an injection $h \cn \ufs{q} \to \ufsr$}
\entry{$\ord{n}_{\emptyset}^{\be_{\emptyset}}$, $\psi_{\emptyset}$}{ordn_emptyFG}{$\ord{1}$ and $1_{\ord{1}}$}
\entry{$\GG$}{def:GG}{pointed $G$-category of finite tuples of pointed finite $G$-sets}
\entry{$\vstar$}{GG_objects}{initial-terminal basepoint of $\GG$}
\entry{$\angordmal$}{GG_obj_gaction}{object in $\GG$}
\entry{$(f, \ang{\psi})$}{fangpsiGG}{morphism in $\GG$}
\entry{$\ifg$}{ifgGG}{length-1 inclusion $\FG \to \GG$}
\entry{$(\oplus, \ang{}, \xi)$}{def:GG_permutative}{naive permutative $G$-category structure on $\GG$}
\entry{$\sma$}{def:smashFGGG}{smash functor $\GG \to \FG$}
\entry{$\GGCatg$}{def:ggcatg}{2-category of $\GGG$-categories}
\entry{$\ig$}{ig}{pointed full subcategory inclusion $\Gsk \to \GG$}
\entry{$\igst$}{igst_iifunctor}{pullback 2-functor $\GGCatg \to \GGCatii$ along $\ig$}
\entry{$\Lg$}{thm:ggcat_ggcatg_iieq}{adjoint 2-equivalence $\GGCatii \to \GGCatg$}
\entry{$\GGpunc(-;-)$}{Lg_f_nbe}{nonzero morphisms in $\GG$}
\entry{$\ug$}{ug}{unit 2-natural isomorphism $1 \to \igst\Lg$}
\entry{$\vg$}{vg}{counit 2-natural isomorphism $\Lg\igst \to 1$}
\entry{$\jin$}{jin_nbe}{isomorphism $\ig\angordn \fiso \nbe$ in $\GG$}
\entry{$\Xnbeta$}{Xnbeta}{pointed category $X\angordn$ with twisted $G$-action}

\newchHeader{ch:hgo}
\entry{$\Aangordnbe$}{def:nbeta_gcat}{pointed $G$-category of $\angordnbe$-systems}
\entry{$g \cdot (a,\glu)$, $(ga, g\glu)$}{nsystem_gactionGG}{$g$-action on $\angordnbe$-system $(a,\glu)$}
\entry{$\Asgangordnbe$}{Asgangordnbe}{pointed $G$-category of strong $\angordnbe$-systems}
\entry{$\Aupom$, $\Aupomsg$}{Aupom}{pointed functors $\Aangordmal \to \Aangordnbe$ and $\Asgangordmal \to \Asgangordnbe$}
\entry{$\Hgo\A$, $\Adash$}{A_ptfunctorGG}{$H$-theory of an $\Op$-pseudoalgebra $\A$}
\entry{$\Hgosg\A$, $\Asgdash$}{A_ptfunctorGG}{strong $H$-theory of an $\Op$-pseudoalgebra $\A$}
\entry{$\Hgo f$}{hgo_f}{$H$-theory of a lax $\Op$-morphism}
\entry{$\Hgosg f$}{hgosg_f}{strong $H$-theory of an $\Op$-pseudomorphism}
\entry{$\Hgo \omega$}{hgo_omega}{$H$-theory of an $\Op$-transformation}
\entry{$\Hgosg \omega$}{hgosg_omega}{strong $H$-theory of an $\Op$-transformation}
\entry{$\Hgo$}{Hgo_twofunctor}{$H$-theory 2-functor $\AlglaxO \to \GGCatg$}
\entry{$\Hgosg$}{Hgo_twofunctor}{strong $H$-theory 2-functor $\AlgpspsO \to \GGCatg$}
\entry{$(\Lg,\igst)$}{thm:ggtop_ggtopg_iieq}{adjoint equivalence $\GGTopii \to \GGTopg$}
\entry{$\Kgg X$}{Kgg_object}{orthogonal $G$-spectrum associated to a $\GGG$-space $X$}
\entry{$(S^V)^{\sma\malp}$}{SVmalp}{pointed $G$-space $\Topgst(\sma\malp, S^V)$}
\entry{$\asm_{\malp}$}{assembly_sph}{pointed $G$-morphism $(S^V)^{\sma\malp} \sma S^W \to (S^{V \oplus W})^{\sma\malp}$}
\entry{$\Kgg$}{def:Kgg_functor}{prolongation functor $\GGTopg \to \Gspec$}
\entry{$\kiso$}{KgKgg}{natural isomorphism $\Kgg\Lg \fiso \Kg$}
\entry{$\Khgo$}{Khgo_functors}{equivariant $K$-theory $\Kgg\clast\Hgo \cn \AlglaxO \to \Gspec$}
\entry{$\Khgosg$}{Khgo_functors}{equivariant $K$-theory $\Kgg\clast\Hgosg \cn \AlgpspsO \to \Gspec$}
\entry{$\Kiso$, $\Kisosg$}{KgoKhgo}{natural isomorphisms $\Khgo \fiso \Kgo$ and $\Khgosg \fiso \Kgosg$}

\newchHeader{ch:shih}
\entry{$\FGCat$}{def:fgcat}{2-category of $\Fskg$-categories}
\entry{$\FGCatg$}{def:fgcatg}{2-category of $\FGG$-categories}
\entry{$\ig$}{ig_FG}{pointed full subcategory inclusion $\Fsk \to \FG$}
\entry{$(\Lg,\igst)$}{thm:fgcat_fgcatg_iieq}{adjoint 2-equivalence $\FGCat \to \FGCatg$}
\entry{$(X\ordn)_\be$}{def:gnbeta}{pointed category $X\ordn$ with twisted $G$-action}
\entry{$\,\jin$}{Xjin}{isomorphism $\ordn \fiso \nbeta$ in $\FG$}
\entry{$(a,\gl)$}{nsys}{$\ordn$-system}
\entry{$a_s$}{nsys_s}{$s$-component object $a_s$}
\entry{$\gl_{x;\, s, \ang{s_i}_{i \in \ufs{r}}}$}{gl-morphism}{gluing morphism}
\entry{$\theta_s$}{theta_s}{$s$-component morphism of $\theta$}
\entry{$\Aordn$}{A_ordn}{pointed category of $\ordn$-systems in $\A$}
\entry{$\Asgordn$}{Asgordn}{pointed category of strong $\ordn$-systems in $\A$}
\entry{$\Aordnbe$}{def:nsys_gcat}{pointed $G$-category of $\ordnbe$-systems in $\A$}
\entry{$g \cdot (a,\gl)$, $(ga, g\gl)$}{nsys_gaction}{$g$-action on an $\ordnbe$-system $(a,\gl)$}
\entry{$g\tha$}{g_tha}{$g$-action on a morphism $\tha$ of $\ordnbe$-systems}
\entry{$\Asgordnbe$}{Asgordnbe}{pointed $G$-category of strong $\ordnbe$-systems in $\A$}
\entry{$\Apsi$}{psitil_f}{pointed functor $\Aordmal \to \Aordnbe$}
\entry{$\big(a^{\psitil}, \gl^{\psitil}\big)$}{apsitil_s}{$\ordnbe$-system $(\Apsi) (a,\gl)$}
\entry{$\tha^{\psitil}$}{thapsitil_comp}{morphism $(\Apsi)(\theta)$}
\entry{$\Asgpsi$}{psitilsg_f}{pointed functor $\Asgordmal \to \Asgordnbe$}
\entry{$\Sgo\A$, $\Adash$}{sys_FGcat}{Shimakawa $H$-theory of an $\Op$-pseudoalgebra $\A$}
\entry{$\Sgosg\A$, $\Asgdash$}{sys_FGcat}{Shimakawa strong $H$-theory of an $\Op$-pseudoalgebra $\A$}
\entry{$\Sgo$}{def:sgo}{Shimakawa $H$-theory $\AlglaxO \to \FGCatg$}
\entry{$\Sgosg$}{def:sgosg}{Shimakawa strong $H$-theory $\AlgpspsO \to \FGCatg$}
\entry{$\Jgos$}{jgos_jgossg}{Shimakawa $J$-theory $\AlglaxO \to \FGCat$}
\entry{$\Jgossg$}{jgos_jgossg}{Shimakawa strong $J$-theory $\AlgpspsO \to \FGCat$}

\newchHeader{ch:shimakawa_K}
\entry{$\FGTopg$}{def:ggtopg}{category of $\FGG$-spaces}
\entry{$\Kfg X$}{Kfgx}{orthogonal $G$-spectrum associated to an $\FGG$-space $X$}
\entry{$(S^V)^{\mal}$}{SVmal}{pointed $G$-space $\Topgst(\mal, S^V)$}
\entry{$\Kfg$}{def:Kfg_functor}{prolongation functor $\FGTopg \to \Gspec$}
\entry{$\WG$}{PwgUgs}{pointed $G$-category of pointed finite $G$-CW complexes and pointed morphisms}
\entry{$\Pwg$}{Pwgxa}{prolongation functor $\FGTopg \to \WGTopg$}
\entry{$\Ksho$}{ksho_kshosg}{Shimakawa $K$-theory $\Kfg\clast\Sgo \cn \AlglaxO \to \Gspec$}
\entry{$\Kshosg$}{ksho_kshosg}{Shimakawa strong $K$-theory $\Kfg\clast\Sgosg \cn \AlgpspsO \to \Gspec$}
\entry{$\Delta$}{not:Delta_cat}{category of $\ordr$ for $r \geq 0$ and weakly order-preserving functions}
\entry{$d^i$, $s^i$}{not:face_degen}{coface and codegeneracy}
\entry{$X_\crdot$}{def:simp_obj_iii}{simplicial object $\Deltaop \to \V$}
\entry{$X_r$, $d_i$, $s_i$}{def:simp_obj_iii}{$r$-simplex object, $i$-th face, and $i$-th degeneracy}
\entry{$\V^{\Deltaop}$}{def:simp_obj_iv}{category of simplicial objects and simplicial morphisms in $\V$}
\entry{$\htpy \cn \fun \simeq \fun'$}{def:simp_obj_v}{homotopy of simplicial morphisms}
\entry{$\Vse$}{def:cat_bar}{$\V$-category given by the objects and internal hom of $\V$}
\entry{$\Bcdot(\hun,\C,\fun)$}{bcdot_hcf}{simplicial bar construction}
\entry{$\Bc_r(\hun,\C,\fun)$}{bc_r}{$r$-simplex object of the simplicial bar construction}
\entry{$\simpdot$}{not:simpdot}{cosimplicial object $\Delta \to \V$}
\entry{$|X_\crdot|$}{x_dot}{realization $\int^{\ordr \in \Delta} X_{\ordr} \otimes \simp^{\ordr}$}
\entry{$\Bc(\hun,\C,\fun)$}{bar_realization}{realized bar construction $\int^{\ordr \in \Delta} \Bc_r(\hun,\C,\fun) \otimes \simp^{\ordr}$}
\entry{$\C_c$}{representable_c}{representable $\V$-functor $\C(-,c) \cn \C^{\op} \to \Vse$}
\entry{$\Bc(\C,\C,\fun)$}{Bc_CCf}{bar construction of $\fun \cn \C \to \Vse$}
\entry{$\Bc(\C_c,\C,\fun)$}{real_Cc}{value of $\Bc(\C,\C,\fun)$ at $c$}
\entry{$(\fun_c)_\crdot$}{not:funcdot}{constant simplicial object at $\fun_c$}
\entry{$\retn_c$}{retn_c}{simplicial retraction $\Bcdot(\C_c,\C,\fun) \to (\fun_c)_\crdot$ at $c$}
\entry{$\retn$}{retn}{retraction $\Bc(\C,\C,\fun) \to \fun$ of $\fun$}
\entry{$\secn_c$}{secn_c}{section $(\fun_c)_\crdot \to \Bcdot(\C_c,\C,\fun)$ at $c$}
\entry{$\CVse$}{not:CVse}{category of $\V$-functors $\C \to \Vse$ and $\V$-natural transformations}
\entry{$\Bc(\C,\C,-)$}{BCC_functor}{bar functor on $\CVse$}
\entry{$\retn$}{retn_bar_id}{retraction $\Bc(\C,\C,-) \to 1$}
\entry{$\simpdotpl$}{not:simpdotpl}{cosimplicial pointed $G$-space $\Delta \to \Gtopst$}
\entry{$\simppl^r$}{not:simpplr}{pointed $G$-space $\simp^r \sqcup {*}$}
\entry{$\simp^r$}{top_simplex}{topological $r$-simplex with trivial $G$-action}
\entry{$|X_\crdot|$}{real_simpgspace}{realization $\int^{\ordr \in \Delta} X_{\ordr} \sma \simppl^r$ of a simplicial pointed $G$-space $X_\crdot$} 
\entry{$\Bc(\FG,\FG,X)$, $\Bc X$}{bar_FG}{bar construction of an $\FGG$-space $X$}
\entry{$\Bc\big(\FG(-,\mal),\FG,X\big)$}{real_FG}{value of $\Bc X$ at $\mal \in \FG$}
\entry{$\BcFG$, $\Bc$}{bar_functor_FG}{bar functor on $\FGTopg$}
\entry{$\retn_X$}{ret_FG}{retraction of an $\FGG$-space $X$}
\entry{$\retn$}{retn_barFG_id}{retraction $\Bc \to 1$ on $\FGTopg$}
\entry{$\Khsho$}{khsho_khshosg}{homotopical Shimakawa $K$-theory $\Kfg\Bc\clast\Sgo \cn \AlglaxO \to \Gspec$}
\entry{$\Khshosg$}{khsho_khshosg}{homotopical Shimakawa strong $K$-theory $\Kfg\Bc\clast\Sgosg \cn \AlgpspsO \to \Gspec$}
\entry{$\Lat_r X_\crdot$}{def:reedy_cof}{$r$-th latching $G$-space of a simplicial $G$-space $X_\crdot$}

\newchHeader{ch:shim_top}
\entry{$\Cpunc(c,c')$}{Cpunc}{set of nonzero morphisms $\C(c,c') \setminus \{0\}$}
\entry{$\funst$}{funst_pullback}{pullback functor $\DTopgst \to \CTopgst$ along $\fun \cn \C \to \D$}
\entry{$\funl$}{fun_adj}{left adjoint of $\funst$}
\entry{$\funu$}{funu}{unit $1 \to \funst\funl$}
\entry{$\funv$}{funv}{counit $\funl\funst \to 1$}
\entry{$\KC$}{def:KC}{prolongation functor $\CTopgst \to \Gspec$}
\entry{$\KC X$}{KC_object}{orthogonal $G$-spectrum associated to a pointed $G$-functor $X \cn \C \to \Topgst$}
\entry{$(S^V)^{\tF c}$}{SVFc}{pointed $G$-space $\Topgst(\tF c, S^V)$}
\entry{$\asm_{c}$}{KC_assembly}{pointed $G$-morphism $(S^V)^{\tF c} \sma S^W \to (S^{V \oplus W})^{\tF c}$}
\entry{$\kfgiso$}{not:kfgiso}{natural isomorphism $\KD\funl \fiso \KC$}
\entry{$\kggik$}{KggKfg_isos}{natural isomorphism $\Kgg\ifgl \fiso \Kfg$}
\entry{$\kfgsk$}{KggKfg_isos}{natural isomorphism $\Kfg\!\smal\! \fiso \Kgg$}
\entry{$\kbgg$}{not:kbgg}{natural transformation $\Kfg\Bc \to \Kgg\ifgl$}
\entry{$\GGssetgst$}{def:GGssetgst}{category of $\GG$-simplicial $G$-sets}
\entry{$\kgbk$}{not:kgbk}{natural transformation $\Kfg\Bc\!\smal \to \Kgg$}

\newchHeader{ch:h_comparison}
\entry{$\Pist$}{Pistar_twonat}{$H$-theory comparison $\smast\Sgo \to \Hgo$}
\entry{$\Pist_\A$}{PistA}{$\A$-component $G$-natural transformation $\smast\Sgo\A \to \Hgo\A$}
\entry{$\Pist_{\A,\angordnbe}$}{PistAnbe}{$\angordnbe$-component pointed $G$-functor of $\Pist_\A$}
\entry{$(a^\stimes,\gl^\stimes)$}{agltimes}{$\angordnbe$-system $\Pist_{\A,\angordnbe} (a,\gl)$}
\entry{$\angstimes$}{angstimes}{product subset $\txprod_{j\in \ufs{q}}\, s_j \subseteq 
\txprod_{j\in \ufs{q}}\, \ufs{n}_j$}
\entry{$\Pistsg$}{Pistarsg_twonat}{strong $H$-theory comparison $\smast\Sgosg \to \Hgosg$}
\entry{$\Pistsg_\A$}{PistsgA}{$\A$-component $G$-natural transformation $\smast\Sgosg\A \to \Hgosg\A$}
\entry{$\Pistsg_{\A,\angordnbe}$}{PistsgAnbe}{$\angordnbe$-component pointed $G$-functor of $\Pistsg_\A$}
\entry{$\proCnbe$}{not:proCnbe}{$\angordnbe$-twisted product of a pointed $G$-category $\C$}
\entry{$\ang{c_{\bdi} \in \C}_{\bdi \sins \ufs{n_1n_2\cdots n_q}}$}{proCnbe_object}{object or morphism in $\proCnbe$}
\entry{$\,\bdi$}{bdi}{$q$-tuple $\ang{i_j}_{j \in \ufs{q}} \in \txprod_{j \in \ufs{q}}\, \ufs{n}_j$}
\entry{$g\bdi$}{gbdi}{diagonal $g$-action $\ang{g i_j}_{j \in \ufs{q}}$}
\entry{$gc$}{proCnbe_gaction}{$g$-action $\ang{g c_{\ginv\bdi}}_{\bdi \sins \ufs{n_1n_2\cdots n_q}} \in \proCnbe$ on $c = \ang{c_{\bdi}}_{\bdi \sins \ufs{n_1n_2\cdots n_q}}$}
\entry{$\zb$}{sgotoprod}{pointed $G$-functor $\Asmaangordnbe \to \proAnbe$}
\entry{$\zb(a,\gl)$}{sgotoprod_def}{$\ang{a_{\{\bdi\}}}_{\bdi \sins \ufs{n_1 n_2 \cdots n_q}}$}
\entry{$\zb\theta$}{sgotoprod_def}{$\ang{\theta_{\{\bdi\}}}_{\bdi \sins \ufs{n_1 n_2 \cdots n_q}}$}
\entry{$\{\bdi\}$}{not:bdiangij}{one-element subset $\{\ang{i_j}_{j \in \ufs{q}}\} \subseteq \txprod_{j \in \ufs{q}}\, \ufs{n}_j$}
\entry{$\zbsg$}{sgosgtoprod}{restriction of $\zb$ to $\Asgsmaangordnbe$}
\entry{$\zd$}{hgotoprod}{pointed $G$-functor $\Aangordnbe \to \proAnbe$}
\entry{$\zd(a,\glu)$}{hgotoprod_def}{$\ang{a_{\ang{\{i_j\}}_{j \in \ufsq}}}_{\bdi \sins \ufs{n_1 n_2 \cdots n_q}}$}
\entry{$\zd\theta$}{hgotoprod_def}{$\ang{\theta_{\ang{\{i_j\}}_{j \in \ufsq}}}_{\bdi \sins \ufs{n_1 n_2 \cdots n_q}}$}
\entry{$\ang{\{i_j\}}_{j \in \ufsq}$}{not:angij}{$q$-tuple $(\{i_1\}, \{i_2\}, \ldots, \{i_q\})$}
\entry{$\zdsg$}{hgosgtoprod}{restriction of $\zd$ to $\Asgangordnbe$}
\entry{$\bdk^t$}{not:bdkt}{least element in its $G$-orbit $G\bdk^t$}
\entry{$G_t$}{not:Gt}{stabilizer of $\bdk^t$}
\entry{$N_t$}{not:Gt}{cardinality of $G\bdk^t$}
\entry{$g_{t,d}$}{gtd}{element in $G$ such that $g_{t,d}(1) = d$ and $g_{t,1} = e$}
\entry{$\uph$}{uph}{functor $\proCnbe \to \txprod_{t \in \ufsr}\, \C$ sending $\ang{c_{\bdi}}$ to $\ang{c_{\bdk^t}}_{t \in \ufsr}$}
\entry{$\uphinv$}{uph}{functor $\txprod_{t \in \ufsr}\, \C \to \proCnbe$ sending $\ang{c_t}_{t \in \ufsr}$ to $\ang{g_{t,d} c_t}_{d \in \ufsN_t,\, t \in \ufsr}$}

\newchHeader{ch:sgoprod}
\entry{$\ep_n$}{Opn_object}{arbitrarily chosen object in $\Op(n)$}
\entry{$\zbsgad$}{zbsgad}{pointed functor $\proAnbe \to \Asgsmaangordnbe$}
\entry{$|s|$}{subset_s}{cardinality of a set $s$}
\entry{$\sig_{s,\ang{s_\ell}_{\ell \in \ufs{r}}}$}{sigma_angsi}{permutation from $s$ to $\coprod_{\ell \in \ufs{r}}\, s_\ell$}
\entry{$\alp_{x;\, s, \ang{s_\ell}_{\ell \in \ufs{r}}}$}{al_Op}{unique isomorphism $\ga(x \sscs \ang{\ep_{\cardsl}}_{\ell \in \ufs{r}}) \sig_{s,\ang{s_\ell}_{\ell \in \ufs{r}}} \fiso \ep_{\cards}$ in $\Op(\cards)$}
\entry{$\unis$}{unit_zbsg}{unit $1_1 \cn 1 \to \zbsg\zbsgad$}
\entry{$\cous$}{counit_zbsg}{counit $\zbsgad\zbsg \fiso 1$}
\entry{$(\func,\pse)$}{def:pseudoG}{pseudo $G$-equivariant functor with constraint $\pse$}
\entry{$\zbpg$}{zbsgad_constraint}{pseudo $G$-equivariant constraint for $\zbsgad$}
\entry{$\si^{g,s}$}{si_g_s}{permutation $\ginv s \to s$ induced by $g \in G$}
\entry{$\al^{g,s}$}{al_g_s}{unique isomorphism $\ginv \ep_{\cards} \si^{g,s} \fiso \ep_{\cards}$ in $\Op(\cards)$}

\newchHeader{ch:hgoprod}
\entry{$\zdsgad$}{zdsgad}{pointed functor $\proAnbe \to \Asgangordnbe$}
\entry{$\bskl$}{bskl}{marker $\angscompkskell$}
\entry{$\usi_{\angs, k, \ang{s_{k,\ell}}_{\ell \in \ufs{r}}}$}{usi_perm}{permutation $\angstimes \to \txcoprod_{\ell \in \ufs{r}}\, \bskltimes$}
\entry{$\ual_{x; \angs, k, \ang{s_{k,\ell}}_{\ell \in \ufs{r}}}$}{ual_iso}{unique isomorphism $\ga\big(x \sscs \ang{\ep_{\bskltimesc}}_{\ell \in \ufs{r}}\big) \usi_{\angs, k, \ang{s_{k,\ell}}_{\ell \in \ufs{r}}} \fiso \ep_{\angstimesc}$ in $\Op(\angstimesc)$}
\entry{$\unit $}{unit_zdsg}{unit $1_1 \cn 1 \to \zdsg\zdsgad$}
\entry{$\counit$}{counit_zdsg}{counit $\zdsgad\zdsg \to 1$}
\entry{$\bdip$}{bdip}{$q$-tuple $(\{i_1\}, \{i_2\}, \ldots, \{i_q\})$}
\entry{$\nbehat$}{nbe_prime}{object $\bordi{n}{\be}{j}_{j=2}^q \in \GG$}
\entry{$\angshat$}{s_prime}{marker $\ang{s_j \subseteq \ufs{n}_j}_{j=2}^q$}
\entry{$\angshattimes$}{sone_angshat}{product subset $\txprod_{j=2}^q\, s_j \subseteq \txprod_{j=2}^q\, \ufs{n}_j$}
\entry{$(a_{\{i_1\}}, \glu_{\{i_1\}})$}{aglu_ione}{$i_1$-restriction of $(a,\glu)$}
\entry{$\zdpg$}{zdsgad_constraint}{pseudo $G$-equivariant constraint for $\zdsgad$}
\entry{$\usi^{g,\angs}$}{usi_g_angs}{permutation $\ginvangstimes \to \angstimes$}
\entry{$\ual^{g,\angs}$}{ual_g_angs}{unique isomorphism $\ginv \ep_{\angstimesc} \usi^{g,\angs} \fiso \ep_{\angstimesc}$ in $\Op(\angstimesc)$}

\newchHeader{ch:compgen}
\entry{$\inc$, $\incC$}{incC}{inclusion $G$-functor $\C \to \Catg(\EG,\C)$}
\entry{$\pn$}{proC}{right adjoint inverse $\Catg(\EG,\C) \to \C$ of $\inc$}
\entry{$\ppng$}{ppng}{pseudo $G$-equivariant constraint of $\pn$}
\entry{$(f_g, \gl^{f_g})$}{f_sub_g}{image of $g \in G$ under a functor $f \cn \EG \to \Ahmal$}
\entry{$f_{g,s}$}{fgs}{$s$-component object of $f_g$}
\entry{$f_{g,s,h}$}{fgsh}{image of $h \in G$ under the functor $f_{g,s} \cn \EG \to \A$}
\entry{$\cni$}{cni_functor}{$G$-equivariant inverse $\Catg(\EG, \Ahmal) \to \Ahmal$ of the inclusion $G$-functor for $\Ahmal$}
\entry{$\cnisg$}{cnisg_functor}{strong variant $\Catg(\EG, \sgAhmal) \to \sgAhmal$ of $\cni$}
\entry{$\ucni$}{ucni}{$G$-equivariant unit $1_1 \cn 1 \to \cni\inc^{\Ahmal}$}
\entry{$\ucnisg$}{ucnisg}{$G$-equivariant unit $1_1 \cn 1 \to \cnisg\inc^{\sgAhmal}$}
\entry{$\ccni$}{ccni}{$G$-equivariant counit $\inc^{\Ahmal}\cni \to 1$}
\entry{$\ccnisg$}{ccnisg}{$G$-equivariant counit $\inc^{\sgAhmal}\cnisg \to 1$}
\entry{$(\ion,\pon,\uon,\von)$}{EHEG}{adjoint $H$-equivalence $\EH \to \EG$ for a subgroup $H \subseteq G$}
\entry{$(f_g, \glu^{f_g})$}{f_g}{image of $g \in G$ under a functor $f \cn \EG \to \Ahnbe$}
\entry{$f_{g,\angs}$}{fgangs}{$\angs$-component object of $f_g$}
\entry{$f_{g,\angs,h}$}{fgangsh}{image of $h \in G$ under the functor $f_{g,\angs} \cn \EG \to \A$}
\entry{$\ci$}{ci_functor}{$G$-equivariant inverse $\Catg(\EG, \Ahnbe) \to \Ahnbe$ of the inclusion $G$-functor for $\Ahnbe$}
\entry{$\cisg$}{cisg_functor}{strong variant $\Catg(\EG, \sgAhnbe) \to \sgAhnbe$ of $\ci$}
\entry{$\uci$}{uci}{$G$-equivariant unit $1_1 \cn 1 \to \ci\inc^{\Ahnbe}$}
\entry{$\ucisg$}{ucisg}{$G$-equivariant unit $1_1 \cn 1 \to \cisg\inc^{\sgAhnbe}$}
\entry{$\cci$}{cci}{$G$-equivariant counit $\inc^{\Ahnbe}\ci \to 1$}
\entry{$\ccisg$}{ccisg}{$G$-equivariant counit $\inc^{\sgAhnbe}\cisg \to 1$}
\entry{$\fun^H$}{def:cat_weakg}{$H$-fixed subfunctor of a $G$-functor $\fun$}
\entry{$\eqg$}{not:cateqg}{categorical weak $G$-equivalence}

\newchHeader{ch:special}
\entry{$\angordone$, $\ang{\ordone}_{j \in \ufsq}$}{angordone}{object in $\Gsk$ consisting of copies of $\ordone = \{0,1\}$}
\entry{$\segbdi$}{i-char}{$\bdi$-characteristic morphism $\nbe \to \angordone$}
\entry{$\seg_{i_j}$}{ij-char}{$i_j$-characteristic function $\ordn_j \to \ordone$}
\entry{$\segnbe$, $\segnbex$}{nbe_segal}{$\nbe$-Segal functor $X\nbe \to \Xonenbe$ for a $\GGG$-category $X$}
\entry{$\Siangordn$}{Siangordn}{product group $\txprod_{j=1}^q \Si_{n_j}$}
\entry{$\FGSin$}{not:FGSin}{family of graph subgroups of $\GSin$}
\entry{$\bsi$, $\ang{\si_j \in \Si_{n_j}}_{j \in \ufsq}$}{Siangord_obj}{$q$-tuple of permutations}
\entry{$\segn$, $\segnx$}{nsegal}{$\angordn$-Segal functor $X\angordn \to X\angordone^{n_1\cdots n_q}$ for a $\Gskg$-category $X$}
\entry{$\Xh$}{Xh}{$H$-restriction of a $\Gskg$-category $X$}
\entry{$\Gpl$}{not:Gpl}{$G \sqcup {*}$}
\entry{$\,\bbr$}{i_bbr}{$G$-inverse $\Catg(\EG,\Chnbe) \to \Chnbe$ of the inclusion $G$-functor for $\Chnbe$}

\newchHeader{ch:kgmmo}
\entry{$\Op_n$}{not:Opn}{$\Op(n)$ for an operad $\Op$}
\entry{$\DG$}{def:dgo}{$\Gcatst$-category associated to a reduced $\Gcat$-operad $\Op$}
\entry{$(\psi; \ang{x_j}_{j \in \ufsn})$}{dgmn_obj}{object or morphism in $\DG(\mal,\nbeta)$}
\entry{$\objzero$}{dgmn_basept}{basepoint $(0; \ang{*}_{j \in \ufsn})$ of $\DG(\mal,\nbeta)$}
\entry{$\objone_{\objx}$}{dg_idtwocell}{identity morphism $(\psi; \ang{1_{x_j}}_{j \in \ufsn})$}
\entry{$g (\psi; \ang{x_j}_{j \in \ufsn})$}{dgmn_gact}{$g$-action on $\DG(\mal,\nbeta)$}
\entry{$\tau_{\ginv}$}{tauginv}{permutation $(g\psi\ginv)^{-1}j \to (g\psi)^{-1}j$ given by $\ginv$}
\entry{$\objone_{\mal}$}{dg_idonecell}{identity 1-cell $(1_{\mal}; \ang{\opu}_{i \in \ufsm})$ of $\mal \in \DG$} 
\entry{$\comp$}{dg_comp}{composition of $\DG$}
\entry{$\tau^j_{\psi,\phi}$}{tauj_psiphi}{shuffle $(\psi\phi)^{-1}(j) \to \coprod_{i \in \psiinv j} \phiinv i$}
\entry{$\Pig$}{def:Pig}{pointed $G$-subcategory of $\FG$ consisting of injections}
\entry{$\gio$}{dgcoop}{$\Gcatst$-inclusion $\Pig \to \DG$}
\entry{$\gxi$}{dgcoop}{$\Gcatst$-projection $\DG \to \FG$}
\entry{$\DGCatg$}{def:dgcatg_iicat_i}{2-category with $\Pig$-strict $\Gcatst$-pseudotransformations as 1-cells}
\entry{$\dgcatg$}{def:dgcatg_iicat_ii}{sub-2-category of $\DGCatg$ with $\Gcatst$-natural transformations as 1-cells}
\entry{$\gxist$}{gxist}{pullback 2-functors $\FGCatg \to \dgcatg$ and $\FGCatg \to \DGCatg$ along $\gxi$}
\entry{$\AlgO$}{def:algstpsO}{2-category of $\Op$-algebras, morphisms, and $\Op$-transformations}
\entry{$\AlgstpsO$}{def:algstpsO}{2-category of $\Op$-algebras, $\Op$-pseudomorphisms, and $\Op$-transformations}
\entry{$\Rg\A$}{rga}{$\DGG$-category associated to an $\Op$-algebra $\A$}
\entry{$\Rg f$}{rgf}{$\Pig$-strict $\Gcatst$-pseudotransformation associated to an $\Op$-pseudomorphism $f$}
\entry{$\Rg\omega$}{rgom}{$\Gcatst$-modification associated to an $\Op$-transformation $\omega$}
\entry{$\Rg$}{rg_twofunctor}{2-functors $\AlgO \to \dgcatg$ or $\AlgstpsO \to \DGCatg$}
\entry{$X^2$}{psfggcat_xtwo}{composition constraint of a pseudo $\FGG$-category $X$}
\entry{$\FGCatgps$}{def:psfggcat_iicat}{2-category of pseudo $\FGG$-categories}
\entry{$\Incj$}{fgcatg_psfgcatg}{inclusion 2-functor $\FGCatg \to \FGCatgps$}
\entry{$\gze$}{gze}{strictly unital $\Gcatst$-pseudofunctor $\FG \to \DG$}
\entry{$\gze^2$}{gze_two}{composition constraint of $\gze$}
\entry{$\gzest$}{gzest}{pullback 2-functor $\DGCatg \to \FGCatgps$ along $\gze$}
\entry{$\str$}{str_iifunctor}{strictification left 2-adjoint $\FGCatgps \to \FGCatg$ of $\Incj$}
\entry{$\stru$}{stru}{unit $1 \to \Incj\str$}
\entry{$\strv$}{strv}{counit $\str\Incj \to 1$}
\entry{$\Kgmmo$}{kgmmo_diag}{GMMO $K$-theory $\AlgstpsO \to \Gspec$ for $\Op$}
\entry{$\cgs_\A$}{cgsa}{$\Gcatst$-pseudotransformation $\gxist\Sgosg\algi\A \to \Rg\A$ for an $\Op$-algebra $\A$}
\entry{$\cgs$}{not:cgs}{2-natural transformation $\gxist\Sgosg\algi \to \Rg$}

\newchHeader{ch:kgl}
\entry{$G/H$}{ex:univ_gset}{orbit $G$-set $\{gH \tmid g \in G\}$}
\entry{$\ome$}{not:ome}{set of natural numbers $\{0,1,2,\ldots\}$}
\entry{$\ugsetw$}{ugsetw}{universal $G$-set $\coprod_{H \subseteq G} \, \ome \times G/H$}
\entry{$\omea$}{def:omea}{countable $G$-set of functions $A \to \ome$}
\entry{$\omeg$}{ex:omea}{universal $G$-set of functions $G \to \ome$}
\entry{$\Injm$}{def:mcat_i}{monoid of injections $\ome \to \ome$}
\entry{$\schm$}{def:mcat_ii}{translation category $\tn\Injm$}
\entry{$[v,u]$}{not:mor_vu}{unique isomorphism $u \to v$ in $\schm$}
\entry{$(-)_\mstar$}{not:ustar}{action functor or natural isomorphism on an $\schm$-category}
\entry{$\C^{\suppempty}$}{not:Csuppempty}{full subcategory of $\C$ consisting of $\schm$-fixed objects}
\entry{$\supp(x)$}{not:suppx}{support of an object $x$ in an $\schm$-category}
\entry{$\Mcatt$}{not:Mcatt}{category of tame $\schm$-categories and $\schm$-functors}
\entry{$[v,u]_{\mstar}^x$}{vustar_x}{$x$-component isomorphism $u_{\mstar}(x) \to v_{\mstar}(x)$ of $[v,u]_{\mstar}$}
\entry{$\Fome$}{ex:sch2.14}{tame $\schm$-category of finite subsets of $\ome$ and bijections}
\entry{$\rep_S$}{rep_S}{chosen bijection $\ome \fiso S$ for a countably infinite set $S$}
\entry{$\schj$}{def:repar_ii}{2-category of countably infinite sets, injections, and unique 2-cells}
\entry{$\Cbrac$}{Cbrac}{reparametrization 2-functor $\schj \to \Cat$ for an $\schm$-category $\C$}
\entry{$p^\ome$}{repU_p_repV}{injection $\rep_V^{-1} p \rep_U \cn \ome \to \ome$ for an injection $p \cn U \to V$}
\entry{$[q^\ome,p^\ome]_{\mstar}$}{Cbrac_iicell}{$\schm$-action natural isomorphism $p^\ome_{\mstar} \to q^\ome_{\mstar}$}
\entry{$\Comeg$}{def:sch2.21_i}{$\omeg$-reparametrization of $\C$}
\entry{$g^\ome$}{rep_g_rep}{bijection $\rep_{\omeg}^{-1} g \rep_{\omeg} \cn \ome \to \ome$}
\entry{$u^\ome$}{rep_u_rep}{injection $\rep_{\omeg}^{-1} u \rep_{\omeg} \cn \ome \to \ome$}
\entry{$\Fixg\C$}{fixgc}{$G$-fixed $\schm$-category $\Comeg^G$}
\entry{$\fomeg$}{fomeg}{$\omeg$-reparametrization of $\fun$}
\entry{$\Bxtimes$}{def:bxtimes_ii}{box product of $\schm$-categories}
\entry{$(\Mcatt, \Bxtimes, \bone, \bxi)$}{def:bxtimes_iv}{symmetric monoidal category of tame $\schm$-categories and $\schm$-functors}
\entry{$\iota$}{not:bxinc}{inclusion $\schm$-functor $\C\bxtimes\D \to \C\ttimes\D$}
\entry{$\Parcat$}{def:parcat}{category of parsummable categories and parsummable functors}
\entry{$(\C,\psum,\pzero)$}{expl:parcat}{parsummable category}
\entry{$\FMCat$}{def:fmcat}{category of $\FM$-categories}
\entry{$\Jsc$}{Jsc}{global $J$-theory $\Parcat \to \FMCat$}
\entry{$\bdI$}{def:sch8.2_i}{category of finite sets and injections}
\entry{$\ICat$}{def:sch8.2_ii}{category of $\bdI$-categories and morphisms}
\entry{$\bone$}{def:sch8.2_ii}{constant functor $\bdI \to \Cat$ at the terminal category $\bone$}
\entry{$\ITop$}{def:sch8.2_iii}{category of $\bdI$-spaces and morphisms}
\entry{$*$}{def:sch8.2_iii}{constant functor $\bdI \to \Top$ at the one-point space $*$}
\entry{$\upr$}{upr_functor}{functor $\Mcat \to \ICat$}
\entry{$i_!$}{extzero}{injection $\omea \to \omeb$ extending an injection $i \cn A \to B$}
\entry{$i_!^\ome$}{iome_ext}{injection $\rep_{\omeb}^{-1} i_! \rep_{\omea} \cn \ome \to \ome$}
\entry{$\FICat$}{def:ficat_i}{category of $\FI$-categories}
\entry{$\FITop$}{def:ficat_ii}{category of $\FI$-spaces}
\entry{$\uprst$}{fmcat_ficat}{functor $\FMCat \to \FICat$ induced by $\upr$}
\entry{$\clast$}{ficat_fitop}{functor $\FICat \to \FITop$ induced by $\cla$}
\entry{$\bR[A]$}{def:asphere}{real vector space of functions $A \to \bR$}
\entry{$S^A$}{not:asphere}{$A$-sphere $\bR[A] \sqcup \{\infty\}$}
\entry{$S^\emptyset$}{not:emptysp}{$\emptyset$-sphere $\{*,\infty\}$}
\entry{$X_A$}{def:sptop_i}{value of a symmetric spectrum $X$ at a finite set $A$}
\entry{$i_*$}{not:istarab}{structure morphism $X_A \sma S^{B \setminus i(A)} \to X_B$ for an injection $i \cn A \to B$}
\entry{$f_A$}{def:sptop_ii}{value of a symmetric spectrum morphism $f$ at a finite set $A$}
\entry{$\Sptop$}{sptop_mor}{category of symmetric spectra}
\entry{$\Sptopg$}{def:sptop_iii}{category of $G$-symmetric spectra}
\entry{$X_G$}{def:sptop_iv}{underlying $G$-spectrum of a symmetric spectrum $X$}
\entry{$\pinhu(X)$}{not:pinhu}{$n$-th $H$-equivariant homotopy group of a $G$-symmetric spectrum $X$}
\entry{$\Gspu$}{gsp_forget}{forgetful functor $\Gspec \to \Sptopg$}
\entry{$\Kfi$}{kfi_functor}{prolongation functor $\FITop \to \Sptop$}
\entry{$\assm_n$}{not:assmn}{pointed morphism $(S^A)^n \sma S^{B \setminus i(A)} \to (S^B)^n$}
\entry{$\Ksc$}{def:sch3.3}{global $K$-theory $\Kfi\clast\uprst\Jsc \cn \Parcat \to \Sptop$}
\entry{$\kcg$}{kcga}{$G$-symmetric spectrum associated to a parsummable category $\C$}
\entry{$\ascgc$, $\bscgc$}{sch4.15}{natural $\pistu$-isomorphisms $(\Ksc\C)_G \to \gacomegs \fot{} \kcg$}
\entry{$\Jscg$}{Jscg}{Schwede $J$-theory $\Parcat \to \FGCat$}
\entry{$\Kfgsi$}{kfgsi_functor}{prolongation functor $\FGTop \to \Sptopg$}
\entry{$\Kscgb$}{def:kscgb}{Schwede $K$-theory $\Kfgsi\clast\Jscg \cn \Parcat \to \Sptopg$}
\entry{$\Isc$}{Isc_nt}{natural componentwise $\pistu$-isomorphism $\Kscg \to \Kscgb$}

\newchHeader{ch:kgl_gmmo}
\entry{$\Inop$}{def:inop}{injection operad}
\entry{$\iphi^\si$}{iop_symmetry}{$\si$-action on an injection $\iphi \cn \ufsn \itimes \ome \to \ome$}
\entry{$\Iop$}{iopn}{categorical injection operad with $\Iop(n) = \tn\Inop(n)$}
\entry{$\iphi_*$}{iphist_caction}{$\iphi$-action functor $\C^n \to \C$}
\entry{$[\ipsi, \iphi]_*$}{ipsiiphi_st}{natural isomorphism $\iphi_* \to \ipsi_*$}
\entry{$\Inopg$}{def:global_einf}{injection $G$-operad}
\entry{$\Iopg$}{iopgn}{categorical injection $G$-operad with $\Iopg(n) = \tn\Inopg(n)$}
\entry{$\gacomeg$}{gacomeg}{$\Iopg$-action $G$-functor $\Iopg(n) \times \Comeg^n \to \Comeg$}
\entry{$\iphi^j$}{not:iphij}{injection $\iphi(j,-) \cn \omeg \to \omeg$ for an injection $\iphi \cn \ufsn \itimes \omeg \to \omeg$}
\entry{$\iphijome$}{iphijome}{injection $\rep_{\omeg}^{-1} \iphi^j \rep_{\omeg} \cn \ome \to \ome$}
\entry{$\iphijomest$}{iphijomest}{$\iphijome$-action functor $\C \to \C$}
\entry{$[\ipsi^{j,\ome}, \iphi^{j,\ome}]_{\mstar}$}{not:ipsiphijome}{$\schm$-action natural isomorphism $\iphijomest \to \ipsijomest$}
\entry{$\Ig$}{Ig_functor}{functor $\Parcat \to \Algiopg$}
\entry{$\Ig\C$}{not:IgC}{$\Iopg$-algebra $(\Comeg,\gacomeg)$}
\entry{$\Ig\fun$}{not:IgC}{$\Iopg$-algebra morphism $\fomeg$}
\entry{$\DG(\mal,\nbeta)$}{dgmn_iopg}{pointed $G$-category $\coprod_{\psi \in \FG(\mal,\,\nbeta)} \,\prod_{j \in \ufsn}\, \Iopg(|\psiinv j|)$}
\entry{$(\psi; \ang{\iphi_j}_{j \in \ufsn})$}{dgmn_iopg_obj}{object in $\DG(\mal,\nbeta)$}
\entry{$\objzero$}{dgiopg_bp}{basepoint $(0; \ang{*}_{j \in \ufsn})$ of $\DG(\mal,\nbeta)$}
\entry{$\objone_{\mal}$}{dgiopg_idonecell}{identity 1-cell $(1_{\mal}; \ang{1_{\omeg}}_{i \in \ufsm})$ of an object $\mal \in \DG$} 
\entry{$(\psi; \ang{[\ivphi_j,\iphi_j]}_{j \in \ufsn})$}{dgmn_iopg_mor}{morphism $(\psi; \ang{\iphi_j}_{j \in \ufsn}) \to (\psi; \ang{\ivphi_j}_{j \in \ufsn})$}
\entry{$\objone_{\objx}$}{dgiopg_idtwo}{identity $(\psi; \ang{1_{\iphi_j}}_{j \in \ufsn})$ of $\objx = (\psi; \ang{\iphi_j}_{j \in \ufsn})$}
\entry{$\obij_s$}{ord_bij}{order-preserving bijection $\{1,2,\ldots,|s|\} \fiso s  \subseteq \ufsm$}
\entry{$\comp$}{not:DGcomp}{composition pointed $G$-functor of $\DG$}
\entry{$\cglgc$}{cglgc}{$\Gcatst$-pseudotransformation $\gxist\Lg\Jscg\C \to \Rg\Ig\C$}
\entry{$\cglgcnbeta$}{cglgcn}{$\nbeta$-component pointed $G$-functor $(\gxist\Lg\Jscg\C)\nbeta = \Cboxn \to \Comegnbeta = (\Rg\Ig\C)\nbeta$}
\entry{$\cglgcx$}{cglgcx}{pseudonaturality constraint of $\cglgc$}
\entry{$(-)'$}{t_prime}{$\obij_s(-)$}
\entry{$\cglg$}{cglg_nat}{comparison natural transformation $\gxist\Lg\Jscg \to \Rg\Ig$}
\entry{$\cgglcnbeta$}{cgglcnbeta_functor}{pointed $G$-functor $(\Rg\Ig\C)\nbeta \to (\gxist\Lg\Jscg\C)\nbeta$}
\entry{$\cglgu$}{cglgu_nt}{unit $G$-natural isomorphism $1 \to \cgglcnbeta \cglgcnbeta$}
\entry{$\cglgv$}{cglgv_nt}{counit $G$-natural isomorphism $\cglgcnbeta\cgglcnbeta \to 1$}
\entry{$\Kiopgb$}{kiopgb_def}{GMMO $K$-theory $\Gspu\Kiopg\algi \cn \Algiopg \to \Sptopg$ for $\Iopg$}

\newchHeader{ch:mon_cat}
\entry{$\Ob\C$}{not:ObC}{class of objects in a category $\C$}
\entry{$\C(a,b)$}{not:Cab}{hom set of morphisms $a \to b$}
\entry{$\mcomp_{a,b,c}(h,f)$, $h \comp f$, $hf$}{not:mcompabc}{composition of morphisms}
\entry{$1_a$}{not:onea}{identity morphism of $a$}
\entry{$\Cop$}{def:opposites}{opposite category of $\C$}
\entry{$\C \times \C'$}{def:producat_cat}{Cartesian product}
\entry{$\Rightarrow$}{twocellnotation}{2-cell notation}
\entry{$\phi'\phi$}{not:vercomp}{vertical composition of natural transformations}
\entry{$\psi \ast \phi$}{not:horcomp}{horizontal composition  of natural transformations}
\entry{$(L,R,u,v)$}{def:adjunction}{adjunction with left adjoint $L$, right adjoint $R$, unit $u$, and counit $v$}
\entry{$\int^{a \in \C} F(a,a)$}{not:coend}{coend of a functor $F \cn \Cop \times \C \to \D$}
\entry{$(\C,\otimes,\tu,\alpha,\lambda,\rho)$}{def:monoidalcategory}{monoidal category}
\entry{$(\C,\xi)$}{def:braidedmoncat}{braided monoidal category with braiding $\xi$}
\entry{$[a,-]$}{notation:internal-hom}{internal hom}
\entry{$(F, F^2, F^0)$}{def:monoidalfunctor}{monoidal functor with monoidal constraint $F^2$ and unit constraint $F^0$}
\entry{$(a,\mu,\eta)$}{notation:monoid}{monoid with multiplication $\mu$ and unit $\eta$}
\entry{$(x,\umu)$}{def:modules_i}{$a$-module}
\entry{$(\C,\mcomp,\cone)$}{def:enriched-category}{$\V$-category with composition $\mcomp$ and identity $\cone$} 
\entry{$\A_0$}{not:Azero}{class of objects of a 2-category $\A$}
\entry{$\A_1(a,b)$}{not:Aoneab}{class of 1-cells $a \to b$}
\entry{$\A_2(f,f')$}{not:Atwo}{set of 2-cells $f \to f'$}
\entry{$1_a$}{not:idonecell}{identity 1-cell of an object $a$}
\entry{$1_f$}{not:idtwocell}{identity 2-cell of a 1-cell $f$}
\entry{$\al'\al$}{not:vcompiicell}{vertical composition of 2-cells}
\entry{$gf$}{not:hcompicell}{horizontal composition of 1-cells}
\entry{$\beta * \alpha$}{not:hcompiicell}{horizontal composition of 2-cells}
\entry{$\A(a,b)$}{not:homcat}{hom category with 1-cells $a \to b$ as objects and 2-cells as morphisms}
\entry{$(F_0,F_1,F_2)$}{def:twofunctor}{2-functor}
\entry{$\varphi\phi$}{def:twonatcomposition}{horizontal composition of 2-natural transformations}
\entry{$(L,R,u,v)$}{def:twoadjunction}{2-adjunction}
\entry{$\Psi \Phi$}{not:modvcomp}{vertical composition of modifications}
\entry{$\Phi' * \Phi$}{not:modhcomp}{horizontal composition of modifications}
\entry{$(\Op,\ga,\opu)$}{def:voperad}{$\V$-operad with operadic composition $\ga$ and operadic unit $\opu$}
\entry{$(\A,\gaA)$}{def:operadalg}{$\Op$-algebra with $\Op$-action $\gaA$}
\entry{$\AlgO$}{not:AlgO}{category of $\Op$-algebras and $\Op$-algebra morphisms}
\end{tabbing}

%% file: chap/kgo.tex
The main result in \cite[\namecref{EqK:thm:Kgo_multi} \ref*{EqK:thm:Kgo_multi}]{yau-eqk} constructs the equivariant $K$-theory $\Gcat$-multifunctors and $\Gtop$-multifunctors\index{multifunctorial K-theory@multifunctorial $K$-theory}\index{K-theory@$K$-theory!multifunctorial}
\begin{equation}\label{Kgo_multifunctor}
\begin{tikzpicture}[baseline={(a.base)}]
\def\h{2.5} \def\u{1.1} \def\v{1}
\draw[0cell]
(0,0) node (a) {\phantom{\MultpsO}}
(a)++(0,\v) node (a1) {\MultpsO}
(a)++(0,\u) node (a1'') {\phantom{\MultpsO}}
(a)++(0,\v) node (a1') {\phantom{(\Op)}} 
(a)++(0,-\v) node (a2) {\MultpspsO}
(a)++(0,-\v) node (a2') {\phantom{(\Op)}} 
(a)++(.8*\h,0) node (b) {\GGCat}
(b)++(\h,0) node (c) {\GGTop}
(c)++(.9*\h,0) node (d) {\GSp}
;
\draw[1cell=.9]
(a1') edge node[inner sep=1pt,swap] {\Jgo} (b)
(a2') edge[shorten <=.7ex] node[inner sep=1pt] {\Jgosg} (b)
(b) edge node {\clast} (c)
(c) edge node {\Kg} (d)
;
\draw[1cell=1]
(a1'') [rounded corners=2pt] -- node[pos=.7] {\Kgo} ($(c)+(0,\u)$) -| (d)
;
\draw[1cell=1]
(a2) [rounded corners=2pt] -- node[pos=.7] {\Kgosg} ($(c)+(0,-\v)$) -| (d)
;
\end{tikzpicture}
\end{equation}
for a compact Lie group $G$ and a $\Tinf$-operad $\Op$.  Each of the two composite enriched multifunctors $\Kgo$ and $\Kgosg$ sends $\Op$-pseudoalgebras to orthogonal $G$-spectra, in a way that respects the $\Gcat$-multicategory structure of its domain and the $\Gtop$-multicategory structure of $\GSp$.  \cref{ch:shim_top,part:kgo_shi_comp} compare $\Kgo$ and $\Kgosg$ with Shimakawa's equivariant $K$-theory \cite{shimakawa89,shimakawa91,mmo}.  To facilitate that comparison, this chapter reviews $\Kgo$ and $\Kgosg$ from \cite{yau-eqk}.

Since Shimakawa's equivariant $K$-theory is a functor and not a multifunctor, the comparison with our equivariant $K$-theory occurs at the functor level.  Thus, it suffices to review the underlying functors of $\Kgo$ and $\Kgosg$, as displayed in \cref{KgoKgosg_functors}.
\begin{equation}\label{KgoKgosg_functors}
\begin{tikzpicture}[baseline={(a.base)}]
\def\h{2.5} \def\u{1.1} \def\v{1}
\draw[0cell]
(0,0) node (a) {\phantom{\AlglaxO}}
(a)++(0,\v) node (a1) {\AlglaxO}
(a)++(0,\u) node (a1'') {\phantom{\AlglaxO}}
(a)++(0,\v) node (a1') {\phantom{(\Op)}} 
(a)++(0,-\v) node (a2) {\AlgpspsO}
(a)++(0,-\v) node (a2') {\phantom{(\Op)}} 
(a)++(.8*\h,0) node (b) {\GGCatii}
(b)++(\h,0) node (c) {\GGTopii}
(c)++(.9*\h,0) node (d) {\Gspec}
;
\draw[1cell=.9]
(a1') edge node[inner sep=1pt,swap] {\Jgo} (b)
(a2') edge[shorten <=.7ex] node[inner sep=1pt] {\Jgosg} (b)
(b) edge node {\clast} (c)
(c) edge node {\Kg} (d)
;
\draw[1cell=1]
(a1'') [rounded corners=2pt] -- node[pos=.7] {\Kgo} ($(c)+(0,\u)$) -| (d)
;
\draw[1cell=1]
(a2) [rounded corners=2pt] -- node[pos=.7] {\Kgosg} ($(c)+(0,-\v)$) -| (d)
;
\end{tikzpicture}
\end{equation}
The functors in \cref{KgoKgosg_functors} are obtained from the $\Gcat$-multifunctors and $\Gtop$-multifunctors in \cref{Kgo_multifunctor} by keeping the same object assignments and restricting to $G$-fixed 1-ary 1-cells and 2-cells in the $\Gcat$-enrichment or to 1-ary $G$-fixed subspaces in the $\Gtop$-enrichment.  See \cref{Kgo_functors}.

\organization
This chapter consists of the following sections.

\secname{sec:psdocom-operad} 
This section reviews the intrinsic pairing of an operad, $G$-categorical pseudo-commutative operads, the Barratt-Eccles operad $\BE$, and the $G$-Barratt-Eccles operad $\GBE$.

\secname{sec:psalg} 
This section reviews pseudoalgebras over a reduced $\Gcat$-operad; lax, pseudo, and strict $\Op$-morphisms;  $\Op$-transformations; and the 2-categories $\AlglaxO$, $\AlgpspsO$, and $\AlgstO$ with $\Op$-pseudoalgebras as objects.

\secname{sec:smGcat} 
This section reviews the 2-categories $\smgcat$, $\smgcatsg$, and $\smgcatst$ with naive symmetric monoidal $G$-categories as objects, and the 2-equivalences between them and the 2-categories of $\BE$-pseudoalgebras in \cref{sec:psalg} for the Barratt-Eccles $\Gcat$-operad $\BE$.

\secname{sec:ggcat} 
This section reviews the indexing category $\Gsk$ and the 2-category $\GGCatii$ with pointed functors $\Gsk \to \Gcatst$, called $\Gskg$-categories, as objects.

\secname{sec:jemg} 
This section begins the construction of the 2-functors $\Jgo$ and $\Jgosg$ in \cref{KgoKgosg_functors} by describing the object assignments of the $\Gskg$-categories $\Jgo\A$ and $\Jgosg\A$ for an $\Op$-pseudoalgebra $\A$.

\secname{sec:jemg_morphisms} 
This section finishes the construction of the $\Gskg$-categories $\Jgo\A$ and $\Jgosg\A$ by describing their morphism assignments.

\secname{sec:jgosg_onecells} 
This section reviews the 1-cell assignments of the 2-functors $\Jgo$ and $\Jgosg$, which send lax $\Op$-morphisms and $\Op$-pseudomorphisms to natural transformations between $\Gskg$-categories.

\secname{sec:jgosg_twocells} 
This section reviews the 2-cell assignments of the 2-functors $\Jgo$ and $\Jgosg$, which send $\Op$-transformations to modifications.

\secname{sec:ggtop} 
This section defines the category $\GGTopii$ of $\Gskg$-spaces and the functor $\clast$ in \cref{KgoKgosg_functors}, which is induced by the classifying space functor.

\secname{sec:spectra} 
This section reviews the categories $\GSp$ and $\Gspec$ of orthogonal $G$-spectra.

\secname{sec:semg} 
This section reviews the functor $\Kg$ in \cref{KgoKgosg_functors}, which sends $\Gskg$-spaces and natural transformations to orthogonal $G$-spectra and $G$-morphisms.

\section{Pseudo-Commutative Operads}\label{sec:psdocom-operad}

This section recalls the following concepts from \cite[\namecref{EqK:ch:psalg} \ref*{EqK:ch:psalg}]{yau-eqk}.
\begin{itemize}
\item \desref{Intrinsic pairing}{def:intrinsic_pairing}
\item \desref{The 2-category $\Gcat$ of small $G$-categories, $G$-functors, and $G$-natural transformations}{def:GCat,def:Catg}
\item \desref{Pseudo-commutative operads}{def:pseudocom_operad}
\item \desref{The Barratt-Eccles operad $\BE$ and naive permutative $G$-categories}{def:BE}
\item \desref{The $G$-Barratt-Eccles operad $\GBE$ and genuine permutative $G$-categories}{def:GBE}
\end{itemize}
For brief reviews of symmetric monoidal categories, 2-categories, and enriched operads, see \cref{sec:sym_mon_cat,sec:twocategories,sec:operads}.

\subsection*{Intrinsic Pairing}

The following definition is \cite[Def.\ 3.1]{gmmo23}.

\begin{definition}[Intrinsic Pairing]\label{def:intrinsic_pairing}
Suppose $(\V,\times,\bone)$ is a Cartesian closed category, and $(\Op,\ga,\opu)$ is a reduced $\V$-operad, where \dindex{reduced}{operad}\emph{reduced} means $\Op(0) = \bone$.  The \emph{intrinsic pairing}\index{intrinsic pairing}
\[(\Op,\Op) \fto{\intr} \Op\]
is the family of composites in $\V$
\begin{equation}\label{intr_jk}
\begin{tikzpicture}[vcenter]
\def\u{.65}
\draw[0cell]
(0,0) node (a1) {\Op(j) \times \Op(k)}
(a1)++(3.5,0) node (a2) {\phantom{{\Op(j) \times \Op(k)^j}}}
(a2)++(0,.03) node (a2') {\Op(j) \times \Op(k)^j}
(a2)++(2.5,0) node (a3) {\Op(jk)}
;
\draw[1cell=.9]
(a1) edge node {1 \times \Delta^j} (a2)
(a2) edge node {\ga} (a3)
;
\draw[1cell=1]
(a1) [rounded corners=2pt] |- ($(a2)+(-.5,\u)$)
-- node[pos=.0] {\intr_{j,k}} ($(a2)+(0,\u)$) -| (a3)
;
\end{tikzpicture}
\end{equation}
for $j,k \geq 0$, where $\Delta^j$ is the $j$-fold diagonal.  We often abbreviate $\intr_{j,k}$ to $\intr$.
\end{definition}

\begin{example}[Associative Operad]\label{ex:as_intrinsic}
The \dindex{associative}{operad}associative operad $\As$ \cite[Section 14.2]{yau-operad} is the $\Set$-operad with $\As(n) = \Sigma_n$, the symmetric group on $n$ letters, for each $n \geq 0$.  The right $\Sigma_n$-action on $\As(n)$ is given by the group multiplication of $\Sigma_n$.  The operadic composition
\[\Sigma_n \times \txprod_{i=1}^n \Sigma_{k_i} \fto{\ga} \Sigma_{k_1 + \cdots + k_n}\]
is given by the composite
\begin{equation}\label{as_gamma}
\ga\left(\phi; \phi_1,\ldots,\phi_n\right) 
= \phi\ang{k_1,\ldots,k_n} \circ \left(\phi_1 \oplus \cdots \oplus \phi_n\right)
\end{equation}
of the block sum\index{block sum} $\phi_1 \oplus \cdots \oplus \phi_n$ and the block permutation\index{block permutation} $\phi\ang{k_1,\ldots,k_n}$ induced by $\phi$ that permutes $n$ consecutive blocks of lengths $k_1, \ldots, k_n$.  Algebras over $\As$ are \index{monoid}monoids.

The intrinsic pairing of permutations $\sigma \in \Sigma_j$ and $\tau \in \Sigma_k$ is the composite
\begin{equation}\label{as_intrinsic}
\sigma \intr \tau = \ga\left(\sigma; \tau,\ldots,\tau\right) = \sigma\ang{k,\ldots,k} \circ \tau^j \inspace \Sigma_{jk} 
\end{equation}
of the block sum $\tau^j = \tau \oplus \cdots \oplus \tau$ and the block permutation $\sigma\ang{k,\ldots,k}$ induced by $\sigma$.  The diagram in $\V$
\begin{equation}\label{intrinsic_compatible}
\begin{tikzpicture}[xscale=3,yscale=1.3,vcenter]
\draw[0cell=.9]
(0,0) node (x11) {\Op(j) \times \Op(k)}
(x11)++(1,0) node (x12) {\Op(jk)}
(x11)++(0,-1) node (x21) {\Op(j) \times \Op(k)}
(x12)++(0,-1) node (x22) {\Op(jk)}
;
\draw[1cell=.9]  
(x11) edge node {\intr_{j,k}} (x12)
(x21) edge node {\intr_{j,k}} (x22)
(x11) edge[transform canvas={xshift=1em}] node[swap] {\sigma \times \tau} (x21)
(x12) edge node {\sigma \intr \tau} (x22)
;
\end{tikzpicture}
\end{equation}
commutes, where $\intr_{j,k}$ and $\sigma \intr \tau$ are the intrinsic pairings in \cref{intr_jk,as_intrinsic}.  Each of $\sigma$, $\tau$, and $\sigma \intr \tau$ is the right symmetric group action on $\Op$ for the indicated permutation.
\end{example}

\begin{definition}[Finite Sets]\label{def:transpose_perm}
For each $n \geq 0$, we define the \dindex{unpointed}{finite set}unpointed finite set with $n$ elements, 
\begin{equation}\label{ufsn}
\ufs{n} = \begin{cases} \{1,2,\ldots,n\} & \text{if $n \geq 1$ and}\\
\emptyset & \text{if $n=0$,}
\end{cases}
\end{equation}
equipped with its natural ordering.
\begin{itemize}
\item A sum $\txsum_{j=1}^n$ is denoted by $\txsum_{j \in \ufs{n}}$ and likewise for other operators that involve a running index.  
\item For $j,k \geq 0$, the \emph{$(j,k)$-transpose permutation}\dindex{transpose}{permutation}
\begin{equation}\label{eq:transpose_perm}
\ufs{jk} \fto[\iso]{\twist_{j,k}} \ufs{kj}
\end{equation}
in $\Sigma_{jk}$ is the bijection defined by
\[\twist_{j,k}\big(b + (a-1)k\big) = a + (b-1)j\]
for $(a,b) \in \ufs{j} \times \ufs{k}$.   A \emph{transpose permutation} is a $(j,k)$-transpose permutation for some $j,k \geq 0$.
\item The \emph{lexicographic ordering}\index{lexicographic ordering} of the product $\ufs{j} \times \ufs{k}$, or of $\ufs{jk}$, is the bijection 
\begin{equation}\label{lex_bijection}
\ufs{jk} \fto[\iso]{\lambda_{j,k}} \ufs{j} \times \ufs{k}
\end{equation}
defined by 
\[\lambda_{j,k}\big(b + (a-1)k\big) = (a,b)\]
for $(a,b) \in \ufs{j} \times \ufs{k}$.\defmark
\end{itemize}
\end{definition}

For permutations $\sigma \in \Sigma_j$ and $\tau \in \Sigma_k$, the diagram
\begin{equation}\label{intrinsic_transpose}
\begin{tikzpicture}[xscale=2.5,yscale=1.3,vcenter]
\draw[0cell=1]
(0,0) node (x11) {\ufs{jk}}
(x11)++(1,0) node (x12) {\ufs{jk}}
(x11)++(0,-1) node (x21) {\ufs{kj}}
(x12)++(0,-1) node (x22) {\ufs{kj}}
;
\draw[1cell=1]  
(x11) edge node {\sigma \intr \tau} (x12)
(x21) edge node {\tau \intr \sigma} (x22)
(x11) edge node[swap] {\twist_{j,k}} (x21)
(x12) edge node {\twist_{j,k}} (x22)
;
\end{tikzpicture}
\end{equation}
 commutes, since each composite is given by the assignment
\[\begin{tikzcd}
b + (a-1)k \rar[maps to] & \sigma(a) + \big(\tau(b) - 1\big)j
\end{tikzcd}\]
for $(a,b) \in \ufs{j} \times \ufs{k}$.

\subsection*{Pseudo-Commutative Operads}

\begin{definition}[Equivariant Categories]\label{def:GCat}
$\Cat$ denotes the 2-category of small categories, functors, and natural transformations.  The underlying 1-category of small categories and functors is also denoted by $\Cat$.
\begin{itemize}
\item Each group $G$ is regarded as a 2-category with one object $*$, set of 1-cells $G$, only identity 2-cells, and horizontal composition of 1-cells given by the group multiplication.  
\item $\Gcat$\label{not:Gcat} denotes the 2-category with
\begin{itemize}
\item 2-functors $G \to \Cat$ as objects, called \index{G-category@$G$-category}\index{category!equivariant}\index{equivariant!category}\emph{small $G$-categories};
\item 2-natural transformations as 1-cells, called \index{G-functor@$G$-functor}\index{functor!equivariant}\index{equivariant!functor}\emph{$G$-functors};
\item modifications as 2-cells, called \index{G-natural transformation@$G$-natural transformation}\index{natural transformation!equivariant}\index{equivariant!natural transformation}\emph{$G$-natural transformations}; and
\item horizontal and vertical compositions induced by those of the 2-category $\Cat$.
\end{itemize}
\end{itemize}
The underlying 1-category of $\Gcat$ is equipped with the Cartesian symmetric monoidal structure, diagonal $G$-action, and monoidal unit \label{not:bone}$\boldone$, which is the terminal $G$-category with only one object $*$ and its identity morphism.
\end{definition}

\begin{explanation}[Unpacking $\Gcat$]\label{expl:GCat}
A small \emph{$G$-category} consists of a small category $\C$ and a \emph{$g$-action} isomorphism
\begin{equation}\label{gactioniso}
\C \fto[\iso]{g} \C \foreachspace g \in G
\end{equation}
such that the following two statements hold.
\begin{enumerate}
\item For the identity element $e \in G$, the $e$-action is the identity functor $1_\C$.
\item For $g, h \in G$, the $hg$-action is equal to the composition $h \circ g$.
\end{enumerate}
For $g \in G$ and an object or morphism $x \in \C$, the $g$-action on $x$ is written as either $g \cdot x$ or $gx$.

A \emph{$G$-functor} $F \cn \C \to \D$ between small $G$-categories is a functor such that the following diagram commutes for each $g \in G$.
\begin{equation}\label{Gfunctor}
\begin{tikzpicture}[xscale=2,yscale=1.3,vcenter]
\draw[0cell=.9]
(0,0) node (x11) {\C}
(x11)++(1,0) node (x12) {\D}
(x11)++(0,-1) node (x21) {\C}
(x12)++(0,-1) node (x22) {\D}
;
\draw[1cell=.9]  
(x11) edge node {F} (x12)
(x21) edge node {F} (x22)
(x11) edge node[swap] {g} (x21)
(x12) edge node {g} (x22)
;
\end{tikzpicture}
\end{equation}
For $G$-functors $F,F' \cn \C \to \D$, a \emph{$G$-natural transformation} $\theta \cn F \to F'$ is a natural transformation such that the morphism equality
\begin{equation}\label{Gnattr}
g \cdot Fc = F(gc) \fto{g \cdot \theta_c = \theta_{gc}} g \cdot F'c = F'(gc)
\end{equation}
holds for each $g \in G$ and object $c \in \C$.  The conditions \cref{Gfunctor,Gnattr} are called the \index{G-natural transformation@$G$-natural transformation!G-equivariance@$G$-equivariance}\emph{$G$-equivariance conditions}.  A $G$-functor or a $G$-natural transformation is said to be \emph{$G$-equivariant}.
\end{explanation}

\begin{definition}[Internal Hom]\label{def:Catg}
Suppose $\C$ and $\D$ are small $G$-categories for a group $G$.  The small $G$-category \index{category!internal hom}\index{internal hom}$\Catg(\C,\D)$ has
\begin{itemize}
\item functors $\C \to \D$ as objects,
\item natural transformations as morphisms, and
\item identities and composition defined componentwise in $\D$.
\end{itemize}
Its \emph{conjugation $G$-action}\index{conjugation G-action@conjugation $G$-action} is defined by the composite and whiskering
\begin{equation}\label{conjugation-gaction}
\begin{tikzpicture}[xscale=1.5,yscale=1.5,vcenter]
\draw[0cell=.9]
(0,0) node (a) {\C}
(a)++(1,0) node (b) {\C}
(b)++(1.2,0) node (c) {\D}
(c)++(1,0) node (d) {\D}
;
\draw[1cell=.9]  
(a) edge node {g^{-1}} (b)
(b) edge[bend left] node {F} (c)
(b) edge[bend right] node[swap] {F'} (c)
(c) edge node {g} (d)
;
\draw[2cell]
node[between=b and c at .43, rotate=-90, 2label={above,\theta}] {\Rightarrow}
;
\end{tikzpicture}
\end{equation}
for functors $F,F' \cn \C \to \D$, a natural transformation $\theta \cn F \to F'$, and an element $g \in G$.  The quadruple
\begin{equation}\label{gcat-closed}
\left(\Gcat,\times,\bone,\Catg\right)
\end{equation}
is a Cartesian closed category, with limits and colimits computed in $\Cat$.  Passing to the \index{G-fixed subcategory@$G$-fixed subcategory}$G$-fixed subcategory, there is an equality
\begin{equation}\label{catg_gfixed}
\Catg(\C,\D)^G = \Gcat(\C,\D),
\end{equation}
where $\Gcat(\C,\D)$ is the category of $G$-functors and $G$-natural transformations.
\end{definition}

The following definition is \cite[Def.\ 3.10]{gmmo23}.

\begin{definition}[Pseudo-Commutative Operads]\label{def:pseudocom_operad}
For a group $G$, a \emph{pseudo-commutative operad}\dindex{pseudo-commutative}{operad} $(\Op,\ga,\opu,\pcom)$ in $\Gcat$ consists of the following data.
\begin{description}
\item[Reduced operad] 
$(\Op,\ga,\opu)$ is a reduced $\Gcat$-operad \cite[\ref*{EqK:def:enr-multicategory}]{yau-eqk} for the Cartesian closed category $\Gcat$, where \emph{reduced} means $\Op(0)$ is a terminal $G$-category $\bone$.
\item[Pseudo-commutative structure] 
For $j,k \geq 0$, it is equipped with a $G$-natural isomorphism $\pcom_{j,k}$, called the \emph{$(j,k)$-pseudo-commutativity isomorphism}, as follows.
\begin{equation}\label{pseudocom_isos}
\begin{tikzpicture}[xscale=3,yscale=1.3,vcenter]
\draw[0cell=.9]
(0,0) node (x11) {\Op(j) \times \Op(k)}
(x11)++(1,0) node (x12) {\Op(jk)}
(x11)++(0,-1) node (x21) {\Op(k) \times \Op(j)}
(x12)++(0,-1) node (x22) {\Op(kj)}
;
\draw[1cell=.9]  
(x11) edge node {\intr_{j,k}} (x12)
(x21) edge node[swap] {\intr_{k,j}} (x22)
(x11) edge node[swap] {\twist} (x21)
(x12) edge node {\twist_{k,j}} (x22)
;
\draw[2cell]
node[between=x11 and x22 at .5, rotate=-90, 2label={above,\pcom_{j,k}}, 2label={below,\iso}] {\Rightarrow}
;
\end{tikzpicture}
\end{equation}
The collection $\pcom = \{\pcom_{j,k}\}_{j,k \geq 0}$ is called the \emph{pseudo-commutative structure}.  The arrow $\twist$ is the braiding for $\Gcat$, swapping the two arguments.  Each of $\intr_{j,k}$ and $\intr_{k,j}$ is the intrinsic pairing \cref{intr_jk}.  The arrow $\twist_{k,j}$ is the right symmetric group action $G$-functor on $\Op$ for the $(k,j)$-transpose permutation \cref{eq:transpose_perm}.   For objects $(x,y) \in \Op(j) \times \Op(k)$, the $(x,y)$-component of $\pcom_{j,k}$ is an isomorphism
\begin{equation}\label{pcom_xy}
(x \intr_{j,k} y) \twist_{k,j} \fto[\iso]{\pcom_{j,k; x,y}} y \intr_{k,j} x \inspace \Op(kj).
\end{equation}
The subscripts of $\pcom$ and $\intr$ are sometimes omitted to simplify the notation.
\end{description}
The quadruple $(\Op,\ga,\opu,\pcom)$ is required to satisfy the following four axioms whenever they are defined.
\begin{description}
\item[Unity] 
For the operadic unit $\opu \in \Op(1)$ and an object $y \in \Op(k)$, using $\twist_{k,1} = \id_k \in \Si_k$ and the unity axioms of $\Op$, there is a morphism equality
\begin{equation}\label{pseudocom_unity}
(\opu \intr_{1,k} y) \twist_{k,1} = y \fto{\pcom_{1,k} = 1_y} y \intr_{k,1} \opu = y \inspace \Op(k).
\end{equation}

\item[Symmetry] 
Using $\twist \circ \twist = 1$ and $\twist_{k,j} = \twist_{j,k}^{-1}$, the following composite is equal to the identity $G$-natural transformation of $\intr_{j,k}$.
\begin{equation}\label{pseudocom_sym}

\end{equation}

\item[Equivariance] 
The following pasting diagram equality holds for permutations $\sigma \in \Sigma_j$ and $\tau \in \Sigma_k$.
\begin{equation}\label{pseudocom_equiv}
%
\end{equation}
The two unlabeled regions commute by \cref{intrinsic_compatible}.  The left vertical boundaries are equal by the naturality of the braiding $\twist$ for $\Gcat$.  The right vertical boundaries are equal by \cref{intrinsic_transpose}.

\item[Operadic compatibility] 
Consider objects $x \in \Op(j)$, $y_i \in \Op(k_i)$ for $1 \leq i \leq j$, and $z \in \Op(\ell)$, along with the following notation for $r \in \{1,\ldots,j\}$.
\[\left\{\scalebox{.85}{$\begin{split}
k &= k_1 + \cdots + k_j \qquad y_r^\ell = (\overbracket[.5pt]{y_r,\ldots,y_r}^{\ell}) \in \Op(k_r)^\ell \\
y_\crdot^\ell &= \big(y_1^\ell, \ldots, y_j^\ell\big)\\
\twist_{\ell,k_\centerdot} &= \twist_{\ell,k_1} \times \cdots \times \twist_{\ell,k_j} \in \Sigma_{\ell k} \\
\Twist_{\ell,k_\centerdot} &= \big(\twist_{k_1,\ell} \times \cdots \times \twist_{k_j,\ell}\big) \twist_{\ell,k} \in \Sigma_{\ell k}\\
\boldy &= (y_1,\ldots,y_j) \in \Op(k_1) \times \cdots \times \Op(k_j)\\
\boldy^\ell &= (\overbracket[.5pt]{\boldy,\ldots,\boldy}^{\ell}) \in \left(\Op(k_1) \times \cdots \times \Op(k_j)\right)^\ell\\
z \intr \boldy & = \left(z \intr y_1, \ldots, z \intr y_j\right) \in \Op(\ell k_1) \times \cdots \times \Op(\ell k_j)\\ 
\boldy \intr z &= \left(y_1 \intr z, \ldots, y_j \intr z \right) \in \Op(k_1 \ell) \times \cdots \times \Op(k_j \ell)\\
(\bold y \intr z) \twist_{\ell,\bdot} &= \big((y_1 \intr z)\twist_{\ell,k_1}, \ldots, (y_j \intr z)\twist_{\ell,k_j}\big) \in \Op(\ell k_1) \times \cdots \times \Op(\ell k_j)\\
\pcom_{k_\centerdot,\ell} &= \big( \pcom_{k_1,\ell}, \ldots, \pcom_{k_j,\ell}\big)\\
\end{split}$}\right.\]
Then the following diagram in $\Op(\ell k)$ commutes.
\begin{equation}\label{pseudocom_compat}
\begin{tikzpicture}[xscale=1,yscale=1,vcenter]
\def\h{2.5} \def\f{.85} \def\u{-1.2} \def\v{-1.3}
\draw[0cell=.9]
(0,0) node (x) {\ga\big(x \sscs (\bold y \intr z) \twist_{\ell,\bdot}\big) \Twist_{\ell,k_\centerdot}}
(x)++(-\h,\u) node (x11) {\ga\left(x \sscs z \intr \boldy\right) \Twist_{\ell,k_\centerdot}}
(x11)++(0,\v) node (x21) {\ga\big(x \intr z \sscs y_\crdot^\ell\big) \Twist_{\ell,k_\centerdot}}
(x21)++(0,\v) node (x31) {\ga\big((x \intr z) \twist_{\ell,j} \sscs \boldy^\ell\big)}
(x31)++(\f,\v) node (x41) {\ga\big(z \intr x \sscs \boldy^\ell\big)}
(x)++(\h,\u) node (x12) {\ga\left(x \sscs \boldy \intr z\right) (\twist_{\ell, k_\centerdot} \Twist_{\ell,k_\centerdot})}
(x12)++(0,\v) node (x22) {\ga\left(x \sscs \boldy \intr z\right) \twist_{\ell,k}}
(x22)++(0,\v) node (x32) {\big(\ga\left(x \sscs \boldy\right) \intr z\big) \twist_{\ell,k}}
(x32)++(-\f,\v) node (x42) {z \intr \ga\left(x \sscs \boldy\right)}
;
\draw[1cell=.8]  
(x) edge[transform canvas={xshift={-1.3em}}] node[swap,pos=.5] {\ga\left(1 \sscs \pcom_{k_\centerdot,\ell}\right) \Twist_{\ell,k_\centerdot}} (x11)
(x11) edge[equal] node[swap] {\mathrm{(a)}} (x21)
(x21) edge[equal] node[swap] {\mathrm{(top)}} (x31)
(x31) edge[shorten <=-2pt, shorten >=-2pt] node[swap, pos=.2] {\ga( \pcom_{j,\ell} \sscs 1^{\ell j} )} (x41)
(x41) edge[equal] node {\mathrm{(a)}} (x42)
(x) edge[equal, transform canvas={xshift={1.3em}}] node {\mathrm{(bot)}} (x12)
(x12) edge[equal] node {\dagger} (x22)
(x22) edge[equal] node {\mathrm{(a)}} (x32)
(x32) edge[shorten <=-1pt, shorten >=-1pt] node[pos=.3] {\pcom_{k,\ell}} (x42)
;
\end{tikzpicture}
\end{equation}
The three equalities labeled $\mathrm{(a)}$ hold by the associativity axiom of $\Op$.  The equalities labeled $\mathrm{(top)}$ and $\mathrm{(bot)}$ hold by the top and bottom equivariance axioms of $\Op$.  The equality labeled $\dagger$ holds by the equality
\[\twist_{\ell, k_\centerdot} \Twist_{\ell,k_\centerdot} = \twist_{\ell,k},\]
using the fact that $\twist_{\ell,k_r} = \twist_{k_r,\ell}^{-1}$ for $r \in \{1,\ldots,j\}$.
\end{description}
This finishes the definition of a pseudo-commutative operad in $\Gcat$.  
\end{definition}

\subsection*{Barratt-Eccles Operads and Permutative $G$-Categories}

\begin{definition}[Translation Categories]\label{def:translation_cat}
For a set $S$, the \emph{translation category}\dindex{translation}{category} $\tn S$ is the small category with object set $S$ and each hom set given by the one-element set $*$.  Identity morphisms and composition are uniquely determined by the fact that $*$ is a terminal object in the category $\Set$.  For objects $a, b \in \tn S$, the unique morphism from $a$ to $b$ is usually denoted by either\label{not:b_to_a}
\[ a \fto{[b,a]} b \orspace a \fto{!} b.\]
For a group $G$, the translation category $\EG$ is a $G$-category with $G$ acting by left multiplication on objects.
\end{definition}

In the literature, a translation category is also called an \dindex{indiscrete}{category}\emph{indiscrete category} and a \dindex{chaotic}{category}\emph{chaotic category}.  The following observation from \cite[Cor.\ 4.9]{corner-gurski} and \cite[Lemma 3.12]{gmmo23} is a direct consequence of the terminal property of a one-element set in the category $\Set$.

\begin{proposition}\label{translation_pseudocom}
Suppose $\Op$ is a reduced $\Gcat$-operad for a group $G$ such that each $\Op(n)$ is a translation category.  Then $\Op$ admits a unique pseudo-commutative structure.
\end{proposition}

\begin{definition}[Barratt-Eccles Operad]\label{def:BE}
For a group $G$, the \emph{Barratt-Eccles operad}\dindex{Barratt-Eccles}{operad} $\BE$ is the pseudo-commutative operad in $\Gcat$ given by the translation categories
\[\BE(n) = \tn\As(n) = \tn\Sigma_n \forspace n \geq 0\]
with the trivial $G$-action.  The operad structure on objects is the one on the associative operad \cref{as_gamma}.  On morphisms, the operad structure is uniquely determined by the terminal property of each hom set in $\tn\Sigma_n$.  A $\BE$-algebra in $\Gcat$ is called a \index{naive permutative G-category@naive permutative $G$-category}\index{G-category@$G$-category!naive permutative}\emph{naive permutative $G$-category}.
\end{definition}

\begin{explanation}[Naive Permutative $G$-Categories]\label{expl:naive_perm_Gcat}
A naive permutative $G$-category is, by \cref{def:BE}, a $\Gcat$-multifunctor 
\[\upga \cn \BE \to \Gcat.\]  
The unique object $* \in \BE$ is sent to a small $G$-category $\upga(*) = \C$.  Since $G$ acts trivially on $\BE$, the $G$-functor at level 0
\[\upga \cn \BE(0) = \tn\Sigma_0 \iso \boldone \to \Catg(\bone,\C)\]
corresponds to a $G$-fixed object $\pu \in \C$.  The $G$-functor at level 2
\[\upga \cn \BE(2) = \tn\Sigma_2 \to \Catg(\C \times \C, \C)\]
sends the identity permutation $\id_2 \in \Sigma_2$ to a $G$-functor \cref{Gfunctor}
\[\upga(\id_2) = \oplus \cn  \C \times \C \to \C.\]
For the nonidentity permutation $\tau \in \Sigma_2$, the unique nonidentity isomorphism 
\[ [\tau,\id_2] \cn \id_2 \fto{\iso} \tau \inspace \tn\Sigma_2\]
is sent by $\upga$ to a $G$-natural isomorphism $\xi$ \cref{Gnattr} as follows, with $\twist$ swapping the two arguments.  
\begin{equation}\label{npgc_braiding}
\begin{tikzpicture}[xscale=3, yscale=1.1,vcenter]
\def\v{-1} \def\h{1} \def\m{1} \def\q{15}
\draw[0cell=1] 
(0,0) node (x11) {\C \times \C}
(x11)++(\h,0) node (x12) {\C}
(x11)++(\h/2,\v) node (x2) {\phantom{\C \times \C}}
(x2)++(0,.1) node (x2') {\C \times \C}
;
\draw[1cell] 
(x11) edge node [pos=.45](i) {\oplus} (x12)
(x11) edge[bend right=\q] node[swap,pos=.5] (a) {\twist} (x2)
(x2) edge[bend right=\q] node[swap,pos=.5] {\oplus} (x12)
;
\draw[2cell] 
node[between=i and x2 at .47, shift={(0,0)}, rotate=-90, 2label={above,\xi}] {\Rightarrow}
;
\end{tikzpicture}
\end{equation}
The quadruple $(\C,\oplus,\pu,\xi)$ is a permutative category because the necessary axioms hold in $\BE$.  Coherence properties of $\BE$ \cite[11.4.14]{cerberusIII} imply that the $\Gcat$-multifunctor $\upga$ is completely determined by the data $(\C,\oplus,\pu,\xi)$.  

In summary, a naive permutative $G$-category is precisely a small $G$-category $\C$ together with a permutative category structure $(\oplus,\pu,\xi)$ such that $\oplus$ and the braiding $\xi$ are $G$-equivariant and that the monoidal unit $\pu$ is $G$-fixed.
\end{explanation}

The genuine variant of the Barratt-Eccles operad is introduced by Shimakawa \cite{shimakawa89}.  See also \cite{gm17,gmm17}.

\begin{definition}[$G$-Barratt-Eccles Operad]\label{def:GBE}
Suppose $G$ is a group.  The \emph{$G$-Barratt-Eccles operad}\index{G-Barratt-Eccles operad@$G$-Barratt-Eccles operad}\index{operad!G-Barratt-Eccles@$G$-Barratt-Eccles} $\GBE$ is the pseudo-commutative operad in $\Gcat$ given by the small $G$-categories
\[\GBE(n) = \Catg(\EG, \BE(n)) 
= \Catg(\EG, \tn \Sigma_n) 
\iso \tn[G,\Sigma_n] \forspace n \geq 0.\]
Here $\Catg(-,-)$ is the internal hom of $\Gcat$ \pcref{def:Catg}, and $[G,\Sigma_n]$ is the $G$-set of functions $G \to \Sigma_n$ with $G$ acting trivially on $\Sigma_n$.
\begin{itemize}
\item The $\Gcat$-operadic unit of $\GBE$ is given by the unique $G$-functor 
\[\EG \fto{\opu} \tn\Sigma_1 \iso \boldone.\] 
\item The right $\Sigma_n$-action on $\GBE(n)$ is induced by the one on $\tn\Sigma_n$.
\item The $\Gcat$-operadic composition of $\GBE$ is induced by the $\Cat$-operadic composition of the Barratt-Eccles operad $\BE$ (\cref{def:BE}), using the fact that $\Catg(\EG,-)$ preserves products.
\end{itemize}
A $\GBE$-algebra in $\Gcat$ is called a \index{genuine permutative G-category@genuine permutative $G$-category}\index{G-category@$G$-category!genuine permutative}\emph{genuine permutative $G$-category}.
\end{definition}

The following observation relating naive and genuine permutative $G$-categories is due to Shimakawa \cite[page 256]{shimakawa89}.  See also \cite[Prop.\ 4.6]{gm17}.

\begin{proposition}\label{naive_genuine_pGcat}
For each group $G$, each naive permutative $G$-category $\C$ yields a genuine permutative $G$-category $\Catg(\EG,\C)$.
\end{proposition}

\section{Operadic Pseudoalgebras}\label{sec:psalg}

This section recalls the following concepts from \cite[\namecref{EqK:ch:psalg} \ref*{EqK:ch:psalg}]{yau-eqk}.
\begin{itemize}
\item \desref{Pseudoalgebras of a reduced $\Gcat$-operad}{def:pseudoalgebra}
\item \desref{Lax $\Op$-morphisms and $\Op$-transformations}{def:laxmorphism,def:algtwocells}
\item \desref{The 2-categories $\AlglaxO$, $\AlgpspsO$, and $\AlgstO$ with $\Op$-pseudoalgebras as objects}{oalgps_twocat}
\end{itemize}

\subsection*{Pseudoalgebras}

Operadic pseudoalgebras are introduced in \cite[Def.\ 2.14]{gmmo20}.  Pseudoalgebras generalize algebras \pcref{def:operadalg} in the same way that monoidal categories generalize strict monoidal categories.  Recall that a $\Gcat$-operad $\Op$ is \emph{reduced} if $\Op(0)$ is a terminal $G$-category.

\begin{definition}\label{def:pseudoalgebra}
Suppose $G$ is a group, and $(\Op,\ga,\opu)$ is a reduced $\Gcat$-operad.  An \emph{$\Op$-pseudoalgebra}\index{pseudoalgebra}\index{operad!pseudoalgebra} is a triple $(\A,\gaA,\phiA)$ consisting of the following data.
\begin{description}
\item[Underlying $G$-category] $\A$ is a small $G$-category.
\item[$\Op$-action] It is equipped with a $G$-functor
\begin{equation}\label{gaAn}
\Op(n) \times \A^n \fto{\gaA_n} \A
\end{equation}
for each $n \geq 0$, called the \index{pseudoalgebra!action G-functor@action $G$-functor}\index{action G-functor@action $G$-functor}\emph{$n$-th $\Op$-action $G$-functor}.  When $n=0$, the $G$-functor
\[\Op(0) = \boldone \fto{\gaA_0} \A\]
specifies a $G$-fixed object
\begin{equation}\label{pseudoalg_zero}
\zero = \gaA_0(*) \in \A,
\end{equation} 
called the \emph{basepoint}\index{pseudoalgebra!basepoint}\index{basepoint} of $\A$.  We sometimes abbreviate $\gaA_n$ to $\gaA$.
\item[Associativity constraint] For each sequence of integers
\[\big(n; \ang{m_j}_{j \in \ufs{n}}\big) = \big(n; m_1,\ldots,m_n\big)\] 
with $n > 0$ and each $m_j \geq 0$, it is equipped with a $G$-natural isomorphism, called the \emph{associativity constraint}\index{pseudoalgebra!associativity constraint}\index{associativity constraint} of $\A$, as follows, where $m = \sum_{j=1}^n m_j$.
\begin{equation}\label{phiA}
\begin{tikzpicture}[xscale=1,yscale=1,vcenter]
\def\h{5.5} \def\g{.5} \def\v{-1.3}
\draw[0cell=.9]
(0,0) node (x11) {\Op(n) \times \txprod_{j=1}^n \big( \Op(m_j) \times \A^{m_j} \big)}
(x11)++(\h,0) node (x12) {\Op(n) \times \A^n}
(x12)++(\g,\v) node (x) {\A}
(x11)++(0,2*\v) node (x21) {\big(\Op(n) \times \txprod_{j=1}^n \Op(m_j) \big) \times \A^{m_1 + \cdots + m_n}}
(x21)++(\h,0) node (x22) {\Op(m) \times \A^m}
;
\draw[1cell=.9]  
(x11) edge node {1 \times \txprod_{j=1}^n\, \gaA_{m_j}} (x12)
(x12) edge node[pos=.7] {\gaA_n} (x)
(x11) edge node[swap] {\pi} (x21)
(x21) edge node {\ga \times 1} (x22)
(x22) edge node[swap,pos=.8] {\gaA_m} (x)
;
\draw[2cell]
node [between=x11 and x22 at .5, shift={(0,.2)}, rotate=-90, 2label={below,\iso}, 2label={above,\phiA_{(n;\, m_1,\ldots,m_n)}}] {\Rightarrow}
;
\end{tikzpicture}
\end{equation}
The arrow $\pi$ shuffles the $\Op(m_j)$ factors to the left, keeping their relative order unchanged. 
A typical component of $\phiA_{(n;\, \ang{m_j}_{j \in \ufs{n}})}$ is an isomorphism
\begin{equation}\label{phiA_component}
\begin{tikzpicture}[xscale=5.8,yscale=1,baseline={(a.base)}]
\draw[0cell=.8]
(0,0) node (a) {\gaA_n\Big(x \sscs \bang{ \gaA_{m_j}\big(x_j \sscs \ang{a_{j,i}}_{i\in \ufs{m}_j} \big) }_{j\in \ufs{n}}\Big)}
(a)++(1,0) node (b) {\gaA_m\Big( \ga\big( x \sscs \ang{x_j}_{j\in \ufs{n}} \big) \sscs \bang{\ang{a_{j,i}}_{i\in \ufs{m}_j} }_{j\in \ufs{n}} \Big)}
;
\draw[1cell=.8]
(a) edge node {\phiA_{(n;\, \ang{m_j}_{j \in \ufs{n}})}} node[swap] {\iso} (b)
;
\end{tikzpicture}
\end{equation}
in $\A$ for objects $x \in \Op(n)$, $x_j \in \Op(m_j)$, and $a_{j,i} \in \A$ for $1 \leq j \leq n$ and $1 \leq i \leq m_j$.  By convention, $\phiA_{(0;\ang{})}$ is the identity morphism $1_\zero$ of the basepoint $\zero \in \A$.
\end{description}
The data above are required to satisfy the axioms  \cref{pseudoalg_action_sym,pseudoalg_action_unity,pseudoalg_basept_axiom,pseudoalg_topeq,pseudoalg_boteq,pseudoalg_unity,pseudoalg_comp_axiom} below.
\begin{description}
\item[Action equivariance] 
For each permutation $\sigma \in \Sigma_n$ with $n \geq 2$, the diagram of $G$-functors
\begin{equation}\label{pseudoalg_action_sym}
\begin{tikzpicture}[xscale=3,yscale=1.3,vcenter]
\draw[0cell=.9]
(0,0) node (x11) {\Op(n) \times \A^n}
(x11)++(1,0) node (x12) {\Op(n) \times \A^n}
(x11)++(0,-1) node (x21) {\Op(n) \times \A^n}
(x21)++(1,0) node (x22) {\A}
;
\draw[1cell=.9]  
(x11) edge node {1 \times \sigma} (x12)
(x12) edge node {\gaA_n} (x22)
(x11) edge[transform canvas={xshift=1em}] node[swap] {\sigma \times 1} (x21)
(x21) edge node {\gaA_n} (x22)
;
\end{tikzpicture}
\end{equation}
commutes. The top $\sigma$ permutes the $n$ copies of $\A$ from the left:
\[\sigma(a_1,a_2,\ldots,a_n) = \big(a_{\sigmainv(1)}, a_{\sigmainv(2)}, \ldots, a_{\sigmainv(n)}\big).\]
The left vertical $\sigma$ in \cref{pseudoalg_action_sym} is the right $\sigma$-action on $\Op(n)$.
\item[Action unity] 
The following diagram commutes. 
\begin{equation}\label{pseudoalg_action_unity}
\begin{tikzpicture}[xscale=2.5,yscale=1.3,vcenter]
\draw[0cell=.9]
(0,0) node (x11) {\boldone \times \A}
(x11)++(1,0) node (x12) {\A}
(x11)++(0,-1) node (x21) {\Op(1) \times \A}
(x21)++(1,0) node (x22) {\A}
;
\draw[1cell=.9]  
(x11) edge node {\iso} (x12)
(x12) edge[equal] (x22)
(x11) edge node[swap] {\opu \times 1} (x21)
(x21) edge node {\gaA_1} (x22)
;
\end{tikzpicture}
\end{equation}
\item[Basepoint] 
Consider objects $x \in \Op(n)$,
\[\anga = \big(a_1,\ldots,a_{j-1}, a_{j+1}, \ldots, a_n) \in \A^{n-1}\] 
for any $j \in \{1,\ldots,n\}$,\label{not:dyj}
\[\begin{split}
\dy_j x &= \ga\big(x \sscs \opu^{j-1}, *, \opu^{n-j}\big) \in \Op(n-1), \andspace\\
\dy_j \anga &= \big(a_1,\ldots,a_{j-1}, \zero, a_{j+1}, \ldots, a_n) \in \A^n
\end{split}\]
with $* \in \Op(0) = \boldone$ the unique object, $\opu \in \Op(1)$ the operadic unit, $\opu^k$ the $k$-tuple of copies of $\opu$, and $\zero = \gaA_0(*)$ the basepoint of $\A$ \cref{pseudoalg_zero}.
Then the following component of $\phiA$ \cref{phiA_component} is the identity morphism.
\begin{equation}\label{pseudoalg_basept_axiom} 
\gaA_n\big( x \sscs \dy_j\anga \big) \fto{\phiA_{(n;\, 1^{j-1},0,1^{n-j})} \,=\, 1} 
\gaA_{n-1}\big( \dy_j x \sscs \anga \big)
\end{equation}
This component of $\phiA$ has $m_j = 0$, $x_j = * \in \Op(0)$, $m_r = 1$ for $r \neq j$, and $x_r = \opu \in \Op(1)$.  The domain of $\phiA$ uses the equality
\[\gaA_n\big( x \sscs \dy_j\anga \big) = 
\gaA_n\left( x \sscs \bang{\gaA_1(\opu; a_i)}_{i=1}^{j-1} \scs \gaA_0(*) \scs \bang{\gaA_1(\opu; a_i)}_{i=j+1}^n \right)\]
in $\A$, which holds by the action unity axiom \cref{pseudoalg_action_unity}.

\item[Top equivariance] 
For each permutation $\sigma \in \Sigma_n$ and objects $x \in \Op(n)$, $x_j \in \Op(m_j)$, and $a_{j,i} \in \A$ as in \cref{phiA_component}, the following diagram in $\A$ commutes.
\begin{equation}\label{pseudoalg_topeq}
\begin{tikzpicture}[xscale=1,yscale=1,vcenter]
\def\h{3} \def\g{0} \def\v{1.4} \def\u{-1}
\draw[0cell=.8]
(0,0) node (x11) {\gaA_n\Big(x\sigma \sscs \big\langle \gaA_{m_{\sigma(j)}}\big(x_{\sigma(j)} \sscs \ang{a_{\sigma(j), i}}_{i \in \ufs{m}_{\sigma(j)}} \big) \big\rangle_{j\in \ufs{n}}\Big)}
(x11)++(\h,\v) node (x12) {\gaA_m\left( \ga\big( x\sigma \sscs \ang{x_{\sigma(j)}}_{j\in \ufs{n}} \big) \sscs \bang{\ang{a_{\sigma(j),i}}_{i \in \ufs{m}_{\sigma(j)}} }_{j\in \ufs{n}} \right)}
(x11)++(\h+.5,\u) node (x) {\gaA_m\left( \ga\big( x \sscs \ang{x_j}_{j\in \ufs{n}} \big) \sigmabar \sscs \bang{\ang{a_{\sigma(j),i}}_{i \in \ufs{m}_{\sigma(j)}} }_{j\in \ufs{n}} \right)}
(x11)++(0,2*\u) node (x21) {\gaA_n\Big(x \sscs \big\langle \gaA_{m_j}\big(x_j \sscs \ang{a_{j,i}}_{i \in \ufs{m}_j} \big) \big\rangle_{j\in \ufs{n}}\Big)}
(x21)++(\h,-\v) node (x22) {\gaA_m\left( \ga\big( x \sscs \ang{x_j}_{j\in \ufs{n}} \big) \sscs \bang{\ang{a_{j,i}}_{i \in \ufs{m}_j} }_{j\in \ufs{n}} \right)}
;
\draw[1cell=.85]  
(x11) edge[transform canvas={xshift={-2em}}] node[pos=.2] {\phiA_{(n;\, m_{\sigma(1)}, \ldots, m_{\sigma(n)})}} (x12)
(x12) edge[equal] node[pos=.6] {\spadesuit} (x)
(x) edge[equal] (x22)
(x11) edge[equal, transform canvas={xshift={-1em}}](x21)
(x21) edge[transform canvas={xshift={-2em}}] node[swap,pos=.2] {\phiA_{(n;\, m_1,\ldots,m_n)}} (x22)
;
\end{tikzpicture}
\end{equation}
The equality labeled $\spadesuit$ holds by the top equivariance axiom for $\Op$, where
\[\sigmabar = \sigma\bang{m_{\sigma(1)}, \ldots, m_{\sigma(n)}} \in \Sigma_{m_1+\cdots+m_n}\]
is the block permutation induced by $\sigma$ that permutes $n$ consecutive blocks of lengths $m_{\sigma(1)}, \ldots, m_{\sigma(n)}$.  The other two equalities hold by the action equivariance axiom \cref{pseudoalg_action_sym}.

\item[Bottom equivariance] 
For permutations $\tau_j \in \Sigma_{m_j}$ for $1 \leq j \leq n$, the following diagram in $\A$ commutes.
\begin{equation}\label{pseudoalg_boteq}

\end{equation}
The equality labeled $\clubsuit$ holds by the bottom equivariance axiom for $\Op$, where
\[\tau^\times = \tau_1 \oplus \cdots \oplus \tau_n\in \Sigma_{m_1+\cdots+m_n}\]
is the block sum.  The other two equalities hold by the action equivariance axiom \cref{pseudoalg_action_sym}.

\item[Unity] 
For objects $x \in \Op(n)$ and $\anga = \ang{a_j}_{j\in \ufs{n}} \in \A^n$, the following diagram in $\A$ commutes, where $1^n = \ang{1}_{j \in \ufs{n}}$.
\begin{equation}\label{pseudoalg_unity}
%
\end{equation}
The two equalities along the left boundary hold by the action unity axiom \cref{pseudoalg_action_unity}.  The two equalities along the right boundary hold by the unity axioms for $\Op$.

\item[Composition] 
Consider objects 
\[x \in \Op(n), \quad x_j \in \Op(m_j), \quad x_{j,i} \in \Op(k_{j,i}), \andspace a_{j,i,h} \in \A\]
for $1 \leq j \leq n$, $1 \leq i \leq m_j$, and $1 \leq h \leq k_{j,i}$, along with the following notation, where $\bdot\,$  denotes a running index in a finite sequence.
\[\left\{\scalebox{.9}{$
\begin{aligned}
m_{\crdot} &= \ang{m_j}_{j\in \ufs{n}} & m &= \txsum_{j=1}^n m_j & k_j &= \txsum_{i=1}^{m_j} k_{j,i} & k &= \txsum_{j=1}^n k_j\\ 
k_{\crdot} &= \ang{k_j}_{j\in \ufs{n}} & k_{j,\crdot} &= \ang{k_{j,i}}_{i \in \ufs{m}_j} & k_{\crdot,\crdot} &= \ang{k_{j,\crdot} }_{j\in \ufs{n}} &&\\
x_{\crdot} &= \ang{x_j}_{j\in \ufs{n}} & x_{j,\crdot} &= \ang{x_{j,i}}_{i \in \ufs{m}_j} & x_{\crdot,\crdot} &= \ang{x_{j,\crdot} }_{j\in \ufs{n}} &&\\
a_{j,i,\crdot} &= \ang{a_{j,i,h}}_{h \in \ufs{k}_{j,i}} & a_{j,\crdot,\crdot} &= \ang{a_{j,i,\crdot}}_{i \in \ufs{m}_j}
& a_{\crdot,\crdot,\crdot} &= \ang{a_{j,\crdot,\crdot}}_{j\in \ufs{n}} &&
\end{aligned}
$}
\right.\]
Then the following diagram in $\A$ commutes.
\begin{equation}\label{pseudoalg_comp_axiom}
\begin{tikzpicture}[vcenter]
\def\h{4} \def\g{0} \def\v{1.4} \def\u{-1}
\draw[0cell=.8]
(0,0) node (x11) {\gaA_n\Big(x \sscs \bang{\gaA_{m_j} \big( x_j \sscs \ang{ \gaA_{k_{j,i}} ( x_{j,i} \sscs a_{j,i,\crdot} ) }_{i \in \ufs{m}_j} \big) }_{j\in \ufs{n}}\Big)}
(x11)++(\h,\v) node (x12) {\gaA_n \Big(x \sscs \bang{\gaA_{k_j}\big( \ga\left(x_j \sscs x_{j,\crdot} \right) \sscs a_{j,\crdot,\crdot} \big) }_{j\in \ufs{n}} \Big)}
(x11)++(\h+.5,\u) node (x) {\gaA_k\Big(\ga\big(x \sscs \bang{ \ga\left(x_j \sscs x_{j,\crdot} \right) }_{j\in \ufs{n}} \big) \sscs a_{\crdot,\crdot,\crdot} \Big)}
(x11)++(0,2*\u) node (x21) {\gaA_m\Big(\ga\left(x \sscs x_\crdot \right) \sscs  \bang{ \bang{\gaA_{k_{j,i}}\left( x_{j,i} \sscs a_{j,i,\crdot} \right) }_{i \in \ufs{m}_j} }_{j\in \ufs{n}} \Big)}
(x21)++(\h,-\v) node (x22) {\gaA_k \Big( \ga\big( \ga\left(x \sscs x_\crdot \right) \sscs x_{\crdot,\crdot} \big) \sscs a_{\crdot,\crdot,\crdot} \Big)}
;
\draw[1cell=.85]  
(x11) edge[transform canvas={xshift={-3em}}, shorten >=-3pt] node[pos=.4] {\gaA_n\big(1_x \sscs \ang{ \phiA_{(m_j;\, k_{j,\crdot})} }_{j\in \ufs{n}} \big)} (x12)
(x12) edge node[pos=.6] {\phiA_{(n;\, k_{\crdot})}} (x)
(x) edge[equal] (x22)
(x11) edge[transform canvas={xshift={-1em}}] node[swap] {\phiA_{(n;\, m_{\crdot})}} (x21)
(x21) edge[transform canvas={xshift={-3em}}, shorten >=-1em] node[swap,pos=.5] {\phiA_{(m;\, k_{\crdot,\crdot} )}} (x22)
;
\end{tikzpicture}
\end{equation}
The lower-right equality holds by the associativity axiom for $\Op$.
\end{description}
This finishes the definition of an $\Op$-pseudoalgebra $(\A,\gaA,\phiA)$.  If each component of the associativity constraint $\phiA$ is the identity natural transformation, then $(\A,\gaA)$ is called an \index{algebra}\index{operad!algebra}\emph{$\Op$-algebra}.
\end{definition}

The following observation is \cite[\ref*{EqK:phi_id}]{yau-eqk}; it says that components of the associativity constraint $\phiA$ involving only $* \in \Op(0)$ are identities.

\begin{lemma}\label{phi_id}
For each $\Op$-pseudoalgebra $(\A,\gaA,\phiA)$ and object $x \in \Op(n)$ with $n>0$, the component of the associativity constraint 
\[\gaA_n\big(x ; \ang{\zero}_{j \in \ufs{n}}\big) = \gaA_n\big(x ; \ang{\gaA_0(*)}_{j \in \ufs{n}} \big) \fto{\phiA_{(n;\,0^n)}} \gaA_0(*) = \zero\]
is equal to $1_\zero$, where $0^n = \ang{0}_{j \in \ufs{n}}$.
\end{lemma}

The following $\Op$-pseudoalgebra analogue of \cref{naive_genuine_pGcat} follows from the fact that the functor $\Catg(\EG,-)$ \pcref{def:Catg} preserves finite products, $G$-functors, and $G$-natural transformations. 

\begin{proposition}\label{catgego}
Suppose $\Op$ is a $\Gcat$-operad for a group $G$.  Then the following statements hold.
\begin{enumerate}
\item\label{ng_i} Applying $\Catg(\EG,-)$ levelwise to $\Op$ and to the $\Gcat$-operad structure of $\Op$ yields a $\Gcat$-operad 
\[\Oph = \Catg(\EG,\Op),\] 
which is, furthermore, reduced if $\Op$ is.
\item\label{ng_ii} Applying $\Catg(\EG,-)$ to an $\Op$-pseudoalgebra $(\A,\gaA,\phiA)$ yields an $\Oph$-pseudoalgebra 
\[\Ah = \Catg(\EG,\A).\]
Moreover, if $\A$ is an $\Op$-algebra, then $\Ah$ is an $\Oph$-algebra.
\end{enumerate}
\end{proposition}

The following definition is a lax extension of \cite[Def.\ 2.23]{gmmo20} and \cite[Def.\ 2.4]{corner-gurski}.

\begin{definition}\label{def:laxmorphism}
Suppose $(\A,\gaA,\phiA)$ and $(\B,\gaB,\phiB)$ are $\Op$-pseudoalgebras \pcref{def:pseudoalgebra} for a reduced $\Gcat$-operad $(\Op,\ga,\opu)$ and a group $G$.  A \emph{lax $\Op$-morphism}\index{lax morphism}\index{operad!lax morphism}
\[\big(\A,\gaA,\phiA\big) \fto{(f,\actf)} \big(\B,\gaB,\phiB\big)\]
consists of a $G$-functor 
\begin{equation}\label{fAB}
\A \fto{f} \B
\end{equation} 
and a $G$-natural transformation $\actf_n$ for each $n \geq 0$ as follows.
\begin{equation}\label{laxmorphism_constraint}
\begin{tikzpicture}[xscale=3.5,yscale=1.3,vcenter]
\draw[0cell=.9]
(0,0) node (x11) {\Op(n) \times \A^n}
(x11)++(1,0) node (x12) {\Op(n) \times \B^n}
(x11)++(0,-1) node (x21) {\A}
(x21)++(1,0) node (x22) {\B}
;
\draw[1cell=.9]  
(x11) edge node {1 \times f^n} (x12)
(x12) edge node {\gaB_n} (x22)
(x11) edge node[swap] {\gaA_n} (x21)
(x21) edge node[swap] {f} (x22)
;
\draw[2cell]
node[between=x11 and x22 at .5, shift={(0,0)}, rotate=-90, 2label={above,\actf_n}] {\Rightarrow}
;
\end{tikzpicture}
\end{equation}
We call $\actf_n$ an \emph{action constraint}\index{action constraint}.  A typical component of $\actf_n$ is a morphism 
\begin{equation}\label{actf_component}
\gaB_n\big(x \sscs \ang{f(a_j)}_{j\in \ufs{n}} \big) \fto{\actf_n} 
f\big(\gaA_n (x \sscs \ang{a_j}_{j\in \ufs{n}} ) \big) \inspace \B
\end{equation}
for objects $x \in \Op(n)$ and $\ang{a_j}_{j\in \ufs{n}} \in \A^n$.  The data $(f,\actf)$ are required to satisfy the axioms \cref{laxmorphism_basepoint,laxmorphism_unity,laxmorphism_equiv,laxmorphism_associativity} below.
\begin{description}
\item[Basepoint] $\actf_0$ is the identity morphism:
\begin{equation}\label{laxmorphism_basepoint}
\zero^\B = \gaB_0(*) \fto{\actf_0 \,=\, 1} f\big(\gaA_0(*)\big) = f(\zero^\A).
\end{equation}
In particular, the $G$-functor $f$ preserves the basepoint.

\item[Unity] 
For each object $a \in \A$, the following diagram in $\B$ commutes, where the two equalities hold by the action unity axiom \cref{pseudoalg_action_unity}.
\begin{equation}\label{laxmorphism_unity}
\begin{tikzpicture}[xscale=3.5,yscale=1,vcenter]
\draw[0cell=.9]
(0,0) node (x11) {\gaB_1\big( \opu \sscs f(a) \big)}
(x11)++(1,0) node (x12) {f\big(\gaA_1(\opu \sscs a) \big)}
(x11)++(0,-1) node (x21) {f(a)}
(x21)++(1,0) node (x22) {f(a)}
;
\draw[1cell=.9]  
(x11) edge node {\actf_1} (x12)
(x21) edge node {1} (x22)
(x11) edge[equal] (x21)
(x12) edge[equal] (x22)
;
\end{tikzpicture}
\end{equation}

\item[Equivariance] 
For each permutation $\sigma \in \Sigma_n$ and objects $x \in \Op(n)$ and $\ang{a_j}_{j\in \ufs{n}} \in \A^n$, the following diagram in $\B$ commutes, where the two equalities hold by the action equivariance axiom \cref{pseudoalg_action_sym}.
\begin{equation}\label{laxmorphism_equiv}
\begin{tikzpicture}[xscale=4.5,yscale=1.2,vcenter]
\draw[0cell=.9]
(0,0) node (x11) {\gaB_n\big(x\sigma \sscs \ang{f(a_{\sigma(j)})}_{j\in \ufs{n}} \big)}
(x11)++(1,0) node (x12) {f\big(\gaA_n \big(x\sigma \sscs \ang{a_{\sigma(j)}}_{j\in \ufs{n}}  \big) \big)}
(x11)++(0,-1) node (x21) {\gaB_n\big(x \sscs \bang{f(a_j)}_{j\in \ufs{n}} \big)}
(x21)++(1,0) node (x22) {f\big(\gaA_n (x \sscs \ang{a_j}_{j\in \ufs{n}} ) \big)}
;
\draw[1cell=.9]  
(x11) edge node {\actf_n} (x12)
(x21) edge node {\actf_n} (x22) 
(x11) edge[equal] (x21)
(x12) edge[equal] (x22)
;
\end{tikzpicture}
\end{equation}

\item[Associativity] Consider objects 
\[x \in \Op(n),\quad x_j \in \Op(m_j), \andspace a_{j,i} \in \A\]
for $1 \leq j \leq n$ and $1 \leq i \leq m_j$, along with the following notation, where $\bdot\,$ denotes a running index in a finite sequence.
\[\left\{\scalebox{.9}{$
\begin{aligned}
m_{\crdot} &= \ang{m_j}_{j\in \ufs{n}} & a_{j,\crdot} &= \ang{a_{j,i}}_{i \in \ufs{m}_j} & a_{\crdot,\crdot} &= \ang{a_{j,\crdot}}_{j\in \ufs{n}} & x_{\crdot} &= \ang{x_j}_{j\in \ufs{n}}\\
\actf_{m_{\crdot}} &= \bang{\actf_{m_j}}_{j\in \ufs{n}} & fa_{j,\crdot} &= \bang{f(a_{j,i})}_{i \in \ufs{m}_j} & fa_{\crdot,\crdot} &= \bang{fa_{j,\crdot}}_{j\in \ufs{n}} & m &= \txsum_{j=1}^n m_j
\end{aligned}
$}
\right.\]
Then the following diagram in $\B$ commutes.
\begin{equation}\label{laxmorphism_associativity}
\begin{tikzpicture}[xscale=1,yscale=1,vcenter]
\def\h{3.3} \def\g{0} \def\v{1.2} \def\u{-1}
\draw[0cell=.8]
(0,0) node (x11) {\gaB_n\big(x \sscs \big\langle \gaB_{m_j} ( x_j \sscs fa_{j,\crdot} ) \big\rangle_{j\in \ufs{n}} \big)}
(x11)++(\h,\v) node (x12) {\gaB_m\big( \ga(x\sscs x_\crdot) \sscs fa_{\crdot,\crdot} \big)}
(x11)++(\h+.5,\u) node (x) {f \gaA_m\big( \ga(x \sscs x_\crdot) \sscs a_{\crdot,\crdot} \big)}
(x11)++(0,2*\u) node (x21) {\gaB_n\big(x \sscs \big\langle f\gaA_{m_j} (x_j \sscs a_{j,\crdot} ) \big\rangle_{j\in \ufs{n}} \big)}
(x21)++(\h,-\v) node (x22) {f\gaA_n\big( x \sscs \bang{\gaA_{m_j}(x_j \sscs a_{j,\crdot}) }_{j\in \ufs{n}} \big)}
;
\draw[1cell=.85]  
(x11) edge[transform canvas={xshift={-1em}}, shorten >=0ex] node[pos=.4] {\phiB_{(n;\, m_\crdot)}} (x12)
(x12) edge node[pos=.6] {\actf_m} (x)
(x11) edge[transform canvas={xshift={0em}}] node[swap] {\gaB_n\big( 1_x \sscs \actf_{m_{\crdot}} \big)} (x21)
(x21) edge[transform canvas={xshift={-1em}}, shorten >=0pt] node[swap,pos=.5] {\actf_n} (x22)
(x22) edge node[swap,pos=.6] {f\phiA_{(n;\, m_\crdot)}} (x)
;
\end{tikzpicture}
\end{equation}
\end{description}
This finishes the definition of a lax $\Op$-morphism.  Moreover, $(f,\actf)$ is called
\begin{itemize}
\item an \emph{$\Op$-pseudomorphism}\index{pseudomorphism}\index{operad!pseudomorphism} if each $\actf_n$ is a $G$-natural isomorphism and
\item a \emph{strict $\Op$-morphism}\index{strict morphism}\index{operad!strict morphism} if each $\actf_n$ is the identity.
\end{itemize}
A strict $\Op$-morphism between $\Op$-algebras is called an \index{algebra morphism}\index{operad!algebra morphism}\emph{$\Op$-algebra morphism}.
\end{definition}

The following definition is a lax extension of \cite[Def.\ 2.24]{gmmo20} and \cite[Def.\ 2.7]{corner-gurski}.

\begin{definition}\label{def:algtwocells}
Suppose 
\[(f,\actf), (h,\acth) \cn \big(\A,\gaA,\phiA\big) \to \big(\B,\gaB,\phiB\big)\]
are lax $\Op$-morphisms \pcref{def:laxmorphism} between $\Op$-pseudoalgebras \pcref{def:pseudoalgebra} for a reduced $\Gcat$-operad $(\Op,\ga,\opu)$ and a group $G$.  An \emph{$\Op$-transformation}\index{transformation}\index{operad!transformation}
\begin{equation}\label{Otransform_iicell}
\begin{tikzpicture}[vcenter]
\def\t{30}
\draw[0cell]
(0,0) node (a1) {\phantom{A}}
(a1)++(1.8,0) node (a2) {\phantom{A}}
(a1)++(-.7,0) node (a1') {(\A,\gaA,\phiA)}
(a2)++(.7,0) node (a2') {(\B,\gaB,\phiB)}
;
\draw[1cell=.9]
(a1) edge[bend left=\t] node {(f,\actf)} (a2)
(a1) edge[bend right=\t] node[swap] {(h,\acth)} (a2)
;
\draw[2cell=1]
node[between=a1 and a2 at .4, rotate=-90, 2label={above,\omega}] {\Rightarrow}
;
\end{tikzpicture}
\end{equation}
is a $G$-natural transformation $\omega \cn f \to h$ such that, for objects $x \in \Op(n)$ and $\anga = \ang{a_j}_{j\in \ufs{n}} \in \A^n$, the diagram
\begin{equation}\label{Otransformation_ax}
\begin{tikzpicture}[xscale=3.6,yscale=1.5,vcenter]
\draw[0cell=.9]
(0,0) node (x11) {\gaB_n\big(x \sscs \ang{f(a_j)}_{j\in \ufs{n}} \big)}
(x11)++(1,0) node (x12) {f\gaA_n \big(x \sscs \anga \big)}
(x11)++(0,-1) node (x21) {\gaB_n\big(x \sscs \bang{h(a_j)}_{j\in \ufs{n}} \big)}
(x21)++(1,0) node (x22) {h\gaA_n \big(x \sscs \anga \big)}
;
\draw[1cell=.9]  
(x11) edge node {\actf_n} (x12)
(x21) edge node {\acth_n} (x22) 
(x11) edge[transform canvas={xshift={2em}}] node[swap] {\gaB_n\big(1_x \sscs \ang{\omega_{a_j}}_{j\in \ufs{n}} \big)} (x21)
(x12) edge[transform canvas={xshift={-1em}}, shorten >=1pt] node {\omega_{\gaA_n(x;\, \anga)}} (x22)
;
\end{tikzpicture}
\end{equation}
in $\B$ commutes.
\end{definition}

\begin{explanation}\label{expl:Otr_pointed}
The basepoint axiom \cref{laxmorphism_basepoint} and the $n=0$ case of the commutative diagram \cref{Otransformation_ax} imply that the basepoint-component $\omega_{\zero^\A}$ is the identity morphism of the basepoint $\zero^\B \in \B$.  In other words, each $\Op$-transformation is necessarily pointed as a $G$-natural transformation.
\end{explanation}

The following observation is \cite[\ref*{EqK:oalgps_twocat}]{yau-eqk}.  For a brief review of 2-categories, see \cref{sec:twocategories}.

\begin{proposition}\label{oalgps_twocat}
Suppose $(\Op,\ga,\opu)$ is a reduced $\Gcat$-operad for a group $G$.  Then there is a 2-category $\AlglaxO$\index{pseudoalgebra!2-category}\index{2-category!pseudoalgebra} with 
\begin{itemize}
\item $\Op$-pseudoalgebras (\cref{def:pseudoalgebra}) as objects,
\item lax $\Op$-morphisms (\cref{def:laxmorphism}) as 1-cells, and
\item $\Op$-transformations (\cref{def:algtwocells}) as 2-cells.
\end{itemize}
Moreover, there are sub-2-categories $\AlgpspsO$ and $\AlgstO$ with the same objects and 2-cells as $\AlglaxO$, and with 1-cells given by, respectively, $\Op$-pseudomorphisms and strict $\Op$-morphisms.
\end{proposition}

\begin{explanation}\label{expl:oalgps_twocat}
For each variant $\va \in \{\sflax,\sfps,\sfst\}$, the 2-category $\AlgpsvO$ can also be obtained from the $\Gcat$-multicategory $\MultvO$ in \cite[\ref*{EqK:thm:multpso}]{yau-eqk} by
\begin{itemize}
\item first passing to $G$-fixed subcategories in the $\Gcat$-enrichment and
\item then restricting to 1-ary 1-cells and 1-ary 2-cells.
\end{itemize} 
A subtle point to keep in mind is that, while the 2-category $\AlgpsvO$ exists for any reduced $\Gcat$-operad $\Op$, the $\Gcat$-multicategory $\MultvO$ requires a pseudo-commutative $\Gcat$-operad \pcref{def:pseudocom_operad} with $\Op(1) = \bone$.  See \cite[\ref*{EqK:ex:underlying_iicat}]{yau-eqk} for further discussion.  

In the 2-category $\AlglaxO$, the horizontal composition of two composable lax $\Op$-morphisms
\[\begin{tikzpicture}[xscale=1,yscale=1,vcenter]
\def\h{3}
\draw[0cell=.9]
(0,0) node (a) {\big(\A,\gaA,\phiA\big)}
(a)++(\h,0) node (b) {\big(\B,\gaB,\phiB\big)}
(b)++(\h,0) node (c) {\big(\C,\gaC,\phiC\big)}
;
\draw[1cell=.9]
(a) edge node {(f,\actf)} (b)
(b) edge node {(h,\acth)} (c)
;
\end{tikzpicture}\]
is defined by the composite $G$-functor
\begin{equation}\label{Omorphism_comp}
\A \fto{hf} \C.
\end{equation}
For each $n \geq 0$, the action constraint $\acthf_n$ is defined as the pasting
\begin{equation}\label{Omorphism_paste}
\begin{tikzpicture}[xscale=3,yscale=1.3,vcenter]
\draw[0cell=.9]
(0,0) node (x11) {\Op(n) \times \A^n}
(x11)++(1,0) node (x12) {\Op(n) \times \B^n}
(x12)++(1,0) node (x13) {\Op(n) \times \C^n}
(x11)++(0,-1) node (x21) {\A}
(x21)++(1,0) node (x22) {\B}
(x22)++(1,0) node (x23) {\C}
;
\draw[1cell=.9]  
(x11) edge node {1 \times f^n} (x12)
(x12) edge node {1 \times h^n} (x13)
(x11) edge node[swap] {\gaA_n} (x21)
(x12) edge node[swap] {\gaB_n} (x22)
(x13) edge node {\gaC_n} (x23)
(x21) edge node[swap] {f} (x22)
(x22) edge node[swap] {h} (x23)
;
\draw[2cell]
node[between=x11 and x22 at .5, shift={(-.1,0)}, rotate=-90, 2label={above,\actf_n}] {\Rightarrow}
node[between=x12 and x23 at .5, shift={(0,0)}, rotate=-90, 2label={above,\acth_n}] {\Rightarrow}
;
\end{tikzpicture}
\end{equation}
as a $G$-natural transformation.
\end{explanation}

\section{Symmetric Monoidal $G$-Categories}\label{sec:smGcat}

This section recalls the following concepts from \cite[\namecref{EqK:ch:psalg} \ref*{EqK:ch:psalg}]{yau-eqk}.
\begin{itemize}
\item \desref{The 2-categories $\smgcat$, $\smgcatsg$, and $\smgcatst$ with naive symmetric monoidal $G$-categories as objects}{def:naive_smGcat,def:smGfunctor,def:monGnat,def:smGcat_twocat}
\item \desref{The 2-equivalence between $\AlglaxBE$, where $\BE$ is the Barratt-Eccles $\Gcat$-operad, and $\smgcat$, and their pseudo and strict variants}{thm:BEpseudoalg}
\item \desref{Genuine symmetric monoidal $G$-categories}{def:GBE_pseudoalg}
\end{itemize}

\subsection*{Naive Symmetric Monoidal $G$-Categories}

\begin{definition}\label{def:naive_smGcat}
For a group $G$, a \emph{naive symmetric monoidal $G$-category}\index{naive symmetric monoidal G-category@naive symmetric monoidal $G$-category}\index{G-category@$G$-category!naive symmetric monoidal} 
\[(\A,\otimes,\zero,\alpha,\lambda,\rho,\xi)\] 
consists of
\begin{itemize}
\item a small $G$-category $\A$ and
\item a strictly unital symmetric monoidal structure $(\otimes,\zero,\alpha,\lambda=1,\rho=1,\xi)$
\end{itemize}
such that the following two statements hold.
\begin{enumerate}
\item The monoidal unit $\zero \in \A$ is $G$-fixed.
\item The monoidal product $\otimes$, the associativity isomorphism $\alpha$, and the braiding $\xi$ are $G$-equivariant.\defmark
\end{enumerate}
\end{definition}

\begin{definition}\label{def:smGfunctor}
Suppose $(\A,\otimes,\zero,\alpha,\xi)$ and $(\B,\otimes,\zero,\alpha,\xi)$ are naive symmetric monoidal $G$-categories for a group $G$.  A \emph{symmetric monoidal $G$-functor}\index{symmetric monoidal G-functor@symmetric monoidal $G$-functor}\index{G-functor@$G$-functor!symmetric monoidal}
\[(f,f^2,f^0) \cn \A \to \B\]
is a symmetric monoidal functor such that the following three conditions hold.
\begin{enumerate}
\item $f$ is a $G$-functor.
\item The unit constraint $f^0 \cn \zero^\B \to f(\zero^\A)$ is $G$-fixed.
\item The monoidal constraint $f^2$ is a $G$-natural transformation.
\end{enumerate}
Moreover, $(f,f^2,f^0)$ is called
\begin{itemize}
\item \emph{strictly unital} if $f^0 = 1_{\zero^\B}$;
\item \emph{strictly unital strong} if $f^0 = 1_{\zero^\B}$ and $f^2$ is invertible; and
\item \emph{strict} if $f^0$ and $f^2$ are identities.\defmark
\end{itemize}
\end{definition}

\begin{definition}\label{def:monGnat}
For a group $G$, suppose
\[(f,f^2,f^0 = 1), (h,h^2,h^0 = 1) \cn (\A, \otimes, \zero, \alpha,\xi) \to (\B, \otimes, \zero, \alpha,\xi)\]
are strictly unital symmetric monoidal $G$-functors between naive symmetric monoidal $G$-categories.  A \emph{monoidal $G$-natural transformation}\index{monoidal G-natural transformation@monoidal $G$-natural transformation}\index{G-natural transformation@$G$-natural transformation!monoidal}
\[(f,f^2,f^0 = 1) \fto{\psi} (h,h^2,h^0 = 1)\]
is a monoidal natural transformation that is also $G$-equivariant.
\end{definition}

\begin{definition}\label{def:smGcat_twocat}
For a group $G$, the 2-category $\smgcat$\index{naive symmetric monoidal G-category@naive symmetric monoidal $G$-category!2-category}\index{2-category!naive symmetric monoidal G-category@naive symmetric monoidal $G$-category} is defined as follows.
\begin{itemize}
\item Objects are naive symmetric monoidal $G$-categories (\cref{def:naive_smGcat}).
\item 1-cells are strictly unital symmetric monoidal $G$-functors (\cref{def:smGfunctor}).
\item 2-cells are monoidal $G$-natural transformations (\cref{def:monGnat}).
\item Horizontal composition of 1-cells is composition of symmetric monoidal functors. 
\item Horizontal and vertical compositions of 2-cells are those of natural transformations. 
\end{itemize}
Moreover, the sub-2-categories $\smgcatsg$ and $\smgcatst$ have the same objects and 2-cells as $\smgcat$.
\begin{itemize}
\item 1-cells in $\smgcatsg$ are strictly unital strong symmetric monoidal $G$-functors.
\item 1-cells in $\smgcatst$ are strict symmetric monoidal $G$-functors.\defmark
\end{itemize}
\end{definition}

Considering the 2-categories $\AlglaxBE$, $\AlgpspsBE$, and $\AlgstBE$ \pcref{oalgps_twocat} for the Barratt-Eccles $\Gcat$-operad $\BE$ \pcref{def:BE}, the following observation says that, up to a 2-equivalence, $\BE$-pseudoalgebras and naive symmetric monoidal $G$-categories \pcref{def:pseudoalgebra,def:naive_smGcat} are interchangeable.  This result is proved in detail in  \cite[\ref*{EqK:thm:BEpseudoalg}]{yau-eqk}.

\begin{theorem}\label{thm:BEpseudoalg}
There are \index{2-equivalence}2-equivalences
\[\begin{split}
\AlglaxBE & \fto[\sim]{\Phi} \smgcat,\\
\AlgpspsBE & \fto[\sim]{\Phi} \smgcatsg, \andspace\\
\AlgstBE & \fto[\sim]{\Phi} \smgcatst.
\end{split}\]
\end{theorem}

\subsection*{Genuine Symmetric Monoidal $G$-Categories}

The following definition is \cite[Def.\ 0.3]{gmmo20}.

\begin{definition}\label{def:GBE_pseudoalg}
For the $G$-Barratt-Eccles $\Gcat$-operad $\GBE$, $\GBE$-pseudoalgebras in $\Gcat$ \pcref{def:GBE,def:pseudoalgebra} are called \index{G-Barratt-Eccles operad@$G$-Barratt-Eccles operad!pseudoalgebra}\index{pseudoalgebra!G-Barratt-Eccles operad@$G$-Barratt-Eccles operad}\index{G-category@$G$-category!genuine symmetric monoidal}\index{genuine symmetric monoidal G-category@genuine symmetric monoidal $G$-category}\emph{genuine symmetric monoidal $G$-categories}.
\end{definition}

Applied to the Barratt-Eccles $\Gcat$-operad $\BE$ (\cref{def:BE}), \cref{catgego} \eqref{ng_ii} yields the following observation.

\begin{proposition}\label{naive_genuine_smgcat}
For each group $G$, each $\BE$-pseudoalgebra $\A$ yields a $\GBE$-pseudoalgebra $\Catg(\EG,\A)$.
\end{proposition}

\cref{thm:BEpseudoalg,naive_genuine_smgcat} together imply that $\Catg(\EG,-)$ sends naive symmetric monoidal $G$-categories to genuine symmetric monoidal $G$-categories.

 

%% file: chap/kgo_ii.tex
\section{$\Gskg$-Categories}
\label{sec:ggcat}

This section reviews the 2-category $\GGCatii$, whose objects are pointed functors
\[(\Gsk,\vstar) \fto{f} (\Gcatst,\boldone),\]
called $\Gskg$-categories.  This 2-category is obtained from the symmetric monoidal closed category $\GGCat$ in \cite[\namecref{EqK:ch:ggcat} \ref*{EqK:ch:ggcat}]{yau-eqk} by keeping the same objects and passing to the $G$-fixed subcategories in its $\Gcat$-enrichment.

\secoutline
\begin{itemize}
\item \cref{def:Fsk,def:Fsk_permutative,def:Fsk_smashpower,def:injections} recall the indexing category $\Fsk$ of pointed finite sets and some auxiliary constructions.
\item \cref{def:Gsk,def:Gsk_permutative,def:smashFskGsk} recall the indexing category $\Gsk$ and its comparison with $\Fsk$.
\item The 2-category $\GGCatii$ is defined in \cref{def:ptGcat,def:gcatst,def:GGCat}, with further elaboration given in \cref{expl:GGCat_iicat}.
\end{itemize}

\subsection*{The Indexing Category $\Fsk$}

The category $\Fsk$ in the next definition is the opposite of Segal's category $\Ga$ \cite{segal}.

\begin{definition}\label{def:Fsk}\
\begin{itemize}
\item A \emph{pointed set}\index{pointed set} $(S, \bp)$ is a set $S$ equipped with a distinguished element $\bp$, called the \emph{basepoint}.  A \emph{pointed finite set}\index{pointed finite set}\index{finite set!pointed} is a pointed set whose underlying set is finite.
\item The pointed finite set
\begin{equation}\label{ordn}
\ord{n} = \{0 < 1 < \cdots < n\}
\end{equation}
is equipped with the basepoint $0$ and its natural ordering.  Unless otherwise specified, a \emph{pointed finite set} means $\ordn$ for some $n \geq 0$.
\item A \emph{pointed function}\index{pointed function} between pointed sets is a basepoint-preserving function.
\item The small pointed category $\Fsk$\label{not:Fsk} has objects $\ordn$ for $n \geq 0$, pointed functions as morphisms, and basepoint $\ord{0}$.\defmark
\end{itemize}
\end{definition}

Recall that a permutative category is a strict symmetric monoidal category \pcref{def:symmoncat}. 

\begin{definition}[Permutative Structure]\label{def:Fsk_permutative}
There is a permutative category structure $(\sma, \ord{1}, \xi)$ on $\Fsk$ defined as follows.
\begin{description}
\item[Monoidal product] The monoidal product is the functor
\[\Fsk \times \Fsk \fto{\sma} \Fsk\]
given on objects by
\begin{equation}\label{m-sma-n}
\ord{m} \sma \ord{n} = \ord{mn}.
\end{equation}
It is extended to morphisms using the lexicographic ordering \cref{lex_bijection}.
\item[Monoidal unit] The monoidal unit is the object $\ord{1} = \{0 < 1\} \in \Fsk$.
\item[Braiding] For a pair of objects $(\ord{m}, \ord{n}) \in \Fsk^2$, the braiding
\begin{equation}\label{Fsk_braiding}
\ord{m} \sma \ord{n} \fto[\iso]{\xi_{\ord{m}, \ord{n}}} \ord{n} \sma \ord{m}
\end{equation}
is given by the $(m,n)$-transpose permutation $\twist_{m,n}$ \cref{eq:transpose_perm} away from the basepoint $0 \in \ord{mn}$.
\end{description}
It is proved in \cite[\ref*{EqK:Fsk_permutative}]{yau-eqk} that $(\Fsk,\sma,\ord{1},\xi)$ is a permutative category.
\end{definition}

\begin{definition}[Smash Powers]\label{def:Fsk_smashpower}
For $q \geq 0$, the small pointed category $\Fsk^{(q)}$ is defined as follows.
\begin{description}
\item[$q=0$] The pointed category\label{not:Fsk_smash_zero}
\[\Fsk^{(0)} = \big\{\vstar \rightleftarrows \ang{}\big\}\]
consists of the initial-terminal basepoint $\vstar$, the empty tuple $\ang{}$, the identity morphisms of $\vstar$ and $\ang{}$, and the nonidentity morphisms $\vstar \to \ang{} \to \vstar$.
\item[$q>0$] In this case, the small pointed category\label{not:Fskq}
\[\Fsk^{(q)} = \Fsk^{\sma q}\]
is the $q$-fold smash power of $\Fsk$ defined as follows. 
\begin{description}
\item[Objects]  
It has an initial-terminal basepoint $\vstar$.  A typical object in $\Fsk^{(q)}$ is a $q$-tuple
\begin{equation}\label{angordn}
\ang{\ord{n}} = \ang{\ord{n}_i}_{i \in \ufs{q}} = (\ord{n}_1, \ldots, \ord{n}_q)
\end{equation}
with each $\ord{n}_i$ a pointed finite set \cref{ordn}.  If any $\ord{n}_i = \ord{0}$, then $\ang{\ord{n}} = \vstar$. 
\item[Morphisms] 
A typical morphism in $\Fsk^{(q)}$ is a $q$-tuple
\begin{equation}\label{angpsi}
\ang{\psi} = \ang{\psi_i}_{i \in \ufs{q}} \cn \ang{\ord{m}_i}_{i \in \ufs{q}} \to \ang{\ord{n}_i}_{i \in \ufs{q}}
\end{equation}
with each
\[\ord{m}_i \fto{\psi_i} \ord{n}_i\]
a morphism in $\Fsk$, meaning a pointed function.  If any $\psi_i$ factors through $\ord{0}$, then $\ang{\psi}$ is the \index{0-morphism}\emph{0-morphism}, meaning that it factors through the basepoint $\vstar$.  A morphism that is not the 0-morphism is called a \index{nonzero morphism}\emph{nonzero morphism}.
\end{description}
\end{description}  
A nonbasepoint object or a nonzero morphism in $\Fsk^{(q)}$ is said to have \emph{length}\index{length} $q$.
\end{definition}

\begin{definition}[Reindexing]\label{def:injections}
The category $\Inj$ has the unpointed finite sets $\ufs{n} = \{1, 2,\ldots , n\}$ \cref{ufsn} for $n \geq 0$ as objects and injections as morphisms.

For an injection $h \cn \ufs{q} \to \ufs{r}$, the pointed functor\label{not:reind_func}
\[\Fsk^{(q)} \fto{h_*} \Fsk^{(r)}\]
is defined as follows.
\begin{description}
\item[$q = r = 0$] 
In this case, $h$ is the identity function on $\ufs{0} = \emptyset$, and $h_*$ is the identity functor on $\Fsk^{(0)}$.
\item[$q = 0 < r$] 
In this case, the pointed functor $h_*$ is determined by the following object assignment.
\begin{equation}\label{hangempty}
\left\{\begin{split}
h_* \vstar &= \vstar\\
h_*\ang{} &= \ang{\ord{1}}_{j \in \ufs{r}} = \left(\ord{1}, \ldots, \ord{1}\right) \in \Fsk^{(r)}
\end{split}\right.
\end{equation}
\item[$q>0$] 
Given an object $\ang{\ord{n}_i}_{i \in \ufs{q}}$ \cref{angordn} and a morphism $\ang{\psi_i}_{i \in \ufs{q}}$ \cref{angpsi} in $\Fsk^{(q)}$, we define the object and morphism
\begin{equation}\label{reindexing_functor}
\begin{split}
h_*\ang{\ord{n}_i}_{i \in \ufs{q}} &= \ang{\ord{n}_{\hinv(j)}}_{j \in \ufs{r}} \andspace\\
h_*\ang{\psi_i}_{i \in \ufs{q}} &= \ang{\psi_{\hinv(j)}}_{j \in \ufs{r}}
\end{split}
\end{equation}
in $\Fsk^{(r)}$.  If $\hinv(j) = \emptyset$, then
\begin{equation}\label{ordn_empty}
\ord{n}_{\emptyset} = \ord{1} \andspace \ord{1} \fto{\psi_{\emptyset} = 1_{\ord{1}}} \ord{1}.
\end{equation}
\end{description}
We call $h$ a \emph{reindexing injection}\index{reindexing injection} and $h_*$ a \index{reindexing functor}\emph{reindexing functor}.
\end{definition}

\subsection*{The Indexing Category $\Gsk$}

Next, we recall the category $\Gsk$ from \cite{elmendorf-mandell}, which is discussed thoroughly in \cite[Ch.\ 9--13]{cerberusIII}.  

\begin{definition}\label{def:Gsk}
The small pointed category $\Gsk$ is defined as follows.
\begin{description}
\item[Objects] The pointed set of objects of $\Gsk$ is defined as the wedge
\begin{equation}\label{Gsk_objects}
\Ob(\Gsk) = \bigvee_{q \geq 0} \Ob(\Fsk^{(q)})
\end{equation}
that identifies the basepoints $\vstar \in \Fsk^{(q)}$ for $q \geq 0$.  The identified object is the initial-terminal basepoint $\vstar \in \Gsk$.
\item[Morphisms]
Given an object $\angordm \in \Fsk^{(p)}$ and an object $\angordn \in \Fsk^{(q)}$, the pointed set of morphisms is defined as the following wedge.
\begin{equation}\label{Gsk_morphisms}
\begin{split}
\Gsk\big(\angordm, \angordn\big)
&= \bigvee_{f \in \Inj(\ufs{p},\, \ufs{q})}~ \Fsk^{(q)}\big(f_* \angordm, \angordn \big)\\
&= \bigvee_{f \in \Inj(\ufs{p},\, \ufs{q})}~ \bigwedge_{i \in \ufs{q}}~ \Fsk\big(\ord{m}_{\finv(i)} \scs \ord{n}_i \big)
\end{split}
\end{equation}
The basepoint of each pointed set of morphisms is the \index{0-morphism}\emph{0-morphism}, which is the unique morphism that factors through $\vstar \in \Gsk$.  A morphism that is not the 0-morphism is called a \index{nonzero morphism}\emph{nonzero morphism}. 

In \cref{Gsk_morphisms}, for each reindexing injection $f \cn \ufs{p} \to \ufs{q}$, 
\[\Fsk^{(p)} \fto{f_*} \Fsk^{(q)}\]
is the reindexing functor \pcref{def:injections}. 
\begin{itemize}
\item If $q=0$ in \cref{Gsk_morphisms}, then $p=0$ and $f \cn \ufs{0} \to \ufs{0}$ is $1_\emptyset$.  The pointed set
\begin{equation}\label{Gsk_empty_mor}
\Gsk\big(\ang{}, \ang{}\big) = \Fsk^{(0)}\big(\ang{}, \ang{}\big)
\end{equation}
consists of the identity morphism of $\ang{}$ and the 0-morphism $\ang{} \to \vstar \to \ang{}$.
\item For $q>0$ in \cref{Gsk_morphisms}, a morphism in $\Gsk$ is a pair
\begin{equation}\label{fangpsi}
\angordm \fto{(f, \ang{\psi})} \angordn
\end{equation}
consisting of
\begin{itemize}
\item a reindexing injection $f \cn \ufs{p} \to \ufs{q}$ and
\item a morphism 
\[\ang{\psi} = \ang{\psi_i}_{i \in \ufs{q}} \cn f_* \angordm \to \angordn
\inspace \Fsk^{(q)}.\]
\end{itemize}
If some pointed function
\[\ord{m}_{\finv(i)} \fto{\psi_i} \ord{n}_i\]
factors through $\ord{0} \in \Fsk$, then $(f, \ang{\psi})$ is the 0-morphism, factoring through the basepoint $\vstar$.  
\end{itemize}
A morphism in $\Gsk$ is also denoted by a generic symbol, such as $\upom$.
\item[Identities] The identity morphism 
\[1_{\angordm} = \big(1_{\ufs{p}}, \ang{1_{\ordm_k}}_{k \in \ufs{p}}\big)\]
of an object $\angordm$ of length $p \geq 0$ consists of the identity function on $\ufs{p}$ and the identity function on $\ordm_k$ for each $k \in \ufs{p}$.
\item[Composition] Consider composable morphisms in $\Gsk$
\begin{equation}\label{Gsk_composable}
\angordm \fto{(f,\ang{\psi})} \angordn \fto{(h,\ang{\phi})} \angordl
\end{equation}
for objects $\angordm \in \Fsk^{(p)}$, $\angordn \in \Fsk^{(q)}$, and $\angordl \in \Fsk^{(r)}$.
\begin{itemize}
\item If $q=0$, then $p=0$, and $(f,\ang{\psi})$ is either the identity morphism $1_{\ang{}}$ or the 0-morphism.  Their composites with $(h,\ang{\phi})$ are, respectively, $(h,\ang{\phi})$ and the 0-morphism.
\item If $q>0$, then $r>0$, and the composite is defined as the pair
\begin{equation}\label{Gsk_composite}
(h,\ang{\phi}) \circ (f,\ang{\psi}) = \big(hf, \ang{\phi} \circ h_*\ang{\psi} \big) \cn \angordm \to \angordl.
\end{equation}
\end{itemize}
\end{description}
It is proved in \cite[\ref*{EqK:Gsk_Gcategory}]{yau-eqk} that $\Gsk$ is a well-defined pointed category.  
\end{definition}

\begin{definition}[Length-1 Inclusion]\label{def:ifg}
The pointed full subcategory inclusion\index{length-1 inclusion}
\begin{equation}\label{ifg}
\Fsk \fto{\ifg} \Gsk
\end{equation}
sends a pointed finite set $\ordn \in \Fsk$ for $n \geq 0$ to the length-1 object of $\Gsk$ consisting of $\ordn$.  This is the initial-terminal basepoint $\vstar \in \Gsk$ if $\ordn = \ordz$. We usually abbreviate $\ifg\ordn$ to $\ordn$.  A pointed function $\psi \cn \ordm \to \ordn$ in $\Fsk$ is sent by $\ifg$ to the morphism
\[\ordm \fto{(1_{\ufs{1}}, \psi)} \ordn \inspace \Gsk.\]
This is the unique morphism from or to the basepoint $\vstar \in \Gsk$ if either $\ordm$ or $\ordn$ is $\ordz$.
\end{definition}

\begin{definition}[Permutative Structure]\label{def:Gsk_permutative}
There is a permutative category structure $(\oplus, \ang{}, \xi)$ on $\Gsk$ defined as follows.
\begin{description}
\item[Monoidal product on objects] The monoidal product
\[\Gsk \times \Gsk \fto{\oplus} \Gsk\]
is given on objects by concatenation\index{concatenation}
\begin{equation}\label{Gsk_oplus_obj}
\angordm \oplus \angordn = 
\big(\ord{m}_1, \ldots, \ord{m}_p, \ord{n}_1, \ldots, \ord{n}_q\big)
\end{equation}
for $\angordm \in \Fsk^{(p)}$ and $\angordn \in \Fsk^{(q)}$.  If any $\ord{m}_i$ or $\ord{n}_j$ is $\ord{0}$, then the right-hand side of \cref{Gsk_oplus_obj} has an entry of $\ord{0}$, so
\begin{equation}\label{Gsk_oplus_vstar}
\vstar \oplus \angordn = \vstar = \angordm \oplus \vstar.
\end{equation}
\item[Monoidal product on morphisms]
When one of the two morphisms has either domain or codomain given by the initial-terminal basepoint $\vstar \in \Gsk$, their monoidal product is uniquely defined.

For other cases, we consider nonbasepoint objects $\angordm \in \Fsk^{(p)}$, $\angordn \in \Fsk^{(q)}$, $\angordj \in \Fsk^{(r)}$, and $\angordl \in \Fsk^{(s)}$, and morphisms
\begin{equation}\label{fpsi_hphi}
\angordm \fto{(f,\ang{\psi})} \angordn \andspace \angordj \fto{(h,\ang{\phi})} \angordl
\end{equation}
in $\Gsk$.  Then we define the morphism
\begin{equation}\label{Gsk_oplus_morphism}
\begin{split}
&(f, \ang{\psi}) \oplus (h, \ang{\phi}) \\
&= \big(f \oplus h, \ang{\psi} \oplus \ang{\phi}\big) \cn
\angordm \oplus \angordj \to \angordn \oplus \angordl
\end{split}
\end{equation}
with reindexing injection
\begin{equation}\label{f_oplus_h}
\ufs{p+r} \fto{f \oplus h} \ufs{q+s}
\end{equation}
given by
\[(f \oplus h)(i) = \begin{cases} f(i) & \text{if $1 \leq i \leq p$ and}\\
q + h(i-p) & \text{if $p+1 \leq i \leq p+r$}.
\end{cases}\]
The morphism $\ang{\psi} \oplus \ang{\phi}$ is the concatenation of the $q$-tuple $\ang{\psi}$ and the $s$-tuple $\ang{\phi}$ as follows.
\begin{equation}\label{psi_oplus_phi}
\begin{tikzpicture}[vcenter]
\draw[0cell=.9]
(0,0) node (a) {(f \oplus h)_* (\angordm \oplus \angordj)}
(a)++(0,-1) node (b) {f_* \angordm \oplus h_* \angordj}
(b)++(5,0) node (c) {\angordn \oplus \angordl}
(c)++(0,1) node (d) {\angordn \oplus \angordl}
;
\draw[1cell=.9]
(a) edge[equal] (b)
(c) edge[equal] (d)
(a) edge node {\ang{\psi} \oplus \ang{\phi}} (d)
(b) edge node {\big(\ang{\psi_i}_{i \in \ufs{q}} \, , \ang{\phi_k}_{k \in \ufs{s}} \big)} (c)
;
\end{tikzpicture}
\end{equation}
\item[Monoidal unit]
The monoidal unit is the empty tuple $\ang{} \in \Gsk$.
\item[Braiding] 
For a pair of objects $(\angordm, \angordn) \in \Gsk^2$, the braiding $\xi$ is $1_\vstar$ if either $\angordm$ or $\angordn$ is the basepoint $\vstar \in \Gsk$.  For nonbasepoint objects, the braiding is the isomorphism
\begin{equation}\label{Gsk_braiding}
\xi_{\angordm, \angordn} = (\tau_{p,q} \, , \ang{1})
\cn \angordm \oplus \angordn \fto{\iso} \angordn \oplus \angordm
\end{equation}
with $p$ the length of $\angordm$ and $q$ the length of $\angordn$.  Its reindexing injection is the block permutation 
\[\ufs{p+q} \fto[\iso]{\tau_{p,q}} \ufs{q+p}\]
that interchanges the first $p$ elements with the last $q$ elements:
\[\tau_{p,q}(i) = \begin{cases} 
q+i & \text{if $1 \leq i \leq p$ and}\\
i-p & \text{if $p+1 \leq i \leq p+q$}.
\end{cases}\]
The morphism $\ang{1}$ in \cref{Gsk_braiding} is the $(q+p)$-tuple with each entry given by an identity function of some $\ordn_i$ for $i \in \ufs{q}$ or some $\ordm_k$ for $k \in \ufs{p}$.
\end{description}
It is proved in \cite[\ref*{EqK:Gsk_permutative}]{yau-eqk} that $(\Gsk,\oplus, \ang{}, \xi)$ is a permutative category.
\end{definition}

\begin{definition}[Comparing Permutative Structures]\label{def:smashFskGsk}
We define a functor
\[\Gsk \fto{\sma} \Fsk\]
as follows.
\begin{description}
\item[Objects] The object assignment is defined as follows for $\ang{\ord{m}_k}_{k \in \ufs{p}} \in \Gsk \setminus \{\vstar, \ang{}\}$.  
\begin{equation}\label{smash_Gskobjects}
\left\{
\begin{split}
\sma \vstar &= \ord{0}\\
\sma \ang{} &= \ord{1}\\
\sma \ang{\ord{m}_k}_{k \in \ufs{p}} &= \sma_{k \in \ufs{p}} \,\ord{m}_k
= \ord{m_1 m_2 \cdots m_p}
\end{split}\right.
\end{equation}
In the last case in \cref{smash_Gskobjects}, $\sma_{k \in \ufs{p}}$ is the $p$-fold iterate of the monoidal product of $\Fsk$ \cref{m-sma-n}. 
\item[Morphisms]
The identity morphism and the 0-morphism in $\Gsk\big(\ang{}, \ang{}\big)$, as discussed in \cref{Gsk_empty_mor}, are sent by $\sma$ to, respectively, the identity morphism and the 0-morphism in $\Fsk\big(\ord{1}, \ord{1}\big)$.

For a morphism 
\[\angordm \fto{(f, \ang{\psi})} \angordn \inspace \Gsk\]
as defined in \cref{fangpsi}, the morphism
\[\sma \angordm \fto{\sma (f, \ang{\psi})} \!\smam\angordn \inspace \Fsk\]
is defined as the following composite pointed function.
\begin{equation}\label{smash_fpsi}
\begin{tikzpicture}[vcenter]
\def\t{.6ex}
\draw[0cell=.9]
(0,0) node (a) {\bigwedge_{k \in \ufs{p}} \ord{m}_k}
(a)++(2.8,0) node (b) {\bigwedge_{\finv(i) \neq \emptyset} \ord{m}_{\finv(i)}}
(b)++(3,0) node (c) {\bigwedge_{i \in \ufs{q}} \ord{m}_{\finv(i)}}
(c)++(2.5,0) node (d) {\bigwedge_{i \in \ufs{q}} \ord{n}_i}
;
\draw[1cell=.9]
(a) edge[transform canvas={yshift=\t}] node {f_*} node[swap] {\iso} (b)
(b) edge[transform canvas={yshift=\t}] node {\iso} (c)
(c) edge[transform canvas={yshift=\t}] node {\sma_i\, \psi_i} (d)
;
\end{tikzpicture}
\end{equation}
The three pointed functions in \cref{smash_fpsi} are defined as follows. 
\begin{itemize}
\item $f_*$ permutes the $p$ entries according to the injection $f \cn \ufs{p} \to \ufs{q}$.  The indexing set in the codomain is given by $\big\{ i \in \ufs{q} \cn \finv(i) \neq \emptyset\}$.
\item Using \cref{ordn_empty}, the middle pointed bijection in \cref{smash_fpsi} inserts a copy of the smash unit $\ord{1} = \ord{m}_{\emptyset}$ for each index $i \in \ufs{q}$ not in the image of $f$.
\item $\sma_i\, \psi_i$ is the smash product of the pointed functions $\psi_i$ for $i \in \ufs{q}$.
\end{itemize}
\end{description}
It is proved in \cite[\ref*{EqK:sma_symmon}]{yau-eqk} that $\sma \cn \Gsk \to \Fsk$ is a strict symmetric monoidal pointed functor.
\end{definition}

\subsection*{Pointed $G$-Categories}
Recall the 2-category $\Gcat$ \pcref{def:GCat}.

\begin{definition}\label{def:ptGcat}
For a group $G$, we define the following notions.
\begin{itemize}
\item A \emph{pointed $G$-category}\index{pointed G-category@pointed $G$-category}\index{G-category@$G$-category!pointed} is a $G$-category equipped with a $G$-fixed object, called the \emph{basepoint}.
\item A \emph{pointed $G$-functor}\index{pointed G-functor@pointed $G$-functor}\index{G-functor@$G$-functor!pointed} between pointed $G$-categories is a basepoint-preserving $G$-functor. 
\item A \emph{pointed $G$-natural transformation}\index{pointed G-natural transformation@pointed $G$-natural transformation}\index{G-natural transformation@$G$-natural transformation!pointed} between pointed $G$-functors is a $G$-natural transformation whose basepoint-component is the identity morphism at the basepoint.
\item A \emph{pointed $G$-modification}\index{pointed G-modification@pointed $G$-modification}\index{G-modification@$G$-modification!pointed} between pointed $G$-natural transformations is a $G$-equivariant modification whose basepoint-component is the identity 2-cell of the identity 1-cell at the basepoint.\defmark
\end{itemize}
\end{definition}

The next definition is the pointed analogue of $\Gcat$.

\begin{definition}\label{def:gcatst}
For a group $G$, the 2-category\index{pointed G-category@pointed $G$-category!2-category}\index{2-category!pointed G-category@pointed $G$-category}
\begin{equation}\label{Gcatst}
\Gcatst
\end{equation}
has
\begin{itemize}
\item small pointed $G$-categories as objects,
\item pointed $G$-functors as 1-cells, and
\item pointed $G$-natural transformations as 2-cells.
\end{itemize}
Identity 1-cells and 2-cells, vertical composition of 2-cells, and horizontal composition of 1-cells and 2-cells are given by the corresponding structures for functors and natural transformations.  The underlying 1-category of $\Gcatst$ is denoted by the same notation.  For the trivial group $G$, $\Gcatst$ is denoted by $\Catst$.
\end{definition}

\begin{definition}[Symmetric Monoidal Closed Structure]\label{expl:Gcatst}
The complete and cocomplete Cartesian closed category \cref{gcat-closed}
\[(\Gcat, \times, \boldone, \Catg)\]
yields the complete and cocomplete symmetric monoidal closed category
\begin{equation}\label{Gcatst_smc}
(\Gcatst, \sma, \bonep, \Catgst)
\end{equation}
defined as follows.
\begin{itemize}
\item For small pointed $G$-categories $\C$ and $\D$, the pointed category $\C \sma \D$ is defined in $\Catst$ with $G$ acting diagonally, meaning that
\[g (c,d) = (gc,gd)\]
for $g \in G$, $c \in \C$, and $d \in \D$. 
\item The monoidal unit is the discrete category
\begin{equation}\label{gcatst_unit}
\bonep = \boldone \bincoprod \boldone
\end{equation} 
with two objects and the trivial $G$-action. 
\item The internal hom is given by the small pointed $G$-category 
\begin{equation}\label{catgst_cd}
\Catgst(\C,\D)
\end{equation}
with
\begin{itemize}
\item pointed functors $\C \to \D$ as objects,
\item basepoint given by the constant functor at the basepoint in $\D$,
\item pointed natural transformations as morphisms,
\item identities and composition defined componentwise in $\D$, and
\item the conjugation $G$-action \cref{conjugation-gaction}.
\end{itemize}
\end{itemize}
Passing to the $G$-fixed subcategory yields the category
\begin{equation}\label{catgst_gcatst}
\Catgst(\C,\D)^G = \Gcatst(\C,\D)
\end{equation}
of pointed $G$-functors $\C \to \D$ and pointed $G$-natural transformations.  The notation
\begin{equation}\label{Catgst_iicat}
\Catgst
\end{equation}
also denotes the following categories with small pointed $G$-categories as objects.
\begin{enumerate}
\item The 2-category with hom categories given by the categories in \cref{catgst_cd}.
\item The pointed $G$-category with the terminal $G$-category $\bone$ as the basepoint object, the trivial $G$-action on objects, pointed functors as morphisms, and the conjugation $G$-action \cref{conjugation-gaction} on morphisms.\defmark
\end{enumerate}
\end{definition}

\subsection*{$\Gskg$-Categories}

The 2-category $\GGCatii$ in the next definition is obtained from the symmetric monoidal $\Gcat$-category $\GGCat$ in \cite[\ref*{EqK:def:GGCat}]{yau-eqk} by keeping the same objects and passing to the $G$-fixed subcategories in its $\Gcat$-enrichment.

\begin{definition}[$\Gskg$-Categories]\label{def:GGCat}
For a group $G$, using the 2-category $\Gcatst$ \pcref{def:gcatst}, the 2-category $\GGCatii$ is defined as follows.
\begin{description}
\item[Objects]
An object in $\GGCatii$, called a \index{G-G-category@$\Gskg$-category}\emph{$\Gskg$-category}, is a pointed functor
\begin{equation}\label{ggcat_obj}
(\Gsk, \vstar) \fto{X} (\Gcatst,\boldone).
\end{equation}
\item[1-cells]
A 1-cell $\theta \cn X \to X'$ in $\GGCatii$ is a natural transformation as follows.
\begin{equation}\label{ggcat_mor}
\begin{tikzpicture}[vcenter]
\def\t{28}
\draw[0cell]
(0,0) node (a1) {\Gsk}
(a1)++(1.8,0) node (a2) {\phantom{\Gskel}}
(a2)++(.3,0) node (a2') {\Gcatst}
;
\draw[1cell=.9]
(a1) edge[bend left=\t] node {X} (a2)
(a1) edge[bend right=\t] node[swap] {X'} (a2)
;
\draw[2cell]
node[between=a1 and a2 at .45, rotate=-90, 2label={above,\theta}] {\Rightarrow}
;
\end{tikzpicture}
\end{equation}
\item[2-cells]
A 2-cell $\Theta \cn \theta \to \ups$ in $\GGCatii$ is a modification as follows.
\begin{equation}\label{ggcat_iicell}
\begin{tikzpicture}[vcenter]
\def\t{25}
\draw[0cell]
(0,0) node (a) {\Gsk}
(a)++(3,0) node (b) {\phantom{\GG}}
(b)++(.2,0) node (b') {\Gcatst}
;
\draw[1cell=.8]
(a) edge[bend left=\t] node {X} (b)
(a) edge[bend right=\t] node[swap] {X'} (b)
;
\draw[2cell=.9]
node[between=a and b at .32, rotate=-90, 2label={below,\theta}] {\Rightarrow}
node[between=a and b at .65, rotate=-90, 2label={above,\ups}] {\Rightarrow}
;
\draw[2cell]
node[between=a and b at .5, rotate=0, shift={(0,-.2)}, 2labelmed={above,\!\Theta}] {\Rrightarrow}
;
\end{tikzpicture}
\end{equation}
\item[Other structures]
Identity 1-cells and 2-cells, vertical composition of 2-cells, and horizontal composition of 1-cells and 2-cells are defined componentwise in the 2-category $\Gcatst$.
\end{description}
The underlying 1-category of $\GGCatii$ is denoted by the same notation.
\end{definition}

\begin{explanation}[Unpacking $\GGCatii$]\label{expl:GGCat_iicat}
The 2-category $\GGCatii$ in \cref{def:GGCat} is given explicitly as follows.
\begin{description}
\item[Objects]
A $\Gskg$-category $X \cn \Gsk \to \Gcatst$ \cref{ggcat_obj} consists of the following data.
\begin{itemize}
\item $X$ sends each object $\angordm \in \Gsk$ \cref{Gsk_objects} to a small pointed $G$-category $X\angordm$ such that $X\vstar = \boldone$.  The \emph{canonical basepoint}\dindex{canonical}{basepoint} of $X\angordm$ is given by the $G$-functor
\begin{equation}\label{fm_pointed}
X(\vstar \to \angordm) \cn X\vstar = \boldone \to X\angordm,
\end{equation}
where $\vstar \to \angordm$ is the unique morphism in $\Gsk(\vstar,\angordm)$. 
\item $X$ sends each morphism $\upom \cn \angordm \to \angordn$ in $\Gsk$ \cref{Gsk_morphisms} to a pointed $G$-functor
\begin{equation}\label{f_upom}
X\angordm \fto{X\upom} X\angordn
\end{equation}
such that $X$ preserves identities and composition of morphisms.  
\end{itemize}
\item[1-cells]
A 1-cell $\theta \cn X \to X'$ in $\GGCatii$ \cref{ggcat_mor} consists of, for each object $\angordm \in \Gsk$, an $\angordm$-component pointed $G$-functor
\begin{equation}\label{ggcat_mor_component}
X\angordm \fto{\theta_{\angordm}} X'\angordm
\end{equation}
such that, for each morphism $\upom \cn \angordm \to \angordn$ in $\Gsk$, the following naturality diagram of pointed $G$-functors commutes.
\begin{equation}\label{ggcat_mor_naturality}
\begin{tikzpicture}[vcenter]
\def\v{-1.4}
\draw[0cell]
(0,0) node (a11) {X\angordm}
(a11)++(2.5,0) node (a12) {X'\angordm}
(a11)++(0,\v) node (a21) {X\angordn}
(a12)++(0,\v) node (a22) {X'\angordn}
;
\draw[1cell=.9]
(a11) edge node {\theta_{\angordm}} (a12)
(a12) edge node {X'\upom} (a22)
(a11) edge node[swap] {X\upom} (a21)
(a21) edge node {\theta_{\angordn}} (a22)
;
\end{tikzpicture}
\end{equation}
Identity 1-cells and horizontal composition of 1-cells in $\GGCatii$ are defined componentwise using the components in \cref{ggcat_mor_component}.  A 1-cell is automatically pointed in the sense that its $\vstar$-component
\begin{equation}\label{ggcat_mor_vstar}
X\vstar = \boldone \fto{\theta_\vstar} X'\vstar = \boldone
\end{equation}
is the identity functor on $\bone$.
\item[2-cells]
A 2-cell $\Theta \cn \theta \to \ups$ in $\GGCatii$ \cref{ggcat_iicell} consists of, for each object $\angordm \in \Gsk$, an $\angordm$-component pointed $G$-natural transformation
\begin{equation}\label{ggcat_inthom_Theta}
\begin{tikzpicture}[vcenter]
\def\t{25}
\draw[0cell]
(0,0) node (a1) {\phantom{X'}}
(a1)++(2.5,0) node (a2) {\phantom{X'}}
(a1)++(-.2,0) node (a1') {X\angordm}
(a2)++(.25,0) node (a2') {X'\angordm}
;
\draw[1cell=.85]
(a1) edge[bend left=\t] node {\theta_{\angordm}} (a2)
(a1) edge[bend right=\t] node[swap] {\ups_{\angordm}} (a2)
;
\draw[2cell]
node[between=a1 and a2 at .37, rotate=-90, 2label={above, \Theta_{\angordm}}] {\Rightarrow}
;
\end{tikzpicture}
\end{equation}
such that, for each morphism $\upom \cn \angordm \to \angordn$ in $\Gsk$, the following two whiskered $G$-natural transformations are equal.
\begin{equation}\label{ggcat_inthom_Theta_modax}
\begin{tikzpicture}[vcenter]
\def\t{25} \def\v{-1.6}
\draw[0cell]
(0,0) node (a1) {\phantom{X'}}
(a1)++(2.5,0) node (a2) {\phantom{X'}}
(a1)++(-.2,0) node (a1') {X\angordm}
(a2)++(.25,0) node (a2') {X'\angordm}
(a1)++(0,\v) node (b1) {\phantom{X'}}
(a2)++(0,\v) node (b2) {\phantom{X'}}
(b1)++(-.2,0) node (b1') {X\angordn}
(b2)++(.25,0) node (b2') {X'\angordn}
;
\draw[1cell=.8]
(a1) edge[bend left=\t] node[pos=.4] {\theta_{\angordm}} (a2)
(a1) edge[bend right=\t] node[swap,pos=.6] {\ups_{\angordm}} (a2)
(b1) edge[bend left=\t] node[pos=.4] {\theta_{\angordn}} (b2)
(b1) edge[bend right=\t] node[swap,pos=.6] {\ups_{\angordn}} (b2)
(a1) edge node[swap] {X\upom} (b1)
(a2) edge node {X'\upom} (b2)
;
\draw[2cell=.9]
node[between=a1 and a2 at .37, rotate=-90, 2label={above, \Theta_{\angordm}}] {\Rightarrow}
node[between=b1 and b2 at .37, rotate=-90, 2label={above, \Theta_{\angordn}}] {\Rightarrow}
;
\end{tikzpicture}
\end{equation}
\end{description}
Identities, horizontal composition, and vertical composition of 2-cells are given componentwise using the components in \cref{ggcat_inthom_Theta}.   A 2-cell is automatically pointed in the sense that its $\vstar$-component $\Theta_{\vstar}$ is the identity natural transformation on the identity functor on $\bone$.
\end{explanation}

 

%% file: chap/kgo_iii.tex
\section{$\Gskg$-Categories from Operadic Pseudoalgebras: Objects}
\label{sec:jemg}

The following assumption is in effect for the rest of this chapter.

\begin{assumption}[$\Tinf$-Operads]\label{as:OpA}
$(\Op,\ga,\opu,\pcom)$ is a pseudo-commutative operad in $\Gcat$ for an arbitrary group $G$ \pcref{def:pseudocom_operad} such that $\Op(1)$ is a terminal $G$-category.  Such a $\Gcat$-operad is called a \index{T-infinity operad@$\Tinf$-operad}\index{operad!T-infinity@$\Tinf$}\emph{$\Tinf$-operad}.  
\end{assumption}

\begin{example}\label{ex:Tinf}
The Barratt-Eccles operad $\BE$ and the $G$-Barratt-Eccles operad $\GBE$ \pcref{def:BE,def:GBE} are $\Tinf$-operads.  More generally, if $\Op$ is a $\Tinf$-operad, then so is the $\Gcat$-operad $\Catg(\EG,\Op)$ in \cref{catgego} \eqref{ng_i}.
\end{example}

This section and \cref{sec:jemg_morphisms,sec:jgosg_onecells,sec:jgosg_twocells} recall the (strong) $J$-theory 2-functors
\[\AlglaxO \fto{\Jgo} \GGCatii \andspace \AlgpspsO \fto{\Jgosg} \GGCatii\]
for a $\Tinf$-operad $\Op$.
\begin{itemize}
\item The domain of $\Jgo$ is the 2-category $\AlglaxO$ \pcref{oalgps_twocat} of $\Op$-pseudoalgebras, lax $\Op$-morphisms, and $\Op$-transformations \pcref{def:pseudoalgebra,def:laxmorphism,def:algtwocells}.  The domain $\AlgpspsO$ of $\Jgosg$ has $\Op$-pseudomorphisms as 1-cells.
\item The codomain is the 2-category $\GGCatii$ \pcref{def:GGCat} of $\Gskg$-categories, natural transformations, and modifications.
\end{itemize}  
The (strong) $J$-theory 2-functors $\Jgo$ and $\Jgosg$ are obtained from the (strong) $J$-theory $\Gcat$-multifunctors 
\[\MultpsO \fto{\Jgo} \GGCat \andspace \MultpspsO \fto{\Jgosg} \GGCat\]
in \cite[\namecref{EqK:thm:Jgo_multifunctor} \ref*{EqK:thm:Jgo_multifunctor}]{yau-eqk} by restricting to $G$-fixed 1-ary 1-cells and 2-cells.  There is no change in the object assignments.  Thus, for an $\Op$-pseudoalgebra $\A$, $\Jgosg\A$ is the same $\Gskg$-category regardless of whether $\Jgosg$ is regarded as a $\Gcat$-multifunctor or as a 2-functor, and likewise for $\Jgo$.  

\secoutline
For the rest of this section, $(\A,\gaA,\phiA)$ denotes an $\Op$-pseudoalgebra for a $\Tinf$-operad $(\Op,\ga,\opu,\pcom)$.  The object assignments of $\Jgo$ and $\Jgosg$ send $\A$ to $\Gskg$-categories \cref{ggcat_obj}, which mean pointed functors  
\[\Jgo\A = \Adash \andspace \Jgosg\A = \Asgdash \cn \Gsk \to \Gcatst.\]
This section recalls from \cite[\namecref{EqK:sec:jemg_objects} \ref*{EqK:sec:jemg_objects}]{yau-eqk} the small pointed $G$-categories
\[(\Jgo\A)\angordn = \Aangordn \andspace (\Jgosg\A)\angordn = \Asgangordn\]
for each object $\angordn \in \Gsk$ \cref{Gsk_objects}.  Their objects are called, respectively, \emph{$\angordn$-systems} and \emph{strong $\angordn$-systems}.
\begin{itemize}
\item \cref{not:compk} fixes notation for substitution and partition. 
\item \cref{def:nsystem} defines \emph{$\angordn$-systems in $\A$}.  An $\angordn$-system in $\A$ consists of a collection of objects in $\A$, together with gluing morphisms relating these objects.  Strong $\angordn$-systems are those with invertible gluing morphisms.
\item \cref{def:nsystem_morphism} defines the small pointed categories $\Aangordn$ of $\angordn$-systems and $\Asgangordn$ of strong $\angordn$-systems.
\item \cref{def:Aangordn_gcat} defines the pointed $G$-categories $\Aangordn$ and $\Aangordnsg$.
\item \cref{ex:nsystem_ones} explains $\Aangordn$ and $\Aangordnsg$ when $\angordn$ is a tuple of copies of $\ord{1}$.
\end{itemize}

\subsection*{Systems}
For $m \geq 0$, recall that $\ufs{m} = \{1 < \cdots <m\}$ denotes an unpointed finite set with $m$ elements \cref{ufsn} and that $\ord{m} = \{0 < 1 < \cdots < m\}$ denotes a pointed finite set with basepoint 0 \cref{ordn}. 
The following notation is used to denote substitution of entries in a tuple and partitions of sets.

\begin{notation}[Substitution and Partition]\label{not:compk}
Suppose $\ang{s}  = \ang{s_i}_{i \in \ufs{m}}$ is an $m$-tuple of symbols for some $m > 0$. 
\begin{itemize}
\item For any $k \in \ufs{m}$ and any symbol $t$, the notation\index{substitution}
\begin{equation}\label{compk}
\ang{s} \compk t = \big(s_1, \ldots, s_{k-1}, t, s_{k+1}, \ldots, s_m\big)
\end{equation}
denotes the $m$-tuple obtained from $\ang{s}$ by replacing its $k$-th entry by $t$. 
\item For $k \neq \ell \in \ufs{m}$ and a pair of symbols $(t,u)$, the notation
\begin{equation}\label{compkl}
\ang{s} \compk t \compell u = \big(\ang{s} \compk t\big) \compell u = \big(\ang{s} \compell u\big) \compk t
\end{equation}
denotes the $m$-tuple obtained from $\ang{s}$ by replacing its $k$-th entry by $t$ and its $\ell$-th entry by $u$.
\end{itemize}
A \emph{partition}\index{partition} of a set $S$, denoted
\begin{equation}\label{partition}
S = \coprod_{j \in \ufs{p}}\, S_j,
\end{equation} 
is a $p$-tuple $\ang{S_j}_{j \in \ufs{p}}$ of pairwise disjoint, possibly-empty subsets of $S$ whose union is $S$.  Pairwise disjointness means
\[S_i \cap S_j = \emptyset \forspace i \neq j \in \ufs{p}.\]
The case $p=0$ can only happen if $S$ is empty.
\end{notation}

In \cref{def:nsystem}, note that the pseudo-commutative structure of $\Op$ \cref{pseudocom_isos} is only used in the commutativity axiom \cref{system_commutativity}.

\begin{definition}[$\angordn$-Systems]\label{def:nsystem}
Given a $\Tinf$-operad $(\Op,\ga,\opu,\pcom)$ \pcref{as:OpA}, an $\Op$-pseudoalgebra $(\A,\gaA,\phiA)$ \pcref{def:pseudoalgebra}, and an object \cref{Gsk_objects}
\[\angordn = \ang{\ord{n}_j}_{j \in \ufs{q}} = \big(\ord{n}_1, \ldots, \ord{n}_q\big) \in \Gskel\]
of length $q > 0$ with each $\ord{n}_j \in \Fskel$ \pcref{def:Fsk}, an \emph{$\angordn$-system in $\A$}\index{system} is defined to be a pair
\begin{equation}\label{nsystem}
(a,\glu)
\end{equation}
consisting of the following data.
\begin{description}
\item[Component objects]
For each $q$-tuple of subsets
\begin{equation}\label{marker}
\ang{s} = \ang{s_j}_{j \in \ufs{q}} = \ang{s_j \subseteq \ufs{n}_j}_{j \in \ufs{q}}
\end{equation}
with $\ufs{n}_j = \ord{n}_j \setminus \{0\}$, $(a,\glu)$ consists of an \emph{$\ang{s}$-component object}\index{component object}
\begin{equation}\label{a_angs}
a_{\ang{s}} \in \A.
\end{equation}
The $q$-tuple $\ang{s}$ is called the \emph{marker}\index{marker} of the object $a_{\ang{s}}$.
\item[Gluing]
For each object $x \in \Op(r)$ with $r \geq 0$, marker $\ang{s}$, index $k \in \ufs{q}$, and partition of $s_k$ into $r$ subsets
\begin{equation}\label{subset-partition}
s_k = \coprod_{i \in \ufs{r}} \, s_{k,i} \subseteq \ufs{n}_k,
\end{equation}
$(a,\glu)$ consists of a \emph{gluing morphism}\index{gluing morphism} at $\big(x; \angs, k, \ang{s_{k,i}}_{i \in \ufs{r}}\big)$:
\begin{equation}\label{gluing-morphism}
\gaA_r\big(x\sscs \ang{a_{\ang{s} \compk\, s_{k,i}}}_{i \in \ufs{r}} \big) 
\fto{\glu_{x; \ang{s} \csp k, \ang{s_{k,i}}_{i \in \ufs{r}}}} a_{\ang{s}} \inspace \A.
\end{equation}
In the domain of \cref{gluing-morphism}, each marker 
\begin{equation}\label{angs_compk_ski}
\ang{s} \compk s_{k,i} = \big(s_1, \ldots, s_{k-1}, s_{k,i} \scs s_{k+1}, \ldots, s_q \big)
\end{equation}
is obtained from $\ang{s}$ by replacing the $k$-th subset $s_k$ by $s_{k,i} \subseteq \ufs{n}_k$, as defined in \cref{compk}.  The gluing morphism in \cref{gluing-morphism} is also denoted by $(a,\glu)_{x; \ang{s} \csp k, \ang{s_{k,i}}_{i \in \ufs{r}}}$.
\end{description}
The pair $(a,\glu)$ is required to satisfy the axioms \cref{system_obj_unity,system_naturality,system_unity_i,system_unity_iii,system_equivariance,system_associativity,system_commutativity} whenever they are defined.
\begin{description}
\item[Object unity]
If $s_j = \emptyset \subseteq \ufs{n}_j$ for some $j \in \ufs{q}$, then
\begin{equation}\label{system_obj_unity}
a_{\ang{s}} = \zero = \gaA_0(*) \in \A,
\end{equation}
the basepoint of $\A$ \cref{pseudoalg_zero}, where $* \in \Op(0)$ is the unique object.
\item[Naturality]
For each morphism $h \cn x \to y$ in the $G$-category $\Op(r)$ with $r \geq 0$, the following diagram in $\A$ commutes.
\begin{equation}\label{system_naturality}

\end{equation}
\item[Unity]
The gluing morphism $\glu_{x; \ang{s} \csp k, \ang{s_{k,i}}_{i \in \ufs{r}}}$ in \cref{gluing-morphism} is the identity morphism in each of the following two cases.
\begin{itemize}
\item If $s_j = \emptyset$ for some $j \in \ufs{q}$, then the following diagram in $\A$ commutes.
\begin{equation}\label{system_unity_i}
%
\end{equation}
\begin{itemize}
\item The vertical equalities in \cref{system_unity_i} follow from the object unity axiom \cref{system_obj_unity} and, if $j=k$, the fact that $s_{k,i} \subseteq s_k$.
\item The equality labeled $\mathbf{b}$ follows from \cref{phi_id}.
\end{itemize}
\item If $r=1$---which implies that $x$ is the operadic unit $\opu \in \Op(1)$---then the following diagram in $\A$ commutes.
\begin{equation}\label{system_unity_iii}
%
\end{equation}
The equality in \cref{system_unity_iii} follows from the action unity axiom \cref{pseudoalg_action_unity} for $\A$.
\end{itemize}
\item[Equivariance]
For each permutation $\si \in \Sigma_r$, the following diagram in $\A$ commutes.
\begin{equation}\label{system_equivariance}
%
\end{equation}
In \cref{system_equivariance}, the left vertical equality follows from the action equivariance axiom \cref{pseudoalg_action_sym} for $\A$.  The bottom horizontal gluing morphism uses the partition
\[s_k = \coprod_{i \in \ufs{r}}\, s_{k,\sigmainv(i)}.\]
\item[Associativity]
Suppose we are given objects
\[\begin{split}
\left(x \sscs \ang{x_i}_{i \in \ufs{r}} \right) &\in \Op(r) \times \txprod_{i \in \ufs{r}}\, \Op(t_i) \andspace\\
\boldx = \ga\big(x\sscs \ang{x_i}_{i \in \ufs{r}} \big) &\in \Op(t)
\end{split}\]
with $t = \sum_{i \in \ufs{r}} t_i$, a marker $\ang{s} = \ang{s_j \subseteq \ufs{n}_j}_{j \in \ufs{q}}$ as defined in \cref{marker}, an index $k \in \ufs{q}$, and partitions
\[s_k = \coprod_{i \in \ufs{r}} s_{k,i} \andspace s_{k,i} = \coprod_{\ell \in \ufs{t}_i} s_{k,i,\ell}\]
of, respectively, $s_k$ and $s_{k,i}$ for each $i \in \ufs{r}$.  Then the following diagram in $\A$ commutes.
\begin{equation}\label{system_associativity}
\begin{tikzpicture}[vcenter]
\def\h{4} \def\v{-2} \def\t{2em}
\draw[0cell=.85]
(0,0) node (a) {\gaA_r\big(x \sscs \bang{\gaA_{t_i}\big(x_i \sscs \ang{a_{\ang{s} \compk\, s_{k,i,\ell}}}_{\ell \in \ufs{t}_i}  \big)}_{i \in \ufs{r}} \big)}
(a)++(-\h,\v/2) node (b) {\gaA_r\big(x \sscs \ang{a_{\ang{s} \compk\, s_{k,i}}}_{i \in \ufs{r}} \big)}
(a)++(0,\v) node (c) {\gaA_t\big(\boldx \sscs \ang{\ang{a_{\ang{s} \compk\, s_{k,i,\ell}}}_{\ell \in \ufs{t}_i}}_{i \in \ufs{r}} \big)}
(c)++(-.5em,0) node (c') {\phantom{\gaA_t\Big(\boldx \sscs \bang{\bang{a_{\ang{s} \compk\, s_{k,i,\ell}}}_{\ell \in \ufs{t}_i}}_{i \in \ufs{r}} \Big)}}
(b)++(0,\v) node (d) {\phantom{a_{\ang{s}}}}
(d)++(\t,0) node (d') {a_{\ang{s}}}
;
\draw[1cell=.85]
(a) edge[bend right=8, shorten >=-0ex] node[swap,pos=.6,inner sep=0pt] {\gaA_r\big(x \sscs \ang{\glu_{x_i;\, \ang{s} \compk\, s_{k,i} \scs k, \ang{s_{k,i,\ell}}_{\ell \in \ufs{t}_i}}}_{i \in \ufs{r}} \big)} (b)
(b) edge[bend right=10,transform canvas={xshift=\t}] node[swap] {\glu_{x; \ang{s} \csp k, \ang{s_{k,i}}_{i \in \ufs{r}}}} (d)
(a) edge[bend left=10,transform canvas={xshift=-\t}] node {\phiA_{(r;\, t_1,\ldots,t_r)}} node[swap] {\iso} (c)
(c') edge[bend left=6] node {\glu_{\boldx; \ang{s} \csp k, \ang{\ang{s_{k,i,\ell}}_{\ell \in \ufs{t}_i}}_{i \in \ufs{r}}}} (d')
;
\end{tikzpicture}
\end{equation}
In \cref{system_associativity}, $\phiA_{(r;\, t_1,\ldots,t_r)}$ is the component of the associativity constraint \cref{phiA} of $\A$ at the objects $x$, $\ang{x_i}_{i \in \ufs{r}}$, and $\ang{\ang{a_{\ang{s} \compk\, s_{k,i,\ell}}}_{\ell \in \ufs{t}_i}}_{i \in \ufs{r}}$.  The bottom gluing morphism uses the partition
\[s_k = \coprod_{i \in \ufs{r}} \coprod_{\ell \in \ufs{t}_i}\, s_{k,i,\ell}.\] 
\item[Commutativity]
Suppose we are given a pair of objects
\[(x,y) \in \Op(r) \times \Op(t),\]
a marker $\ang{s} = \ang{s_j \subseteq \ufs{n}_j}_{j \in \ufs{q}}$, two distinct indices $k,\ell \in \ufs{q}$, and partitions
\[s_k = \coprod_{i \in \ufs{r}} s_{k,i} \subseteq \ufs{n}_k \andspace 
s_\ell = \coprod_{p \in \ufs{t}} s_{\ell,p} \subseteq \ufs{n}_\ell\]
of, respectively, $s_k$ and $s_\ell$.
 Then the following diagram in $\A$ commutes.
\begin{equation}\label{system_commutativity}

\end{equation}
\begin{itemize}
\item In \cref{system_commutativity}, the marker
\[\ang{s} \compk\, s_{k,i}\, \compell\, s_{\ell,p}\]
is obtained from $\ang{s}$ by replacing its $k$-th subset by $s_{k,i}$ and its $\ell$-th subset by $s_{\ell,p}$, as defined in \cref{compkl}.
\item The equality labeled $\mathbf{eq}$ follows from the action equivariance axiom \cref{pseudoalg_action_sym} for $\A$, applied to the transpose permutation $\twist_{t,r} \in \Sigma_{tr}$ defined in \cref{eq:transpose_perm}.  The transpose permutation $\twist_{t,r}$ switches the order of indexing from $\ang{\ang{\Cdots}_{i \in \ufs{r}}}_{p \in \ufs{t}}$ to $\ang{\ang{\Cdots}_{p \in \ufs{t}}}_{i \in \ufs{r}}$.
\item Each $\intr$ is an instance of the intrinsic pairing of $\Op$ (\cref{def:intrinsic_pairing}), with
\[\begin{split}
x \intr y &= \ga\big(x \sscs \ang{y}_{i \in \ufs{r}}\big) \in \Op(rt) \andspace\\
y \intr x &= \ga\big(y \sscs \ang{x}_{p \in \ufs{t}}\big) \in \Op(tr).
\end{split}\]
The top boundary arrow in \cref{system_commutativity} involves the $(x,y)$-component of the $(r,t)$-pseudo-commutativity isomorphism \cref{pseudocom_isos} of $\Op$:
\[(x \intr y) \twist_{t,r} \fto[\iso]{\pcom_{r,t}} y \intr x.\]
\item $\phiA_{(t;\, r,\ldots,r)}$ and $\phiA_{(r;\, t,\ldots,t)}$ are instances of the associativity constraint \cref{phiA} of $\A$.  
\end{itemize}
\end{description}
This finishes the definition of an $\angordn$-system $(a,\glu)$.  

Moreover, we define the following.
\begin{itemize}
\item A \emph{strong $\angordn$-system}\dindex{strong}{system} is an $\angordn$-system such that each gluing morphism $\glu_{x; \ang{s} \csp k, \ang{s_{k,i}}_{i \in \ufs{r}}}$ in \cref{gluing-morphism} is an isomorphism.
\item The \emph{base $\angordn$-system}\dindex{base}{system} is the $\angordn$-system $(\zero,1_\zero)$\label{not:basesystem} with
\begin{itemize}
\item each component object given by the basepoint $\zero \in \A$ and
\item each gluing morphism given by the identity morphism $1_\zero$ of $\zero$.\defmark
\end{itemize}
\end{itemize}
\end{definition}

\subsection*{$G$-Categories of $\angordn$-Systems}
Morphisms of $\angordn$-systems are defined next.

\begin{definition}[Morphisms of $\angordn$-Systems]\label{def:nsystem_morphism}
Under the same hypotheses as \cref{def:nsystem}, suppose $(a,\glu^a)$ and $(b,\glu^b)$ are $\angordn$-systems in $\A$.  A \emph{morphism}\index{system!morphism} of $\angordn$-systems
\begin{equation}\label{nsystem_mor}
(a,\glu^a) \fto{\theta} (b,\glu^b)
\end{equation}
consists of, for each marker $\ang{s} = \ang{s_j \subseteq \ufs{n}_j}_{j \in \ufs{q}}$, an \emph{$\ang{s}$-component morphism}\index{component morphism}
\begin{equation}\label{theta_angs}
a_{\ang{s}} \fto{\theta_{\ang{s}}} b_{\ang{s}} \inspace \A
\end{equation}
such that the following two axioms are satisfied.
\begin{description}
\item[Unity] If $s_j = \emptyset \subseteq \ufs{n}_j$ for some $j \in \ufs{q}$, then there is an equality of morphisms
\begin{equation}\label{nsystem_mor_unity}
a_{\ang{s}} = \zero \fto{\theta_{\ang{s}} \,=\, 1_\zero} b_{\ang{s}} = \zero. 
\end{equation}
The object equalities $a_{\ang{s}} = b_{\ang{s}} = \zero$ follow from the object unity axiom \cref{system_obj_unity}.
\item[Compatibility] For each object $x \in \Op(r)$ with $r \geq 0$, marker $\ang{s}$, index $k \in \ufs{q}$, and partition
\[s_k = \coprod_{i \in \ufs{r}} \, s_{k,i} \subseteq \ufs{n}_k,\]
the following diagram in $\A$ commutes.
\begin{equation}\label{nsystem_mor_compat}
\begin{tikzpicture}[vcenter]
\def\v{-1.5} \def\u{-.08}
\draw[0cell]
(0,0) node (a1) {\gaA_r\big(x \sscs \ang{a_{\ang{s} \compk\, s_{k,i}}}_{i \in \ufs{r}} \big)}
(a1)++(0,\v) node (b1) {\gaA_r\big(x \sscs \ang{b_{\ang{s} \compk\, s_{k,i}}}_{i \in \ufs{r}} \big)}
(a1)++(5,0) node (a2) {\phantom{a_{\ang{s}}}}
(a2)++(0,\u) node (a2') {a_{\ang{s}}}
(a2)++(0,\v) node (b2) {\phantom{b_{\ang{s}}}}
(b2)++(0,\u) node (b2') {b_{\ang{s}}}
;
\draw[1cell=.9]
(a1) edge node {\glu^a_{x; \ang{s} \csp k, \ang{s_{k,i}}_{i \in \ufs{r}}}} (a2)
(b1) edge node {\glu^b_{x; \ang{s} \csp k, \ang{s_{k,i}}_{i \in \ufs{r}}}} (b2)
(a1) edge[transform canvas={xshift=3em}] node[swap] {\gaA_r\big(1_x \sscs \ang{\theta_{\ang{s} \compk\, s_{k,i}}}_{i \in \ufs{r}} \big)} (b1)
(a2') edge node {\theta_{\ang{s}}} (b2')
;
\end{tikzpicture}
\end{equation}
\end{description}
This finishes the definition of a morphism of $\angordn$-systems.  

Moreover, we define the following.
\begin{itemize}
\item Composition and identity morphisms of $\angordn$-systems in $\A$ are defined componentwise using \cref{theta_angs}.  Denote by $\Aangordn$ the category of $\angordn$-systems in $\A$ and their morphisms.
\item The pointed category of \emph{$\angordn$-systems in $\A$}\index{system!pointed category} is defined as follows.\label{not:Aangordn} 
\begin{equation}\label{A_angordn}
\Aangordn = \begin{cases}
\boldone & \text{if $\angordn = \vstar$},\\
(\A,\zero) & \text{if $\angordn = \ang{}$, and}\\
\big(\Aangordn, (\zero,1_\zero)\big) & \text{if $\angordn \neq \vstar, \ang{}$.}
\end{cases}
\end{equation}
The pointed category 
\begin{equation}\label{Aangordnsg}
\Aangordnsg
\end{equation} 
of \emph{strong $\angordn$-systems in $\A$} is defined as the pointed full subcategory of $\Aangordn$ with strong $\angordn$-systems in $\A$ as objects if $\angordn \neq \vstar, \ang{}$.  The other two cases are defined as
\[\Asgstar = \bone \andspace \Asgang = (\A,\zero).\]
\end{itemize}
This finishes the definition.
\end{definition}

Recall from \cref{def:ptGcat} that a pointed $G$-category is a $G$-category equipped with a distinguished $G$-fixed object called the basepoint.

\begin{definition}[Pointed $G$-Categories of $\angordn$-Systems]\label{def:Aangordn_gcat}
The pointed category $\Aangordn$ (\cref{def:nsystem_morphism}) is extended to a pointed $G$-category\index{system!pointed G-category@pointed $G$-category} as follows.  If $\angordn = \vstar$ or $\ang{}$, then
\begin{equation}\label{vstar_system}
\sys{\A}{\vstar} = \boldone \andspace \sys{\A}{\ang{}} = (\A,\zero)
\end{equation}
are already pointed $G$-categories.  Suppose $\angordn \in \Gsk \setminus \{\vstar, \ang{}\}$ is an object of length $q > 0$. 
\begin{description}
\item[$G$-action on systems]  
For an element $g \in G$ and an $\angordn$-system $(a,\glu) \in \Aangordn$ (\cref{def:nsystem}), the $\angordn$-system in $\A$ 
\begin{equation}\label{nsystem_gaction}
g \cdot (a,\glu) = (ga, g\glu)
\end{equation}
is defined as follows.
\begin{description}
\item[Component objects] For each marker $\ang{s} = \ang{s_j \subseteq \ufs{n}_j}_{j \in \ufs{q}}$, the $\ang{s}$-component object of $(ga,g\glu)$ is the image of $a_{\angs}$ under the $g$-action on $\A$:
\begin{equation}\label{ga_scomponent}
(ga)_{\ang{s}} = g a_{\ang{s}} \in \A.
\end{equation}
\item[Gluing] Given an object $x \in \Op(r)$ with $r \geq 0$, a marker $\ang{s}$, an index $k \in \ufs{q}$, and a partition
\[s_k = \coprod_{i \in \ufs{r}}\, s_{k,i} \subseteq \ufs{n}_k\]
of $s_k$ into $r$ subsets, the gluing morphism of $(ga,g\glu)$ at $(x; \angs, k, \ang{s_{k,i}}_{i \in \ufs{r}})$ is defined by the following commutative diagram.
\begin{equation}\label{ga_gluing}
\begin{tikzpicture}[vcenter]
\def\v{-1}
\draw[0cell=.8]
(0,0) node (a1) {\gaA_r\big(x \sscs \ang{(ga)_{\ang{s} \compk\, s_{k,i}} }_{i \in \ufs{r}} \big)}
(a1)++(0,\v) node (a2) {\gaA_r\big(x \sscs \ang{g a_{\ang{s} \compk\, s_{k,i}}}_{i \in \ufs{r}}\big)}
(a2)++(0,\v) node (a3) {g \gaA_r \big(\ginv x \sscs \ang{a_{\ang{s} \compk\, s_{k,i}}}_{i \in \ufs{r}} \big)}
(a1)++(5,0) node (b1) {(ga)_{\ang{s}}}
(b1)++(0,2*\v) node (b3) {g a_{\ang{s}}}
;
\draw[1cell=.9]
(a1) edge node {(g\glu)_{x; \ang{s},\, k, \ang{s_{k,i}}_{i \in \ufs{r}}}} (b1)
(a3) edge node {g \glu_{\ginv x; \ang{s},\, k, \ang{s_{k,i}}_{i \in \ufs{r}}}} (b3)
(a1) edge[equal,shorten >=-.5ex] (a2)
(a2) edge[equal,shorten >=-.5ex] node[swap,inner sep=4pt] {\mathbf{f}} (a3)
(b1) edge[equal] (b3)
;
\end{tikzpicture}
\end{equation}
\begin{itemize}
\item The two unlabeled equalities in \cref{ga_gluing} follow from the definition of $(ga)_{\ang{s}}$ in \cref{ga_scomponent}.
\item The equality labeled $\mathbf{f}$ follows from the $G$-functoriality of the $r$-th $\Op$-action $G$-functor $\gaA_r$ \cref{gaAn}.
\item In the bottom horizontal arrow in \cref{ga_gluing}, 
\[\gaA_r \big(\ginv x \sscs \ang{a_{\ang{s} \compk\, s_{k,i}}}_{i \in \ufs{r}} \big)
\fto{\glu_{\ginv x; \ang{s},\, k, \ang{s_{k,i}}_{i \in \ufs{r}}}} a_{\ang{s}}\]
is the gluing morphism of $(a,\glu)$ at $(\ginv x; \ang{s}, k, \ang{s_{k,i}}_{i \in \ufs{r}})$, and $g\glu_{\cdots}$ is its image under the $g$-action on $\A$.
\end{itemize}
\end{description}
This finishes the definition of the $\angordn$-system $(ga,g\glu)$.
\item[$G$-action on morphisms]  
For a morphism of $\angordn$-systems in $\A$ (\cref{def:nsystem_morphism})
\[(a,\glu^a) \fto{\theta} (b,\glu^b),\]
the morphism of $\ordtu{n}$-systems
\begin{equation}\label{gtheta}
(ga,g\glu^a) \fto{g\theta} (gb, g\glu^b)
\end{equation}
is defined by, for each marker $\ang{s} = \ang{s_j \subseteq \ufs{n}_j}_{j \in \ufs{q}}$, the $\ang{s}$-component morphism
\begin{equation}\label{gtheta_angs}
(ga)_{\ang{s}} = ga_{\ang{s}}
\fto{(g\theta)_{\ang{s}} \,=\, g\theta_{\ang{s}}}
(gb)_{\ang{s}} = gb_{\ang{s}}.
\end{equation}
This finishes the definition of the pointed $G$-category $\Aangordn$.
\item[Strong variant]  
The pointed category of strong $\angordn$-systems in $\A$ 
\begin{equation}\label{sgAordnbe}
\Aangordnsg
\end{equation}
is extended to a pointed $G$-category by restricting the $G$-action on $\Aangordn$ defined in \crefrange{vstar_system}{gtheta_angs} to the full subcategory $\Aangordnsg$.\defmark
\end{description}
\end{definition}

\begin{example}[Tuple of $\ord{1}$]\label{ex:nsystem_ones}
Consider \cref{def:nsystem,def:nsystem_morphism,def:Aangordn_gcat} for the object
\[\angordn = \ordtu{1}_{j \in \ufs{q}} = \big(\ord{1},\ldots,\ord{1}\big) \in \Gsk\]
consisting of $q$ copies of $\ord{1}$.  The only subsets of $\ufs{1} = \{1\}$ are $\emptyset$ and $\{1\}$.  There is a canonical isomorphism of pointed $G$-categories
\begin{equation}\label{Aordtuone}
(\A,\zero) \iso \big(\A\ordtu{1}, (\zero,1_\zero)\big) 
= \big(\syssg{\A}{\ordtu{1}}, (\zero,1_\zero)\big).
\end{equation}
Under this pointed isomorphism, an object $b \in \A$ corresponds to the $\ordtu{1}$-system with
\begin{itemize}
\item $\ang{\{1\}}_{j \in \ufs{q}}$-component object given by $b$, 
\item all other component objects given by the basepoint $\zero \in \A$, and 
\item all gluing morphisms given by identities.
\end{itemize}
Such a $\ordtu{1}$-system is necessarily strong because each gluing morphism is invertible. 
\end{example}

 

%% file: chap/kgo_iv.tex
\section{$\Gskg$-Categories from Operadic Pseudoalgebras: Morphisms}
\label{sec:jemg_morphisms}

This section finishes the construction of the object assignments of the (strong) $J$-theory 2-functors
\[\AlglaxO \fto{\Jgo} \GGCatii \andspace \AlgpspsO \fto{\Jgosg} \GGCatii\]
by recalling from \cite[Sections \ref*{EqK:sec:jemg_morphisms} and \ref*{EqK:sec:jemg_morphisms_ii}]{yau-eqk} the morphism assignments of the $\Gskg$-categories 
\begin{equation}\label{JgoJgosg_arrows}
\begin{tikzcd}
\Gsk \arrow[r, "\Jgo\A \,=\, \Adash", ->,shift left=.5ex]
  \arrow[r, "\Jgosg\A \,=\, \Asgdash" below, ->,shift right=.4ex] & [4.5em] \Gcatst
\end{tikzcd}
\end{equation}
for a $\Tinf$-operad $\Op$ \pcref{as:OpA} and an $\Op$-pseudoalgebra $(\A,\gaA,\phiA)$ \pcref{def:pseudoalgebra}.

\secoutline
Similar to \cref{sec:jemg}, we first consider categories of systems and then restrict to the full subcategories of strong systems.  See \cref{AF_sg}.  For each morphism $\upom \cn \angordm \to \angordn$ in $\Gsk$, the pointed $G$-functor 
\[\Aangordm \fto{\Aupom} \Aangordn\]
is constructed in several steps.  When the domain object $\angordm$ has positive length, the morphism $\upom$ has the form $(f,\angpsi)$ \cref{fangpsi}.  The pointed $G$-functor $\Afpsi$ is the composite of two pointed $G$-functors, denoted $\ftil$ and $\psitil$. 
\begin{itemize}
\item \cref{def:ftil_functor} recalls the pointed $G$-functor
\[\Aangordm \fto{\ftil} \Afangordm\]
from $\angordm$-systems to $f_*\angordm$-systems.
\item \cref{def:psitil_functor} recalls the pointed $G$-functor
\[\Afangordm \fto{\psitil} \Aangordn\]
from $f_*\angordm$-systems to $\angordn$-systems.
\item \cref{def:Afangpsi} defines the pointed $G$-functors $\Aupom$ and $\Aupomsg$ in all possible cases.  
\item \cref{A_ptfunctor} states that $\Jgo\A = \Adash$ and $\Jgosg\A = \Asgdash$ are $\Gskg$-categories.
\end{itemize}

We sometimes restrict attention to morphisms in $\Gsk$ of the form stated in \cref{as:fpsi}.  We will state it explicitly whenever this assumption is in effect.

\begin{assumption}\label{as:fpsi}
Denote by
\[\angordm = \ang{\ordm_i}_{i \in \ufs{p}} \fto{(f, \ang{\psi})} 
\angordn = \ang{\ordn_j}_{j \in \ufs{q}}\]
a morphism in $\Gsk$ consisting of a reindexing injection 
\[\ufs{p} \fto{f} \ufs{q} \withspace p>0\] 
and a morphism 
\[f_*\angordm =  \bang{\ordm_{\finv(j)}}_{j \in \ufs{q}} 
\fto{\ang{\psi} = \ang{\psi_j}_{j \in \ufs{q}}} \angordn \inspace \Fsk^{(q)},\]
as defined in \cref{reindexing_functor,fangpsi}.
\end{assumption}

\begin{definition}[The $G$-Functor $\ftil$]\label{def:ftil_functor}
Under \cref{as:OpA,as:fpsi},  the pointed $G$-functor
\begin{equation}\label{ftil_functor}
\Aangordm \fto{\ftil} \Afangordm
\end{equation}
is defined as follows.
\begin{description}
\item[Objects] 
Given an $\angordm$-system $(a,\glu) \in \Aangordm$ \pcref{def:nsystem}, the $f_*\angordm$-system
\begin{equation}\label{ftil_aglu}
\ftil(a,\glu) = \left(\atil, \glutil\right) \in \Afangordm
\end{equation}
is defined as follows.
\begin{description}
\item[Component objects] Given a marker
\begin{equation}\label{angs_finv}
\ang{s} = \bang{s_j \subseteq \ufs{m}_{\finv(j)}}_{j \in \ufs{q}},
\end{equation}
we first define the marker
\begin{equation}\label{ftil_angs}
\ftil_*\ang{s} = \bang{s_{f(i)} \subseteq \ufs{m}_i}_{i \in \ufs{p}}
\end{equation}
obtained from $\ang{s}$ by
\begin{itemize}
\item removing those $s_j$ with $\finv(j) = \emptyset$ and
\item permuting the resulting $p$-tuple according to $\finv$.
\end{itemize}  
The $\ang{s}$-component object of $(\atil,\glutil)$ is defined as
\begin{equation}\label{atil_component}
\atil_{\ang{s}} = \begin{cases}
\zero & \text{if $s_j = \emptyset$ for some $j \in \ufs{q}$, and}\\
a_{\ftil_*\ang{s}} & \text{if $s_j \neq \emptyset$ for each $j \in \ufs{q}$.}
\end{cases}
\end{equation}
\begin{itemize}
\item In \cref{atil_component}, $\zero = \gaA_0(*)$ is the basepoint of $\A$, and the first case is forced by the object unity axiom \cref{system_obj_unity}.
\item In the second case, $a_{\ftil_*\ang{s}}$ is the $\ftil_*\ang{s}$-component object of $(a,\glu)$.  Note that for an index $j \in \ufs{q}$ such that $\finv(j) = \emptyset$---which implies $\ufs{m}_{\finv(j)} = \{1\}$---the condition $s_j \neq \emptyset$ means $s_j = \{1\}$.
\end{itemize}
\item[Gluing] Given an object $x \in \Op(r)$ with $r \geq 0$, a marker $\ang{s}$ as defined in \cref{angs_finv}, an index $k \in \ufs{q}$, and a partition
\[s_k = \coprod_{\ell \in \ufs{r}}\, s_{k,\ell} \subseteq \ufs{m}_{\finv(k)},\]
there are three possible cases of the corresponding gluing morphism of $(\atil,\glutil)$ as follows.
\begin{itemize}
\item If $s_j \neq \emptyset$ for each $j \in \ufs{q}$ and $\finv(k) \neq \emptyset$, then the corresponding gluing morphism of $(\atil,\glutil)$ is defined by the following commutative diagram.
\begin{equation}\label{glutil_component}
\begin{tikzpicture}[vcenter]
\def\v{-1.2}
\draw[0cell=.9]
(0,0) node (a1) {\gaA_r\big(x \sscs \ang{\atil_{\ang{s} \compk\, s_{k,\ell}}}_{\ell \in \ufs{r}}\big)}
(a1)++(5.5,0) node (a2) {\atil_{\ang{s}}}
(a1)++(0,\v) node (b1) {\gaA_r\big(x \sscs \ang{a_{\ftil_*\ang{s} \,\comp_{\finv(k)}\, s_{k,\ell}}}_{\ell \in \ufs{r}}\big)}
(a2)++(0,\v) node (b2) {a_{\ftil_*\ang{s}}}
;
\draw[1cell=.9]
(a1) edge node {\glutil_{x;\, \ang{s},\, k, \ang{s_{k,\ell}}_{\ell \in \ufs{r}}}} (a2)
(b1) edge node {\glu_{x;\, \ftil_*\ang{s},\, \finv(k), \ang{s_{k,\ell}}_{\ell \in \ufs{r}}}} (b2)
(a1) edge[equal,shorten >=-1ex] (b1)
(a2) edge[equal] (b2)
;
\end{tikzpicture}
\end{equation}
The bottom horizontal arrow in \cref{glutil_component} is a gluing morphism of $(a,\glu)$.
\item If $s_j = \emptyset$ for some $j \in \ufs{q}$, then the unity axiom \cref{system_unity_i} forces the definition
\begin{equation}\label{glutil_component_ii}
\glutil_{x;\, \ang{s},\, k, \ang{s_{k,\ell}}_{\ell \in \ufs{r}}} = 1_\zero \cn \zero \to \zero.
\end{equation}
\item If $s_j \neq \emptyset$ for each $j \in \ufs{q}$ and $\finv(k) = \emptyset$, then we define
\begin{equation}\label{glutil_component_iii}
\glutil_{x;\, \ang{s},\, k, \ang{s_{k,\ell}}_{\ell \in \ufs{r}}} 
= 1_{a_{\ftil_*\ang{s}}} \cn a_{\ftil_*\ang{s}} \to a_{\ftil_*\ang{s}}.
\end{equation}
\end{itemize}
\end{description}
\item[Morphisms] 
Given a morphism of $\angordm$-systems \pcref{def:nsystem_morphism}
\[(a,\glu^a) \fto{\theta} (b,\glu^b),\]
the morphism of $f_*\angordm$-systems
\[(\atil,\glutil^a) \fto{\ftil(\theta) = \thatil} (\btil,\glutil^b)\]
has, for each marker $\ang{s}$ as defined in \cref{angs_finv}, $\ang{s}$-component morphism defined as
\begin{equation}\label{thatil_component}
\begin{split}
&\big(\atil_{\ang{s}} \fto{\thatil_{\ang{s}}} \btil_{\ang{s}}\big)\\ 
&= \begin{cases}
\zero \fto{1_\zero} \zero & \text{if $s_j = \emptyset$ for some $j \in \ufs{q}$, and}\\
a_{\ftil_*\ang{s}} \fto{\theta_{\ftil_*\ang{s}}} b_{\ftil_*\ang{s}} & \text{if $s_j \neq \emptyset$ for each $j \in \ufs{q}$.}
\end{cases}
\end{split}
\end{equation}
In \cref{thatil_component}, $\theta_{\ftil_*\ang{s}}$ is the $\ftil_*\ang{s}$-component morphism of $\theta$.  The first case of \cref{thatil_component} is forced by the unity axiom \cref{nsystem_mor_unity}.  
\end{description}
This finishes the definition of the pointed $G$-functor $\ftil$ in \cref{ftil_functor}.
\end{definition}

Next, we recall the second pointed $G$-functor that comprises $\Afpsi$.

\begin{definition}[The $G$-Functor $\psitil$]\label{def:psitil_functor}
Under \cref{as:OpA,as:fpsi}, the pointed $G$-functor
\begin{equation}\label{psitil_functor}
\Afangordm \fto{\psitil} \Aangordn
\end{equation}
is defined as follows.
\begin{description}
\item[Objects] 
Given an $f_*\angordm$-system $(a,\glu) \in \Afangordm$ \pcref{def:nsystem}, the $\angordn$-system
\[\psitil(a,\glu) = (a^{\psitil}, \glu^{\psitil}) \in \Aangordn\]
is defined as follows.
\begin{description}
\item[Component objects] Given a marker
\begin{equation}\label{angs_ordtun}
\ang{s} = \ang{s_j \subseteq \ufs{n}_j}_{j \in \ufs{q}},
\end{equation}
recalling that 
\[\ord{m}_{\finv(j)} \fto{\psi_j} \ord{n}_j\]
is a pointed function for each $j \in \ufs{q}$, we first define the marker
\[\psiinv\ang{s} = \bang{\psiinv_j s_j \subseteq \ufs{m}_{\finv(j)}}_{j \in \ufs{q}}.\]
The $\ang{s}$-component object of $(a^{\psitil}, \glu^{\psitil})$ is defined as the $\psiinv\ang{s}$-component object of $(a,\glu)$:
\begin{equation}\label{apsitil_angs}
a^{\psitil}_{\ang{s}} = a_{\psiinv\ang{s}}.
\end{equation}
\item[Gluing] 
Given an object $x \in \Op(r)$ with $r \geq 0$, a marker $\ang{s}$ as defined in \cref{angs_ordtun}, an index $k \in \ufs{q}$, and a partition
\begin{equation}\label{skpartition}
s_k = \coprod_{\ell \in \ufs{r}}\, s_{k,\ell} \subseteq \ufs{n}_k,
\end{equation}
the corresponding gluing morphism of $(a^{\psitil}, \glu^{\psitil})$ is defined by the following commutative diagram.
\begin{equation}\label{glupsitil}
\begin{tikzpicture}[vcenter]
\def\v{-1.2}
\draw[0cell=.9]
(0,0) node (a1) {\gaA_r\big(x \sscs \ang{a^{\psitil}_{\ang{s} \compk\, s_{k,\ell}}}_{\ell \in \ufs{r}}\big)}
(a1)++(5.5,0) node (a2) {a^{\psitil}_{\ang{s}}}
(a1)++(0,\v) node (b1) {\gaA_r\big(x \sscs \ang{a_{\psiinv\ang{s} \,\comp_{k}\, \psiinv_k s_{k,\ell}}}_{\ell \in \ufs{r}}\big)}
(a2)++(0,\v) node (b2) {a_{\psiinv\ang{s}}}
;
\draw[1cell=.9]
(a1) edge node {\glu^{\psitil}_{x;\, \ang{s},\, k, \ang{s_{k,\ell}}_{\ell \in \ufs{r}}}} (a2)
(b1) edge node {\glu_{x;\, \psiinv\ang{s},\, k, \ang{\psiinv_k s_{k,\ell}}_{\ell \in \ufs{r}}}} (b2)
(a1) edge[equal,shorten >=-1ex] (b1)
(a2) edge[equal] (b2)
;
\end{tikzpicture}
\end{equation}
\end{description}
\item[Morphisms] 
Given a morphism of $f_*\angordm$-systems \pcref{def:nsystem_morphism}
\[(a,\glu^a) \fto{\theta} (b,\glu^b),\]
the morphism of $\angordn$-systems
\[(a^{\psitil}, \glu^{a, \psitil}) \fto{\psitil(\theta) = \tha^{\psitil}} (b^{\psitil}, \glu^{b, \psitil})\]
has, for each marker $\ang{s}$ as defined in \cref{angs_ordtun}, $\ang{s}$-component morphism defined as the $\psiinv\ang{s}$-component morphism of $\theta$:
\begin{equation}\label{thapsitil_component}
a^{\psitil}_{\ang{s}} = a_{\psiinv\ang{s}} \fto{\tha^{\psitil}_{\ang{s}} \,=\, \theta_{\psiinv\ang{s}}} 
b^{\psitil}_{\ang{s}} = b_{\psiinv\ang{s}}.
\end{equation}
\end{description}
This finishes the definition of the pointed $G$-functor $\psitil$.
\end{definition}

Next, we recall the pointed $G$-functor $\Aupom$ for a general morphism $\upom$ in $\Gsk$ \cref{Gsk_morphisms}.

\begin{definition}[The $G$-Functor $\Aupom$]\label{def:Afangpsi}
For a $\Tinf$-operad $\Op$ \pcref{as:OpA}, an $\Op$-pseudoalgebra $\A$ \pcref{def:pseudoalgebra}, and a morphism  $\upom \cn \angordm \to \angordn$ in $\Gsk$, the pointed $G$-functor
\begin{equation}\label{AF}
\Aangordm \fto{\Aupom} \Aangordn
\end{equation}
is defined as follows.  First, suppose $\upom$ has the form
\[\angordm = \ang{\ordm_i}_{i \in \ufs{p}} \fto{\upom = (f,\angpsi)} 
\angordn = \ang{\ordn_j}_{j \in \ufs{q}}\]
with $\angordm, \angordn \in \Gsk \setminus \{\vstar, \ang{}\}$, as defined in \cref{as:fpsi}.  We define $\Aupom$ as the composite pointed $G$-functor
\begin{equation}\label{A_fangpsi}
\begin{tikzpicture}[vcenter]
\def\h{2.5} \def\t{.7}
\draw[0cell]
(0,0) node (a1) {\Aangordm}
(a1)++(\h,0) node (a2) {\Afangordm}
(a2)++(\h,0) node (a3) {\Aangordn}
;
\draw[1cell=.9]
(a1) edge node{\ftil} (a2)
(a2) edge node {\psitil} (a3)
;
\draw[1cell=.9]
(a1) [rounded corners=2pt, shorten <=0ex] |- ($(a2)+(-1,\t)$)
-- node {\Afpsi} ($(a2)+(1,\t)$) -| (a3)
;
\end{tikzpicture}
\end{equation}
where $\ftil$ and $\psitil$ are the pointed $G$-functors in, respectively, \cref{def:ftil_functor,def:psitil_functor}.  
\begin{description}
\item[Marginal cases]  
The remaining cases of $\Aupom$ are defined as follows.
\begin{itemize}
\item If $\angordm = \vstar \in \Gsk$, then 
\begin{equation}\label{A_vstar_n}
\sys{\A}{\vstar} = \boldone \fto{\Aupom} \Aangordn
\end{equation}
sends the unique object of $\boldone$ to the base $\angordn$-system $(\zero,1_\zero)$. 
\item If $\angordn = \vstar \in \Gsk$, then 
\begin{equation}\label{A_m_vstar}
\Aangordm \fto{\Aupom} \sys{\A}{\vstar} = \boldone
\end{equation}
is the unique pointed $G$-functor to the terminal pointed $G$-category.
\item If $\angordm = \angordn = \ang{}$, the empty sequence, then, by \cref{Gsk_empty_mor}, $\upom \cn \ang{} \to \ang{}$ is either the identity morphism $1_{\ang{}}$ or the 0-morphism $\ang{} \to \vstar \to \ang{}$.
In these two cases, we define 
\begin{equation}\label{A_empty_empty}
\sys{\A}{\ang{}} = \A \fto{\Aupom} \sys{\A}{\ang{}} = \A
\end{equation}
as, respectively, the identity functor and the constant functor at the $G$-fixed basepoint $\zero \in \A$.
\item Suppose $\angordm = \ang{}$ and $\angordn \in \Gsk \setminus \{\vstar, \ang{}\}$ has length $q>0$.  Then the morphism $\upom$ factors as the following composite in $\Gsk$ \cref{Gsk_composite}.
\begin{equation}\label{Fiqangpsi}
\begin{tikzpicture}[vcenter]
\def\t{.7} \def\h{3}
\draw[0cell]
(0,0) node (a1) {\ang{}}
(a1)++(\h,0) node (a2) {\ordtu{1}_{j \in \ufs{q}}}
(a2)++(\h,0) node (a3) {\angordn}
;
\draw[1cell=.9]
(a1) edge node {(\im_q,\ang{1_{\ord{1}}}_{j \in \ufs{q}})} (a2)
(a2) edge node {(1_{\ufs{q}}, \ang{\psi_j}_{j \in \ufs{q}})} (a3)
;
\draw[1cell=1]
(a1) [rounded corners=2pt, shorten <=0ex] |- ($(a2)+(-1,\t)$)
-- node {\upom \,=\, (\im_q, \ang{\psi_j}_{j \in \ufs{q}})} ($(a2)+(1,\t)$) -| (a3)
;
\end{tikzpicture}
\end{equation}
In this composite,
\begin{itemize}
\item $\im_q \cn \emptyset \to \ufs{q}$ is the unique reindexing injection, and
\item $\psi_j \cn \ord{1} \to \ord{n}_j$ is a pointed function for each $j \in \ufs{q}$.
\end{itemize}
We define $\Aupom$ as the following composite pointed $G$-functor.
\begin{equation}\label{AF_empty_n}
\begin{tikzpicture}[vcenter]
\def\t{.7}
\draw[0cell]
(0,0) node (a1) {\sys{\A}{\ang{}}}
(a1)++(.9,0) node (a1') {\phantom{\A}}
(a1')++(0,.02) node (a1'') {\A}
(a1')++(2,0) node (a2) {\phantom{\sys{\A}{\ordtu{1}_{j \in \ufs{q}}}}}
(a2)++(0,-.03) node (a2') {\sys{\A}{\ordtu{1}_{j \in \ufs{q}}}}
(a2)++(3.5,0) node (a3) {\Aangordn}
;
\draw[1cell=.9]
(a1) edge[equal] (a1')
(a1') edge node {\iso} (a2)
(a2) edge node {\sys{\A}{(1_{\ufs{q}}, \ang{\psi_j}_{j \in \ufs{q}})}} (a3)
;
\draw[1cell=1]
(a1) [rounded corners=2pt, shorten <=0ex] |- ($(a2)+(0,\t)$)
-- node[pos=.3] {\Aupom} ($(a2)+(1,\t)$) -| (a3)
;
\end{tikzpicture}
\end{equation}
In \cref{AF_empty_n},
\begin{itemize}
\item the functor labeled $\iso$ is the pointed $G$-isomorphism in \cref{Aordtuone}, and
\item $\sys{\A}{(1_{\ufs{q}}, \ang{\psi_j}_{j \in \ufs{q}})}$ is the pointed $G$-functor in \cref{A_fangpsi} with $f = 1_{\ufs{q}}$.
\end{itemize}
\end{itemize}
This finishes the definition of the pointed $G$-functor $\Aupom$.  
\item[Strong variant]  
We define the pointed $G$-functor
\begin{equation}\label{AF_sg}
\Aangordmsg \fto{\Aupomsg} \Aangordnsg
\end{equation}
by restricting the definition of $\Aupom$ to the full subcategory of strong $\angordm$-systems.\defmark
\end{description}
\end{definition}

The following result is \cite[\namecref{EqK:A_ptfunctor} \ref*{EqK:A_ptfunctor}]{yau-eqk}.

\begin{lemma}[$J$-Theory on Objects]\label{A_ptfunctor}\index{J-theory@$J$-theory!object}
For a $\Tinf$-operad $\Op$ \pcref{as:OpA} and an $\Op$-pseudoalgebra $\A$ \pcref{def:pseudoalgebra}, the object and morphism assignments 
\[\angordn \mapsto \Aangordn \andspace \upom \mapsto \Aupom\]
in, respectively,  \cref{def:Aangordn_gcat,def:Afangpsi} define a pointed functor
\[\Gsk \fto{\Jgo\A = \Adash} \Gcatst.\]
Moreover, the strong variant 
\[\Gsk \fto{\Jgosg\A = \Asgdash} \Gcatst\]
defined in \cref{sgAordnbe,AF_sg} is also a pointed functor.
\end{lemma}

 

%% file: chap/kgo_v.tex
\section{$J$-Theory on 1-Cells}
\label{sec:jgosg_onecells}

This section recalls the 1-cell assignments of the (strong) $J$-theory 2-functors
\[\AlglaxO \fto{\Jgo} \GGCatii \andspace \AlgpspsO \fto{\Jgosg} \GGCatii\]
for a $\Tinf$-operad $\Op$ \pcref{as:OpA}.  Their object assignments, $\A \mapsto \Adash$ and $\A \mapsto \Asgdash$, are given in \cref{A_ptfunctor}.  The 1-cell assignments of $\Jgo$ and $\Jgosg$ are obtained from the constructions in \cite[Sections \ref*{EqK:sec:jemg_pos_i} and \ref*{EqK:sec:jemg_pos_i_proof}]{yau-eqk} by restricting to $G$-fixed 1-ary 1-cells.  Recall that
\begin{itemize}
\item 1-cells in $\AlglaxO$ and $\AlgpspsO$ are lax $\Op$-morphisms and $\Op$-pseudomorphisms \pcref{def:laxmorphism}, and 
\item 1-cells in $\GGCatii$ are natural transformations \cref{ggcat_mor}.
\end{itemize}  
Similar to \cref{sec:jemg,sec:jemg_morphisms}, we first consider the lax case \cref{Jgo_f} and then restrict the definitions to the strong case \cref{Jgosg_f}.

\begin{definition}[$J$-Theory on 1-Cells]\label{def:Jgo_pos_obj}\index{J-theory@$J$-theory!1-cell}
Suppose $(\A,\gaA,\phiA)$ and $(\B,\gaB,\phiB)$ are $\Op$-pseudoalgebras \pcref{def:pseudoalgebra} for a $\Tinf$-operad $\Op$ \pcref{as:OpA}, and suppose 
\[(\A,\gaA,\phiA) \fto{(f,\actf)} (\B,\gaB,\phiB)\]
is a lax $\Op$-morphism \pcref{def:laxmorphism}.  We define a natural transformation
\begin{equation}\label{Jgo_f}
\begin{tikzpicture}[vcenter]
\def\t{28}
\draw[0cell]
(0,0) node (a1) {\Gsk}
(a1)++(2,0) node (a2) {\phantom{\Gskel}}
(a2)++(.3,0) node (a2') {\Gcatst}
;
\draw[1cell=.8]
(a1) edge[bend left=\t] node {\Adash} (a2)
(a1) edge[bend right=\t] node[swap] {\Bdash} (a2)
;
\draw[2cell=.9]
node[between=a1 and a2 at .37, rotate=-90, 2label={above,\Jgo f}] {\Rightarrow}
;
\end{tikzpicture}
\end{equation}
as follows, where we use the description of a 1-cell in $\GGCatii$ given in \cref{ggcat_mor_component}.  By \cref{ggcat_mor_vstar,A_angordn}, the $\vstar$-component 
\begin{equation}\label{Jgof_vstar}
\sys{\A}{\vstar} = \bone \fto{(\Jgo f)_{\vstar} = 1} \sys{\B}{\vstar} = \bone
\end{equation}
is the identity functor on the terminal $G$-category $\bone$.  For the empty tuple $\ang{} \in \Gsk$, the $\ang{}$-component
\begin{equation}\label{Jgof_ang}
\sys{\A}{\ang{}} = \A \fto{(\Jgo f)_{\ang{}} = f} \sys{\B}{\ang{}}= \B
\end{equation}
is defined as the given pointed $G$-functor $f$ \cref{laxmorphism_basepoint}.

For a nonbasepoint object $\angordn = \ang{\ordn_j}_{j \in \ufs{q}} \in \Gsk$ \cref{Gsk_objects} of length $q > 0$, the $\angordn$-component pointed $G$-functor
\begin{equation}\label{Jgof_component}
\Aangordn \fto{(\Jgo f)_{\angordn}} \Bangordn
\end{equation}
is defined as follows.
\begin{description}
\item[Component objects]
Suppose $(a,\glu^a) \in \Aangordn$ is an $\angordn$-system \cref{nsystem}.  The $\angordn$-system 
\begin{equation}\label{Jgof_m_objects}
(\Jgo f)_{\angordn}(a,\glu^a) = (b,\glu^b) \in \Bangordn
\end{equation}
has, for each marker $\angs = \ang{s_j \subseteq \ufs{n}_j}_{j \in \ufs{q}}$, $\angs$-component object defined as
\begin{equation}\label{Jgof_m_obj_comp}
b_{\angs} = fa_{\angs} \in \B,
\end{equation}
where $a_{\angs} \in \A$ is the $\angs$-component object of $(a,\glu)$ \cref{a_angs}.
\item[Gluing]
For each object $x \in \Op(r)$ with $r \geq 0$, marker $\ang{s} = \ang{s_j \subseteq \ufs{n}_j}_{j \in \ufs{q}}$, index $k \in \ufs{q}$, and partition of $s_k$ into $r$ subsets
\[s_k = \coprod_{i \in \ufs{r}} \, s_{k,i} \subseteq \ufs{n}_k,\]
the gluing morphism \cref{gluing-morphism} of $(b,\glu^b)$ at $\big(x; \angs, k, \ang{s_{k,i}}_{i \in \ufs{r}}\big)$ is defined as the following composite in $\B$.
\begin{equation}\label{Jgof_m_gluing}
\begin{tikzpicture}[vcenter]
\def\u{-1} \def\h{4.5} \def\a{10} \def\b{.7}
\draw[0cell=.85]
(0,0) node (a1) {\gaB_r\big(x \sscs \ang{b_{\angs \compk\, s_{k,i}}}_{i \in \ufs{r}}\big)}
(a1)++(\h,0) node (a2) {b_{\angs}}
(a1)++(0,\u) node (b1) {\gaB_r\big(x \sscs \ang{f a_{\angs \compk\, s_{k,i}}}_{i \in \ufs{r}}\big)}
(a2)++(0,\u) node (b2) {fa_{\angs}}
(b1)++(\h/2,-1.2) node (c) {f \gaA_r \big(x \sscs \ang{a_{\angs \compk\, s_{k,i}}}_{i \in \ufs{r}} \big)}
;
\draw[1cell=.9]
(a1) edge[equal,shorten <=-.3ex,shorten >=-.5ex] (b1)
(a2) edge[equal] (b2)
(a1) edge node {\glu^b_{x;\, \angs,\, k, \ang{s_{k,i}}_{i \in \ufs{r}}}} (a2)
;
\draw[1cell=.9]
(b1) [rounded corners=2pt, shorten <=-.2ex] |- node[pos=.2,swap] {\actf_r} (c);
\draw[1cell=.9]
(c) [rounded corners=2pt, shorten <=-0ex] -| node[pos=.8] {f \glu^a_{x;\, \angs,\, k, \ang{s_{k,i}}_{i \in \ufs{r}}}} (b2);
\end{tikzpicture}
\end{equation}
The lower-left arrow $\actf_r$ is a component of the action constraint of $f$ \cref{actf_component}.  The lower-right arrow is the image under $f$ of the gluing morphism of $(a,\glu^a)$ at $\big(x; \angs, k, \ang{s_{k,i}}_{i \in \ufs{r}}\big)$.  This finishes the definition of the object assignment of $(\Jgo f)_{\angordn}$ \cref{Jgof_component}.
\item[Morphisms]
Suppose
\[(a,\glu^a) \fto{\theta} (a',\glu^{a'})\]
is a morphism of $\angordn$-systems in $\A$ \pcref{def:nsystem_morphism}.  The morphism of $\angordn$-systems in $\B$
\begin{equation}\label{Jgof_m_theta}
(\Jgo f)_{\angordn} (a,\glu^a) \fto{(\Jgo f)_{\angordn} \theta} 
(\Jgo f)_{\angordn} (a',\glu^{a'}) 
\end{equation}
has, for each marker $\ang{s} = \ang{s_j \subseteq \ufs{n}_j}_{j \in \ufs{q}}$, $\angs$-component morphism \cref{theta_angs} given by the image under $f$ of the $\angs$-component morphism of $\theta$:
\begin{equation}\label{Jgof_m_theta_comp}
f a_{\angs} \fto{((\Jgo f)_{\angordn} \theta)_{\angs} = f\theta_{\angs}} f a'_{\angs}.
\end{equation}
\end{description}
This finishes the definitions of the $\angordn$-component pointed $G$-functor $(\Jgo f)_{\angordn}$ \cref{Jgof_component} and the natural transformation $\Jgo f$ \cref{Jgo_f}.

\parhead{Strong variant}.  Suppose $(f,\actf)$ is an $\Op$-pseudomorphism.  We define a natural transformation
\begin{equation}\label{Jgosg_f}
\begin{tikzpicture}[vcenter]
\def\t{30}
\draw[0cell]
(0,0) node (a1) {\Gsk}
(a1)++(2,0) node (a2) {\phantom{\Gskel}}
(a2)++(.3,0) node (a2') {\Gcatst}
;
\draw[1cell=.8]
(a1) edge[bend left=\t] node {\Asgdash} (a2)
(a1) edge[bend right=\t] node[swap] {\Bsgdash} (a2)
;
\draw[2cell=.9]
node[between=a1 and a2 at .35, rotate=-90, 2label={above,\Jgosg f}] {\Rightarrow}
;
\end{tikzpicture}
\end{equation}
by restricting the definition of $\Jgo f$ to the full subcategories of strong systems.  For each object $\angordn \in \Gsk$, the $\angordn$-component pointed $G$-functor
\begin{equation}\label{Jgosgf_component}
\Asgangordn \fto{(\Jgosg f)_{\angordn}} \Bsgangordn
\end{equation}
is well defined because, in the diagram \cref{Jgof_m_gluing} defining the gluing morphisms, the action constraint $\actf_r$ is an isomorphism, since $f$ is an $\Op$-pseudomorphism.  Moreover, the gluing morphism $\glu^a_{\cdots}$ is an isomorphism because $(a,\glu^a)$ is a strong $\angordn$-system in $\A$.  Thus, $(\Jgosg f)_{\angordn}(a,\glu^a)$ is a strong $\angordn$-system in $\B$.
\end{definition}

\section{$J$-Theory 2-Functors}
\label{sec:jgosg_twocells}

This section defines the 2-cell assignments and finishes the construction of the (strong) $J$-theory 2-functors
\[\AlglaxO \fto{\Jgo} \GGCatii \andspace \AlgpspsO \fto{\Jgosg} \GGCatii\]
for a $\Tinf$-operad $\Op$ \pcref{as:OpA}.  See \cref{thm:Jgo_twofunctor}.  Their object and 1-cell assignments are given in \cref{A_ptfunctor,def:Jgo_pos_obj}.  The 2-cell assignments of $\Jgo$ and $\Jgosg$ are obtained from the constructions in \cite[\namecref{EqK:sec:jemg_pos_ii} \ref*{EqK:sec:jemg_pos_ii}]{yau-eqk} by restricting to $G$-fixed 1-ary 2-cells.  Recall that
\begin{itemize}
\item 2-cells in $\AlglaxO$ and $\AlgpspsO$ are $\Op$-transformations \pcref{def:algtwocells} and 
\item 2-cells in $\GGCatii$ are modifications \cref{ggcat_inthom_Theta}. 
\end{itemize} 
Similar to \cref{sec:jemg,sec:jemg_morphisms,sec:jgosg_onecells}, we first consider the lax case \cref{Jgotheta} and then restrict the definitions to the strong case \cref{Jgosgtheta}.

\begin{definition}[$J$-Theory on 2-Cells]\label{def:Jgo_pos_mor}\index{J-theory@$J$-theory!2-cell}
Suppose
\[(f,\actf) \fto{\tha} (h,\acth)\]
is an $\Op$-transformation between lax $\Op$-morphisms
\[(f,\actf), (h,\acth) \cn \big(\A,\gaA,\phiA\big) \to \big(\B,\gaB,\phiB\big)\]
between $\Op$-pseudoalgebras $(\A,\gaA,\phiA)$ and $(\B,\gaB,\phiB)$ \pcref{def:pseudoalgebra,def:laxmorphism,def:algtwocells} for a $\Tinf$-operad $\Op$ \pcref{as:OpA}.  We define a modification
\begin{equation}\label{Jgotheta}
\begin{tikzpicture}[baseline={(a.base)}]
\def\h{4} \def\u{.7} \def\t{22}
\draw[0cell]
(0,0) node (a1) {\Gsk}
(a1)++(\h,0) node (a2) {\phantom{\Gskel}}
(a2)++(.3,0) node (a2') {\Gcatst}
(a1)++(.5*\h,0) node (a) {}
;
\draw[1cell=.8]
(a1) edge[bend left=\t] node {\Adash} (a2)
(a1) edge[bend right=\t] node[swap] {\Bdash} (a2)
;
\draw[2cell=.9]
node[between=a1 and a at .7, rotate=-90, 2labelmed={below,\Jgo f}] {\Rightarrow}
node[between=a and a2 at .25, rotate=-90, 2label={above,\Jgo h}] {\Rightarrow}
;
\draw[2cell]
node[between=a1 and a2 at .5, shift={(0,-.2)}, 2labelalt={above,\Jgo \tha}] {\Rrightarrow}
;
\end{tikzpicture}
\end{equation}
as follows.  Using \cref{Jgof_vstar}, the $\vstar$-component
\begin{equation}\label{Jgotheta_vstar}
1_\bone \fto{(\Jgo\tha)_\vstar = 1} 1_\bone
\end{equation}
is the identity natural transformation on the identity functor on the terminal $G$-category $\bone$.  Using \cref{expl:Otr_pointed,Jgof_ang}, the $\ang{}$-component 
\begin{equation}\label{Jgotheta_ang}
(\Jgo f)_{\ang{}} = f \fto{(\Jgo\tha)_{\ang{}} = \tha} (\Jgo h)_{\ang{}} = h
\end{equation}
is defined as the given pointed $G$-natural transformation $\tha$.

Suppose $\angordn = \ang{\ordn_j}_{j \in \ufs{q}} \in \Gsk$ \cref{Gsk_objects} is a nonbasepoint object of length $q > 0$.   Using \cref{Jgof_component,Jgof_m_objects,Jgof_m_obj_comp}, the $\angordn$-component pointed $G$-natural transformation
\begin{equation}\label{Jgotheta_m}

\end{equation}
sends an $\angordn$-system $(a,\glu^a) \in \Aangordn$ \cref{nsystem} to the morphism of $\angordn$-systems in $\B$ \pcref{def:nsystem_morphism}
\[(\Jgo f)_{\angordn}(a,\glu^a) \fto{(\Jgo \tha)_{\angordn, (a,\glu^a)}} (\Jgo h)_{\angordn}(a,\glu^a).\]
For each marker $\angs = \ang{s_j \subseteq \ufsn_j}_{j \in \ufsq}$, its $\angs$-component morphism is given by the $a_{\angs}$-component of $\tha$, as displayed in \cref{Jgotheta_m_angs}.
\begin{equation}\label{Jgotheta_m_angs}
%
\end{equation}
This finishes the definitions of the $\angordn$-component pointed $G$-natural transformation $(\Jgo\tha)_{\angordn}$ \cref{Jgotheta_m} and the modification $\Jgo\tha$ \cref{Jgotheta}.

\parhead{Strong variant}.  For $\Op$-pseudomorphisms $(f,\actf)$ and $(h,\acth)$, the modification
\begin{equation}\label{Jgosgtheta}
%
\end{equation}
is defined by reusing the definitions \cref{Jgotheta,Jgotheta_vstar,Jgotheta_ang,Jgotheta_m,Jgotheta_m_angs} and restricting to the full subcategories of strong systems in \cref{Jgotheta_m}.
\end{definition}

The following result is obtained from the $J$-theory $\Gcat$-multifunctors in \cite[\namecref{EqK:thm:Jgo_multifunctor} \ref*{EqK:thm:Jgo_multifunctor}]{yau-eqk} by restricting to $G$-fixed 1-ary 1-cells and 2-cells.

\begin{theorem}\label{thm:Jgo_twofunctor}\index{J-theory@$J$-theory}\index{J-theory@$J$-theory!strong}
For each $\Tinf$-operad $\Op$ \pcref{as:OpA}, the object, 1-cell, and 2-cell assignments in \cref{A_ptfunctor,def:Jgo_pos_obj,def:Jgo_pos_mor} define 2-functors
\[\begin{split}
\AlglaxO & \fto{\Jgo} \GGCatii \andspace\\
\AlgpspsO & \fto{\Jgosg} \GGCatii
\end{split}\]
between the 2-categories in \cref{oalgps_twocat,def:GGCat}.
\end{theorem}

We call $\Jgo$ the \emph{$J$-theory 2-functor} and $\Jgosg$ the \emph{strong $J$-theory 2-functor}.

 

%% file: chap/kgo_vi.tex
\section{$\Gskg$-Spaces from $\Gskg$-Categories}
\label{sec:ggtop}

This section reviews the functor
\[\GGCatii \fto{\clast} \GGTopii\]
from the category $\GGCatii$ of $\Gskg$-categories to the category $\GGTopii$ of $\Gskg$-spaces, induced by the classifying space functor $\cla$.  The functor $\clast$ underlies the symmetric monoidal functor
\[\GGCat \fto{\clast} \GGTop\]
in \cite[\namecref{EqK:thm:ggcat_ggtop} \ref*{EqK:thm:ggcat_ggtop} (\ref*{EqK:ggcat_ggtop_i})]{yau-eqk}.

\secoutline
\begin{itemize}
\item \cref{def:Gtop,def:weakG_top,def:gtopst} review (pointed) $G$-spaces.
\item \cref{def:ggtop_smc} recalls the category $\GGTopii$ of $\Gskg$-spaces, with further elaboration given in \cref{expl:ggtop_obj}.
\item The functor $\clast$ is recorded in \cref{thm:ggcat_ggtop} and discussed further in \cref{expl:clast_functor}.
\end{itemize}

\subsection*{Equivariant Spaces}

\cref{def:Gtop} is the topological analogue of \cref{def:GCat,def:Catg}.

\begin{definition}[$G$-Spaces]\label{def:Gtop}
Denote by $\Top$ the complete and cocomplete category of
\begin{itemize}
\item \index{space}\emph{spaces}, meaning compactly generated weak Hausdorff spaces, and
\item \index{morphism}\emph{morphisms}, meaning continuous maps between spaces.
\end{itemize}
Suppose $G$ is a group, which is also regarded as a category with one object $*$ and morphism set $G$, with composition and identity given by the group multiplication and the group unit.
\begin{itemize}
\item A \emph{$G$-space}\index{G-space@$G$-space} is a functor $G \to \Top$.  In other words, a $G$-space is a space $X$ equipped with a $g$-action homeomorphism $g \cn X \fiso X$ for each $g \in G$ that satisfies the following two conditions.
\begin{enumerate}
\item For the identity element $e \in G$, the $e$-action is equal to $1_X$.
\item For $g,h \in G$, the $hg$-action is equal to the composite $h \circ g$.
\end{enumerate}
\item A \emph{$G$-morphism}\index{G-morphism@$G$-morphism} $f \cn X \to Y$ between $G$-spaces is a natural transformation between functors $G \to \Top$.  In other words, $f$ is a morphism of spaces that is \index{G-equivariant@$G$-equivariant}\emph{$G$-equivariant}:
\[f(gx) = g(fx) \forspace (g,x) \in G \times X.\]
\item $\Gtop$\label{not:Gtop} denotes the complete and cocomplete category of $G$-spaces and $G$-morphisms, with composition and identities defined in $\Top$.
\item For $G$-spaces $X$ and $Y$, $\Topg(X,Y)$\label{not:Topg} denotes the $G$-space of all morphisms $X \to Y$, with the compact-open topology and conjugation $G$-action.  For $g \in G$ and a morphism $h \cn X \to Y$, the \emph{conjugation $G$-action} $g \cdot h$ is defined as the composite
\begin{equation}\label{ginv_h_g}
X \fto{\ginv} X \fto{h} Y \fto{g} Y.
\end{equation}
\end{itemize}
The quadruple
\begin{equation}\label{gtop_smclosed}
(\Gtop, \times, *, \Topg)
\end{equation} 
is a Cartesian closed category.
\end{definition}


\begin{definition}[Weak $G$-Equivalences]\label{def:weakG_top}\
\begin{itemize}
\item For a $G$-space $X$, a \index{G-fixed point@$G$-fixed point}\emph{$G$-fixed point} is a point $x \in X$ such that $gx = x$ for all $g \in G$.  The \emph{$G$-fixed point space} $X^G \subseteq X$ is the subspace consisting of all the $G$-fixed points.
\item For a $G$-morphism $f \cn X \to Y$ between $G$-spaces, the restriction to $G$-fixed point spaces is denoted by $f^G \cn X^G \to Y^G$.
\item A $G$-morphism $f \cn X \to Y$ is called a \index{weak G-equivalence@weak $G$-equivalence}\emph{weak $G$-equivalence} if, for each subgroup $H \subseteq G$, the morphism $f^H \cn X^H \to Y^H$ is a weak homotopy equivalence of spaces.  A weak $G$-equivalence is also denoted by \label{not:eqg}$\eqg$.\defmark
\end{itemize}
\end{definition}

\cref{def:gtopst} is the topological analogue of \cref{def:gcatst,expl:Gcatst}.

\begin{definition}[Pointed $G$-Spaces]\label{def:gtopst}
The complete and cocomplete Cartesian closed category \cref{gtop_smclosed}
\[(\Gtop, \times, *, \Topg)\]
yields the complete and cocomplete symmetric monoidal closed category
\begin{equation}\label{Gtopst_smc}
(\Gtopst, \sma, \stplus, \Topgst).
\end{equation}
If $G$ is the trivial group, then $\Gtopst$ is denoted by $\Topst$, whose objects and morphisms are called, respectively, \index{pointed space}\emph{pointed spaces} and \index{pointed morphism}\emph{pointed morphisms}.
\begin{itemize}
\item An object in $\Gtopst$ is called a \index{pointed G-space@pointed $G$-space}\emph{pointed $G$-space}.  It consists of a $G$-space $X$ \pcref{def:Gtop} and a distinguished $G$-fixed \index{basepoint}\emph{basepoint}.
\item A morphism in $\Gtopst$ is called a \index{pointed G-morphism@pointed $G$-morphism}\emph{pointed $G$-morphism}.  It consists of a $G$-morphism between $G$-spaces that preserves the basepoints. 
\item For pointed $G$-spaces $X$ and $Y$, $\Gtopst(X,Y)$ also denotes the pointed space of pointed $G$-morphisms $X \to Y$, with the compact-open topology and the basepoint given by the constant morphism at the basepoint of $Y$.
\item Composition and identities in $\Gtopst$ are defined in $\Gtop$.
\item $\sma$ is the \index{smash product}smash product.
\item The \index{smash unit}\emph{smash unit} $\stplus$\label{not:stplus} consists of two points, with $G$ acting trivially.
\item For pointed $G$-spaces $X$ and $Y$, the internal hom $\Topgst(X,Y)$\label{not:Topgst} is the $G$-subspace of $\Topg(X,Y)$ consisting of pointed morphisms $X \to Y$.
\begin{itemize}
\item The basepoint of $\Topgst(X,Y)$ is the constant morphism at the basepoint of $Y$.
\item $G$ acts on $\Topgst(X,Y)$ by conjugation \cref{ginv_h_g}.
\end{itemize}
The $G$-fixed point subspace of the pointed $G$-space $\Topgst(X,Y)$ is the pointed space
\[\Gtopst(X,Y) = \Topgst(X,Y)^G\]
of pointed $G$-morphisms.
\end{itemize}
The notation 
\begin{equation}\label{topgst_gtopst_enr}
\Topgst
\end{equation}
also denotes the following categories with pointed $G$-spaces as objects.
\begin{enumerate}
\item The pointed $G$-category with pointed morphisms, the one-point space $*$ as the basepoint, and $G$ acting trivially on objects and by conjugation on morphisms.
\item The $\Gtopst$-category (in the sense of \cref{def:enriched-category}) with hom $\Gtopst$-object given by the pointed $G$-space $\Topgst(X,Y)$ for any pair $(X,Y)$ of pointed $G$-spaces, and composition given by that of pointed morphisms.\defmark
\end{enumerate}
\end{definition}

\subsection*{$\Gskg$-Spaces}

\cref{def:ggtop_smc} is the topological analogue of \cref{def:GGCat}.  The category $\GGTopii$ is obtained from the symmetric monoidal $\Gtop$-category $\GGTop$ in \cite[\ref*{EqK:def:ggtop_smc}]{yau-eqk} by passing to the $G$-fixed subspaces in its $\Gtop$-enrichment.

\begin{definition}[$\Gskg$-Spaces]\label{def:ggtop_smc}
For a group $G$, the category $\GGTopii$ is defined as follows.
\begin{description}
\item[Objects] An object in $\GGTopii$, called a \index{G-G-space@$\Gskg$-space}\emph{$\Gskg$-space}, is a pointed functor
\begin{equation}\label{ggtop_obj}
(\Gsk, \vstar) \fto{f} (\Gtopst,*).
\end{equation}
\item[Morphisms] A morphism $\theta \cn f \to f'$ in $\GGTopii$ is a natural transformation as follows.
\begin{equation}\label{ggtop_mor}
\begin{tikzpicture}[vcenter]
\def\t{28}
\draw[0cell]
(0,0) node (a1) {\Gsk}
(a1)++(1.8,0) node (a2) {\phantom{\Gskel}}
(a2)++(.33,0) node (a2') {\Gtopst}
;
\draw[1cell=.9]
(a1) edge[bend left=\t] node {f} (a2)
(a1) edge[bend right=\t] node[swap] {f'} (a2)
;
\draw[2cell]
node[between=a1 and a2 at .43, rotate=-90, 2label={above,\theta}] {\Rightarrow}
;
\end{tikzpicture}
\end{equation}
Identities and composition are defined componentwise in $\Gtopst$.
\end{description}
The category $\FGTop$\label{not:FGTop} is defined in the same way by replacing the pointed category $(\Gsk,\vstar)$ with $(\Fsk,\ordz)$ \pcref{def:Fsk}.  An object in $\FGTop$ is called an \index{F-G-space@$\Fskg$-space}\emph{$\Fskg$-space}.
\end{definition}

\begin{explanation}[Unpacking $\GGTopii$]\label{expl:ggtop_obj}
A $\Gskg$-space $f \cn \Gsk \to \Gtopst$ consists of the following data.
\begin{itemize}
\item $f$ sends each object $\angordm \in \Gsk$ \cref{Gsk_objects} to a pointed $G$-space $f\angordm$ such that $f\vstar = *$.  The \index{canonical basepoint}\emph{canonical basepoint} of $f\angordm$ is given by the $G$-morphism
\begin{equation}\label{fm_pointed_ggtop}
f(\vstar \to \angordm) \cn f\vstar = * \to f\angordm,
\end{equation}
where $\vstar \to \angordm$ is the unique morphism in $\Gsk(\vstar,\angordm)$. 
\item $f$ sends each morphism $\upom \cn \angordm \to \angordn$ in $\Gsk$ \cref{Gsk_morphisms} to a pointed $G$-morphism between pointed $G$-spaces
\begin{equation}\label{f_upom_ggtop}
f\angordm \fto{f\upom} f\angordn
\end{equation}
such that $f$ preserves identities and composition of morphisms.
\end{itemize}
A morphism $\theta \cn f \to f'$ in $\GGTopii$ consists of, for each object $\angordm \in \Gsk$, an $\angordm$-component pointed $G$-morphism between pointed $G$-spaces
\begin{equation}\label{ggtop_mor_component}
f\angordm \fto{\theta_{\angordm}} f'\angordm
\end{equation}
such that, for each morphism $\upom \cn \angordm \to \angordn$ in $\Gsk$, the following naturality diagram of pointed $G$-morphisms commutes.
\begin{equation}\label{ggtop_mor_naturality}
\begin{tikzpicture}[vcenter]
\def\v{-1.4}
\draw[0cell]
(0,0) node (a11) {f\angordm}
(a11)++(2.5,0) node (a12) {f'\angordm}
(a11)++(0,\v) node (a21) {f\angordn}
(a12)++(0,\v) node (a22) {f'\angordn}
;
\draw[1cell=.9]
(a11) edge node {\theta_{\angordm}} (a12)
(a12) edge node {f'\upom} (a22)
(a11) edge node[swap] {f\upom} (a21)
(a21) edge node {\theta_{\angordn}} (a22)
;
\end{tikzpicture}
\end{equation}
Identity morphisms and composition in $\GGTopii$ are defined componentwise in $\Gtopst$ using the components in \cref{ggtop_mor_component}.  As a natural transformation, $\theta$ is automatically pointed, meaning that $\theta_\vstar = 1_*$.  The category $\FGTop$ admits the same description as $\GGTopii$ with $(\Gsk,\vstar)$ replaced by $(\Fsk,\ordz)$.
\end{explanation}

\subsection*{Classifying Space}

The \index{classifying space}classifying space functor $\cla$  \cite[\namecref{EqK:ch:nerve} \ref*{EqK:ch:nerve}]{yau-eqk} is the composite
\begin{equation}\label{classifying_space}
\begin{tikzpicture}[baseline={(a1.base)}]
\def\h{2.2}
\draw[0cell]
(0,0) node (a1) {\Cat}
(a1)++(\h,0) node (a2) {\sset}
(a2)++(\h,0) node (a3) {\Top}
;
\draw[1cell=.9]
(a1) edge node {\Ner} (a2)
(a2) edge node {\Rea} (a3)
;
\draw[1cell=.9]
(a1) [rounded corners=2pt] |- ($(a2)+(-1,.6)$) -- node {\cla} ($(a2)+(1,.6)$) -| (a3)
 ;
\end{tikzpicture}
\end{equation}
of the nerve $\Ner$ and realization $\Rea$ \cref{x_dot}, where $\sset$ is the category of simplicial sets (\cref{def:simp_obj} \eqref{def:simp_obj_vi}).  It preserves finite products and yields a strong symmetric monoidal functor 
\begin{equation}\label{cla_gcat_gtop}
(\Gcat,\times,\boldone) \fto{(\cla,\clatwo,\clazero)} (\Gtop, \times, *).
\end{equation}
The following observation records the underlying functor in \cite[\namecref{EqK:thm:ggcat_ggtop} \ref*{EqK:thm:ggcat_ggtop} (\ref*{EqK:ggcat_ggtop_i})]{yau-eqk}.

\begin{theorem}\label{thm:ggcat_ggtop}
For a group $G$, composing and whiskering with the classifying space functor $\cla$ in \cref{cla_gcat_gtop} induce a functor\index{classifying space}\index{G-G-category@$\Gskg$-category}
\begin{equation}\label{clast}
\GGCatii \fto{\clast} \GGTopii
\end{equation}
between the categories in \cref{def:GGCat,def:ggtop_smc}.
\end{theorem}

\begin{explanation}[Unpacking $\clast$]\label{expl:clast_functor}
The functor $\clast$ in \cref{clast} sends a pointed functor $f \cn \Gsk \to \Gcatst$ to the composite pointed functor
\begin{equation}\label{clast_obj}
\begin{tikzpicture}[baseline={(a1.base)}]
\def\b{.6}
\draw[0cell]
(0,0) node (a1) {(\Gsk,\vstar)}
(a1)++(2.4,0) node (a2) {(\Gcatst,\boldone)}
(a2)++(2.9,0) node (a3) {(\Gtopst,*).}
;
\draw[1cell=.9]
(a1) edge node {f} (a2)
(a2) edge node {\cla} (a3)
;
\draw[1cell=.9]
(a1) [rounded corners=2pt, shorten <=0ex] |- ($(a2)+(-1,\b)$)
-- node {\clast f} ($(a2)+(1.25,\b)$) -| (a3)
;
\end{tikzpicture}
\end{equation}
This is well defined because $\cla\boldone = *$.  The functor $\clast$ sends a natural transformation $\tha$ to the whiskering
\begin{equation}\label{clast_mor}
\begin{tikzpicture}[baseline={(a1.base)}]
\def\t{28}
\draw[0cell]
(0,0) node (a1) {\Gsk}
(a1)++(1.8,0) node (a2) {\phantom{\Gsk}}
(a2)++(.3,0) node (a2') {\Gcatst}
(a2')++(.2,0) node (a2'') {\phantom{Gsk}}
(a2'')++(2.3,0) node (a3) {\Gtopst.}
;
\draw[1cell=.9]
(a1) edge[bend left=\t] node {f} (a2)
(a1) edge[bend right=\t] node[swap] {f'} (a2)
(a2'') edge node {\cla} (a3)
;
\draw[2cell]
node[between=a1 and a2 at .43, rotate=-90, 2label={above,\theta}] {\Rightarrow}
;
\end{tikzpicture}
\end{equation}
For each object $\angordm \in \Gsk$ \cref{Gsk_objects}, the $\angordm$-component of $\clast\theta$ is the pointed $G$-morphism between pointed $G$-spaces
\[\cla f\angordm \fto{\cla\theta_{\angordm}} \cla f'\angordm\]
obtained from $\theta_{\angordm}$ \cref{ggcat_mor_component} by applying $\cla$.
\end{explanation}

 

%% file: chap/kgo_vii.tex
\section{Orthogonal $G$-Spectra}
\label{sec:spectra}

For a compact Lie group $G$, this section reviews the categories of orthogonal $G$-spectra in \cite{mandell_may} by adapting \cite[\namecref{EqK:ch:spectra} \ref*{EqK:ch:spectra}]{yau-eqk}.

\secoutline
\begin{itemize}
\item \cref{def:g_universe,def:indexing_gspace,def:IU_spaces} recall the symmetric monoidal category $\IU$ for a complete $G$-universe $\univ$.
\item \cref{def:iu_space,def:iu_morphism,def:iut_gtop_enr,def:giut_top_enr} review the categories $\IUT$ and $\GIUT$ of $\IU$-spaces.
\item \cref{def:g_sphere,def:gsp_module,def:gsp_morphism} review the categories $\GSp$ and $\Gspec$ of orthogonal $G$-spectra.
\end{itemize}

\subsection*{The Category $\IU$}

Recall the notions of (pointed) $G$-spaces in \cref{def:Gtop,def:gtopst}.

\begin{definition}[Complete $G$-Universe]\label{def:g_universe}
Suppose $G$ is a compact Lie group.
\begin{enumerate}
\item A \emph{real $G$-inner product space}\index{G-inner product space@$G$-inner product space} is a pair $(X,\mu)$ consisting of
\begin{itemize}
\item a real inner product space $X$ and
\item a $G$-action $\mu \cn G \times X \to X$
\end{itemize}
such that, for each $g \in G$, $\mu(g,-)$ is a linear isometric isomorphism. 
\item A \emph{complete $G$-universe}\index{complete $G$-universe} $\univ$\label{not:univ} is a real $G$-inner product space that contains countably many copies of each irreducible $G$-representation.
\end{enumerate}
Fix a complete $G$-universe $\univ$ for the rest of this chapter.  One possible choice of $\univ$ is the direct sum of countably many copies of the regular representation of $G$.
\end{definition}

\begin{definition}[Indexing $G$-Spaces]\label{def:indexing_gspace}\
\begin{itemize}
\item An \emph{indexing $G$-space}\index{indexing G-space@indexing $G$-space} in $\univ$ is a finite dimensional real $G$-inner product subspace $V \subset \univ$.
\item $\IU$\label{not:IU} denotes the collection of all real $G$-inner product spaces that are isomorphic to some indexing $G$-spaces in $\univ$ via $G$-equivariant linear isometric isomorphisms.\defmark
\end{itemize}
\end{definition}

Basic concepts of enriched category theory are reviewed in \cref{sec:enrichedcat}.

\begin{convention}[Disjoint Basepoint]\label{conv:disjoint_gbasept}
For each $G$-space $X$, denote by $X_\splus$ the pointed $G$-space obtained from $X$ by adjoining a disjoint $G$-fixed basepoint\dindex{disjoint}{basepoint} $*$.  Applying this procedure to each hom $G$-space, a $\Gtop$-category becomes a $\Gtopst$-category.
\end{convention}

\begin{definition}[Indexing Category $\IU$]\label{def:IU_spaces}
Extend the collection $\IU$ \pcref{def:indexing_gspace} to a category whose morphisms are linear isometric isomorphisms between real $G$-inner product spaces.  Note that morphisms in $\IU$ are \emph{not} required to be $G$-equivariant.  Furthermore, $\IU$ is equipped with the following structures.
\begin{enumerate}
\item\label{def:iu_spaces_i} $(\IU,\oplus,0,\xi)$ is a symmetric monoidal category.
\begin{itemize}
\item The monoidal product $\oplus$ is the direct sum for real $G$-inner product spaces and linear isometric isomorphisms.
\item The monoidal unit is the one-point space 0.
\item The braiding 
\[V \oplus W \fto[\iso]{\xi_{V,W}} W \oplus V\]
permutes the two arguments.
\end{itemize}
\item\label{def:iu_spaces_ii} For objects $V,W \in \IU$, the set $\IU(V,W)$ of linear isometric isomorphisms is topologized as a $G$-subspace of $\Topg(V,W)$, so the group $G$ acts on $\IU(V,W)$ by conjugation \cref{ginv_h_g}.  With these hom $G$-spaces, $\IU$ becomes a $\Gtop$-category.
\item\label{def:iu_spaces_iii} By \cref{conv:disjoint_gbasept}, $\IU$ is also regarded as a $\Gtopst$-category.
\item\label{def:iu_spaces_iv} Denote by $\IUsk$ the small full subcategory of $\IU$ whose objects are the indexing $G$-spaces in $\univ$.  The inclusion $\IUsk \to \IU$ is an equivalence of categories, so $\IUsk$ is a small skeleton of $\IU$. 
\item\label{def:iu_spaces_v} The small skeleton $\IUsk$ inherits a symmetric monoidal structure from $(\IU,\oplus)$ such that the inclusion functor 
\[(\IUsk,\oplus) \to (\IU,\oplus)\]
is strong symmetric monoidal.
\item\label{def:iu_spaces_vi} Using the Axiom of Choice, for each object $V \in \IU$, we choose a $G$-linear isometric isomorphism
\begin{equation}\label{VV'}
V \fto[\iso]{\upphi_V} V' \withspace V' \in \IUsk
\end{equation}
such that 
\begin{equation}\label{upphi_one}
\upphi_V = 1_V \ifspace V \in \IUsk.
\end{equation}
\end{enumerate}
This finishes the definition.
\end{definition}

\subsection*{$\IU$-Spaces}

Denote by $(\pSet, \sma, \stplus)$\label{not:pSet} the symmetric monoidal closed category of pointed sets and pointed functions, with the monoidal product given by the smash product $\sma$.  The underlying pointed set functor
\begin{equation}\label{gtopst_pset}
\Gtopst \fto{\und} \pSet
\end{equation}
is strict symmetric monoidal and faithful.

\begin{definition}\label{def:iu_space}
An \emph{$\IU$-space}\index{IU-space@$\IU$-space} is a $\Gtopst$-functor 
\[\IU \fto{X} \Topgst\]
from the $\Gtopst$-category $\IU$ \pcref{def:IU_spaces} to the $\Gtopst$-category $\Topgst$ \cref{topgst_gtopst_enr}.  Moreover, the $\pSet$-functor obtained from an $\IU$-space $X$ by changing enrichment along the functor $\und$ \cref{gtopst_pset} is called the \emph{underlying functor}\index{underlying functor} of $X$ and is denoted by the same notation.
\end{definition}

\begin{explanation}[Unraveling $\IU$-Spaces]\label{expl:iu_space}
An $\IU$-space $X \cn \IU \to \Topgst$ is determined by
\begin{itemize}
\item a pointed $G$-space $X_V$ for each object $V \in \IU$ and
\item a component pointed $G$-morphism between pointed $G$-spaces
\begin{equation}\label{iu_space_comp_mor}
\IU(V,W)_\splus \fto{X} \Topgst(X_V, X_W)
\end{equation}
for each pair $(V,W)$ of objects in $\IU$
\end{itemize}
such that the following diagrams in $\Gtopst$ commute for objects $U,V,W \in \IU$, where $\mcomp$ denotes composition.
\begin{equation}\label{iu_space_axioms}
\begin{tikzpicture}[vcenter]
\def\v{-1.4}
\draw[0cell=.8]
(0,0) node (a11) {\IU(V,W)_\splus \sma \IU(U,V)_\splus}
(a11)++(3.7,0) node (a12) {\IU(U,W)_\splus}
(a11)++(0,\v) node (a21) {\Topgst(X_V,X_W) \sma \Topgst(X_U,X_V)}
(a12)++(0,\v) node (a22) {\Topgst(X_U,X_W)}
(a12)++(1.7,0) node (b11) {\stplus}
(b11)++(1.6,0) node (b12) {\IU(V,V)_\splus}
(b11)++(.9,\v) node (b22) {\Topgst(X_V,X_V)}
;
\draw[1cell=.8]
(a11) edge node {\mcomp} (a12)
(a12) edge node {X} (a22)
(a11) edge[transform canvas={xshift=1.2em}] node[swap] {X \sma X} (a21)
(a21) edge node {\mcomp} (a22)
(b11) edge node {1_V} (b12)
(b12) edge node[pos=.3] {X} (b22)
(b11) edge node[swap,pos=.3] {1_{X_V}} (b22)
;
\end{tikzpicture}
\end{equation}
The pointed $G$-space $X_V$ is also denoted by $X(V)$.
\begin{description}
\item[Components] 
The component pointed $G$-morphism in \cref{iu_space_comp_mor} is equivalent to a $G$-morphism between $G$-spaces
\begin{equation}\label{iu_space_comp'}
\IU(V,W) \fto{X} \Topgst(X_V, X_W),
\end{equation}
which sends each linear isometric isomorphism $f \cn V \fiso W$ in $\IU$ to a pointed homeomorphism between pointed $G$-spaces
\begin{equation}\label{iu_space_xf}
X_V \fto[\iso]{X_f} X_W.
\end{equation}
Note that $X_f$ is \emph{not} required to be $G$-equivariant.
\item[Equivariance] 
The $G$-equivariance of \cref{iu_space_comp_mor} means that, for each element $g \in G$ and each linear isometric isomorphism $f \cn V \fiso W$ in $\IU$, the following diagram of pointed homeomorphisms commutes.
\begin{equation}\label{x_gfginv}
\begin{tikzpicture}[vcenter]
\def\v{-1.4}
\draw[0cell]
(0,0) node (a11) {X_V}
(a11)++(2.5,0) node (a12) {X_W}
(a11)++(0,\v) node (a21) {X_V}
(a12)++(0,\v) node (a22) {X_W}
;
\draw[1cell=.9]
(a11) edge node {X_{gf\ginv}} (a12)
(a11) edge node[swap] {\ginv} (a21)
(a21) edge node {X_f} (a22)
(a22) edge node[swap] {g} (a12)
;
\end{tikzpicture}
\end{equation}
In particular, if $f \cn V \fiso W$ is $G$-equivariant, then so is $X_f = X_{gf\ginv}$. 
\item[Underlying functor] 
The $\pSet$-functor $X$ sends each object $V \in \IU$ to the pointed $G$-space $X_V$ and each linear isometric isomorphism $f \cn V \fiso W$ in $\IU$ to the pointed homeomorphism $X_f$.\defmark
\end{description}
\end{explanation}

\subsection*{$\IU$-Morphisms}

\begin{definition}\label{def:iu_morphism}\
\begin{enumerate}
\item\label{def:iu_morphism_i}
For $\IU$-spaces $X,Y \cn \IU \to \Topgst$ \pcref{def:iu_space}, an \emph{$\IU$-morphism}\index{IU-morphism@$\IU$-morphism} $\theta \cn X \to Y$ is a $\pSet$-natural transformation
\begin{equation}\label{iu_mor}
\begin{tikzpicture}[vcenter]
\def\t{28}
\draw[0cell]
(0,0) node (a1) {\phantom{A}}
(a1)++(1.8,0) node (a2) {\phantom{A}}
(a1)++(-.1,0) node (a1') {\IU}
(a2)++(.21,-.04) node (a2') {\Topgst}
;
\draw[1cell=.9]
(a1) edge[bend left=\t] node {X} (a2)
(a1) edge[bend right=\t] node[swap] {Y} (a2)
;
\draw[2cell]
node[between=a1 and a2 at .45, rotate=-90, 2label={above,\theta}] {\Rightarrow}
;
\end{tikzpicture}
\end{equation}
between the underlying functors of $X$ and $Y$.
\item\label{def:iu_morphism_ii}  
Denote by
\begin{equation}\label{IUT}
\IUT
\end{equation}
the category of $\IU$-spaces \pcref{def:iu_space} and $\IU$-morphisms \cref{iu_mor}.  Identities and composition are given by those of $\pSet$-natural transformations.
\item\label{def:iu_morphism_iii}  
A \emph{$G$-equivariant $\IU$-morphism}\index{G-equivariant IU-morphism@$G$-equivariant $\IU$-morphism}\index{IU-morphism@$\IU$-morphism!G-equivariant@$G$-equivariant} $\psi \cn X \to Y$ is a $\Gtopst$-natural transformation.
\item\label{def:iu_morphism_iv}
Denote by
\begin{equation}\label{GIUT}
\GIUT
\end{equation}
the category of $\IU$-spaces \pcref{def:iu_space} and $G$-equivariant $\IU$-morphisms.  Identities and composition are given by those of $\Gtopst$-natural transformations.\defmark
\end{enumerate}
\end{definition}

\begin{explanation}[$\IU$-Morphisms]\label{expl:iu_morphism}
An $\IU$-morphism $\theta \cn X \to Y$ \cref{iu_mor} is determined by a pointed function
\begin{equation}\label{iu_morphism_comp}
\stplus \fto{\theta_V} \Topgst(X_V,Y_V) \foreachspace V \in \IU,
\end{equation}
meaning a pointed morphism between pointed $G$-spaces
\begin{equation}\label{iu_mor_component}
X_V \fto{\theta_V} Y_V,
\end{equation}
such that the following naturality diagram commutes for each linear isometric isomorphism $f \cn V \fiso W$ in $\IU$.
\begin{equation}\label{iu_mor_natural}
\begin{tikzpicture}[vcenter]
\def\v{-1.4}
\draw[0cell]
(0,0) node (a11) {X_V}
(a11)++(2.3,0) node (a12) {Y_V}
(a11)++(0,\v) node (a21) {X_W}
(a12)++(0,\v) node (a22) {Y_W}
;
\draw[1cell=.9]
(a11) edge node {\theta_V} (a12)
(a12) edge node {Y_f} node[swap] {\iso} (a22)
(a11) edge node[swap] {X_f} node {\iso} (a21)
(a21) edge node {\theta_W} (a22)
;
\end{tikzpicture}
\end{equation}
The component pointed morphisms $\theta_V$ \cref{iu_mor_component} are \emph{not} required to be $G$-equivariant.  Identities and composition are defined using these components.

\parhead{Small hom sets}.  For each object $V \in \IU$, the naturality diagram \cref{iu_mor_natural} implies that $\theta_V$ is equal to the following composite pointed morphism, where $\upphi_V \cn V \fiso V'$ is the chosen $G$-linear isometric isomorphism in \cref{VV'} with $V' \in \IUsk$.
\begin{equation}\label{iu_mor_extend}
\begin{tikzpicture}[vcenter]
\def\v{-1.4}
\draw[0cell]
(0,0) node (a11) {X_V}
(a11)++(2.3,0) node (a12) {Y_V}
(a11)++(0,\v) node (a21) {X_{V'}}
(a12)++(0,\v) node (a22) {Y_{V'}}
;
\draw[1cell=.9]
(a11) edge node {\theta_V} (a12)
(a11) edge node[swap] {X_{\upphi_V}} node {\iso} (a21)
(a21) edge node {\theta_{V'}} (a22)
(a22) edge node[swap,pos=.45] {Y_{\upphi_V^{-1}}} node {\iso} (a12)
;
\end{tikzpicture}
\end{equation}
Thus, each $\IU$-morphism $\theta \cn X \to Y$ is determined by the \emph{set} of components
\begin{equation}\label{theta_iusk}
\Big\{X_V \fto{\theta_{V}} Y_V \cn V \in \IUsk\Big\}.
\end{equation}

\parhead{Determination on $\IUsk$}.  Conversely, to define an $\IU$-morphism $\theta \cn X \to Y$, it is enough to define the set of components in \cref{theta_iusk} such that the naturality diagram \cref{iu_mor_natural} commutes for each morphism $f \cn V \fiso W$ in $\IUsk$.  Indeed, a partially-defined $\theta$ \cref{theta_iusk} must be extended to all of $\IU$ using \cref{iu_mor_extend}.  For a general isomorphism $f \cn V \fiso W$ in $\IU$, the naturality diagram
\begin{equation}\label{iusk_mor}
\begin{tikzpicture}[vcenter]
\def\h{2} \def\v{-1.4}
\draw[0cell]
(0,0) node (a11) {X_V}
(a11)++(\h,0) node (a12) {X_{V'}}
(a12)++(\h,0) node (a13) {Y_{V'}}
(a13)++(\h,0) node (a14) {Y_V}
(a11)++(0,\v) node (a21) {X_W}
(a12)++(0,\v) node (a22) {X_{W'}}
(a13)++(0,\v) node (a23) {Y_{W'}}
(a14)++(0,\v) node (a24) {Y_W}
;
\draw[1cell=.9]
(a11) edge node {X_{\upphi_V}} node[swap] {\iso} (a12)
(a12) edge node {\theta_{V'}} (a13)
(a13) edge node {Y_{\upphi_V^{-1}}} node[swap] {\iso} (a14)
(a21) edge node {X_{\upphi_W}} node[swap] {\iso} (a22)
(a22) edge node {\theta_{W'}} (a23)
(a23) edge node {Y_{\upphi_W^{-1}}} node[swap] {\iso} (a24)
(a11) edge node[swap] {X_f} (a21)
(a14) edge node {Y_f} (a24)
;
\end{tikzpicture}
\end{equation}
commutes because, using the functoriality of $X$ and $Y$ \cref{iu_space_axioms}, it is the naturality diagram \cref{iu_mor_natural} of $\theta$ for the composite
\[V' \fto[\iso]{\upphi_V^{-1}} V \fto{f} W \fto[\iso]{\upphi_W} W' \inspace \IUsk.\]
In summary, specifying an $\IU$-morphism is equivalent to specifying the set of components in \cref{theta_iusk} such that the naturality diagram \cref{iu_mor_natural} commutes for each morphism $f$ in $\IUsk$.
\end{explanation}

\begin{explanation}[Equivariant $\IU$-Morphisms]\label{expl:eqiu_morphism}
A $G$-equivariant $\IU$-morphism $\psi \cn X \to Y$ is determined by a pointed $G$-morphism
\begin{equation}\label{eqiu_morphism_comp}
\stplus \fto{\psi_V} \Topgst(X_V,Y_V) \foreachspace V \in \IU,
\end{equation}
meaning a pointed $G$-morphism between pointed $G$-spaces
\begin{equation}\label{eqiu_mor_component}
X_V \fto{\psi_V} Y_V,
\end{equation}
such that the naturality diagram \cref{iu_mor_natural} commutes. 
Identities and composition are defined using the components $\psi_V$.  The discussion in the paragraphs of \cref{theta_iusk,iusk_mor} also applies to $G$-equivariant $\IU$-morphisms.  Thus, specifying a $G$-equivariant $\IU$-morphism $\psi$ is equivalent to specifying the \emph{set} of component pointed $G$-morphisms
\begin{equation}\label{psi_iusk}
\Big\{X_V \fto{\psi_{V}} Y_V \cn V \in \IUsk\Big\}
\end{equation}
such that the naturality diagram \cref{iu_mor_natural} commutes for each morphism $f$ in $\IUsk$.  The components of $\psi$ for $V \in \IU \setminus \IUsk$ are given by the diagrams \cref{iu_mor_extend} using the chosen $G$-linear isometric isomorphisms $\upphi_V$ in \cref{VV'}.
\end{explanation}

\subsection*{Topological Enrichment}

\begin{definition}[$\Gtop$-Enrichment of $\IUT$]\label{def:iut_gtop_enr}
For $\IU$-spaces $X,Y \cn \IU \to \Topgst$ \pcref{def:iu_space}, the morphism set
\begin{equation}\label{IUTXY_topology}
\IUT(X,Y) \bigsubset \prod_{V \in \IUsk} \Topgst(X_V,Y_V)
\end{equation}
is given the $G$-subspace topology, so $G$ acts componentwise by conjugation.  Equipped with these hom $G$-spaces, the category $\IUT$ in \cref{IUT} becomes a $\Gtop$-category.
\end{definition}

\begin{explanation}[Conjugation $G$-Action on $\IU$-Morphisms]\label{expl:iut_gtop_enr}
For an $\IU$-morphism $\theta \cn X \to Y$ in the $G$-space $\IUT(X,Y)$ in \cref{IUTXY_topology} and an element $g \in G$, the $g$-action on $\theta$ yields the $\IU$-morphism 
\[X \fto{g \cdot \theta} Y\]
whose $V$-component pointed morphism, for $V \in \IU$, is given by conjugating $\theta_V$ as follows.
\begin{equation}\label{iu_mor_gaction}
\begin{tikzpicture}[vcenter]
\def\v{-1.4}
\draw[0cell]
(0,0) node (a11) {X_V}
(a11)++(2.5,0) node (a12) {Y_V}
(a11)++(0,\v) node (a21) {X_{V}}
(a12)++(0,\v) node (a22) {Y_{V}}
;
\draw[1cell=.9]
(a11) edge node {(g \cdot \theta)_V} (a12)
(a11) edge node[swap] {\ginv} node {\iso} (a21)
(a21) edge node {\theta_{V}} (a22)
(a22) edge node[swap,pos=.5] {g} node {\iso} (a12)
;
\end{tikzpicture}
\end{equation}
Recall that $G$-equivariant $\IU$-morphisms are componentwise pointed $G$-morphisms \cref{eqiu_mor_component}.  Passing from the hom $G$-space $\IUT(X,Y)$ to its $G$-fixed point subspace yields the equality
\begin{equation}\label{giutxy_fixedpt}
\GIUT(X,Y) = \IUT(X,Y)^G,
\end{equation}
where $\GIUT(X,Y)$ is the hom set in the category $\GIUT$ in \cref{GIUT}.
\end{explanation}

\begin{definition}[$\Top$-Enrichment of $\GIUT$]\label{def:giut_top_enr}
For $\IU$-spaces $X,Y \cn \IU \to \Topgst$ \pcref{def:iu_space}, the morphism set $\GIUT(X,Y)$ is given the subspace topology using \cref{IUTXY_topology,giutxy_fixedpt}.  Equipped with these hom spaces, the category $\GIUT$ in \cref{GIUT} becomes a $\Top$-category.
\end{definition}

\subsection*{Equivariant Orthogonal Spectra}

\begin{definition}[$G$-Sphere]\label{def:g_sphere}\
\begin{itemize}
\item For each object $V \in \IU$ \pcref{def:indexing_gspace}, the \emph{$V$-sphere}\index{sphere} $S^V$ is the pointed $G$-space given by the one-point compactification $V \sqcup \{\infty\}$ of $V$ with $G$-fixed basepoint $\infty$. 
\item For each pair of objects $(V,W) \in (\IUsk)^2$, the $(V,W)$-component pointed $G$-homeomorphism
\begin{equation}\label{gsphere_multiplication}
S^V \sma S^W \fto[\iso]{\mu_{V,W}} S^{V \oplus W}
\end{equation}
is defined by
\begin{equation}\label{gsp_mult}
\mu_{V,W}(x; y) = \begin{cases}
\infty & \text{if $x = \infty \in S^V$ or $y = \infty \in S^W$, and}\\
x \oplus y & \text{if $(x;y) \in V \sma W$.}
\end{cases}
\end{equation}
\item For a morphism $f \cn V \fiso W$ in $\IU$, its pointed extension $S^V \to S^W$, sending $\infty$ to $\infty$, is also denoted by $f$.\defmark
\end{itemize}
\end{definition}

\begin{definition}[Orthogonal $G$-Spectra]\label{def:gsp_module}
Suppose $G$ is a compact Lie group.  An \emph{orthogonal $G$-spectrum}\index{orthogonal G-spectrum@orthogonal $G$-spectrum}\index{G-spectrum@$G$-spectrum!orthogonal} is a pair $(X,\umu)$ consisting of
\begin{itemize}
\item an $\IU$-space $X \cn \IU \to \Topgst$ \pcref{def:iu_space} and
\item $(V,W)$-component pointed $G$-morphisms
\begin{equation}\label{gsp_action_vw}
X_V \sma S^W \fto{\umu_{V,W}} X_{V \oplus W} \forspace (V,W) \in (\IUsk)^2
\end{equation}
\end{itemize}
such that the following naturality, unity, and associativity axioms hold.
\begin{description}
\item[Naturality]
For each pair of linear isometric isomorphisms
\[V \fto[\iso]{f} V' \andspace W \fto[\iso]{h} W' \inspace \IUsk,\]
the following naturality diagram of pointed morphisms commutes.
\begin{equation}\label{gsp_action_nat}

\end{equation}
\item[Unity]
The following unity diagram in $\Gtopst$ commutes for each object $V \in \IUsk$, where each instance of $\rho$ denotes a right unit isomorphism.
\begin{equation}\label{gsp_unity}
%
\end{equation}
\item[Associativity]  
The following associativity diagram in $\Gtopst$ commutes for each triple of objects $(U,V,W) \in (\IUsk)^3$, where each instance of $\al$ denotes an associativity isomorphism.
\begin{equation}\label{gsp_assoc}
%
\end{equation}
\end{description}
This finishes the definition of an orthogonal $G$-spectrum $(X,\umu)$.  We call $\umu$ the \index{sphere action}\emph{sphere action}.
\end{definition}

\begin{definition}[Categories of Orthogonal $G$-Spectra]\label{def:gsp_morphism}
Suppose $G$ is a compact Lie group.
\begin{enumerate}
\item\label{def:gsp_morphism_i} 
For orthogonal $G$-spectra $(X,\umu^X)$ and $(Y,\umu^Y)$, a \emph{morphism} of orthogonal $G$-spectra
\[(X,\umu^X) \fto{\tha} (Y,\umu^Y)\]
is an $\IU$-morphism \cref{iu_mor} such that the following diagram of pointed morphisms commutes for each pair of objects $(V,W) \in (\IUsk)^2$.
\begin{equation}\label{gsp_mor_axiom}
\begin{tikzpicture}[vcenter]
\def\v{-1.4}
\draw[0cell]
(0,0) node (a11) {X_V \sma S^W}
(a11)++(3,0) node (a12) {\phantom{X_{V \oplus W}}}
(a12)++(0,-.06) node (a12') {X_{V \oplus W}}
(a11)++(0,\v) node (a21) {Y_{V} \sma S^{W}}
(a12)++(0,\v) node (a22) {\phantom{Y_{V \oplus W}}}
(a22)++(0,-.06) node (a22') {Y_{V \oplus W}}
;
\draw[1cell=.9]
(a11) edge node {\umu^X_{V,W}} (a12)
(a12') edge node {\theta_{V \oplus W}} (a22')
(a11) edge node[swap] {\theta_V \sma 1} (a21)
(a21) edge node {\umu^Y_{V,W}} (a22)
;
\end{tikzpicture}
\end{equation}
Identities and composition of morphisms are defined using the components $\theta_V$ \cref{iu_mor_component}.  The $\Gtop$-category of orthogonal $G$-spectra and morphisms is denoted by $\GSp$.\label{not:GSp}
\item\label{def:gsp_morphism_ii} 
A \emph{$G$-morphism} of orthogonal $G$-spectra is a $G$-equivariant $\IU$-morphism (\cref{def:iu_morphism} \cref{def:iu_morphism_iii}) that makes the diagrams \cref{gsp_mor_axiom} commute.  The $\Top$-category of orthogonal $G$-spectra and $G$-morphisms is denoted by \label{not:Gspec}$\Gspec$.\defmark
\end{enumerate}
\end{definition}

In \cite[Def.\ II.2.6]{mandell_may}, the categories $\GSp$ and $\Gspec$ are denoted by, respectively, $\mathscr{I}_G\mathscr{S}$ and $G\mathscr{I}\mathscr{S}$.

\section{Orthogonal $G$-Spectra from $\Gskg$-Spaces}
\label{sec:semg}

For a compact Lie group $G$, this section reviews the prolongation functor
\[\GGTopii \fto{\Kg} \Gspec\]
from the category $\GGTopii$ of $\Gskg$-spaces and natural transformations \pcref{def:ggtop_smc} to the category $\Gspec$ of orthogonal $G$-spectra and $G$-morphisms (\cref{def:gsp_morphism} \cref{def:gsp_morphism_ii}).  The functor $\Kg$ is obtained from the unital symmetric monoidal $\Gtop$-functor
\[\GGTop \fto{(\Kg,\Kgtwo,\Kgzero)} \GSp\]
in \cite[\namecref{EqK:ch:semg} \ref*{EqK:ch:semg}]{yau-eqk} by keeping the same object assignment and passing to $G$-fixed subspaces in its $\Gtop$-enrichment. 
Thus, for a $\Gskg$-space $X \cn \Gsk \to \Gtopst$, $\Kg X$ is the same orthogonal $G$-spectrum regardless of whether $\Kg$ is regarded as a unital symmetric monoidal $\Gtop$-functor or as a functor.

\secoutline
\begin{itemize}
\item \cref{def:ggspace_gspectra} defines the object assignment of $\Kg$, which sends $\Gskg$-spaces to orthogonal $G$-spectra.
\item \cref{def:ggtop_gsp_mor} defines the morphism assignment of $\Kg$, which sends natural transformations to $G$-morphisms of orthogonal $G$-spectra.
\item The diagram \cref{Kgo_functors} summarizes the equivariant $K$-theory functors $\Kgo$ and $\Kgosg$ reviewed in this chapter.
\end{itemize}

\begin{definition}[$\Kg$ on Objects]\label{def:ggspace_gspectra}
Given a compact Lie group $G$ and a $\Gskg$-space \cref{ggtop_obj}
\[(\Gsk,\vstar) \fto{X} (\Gtopst,*),\]
the \index{orthogonal G-spectrum@orthogonal $G$-spectrum!multifunctorial K-theory@multifunctorial $K$-theory}\index{multifunctorial K-theory@multifunctorial $K$-theory!orthogonal G-spectrum@orthogonal $G$-spectrum}orthogonal $G$-spectrum \pcref{def:gsp_module}
\begin{equation}\label{Kg_object}
(\Kg X, \umu) \in \Gspec
\end{equation}
is defined as follows.
\begin{description}
\item[Object assignment of $\IU$-space]
The $\IU$-space \pcref{def:iu_space}
\begin{equation}\label{Kgx_iuspace}
\Kg X \cn \IU \to \Topgst
\end{equation}
sends each object $V \in \IU$ to the coend \pcref{def:coends}
\begin{equation}\label{Kgxv}
(\Kg X)_V = \int^{\angordn \in \Gsk} (S^V)^{\sma\angordn} \sma X\angordn
\end{equation}
 taken in $\Gtopst$.
\begin{itemize}
\item The pointed finite set $\sma\angordn \in \Fsk$ \pcref{def:smashFskGsk} is regarded as a discrete pointed $G$-space with the trivial $G$-action.
\item $S^V$ is the $V$-sphere \pcref{def:g_sphere}.
\item The pointed $G$-space 
\begin{equation}\label{SVsman}
(S^V)^{\sma\angordn} = \Topgst(\sma\angordn, S^V)
\end{equation}
consists of pointed morphisms $\sma\angordn \to S^V$ \cref{Gtopst_smc}, with $G$ acting by conjugation.  Since $G$ acts trivially on $\sma\angordn$, the group $G$ acts on $(S^V)^{\sma\angordn}$ by postcomposing with the $G$-action on $S^V$.
\end{itemize}
\item[Morphism assignment of $\IU$-space]
For a linear isometric isomorphism $f \cn V \fiso W$ in $\IU$, the pointed homeomorphism \cref{iu_space_xf}
\begin{equation}\label{Kgxf}
(\Kg X)_V \fto[\iso]{(\Kg X)_f} (\Kg X)_W
\end{equation}
is induced by the pointed homeomorphisms 
\[(S^V)^{\sma\angordn} \fto[\iso]{f \circ -} (S^W)^{\sma\angordn} \forspace \angordn \in \Gsk\]
that postcompose with the pointed homeomorphism $f \cn S^V \fiso S^W$.
\item[Sphere action]
For each pair of objects $(V,W) \in (\IUsk)^2$, the $(V,W)$-component pointed $G$-morphism \cref{gsp_action_vw} is defined by the following commutative diagram in $\Gtopst$.
\begin{equation}\label{Kgx_action_vw}
\begin{tikzpicture}[vcenter]
\def\h{5} \def\u{-1} \def\v{-1.4}
\draw[0cell=.9]
(0,0) node (a11) {(\Kg X)_V \sma S^W}
(a11)++(\h,0) node (a12) {(\Kg X)_{V \oplus W}}
(a11)++(0,\u) node (a21) {\big( \txint^{\angordn \in \Gsk} (S^V)^{\sma\angordn} \sma X\angordn \big) \sma S^W}
(a12)++(0,\u) node (a22) {\txint^{\angordn \in \Gsk} (S^{V \oplus W})^{\sma\angordn} \sma X\angordn}
(a21)++(0,\v) node (a3) {\txint^{\angordn \in \Gsk} \big( (S^V)^{\sma\angordn} \sma S^W \big) \sma X\angordn}
;
\draw[1cell=.8]
(a11) edge node {\umu_{V,W}} (a12)
(a11) edge[equal] (a21)
(a12) edge[equal] (a22)
(a21) edge node[swap] {\iso} (a3)
(a3) [rounded corners=2pt] -| node[pos=.25] {\assm = \txint\! \assm_{\angordn} \sma 1} (a22)
;
\end{tikzpicture}
\end{equation}
\begin{itemize}
\item The pointed $G$-homeomorphism denoted by $\iso$ uses
\begin{itemize}
\item the fact that $- \sma S^W$ commutes with coends and
\item the associativity isomorphism and braiding for the symmetric monoidal category $(\Gtopst,\sma)$ \cref{Gtopst_smc} to move $S^W$ to the left of $X\angordn$.
\end{itemize}
\item The pointed $G$-morphism $\assm$ is induced by the pointed $G$-morphisms
\begin{equation}\label{assembly_sph}
(S^V)^{\sma\angordn} \sma S^W \fto{\assm_{\angordn}} (S^{V \oplus W})^{\sma\angordn}
\end{equation}
for $\angordn \in \Gsk$ defined by the assignment
\[\begin{split}
& \big(\mq \smam\angordn \fto{t} S^V ; w \big) \in (S^V)^{\sma\angordn} \sma S^W\\
&\mapsto \big(\mq \smam\angordn \fto{t} S^V \fto{- \oplus w} S^{V \oplus W} \big) \in (S^{V \oplus W})^{\sma\angordn}.
\end{split}\]
\end{itemize}
\end{description}
It is proved in \cite[\namecref{EqK:Kgx_welldef} \ref*{EqK:Kgx_welldef}]{yau-eqk} that $(\Kg X, \umu)$ is an orthogonal $G$-spectrum.
\end{definition}

\begin{definition}[$\Kg$ on Morphisms]\label{def:ggtop_gsp_mor}
Given a compact Lie group $G$ and a natural transformation \cref{ggtop_mor}
\begin{equation}\label{ggspace_mor_theta}

\end{equation}
between $\Gskg$-spaces $X$ and $X'$ \cref{ggtop_obj}, the $G$-morphism of orthogonal $G$-spectra (\cref{def:gsp_morphism} \cref{def:gsp_morphism_ii})
\begin{equation}\label{Kg_theta}
(\Kg X,\umu) \fto{\Kg\theta} (\Kg X',\umu')
\end{equation}
has, for each object $V \in \IU$, $V$-component pointed $G$-morphism defined by the following commutative diagram in $\Gtopst$.
\begin{equation}\label{Kg_theta_v}
%
\end{equation}
It follows from \cite[\namecref{EqK:Kg_gtop_functor} \ref*{EqK:Kg_gtop_functor}]{yau-eqk}, by passing to $G$-fixed subspaces in the $\Gtop$-enrichment, that 
\begin{equation}\label{Kg_functor}
\GGTopii \fto{\Kg} \Gspec
\end{equation}
is a functor with
\begin{itemize}
\item object assignment $X \mapsto \Kg X$ \cref{Kg_object} and
\item morphism assignment $\theta \mapsto \Kg\theta$ \cref{Kg_theta}.\defmark
\end{itemize}
\end{definition}

\subsection*{Summary}

The diagram 
\begin{equation}\label{Kgo_functors}
%
\end{equation}
summarizes the two equivariant $K$-theory functors\index{multifunctorial K-theory@multifunctorial $K$-theory!underlying functor}\index{K-theory@$K$-theory!multifunctorial - underlying functor}\index{K-theory@$K$-theory!via J-theory@via $J$-theory}
\[\begin{split}
\Kgo &= \Kg \circ \clast \circ \Jgo \andspace\\ 
\Kgosg &= \Kg \circ \clast \circ \Jgosg
\end{split}\] 
constructed in this chapter.  The functors $\Kgo$ and $\Kgosg$ are obtained from the $\Gcat$-multifunctors and $\Gtop$-multifunctors in \cref{Kgo_multifunctor} by keeping the same object assignments and restricting to $G$-fixed 1-ary 1-cells and 2-cells in the $\Gcat$-enrichment or to 1-ary $G$-fixed subspaces in the $\Gtop$-enrichment.  Each of $\Kgo$ and $\Kgosg$ is a composite of three functors.
\begin{description}
\item[$J$-theory] For a group $G$ and a $\Tinf$-operad $\Op$ \pcref{as:OpA}, $\Jgo$ and $\Jgosg$ are the (strong) $J$-theory 2-functors in \cref{thm:Jgo_twofunctor}.  At the object level, they send $\Op$-pseudoalgebras to $\Gskg$-categories \pcref{A_ptfunctor}.  At the 1-cell level, they send lax $\Op$-morphisms and $\Op$-pseudomorphisms to natural transformations \pcref{def:Jgo_pos_obj}.
\item[Classifying space] For a group $G$, the functor $\clast$ \cref{clast} sends $\Gskg$-categories and natural transformations to $\Gskg$-spaces and natural transformations by composing and whiskering with the classifying space functor \pcref{expl:clast_functor}.
\item[Prolongation] For a compact Lie group $G$, the functor $\Kg$ \cref{Kg_functor} sends $\Gskg$-spaces and natural transformations to orthogonal $G$-spectra and $G$-morphisms \pcref{def:ggspace_gspectra,def:ggtop_gsp_mor}.
\end{description}
For an $\Op$-pseudoalgebra $\A$ \pcref{def:pseudoalgebra} and an object $V \in \IU$ \pcref{def:indexing_gspace}, $\Kgo$ and $\Kgosg$ yield the pointed $G$-spaces
\begin{equation}\label{kgoav}
\begin{split}
(\Kgo\A)_V &= \int^{\angordn \in \Gsk} (S^V)^{\sma\angordn} \sma \cla(\Aangordn) \andspace\\
(\Kgosg\A)_V &= \int^{\angordn \in \Gsk} (S^V)^{\sma\angordn} \sma \cla(\Aangordnsg)\\
\end{split}
\end{equation}
with the sphere action defined in \cref{Kgx_action_vw}.  For each object $\angordn \in \Gsk$ \cref{Gsk_objects}, $\Aangordn$ is the pointed $G$-category of $\angordn$-systems in $\A$ \pcref{def:nsystem,def:nsystem_morphism,def:Aangordn_gcat}, and $\Aangordnsg$ is the pointed full $G$-subcategory of strong $\angordn$-systems defined in \cref{Aangordnsg,sgAordnbe}.

 

%% file: chap/ggcatg.tex
For an arbitrary group $G$, this chapter constructs a 2-category $\GGCatg$ \pcref{def:ggcatg} and a 2-equivalence \pcref{thm:ggcat_ggcatg_iieq}
\[\GGCatg \fto[\sim]{\igst} \GGCatii\]
to the codomain $\GGCatii$ of the $J$-theory 2-functors $\Jgo$ and $\Jgosg$ \pcref{thm:Jgo_twofunctor}.  \cref{ch:hgo} shows that each $J$-theory 2-functor factors as an $H$-theory 2-functor with codomain $\GGCatg$, followed by the 2-equivalence $\igst$.  These factorizations of the $J$-theory 2-functors are used in \cref{part:kgo_shi_comp} to categorically compare our equivariant $K$-theory with Shimakawa $K$-theory.  \Cref{table.ggcat} summarizes the 2-categories $\GGCatii$ and $\GGCatg$.
\begin{figure}[H] 
\centering
\resizebox{.9\linewidth}{!}{%
{\renewcommand{\arraystretch}{1.3}%
{\setlength{\tabcolsep}{1ex}
\begin{tabular}{c|cr|cr}
& $\GGCatii$ & \eqref{def:GGCat} & $\GGCatg$ & \eqref{def:ggcatg} \\ \hline
indexing categories & $\Gsk$ & \eqref{def:Gsk} & $\GG$ & \eqref{def:GG} \\
base categories & $\Gcatst$ & \eqref{def:gcatst} & $\Catgst$ & \eqref{Catgst_iicat} \\
objects & pointed functors $\Gsk \to \Gcatst$ & \eqref{ggcat_obj} & pointed $G$-functors $\GG \to \Catgst$ & \eqref{ggcatg_obj} \\
1-cells & natural transformations & \eqref{ggcat_mor} & $G$-natural transformations & \eqref{ggcatg_icell} \\
2-cells & modifications & \eqref{ggcat_iicell} & $G$-modifications & \eqref{ggcatg_iicell} \\
\end{tabular}}}}
\caption{A comparison of two 2-categories.}
\label{table.ggcat}
\end{figure}

\organization
This chapter consists of the following sections.

\secname{sec:FG}
This section constructs the small pointed $G$-category $\FG$, which is the $G$-equivariant analogue of the category $\Fsk$ of pointed finite sets and pointed functions.  The objects of $\FG$ are pointed finite $G$-sets.  Its morphisms are pointed functions between underlying pointed finite sets.  The group $G$ acts on morphisms by conjugation.  When it is equipped with the smash product of pointed finite $G$-sets, $\FG$ becomes a naive permutative $G$-category.

\secname{sec:GG}
This section constructs the small pointed $G$-category $\GG$, which is the $G$-equivariant analogue of the category $\Gsk$ \pcref{def:Gsk}.  The objects of $\GG$ are finite tuples of pointed finite $G$-sets.  Its morphisms are defined like those for $\Gsk$.  The $G$-action is inherited from the one on $\FG$.  When it is equipped with the concatenation product, $\GG$ becomes a naive permutative $G$-category.  There is a strict symmetric monoidal pointed $G$-functor $\GG \to \FG$ given by the smash product on objects.

\secname{sec:ggcatg}
Using $\GG$ as the indexing $G$-category and the 2-category $\Catgst$ \cref{Catgst_iicat}, this section constructs the 2-category $\GGCatg$, whose objects, called $\GGG$-categories, are pointed $G$-functors $\GG \to \Catgst$.  There is a full subcategory inclusion $\ig \cn \Gsk \to \GG$ that equips each pointed finite set with the trivial $G$-action.  Precomposition with the functor $\ig$ induces a comparison 2-functor 
\[\GGCatg \fto{\igst} \GGCatii.\]

\secname{sec:ggcat_ggcatg_iieq}
This section proves that the 2-functor $\igst$ is a 2-equivalence by constructing an explicit inverse 2-functor
\[\GGCatii \fto{\Lg} \GGCatg\]
along with unit and counit 2-natural isomorphisms.

\section{The Indexing $G$-Category $\FG$}
\label{sec:FG}

For an arbitrary group $G$, this section defines the small pointed $G$-category $\FG$, along with some auxiliary constructions to be used in \cref{sec:GG} to construct the indexing $G$-category $\GG$.   The $G$-category $\FG$ consists of  pointed finite $G$-sets and pointed functions with the conjugation $G$-action.  It is the $G$-equivariant analogue of the category $\Fskel$ \pcref{def:Fsk}.

\secoutline
\begin{itemize}
\item \cref{def:ptGset,def:FG} define the small pointed $G$-category $\FG$.  
\item \cref{def:FG_permutative,FG_permutative} define a naive permutative $G$-category structure on $\FG$.
\item \cref{def:FG_smashpower,def:injectionsFG} define smash powers of $\FG$ and reindexing functors between them.
\end{itemize}

\subsection*{Pointed $G$-Sets}

Recall from \cref{ordn} that $\ordn$ denotes the pointed finite set $\{0 < 1 < \cdots < n\}$ with basepoint 0.  For a set $S$, $\Aut(S)$\label{not:AutS} denotes the group of self-bijections $S \fiso S$.

\begin{definition}[Pointed Finite $G$-Sets]\label{def:ptGset}
Suppose $G$ is a group and $n \geq 0$.
\begin{itemize}
\item A \emph{$G$-set}\index{G-set@$G$-set} is a pair $(S,\mu)$ consisting of a set $S$, called the \emph{underlying set}, and a group homomorphism $\mu \cn G \to \Aut(S)$, called the \index{G-action@$G$-action}\emph{$G$-action}.  An element $\mu(g)(a) \in S$ is often denoted by $g \cdot a$ or $ga$ for $(g,a) \in G \times S$.  A $G$-set $(S,\mu)$ is \emph{trivial}\index{G-action@$G$-action!trivial}\index{trivial $G$-action} if $\mu$ is the constant homomorphism at the identity $1_S$.  A $G$-set $(S,\mu)$ is \emph{finite}\index{G-set@$G$-set!finite}\index{finite G-set@finite $G$-set} if $S$ is a finite set. 
\item A \emph{pointed $G$-set} is a $G$-set $(S,\mu)$ equipped with a $G$-fixed basepoint $\bp \in S$.  Such a pointed $G$-set is also denoted by $S^\mu$.  A \emph{pointed finite $G$-set} is a pointed $G$-set whose underlying set is finite.
\item Each group homomorphism $\be \cn G \to \Sigma_n$ is extended to a pointed finite $G$-set 
\begin{equation}\label{ordn_be}
(\ord{n}, \be) = \ord{n}^\be
\end{equation}
with the basepoint $0 \in \ord{n}$ fixed by $\be(g)$ for each $g \in G$.  The pointed set $\ordn$ is called the \emph{underlying pointed finite set} of $\ordn^\be$.
\end{itemize}
Unless otherwise specified, a pointed finite $G$-set means one of the form $\ord{n}^\be$ \cref{ordn_be}.
\end{definition}

Recall that a pointed $G$-category is a $G$-category equipped with a $G$-fixed object, called the basepoint \pcref{def:GCat,def:ptGcat}.

\begin{definition}[$G$-Category of Pointed Finite $G$-Sets]\label{def:FG}
For a group $G$, the small pointed $G$-category $\FG$\index{FG@$\FG$} is defined as follows.
\begin{description}
\item[Objects] The objects of $\FG$ are pointed finite $G$-sets of the form $\ord{n}^\be$ \cref{ordn_be}.  The group $G$ acts trivially on objects, meaning $g \cdot \ordn^\be = \ordn^\be$ for each $g \in G$.
\item[Morphisms] The morphisms of $\FG$ are pointed functions between underlying pointed finite sets, with $G$ acting by conjugation.  A morphism that factors through the object $\ord{0}$ is called a \index{0-morphism}\emph{0-morphism}.
\item[Composition and identity morphisms] These are defined as those for pointed functions.
\item[Basepoint] The basepoint of $\FG$ is the initial-terminal object $\ord{0}$.
\end{description}
If $G$ is the trivial group, then $\FG$ is identified with the pointed category $\Fsk$ \pcref{def:Fsk}.
\end{definition}

\begin{explanation}[Unpacking $\FG$]\label{expl:FG}
Each hom set of $\FG$ is a pointed $G$-set, with basepoint given by the 0-morphism.  Morphisms in $\FG$ are \emph{not} required to be $G$-equivariant.  For a pointed function between pointed finite $G$-sets 
\[\ord{m}^\al \fto{\psi} \ord{n}^\be\]
and $g \in G$, the \emph{conjugation $g$-action}\index{conjugation action}\index{G-action@$G$-action!conjugation} on $\psi$ is given by the following composite pointed function.
\begin{equation}\label{gpsi}
\begin{tikzpicture}[baseline={(a.base)}]
\def\h{2} \def\u{.7}
\draw[0cell]
(0,0) node (a) {\ord{m}^\al}
(a)++(1.1*\h,0) node (b) {\ord{m}^\al}
(b)++(.9*\h,0) node (c) {\ord{n}^\be}
(c)++(\h,0) node (d) {\ord{n}^\be}
;
\draw[1cell=.9]
(a) edge node {\al(g)^{-1}} node[swap] {\iso} (b)
(b) edge node {\psi} (c)
(c) edge node {\be(g)} node[swap] {\iso} (d)
;
\draw[1cell=.9]
(a) [rounded corners=2pt] |- ($(b)+(0,\u)$)
-- node {g \cdot \psi} ($(c)+(0,\u)$) -| (d)
;
\end{tikzpicture}
\end{equation}
Abbreviating $\al(g)$ and $\be(g)$ to $g$, we have
\[g \cdot \psi = g \circ \psi \circ \ginv.\]
In particular, $\psi$ is $G$-equivariant if and only if it is $G$-fixed, meaning
\[g \cdot \psi = \psi \foreachspace g \in G.\]
The pointed $G$-category $\FG$ is denoted by $\GaG$ in \cite{shimakawa89,shimakawa91}.  The notation $\FG$ follows the usage in \cite{gmmo19,gmmo23,mmo}. 
\end{explanation}

\subsection*{Naive Permutative $G$-Category}
The next definition extends the permutative category $(\Fsk,\sma,\ord{1},\xi)$ \pcref{def:Fsk_permutative} to a naive permutative $G$-category structure \pcref{expl:naive_perm_Gcat} on $\FG$.  Recall from \cref{ufsn} that $\ufs{n}$ denotes the unpointed finite set $\{1 < 2 < \cdots < n\}$ with $\ufs{0} = \emptyset$.

\begin{definition}\label{def:FG_permutative}
The small pointed $G$-category $\FG$ \pcref{def:FG} is given the following naive permutative $G$-category structure $(\sma, \ord{1}, \xi)$.
\begin{description}
\item[Monoidal product on objects] The monoidal product
\[\FG \times \FG \fto{\sma} \FG\]
has object assignment
\begin{equation}\label{mal-sma-nbe}
\ord{m}^\al \sma \ord{n}^\be = \ord{mn}^{\al \sma \be}
\end{equation}
with $\ord{mn}$ ordered lexicographically \cref{lex_bijection} away from the basepoint $0 \in \ord{mn}$.  The $G$-action $\al \sma \be$ fixes 0 and is given diagonally by
\begin{equation}\label{al-sma-be}
(\al \sma \be)(g)(a,b) = \big(\al(g)(a), \be(g)(b)\big)
\end{equation}
for $g \in G$ and $(a,b) \in \ufs{m} \times \ufs{n}$.  This $g$-action is usually denoted by
\[g \cdot (a,b) = (g \cdot a, g \cdot b)\]
with $\cdot$ omitted if there is no danger of confusion.
\item[Monoidal product on morphisms] 
The smash product $\sma$ is extended to morphisms using the lexicographic ordering.  In other words, for pointed functions
\begin{equation}\label{fmp-hnq}
\ord{m}^\al \fto{f} \ord{p}^\chi \andspace \ord{n}^\be \fto{h} \ord{q}^\de,
\end{equation}
the pointed function
\[\ord{m}^\al \sma \ord{n}^\be \fto{f \sma h} \ord{p}^\chi \sma \ord{q}^\de\]
is defined by
\begin{equation}\label{fsmah-ab}
(f \sma h)(0) = 0 \andspace (f \sma h)(a,b) = \big(f(a), h(b)\big)
\end{equation}
for $(a,b) \in \ufs{m} \times \ufs{n}$. 
\item[Monoidal unit] The monoidal unit is the pointed finite set $\ord{1}$ with the trivial $G$-action.
\item[Braiding] The component of the braiding $\xi$ at a pair of objects $(\ord{m}^\al, \ord{n}^\be)$ is the pointed bijection
\begin{equation}\label{FG_braiding}
\ord{m}^\al \sma \ord{n}^\be \fto[\iso]{\xi_{\ord{m}^\al, \ord{n}^\be}} \ord{n}^\be \sma \ord{m}^\al
\end{equation}
given by the $(m,n)$-transpose permutation $\twist_{m,n}$ \cref{eq:transpose_perm} away from the basepoint $0 \in \ord{mn}$.  In other words, using the lexicographic ordering, the braiding is given by
\[\xi_{\ord{m}^\al, \ord{n}^\be} (a,b) = (b,a)\]
\end{description}
for $(a,b) \in \ufs{m} \times \ufs{n}$. 
\end{definition}

\begin{lemma}\label{FG_permutative}
The quadruple $(\FG,\sma,\ord{1},\xi)$ in \cref{def:FG,def:FG_permutative} is a naive permutative $G$-category.
\end{lemma}

\begin{proof}
The proof given in \cite[\namecref{EqK:Fsk_permutative} \ref*{EqK:Fsk_permutative}]{yau-eqk} that the quadruple $(\Fsk,\sma,\ord{1},\xi)$ is a permutative category also applies to $(\FG,\sma,\ord{1},\xi)$.

\parhead{$G$-equivariance}.  To check that the structure $(\sma, \ord{1}, \xi)$ is $G$-equivariant, first observe that the monoidal unit $\ord{1}$ is $G$-fixed, since each object of $\FG$ is $G$-fixed.  

Since $G$ acts trivially on objects of $\FG$ and $\FG \times \FG$, the functor $\sma$ is $G$-equivariant if and only if it is $G$-equivariant on morphisms.  Using \cref{gpsi,al-sma-be,fsmah-ab}, the following computation proves that $\sma$ is $G$-equivariant for pointed functions $f$ and $h$ as defined in \cref{fmp-hnq}, where $g \in G$ and $(a,b) \in \ufs{m} \times \ufs{n}$.
\[\begin{split}
\big(g \cdot (f \sma h) \big)(a,b) 
&= g \big((f \sma h)(\ginv a, \ginv b) \big)\\
&= g \big(f(\ginv a), h(\ginv b)\big)\\
&= \big(g f(\ginv a), g h(\ginv b)\big)\\
&= \big( (g \cdot f)(a), (g \cdot h)(b\big))\\
&= \big((g \cdot f) \sma (g \cdot h)\big)(a,b)
\end{split}\]
This proves that $\sma$ is a $G$-functor.

Since $G$ acts trivially on objects of $\FG \times \FG$, the braiding $\xi$ is $G$-equivariant if and only if 
\[g \cdot \xi_{\ord{m}^\al, \ord{n}^\be} = \xi_{\ord{m}^\al, \ord{n}^\be} \forspace g \in G.\]
This equality is proved using \cref{gpsi,al-sma-be,FG_braiding} as follows, where $(a,b) \in \ufs{m} \times \ufs{n}$.
\[\begin{split}
\big(g \cdot \xi_{\ord{m}^\al, \ord{n}^\be}\big)(a,b)
&= g \xi_{\ord{m}^\al, \ord{n}^\be}  (\ginv a, \ginv b)\\
&= g (\ginv b, \ginv a)\\
&=\big(g \ginv b, g \ginv a\big)\\
&= (b,a)\\
&= \xi_{\ord{m}^\al, \ord{n}^\be}(a,b)
\end{split}\]
This proves that the braiding $\xi$ is a $G$-natural isomorphism.
\end{proof}

\subsection*{Smash Powers}
\cref{def:FG_smashpower} extends the smash powers of $\Fsk$ \pcref{def:Fsk_smashpower} to $\FG$.

\begin{definition}\label{def:FG_smashpower}
For a group $G$ and $q \geq 0$, the small pointed $G$-category $\FGsma{q}$\index{FGq@$\FGsma{q}$}\index{smash power} is defined as follows.
\begin{description}
\item[$q=0$] The pointed $G$-category
\[\FGsma{0} = \big\{\vstar \rightleftarrows \ang{}\big\}\]
consists of the initial-terminal basepoint $\vstar$, the empty tuple $\ang{}$, the identity morphisms of $\vstar$ and $\ang{}$, and the nonidentity morphisms $\vstar \to \ang{} \to \vstar$.  The group $G$ acts trivially on $\FGsma{0}$.
\item[$q>0$] The small pointed $G$-category $\FGsma{q}$ is the $q$-fold smash power of $\FG$ defined as follows. 
\begin{description}
\item[Objects] It has an initial-terminal basepoint $\vstar$.  A typical object in $\FGsma{q}$ is a $q$-tuple
\begin{equation}\label{ordnbeta}
\ang{\ord{n}^\be} = \ang{\ord{n}_i^{\be_i}}_{i \in \ufs{q}} 
= \big( \ord{n}_1^{\be_1}, \ldots, \ord{n}_q^{\be_q}\big)
\end{equation}
with each $\ord{n}_i^{\be_i}$ a pointed finite $G$-set \cref{ordn_be}.  If any $\ord{n}_i = \ord{0}$, then $\ang{\ord{n}^\be} = \vstar$. 
\item[Morphisms] A typical morphism in $\FGsma{q}$ is a $q$-tuple
\begin{equation}\label{angpsiFG}
\ang{\psi} = \ang{\psi_i}_{i \in \ufs{q}} \cn 
\ang{\ord{m}_i^{\al_i}}_{i \in \ufs{q}} \to \ang{\ord{n}_i^{\be_i}}_{i \in \ufs{q}}
\end{equation}
with each
\[\ord{m}_i^{\al_i} \fto{\psi_i} \ord{n}_i^{\be_i}\]
a morphism in $\FG$, meaning a pointed function between underlying pointed finite sets.  If any $\psi_i$ factors through $\ord{0}$, then $\ang{\psi}$ is the \index{0-morphism}\emph{0-morphism}, meaning that it factors through the basepoint $\vstar$.  A \index{nonzero morphism}\emph{nonzero morphism} is a morphism that does not factor through $\vstar$.
\item[$G$-action] The group $G$ acts trivially on objects of $\FGsma{q}$ and by conjugation on morphisms.
\end{description}
\end{description}
A nonbasepoint object or a nonzero morphism in $\FGsma{q}$ is said to have \emph{length} $q$.
\end{definition}

\begin{explanation}\label{expl:gangpsi}
The trivial $G$-action on objects of $\FGsma{q}$ means that
\[g \cdot \angordnbe = \angordnbe\]
for each $g \in G$ and object $\angordnbe \in \FGsma{q}$.  The conjugation $g$-action on a morphism $\angpsi$ yields the morphism
\begin{equation}\label{gangpsi}
g \cdot \ang{\psi} = \ang{g \cdot \psi_i}_{i \in \ufs{q}} \cn 
\ang{\ord{m}_i^{\al_i}}_{i \in \ufs{q}} \to \ang{\ord{n}_i^{\be_i}}_{i \in \ufs{q}}
\end{equation}
with each pointed function
\begin{equation}\label{gpsii}
\begin{tikzpicture}[vcenter]
\def\h{2} \def\u{.7}
\draw[0cell]
(0,0) node (a) {\ord{m}_i^{\al_i}}
(a)++(1.2*\h,0) node (b) {\ord{m}_i^{\al_i}}
(b)++(.9*\h,0) node (c) {\ord{n}_i^{\be_i}}
(c)++(\h,0) node (d) {\ord{n}_i^{\be_i}}
;
\draw[1cell=.9]
(a) edge node {\al_i(g)^{-1}} node[swap] {\iso} (b)
(b) edge node {\psi_i} (c)
(c) edge node {\be_i(g)} node[swap] {\iso} (d)
;
\draw[1cell=.9]
(a) [rounded corners=2pt] |- ($(b)+(0,\u)$)
-- node {g \cdot \psi_i} ($(c)+(0,\u)$) -| (d)
;
\end{tikzpicture}
\end{equation}
defined as in \cref{gpsi}.  Each morphism set of $\FGsma{q}$ is a pointed $G$-set, with the 0-morphism as the basepoint.
\end{explanation}

\subsection*{Reindexing}
General morphisms in $\GG$ involve the reindexing device in the following definition.

\begin{definition}\label{def:injectionsFG}
For an injection $h \cn \ufs{q} \to \ufs{r}$, the pointed $G$-functor  
\[\FGsma{q} \fto{h_*} \FGsma{r}\]
is defined as follows.
\begin{description}
\item[$q = r = 0$] 
In this case, $h$ is the identity function on $\ufs{0} = \emptyset$, and $h_*$ is defined as the identity functor on $\FGsma{0}$.
\item[$q = 0 < r$] 
In this case, the pointed functor $h_*$ is determined by the following object assignment.
\begin{equation}\label{hangemptyFG}
\left\{\begin{split}
h_* \vstar &= \vstar\\
h_*\ang{} &= \ang{\ord{1}}_{j \in \ufs{r}} = \left(\ord{1}, \ldots, \ord{1}\right) \in \FGsma{r}
\end{split}\right.
\end{equation}
In the object $\ang{\ord{1}}_{j \in \ufs{r}}$, each copy of $\ord{1} \in \FG$ has the trivial $G$-action.
\item[$q>0$] 
Given an object $\ang{\ord{n}_i^{\be_i}}_{i \in \ufs{q}}$ \cref{ordnbeta} and a morphism $\ang{\psi_i}_{i \in \ufs{q}}$ \cref{angpsiFG}, we define the object and morphism
\begin{equation}\label{reindexing_functorFG}
\begin{split}
h_*\ang{\ord{n}_i^{\be_i}}_{i \in \ufs{q}} &= \bang{\ord{n}_{\hinv(j)}^{\be_{\hinv(j)}}}_{j \in \ufs{r}} \andspace\\
h_*\ang{\psi_i}_{i \in \ufs{q}} &= \ang{\psi_{\hinv(j)}}_{j \in \ufs{r}}
\end{split}
\end{equation}
in $\FGsma{r}$.  If $\hinv(j) = \emptyset$, then
\begin{equation}\label{ordn_emptyFG}
\ord{n}_{\emptyset}^{\be_{\emptyset}} = \ord{1} \andspace 
\ord{1} \fto{\psi_{\emptyset} = 1_{\ord{1}}} \ord{1}.
\end{equation}
\end{description}
We call $h$ a \index{reindexing injection}\emph{reindexing injection} and $h_*$ a \index{reindexing G-functor@reindexing $G$-functor}\emph{reindexing $G$-functor}.
\end{definition}

\begin{explanation}\label{expl:reindexingFG}
The pointed $G$-functor $h_*$ in \cref{def:injectionsFG} is well defined for the following reasons.
\begin{itemize}
\item Suppose $\ord{n}_i = \ord{0}$ for some $i \in \ufs{q}$. Since $h \cn \ufs{q} \to \ufs{r}$ is an injection, at least one entry in $h_*\ang{\ord{n}_i^{\be_i}}_{i \in \ufs{q}}$ is $\ord{0}$.  Thus, we have that
\[\ang{\ord{n}_i^{\be_i}}_{i \in \ufs{q}} = \vstar \in \FGsma{q} 
\impliespace h_*\ang{\ord{n}_i^{\be_i}}_{i \in \ufs{q}} = \vstar \in \FGsma{r}.\]
\item If some $\psi_i$ factors through $\ord{0}$, then at least one entry in $h_*\ang{\psi_i}_{i \in \ufs{q}}$ factors through $\ord{0}$.  Thus, if $\ang{\psi}$ is the 0-morphism, then so is $h_*\ang{\psi_i}_{i \in \ufs{q}}$.
\item $h_*$ is $G$-equivariant because $G$ acts trivially on objects, and the identity morphism $1_{\ord{1}}$ is $G$-fixed.
\end{itemize}  
The assignment $h \mapsto h_*$ is functorial in the sense that it preserves identities and composition.
\end{explanation}

\section{The Indexing $G$-Category $\GG$}
\label{sec:GG}

This section defines the $G$-equivariant extension of the category $\Gsk$ \pcref{def:Gsk}, where the role of $\Fsk$ is now played by the pointed $G$-category $\FG$ \pcref{def:FG}.

\secoutline
\begin{itemize}
\item \cref{def:GG,GG_Gcategory} construct the small pointed $G$-category $\GG$.
\item \cref{def:GG_permutative,GG_permutative} construct a naive permutative $G$-category structure on $\GG$.
\item \cref{def:smashFGGG,sma_symmon} construct the strict symmetric monoidal pointed $G$-functor $\sma \cn \GG \to \FG$.
\end{itemize}

\subsection*{The $G$-Category $\GG$}

Recall from \cref{def:injections} that $\Inj$ is the category with the unpointed finite sets $\ufs{n} = \{1, 2,\ldots , n\}$ \cref{ufsn} for $n \geq 0$ as objects and injections as morphisms.  

\begin{definition}\label{def:GG}
For a group $G$, the small pointed $G$-category $\GG$\index{GG@$\GG$} is defined as follows.
\begin{description}
\item[Objects] The object pointed $G$-set is the wedge
\begin{equation}\label{GG_objects}
\Ob(\GG) = \bigvee_{q \geq 0} \Ob(\FGsma{q})
\end{equation}
that identifies the basepoints $\vstar \in \FGsma{q}$ \pcref{def:FG_smashpower} for $q \geq 0$.  The identified object is the initial-terminal basepoint $\vstar \in \GG$.  The group $G$ acts trivially on the objects of $\GG$, meaning 
\begin{equation}\label{GG_obj_gaction}
g \cdot \angordmal = \angordmal
\end{equation}
for $g \in G$ and $\angordmal \in \GG$.
\item[Morphisms]
Given an object $\ang{\ord{m}^\al} \in \FGsma{p}$ and an object $\ang{\ord{n}^\be} \in \FGsma{q}$, the pointed $G$-set of morphisms is defined as the following wedge.
\begin{equation}\label{GG_morphisms}
\begin{split}
\GG\big(\ang{\ord{m}^\al}, \ang{\ord{n}^\be}\big) 
&= \bigvee_{f \in \Inj(\ufs{p},\, \ufs{q})}~ \FGsma{q}\big(f_* \ang{\ord{m}^\al}, \bang{\ord{n}^\be} \big)\\
&= \bigvee_{f \in \Inj(\ufs{p},\, \ufs{q})}~ \bigwedge_{i \in \ufs{q}}~ \FG\big(\ord{m}_{\finv(i)}^{\al_{\finv(i)}} \scs \ord{n}_i^{\be_i} \big)
\end{split}
\end{equation}
The basepoint of each pointed $G$-set of morphisms is the \index{0-morphism}\emph{0-morphism}, which is the unique morphism that factors through $\vstar \in \GG$.  A morphism that is not the 0-morphism is called a \index{nonzero morphism}\emph{nonzero morphism}. 

In \cref{GG_morphisms}, for each reindexing injection $f \cn \ufs{p} \to \ufs{q}$, 
\[\FGsma{p} \fto{f_*} \FGsma{q}\]
is the reindexing $G$-functor in \cref{def:injectionsFG}. 
\begin{itemize}
\item If $q=0$ in \cref{GG_morphisms}, then $p=0$ and $f \cn \ufs{0} \to \ufs{0}$ is $1_\emptyset$.  The pointed set
\begin{equation}\label{GG_empty_mor}
\GG\big(\ang{}, \ang{}\big) = \FGsma{0}\big(\ang{}, \ang{}\big)
\end{equation}
consists of the identity morphism of $\ang{}$ and the 0-morphism $\ang{} \to \vstar \to \ang{}$.  The group $G$ acts trivially on this set.
\item If $q>0$ in \cref{GG_morphisms}, a morphism in $\GG$ is a pair
\begin{equation}\label{fangpsiGG}
\ang{\ord{m}^\al} \fto{(f, \ang{\psi})} \ang{\ord{n}^\be}
\end{equation}
consisting of
\begin{itemize}
\item a reindexing injection $f \cn \ufs{p} \to \ufs{q}$ and
\item a morphism 
\[\ang{\psi} = \ang{\psi_i}_{i \in \ufs{q}} \cn f_* \ang{\ord{m}^\al} \to \ang{\ord{n}^\be} \inspace \FGsma{q}.\]
\end{itemize}
If some pointed function
\[\ord{m}_{\finv(i)}^{\al_{\finv(i)}} \fto{\psi_i} \ord{n}_i^{\be_i}\]
factors through $\ord{0} \in \FG$, then $(f, \ang{\psi})$ is the 0-morphism, factoring through the basepoint $\vstar$.  
\item For $g \in G$, the $g$-action is given by 
\begin{equation}\label{GG_Gaction}
g \cdot (f, \ang{\psi}) = \big(f, g \cdot \ang{\psi}\big) \cn \ang{\ord{m}^\al} \to \ang{\ord{n}^\be},
\end{equation}
where 
\[g \cdot \ang{\psi} = \ang{g\psi_i \ginv}_{i \in \ufs{q}} \cn f_* \ang{\ord{m}^\al} \to \ang{\ord{n}^\be}\]
is the conjugation $g$-action on $\ang{\psi}$ defined in \cref{gangpsi}.
\end{itemize}
A morphism in $\GG$ is also denoted by a generic symbol, such as $\upom$.
\item[Identities]
For an object $\angordmal$ of length $p \geq 0$, the identity morphism
\[\Big(1_{\ufs{p}}, \ang{1_{\ordm_k^{\al_k}}}_{k \in \ufs{p}}\Big)\]
consists of the identity function on $\ufs{p}$ and the identity morphism on $\ordm_k^{\al_k}$ for each $k \in \ufs{p}$.
\item[Composition] 
Consider composable morphisms in $\GG$
\[\angordmal \fto{(f,\ang{\psi})} \angordnbe \fto{(h,\ang{\phi})} \angordellde\]
for objects $\angordmal \in \FGsma{p}$, $\angordnbe \in \FGsma{q}$, and $\angordellde \in \FGsma{r}$.
\begin{itemize}
\item If $q=0$, then $p=0$, and $(f,\ang{\psi})$ is either the identity morphism $1_{\ang{}}$ or the 0-morphism.  Their composites with $(h,\ang{\phi})$ are, respectively, $(h,\ang{\phi})$ and the 0-morphism.
\item If $q>0$, then $r>0$, and the composite is defined as the pair
\begin{equation}\label{GG_composite}
(h,\ang{\phi}) \circ (f,\ang{\psi}) = \big(hf, \ang{\phi} \circ h_*\ang{\psi} \big) \cn \angordmal \to \angordellde.
\end{equation}
\end{itemize}
\end{description}
If $G$ is the trivial group, then $\GG$ is identified with the pointed category $\Gsk$ \pcref{def:Gsk}.  \cref{GG_Gcategory} verifies that $\GG$ is a pointed $G$-category.  
\end{definition}

\begin{explanation}[$G$-Action]\label{expl:GG_Gaction}
The $G$-action on morphisms of $\GG$ \cref{GG_Gaction} has the following properties.
\begin{itemize}
\item The $G$-action does not change the reindexing injection.
\item Each 0-morphism, which factors through the initial-terminal basepoint $\vstar \in \GG$, is $G$-fixed.
\item The $G$-action fixes each morphism whose domain and codomain consist of only trivial $G$-sets.\defmark
\end{itemize}
\end{explanation}

\begin{explanation}[Composition]\label{expl:GG_composite}
On the right-hand side of the equality in \cref{GG_composite}, the first entry is the composite injection
\[\ufs{p} \fto{f} \ufs{q} \fto{h} \ufs{r}.\]
The second entry is the $r$-tuple
\begin{equation}\label{GG_comp}
\ang{\phi} \circ h_*\ang{\psi} = 
\Bang{\ord{m}_{(hf)^{-1}(j)}^{\al_{(hf)^{-1}(j)}} \fto{\psi_{\hinv(j)}} \ord{n}_{\hinv(j)}^{\be_{\hinv(j)}} \fto{\phi_j} \ord{\ell}_j^{\de_j} }_{j \in \ufs{r}}
\end{equation}
of pointed functions.
\end{explanation}

\begin{lemma}\label{GG_Gcategory}
$\GG$ in \cref{def:GG} is a small pointed $G$-category .
\end{lemma}

\begin{proof}
The proof given in \cite[\namecref{EqK:Gsk_Gcategory} \ref*{EqK:Gsk_Gcategory}]{yau-eqk} that $\Gsk$ is a pointed category also applies to $\GG$.  Next, we check that composition is $G$-equivariant.  By \cref{GG_Gaction}, the $G$-action does not change the reindexing injection, so we concentrate on the second entry.  Using \cref{GG_comp}, the following computation for $g \in G$ proves that composition is $G$-equivariant.
\begin{equation}\label{GG_comp_geq}
\begin{split}
\big(g \cdot \ang{\phi}\big) \circ h_*\big(g \cdot \ang{\psi}\big)
&= \bang{g \cdot \phi_j}_{j \in \ufs{r}} \circ h_*\bang{g \cdot \psi_i}_{i \in \ufs{q}}\\
&= \bang{\big(g \phi_j \ginv) (g \psi_{\hinv(j)} \ginv) }_{j \in \ufs{r}}\\
&= \bang{ g \phi_j \psi_{\hinv(j)} \ginv}_{j \in \ufs{r}}\\
&= g \cdot \bang{\phi_j \psi_{\hinv(j)}}_{j \in \ufs{r}}\\
&= g \cdot \big( \ang{\phi} \circ h_* \ang{\psi}\big)
\end{split}
\end{equation}
This proves that $\GG$ is a small pointed $G$-category.
\end{proof}

\cref{def:ifgGG} defines the $G$-equivariant analogue of the inclusion $\Fsk \to \Gsk$ \pcref{def:ifg}.

\begin{definition}[Length-1 Inclusion]\label{def:ifgGG}
The pointed full $G$-subcategory inclusion\index{length-1 inclusion}
\begin{equation}\label{ifgGG}
\FG \fto{\ifg} \GG
\end{equation}
sends a pointed finite $G$-set $\ordn^\be \in \FG$ to the length-1 object of $\GG$ consisting of $\ordn^\be$.  This is the initial-terminal basepoint $\vstar \in \GG$ if $\nbeta = \ordz$.  We usually abbreviate $\ifg\ordn^\be$ to $\ordn^\be$.  A pointed function $\psi \cn \ordm^\al \to \ordn^\be$ in $\FG$ is sent by $\ifg$ to the morphism
\[\ordm^\al \fto{(1_{\ufs{1}}, \psi)} \ordn^\be \inspace \GG.\]
This is the unique morphism from or to the basepoint $\ifg\ordz = \vstar \in \GG$ if either $\ordm^\al$ or $\ordn^\be$ is $\ordz$.
\end{definition}

\subsection*{Naive Permutative $G$-Category}

\cref{def:GG_permutative} extends the permutative category structure on $\Gsk$ \pcref{def:Gsk_permutative} to a naive permutative $G$-category structure (\cref{expl:naive_perm_Gcat}) on $\GG$.

\begin{definition}\label{def:GG_permutative}
The small pointed $G$-category $\GG$ \pcref{def:GG} is given a naive permutative $G$-category structure $(\oplus, \ang{}, \xi)$ as follows.
\begin{description}
\item[Monoidal product on objects] 
The monoidal product 
\[\GG \times \GG \fto{\oplus} \GG\]
is defined on objects by concatenation
\begin{equation}\label{GG_oplus_obj}
\angordmal \oplus \angordnbe = 
\big(\ord{m}_1^{\al_1}, \ldots, \ord{m}_p^{\al_p}, \ord{n}_1^{\be_1}, \ldots, \ord{n}_q^{\be_q} \big)
\end{equation}
for $\angordmal \in \FGsma{p}$ and $\angordnbe \in \FGsma{q}$.  If any $\ord{m}_i$ or $\ord{n}_j$ is $\ord{0}$, then the right-hand side of \cref{GG_oplus_obj} has at least one entry of $\ord{0}$.  Thus, we have
\begin{equation}\label{GG_oplus_vstar}
\vstar \oplus \angordnbe = \vstar = \angordmal \oplus \vstar.
\end{equation}
\item[Monoidal product on morphisms]
By \cref{GG_oplus_vstar} and the fact that the basepoint $\vstar$ is initial and terminal in $\GG$, $\oplus$ is uniquely defined when one morphism has either domain or codomain $\vstar$.

To define $\oplus$ for other morphisms, we consider morphisms
\[\sord{m}{\al} \fto{(f,\ang{\psi})} \sord{n}{\be} \andspace 
\sord{\jmath}{\chi} \fto{(h,\ang{\phi})} \sord{\ell}{\de}\]
in $\GG$ for objects $\sord{m}{\al} \in \FGsma{p}$, $\sord{n}{\be} \in \FGsma{q}$, $\sord{\jmath}{\chi} \in \FGsma{r}$, and $\sord{\ell}{\de} \in \FGsma{s}$.  Then we define the morphism
\begin{equation}\label{GG_oplus_morphism}
\begin{split}
&\big(f, \ang{\psi}\big) \oplus \big(h, \ang{\phi}\big) \\
&= \big(f \oplus h, \ang{\psi} \oplus \ang{\phi}\big) \cn
\sord{m}{\al} \oplus \sord{\jmath}{\chi} \to \sord{n}{\be} \oplus \sord{\ell}{\de}
\end{split}
\end{equation}
with reindexing injection
\[\ufs{p+r} \fto{f \oplus h} \ufs{q+s}\]
defined by
\[(f \oplus h)(i) = \begin{cases} f(i) & \text{if $1 \leq i \leq p$ and}\\
q + h(i-p) & \text{if $p+1 \leq i \leq p+r$}.
\end{cases}\]
The morphism $\ang{\psi} \oplus \ang{\phi}$ is the concatenation of the $q$-tuple $\ang{\psi}$ and the $s$-tuple $\ang{\phi}$, as displayed in the next diagram.
\begin{equation}\label{psiplusphi}
\begin{tikzpicture}[vcenter]
\draw[0cell=.9]
(0,0) node (a) {(f \oplus h)_*\big(\sord{m}{\al} \oplus \sord{\jmath}{\chi} \big)}
(a)++(0,-1) node (b) {f_* \sord{m}{\al} \oplus h_* \sord{\jmath}{\chi}}
(b)++(5.5,0) node (c) {\sord{n}{\be} \oplus \sord{\ell}{\de}}
(c)++(0,1) node (d) {\sord{n}{\be} \oplus \sord{\ell}{\de}}
;
\draw[1cell=.9]
(a) edge[equal] (b)
(c) edge[equal] (d)
(a) edge node {\ang{\psi} \oplus \ang{\phi}} (d)
(b) edge node {\big(\ang{\psi_i}_{i \in \ufs{q}} \scs \ang{\phi_k}_{k \in \ufs{s}} \big)} (c)
;
\end{tikzpicture}
\end{equation}
\item[Monoidal unit]
The monoidal unit is the empty tuple $\ang{} \in \GG$.
\item[Braiding] 
The component of the braiding $\xi$ at a pair of objects $(\sord{m}{\al}, \sord{n}{\be})$ is $1_\vstar$ if either $\sord{m}{\al}$ or $\sord{n}{\be}$ is the basepoint $\vstar \in \GG$.  For nonbasepoint objects, it is the isomorphism
\begin{equation}\label{GG_braiding}
\xi_{(\sord{m}{\al}, \sord{n}{\be})} = \big(\tau_{p,q} \, , \ang{1} \big)
\cn \sord{m}{\al} \oplus \sord{n}{\be} \fto{\iso} \sord{n}{\be} \oplus \sord{m}{\al}
\end{equation}
with reindexing bijection
\[\ufs{p+q} \fto[\iso]{\tau_{p,q}} \ufs{q+p}\]
given by the block permutation that interchanges the first $p$ elements with the last $q$ elements:
\[\tau_{p,q}(i) = \begin{cases} 
q+i & \text{if $1 \leq i \leq p$ and}\\
i-p & \text{if $p+1 \leq i \leq p+q$}.
\end{cases}\]
The morphism $\ang{1}$ in \cref{GG_braiding} is the $(q+p)$-tuple with each entry given by an identity morphism of some $\ordn_i^{\be_i}$ or some $\ordm_k^{\al_k}$.
\end{description}
\cref{GG_permutative} verifies that this triple is a naive permutative $G$-category structure on $\GG$.
\end{definition}

\begin{lemma}\label{GG_permutative}
The quadruple $(\GG,\oplus,\ang{},\xi)$ in \cref{def:GG,def:GG_permutative} is a naive permutative $G$-category.
\end{lemma}

\begin{proof}
The proof given in \cite[\namecref{EqK:Gsk_permutative} \ref*{EqK:Gsk_permutative}]{yau-eqk} that the quadruple $(\Gsk,\oplus,\ang{},\xi)$ is a permutative category also applies to $(\GG,\oplus,\ang{},\xi)$.

\parhead{$G$-equivariance}.  To check that the structure $(\oplus,\ang{},\xi)$ is $G$-equivariant, first observe that the monoidal unit $\ang{}$ is $G$-fixed, since each object of $\GG$ is $G$-fixed.  The functor $\oplus$ is $G$-equivariant by the following computation for each $g \in G$, using \cref{GG_Gaction,GG_oplus_morphism}.
\begin{equation}\label{GG_mor_geq}
\begin{split}
& g \cdot \big(f \oplus h, \ang{\psi} \oplus \ang{\phi}\big) \\
&= \big(f \oplus h, g\cdot (\ang{\psi} \oplus \ang{\phi})\big)\\
&= \big(f \oplus h, (g \cdot \ang{\psi}) \oplus (g \cdot \ang{\phi})\big)\\
&= \big(f, g \cdot \ang{\psi}\big) \oplus \big(h, g \cdot \ang{\phi}\big)
\end{split}
\end{equation}
The braiding $\xi$ is $G$-equivariant by the following computation, using \cref{GG_Gaction,GG_braiding}.
\[\begin{split}
g \cdot \xi_{(\sord{m}{\al}, \sord{n}{\be})} 
&= g \cdot \big(\tau_{p,q} \spc \ang{1} \big)\\
&= \big(\tau_{p,q} \spc \ang{g \cdot 1} \big)\\
&= \big(\tau_{p,q} \spc \ang{1} \big)\\
& = \xi_{(\sord{m}{\al}, \sord{n}{\be})}
\end{split}\]
This proves that the permutative structure $(\oplus, \ang{}, \xi)$ is $G$-equivariant.
\end{proof}

\subsection*{Comparing $\GG$ and $\FG$}

Next, we compare the naive permutative $G$-categories $\FG$ and $\GG$ in \cref{FG_permutative,GG_permutative} via the functor in \cref{def:smashFGGG}, which is the $G$-equivariant analogue of \cref{def:smashFskGsk}.

\begin{definition}\label{def:smashFGGG}
We define a functor
\[\GG \fto{\sma} \FG\]
as follows.
\begin{description}
\item[Objects] 
The object assignment of $\sma$ is defined as follows for $\ang{\ord{m}_k^{\al_k}}_{k \in \ufs{p}} \in \GG \setminus \{\vstar,\ang{}\}$.  
\begin{equation}\label{smash_GGobjects}
\left\{
\begin{split}
\sma \vstar &= \ord{0}\\
\sma \ang{} &= \ord{1}\\
\sma \ang{\ord{m}_k^{\al_k}}_{k \in \ufs{p}} &= \sma_{k \in \ufs{p}} \,\ord{m}_k^{\al_k} 
= \ord{m_1 \cdots m_p}^{\sma_{k \in \ufs{p}}\, \al_k}
\end{split}\right.
\end{equation}
Here $\sma_{k \in \ufs{p}}$ is the $p$-fold iterate of the monoidal product of $\FG$ \cref{mal-sma-nbe}.  This is well defined because if some $\ord{m}_k = \ord{0}$, then $\ord{m_1 \cdots m_p} = \ord{0}$.
\item[Morphisms]
The two morphisms in $\GG\big(\ang{}, \ang{}\big)$ are the identity morphism and the 0-morphism \cref{GG_empty_mor}.  They are sent by $\sma$ to, respectively, the identity morphism and the 0-morphism in $\FG(\ord{1}, \ord{1})$.

For a morphism \cref{fangpsiGG}
\[\sord{m}{\al} \fto{(f, \ang{\psi})} \sord{n}{\be} \inspace \GG,\]
the morphism
\[\smad\sord{m}{\al} \fto{\smad (f, \ang{\psi})} \smad\sord{n}{\be} \inspace \FG\]
is defined as the following composite pointed function.
\begin{equation}\label{smash_fpsiFG}
\begin{tikzpicture}[vcenter]
\def\t{.7ex}
\draw[0cell=.9]
(0,0) node (a) {\bigwedge_{k \in \ufs{p}} \ord{m}_k^{\al_k}}
(a)++(2.8,0) node (b) {\bigwedge_{\finv(i) \neq \emptyset} \ord{m}_{\finv(i)}^{\al_{\finv(i)}}}
(b)++(3,0) node (c) {\bigwedge_{i \in \ufs{q}} \ord{m}_{\finv(i)}^{\al_{\finv(i)}}}
(c)++(2.5,0) node (d) {\bigwedge_{i \in \ufs{q}} \ord{n}_i^{\be_i}}
;
\draw[1cell=.9]
(a) edge[transform canvas={yshift=\t}] node {f_*} node[swap] {\iso} (b)
(b) edge[transform canvas={yshift=\t}] node {\iso} (c)
(c) edge[transform canvas={yshift=\t}] node {\sma_i\, \psi_i} (d)
;
\end{tikzpicture}
\end{equation}
The three pointed functions in \cref{smash_fpsiFG} are defined as follows. 
\begin{itemize}
\item $f_*$ permutes the $p$ entries according to the reindexing injection $f \cn \ufs{p} \to \ufs{q}$.  The indexing set in the codomain is given by $\big\{ i \in \ufs{q} \cn \finv(i) \neq \emptyset\}$.
\item Using \cref{ordn_emptyFG}, the middle pointed bijection in \cref{smash_fpsiFG} inserts a copy of the smash unit 
\[\ord{1} = \ord{m}_{\emptyset}^{\al_{\emptyset}}\]
for each index $i \in \ufs{q}$ not in the image of $f$.
\item $\sma_i\, \psi_i$ is the smash product of the pointed functions $\psi_i$ for $i \in \ufs{q}$.
\end{itemize}
This yields a well-defined morphism $\smad(f,\ang{\psi})$ because, if any $\psi_i$ factors through $\ord{0}$, then the composite in \cref{smash_fpsiFG} also factors through $\ord{0}$.
\end{description}
\cref{sma_symmon} proves that $\sma$ is a strict symmetric monoidal pointed $G$-functor.
\end{definition}

\begin{explanation}\label{expl:smash_fpsiFG}
In \cref{smash_fpsiFG}, the domain $\sma_{k \in \ufs{p}}\, \ord{m}_k^{\al_k}$ is the smash product defined in \cref{mal-sma-nbe}.  It uses the lexicographic ordering \cref{lex_bijection} to make the identification
\[\ord{m_1 \cdots m_p}^{\sma_{k \in \ufs{p}}\, \al_k} = \bigwedge_{k \in \ufs{p}} \ord{m}_k^{\al_k}.\]
Using this identification, the image of an element $\ang{a_k}_{k \in \ufs{p}} \in \sma_{k \in \ufs{p}}\, \ord{m_k}^{\al_k}$ under the composite \cref{smash_fpsiFG} is
\begin{equation}\label{smafangpsiangs}
\smad(f,\ang{\psi}) \ang{a_k}_{k \in \ufs{p}} 
= \bang{\psi_i a_{\finv(i)}}_{i \in \ufs{q}} \in \bigwedge_{i \in \ufs{q}} \ord{n}_i^{\be_i},
\end{equation}
where
\[a_{\emptyset} = 1 \in \ord{1} = \ord{m}_\emptyset^{\al_\emptyset}\]
for each index $i \in \ufs{q}$ not in the image of $f$.
\end{explanation}

\begin{lemma}\label{sma_symmon}
The assignments in \cref{def:smashFGGG} define a strict symmetric monoidal pointed $G$-functor
\[\big(\GG,\oplus,\ang{},\xi\big) \fto{\sma} \big(\FG,\sma,\ord{1},\xi\big).\]
\end{lemma}

\begin{proof}
The proof given in \cite[\namecref{EqK:sma_symmon} \ref*{EqK:sma_symmon}]{yau-eqk} that $\sma \cn \Gsk \to \Fsk$ in \cref{def:smashFskGsk} is a strict symmetric monoidal pointed functor also applies to $\sma \cn \GG \to \FG$.

\parhead{$G$-equivariance}.  The functor $\sma \cn \GG \to \FG$ is $G$-equivariant on objects because $G$ acts trivially on the objects of $\GG$ and $\FG$.  By \cref{gpsi,GG_Gaction}, the $G$-equivariance of $\sma$ on morphisms means that, for each $g \in G$ and morphism $(f,\ang{\psi})$ \cref{fangpsiGG}, the morphism equality 
\begin{equation}\label{sma-fgpsi}
\smad \big(f, \ang{g \cdot \psi_i}_{i \in \ufs{q}} \big) 
= g \cdot \smad(f, \ang{\psi}) \cn 
\smad\sord{m}{\al} \to \smad\sord{n}{\be}
\end{equation}
holds in $\FG$.  By \cref{al-sma-be,smash_fpsiFG}, the left-hand and right-hand sides of \cref{sma-fgpsi} are, respectively, the top and left-bottom-right boundary composites in the diagram \cref{sma-fgpsi-diagram}.
\begin{equation}\label{sma-fgpsi-diagram}
\begin{tikzpicture}[vcenter]
\def\t{.7ex}
\def\v{-1.75}
\draw[0cell=.85]
(0,0) node (a) {\bigwedge_{k \in \ufs{p}} \ord{m}_k^{\al_k}}
(a)++(3,0) node (b) {\bigwedge_{\finv(i) \neq \emptyset} \ord{m}_{\finv(i)}^{\al_{\finv(i)}}}
(b)++(3,0) node (c) {\bigwedge_{i \in \ufs{q}} \ord{m}_{\finv(i)}^{\al_{\finv(i)}}}
(c)++(3,0) node (d) {\bigwedge_{i \in \ufs{q}} \ord{n}_i^{\be_i}}
(a)++(0,\v) node (a2) {\bigwedge_{k \in \ufs{p}} \ord{m}_k^{\al_k}}
(b)++(0,\v) node (b2) {\bigwedge_{\finv(i) \neq \emptyset} \ord{m}_{\finv(i)}^{\al_{\finv(i)}}}
(c)++(0,\v) node (c2) {\bigwedge_{i \in \ufs{q}} \ord{m}_{\finv(i)}^{\al_{\finv(i)}}}
(d)++(0,\v) node (d2) {\bigwedge_{i \in \ufs{q}} \ord{n}_i^{\be_i}}
;
\draw[1cell=.85]
(a) edge[transform canvas={yshift=\t}] node {f_*} node[swap] {\iso} (b)
(b) edge[transform canvas={yshift=\t}] node {\iso} (c)
(c) edge[transform canvas={yshift=\t}] node {\sma_i\, g\psi_i \ginv} (d)
(a2) edge[transform canvas={yshift=\t}] node {f_*} node[swap] {\iso} (b2)
(b2) edge[transform canvas={yshift=\t}] node {\iso} (c2)
(c2) edge[transform canvas={yshift=\t}] node {\sma_i\, \psi_i} (d2)
(a) edge node[swap] {\sma_k\, \ginv} (a2)
(b) edge node[swap] {\sma \ginv} (b2)
(c) edge[transform canvas={xshift=-1ex}] node[swap] {\sma_i\, \ginv} (c2)
(d2) edge[transform canvas={xshift=-1ex}] node[swap] {\sma_i\, g} (d) 
;
\end{tikzpicture}
\end{equation}
The three regions in the diagram \cref{sma-fgpsi-diagram} commute for the following reasons.
\begin{itemize}
\item The left region commutes by the naturality of the braiding for the smash product.
\item The middle region commutes because $G$ acts trivially on $\ord{1} = \ord{m}_\emptyset^{\alpha_\emptyset}$.
\item The right region commutes by the functoriality of $\sma \cn \FG \times \FG \to \FG$.
\end{itemize}  
This proves that the strict symmetric monoidal pointed functor $\sma \cn \GG \to \FG$ is $G$-equivariant.
\end{proof}

\section{The 2-Category of $\GGG$-Categories}
\label{sec:ggcatg}

This section constructs the 2-category $\GGCatg$ and a 2-functor $\igst$ from it to the 2-category $\GGCatii$.  \cref{thm:ggcat_ggcatg_iieq} proves that $\igst$ is a 2-equivalence.

\secoutline
\begin{itemize}
\item \cref{def:ggcatg} defines the 2-category $\GGCatg$, whose objects, called $\GGG$-categories, are pointed $G$-functors $\GG \to \Catgst$. 
\item \cref{expl:GGcategory} unpacks the 2-category $\GGCatg$.
\item \cref{def:gskel_gg} defines the full subcategory inclusion $\ig \cn \Gsk \to \GG$.
\item \cref{igst_iifunctor} constructs the induced 2-functor $\igst \cn \GGCatg \to \GGCatii$.
\end{itemize}

\subsection*{$\GGG$-Categories}

\cref{def:ggcatg} extends the 2-category $\GGCatii$ \pcref{def:GGCat} by replacing $\Gsk$ and $\Gcatst$ with the pointed $G$-category $\GG$ \pcref{def:GG} and $\Catgst$ \cref{Catgst_iicat}.  Basic concepts of 2-category theory are reviewed in \cref{sec:twocategories}.

\begin{definition}\label{def:ggcatg}
For a group $G$, the 2-category $\GGCatg$ is defined as follows.
\begin{description}
\item[Objects] An object in $\GGCatg$, called a \index{GGG-category@$\GGG$-category}\emph{$\GGG$-category}, is a pointed $G$-functor
\begin{equation}\label{ggcatg_obj}
(\GG,\vstar) \fto{X} (\Catgst,\bone).
\end{equation}
\item[1-cells] A 1-cell $\tha \cn X \to X'$ in $\GGCatg$ is a $G$-natural transformation as follows.
\begin{equation}\label{ggcatg_icell}
\begin{tikzpicture}[vcenter]
\def\t{28}
\draw[0cell]
(0,0) node (a1) {\phantom{\Gsk}}
(a1)++(1.8,0) node (a2) {\phantom{\Gsk}}
(a1)++(-.08,0) node (a1') {\GG}
(a2)++(.2,0) node (a2') {\Catgst}
;
\draw[1cell=.9]
(a1) edge[bend left=\t] node {X} (a2)
(a1) edge[bend right=\t] node[swap] {X'} (a2)
;
\draw[2cell]
node[between=a1 and a2 at .45, rotate=-90, 2label={above,\theta}] {\Rightarrow}
;
\end{tikzpicture}
\end{equation}
\item[2-cells] A 2-cell $\Theta \cn \tha \to \ups$ in $\GGCatg$ is a $G$-modification as follows.
\begin{equation}\label{ggcatg_iicell}
\begin{tikzpicture}[vcenter]
\def\t{25}
\draw[0cell]
(0,0) node (a) {\GG}
(a)++(3,0) node (b) {\phantom{\GG}}
(b)++(.15,0) node (b') {\Catgst}
;
\draw[1cell=.8]
(a) edge[bend left=\t] node {X} (b)
(a) edge[bend right=\t] node[swap] {X'} (b)
;
\draw[2cell=.9]
node[between=a and b at .32, rotate=-90, 2label={below,\theta}] {\Rightarrow}
node[between=a and b at .65, rotate=-90, 2label={above,\ups}] {\Rightarrow}
;
\draw[2cell]
node[between=a and b at .5, rotate=0, shift={(0,-.2)}, 2labelmed={above,\Theta}] {\Rrightarrow}
;
\end{tikzpicture}
\end{equation}
\item[Other structures] Identity 1-cells and 2-cells, vertical composition of 2-cells, and horizontal composition of 1-cells and 2-cells are defined componentwise in the 2-category $\Catgst$.
\end{description}
The underlying 1-category of $\GGCatg$ is denoted by the same notation.
\end{definition}

\begin{explanation}[Unpacking $\GGCatg$]\label{expl:GGcategory}
The 2-category $\GGCatg$ in \cref{def:ggcatg} is given explicitly as follows.
\begin{description}
\item[Objects]
A $\GGG$-category $X \cn \GG \to \Catgst$ \cref{ggcatg_obj} consists of the following data.
\begin{itemize}
\item $X$ sends each object $\angordmal \in \GG$ \cref{GG_objects} to a small pointed $G$-category $X\angordmal$ such that $X\vstar = \bone$.  Its $G$-fixed basepoint is given by the pointed functor
\[X(\vstar \to \angordmal) \cn X\vstar = \bone \to X\angordmal\]
for the $G$-fixed unique morphism $\vstar \to \angordmal$ in $\GG$.  
\item $X$ sends each morphism $\upom \cn \angordmal \to \angordnbe$ in $\GG$ \cref{GG_morphisms} to a pointed functor
\begin{equation}\label{X_upom_GG}
X\angordmal \fto{X\upom} X\angordnbe
\end{equation}
such that $X$ preserves identity morphisms and composition.  The functor $X\upom$ is \emph{not} generally $G$-equivariant, in contrast to $\Gskg$-categories \cref{f_upom}.
\item The $G$-equivariance of $X$ means the equality of functors
\begin{equation}\label{GGGcat_Gequiv}
X(g \cdot \upom) = g (X\upom) \ginv
\end{equation}
for each $g \in G$ and each morphism $\upom$ in $\GG$.  The morphism $g \cdot \upom$ is defined in \cref{GG_Gaction}, and $g (X\upom) \ginv$ is the conjugation $g$-action \cref{conjugation-gaction} on $X\upom$.  In particular, if the morphism $\upom$ is $G$-fixed, then $X\upom$ is a pointed $G$-functor.
\end{itemize}
\item[1-cells]
A 1-cell $\tha \cn X \to X'$ in $\GGCatg$ \cref{ggcatg_icell} consists of, for each object $\angordmal \in \GG$, an $\angordmal$-component pointed functor
\begin{equation}\label{ggcatg_icell_comp}
X\angordmal \fto{\tha_{\angordmal}} X'\angordmal
\end{equation}
such that, for each morphism $\upom \cn \angordmal \to \angordnbe$ in $\GG$, the following naturality diagram of pointed functors commutes.
\begin{equation}\label{ggcatg_icell_nat}
\begin{tikzpicture}[vcenter]
\def\v{-1.4}
\draw[0cell]
(0,0) node (a11) {X\angordmal}
(a11)++(2.75,0) node (a12) {X'\angordmal}
(a11)++(0,\v) node (a21) {X\angordnbe}
(a12)++(0,\v) node (a22) {X'\angordnbe}
;
\draw[1cell=.9]
(a11) edge node {\theta_{\angordmal}} (a12)
(a12) edge node {X'\upom} (a22)
(a11) edge node[swap] {X\upom} (a21)
(a21) edge node {\theta_{\angordnbe}} (a22)
;
\end{tikzpicture}
\end{equation}
Since $G$ acts trivially on the objects of $\GG$ \cref{GG_obj_gaction}, the $G$-equivariance of $\tha$ means the equality of functors
\begin{equation}\label{ggcatg_icell_geq}
\tha_{\angordmal} = g \theta_{\angordmal} \ginv
\end{equation}
for each $g \in G$ and object $\angordmal \in \GG$.   In other words, each component of a 1-cell in $\GGCatg$ is a pointed $G$-functor, similar to 1-cells in $\GGCatii$ \cref{ggcat_mor_component}.  

Identity 1-cells and horizontal composition of 1-cells in $\GGCatg$ are defined componentwise using the components in \cref{ggcatg_icell_comp}.  A 1-cell is automatically pointed, meaning $\tha_{\vstar} = 1_{\bone}$.
\item[2-cells]
A 2-cell $\Theta \cn \theta \to \ups$ in $\GGCatg$ \cref{ggcatg_iicell} consists of, for each object $\angordmal \in \GG$, an $\angordmal$-component pointed natural transformation
\begin{equation}\label{ggcatg_iicell_comp}
\begin{tikzpicture}[vcenter]
\def\t{25}
\draw[0cell]
(0,0) node (a1) {\phantom{X'}}
(a1)++(2.5,0) node (a2) {\phantom{X'}}
(a1)++(-.25,0) node (a1') {X\angordmal}
(a2)++(.33,0) node (a2') {X'\angordmal}
;
\draw[1cell=.85]
(a1) edge[bend left=\t] node {\theta_{\angordmal}} (a2)
(a1) edge[bend right=\t] node[swap] {\ups_{\angordmal}} (a2)
;
\draw[2cell]
node[between=a1 and a2 at .37, rotate=-90, 2label={above, \Theta_{\angordmal}}] {\Rightarrow}
;
\end{tikzpicture}
\end{equation}
such that, for each morphism $\upom \cn \angordmal \to \angordnbe$ in $\GG$, the following two whiskered natural transformations are equal.
\begin{equation}\label{ggcatg_iicell_modax}
\begin{tikzpicture}[vcenter]
\def\t{25} \def\v{-1.6}
\draw[0cell]
(0,0) node (a1) {\phantom{X'}}
(a1)++(2.5,0) node (a2) {\phantom{X'}}
(a1)++(-.25,0) node (a1') {X\angordmal}
(a2)++(.3,0) node (a2') {X'\angordmal}
(a1)++(0,\v) node (b1) {\phantom{X'}}
(a2)++(0,\v) node (b2) {\phantom{X'}}
(b1)++(-.25,0) node (b1') {X\angordnbe}
(b2)++(.3,0) node (b2') {X'\angordnbe}
;
\draw[1cell=.8]
(a1) edge[bend left=\t] node[pos=.4] {\theta_{\angordmal}} (a2)
(a1) edge[bend right=\t] node[swap,pos=.6] {\ups_{\angordmal}} (a2)
(b1) edge[bend left=\t] node[pos=.4] {\theta_{\angordnbe}} (b2)
(b1) edge[bend right=\t] node[swap,pos=.6] {\ups_{\angordnbe}} (b2)
(a1) edge[transform canvas={xshift=-.7ex}] node[swap] {X\upom} (b1)
(a2) edge node {X'\upom} (b2)
;
\draw[2cell=.9]
node[between=a1 and a2 at .37, rotate=-90, 2label={above, \Theta_{\angordmal}}] {\Rightarrow}
node[between=b1 and b2 at .37, rotate=-90, 2label={above, \Theta_{\angordnbe}}] {\Rightarrow}
;
\end{tikzpicture}
\end{equation}
Since $G$ acts trivially on the objects of $\GG$ \cref{GG_obj_gaction}, the $G$-equivariance of $\Theta$ means the equality of natural transformations
\begin{equation}\label{ggcatg_iicell_geq}
\Theta_{\angordmal} = g * \Theta_{\angordmal} * \ginv
\end{equation}
for each $g \in G$ and object $\angordmal \in \GG$, where $*$ denotes horizontal composition of natural transformations.   In other words, each component of a 2-cell in $\GGCatg$ is a pointed $G$-natural transformation, similar to 2-cells in $\GGCatii$ \cref{ggcat_inthom_Theta}.  

Identities, horizontal composition, and vertical composition of 2-cells are given componentwise using the components in \cref{ggcatg_iicell_comp}.   A 2-cell is automatically pointed, meaning $\Theta_\vstar = 1_{1_{\bone}}$.\defmark
\end{description}
\end{explanation}

\subsection*{$\Gskg$-Categories from $\GGG$-Categories}

\cref{def:gskel_gg} compares the indexing categories $\Gsk$ and $\GG$ \pcref{def:Gsk,def:GG}.

\begin{definition}[From $\Gsk$ to $\GG$]\label{def:gskel_gg}
For a group $G$, define the pointed full subcategory inclusion
\begin{equation}\label{ig}
\Gsk \fto{\ig} \GG
\end{equation}
by sending each object $\angordn = \ang{\ordn_i}_{i \in \ufs{q}} \in \Gsk$ to the object 
\begin{equation}\label{ig_angordn}
\ig\angordn = \ang{\ordn_i^{\eps_i}}_{i \in \ufs{q}} \in \GG,
\end{equation}
where each $\eps_i$ is the trivial $G$-action on $\ordn_i$.  The functor $\ig$ is well defined by \cref{Gsk_morphisms}, \cref{GG_morphisms}, and the fact that morphisms in $\FG$ \pcref{def:FG} are pointed functions between underlying pointed finite sets.  If there is no danger of confusion, then $\ig\angordn$ is abbreviated to $\angordn$ and similarly for morphisms.
\end{definition}

Using the full subcategory inclusion $\ig$, \cref{igst_iifunctor} shows that $\GGG$-categories yield $\Gskg$-categories and similarly for 1-cells and 2-cells \pcref{def:GGCat,def:ggcatg}.

\begin{lemma}\label{igst_iifunctor}
Precomposition with the functor $\ig \cn \Gsk \to \GG$ in \cref{ig} induces a 2-functor
\[\GGCatg \fto{\igst} \GGCatii.\]
\end{lemma}

\begin{proof}
We first verify that $\igst$ is well defined on objects, 1-cells, and 2-cells.  
\begin{description}
\item[Objects]  
To see that $\igst$ is well defined on objects, we consider a pointed $G$-functor \cref{ggcatg_obj}
\[(\GG,\vstar) \fto{X} (\Catgst,\boldone).\]
The pointed functor \cref{ggcat_obj}
\[(\Gsk,\vstar) \fto{\igst X} (\Gcatst,\boldone)\]
is defined on objects by
\[(\igst X)\angordn = X(\ig\angordn) = X\angordn \forspace \angordn \in \Gsk\]
and similarly for morphisms, where $\ig\angordn$ is defined in \cref{ig_angordn}.  For each morphism $\upom$ in $\Gsk$, the pointed functor 
\[(\igst X)\upom = X(\ig\upom)\] 
is $G$-equivariant by the $G$-equivariance of $X$ \cref{GGGcat_Gequiv} and the fact that $G$ acts trivially on $\ig\upom$ \pcref{expl:GG_Gaction}.  Identity morphisms and composition in both categories $\Catgst$ and $\Gcatst$ are defined in $\Cat$.  Thus, $\igst X$ is a well-defined $\Gskg$-category.
\item[1-cells]  
The 1-cell assignment of $\igst$ is well defined because each component of each 1-cell in $\GGCatg$ is a pointed $G$-functor \cref{ggcatg_icell_geq}, as required for components of 1-cells in $\GGCatii$ \cref{ggcat_mor_component}.  The naturality of 1-cells in $\GGCatii$ \cref{ggcat_mor_naturality} is covered by the naturality of 1-cells in $\GGCatg$ \cref{ggcatg_icell_nat}.
\item[2-cells]  
Along the same lines, the 2-cell assignment of $\igst$ is well defined by the description of 2-cells in $\GGCatii$ given in \cref{ggcat_inthom_Theta,ggcat_inthom_Theta_modax} and the description of 2-cells in $\GGCatg$ given in \cref{ggcatg_iicell_modax,ggcatg_iicell_geq}.
\item[2-functoriality]  
The assignment $\igst$ preserves the other 2-categorical structures---identity 1-cells and 2-cells, vertical composition, and horizontal composition---because these structures are defined componentwise in $\Catgst$ for $\GGCatg$ and in $\Gcatst$ for $\GGCatii$.
\end{description}
This proves that $\igst$ is a 2-functor.
\end{proof}

\section{2-Equivalence between $\Gskg$-Categories and $\GGG$-Categories}
\label{sec:ggcat_ggcatg_iieq}

This section first proves \cref{thm:ggcat_ggcatg_iieq}, which states that the 2-functor 
\[\GGCatg \fto{\igst} \GGCatii\]
in \cref{igst_iifunctor} is a 2-equivalence \pcref{def:twoequivalence}.
\cref{thm:ggcat_ggcatg_iieq} also holds if $\Gsk$ and $\GG$ are replaced by $\Fsk$ and $\FG$; see \cref{thm:fgcat_fgcatg_iieq}.  \cref{Xjininv} shows how, for a $\GGG$-category $X$ and an object $\nbe \in \GG$, the pointed $G$-category $X\nbe$ can be described in terms of $X\angordn$ with a twisted $G$-action.

\begin{theorem}\label{thm:ggcat_ggcatg_iieq}
For each group $G$, there is an adjoint 2-equivalence\index{adjoint 2-equivalence}
\begin{equation}\label{List_arrows}
\begin{tikzpicture}[vcenter]
\draw[0cell]
(0,0) node (a1) {\GGCatii}
(a1)++(2.5,0) node (a2) {\GGCatg}
;
\draw[1cell=.9]
(a1) edge[transform canvas={yshift=.5ex}] node {\Lg} (a2)
(a2) edge[transform canvas={yshift=-.4ex}] node {\igst} (a1)
;
\end{tikzpicture}
\end{equation}
between the 2-categories in \cref{def:GGCat,def:ggcatg}.
\end{theorem}

\begin{proof}
We first construct
\begin{itemize}
\item the inverse 2-functor $\Lg$ \cref{Lg} of $\igst$,
\item the unit 2-natural isomorphism $\ug$ \cref{ug}, and
\item the counit 2-natural isomorphism $\vg$ \cref{vg}.
\end{itemize}  
In \cref{Lgigst_lefttriangle,Lgigst_righttriangle}, we verify that the quadruple $(\Lg,\igst,\ug,\vg)$ satisfies the two triangle identities for a 2-adjunction \pcref{def:twoadjunction}. 
\begin{description}
\item[Inverse $\Lg$]  
The inverse 2-functor 
\begin{equation}\label{Lg}
\GGCatii \fto{\Lg} \GGCatg
\end{equation}
of $\igst$ sends a pointed functor $X \cn \Gsk \to \Gcatst$ to the pointed $G$-functor
\[(\GG,\vstar) \fto{\Lg X} (\Catgst,\bone)\]
whose value at an object $\nbe \in \GG$ is defined as the coend 
\begin{equation}\label{Lg_f_nbe}
(\Lg X)\nbe = 
\ecint^{\angordm \in \Gskel} \bigvee_{\GGpunc(\ig\angordm; \nbe)} X\angordm
\end{equation}
taken in $\Catst$.  In the wedge index in \cref{Lg_f_nbe}, $\GGpunc(-;-)$ denotes the set of nonzero morphisms in $\GG$ \cref{GG_morphisms}.  An empty wedge---which happens, for example, if $\angordm = \vstar \in \Gsk$ or $\nbe = \vstar \in \GG$---is defined as the terminal category $\boldone$.  The coend in \cref{Lg_f_nbe} is a quotient of the pointed category
\[\bigvee_{\angordm \in \Gsk} \bigvee_{\GGpunc(\ig\angordm; \nbe)} X\angordm.\]
Each object or morphism is represented by a pair
\begin{equation}\label{LX_reps}
\big(\ig\angordm \fto{\upom} \nbe ; x \big) \in \GG(\ig\angordm; \nbe) \ttimes X\angordm
\end{equation}
with $\angordm \in \Gsk$ and $\ig\angordm \in \GG$ \cref{ig_angordn}.  The pair $(\upom;x)$ represents the basepoint if either $x \in X\angordm$ is the basepoint or $\upom$ is the 0-morphism.  The defining relation of the coend $(\Lg X)\nbe$ \cref{Lg_f_nbe} identifies, for each triple
\[\begin{split}
& \big(\ig\angordm \fto{\upom} \nbe ; \angordl \fto{\upom'} \angordm ; x \big) \\
& \in \GG(\ig\angordm; \nbe) \ttimes \Gsk(\angordl; \angordm) \ttimes X\angordl,
\end{split}\]
the pairs
\begin{equation}\label{LX_relations}
\begin{split}
(\upom\upom'; x) & \in \GG(\ig\angordl; \nbe) \ttimes X\angordl \andspace\\
(\upom; (X\upom')(x)) & \in \GG(\ig\angordm; \nbe) \ttimes X\angordm.
\end{split}
\end{equation}
\begin{itemize}
\item The pointed functor $\Lg X$ is defined on morphisms of $\GG$ using the variable $\angordnbe$ in the wedge index $\GGpunc(\ig\angordm; \nbe)$ in \cref{Lg_f_nbe}.
\item The 1-cell and 2-cell assignments of $\Lg$ are defined componentwise using the term $X\angordm$ in the coend in \cref{Lg_f_nbe}.
\end{itemize}
\item[$G$-action] 
Using the $G$-action on $\GG$ \cref{GG_Gaction} and the trivial $G$-action on the entries of $\ig\angordm \in \GG$, the group $G$ acts diagonally on representatives \cref{LX_reps} of the pointed category $(\Lg X)\nbe$:
\begin{equation}\label{LXnbe_gaction}
g \cdot ((\gi, \angpsi) ; x) = ((\gi, \ang{g\psi}) ; gx)
\end{equation}
for $g \in G$, $(\gi, \angpsi) \in \GG(\ig\angordm; \nbe)$, and $x \in X\angordm$.  This $G$-action is well defined because, in the context of \cref{LX_relations}, $X\upom'$ is a $G$-functor \cref{f_upom} for each morphism $\upom'$ in $\Gsk$.
\item[$G$-equivariance of $\Lg X$]  
The $G$-equivariance of $\Lg X$ \cref{GGGcat_Gequiv} means that, for each $g \in G$ and each morphism
\[\lde \fto{(\gj,\angphi)} \nbe = \sordi{n}{\be}{j}_{j \in \ufs{q}} \inspace \GG,\]
there is an equality of functors
\begin{equation}\label{Lgf_gequiv}
\begin{split}
& (\Lg X)(\gj, \ang{g \phi\ginv}) \\
&= g \circ (\Lg X)(\gj,\angphi) \circ \ginv 
\cn (\Lg X)\lde \to (\Lg X)\nbe,
\end{split}
\end{equation}
where
\[\ang{g\phi\ginv} = 
\Bang{\ordi{\ell}{\de}{\gjinv(j)} \fto{\ginv} \ordi{\ell}{\de}{\gjinv(j)} \fto{\phi_j} \ordi{n}{\be}{j} \fto{g} \ordi{n}{\be}{j}}_{j \in \ufs{q}}.\]
The equality \cref{Lgf_gequiv} holds if either $\lde$ or $\nbe$ is the basepoint $\vstar \in \GG$, since $(\Lg X)\vstar = \boldone$.  

For nonbasepoint objects $\lde$ and $\nbe$, we consider a representing pair
\[((\gi, \angpsi) ; x) \in \GGpunc(\ig\angordm; \lde) \times X\angordm\]
in which $x \in X\angordm$ is either an object or a morphism.  The following computation proves that the images of $((\gi, \angpsi) ; x)$ under the two functors in \cref{Lgf_gequiv} are equal.
\[\begin{split}
& (\Lg X)(\gj, \ang{g \phi\ginv}) ((\gi, \angpsi) ; x) \\
&= \big((\gj, \ang{g \phi\ginv}) \circ (\gi, \angpsi) ; x \big) \\
&= \big( \big(\gj \gi, \ang{g \phi_j \ginv \psi_{\gjinv(j)}}_{j \in \ufs{q}} \big) ; x\big) \\
&= \big( \big(\gj \gi, \ang{g\phi_j}_{j \in \ufs{q}} \circ \gj_*\ang{\ginv\psi} \big) ; g\ginv x\big) \\
&= g \cdot \big( \big(\gj \gi, \angphi \circ \gj_*\ang{\ginv\psi} \big) ; \ginv x \big) \\
&= g \cdot \big( (\gj, \angphi) \circ (\gi, \ang{\ginv\psi}) ; \ginv x \big) \\
&= g \cdot \big((\Lg X)(\gj,\angphi) ((\gi, \ang{\ginv \psi}) ; \ginv x) \big) \\
&= (g \circ (\Lg X)(\gj,\angphi) \circ \ginv) ((\gi, \angpsi) ; x)
\end{split}\]
This proves that $\Lg X$ is a pointed $G$-functor.  The 2-functoriality of $\Lg$ \cref{Lg} follows from the fact that the 2-category structures of $\GGCatii$ and $\GGCatg$ are defined componentwise in, respectively, $\Gcatst$ and $\Catgst$.
\item[Unit]  
The unit 2-natural isomorphism
\begin{equation}\label{ug}

\end{equation}
sends a pointed functor $X \cn \Gsk \to \Gcatst$ to the natural isomorphism
\begin{equation}\label{ugx_iicell}
%
\end{equation}
whose value at an object $\angordm \in \Gsk$ is the pointed $G$-isomorphism
\[X\angordm \fto[\iso]{\ug_{X,\angordm}} (\igst\Lg X)\angordm = (\Lg X)(\ig\angordm)\]
defined as follows
\begin{itemize}
\item If $\angordm = \vstar \in \Gskel$, then $\ug_{X,\angordm} = 1_{\boldone} \cn \boldone \to \boldone$.
\item If $\angordm \neq \vstar$, then $\ug_{X,\angordm}$ is defined as the following composite.
\begin{equation}\label{ugf_angordm}
%
\end{equation}
\begin{itemize}
\item $\ug_{X,\angordm}^1$ is the inclusion into the wedge summand corresponding to the identity morphism of $\ig\angordm \in \GG$. 
\item $\ug_{X,\angordm}^2$ is the universal functor into the coend corresponding to the object $\angordt = \angordm \in \Gskel$. 
\end{itemize} 
\end{itemize}
Thus, the functor $\ug_{X,\angordm}$ \cref{ugf_angordm} sends an object or a morphism $x \in X\angordm$ to the representing pair
\begin{equation}\label{ugx_angordmx}
\ug_{X,\angordm}(x) = (1_{\angordm}; x) \in \GGpunc(\ig\angordm; \ig\angordm) \times X\angordm.
\end{equation}
Neither $\ug_{X,\angordm}^1$ nor $\ug_{X,\angordm}^2$ is an isomorphism in general, but the composite $\ug_{X,\angordm}$ is an isomorphism.  Its inverse is a version of the \emph{density isomorphism} \cite[III.3.7.8]{cerberusIII}.  The 2-naturality of $\ug$ follows from the description \cref{ugx_angordmx} of $\ug_{X,\angordm}(x)$.
\item[Counit]  
The counit 2-natural isomorphism
\begin{equation}\label{vg}
\begin{tikzpicture}[vcenter]
\def\t{20} \def\h{1.6} \def\v{1}
\draw[0cell]
(0,0) node (a1) {\phantom{A}}
(a1)++(\h,0) node (a2) {\phantom{A}}
(a1)++(-.5,0) node (a1') {\GGCatg}
(a2)++(.5,0) node (a2') {\GGCatg}
(a1)++(\h/2,\v) node (b) {\GGCatii}
(a1)++(-.2,0) node (a1'') {\phantom{\GG}}
(a2)++(.2,0) node (a2'') {\phantom{\GG}}
;
\draw[1cell=.9]
(a1) edge[bend right=\t] node[swap] {1} (a2)
(a1'') edge[bend left=\t] node[pos=.2] {\igst} (b)
(b) edge[bend left=\t] node[pos=.8] {\Lg} (a2'') 
;
\draw[2cell]
node[between=a1 and a2 at .42, shift={(0,\v/3)}, rotate=-90, 2label={above,\vg}] {\Rightarrow}
;
\end{tikzpicture}
\end{equation}
sends a pointed $G$-functor $X \cn \GG \to \Catgst$ to the $G$-natural isomorphism
\begin{equation}\label{vgf}
\begin{tikzpicture}[vcenter]
\def\t{28}
\draw[0cell]
(0,0) node (a1) {\phantom{\Gskel}}
(a1)++(2,0) node (a2) {\phantom{\Gskel}}
(a1)++(-.05,0) node (a1') {\GG}
(a2)++(.2,0) node (a2') {\Catgst}
;
\draw[1cell=.9]
(a1) edge[bend left=\t] node {\Lg\igst X} (a2)
(a1) edge[bend right=\t] node[swap] {X} (a2)
;
\draw[2cell]
node[between=a1 and a2 at .4, rotate=-90, 2label={above,\vg_X}] {\Rightarrow}
;
\end{tikzpicture}
\end{equation}
whose value at an object $\nbe \in \GG$ is the pointed isomorphism
\begin{equation}\label{vgf_nbe}
(\Lg\igst X)\nbe \fto[\iso]{\vg_{X,\nbe}} X\nbe
\end{equation}
defined as follows.
\begin{itemize}
\item If $\nbe = \vstar$, then $\vg_{X,\nbe} = 1_{\boldone} \cn \boldone \to \boldone$.
\item If $\nbe \neq \vstar$, then the pointed functor
\[(\Lg\igst X)\nbe = \ecint^{\angordm \in \Gskel}\mquad \bigvee_{\GGpunc(\ig\angordm; \nbe)} \mquad X\angordm \fto[\iso]{\vg_{X,\nbe}} X\nbe\]
sends a representing pair
\begin{equation}\label{upom_x}
(\upom; x) \in \GGpunc(\ig\angordm; \nbe) \times X\angordm,
\end{equation}
where $x \in X\angordm$ is either an object or a morphism, to
\begin{equation}\label{vgxnbe}
\vg_{X,\nbe} (\upom; x) = (X\upom)(x) \in X\nbe.
\end{equation}
\end{itemize}
\item[Invertibility of $\vg_{X,\nbe}$] 
To see that $\vg_{X,\nbe}$ \cref{vgf_nbe} is an isomorphism,  we use the object $\angordn \in \Gsk$ obtained from $\nbe \in \GG$ by forgetting the $G$-action $\be_j$ on $\ordn_j$ for each $j \in \ufs{q}$, where $q$ is the length of $\nbe$.  We define the isomorphism 
\begin{equation}\label{jin_nbe}
\ig\angordn \fto[\iso]{\jin = (1_{\ufs{q}}, \ang{1})} \nbe 
\end{equation}
in $\GG$ \cref{GG_morphisms} given by 
\begin{itemize}
\item the identity function on $\ufs{q}$ and
\item the identity pointed functions $1 \cn \ordn_j \to \ordi{n}{\be}{j}$ for $j \in \ufs{q}$.
\end{itemize}
Its inverse isomorphism in $\GG$ is defined in the same way:
\begin{equation}\label{jininv}
\nbe \fto[\iso]{\jininv = (1_{\ufs{q}}, \ang{1})} \ig\angordn.
\end{equation}
The inverse of $\vg_{X,\nbe}$ is defined as the following composite pointed functor.
\begin{equation}\label{vgfn_inv}
\begin{tikzpicture}[vcenter]
\def\v{-1.4}
\draw[0cell=.9]
(0,0) node (a11) {X\nbe}
(a11)++(3,0) node (a12) {(\Lg\igst X)\nbe}
(a12)++(3.2,0) node (a12') {\txint^{\angordm \in \Gskel} \txwedge_{\GGpunc(\ig\angordm; \nbe)} X\angordm}
(a11)++(0,\v) node (a21) {X\angordn}
(a12)++(0,\v) node (a22) {\txwedge_{\GGpunc(\ig\angordn; \nbe)} X\angordn}
;
\draw[1cell=.9]
(a11) edge node {\vginv_{X,\nbe}} (a12)
(a12) edge[equal] (a12')
(a11) edge node[pos=.45] {\iso} node[swap,pos=.4] {X\jininv} (a21)
(a21) edge node {\wg_1} (a22)
(a22) [rounded corners=2pt] -- node[pos=.55] {\wg_2} ($(a12')+(0,\v)$) -| (a12')
;
\end{tikzpicture}
\end{equation}
\begin{itemize}
\item Since $\jininv$ is a morphism in $\GG$, $X\jininv$ is a 1-cell in $\Catgst$, which means a pointed functor.
\item $\wg_1$ is the inclusion into the wedge summand corresponding to the isomorphism $\jin \cn \ig\angordn \to \nbe$ in $\GG$.
\item $\wg_2$ is the universal pointed functor corresponding to the object $\angordm = \angordn \in \Gskel$.
\end{itemize}
Thus, $\vginv_{X,\nbe}$ \cref{vgfn_inv} sends an object or a morphism $x \in X\nbe$ to the representing pair
\begin{equation}\label{vgfninv_x}
\vginv_{X,\nbe}(x) = \big(\jin; (X\jininv)(x) \big) \in \GGpunc(\ig\angordn; \nbe) \times X\angordn.
\end{equation}
The following equalities in $X\nbe$ prove that $\vg_{X,\nbe} \vginv_{X,\nbe}$ is the identity functor.
\[\begin{split}
& \vg_{X,\nbe} \vginv_{X,\nbe} (x) \\
&= \vg_{X,\nbe} \big(\jin; (X\jininv)(x) \big) \\
&= (X\jin) (X\jininv)(x) \\
&= (X 1_{\nbe}) x \\
&= 1_{X\nbe} x = x
\end{split}\]
Using the universal properties of coends, the following equalities in $(\Lg\igst X)\nbe$ prove that $\vginv_{X,\nbe} \vg_{X,\nbe}$ is the identity functor, where $(\upom;x)$ is a representing pair \cref{upom_x}.
\[\begin{split}
& \vginv_{X,\nbe} \vg_{X,\nbe} (\upom; x) \\
&=  \vginv_{X,\nbe} \big((X\upom)(x)\big) \\
&= \big(\jin ; (X\jininv) (X\upom)(x) \big) \\
&= \big(\jin ; X(\jininv\upom) (x) \big) \\
&= \big(\jin \circ \jininv \circ \upom ; x \big) \\
&= (\upom; x) \\
\end{split}\]
This proves that $\vg_{X,\nbe}$ is an isomorphism.
\item[$G$-equivariance of $\vg_X$] 
By \cref{ggcatg_icell_geq}, the $G$-equivariance of $\vg_X$ means that its component pointed functors are $G$-equivariant.  Using \cref{GGGcat_Gequiv,vgxnbe}, the following equalities in $X\nbe$ prove that the isomorphism $\vg_{X,\nbe}$ is $G$-equivariant.
\begin{equation}\label{vg_gnatural}
\begin{split}
& \vg_{X,\nbe}\big(g \cdot (\upom; x) \big) \\
&= \vg_{X,\nbe}\big( g\cdot \upom; gx \big) \\
&= X (g\cdot \upom) (gx) \\
&= g \big((X\upom)(\ginv gx)\big) \\
&= g \big((X\upom)(x)\big) \\
&= g \big(\vg_{X,\nbe} (\upom; x) \big)
\end{split}
\end{equation}
This proves that $\vg_X$ is a $G$-natural isomorphism.  The 2-naturality of $\vg$ follows from the description \cref{vgxnbe} of $\vg_{X,\nbe} (\upom; x)$.  Note that $\vginv_{X,\nbe}$ \cref{vgfn_inv} is also a $G$-isomorphism because it is inverse to the $G$-isomorphism $\vg_{X,\nbe}$.
\item[Left triangle identity]  
This triangle identity states that, for each pointed functor $X \cn \Gsk \to \Gcatst$, the following composite is the identity $G$-natural transformation.
\begin{equation}\label{Lgigst_lefttriangle}
\Lg X \fto{\Lg\ug_X} \Lg\igst\Lg X \fto{\vg_{\Lg X}} \Lg X
\end{equation}
It suffices to check that each nonbasepoint component of this composite is the identity functor.  For an object $\nbe \in \GG \setminus \{\vstar\}$ and a representing pair $(\upom; x)$ \cref{upom_x}, the following equalities in $(\Lg X)\nbe$ prove that the $\nbe$-component of the composite in \cref{Lgigst_lefttriangle} is the identity functor.
\[\begin{split}
& (\vg_{\Lg X, \nbe}) (\Lg\ug_X)_{\nbe} (\upom; x) \\
&= \vg_{\Lg X, \nbe} \big(\upom; \ug_{X,\angordm} x \big) \\
&= \vg_{\Lg X, \nbe} \big(\upom; (1_{\angordm}; x) \big) \\
&= \big((\Lg X)\upom\big) (1_{\angordm}; x) \\
&= \big(\upom\circ 1_{\angordm}; x \big) \\
&= (\upom; x)
\end{split}\]
This proves the left triangle identity.
\item[Right triangle identity]  
This triangle identity states that, for each pointed $G$-functor $X \cn \GG \to \Catgst$, the following composite is the identity natural transformation.
\begin{equation}\label{Lgigst_righttriangle}
\igst X \fto{\ug_{\igst X}} \igst\Lg\igst X \fto{\igst\vg_X} \igst X
\end{equation}
It suffices to check that each nonbasepoint component of this composite is the identity functor.  For an object $\angordm \in \Gsk \setminus \{\vstar\}$ and an object or a morphism $x \in (\igst X)\angordm = X\angordm$, the following equalities in $X\angordm$ prove that the $\angordm$-component of the composite in \cref{Lgigst_righttriangle} is the identity functor.  
\[\begin{split}
& (\igst\vg_X)_{\angordm} (\ug_{\igst X, \angordm}) (x) \\
&= \vg_{X,\angordm} (1_{\angordm}; x) \\
&= (X 1_{\angordm}) x \\
&= 1_{X\angordm} x = x
\end{split}\]
This proves the right triangle identity.  
\end{description}
In summary, the quadruple $(\Lg,\igst,\ug,\vg)$ is an adjoint 2-equivalence.
\end{proof}

\subsection*{Reconstructing $\GGG$-Categories}

By \cref{thm:ggcat_ggcatg_iieq}, for each pointed $G$-functor $X \cn \GG \to \Catgst$ and object $\nbe \in \GG \setminus \{\vstar,\ang{}\}$, the pointed $G$-category $X\nbe$ can be reconstructed from the pointed $G$-categories $X\angordm$ for $\angordm \in \Gsk$ via the counit $\vg \cn \Lg\igst \fiso 1$ \cref{vgf_nbe}.  The rest of this section makes this reconstruction explicit using the following definitions.  Recall from \cref{jin_nbe} that $\angordn \in \Gsk$ is the object obtained from $\nbe$ by forgetting the $G$-action $\be_j$ on $\ordn_j$ for each $j$, and $\ig\angordn \in \GG$ equips each $\ordn_j$ with the trivial $G$-action.  We sometimes abbreviate $\ig\angordn$ to $\angordn$.

\begin{definition}\label{def:gnbe}
Given an object $\nbe \in \GG \setminus \{\vstar,\ang{}\}$ of length $q>0$ and an element $g \in G$, we define the isomorphisms 
\begin{equation}\label{betag}
\begin{split}
& \angordn \fto[\iso]{\betag = (1_{\ufsq}, \ang{\be_j g}_{j \in \ufsq})} \angordn \andspace\\
& \nbe \fto[\iso]{\betag = (1_{\ufsq}, \ang{\be_j g}_{j \in \ufsq})} \nbe
\end{split}
\end{equation}
in $\GG$ \cref{GG_morphisms}, each consisting of
\begin{itemize}
\item the identity function on $\ufsq$ and
\item for each $j \in \ufsq$, the pointed bijection $\be_j g \cn \ordn_j \fiso \ordn_j$ given by the $G$-action $\be_j$ on $\ordn_j^{\be_j}$.
\end{itemize}
Given a pointed $G$-functor $X \cn \GG \to \Catgst$, the pointed isomorphism
\begin{equation}\label{Xbetag}
X\angordn \fto[\iso]{X\betag} X\angordn
\end{equation}
is $G$-equivariant by \cref{GG_Gaction}, \cref{GGGcat_Gequiv}, and the trivial $G$-action on the entries of $\angordn$.  We define the pointed $G$-category
\begin{equation}\label{Xnbeta}
\Xnbeta
\end{equation}
with
\begin{itemize}
\item underlying pointed category given by $X\angordn$ and
\item $g$-action functor given by the composite
\begin{equation}\label{Xnbeta_gaction}
\begin{tikzpicture}[vcenter]
\def\u{.65}
\draw[0cell]
(0,0) node (a1) {X\angordn}
(a1)++(1.8,0) node (a2) {X\angordn}
(a2)++(2,0) node (a3) {X\angordn}
;
\draw[1cell=.9]
(a1) edge node {g} (a2)
(a2) edge node {X\betag} (a3)
(a1) [rounded corners=2pt] |- ($(a2)+(-1,\u)$) -- node {g \cdot -} ($(a2)+(1,\u)$) -| (a3)
;
\end{tikzpicture}
\end{equation}
for $g \in G$.
\end{itemize}  
In other words, the $g$-action on $\Xnbeta$ is obtained from the original $g$-action on $X\angordn$ by twisting it by the $G$-functor $X\betag$.
\end{definition}

Recall from \cref{jininv} the isomorphism $\jininv \cn \nbe \fiso \angordn$ given by $(1_{\ufsq}, \ang{1})$.  \cref{Xjininv} proves that the pointed $G$-category $X\nbe$ is given by $\Xnbeta$ via $X\jininv$.

\begin{lemma}\label{Xjininv}
In the context of \cref{def:gnbe}, the pointed isomorphism
\[X\nbe \fto[\iso]{X\jininv} \Xnbeta\]
is $G$-equivariant.
\end{lemma}

\begin{proof}
The pointed functor $X\jininv$ is well defined because the underlying pointed category of $\Xnbeta$ is $X\angordn$.  It is an isomorphism of categories because $\jininv$ is the inverse of the morphism $\jin$ \cref{jin_nbe}.

For each $g \in G$, by \cref{GG_composite,jininv,betag}, there is a commutative diagram of isomorphisms
\begin{equation}\label{jininv_betag}
\begin{tikzpicture}[vcenter]
\def\v{-1.3}
\draw[0cell]
(0,0) node (a11) {\nbe}
(a11)++(2,0) node (a12) {\angordn}
(a11)++(0,\v) node (a21) {\nbe}
(a12)++(0,\v) node (a22) {\angordn}
;
\draw[1cell=.9]
(a11) edge node {\jininv} (a12)
(a12) edge node {\betag} (a22)
(a11) edge node[swap] {\betag} (a21)
(a21) edge node {\jininv} (a22)
;
\end{tikzpicture}
\end{equation}
in $\GG$.  By \cref{GG_Gaction}, the $g$-action on $\jininv$ yields the isomorphism
\begin{equation}\label{g_jininv}
\nbe \fto{g \cdot \jininv = \jininv \betaginv} \angordn
\end{equation}
because $G$ acts trivially on the entries of $\angordn = \ig\angordn$.  For an object or a morphism $x \in X\nbe$, the following equalities in $\Xnbeta$ prove that $X\jininv$ is $G$-equivariant.
\[\begin{aligned}
& (X\jininv)(gx) && \\
&= [X(\betag \jininv \betaginv)] (gx) && \text{by \cref{jininv_betag}}\\
&= (X\betag) [X(\jininv \betaginv) (gx)] && \text{by functoriality of $X$}\\
&= (X\betag) [X(g \cdot \jininv) (gx)] && \text{by \cref{g_jininv}}\\
&= (X\betag) [g (X\jininv) (\ginv gx)] && \text{by \cref{GGGcat_Gequiv}}\\
&= (X\betag) [g ((X\jininv)x) ] && \text{by $\ginv g = 1$}\\
&= g \cdot [(X\jininv)x] && \text{by \cref{Xnbeta_gaction}}
\end{aligned}\]
This proves that $X\jininv$ is a $G$-isomorphism.
\end{proof}

\begin{remark}\label{rk:shi_prop2}
\cref{Xjininv} is analogous to \cite[Prop.\ 2]{shimakawa91}, which deals with $\FGG$-spaces instead of $\GGG$-categories.
\end{remark}

%% file: chap/hgo.tex
To effectively compare Shimakawa $K$-theory \pcref{ch:shimakawa_K}
\begin{equation}\label{Ksho_chintro}
\begin{tikzpicture}[vcenter]
\def\h{2.5} \def\u{.6}
\draw[0cell]
(0,0) node (a1) {\phantom{\AlglaxO}}
(a1)++(0,-.04) node (a1') {\AlglaxO}
(a1)++(\h,0) node (a2) {\FGCatg}
(a2)++(\h,0) node (a3) {\FGTopg}
(a3)++(.9*\h,0) node (a4) {\Gspec} 
;
\draw[1cell=.9]
(a1) edge node {\Sgo} (a2)
(a2) edge node {\clast} (a3)
(a3) edge node {\Kfg} (a4)
(a1') [rounded corners=2pt, shorten <=-.2ex] |- ($(a2)+(0,\u)$) -- node {\Ksho} ($(a3)+(0,\u)$) -| (a4)
;
\end{tikzpicture}
\end{equation}
and its homotopical variant with our equivariant $K$-theory functor $\Kgo$ \cref{Kgo_functors}, this chapter constructs an equivalent variant of $\Kgo$, denoted by $\Khgo$.  The following diagram summarizes $\Khgo$ and $\Kgo$ in the top and bottom halves. 
\begin{equation}\label{KhgoKgo_chintro}
\begin{tikzpicture}[vcenter]
\def\c{2} \def\u{.7} \def\h{2.5} \def\d{.7} \def\t{10}
\draw[0cell=.9]
(0,0) node (a1) {\AlglaxO}
(a1)++(\c,\u) node (a2) {\GGCatg}
(a2)++(\h,0) node (a3) {\GGTopg}
(a3)++(\c,-\u) node (a4) {\Gspec}
(a1)++(\c,-\u) node (b2) {\GGCatii}
(b2)++(\h,0) node (b3) {\GGTopii}
;
\draw[1cell=.8]
(a1) edge[bend left=\t] node[pos=.7] {\Hgo} (a2)
(a2) edge node {\clast} (a3)
(a3) edge[bend left=\t] node[pos=.3] {\Kgg} (a4)
(a1) edge[bend right=\t] node[swap,pos=.7] {\Jgo} (b2)
(b2) edge node {\clast} (b3)
(b3) edge[bend right=\t] node[swap,pos=.3] {\Kg} (a4)
(a2) edge node {\igst} (b2)
(a3) edge[transform canvas={xshift=-.7em}] node[swap] {\igst} (b3)
(b3) edge node[swap] {\Lg} (a3)
(a1) [rounded corners=2pt] |- ($(a2)+(0,\d)$) -- node {\Khgo} ($(a3)+(0,\d)$) -| (a4)
;
\draw[1cell=.8]
(a1) [rounded corners=2pt] |- ($(b2)+(0,-\d)$) -- node {\Kgo} node[swap] {\phantom{x}} ($(b3)+(0,-\d)$) -| (a4)
;
\draw[2cell=.8]
node[between=a3 and b3 at .6, shift={(.9,0)}, rotate=-60, 2label={above,\kiso}] {\Rightarrow}
;
\end{tikzpicture}
\end{equation}
The functor $\Khgo$ factors through the category $\GGCatg$ of $\GGG$-categories, where $\GG$ is the $G$-equivariant version of $\Gsk$ \pcref{def:GG}.  The $J$-theory functor $\Jgo$ factors as $\igst\Hgo$ \pcref{jgohgoigst}, where $\igst \cn \GGCatg \to \GGCatii$ is the 2-equivalence in \cref{thm:ggcat_ggcatg_iieq}.  The middle square commutes, and the vertical pair $(\Lg,\igst)$ is an adjoint equivalence \pcref{thm:ggtop_ggtopg_iieq}.  In the right region, the prolongation functor $\Kg$ factors as the composite $\Kgg\Lg$ up to a natural isomorphism $\kiso$ \pcref{KgKgg}.  These facts imply that the functors $\Kgo$ and $\Khgo$ are naturally isomorphic, and the strong variant involving $\Kgosg$ and $\Khgosg$ is also true.  See \cref{KgoKhgo}.  \cref{ch:shim_top} compares $\Kgg$ with the last step of the homotopical Shimakawa $K$-theory.  \cref{part:kgo_shi_comp} compares $H$-theory $\Hgo$ with Shimakawa $H$-theory.  

\organization
This chapter consists of the following sections.

\secname{sec:hgo_objects}
This section constructs the object assignments of $\Hgo$ and $\Hgosg$, which send $\Op$-pseudoalgebras to pointed $G$-functors $\GG \to \Catgst$.

\secname{sec:hgo_onecells}
This section constructs the 1-cell assignments of $\Hgo$ and $\Hgosg$, which send lax $\Op$-morphisms and $\Op$-pseudomorphisms to $G$-natural transformations.

\secname{sec:hgo_twocells}
This section constructs the 2-cell assignments of $\Hgo$ and $\Hgosg$, which send $\Op$-transformations to $G$-modifications.  \cref{Hgo_twofunctor} records the fact that $\Hgo$ and $\Hgosg$ are 2-functors.  \cref{jgohgoigst} records the factorization $\Jgo = \igst\Hgo$ and its strong variant.

\secname{sec:GGTopg}
This section constructs the category $\GGTopg$ of $\GGG$-spaces and the functor $\clast$ from $\GGCatg$ to $\GGTopg$ induced by the classifying space functor $\cla$.  

\secname{sec:Kgg}
This section constructs the prolongation functor $\Kgg$ from $\GGTopg$ to the category $\Gspec$ of orthogonal $G$-spectra.  \cref{KgKgg} proves that the functors $\Kgg\Lg$ and $\Kg$ are naturally isomorphic.

\secname{sec:kgo_hthy}
This section defines the equivariant $K$-theory functor $\Khgo$ and its strong variant $\Khgosg$.  \cref{KgoKhgo} proves that $\Kgo$ and $\Kgosg$ are naturally isomorphic to, respectively, $\Khgo$ and $\Khgosg$.

\section{$H$-Theory on Objects}
\label{sec:hgo_objects}

This section defines, for each $\Tinf$-operad $\Op$, the object assignments of the (strong) $H$-theory 2-functors
\[\AlglaxO \fto{\Hgo} \GGCatg \andspace \AlgpspsO \fto{\Hgosg} \GGCatg,\]
which send $\Op$-pseudoalgebras to $\GGG$-categories, meaning pointed $G$-functors $\GG \to \Catgst$.

\secoutline
\begin{itemize}
\item \cref{def:nbeta_gcat} defines, for each $\Op$-pseudoalgebra $\A$, the object assignments of $\Hgo\A = \Adash$ and $\Hgosg\A = \Asgdash$, which send objects of $\GG$ to small pointed $G$-categories.
\item \cref{expl:nbeta_gcat} discusses the fact that \cref{def:nbeta_gcat} extends the construction of $\angordn$-systems to $\angordnbe \in \GG$.
\item \cref{def:AfangpsiGG} defines the morphism assignments of $\Adash$ and $\Asgdash$, which send morphisms of $\GG$ to pointed functors.
\item \cref{A_ptfunctorGG} proves that $\Adash$ and $\Asgdash$ are pointed $G$-functors.
\end{itemize}

\subsection*{$H$-Theory on Objects: Object Assignment}

Recall that, for $n \geq 0$, $\ufs{n} = \{1,2,\ldots,n\}$ denotes an unpointed finite set \cref{ufsn} and that $\ord{n} = \{0 < 1 < \cdots < n\}$ denotes a pointed finite set \cref{ordn} with basepoint 0.  Also recall \cref{not:compk} for substitution and partition.  \cref{def:nbeta_gcat} extends \cref{def:Aangordn_gcat} with the object $\angordn \in \Gsk$ replaced by $\angordnbe \in \GG$.

\begin{definition}[Pointed $G$-Categories of $\nbe$-Systems]\label{def:nbeta_gcat}
Given a $\Tinf$-operad $(\Op,\ga,\opu,\pcom)$ \pcref{as:OpA}, an $\Op$-pseudoalgebra $(\A,\gaA,\phiA)$ \pcref{def:pseudoalgebra}, and an object \cref{GG_objects}
\[\angordnbe = \sordi{n}{\be}{j}_{j \in \ufsq} 
= \big(\ordn_1^{\be_1}, \ldots, \ordn_q^{\be_q}\big) \in \GG,\]
we define the small pointed $G$-category $\Aangordnbe$ as follows.  A strong variant is defined in \cref{Asgangordnbe}.
\begin{description}
\item[Base cases] 
If $\angordnbe$ is either the basepoint $\vstar$ or the empty tuple $\ang{}$, then we define the pointed $G$-categories
\begin{equation}\label{vstar_systemGG}
\sys{\A}{\vstar} = \boldone \andspace \sys{\A}{\ang{}} = (\A,\zero)
\end{equation}
as in \cref{vstar_system}.
\item[Underlying pointed categories]
For $\angordnbe \in \GG \setminus \{\vstar,\ang{}\}$, we denote by $\angordn \in \Gsk$ the object obtained by forgetting the $G$-action $\be_j$ on $\ordn_j$ for each $j \in \ufsq$.  The underlying pointed category of $\Aangordnbe$ is defined as
\begin{equation}\label{Aangordnbe_pt}
\Aangordnbe = \big(\Aangordn, (\zero, 1_\zero)\big)
\end{equation}
as in \cref{A_angordn}.  An object in $\Aangordnbe$ is called an \index{system}\emph{$\angordnbe$-system}, which consists of the same data as an $\angordn$-system \pcref{def:nsystem}, and similarly for morphisms \pcref{def:nsystem_morphism}.
\item[$G$-action on $\nbe$-systems]
To define the $G$-action on the pointed category $\Aangordnbe$, suppose $g \in G$ and $(a,\glu) \in \Aangordnbe$ is an $\angordnbe$-system in $\A$ \cref{nsystem}.  We define the $\angordnbe$-system in $\A$ 
\begin{equation}\label{nsystem_gactionGG}
g \cdot (a,\glu) = (ga, g\glu)
\end{equation}
as follows.
\begin{description}
\item[Component objects] 
For each marker $\ang{s} = \ang{s_j \subseteq \ufs{n}_j}_{j \in \ufs{q}}$, the $\ang{s}$-component object of $(ga,g\glu)$ is defined as
\begin{equation}\label{ga_scomponentGG}
(ga)_{\ang{s}} = g a_{\ginv\ang{s}} \in \A,
\end{equation}
where the marker on the right-hand side is
\begin{equation}\label{ginvs}
\ginv\ang{s} = \ang{(\be_j g)^{\inv} s_j \subseteq \ufs{n}_j}_{j \in \ufs{q}}
= \ang{\ginv s_j}_{j \in \ufs{q}}.
\end{equation}
The object $a_{\ginv\ang{s}}$ is the $\ginv\ang{s}$-component object of $(a,\glu)$ \cref{a_angs}, and $g a_{\ginv\ang{s}}$ is its image under the $g$-action on $\A$.
\item[Gluing] 
Given an object $x \in \Op(r)$ with $r \geq 0$, a marker $\ang{s}$, an index $k \in \ufs{q}$, and a partition
\[s_k = \coprod_{i \in \ufs{r}}\, s_{k,i} \subseteq \ufs{n}_k\]
of $s_k$ into $r$ subsets, the gluing morphism \cref{gluing-morphism} of $(ga,g\glu)$ at $(x; \angs, k, \ang{s_{k,i}}_{i \in \ufs{r}})$ is defined by the following commutative diagram in $\A$.
\begin{equation}\label{ga_gluingGG}
\begin{tikzpicture}[vcenter]
\def\v{-1}
\draw[0cell=.85]
(0,0) node (a1) {\gaA_r\big(x; \ang{(ga)_{\ang{s} \compk\, s_{k,i}}}_{i \in \ufs{r}} \big)}
(a1)++(0,\v) node (a2) {\gaA_r\big(x; \ang{g a_{\ginv\ang{s} \,\compk\, (\ginv s_{k,i})}}_{i \in \ufs{r}}\big)}
(a2)++(0,\v) node (a3) {g \gaA_r \big(\ginv x; \ang{a_{\ginv\ang{s} \,\compk\, (\ginv s_{k,i})}}_{i \in \ufs{r}} \big)}
(a1)++(6.5,0) node (b1) {(ga)_{\ang{s}}}
(b1)++(0,2*\v) node (b3) {\phantom{g a_{\ginv\ang{s}}}}
(b3)++(0,-.05) node (b3') {g a_{\ginv\ang{s}}}
;
\draw[1cell=.85]
(a1) edge node {(g\glu)_{x;\, \ang{s},\, k, \ang{s_{k,i}}_{i \in \ufs{r}}}} (b1)
(a3) edge node {g \glu_{\ginv x;\, \ginv\ang{s},\, k,\, \ang{\ginv s_{k,i}}_{i \in \ufs{r}}}} (b3)
(a1) edge[equal,shorten >=-.5ex] (a2)
(a2) edge[equal,shorten >=-.5ex] node[swap] {(\mathbf{f})} (a3)
(b1) edge[equal] (b3')
;
\end{tikzpicture}
\end{equation}
\begin{itemize}
\item The two unlabeled equalities in \cref{ga_gluingGG} follow from the definition of $(ga)_{\ang{s}}$ in \cref{ga_scomponentGG} and the partition
\[\ginv s_k = \coprod_{i \in \ufs{r}}\, \ginv s_{k,i}.\]
\item The equality labeled $(\mathbf{f})$ follows from the $G$-functoriality of $\gaA_r$ \cref{gaAn}.
\item In the bottom horizontal arrow in \cref{ga_gluingGG}, 
\[\glu_{\ginv x;\, \ginv\ang{s},\, k,\, \ang{\ginv s_{k,i}}_{i \in \ufs{r}}} \cn 
\gaA_r \big(\ginv x ; \ang{a_{\ginv\ang{s} \,\compk\, (\ginv s_{k,i})}}_{i \in \ufs{r}} \big)
\to a_{\ginv\ang{s}}\]
is the indicated gluing morphism of $(a,\glu)$, and $g\glu_{\cdots}$ is its image under the $g$-action on $\A$.
\end{itemize}
\item[Axioms] Each of the axioms of an $\angordn$-system,  \cref{system_obj_unity,system_naturality,system_unity_i,system_unity_iii,system_equivariance,system_associativity,system_commutativity}, for $(ga,g\glu)$ follows from the corresponding axiom for $(a,\glu)$ and the following facts.
\begin{itemize}
\item The $g$-action on $\A$ is a functor.
\item The basepoint $\zero \in \A$ and its identity morphism $1_\zero$ are $G$-fixed.  In particular, the base $\angordnbe$-system $(\zero,1_\zero)$ is $G$-fixed. 
\item The object $* \in \Op(0)$ and the operadic unit $\opu \in \Op(1)$ are $G$-fixed.
\item The right symmetric group action on $\Op$, the pseudo-commutative structure $\pcom$ on $\Op$ \cref{pseudocom_isos}, and the associativity constraint $\phiA$ of $\A$ \cref{phiA} are $G$-equivariant.
\end{itemize}
This finishes the definition of the $\angordnbe$-system $g \cdot (a,\glu) = (ga,g\glu)$.
\end{description}
\item[$G$-action on morphisms]  
For a morphism of $\angordnbe$-systems in $\A$ \pcref{def:nsystem_morphism}
\[(a,\glu^a) \fto{\tha} (b,\glu^b),\]
the morphism of $\angordnbe$-systems
\begin{equation}\label{gthetaGG}
(ga,g\glu^a) \fto{g\theta} (gb, g\glu^b)
\end{equation}
is defined by, for each marker $\ang{s} = \bang{s_j \subseteq \ufs{n}_j}_{j \in \ufs{q}}$, the $\ang{s}$-component morphism
\begin{equation}\label{gtheta_angsGG}
(ga)_{\ang{s}} = ga_{\ginv\ang{s}}
\fto{(g\theta)_{\ang{s}} \,=\, g\theta_{\ginv\ang{s}}}
(gb)_{\ang{s}} = gb_{\ginv\ang{s}}.
\end{equation}
\begin{itemize}
\item The unity axiom \cref{nsystem_mor_unity} holds for $g\theta$ because the identity morphism $1_\zero$ is $G$-fixed. 
\item The compatibility axiom \cref{nsystem_mor_compat} holds for $g\theta$ by the compatibility axiom for $\theta$ and the functoriality of the $g$-action on $\A$.
\end{itemize}
This finishes the definition of the morphism $g\theta$ of $\angordnbe$-systems.  
\item[Functoriality of $G$-action] 
For each $g \in G$, the functoriality of the $g$-action on $\Aangordnbe$, as defined in \cref{nsystem_gactionGG,gthetaGG}, follows from
\begin{itemize}
\item the definition \cref{gtheta_angsGG} of $(g\theta)_{\ang{s}}$ and
\item the functoriality of the $g$-action on $\A$.
\end{itemize}
As $g \in G$ varies, this defines a $G$-action on the pointed category $\Aangordnbe$ by
\begin{itemize}
\item the definitions \cref{ga_scomponentGG,ginvs,ga_gluingGG,gtheta_angsGG}, and 
\item the $G$-action axioms for $\A$, $\Op(r)$, and each $\ord{n}_j^{\be_j}$ for $j \in \ufs{q}$.
\end{itemize}
\end{description}  
This finishes the definition of the pointed $G$-category $\Aangordnbe$.

\parhead{Strong variant}.
The small pointed $G$-category 
\begin{equation}\label{Asgangordnbe}
\Asgangordnbe
\end{equation}
has underlying pointed category $\Asgangordn$ \cref{Aangordnsg}.
\begin{itemize}
\item An object in $\Asgangordnbe$ is called a \index{system!strong}\index{strong system}\emph{strong $\angordnbe$-system}, which consists of the same data as a strong $\angordn$-system.  Recall from the last paragraph of \cref{def:nsystem} that an $\angordn$-system is \emph{strong} if its gluing morphisms are isomorphisms. 
\item The $G$-action on $\Asgangordnbe$ is the restriction of the $G$-action on $\Aangordnbe$ to the full subcategory of strong $\angordnbe$-systems.  This $G$-action is well defined because, if each gluing morphism $\glu_{\cdots}$ is an isomorphism in \cref{ga_gluingGG}, then so is its image $g\glu_{\cdots}$ under the $g$-action.\defmark
\end{itemize}
\end{definition}

\begin{explanation}[Compatibility of Definitions]\label{expl:nbeta_gcat}
Each object $\angordn = \ang{\ordn_j}_{j \in \ufsq} \in \Gsk$ is also regarded as the object $\sordi{n}{\eps}{j}_{j \in \ufsq} \in \GG$ via the full subcategory inclusion $\ig \cn \Gsk \to \GG$ \cref{ig} that equips each pointed finite set $\ordn_j$ with the trivial $G$-action.  The pointed $G$-categories $\Aangordn$ and $\sys{\A}{\sordi{n}{\eps}{j}_{j \in \ufsq}}$, as defined in \cref{def:Aangordn_gcat,def:nbeta_gcat}, are the same for the following reasons.
\begin{itemize}
\item By \cref{Aangordnbe_pt}, their underlying pointed categories are the same. 
\item Their $G$-actions are the same because the marker $\ginv\angs$ \cref{ginvs} is equal to $\angs$ when each $\ordn_j$ has the trivial $G$-action.  Thus, \cref{ga_scomponentGG,ga_gluingGG,gtheta_angsGG} reduce to, respectively, \cref{ga_scomponent,ga_gluing,gtheta_angs}.
\end{itemize}
In summary, for each object $\angordn \in \Gsk$, \cref{def:Aangordn_gcat,def:nbeta_gcat} yield the same pointed $G$-category $\Aangordn$.  The same remark also applies to the strong variant $\Asgangordn$ in \cref{sgAordnbe,Asgangordnbe}.
\end{explanation}

\subsection*{$H$-Theory on Objects: Morphism Assignment}

For a morphism $\upom \cn \angordmal \to \angordnbe$ in $\GG$ \cref{GG_morphisms}, we denote by $\upom \cn \angordm \to \angordn$ the morphism in $\Gsk$ \cref{Gsk_morphisms} obtained by forgetting the $G$-actions $\al_i$ on $\ordm_i$ and $\be_j$ on $\ordn_j$.

\begin{definition}\label{def:AfangpsiGG}
For a $\Tinf$-operad $\Op$ \pcref{as:OpA}, an $\Op$-pseudoalgebra $\A$ \pcref{def:pseudoalgebra}, and a morphism  $\upom \cn \angordmal \to \angordnbe$ in $\GG$ \cref{GG_morphisms}, we define the pointed functor
\begin{equation}\label{Aupom}
\Aangordmal \fto{\Aupom} \Aangordnbe
\end{equation}
as the underlying pointed functor of $\Aupom \cn \Aangordm \to \Aangordn$ in \cref{AF}.

\parhead{Strong variant}.  We define the pointed functor
\begin{equation}\label{Aupomsg}
\Asgangordmal \fto{\Aupomsg} \Asgangordnbe
\end{equation}
as the underlying pointed functor of $\Aupomsg \cn \Asgangordm \to \Asgangordn$ in \cref{AF_sg}.
\end{definition}

\begin{lemma}[$H$-Theory on Objects]\label{A_ptfunctorGG}\index{H-theory@$H$-theory!object}
For a $\Tinf$-operad $\Op$ and an $\Op$-pseudoalgebra $\A$, the object and morphism assignments 
\[\angordnbe \mapsto \Aangordnbe \andspace \upom \mapsto \Aupom\]
in, respectively, \cref{def:nbeta_gcat,def:AfangpsiGG} define a pointed $G$-functor
\[\GG \fto{\Hgo\A = \Adash} \Catgst.\]
Moreover, the strong variant 
\[\GG \fto{\Hgosg\A = \Asgdash} \Catgst\]
defined in \cref{Asgangordnbe,Aupomsg} is also a pointed $G$-functor.
\end{lemma}

\begin{proof}
We prove that $\Hgo\A = \Adash$ is a pointed $G$-functor.  The same argument applies to the strong variant by restricting to strong $\angordnbe$-systems.  
\begin{description}
\item[Pointed functoriality] 
By \cref{Aangordnbe_pt}, $\Aangordnbe$ and $\Aangordn$ are the same as pointed categories.  By \cref{def:FG,GG_morphisms}, each morphism $\upom \cn \angordmal \to \angordnbe$ in $\GG$ is uniquely determined by the morphism $\upom \cn \angordm \to \angordn$ in $\Gsk$ obtained by forgetting the $G$-action on each $\ordm_i$ and $\ordn_j$.  Thus, the proof of \cref{A_ptfunctor} given in \cite[\namecref{EqK:A_ptfunctor} \ref*{EqK:A_ptfunctor}]{yau-eqk} proves that $\Hgo\A = \Adash$ is a pointed functor.
\item[$G$-equivariance] 
On objects, $\Hgo\A = \Adash$ is $G$-equivariant because $G$ acts trivially on the objects of $\GG$ \cref{GG_objects} and $\Catgst$ \cref{Catgst_iicat}.  To show that $\Adash$ is $G$-equivariant on morphisms \cref{GGGcat_Gequiv}, we need to show that, for each $g \in G$ and each morphism $\upom \cn \angordmal \to \angordnbe$ in $\GG$, there is an equality of functors
\begin{equation}\label{AgF}
\sys{\A}{(g \cdot \upom)} = g (\Aupom) \ginv\cn \Aangordmal \to \Aangordnbe.
\end{equation}
There are several cases.
\begin{description}
\item[Base case]  
If either $\angordmal$ or $\angordnbe$ is the initial-terminal basepoint $\vstar \in \GG$, then $g \cdot \upom = \upom$ and
\[\sys{\A}{(g \cdot \upom)} = \Aupom\]
is the constant functor at the basepoint.  There is also an equality 
\[g (\Aupom) \ginv = \Aupom\]
because the basepoint of each pointed $G$-category, including $\Aangordnbe$, is $G$-fixed.
\item[Empty domain and codomain]  
Suppose $\sord{m}{\al} = \ang{} = \sord{n}{\be}$.  Then $\upom \cn \ang{} \to \ang{}$ is either the identity morphism, $1_{\ang{}}$, or the 0-morphism, $\ang{} \to \vstar \to \ang{}$.  In either case, $g \cdot \upom = \upom$ and
\[\sys{\A}{(g \cdot \upom)} = \Aupom.\]
\begin{itemize}
\item If $\upom$ is the identity morphism, then $\Aupom$ is the identity functor, which is fixed by the conjugation $G$-action.
\item If $\upom$ is the 0-morphism, then $\Aupom$ is the constant functor at the basepoint $\zero \in \A\ang{} = \A$ by definition \cref{A_empty_empty}.  There is an equality
\[g (\Aupom) \ginv = \Aupom\]
because the basepoint $\zero \in \A$ is $G$-fixed by \cref{pseudoalg_zero}.
\end{itemize}
In either case, the desired equality \cref{AgF} holds.
\item[Empty domain]  
Suppose $\sord{m}{\al} = \ang{}$ and $\vstar \neq \sord{n}{\be}$ has length $q > 0$.  Then $\upom$ factors as a composite
\begin{equation}\label{emptydom_diag}
\begin{tikzpicture}[vcenter]
\def\t{.7} \def\h{3}
\draw[0cell]
(0,0) node (a1) {\ang{}}
(a1)++(\h,0) node (a2) {\ordtu{1}_{j \in \ufs{q}}}
(a2)++(\h,0) node (a3) {\phantom{\sord{n}{\be}}}
(a3)++(0,.05) node (a3') {\sord{n}{\be}}
;
\draw[1cell=.9]
(a1) edge node {(\im_q,\ang{1_{\ord{1}}}_{j \in \ufs{q}})} (a2)
(a2) edge node {(1_{\ufs{q}}, \ang{\psi_j}_{j \in \ufs{q}})} (a3)
;
\draw[1cell=.9]
(a1) [rounded corners=3pt, shorten <=0ex] |- ($(a2)+(-1,\t)$)
-- node {\upom \,=\, (\im_q, \ang{\psi_j}_{j \in \ufs{q}})} ($(a2)+(1,\t)$) -| (a3')
;
\end{tikzpicture}
\end{equation}
with $\im_q \cn \emptyset \to \ufsq$ the unique function, and $\Aupom$ is the following composite.
\begin{equation}\label{Aemptydom_diag}
\begin{tikzpicture}[vcenter]
\def\t{.7}
\draw[0cell]
(0,0) node (a1) {\sys{\A}{\ang{}}}
(a1)++(.9,0) node (a1') {\phantom{\A}}
(a1')++(0,.03) node (a1'') {\A}
(a1')++(1.6,0) node (a2) {\phantom{\sys{\A}{\ordtu{1}_{j \in \ufs{q}}}}}
(a2)++(0,-.03) node (a2') {\sys{\A}{\ordtu{1}_{j \in \ufs{q}}}}
(a2)++(3.7,0) node (a3) {\sys{\A}{\sord{n}{\be}}}
;
\draw[1cell=.9]
(a1) edge[equal, shorten <=-.4ex] (a1')
(a1') edge node {\iso} (a2)
(a2) edge node {\sys{\A}{(1_{\ufs{q}}, \ang{\psi_j}_{j \in \ufs{q}})}} (a3)
;
\draw[1cell=.9]
(a1) [rounded corners=3pt, shorten <=0ex] |- ($(a2)+(0,\t)$)
-- node {\Aupom} ($(a2)+(1,\t)$) -| (a3)
;
\end{tikzpicture}
\end{equation}
The pointed isomorphism $\A \iso \sys{\A}{\ordtu{1}}$, defined in \cref{Aordtuone}, is $G$-equivariant by \cref{nsystem_gactionGG,ga_scomponentGG,ginvs,ga_gluingGG} and the following facts.
\begin{itemize}
\item The $G$-action on $\ord{1}$ is trivial.
\item For each $g \in G$, the $g$-action on $\A$ is a functor, so it preserves identity morphisms.
\end{itemize}
Thus, to show that $\Aupom$ is $G$-equivariant, it suffices to show that the functor $\sys{\A}{(1_{\ufs{q}}, \ang{\psi_j}_{j \in \ufs{q}})}$ is $G$-equivariant.  This is included in the following case as $f = 1_{\ufs{q}}$.
\item[Nonempty domain]  
For the final case, suppose $\upom$ has the form \cref{fangpsiGG}
\[\upom = (f,\ang{\psi_j}_{j \in \ufs{q}}) \cn \sord{m}{\al} = \sordi{m}{\al}{i}_{i \in \ufs{p}} \to 
\sord{n}{\be} = \sordi{n}{\be}{j}_{j \in \ufs{q}}\]
with $\angordmal, \angordnbe \in \GG \setminus \{\vstar,\ang{}\}$.  By \cref{A_fangpsi}, $\Aupom$ is the composite 
\begin{equation}\label{Aupom_diag}
\begin{tikzpicture}[vcenter]
\def\h{2.5} \def\t{.7}
\draw[0cell]
(0,0) node (a1) {\Aangordmal}
(a1)++(\h,0) node (a2) {\sys{\A}{(f_*\sord{m}{\al})}}
(a2)++(\h,0) node (a3) {\phantom{\Aangordnbe}}
(a3)++(0,.05) node (a3') {\Aangordnbe}
;
\draw[1cell=.9]
(a1) edge node{\ftil} (a2)
(a2) edge node {\psitil} (a3)
;
\draw[1cell=1]
(a1) [rounded corners=3pt, shorten <=0ex] |- ($(a2)+(-1,\t)$)
-- node {\sys{\A}{(f,\angpsi)}} ($(a2)+(1,\t)$) -| (a3')
;
\end{tikzpicture}
\end{equation}
of the pointed functors $\ftil$ and $\psitil$ in, respectively, \cref{def:ftil_functor,def:psitil_functor}.  First, we show that the two functors in \cref{AgF} are equal on objects. 
\begin{description}
\item[Component objects]  
For an $\angordmal$-system $(a,\glu) \in \sys{\A}{\sord{m}{\al}}$ and a marker $\ang{s} = \bang{s_j \subseteq \ufs{n}_j}_{j \in \ufs{q}}$, using \cref{atil_component,apsitil_angs}, the left-hand side of \cref{AgF} yields the following $\ang{s}$-component object.
\begin{equation}\label{A_fgpsiginv_aglu}
\begin{split}
&\big(\sys{\A}{(g \cdot (f,\angpsi))} (a,\glu)\big)_{\ang{s}}\\
&= \big(\sys{\A}{(f,\ang{g\psi_j\ginv}_{j \in \ufs{q}})} (a,\glu)\big)_{\ang{s}}\\
&= \atil_{\ang{(g\psi_j\ginv)^{-1} s_j}_{j \in \ufs{q}}}\\
&= \begin{cases}
\zero & \text{if $(g\psi_j \ginv)^{-1} s_j = \emptyset$ for some $j \in \ufs{q}$}\\
a_{\ang{(g \psi_{f(i)} \ginv)^{-1} s_{f(i)}}_{i \in \ufs{p}}} & \text{if $(g\psi_j \ginv)^{-1} s_j \neq \emptyset$ for each $j \in \ufs{q}$}
\end{cases}\\
&= \begin{cases}
\zero & \text{if $(g\psi_j)^{-1} s_j = \emptyset$ for some $j \in \ufs{q}$}\\
a_{\ang{(g \psi_{f(i)} \ginv)^{-1} s_{f(i)}}_{i \in \ufs{p}}} & \text{if $(g\psi_j)^{-1} s_j \neq \emptyset$ for each $j \in \ufs{q}$}
\end{cases}
\end{split}
\end{equation}
The last step in \cref{A_fgpsiginv_aglu} uses the fact that the $\ginv$-action is a bijection 
\begin{equation}\label{ginv_action}
\ord{m}_{\finv(j)} \fto[\iso]{\ginv} \ord{m}_{\finv(j)} \quad \text{for each $j \in \ufs{q}$.}
\end{equation}
This implies that 
\[(g\psi_j \ginv)^{-1} s_j \neq \emptyset \iffspace (g\psi_j)^{-1} s_j \neq \emptyset.\]

Using \cref{atil_component,apsitil_angs,ga_scomponentGG}, the right-hand side of \cref{AgF} yields the following $\ang{s}$-component object.
\begin{equation}\label{gAfpsiginv_aglu}
\begin{split}
&\big((g \sys{\A}{(f,\angpsi)} \ginv) (a,\glu)\big)_{\ang{s}}\\
&= g \big(\sys{\A}{(f,\angpsi)} (\ginv (a,\glu))\big)_{\ang{\ginv s_j}_{j \in \ufs{q}}}\\
&= \scalebox{.9}{$\begin{cases}
g\zero & \text{if $\psiinv_j \ginv s_j = \emptyset$ for some $j \in \ufs{q}$}\\
g \big((\ginv (a,\glu))_{\ang{\psiinv_{f(i)} \ginv s_{f(i)}}_{i \in \ufs{p}}}\big)
& \text{if $\psiinv_j \ginv s_j \neq \emptyset$ for each $j \in \ufs{q}$}
\end{cases}$}\\
&= \scalebox{.9}{$\begin{cases}
\zero & \text{if $\psiinv_j \ginv s_j = \emptyset$ for some $j \in \ufs{q}$}\\
g \ginv a_{\ang{(\ginv)^{-1} \psiinv_{f(i)} \ginv s_{f(i)}}_{i \in \ufs{p}}}
& \text{if $\psiinv_j \ginv s_j \neq \emptyset$ for each $j \in \ufs{q}$}
\end{cases}$}\\
&= \begin{cases}
\zero & \text{if $(g\psi_j)^{-1} s_j = \emptyset$ for some $j \in \ufs{q}$}\\
a_{\ang{(g \psi_{f(i)} \ginv)^{-1} s_{f(i)}}_{i \in \ufs{p}}}
& \text{if $(g\psi_j)^{-1} s_j \neq \emptyset$ for each $j \in \ufs{q}$}
\end{cases}
\end{split}
\end{equation}
Thus, by \cref{A_fgpsiginv_aglu,gAfpsiginv_aglu}, the $\angordnbe$-systems
\begin{equation}\label{Agupomaglu}
\sys{\A}{(g \cdot (f,\angpsi))} (a,\glu) \andspace
(g \sys{\A}{(f,\angpsi)} \ginv) (a,\glu)
\end{equation}
have the same $\ang{s}$-component objects.  Next, we show that these two $\angordnbe$-systems have the same gluing morphisms.
\item[Gluing] 
In the context of \cref{skpartition}, first suppose $(g\psi_j\ginv)^{-1} s_j \neq \emptyset$---which is equivalent to $(g\psi_j)^{-1} s_j \neq \emptyset$, as we explain in \cref{ginv_action}---for each $j \in \ufs{q}$ and $\finv(k) \neq \emptyset$. 
\begin{itemize}
\item By \cref{glutil_component,glupsitil}, the gluing morphism of the $\angordnbe$-system $\sys{\A}{(g \cdot (f,\angpsi))} (a,\glu)$ is given by
\begin{equation}\label{AgF_left_gluing}
\glu_{x;\, \ang{(g\psi_{f(i)} \ginv)^{-1} s_{f(i)}}_{i \in \ufs{p}} \scs \finv(k) \scs \ang{(g\psi_k \ginv)^{-1} s_{k,\ell}}_{\ell \in \ufs{r}}}.
\end{equation}
\item By \cref{glutil_component,glupsitil,ga_gluingGG}, the gluing morphism of the $\angordnbe$-system $(g \sys{\A}{(f,\angpsi)} \ginv) (a,\glu)$ is given by
\begin{equation}\label{AgF_right_gluing}
g\ginv \glu_{(\ginv)^{-1} \ginv x;\, \ang{(\ginv)^{-1} \psiinv_{f(i)} \ginv s_{f(i)}}_{i \in \ufs{p}} \scs \finv(k) \scs \ang{(\ginv)^{-1} \psiinv_k \ginv s_{k,\ell}}_{\ell \in \ufs{r}}}.
\end{equation}
\end{itemize}
Since $g\ginv = 1$, the gluing morphisms in \cref{AgF_left_gluing,AgF_right_gluing} are equal.  All other gluing morphisms of the two $\angordnbe$-systems in \cref{Agupomaglu} are identities.  Thus, the two $\angordnbe$-systems in \cref{Agupomaglu} have the same gluing morphisms.  This proves the equality \cref{AgF} on objects.
\item[Morphisms]  
To prove the equality \cref{AgF} on morphisms of $\angordmal$-systems, we reuse the computation in \cref{A_fgpsiginv_aglu,gAfpsiginv_aglu}.  Instead of \cref{atil_component,apsitil_angs,ga_scomponentGG} for component objects, here we use \cref{thatil_component,thapsitil_component,gtheta_angsGG} for component morphisms.  
\end{description}
\end{description}
\end{description}
This completes the proof of the desired equality \cref{AgF}, proving that $\Hgo\A = \Adash$ is a pointed $G$-functor.
\end{proof}

\section{$H$-Theory on 1-Cells}
\label{sec:hgo_onecells}

This section constructs the 1-cell assignments of the (strong) $H$-theory 2-functors
\[\AlglaxO \fto{\Hgo} \GGCatg \andspace \AlgpspsO \fto{\Hgosg} \GGCatg\]
for a $\Tinf$-operad $\Op$ \pcref{as:OpA}.  Their object assignments, 
\[\A \mapsto \Hgo\A = \Adash \andspace \A \mapsto \Hgosg\A = \Asgdash,\]
are given in \cref{A_ptfunctorGG}.  The 1-cell assignments of $\Hgo$ and $\Hgosg$ send, respectively, lax $\Op$-morphisms and $\Op$-pseudomorphisms \pcref{def:laxmorphism} to 1-cells in $\GGCatg$, which are $G$-natural transformations \cref{ggcatg_icell}.  Components of 1-cells in $\GGCatg$ are pointed $G$-functors \cref{ggcatg_icell_geq}.

\secoutline
\begin{itemize}
\item \cref{def:hgo_icell} constructs a $G$-natural transformation $\Hgo f$ for each lax $\Op$-morphism $f$ and a $G$-natural transformation $\Hgosg f$ for each $\Op$-pseudomorphism $f$.
\item \cref{hgo_icell_welldef} proves that $\Hgo f$ and $\Hgosg f$ are well defined.
\end{itemize}

\begin{definition}\label{def:hgo_icell}
Suppose
\[\big(\A,\gaA,\phiA\big) \fto{(f,\actf)} \big(\B,\gaB,\phiB\big)\]
is a lax $\Op$-morphism \pcref{def:laxmorphism} between $\Op$-pseudoalgebras $(\A,\gaA,\phiA)$ and $(\B,\gaB,\phiB)$ \pcref{def:pseudoalgebra} for a $\Tinf$-operad $(\Op,\ga,\opu,\pcom)$ \pcref{as:OpA}.  For the pointed $G$-functors $\Hgo\A = \Adash$ and $\Hgo\B = \Bdash$ in \cref{A_ptfunctorGG}, we define a $G$-natural transformation
\begin{equation}\label{hgo_f}
\begin{tikzpicture}[vcenter]
\def\t{30}
\draw[0cell]
(0,0) node (a1) {\GG}
(a1)++(2,0) node (a2) {\phantom{\GG}}
(a2)++(.12,0) node (a2') {\Catgst}
;
\draw[1cell=.8]
(a1) edge[bend left=\t] node {\Adash} (a2)
(a1) edge[bend right=\t] node[swap] {\Bdash} (a2)
;
\draw[2cell=.9]
node[between=a1 and a2 at .33, rotate=-90, 2label={above,\Hgo f}] {\Rightarrow}
;
\end{tikzpicture}
\end{equation}
as follows.  A strong variant is defined in \cref{hgosg_f}.  For an object \cref{GG_objects}
\[\angordnbe = \sordi{n}{\be}{j}_{j \in \ufsq} 
= \big(\ordn_1^{\be_1}, \ldots, \ordn_q^{\be_q}\big) \in \GG,\]
the $\angordnbe$-component pointed $G$-functor
\begin{equation}\label{hgo_f_nbe}
\Aangordnbe \fto{(\Hgo f)_{\angordnbe}} \Bangordnbe
\end{equation}
is defined as follows.
\begin{description}
\item[Base cases]
If $\angordnbe$ is either the basepoint $\vstar$ or the empty tuple $\ang{}$, then $(\Hgo f)_{\angordnbe}$ is defined as, respectively, $1_{\bone}$ and $f \cn \A \to \B$.  This is well defined by \cref{laxmorphism_basepoint,vstar_systemGG}.
\item[Component objects]
Suppose $\angordnbe \in \GG \setminus \{\vstar,\ang{}\}$ and $(a,\glu^a) \in \Aangordnbe$ is an $\angordnbe$-system in $\A$ \cref{Aangordnbe_pt}, which consists of the same data as an $\angordn$-system \pcref{def:nsystem}.  The $\angordnbe$-system
\begin{equation}\label{hgof_aglu}
(\Hgo f)_{\angordnbe} (a,\glu^a) \in \Bangordnbe
\end{equation}
has, for each marker $\angs = \ang{s_j \subseteq \ufsn_j}_{j \in \ufsq}$, $\angs$-component object \cref{a_angs} defined as the image under $f$ of the $\angs$-component object of $(a,\glu^a)$:
\begin{equation}\label{hgof_aglu_comp}
\big((\Hgo f)_{\angordnbe} (a,\glu^a)\big)_{\angs} = fa_{\angs} \in \B.
\end{equation}
\item[Gluing]
For each object $x \in \Op(r)$ with $r \geq 0$, marker $\ang{s}$, index $k \in \ufs{q}$, and partition of $s_k$ into $r$ subsets
\[s_k = \coprod_{i \in \ufs{r}} \, s_{k,i} \subseteq \ufs{n}_k,\]
the gluing morphism \cref{gluing-morphism} of the $\angordnbe$-system $(\Hgo f)_{\angordnbe} (a,\glu^a)$ at $(x; \angs, k, \ang{s_{k,i}}_{i \in \ufs{r}})$ is defined as the following composite in $\B$.
\begin{equation}\label{hgof_gluing}
\begin{tikzpicture}[vcenter]
\def\u{-1} \def\h{4.5} \def\a{10} \def\b{.7}
\draw[0cell=.85]
(0,0) node (a1) {\gaB_r\big(x \sscs \ang{((\Hgo f)_{\angordnbe} (a,\glu^a))_{\angs \,\compk\, s_{k,i}}}_{i \in \ufsr}\big)}
(a1)++(\h,0) node (a2) {\big((\Hgo f)_{\angordnbe} (a,\glu^a)\big)_{\angs}}
(a1)++(0,\u) node (b1) {\gaB_r\big(x \sscs \ang{f a_{\angs \,\compk\, s_{k,i}}}_{i \in \ufsr}\big)}
(a2)++(0,\u) node (b2) {f a_{\angs}}
(b1)++(\h/2,-1.2) node (c) {f \gaA_r \big(x \sscs \ang{a_{\angs \,\compk\, s_{k,i}}}_{i \in \ufsr} \big)}
;
\draw[1cell=.9]
(a1) edge[equal,shorten <=-0ex,shorten >=-.5ex] (b1)
(a2) edge[equal] (b2)
;
\draw[1cell=.9]
(a1) [rounded corners=2pt, shorten <=-.2ex] |- ($(a1)+(1,\b)$)
-- node {\glu_{x; \angs,\, k, \ang{s_{k,i}}_{i \in \ufs{r}}}} ($(a2)+(-1,\b)$) -| (a2)
;
\draw[1cell=.9]
(b1) [rounded corners=2pt, shorten <=-.2ex] |- node[pos=.2,swap] {\actf_r} (c);
\draw[1cell=.9]
(c) [rounded corners=2pt, shorten <=-0ex] -| node[pos=.8] {f \glu^a_{x; \angs,\, k, \ang{s_{k,i}}_{i \in \ufs{r}}}} (b2);
\end{tikzpicture}
\end{equation}
The lower-left arrow $\actf_r$ is a component of the $r$-th action constraint of $f$ \cref{actf_component}.  The lower-right arrow is the image under $f$ of the indicated gluing morphism of $(a,\glu^a)$.
\item[Morphisms]
Suppose we are given a morphism of $\angordnbe$-systems in $\A$ \pcref{def:nsystem_morphism}
\[(a,\glu^a) \fto{\tha} (b,\glu^b).\]
The morphism of $\angordnbe$-systems in $\B$
\begin{equation}\label{hgof_tha}
(\Hgo f)_{\angordnbe} (a,\glu^a) \fto{(\Hgo f)_{\angordnbe} \tha}
(\Hgo f)_{\angordnbe} (b,\glu^b)
\end{equation}
has, for each marker $\ang{s} = \ang{s_j \subseteq \ufs{n}_j}_{j \in \ufs{q}}$, $\ang{s}$-component morphism \cref{theta_angs} defined as the image under $f$ of the $\angs$-component of $\tha$:
\begin{equation}\label{hgof_tha_angs}
fa_{\angs} \fto{((\Hgo f)_{\angordnbe} \tha)_{\ang{s}} \,=\, f\tha_{\angs}} fb_{\angs}.
\end{equation}
\end{description}
This finishes the construction of $\Hgo f$ \cref{hgo_f}.  \cref{hgo_icell_welldef} proves that $\Hgo f$ is a well-defined $G$-natural transformation.

\parhead{Strong variant}.
Suppose $(f,\actf)$ is an $\Op$-pseudomorphism, meaning that each action constraint $\actf_n$ is a $G$-natural isomorphism.  For the pointed $G$-functors $\Hgosg\A = \Asgdash$ and $\Hgosg\B = \Bsgdash$ in \cref{A_ptfunctorGG}, we define a $G$-natural transformation
\begin{equation}\label{hgosg_f}
\begin{tikzpicture}[vcenter]
\def\t{30}
\draw[0cell]
(0,0) node (a1) {\GG}
(a1)++(2,0) node (a2) {\phantom{\GG}}
(a2)++(.12,0) node (a2') {\Catgst}
;
\draw[1cell=.8]
(a1) edge[bend left=\t] node {\Asgdash} (a2)
(a1) edge[bend right=\t] node[swap] {\Bsgdash} (a2)
;
\draw[2cell=.9]
node[between=a1 and a2 at .33, rotate=-90, 2label={above,\Hgosg f}] {\Rightarrow}
;
\end{tikzpicture}
\end{equation}
by restricting \cref{hgo_f_nbe,hgof_aglu,hgof_aglu_comp,hgof_gluing,hgof_tha,hgof_tha_angs} to strong $\angordnbe$-systems for $\angordnbe \in \GG$.  This is well defined because, in the diagram \cref{hgof_gluing}, both $\actf_r$ and $\glu^a_{\cdots}$, and hence also $\glu_{\cdots}$, are now isomorphisms.
\end{definition}

\begin{lemma}[$H$-Theory on 1-Cells]\label{hgo_icell_welldef}\index{H-theory@$H$-theory!1-cell}
For each lax $\Op$-morphism $f \cn \A \to \B$ between $\Op$-pseudoalgebras, the assignment
\begin{equation}\label{Hgof_lemma}
\begin{tikzpicture}[vcenter]
\def\t{30}
\draw[0cell]
(0,0) node (a1) {\GG}
(a1)++(2,0) node (a2) {\phantom{\GG}}
(a2)++(.12,0) node (a2') {\Catgst}
;
\draw[1cell=.8]
(a1) edge[bend left=\t] node {\Adash} (a2)
(a1) edge[bend right=\t] node[swap] {\Bdash} (a2)
;
\draw[2cell=.9]
node[between=a1 and a2 at .33, rotate=-90, 2label={above,\Hgo f}] {\Rightarrow}
;
\end{tikzpicture}
\end{equation}
in \cref{hgo_f} is a $G$-natural transformation.  Moreover, if $f$ is an $\Op$-pseudomorphism, then the assignment
\begin{equation}\label{Hgosgf_lemma}
\begin{tikzpicture}[vcenter]
\def\t{30}
\draw[0cell]
(0,0) node (a1) {\GG}
(a1)++(2,0) node (a2) {\phantom{\GG}}
(a2)++(.12,0) node (a2') {\Catgst}
;
\draw[1cell=.8]
(a1) edge[bend left=\t] node {\Asgdash} (a2)
(a1) edge[bend right=\t] node[swap] {\Bsgdash} (a2)
;
\draw[2cell=.9]
node[between=a1 and a2 at .33, rotate=-90, 2label={above,\Hgosg f}] {\Rightarrow}
;
\end{tikzpicture}
\end{equation}
in \cref{hgosg_f} is a $G$-natural transformation.
\end{lemma}

\begin{proof}
The $G$-naturality of $\Hgosg f$ follows from the $G$-naturality of $\Hgo f$ by restricting to strong $\angordnbe$-systems.  The explanation in \cite[\ref*{EqK:def:Jgo_pos_obj}, \ref*{EqK:Jgof_obj_welldef}, and \ref*{EqK:Jgof_natural}]{yau-eqk}, restricted to the 1-ary case, proves that $\Hgo f$ is a natural transformation.  It remains to show that $\Hgo f$ is $G$-equivariant.  By \cref{ggcatg_icell_geq}, the $G$-equivariance of $\Hgo f$ means that each of its component pointed functors is $G$-equivariant.  For $\angordnbe = \vstar$ and $\ang{}$, $(\Hgo f)_{\angordnbe}$ is given by, respectively, $1_{\bone}$ and $f$, which are $G$-functors.  Suppose $\angordnbe \in \GG \setminus \{\vstar,\ang{}\}$. 
\begin{description}
\item[Component objects]
To prove that the pointed functor $(\Hgo f)_{\angordnbe}$ \cref{hgo_f_nbe} is $G$-equivariant on objects, we consider $g \in G$ and an $\angordnbe$-system $(a,\glu^a)$ in $\A$ \cref{nsystem}.  Using \cref{ga_scomponentGG}, \cref{hgof_aglu_comp}, and the $G$-equivariance of $f$, the computation \cref{hgof_gaglu_comp} shows that the $\angordnbe$-systems in $\B$
\begin{equation}\label{hgof_gaglu}
(\Hgo f)_{\angordnbe} (g \cdot (a,\glu^a)) \andspace 
g \cdot \big((\Hgo f)_{\angordnbe} (a,\glu^a) \big)
\end{equation}
have the same $\angs$-component object for each marker $\angs = \ang{s_j \subseteq \ufsn_j}_{j \in \ufsq}$.
\begin{equation}\label{hgof_gaglu_comp}
\begin{split}
& [(\Hgo f)_{\angordnbe} (g \cdot (a,\glu^a))]_{\angs} \\
&= f (g \cdot (a,\glu^a))_{\angs} \\
&=  f(ga_{\ginv\angs}) \\
&= gf(a_{\ginv\angs}) \\
&= g [(\Hgo f)_{\angordnbe} (a,\glu^a)]_{\ginv \angs} \\
&= \big[g \cdot \big((\Hgo f)_{\angordnbe} (a,\glu^a) \big)\big]_{\angs}
\end{split}
\end{equation}
\item[Gluing]
Using \cref{ga_gluingGG,hgof_gluing}, the gluing morphisms of the two $\angordnbe$-systems in \cref{hgof_gaglu} at $(x; \angs, k, \ang{s_{k,i}}_{i \in \ufsr})$ are, respectively, the top and bottom boundary composites in the diagram \cref{hgof_gaglu_gluing} in $\B$.
\begin{equation}\label{hgof_gaglu_gluing}
\begin{tikzpicture}[vcenter]
\def\h{5} \def\v{-1.4}
\draw[0cell=.85]
(0,0) node (a11) {f\gaA_r(x; \ang{ga_{\ginv\angs \,\compk\, (\ginv s_{k,i})}}_{i \in \ufsr})}
(a11)++(\h,0) node (a12) {fg\gaA_r(\ginv x; \ang{a_{\ginv\angs \,\compk\, (\ginv s_{k,i})}}_{i \in \ufsr})}
(a11)++(0,\v) node (a21) {\gaB_r(x; \ang{fga_{\ginv\angs \,\compk\, (\ginv s_{k,i})}}_{i \in \ufsr})}
(a12)++(0,\v) node (a22) {fga_{g\inv\angs}}
(a21)++(0,\v) node (a31) {g\gaB_r(\ginv x; \ang{f a_{\ginv\angs \,\compk\, (\ginv s_{k,i})}}_{i \in \ufsr})}
(a22)++(0,\v) node (a32) {gfa_{\ginv\angs}}
(a31)++(0,\v) node (a41) {\phantom{\gaA_r}}
(a41)++(.3*\h,0) node (a41') {\phantom{\gaA_r}}
(a41)++(.3,0) node (a41'') {gf\gaA_r(\ginv x; \ang{a_{\ginv\angs \,\compk\, (\ginv s_{k,i})}}_{i \in \ufsr})}
;
\draw[1cell=.8]
(a12) edge node[sloped] {=} (a41')
(a21) edge node {\actf_r} (a11)
(a11) edge[equal] (a12)
(a12) edge node {fg\glu^a_{\ginv x;\, \ginv\angs,\, k, \ang{\ginv s_{k,i}}_{i \in \ufsr}}} (a22)
(a22) edge[equal] (a32)
(a21) edge[equal] (a31)
(a31) edge node[swap] {g\actf_r} (a41)
(a41'') [rounded corners=2pt] -| node[pos=.8,swap] {gf \glu^a_{\ginv x;\, \ginv\angs,\, k, \ang{\ginv s_{k,i}}_{i \in \ufs{r}}}} (a32)
;
\end{tikzpicture}
\end{equation}
The four equalities in \cref{hgof_gaglu_gluing} follow from the $G$-equivariance of the $G$-functor $f$ \cref{fAB} and the $\Op$-action $G$-functors $\gaA_r$ and $\gaB_r$ \cref{gaAn}.  The upper-left and lower-right triangles commute by the $G$-equivariance of, respectively, the action constraint $\actf_r$ \cref{laxmorphism_constraint} and $f$.  This proves that the two $\angordnbe$-systems in \cref{hgof_gaglu} are equal, so $(\Hgo f)_{\angordnbe}$ is $G$-equivariant on objects. 
\item[Morphisms]
To prove the $G$-equivariance of $(\Hgo f)_{\angordnbe}$ on morphisms of $\angordnbe$-systems in $\A$ \pcref{def:nsystem_morphism}, we reuse the computation \cref{hgof_gaglu_comp} by replacing \cref{ga_scomponentGG,hgof_aglu_comp} with \cref{gtheta_angsGG,hgof_tha_angs}.\qedhere
\end{description}
\end{proof}

\section{$H$-Theory 2-Functors}
\label{sec:hgo_twocells}

This section constructs the (strong) $H$-theory 2-functors
\[\AlglaxO \fto{\Hgo} \GGCatg \andspace \AlgpspsO \fto{\Hgosg} \GGCatg\]
for a $\Tinf$-operad $\Op$ \pcref{as:OpA} and observes that they factor the (strong) $J$-theory 2-functors $\Jgo$ and $\Jgosg$ \pcref{thm:Jgo_twofunctor}.  The object and 1-cell assignments of $\Hgo$ and $\Hgosg$ are given in \cref{A_ptfunctorGG,hgo_icell_welldef}.  

\secoutline
\begin{itemize}
\item \cref{def:hgo_iicell} constructs the 2-cell assignments of $\Hgo$ and $\Hgosg$.
\item \cref{hgo_iicell_welldef} proves that the 2-cell assignments of $\Hgo$ and $\Hgosg$ are well defined.
\item \cref{Hgo_twofunctor} records the 2-functoriality of $\Hgo$ and $\Hgosg$.
\item \cref{jgohgoigst} proves that $\Jgo$ and $\Jgosg$ factor through, respectively, $\Hgo$ and $\Hgosg$.
\end{itemize}

For \cref{def:hgo_iicell}, recall that the 2-cells in $\AlglaxO$ and $\AlgpspsO$ \pcref{oalgps_twocat} are $\Op$-transformations \pcref{def:algtwocells}.  The 2-cells in $\GGCatg$ are $G$-modifications \pcref{def:ggcatg}.  Their components are pointed $G$-natural transformations by \cref{ggcatg_iicell_comp,ggcatg_iicell_modax,ggcatg_iicell_geq}.

\begin{definition}\label{def:hgo_iicell}
Suppose 
\begin{equation}\label{Otr_iicell}

\end{equation}
as follows.   A strong variant is defined in \cref{hgosg_omega}.  For an object \cref{GG_objects}
\[\angordnbe = \sordi{n}{\be}{j}_{j \in \ufsq} 
= \big(\ordn_1^{\be_1}, \ldots, \ordn_q^{\be_q}\big) \in \GG,\]
the $\angordnbe$-component pointed $G$-natural transformation
\begin{equation}\label{hgo_om_nbe}
%
\end{equation}
is defined as follows. 
\begin{description}
\item[Base cases]
If $\angordnbe$ is either the basepoint $\vstar$ or the empty tuple $\ang{}$, then $(\Hgo\omega)_{\angordnbe}$ is defined as, respectively, $1_{1_{\bone}}$ and $\omega \cn f \to h$.  This is well defined by \cref{expl:Otr_pointed} and the base cases in \cref{def:hgo_icell}.
\item[Nonbase cases]
For $\angordnbe \in \GG \setminus \{\vstar,\ang{}\}$ and an $\angordnbe$-system $(a,\glu^a) \in \Aangordnbe$ \cref{Aangordnbe_pt}, the $(a,\glu^a)$-component morphism
\begin{equation}\label{hgoom_a}
(\Hgo f)_{\angordnbe} (a,\glu^a)
\fto{(\Hgo\omega)_{\angordnbe,\, (a,\glu^a)}}
(\Hgo h)_{\angordnbe} (a,\glu^a) \inspace \Bangordnbe
\end{equation}
has, for each marker $\angs = \ang{s_j \subseteq \ufsn_j}_{j \in \ufsq}$, $\angs$-component morphism defined as the $a_{\angs}$-component of $\omega$:
\begin{equation}\label{omega_aangs}
fa_{\angs} \fto{(\Hgo\omega)_{\angordnbe,\, (a,\glu^a), \angs} \,=\, \omega_{a_{\angs}}} ha_{\angs}.
\end{equation}
This morphism in $\B$ is well defined by \cref{hgof_aglu_comp}.
\end{description}

\parhead{Strong variant}.
Suppose $(f,\actf)$ and $(h,\acth)$ are $\Op$-pseudomorphisms.  For the $G$-natural transformations $\Hgosg f$ and $\Hgosg h$ in \cref{hgo_icell_welldef}, we define a $G$-modification
\begin{equation}\label{hgosg_omega}
\begin{tikzpicture}[vcenter]
\def\h{4} \def\u{.5} \def\t{23}
\draw[0cell]
(0,0) node (a1) {\GG}
(a1)++(\h/2,0) node (a2) {\phantom{A}}
(a1)++(\h,0) node (a3) {\phantom{\GG}}
(a3)++(.1,0) node (a3') {\Catgst}
;
\draw[1cell=.8]
(a1) edge[bend left=\t] node (A) {\Asgdash} (a3)
(a1) edge[bend right=\t] node[swap] (B) {\Bsgdash} (a3)
;
\draw[2cell=.8]
node[between=a1 and a2 at .7, rotate=-90, 2labelmed={below,\Hgosg f}] {\Rightarrow}
node[between=a2 and a3 at .25, rotate=-90, 2label={above,\Hgosg h}] {\Rightarrow}
;
\draw[2cell=.9]
node[between=a1 and a3 at .5, shift={(0,-.25)},  2labelmed={above,\Hgosg\omega}] {\Rrightarrow}
;
\end{tikzpicture}
\end{equation}
by restricting \cref{hgo_om_nbe,hgoom_a,omega_aangs} to strong $\angordnbe$-systems for $\angordnbe \in \GG$.  This is well defined because, by \cref{Aangordnsg,Aangordnbe_pt,Asgangordnbe}, $\Bsgangordnbe$ is a full subcategory of $\Bangordnbe$.
\end{definition}

\begin{lemma}[$H$-Theory on 2-Cells]\label{hgo_iicell_welldef}\index{H-theory@$H$-theory!2-cell}
For each $\Op$-transformation $\omega \cn f \to h$ between lax $\Op$-morphisms $f,h \cn \A \to \B$ between $\Op$-pseudoalgebras, the assignment
\[\Hgo f \fto{\Hgo \omega} \Hgo h\]
in \cref{hgo_omega} is a $G$-modification.  Moreover, if $f$ and $h$ are $\Op$-pseudomorphisms, then the assignment
\[\Hgosg f \fto{\Hgosg \omega} \Hgosg h\]
in \cref{hgosg_omega} is a $G$-modification.
\end{lemma}

\begin{proof}
It suffices to prove the first assertion for $\Hgo\omega$ because $\Hgosg\omega$ is defined by restricting to the full subcategories of strong $\angordnbe$-systems.   The explanation in \cite[\ref*{EqK:def:Jgo_pos_mor}, \ref*{EqK:Jgotheta_system_mor}, \ref*{EqK:Jgotheta_m_natural}, and \ref*{EqK:Jgotheta_modification}]{yau-eqk}, restricted to the 1-ary case, proves that $\Hgo \omega$ is a modification.  It remains to show that $\Hgo\omega$ is $G$-equivariant.

\parhead{$G$-equivariance}. 
By \cref{ggcatg_iicell_geq}, the $G$-equivariance of $\Hgo\omega$ means that its component pointed natural transformations are $G$-equivariant.  
\begin{itemize}
\item For $\angordnbe = \vstar$ and $\ang{}$, $(\Hgo\omega)_{\angordnbe}$ is given by, respectively, $1_{1_{\bone}}$ and $\omega$, which are $G$-natural transformations.  
\item For $\angordnbe \in \GG \setminus \{\vstar,\ang{}\}$, to prove that $(\Hgo\omega)_{\angordnbe}$ \cref{hgo_om_nbe} is $G$-equivariant, we reuse the computation \cref{hgof_gaglu_comp} with $\omega$ in place of $f$, along with \cref{ga_scomponentGG,gtheta_angsGG,omega_aangs}.\qedhere
\end{itemize}
\end{proof}

\begin{proposition}\label{Hgo_twofunctor}\index{H-theory@$H$-theory}
For each $\Tinf$-operad $\Op$ \pcref{as:OpA}, the object, 1-cell, and 2-cell assignments in \cref{A_ptfunctorGG,hgo_icell_welldef,hgo_iicell_welldef} define 2-functors
\[\begin{split}
\AlglaxO & \fto{\Hgo} \GGCatg \andspace\\
\AlgpspsO & \fto{\Hgosg} \GGCatg
\end{split}\]
between the 2-categories in \cref{oalgps_twocat,def:ggcatg}.
\end{proposition}

\begin{proof}
The explanation in \cite[\ref*{EqK:Jgo_gamma} and \ref*{EqK:thm:Jgo_multifunctor}]{yau-eqk}, restricted to the 1-ary case, proves that $\Hgo$ and $\Hgosg$ are 2-functors.
\end{proof}

The 2-functors $\Hgo$ and $\Hgosg$ are called \emph{$H$-theory for $\Op$} and \emph{strong $H$-theory for $\Op$}.

\subsection*{$J$-Theory Factors through $H$-Theory}

Recall the 2-equivalence \pcref{thm:ggcat_ggcatg_iieq}
\[\GGCatg \fto[\sim]{\igst} \GGCatii\]
induced by precomposition with the pointed full subcategory inclusion $\ig \cn \Gsk \to \GG$ \cref{ig}.  \cref{jgohgoigst} states that the $J$-theory 2-functor $\Jgo$ \pcref{thm:Jgo_twofunctor} factors as the composite of $\Hgo$ \pcref{Hgo_twofunctor} and $\igst$.  The strong variant is also true.

\begin{proposition}\label{jgohgoigst}
For each $\Tinf$-operad $\Op$ \pcref{as:OpA}, the following diagrams of 2-functors commute.
\begin{equation}\label{HistJ}
\begin{tikzpicture}[vcenter]
\def\h{2} \def\v{1.4}
\draw[0cell]
(0,0) node (a1) {\phantom{A^l}}
(-.5,-.1) node (a1') {\AlglaxO}
(a1)++(\h,\v/2) node (a2) {\GGCatg}
(a2)++(-.4,0) node (a2') {\phantom{A}}
(a1)++(\h,-\v/2) node (a3) {\GGCatii}
(a3)++(-.4,0) node (a3') {\phantom{A}}
(a1)++(2.5*\h,0) node (b1) {\phantom{A^l}}
(b1)++(-.45,0) node (b1') {\AlgpspsO}
(b1)++(\h,\v/2) node (b2) {\GGCatg}
(b2)++(-.4,0) node (b2') {\phantom{A}}
(b1)++(\h,-\v/2) node (b3) {\GGCatii}
(b3)++(-.4,0) node (b3') {\phantom{A}}
;
\draw[1cell]
(a1) edge node {\Hgo} (a2')
(a2) edge node {\igst} node[swap,sloped] {\sim} (a3)
(a1) edge node[swap] {\Jgo} (a3')
(b1) edge node {\Hgosg} (b2')
(b2) edge node {\igst} node[swap,sloped] {\sim} (b3)
(b1) edge node[swap] {\Jgosg} (b3')
;
\end{tikzpicture}
\end{equation}
\end{proposition}

\begin{proof}
The 2-functors
\begin{equation}\label{Jgo_igst_Hgo}
\Jgo \andspace \igst \Hgo \cn \AlglaxO \to \GGCatii
\end{equation}
are equal on objects, 1-cells, and 2-cells of $\AlglaxO$ \pcref{oalgps_twocat} for the following reasons.
\begin{description}
\item[Objects] 
For an $\Op$-pseudoalgebra $\A$ \pcref{def:pseudoalgebra}, the pointed functors
\[\Jgo\A \andspace \igst\Hgo\A \cn \Gsk \to \Gcatst\]
are equal on objects of $\Gsk$ by \cref{expl:nbeta_gcat}.  These two pointed functors are equal on morphisms of $\Gsk$ by \cref{Aupom}.
\item[1-cells]
The two 2-functors in \cref{Jgo_igst_Hgo} are equal on each lax $\Op$-morphism $f$ between $\Op$-pseudoalgebras \pcref{def:laxmorphism} by
\begin{itemize}
\item the definition of $\Jgo f$ in \cref{Jgof_m_obj_comp,Jgof_m_gluing,Jgof_m_theta_comp} and
\item the definition of $\Hgo f$ in \cref{hgof_aglu_comp,hgof_gluing,hgof_tha_angs}.
\end{itemize}
\item[2-cells]
The two 2-functors in \cref{Jgo_igst_Hgo} are equal on each $\Op$-transformation \pcref{def:algtwocells} between lax $\Op$-morphisms between $\Op$-pseudoalgebras by \cref{Jgotheta_m_angs,omega_aangs}.
\end{description}
The same argument also proves the equality
\[\Jgosg = \igst \Hgosg\]
by restricting to strong systems.
\end{proof}

\section{$\GGG$-Spaces from $\GGG$-Categories}
\label{sec:GGTopg}

This section constructs the passage from $\GGG$-categories to $\GGG$-spaces.  Throughout this section, $G$ denotes an arbitrary group.

\secoutline
\begin{itemize}
\item \cref{def:GGTopg} defines the categories $\GGTopg$ of $\GGG$-spaces, with further elaboration given in \cref{expl:GGTopg}.
\item \cref{GGCatg_GGTopg} constructs the functor from $\GGCatg$ to $\GGTopg$ induced by the classifying space functor $\cla$, with further elaboration given in \cref{expl:GGCatg_GGTopg}.
\item \cref{thm:ggtop_ggtopg_iieq} proves that the full subcategory inclusion $\ig \cn \Gsk \to \GG$ \cref{ig} induces an equivalence between the categories $\GGTopg$ and $\GGTopii$.
\end{itemize}

\subsection*{$\GGG$-Spaces}
Recall the complete and cocomplete symmetric monoidal closed category \cref{Gtopst_smc}
\[(\Gtopst, \sma, \stplus, \Topgst)\]
of pointed $G$-spaces, with internal hom pointed $G$-spaces $\Topgst(X,Y)$.  The notation $\Topgst$ \cref{topgst_gtopst_enr} also denotes the pointed $G$-category of pointed $G$-spaces and pointed morphisms with the conjugation $G$-action \cref{ginv_h_g}.  Recall the pointed $G$-category $\GG$ \pcref{def:GG}.  \cref{def:GGTopg} defines the topological analogue of the 2-category $\GGCatg$ of $\GGG$-categories \pcref{def:ggcatg}.

\begin{definition}\label{def:GGTopg}
Suppose $G$ is a group.  The category $\GGTopg$ has pointed $G$-functors
\begin{equation}\label{GGTopg_obj}
(\GG,\vstar) \fto{X} (\Topgst,*)
\end{equation}
as objects, called \index{GGG-space@$\GGG$-space}\emph{$\GGG$-spaces}, and $G$-natural transformations
\begin{equation}\label{GGTopg_icell}
\begin{tikzpicture}[vcenter]
\def\t{28}
\draw[0cell]
(0,0) node (a1) {\phantom{\Gsk}}
(a1)++(1.8,0) node (a2) {\phantom{\Gsk}}
(a1)++(-.08,0) node (a1') {\GG}
(a2)++(.23,0) node (a2') {\Topgst}
;
\draw[1cell=.9]
(a1) edge[bend left=\t] node {X} (a2)
(a1) edge[bend right=\t] node[swap] {X'} (a2)
;
\draw[2cell]
node[between=a1 and a2 at .45, rotate=-90, 2label={above,\theta}] {\Rightarrow}
;
\end{tikzpicture}
\end{equation}
as morphisms.  Identities and composition are those of natural transformations.
\end{definition}

\begin{explanation}[Unpacking $\GGTopg$]\label{expl:GGTopg}
The category $\GGTopg$ in \cref{def:GGTopg} is given explicitly as follows.
\begin{description}
\item[Objects]
A $\GGG$-space $X \cn \GG \to \Topgst$ \cref{GGTopg_obj} consists of the following data.
\begin{description}
\item[Pointed $G$-spaces] 
$X$ sends each object $\malp \in \GG$ \cref{GG_objects} to a pointed $G$-space $X\malp$ such that $X\vstar = *$.  Its $G$-fixed basepoint is given by the pointed morphism
\[X(\vstar \to \malp) \cn X\vstar = * \to X\malp\]
for the $G$-fixed unique morphism $\vstar \to \malp$ in $\GG$.  
\item[Pointed morphisms] 
$X$ sends each morphism $\upom \cn \malp \to \nbe$ in $\GG$ \cref{GG_morphisms} to a pointed morphism
\begin{equation}\label{X_upom_GGTopg}
X\malp \fto{X\upom} X\nbe
\end{equation}
such that $X$ preserves identity morphisms and composition.  The morphism $X\upom$ is \emph{not} generally $G$-equivariant.
\item[Equivariance] 
The $G$-equivariance of $X$ means the equality of morphisms
\begin{equation}\label{GGGspace_equiv}
X(g \cdot \upom) = g (X\upom) \ginv
\end{equation}
for each $g \in G$ and each morphism $\upom$ in $\GG$.  The morphism $g \cdot \upom$ is the $g$-action on $\upom$ \cref{GG_Gaction}, and $g (X\upom) \ginv$ is the conjugation $g$-action \cref{ginv_h_g} on $X\upom$.  Thus, if the morphism $\upom$ is $G$-fixed, then $X\upom$ is a pointed $G$-morphism.
\end{description}
\item[Morphisms] 
A morphism $\tha \cn X \to X'$ in $\GGTopg$ \cref{GGTopg_icell} consists of, for each object $\malp \in \GG$, an $\malp$-component pointed morphism
\begin{equation}\label{GGTopg_icell_comp}
X\malp \fto{\tha_{\malp}} X'\malp
\end{equation}
such that, for each morphism $\upom \cn \malp \to \nbe$ in $\GG$, the following naturality diagram of pointed morphisms commutes.
\begin{equation}\label{GGTopg_icell_nat}
\begin{tikzpicture}[vcenter]
\def\v{-1.4}
\draw[0cell]
(0,0) node (a11) {X\malp}
(a11)++(2.5,0) node (a12) {X'\malp}
(a11)++(0,\v) node (a21) {X\nbe}
(a12)++(0,\v) node (a22) {X'\nbe}
;
\draw[1cell=.9]
(a11) edge node {\theta_{\malp}} (a12)
(a12) edge node {X'\upom} (a22)
(a11) edge node[swap] {X\upom} (a21)
(a21) edge node {\theta_{\nbe}} (a22)
;
\end{tikzpicture}
\end{equation}
The $G$-equivariance of $\tha$ means the equality of morphisms
\begin{equation}\label{GGTopg_icell_geq}
\tha_{\malp} = g \tha_{\malp} \ginv
\end{equation}
for each $g \in G$ and object $\malp \in \GG$.   In other words, each component of a morphism in $\GGTopg$ is a pointed $G$-morphism.  A morphism $\tha$ is automatically pointed, meaning $\tha_{\vstar} = 1_{*}$.  
\end{description}
Identity morphisms and composition are defined componentwise using the components in \cref{GGTopg_icell_comp}.
\end{explanation}

\subsection*{$\GGG$-Spaces from $\GGG$-Categories}

\cref{thm:ggcat_ggtop} describes the passage from $\Gskg$-categories to $\Gskg$-spaces via the classifying space functor.  \cref{GGCatg_GGTopg} is the analogue involving $\GG$.

\begin{lemma}\label{GGCatg_GGTopg}
For a group $G$, composing and whiskering with the classifying space functor $\cla$ induce a functor\index{classifying space}
\begin{equation}\label{clast_GGCatg}
\GGCatg \fto{\clast} \GGTopg 
\end{equation}
between the categories in \cref{def:ggcatg,def:GGTopg}.
\end{lemma}

\begin{proof}
By functoriality, the classifying space functor $\cla \cn \Cat \to \Top$ \cref{classifying_space}, from small categories to spaces, preserves the conjugation $G$-actions on morphisms.  Thus, it yields a pointed $G$-functor
\begin{equation}\label{cla_Catgst_Topgst}
\Catgst \fto{\cla} \Topgst
\end{equation}
between the pointed $G$-categories in \cref{Catgst_iicat,topgst_gtopst_enr}.  The functoriality of $\clast$ in \cref{clast_GGCatg} follows from the fact that whiskering with $\cla$ \cref{cla_Catgst_Topgst} preserves identities and composition of natural transformations.
\end{proof}

\begin{explanation}[Unpacking $\clast$]\label{expl:GGCatg_GGTopg}
The functor $\clast$ in \cref{clast_GGCatg} sends a pointed $G$-functor $X \cn \GG \to \Catgst$ and a $G$-natural transformation $\tha$ to the following composite pointed $G$-functor and whiskered $G$-natural transformation.
\begin{equation}\label{clast_GGCatg_mor}
\begin{tikzpicture}[baseline={(a1.base)}]
\def\t{28}
\draw[0cell]
(0,0) node (a1) {\GG}
(a1)++(1.8,0) node (a2) {\phantom{\GG}}
(a2)++(.15,0) node (a2') {\Catgst}
(a2')++(2,0) node (a3) {\Topgst}
;
\draw[1cell=.9]
(a1) edge[bend left=\t] node {X} (a2)
(a1) edge[bend right=\t] node[swap] {X'} (a2)
(a2') edge node {\cla} (a3)
;
\draw[2cell]
node[between=a1 and a2 at .43, rotate=-90, 2label={above,\tha}] {\Rightarrow}
;
\end{tikzpicture}
\end{equation}
For each object $\malp \in \GG$, the $\malp$-component of $\clast\tha$ is the pointed $G$-morphism between pointed $G$-spaces
\[\cla X\malp \fto{\cla\tha_{\malp}} \cla X'\malp\]
obtained from $\tha_{\malp}$ \cref{GGTopg_icell_comp} by applying $\cla$.
\end{explanation}

\subsection*{$\Gskg$-Spaces and $\GGG$-Spaces}

Recall the full subcategory inclusion $\ig \cn \Gsk \to \GG$ \cref{ig}.  \cref{thm:ggtop_ggtopg_iieq} is used in \cref{KgKgg} to compare the functor $\Kg \cn \GGTopii \to \Gspec$ \cref{Kg_functor} and a variant $\Kgg \cn \GGTopg \to \Gspec$ constructed in \cref{sec:Kgg}.

\begin{lemma}\label{thm:ggtop_ggtopg_iieq}
For each group $G$, there is an adjoint equivalence\index{adjoint equivalence}
\begin{equation}\label{Ligst_adeq}
\begin{tikzpicture}[vcenter]
\draw[0cell]
(0,0) node (a1) {\GGTopii}
(a1)++(2.5,0) node (a2) {\GGTopg}
;
\draw[1cell=.9]
(a1) edge[transform canvas={yshift=.5ex}] node {\Lg} (a2)
(a2) edge[transform canvas={yshift=-.4ex}] node {\igst} (a1)
;
\end{tikzpicture}
\end{equation}
between the categories in \cref{def:ggtop_smc,def:GGTopg}, in which the right adjoint $\igst$ is induced by $\ig$.
\end{lemma}

\begin{proof}
We reuse the proofs of \cref{igst_iifunctor,thm:ggcat_ggcatg_iieq} by replacing \cref{Gcatst_smc} 
\[(\Gcatst, \sma, \bonep, \Catgst)\] 
with the topological analogue 
\[(\Gtopst, \sma, \stplus, \Topgst)\] 
defined in \cref{Gtopst_smc}.
\end{proof}

\section{Orthogonal $G$-Spectra from $\GGG$-Spaces}
\label{sec:Kgg}

This section constructs the prolongation functor
\[\GGTopg \fto{\Kgg} \Gspec\]
that sends $\GGG$-spaces to orthogonal $G$-spectra and compares it with the prolongation functor $\Kg \cn \GGTopii \to \Gspec$.  The functor $\Kgg$ is used in \cref{sec:kgo_hthy} to construct equivalent variants of the functors $\Kgo$ and $\Kgosg$ \cref{Kgo_functors}.

\secoutline
\begin{itemize}
\item \cref{def:gggspace_gspectra} constructs the orthogonal $G$-spectrum $\Kgg X$ associated to a $\GGG$-space $X$.
\item \cref{kggx_welldef} proves that $\Kgg X$ is well defined.
\item \cref{def:Kgg_mor} constructs the $G$-morphism $\Kgg\tha$ between orthogonal $G$-spectra associated to a $G$-natural transformation $\tha$ between $\GGG$-spaces.
\item \cref{kggpsi_welldef} proves that $\Kgg\tha$ is well defined.
\item \cref{def:Kgg_functor} defines the functor $\Kgg$ from the category $\GGTopg$ of $\GGG$-spaces to the category $\Gspec$ of orthogonal $G$-spectra.
\item \cref{KgKgg} proves that the functors $\Kg$ and $\Kgg$ agree up to a natural isomorphism.
\end{itemize}

\subsection*{Object Assignment of $\Kgg$}

\cref{def:gggspace_gspectra} is the $\GG$-analogue of \cref{def:ggspace_gspectra}, which defines the orthogonal $G$-spectra associated to $\Gskg$-spaces.

\begin{definition}[$\Kgg$ on Objects]\label{def:gggspace_gspectra}
Given a compact Lie group $G$ and a $\GGG$-space \cref{GGTopg_obj}
\[(\GG,\vstar) \fto{X} (\Topgst,*),\]
the \index{orthogonal G-spectrum@orthogonal $G$-spectrum!multifunctorial K-theory@multifunctorial $K$-theory}\index{multifunctorial K-theory@multifunctorial $K$-theory!orthogonal G-spectrum@orthogonal $G$-spectrum}orthogonal $G$-spectrum \pcref{def:gsp_module}
\begin{equation}\label{Kgg_object}
(\Kgg X, \umu) \in \Gspec
\end{equation}
is defined as follows.
\begin{description}
\item[Object assignment of $\IU$-space]
The $\IU$-space \pcref{def:iu_space} $\Kgg X$ sends each object $V \in \IU$ to the coend \pcref{def:coends}
\begin{equation}\label{Kggxv}
(\Kgg X)_V = \int^{\malp \in \GG} (S^V)^{\sma\malp} \sma X\malp.
\end{equation}
\begin{description}
\item[Coend] 
The coend in \cref{Kggxv} is taken in the category $\Topst$ of pointed spaces and pointed morphisms.  The pointed finite $G$-set $\sma\malp\in \FG$ \cref{smash_GGobjects} is regarded as a discrete pointed $G$-space, and $S^V$ is the $V$-sphere \pcref{def:g_sphere}.  The pointed $G$-space 
\begin{equation}\label{SVmalp}
(S^V)^{\sma\malp} = \Topgst(\sma\malp, S^V)
\end{equation}
consists of pointed morphisms $\sma\malp \to S^V$ \cref{Gtopst_smc}, with $G$ acting by conjugation \cref{ginv_h_g}.
\item[$G$-action] 
The group $G$ acts diagonally on representatives.  This means that, for an element $g \in G$ and a representative pair
\begin{equation}\label{Kggxv_rep}
\big(\smam\malp \fto{\upom} S^V ; x \in X\malp \big) \in (S^V)^{\sma\malp} \ttimes X\malp
\end{equation}
in $(\Kgg X)_V$, the diagonal $g$-action is given by
\begin{equation}\label{Kggxv_rep_gact}
g \cdot (\upom; x) = (g\upom\ginv ; gx),
\end{equation}
where $g\upom\ginv$ means the composite pointed morphism
\[\sma\malp \fto{\ginv} \smam\malp \fto{\upom} S^V \fto{g} S^V.\]
\end{description}
\item[Morphism assignment of $\IU$-space]
For a linear isometric isomorphism $f \cn V \fiso W$ in $\IU$, the pointed homeomorphism \cref{iu_space_xf}
\begin{equation}\label{Kggxf}
(\Kgg X)_V \fto[\iso]{(\Kgg X)_f} (\Kgg X)_W
\end{equation}
is induced by the pointed homeomorphisms 
\[(S^V)^{\sma\malp} \fto[\iso]{f \circ -} (S^W)^{\sma\malp} \forspace \malp \in \GG\]
that postcompose with the pointed homeomorphism $f \cn S^V \fiso S^W$.    In terms of representatives \cref{Kggxv_rep}, it is given by
\begin{equation}\label{Kggxf_rep}
(\Kgg X)_f(\upom; x) = (f\upom; x).
\end{equation}
\item[Sphere action]
For each pair of objects $(V,W) \in (\IUsk)^2$, the $(V,W)$-component pointed $G$-morphism \cref{gsp_action_vw} is defined by the following commutative diagram in $\Gtopst$.
\begin{equation}\label{Kggx_action_vw}
\begin{tikzpicture}[vcenter]
\def\h{5} \def\u{-1} \def\v{-1.4}
\draw[0cell=.9]
(0,0) node (a11) {(\Kgg X)_V \sma S^W}
(a11)++(\h,0) node (a12) {(\Kgg X)_{V \oplus W}}
(a11)++(0,\u) node (a21) {\big( \txint^{\malp \in \GG} (S^V)^{\sma\malp} \sma X\malp \big) \sma S^W}
(a12)++(0,\u) node (a22) {\txint^{\malp \in \GG} (S^{V \oplus W})^{\sma\malp} \sma X\malp}
(a21)++(0,\v) node (a3) {\txint^{\malp \in \GG} \big( (S^V)^{\sma\malp} \sma S^W \big) \sma X\malp}
;
\draw[1cell=.8]
(a11) edge node {\umu_{V,W}} (a12)
(a11) edge[equal] (a21)
(a12) edge[equal] (a22)
(a21) edge node[swap] {\iso} (a3)
(a3) [rounded corners=2pt] -| node[pos=.25] {\asm = \txint\! \asm_{\malp} \sma 1} (a22)
;
\end{tikzpicture}
\end{equation}
\begin{itemize}
\item The pointed $G$-homeomorphism denoted by $\iso$ first commutes $- \sma S^W$ with the coend.  Then it moves $S^W$ to the left of $X\malp$ using the associativity isomorphism and braiding for the symmetric monoidal category $(\Gtopst,\sma)$ \cref{Gtopst_smc}.
\item The pointed $G$-morphism $\asm$ is induced by the pointed $G$-morphisms
\begin{equation}\label{assembly_sph}
(S^V)^{\sma\malp} \sma S^W \fto{\asm_{\malp}} (S^{V \oplus W})^{\sma\malp}
\end{equation}
for $\malp \in \GG$ defined by the assignment
\[\begin{split}
& \big(\smam\malp \fto{\upom} S^V ; y \big) \in (S^V)^{\sma\malp} \sma S^W\\
&\mapsto \big(\smam\malp \fto{\upom} S^V \fto{- \oplus y} S^{V \oplus W} \big) \in (S^{V \oplus W})^{\sma\malp}.
\end{split}\]
In other words, 
\[\asm_{\malp}(\upom; y) = \upom \oplus y\]
sends an element $\bdi \in \sma\malp$ to the point 
\begin{equation}\label{assembly_explicit}
\big( \asm_{\malp}(\upom; y)\big)(\bdi) 
= (\upom \bdi) \oplus y \in S^{V \oplus W}.
\end{equation}
For a representative pair $(\upom; x)$ of $(\Kgg X)_V$ \cref{Kggxv_rep} and a point $y \in S^W$, $\umu_{V,W}$ is given by
\begin{equation}\label{Kggx_actrep}
\umu_{V,W}\big((\upom; x); y\big) = (\upom \oplus y; x).
\end{equation}
\end{itemize}
\end{description}
\cref{kggx_welldef} proves that $(\Kfg X,\umu)$ is well defined.
\end{definition}

\begin{lemma}\label{kggx_welldef}
The pair $(\Kgg X,\umu)$ in \cref{Kgg_object} is an orthogonal $G$-spectrum.
\end{lemma}

\begin{proof}
Using \cref{expl:iu_space}, we first prove that $\Kgg X$, defined in \crefrange{Kggxv}{Kggxf}, is an $\IU$-space.
\begin{description}
\item[$G$-action]
To see that the pointed space $(\Kgg X)_V$ \cref{Kggxv} is a pointed $G$-space, we observe that the $g$-action \cref{Kggxv_rep_gact} is well defined as follows.  Recall that the defining relation of the coend $(\Kgg X)_V$ identifies the two representative pairs
\begin{equation}\label{Kggxv_relation}
\big(\upom(\sma\upla) ; x\big) \sim \big(\upom ; (X\upla)x\big)
\end{equation}
for 
\begin{itemize}
\item a morphism $\nbe \fto{\upla} \malp$ in $\GG$,
\item a pointed morphism $\sma\malp \fto{\upom} S^V$, and
\item a point $x \in X\nbe$.
\end{itemize}
The following computation proves that the diagonal $g$-actions on the two sides of \cref{Kggxv_relation} are equal in $(\Kgg X)_V$.
\begin{equation}\label{Kggxv_gactwell}
\begin{aligned}
& g \cdot (\upom(\sma\upla) ; x) &&\\
&= (g\upom (\sma\upla) \ginv ; gx) && \text{by \cref{Kggxv_rep_gact}}\\
&= \big((g\upom\ginv)(g (\sma\upla) \ginv) ; gx \big) && \text{by $\ginv g = 1$}\\
&= \big((g\upom\ginv) \smam(g\cdot\upla) ; gx \big) && \text{by \cref{sma-fgpsi}}\\
&= \big((g\upom\ginv) ; (X(g\cdot\upla))(gx) \big) && \text{by \cref{Kggxv_relation}}\\
&= \big(g\upom\ginv ; (g(X\upla)\ginv)(gx) \big) && \text{by \cref{GGGspace_equiv}}\\
&= (g\upom\ginv ; g(X\upla)x) && \text{by $\ginv g = 1$}\\
&= g \cdot (\upom ; (X\upla)x) && \text{by \cref{Kggxv_rep_gact}}
\end{aligned}
\end{equation}
This proves that $(\Kgg X)_V$ \cref{Kggxv} is a pointed $G$-space.
\item[Functoriality]
The assignment  \cref{Kggxf} 
\[f \mapsto (\Kgg X)_f\]
preserves identities and composition in the sense of \cref{iu_space_axioms} because $(\Kgg X)_f$ is induced by postcomposition with $f$ \cref{Kggxf_rep}.
\item[Equivariance]
For an element $g \in G$, a linear isometric isomorphism $f \cn V \fiso W$ in $\IU$, and a representative pair $(\upom; x)$ of $(\Kgg X)_V$ \cref{Kggxv_rep}, the equivariance diagram \cref{x_gfginv} for $\Kgg X$ commutes by the following equalities in $(\Kgg X)_W$.
\begin{equation}\label{Kggx_equivar}
\begin{aligned}
& g \cdot (\Kgg X)_f \big( \ginv \cdot (\upom; x)\big) &&\\
&= g \cdot (\Kgg X)_f (\ginv\upom g; \ginv x) && \text{by \cref{Kggxv_rep_gact}}\\
&= g \cdot (f\ginv\upom g; \ginv x) && \text{by \cref{Kggxf_rep}}\\
&= (gf\ginv\upom g \ginv; g \ginv x) && \text{by \cref{Kggxv_rep_gact}}\\
&= (gf\ginv\upom; x) && \text{by $g\ginv = 1$}\\
&= (\Kgg X)_{gf\ginv} (\upom; x) && \text{by \cref{Kggxf_rep}}
\end{aligned}
\end{equation}
This proves that $\Kgg X$ is an $\IU$-space.
\end{description}

Next, we prove that $(\Kgg X,\umu)$ satisfies the defining conditions \crefrange{gsp_action_vw}{gsp_assoc} for an orthogonal $G$-spectrum.
\begin{description}
\item[$G$-equivariance of sphere action]
To show that the $(V,W)$-component pointed morphism $\umu_{V,W}$ \cref{Kggx_action_vw} is $G$-equivariant, it suffices to prove that the pointed morphism $\asm_{\malp}$ \cref{assembly_sph} is $G$-equivariant for each object $\malp \in \GG$.  The group $G$ acts diagonally on the smash product $(S^V)^{\sma\malp} \sma S^W$ and by conjugation on $(S^V)^{\sma\malp}$ and $(S^{V \oplus W})^{\sma\malp}$.  For an element $g \in G$, an element $\bdi \in \sma\malp$, a pointed morphism $\upom \cn \msmam\malp \to S^V$, and a point $y \in S^W$, the following equalities in $S^{V \oplus W}$ prove that $\asm_{\malp}$ is $G$-equivariant.
\begin{equation}\label{Kggx_sph_eq}
\begin{aligned}
& [g \cdot (\asm_{\malp} (\upom; y))] (\bdi) &&\\
&= g [(\asm_{\malp} (\upom; y) ) (\ginv \bdi) ] && \text{by \cref{SVmalp}}\\
&= g (\upom(\ginv \bdi) \oplus y) && \text{by \cref{assembly_explicit}}\\
&= g\upom(\ginv \bdi) \oplus gy && \text{by diagonal action}\\
&= [\asm_{\malp} (g\upom\ginv; gy)] (\bdi) && \text{by \cref{assembly_explicit}}\\
&= [\asm_{\malp} (g \cdot (\upom; y))] (\bdi) && \text{by diagonal action}
\end{aligned}
\end{equation}
This proves that $\umu_{V,W}$ \cref{Kggx_action_vw} is a pointed $G$-morphism.
\item[Naturality]
For linear isometric isomorphisms $f \cn V \fiso V'$ and $h \cn W \fiso W'$ in $\IUsk$, a representative pair $(\upom; x)$ of $(\Kgg X)_V$ \cref{Kggxv_rep}, and a point $y \in S^W$, the following equalities in $(\Kgg X)_{V' \oplus W'}$ prove that the naturality diagram \cref{gsp_action_nat} for $\Kgg X$ commutes.
\[\begin{aligned}
& (\umu_{V',W'}) \big((\Kgg X)_f \sma h \big) ((\upom; x); y) &&\\
&= (\umu_{V',W'}) \big((f\upom; x) ; hy\big) && \text{by \cref{Kggxf_rep}} \\
&= (f\upom \oplus hy; x) && \text{by \cref{Kggx_actrep}}\\
&= \big((f \oplus h)(\upom \oplus y) ; x\big) &&\\
&= (\Kgg X)_{f \oplus h} (\upom \oplus y ; x) && \text{by \cref{Kggxf_rep}}\\
&= (\Kgg X)_{f \oplus h} (\umu_{V,W}) ((\upom; x) ; y) && \text{by \cref{Kggx_actrep}}
\end{aligned}\]
\item[Unity]
For a representative pair $(\upom; x)$ of $(\Kgg X)_V$ \cref{Kggxv_rep} and the nonbasepoint $0 \in S^0$, the following equalities in $(\Kgg X)_V$ prove that the unity diagram \cref{gsp_unity} for $\Kgg X$ commutes, where $\rho \cn V \oplus 0 \fiso V$ is the right unit isomorphism.
\[\begin{aligned}
& (\Kgg X)_{\rho} (\umu_{V,0}) ((\upom; x); 0) &&\\
&= (\Kgg X)_{\rho} (\upom \oplus 0; x) && \text{by \cref{Kggx_actrep}}\\
&= (\rho(\upom \oplus 0); x) && \text{by \cref{Kggxf_rep}}\\
&= (\upom; x)
\end{aligned}\]
\item[Associativity]
For a representative pair $(\upom; x)$ of $(\Kgg X)_U$ \cref{Kggxv_rep}, $v \in S^V$, and $y \in S^W$, the following equalities in $(\Kgg X)_{U \oplus (V \oplus W)}$ prove that the associativity diagram \cref{gsp_assoc} for $\Kgg X$ commutes, where 
\[(U \oplus V) \oplus W \fto[\iso]{\al} U \oplus (V \oplus W)\]
is the associativity isomorphism.
\[\begin{aligned}
& (\Kgg X)_\al (\umu_{U \oplus V,W}) (\umu_{U,V} \sma 1) \big((\upom; x) ; v ; y \big) &&\\
&= (\Kgg X)_\al (\umu_{U \oplus V,W}) \big((\upom \oplus v; x) ; y \big) && \text{by \cref{Kggx_actrep}}\\
&= (\Kgg X)_\al \big((\upom \oplus v) \oplus y ; x\big) && \text{by \cref{Kggx_actrep}}\\
&= \big(\upom \oplus (v \oplus y) ; x\big) && \text{by \cref{Kggxf_rep}}\\
&= (\umu_{U, V \oplus W}) \big((\upom; x) ; v \oplus y \big) && \text{by \cref{Kggx_actrep}}\\
&= (\umu_{U, V \oplus W}) (1 \sma \mu_{V,W}) \big((\upom; x) ; v ; y \big) && \text{by \cref{gsp_mult}}
\end{aligned}\]
\end{description}
This proves that $(\Kgg X,\umu)$ is an orthogonal $G$-spectrum.
\end{proof}

\begin{explanation}[$\Kg$ vs. $\Kgg$]\label{expl:Kggxv}
There is a subtle difference between the pointed $G$-spaces
\[\begin{split}
(\Kg X)_V &= \int^{\angordn \in \Gsk} (S^V)^{\sma\angordn} \sma X\angordn \andspace\\
(\Kgg X)_V &= \int^{\malp \in \GG} (S^V)^{\sma\malp} \sma X\malp
\end{split}\]
constructed in \cref{def:ggspace_gspectra,def:gggspace_gspectra}. 
\begin{itemize}
\item For a $\Gskg$-space $X$ \cref{ggtop_obj}, the coend $(\Kg X)_V$ is taken in $\Gtopst$.  This is well defined because $X$ sends morphisms in $\Gsk$ to pointed $G$-morphisms \cref{f_upom_ggtop}. 
\item For a $\GGG$-space $X$ \cref{GGTopg_obj}, the coend $(\Kgg X)_V$ is taken in $\Topst$ instead of $\Gtopst$.  The reason is that $X$ sends morphisms in $\GG$ to pointed morphisms \cref{X_upom_GGTopg} that are \emph{not} necessarily $G$-equivariant.  The group $G$ acts diagonally on representative pairs of $(\Kgg X)_V$, as defined in \cref{Kggxv_rep_gact}.  The computation \cref{Kggxv_gactwell} shows that this $G$-action on $(\Kgg X)_V$ is well defined.
\end{itemize}
To reconcile the difference between $\Kg$ and $\Kgg$, observe that the coend $(\Kg X)_V$ can also be taken in $\Topst$, with $G$ acting diagonally on representative pairs.  The functors $\Kg$ and $\Kgg$ differ only by an equivalence and a natural isomorphism; see \cref{KgKgg}.
\end{explanation}

\subsection*{Morphism Assignment of $\Kgg$}

Next, we define the morphism assignment of $\Kgg$.  A morphism between $\GGG$-spaces is a $G$-natural transformation \pcref{def:GGTopg}.  A $G$-morphism between orthogonal $G$-spectra is a $G$-equivariant $\IU$-morphism that is compatible with the sphere actions \pcref{def:iu_morphism,def:gsp_morphism}.  \cref{def:Kgg_mor} uses \cref{def:gggspace_gspectra}.

\begin{definition}[$\Kgg$ on Morphisms]\label{def:Kgg_mor}
For a compact Lie group $G$ and a $G$-natural transformation 
\begin{equation}\label{Kgg_theta_iicell}
\begin{tikzpicture}[vcenter]
\def\t{28}
\draw[0cell]
(0,0) node (a1) {\phantom{\Gsk}}
(a1)++(1.8,0) node (a2) {\phantom{\Gsk}}
(a1)++(-.08,0) node (a1') {\GG}
(a2)++(.23,0) node (a2') {\Topgst}
;
\draw[1cell=.9]
(a1) edge[bend left=\t] node {X} (a2)
(a1) edge[bend right=\t] node[swap] {Y} (a2)
;
\draw[2cell]
node[between=a1 and a2 at .42, rotate=-90, 2label={above,\tha}] {\Rightarrow}
;
\end{tikzpicture}
\end{equation}
between $\GGG$-spaces $X$ and $Y$, the $G$-morphism between orthogonal $G$-spectra
\begin{equation}\label{Kgg_psi}
(\Kgg X, \umu) \fto{\Kgg\tha} (\Kgg Y, \umu)
\end{equation}
has, for each object $V \in \IU$, $V$-component pointed $G$-morphism \cref{eqiu_mor_component} defined by the commutative diagram
\begin{equation}\label{Kgg_psiv}
\begin{tikzpicture}[vcenter]
\def\v{-1.4}
\draw[0cell=.9]
(0,0) node (a11) {(\Kgg X)_V}
(a11)++(3.4,0) node (a12) {\txint^{\malp \sins \GG} (S^V)^{\sma\malp} \sma X\malp}
(a11)++(0,\v) node (a21) {(\Kgg Y)_V}
(a12)++(0,\v) node (a22) {\txint^{\malp \sins \GG} (S^V)^{\sma\malp} \sma Y\malp}
;
\draw[1cell=.8]
(a11) edge[equal] (a12)
(a21) edge[equal] (a22)
(a11) edge[transform canvas={xshift=1em}] node[swap] {(\Kgg\tha)_V} (a21)
(a12) edge[transform canvas={xshift=-2em}, shorten <=-.5ex] node {\txint^{\malp} 1 \sma \tha_{\malp}} (a22)
;
\end{tikzpicture}
\end{equation}
in $\Gtopst$ \cref{Gtopst_smc}.  
\end{definition}

\begin{lemma}\label{kggpsi_welldef}
The assignment $\Kgg\tha$ \cref{Kgg_psi} is a $G$-morphism between orthogonal $G$-spectra.
\end{lemma}

\begin{proof}
The pointed morphism $(\Kgg\tha)_V$ \cref{Kgg_psiv} sends a representative pair $(\upom; x)$ of $(\Kgg X)_V$ \cref{Kggxv_rep} to the representative pair
\begin{equation}\label{Kgg_psiv_rep}
(\Kgg\tha)_V (\upom; x) = \big(\upom; \tha_{\malp} x\big)
\end{equation}
of $(\Kgg Y)_V$. 
\begin{description}
\item[Well definedness]
To see that $(\Kgg\tha)_V$ is well defined, we consider the two representative pairs 
\[(\upom(\sma\upla); x) \andspace (\upom; (X\upla)x)\]
in \cref{Kggxv_relation}, which are identified in the coend $(\Kgg X)_V$.  These representative pairs are sent by $(\Kgg\tha)_V$ to
\[\big(\upom(\sma\upla); \tha_{\nbe} x\big) \andspace \big(\upom; \tha_{\malp} (X\upla)x\big).\]
By the defining relation \cref{Kggxv_relation} of the coend $(\Kgg Y)_V$ and the naturality of $\tha$ \cref{GGTopg_icell_nat}, each of the preceding two representative pairs is equal to 
\[\big(\upom; (Y\upla) \tha_{\nbe} x\big).\]
Thus, $(\Kgg\tha)_V$ is a pointed morphism.
\item[Naturality] 
The naturality of $\Kgg\tha$ with respect to morphisms in $\IU$ follows from \cref{Kggxf_rep,Kgg_psiv_rep}.  Indeed, each composite in the naturality diagram \cref{iu_mor_natural} for $\Kgg\tha$ sends a representative pair $(\upom; x)$ to $(f\upom; \tha_{\malp} x)$.  Thus, $\Kgg\tha$ is an $\IU$-morphism.
\item[Equivariance] 
The $\IU$-morphism $\Kgg\tha$ is $G$-equivariant by \cref{Kggxv_rep_gact}, \cref{Kgg_psiv_rep}, and the $G$-equivariance \cref{GGTopg_icell_geq} of $\tha_{\malp}$.
\item[Compatibility] 
$\Kgg\tha$ preserves the sphere actions on $\Kgg X$ and $\Kgg Y$ by \cref{Kggx_actrep,Kgg_psiv_rep}.  Indeed, each composite in the compatibility diagram \cref{gsp_mor_axiom} for $\Kgg\tha$ sends a representative $((\upom; x); y)$ of $(\Kgg X)_V \sma S^W$ to the representative $(\upom \oplus y; \tha_{\malp} x)$ of $(\Kgg Y)_{V \oplus W}$.
\end{description}
This proves that $\Kgg\tha$ is a $G$-morphism between orthogonal $G$-spectra.
\end{proof}

\begin{definition}\label{def:Kgg_functor}
For a compact Lie group $G$, the functor
\[\GGTopg \fto{\Kgg} \Gspec\]
is defined by
\begin{itemize}
\item the object assignment $X \mapsto (\Kgg X,\umu)$ \pcref{def:gggspace_gspectra} and
\item the morphism assignment $\tha \mapsto \Kgg\tha$ \pcref{def:Kgg_mor}.  
\end{itemize}
Its functoriality follows from \cref{Kgg_psiv_rep} together with the fact that identities and composition are defined componentwise in $\GGTopg$ and $\Gspec$ \pcref{def:gsp_morphism,def:GGTopg}.
\end{definition}

\subsection*{Comparison of $\Kg$ and $\Kgg$}

Recall the adjoint equivalence \pcref{thm:ggtop_ggtopg_iieq} 
\begin{equation}\label{Ligst_display}
\begin{tikzpicture}[vcenter]
\draw[0cell]
(0,0) node (a1) {\GGTopii}
(a1)++(2.5,0) node (a2) {\GGTopg}
;
\draw[1cell=.9]
(a1) edge[transform canvas={yshift=.5ex}] node {\Lg} (a2)
(a2) edge[transform canvas={yshift=-.4ex}] node {\igst} (a1)
;
\end{tikzpicture}
\end{equation}
between the categories in \cref{def:ggtop_smc,def:GGTopg}.
\begin{itemize}
\item The right adjoint $\igst$ is induced by the full subcategory inclusion $\ig \cn \Gsk \to \GG$ \cref{ig}. 
\item The left adjoint $\Lg$ is defined objectwise as a coend \cref{Lg_f_nbe} in $\Topst$, with $G$ acting diagonally on representatives \cref{LXnbe_gaction}.
\end{itemize}  
\cref{KgKgg} proves that the prolongation functors $\Kg$ and $\Kgg$ \pcref{def:ggtop_gsp_mor,def:Kgg_functor} correspond to each other via the equivalence $\Lg$.

\begin{lemma}\label{KgKgg}
For a compact Lie group $G$, there is a natural isomorphism
\begin{equation}\label{KggLKg}
\begin{tikzpicture}[vcenter]
\def\h{2} \def\v{-1.3} \def\t{15}
\draw[0cell]
(0,0) node (a1) {\GGTopg}
(a1)++(0,\v) node (a2) {\GGTopii}
(a1)++(\h,\v/2) node (a3) {\phantom{\Gspec}}
(a3)++(-.1*\h,0) node (a3') {\Gspec}
;
\draw[1cell=.9]
(a1) edge[bend left=\t] node {\Kgg} (a3)
(a2) edge[transform canvas={xshift=-1ex}] node {\Lg} (a1)
(a2) edge[bend right=\t] node[swap] {\Kg} (a3)
;
\draw[2cell=.9]
node[between=a1 and a2 at .6, shift={(.35*\h,0)}, rotate=-50, 2label={above,\kiso}] {\Rightarrow}
;
\end{tikzpicture}
\end{equation}
between the functors $\Kgg\Lg$ and $\Kg$.
\end{lemma}

\begin{proof}
Recall from \cref{def:gsp_morphism} \eqref{def:gsp_morphism_ii} that a $G$-morphism between orthogonal $G$-spectra in $\Gspec$ is a $G$-equivariant $\IU$-morphism (\cref{def:iu_morphism} \eqref{def:iu_morphism_iii}) that is compatible with the sphere actions \cref{gsp_mor_axiom}.  The natural isomorphism $\kiso$ is defined by the pointed $G$-homeomorphisms in \cref{kgkggxv} for a pointed functor $X \cn \Gsk \to \Gtopst$ \cref{ggtop_obj} and an object $V \in \IU$ \pcref{def:indexing_gspace}.  The wedge $\txwedge$ is indexed by the $G$-set $\GGpunc(\ig\angordn; \malp)$ of nonzero morphisms $\ig\angordn \to \malp$ in $\GG$.
\begin{equation}\label{kgkggxv}
\begin{aligned}
& (\Kgg \Lg X)_V &&\\
&= \txint^{\malp \in \GG} (S^V)^{\sma\malp} \sma (\Lg X)\malp && \text{by \cref{Kggxv}}\\
&= \txint^{\malp \in \GG} (S^V)^{\sma\malp} \sma \big[\txint^{\angordn \in \Gsk} \txwedge X\angordn \big] && \text{by \cref{Lg_f_nbe}}\\
&\iso \txint^{\angordn \in \Gsk} \big[\txint^{\malp \in \GG} \txwedge (S^V)^{\sma\malp}\big] \sma X\angordn && \text{by commutation}\\
&\iso \txint^{\angordn \in \Gsk} (S^V)^{\sma\angordn} \sma X\angordn && \text{by Yoneda}\\
&= (\Kg X)_V && \text{by \cref{Kgxv}}
\end{aligned}
\end{equation}
\begin{itemize}
\item The first natural $G$-homeomorphism in \cref{kgkggxv} uses
\begin{itemize}
\item the commutation of $(S^V)^{\sma\malp} \msmam -$ with the coend $\txint^{\angordn \in \Gsk}$ and the wedge $\txwedge$; 
\item the commutation of the coends $\txint^{\malp \in \GG}$ and $\txint^{\angordn \in \Gsk}$; and 
\item the commutation of $- \msmam X\angordn$ with the wedge $\txwedge$ and the coend $\txint^{\malp \in \GG}$.
\end{itemize}
\item The second natural $G$-homeomorphism in \cref{kgkggxv} uses the Enriched Yoneda Density Theorem \cite[3.7.8]{cerberusIII}.  A set-theoretic version of the Yoneda Density Theorem is \cite[Prop.\ 2.2.1]{loregian}.
\end{itemize}
The pointed $G$-homeomorphisms in \cref{kgkggxv} are natural in $V \in \IU$ and $X \in \GGTopii$.  They are compatible with the sphere actions by \cref{Kgx_action_vw,Kggx_action_vw}.
\end{proof}

\begin{explanation}[Unpacking $\kiso$]\label{expl:kgkggxv}
The pointed $G$-homeomorphism \cref{kgkggxv}
\[(\Kgg\Lg X)_V \fto[\iso]{\kiso_{X,V}} (\Kg X)_V\]
sends a representative
\[\begin{split}
&\big(\smam\malp \fto{\upom} S^V ; \ig\angordn \fto{\upla} \malp ; x \big)\\
&\in (S^V)^{\sma\malp} \times \GGpunc(\ig\angordn; \malp) \times X\angordn
\end{split}\]
of $(\Kgg\Lg X)_V$ to the representative
\begin{equation}\label{kiso_xv}
\begin{split}
&\big(\smam \angordn \fto{\sma\upla} \smam\malp \fto{\upom} S^V ; x \big)\\
&\in (S^V)^{\sma\angordn} \times X\angordn
\end{split}
\end{equation}
of $(\Kg X)_V$.  The inverse pointed $G$-homeomorphism
\[(\Kg X)_V \fto[\iso]{\kiso^{-1}_{X,V}} (\Kgg\Lg X)_V\]
sends a representative
\[(\upom; x) \in (S^V)^{\sma\angordn} \times X\angordn\]
of $(\Kg X)_V$ to the representative
\begin{equation}\label{kisoinv_xv}
\begin{split}
&\big(\smam\angordn\fto{\upom} S^V ; \ig\angordn \fto{1} \ig\angordn ; x \big)\\
&\in (S^V)^{\sma\angordn} \times \GGpunc(\ig\angordn; \ig\angordn) \times X\angordn
\end{split}
\end{equation}
of $(\Kgg\Lg X)_V$.
\end{explanation}

\section{Equivariant $K$-Theory via $\GGG$-Categories}
\label{sec:kgo_hthy}

This section defines an equivalent variant of our equivariant $K$-theory functor $\Kgo$ \cref{Kgo_functors}, denoted by $\Khgo$, using $H$-theory $\Hgo$ \pcref{Hgo_twofunctor} instead of $J$-theory $\Jgo$ \pcref{thm:Jgo_twofunctor}.  There is also a strong variant, denoted by $\Khgosg$, involving strong $H$-theory $\Hgosg$.

\secoutline
\begin{itemize}
\item \cref{def:Khgo} defines the equivariant $K$-theory functors $\Khgo$ and $\Khgosg$, with further elaboration given in \cref{expl:Khgo}.
\item \cref{KgoKhgo} proves that $\Kgo$ and $\Khgo$ are naturally isomorphic.  The strong variant involving $\Kgosg$ and $\Khgosg$ is also true.
\item \cref{expl:Kiso} unpacks these natural isomorphisms.
\end{itemize}

\begin{definition}[$\Khgo$ and $\Khgosg$]\label{def:Khgo}
For a compact Lie group $G$ and a $\Tinf$-operad $\Op$ \pcref{as:OpA}, the functors\index{K-theory@$K$-theory!via H-theory@via $H$-theory}
\[\begin{split}
\Khgo &= \Kgg \circ \clast \circ \Hgo \andspace\\ 
\Khgosg &= \Kgg \circ \clast \circ \Hgosg
\end{split}\] 
are defined as the following composites.
\begin{equation}\label{Khgo_functors}
\begin{tikzpicture}[baseline={(a.base)}]
\def\h{2.5} \def\u{1.1} \def\v{1}
\draw[0cell]
(0,0) node (a) {\phantom{\AlglaxO}}
(a)++(0,\v) node (a1) {\AlglaxO}
(a)++(0,\u) node (a1'') {\phantom{\AlglaxO}}
(a)++(0,\v) node (a1') {\phantom{(\Op)}} 
(a)++(0,-\v) node (a2) {\AlgpspsO}
(a)++(0,-\v) node (a2') {\phantom{(\Op)}} 
(a)++(.8*\h,0) node (b) {\GGCatg}
(b)++(\h,0) node (c) {\GGTopg}
(c)++(.9*\h,0) node (d) {\Gspec}
;
\draw[1cell=.9]
(a1') edge node[inner sep=1pt,swap] {\Hgo} (b)
(a2') edge[shorten <=.7ex] node[inner sep=1pt] {\Hgosg} (b)
(b) edge node {\clast} (c)
(c) edge node {\Kgg} (d)
;
\draw[1cell=1]
(a1'') [rounded corners=2pt] -- node[pos=.7] {\Khgo} ($(c)+(0,\u)$) -| (d)
;
\draw[1cell=1]
(a2) [rounded corners=2pt] -- node[pos=.7] {\Khgosg} ($(c)+(0,-\v)$) -| (d)
;
\end{tikzpicture}
\end{equation}
The categories in the diagram \cref{Khgo_functors} are defined in \cref{oalgps_twocat,def:ggcatg,def:GGTopg,def:gsp_module,def:gsp_morphism}.
\begin{description}
\item[$H$-theory] For a group $G$ and a $\Tinf$-operad $\Op$, $\Hgo$ and $\Hgosg$ are the (strong) $H$-theory 2-functors in \cref{Hgo_twofunctor}.  At the object level, they send $\Op$-pseudoalgebras to $\GGG$-categories \pcref{A_ptfunctorGG}.  At the 1-cell level, they send lax $\Op$-morphisms and $\Op$-pseudomorphisms to $G$-natural transformations \pcref{hgo_icell_welldef}.
\item[Classifying space] For a group $G$, the functor $\clast$ \cref{clast_GGCatg} sends $\GGG$-categories and $G$-natural transformations to $\GGG$-spaces and $G$-natural transformations by composing and whiskering with the classifying space functor \pcref{expl:GGCatg_GGTopg}.
\item[Prolongation] For a compact Lie group $G$, the functor $\Kgg$ \pcref{def:gggspace_gspectra,def:Kgg_mor} sends $\GGG$-spaces and $G$-natural transformations to orthogonal $G$-spectra and $G$-morphisms.\defmark
\end{description}
\end{definition}

\begin{explanation}[Unpacking]\label{expl:Khgo}
For an $\Op$-pseudoalgebra $\A$ \pcref{def:pseudoalgebra} and an object $V \in \IU$ \pcref{def:indexing_gspace}, $\Khgo$ and $\Khgosg$ yield the following pointed $G$-spaces.
\begin{equation}\label{khgoav}
\begin{split}
(\Khgo\A)_V &= \int^{\malp \in \GG} (S^V)^{\sma\malp} \sma \cla(\Amalp)\\
(\Khgosg\A)_V &= \int^{\malp \in \GG} (S^V)^{\sma\malp} \sma \cla(\Asgmalp)\\
\end{split}
\end{equation}
The sphere action is defined in \cref{Kggx_action_vw}.  For each object $\malp \in \GG$ \cref{GG_objects}, $\Amalp$ is the pointed $G$-category of $\malp$-systems in $\A$ \pcref{def:nbeta_gcat}, and $\Asgmalp$ is the pointed full $G$-subcategory of strong $\malp$-systems \cref{Asgangordnbe}.
\end{explanation}

\cref{KgoKhgo} proves that the equivariant $K$-theory functors $\Kgo$ \cref{Kgo_functors} and $\Khgo$ \cref{Khgo_functors} are naturally isomorphic, and so are the strong variants $\Kgosg$ and $\Khgosg$.

\begin{theorem}\label{KgoKhgo}
For a compact Lie group $G$ and a $\Tinf$-operad $\Op$ \pcref{as:OpA}, there are natural isomorphisms
\begin{equation}\label{KisoKisosg}

\end{equation}
between the functors in \cref{Kgo_functors,Khgo_functors}.
\end{theorem}

\begin{proof}
To construct the natural isomorphism $\Kiso$, we consider the following diagram of functors.
\begin{equation}\label{kgo_khgo_diag}
%
\end{equation}
\begin{itemize}
\item The top and bottom regions are the definitions of $\Khgo$ \cref{Khgo_functors} and $\Kgo$ \cref{Kgo_functors}.
\item The left region involving $\Hgo$, $\Jgo$, and $\igst$ commutes by \cref{jgohgoigst}.
\item The middle region involving $\clast$ and $\igst$ commutes because each $\igst$ precomposes with the full subcategory inclusion $\ig \cn \Gsk \to \GG$, while each $\clast$ postcomposes with the classifying space functor $\cla$.
\item The vertical pair $(\Lg,\igst)$ is an adjoint equivalence \pcref{thm:ggtop_ggtopg_iieq}.  Its counit natural isomorphism $\vg \cn \Lg\igst \fiso 1$ is defined in \crefrange{vg}{vgxnbe}, with $(\Gcatst,\Catgst)$ replaced by $(\Gtopst,\Topgst)$.
\item The natural isomorphism $\kiso \cn \Kgg\Lg \fiso \Kg$ is given by \cref{KgKgg}.
\end{itemize}
Using the diagram \cref{kgo_khgo_diag}, the natural isomorphism $\Kiso$ is given by the following composite.
\begin{equation}\label{Kiso_diag}
\begin{tikzpicture}[vcenter]
\def\h{4.2} \def\v{1.3}
\draw[0cell]
(0,0) node (a1) {\Kgg \clast \Hgo}
(a1)++(1.6,0) node (a0) {\Khgo}
(a1)++(0,-\v) node (a2) {\Kgg \Lg \igst \clast \Hgo}
(a2)++(\h,0) node (a3) {\Kg \igst \clast \Hgo}
(a3)++(.57*\h,0) node (a4) {\Kg \clast \igst \Hgo}
(a4)++(0,\v) node (a5) {\Kg \clast \Jgo}
(a5)++(-1.45,0) node (a6) {\Kgo}
;
\draw[1cell=.9]
(a0) edge node {\Kiso} node[swap] {\iso} (a6)
(a0) edge[equal] (a1)
(a1) edge[transform canvas={xshift=1em}] node {\iso} node[swap] {\Kgg \vginv_{\clast \Hgo}} (a2)
(a2) edge node {\kiso_{\igst \clast \Hgo}} node[swap] {\iso} (a3)
(a3) edge[equal] (a4)
(a4) edge[equal] (a5)
(a5) edge[equal] (a6)
;
\end{tikzpicture}
\end{equation}
The natural isomorphism $\Kisosg \cn \Khgosg \fiso \Kgosg$ is constructed in the same way by replacing $J$-theory $\Jgo$ and $H$-theory $\Hgo$ in the diagram \cref{kgo_khgo_diag} with their strong variants $\Jgosg$ \pcref{thm:Jgo_twofunctor} and $\Hgosg$ \pcref{Hgo_twofunctor}.
\end{proof}

\begin{explanation}[Unpacking $\Kiso$ and $\Kisosg$]\label{expl:Kiso}
The natural isomorphism $\Kiso \cn \Khgo \fiso \Kgo$ \cref{Kiso_diag} is given explicitly as follows for an $\Op$-pseudoalgebra $\A$ \pcref{def:pseudoalgebra} and an object $V \in \IU$ \pcref{def:indexing_gspace}.  Using \cref{kgoav,khgoav,vgfninv_x,kiso_xv}, the pointed $G$-homeomorphism 
\begin{equation}\label{KisoAV}
\begin{tikzpicture}[vcenter]
\def\h{3.5}
\draw[0cell]
(0,0) node (a1) {(\Khgo\A)_V}
(a1)++(\h,0) node (a2) {\int^{\malp \in \GG} (S^V)^{\sma\malp} \sma \cla\Amalp}
(a1)++(0,-1.3) node (b1) {(\Kgo\A)_V}
(b1)++(.9*\h,0) node (b2) {\int^{\angordn \in \Gsk} (S^V)^{\sma\angordn} \sma \cla\Aangordn}
;
\draw[1cell=.9]
(a1) edge[equal] (a2)
(b1) edge[equal] (b2)
(a1) edge node {\iso} node[swap] {\Kiso_{\A,V}} (b1)
;
\end{tikzpicture}
\end{equation}
sends a representative \cref{Kggxv_rep}
\[(\upom; x) \in (S^V)^{\sma\malp} \times \cla\Amalp\]
of $(\Khgo\A)_V$ to the representative
\[\big(\upom(\sma\jin) ; (\cla\A\jininv)(x) \big) \in (S^V)^{\sma\angordm} \times \cla\Aangordm\]
of $(\Kgo\A)_V$.  
\begin{description}
\item[First component]
For each object $\malp \in \GG$ \cref{GG_objects}, the isomorphism \cref{jin_nbe}
\[\ig\angordm \fto[\iso]{\jin} \malp \inspace \GG\]
is defined by identity functions, and so is its inverse \cref{jininv}
\[\malp \fto[\iso]{\jininv} \ig\angordm.\]
By \cref{smash_fpsiFG}, the isomorphism
\[\sma\angordm \fto[\iso]{\sma\jin} \msmam\malp \inspace \FG\]
is given by the identity function.  The composite
\[\sma\angordm \fto{\sma\jin} \msmam\malp \fto{\upom} S^V\]
is given by the same pointed function as $\upom$.
\item[Second component]
By \cref{def:Aangordn_gcat,def:nbeta_gcat}, the pointed $G$-category $\Aangordm$ of $\angordm$-systems and the pointed $G$-category $\Amalp$ of $\malp$-systems are equal as pointed categories, but \emph{not} generally as $G$-categories.  By \cref{def:ftil_functor,def:psitil_functor,def:Afangpsi,def:AfangpsiGG}, the pointed isomorphism
\[\Amalp \fto[\iso]{\A\jininv} \Aangordm\]
is given by the identity functor.  Passing to classifying spaces, the pointed $G$-spaces $\cla\Amalp$ and $\cla\Aangordm$ are equal as pointed spaces, but not generally as $G$-spaces.  The pointed isomorphism
\[\cla\Amalp \fto[\iso]{\cla\A\jininv} \cla\Aangordm\]
is given by the identity function.
\item[Inverse]
By \cref{vgxnbe,kisoinv_xv}, the inverse pointed $G$-homeomorphism 
\[(\Kgo\A)_V \fto[\iso]{\Kiso^{-1}_{\A,V}} (\Khgo\A)_V\]
sends a representative
\[(\upom; x) \in (S^V)^{\sma\angordn} \times \cla\Aangordn\]
of $(\Kgo\A)_V$ to the representative $(\upom; x)$ of $(\Khgo\A)_V$.
\item[Strong variant]
The natural isomorphism 
\[(\Khgosg\A)_V \fto[\iso]{\Kisosg_{\A,V}} (\Kgosg\A)_V\]
admits the same description as $\Kiso_{\A,V}$ with $\Amalp$ and $\Aangordm$ replaced by their strong variants $\Asgmalp$ \cref{Asgangordnbe} and $\Aangordmsg$ \cref{sgAordnbe}.\defmark
\end{description}
\end{explanation}

%% file: chap/shih.tex
For a group $G$ and a 1-connected $\Gcat$-operad $\Op$, this chapter constructs Shimakawa (strong) $H$-theory 2-functors
\[\AlglaxO \fto{\Sgo} \FGCatg \andspace \AlgpspsO \fto{\Sgosg} \FGCatg\]
between the 2-categories in \cref{oalgps_twocat,def:fgcatg}, along with Shimakawa (strong) $J$-theory 2-functors
\[\AlglaxO \fto{\Jgos} \FGCat \andspace \AlgpspsO \fto{\Jgossg} \FGCat.\]
The 2-functors $\Sgo$ and $\Sgosg$ are based on the indexing $G$-category $\FG$ of pointed finite $G$-sets and pointed functions, with the conjugation $G$-action \pcref{def:FG}.  Using $\Sgo$ and $\Sgosg$, \cref{ch:shimakawa_K} discusses Shimakawa equivariant $K$-theory machines that send $\Op$-pseudoalgebras to orthogonal $G$-spectra. \cref{ch:shim_top,part:kgo_shi_comp} compare Shimakawa $K$-theory with our equivariant $K$-theory \cref{Khgo_functors}. 

\organization
This chapter consists of the following sections.

\secname{sec:fgcatg}
This section defines the 2-categories $\FGCat$ and $\FGCatg$.  There is an adjoint 2-equivalence 
\begin{equation}\label{LgIgst_chi}
\begin{tikzpicture}[vcenter]
\draw[0cell]
(0,0) node (a1) {\FGCat}
(a1)++(2.5,0) node (a2) {\phantom{\FGCatg}}
(a2)++(0,-.04) node (a2') {\FGCatg}
;
\draw[1cell=.9]
(a1) edge[transform canvas={yshift=.5ex}] node {\Lg} (a2)
(a2) edge[transform canvas={yshift=-.4ex}] node {\igst} (a1)
;
\end{tikzpicture}
\end{equation}
in which $\igst$ is induced by the full subcategory inclusion $\ig \cn \Fsk \to \FG$.

\secname{sec:sgo_obj}
This section constructs the object assignments of $\Sgo$ and $\Sgosg$, which send an $\Op$-pseudoalgebra $\A$ to pointed $G$-functors
\begin{equation}\label{SgoSgosg_chi}
\begin{tikzpicture}[vcenter]
\draw[0cell]
(0,0) node (a1) {\FG}
(a1)++(3.5,0) node (a2) {\Catgst.}
;
\draw[1cell=.9]
(a1) edge[transform canvas={yshift=.5ex}] node {\Sgo\A \,=\, \Adash} (a2)
(a1) edge[transform canvas={yshift=-.4ex}] node[swap] {\Sgosg\A \,=\, \Asgdash} (a2)
;
\end{tikzpicture}
\end{equation}

\secname{sec:sgo_iifunctor}
This section finishes the construction of the 2-functors $\Sgo$ and $\Sgosg$ by defining their 1-cell and 2-cell assignments.  The 1-cell assignments of $\Sgo$ and $\Sgosg$ send lax $\Op$-morphisms and $\Op$-pseudomorphisms to $G$-natural transformations.  Their 2-cell assignments send $\Op$-transformations to $G$-modifications.  For the Barratt-Eccles operad $\BE$, $\Sgosgbe$ recovers Shimakawa's original construction in \cite{shimakawa89}; see \cref{expl:shi89_section2}.  The 2-functor $\Jgos$ is the composite of $\Sgo$ and $\igst$.  The 2-functor $\Jgossg$ is the composite of $\Sgosg$ and $\igst$.

\section{$\Fskg$-Categories and $\FGG$-Categories}
\label{sec:fgcatg}

This section defines the codomain 2-category $\FGCatg$ of Shimakawa $H$-theory and the 2-equivalence between the 2-categories $\FGCat$ and $\FGCatg$. 

\secoutline
\begin{itemize}
\item \cref{def:fgcat} defines the 2-category $\FGCat$ of $\Fskg$-categories, with further description given in \cref{expl:fgcat_iicat}.
\item \cref{def:fgcatg} defines the 2-category $\FGCatg$ of $\FGG$-categories, with further description given in \cref{expl:fgcatg_iicat}.
\item \cref{thm:fgcat_fgcatg_iieq} records the 2-equivalence between these 2-categories induced by the inclusion $\Fsk \to \FG$.  This observation lifts a result of Shimakawa \cite{shimakawa91} about $\GaG$-spaces and $\Ga\mh G$-spaces.  See \cref{rk:shi91}.
\item \cref{Xjin} shows how, for an $\FGG$-category $X$ and an object $\nbeta \in \FG$, the pointed $G$-category $X\nbeta$ can be described in terms of $X\ordn$ with a twisted $G$-action.
\end{itemize}

\subsection*{$\Fskg$-Categories}

For \cref{def:fgcat}, recall from \cref{def:Fsk} the small pointed category $\Fsk$ with pointed finite sets $\ordn = \{0,1,\ldots,n\}$ \cref{ordn} as objects, pointed functions as morphisms, and basepoint $\ordz$.  Also recall the 2-category $\Gcatst$ \cref{Gcatst}.  The 2-category $\FGCat$ in \cref{def:fgcat} is the $\Fsk$-analogue of the 2-category $\GGCatii$ \pcref{def:GGCat}.

\begin{definition}\label{def:fgcat}
For a group $G$, the 2-category $\FGCat$ is defined as follows.
\begin{description}
\item[Objects]
An object in $\FGCat$, called an \index{F-G-category@$\Fskg$-category}\emph{$\Fskg$-category}, is a pointed functor
\begin{equation}\label{fgcat_obj}
(\Fsk, \ordz) \fto{X} (\Gcatst,\boldone).
\end{equation}
\item[1-cells]
A 1-cell $\theta \cn X \to X'$ in $\FGCat$ is a natural transformation as follows.
\begin{equation}\label{fgcat_mor}
\begin{tikzpicture}[vcenter]
\def\t{28}
\draw[0cell]
(0,0) node (a1) {\Fsk}
(a1)++(1.8,0) node (a2) {\phantom{\Fsk}}
(a2)++(.3,0) node (a2') {\Gcatst}
;
\draw[1cell=.9]
(a1) edge[bend left=\t] node {X} (a2)
(a1) edge[bend right=\t] node[swap] {X'} (a2)
;
\draw[2cell]
node[between=a1 and a2 at .45, rotate=-90, 2label={above,\theta}] {\Rightarrow}
;
\end{tikzpicture}
\end{equation}
\item[2-cells]
A 2-cell $\Theta \cn \theta \to \ups$ in $\FGCat$ is a modification as follows.
\begin{equation}\label{fgcat_iicell}
\begin{tikzpicture}[vcenter]
\def\t{25}
\draw[0cell]
(0,0) node (a) {\Fsk}
(a)++(3,0) node (b) {\phantom{\Fsk}}
(b)++(.3,0) node (b') {\Gcatst}
;
\draw[1cell=.8]
(a) edge[bend left=\t] node {X} (b)
(a) edge[bend right=\t] node[swap] {X'} (b)
;
\draw[2cell=.9]
node[between=a and b at .34, rotate=-90, 2label={below,\theta}] {\Rightarrow}
node[between=a and b at .65, rotate=-90, 2label={above,\ups}] {\Rightarrow}
;
\draw[2cell]
node[between=a and b at .52, rotate=0, shift={(0,-.16)}, 2labelalt={above,\!\!\Theta}] {\Rrightarrow}
;
\end{tikzpicture}
\end{equation}
\item[Other structures]
Identity 1-cells and 2-cells, vertical composition of 2-cells, and horizontal composition of 1-cells and 2-cells are defined componentwise in the 2-category $\Gcatst$.
\end{description}
The underlying 1-category of $\FGCat$ is denoted by the same notation.
\end{definition}

\begin{explanation}[Unpacking $\FGCat$]\label{expl:fgcat_iicat}
The 2-category $\FGCat$ in \cref{def:fgcat} is given explicitly as follows.
\begin{description}
\item[Objects]
An $\Fskg$-category $X \cn \Fsk \to \Gcatst$ \cref{fgcat_obj} consists of the following data.
\begin{itemize}
\item $X$ sends each pointed finite set $\ordm \in \Fsk$ to a small pointed $G$-category $X\ordm$ such that $X\ordz = \boldone$.  The \emph{canonical basepoint}\dindex{canonical}{basepoint} of $X\ordm$ is given by the $G$-functor
\begin{equation}\label{fgcat_obj_obj}
X(\ordz \to \ordm) \cn X\ordz = \boldone \to X\ordm,
\end{equation}
where $\ordz \to \ordm$ is the unique pointed function. 
\item $X$ sends each pointed function $\psi \cn \ordm \to \ordn$ in $\Fsk$ to a pointed $G$-functor
\begin{equation}\label{fgcat_obj_mor}
X\ordm \fto{X\psi} X\ordn
\end{equation}
such that $X$ preserves identities and composition of morphisms.  
\end{itemize}
\item[1-cells]
A 1-cell $\theta \cn X \to X'$ in $\FGCat$ \cref{fgcat_mor} consists of, for each pointed finite set $\ordm \in \Fsk$, an $\ordm$-component pointed $G$-functor
\begin{equation}\label{fgcat_mor_component}
X\ordm \fto{\theta_{\ordm}} X'\ordm
\end{equation}
such that, for each pointed function $\psi \cn \ordm \to \ordn$ in $\Fsk$, the following naturality diagram of pointed $G$-functors commutes.
\begin{equation}\label{fgcat_mor_naturality}

\end{equation}
Identity 1-cells and horizontal composition of 1-cells in $\FGCat$ are defined componentwise using the components in \cref{fgcat_mor_component}.  A 1-cell is automatically pointed in the sense that $\theta_{\ordz} = 1_\bone$.
\item[2-cells]
A 2-cell $\Theta \cn \theta \to \ups$ in $\FGCat$ \cref{fgcat_iicell} consists of, for each pointed finite set $\ordm \in \Fsk$, an $\ordm$-component pointed $G$-natural transformation
\begin{equation}\label{fgcat_inthom_Theta}
%
\end{equation}
such that, for each pointed function $\psi \cn \ordm \to \ordn$ in $\Fsk$, the following two whiskered $G$-natural transformations are equal.
\begin{equation}\label{fgcat_inthom_Theta_modax}
%
\end{equation}
Identities, horizontal composition, and vertical composition of 2-cells are given componentwise using the components in \cref{fgcat_inthom_Theta}.   A 2-cell is automatically pointed in the sense that $\Theta_{\ordz} = 1_{1_{\bone}}$.\defmark
\end{description}
\end{explanation}

\subsection*{$\FGG$-Categories}
For \cref{def:fgcatg}, recall from \cref{def:FG} the small pointed $G$-category $\FG$ with pointed finite $G$-sets $\ordn^\be$ \cref{ordn_be} as objects, pointed functions as morphisms, and the conjugation $G$-action on morphisms.  Also recall the 2-category $\Catgst$ \cref{Catgst_iicat}.  The 2-category $\FGCatg$ in \cref{def:fgcatg} is the $\FG$-analogue of the 2-category $\GGCatg$ \pcref{def:ggcatg}.

\begin{definition}\label{def:fgcatg}
For a group $G$, the 2-category $\FGCatg$ is defined as follows.
\begin{description}
\item[Objects] An object in $\FGCatg$, called an \index{FGG-category@$\FGG$-category}\emph{$\FGG$-category}, is a pointed $G$-functor
\begin{equation}\label{fgcatg_obj}
(\FG,\ordz) \fto{X} (\Catgst,\boldone).
\end{equation}
\item[1-cells] A 1-cell $\theta \cn X \to X'$ in $\FGCatg$ is a $G$-natural transformation as follows.  
\begin{equation}\label{fgcatg_icell}
\begin{tikzpicture}[vcenter]
\draw[0cell]
(0,0) node (a) {\FG}
(a)++(2,0) node (b) {\phantom{\FG}}
(b)++(.15,0) node (b') {\Catgst}
;
\draw[1cell=.9]
(a) edge[bend left=25] node {X} (b)
(a) edge[bend right=25] node[swap] {X'} (b)
;
\draw[2cell]
node[between=a and b at .45, rotate=-90, 2label={above,\theta}] {\Rightarrow}
;
\end{tikzpicture}
\end{equation}
\item[2-cells]  A 2-cell $\Theta \cn \theta \to \ups$ in $\FGCatg$ is a $G$-modification as follows.  
\begin{equation}\label{fgcatg_iicell}
\begin{tikzpicture}[vcenter]
\def\t{25}
\draw[0cell]
(0,0) node (a) {\FG}
(a)++(3,0) node (b) {\phantom{\FG}}
(b)++(.15,0) node (b') {\Catgst}
;
\draw[1cell=.8]
(a) edge[bend left=\t] node {X} (b)
(a) edge[bend right=\t] node[swap] {X'} (b)
;
\draw[2cell=.9]
node[between=a and b at .34, rotate=-90, 2label={below,\theta}] {\Rightarrow}
node[between=a and b at .65, rotate=-90, 2label={above,\ups}] {\Rightarrow}
;
\draw[2cell]
node[between=a and b at .52, rotate=0, shift={(0,-.16)}, 2labelalt={above,\!\!\Theta}] {\Rrightarrow}
;
\end{tikzpicture}
\end{equation}
\item[Other structures] Identity 1-cells and 2-cells, vertical composition of 2-cells, and horizontal composition of 1-cells and 2-cells are defined componentwise in the 2-category $\Catgst$.
\end{description}
The underlying 1-category of $\FGCatg$ is denoted by the same notation.
\end{definition}

\begin{explanation}[Unpacking $\FGCatg$]\label{expl:fgcatg_iicat}
The 2-category $\FGCatg$ in \cref{def:fgcatg} is given explicitly as follows.
\begin{description}
\item[Objects]
An $\FGG$-category $X \cn \FG \to \Catgst$ \cref{fgcatg_obj} consists of the following data.
\begin{itemize}
\item $X$ sends each pointed finite $G$-set $\ordm^\al \in \FG$ \cref{ordn_be} to a small pointed $G$-category $X\ordm^\al$ such that $X\ordz = \bone$.  Its $G$-fixed basepoint is given by the pointed functor
\[X(\ordz \to \ordm^\al) \cn X\ordz = \bone \to X\ordm^\al\]
for the unique pointed function $\ordz \to \ordm^\al$ in $\FG$.  
\item $X$ sends each pointed function $\psi \cn \ordm^\al \to \ordn^\be$ in $\FG$ to a pointed functor
\begin{equation}\label{X_psi_FG}
X\ordm^\al \fto{X\psi} X\ordn^\be
\end{equation}
such that $X$ preserves identity morphisms and composition.  The functor $X\psi$ is \emph{not} generally $G$-equivariant.
\item The $G$-equivariance of $X$ means the equality of functors
\begin{equation}\label{FGGcat_Gequiv}
X(g \cdot \psi) = g (X\psi) \ginv
\end{equation}
for each $g \in G$ and each morphism $\psi$ in $\FG$.  The morphism $g \cdot \psi$ is the conjugation $g$-action on $\psi$ \cref{gpsi}, and $g (X\psi) \ginv$ is the conjugation $g$-action \cref{conjugation-gaction} on $X\psi$.  In particular, if the morphism $\psi$ is $G$-equivariant, then $X\psi$ is a pointed $G$-functor.
\end{itemize}
\item[1-cells]
A 1-cell $\tha \cn X \to X'$ in $\FGCatg$ \cref{fgcatg_icell} consists of, for each object $\ordm^\al \in \FG$, an $\ordm^\al$-component pointed functor
\begin{equation}\label{fgcatg_icell_comp}
X\ordm^\al \fto{\tha_{\ordm^\al}} X'\ordm^\al
\end{equation}
such that, for each morphism $\psi \cn \ordm^\al \to \ordn^\be$ in $\FG$, the following naturality diagram of pointed functors commutes.
\begin{equation}\label{fgcatg_icell_nat}
\begin{tikzpicture}[vcenter]
\def\v{-1.4}
\draw[0cell]
(0,0) node (a11) {X\ordm^\al}
(a11)++(2.2,0) node (a12) {X'\ordm^\al}
(a11)++(0,\v) node (a21) {X\ordn^\be}
(a12)++(0,\v) node (a22) {X'\ordn^\be}
;
\draw[1cell=.9]
(a11) edge node {\theta_{\ordm^\al}} (a12)
(a12) edge node {X'\psi} (a22)
(a11) edge node[swap] {X\psi} (a21)
(a21) edge node {\theta_{\ordn^\be}} (a22)
;
\end{tikzpicture}
\end{equation}
Since $G$ acts trivially on the objects of $\FG$, the $G$-equivariance of $\tha$ means the equality of functors
\begin{equation}\label{fgcatg_icell_geq}
\tha_{\ordm^\al} = g \tha_{\ordm^\al} \ginv
\end{equation}
for each $g \in G$ and object $\ordm^\al \in \FG$.   In other words, each component of a 1-cell in $\FGCatg$ is a pointed $G$-functor.  

Identity 1-cells and horizontal composition of 1-cells in $\FGCatg$ are defined componentwise using the components in \cref{fgcatg_icell_comp}.  A 1-cell is automatically pointed, meaning that $\tha_{\ordz} = 1_{\bone}$.
\item[2-cells]
A 2-cell $\Theta \cn \theta \to \ups$ in $\FGCatg$ \cref{fgcatg_iicell} consists of, for each object $\ordm^\al \in \FG$, an $\ordm^\al$-component pointed natural transformation
\begin{equation}\label{fgcatg_iicell_comp}
\begin{tikzpicture}[vcenter]
\def\t{25}
\draw[0cell]
(0,0) node (a1) {\phantom{X'}}
(a1)++(2.2,0) node (a2) {\phantom{X'}}
(a1)++(-.15,0) node (a1') {X\ordm^\al}
(a2)++(.23,0) node (a2') {X'\ordm^\al}
;
\draw[1cell=.85]
(a1) edge[bend left=\t] node {\theta_{\ordm^\al}} (a2)
(a1) edge[bend right=\t] node[swap] {\ups_{\ordm^\al}} (a2)
;
\draw[2cell]
node[between=a1 and a2 at .37, rotate=-90, 2label={above, \Theta_{\ordm^\al}}] {\Rightarrow}
;
\end{tikzpicture}
\end{equation}
such that, for each morphism $\psi \cn \ordm^\al \to \ordn^\be$ in $\FG$, the following two whiskered natural transformations are equal.
\begin{equation}\label{fgcatg_iicell_modax}
\begin{tikzpicture}[vcenter]
\def\t{25} \def\v{-1.6}
\draw[0cell]
(0,0) node (a1) {\phantom{X'}}
(a1)++(2.5,0) node (a2) {\phantom{X'}}
(a1)++(-.15,0) node (a1') {X\ordm^\al}
(a2)++(.2,0) node (a2') {X'\ordm^\al}
(a1)++(0,\v) node (b1) {\phantom{X'}}
(a2)++(0,\v) node (b2) {\phantom{X'}}
(b1)++(-.15,0) node (b1') {X\ordn^\be}
(b2)++(.2,0) node (b2') {X'\ordn^\be}
;
\draw[1cell=.8]
(a1) edge[bend left=\t] node[pos=.4] {\theta_{\ordm^\al}} (a2)
(a1) edge[bend right=\t] node[swap,pos=.6] {\ups_{\ordm^\al}} (a2)
(b1) edge[bend left=\t] node[pos=.4] {\theta_{\ordn^\be}} (b2)
(b1) edge[bend right=\t] node[swap,pos=.6] {\ups_{\ordn^\be}} (b2)
(a1) edge[transform canvas={xshift=-.7ex}] node[swap] {X\psi} (b1)
(a2) edge node {X'\psi} (b2)
;
\draw[2cell=.9]
node[between=a1 and a2 at .37, rotate=-90, 2label={above, \Theta_{\ordm^\al}}] {\Rightarrow}
node[between=b1 and b2 at .37, rotate=-90, 2label={above, \Theta_{\ordn^\be}}] {\Rightarrow}
;
\end{tikzpicture}
\end{equation}
Since $G$ acts trivially on the objects of $\FG$, the $G$-equivariance of $\Theta$ means the equality of natural transformations
\begin{equation}\label{fgcatg_iicell_geq}
\Theta_{\ordm^\al} = g * \Theta_{\ordm^\al} * \ginv
\end{equation}
for each $g \in G$ and object $\ordm^\al \in \FG$, where $*$ denotes horizontal composition of natural transformations.   In other words, each component of a 2-cell in $\FGCatg$ is a pointed $G$-natural transformation.  

Identities, horizontal composition, and vertical composition of 2-cells are given componentwise using the components in \cref{fgcatg_iicell_comp}.   A 2-cell is automatically pointed, meaning that $\Theta_{\ordz} = 1_{1_{\bone}}$.\defmark
\end{description}
\end{explanation}

\begin{remark}[$\GaG$-Categories]\label{rk:GaGcat}
In Shimakawa's work \cite{shimakawa89,shimakawa91}, the pointed $G$-category $\FG$ is denoted by $\GaG$.  In \cite[p.\ 251]{shimakawa89}, a \emph{$\GaG$-category} is defined as a $G$-functor $\FG \to \Catgst$, \emph{without} requiring it to be a pointed functor.  The pointed condition is included in \cite[Def.\ 2.1 (a)]{shimakawa89} for a \emph{special $\GaG$-category}.  Thus, an $\FGG$-category \cref{fgcatg_obj} is the same as a $\GaG$-category that is pointed, sending $\ordz$ to $\bone$.  We include the pointed condition in $\FGG$-categories because $G$-functors $\FG \to \Catgst$ of practical interest, including those produced by Shimakawa $H$-theory, are pointed.
\end{remark}

\subsection*{2-Equivalence Between $\Fskg$-Categories and $\FGG$-Categories}

Recall the 2-equivalence \pcref{thm:ggcat_ggcatg_iieq}
\[\GGCatg \fto[\sim]{\igst} \GGCatii\]
induced by the pointed full subcategory inclusion $\ig \cn \Gsk \to \GG$ \cref{ig}.  \cref{thm:fgcat_fgcatg_iieq} is the analogue involving $\FGCatg$, $\FGCat$, and the inclusion $\Fsk \to \FG$ in \cref{ig_FG}.

\begin{definition}[From $\Fsk$ to $\FG$]\label{def:fsk_fg}
For a group $G$, define the pointed full subcategory inclusion
\begin{equation}\label{ig_FG}
\Fsk \fto{\ig} \FG
\end{equation}
by sending each pointed finite set $\ordn \in \Fsk$ \cref{ordn} to the pointed finite $G$-set $\ordn^\eps$ \cref{ordn_be} equipped with the trivial $G$-action $\eps \cn G \to \Si_n$.  The functor $\ig$ is well defined by the fact that morphisms in $\FG$ \pcref{def:FG} are pointed functions between underlying pointed finite sets.  If there is no danger of confusion, then $\ig\ordn$ is abbreviated to $\ordn$ and similarly for morphisms.
\end{definition}

\begin{lemma}\label{thm:fgcat_fgcatg_iieq}
For each group $G$, there is an adjoint 2-equivalence\index{adjoint 2-equivalence}
\begin{equation}\label{Lgigst_gcat}
\begin{tikzpicture}[vcenter]
\draw[0cell]
(0,0) node (a1) {\FGCat}
(a1)++(2.5,0) node (a2) {\phantom{\FGCatg}}
(a2)++(0,-.04) node (a2') {\FGCatg}
;
\draw[1cell=.9]
(a1) edge[transform canvas={yshift=.5ex}] node {\Lg} (a2)
(a2) edge[transform canvas={yshift=-.4ex}] node {\igst} (a1)
;
\end{tikzpicture}
\end{equation}
between the 2-categories in \cref{def:fgcat,def:fgcatg}, where the right 2-adjoint $\igst$ is given by precomposition with $\ig \cn \Fsk \to \FG$ in \cref{ig_FG}.
\end{lemma}

\begin{proof}
The proofs of \cref{igst_iifunctor,thm:ggcat_ggcatg_iieq} are applicable by replacing the indexing categories $(\Gsk,\vstar)$ and $(\GG,\vstar)$ with, respectively, $(\Fsk,\ordz)$ and $(\FG,\ordz)$.  Thus, for a pointed functor $X \cn \Fsk \to \Gcatst$, the pointed $G$-functor
\[(\FG,\ordz) \fto{\Lg X} (\Catgst,\bone)\]
sends a pointed finite $G$-set $\nbeta \in \FG$ to the coend 
\begin{equation}\label{LX_nbeta}
(\Lg X)\nbeta = \ecint^{\ordm \in \Fsk} \bigvee_{\FGpunc(\ig\ordm \,;\, \nbeta)} X\ordm
\end{equation}
taken in $\Catst$.  The wedge is indexed by the set $\FGpunc(\ig\ordm ; \nbeta)$ of nonzero morphisms in $\FG(\ig\ordm ; \nbeta)$.  The group $G$ acts diagonally on representatives, meaning
\begin{equation}\label{LXnbeta_gaction}
g(\psi ; x) = (g\psi ; gx)
\end{equation}
for $g \in G$, $\psi \in \FGpunc(\ig\ordm; \nbeta)$, and $x \in X\ordm$.  
\end{proof}

\begin{remark}\label{rk:shi91}
\cref{thm:fgcat_fgcatg_iieq} is a 2-categorical analogue of \cite[Theorem 1]{shimakawa91}.  That result states that there is an equivalence between the category of $\GaG$-spaces and the category of $\Ga\mh G$-spaces.  See \cref{fgtop_fgtopg_eq}.
\end{remark}

For each pointed $G$-functor $X \cn \FG \to \Catgst$ and an object $\nbeta \in \FG$, the pointed $G$-category $X\nbeta$ can be reconstructed from $X\ordn$ using the following definition.

\begin{definition}\label{def:gnbeta}
Given a pointed $G$-functor $X \cn \FG \to \Catgst$ and a pointed finite $G$-set $\nbeta \in \FG$, define the pointed $G$-category $(X\ordn)_\be$ as follows.
\begin{itemize}
\item The underlying pointed category of $(X\ordn)_\be$ is the underlying pointed category of $X\ordn$, where $\ordn \in \FG$ is equipped with the trivial $G$-action.
\item For each $g \in G$, the $g$-action functor on $(X\ordn)_\be$ is the composite
\begin{equation}\label{Xordnbe_g}
\begin{tikzpicture}[vcenter]
\def\u{.65}
\draw[0cell]
(0,0) node (a1) {X\ordn}
(a1)++(1.8,0) node (a2) {X\ordn}
(a2)++(2,0) node (a3) {X\ordn}
;
\draw[1cell=.9]
(a1) edge node {g} (a2)
(a2) edge node {X(\be g)} (a3)
(a1) [rounded corners=2pt] |- ($(a2)+(-1,\u)$) -- node {g \cdot -} ($(a2)+(1,\u)$) -| (a3)
;
\end{tikzpicture}
\end{equation}
of the $g$-action on $X\ordn$ and the image under $X$ of the pointed bijection $\be g \cn \ordn \fiso \ordn$.
\end{itemize}
Moreover, denote by $\jin \cn \ordn \fiso \nbeta$ the isomorphism in $\FG$ given by the identity morphism on $\ordn$.
\end{definition}

\begin{lemma}\label{Xjin}
In the context of \cref{def:gnbeta}, the pointed isomorphism
\begin{equation}\label{Xjin_giso}
(X\ordn)_\be \fto[\iso]{X\jin} X\nbeta
\end{equation}
is $G$-equivariant.
\end{lemma}

\begin{proof}
The proof of \cref{Xjininv} is applicable by restricting to length-1 objects in $\GG$.
\end{proof}

\section{$\FGG$-Categories from Operadic Pseudoalgebras}
\label{sec:sgo_obj}

This section constructs the object assignments of Shimakawa (strong) $H$-theory 2-functors
\[\AlglaxO \fto{\Sgo} \FGCatg \andspace \AlgpspsO \fto{\Sgosg} \FGCatg.\]
The objects of the domain 2-categories $\AlglaxO$ and $\AlgpspsO$ \pcref{oalgps_twocat} are $\Op$-pseudoalgebras.  The objects of the codomain 2-category $\FGCatg$ \pcref{def:fgcatg} are $\FGG$-categories, meaning pointed $G$-functors $\FG \to \Catgst$. 
The definitions in this section are analogous to those in \cref{sec:jemg,sec:jemg_morphisms,sec:hgo_objects}, where the roles of the indexing categories $(\Gsk,\vstar)$ and $(\GG,\vstar)$ \pcref{def:Gsk,def:GG} are now played by $(\Fsk,\ordz)$ and $(\FG,\ordz)$ \pcref{def:Fsk,def:FG}.

Recall the Cartesian closed category $\Gcat$ \pcref{def:GCat,def:Catg}.  \cref{as:Op_iconn} is in effect throughout this section.

\begin{assumption}\label{as:Op_iconn}
We assume that $\Op$ is a $\Gcat$-operad that is \emph{1-connected}\index{1-connected}\index{operad!1-connected} in the sense that 
\begin{equation}\label{i_connected}
\Op(0) = \{*\} \andspace \Op(1) = \{\opu\}
\end{equation} 
are terminal $G$-categories.  We assume that $(\A,\gaA,\phiA)$ is an $\Op$-pseudoalgebra \pcref{def:pseudoalgebra}.
\end{assumption}

For example, $\Tinf$-operads \pcref{as:OpA}, including the Barratt-Eccles operad $\BE$ and the $G$-Barratt-Eccles operad $\GBE$ \pcref{def:BE,def:GBE}, are 1-connected.  

\secoutline
\begin{itemize}
\item \cref{def:nsys} defines (strong) $\ordn$-systems in $\A$ for each pointed finite set $\ordn$.
\item \cref{def:nsys_morphism} defines the pointed categories $\Aordn$ and $\Asgordn$ of (strong) $\ordn$-systems.
\item \cref{def:nsys_gcat} defines the pointed $G$-categories $\Aordnbe$ and $\Asgordnbe$ of (strong) $\ordnbe$-systems for each pointed finite $G$-set $\ordnbe$. 
\item \cref{def:Apsi} defines the pointed functors $\Apsi$ and $\Asgpsi$ for each pointed function $\psi$ between pointed finite $G$-sets.
\item \cref{sys_FGcat} proves that each assignment
\begin{equation}\label{SgoSgosg_seci}
\begin{tikzpicture}[vcenter]
\draw[0cell]
(0,0) node (a1) {(\FG,\ordz)}
(a1)++(4,0) node (a2) {(\Catgst,\bone)}
;
\draw[1cell=.9]
(a1) edge[transform canvas={yshift=.5ex}] node {\Sgo\A \,=\, \Adash} (a2)
(a1) edge[transform canvas={yshift=-.4ex}] node[swap] {\Sgosg\A \,=\, \Asgdash} (a2)
;
\end{tikzpicture}
\end{equation}
is a pointed $G$-functor.
\end{itemize}

\subsection*{Categories of $\ordn$-Systems}
Recall the notation for unpointed finite sets \cref{ufsn}:
\[\ufs{n} = \{1,2,\ldots,n\} = \ord{n} \setminus \{0\}.\]
\cref{def:nsys} is the $\Fsk$-analogue of \cref{def:nsystem}.

\begin{definition}[$\ordn$-Systems]\label{def:nsys}
Given a pointed finite set $\ord{n} \in \Fsk$ \cref{ordn}, an \index{system!Shimakawa H-theory@Shimakawa $H$-theory}\index{Shimakawa H-theory@Shimakawa $H$-theory!system}\emph{$\ord{n}$-system in $\A$} is defined to be a pair
\begin{equation}\label{nsys}
(a,\gl)
\end{equation}
consisting of the following data.
\begin{description}
\item[Component objects] 
For each subset $s \subseteq \ufs{n}$, $(a,\gl)$ consists of an \emph{$s$-component object}\index{Shimakawa H-theory@Shimakawa $H$-theory!component object} 
\begin{equation}\label{nsys_s}
a_s \in \A.
\end{equation}
\item[Gluing]
For each object $x \in \Op(r)$ with $r \geq 0$, subset $s \subseteq \ufsn$, and partition
\[s = \coprod_{i \in \ufs{r}}\, s_i \subseteq \ufs{n},\]
$(a,\gl)$ consists of a \emph{gluing morphism}\index{Shimakawa H-theory@Shimakawa $H$-theory!gluing morphism} at $(x; s, \ang{s_i}_{i \in \ufs{r}})$:
\begin{equation}\label{gl-morphism}
\gaA_r\big(x; \ang{a_{s_i}}_{i \in \ufs{r}}\big) 
\fto{\gl_{x;\, s, \ang{s_i}_{i \in \ufs{r}}}} a_s \inspace \A.
\end{equation}
This gluing morphism is also denoted by $(a,\gl)_{x;\, s, \ang{s_i}_{i \in \ufs{r}}}$.
\end{description}
The pair $(a,\gl)$ is required to satisfy the axioms \cref{nsys_obj_unity,nsys_naturality,nsys_unity_empty,nsys_unity_one,nsys_equivariance,nsys_associativity}
whenever they are defined.
\begin{description}
\item[Object unity] 
There is an equality
\begin{equation}\label{nsys_obj_unity}
a_\emptyset = \zero = \gaA_0(*) \in \A,
\end{equation}
the basepoint of $\A$ \cref{pseudoalg_zero}.
\item[Naturality]
For each morphism $h \cn x \to y$ in the $G$-category $\Op(r)$ with $r \geq 0$, the following diagram in $\A$ commutes.
\begin{equation}\label{nsys_naturality}
\begin{tikzpicture}[vcenter]
\def\v{-1.5} \def\u{-.08}
\draw[0cell]
(0,0) node (a1) {\gaA_r\big(x; \ang{a_{s_i}}_{i \in \ufs{r}} \big)}
(a1)++(0,\v) node (b1) {\gaA_r\big(y; \ang{a_{s_i}}_{i \in \ufs{r}} \big)}
(a1)++(3.8,0) node (a2) {a_s}
(a2)++(0,\v) node (b2) {a_s}
;
\draw[1cell=.9]
(a1) edge node {\gl_{x;\, s, \ang{s_{i}}_{i \in \ufs{r}}}} (a2)
(b1) edge node {\gl_{y;\, s, \ang{s_{i}}_{i \in \ufs{r}}}} (b2)
(a1) edge[transform canvas={xshift=2em}] node[swap] {\gaA_r(h; \ang{1}_{i \in \ufs{r}})} (b1)
(a2) edge[equal] (b2)
;
\end{tikzpicture}
\end{equation}
\item[Unity]
The gluing morphism \cref{gl-morphism} is the identity morphism in each of the following two cases.
\begin{itemize}
\item If $s = \emptyset \subseteq \ufsn$---which implies $s_i = \emptyset$ for each $i \in \ufsr$---then the following diagram in $\A$ commutes.
\begin{equation}\label{nsys_unity_empty}
\begin{tikzpicture}[vcenter]
\def\v{-1.2} \def\u{-.08}
\draw[0cell]
(0,0) node (a1) {\gaA_r\big(x \sscs \ang{a_{s_i}}_{i \in \ufs{r}} \big)}
(a1)++(0,\v) node (b1) {\gaA_r\big(x \sscs \ang{\zero}_{i \in \ufs{r}} \big)}
(b1)++(2.6,0) node (b) {\gaA_0(*) = \zero}
(a1)++(5,0) node (a2) {\phantom{a_s}}
(a2)++(0,\u) node (a2') {a_s}
(a2)++(0,\v) node (b2) {\zero}
;
\draw[1cell]
(a1) edge node {\gl_{x;\, s, \ang{s_i}_{i \in \ufs{r}}}} (a2)
(b) edge node {1} (b2)
(a1) edge[equal] (b1)
(a2') edge[equal] (b2)
(b1) edge[equal] node {\mathbf{b}} (b)
;
\end{tikzpicture}
\end{equation}
The vertical equalities in \cref{nsys_unity_empty} follow from the object unity axiom \cref{nsys_obj_unity}.  The equality labeled $\mathbf{b}$ follows from \cref{phi_id}.
\item If $r=1$---which implies $x = \opu \in \Op(1)$---then $\gl_{\opu;\, s,s}$ is equal to $1_{a_s}$:
\begin{equation}\label{nsys_unity_one}
\gaA_1(\opu; a_s) = a_s \fto{\gl_{\opu;\, s,s} = 1_{a_s}} a_s.
\end{equation}
The equality in the domain follows from the action unity axiom \cref{pseudoalg_action_unity} for $\A$.
\end{itemize}
\item[Equivariance]
For each permutation $\si \in \Sigma_r$, the following diagram in $\A$ commutes.
\begin{equation}\label{nsys_equivariance}
\begin{tikzpicture}[vcenter]
\def\v{-1.2} \def\u{-.08}
\draw[0cell]
(0,0) node (a1) {\gaA_r\big(x\si ; \ang{a_{s_i}}_{i \in \ufs{r}} \big)}
(a1)++(0,\v) node (b1) {\gaA_r\big(x ; \ang{a_{s_{\sigmainv(i)}}}_{i \in \ufs{r}} \big)}
(a1)++(4.2,0) node (a2) {a_s}
(a2)++(0,\v) node (b2) {a_s}
;
\draw[1cell=.9]
(a1) edge node {\gl_{x\si ;\, s, \ang{s_i}_{i \in \ufs{r}}}} (a2)
(b1) edge node {\gl_{x;\, s, \ang{s_{\sigmainv(i)}}_{i \in \ufs{r}}}} (b2)
(a1) edge[equal,shorten <=-.5ex, shorten >=-.5ex] (b1)
(a2) edge[equal] (b2)
;
\end{tikzpicture}
\end{equation}
In \cref{nsys_equivariance}, the left vertical equality follows from the action equivariance axiom \cref{pseudoalg_action_sym} for $\A$.  The bottom horizontal gluing morphism uses the partition
\[s = \coprod_{i \in \ufs{r}}\, s_{\sigmainv(i)}.\]
\item[Associativity]
Suppose we are given objects
\[\begin{split}
\left(x \sscs \ang{x_i}_{i \in \ufs{r}} \right) &\in \Op(r) \times \txprod_{i \in \ufs{r}}\, \Op(t_i) \andspace\\
\boldx = \ga\big(x\sscs \ang{x_i}_{i \in \ufs{r}} \big) &\in \Op(t)
\end{split}\]
with $t = \sum_{i \in \ufs{r}} t_i$, and partitions
\[s = \coprod_{i \in \ufs{r}} s_i \andspace s_i = \coprod_{\ell \in \ufs{t}_i} s_{i,\ell}\]
for each $i \in \ufs{r}$.  Then the following diagram in $\A$ commutes.
\begin{equation}\label{nsys_associativity}
\begin{tikzpicture}[vcenter]
\def\h{2} \def\v{-1.5} \def\t{0em} \def\c{15} \def\b{7}
\draw[0cell=.85]
(0,0) node (a) {\gaA_r\big(x \sscs \bang{\gaA_{t_i}\big(x_i \sscs \ang{a_{s_{i,\ell}}}_{\ell \in \ufs{t}_i} \big)}_{i \in \ufs{r}} \big)}
(a)++(-\h,\v) node (b) {\gaA_r\big(x \sscs \ang{a_{s_i}}_{i \in \ufs{r}} \big)}
(a)++(\h,\v) node (c) {\gaA_t\big(\boldx \sscs \ang{\ang{a_{s_{i,\ell}}}_{\ell \in \ufs{t}_i}}_{i \in \ufs{r}} \big)}
(a)++(0,\v-1) node (d) {a_s}
;
\draw[1cell=.85]
(a) edge[bend right=\c, shorten >=0ex] node[swap,pos=.8] {\gaA_r\big(x \sscs \bang{\gl_{x_i;\, s_i, \ang{s_{i,\ell}}_{\ell \in \ufs{t}_i}}}_{i \in \ufs{r}} \big)} (b)
(b) edge[bend right=\b,transform canvas={xshift=\t}] node[swap,pos=.4] {\gl_{x;\, s, \ang{s_i}_{i \in \ufs{r}}}} (d)
(a) edge[bend left=\c,transform canvas={xshift=-\t}] node[pos=.8] {\phiA_{(r;\, t_1,\ldots,t_r)}} node[swap,pos=.35] {\iso} (c)
(c) edge[bend left=\b] node[pos=.4] {\gl_{\boldx;\, s, \ang{\ang{s_{i,\ell}}_{\ell \in \ufs{t}_i}}_{i \in \ufs{r}}}} (d)
;
\end{tikzpicture}
\end{equation}
In \cref{nsys_associativity}, $\phiA_{(r;\, t_1,\ldots,t_r)}$ is the component of the associativity constraint \cref{phiA} of $\A$ at the objects $x$, $\ang{x_i}_{i \in \ufs{r}}$, and $\ang{\ang{a_{s_{i,\ell}}}_{\ell \in \ufs{t}_i}}_{i \in \ufs{r}}$.  The lower right gluing morphism uses the partition
\[s = \coprod_{i \in \ufs{r}} \coprod_{\ell \in \ufs{t}_i}\, s_{i,\ell}.\] 
\end{description}
This finishes the definition of an $\ord{n}$-system $(a,\gl)$ in $\A$.  Moreover, we define the following.
\begin{itemize}
\item A \emph{strong $\ord{n}$-system}\index{Shimakawa H-theory@Shimakawa $H$-theory!strong system} is an $\ord{n}$-system such that each gluing morphism $\gl_{x;\, s, \ang{s_i}_{i \in \ufs{r}}}$ in \cref{gl-morphism} is an isomorphism.
\item The \emph{base $\ord{n}$-system}\index{Shimakawa H-theory@Shimakawa $H$-theory!base system} is the $\ord{n}$-system $(\zero,1_\zero)$ with
\begin{itemize}
\item each component object given by the basepoint $\zero \in \A$ and
\item each gluing morphism given by the identity morphism $1_\zero$ of $\zero$.\defmark
\end{itemize}
\end{itemize}
\end{definition}

\begin{explanation}\label{expl:nsys}
Consider \cref{def:nsys} of an $\ordn$-system.
\begin{enumerate}
\item\label{expl:nsys_i} The $r=0$ case of the unity axiom \cref{nsys_unity_empty} means the morphism equality
\begin{equation}\label{nsys_unity_zero}
\gaA_0(*) = \zero \fto{\gl_{*;\, \emptyset,\ang{}} = 1_\zero} a_\emptyset = \zero.
\end{equation}
\item\label{expl:nsys_ii} 
In \cref{def:nsystem} of an $\angordn$-system, $\Op$ is assumed to be a $\Tinf$-operad \pcref{as:OpA}.  The pseudo-commutative structure of $\Op$ is only used in the commutativity axiom \cref{system_commutativity}, which only applies when $\angordn \in \Gsk$ has length $q > 1$.  Thus, for an object $\angordn \in \Gsk$ of length 1 \pcref{def:ifg}, \cref{def:nsystem} reduces to \cref{def:nsys}.\defmark
\end{enumerate}
\end{explanation}

\cref{def:nsys_morphism} is the $\Fsk$-analogue of \cref{def:nsystem_morphism}.

\begin{definition}[Morphisms of $\ord{n}$-Systems]\label{def:nsys_morphism}
Given $\ord{n}$-systems $(a,\gl^a)$ and $(b,\gl^b)$ in $\A$ \cref{nsys} for a pointed finite set $\ord{n} \in \Fsk$ \cref{ordn}, a \emph{morphism} of $\ord{n}$-systems
\begin{equation}\label{nsys_mor}
(a,\gl^a) \fto{\tha} (b,\gl^b)
\end{equation}
consists of, for each subset $s \subseteq \ufs{n}$, an \emph{$s$-component morphism}
\begin{equation}\label{theta_s}
a_s \fto{\theta_s} b_s \inspace \A
\end{equation}
such that the axioms \cref{nsys_mor_unity,nsys_mor_compat} are satisfied.
\begin{description}
\item[Unity] If $s = \emptyset$, then there is an equality of morphisms
\begin{equation}\label{nsys_mor_unity}
a_\emptyset = \zero \fto{\theta_{\emptyset} \,=\, 1_\zero} b_\emptyset = \zero. 
\end{equation}
The equalities $a_\emptyset = b_\emptyset = \zero$ follow from the object unity axiom \cref{nsys_obj_unity}.
\item[Compatibility] For each object $x \in \Op(r)$ with $r \geq 0$, subset $s \subseteq \ufsn$, and partition
\[s = \coprod_{i \in \ufs{r}} \, s_i \subseteq \ufs{n},\]
the following diagram in $\A$ commutes.
\begin{equation}\label{nsys_mor_compat}
\begin{tikzpicture}[vcenter]
\def\v{-1.5} \def\u{-.08}
\draw[0cell]
(0,0) node (a1) {\gaA_r\big(x \sscs \ang{a_{s_i}}_{i \in \ufs{r}} \big)}
(a1)++(0,\v) node (b1) {\gaA_r\big(x \sscs \ang{b_{s_i}}_{i \in \ufs{r}} \big)}
(a1)++(3.6,0) node (a2) {a_s}
(a2)++(0,\v) node (b2) {b_s}
;
\draw[1cell=.9]
(a1) edge node {\gl^a_{x;\, s, \ang{s_i}_{i \in \ufs{r}}}} (a2)
(b1) edge node {\gl^b_{x;\, s, \ang{s_i}_{i \in \ufs{r}}}} (b2)
(a1) edge[transform canvas={xshift=2em}] node[swap] {\gaA_r\big(1_x \sscs \ang{\theta_{s_i}}_{i \in \ufs{r}} \big)} (b1)
(a2) edge node {\theta_s} (b2)
;
\end{tikzpicture}
\end{equation}
\end{description}
This finishes the definition of a morphism of $\ord{n}$-systems.  Moreover, we define the following.
\begin{itemize}
\item Composition and identity morphisms of $\ord{n}$-systems are defined componentwise using \cref{theta_s}. 
\item Define the pointed category 
\begin{equation}\label{A_ordn}
\Aordn = \big(\Aordn, (\zero,1_\zero)\big)
\end{equation}
of $\ord{n}$-systems in $\A$\index{system!pointed category} \pcref{def:nsys} and their morphisms, with the basepoint given by the base $\ordn$-system $(\zero,1_\zero)$.
\item Define the pointed category 
\begin{equation}\label{Asgordn}
\Asgordn
\end{equation} 
of \emph{strong $\ordn$-systems in $\A$} as the pointed full subcategory of $\Aordn$ with strong $\ordn$-systems in $\A$ as objects.\defmark
\end{itemize}
\end{definition}

\begin{example}\label{ex:nsys_zero_one}
We consider \cref{def:nsys,def:nsys_morphism}.
\begin{enumerate}
\item\label{ex:nsys_zero} For the object $\ord{0} = \{0\} \in \Fsk$, there are equalities
\[\Aordz = \Asgordz = \big\{(\zero,1_\zero)\big\}.\]
Since $\Aordz$ and $\Asgordz$ are terminal categories, we identify them with $\bone$
\item\label{ex:nsys_one} For the object $\ord{1} = \{0,1\} \in \Fsk$, there is a canonical isomorphism of pointed categories
\[(\A,\zero) \iso \Aordone = \Asgordone.\]
Under this isomorphism, an object $b \in \A$ corresponds to the strong $\ord{1}$-system with
\begin{itemize}
\item $\{1\}$-component object given by $b$,
\item $\emptyset$-component object given by the basepoint $\zero \in \A$, and
\item each gluing morphism given by the identity.
\end{itemize}
\item\label{ex:unity} The unity axiom \cref{nsys_mor_unity} is actually the $r=0$ case of the compatibility axiom \cref{nsys_mor_compat}.  The unity axiom is stated separately because, in practice, this case is checked separately.\defmark
\end{enumerate}
\end{example}

\subsection*{$G$-Categories of $\ordnbe$-Systems}

\cref{def:nsys_gcat} is the $\FG$-analogue of \cref{def:nbeta_gcat}.

\begin{definition}\label{def:nsys_gcat}
Given a pointed finite $G$-set $\ordnbe \in \FG$ \cref{ordn_be}, we define the small pointed $G$-category $\Aordnbe$ as follows.  A strong variant is defined in \cref{Asgordnbe}.
\begin{description}
\item[Pointed category]
The underlying pointed category of $\Aordnbe$ is defined as
\begin{equation}\label{A_nbeta}
\Aordnbe = \Aordn,
\end{equation}
the pointed category of $\ordn$-systems in $\A$ \cref{A_ordn}.  An object in $\Aordnbe$, called an \emph{$\ordnbe$-system}, consists of the same data as an $\ordn$-system \cref{nsys}, and similarly for morphisms \pcref{def:nsys_morphism}.
\item[$G$-action on $\ordnbe$-systems] 
 For each $g \in G$ and $\ordnbe$-system $(a,\gl) \in \Aordnbe$, the $\ordnbe$-system in $\A$ 
\begin{equation}\label{nsys_gaction}
g \cdot (a,\gl) = (ga, g\gl)
\end{equation}
is defined as follows.
\begin{description}
\item[Component objects] 
For each subset $s \subseteq \ufs{n}$, its $s$-component object is given by
\begin{equation}\label{ga_scomp}
(ga)_{s} = g a_{\ginv s} \in \A.
\end{equation}
On the right-hand side of \cref{ga_scomp}, the subset is given by
\begin{equation}\label{ginv_s}
\ginv s = (\be g)^{-1} s \subseteq \ufs{n}.
\end{equation}
The object $a_{\ginv s}$ is the $\ginv s$-component object of $(a,\gl)$, and $g a_{\ginv s}$ is its image under the $g$-action on $\A$.  With the $G$-action $\beta$ on $\ord{n}^\be$ understood, we abbreviate $\be g$ to $g$.  
\item[Gluing] 
For an object $x \in \Op(r)$ with $r \geq 0$, a subset $s \subseteq \ufsn$, and a partition
\[s = \coprod_{i \in \ufs{r}}\, s_{i} \subseteq \ufs{n},\]
the gluing morphism of $(ga,g\gl)$ at $(x; s, \ang{s_i}_{i \in \ufsr})$ is defined by the following commutative diagram in $\A$.
\begin{equation}\label{ga_gl}
\begin{tikzpicture}[vcenter]
\def\v{-1}
\draw[0cell=.9]
(0,0) node (a1) {\gaA_r\big(x \sscs \ang{(ga)_{s_i}}_{i \in \ufs{r}} \big)}
(a1)++(0,\v) node (a2) {\gaA_r\big(x \sscs \ang{g a_{\ginv s_i}}_{i \in \ufs{r}}\big)}
(a2)++(0,\v) node (a3) {g \gaA_r \big(\ginv x \sscs \ang{a_{\ginv s_i}}_{i \in \ufs{r}} \big)}
(a1)++(5.5,0) node (b1) {(ga)_s}
(b1)++(0,2*\v) node (b3) {g a_{\ginv s}}
;
\draw[1cell=.85]
(a1) edge node {(g\gl)_{x;\, s, \ang{s_i}_{i \in \ufs{r}}}} (b1)
(a3) edge node {g \gl_{\ginv x;\, \ginv s, \ang{\ginv s_i}_{i \in \ufs{r}}}} (b3)
(a1) edge[equal,shorten >=-.5ex] (a2)
(a2) edge[equal,shorten >=-.5ex] node[swap] {(\mathbf{f})} (a3)
(b1) edge[equal] (b3)
;
\end{tikzpicture}
\end{equation}
\begin{itemize}
\item The two unlabeled equalities in \cref{ga_gl} follow from the definition \cref{ga_scomp} of $(ga)_{s}$ and the partition
\begin{equation}\label{ginv_s_partition}
\ginv s = \coprod_{i \in \ufs{r}}\, \ginv s_i \subseteq \ufs{n}.
\end{equation}
\item The equality labeled $(\mathbf{f})$ follows from the $G$-functoriality of $\gaA_r$ \cref{gaAn}.
\item In the bottom horizontal arrow in \cref{ga_gl}, 
\[\gl_{\ginv x;\, \ginv s, \ang{\ginv s_i}_{i \in \ufs{r}}} \cn 
\gaA_r \big(\ginv x \sscs \ang{a_{\ginv s_i}}_{i \in \ufs{r}} \big)
\to a_{\ginv s}\]
is the indicated gluing morphism of $(a,\gl)$, and $g\gl_{\cdots}$ is its image under the $g$-action on $\A$.
\end{itemize}
\item[Axioms]
Each of the axioms of an $\ord{n}$-system,  \cref{nsys_obj_unity,nsys_naturality,nsys_unity_empty,nsys_unity_one,nsys_equivariance,nsys_associativity}, for $(ga,g\gl)$ follows from the corresponding axiom for $(a,\gl)$ and the following facts.
\begin{itemize}
\item The $g$-action on $\A$ is a functor.
\item The basepoint $\zero \in \A$ and its identity morphism $1_\zero$ are both $G$-fixed.  In particular, the base $\ordnbe$-system $(\zero,1_\zero)$ is $G$-fixed. 
\item The object $* \in \Op(0)$ and the operadic unit $\opu \in \Op(1)$ are both $G$-fixed.
\item The right symmetric group action on $\Op$ and the associativity constraint $\phiA$ of $\A$ \cref{phiA} are $G$-equivariant.
\end{itemize}
\end{description}
This finishes the definition of the $\ordnbe$-system $g \cdot (a,\gl) = (ga,g\gl)$.
\item[$G$-action on morphisms] 
For a morphism of $\ordnbe$-systems in $\A$ \pcref{def:nsys_morphism}
\[(a,\gl^a) \fto{\tha} (b,\gl^b),\]
the morphism of $\ordnbe$-systems
\begin{equation}\label{g_tha}
(ga,g\gl^a) \fto{g\theta} (gb, g\gl^b)
\end{equation}
is defined by, for each subset $s \subseteq \ufs{n}$, the $s$-component morphism
\begin{equation}\label{gtha_s}
(ga)_{s} = ga_{\ginv s}
\fto{(g\theta)_{s} \,=\, g\theta_{\ginv s}}
(gb)_{s} = gb_{\ginv s}.
\end{equation}
\begin{itemize}
\item The unity axiom \cref{nsys_mor_unity} holds for $g\theta$ because the identity morphism $1_\zero$ is $G$-fixed. 
\item The compatibility axiom \cref{nsys_mor_compat} holds for $g\theta$ by the compatibility axiom for $\theta$ and the functoriality of the $g$-action on $\A$.
\end{itemize}
This finishes the definition of the morphism $g\theta$ of $\ordnbe$-systems.  
\item[Functoriality of $G$-action]
For each $g \in G$, the functoriality of the $g$-action on $\Aordnbe$, as defined in \cref{nsys_gaction,g_tha}, follows from
\begin{itemize}
\item the definition \cref{gtha_s} of $(g\theta)_s$ and
\item the functoriality of the $g$-action on $\A$.
\end{itemize}
As $g \in G$ varies, this defines a $G$-action on the pointed category $\Aordnbe$ by
\begin{itemize}
\item the definitions \cref{ga_scomp,ginv_s,ga_gl,gtha_s}, and 
\item the $G$-action axioms for $\A$, $\Op(r)$, and $\ordnbe$.
\end{itemize}
\end{description}
This finishes the definition of the pointed $G$-category $\Aordnbe$.

\parhead{Strong variant}.  
The small pointed $G$-category 
\begin{equation}\label{Asgordnbe}
\Asgordnbe
\end{equation}
has underlying pointed category $\Asgordn$ \cref{Asgordn}.
\begin{itemize}
\item An object in $\Asgordnbe$ is called a \emph{strong $\ordnbe$-system}, which consists of the same data as a strong $\ordn$-system.  Recall from the last paragraph of \cref{def:nsys} that an $\ordn$-system is \emph{strong} if its gluing morphisms are isomorphisms. 
\item The $G$-action on $\Asgordnbe$ is the restriction of the $G$-action on $\Aordnbe$ to the full subcategory of strong $\ordnbe$-systems.  This $G$-action is well defined because, if each gluing morphism $\gl_{\cdots}$ is an isomorphism in \cref{ga_gl}, then so is its image $g\gl_{\cdots}$ under the $g$-action.\defmark
\end{itemize}
\end{definition}

\subsection*{Functors Between $G$-Categories of Systems}

Next, we define functors between $G$-categories of systems associated to morphisms of $\FG$ \pcref{def:FG}, which are pointed functions between underlying pointed finite sets.

\begin{definition}\label{def:Apsi}
Given a pointed function $\psi \cn \ordmal \to \ordnbe$ in $\FG$, we define a pointed functor
\begin{equation}\label{psitil_f}
\Aordmal \fto{\Apsi} \Aordnbe
\end{equation}
between the pointed categories in \cref{A_nbeta} as follows.  A strong variant is defined in \cref{psitilsg_f}.
\begin{description}
\item[Component objects] 
Given an $\ordmal$-system $(a,\gl) \in \Aordmal$ \cref{nsys}, the $\ordnbe$-system
\[(\Apsi) (a,\gl) = \big(a^{\psitil}, \gl^{\psitil}\big) \in \Aordnbe\]
has, for each subset $s \subseteq \ufsn$, $s$-component object defined as the $\psiinv s$-component object of $(a,\gl)$:
\begin{equation}\label{apsitil_s}
a^{\psitil}_s = a_{\psiinv s}.
\end{equation}
\item[Gluing] 
For an object $x \in \Op(r)$ with $r \geq 0$, a subset $s \subseteq \ufsn$, and a partition
\[s = \coprod_{i \in \ufs{r}}\, s_i \subseteq \ufs{n},\]
the gluing morphism \cref{gl-morphism} of $(\Apsi)(a,\gl)$ at $(x; s, \ang{s_i}_{i \in \ufsr})$ is defined by the following commutative diagram in $\A$.
\begin{equation}\label{glpsitil}
\begin{tikzpicture}[vcenter]
\def\v{-1.2}
\draw[0cell=1]
(0,0) node (a1) {\gaA_r\big(x \sscs \ang{a^{\psitil}_{s_i}}_{i \in \ufs{r}}\big)}
(a1)++(4.7,0) node (a2) {a^{\psitil}_s}
(a1)++(0,\v) node (b1) {\gaA_r\big(x \sscs \ang{a_{\psiinv s_i}}_{i \in \ufs{r}}\big)}
(a2)++(0,\v) node (b2) {a_{\psiinv s}}
;
\draw[1cell=.9]
(a1) edge node {\gl^{\psitil}_{x;\, s, \ang{s_i}_{i \in \ufs{r}}}} (a2)
(b1) edge node {\gl_{x;\, \psiinv s, \ang{\psiinv s_i}_{i \in \ufs{r}}}} (b2)
(a1) edge[equal,shorten >=-.3ex] (b1)
(a2) edge[equal] (b2)
;
\end{tikzpicture}
\end{equation}
Each of the axioms of an $\ord{n}$-system, \cref{nsys_obj_unity,nsys_naturality,nsys_unity_empty,nsys_unity_one,nsys_equivariance,nsys_associativity}, holds for $\big(a^{\psitil}, \gl^{\psitil}\big)$ by the corresponding axiom for $(a,\gl)$ and the fact that taking preimages preserves the empty set and partitions.
\item[Morphisms]
For a morphism of $\ordmal$-systems in $\A$ \pcref{def:nsys_morphism}
\[(a,\gl^a) \fto{\theta} (b,\gl^b),\]
the morphism of $\ordnbe$-systems
\[\big(a^{\psitil}, \gl^{a, \psitil}\big) \fto{(\Apsi)(\theta) = \tha^{\psitil}} \big(b^{\psitil}, \gl^{b, \psitil}\big)\]
has, for each subset $s \subseteq \ufs{n}$, $s$-component morphism defined as the $\psiinv s$-component morphism of $\theta$:
\begin{equation}\label{thapsitil_comp}
a^{\psitil}_s = a_{\psiinv s} \fto{\tha^{\psitil}_s \,=\, \theta_{\psiinv s}} 
b^{\psitil}_s = b_{\psiinv s}.
\end{equation}
Each of the unity axiom \cref{nsys_mor_unity} and the compatibility axiom \cref{nsys_mor_compat} for a morphism of $\ordnbe$-systems holds for $\tha^{\psitil}$ by the corresponding axiom for $\theta$ and the fact that taking preimages preserves the empty set and partitions.
\item[Pointed Functoriality]
The functoriality of $\Apsi$ follows from the definition \cref{thapsitil_comp} of $\tha^{\psitil}_s$.  The functor $\Apsi$ is pointed, sending the base $\ordmal$-system $(\zero,1_\zero)$ to the base $\ordnbe$-system, by \cref{apsitil_s,glpsitil}.
\end{description}
This finishes the definition of the pointed functor $\Apsi$.

\parhead{Strong variant}.  For the pointed categories $\Asgordmal$ and $\Asgordnbe$ \cref{Asgordnbe}, we define the pointed subfunctor
\begin{equation}\label{psitilsg_f}
\Asgordmal \fto{\Asgpsi} \Asgordnbe
\end{equation}
of $\Apsi$ by restricting to strong $\ordmal$-systems.  This yields a well-defined pointed functor $\Asgpsi$ because, if each component of the gluing morphism $\gl$ is an isomorphism, then each component of $\gl^{\psitil}$ in \cref{glpsitil} is an isomorphism.
\end{definition}

While \cref{def:Apsi} is analogous to \cref{def:psitil_functor}, we emphasize that the pointed functors $\Apsi$ \cref{psitil_f} and $\Asgpsi$ \cref{psitilsg_f} are \emph{not} generally $G$-equivariant.  The $G$-equivariance of $\Adash$ takes the form of the commutative diagram \cref{Adash_Gequiv}.  \cref{sys_FGcat} is the $\FG$-analogue of \cref{A_ptfunctorGG}.  

\begin{lemma}[Shimakawa $H$-Theory on Objects]\label{sys_FGcat}\index{Shimakawa H-theory@Shimakawa $H$-theory!object}
Under \cref{as:Op_iconn}, the object and morphism assignments 
\[\ordnbe \mapsto \Aordnbe \andspace \psi \mapsto \Apsi\]
in, respectively, \cref{def:nsys_gcat,def:Apsi} define a pointed $G$-functor
\[(\FG,\ordz) \fto{\Sgo\A \,=\, \Adash} (\Catgst,\bone).\]
Moreover, the strong variant 
\[(\FG,\ordz) \fto{\Sgosg\A \,=\, \Asgdash} (\Catgst,\bone).\]
defined in \cref{Asgordnbe,psitilsg_f} is also a pointed $G$-functor.
\end{lemma}

\begin{proof}
We prove that $\Sgo\A = \Adash$ is a pointed $G$-functor.  The proof for the strong variant $\Asgdash$ is the same after restricting to strong systems.  
\begin{description}
\item [Pointedness]
By \cref{A_nbeta} and \cref{ex:nsys_zero_one} \eqref{ex:nsys_zero}, $\Aordz = \bone$, so $\Adash$ is pointed.  
\item[Identities]
For an object $\ordmal \in \FG$ with identity morphism $1_{\ord{m}^\al}$, the pointed functor
\[\sys{\A}{\ord{m}^\al} \fto{\sys{\A}{1_{\ord{m}^\al}}} \sys{\A}{\ord{m}^\al}\]
is the identity functor by \cref{apsitil_s,glpsitil,thapsitil_comp}, since $1_{\ord{m}^\al}^{-1}s = s$ for each subset $s \subseteq \ufsn$.
\item[Composition]
For morphisms
\[\ordmal \fto{\psi_1} \ordnbe \fto{\psi_2} \ordlde \inspace \FG,\]
the composite functor 
\[\Aordmal \fto{\sys{\A}{\psi_1}} \Aordnbe \fto{\sys{\A}{\psi_2}} \Aordlde\]
is equal to $\sys{\A}{(\psi_2 \psi_1)}$ by \cref{apsitil_s,glpsitil,thapsitil_comp}, since 
\[\psiinv_1 \psiinv_2 s = (\psi_2 \psi_1)^{-1} s \forspace s \subseteq \ufsl.\]
This proves that $\Adash$ is a pointed functor.
\item[$G$-equivariance]  
By \cref{FGGcat_Gequiv}, the $G$-equivariance of $\Adash$ means that the following diagram of pointed functors commutes for each $g \in G$ and morphism $\psi \cn \ordmal \to \ordnbe$ in $\FG$.
\begin{equation}\label{Adash_Gequiv}
\begin{tikzpicture}[vcenter]
\def\v{-1.3}
\draw[0cell]
(0,0) node (a1) {\Aordmal}
(a1)++(3,0) node (a2) {\Aordnbe}
(a1)++(0,\v) node (b1) {\Aordmal}
(a2)++(0,\v) node (b2) {\Aordnbe}
;
\draw[1cell=.9]
(a1) edge node {\sys{\A}(g\psi\ginv)} (a2)
(a1) edge node[swap] {\ginv} (b1)
(b1) edge node {\Apsi} (b2)
(b2) edge[shorten <=-.3ex] node[swap] {g} (a2)
;
\end{tikzpicture}
\end{equation}
\begin{description}
\item[Component objects]
Using \cref{ga_scomp}, \cref{apsitil_s}, and $g\ginv = 1$, the following computation shows that the two composites in \cref{Adash_Gequiv} agree on each $\ordmal$-system $(a,\gl)$ and subset $s \subseteq \ufs{n}$.
\begin{equation}\label{Adash_Geq_objcomp}
\begin{split}
&\big(\sys{\A}(g\psi\ginv) (a,\gl) \big)_s\\
&= a_{(g\psi\ginv)^{-1} s}\\
&= g\big(\ginv a_{(\ginv)^{-1} \psiinv \ginv s}\big)\\
&= g\big( (\ginv a)_{\psiinv \ginv s} \big)\\
&= g\big((\sys{\A}{\psi}) (\ginv a) \big)_{\ginv s}\\ 
&= \big((g (\sys{\A}{\psi}) \ginv) (a,\gl) \big)_s
\end{split}
\end{equation}
\item[Gluing]
Using \cref{ga_gl}, \cref{glpsitil}, and $g\ginv = 1$, the following computation shows that the two composites in \cref{Adash_Gequiv} yield the same gluing morphism at $(x; s, \ang{s_i}_{i \in \ufs{r}})$ for each object $x \in \Op(r)$ with $r \geq 0$, subset $s \subseteq \ufsn$, and partition $s = \coprod_{i \in \ufs{r}}\, s_i$.
\begin{equation}\label{Adash_geq_gluing}
\begin{split}
&\big(\sys{\A}(g\psi\ginv) (a,\gl) \big)_{x;\, s, \ang{s_i}_{i \in \ufs{r}}}\\
&= \gl_{x;\, (g\psi\ginv)^{-1} s,\, \ang{(g\psi\ginv)^{-1} s_i}_{i \in \ufs{r}}}\\
&= g \big(\ginv \gl_{(\ginv)^{-1}\ginv x;\, (\ginv)^{-1}\psiinv\ginv s,\, \ang{(\ginv)^{-1}\psiinv\ginv s_i}_{i \in \ufs{r}}}  \big)\\
&= g \big( (\ginv \gl)_{\ginv x;\, \psiinv\ginv s,\, \ang{\psiinv\ginv s_i}_{i \in \ufs{r}}} \big)\\
&= g \big( ((\sys{\A}{\psi}) (\ginv \gl))_{\ginv x;\, \ginv s,\, \ang{\ginv s_i}_{i \in \ufs{r}}} \big)\\
&= \big((g (\sys{\A}{\psi}) \ginv) (a,\gl) \big)_{x;\, s, \ang{s_i}_{i \in \ufs{r}}}
\end{split}
\end{equation}
This proves that the diagram \cref{Adash_Gequiv} commutes on objects.  
\item[Morphisms]
The proof that the diagram \cref{Adash_Gequiv} commutes on morphisms of $\ordmal$-systems is the same as \cref{Adash_Geq_objcomp}, using \cref{gtha_s,thapsitil_comp} instead of \cref{ga_scomp,apsitil_s}.\qedhere
\end{description}
\end{description}
\end{proof}

\section{$H$-Theory and $J$-Theory 2-Functors}
\label{sec:sgo_iifunctor}

Extending the object assignments in \cref{sys_FGcat}, this section constructs Shimakawa (strong) $H$-theory 2-functors
\[\AlglaxO \fto{\Sgo} \FGCatg \andspace \AlgpspsO \fto{\Sgosg} \FGCatg\]
between the 2-categories in \cref{oalgps_twocat,def:fgcatg}.  Shimakawa (strong) $J$-theory 2-functors are the composites
\[\AlglaxO \fto{\Jgos = \igst\Sgo} \FGCat \andspace \AlgpspsO \fto{\Jgossg = \igst\Sgosg} \FGCat\]
of Shimakawa (strong) $H$-theory with the 2-equivalence $\igst \cn \FGCatg \to \FGCat$ \pcref{thm:fgcat_fgcatg_iieq}.

\secoutline
\begin{itemize}
\item \cref{def:sgo} defines the 1-cell and 2-cell assignments of $\Sgo$ and $\Sgosg$.
\item \cref{Sgo_twofunctor} records the fact that $\Sgo$ and $\Sgosg$ are 2-functors.
\item \cref{expl:shi89_section2} discusses the fact that, for the Barratt-Eccles operad $\BE$, $\Sgosgbe$ recovers Shimakawa's original construction of a $\GaG$-category from a naive symmetric monoidal $G$-category.
\item \cref{def:jgos} defines Shimakawa (strong) $J$-theory $\Jgos$ and $\Jgossg$.
\end{itemize}

\subsection*{Shimakawa $H$-Theory}

\begin{definition}[Shimakawa $H$-Theory]\label{def:sgo}
For a 1-connected $\Gcat$-operad $\Op$ \pcref{as:Op_iconn}, we define the 2-functor
\[\AlglaxO \fto{\Sgo} \FGCatg,\]
called \emph{Shimakawa $H$-theory}\index{Shimakawa H-theory@Shimakawa $H$-theory} for $\Op$, using the assignments on objects, 1-cells, and 2-cells in, respectively, \cref{sgo_object,sgo_1cell,sgo_2cell}.  A strong variant is defined in \cref{def:sgosg}.
\begin{description}
\item[Objects] 
$\Sgo$ sends each $\Op$-pseudoalgebra $\A$ \pcref{def:pseudoalgebra} to the pointed $G$-functor in \cref{sys_FGcat}:
\begin{equation}\label{sgo_object}
(\FG,\ordz) \fto{\Sgo\A = \Adash} (\Catgst,\bone).
\end{equation}
Its object and morphism assignments are given in, respectively, \cref{def:nsys_gcat,def:Apsi}.
\item[1-cells] 
For a lax $\Op$-morphism \pcref{def:laxmorphism}
\[\big(\A,\gaA,\phiA\big) \fto{(f,\actf)} \big(\B,\gaB,\phiB\big)\]
between $\Op$-pseudoalgebras, the $G$-natural transformation
\begin{equation}\label{sgo_1cell}
\def\t{25}
\begin{tikzpicture}[baseline={(a.base)}]
\draw[0cell]
(0,0) node (a) {\FG}
(a)++(2.5,0) node (b) {\phantom{\FG}}
(b)++(.15,0) node (b') {\Catgst}
;
\draw[1cell=.9]
(a) edge[bend left=\t] node {\Adash} (b)
(a) edge[bend right=\t] node[swap] {\Bdash} (b)
;
\draw[2cell]
node[between=a and b at .35, rotate=-90, 2label={above,\Sgo f}] {\Rightarrow}
;
\end{tikzpicture}
\end{equation} 
has, for each pointed finite $G$-set $\ordnbe \in \FG$, $\ordnbe$-component pointed $G$-functor
\begin{equation}\label{sgosf_nbe}
\Aordnbe \fto{(\Sgo f)_{\ordnbe}} \Bordnbe
\end{equation}
defined as follows.
\begin{description}
\item[Component objects]
For an $\ordnbe$-system $(a,\gl^a) \in \Aordnbe$ \cref{nsys} and a subset $s \subseteq \ufs{n}$, the $s$-component object of the $\ordnbe$-system $(\Sgo f)_{\ordnbe} (a,\gl^a) \in \Bordnbe$ is defined as
\begin{equation}\label{sgof_nbe_a}
\big((\Sgo f)_{\ordnbe} (a,\gl^a)\big)_s = f a_s \in \B.
\end{equation}
\item[Gluing] 
For each object $x \in \Op(r)$ with $r \geq 0$, subset $s \subseteq \ufsn$, and partition $s = \coprod_{i \in \ufs{r}}\, s_i$, the gluing morphism of the $\ordnbe$-system $(\Sgo f)_{\ordnbe} (a,\gl^a)$ at $(x; s, \ang{s_i}_{i \in \ufs{r}})$ is defined as the following composite in $\B$.
\begin{equation}\label{sgof_nbe_gl}
\begin{tikzpicture}[vcenter]
\def\u{-1} \def\h{4.5} \def\a{10} \def\b{.7}
\draw[0cell=.85]
(0,0) node (a1) {\gaB_r\big(x \sscs \bang{\big((\Sgo f)_{\ordnbe} (a,\gl)\big)_{s_i}}_{i \in \ufs{r}}\big)}
(a1)++(\h,0) node (a2) {\big((\Sgo f)_{\ordnbe} (a,\gl)\big)_s}
(a1)++(0,\u) node (b1) {\gaB_r\big(x \sscs \ang{fa_{s_i}}_{i \in \ufs{r}}\big)}
(a2)++(0,\u) node (b2) {fa_s}
(b1)++(\h/2,-1.2) node (c) {f \gaA_r \big(x \sscs \ang{a_{s_i}}_{i \in \ufs{r}} \big)}
;
\draw[1cell=.9]
(a1) edge[equal,shorten <=-0ex,shorten >=-.5ex] (b1)
(a2) edge[equal] (b2)
;
\draw[1cell=.9]
(a1) [rounded corners=2pt, shorten <=-.2ex] |- ($(a1)+(1,\b)$)
-- node {\gl_{x;\, s, \ang{s_{i}}_{i \in \ufs{r}}}} ($(a2)+(-1,\b)$) -| (a2)
;
\draw[1cell=.9]
(b1) [rounded corners=2pt, shorten <=-.2ex] |- node[pos=.2,swap] {\actf_r} (c);
\draw[1cell=.9]
(c) [rounded corners=2pt, shorten <=-0ex] -| node[pos=.8] {f \gl^a_{x;\, s, \ang{s_{i}}_{i \in \ufs{r}}}} (b2);
\end{tikzpicture}
\end{equation}
The lower-left arrow $\actf_r$ is a component of the $r$-th action constraint of $f$ \cref{actf_component}.  The lower-right arrow is the image under $f$ of the indicated gluing morphism of $(a,\glu^a)$.
\item[Morphisms]
For a morphism \pcref{def:nsys_morphism}
\[(a,\gl^a) \fto{\theta} (b,\gl^b) \inspace \Aordnbe,\]
the morphism of $\ordnbe$-systems in $\B$
\begin{equation}\label{sgof_tha}
(\Sgo f)_{\ordnbe} (a,\gl^a) \fto{(\Sgo f)_{\ordnbe} \theta} (\Sgo f)_{\ordnbe} (b,\gl^b)
\end{equation}
has, for each subset $s \subseteq \ufsn$, $s$-component morphism defined as
\begin{equation}\label{sgof_tha_s}
f a_s \fto{\big((\Sgo f)_{\ordnbe} \theta\big)_s \,=\, f\theta_s} f b_s \inspace \B.
\end{equation}
\end{description}
\item[2-cells]
Suppose we are given an $\Op$-transformation $\omega$ \pcref{def:algtwocells} between lax $\Op$-morphisms between $\Op$-pseudoalgebras as follows.
\begin{equation}\label{Otransf_iicell}

\end{equation}
has, for each pointed finite $G$-set $\ordnbe \in \FG$, $\ordnbe$-component pointed $G$-natural transformation
\begin{equation}\label{sgo_omega_nbe}
%
\end{equation} 
defined by the following $s$-component morphism for each $\ordnbe$-system $(a,\gl^a) \in \Aordnbe$ and subset $s \subseteq \ufs{n}$.
\begin{equation}\label{sgo_omega_nbe_s}
%
\end{equation}
\end{description}
This finishes the definition of the 2-functor $\Sgo$.

\parhead{Strong variant}.  The 2-functor
\begin{equation}\label{def:sgosg}
\AlgpspsO \fto{\Sgosg} \FGCatg,
\end{equation}
called \emph{Shimakawa strong $H$-theory}\index{Shimakawa H-theory@Shimakawa $H$-theory!strong} for $\Op$, is defined as follows.
\begin{description}
\item[Domain] $\AlgpspsO$ is the sub-2-category of $\AlglaxO$ with the same objects and 2-cells, and with 1-cells given by $\Op$-pseudomorphisms \pcref{oalgps_twocat}.  
\item[Objects] $\Sgosg$ sends an $\Op$-pseudoalgebra $\A$ to the pointed $G$-functor in \cref{sys_FGcat}:
\[(\FG,\ordz) \fto{\Sgosg\A = \Asgdash} (\Catgst,\bone).\]
\item[1-cells] For an $\Op$-pseudomorphism $(f,\actf)$, the $G$-natural transformation
\begin{equation}\label{Sgosgf_iicell}
\begin{tikzpicture}[baseline={(a.base)}]
\def\t{25}
\draw[0cell]
(0,0) node (a) {\FG}
(a)++(2.5,0) node (b) {\phantom{\FG}}
(b)++(.15,0) node (b') {\Catgst}
;
\draw[1cell=.9]
(a) edge[bend left=\t] node {\Asgdash} (b)
(a) edge[bend right=\t] node[swap] {\Bsgdash} (b)
;
\draw[2cell]
node[between=a and b at .35, rotate=-90, 2label={above,\Sgosg f}] {\Rightarrow}
;
\end{tikzpicture}
\end{equation}
is defined by \cref{sgof_nbe_a,sgof_nbe_gl,sgof_tha_s}, applied to strong systems.  This is well defined because, in the diagram \cref{sgof_nbe_gl}, the morphisms $\actf_r$ and $\gl^a_{x;\, s, \ang{s_i}_{i \in \ufs{r}}}$, and hence also $\gl_{x;\, s, \ang{s_i}_{i \in \ufs{r}}}$, are now isomorphisms.
\item[2-cells] The 2-cell assignment of $\Sgosg$ is defined by \cref{sgo_omega_nbe_s}, applied to strong systems and $\Op$-pseudomorphisms $(f,\actf)$ and $(h,\acth)$.
\end{description}
This completes the definition of the 2-functor $\Sgosg$.
\end{definition}

\begin{proposition}[Shimakawa $H$-Theory]\label{Sgo_twofunctor}
For each 1-connected $\Gcat$-operad $\Op$ \pcref{as:Op_iconn}, the object, 1-cell, and 2-cell assignments in \cref{def:sgo} define 2-functors\index{Shimakawa H-theory@Shimakawa $H$-theory}
\[\begin{split}
\AlglaxO & \fto{\Sgo} \FGCatg \andspace\\
\AlgpspsO & \fto{\Sgosg} \FGCatg
\end{split}\]
between the 2-categories in \cref{oalgps_twocat,def:fgcatg}.
\end{proposition}

\begin{proof}
We reuse the proofs of \cref{hgo_icell_welldef,hgo_iicell_welldef,Hgo_twofunctor} by restricting to objects in $\GG$ of length 1 \pcref{def:ifgGG}.
\end{proof}

\begin{explanation}[Shimakawa's Construction]\label{expl:shi89_section2}
For the Barratt-Eccles $\Gcat$-operad $\BE$ \pcref{def:BE}, the object assignment of Shimakawa strong $H$-theory \pcref{sys_FGcat}
\[\AlgpspsBE \fto{\Sgosgbe} \FGCatg,\]
sending $\BE$-pseudoalgebras to $\FGG$-categories, recovers Shimakawa's original construction in \cite[p.\ 251--252]{shimakawa89}, sending naive symmetric monoidal $G$-categories \pcref{def:naive_smGcat} to $\GaG$-categories.  As we prove in detail in \cite[\ref*{EqK:thm:BEpseudoalg}]{yau-eqk} and briefly recall in \cref{thm:BEpseudoalg}, $\BE$-pseudoalgebras correspond to naive symmetric monoidal $G$-categories under certain 2-equivalences.  Shimakawa (strong) $H$-theory $\Sgo$ and $\Sgosg$ \pcref{Sgo_twofunctor} extend Shimakawa's original construction in the following two ways.
\begin{enumerate}
\item $\Sgo$ and $\Sgosg$ are 2-functors, taking into account 1-cells and 2-cells in their domains and codomains.
\item $\Sgo$ and $\Sgosg$ apply to any 1-connected $\Gcat$-operad $\Op$ \pcref{as:Op_iconn}, not just the Barratt-Eccles $\Gcat$-operad $\BE$.  
\end{enumerate}
For example, $\Sgo$ and $\Sgosg$ can be applied to the $G$-Barratt-Eccles operad $\GBE$, whose pseudoalgebras are \emph{genuine} symmetric monoidal $G$-categories \pcref{def:GBE,def:GBE_pseudoalg}.  See \cref{expl:c_shimakawa,expl:khsho_shi,expl:k_shimakawa} for more discussion of Shimakawa's original construction.
\end{explanation}

\subsection*{Shimakawa $J$-Theory}

Recall from \cref{jgohgoigst} that $J$-theory $\Jgo$ factors as $H$-theory $\Hgo$ followed by $\igst \cn \GGCatg \to \GGCatii$, and similarly for the strong variant.  The next definitions are the Shimakawa analogues of those factorizations.  \cref{sec:weq_app} compares strong $J$-theory and Shimakawa strong $J$-theory.

\begin{definition}[Shimakawa $J$-Theory]\label{def:jgos}
For a 1-connected $\Gcat$-operad $\Op$ \pcref{as:Op_iconn}, we define \emph{Shimakawa $J$-theory}\index{Shimakawa J-theory@Shimakawa $J$-theory}\index{J-theory@$J$-theory!Shimakawa} $\Jgos$ and \emph{Shimakawa strong $J$-theory} $\Jgossg$ for $\Op$ as the following composite 2-functors.
\begin{equation}\label{jgos_jgossg}
\begin{tikzpicture}[vcenter]
\def\h{2.5} \def\u{.7} \def\v{-1.8}
\draw[0cell]
(0,0) node (a1) {\AlglaxO}
(a1)++(.12*\h,.1*\v) node (a1') {\phantom{A}}
(a1)++(0,\v) node (a2) {\AlgpspsO}
(a2)++(.12*\h,-.1*\v) node (a2') {\phantom{A}}
(a1)++(\h,\v/2) node (a3) {\phantom{\FGCatg}}
(a3)++(0,-.02) node (a3') {\FGCatg}
(a3)++(\h,0) node (a4) {\FGCat}
;
\draw[1cell=.9]
(a1') edge[shorten >=-.4ex] node[swap, inner sep=1pt] {\Sgo} (a3)
(a2') edge[shorten >=-.4ex] node[inner sep=1pt] {\Sgosg} (a3)
(a3) edge node {\igst} (a4)
(a1) [rounded corners=2pt] -| node[pos=.25] {\Jgos}  (a4)
;
\draw[1cell=.9]
(a2) [rounded corners=2pt] -| node[pos=.25] {\Jgossg}  (a4)
;
\end{tikzpicture}
\end{equation}
The 2-functors $\Sgo$ and $\Sgosg$ are Shimakawa (strong) $H$-theory \pcref{Sgo_twofunctor}, and $\igst$ is the 2-equivalence in \cref{thm:fgcat_fgcatg_iieq}.
\end{definition}

\begin{explanation}\label{expl:jgos}
The 2-equivalence \pcref{thm:fgcat_fgcatg_iieq}
\[\FGCatg \fto{\igst} \FGCat\] 
is given by precomposition with the full subcategory inclusion $\ig \cn \Fsk \to \FG$ \cref{ig_FG}.  Thus, $\Jgos$ and $\Jgossg$ are given explicitly as follows.
\begin{description}
\item[Objects]
Shimakawa $J$-theory 
\[\AlglaxO \fto{\Jgos = \igst\Sgo} \FGCat\] 
sends an $\Op$-pseudoalgebra $\A$ \pcref{def:pseudoalgebra} to the $\Fskg$-category \cref{fgcat_obj}
\begin{equation}\label{jgosan}
(\Fsk,\ord{0}) \fto{\Jgos\A = \Adash} (\Gcatst,\bone)
\end{equation}
whose value at a pointed finite set $\ordn \in \Fsk$ is the pointed category 
\[(\Jgos\A)\ordn = \Aordn\]
of $\ordn$-systems in $\A$ \cref{A_ordn}.  The value of $\Jgos\A$ at a morphism in $\Fsk$ is given in \cref{def:Apsi}, restricted to pointed finite sets $\ordn \in \Fsk$ with trivial $G$-actions.

\parhead{$G$-action}.
The $G$-action on $(\Jgos\A)\ordn$ is constructed in \crefrange{nsys_gaction}{gtha_s}, with the trivial $G$-action on $\ordn$.  In other words, for each $g \in G$ and $\ordn$-system $(a,\gl) \in \Aordn$, the $\ordn$-system in $\A$ 
\begin{equation}\label{g_Aordn}
g \cdot (a,\gl) = (ga,g\gl)
\end{equation}
has, for each subset $s \subseteq \ufsn$, $s$-component object
\begin{equation}\label{gas_Aordn}
(ga)_s = g a_s \in \A.
\end{equation}
For an object $x \in \Op(r)$ with $r \geq 0$, a subset $s \subseteq \ufsn$, and a partition $s = \txcoprod_{i \in \ufs{r}}\, s_{i} \subseteq \ufs{n}$, the gluing morphism of $(ga,g\gl)$ at $(x; s, \ang{s_i}_{i \in \ufsr})$ is defined by the following commutative diagram in $\A$.
\begin{equation}\label{ggl_Aordn}

\end{equation}
For a morphism $\tha$ of $\ordn$-systems, the morphism $g \tha$ has $s$-component morphism
\begin{equation}\label{gtha_scomp}
(g\tha)_s = g \tha_s.
\end{equation}
The pointed $G$-category $\Aordn$ coincides with $J$-theory $(\Jgo\A)\ordn$ \pcref{def:Aangordn_gcat}, regarding $\ordn$ as an object of $\Gsk$ of length 1.
\item[1-cells]
Shimakawa $J$-theory sends a lax $\Op$-morphism $f \cn \A \to \B$ between $\Op$-pseudoalgebras \pcref{def:laxmorphism} to the natural transformation \cref{fgcat_mor}
\begin{equation}\label{jgosf}
\def\t{25}
%
\end{equation} 
as defined in \crefrange{sgosf_nbe}{sgof_tha_s}, restricted to pointed finite sets $\ordn \in \Fsk$ with trivial $G$-actions.
\item[2-cells]
Shimakawa $J$-theory sends an $\Op$-transformation $\omega$ \pcref{def:algtwocells} between lax $\Op$-morphisms between $\Op$-pseudoalgebras 
\begin{equation}\label{Otr_omega}
%
\end{equation}
to the modification \cref{fgcat_iicell}
\begin{equation}\label{jgosomega}
\def\t{25}
%
\end{equation}
as defined in \cref{sgo_omega_nbe,sgo_omega_nbe_s}, restricted to pointed finite sets $\ordn \in \Fsk$ with trivial $G$-actions.
\item[Strong variant]
Shimakawa strong $J$-theory 
\[\AlgpspsO \fto{\Jgossg = \igst\Sgosg} \FGCat\] 
is obtained in the same way from Shimakawa strong $H$-theory $\Sgosg$ \cref{def:sgosg}, by restricting to pointed finite sets $\ordn \in \Fsk$ with trivial $G$-actions.  The preceding discussion, \crefrange{jgosan}{jgosomega}, applies to $\Jgossg$ by restricting to strong $\ordn$-systems \cref{Asgordn}, with invertible gluing morphisms, and $\Op$-pseudomorphisms \pcref{def:laxmorphism}, with invertible action constraints.\defmark
\end{description}
\end{explanation}

%% file: chap/shimakawa_K.tex
This chapter constructs and compares Shimakawa $K$-theory \cref{ksho_kshosg} 
\begin{equation}\label{Ksho_chi}
\begin{tikzpicture}[vcenter]
\def\h{2} \def\u{.6}
\draw[0cell=.9]
(0,0) node (a1) {\phantom{\AlglaxO}}
(a1)++(0,-.04) node (a1') {\AlglaxO}
(a1)++(\h,0) node (a2) {\FGCatg}
(a2)++(\h,0) node (a3) {\FGTopg}
(a3)++(\h,0) node (a4) {\Gspec} 
;
\draw[1cell=.85]
(a1) edge node {\Sgo} (a2)
(a2) edge node {\clast} (a3)
(a3) edge node {\Kfg} (a4)
(a1) [rounded corners=2pt, shorten <=-.2ex] |- ($(a2)+(0,\u)$) -- node {\Ksho} ($(a3)+(0,\u)$) -| (a4)
;
\end{tikzpicture}
\end{equation}
and its homotopical variant \cref{khsho_khshosg}
\begin{equation}\label{Khsho_chi}
\begin{tikzpicture}[vcenter]
\def\h{2} \def\u{.6}
\draw[0cell=.9]
(0,0) node (a1) {\phantom{\AlglaxO}}
(a1)++(0,-.04) node (a1') {\AlglaxO}
(a1)++(\h,0) node (a2) {\FGCatg}
(a2)++(\h,0) node (a3) {\FGTopg}
(a3)++(\h,0) node (a4) {\FGTopg}
(a4)++(\h,0) node (a5) {\Gspec} 
;
\draw[1cell=.85]
(a1) edge node {\Sgo} (a2)
(a2) edge node {\clast} (a3)
(a3) edge node {\Bc} (a4)
(a4) edge node {\Kfg} (a5)
(a1) [rounded corners=2pt, shorten <=-.2ex] |- ($(a2)+(0,\u)$) -- node {\Khsho} ($(a4)+(0,\u)$) -| (a5)
;
\end{tikzpicture}
\end{equation}
for a 1-connected $\Gcat$-operad $\Op$ and a compact Lie group $G$.  At the object level, each of the two functors $\Ksho$ and $\Khsho$ sends $\Op$-pseudoalgebras \pcref{def:pseudoalgebra} to orthogonal $G$-spectra \pcref{def:gsp_module}.  
\begin{itemize}
\item The functor $\Sgo$ is Shimakawa $H$-theory \pcref{def:sgo}, which sends $\Op$-pseudoalgebras to $\FGG$-categories \pcref{sys_FGcat}.
\item The functor $\clast$ \cref{clast_ggcatg} is induced by the classifying space functor $\cla \cn \Cat \to \Top$ \cref{classifying_space}.  Its sends $\FGG$-categories to $\FGG$-spaces.
\item The prolongation functor $\Kfg$ extends each $\FGG$-space to an orthogonal $G$-spectrum \cref{Kfgx}. 
\item In the homotopical Shimakawa $K$-theory $\Khsho$, $\Bc$ denotes the bar functor $\BcFG$ \cref{bar_functor_FG}. 
\end{itemize} 
The strong variants of $\Ksho$ and $\Khsho$, denoted by $\Kshosg$ and $\Khshosg$, are obtained by replacing Shimakawa $H$-theory $\Sgo$ with its strong variant $\Sgosg$ \cref{def:sgosg}:
\begin{equation}\label{KshosgKhshosg_chi}
\begin{tikzpicture}[vcenter]
\draw[0cell]
(0,0) node (a1) {\AlgpspsO}
(a1)++(4.3,0) node (a2) {\Gspec.}
;
\draw[1cell=.9]
(a1) edge[transform canvas={yshift=.5ex}] node {\Kshosg = \Kfg \clast \Sgosg} (a2)
(a1) edge[transform canvas={yshift=-.4ex}] node[swap] {\Khshosg = \Kfg \Bc \clast \Sgosg} (a2)
;
\end{tikzpicture}
\end{equation}
\cref{ch:shim_top,part:kgo_shi_comp} compare Shimakawa (strong) $K$-theory with our equivariant $K$-theory \cref{Khgo_functors}. 

\subsection*{Comparison of Shimakawa $K$-Theories}

The main observation of this chapter is \cref{thm:ksho_khsho}.  It proves that, for a finite group $G$, the functors $\Ksho$ and $\Khsho$ are naturally weakly $G$-equivalent, and the functors $\Kshosg$ and $\Khshosg$ are also naturally weakly $G$-equivalent.  Thus, starting from $\Op$-pseudoalgebras, Shimakawa $K$-theory $\Ksho$ and its homotopical variant $\Khsho$ yield orthogonal $G$-spectra that are naturally levelwise weakly $G$-equivalent.  The same statement also holds for the strong variants $\Kshosg$ and $\Khshosg$.

\subsection*{Shimakawa's Construction}

Shimakawa's original equivariant $K$-theory machine, constructed in \cite{shimakawa89}, is naturally weakly $G$-equivalent to the homotopical Shimakawa strong $K$-theory 
\begin{equation}\label{Khshbesg_chi}
\begin{tikzpicture}[vcenter]
\def\h{2} \def\u{.6}
\draw[0cell=.9]
(0,0) node (a1) {\AlgpspsBE}
(a1)++(\h,0) node (a2) {\FGCatg}
(a2)++(\h,0) node (a3) {\FGTopg}
(a3)++(\h,0) node (a4) {\FGTopg}
(a4)++(\h,0) node (a5) {\Gspec} 
;
\draw[1cell=.85]
(a1) edge node {\Sgosgbe} (a2)
(a2) edge node {\clast} (a3)
(a3) edge node {\Bc} (a4)
(a4) edge node {\Kfg} (a5)
(a1) [rounded corners=2pt] |- ($(a2)+(0,\u)$) -- node {\Khshbesg} ($(a4)+(0,\u)$) -| (a5)
;
\end{tikzpicture}
\end{equation}
for the Barratt-Eccles operad $\BE$ \pcref{def:BE}.  By \cref{thm:ksho_khsho}, it is also naturally weakly $G$-equivalent to Shimakawa strong $K$-theory 
\begin{equation}\label{Kshbesg_chi}
\begin{tikzpicture}[vcenter]
\def\h{2} \def\u{.6}
\draw[0cell=.9]
(0,0) node (a1) {\AlgpspsBE}
(a1)++(\h,0) node (a2) {\FGCatg}
(a2)++(\h,0) node (a3) {\FGTopg}
(a3)++(\h,0) node (a4) {\Gspec} 
;
\draw[1cell=.85]
(a1) edge node {\Sgosgbe} (a2)
(a2) edge node {\clast} (a3)
(a3) edge node {\Kfg} (a4)
(a1) [rounded corners=2pt] |- ($(a2)+(0,\u)$) -- node {\Kshbesg} ($(a3)+(0,\u)$) -| (a4)
;
\end{tikzpicture}
\end{equation}
for $\BE$.  See \cref{expl:shi89_section2,expl:khsho_shi,expl:k_shimakawa}.

\organization
This chapter consists of the following sections.

\secname{sec:fgtopg}
This section defines the category $\FGTopg$ of $\FGG$-spaces and constructs the functor $\clast$ from $\FGCatg$ to $\FGTopg$.

\secname{sec:c_shim}
This section constructs Shimakawa $K$-theory $\Ksho$ and its strong variant $\Kshosg$ by constructing the functor $\Kfg$ that sends $\FGG$-spaces to orthogonal $G$-spectra.

\secname{sec:bar_const}
To prepare for \cref{sec:h_shim}, this section recalls the bar functor in the general setting of enriched categories.  There is a natural transformation $\retn$ \cref{retn_bar_id}, called the retraction, that relates the bar functor and the identity functor.

\secname{sec:h_shim}
With the bar functor in hand, this section defines the homotopical Shimakawa $K$-theory $\Khsho$ and its strong variant $\Khshosg$.

\secname{sec:invariance_thm}
This section proves that, for a finite group $G$, prolonging $\FGG$-spaces to pointed finite $G$-CW complexes preserves componentwise weak $G$-equivalences between \emph{proper} $\FGG$-spaces \pcref{thm:invariance}.  Properness is defined in terms of Reedy cofibrant simplicial $G$-spaces.  It ensures that prolongation is homotopically well-behaved.  \cref{thm:invariance} is a variant of the Invariance Theorem from \cite{gmmo19}.  See \cref{rk:invariance}.

\secname{sec:ch_shim}
Using \cref{thm:invariance}, this section proves that, for a finite group $G$, the retraction $\retn$ induces a natural weak $G$-equivalence between Shimakawa $K$-theory $\Ksho$ and its homotopical variant $\Khsho$.  Replacing Shimakawa $H$-theory $\Sgo$ with the strong variant $\Sgosg$, the same proof shows that the retraction $\retn$ induces a natural weak $G$-equivalence between Shimakawa strong $K$-theory $\Kshosg$ and its homotopical variant $\Khshosg$.

\section{$\FGG$-Spaces from $\FGG$-Categories}
\label{sec:fgtopg}

This section constructs the passage from $\FGG$-categories to $\FGG$-spaces.  Throughout this section, $G$ denotes an arbitrary group.

\secoutline
\begin{itemize}
\item \cref{def:ggtopg} defines the categories $\FGTopg$ of $\FGG$-spaces, with further elaboration given in \cref{expl:ggtopg}.
\item \cref{ggcatg_ggtopg} constructs the functor from $\FGCatg$ to $\FGTopg$ induced by the classifying space functor $\cla$, with further elaboration given in \cref{expl:ggcatg_ggtopg}.
\end{itemize}

\subsection*{$\FGG$-Spaces}
Recall the pointed $G$-category $\Topgst$ \cref{topgst_gtopst_enr} of pointed $G$-spaces and pointed morphisms with the conjugation $G$-action \cref{ginv_h_g}.  Recall the pointed $G$-category $\FG$ of pointed finite $G$-sets and pointed morphisms with the conjugation $G$-action \pcref{def:FG}.  \cref{def:ggtopg} defines the topological analogue of the 2-category $\FGCatg$ of $\FGG$-categories \pcref{def:fgcatg}.

\begin{definition}\label{def:ggtopg}
Suppose $G$ is a group.  The category $\FGTopg$ has pointed $G$-functors
\begin{equation}\label{ggtopg_obj}
(\FG,\ord{0}) \fto{X} (\Topgst,*)
\end{equation}
as objects, called \index{FGG-space@$\FGG$-space}\emph{$\FGG$-spaces}, and $G$-natural transformations
\begin{equation}\label{ggtopg_icell}
\begin{tikzpicture}[vcenter]
\def\t{28}
\draw[0cell]
(0,0) node (a1) {\phantom{\Gsk}}
(a1)++(1.8,0) node (a2) {\phantom{\Gsk}}
(a1)++(-.08,0) node (a1') {\FG}
(a2)++(.23,0) node (a2') {\Topgst}
;
\draw[1cell=.9]
(a1) edge[bend left=\t] node {X} (a2)
(a1) edge[bend right=\t] node[swap] {X'} (a2)
;
\draw[2cell]
node[between=a1 and a2 at .45, rotate=-90, 2label={above,\theta}] {\Rightarrow}
;
\end{tikzpicture}
\end{equation}
as morphisms.  Identities and composition are those of natural transformations.
\end{definition}

\begin{explanation}[Unpacking $\FGTopg$]\label{expl:ggtopg}
The category $\FGTopg$ in \cref{def:ggtopg} is given explicitly as follows.
\begin{description}
\item[Objects]
An $\FGG$-space $X \cn \FG \to \Topgst$ \cref{ggtopg_obj} consists of the following data.
\begin{description}
\item[Pointed $G$-spaces] $X$ sends each pointed finite $G$-set $\mal \in \FG$ \cref{ordn_be} to a pointed $G$-space $X\mal$ such that $X\ord{0} = *$.  Its $G$-fixed basepoint is given by the pointed morphism
\[X(\ord{0} \to \mal) \cn X\ord{0} = * \to X\mal\]
for the $G$-fixed unique morphism $\ord{0} \to \mal$ in $\FG$.  
\item[Pointed morphisms] $X$ sends each pointed function $\psi \cn \mal \to \nbeta$ in $\FG$ to a pointed morphism
\begin{equation}\label{X_upom_ggtopg}
X\mal \fto{X\psi} X\nbeta
\end{equation}
such that $X$ preserves identity morphisms and composition.  The morphism $X\psi$ is \emph{not} generally $G$-equivariant.
\item[Equivariance] The $G$-equivariance of $X$ means the equality of morphisms
\begin{equation}\label{gggspace_equiv}
X(g \cdot \psi) = g (X\psi) \ginv
\end{equation}
for each $g \in G$ and each morphism $\psi$ in $\FG$.  The morphism $g \cdot \psi$ is the conjugation $g$-action on $\psi$ \cref{gpsi}, and $g (X\psi) \ginv$ is the conjugation $g$-action \cref{ginv_h_g} on $X\psi$.  Thus, if the morphism $\psi$ is $G$-fixed, then $X\psi$ is a pointed $G$-morphism.
\end{description}
\item[Morphisms] 
A morphism $\tha \cn X \to X'$ in $\FGTopg$ \cref{ggtopg_icell} consists of, for each object $\mal \in \FG$, an $\mal$-component pointed morphism
\begin{equation}\label{ggtopg_icell_comp}
X\mal \fto{\tha_{\mal}} X'\mal
\end{equation}
such that, for each morphism $\psi \cn \mal \to \nbeta$ in $\FG$, the following naturality diagram of pointed morphisms commutes.
\begin{equation}\label{ggtopg_icell_nat}
\begin{tikzpicture}[vcenter]
\def\v{-1.4}
\draw[0cell]
(0,0) node (a11) {X\mal}
(a11)++(2.5,0) node (a12) {X'\mal}
(a11)++(0,\v) node (a21) {X\nbeta}
(a12)++(0,\v) node (a22) {X'\nbeta}
;
\draw[1cell=.9]
(a11) edge node {\theta_{\mal}} (a12)
(a12) edge node {X'\psi} (a22)
(a11) edge node[swap] {X\psi} (a21)
(a21) edge node {\theta_{\nbeta}} (a22)
;
\end{tikzpicture}
\end{equation}
The $G$-equivariance of $\tha$ means the equality of morphisms
\begin{equation}\label{ggtopg_icell_geq}
\tha_{\mal} = g \theta_{\mal} \ginv
\end{equation}
for each $g \in G$ and pointed finite $G$-set $\mal \in \FG$.   In other words, each component of a morphism in $\FGTopg$ is a pointed $G$-morphism.  A morphism $\tha$ is automatically pointed, meaning $\tha_{\ord{0}} = 1_{*}$.  Identity morphisms and composition are defined componentwise using the components in \cref{ggtopg_icell_comp}.\defmark
\end{description}
\end{explanation}

\subsection*{$\FGG$-Spaces from $\FGG$-Categories}

\cref{ggcatg_ggtopg} is the $\FG$-analogue of \cref{GGCatg_GGTopg}.  The proof is the same after replacing $\GG$ by $\FG$.

\begin{lemma}\label{ggcatg_ggtopg}
For a group $G$, composing and whiskering with the classifying space functor $\cla$ induce a functor\index{classifying space}
\begin{equation}\label{clast_ggcatg}
\FGCatg \fto{\clast} \FGTopg
\end{equation}
between the categories in \cref{def:fgcatg,def:ggtopg}.
\end{lemma}

\begin{explanation}[Unpacking $\clast$]\label{expl:ggcatg_ggtopg}
The functor $\clast$ in \cref{clast_ggcatg} sends a pointed $G$-functor $X \cn \FG \to \Catgst$ and a $G$-natural transformation $\tha$ to the following composite pointed $G$-functor and whiskered $G$-natural transformation.
\begin{equation}\label{clast_ggcatg_mor}
\begin{tikzpicture}[baseline={(a1.base)}]
\def\t{28}
\draw[0cell]
(0,0) node (a1) {\FG}
(a1)++(1.8,0) node (a2) {\phantom{\GG}}
(a2)++(.15,0) node (a2') {\Catgst}
(a2')++(2,0) node (a3) {\Topgst}
;
\draw[1cell=.9]
(a1) edge[bend left=\t] node {X} (a2)
(a1) edge[bend right=\t] node[swap] {X'} (a2)
(a2') edge node {\cla} (a3)
;
\draw[2cell]
node[between=a1 and a2 at .43, rotate=-90, 2label={above,\tha}] {\Rightarrow}
;
\end{tikzpicture}
\end{equation}
For each pointed finite $G$-set $\mal \in \FG$, the $\mal$-component of $\clast\tha$ is the pointed $G$-morphism between pointed $G$-spaces
\[\cla X\mal \fto{\cla\tha_{\mal}} \cla X'\mal\]
obtained from $\tha_{\mal}$ \cref{ggtopg_icell_comp} by applying $\cla$.
\end{explanation}

\section{Shimakawa $K$-Theory}
\label{sec:c_shim}

This section constructs Shimakawa (strong) $K$-theory functors
\[\AlglaxO \fto{\Ksho} \Gspec \andspace \AlgpspsO \fto{\Kshosg} \Gspec\]
from $\Op$-pseudoalgebras to orthogonal $G$-spectra for a compact Lie group $G$ and a 1-connected $\Gcat$-operad $\Op$.

\secoutline
\begin{itemize}
\item \cref{def:Kfg_obj} constructs the orthogonal $G$-spectrum $\Kfg X$ associated to an $\FGG$-space $X$.
\item \cref{def:Kfg_mor} constructs the $G$-morphism $\Kfg\tha$ between orthogonal $G$-spectra associated to a $G$-natural transformation $\tha$ between $\FGG$-spaces.
\item \cref{def:Kfg_functor} defines the functor $\Kfg$ from the category $\FGTopg$ of $\FGG$-spaces to the category $\Gspec$ of orthogonal $G$-spectra.
\item \cref{def:ksho} defines Shimakawa $K$-theory $\Ksho$ and its strong variant $\Kshosg$, with further elaboration given in \cref{expl:ksho,expl:c_shimakawa}.
\end{itemize}

\subsection*{Orthogonal $G$-Spectra from $\FGG$-Spaces: Objects}

\cref{def:Kfg_obj} defines the orthogonal $G$-spectrum associated to an $\FGG$-space.  It is obtained from \cref{def:gggspace_gspectra} by replacing the indexing $G$-category $\GG$ with $\FG$ \pcref{def:FG}.

\begin{definition}[$\Kfg$ on Objects]\label{def:Kfg_obj}
Given a compact Lie group $G$ and an $\FGG$-space \pcref{def:ggtopg}
\[(\FG,\ord{0}) \fto{X} (\Topgst,*),\]
the orthogonal $G$-spectrum \pcref{def:gsp_module}
\begin{equation}\label{Kfgx}
(\Kfg X, \umu) \in \Gspec
\end{equation}
is defined as follows.
\begin{description}
\item[Object assignment of $\IU$-space]
The $\IU$-space \pcref{def:iu_space} $\Kfg X$ sends each object $V \in \IU$ \pcref{def:IU_spaces} to the pointed $G$-space
\begin{equation}\label{Kfgxv}
(\Kfg X)_V = \int^{\mal \sins \FG} (S^V)^{\mal} \sma X\mal.
\end{equation}
\begin{description}
\item[Coend] 
The coend in \cref{Kfgxv} is taken in the category $\Topst$ of pointed spaces and pointed morphisms.  Each pointed finite $G$-set $\mal \in \FG$ \pcref{def:FG} is regarded as a discrete pointed $G$-space, and $S^V$ is the $V$-sphere \pcref{def:g_sphere}.  The pointed $G$-space 
\begin{equation}\label{SVmal}
(S^V)^{\mal} = \Topgst(\mal, S^V)
\end{equation}
consists of pointed morphisms $\mal \to S^V$ \cref{Gtopst_smc}, with $G$ acting by conjugation \cref{ginv_h_g}.
\item[$G$-action] 
The group $G$ acts diagonally on representatives.  This means that, for an element $g \in G$ and a representative pair
\begin{equation}\label{Kfgxv_rep}
\big(\mal \fto{\upom} S^V ; x \in X\mal \big) \in (S^V)^{\mal} \ttimes X\mal
\end{equation}
in $(\Kfg X)_V$, the diagonal $g$-action is given by
\begin{equation}\label{Kfgxv_rep_gact}
g \cdot (\upom; x) = (g\upom\ginv ; gx),
\end{equation}
where $g\upom\ginv$ means the composite pointed morphism
\[\mal \fto{\ginv} \mal \fto{\upom} S^V \fto{g} S^V.\]
\end{description}
\item[Morphism assignment of $\IU$-space]
For a linear isometric isomorphism $f \cn V \fiso W$ in $\IU$, the pointed homeomorphism \cref{iu_space_xf}
\begin{equation}\label{Kfgxf}
(\Kfg X)_V \fto[\iso]{(\Kfg X)_f} (\Kfg X)_W
\end{equation}
is induced by the pointed homeomorphisms 
\[(S^V)^{\mal} \fto[\iso]{f \circ -} (S^W)^{\mal} \forspace \mal \in \FG\]
that postcompose with the pointed homeomorphism $f \cn S^V \fiso S^W$.  In terms of representatives \cref{Kfgxv_rep}, it is given by
\begin{equation}\label{Kfgxf_rep}
(\Kfg X)_f(\upom; x) = (f\upom; x).
\end{equation}
\item[Sphere action]
For each pair of objects $(V,W) \in (\IUsk)^2$, the $(V,W)$-component pointed $G$-morphism \cref{gsp_action_vw} is defined by the following commutative diagram in $\Gtopst$.
\begin{equation}\label{Kfgx_action_vw}
\begin{tikzpicture}[vcenter]
\def\h{4.7} \def\u{-1} \def\v{-1.4}
\draw[0cell=.9]
(0,0) node (a11) {(\Kfg X)_V \sma S^W}
(a11)++(\h,0) node (a12) {(\Kfg X)_{V \oplus W}}
(a11)++(0,\u) node (a21) {\big( \txint^{\mal \in \FG} (S^V)^{\mal} \sma X\mal \big) \sma S^W}
(a12)++(0,\u) node (a22) {\txint^{\mal \in \FG} (S^{V \oplus W})^{\mal} \sma X\mal}
(a21)++(0,\v) node (a3) {\txint^{\mal \in \FG} \big( (S^V)^{\mal} \sma S^W \big) \sma X\mal}
;
\draw[1cell=.8]
(a11) edge node {\umu_{V,W}} (a12)
(a11) edge[equal] (a21)
(a12) edge[equal] (a22)
(a21) edge node[swap] {\iso} (a3)
(a3) [rounded corners=2pt] -| node[pos=.25] {\asm = \txint\! \asm_{\mal} \sma 1} (a22)
;
\end{tikzpicture}
\end{equation}
\begin{itemize}
\item The pointed $G$-homeomorphism denoted by $\iso$ first commutes $- \sma S^W$ with the coend.  Then it moves $S^W$ to the left of $X\mal$ using the associativity isomorphism and braiding for the symmetric monoidal category $(\Gtopst,\sma)$ \cref{Gtopst_smc}.
\item The pointed $G$-morphism $\asm$ is induced by the pointed $G$-morphisms
\begin{equation}\label{assem_sph}
(S^V)^{\mal} \sma S^W \fto{\asm_{\mal}} (S^{V \oplus W})^{\mal}
\end{equation}
for $\mal \in \FG$ defined by the assignment
\[\begin{split}
& \big(\mal \fto{\upom} S^V ; y \big) \in (S^V)^{\mal} \sma S^W\\
&\mapsto \big(\mal \fto{\upom} S^V \fto{- \oplus y} S^{V \oplus W} \big) \in (S^{V \oplus W})^{\mal}.
\end{split}\]
In other words, $\asm_{\mal}(\upom; y) = \upom \oplus y$ sends an element $i \in \mal$ to the point 
\begin{equation}\label{asm_explicit}
\big( \asm_{\mal}(\upom; y)\big)(i) = (\upom i) \oplus y \in S^{V \oplus W}.
\end{equation}
For a representative pair $(\upom; x)$ of $(\Kfg X)_V$ \cref{Kfgxv_rep} and a point $y \in S^W$, $\umu_{V,W}$ is given by
\begin{equation}\label{Kfgx_actrep}
\umu_{V,W}\big((\upom; x); y\big) = (\upom \oplus y; x).
\end{equation}
\end{itemize}
\end{description}
This finishes the definition of the orthogonal $G$-spectrum $(\Kfg X,\umu)$.
\end{definition}

\begin{lemma}\label{kfgx_welldef}
The pair $(\Kfg X,\umu)$ in \cref{Kfgx} is an orthogonal $G$-spectrum.
\end{lemma}

\begin{proof}
We reuse the proof of \cref{kggx_welldef} by replacing $\GG$ with $\FG$.
\end{proof}

\subsection*{Orthogonal $G$-Spectra from $\FGG$-Spaces: Morphisms}

Next, we define the morphism assignment of $\Kfg$.  A morphism between $\FGG$-spaces is a $G$-natural transformation \pcref{def:ggtopg}.  A $G$-morphism between orthogonal $G$-spectra is a $G$-equivariant $\IU$-morphism that is compatible with the sphere actions \pcref{def:iu_morphism,def:gsp_morphism}.  \cref{def:Kfg_mor} is obtained from \cref{def:Kgg_mor} by replacing the indexing $G$-category $\GG$ with $\FG$.

\begin{definition}[$\Kfg$ on Morphisms]\label{def:Kfg_mor}
For a compact Lie group $G$ and a $G$-natural transformation 
\begin{equation}\label{mortheta_iicell}
\begin{tikzpicture}[vcenter]
\def\t{28}
\draw[0cell]
(0,0) node (a1) {\phantom{\Gsk}}
(a1)++(1.8,0) node (a2) {\phantom{\Gsk}}
(a1)++(-.08,0) node (a1') {\FG}
(a2)++(.23,0) node (a2') {\Topgst}
;
\draw[1cell=.9]
(a1) edge[bend left=\t] node {X} (a2)
(a1) edge[bend right=\t] node[swap] {Y} (a2)
;
\draw[2cell]
node[between=a1 and a2 at .42, rotate=-90, 2label={above,\tha}] {\Rightarrow}
;
\end{tikzpicture}
\end{equation}
between $\FGG$-spaces $X$ and $Y$, the $G$-morphism between orthogonal $G$-spectra
\begin{equation}\label{Kfg_psi}
(\Kfg X, \umu) \fto{\Kfg\tha} (\Kfg Y, \umu)
\end{equation}
has, for each object $V \in \IU$, $V$-component pointed $G$-morphism \cref{eqiu_mor_component} defined by the commutative diagram
\begin{equation}\label{Kfg_psiv}
\begin{tikzpicture}[vcenter]
\def\v{-1.4}
\draw[0cell=.9]
(0,0) node (a11) {(\Kfg X)_V}
(a11)++(3,0) node (a12) {\txint^{\mal \sins \FG} (S^V)^{\mal} \sma X\mal}
(a11)++(0,\v) node (a21) {(\Kfg Y)_V}
(a12)++(0,\v) node (a22) {\txint^{\mal \sins \FG} (S^V)^{\mal} \sma Y\mal}
;
\draw[1cell=.8]
(a11) edge[equal] (a12)
(a21) edge[equal] (a22)
(a11) edge[transform canvas={xshift=1em}] node[swap] {(\Kfg\tha)_V} (a21)
(a12) edge[transform canvas={xshift=-2em}, shorten <=-.5ex] node {\txint^{\mal} 1 \sma \tha_{\mal}} (a22)
;
\end{tikzpicture}
\end{equation}
in $\Gtopst$.  More explicitly, $(\Kfg\tha)_V$ sends a representative pair $(\upom; x)$ of $(\Kfg X)_V$ \cref{Kfgxv_rep} to the representative pair
\begin{equation}\label{Kfg_psiv_rep}
(\Kfg\tha)_V (\upom; x) = \big(\upom; \tha_{\mal} x\big)
\end{equation}
of $(\Kfg Y)_V$. 
\end{definition}

\begin{lemma}\label{kfgpsi_welldef}
The assignment $\Kfg\tha$ \cref{Kfg_psi} is a $G$-morphism between orthogonal $G$-spectra.
\end{lemma}

\begin{proof}
We reuse the proof of \cref{kggpsi_welldef} by replacing $\GG$ with $\FG$.
\end{proof}

\begin{definition}\label{def:Kfg_functor}
For a compact Lie group $G$, the functor
\[\FGTopg \fto{\Kfg} \Gspec\]
is defined by
\begin{itemize}
\item the object assignment $X \mapsto (\Kfg X,\umu)$ \pcref{def:Kfg_obj} and
\item the morphism assignment $\tha \mapsto \Kfg\tha$ \pcref{def:Kfg_mor}.  
\end{itemize}
Its functoriality follows from \cref{Kfg_psiv_rep} together with the fact that identities and composition are defined componentwise in $\FGTopg$ and $\Gspec$ \pcref{def:gsp_morphism,def:ggtopg}.
\end{definition}

\begin{explanation}[Conceptual Machine]\label{expl:Kfg_functor}
In \cite[Section 2]{mmo}, the functor $\Kfg$ \pcref{def:Kfg_functor} is presented as a composite
\begin{equation}\label{PwgUgs}
\FGTopg \fto{\Pwg} \WGTopg \fto{\Ugs} \Gspec
\end{equation}
and called the \index{conceptual Segal machine}\emph{conceptual Segal machine} on $\FGG$-spaces.
\begin{itemize}
\item $\WG$ is the pointed $G$-category of pointed finite $G$-CW complexes and pointed morphisms, with $G$ acting by conjugation \cref{ginv_h_g}.  The category $\WGTopg$ is defined like $\FGTopg$ \pcref{def:ggtopg} with the pointed $G$-category $(\FG,\ord{0})$ replaced by $(\WG,*)$.  Thus, an object in $\WGTopg$, called a \emph{$\WGG$-space}, is a pointed $G$-functor 
\[(\WG,*) \to (\Topgst,*),\] 
and a morphism in $\WGTopg$ is a $G$-natural transformation.  Replacing $(\FG,\ord{0})$ by $(\WG,*)$, the description of $\FGTopg$ in \cref{expl:ggtopg} also applies to $\WGTopg$.
\item The prolongation functor $\Pwg$ is defined objectwise as a coend like \cref{Kfgxv}.  Thus, for an $\FGG$-space $X \cn \FG \to \Topgst$, the $\WGG$-space $\Pwg X$ sends a pointed finite $G$-CW complex $A \in \WG$ to the pointed $G$-space
\begin{equation}\label{Pwgxa}
(\Pwg X)_A = \int^{\mal \sins \FG} A^{\mal} \sma X\mal.
\end{equation}
The pointed $G$-space 
\begin{equation}\label{A_to_mal}
A^{\mal} = \Topgst(\mal,A)
\end{equation}
consists of all pointed morphisms $\mal \to A$, with $G$ acting by conjugation.
\item The functor $\Ugs$ restricts to the $V$-spheres $S^V$ for $V \in \IU$ \pcref{def:g_sphere} and the resulting orthogonal $G$-spectra \cref{Kfgx_action_vw}.\defmark
\end{itemize}
\end{explanation}

\subsection*{Shimakawa $K$-Theory}

Next, we define the first Shimakawa equivariant $K$-theory machine from $\Op$-pseudoalgebras to orthogonal $G$-spectra \pcref{def:pseudoalgebra,def:gsp_module}.  Recall that a $\Gcat$-operad $\Op$ is \emph{1-connected} if $\Op(0)$ and $\Op(1)$ are terminal $G$-categories \cref{i_connected}.

\begin{definition}\label{def:ksho}
For a compact Lie group $G$ and a 1-connected $\Gcat$-operad $\Op$, we define \emph{Shimakawa $K$-theory} $\Ksho$ and \emph{Shimakawa strong $K$-theory} $\Kshosg$ for $\Op$ as the following composite functors.\index{Shimakawa K-theory@Shimakawa $K$-theory}\index{Shimakawa K-theory@Shimakawa $K$-theory!strong}\index{K-theory@$K$-theory!Shimakawa}
\begin{equation}\label{ksho_kshosg}
\begin{tikzpicture}[vcenter]
\def\h{2.2} \def\u{.7} \def\v{-1.8}
\draw[0cell]
(0,0) node (a1) {\AlglaxO}
(a1)++(.12*\h,.1*\v) node (a1') {\phantom{A}}
(a1)++(0,\v) node (a2) {\AlgpspsO}
(a2)++(.12*\h,-.1*\v) node (a2') {\phantom{A}}
(a1)++(\h,\v/2) node (a3) {\phantom{\FGCatg}}
(a3)++(0,-.01) node (a3') {\FGCatg}
(a3)++(\h,0) node (a4) {\FGTopg}
(a4)++(\h,0) node (a5) {\Gspec}
;
\draw[1cell=.9]
(a1') edge[shorten >=-.4ex] node[swap, inner sep=1pt] {\Sgo} (a3)
(a2') edge[shorten >=-.4ex] node[inner sep=1pt] {\Sgosg} (a3)
(a3) edge node {\clast} (a4)
(a4) edge node {\Kfg} (a5)
(a1) [rounded corners=2pt] -| node[pos=.21] {\Ksho}  (a5)
;
\draw[1cell=.9]
(a2) [rounded corners=2pt] -| node[pos=.21] {\Kshosg}  (a5)
;
\end{tikzpicture}
\end{equation}
\begin{description}
\item[$\Op$-pseudoalgebras to $\FGG$-categories] 
$\Sgo$ and $\Sgosg$ are Shimakawa (strong) $H$-theory \pcref{Sgo_twofunctor} between the categories in \cref{oalgps_twocat,def:fgcatg}. 
\item[$\FGG$-categories to $\FGG$-spaces] 
$\clast$ is induced by the classifying space functor $\cla \cn \Cat \to \Top$ \pcref{ggcatg_ggtopg}.  The functors $\Sgo$, $\Sgosg$, and $\clast$ are defined for arbitrary groups $G$. 
\item[$\FGG$-spaces to orthogonal $G$-spectra] 
The functor $\Kfg$ \pcref{def:Kfg_functor}, going between the categories in \cref{def:gsp_morphism,def:ggtopg}, requires $G$ to be a compact Lie group.\defmark
\end{description}
\end{definition}

\begin{explanation}[Unpacking]\label{expl:ksho}
For an $\Op$-pseudoalgebra $\A$ and an object $V \in \IU$ \pcref{def:pseudoalgebra,def:indexing_gspace}, the pointed $G$-spaces $(\Ksho \A)_V$ and $(\Kshosg\A)_V$ are given by the coends \cref{Kfgxv}
\[\begin{split}
(\Ksho \A)_V &= \int^{\mal \sins \FG} (S^V)^{\mal} \sma \cla (\Amal) \andspace\\
(\Kshosg \A)_V &= \int^{\mal \sins \FG} (S^V)^{\mal} \sma \cla (\Asgmal)
\end{split}\]
with the diagonal $G$-action and the sphere action defined in \cref{Kfgxv_rep_gact} and \cref{Kfgx_action_vw}.  The pointed $G$-categories $\Amal$ and $\Asgmal$ are those of (strong) $\mal$-systems in $\A$ \pcref{def:nsys_gcat}.  See \cref{expl:khsho} for the homotopical variants.
\end{explanation}

\begin{explanation}[Shimakawa's Construction]\label{expl:c_shimakawa}
Unlike our $\Gtop$-multifunctors $\Kgo$ and $\Kgosg$ \cref{Kgo_multifunctor}, Shimakawa (strong) $K$-theory $\Ksho$ and $\Kshosg$ are \emph{not} multifunctors and do not generally preserve multiplicative structures.  For the Barratt-Eccles $\Gcat$-operad $\BE$, Shimakawa strong $H$-theory $\Sgosgbe$ is Shimakawa's original construction \pcref{expl:shi89_section2}.  On the other hand, the last step $\Kfg$ is not identical to Shimakawa's original construction, which uses a homotopical variant of the coend \cref{Kfgxv}.  We recall Shimakawa's homotopical construction in \cref{sec:h_shim}.  See \cref{expl:khsho_shi,expl:k_shimakawa} for further discussion of Shimakawa's original construction.  
\end{explanation}

\section{Bar Constructions}
\label{sec:bar_const}

To prepare for \cref{sec:h_shim}, this section recalls the bar construction used in the homotopical version of Shimakawa's construction from $\FGG$-spaces to orthogonal $G$-spectra, following \cite[Section 3.1]{mmo}.  Another reference for the bar construction is \cite[Section 12]{shulman_holim}.  A reference for simplicial objects is \cite[Ch.\ 7]{cerberusIII}.  Monoidal categories and enriched categories are reviewed in \cref{sec:sym_mon_cat,sec:enrichedcat}.

\secoutline
\begin{itemize}
\item \cref{def:simp_obj} recalls simplicial objects, simplicial morphisms, and homotopy.
\item \cref{def:cat_bar} recalls the bar construction, realization, and two auxiliary constructions called retractions and sections.
\item \cref{expl:retn,expl:secn} further elaborate retractions and sections. 
\item \cref{def:bar_functor} defines the morphism assignment of the bar construction and the resulting bar functor.
\item \cref{expl:Bcdot_tha} discusses further naturality properties of retractions and sections.
\end{itemize}  

\begin{definition}[Simplicial Objects]\label{def:simp_obj}\
\begin{enumerate}
\item\label{def:simp_obj_i} 
The category $\Delta$\label{not:Delta_cat}\index{Delta@$\Delta$}  has
\begin{itemize}
\item the totally ordered sets $\ordr = \{0 < 1 < \Cdots < r\}$ for $r \geq 0$ as objects and 
\item weakly order-preserving functions $f \cn \ordr \to \ordt$ as morphisms.  
\end{itemize}
Weakly order-preserving means $f(k) \leq f(\ell)$ for $k < \ell$.
\item\label{def:simp_obj_ii} 
Morphisms of $\Delta$ are generated by the \index{coface}\index{face!co-}\emph{coface} and \index{degeneracy!co-}\index{codegeneracy}\emph{codegeneracy} morphisms for $0 \leq i \leq r$:\label{not:face_degen} 
\begin{equation}\label{coface_codeg}
\ord{r-1} \fto{d^i} \ord{r} \fot{s^i} \ord{r+1}.
\end{equation}
\begin{itemize}
\item The coface $d^i$ is the order-preserving injection whose image does not contain $i \in \ord{r}$. 
\item The codegeneracy $s^i$ is the weakly order-preserving surjection that sends $i,i+1 \in \ord{r+1}$ to $i \in \ord{r}$.
\end{itemize}
The coface and codegeneracy morphisms satisfy the following \index{cosimplicial identities}\index{identities!cosimplicial}\index{simplicial identities!co-}\emph{cosimplicial identities}.
\[\begin{aligned}
d^jd^i & = d^i d^{j-1} \qquad \!\text{if $i < j$}\\
s^j d^i & = \begin{cases}
d^i s^{j-1} & \text{if $i<j$}\\ 
1 & \text{if $i=j$ or $i = j+1$} \\
d^{i-1} s^j  & \text{if $i > j+1$} 
\end{cases} \\
s^j s^i & = s^i s^{j+1} \qquad \text{if $i \leq j$}
\end{aligned}\]
The coface morphisms, codegeneracy morphisms, and cosimplicial identities form a generating set of morphisms and relations for $\Delta$.
\item\label{def:simp_obj_iii} 
A \emph{simplicial object}\index{simplicial object} in a category $\V$ is a functor $X_\crdot \cn \Deltaop \to \V$, where $\Deltaop$ is the opposite category of $\Delta$.  A simplicial object $X_\crdot$ is completely determined by
\begin{itemize}
\item the \index{simplex}\emph{$r$-simplex objects} $X_r = X_\crdot \ord{r}$ for $r \geq 0$, 
\item the \index{face}\emph{$i$-th face morphisms} $d_i = X_\crdot d^i \cn X_r \to X_{r-1}$, and
\item the \index{degeneracy}\emph{$i$-th degeneracy morphisms} $s_i = X_\crdot s^i \cn X_r \to X_{r+1}$
\end{itemize}
that satisfy the following \index{simplicial identities}\emph{simplicial identities}.
\begin{equation}\label{simplicial_id}
\begin{split}
d_i d_j & = d_{j-1} d_i \qquad \!\text{if $i < j$}\\
d_i s_j & = \begin{cases}
s_{j-1} d_i & \text{if $i<j$}\\ 
1 & \text{if $i=j$ or $i = j+1$} \\
s_j d_{i-1}  & \text{if $i > j+1$} 
\end{cases} \\
s_i s_j & = s_{j+1} s_i \qquad \text{if $i \leq j$}
\end{split}
\end{equation}
\item\label{def:simp_obj_iv} 
A \emph{simplicial morphism}\index{simplicial morphism} between simplicial objects is a natural transformation.  In other words, a simplicial morphism $\fun \cn X_\crdot \to Y_\crdot$ consists of morphisms $\fun_r \cn X_r \to Y_r$ for $r \geq 0$ that commute with each face morphism and each degeneracy morphism.  The category of simplicial objects and simplicial morphisms in a category $\V$ is denoted by $\V^{\Deltaop}$.
\item\label{def:simp_obj_v} 
For two simplicial morphisms $\fun,\fun' \cn X_\crdot \to Y_\crdot$, a \emph{homotopy}\index{homotopy} $\htpy \cn \fun \simeq \fun'$ consists of morphisms
\[X_r \fto{\htpy_i} Y_{r+1} \forspace 0 \leq i \leq r\]
that satisfy the following equalities in $\V$.
\begin{equation}\label{htpy_id}
\begin{split}
d_0 \htpy_0 &= \fun_r \qquad d_{r+1} \htpy_r = \fun'_r\\
d_i \htpy_j &= \begin{cases}
\htpy_{j-1} d_i & \text{if $i<j$}\\
d_j \htpy_{j-1} & \text{if $i=j>0$}\\
\htpy_j d_{i-1} & \text{if $i>j+1$}
\end{cases}\\
s_i \htpy_j &= \begin{cases}
\htpy_{j+1} s_i & \text{if $i \leq j$}\\
\htpy_j s_{i-1} & \text{if $i>j$}
\end{cases}
\end{split}
\end{equation}
\item\label{def:simp_obj_vi}
The category of simplicial sets is denoted by $\sset$.  For a simplicial set $X_\crdot$, an element $x \in X_r$ is called an \emph{$r$-simplex}.  An $r$-simplex is \emph{degenerate}\index{simplex!degenerate} if it has the form $s_i x$ for some $x \in X_{r-1}$ and $0 \leq i \leq r-1$.  An $r$-simplex is \emph{nondegenerate}\index{simplex!nondegenerate} if it does not have the form $s_i x$.
\end{enumerate}
This finishes the definition.
\end{definition}

We use \cref{conv:left_normal} for iterated monoidal products, so $a \otimes b \otimes c$ means $(a \otimes b) \otimes c$.

\begin{definition}[Bar Construction]\label{def:cat_bar}
For a complete and cocomplete symmetric monoidal closed category  $(\V,\otimes,\tu)$, we denote by $\Vse$ the $\V$-category given by the objects and internal hom of $\V$.  We consider a triple $(\hun,\C,\fun)$ consisting of a small $\V$-category $(\C,\mcomp,\cone)$ and $\V$-functors
\[\C \fto{\fun} \Vse \andspace \C^{\op} \fto{\hun} \Vse.\]
\begin{description}
\item[Simplicial bar construction]
The \emph{simplicial bar construction}\index{simplicial bar construction}\index{bar construction!simplicial} is the simplicial object in $\V$ 
\begin{equation}\label{bcdot_hcf}
\Deltaop \fto{\Bcdot(\hun,\C,\fun)} \V
\end{equation} 
with $r$-simplex object given by the coproduct
\begin{equation}\label{bc_r}
\Bc_r(\hun,\C,\fun) = 
\coprod_{\ang{c_\ell}_{\ell=0}^r} \hun_{c_r} \otimes \C(c_{r-1}, c_r) \otimes \Cdots \otimes \C(c_0, c_1) \otimes \fun_{c_0}
\end{equation}
in $\V$ for each $r \geq 0$.  The coproduct \cref{bc_r} is indexed by the set of $(r+1)$-tuples $\ang{c_\ell}_{\ell=0}^r$ of objects in $\C$.
\begin{description}
\item[Faces] For $0 \leq i \leq r$, the $i$-th face morphism is induced by the composition of $\C$ and the evaluations of $\fun$ and $\hun$:
\[\begin{aligned}
\C(c_0, c_1) \otimes \fun_{c_0} & \fto{\ev} \fun_{c_1} && \text{if $i=0$,}\\
\C(c_i, c_{i+1}) \otimes \C(c_{i-1}, c_i) & \fto{\mcomp} \C(c_{i-1},c_{i+1}) && \text{if $0 < i < r$, and}\\
\hun_{c_r} \otimes \C(c_{r-1}, c_r) & \fto{\ev} \hun_{c_{r-1}} && \text{if $i=r$.}
\end{aligned}\]
\item[Degeneracies]
The $i$-th degeneracy morphism inserts the identity
\[\tu \fto{\cone_{c_i}} \C(c_i,c_i)\]
of $c_i$
\begin{itemize}
\item to the right of $\C(c_i,c_{i+1})$ if $0 \leq i < r$ and
\item to the left of $\C(c_{r-1},c_r)$ if $i=r$.
\end{itemize} 
\end{description}
The simplicial identities \cref{simplicial_id} for $\Bcdot(\hun,\C,\fun)$ follow from the associativity and unity of $\C$, the $\V$-functoriality of $\fun$ and $\hun$, and the coherence of $\V$.
\item[Realization]
Given a functor \label{not:simpdot}$\simpdot \cn \Delta \to \V$, the \emph{realization}\index{realization} of a simplicial object $X_\crdot \cn \Deltaop \to \V$ is the coend
\begin{equation}\label{x_dot}
|X_\crdot| = \int^{\ordr \in \Delta} X_{\ordr} \otimes \simp^{\ordr} \inspace \V.
\end{equation}
The \emph{realized bar construction}\index{realized bar construction}\index{bar construction!realized} is the realization of the simplicial bar construction $\Bcdot(\hun,\C,\fun)$:
\begin{equation}\label{bar_realization}
\Bc(\hun,\C,\fun) = \int^{\ordr \in \Delta} \Bc_r(\hun,\C,\fun) \otimes \simp^{\ordr}.
\end{equation}
\item[Bar construction]
For an object $c \in \C$, we consider the representable $\V$-functor 
\begin{equation}\label{representable_c}
\C_c = \C(-,c) \cn \C^{\op} \to \Vse.
\end{equation} 
The \index{bar construction}\emph{bar construction of $\fun$} is the $\V$-functor 
\begin{equation}\label{Bc_CCf}
\C \fto{\Bc(\C,\C,\fun)} \Vse
\end{equation}
that sends an object $c \in \C$ to the realized bar construction \cref{bar_realization} with $\hun = \C_c$:
\begin{equation}\label{real_Cc}
\Bc(\C_c,\C,\fun) = \int^{\ordr \in \Delta} \Bc_r(\C_c,\C,\fun) \otimes \simp^{\ordr}.
\end{equation}
\item[Retractions] 
For an object $c \in \C$, $(\fun_c)_\crdot \cn \Deltaop \to \V$\label{not:funcdot} denotes the constant simplicial object at $\fun_c \in \V$.  The \index{simplicial retraction}\index{retraction!simplicial}\emph{simplicial retraction at $c$} is the simplicial morphism
\begin{equation}\label{retn_c}
\Bcdot(\C_c,\C,\fun) \fto{\retn_c} (\fun_c)_\crdot
\end{equation}
defined at each simplicial level by the composition of $\C$ and the evaluation of $\fun$.  The simplicial retraction $\retn_c$ is $\V$-natural in $c$.  Via realization \cref{x_dot}, the simplicial retractions \cref{retn_c} yield a $\V$-natural transformation
\begin{equation}\label{retn}
\begin{tikzpicture}[vcenter]
\def\t{28}
\draw[0cell]
(0,0) node (a1) {\C}
(a1)++(1.8,0) node (a2) {\Vse}
;
\draw[1cell=.9]
(a1) edge[bend left=\t] node {\Bc(\C,\C,\fun)} (a2)
(a1) edge[bend right=\t] node[swap] {\fun} (a2)
;
\draw[2cell]
node[between=a1 and a2 at .45, rotate=-90, 2label={above,\retn}] {\Rightarrow}
;
\end{tikzpicture}
\end{equation}
called the \emph{retraction}\index{retraction} of $\fun$.  See \cref{expl:retn} for further details.
\item[Sections]
The \index{section}\emph{section at $c \in \C$} is the simplicial morphism
\begin{equation}\label{secn_c}
(\fun_c)_\crdot \fto{\secn_c} \Bcdot(\C_c,\C,\fun) 
\end{equation}
defined at each simplicial level $r$ by the tuple of objects $\ang{c}_{\ell=0}^r$ and copies of the identity $\cone_c \cn \tu \to \C(c,c)$.  There are an equality and a homotopy:
\begin{equation}\label{retn_secn}
\begin{split}
\retn_c \secn_c &= 1_{(\fun_c)_\crdot} \andspace\\
\htpy \cn \secn_c \retn_c &\simeq 1_{\Bcdot(\C_c,\C,\fun)}.
\end{split}
\end{equation}
We emphasize that $\secn_c$ is \emph{not} $\V$-natural in $c$.  See \cref{expl:secn} for further details.\defmark
\end{description}
\end{definition}

\begin{explanation}[Retractions]\label{expl:retn}
Consider the retraction 
\[\Bc(\C,\C,\fun) \fto{\retn} \fun\]
of the $\V$-functor $\fun \cn \C \to \Vse$ \cref{retn}.
\begin{description}
\item[$\V$-functoriality of the bar construction] 
For objects $c,d \in \C$, the $\V$-functor \cref{Bc_CCf} 
\[\C \fto{\Bc(\C,\C,\fun)} \Vse\]
has $(c,d)$-component morphism adjoint to the morphism
\begin{equation}\label{BCCf_functor}
\C(c,d) \otimes \Bc(\C_c,\C,\fun) \to \Bc(\C_d,\C,\fun) \inspace \V.
\end{equation}
Using the commutation of $\C(c,d) \Otimes -$ with coends \cref{real_Cc} and coproducts \cref{bc_r}, the preceding morphism is induced by the composition
\[\C(c,d) \otimes \C(c_r,c) \fto{\mcomp} \C(c_r,d)\]
of $\C$.  The $\V$-functoriality of $\Bc(\C,\C,\fun)$ follows from the associativity and unity of $\C$.
\item[Simplicial morphism] 
The simplicial retraction $\retn_c$ at $c \in \C$ \cref{retn_c} is defined at simplicial level $r$ by the following commutative diagrams in $\V$ for $c_0,\ldots,c_r \in \C$.  
\begin{equation}\label{retn_cr}
\begin{tikzpicture}[vcenter]
\def\v{-1.3}
\draw[0cell]
(0,0) node (a11) {\C(c_r,c) \otimes \C(c_{r-1},c_r) \otimes \Cdots \otimes \C(c_0,c_1) \otimes \fun_{c_0}}
(a11)++(2.1,0) node (a11') {\phantom{\C(c_0,c) \otimes \fun_{c_0}}}
(a11)++(5.1,0) node (a12) {\Bc_r(\C_c,\C,\fun)}
(a11')++(0,\v) node (a21) {\C(c_0,c) \otimes \fun_{c_0}}
(a12)++(0,\v) node (a22) {\fun_c = (\fun_c)_r}
;
\draw[1cell=.9]
(a11) edge node {\incn} (a12)
(a12) edge node {(\retn_c)_r} (a22)
(a11') edge node[swap] {\mcomp \otimes 1} (a21)
(a21) edge node {\ev} (a22)
;
\end{tikzpicture}
\end{equation}
The arrow $\incn$ is the inclusion into the coproduct \cref{bc_r}.  The iterated composition $\mcomp$ of $\C$ is the identity if $r=0$.  By the $\V$-functoriality of $\fun$, the composite $\ev(\mcomp \otimes 1)$ is equal to the $(r+1)$-fold iterated evaluations of $\fun$.  The fact that $\retn_c$ is compatible with the face and degeneracy morphisms follows from the $\V$-functoriality of $\fun$ and the associativity and unity of $\C$.
\item[$\V$-naturality] 
The $\V$-naturality of the simplicial retraction $\retn_c$ \cref{retn_c} in $c \in \C$ means the commutativity of the boundary of the following diagram in $\V$ for $c,d \in \C$.
\begin{equation}\label{Vnat_retraction}

\end{equation}
In the preceding diagram, the left region commutes by the associativity of $\C$, and the right region commutes by the $\V$-functoriality of $\fun$.\defmark
\end{description}
\end{explanation}

\begin{explanation}[Sections]\label{expl:secn}
The section 
\[(\fun_c)_\crdot \fto{\secn_c} \Bcdot(\C_c,\C,\fun)\]
at $c \in \C$ \cref{secn_c} is defined at simplicial level $r$ as the following composite in $\V$.
\begin{equation}\label{secn_cr}
%
\end{equation}
The fact that $\secn_c$ is compatible with the face and degeneracy morphisms follows from the unity of $\fun$ and $\C$.
\begin{description}
\item[Nonnaturality] 
Unlike the retraction $\retn$ \cref{retn}, the section $\secn_c$ is \emph{not} $\V$-natural in $c$.  In fact, $\secn_c$ is $\V$-natural in $c$ if and only if the following diagram in $\V$ commutes for $c,d \in \C$.  
\begin{equation}\label{secn_nonnat}
%
\end{equation}
The preceding diagram does not commute in general because the domains of the two inclusions $\iota$ are generally different summands of the coproduct $\Bc_r(\C_d,\C,\fun)$ \cref{bc_r}.
\item[Homotopy] 
For each $c \in \C$, the homotopy \cref{retn_secn} 
\[\htpy \cn \secn_c \retn_c \simeq 1_{\Bcdot(\C_c,\C,\fun)}\]
exists by the \index{extra degeneracy argument}\emph{extra degeneracy argument}.  More precisely, for $0 \leq i \leq r$, the morphism
\[\Bc_r(\C_c,\C,\fun) \fto{\htpy_i} \Bc_{r+1}(\C_c,\C,\fun)\]
is defined by the following commutative diagrams in $\V$ for $c_0, \ldots, c_r \in \C$.  
\begin{equation}\label{extradegen}
\begin{tikzpicture}[vcenter]
\def\v{-1.4} \def\h{4.2} 
\draw[0cell=.8]
(0,0) node (a1) {\C(c_r,c) \otimes [\txotimes_{\ell=1}^{r-i} \C(c_{r-\ell}, c_{r-\ell+1})] \otimes [\txotimes_{\ell=1}^i \C(c_{i-\ell}, c_{i-\ell+1})] \otimes \fun_{c_0}}
(a1)++(0,\v) node (a2) {\C(c_i,c) \otimes [\txotimes_{\ell=1}^i \C(c_{i-\ell}, c_{i-\ell+1})] \otimes \fun_{c_0}}
(a2)++(0,\v) node (a3) {\tu^{\otimes r-i+1} \otimes \C(c_i,c) \otimes [\txotimes_{\ell=1}^i \C(c_{i-\ell}, c_{i-\ell+1})] \otimes \fun_{c_0}}
(a3)++(0,\v) node (a4) {\C(c,c)^{\otimes r-i+1} \otimes \C(c_i,c) \otimes [\txotimes_{\ell=1}^i \C(c_{i-\ell}, c_{i-\ell+1})] \otimes \fun_{c_0}}
(a1)++(\h,.7*\v) node (b1) {\Bc_r(\C_c,\C,\fun)}
(a4)++(\h,-.7*\v) node (b2) {\Bc_{r+1}(\C_c,\C,\fun)}
;
\draw[1cell=.8]
(a1) edge node [swap] {\mcomp \otimes 1^{\otimes i+1}} (a2)
(a2) edge node [swap] {\iso} (a3)
(a3) edge node [swap] {\cone_c^{\otimes r-i+1} \otimes 1^{\otimes i+2}} (a4)
(b1) edge node {\htpy_i} (b2)
(a1) [rounded corners=2pt] -| node[pos=.7] {\iota} (b1)
;
\draw[1cell=.8]
(a4) [rounded corners=2pt] -| node[swap,pos=.7] {\iota} (b2)
;
\end{tikzpicture}
\end{equation}
In the preceding diagram, the composition $\mcomp$ means the identity if $i=r$.  The axioms \cref{htpy_id} for a homotopy hold by \cref{retn_cr}, \cref{secn_cr}, the $\V$-functoriality of $\fun$, and the associativity and unity of $\C$.\defmark
\end{description}
\end{explanation}

Next, we discuss functoriality of the bar construction.

\begin{definition}[Bar Functor]\label{def:bar_functor}
In the context of \cref{def:cat_bar}, suppose 
\begin{equation}\label{Vnt_theta}
\begin{tikzpicture}[vcenter]
\def\t{29}
\draw[0cell]
(0,0) node (a1) {\C}
(a1)++(1.7,0) node (a2) {\Vse}
;
\draw[1cell=.9]
(a1) edge[bend left=\t] node {\fun} (a2)
(a1) edge[bend right=\t] node[swap] {\fun'} (a2)
;
\draw[2cell]
node[between=a1 and a2 at .42, rotate=-90, 2label={above,\tha}] {\Rightarrow}
;
\end{tikzpicture}
\end{equation}
is a $\V$-natural transformation between $\V$-functors $\fun$ and $\fun'$.  For each object $c \in \C$, the $c$-component of $\tha$ is a morphism
\[\tu \fto{\tha_c} \Vse(\fun_c, \fun'_c) \inspace \V.\]
By adjunction and the left unit isomorphism in $\V$, this yields a morphism in $\V$ that we also denote by
\begin{equation}\label{tha_c}
\fun_c \fto{\tha_c} \fun'_c.
\end{equation}
\begin{description}
\item[Simplicial morphism]
The simplicial morphism
\begin{equation}\label{bcdot_hctha}
\Bcdot(\hun,\C,\fun) \fto{\Bcdot(\hun,\C,\tha)} \Bcdot(\hun,\C,\fun')
\end{equation}
is defined at simplicial level $r \geq 0$ by the following commutative diagrams in $\V$ for $c_0,\ldots,c_r \in \C$. 
\begin{equation}\label{Brhctha}
\begin{tikzpicture}[vcenter]
\def\v{-1.4}
\draw[0cell]
(0,0) node (a11) {\hun_{c_r} \otimes [\txotimes_{\ell=1}^r \C(c_{r-\ell}, c_{r-\ell+1})] \otimes \fun_{c_0}}
(a11)++(4.5,0) node (a12) {\Bc_r(\hun,\C,\fun)}
(a11)++(0,\v) node (a21) {\hun_{c_r} \otimes [\txotimes_{\ell=1}^r \C(c_{r-\ell}, c_{r-\ell+1})] \otimes \fun'_{c_0}}
(a12)++(0,\v) node (a22) {\Bc_r(\hun,\C,\fun')}
;
\draw[1cell=.9]
(a11) edge node {\iota} (a12)
(a12) edge[transform canvas={xshift=-1.5em}] node {\Bc_r(\hun,\C,\tha)} (a22)
(a11) edge[transform canvas={xshift=6em}] node[swap] {1^{\otimes r+1} \otimes \tha_{c_0}} (a21)
(a21) edge node {\iota} (a22)
;
\end{tikzpicture}
\end{equation}
\begin{itemize}
\item The compatibility of $\Bcdot(\hun,\C,\tha)$ with the 0-th face morphism $d_0$ follows from the $\V$-naturality of $\tha$. 
\item The compatibility of $\Bcdot(\hun,\C,\tha)$ with other face morphisms and the degeneracy morphisms follows from the functoriality of the monoidal product $\otimes$ of $\V$.
\end{itemize} 
\item[Realization]
Passing the simplicial morphism \cref{bcdot_hctha} to the realized bar constructions \cref{bar_realization} yields a morphism
\begin{equation}\label{bar_real_mor}
\Bc(\hun,\C,\fun) \fto{\Bc(\hun,\C,\tha)} \Bc(\hun,\C,\fun') \inspace \V.
\end{equation}
\item[Bar construction]
Applying \cref{bar_real_mor} with $\hun = \C_c = \C(-,c)$ \cref{representable_c}, the $c$-component morphisms 
\[\Bc(\C_c,\C,\fun) \fto{\Bc(\C_c,\C,\tha)} \Bc(\C_c,\C,\fun') \forspace c\in \C\]
define a $\V$-natural transformation
\begin{equation}\label{BCC_tha}
\begin{tikzpicture}[vcenter]
\def\u{.5}
\draw[0cell]
(0,0) node (a1) {\C}
(a1)++(2.5,0) node (a2) {\Vse}
;
\draw[1cell=.9]
(a1) [rounded corners=3pt] |- ($(a1)+(1,\u)$) -- node {\Bc(\C,\C,\fun)} ($(a2)+(-1,\u)$) -| (a2)
;
\draw[1cell=.9]
(a1) [rounded corners=3pt] |- ($(a1)+(1,-\u)$) -- node[swap] {\Bc(\C,\C,\fun')} ($(a2)+(-1,-\u)$) -| (a2)
;
\draw[2cell=.9]
node[between=a1 and a2 at .25, rotate=-90, 2label={above, \Bc(\C,\C,\tha)}] {\Rightarrow}
;
\end{tikzpicture}
\end{equation}
between the bar constructions of $\fun$ and $\fun'$ \cref{Bc_CCf}, called the \index{bar construction}\emph{bar construction of $\tha$}.  The $\V$-naturality of $\Bc(\C,\C,\tha)$ follows from \cref{BCCf_functor}, \cref{Brhctha}, and the functoriality of the monoidal product $\otimes$ of $\V$.
\item[Bar functor]
Denote by $\CVse$\label{not:CVse} the category with
\begin{itemize}
\item $\V$-functors $\C \to \Vse$ as objects and
\item $\V$-natural transformations as morphisms.
\end{itemize}  
The \emph{bar functor}\index{bar functor}
\begin{equation}\label{BCC_functor}
\CVse \fto{\Bc(\C,\C,-)} \CVse
\end{equation}
is defined by 
\begin{itemize}
\item the object assignment $\fun \mapsto \Bc(\C,\C,\fun)$ \cref{Bc_CCf} and
\item the morphism assignment $\tha \mapsto \Bc(\C,\C,\tha)$ \cref{BCC_tha}.\defmark
\end{itemize}
\end{description}
\end{definition}

\begin{explanation}[Naturality of Retractions and Sections]\label{expl:Bcdot_tha}
The simplicial retraction $\retn_c$ \cref{retn_c} and the section $\secn_c$ at an object $c \in \C$ \cref{secn_c} are natural in the $\V$-functor $\fun \cn \C \to \Vse$.  This means that, for each $\V$-natural transformation $\tha \cn \fun \to \fun'$, the diagram
\begin{equation}\label{ret_sec_nat}
\begin{tikzpicture}[vcenter]
\def\h{2.6} \def\v{-1.4}
\draw[0cell]
(0,0) node (a11) {\Bcdot(\C_c,\C,\fun)}
(a11)++(-\h,0) node (a10) {(\fun_c)_\crdot}
(a11)++(\h,0) node (a12) {(\fun_c)_\crdot}
(a11)++(0,\v) node (a21) {\Bcdot(\C_c,\C,\fun')}
(a10)++(0,\v) node (a20) {(\fun'_c)_\crdot}
(a12)++(0,\v) node (a22) {(\fun'_c)_\crdot}
;
\draw[1cell=.9]
(a10) edge node {\secn_c} (a11)
(a11) edge node {\retn_c} (a12)
(a20) edge node {\secn_c} (a21)
(a21) edge node {\retn_c} (a22)
(a11) edge[transform canvas={xshift=1.5em}] node[swap] {\Bcdot(\C_c,\C,\tha)} (a21)
(a10) edge node[swap] {(\tha_c)_\crdot} (a20)
(a12) edge node {(\tha_c)_\crdot} (a22)
;
\end{tikzpicture}
\end{equation}
of simplicial morphisms commutes, where $(\tha_c)_\crdot$ is the constant simplicial morphism at $\tha_c$ \cref{tha_c}. 
\begin{description}
\item[Sections] The left rectangle commutes by \cref{secn_cr}, \cref{Brhctha}, the naturality of the left unit isomorphism of $\V$, and the functoriality of $\otimes$.  However, $\secn_c$ is not $\V$-natural in $c \in \C$ \cref{secn_nonnat}.
\item[Retractions] The right rectangle commutes by \cref{retn_cr}, \cref{Brhctha}, and the $\V$-naturality of $\tha$.  Passing the right rectangle to realizations yields a commutative diagram
\begin{equation}\label{retract_theta}
\begin{tikzpicture}[vcenter]
\def\v{-1.3}
\draw[0cell]
(0,0) node (a11) {\Bc(\C,\C,\fun)}
(a11)++(2,0) node (a12) {\fun}
(a11)++(0,\v) node (a21) {\Bc(\C,\C,\fun')}
(a12)++(0,\v) node (a22) {\fun'}
;
\draw[1cell=.9]
(a11) edge node {\retn} (a12)
(a12) edge node {\tha} (a22)
(a11) edge[transform canvas={xshift=1.3em}] node[swap] {\Bc(\C,\C,\tha)} (a21)
(a21) edge node {\retn} (a22)
;
\end{tikzpicture}
\end{equation}
in $\CVse$.  In other words, the retraction $\retn$ \cref{retn} defines a natural transformation
\begin{equation}\label{retn_bar_id}
\begin{tikzpicture}[vcenter]
\def\t{24}
\draw[0cell]
(0,0) node (a1) {\phantom{\CVse}}
(a1)++(-.04,0) node (a1') {\CVse}
(a1)++(2.2,0) node (a2) {\CVse}
;
\draw[1cell=.9]
(a1) edge[bend left=\t] node {\Bc(\C,\C,-)} (a2)
(a1) edge[bend right=\t] node[swap] {1} (a2)
;
\draw[2cell]
node[between=a1 and a2 at .43, rotate=-90, 2label={above,\retn}] {\Rightarrow}
;
\end{tikzpicture}
\end{equation}
from the bar functor \cref{BCC_functor} to the identity functor.\defmark
\end{description}
\end{explanation}

\section{Homotopical Shimakawa $K$-Theory}
\label{sec:h_shim}

This section defines the homotopical Shimakawa (strong) $K$-theory functors
\[\AlglaxO \fto{\Khsho} \Gspec \andspace \AlgpspsO \fto{\Khshosg} \Gspec\]
from $\Op$-pseudoalgebras to orthogonal $G$-spectra for a compact Lie group $G$ and a 1-connected $\Gcat$-operad $\Op$.

\secoutline
\begin{itemize}
\item \cref{def:gtop_retn_secn} defines the bar functor $\BcFG$ on the category $\FGTopg$ of $\FGG$-spaces, along with the retraction $\retn \cn \BcFG \to 1$.
\item \cref{def:khsho} defines the homotopical Shimakawa (strong) $K$-theory functors $\Khsho$ and $\Khshosg$ from $\Op$-pseudoalgebras to orthogonal $G$-spectra, with further elaboration given in \cref{expl:khsho}. 
\item \cref{expl:khsho_shi}  discusses the fact that the homotopical Shimakawa strong $K$-theory $\Khshbesg$ for the Barratt-Eccles operad $\BE$ is objectwise level $G$-equivalent to Shimakawa's original equivariant $K$-theory.
\end{itemize}

\subsection*{Bar Construction}

\cref{def:gtop_retn_secn} applies the general constructions in \cref{sec:bar_const} to $\FGG$-spaces \pcref{def:ggtopg}.  In addition to restricting the constructions, we also need to take into account the basepoints throughout.

\begin{definition}[Bar Construction for $\FGG$-Spaces]\label{def:gtop_retn_secn}\index{bar construction!FGG-space@$\FGG$-spaces}
We apply \cref{def:cat_bar,def:bar_functor} in the following context for an arbitrary group $G$.
\begin{itemize}
\item $\V$ is the complete and cocomplete symmetric monoidal closed category 
\[\V = (\Gtopst, \sma, \stplus, \Topgst)\]
of pointed $G$-spaces \cref{Gtopst_smc}.
\item $\C$ is the small pointed $G$-category $\FG$ of pointed finite $G$-sets and pointed functions \pcref{def:FG}.  Each pointed $G$-set $\FG(\mal,\nbeta)$ is regarded as a discrete pointed $G$-space.
\end{itemize}
\begin{description}
\item[Cosimplicial pointed $G$-space]
Define the functor\label{not:simpdotpl}
\[\Delta \fto{\simpdotpl} \Gtopst\] 
that sends $\ordr \in \Delta$ to the pointed $G$-space\label{not:simpplr}
\[\simppl^r = \simp^r \sqcup {*}.\]
This pointed $G$-space is obtained from the topological $r$-simplex 
\begin{equation}\label{top_simplex}
\simp^r = \big\{ (t_0,\ldots,t_{r}) \in \mathbbm{R}^{r+1} \tmid t_i \geq 0,\, \txsum_{i=0}^r t_i = 1 \big\},
\end{equation}
equipped with the trivial $G$-action, by adjoining a disjoint $G$-fixed basepoint $*$.  
\begin{description}
\item[Cofaces] For the $i$-th coface morphism $d^i$ \cref{coface_codeg}, the pointed morphism 
\[\simppl^{r-1} \fto{d^i} \simppl^r\] 
inserts 0
\begin{itemize}
\item to the left of $t_i$ if $0 \leq i < r$ and
\item to the right of $t_{r-1}$ if $i=r$.
\end{itemize} 
\item[Codegeneracies] For the $i$-th codegeneracy morphism $s^i$, the pointed morphism
\[\simppl^{r+1} \fto{s^i} \simppl^r\]
replaces $(t_{i}, t_{i+1})$ by the sum $t_{i}+t_{i+1}$.
\end{description}
\item[Realizations]
The realization of a simplicial pointed $G$-space $X_\crdot \cn \Deltaop \to \Gtopst$ is the pointed $G$-space \cref{x_dot}
\begin{equation}\label{real_simpgspace}
|X_\crdot| = \int^{\ordr \in \Delta} X_{\ordr} \sma \simppl^r.
\end{equation}
\item[Bar construction]
For an $\FGG$-space $X \cn \FG \to \Topgst$ \cref{ggtopg_obj}, the $\FGG$-space 
\begin{equation}\label{bar_FG}
\FG \fto{\Bc(\FG,\FG,X)} \Topgst
\end{equation}
is the bar construction of $X$ \cref{Bc_CCf}.  It is also denoted by $\Bc X$.  It sends each pointed finite $G$-set $\mal \in \FG$ to the pointed $G$-space given by the realized bar construction \cref{real_Cc}
\begin{equation}\label{real_FG}
\Bc\big(\FG(-,\mal),\FG,X\big) = \int^{\ordr \in \Delta} \Bc_r\big(\FG(-,\mal),\FG,X\big) \sma \simppl^r.
\end{equation}
In the coend in \cref{real_FG}, the first term is the pointed $G$-space \cref{bc_r}
\begin{equation}\label{Br_FGX}
\begin{split}
&\Bc_r\big(\FG(-,\mal),\FG,X\big) \\
&= \bigvee_{\ang{\ordm_\ell^{\al_\ell}}_{\ell=0}^r} \mabit
\FG(\ordm_r^{\al_r}, \mal) \sma \FG(\ordm_{r-1}^{\al_{r-1}}, \ordm_r^{\al_r}) \sma \Cdots \sma \FG(\ordm_0^{\al_0}, \ordm_1^{\al_1}) \sma X_{\ordm_0^{\al_0}}
\end{split}
\end{equation}
with the wedge indexed by the set of $(r+1)$-tuples $\ang{\ordm_\ell^{\al_\ell}}_{\ell=0}^r$ of pointed finite $G$-sets.  The group $G$ acts diagonally on the iterated smash product in \cref{Br_FGX} and by conjugation \cref{gpsi} on each pointed $G$-set $\FG(-,-)$.
\item[Bar functor]
Associated to each morphism $\tha \cn X \to X'$ in $\FGTopg$ \cref{ggtopg_icell} is its bar construction \cref{BCC_tha}, which is the $G$-natural transformation
\begin{equation}\label{bar_FG_tha}
\begin{tikzpicture}[vcenter]
\def\u{.5}
\draw[0cell]
(0,0) node (a1) {\FG}
(a1)++(3.3,0) node (a2) {\Topgst}
;
\draw[1cell=.9]
(a1) [rounded corners=3pt, shorten >=-.3ex] |- ($(a1)+(1,\u)$) -- node {\Bc(\FG,\FG,X)} ($(a2)+(-1,\u)$) -| (a2)
;
\draw[1cell=.9]
(a1) [rounded corners=3pt, shorten >=-.3ex] |- ($(a1)+(1,-\u)$) -- node[swap] {\Bc(\FG,\FG,X')} ($(a2)+(-1,-\u)$) -| (a2)
;
\draw[2cell=1]
node[between=a1 and a2 at .22, rotate=-90, 2label={above, \Bc(\FG,\FG,\tha)}] {\Rightarrow}
;
\end{tikzpicture}
\end{equation}
between the bar constructions of $X$ and $X'$.  The bar functor \cref{BCC_functor}
\begin{equation}\label{bar_functor_FG}
\FGTopg \fto{\BcFG} \FGTopg
\end{equation}
is defined by 
\begin{itemize}
\item the object assignment $X \mapsto \Bc(\FG,\FG,X)$ \cref{bar_FG} and
\item the morphism assignment $\tha \mapsto \Bc(\FG,\FG,\tha)$ \cref{bar_FG_tha}.
\end{itemize}
The bar functor $\BcFG$ is also denoted by $\Bc$.
\item[Retractions] 
The retraction of an $\FGG$-space $X$ is the $G$-natural transformation \cref{retn}
\begin{equation}\label{ret_FG}
\begin{tikzpicture}[vcenter]
\def\t{28}
\draw[0cell]
(0,0) node (a1) {\phantom{\Gsk}}
(a1)++(1.8,0) node (a2) {\phantom{\Gsk}}
(a1)++(-.08,0) node (a1') {\FG}
(a2)++(.23,0) node (a2') {\Topgst}
;
\draw[1cell=.9]
(a1) edge[bend left=\t] node {\Bc X} (a2)
(a1) edge[bend right=\t] node[swap] {X} (a2)
;
\draw[2cell]
node[between=a1 and a2 at .37, rotate=-90, 2label={above,\retn_X}] {\Rightarrow}
;
\end{tikzpicture}
\end{equation}
between $\FGG$-spaces, that is, a morphism in $\FGTopg$ \cref{ggtopg_icell}. 
\begin{description}
\item[$G$-homotopy equivalence]
Passing the equality and homotopy in \cref{retn_secn} to realizations implies that each component of the retraction $\retn_X$ is a pointed $G$-homotopy equivalence between pointed $G$-spaces.  In particular, $\retn_X$ is componentwise a weak $G$-equivalence \pcref{def:weakG_top}.
\item[Naturality]
As $X \in \FGTopg$ varies, the retraction $\retn_X$ defines a natural transformation \cref{retn_bar_id}
\begin{equation}\label{retn_barFG_id}
\begin{tikzpicture}[vcenter]
\def\t{28}
\draw[0cell]
(0,0) node (a1) {\phantom{\Fsk}}
(a1)++(1.8,0) node (a2) {\phantom{\Fsk}}
(a1)++(-.43,0) node (a1') {\FGTopg}
(a2)++(.43,0) node (a2') {\FGTopg}
;
\draw[1cell=.9]
(a1) edge[bend left=\t] node {\Bc} (a2)
(a1) edge[bend right=\t] node[swap] {1} (a2)
;
\draw[2cell]
node[between=a1 and a2 at .43, rotate=-90, 2label={above,\retn}] {\Rightarrow}
;
\end{tikzpicture}
\end{equation}
from the bar functor \cref{bar_functor_FG} to the identity functor.\defmark
\end{description}
\end{description}
\end{definition}

\subsection*{Shimakawa Equivariant $K$-Theory Machines}

Next, we define the homotopical versions of Shimakawa (strong) $K$-theory, going from $\Op$-pseudoalgebras to orthogonal $G$-spectra \pcref{def:pseudoalgebra,def:gsp_module}.  For the reader's convenience, \cref{def:khsho} repeats some of the description of the functors from \cref{def:ksho}.

\begin{definition}\label{def:khsho}
For a compact Lie group $G$ and a 1-connected $\Gcat$-operad $\Op$ \cref{i_connected}, we define the \emph{homotopical Shimakawa $K$-theory}\index{homotopical Shimakawa $K$-theory}\index{Shimakawa K-theory@Shimakawa $K$-theory!homotopical} $\Khsho$ and the \emph{homotopical Shimakawa strong $K$-theory} $\Khshosg$ for $\Op$ as the following composite functors.
\begin{equation}\label{khsho_khshosg}
\begin{tikzpicture}[vcenter]
\def\h{2.2} \def\u{.7} \def\v{-1.8}
\draw[0cell]
(0,0) node (a1) {\AlglaxO}
(a1)++(0,.1*\v) node (a1') {\phantom{A}}
(a1)++(0,\v) node (a2) {\AlgpspsO}
(a2)++(0,-.1*\v) node (a2') {\phantom{A}}
(a1)++(.75*\h,\v/2) node (a3) {\phantom{\FGCatg}}
(a3)++(0,-.01) node (a3') {\FGCatg}
(a3)++(\h,0) node (a4) {\FGTopg}
(a4)++(\h,0) node (a4') {\FGTopg}
(a4')++(.9*\h,0) node (a5) {\Gspec}
;
\draw[1cell=.8]
(a1') edge[shorten >=-.4ex] node[swap, inner sep=1pt] {\Sgo} (a3)
(a2') edge[shorten >=-.4ex] node[inner sep=1pt] {\Sgosg} (a3)
(a3) edge node {\clast} (a4)
(a4) edge node {\Bc} (a4')
(a4') edge node {\Kfg} (a5)
;
\draw[1cell=.9]
(a1) [rounded corners=2pt] -| node[pos=.21] {\Khsho}  (a5)
;
\draw[1cell=.9]
(a2) [rounded corners=2pt] -| node[pos=.21] {\Khshosg}  (a5)
;
\end{tikzpicture}
\end{equation}
\begin{description}
\item[$\Op$-pseudoalgebras to $\FGG$-categories] 
$\Sgo$ and $\Sgosg$ are Shimakawa (strong) $H$-theory \pcref{Sgo_twofunctor} between the categories in \cref{oalgps_twocat,def:fgcatg}. 
\item[$\FGG$-categories to $\FGG$-spaces] 
$\clast$ is induced by the classifying space functor $\cla \cn \Cat \to \Top$ \pcref{ggcatg_ggtopg}.
\item[Bar construction]
$\Bc = \BcFG$ is the bar functor on $\FGTopg$ \cref{bar_functor_FG}.  The functors $\Sgo$, $\Sgosg$, $\clast$, and $\Bc$ are defined for arbitrary groups $G$. 
\item[$\FGG$-spaces to orthogonal $G$-spectra] 
The functor $\Kfg$ \pcref{def:Kfg_functor}, going between the categories in \cref{def:gsp_morphism,def:ggtopg}, requires $G$ to be a compact Lie group.\defmark
\end{description}
\end{definition}

\begin{explanation}[Unpacking]\label{expl:khsho}
For an $\Op$-pseudoalgebra $\A$ and an object $V \in \IU$ \pcref{def:pseudoalgebra,def:indexing_gspace}, the pointed $G$-spaces $(\Khsho \A)_V$ and $(\Khshosg\A)_V$ are given by the coends \cref{Kfgxv}
\[\begin{split}
(\Khsho \A)_V &= \int^{\mal \sins \FG} (S^V)^{\mal} \sma \Bc\big(\FG(-,\mal),\FG,\cla \Adash\big) \andspace\\
(\Khshosg \A)_V &= \int^{\mal \sins \FG} (S^V)^{\mal} \sma \Bc\big(\FG(-,\mal),\FG,\cla \Asgdash\big).
\end{split}\]
\begin{itemize}
\item The sphere action is defined in \cref{Kfgx_action_vw}. 
\item The diagonal $G$-action is defined in \cref{Kfgxv_rep_gact,Br_FGX}.  
\item $\Adash = \Sgo\A$ and $\Asgdash = \Sgosg\A$ are Shimakawa (strong) $H$-theory of $\A$ \pcref{sys_FGcat}.
\item Each pointed $G$-space $\Bc(\Cdots)$ is the realized bar construction \cref{real_FG}.  The first one is the coend
\[\Bc\big(\FG(-,\mal),\FG,\cla\Adash\big) 
= \int^{\ordr \in \Delta} \Bc_r\big(\FG(-,\mal),\FG,\cla\Adash\big) \sma \simppl^r\]
where
\[\begin{split}
&\Bc_r\big(\FG(-,\mal),\FG,\cla\Adash\big) \\
&= \bigvee_{\ang{\ordm_\ell^{\al_\ell}}_{\ell=0}^r} \mabit
\FG(\ordm_r^{\al_r}, \mal) \sma \FG(\ordm_{r-1}^{\al_{r-1}}, \ordm_r^{\al_r}) \sma \Cdots \sma \FG(\ordm_0^{\al_0}, \ordm_1^{\al_1}) \sma \cla (\sys{\A}{\ordm_0^{\al_0}}),
\end{split}\]
and likewise for the strong variant involving $\cla\Asgdash$.  The pointed $G$-categories $\sys{\A}{\ordm_0^{\al_0}}$ and $\syssg{\A}{\ordm_0^{\al_0}}$ are those of (strong) $\ordm_0^{\al_0}$-systems in $\A$ \pcref{def:nsys_gcat}.
\end{itemize}
The difference between Shimakawa (strong) $K$-theory and the homotopical Shimakawa (strong) $K$-theory \pcref{def:ksho,def:khsho} is that the latter involves the bar functor $\Bc$ before the final step $\Kfg$.  The bar functor is connected to the identity functor on $\FGTopg$ via the retraction $\retn \cn \Bc \to 1$ \cref{retn_barFG_id}.  The retraction is componentwise a pointed $G$-homotopy equivalence \cref{ret_FG}.  
\end{explanation}

\begin{explanation}[Shimakawa's Construction]\label{expl:khsho_shi}
Continuing \cref{expl:shi89_section2,expl:c_shimakawa}, Shimakawa's original equivariant $K$-theory machine \cite{shimakawa89} corresponds to the functor $\Khshbesg$ \cref{khsho_khshosg} for the Barratt-Eccles operad $\BE$ \pcref{def:BE}.  There are some technical differences, due to more recent advances, that we discuss next.
\begin{description}
\item[$G$-spectra]
Orthogonal $G$-spectra \cite{mandell_may} were invented after Shimakawa's work \cite{shimakawa89}, which uses Lewis-May $G$-spectra \cite[Def.\ I.2.1]{lewis_may_steinberger}.  All of Shimakawa's constructions in \cite{shimakawa89,shimakawa91} still work for orthogonal $G$-spectra.
\item[Smash vs. Cartesian products] 
For an $\FGG$-space $X$, the bar construction $\Bc X$ \cref{bar_FG} uses the smash product.  In contrast, Shimakawa's construction of almost $\Omega$-$G$-spectra in \cite[page 246]{shimakawa89}, denoted by $\Shimg$ there, involves the Cartesian product and a quotient by a $G$-contractible subspace.  This difference can be explained as follows.  Since $\FG(-,\ord{0})$ is a point, so is the realized bar construction \cref{real_FG}
\[\Bc\big(\FG(-,\ord{0}), \FG,X),\]
as required for an $\FGG$-space.  However, the Cartesian variant of the realized bar construction is \emph{not} a point in general.  Thus, taking a suitable quotient is necessary to get a single point when evaluated at $\ord{0}$.  We follow \cite[Section 3]{mmo} in using the smash product in the bar construction $\Bc X$.  By \cite[Theorem 3.19]{gmmo19}, for each $\FGG$-space $X$, the orthogonal $G$-spectra $\Shimg X$ and $\Kfg\Bc X$ are naturally componentwise weakly $G$-equivalent.
\end{description}
The upshot is that Shimakawa's original equivariant $K$-theory machine  is naturally componentwise weakly $G$-equivalent to the homotopical Shimakawa strong $K$-theory $\Khshbesg$ for the Barratt-Eccles operad $\BE$ \cref{khsho_khshosg}.  Moreover, $\Khshbesg$ can be replaced by Shimakawa strong $K$-theory $\Kshbesg$ \cref{ksho_kshosg}; see \cref{expl:k_shimakawa}.
\end{explanation}

\begin{explanation}[Homotopical Machine]\label{expl:gen_shimakawa}
Continuing \cref{expl:Kfg_functor}, in \cite[Def.\ 3.24]{mmo}, the bar functor $\Bc$ \cref{bar_functor_FG} is denoted by $\mathbbm{I}$, and the composite \cref{khsho_khshosg} 
\[\FGTopg \fto{\Bc} \FGTopg \fto{\Kfg} \Gspec\] 
is called the \emph{genuine homotopical Segal machine}\index{genuine homotopical Segal machine} on $\FGG$-spaces.  The latter is denoted by $\mathbbm{S}_{\scriptscriptstyle G}^{\mathbbm{N}_G}$ in \cite[Def.\ 3.17]{gmmo19}.
\end{explanation}

\section{Preservation of Weak $G$-Equivalences}
\label{sec:invariance_thm}

This section proves that, for a componentwise weak $G$-equivalence $\tha$ between $\FGG$-spaces that satisfy a cofibrancy condition called \emph{properness}, the prolongation of $\tha$ to pointed finite $G$-CW complexes is also componentwise a weak $G$-equivalence \pcref{thm:invariance}.  This result is used in \cref{sec:ch_shim} to show that Shimakawa (strong) $K$-theory and its homotopical variant are naturally weakly $G$-equivalent.  Throughout this section, $G$ denotes a finite group.

\secoutline
\begin{itemize}
\item \cref{def:gcof,def:reedy_cof} recall Reedy cofibrancy of simplicial $G$-spaces.
\item \cref{reedy_cofibrant} records two criteria for Reedy cofibrancy.
\item \cref{reedy_realization} proves that realization sends a componentwise weak $G$-equivalence between Reedy cofibrant simplicial $G$-spaces to a weak $G$-equivalence.
\item \cref{def:gssetst,def:proper_fgg} define \emph{proper} $\FGG$-spaces.
\item \cref{thm:invariance} proves that prolongation to pointed finite $G$-CW complexes preserves componentwise weak $G$-equivalences between proper $\FGG$-spaces.
\item \cref{bar_proper} records the fact that the bar construction of an $\FGG$-space is always proper. 
\end{itemize}

\subsection*{Reedy Cofibrancy}
We first recall some definitions and facts about Reedy cofibrancy from \cite[Section 1.2]{mmo}.  Other general references include \cite[App.\ 2]{boardman-vogt}, \cite[Ch.\ 8]{hill_hopkins_ravenel}, \cite[I.1]{lewis_may_steinberger}, \cite[Ch.\ 6]{mayconcise}, and \cite[App.\ B]{schwede_ght}.  

\begin{definition}\label{def:gcof}
A $G$-morphism $f \cn X \to Y$ between $G$-spaces \pcref{def:Gtop} is a \index{G-cofibration@$G$-cofibration}\emph{$G$-cofibration} if, given each solid-arrow commutative diagram of $G$-morphisms 
\begin{equation}\label{Gcofibration}
\begin{tikzpicture}[vcenter]
\def\h{3} \def\v{-2}
\draw[0cell]
(0,0) node (a11) {X}
(a11)++(\h,0) node (a12) {X \times [0,1]}
(a11)++(0,\v) node (a21) {Y}
(a12)++(0,\v) node (a22) {Y \times [0,1]}
(a11)++(\h/2,\v/2) node (a0) {Z}
;
\draw[1cell=.85]
(a11) edge node {i_0} (a12)
(a12) edge node {f \times \Id} (a22)
(a11) edge node[swap] {f} (a21)
(a21) edge node {i_0} (a22)
(a12) edge node[swap,pos=.6] {h} (a0)
(a21) edge node[pos=.6] {j} (a0)
(a22) edge[densely dashed] node[swap,pos=.6] {k} (a0)
;
\end{tikzpicture}
\end{equation}
with each $i_0(x) = (x,0)$, there exists a $G$-morphism $k$ that makes the diagram commute.
\end{definition}

\begin{definition}\label{def:reedy_cof}
For a simplicial $G$-space\index{simplicial G-space@simplicial $G$-space} $X_\crdot \in \SGtop$ \pcref{def:simp_obj} and $r \geq 0$, the \emph{$r$-th latching $G$-space}\index{latching G-space@latching $G$-space} is defined as the union
\[\Lat_r X_\crdot = \bigcup_{i=0}^{r-1} s_i X_{r-1} \bigsubset X_r.\]
A simplicial $G$-space $X_\crdot$ is \emph{Reedy cofibrant}\index{Reedy cofibrant} if the inclusion 
\[\Lat_r X_\crdot \fto{\incn_r} X_r\]
is a $G$-cofibration for each $r \geq 0$.
\end{definition}

The following criteria for Reedy cofibrancy are \cite[Remark 1.9 and Lemma 1.11]{mmo}.
 
\begin{lemma}\label{reedy_cofibrant}
Suppose $X_\crdot$ is a simplicial $G$-space. 
\begin{enumerate}
\item\label{reedy_cof_i}
Suppose $X_\crdot$ is Reedy cofibrant.  For each subgroup $H \subseteq G$, passing to $H$-fixed point spaces levelwise yields a Reedy cofibrant simplicial space $X_\crdot^H$.
\item\label{reedy_cof_ii} 
If each degeneracy morphism 
\[X_r \fto{s_i} X_{r+1} \forspace 0 \leq i \leq r\]
is a $G$-cofibration, then $X_\crdot$ is a Reedy cofibrant simplicial $G$-space.
\end{enumerate}
\end{lemma}

For \cref{reedy_realization}, recall that a $G$-morphism $f$ between $G$-spaces is a \emph{weak $G$-equivalence} if the restriction $f^H$ to $H$-fixed point spaces is a weak equivalence of spaces for each subgroup $H \subseteq G$ \pcref{def:weakG_top}.  The realization $|X_\crdot|$ of a simplicial $G$-space $X_\crdot$ is defined in \cref{x_dot}, using the topological $r$-simplex $\simp^r$ with the trivial $G$-action \cref{top_simplex} for $r \geq 0$.  \cref{reedy_realization} is \cite[Theorem 1.12]{mmo}, which states that realization preserves weak $G$-equivalences between Reedy cofibrant simplicial $G$-spaces.

\begin{theorem}\label{reedy_realization}
Suppose $\fun \cn X_\crdot \to Y_\crdot$ is a simplicial morphism between simplicial $G$-spaces such that the following two conditions hold.
\begin{itemize}
\item Each of $X_\crdot$ and $Y_\crdot$ is Reedy cofibrant. 
\item $\fun_r \cn X_r \to Y_r$ is a weak $G$-equivalence for each $r \geq 0$.  
\end{itemize}
Then the realization
\[|X_\crdot| \fto{|\fun|} |Y_\crdot|\]
is a weak $G$-equivalence.
\end{theorem}

\begin{proof}
For each subgroup $H \subseteq G$, the $H$-fixed point morphism $|\fun|^H$ is isomorphic to the realization $|\fun^H|$ because realization preserves finite limits, including taking $H$-fixed point spaces.  By the assumed Reedy cofibrancy of $X_\crdot$ and $Y_\crdot$, \cref{reedy_cofibrant} \eqref{reedy_cof_i} implies that $\fun^H \cn X_\crdot^H \to Y_\crdot^H$ is a simplicial morphism between Reedy cofibrant simplicial spaces.  It is levelwise a weak equivalence by assumption.  Thus, the \namecref{reedy_realization} follows from the nonequivariant version, which is proved in \cite[A.4]{may-groupcompletion} and \cite[A.44]{schwede_ght}.
\end{proof}

\subsection*{Properness}

\begin{definition}\label{def:gssetst}\
\begin{enumerate}
\item\label{def:gssetst_ssetst} 
Denote by $\pSet$ the category of pointed sets and pointed functions.  Denote by $\ssetst$ the category of simplicial objects in $\pSet$ \pcref{def:simp_obj}.  Its objects and morphisms are called \emph{pointed simplicial sets}\index{pointed simplicial!set} and \emph{pointed simplicial morphisms}.
\item\label{def:gssetst_gsetst}
For a group $G$, denote by $\GSetst$ the category of pointed $G$-sets \pcref{def:ptGset} and pointed $G$-equivariant functions.  Denote by $\Gssetst$ the category of simplicial objects in $\GSetst$.  Its objects are called \index{pointed simplicial!G-set@$G$-set}\emph{pointed simplicial $G$-sets}.
\item\label{def:gssetst_degen}
For a pointed simplicial $G$-set $X_\crdot$, an element $x \in X_r$ is called an \index{simplex}\emph{$r$-simplex}.  An $r$-simplex is \emph{degenerate}\index{simplex!degenerate} if it has the form $s_i x$ for some $(r-1)$-simplex $x \in X_{r-1}$ and some degeneracy morphism $s_i$ with $0 \leq i \leq r-1$.  An $r$-simplex is \emph{nondegenerate}\index{simplex!nondegenerate} if it does not have the form $s_i x$.
\item\label{def:gssetst_finite} 
A \index{pointed simplicial!G-set@$G$-set!finite}\emph{pointed finite simplicial $G$-set} is a pointed simplicial $G$-set with only a finite number of nondegenerate simplices.  Each simplicial level of a pointed finite simplicial $G$-set is a pointed finite $G$-set.  Each pointed finite $G$-set is also regarded as a pointed finite $G$-CW complex with only 0-cells.\defmark
\end{enumerate}
\end{definition}

For \cref{def:proper_fgg}, we recall the following.
\begin{itemize}
\item The prolongation functor \cref{PwgUgs}
\[\FGTopg \fto{\Pwg} \WGTopg\]
sends each $\FGG$-space \pcref{def:ggtopg} to a $\WGG$-space.  Here, $\WG$ is the pointed $G$-category of pointed finite $G$-CW complexes and pointed morphisms with the conjugation $G$-action.  The functor $\Pwg$ is defined objectwise as a coend \cref{Pwgxa}. 
\item Each $\FGG$-space $X$ satisfies the $G$-equivariance property \cref{gggspace_equiv} 
\[X(g \cdot \psi) = g(X\psi)\ginv\]
for pointed morphisms $\psi \in \FG$.  In particular, $\FGG$-spaces preserve $G$-equivariant morphisms, and the same is true for $\WGG$-spaces.  Thus, evaluating a $\WGG$-space at a pointed finite simplicial $G$-set (\cref{def:gssetst} \eqref{def:gssetst_finite}) levelwise yields a simplicial $G$-space.
\end{itemize}

\begin{definition}\label{def:proper_fgg}
An $\FGG$-space $X \cn \FG \to \Topgst$ is \emph{proper}\index{FGG-space@$\FGG$-space!proper}\index{proper!FGG-space@$\FGG$-space} if, for each pointed finite simplicial $G$-set $B_\crdot$, the simplicial $G$-space $(\Pwg X)_{B_\crdot}$ is Reedy cofibrant \pcref{def:reedy_cof}.
\end{definition}

A morphism in each of the categories $\FGTopg$ and $\WGTopg$ \cref{ggtopg_icell} is componentwise a pointed $G$-morphism \cref{ggtopg_icell_geq}.  \cref{thm:invariance} proves that the prolongation functor $\Pwg$ \cref{PwgUgs} preserves componentwise weak $G$-equivalences between proper $\FGG$-spaces.

\begin{theorem}\label{thm:invariance}
Suppose $\tha \cn X \to X'$ is a morphism of $\FGG$-spaces such that the following two conditions hold.
\begin{itemize}
\item Each of $X$ and $X'$ is a proper $\FGG$-space.
\item $\tha_{\mal} \cn X\mal \to X'\mal$ is a weak $G$-equivalence for each pointed finite $G$-set $\mal \in \FG$.
\end{itemize}
Then the $G$-natural transformation
\begin{equation}\label{Ptheta}
\begin{tikzpicture}[vcenter]
\def\t{27}
\draw[0cell]
(0,0) node (a1) {\phantom{\Gsk}}
(a1)++(2.1,0) node (a2) {\phantom{\Gsk}}
(a1)++(-.08,0) node (a1') {\WG}
(a2)++(.23,0) node (a2') {\Topgst}
;
\draw[1cell=.9]
(a1) edge[bend left=\t] node {\Pwg X} (a2)
(a1) edge[bend right=\t] node[swap] {\Pwg X'} (a2)
;
\draw[2cell=.9]
node[between=a1 and a2 at .3, rotate=-90, 2label={above,\Pwg\tha}] {\Rightarrow}
;
\end{tikzpicture}
\end{equation}
is componentwise a weak $G$-equivalence.
\end{theorem}

\begin{proof}
We need to prove that $(\Pwg\tha)_A$ is a weak $G$-equivalence for each pointed finite $G$-CW complex $A$.  We choose a pointed $G$-homotopy equivalence 
\[|B_\crdot| \fto[\sim]{h} A\] 
from the realization of a pointed finite simplicial $G$-set $B_\crdot$, which exists by \cite[B.46 (ii)]{schwede_ght}.  By the naturality of $\Pwg\tha$, the diagram of pointed $G$-morphisms
\begin{equation}\label{Ptheta_nat}
\begin{tikzpicture}[vcenter]
\def\v{-1.3}
\draw[0cell=.9]
(0,0) node (a11) {(\Pwg X)_{|B_\crdot|}}
(a11)++(3.75,0) node (a12) {(\Pwg X')_{|B_\crdot|}}
(a11)++(0,\v) node (a21) {(\Pwg X)_A}
(a12)++(0,\v) node (a22) {(\Pwg X')_A}
;
\draw[1cell=.8]
(a11) edge node {(\Pwg\tha)_{|B_\crdot|}} (a12)
(a12) edge[transform canvas={xshift=-1em}] node {(\Pwg X')_h} (a22)
(a11) edge[transform canvas={xshift=1em}] node[swap] {(\Pwg X)_h} (a21)
(a21) edge node {(\Pwg\tha)_A} (a22)
;
\end{tikzpicture}
\end{equation}
commutes.  By the coend definition \cref{Pwgxa} of $(\Pwg X)_A$, it preserves pointed $G$-homotopies in the variable $A$.  Thus, each of $(\Pwg X)_h$ and $(\Pwg X')_h$ is a $G$-homotopy equivalence.  By the preceding commutative diagram, it suffices to prove that $(\Pwg\tha)_{|B_\crdot|}$ is a weak $G$-equivalence.  

Realization \cref{x_dot} commutes with finite products and $\Pwg$, which is also a coend.  Thus, it suffices to prove that the realization
\[\left|(\Pwg X)_{B_\crdot}\right| \fto{|(\Pwg\tha)_{B_\crdot}|} \left|(\Pwg X')_{B_\crdot}\right|\]
is a weak $G$-equivalence.  By the assumed properness of $X$ and $X'$, the simplicial $G$-spaces $(\Pwg X)_{B_\crdot}$ and $(\Pwg X')_{B_\crdot}$ are Reedy cofibrant.  By \cref{reedy_realization}, it suffices to prove that the $G$-morphism 
\[(\Pwg X)_B \fto{(\Pwg\tha)_B} (\Pwg X')_B\]
is a weak $G$-equivalence for each pointed finite $G$-set $B$.  Observe that $B$ is $G$-isomorphic to some object $\mal \in \FG$ and that $\FG$ is a full subcategory of $\WG$.  Thus, $(\Pwg\tha)_B$ is $G$-isomorphic to $\tha_{\mal}$, which is a weak $G$-equivalence by assumption.
\end{proof}

\begin{remark}[Related Literature]\label{rk:invariance}
\cref{thm:invariance} is a variant of \cite[Theorem 2.6]{gmmo19}.  The only difference between the two Theorems is that the latter does \emph{not} have the finiteness condition on the objects of $\WG$ and in the definition of properness.  In particular, an $\FGG$-space that is proper in the sense of \cite[Def.\ 2.3]{gmmo19} is also proper in the sense of \cref{def:proper_fgg}.  However, the converse is not necessarily true.  A variant of \cref{thm:invariance} for $\Fskg$-spaces is \cite[B.48]{schwede_ght}.
\end{remark}

\begin{lemma}\label{bar_proper}
For each $\FGG$-space $X$, the bar construction $\Bc X$ \cref{bar_FG} is a proper $\FGG$-space.
\end{lemma}

\begin{proof}
By \cite[Lemma 3.18]{gmmo19}, $\Bc X$ is proper in the sense of \cite[Def.\ 2.3]{gmmo19}.  Thus, it is also proper in the sense of  \cref{def:proper_fgg}; see \cref{rk:invariance}.
\end{proof}

\section{Comparison of Shimakawa $K$-Theories}
\label{sec:ch_shim}

This section proves that Shimakawa $K$-theory \cref{ksho_kshosg} and its homotopical variant \cref{khsho_khshosg} 
\begin{equation}\label{KshoKhsho_seci}
\begin{tikzpicture}[vcenter]
\draw[0cell]
(0,0) node (a1) {\AlglaxO}
(a1)++(4.3,0) node (a2) {\phantom{\Gspec}}
(a2)++(0,.05) node (a2') {\Gspec}
;
\draw[1cell=.9]
(a1) edge[transform canvas={yshift=.8ex}] node {\Ksho = \Kfg \clast \Sgo} (a2)
(a1) edge[transform canvas={yshift=-.1ex}] node[swap] {\Khsho = \Kfg \Bc \clast \Sgo} (a2)
;
\end{tikzpicture}
\end{equation}
are naturally componentwise weakly $G$-equivalent.  The strong variant is also true.  Thus, $\Ksho$ and $\Khsho$ are interchangeable, and the strong variants $\Kshosg$ and $\Khshosg$ are also interchangeable.  Unless otherwise specified, $G$ denotes a finite group in this section.  

\secoutline
\begin{itemize}
\item \cref{thm:Kfg_inv} proves that the functor $\Kfg$ sends a componentwise weak $G$-equivalence between proper $\FGG$-spaces to a componentwise weak $G$-equivalence between orthogonal $G$-spectra.
\item \cref{def:Gssetst} defines the category $\FGssetgst$ of $\FG$-simplicial $G$-sets.
\item \cref{reast_proper} proves that realization sends each $\FG$-simplicial $G$-set to a proper $\FGG$-space.
\item \cref{thm:ksho_khsho} proves that $\Ksho$ and $\Khsho$ are naturally componentwise weakly $G$-equivalent, and likewise for the strong variant.
\item \cref{expl:k_shimakawa} discusses the fact that Shimakawa's original equivariant $K$-theory machine \cite{shimakawa89} is naturally componentwise weakly $G$-equivalent to Shimakawa strong $K$-theory \cref{ksho_kshosg}
\[\AlgpspsBE \fto{\Kshbesg} \Gspec\]
for the Barratt-Eccles operad $\BE$.
\end{itemize}

\subsection*{$\Kfg$ Preserves Weak $G$-Equivalences}

Recall the functor \pcref{def:Kfg_functor}
\[\FGTopg \fto{\Kfg} \Gspec\]
that sends $\FGG$-spaces \pcref{def:ggtopg} to orthogonal $G$-spectra \pcref{def:gsp_morphism}.  The functor $\Kfg$ is the last step of (homotopical) Shimakawa (strong) $K$-theory \pcref{def:ksho,def:khsho}.

\begin{theorem}\label{thm:Kfg_inv}
Suppose $\tha \cn X \to X'$ is a morphism of $\FGG$-spaces such that the following two conditions hold.
\begin{itemize}
\item Each of $X$ and $X'$ is a proper $\FGG$-space \pcref{def:proper_fgg}.
\item $\tha_{\mal} \cn X\mal \to X'\mal$ is a weak $G$-equivalence for each pointed finite $G$-set $\mal \in \FG$.
\end{itemize}
Then for each object $V \in \IU$ \pcref{def:indexing_gspace}, the $V$-component pointed $G$-morphism \cref{Kfg_psiv}
\[(\Kfg X)_V \fto{(\Kfg\tha)_V} (\Kfg X')_V\]
is a weak $G$-equivalence.
\end{theorem}

\begin{proof}
With $S^V$ denoting the $V$-sphere \pcref{def:g_sphere}, by \cref{Kfg_psiv,Pwgxa}, $(\Kfg\tha)_V$ is equal to $(\Pwg\tha)_{S^V}$.  The latter is a weak $G$-equivalence by \cref{thm:invariance}.
\end{proof}

\subsection*{Realization is Proper}

\cref{def:Gssetst} is the simplicial set analogue of \cref{def:ggtopg}.

\begin{definition}\label{def:Gssetst}
Suppose $G$ is a group.
\begin{enumerate}
\item\label{def:Gssetst_ssetgst} 
Denote by $\ssetgst$ the pointed $G$-category with pointed simplicial $G$-sets as objects and all pointed simplicial morphisms as morphisms \pcref{def:gssetst}, on which $G$ acts by conjugation \cref{ginv_h_g}.  
\item\label{def:Gssetst_fgssetgst}
The category $\FGssetgst$ has pointed $G$-functors
\begin{equation}\label{fgssetgst_obj}
(\FG,\ord{0}) \to (\ssetgst,*)
\end{equation}
as objects, called \index{FG-simplicial G-set@$\FG$-simplicial $G$-set}\emph{$\FG$-simplicial $G$-sets}, and $G$-natural transformations as morphisms.\defmark
\end{enumerate}
\end{definition}

With $\Topgst$ replaced by $\ssetgst$, the description of $\FGTopg$ in \cref{expl:ggtopg} also applies to $\FGssetgst$.  For \cref{reast_proper}, recall that the classifying space functor $\cla$ \cref{classifying_space} is the composite of the nerve functor $\Ner$ and realization $\Rea$ \cref{x_dot}.  Thus, the induced functor $\clast$ \cref{clast_ggcatg} factors as the composite
\begin{equation}\label{clast_factors}
\begin{tikzpicture}[vcenter]
\def\h{2.5}
\draw[0cell]
(0,0) node (a1) {\FGCatg}
(a1)++(\h,0) node (a2) {\FGssetgst} 
(a2)++(\h,0) node (a3) {\FGTopg}
;
\draw[1cell=.9]
(a1) edge node {\Nerst} (a2)
(a2) edge node {\Reast} (a3)
;
\draw[1cell=.9]
(a1) [rounded corners=2pt] |- ($(a2)+(-1,.6)$) -- node {\clast} ($(a2)+(1,.6)$) -| (a3)
 ;
\end{tikzpicture}
\end{equation}
where $\Nerst$ and $\Reast$ postcompose with, respectively, $\Ner$ and $\Rea$ \pcref{expl:ggcatg_ggtopg}.  \cref{reast_proper} proves that the image of the functor $\Reast$ is always proper \pcref{def:proper_fgg}.

\begin{lemma}\label{reast_proper}
For each $\FG$-simplicial $G$-set $X$ \cref{fgssetgst_obj}, the $\FGG$-space $|X|_*$ is proper.
\end{lemma}

\begin{proof}
We need to prove that, for each pointed finite simplicial $G$-set $B_\crdot$ (\cref{def:gssetst} \eqref{def:gssetst_finite}), the simplicial $G$-space $(\Pwg|X|_*)_{B_\crdot}$ is Reedy cofibrant \pcref{def:reedy_cof}.  By \cref{reedy_cofibrant} \eqref{reedy_cof_ii}, it suffices to prove that each of its degeneracy morphisms is a $G$-cofibration.  A typical degeneracy morphism $s_i$ for $0 \leq i \leq r$ is displayed along the top of the commutative diagram \cref{reast_degen} of pointed $G$-morphisms.
\begin{equation}\label{reast_degen}
\begin{tikzpicture}[vcenter]
\def\h{4} \def\u{-1} \def\v{-1.3} 
\draw[0cell]
(0,0) node (a11) {(\Pwg |X|_*)_{B_r}}
(a11)++(\h,0) node (a12) {(\Pwg |X|_*)_{B_{r+1}}}
(a11)++(0,\u) node (a21) {\txint^{\mal \in \FG} B_r^{\mal} \sma |X\mal|}
(a12)++(0,\u) node (a22) {\txint^{\mal \in \FG} B_{r+1}^{\mal} \sma |X\mal|}
(a21)++(0,\v) node (a31) {\big| \txint^{\mal \in \FG} B_r^{\mal} \sma X\mal\big|}
(a22)++(0,\v) node (a32) {\big| \txint^{\mal \in \FG} B_{r+1}^{\mal} \sma X\mal\big|}
(a31)++(\h/2,\u) node (a4) {\phantom{A}}
;
\draw[1cell=.9]
(a11) edge node {s_i} (a12)
(a11) edge[equal] (a21)
(a12) edge[equal] (a22)
(a21) edge node[swap] {\iso} (a31)
(a22) edge node {\iso} (a32)
(a31) [rounded corners=2pt] |- ($(a4)+(-1,0)$) -- node {\big|\txint s_i^{\mal} \sma 1 \big|} ($(a4)+(1,0)$) -| (a32)
;
\end{tikzpicture}
\end{equation}
\begin{itemize}
\item The two equalities in \cref{reast_degen} and the notation \cref{A_to_mal}
\[(-)^{\mal} = \Topgst(\mal,-)\]
are from the definition of $\Pwg$ \cref{Pwgxa}.  Along the bottom, $s_i \cn B_r \to B_{r+1}$ is the $i$-th degeneracy morphism of $B_\crdot$.
\item The two natural isomorphisms follow from the fact that realization $\Rea$ \cref{x_dot} commutes with each smash product $B_k^{\mal} \sma -$ for $k \in \{r,r+1\}$, which is a left adjoint, and with the coend $\txint^{\mal \in \FG}$.
\end{itemize}

Realization $\Rea$ sends each injective morphism between pointed simplicial $G$-sets to the inclusion of a relative $G$-CW complex.  A relative $G$-CW inclusion is a $G$-cofibration, just as it is the case nonequivariantly \cite[Ch.\ 10.3]{mayconcise}.  Thus, to prove that the top $s_i$ in \cref{reast_degen} is a $G$-cofibration, it suffices to prove that the simplicial morphism
\[\txint^{\mal \in \FG} B_r^{\mal} \sma X\mal \fto{\txint s_i^{\mal} \sma 1} 
\txint^{\mal \in \FG} B_{r+1}^{\mal} \sma X\mal\]
is injective.  The preceding simplicial morphism is injective because the composite
\[B_r \fto{s_i} B_{r+1} \fto{d_i} B_r\]
is the identity morphism by one of the simplicial identities \cref{simplicial_id}.  More explicitly, the composite
\[\txint^{\mal \in \FG} B_r^{\mal} \sma X\mal 
\fto{(\txint d_i^{\mal} \sma 1)(\txint s_i^{\mal} \sma 1)}
\txint^{\mal \in \FG} B_r^{\mal} \sma X\mal \]
is the identity morphism, so $\txint s_i^{\mal} \sma 1$ is injective.
\end{proof}

\subsection*{Comparison Theorem}

\cref{thm:ksho_khsho} proves that Shimakawa $K$-theory \cref{ksho_kshosg} 
\begin{equation}\label{Ksho_compdiag}
\begin{tikzpicture}[vcenter]
\def\h{2} \def\u{.6}
\draw[0cell=.9]
(0,0) node (a1) {\phantom{\AlglaxO}}
(a1)++(0,-.04) node (a1') {\AlglaxO}
(a1)++(\h,0) node (a2) {\FGCatg}
(a2)++(\h,0) node (a3) {\FGTopg}
(a3)++(\h,0) node (a4) {\Gspec} 
;
\draw[1cell=.85]
(a1) edge node {\Sgo} (a2)
(a2) edge node {\clast} (a3)
(a3) edge node {\Kfg} (a4)
(a1) [rounded corners=2pt, shorten <=-.2ex] |- ($(a2)+(0,\u)$) -- node {\Ksho} ($(a3)+(0,\u)$) -| (a4)
;
\end{tikzpicture}
\end{equation}
and its homotopical variant \cref{khsho_khshosg}
\begin{equation}\label{Khsho_compdiag}
\begin{tikzpicture}[vcenter]
\def\h{2} \def\u{.6}
\draw[0cell=.9]
(0,0) node (a1) {\phantom{\AlglaxO}}
(a1)++(0,-.04) node (a1') {\AlglaxO}
(a1)++(\h,0) node (a2) {\FGCatg}
(a2)++(\h,0) node (a3) {\FGTopg}
(a3)++(\h,0) node (a4) {\FGTopg}
(a4)++(\h,0) node (a5) {\Gspec} 
;
\draw[1cell=.85]
(a1) edge node {\Sgo} (a2)
(a2) edge node {\clast} (a3)
(a3) edge node {\Bc} (a4)
(a4) edge node {\Kfg} (a5)
(a1) [rounded corners=2pt, shorten <=-.2ex] |- ($(a2)+(0,\u)$) -- node {\Khsho} ($(a4)+(0,\u)$) -| (a5)
;
\end{tikzpicture}
\end{equation}
are naturally componentwise weakly $G$-equivalent.  The strong variant involving the functors
\begin{equation}\label{KshosgKhshosg_diag}
\begin{tikzpicture}[vcenter]
\draw[0cell]
(0,0) node (a1) {\AlgpspsO}
(a1)++(4.3,0) node (a2) {\Gspec}
;
\draw[1cell=.9]
(a1) edge[transform canvas={yshift=.5ex}] node {\Kshosg = \Kfg \clast \Sgosg} (a2)
(a1) edge[transform canvas={yshift=-.4ex}] node[swap] {\Khshosg = \Kfg \Bc \clast \Sgosg} (a2)
;
\end{tikzpicture}
\end{equation}
is also true.  The only difference between $\Ksho$ and $\Khsho$ is that $\Khsho$ involves the bar functor $\Bc$, while $\Ksho$ does not \pcref{expl:khsho}.  The natural componentwise weak $G$-equivalence connecting $\Ksho$ and $\Khsho$ is defined by the retraction \cref{ret_FG}
\[\Bc X \fto{\retn_X} X \inspace \FGTopg.\]

\begin{theorem}\label{thm:ksho_khsho}
Suppose $\Op$ is a 1-connected $\Gcat$-operad \cref{i_connected} for a finite group $G$, and $\A$ is an $\Op$-pseudoalgebra \pcref{def:pseudoalgebra}.  Then the $G$-morphisms in $\Gspec$ \pcref{def:gsp_morphism}
\begin{equation}\label{thm:ksho_khsho_i}
\begin{split}
\Khsho\A & \fto{\Kfg\retn_{\clast\Sgo\A}} \Ksho\A \andspace\\
\Khshosg\A & \fto{\Kfg\retn_{\clast\Sgosg\A}} \Kshosg\A
\end{split}
\end{equation}
are componentwise weak $G$-equivalences, and they are natural in $\A$.
\end{theorem}

\begin{proof}
We consider $\Kfg\retn_{\clast\Sgo\A}$.  The proof for the strong variant is the same after replacing Shimakawa $H$-theory $\Sgo$ by the strong variant $\Sgosg$ \pcref{def:sgo}.  In the rest of this proof, $\retn_{\clast\Sgo\A}$ is abbreviated to $\retn$.
\begin{description}
\item[Naturality]  
The naturality of $\Kfg\retn$ in $\A$ follows from the naturality of the retraction $\retn$ \cref{retn_barFG_id} and the functoriality of $\Sgo$ \pcref{Sgo_twofunctor}, $\clast$ \pcref{ggcatg_ggtopg}, and $\Kfg$ \pcref{def:Kfg_functor}.
\item[Weak $G$-equivalence]  
For each $\FGG$-space $X$, the retraction $\retn_X$ \cref{ret_FG} is componentwise a weak $G$-equivalence.  By \cref{bar_proper,reast_proper}, the domain and the codomain of the retraction
\[\Bc\clast\Sgo\A \fto{\retn} \clast\Sgo\A \inspace \FGTopg\]
are proper $\FGG$-spaces.  Applying \cref{thm:Kfg_inv} to $\retn$ implies that $\Kfg\retn$ is componentwise a weak $G$-equivalence.\qedhere
\end{description}
\end{proof}

\begin{example}\label{ex:kshgbe_khshgbe}
\cref{thm:ksho_khsho} applies to the $G$-Barratt-Eccles operad $\GBE$ for a finite group $G$ \pcref{def:GBE}.  Its pseudoalgebras are genuine symmetric monoidal $G$-categories \pcref{def:GBE_pseudoalg}.
\end{example}

\begin{explanation}[Shimakawa's Construction]\label{expl:k_shimakawa}
Recall from \cref{expl:khsho_shi} that Shimakawa's original equivariant $K$-theory machine \cite{shimakawa89} is naturally weakly $G$-equivalent to the homotopical Shimakawa strong $K$-theory \cref{khsho_khshosg} 
\begin{equation}\label{Khshbesg_compdiag}
\begin{tikzpicture}[vcenter]
\def\h{2} \def\u{.6}
\draw[0cell=.9]
(0,0) node (a1) {\AlgpspsBE}
(a1)++(\h,0) node (a2) {\FGCatg}
(a2)++(\h,0) node (a3) {\FGTopg}
(a3)++(\h,0) node (a4) {\FGTopg}
(a4)++(\h,0) node (a5) {\Gspec} 
;
\draw[1cell=.85]
(a1) edge node {\Sgosgbe} (a2)
(a2) edge node {\clast} (a3)
(a3) edge node {\Bc} (a4)
(a4) edge node {\Kfg} (a5)
(a1) [rounded corners=2pt] |- ($(a2)+(0,\u)$) -- node {\Khshbesg} ($(a4)+(0,\u)$) -| (a5)
;
\end{tikzpicture}
\end{equation}
for the Barratt-Eccles operad $\BE$ \pcref{def:BE}.  \cref{thm:ksho_khsho} with $\Op = \BE$ implies that Shimakawa's original equivariant $K$-theory machine is naturally weakly $G$-equivalent to Shimakawa strong $K$-theory \cref{ksho_kshosg}
\begin{equation}\label{Kshbesg_compdiag}
\begin{tikzpicture}[vcenter]
\def\h{2} \def\u{.6}
\draw[0cell=.9]
(0,0) node (a1) {\AlgpspsBE}
(a1)++(\h,0) node (a2) {\FGCatg}
(a2)++(\h,0) node (a3) {\FGTopg}
(a3)++(\h,0) node (a4) {\Gspec} 
;
\draw[1cell=.85]
(a1) edge node {\Sgosgbe} (a2)
(a2) edge node {\clast} (a3)
(a3) edge node {\Kfg} (a4)
(a1) [rounded corners=2pt] |- ($(a2)+(0,\u)$) -- node {\Kshbesg} ($(a3)+(0,\u)$) -| (a4)
;
\end{tikzpicture}
\end{equation}
for $\BE$.
\end{explanation}

\begin{remark}\label{rk:ksho_khsho}
The retraction $\retn_X$ is componentwise a pointed $G$-homotopy equivalence, as discussed under \cref{ret_FG}.  However, \cref{thm:ksho_khsho} does \emph{not} claim that either $\Kfg\retn$ is componentwise a pointed $G$-homotopy equivalence.  This stronger statement involving pointed $G$-homotopy equivalences is not known to be true.  The reason is that the componentwise $G$-homotopy inverse of $\retn$, given by the section $\secn_{\mal}$ \cref{secn_c}, is \emph{not} natural in the object $\mal \in \FG$, as explained in \cref{secn_nonnat}.
\end{remark}

%% file: chap/shim_top.tex
This chapter compares Shimakawa's homotopical functor \pcref{def:Kfg_functor,def:gtop_retn_secn}
\[\FGTopg \fto{\Bc} \FGTopg \fto{\Kfg} \Gspec\]
and the prolongation functor \pcref{def:Kgg_functor}
\[\GGTopg \fto{\Kgg} \Gspec.\]
The composite $\Kfg\Bc$ is the last part of the homotopical Shimakawa (strong) $K$-theory \cref{khsho_khshosg}.  The functor $\Kgg$ is the last step of our equivariant $K$-theory \cref{Khgo_functors}.  The main comparison results of this chapter are as follows.
\begin{enumerate}
\item For a finite group $G$ and proper objects, $\Kfg\Bc$ and $\Kgg$ yield naturally componentwise weakly $G$-equivalent orthogonal $G$-spectra.  See \cref{kfgbkgg_compare,kggkfgb_compare}. 
\item For a compact Lie group $G$ and without the bar functor $\Bc$, the functors $\Kfg$ and $\Kgg$ factor through each other up to natural isomorphisms.  See \cref{kfgkgg_compare}.  Thus, any orthogonal $G$-spectrum produced by Shimakawa's functor $\Kfg$ can also be produced by $\Kgg$ and vice versa.
\end{enumerate}

\summary
The following diagram summarizes the comparison between $\Kfg\Bc$ and $\Kgg$.
\begin{equation}\label{KfgBKgg_chi}
\begin{tikzpicture}[vcenter]
\def\t{25} \def\h{2.6} \def\v{-1.4} \def\u{.7} \def\c{.7em}
\draw[0cell]
(0,0) node (a11) {\FGTopg}
(a11)++(\h,0) node (a12) {\GGTopg}
(a12)++(\h,0) node (a13) {\FGTopg}
(a11)++(0,\v) node (a21) {\FGTopg}
(a12)++(0,\v) node (a22) {\Gspec}
(a13)++(0,\v) node (a23) {\FGTopg}
;
\draw[1cell=.9]
(a11) edge node {\ifgl} (a12)
(a12) edge node {\Kgg} (a22)
(a21) edge node[swap] {\Kfg} (a22)
(a11) edge[transform canvas={xshift=-\c}, bend right=\t] node[swap] {\Bc} (a21)
(a11) edge[transform canvas={xshift=\c}, bend left=\t] node {1} (a21)
(a11) [rounded corners=2pt] |- ($(a12)+(-1,\u)$) -- node {1} node[swap, inner sep=5pt] {\iso} ($(a12)+(1,\u)$) -| (a13)
;
\draw[2cell=1]
node[between=a11 and a21 at .55, 2label={above,\retn\phantom{x}}] {\Rightarrow}
node[between=a11 and a21 at .55, shift={(.55*\h,0)}] {\iso}
;
\draw[1cell=.9]
(a12) edge node {\smal} (a13)
(a13) edge[transform canvas={xshift=\c}, bend left=\t] node {\Bc} (a23)
(a13) edge[transform canvas={xshift=-\c}, bend right=\t] node[swap] {1} (a23)
(a23) edge node[pos=.4] {\Kfg} (a22)
;
\draw[2cell=1]
node[between=a13 and a23 at .51, rotate=180, 2labelmed={below,\retn}] {\Rightarrow}
node[between=a13 and a23 at .55, shift={(-.55*\h,0)}] {\iso}
;
\end{tikzpicture}
\end{equation}
This diagram exists for any compact Lie group $G$.  In the left and right regions, for each $\FGG$-space, the retraction $\retn \cn \Bc \to 1$ \cref{retn_barFG_id} is componentwise a pointed $G$-homotopy equivalence.  The functors $\ifgl$ and $\smal$ are the left adjoints induced by the length-1 inclusion functor $\ifg$ and the smash functor $\sma$ \pcref{def:ifgGG,def:smashFGGG}:
\[\FG \fto{\ifg} \GG \fto{\sma} \FG.\]
The other three regions commute up to natural isomorphisms:
\[\smal\ifgl \iso 1, \quad \Kfg \iso \Kgg\ifgl, \andspace \Kgg \iso \Kfg\!\smal\!.\]
For a \emph{finite} group $G$, a proper $\FGG$-space $X$ \pcref{def:proper_fgg}, and a proper $\GGG$-space $Y$ \pcref{def:proper_ggg}, the $G$-morphisms of orthogonal $G$-spectra
\[\begin{split}
& \Kfg\Bc X \fto{\Kfg\retn_X} \Kfg X \fiso \Kgg\ifgl X \andspace \\
& \Kfg\Bc\!\smal\! Y \fto{\Kfg\retn_{\smal Y}} \Kfg\!\smal\! Y \fiso \Kgg Y
\end{split}\]
are componentwise weak $G$-equivalences.  Thus, $\Kfg\Bc X$ can be computed as $\Kgg \ifgl X$.  Conversely, $\Kgg Y$ can be computed as $\Kfg\Bc\!\smal\! Y$.

\organization
This chapter consists of the following sections.

\secname{sec:adj_ptd}
This section constructs in detail the adjunction
\begin{equation}\label{fadjunction}
\begin{tikzpicture}[vcenter]
\draw[0cell]
(0,0) node (a1) {\CTopgst}
(a1)++(2.4,0) node (a2) {\DTopgst}
;
\draw[1cell=.9]
(a1) edge[transform canvas={yshift=.6ex}] node {\funl} (a2)
(a2) edge[transform canvas={yshift=-.4ex}] node {\funst} (a1)
;
\end{tikzpicture}
\end{equation}
induced by a pointed $G$-functor $\fun \cn \C \to \D$ between indexing $G$-categories.  The right adjoint $\funst$ is the pullback functor along $\fun$.  The left adjoint $\funl$ requires some care related to basepoints and $G$-equivariance, since morphisms in $\Topgst$ are not required to be $G$-equivariant.  The functors $\ifgl$ and $\smal$ are examples of $\funl$.

\secname{sec:gsp_ptdiag}
This section constructs the prolongation functor
\[\CTopgst \fto{\KC} \Gspec,\]
from the category of $\CG$-spaces to the category of orthogonal $G$-spectra, for each indexing $G$-category $\C$ equipped with a pointed $G$-functor to $\FG$.  Shimakawa's functor $\Kfg$ and the prolongation functor $\Kgg$ are examples of $\KC$.  The functor $\KC$ is natural in the indexing $G$-category $\C$ \pcref{KC_funl}.

\secname{sec:KfgKgg_comp}
Using the results in \cref{sec:gsp_ptdiag}, this section compares Shimakawa's functors $\Kfg$ and $\Kfg\Bc$ with the functor $\Kgg$, as discussed in the summary.  See \cref{kfgkgg_compare,kfgbkgg_compare,kggkfgb_compare}.

\section{Pointed Equivariant Diagrams}
\label{sec:adj_ptd}

For an arbitrary group $G$, this section constructs the change-of-shape adjunction
\begin{equation}\label{fadj_seci}
\begin{tikzpicture}[vcenter]
\draw[0cell]
(0,0) node (a1) {\CTopgst}
(a1)++(2.4,0) node (a2) {\DTopgst}
;
\draw[1cell=.9]
(a1) edge[transform canvas={yshift=.6ex}] node {\funl} (a2)
(a2) edge[transform canvas={yshift=-.4ex}] node {\funst} (a1)
;
\end{tikzpicture}
\end{equation}
associated to each pointed $G$-functor $\fun \cn \C \to \D$ between indexing $G$-categories.  \cref{sec:gsp_ptdiag} constructs the prolongation functor
\[\CTopgst \fto{\KC} \Gspec\]
that produces orthogonal $G$-spectra.  The left adjoint $\funl$ is used in \cref{KC_funl} to prove that $\KC$ is natural in the indexing $G$-category $\C$.

\secoutline
\begin{itemize}
\item \cref{def:ctopgst} defines indexing $G$-categories and the pointed equivariant diagram category $\CTopgst$.  
\item \cref{ex:ctopgst} observes that the categories $\GGTopg$ and $\FGTopg$ are examples of $\CTopgst$.
\item \cref{fun_adj} constructs the adjunction $(\funl,\funst)$.
\item \cref{ex:smai_adj} applies \cref{fun_adj} to the length-1 inclusion functor $\ifg \cn \FG \to \GG$ and the smash functor $\sma \cn \GG \to \FG$. 
\end{itemize}

Recall the pointed $G$-category $\Topgst$ of pointed $G$-spaces and pointed morphisms with the conjugation $G$-action \pcref{def:GCat,def:ptGcat,def:gtopst}.  

\begin{definition}\label{def:ctopgst}
Suppose $G$ is a group.
\begin{enumerate}
\item\label{def:ctopgst_zero} An \emph{indexing $G$-category}\index{indexing G-category@indexing $G$-category}\index{G-category@$G$-category!indexing} is a small pointed $G$-category $(\C,*)$ with an initial-terminal basepoint $*$ such that $G$ acts trivially on objects, meaning $gc = c$ for each $g \in G$ and object $c \in \C$. 
\item\label{def:ctopgst_i} For an indexing $G$-category $\C$ and objects $c,c' \in \C$, the \emph{0-morphism}\index{0-morphism} $c \to c'$ is the unique morphism 
\[0 \cn c \to * \to c'\] 
that factors through the initial-terminal basepoint $*$.  The set of nonzero morphisms is denoted by
\begin{equation}\label{Cpunc}
\Cpunc(c,c') = \C(c,c') \setminus \{0\}.
\end{equation}
\item\label{def:ctopgst_ii} For an indexing $G$-category $(\C,*)$, define the category $\CTopgst$ with pointed $G$-functors
\begin{equation}\label{ctopgst_obj}
(\C,*) \fto{X} (\Topgst,*)
\end{equation}
as objects, called \index{G-space@$G$-space!diagram}\emph{$\CG$-spaces}, and $G$-natural transformations
\begin{equation}\label{ctopgst_mor}
\begin{tikzpicture}[vcenter]
\def\t{28}
\draw[0cell]
(0,0) node (a1) {\C}
(a1)++(1.8,0) node (a2) {\phantom{\Gsk}}
(a2)++(.2,-.03) node (a2') {\Topgst}
;
\draw[1cell=.9]
(a1) edge[bend left=\t] node {X} (a2)
(a1) edge[bend right=\t] node[swap] {X'} (a2)
;
\draw[2cell]
node[between=a1 and a2 at .43, rotate=-90, 2label={above,\tha}] {\Rightarrow}
;
\end{tikzpicture}
\end{equation}
as morphisms.  Identities and composition are defined componentwise in $\Topgst$.\defmark
\end{enumerate}
\end{definition}

\begin{example}\label{ex:ctopgst}\
\begin{enumerate}
\item The categories $\FG$ and $\GG$ \pcref{def:FG,def:GG} are indexing $G$-categories in the sense of \cref{def:ctopgst} \cref{def:ctopgst_zero}. 
\item The categories $\GGTopg$ and $\FGTopg$ \pcref{def:GGTopg,def:ggtopg} are examples of $\CTopgst$ in the sense of \cref{def:ctopgst} \cref{def:ctopgst_ii}.  The explicit description of $\GGTopg$ in \cref{expl:GGTopg} also applies to $\CTopgst$, after replacing the indexing $G$-category $(\GG,\vstar)$ with $(\C,*)$.\defmark
\end{enumerate}
\end{example}

Each pointed $G$-functor $\fun \cn \C \to \D$ between indexing $G$-categories induces a pullback functor
\begin{equation}\label{funst_pullback}
\CTopgst \fot{\funst} \DTopgst
\end{equation}
given by precomposition with $\fun$:
\begin{equation}\label{ftheta_comp}
\begin{tikzpicture}[vcenter]
\def\t{28}
\draw[0cell]
(0,0) node (a1) {\D}
(a1)++(-1.5,0) node (a0) {\C}
(a1)++(1.8,0) node (a2) {\phantom{\Gsk}}
(a2)++(.25,-.03) node (a2') {\Topgst.}
;
\draw[1cell=.9]
(a0) edge node {\fun} (a1)
(a1) edge[bend left=\t] node {X} (a2)
(a1) edge[bend right=\t] node[swap] {X'} (a2)
;
\draw[2cell]
node[between=a1 and a2 at .43, rotate=-90, 2label={above,\tha}] {\Rightarrow}
;
\end{tikzpicture}
\end{equation}
\cref{fun_adj} constructs the left adjoint $\funl$ of $\funst$, which is essentially a left Kan extension.  However, there are some subtleties related to basepoints and $G$-equivariance.  In particular, $\funl$ involves coends, but morphisms in $\Topgst$ are \emph{not} necessarily $G$-equivariant. 

\begin{lemma}\label{fun_adj}
For a group $G$, suppose $\fun \cn \C \to \D$ is a pointed $G$-functor between indexing $G$-categories.  Then there is an adjunction
\begin{equation}\label{fadj_statement}
\begin{tikzpicture}[vcenter]
\draw[0cell]
(0,0) node (a1) {\CTopgst}
(a1)++(2.4,0) node (a2) {\DTopgst}
;
\draw[1cell=.9]
(a1) edge[transform canvas={yshift=.6ex}] node {\funl} (a2)
(a2) edge[transform canvas={yshift=-.4ex}] node {\funst} (a1)
;
\end{tikzpicture}
\end{equation}
with right adjoint given by the pullback functor $\funst$ \cref{funst_pullback}. 
\end{lemma}

\begin{proof}
We first construct 
\begin{itemize}
\item the left adjoint $\funl$ \cref{funl},
\item the unit $\funu \cn 1 \to \funst\funl$ \cref{funu}, and
\item the counit $\funv \cn \funl\funst \to 1$ \cref{funv}.
\end{itemize}
We verify the two triangle identities for an adjunction in \cref{fun_lefttri,fun_righttri}.
\begin{description}
\item[Left adjoint]  
The functor
\begin{equation}\label{funl}
\CTopgst \fto{\funl} \DTopgst
\end{equation}
sends a pointed $G$-functor $X \cn \C \to \Topgst$ \cref{ctopgst_obj} to the pointed $G$-functor
\[\D \fto{\funl X} \Topgst\]
whose value at an object $d \in \D$ is the coend 
\begin{equation}\label{funlx_d}
(\funl X)d = \ecint^{c \in \C} \bigvee_{\Dpun(\fun c,\, d)} Xc
\end{equation}
taken in $\Topst$.  Its $G$-action is defined in \cref{funlx_gaction}.
\begin{description}
\item[Representatives]
In the wedge index in \cref{funlx_d}, $\Dpun(\fun c, d)$ is the set of nonzero morphisms $\fun c \to d$ \cref{Cpunc}.  An empty wedge---which happens, for example, if $d = * \in \D$---is defined as a point $*$.  The coend in \cref{funlx_d} is a quotient of the wedge
\[\bigvee_{c \in \C} \bigvee_{\Dpun(\fun c,\, d)} Xc.\]
Each point of $(\funl X)d$ is represented by a pair
\begin{equation}\label{funlx_reps}
(\upom; x) \in \D(\fun c, d) \ttimes Xc.
\end{equation}
The pair $(\upom; x)$ represents the basepoint if $x \in Xc$ is the basepoint or if $\upom$ is the 0-morphism, factoring through the basepoint $* \in \D$. 
\item[Relations] 
The defining relation of the coend $(\funl X)d$ identifies, for each triple
\[(\upom; \uppsi; x) \in \D(\fun c, d) \ttimes \C(c', c) \ttimes Xc',\]
the two pairs
\begin{equation}\label{funlx_relations}
\begin{split}
\big(\upom(\fun\uppsi); x\big) & \in \D(\fun c', d) \ttimes Xc' \andspace\\
\big(\upom; (X\uppsi)x\big) & \in \D(\fun c, d) \ttimes Xc.
\end{split}
\end{equation}
\item[Morphism assignment of $\funl X$] 
For a morphism $\upla \cn d \to d'$ in $\D$, the pointed morphism 
\[(\funl X)d \fto{(\funl X)\upla} (\funl X)d'\]
sends a representative pair $(\upom; x)$ in $(\funl X)d$ \cref{funlx_reps} to the representative pair
\begin{equation}\label{funlx_mor}
\big((\funl X)\upla\big) (\upom; x) = (\upla\upom ; x) \in \D(\fun c, d') \ttimes Xc
\end{equation}
in $(\funl X)d'$.  The pointed morphism $(\funl X)\upla$ is well defined because it sends the two representative pairs in \cref{funlx_relations} to the same point in $(\funl X)d'$.  The pointed functoriality of $\funl X$ follows from \cref{funlx_mor}.  Its $G$-equivariance is proved in \cref{funlx_geq}.
\item[$G$-action] 
The group $G$ acts diagonally on representatives \cref{funlx_reps} of the pointed space $(\funl X)d$:
\begin{equation}\label{funlx_gaction}
g (\upom ; x) = (g \upom ; gx) \in \D(\fun c, d) \ttimes Xc
\end{equation}
for $g \in G$, $\upom \in \D(\fun c, d)$, and $x \in Xc$.  Note that $g\upom$ is a morphism $\fun c \to d$ because $\D$ is an indexing $G$-category (\cref{def:ctopgst} \cref{def:ctopgst_zero}), where $G$ acts trivially on objects.  The following equalities in $(\funl X)d$ prove that the $G$-actions on the two representative pairs in \cref{funlx_relations} are the same.
\[\begin{aligned}
& g \big(\upom(\fun\uppsi); x\big) &&\\
&= \big(g (\upom(\fun\uppsi)) ; gx \big) && \text{by \cref{funlx_gaction}}\\
&= \big( (g \upom) (g \fun\uppsi) ; gx \big) && \text{by functoriality of $G$-action}\\
&= \big( (g \upom) \fun(g\uppsi) ; gx \big) && \text{by $G$-equivariance of $\fun$}\\
&= \big(g \upom ; X (g\uppsi) (gx) \big) && \text{by \cref{funlx_relations}}\\
&= \big(g \upom ; (g (X\uppsi) \ginv) (gx) \big) && \text{by $G$-equivariance of $X$}\\
&= \big(g \upom ; g ((X\uppsi)x) \big) && \text{by $\ginv g = 1$} \\
&= g \big(\upom; (X\uppsi)x\big) && \text{by \cref{funlx_gaction}}
\end{aligned}\]
Thus, the $G$-action on the pointed space $(\funl X)d$ is well defined.
\item[$G$-equivariance of $\funl X$]  
The $G$-equivariance of $\funl X$ means that, for each $g \in G$ and each morphism $\upla \cn d \to d'$ in $\D$, there is an equality of pointed morphisms
\begin{equation}\label{funlx_gequiv}
(\funl X)d \fto[= g \big( (\funl X)\upla \big) \ginv]{(\funl X)(g \upla)} (\funl X)d'.
\end{equation}
The following equalities in $(\funl X)d'$ prove that the two morphisms in \cref{funlx_gequiv} are the same on each representative pair $(\upom; x)$ in $(\funl X)d$ \cref{funlx_reps}.
\begin{equation}\label{funlx_geq}
\begin{aligned} 
& \big( (\funl X)(g \upla) \big) (\upom; x) &&\\
&= \big( (g \upla)\upom ; x\big) && \text{by \cref{funlx_mor}} \\
&= \big( g (\upla(\ginv \upom)) ; g \ginv x\big) && \text{by functoriality of $G$-action}\\
&= g \big(\upla(\ginv \upom) ; \ginv x \big) && \text{by \cref{funlx_gaction}}\\
&= \big( g ((\funl X)\upla)\big) (\ginv \upom ; \ginv x) && \text{by \cref{funlx_mor}}\\
&= \big( g ( (\funl X)\upla ) \ginv \big) (\upom; x) && \text{by \cref{funlx_gaction}}
\end{aligned}
\end{equation}
Thus, $\funl X \cn \D \to \Topgst$ is a pointed $G$-functor.
\item[Morphism assignment of $\funl$] 
For a $G$-natural transformation $\tha \cn X \to X'$ in $\CTopgst$ \cref{ctopgst_mor}, the $G$-natural transformation
\begin{equation}\label{ftha_iicell}
\begin{tikzpicture}[vcenter]
\def\t{28}
\draw[0cell]
(0,0) node (a1) {\D}
(a1)++(2,0) node (a2) {\phantom{\Gsk}}
(a2)++(.25,0) node (a2') {\Topgst}
;
\draw[1cell=.9]
(a1) edge[bend left=\t] node {\funl X} (a2)
(a1) edge[bend right=\t] node[swap] {\funl X'} (a2)
;
\draw[2cell]
node[between=a1 and a2 at .4, rotate=-90, 2label={above,\funl\tha}] {\Rightarrow}
;
\end{tikzpicture}
\end{equation}
has, for each object $d \in \D$, $d$-component pointed morphism
\[(\funl X)d \fto{(\funl\tha)_d} (\funl X')d\]
that sends a representative pair $(\upom; x)$ in $(\funl X)d$ \cref{funlx_reps} to the representative pair
\begin{equation}\label{funl_tha_d}
(\funl\tha)_d (\upom; x) = (\upom; \tha_c x) \in \D(\fun c, d) \ttimes X'c
\end{equation}
in $(\funl X')d$. 
\begin{itemize}
\item The pointed morphism $(\funl \tha)_d$ is well defined because it sends the two representative pairs in \cref{funlx_relations} to the same point in $(\funl X')d$ by the naturality of $\tha$. 
\item The naturality of $(\funl\tha)_d$ in $d \in \D$ follows from \cref{funlx_mor,funl_tha_d}. 
\item The $G$-equivariance of $\funl\tha$ means the $G$-equivariance of each $(\funl\tha)_d$, which follows from \cref{funlx_gaction}, \cref{funl_tha_d}, and the $G$-equivariance of the $c$-component morphism $\tha_c \cn Xc \to X'c$.
\item The functoriality of $\funl$ follows from \cref{funl_tha_d} and the fact that identities and composition in $\CTopgst$ and $\DTopgst$ are defined componentwise.
\end{itemize}
\end{description}
This finishes the construction of the functor $\funl$ \cref{funl}.
\item[Unit]
The unit of the adjunction $(\funl,\funst)$ is the natural transformation
\begin{equation}\label{funu}
1_{\CTopgst} \fto{\funu} \funst\funl
\end{equation}
that sends a pointed $G$-functor $X \cn \C \to \Topgst$ to the $G$-natural transformation
\begin{equation}\label{funux}
\begin{tikzpicture}[vcenter]
\def\t{28}
\draw[0cell]
(0,0) node (a1) {\C}
(a1)++(2,0) node (a2) {\phantom{\Gsk}}
(a2)++(.25,0) node (a2') {\Topgst.}
;
\draw[1cell=.9]
(a1) edge[bend left=\t] node {X} (a2)
(a1) edge[bend right=\t] node[swap] {\funst\funl X} (a2)
;
\draw[2cell]
node[between=a1 and a2 at .4, rotate=-90, 2label={above,\funu_X}] {\Rightarrow}
;
\end{tikzpicture}
\end{equation}
Its component at an object $c' \in \C$ is the pointed morphism
\begin{equation}\label{funux_c}
Xc' \fto{\funu_{X,\, c'}} (\funst\funl X)c'
= \ecint^{c\in \C} \mhquad \bigvee_{\Dpun(\fun c,\, \fun c')} \mhquad Xc\\
\end{equation}
that sends a point $x \in Xc'$ to the representative pair \cref{funlx_reps}
\begin{equation}\label{funuxc_x}
\funu_{X,\, c'}(x) = \big(1_{\fun c'} ; x \big) \in \D(\fun c', \fun c') \ttimes Xc'
\end{equation}
in $(\funst\funl X)c' = (\funl X)(\fun c')$.
\begin{itemize}
\item The $G$-equivariance of $\funu_X$ \cref{funux} means the $G$-equivariance of each $\funu_{X,\, c'}$ \cref{funux_c}, which follows from \cref{funlx_gaction,funuxc_x}. 
\item The naturality of $\funu_{X,\, c'}$ in $c' \in \C$ follows from \cref{funlx_relations,funlx_mor,funuxc_x}. 
\item The naturality of $\funu_X$ in $X \in \CTopgst$ follows from \cref{funl_tha_d,funuxc_x}.
\end{itemize}
This finishes the construction of the unit $\funu$ \cref{funu}.
\item[Counit]
The counit of the adjunction $(\funl,\funst)$ is the natural transformation
\begin{equation}\label{funv}
\funl\funst \fto{\funv} 1_{\DTopgst}
\end{equation}
that sends a pointed $G$-functor $X \cn \D \to \Topgst$ to the $G$-natural transformation
\begin{equation}\label{funvx}
\begin{tikzpicture}[vcenter]
\def\t{28}
\draw[0cell]
(0,0) node (a1) {\D}
(a1)++(2,0) node (a2) {\phantom{\Gsk}}
(a2)++(.25,0) node (a2') {\Topgst.}
;
\draw[1cell=.9]
(a1) edge[bend left=\t] node {\funl\funst X} (a2)
(a1) edge[bend right=\t] node[swap] {X} (a2)
;
\draw[2cell]
node[between=a1 and a2 at .4, rotate=-90, 2label={above,\funv_X}] {\Rightarrow}
;
\end{tikzpicture}
\end{equation}
Its component at an object $d \in \D$ is the pointed morphism
\begin{equation}\label{funvx_d}
(\funl\funst X)d = \ecint^{c \in \C} \mhquad \bigvee_{\Dpun(\fun c,\, d)} \mhquad X(\fun c)
\fto{\funv_{X,\, d}} Xd
\end{equation}
that sends a representative pair 
\[(\upom; x) \in \D(\fun c, d) \ttimes X(\fun c)\] 
in $(\funl\funst X)d$ to the point
\begin{equation}\label{funvx_d_rep}
\funv_{X,\, d} (\upom; x) = (X\upom)x \in Xd.
\end{equation}
\begin{itemize}
\item The pointed morphism $\funv_{X,\, d}$ is well defined by \cref{funlx_relations}, \cref{funvx_d_rep}, and the functoriality of $X$. 
\item The naturality of $\funv_{X,\, d}$ in $d \in \D$ follows from \cref{funlx_mor}, \cref{funvx_d_rep}, and the functoriality of $X$.
\item The $G$-equivariance of $\funv_X$ means the $G$-equivariance of each $\funv_{X,\, d}$.  The following equalities in $Xd$ prove that $\funv_{X,\, d}$ is $G$-equivariant.
\[\begin{aligned}
& \funv_{X,\, d} \big(g (\upom; x) \big) &&\\
&= \funv_{X,\, d} (g \upom ; gx) && \text{by \cref{funlx_gaction}}\\
&= \big(X(g \upom)\big) (gx) && \text{by \cref{funvx_d_rep}}\\
&= (g (X\upom)) (\ginv gx) && \text{by $G$-equivariance of $X$}\\
&= g\big( (X\upom)x\big) && \text{by $\ginv g = 1$}\\
&= g\big( \funv_{X,\, d} (\upom; x) \big) && \text{by \cref{funvx_d_rep}}
\end{aligned}\]
Thus, $\funv_X$ \cref{funvx} is a $G$-natural transformation.
\item The naturality of $\funv_X$ in $X \in \DTopgst$ follows from \cref{funl_tha_d}, \cref{funvx_d_rep}, and the naturality of morphisms in $\DTopgst$ \cref{ctopgst_mor}.
\end{itemize}
This finishes the construction of the counit $\funv$ \cref{funv}.
\item[Left triangle identity]
This triangle identity states that, for each pointed $G$-functor $X \cn \C \to \Topgst$, the composite 
\begin{equation}\label{fun_lefttri}
\funl X \fto{\funl\funu_X} \funl\funst\funl X \fto{\funv_{\funl X}} \funl X
\end{equation}
is the identity $G$-natural transformation.  For an object $d \in \D$ and a representative pair $(\upom; x)$ in $(\funl X)d$ \cref{funlx_reps}, the following equalities in $(\funl X)d$ prove that the $d$-component of the composite in \cref{fun_lefttri} is the identity morphism.
\[\begin{aligned}
& (\funv_{\funl X,\, d}) (\funl\funu_X)_d (\upom; x) &&\\
&= (\funv_{\funl X,\, d}) \big(\upom; \funu_{X,\, c} x \big) && \text{by \cref{funl_tha_d}}\\
&= (\funv_{\funl X,\, d}) \big(\upom; (1_{\fun c}; x) \big) && \text{by \cref{funuxc_x}} \\
&= \big((\funl X)\upom\big) (1_{\fun c}; x) && \text{by \cref{funvx_d_rep}}\\
&= \big(\upom 1_{\fun c}; x \big) && \text{by \cref{funlx_mor}} \\
&= (\upom; x)
\end{aligned}\]
This proves the left triangle identity.
\item[Right triangle identity]
This triangle identity states that, for each pointed $G$-functor $X \cn \D \to \Topgst$, the composite 
\begin{equation}\label{fun_righttri}
\funst X \fto{\funu_{\funst X}} \funst\funl\funst X \fto{\funst\funv_X} \funst X
\end{equation}
is the identity natural transformation.  For an object $c \in \C$ and a point $x \in (\funst X)c = X(\fun c)$, the following equalities in $X(\fun c)$ prove that the $c$-component of the composite in \cref{fun_righttri} is the identity morphism.  
\[\begin{aligned}
& (\funst\funv_X)_c (\funu_{\funst X,\, c}) (x) && \\
&= (\funv_{X,\,\fun c}) \big(1_{\fun c}; x\big) && \text{by \cref{funuxc_x}} \\
&= (X 1_{\fun c}) x && \text{by \cref{funvx_d_rep}} \\
&= 1_{X(\fun c)} x && \text{by functoriality} \\
&= x &&
\end{aligned}\]
\end{description}
In summary, the quadruple $(\funl,\funst,\funu,\funv)$ is an adjunction.
\end{proof}

\begin{example}\label{ex:smai_adj}
For the indexing $G$-categories $\FG$ and $\GG$ \pcref{def:FG,def:GG}, the length-1 inclusion functor $\ifg$ and the smash functor $\sma$ \pcref{def:ifgGG,def:smashFGGG},
\begin{equation}\label{ismaone}
\begin{tikzpicture}[vcenter]
\def\u{.6} \def\h{1.7}
\draw[0cell]
(0,0) node (a1) {\FG}
(a1)++(\h,0) node (a2) {\GG}
(a2)++(\h,0) node (a3) {\FG,}
;
\draw[1cell=.9]
(a1) edge node {\ifg} (a2)
(a2) edge node {\sma} (a3)
(a1) [rounded corners=2pt] |- ($(a2)+(-1,\u)$) -- node {1} ($(a2)+(1,\u)$) -| (a3)
;
\end{tikzpicture}
\end{equation}
are pointed $G$-functors such that $\sma\ifg = 1_{\FG}$.  By \cref{fun_adj}, there are adjunctions $(\ifgl,\ifgst)$ and $(\smal,\smast)$, as displayed in the diagram
\begin{equation}\label{isma_adjoints}
\begin{tikzpicture}[vcenter]
\def\u{.7} \def\h{2.6} \def\b{.6ex} \def\c{-.4ex}
\draw[0cell]
(0,0) node (a1) {\FGTopg}
(a1)++(\h,0) node (a2) {\GGTopg}
(a2)++(\h,0) node (a3) {\FGTopg,}
;
\draw[1cell=.9]
(a1) edge[transform canvas={yshift=\b}] node {\ifgl} (a2)
(a2) edge[transform canvas={yshift=\b}] node {\smal} (a3)
(a3) edge[transform canvas={yshift=\c}] node {\smast} (a2)
(a2) edge[transform canvas={yshift=\c}] node {\ifgst} (a1)
(a3) [rounded corners=2pt] |- ($(a2)+(1,-\u)$) -- node[swap] {1} ($(a2)+(-1,-\u)$) -| (a1)
;
\draw[1cell=.9]
(a1) [rounded corners=2pt] |- ($(a2)+(-1,\u)$) -- node {1} ($(a2)+(1,\u)$) -| (a3)
;
\end{tikzpicture}
\end{equation}
such that
\[\begin{split}
\ifgst\smast &= (\sma\ifg)^* = 1_{\FGTopg} \andspace\\ 
\smal\ifgl &\iso 1_{\FGTopg}.
\end{split}\]
The previous natural isomorphism follows from the fact that the composite $\smal\ifgl$ is left adjoint to the identity functor on $\FGTopg$.
\begin{description}
\item[Left adjoint $\ifgl$]
Applying \cref{funlx_d} to $\ifg$, for a pointed $G$-functor $X \cn \FG \to \Topgst$, the pointed $G$-functor
\[\GG \fto{\ifgl X} \Topgst\]
sends an object $\nbe \in \GG$ to the pointed space
\begin{equation}\label{ifglx_nbe}
(\ifgl X)\nbe = \int^{\mal \in \FG} \mhquad \bigvee_{\GGpunc(\ifg\mal;\, \nbe)} \mhquad X\mal.
\end{equation}
\item[Left adjoint $\smal$]
Applying \cref{funlx_d} to $\sma$, for a pointed $G$-functor $Y \cn \GG \to \Topgst$, the pointed $G$-functor
\[\FG \fto{\smal Y} \Topgst\]
sends a pointed finite $G$-set $\mal \in \FG$ to the pointed space
\begin{equation}\label{smaly_mal}
(\smal Y)\mal = \int^{\nbe \in \GG} \mhquad \bigvee_{\FGpunc(\sma\nbe;\, \mal)} \mhquad Y\nbe.
\end{equation}
\end{description}
In \cref{ifglx_nbe,smaly_mal}, the group $G$ acts diagonally on representatives, as defined in \cref{funlx_gaction}.
\end{example}

\section{Orthogonal $G$-Spectra from Equivariant Diagrams}
\label{sec:gsp_ptdiag}

For a compact Lie group $G$, this section constructs the prolongation functor
\[\CTopgst \fto{\KC} \Gspec,\]
sending $\CG$-spaces \cref{ctopgst_obj} to orthogonal $G$-spectra (\cref{def:gsp_morphism} \cref{def:gsp_morphism_ii}), and observes that it is natural in the indexing $G$-category $\C$.  These results are used in \cref{sec:KfgKgg_comp} to compare Shimakawa's (homotopical) functors $\Kfg$ and $\Kfg\Bc$ with the prolongation functor $\Kgg$.

\secoutline
\begin{itemize}
\item \cref{def:KC} defines the functor $\KC$ for an indexing $G$-category $\C$ \pcref{def:ctopgst} equipped with a pointed $G$-functor to $\FG$.  Recall the indexing $G$-category $\FG$ of pointed finite $G$-sets and pointed functions with the conjugation $G$-action \pcref{def:FG}.
\item \cref{ex:KC} observes that $\KC$ generalizes both $\Kgg$ and $\Kfg$ \pcref{def:Kgg_functor,def:Kfg_functor}.
\item \cref{KC_funl} proves that $\KC$ is natural in the indexing $G$-category $\C$.  It is further elaborated in \cref{expl:kfgisoxv}.
\end{itemize}

\begin{definition}[The Functor $\KC$]\label{def:KC}
Suppose $(\C,*)$ is an indexing $G$-category for a compact Lie group $G$, and  $\tF \cn \C \to \FG$ is a pointed $G$-functor.  Define the functor
\[\CTopgst \fto{\KC} \Gspec\]
as follows.  Note that $\tF$ is used in the construction of $\KC$, but it is suppressed from the notation $\KC$ to simplify the typography.
\begin{description}
\item[Object assignment of $\KC$] 
Given a pointed $G$-functor $X \cn \C \to \Topgst$ \cref{ctopgst_obj}, the orthogonal $G$-spectrum \pcref{def:gsp_module}
\begin{equation}\label{KC_object}
(\KC X, \umu) \in \Gspec
\end{equation}
sends an object $V \in \IU$ \pcref{def:indexing_gspace}  to the coend
\begin{equation}\label{KCxv}
(\KC X)_V = \int^{c \in \C} (S^V)^{\tF c} \sma Xc
\end{equation}
taken in $\Topst$.  Its $G$-action is defined in \cref{KCxv_rep_gact}.
\begin{description}
\item[Representatives] 
For each object $c \in \C$, the pointed finite $G$-set $\tF c \in \FG$ is regarded as a discrete pointed $G$-space, and $S^V$ is the $V$-sphere \pcref{def:g_sphere}.  The pointed $G$-space 
\begin{equation}\label{SVFc}
(S^V)^{\tF c} = \Topgst(\tF c, S^V)
\end{equation}
consists of pointed morphisms $\tF c \to S^V$ \cref{Gtopst_smc}, with $G$ acting by conjugation \cref{ginv_h_g}.  The pointed space $(\KC X)_V$ is a quotient of the wedge
\[\bigvee_{c \in \C} (S^V)^{\tF c} \sma Xc.\]
Each point in $(\KC X)_V$ is represented by a pair
\begin{equation}\label{KCxv_rep}
(\upom; x) \in (S^V)^{\tF c} \ttimes Xc.
\end{equation}
The pair $(\upom; x)$ represents the basepoint if $x \in Xc$ is the basepoint or if $\upom$ is  constant at the basepoint of $S^V$.
\item[Relations]
The defining relation of the coend \cref{KCxv} identities, for each triple
\[(\upom; \uppsi; x) \in (S^V)^{\tF c} \ttimes \C(c',c) \ttimes Xc',\]
the pairs
\begin{equation}\label{KCxv_relations}
\begin{split}
\big(\upom(\tF \uppsi) ; x\big) & \in (S^V)^{\tF c'} \ttimes Xc' \andspace\\
\big(\upom; (X\uppsi)x \big) & \in (S^V)^{\tF c} \ttimes Xc.
\end{split}
\end{equation}
\item[$G$-action] 
The group $G$ acts diagonally on representatives.  This means that, for an element $g \in G$ and a representative pair $(\upom; x)$ in $(\KC X)_V$ \cref{KCxv_rep}, the diagonal $g$-action is given by
\begin{equation}\label{KCxv_rep_gact}
g \cdot (\upom; x) = (g\upom\ginv ; gx),
\end{equation}
where $g\upom\ginv$ means the composite pointed morphism
\[\tF c \fto{\ginv} \tF c \fto{\upom} S^V \fto{g} S^V.\]
\item[Morphism assignment of $\KC X$]
For a linear isometric isomorphism $f \cn V \fiso W$ in $\IU$, the pointed homeomorphism \cref{iu_space_xf}
\begin{equation}\label{KCxf}
(\KC X)_V \fto[\iso]{(\KC X)_f} (\KC X)_W
\end{equation}
is induced by the pointed homeomorphisms 
\[(S^V)^{\tF c} \fto[\iso]{f \circ -} (S^W)^{\tF c} \forspace c \in \C\]
that postcompose with the pointed homeomorphism $f \cn S^V \fiso S^W$.    In terms of representatives \cref{KCxv_rep}, it is given by
\begin{equation}\label{KCxf_rep}
(\KC X)_f(\upom; x) = (f\upom; x).
\end{equation}
\item[Sphere action on $\KC X$]
For each pair of objects $(V,W) \in (\IUsk)^2$, the $(V,W)$-component pointed $G$-morphism of $\umu$ \cref{gsp_action_vw} is defined by the following commutative diagram in $\Gtopst$.
\begin{equation}\label{KCx_action_vw}
\begin{tikzpicture}[vcenter]
\def\h{4.2} \def\u{-1} \def\v{-1.4}
\draw[0cell=.9]
(0,0) node (a11) {(\KC X)_V \sma S^W}
(a11)++(\h,0) node (a12) {(\KC X)_{V \oplus W}}
(a11)++(0,\u) node (a21) {\big( \txint^{c \in \C} (S^V)^{\tF c} \sma Xc \big) \sma S^W}
(a12)++(0,\u) node (a22) {\txint^{c \in \C} (S^{V \oplus W})^{\tF c} \sma Xc}
(a21)++(0,\v) node (a3) {\txint^{c \in \C} \big( (S^V)^{\tF c} \sma S^W \big) \sma Xc}
;
\draw[1cell=.8]
(a11) edge node {\umu_{V,W}} (a12)
(a11) edge[equal] (a21)
(a12) edge[equal] (a22)
(a21) edge node[swap] {\iso} (a3)
(a3) [rounded corners=2pt] -| node[pos=.25] {\asm = \txint^{c \in \C}\! \asm_{c} \sma 1} (a22)
;
\end{tikzpicture}
\end{equation}
\begin{itemize}
\item The pointed $G$-homeomorphism denoted by $\iso$ first commutes $- \sma S^W$ with the coend.  Then it moves $S^W$ to the left of $Xc$ using the associativity isomorphism and braiding for the symmetric monoidal category $(\Gtopst,\sma)$ \cref{Gtopst_smc}.
\item The pointed $G$-morphism $\asm$ is induced by the pointed $G$-morphisms
\begin{equation}\label{KC_assembly}
(S^V)^{\tF c} \sma S^W \fto{\asm_{c}} (S^{V \oplus W})^{\tF c}
\end{equation}
for $c \in \C$ defined by the assignment
\[\begin{split}
& (\upom; y) \in (S^V)^{\tF c} \sma S^W\\
&\mapsto \big(\tF c \fto{\upom} S^V \fto{- \oplus y} S^{V \oplus W} \big) \in (S^{V \oplus W})^{\tF c}.
\end{split}\]
In other words, 
\[\asm_{c}(\upom; y) = \upom \oplus y\]
sends an element $i \in \tF c$ to the point 
\begin{equation}\label{KC_assembly_exp}
\big( \asm_{c}(\upom; y)\big)(i) 
= (\upom i) \oplus y \in S^{V \oplus W}.
\end{equation}
For a representative pair $(\upom; x)$ in $(\KC X)_V$ \cref{KCxv_rep} and a point $y \in S^W$, $\umu_{V,W}$ is given by
\begin{equation}\label{KCx_actrep}
\umu_{V,W}\big((\upom; x); y\big) = (\upom \oplus y; x).
\end{equation}
\end{itemize}
Replacing the indexing $G$-category $(\GG,\vstar)$ \pcref{def:GG} and the pointed $G$-functor $\sma \cn \GG \to \FG$ \pcref{sma_symmon} by, respectively, $(\C,*)$ and $\tF \cn \C \to \FG$, the proof of \cref{kggx_welldef} proves that $(\KC X, \umu)$ is an orthogonal $G$-spectrum.
\end{description}
\item[Morphism assignment of $\KC$]
The functor $\KC$ sends a $G$-natural transformation between $\CG$-spaces $X$ and $X'$ \cref{ctopgst_mor}
\begin{equation}\label{tha_cgspace}
\begin{tikzpicture}[vcenter]
\def\t{28}
\draw[0cell]
(0,0) node (a1) {\C}
(a1)++(1.8,0) node (a2) {\phantom{\Gsk}}
(a2)++(.23,-.03) node (a2') {\Topgst}
;
\draw[1cell=.9]
(a1) edge[bend left=\t] node {X} (a2)
(a1) edge[bend right=\t] node[swap] {X'} (a2)
;
\draw[2cell]
node[between=a1 and a2 at .42, rotate=-90, 2label={above,\tha}] {\Rightarrow}
;
\end{tikzpicture}
\end{equation}
to the $G$-morphism between orthogonal $G$-spectra (\cref{def:gsp_morphism} \cref{def:gsp_morphism_ii})
\begin{equation}\label{KC_tha}
(\KC X, \umu) \fto{\KC\tha} (\KC X', \umu).
\end{equation}
Its component at an object $V \in \IU$ \pcref{def:indexing_gspace}  is the pointed $G$-morphism defined by the following commutative diagram in $\Gtopst$ \cref{Gtopst_smc}.
\begin{equation}\label{KC_tha_v}
\begin{tikzpicture}[vcenter]
\def\v{-1.4}
\draw[0cell=1]
(0,0) node (a11) {(\KC X)_V}
(a11)++(2.7,0) node (a12) {\txint^{c \in \C} (S^V)^{\tF c} \sma Xc}
(a11)++(0,\v) node (a21) {(\KC X')_V}
(a12)++(0,\v) node (a22) {\txint^{c \in \C} (S^V)^{\tF c} \sma X'c}
;
\draw[1cell=.9]
(a11) edge[equal] (a12)
(a21) edge[equal] (a22)
(a11) edge[transform canvas={xshift=1em}] node[swap] {(\KC\tha)_V} (a21)
(a12) edge[transform canvas={xshift=-2em}, shorten <=-.5ex] node {\txint^{c} 1 \sma \tha_{c}} (a22)
;
\end{tikzpicture}
\end{equation}
It sends a representative pair $(\upom; x)$ in $(\KC X)_V$ \cref{KCxv_rep} to the representative pair
\begin{equation}\label{KC_tha_rep}
(\KC\tha)_V (\upom; x) = \big(\upom; \tha_{c} x\big)
\end{equation}
of $(\KC X')_V$.  Replacing $(\GG,\vstar)$ and $\sma \cn \GG \to \FG$ by, respectively, $(\C,*)$ and $\tF \cn \C \to \FG$, the proof of \cref{kggpsi_welldef} proves that $(\KC\tha)$ is a $G$-morphism between orthogonal $G$-spectra.
\item[Functoriality] 
The functoriality of $\KC$ follows from \cref{KC_tha_rep} and the fact that identities and composition are defined componentwise in $\Gspec$ and $\CTopgst$ \pcref{def:gsp_morphism,def:ctopgst}.\defmark
\end{description}
\end{definition}

\begin{example}\label{ex:KC}\
\begin{enumerate}
\item The prolongation functor \pcref{def:Kgg_functor}
\[\GGTopg \fto{\Kgg} \Gspec \] 
is the instance of the functor $\KC$ for the indexing $G$-category $(\GG,\vstar)$ and the pointed $G$-functor $\sma \cn \GG \to \FG$ \pcref{def:GG,def:smashFGGG}. 
\item The prolongation functor \pcref{def:Kfg_functor}
\[\FGTopg \fto{\Kfg} \Gspec \] 
is the instance of the functor $\KC$ for the indexing $G$-category $(\FG,\ordz)$ \pcref{def:FG} and the identity functor $1 \cn \FG \to \FG$.\defmark
\end{enumerate}
\end{example}

\subsection*{Naturality of $\KC$}

\cref{KC_funl} proves that the prolongation functor $\KC$ \pcref{def:KC} is compatible with changing the indexing $G$-category $\C$ \pcref{def:ctopgst} via the left adjoint in \cref{fun_adj}.

\begin{theorem}\label{KC_funl}
Suppose $\C$ and $\D$ are indexing $G$-categories for a compact Lie group $G$, and
\begin{equation}\label{fHF_diag}
\begin{tikzpicture}[vcenter]
\def\u{.6} \def\h{1.6}
\draw[0cell]
(0,0) node (a1) {\C}
(a1)++(\h,0) node (a2) {\D}
(a2)++(\h,0) node (a3) {\FG}
;
\draw[1cell=.9]
(a1) edge node {\fun} (a2)
(a2) edge node {\tH} (a3)
(a1) [rounded corners=2pt] |- ($(a2)+(-1,\u)$) -- node {\tF} ($(a2)+(1,\u)$) -| (a3)
;
\end{tikzpicture}
\end{equation}
is a commutative diagram of pointed $G$-functors.  Then there is a natural isomorphism\label{not:kfgiso}
\begin{equation}\label{fKCKD}
\begin{tikzpicture}[vcenter]
\def\h{2.5} \def\v{-1.3} \def\t{20}
\draw[0cell]
(0,0) node (a1) {\CTopgst}
(a1)++(\h,0) node (a2) {\DTopgst}
(a1)++(\h/2,\v) node (a3) {\phantom{\Gspec}}
(a3)++(0,-.15*\v) node (a3') {\Gspec}
;
\draw[1cell=.9]
(a1) edge[bend right=\t] node[swap,pos=.3] {\KC} (a3)
(a1) edge node {\funl} (a2)
(a2) edge[bend left=\t] node[pos=.3] {\KD} (a3)
;
\draw[2cell=.9]
node[between=a1 and a2 at .5, shift={(0,.45*\v)}, rotate=200, 2labelmed={below,\kfgiso}] {\Rightarrow}
;
\end{tikzpicture}
\end{equation}
between the functors $\KD\funl$ and $\KC$.
\end{theorem}

\begin{proof}
The natural isomorphism 
\[\KD\funl \fto[\iso]{\kfgiso} \KC\] 
is defined by the pointed $G$-homeomorphisms in \cref{kfgisoxv} for a pointed $G$-functor $X \cn \C \to \Topgst$ \cref{ctopgst_obj} and an object $V \in \IU$ \pcref{def:indexing_gspace}.  The wedge $\txwedge$ is indexed by the $G$-set $\Dpun(\fun c, d)$ of nonzero morphisms $\fun c \to d$ in $\D$ \cref{Cpunc}.
\begin{equation}\label{kfgisoxv}
\begin{aligned}
& (\KD \funl X)_V &&\\
&= \txint^{d \in \D} (S^V)^{\tH d} \sma (\funl X)d && \text{by \cref{KCxv} for $(\D,\tH)$}\\
&= \txint^{d \in \D} (S^V)^{\tH d} \sma \big[\txint^{c \in \C} \txwedge Xc\big] && \text{by \cref{funlx_d}}\\
&\iso \txint^{c \in \C} \big[\txint^{d \in \D} \txwedge (S^V)^{\tH d}\big] \sma Xc && \text{by commutation}\\
&\iso \txint^{c \in \C} (S^V)^{\tH\fun c} \sma Xc && \text{by Yoneda}\\
&=  \txint^{c \in \C} (S^V)^{\tF c} \sma Xc && \text{by $\tH\fun = \tF$}\\
&= (\KC X)_V && \text{by \cref{KCxv} for $(\C,\tF)$}
\end{aligned}
\end{equation}
\begin{itemize}
\item The first pointed homeomorphism in \cref{kfgisoxv} uses
\begin{itemize}
\item the commutation of $(S^V)^{\tH d} \msmam -$ with the coend $\txint^{c \in \C}$ and the wedge $\txwedge$; 
\item the commutation of the coends $\txint^{d \in \D}$ and $\txint^{c \in \C}$; and 
\item the commutation of $- \msmam Xc$ with the wedge $\txwedge$ and the coend $\txint^{d \in \D}$.
\end{itemize}
\item The second pointed homeomorphism in \cref{kfgisoxv} is induced by the Enriched Yoneda Density Theorem \cite[3.7.8]{cerberusIII}, which yields the pointed homeomorphism
\[\int^{d \in \D} \bigvee_{\Dpun(\fun c, d)} (S^V)^{\tH d} \iso (S^V)^{\tH\fun c}\]
that is natural in $c \in \C$.
\item The pointed homeomorphisms in \cref{kfgisoxv} are $G$-equivariant by \cref{funlx_gaction,KCxv_rep_gact}.  They are natural in $V \in \IU$ by \cref{KCxf_rep}.  They are compatible with the sphere actions by \cref{KCx_actrep}.  They are natural in $X \in \CTopgst$ by \cref{funl_tha_d,KC_tha_v}.  These properties can also be read off from the explicit description of $\kfgiso$ given in \cref{kfgisoxv_exp}.
\end{itemize}
Thus, the pointed $G$-homeomorphisms in \cref{kfgisoxv} assemble into a natural isomorphism $\kfgiso$.
\end{proof}

\begin{explanation}[Unpacking $\kfgiso$]\label{expl:kfgisoxv}
The natural pointed $G$-homeomorphism \cref{kfgisoxv}
\[(\KD\funl X)_V \fto[\iso]{\kfgiso_{X,V}} (\KC X)_V\]
sends a representative
\[(\upom; \upla; x) \in (S^V)^{\tH d} \times \Dpun(\fun c, d) \times Xc\]
in $(\KD\funl X)_V$ to the representative
\begin{equation}\label{kfgisoxv_exp}
\big(\tF c = \tH\fun c \fto{\tH\upla} \tH d \fto{\upom} S^V ; x \big)
\in (S^V)^{\tF c} \times Xc
\end{equation}
in $(\KC X)_V$.  The inverse pointed $G$-homeomorphism
\[(\KC X)_V \fto[\iso]{\kfgiso^{-1}_{X,V}} (\KD\funl X)_V\]
sends a representative
\[(\upom; x) \in (S^V)^{\tF c} \times Xc\]
in $(\KC X)_V$ to the representative
\[\big(\tH\fun c = \tF c \fto{\upom} S^V ; \fun c \fto{1} \fun c ; x \big)
\in (S^V)^{\tH\fun c} \times \Dpun(\fun c, \fun c) \times Xc\]
in $(\KD\funl X)_V$.
\end{explanation}

\section{Shimakawa $G$-Spectra and $\GGG$-Spaces}
\label{sec:KfgKgg_comp}

This section compares Shimakawa's (homotopical) functor (\cref{def:Kfg_functor} and \cref{bar_functor_FG})
\[\FGTopg \fto{\Bc} \FGTopg \fto{\Kfg} \Gspec\]
and the prolongation functor \pcref{def:Kgg_functor}
\[\GGTopg \fto{\Kgg} \Gspec.\]
The conclusion is that $\Kfg$ and $\Kgg$ compute each other.   In the homotopical case, $\Kfg\Bc$ and $\Kgg$ compute each other on proper objects up to natural componentwise weak $G$-equivalences.

\secoutline
\begin{itemize}
\item \cref{kfgkgg_compare} establishes the natural isomorphisms
\[\Kgg\ifgl \fto[\iso]{\kggik} \Kfg \andspace \Kfg\!\smal\! \fto[\iso]{\kfgsk} \Kgg.\]
The functors $\ifgl$ and $\smal$ are the left adjoints induced by the length-1 inclusion functor $\ifg \cn \FG \to \GG$ and the smash functor $\sma \cn \GG \to \FG$.  \cref{expl:kfgkgg_compare} describes these natural isomorphisms explicitly.  Thus, each orthogonal $G$-spectrum produced by Shimakawa's functor $\Kfg$ can also be obtained from a $\GGG$-space by applying $\Kgg$, and the converse is also true.
\item \cref{kfgbkgg_compare} establishes the natural transformation
\begin{equation}\label{kbgg_seci}
\begin{tikzpicture}[vcenter]
\def\h{2.2}
\draw[0cell]
(0,0) node (a1) {\Kfg\Bc}
(a1)++(\h,0) node (a2) {\Kfg} 
(a2)++(.9*\h,0) node (a3) {\phantom{\Kgg\ifgl}} 
(a3)++(0,-.04) node (a3') {\Kgg\ifgl.}
;
\draw[1cell=.9]
(a1) edge node {\Kfg\retn} (a2)
(a2) edge node {\kggik^{-1}} node[swap] {\iso} (a3)
;
\draw[1cell=.9]
(a1) [rounded corners=2pt] |- ($(a2)+(-1,.6)$) -- node {\kbgg} ($(a2)+(1,.6)$) -| (a3')
 ;
\end{tikzpicture}
\end{equation}
Its component at each proper $\FGG$-space is componentwise a weak $G$-equivalence.  Thus, for each proper $\FGG$-space $X$---including all componentwise realizations of $\FG$-simplicial $G$-sets---Shimakawa's orthogonal $G$-spectrum $\Kfg\Bc X$ can also be obtained from a $\GGG$-space by applying $\Kgg$, up to a natural componentwise weak $G$-equivalence.
\item \cref{def:proper_ggg} defines \emph{proper} $\GGG$-spaces.  \cref{reast_properGG} proves that the componentwise realization of each $\GG$-simplicial $G$-set is a proper $\GGG$-space.
\item \cref{kggkfgb_compare} establishes the natural transformation
\begin{equation}\label{kgbk_seci}
\begin{tikzpicture}[vcenter]
\def\h{2.7} \def\u{.65}
\draw[0cell]
(0,0) node (a1) {\Kfg\Bc\smal}
(a1)++(\h,0) node (a2) {\Kfg\smal} 
(a2)++(.7*\h,0) node (a3) {\phantom{\Kgg}}
(a3)++(0,.04) node (a3') {\Kgg.}
;
\draw[1cell=.9]
(a1) edge node {\Kfg\retn_{\smal}} (a2)
(a2) edge node {\kfgsk} node[swap] {\iso} (a3)
;
\draw[1cell=.9]
(a1) [rounded corners=2pt] |- ($(a2)+(-.3*\h,\u)$) -- node {\kgbk} ($(a2)+(0,\u)$) -| (a3')
 ;
\end{tikzpicture}
\end{equation}
Its component at each proper $\GGG$-space is componentwise a weak $G$-equivalence.  Thus, for each proper $\GGG$-space $Y$, the orthogonal $G$-spectrum $\Kgg Y$ can also be obtained from an $\FGG$-space by applying Shimakawa's homotopical functor $\Kfg\Bc$, up to a natural componentwise weak $G$-equivalence.
\end{itemize}

\subsection*{Comparison of $\Kfg$ and $\Kgg$}

Using \cref{KC_funl} along with the length-1 inclusion functor $\ifg \cn \FG \to \GG$ and the smash functor $\sma \cn \GG \to \FG$ \pcref{def:ifgGG,def:smashFGGG}, \cref{kfgkgg_compare} compares the prolongation functors $\Kgg$ and $\Kfg$ \pcref{def:Kgg_functor,def:Kfg_functor}.

\begin{theorem}\label{kfgkgg_compare}\index{Shimakawa K-theory@Shimakawa $K$-theory!comparison with multifunctorial $K$-theory}\index{multifunctorial K-theory@multifunctorial $K$-theory!comparison with Shimakawa}
For each compact Lie group $G$, there is a diagram
\begin{equation}\label{kfgkgg_diag}
\begin{tikzpicture}[vcenter]
\def\g{.25} \def\h{2.5} \def\u{.6} \def\t{27}
\draw[0cell]
(0,0) node (a1) {\FGTopg}
(a1)++(\h,0) node (a2) {\GGTopg}
(a2)++(\h,0) node (a3) {\FGTopg}
(a2)++(0,-1.4) node (a4) {\Gspec}
(a4)++(0,-.15) node (a4') {\phantom{\Gspec}}
;
\draw[1cell=.9]
(a1) edge node {\ifgl} (a2)
(a2) edge node {\smal} (a3)
(a2) edge node[pos=.4] {\Kgg} (a4)
(a1) edge[bend right=\t] node[swap,pos=.2] {\Kfg} (a4)
(a3) edge[bend left=\t] node[pos=.2] {\Kfg} (a4)
(a1) [rounded corners=2pt] |- ($(a2)+(-1,\u)$) -- node {1} ($(a2)+(1,\u)$) -| (a3)
;
\draw[2cell=.9]
node[between=a1 and a4 at .5, shift={(\g,0)}, rotate=180, 2labelmed={below,\kggik}] {\Rightarrow}
node[between=a3 and a4 at .5, shift={(-\g,0)}, rotate=180, 2labelmed={below,\kfgsk}] {\Rightarrow}
;
\end{tikzpicture}
\end{equation}
such that each of the three regions commutes up to a natural isomorphism.
\end{theorem}

\begin{proof}
The functors $\ifgl$ and $\smal$ are the left adjoints of, respectively, the pullback functors $\ifgst$ and $\smast$ \pcref{fun_adj}.  The natural isomorphism 
\[\smal\ifgl \iso 1\] 
is discussed in \cref{ex:smai_adj}.  

As we discuss in \cref{ex:KC}, $\Kfg$ is the instance of $\KC$ \pcref{def:KC} for the indexing $G$-category $(\FG,\ordz)$ and the identity functor $1_{\FG}$, while $\Kgg$ is the instance of $\KC$ for $(\GG,\vstar)$ and the smash functor $\sma \cn \GG \to \FG$.  The natural isomorphisms 
\begin{equation}\label{KggKfg_isos}
\Kgg\ifgl \fto[\iso]{\kggik} \Kfg \andspace \Kfg\!\smal\! \fto[\iso]{\kfgsk} \Kgg
\end{equation}
are the instances of \cref{KC_funl} for the following commutative diagrams of pointed $G$-functors.
\begin{equation}\label{ismai}
\begin{tikzpicture}[baseline={(a1.base)}]
\def\u{.6} \def\h{1.5}
\draw[0cell]
(0,0) node (a1) {\FG}
(a1)++(\h,0) node (a2) {\GG}
(a2)++(\h,0) node (a3) {\FG}

;
\draw[1cell=.9]
(a1) edge node {\ifg} (a2)
(a2) edge node {\sma} (a3)
(a1) [rounded corners=2pt] |- ($(a2)+(-1,\u)$) -- node {1} ($(a2)+(1,\u)$) -| (a3)
;
\begin{scope}[shift={(3*\h,0)}]
\draw[0cell]
(0,0) node (a1) {\GG}
(a1)++(\h,0) node (a2) {\FG}
(a2)++(\h,0) node (a3) {\FG}
;
\draw[1cell=.9]
(a1) edge node {\sma} (a2)
(a2) edge node {1} (a3)
(a1) [rounded corners=2pt] |- ($(a2)+(-1,\u)$) -- node {\sma} ($(a2)+(1,\u)$) -| (a3)
;
\end{scope}
\end{tikzpicture}
\end{equation}
This finishes the proof.
\end{proof}

\begin{explanation}[Unpacking $\kggik$ and $\kfgsk$]\label{expl:kfgkgg_compare}
We consider a pointed $G$-functor $X \cn \FG \to \Topgst$, a pointed $G$-functor $Y \cn \GG \to \Topgst$, and an object $V \in \IU$ \pcref{def:indexing_gspace}.  By \cref{Kggxv,Kfgxv,ifglx_nbe,smaly_mal}, there are pointed $G$-spaces as follows.
\[\begin{split}
(\Kgg\ifgl X)_V 
&= \int^{\nbe \in \GG} (S^V)^{\sma\nbe} \sma 
\Big[\int^{\mal \in \FG} \mquad \bigvee_{\GGpunc(\ifg\mal;\, \nbe)} \mquad X\mal \Big]\\
&\iso \int^{\nbe \in \GG} \int^{\mal \in \FG} (S^V)^{\sma\nbe} \sma 
\Big[\bigvee_{\GGpunc(\ifg\mal;\, \nbe)} \mquad X\mal \Big]\\
(\Kfg \!\smal\! Y)_V 
&= \int^{\mal \in \FG} (S^V)^{\mal} \sma 
\Big[\int^{\nbe \in \GG} \mquad \bigvee_{\FGpunc(\sma\nbe;\, \mal)} \mquad Y\nbe \Big]\\
&\iso \int^{\mal \in \FG} \int^{\nbe \in \GG} (S^V)^{\mal} \sma 
\Big[\bigvee_{\FGpunc(\sma\nbe;\, \mal)} \mquad Y\nbe \Big]\\
\end{split}\]
\cref{expl:kfgisoxv} specializes to the natural isomorphisms in \cref{KggKfg_isos} as follows.  
\begin{description} 
\item[Unpacking $\kggik$] 
The natural pointed $G$-homeomorphism
\[(\Kgg\ifgl X)_V \fto[\iso]{\kggik_{X,V}} (\Kfg X)_V\]
sends a representative
\[(\upom; \upla; x) \in (S^V)^{\sma\nbe} \times \GGpunc(\ifg \mal; \nbe) \times X\mal\]
in $(\Kgg\ifgl X)_V$ to the representative
\[\big(\mal = \sma\ifg\mal \fto{\sma\upla} \msmam\nbe \fto{\upom} S^V ; x \big)
\in (S^V)^{\mal} \times X\mal\]
in $(\Kfg X)_V$.  The inverse pointed $G$-homeomorphism
\begin{equation}\label{kfgkgg_isoxv}
(\Kfg X)_V \fto[\iso]{\kggik^{-1}_{X,V}} (\Kgg\ifgl X)_V
\end{equation}
sends a representative
\[(\upom; x) \in (S^V)^{\mal} \times X\mal\]
in $(\Kfg X)_V$ to the representative
\[\big(\!\smam\ifg\mal = \mal \fto{\upom} S^V ; 1 ; x \big) 
\in (S^V)^{\sma\ifg\mal} \times \GGpunc(\ifg\mal; \ifg\mal) \times X\mal\]
in $(\Kgg\ifgl X)_V$.
\item[Unpacking $\kfgsk$] 
The natural pointed $G$-homeomorphism
\begin{equation}\label{kfgsmalkgg_iso}
(\Kfg \!\smal\! Y)_V \fto[\iso]{\kfgsk_{Y,V}} (\Kgg Y)_V
\end{equation}
sends a representative
\[(\upom; \upla; y) \in (S^V)^{\mal} \times \FGpunc(\sma\nbe; \mal) \times Y\nbe\]
in $(\Kfg \!\smal\! Y)_V$ to the representative
\[\big(\msmam\nbe \fto{\upla} \mal \fto{\upom} S^V ; y \big)
\in (S^V)^{\sma\nbe} \times Y\nbe\]
in $(\Kgg Y)_V$.  The inverse pointed $G$-homeomorphism
\[ (\Kgg Y)_V \fto[\iso]{\kfgsk^{-1}_{Y,V}} (\Kfg \!\smal\! Y)_V\]
sends a representative
\[(\upom; y) \in (S^V)^{\sma\nbe} \times Y\nbe\]
in $(\Kgg Y)_V$ to the representative
\[(\upom; 1; y) \in (S^V)^{\sma\nbe} \times \FGpunc(\sma\nbe; \sma\nbe) \times Y\nbe\]
in $(\Kfg \!\smal\! Y)_V$.\defmark
\end{description}
\end{explanation}

\subsection*{Comparison of $\Kfg\Bc$ and $\Kgg$}

Shimakawa's homotopical functor \pcref{def:Kfg_functor,def:gtop_retn_secn}
\[\FGTopg \fto{\Bc} \FGTopg \fto{\Kfg} \Gspec\]
is the last part of the homotopical Shimakawa (strong) $K$-theory \cref{khsho_khshosg}.  The bar functor \cref{bar_functor_FG}
\[\Bc = \BcFG \cn \FGTopg \to \FGTopg\] 
is connected to the identity functor via the retraction \cref{retn_barFG_id}
\[\Bc \fto{\retn} 1_{\FGTopg}.\]
An explicit, point-set level description of $\retn$ is given in \cref{retn_cr,real_FG,Br_FGX}.  For a finite group $G$, an $\FGG$-space $X$ is \emph{proper} if the simplicial $G$-space $(\Pwg X)_{B_\crdot}$ is Reedy cofibrant for each pointed finite simplicial $G$-set $B_\crdot$ \pcref{def:proper_fgg}.  \cref{kfgbkgg_compare} describes Shimakawa's homotopical functor $\Kfg\Bc$ in terms of the prolongation functor $\Kgg$ \pcref{def:Kgg_functor}, the retraction $\retn$, and the natural isomorphism $\kggik \cn \Kgg\ifgl \fiso \Kfg$ in \cref{kfgkgg_compare}.

\begin{theorem}\label{kfgbkgg_compare}\index{Shimakawa K-theory@Shimakawa $K$-theory!comparison with multifunctorial $K$-theory}\index{multifunctorial K-theory@multifunctorial $K$-theory!comparison with Shimakawa}
Suppose $G$ is a compact Lie group in \cref{kfgbkgg_compare_i} and a finite group in \cref{kfgbkgg_compare_ii,kfgbkgg_compare_iii}.
\begin{enumerate}
\item\label{kfgbkgg_compare_i} 
The retraction $\retn$ and the natural isomorphism $\kggik$ induce a natural transformation $\kbgg$ as follows.
\begin{equation}\label{kbgg_thm}
\begin{tikzpicture}[vcenter]
\def\t{25} \def\h{2.6} \def\v{-1.4} \def\c{1em}
\draw[0cell]
(0,0) node (a11) {\FGTopg}
(a11)++(\h,0) node (a12) {\GGTopg}
(a11)++(0,\v) node (a21) {\FGTopg}
(a12)++(0,\v) node (a22) {\Gspec}
;
\draw[1cell=.9]
(a11) edge node {\ifgl} (a12)
(a12) edge node {\Kgg} (a22)
(a21) edge node[swap] {\Kfg} (a22)
(a11) edge[transform canvas={xshift=-\c}, bend right=\t] node[swap] {\Bc} (a21)
(a11) edge[transform canvas={xshift=\c}, bend left=\t] node {1} (a21)
;
\draw[2cell=1]
node[between=a11 and a21 at .6, 2labelmed={above,\retn\phantom{x}}] {\Rightarrow}
node[between=a11 and a21 at .6, shift={(.63*\h,0)}, 2labelmed={above,\kggik^{-1}}] {\Rightarrow}
;
\end{tikzpicture}
\end{equation}
\item\label{kfgbkgg_compare_ii}
The component of $\kbgg$ at each proper $\FGG$-space is componentwise a weak $G$-equivalence.
\item\label{kfgbkgg_compare_iii}
For each $\FG$-simplicial $G$-set $X$, the component of $\kbgg$ at $|X|_* \in \FGTopg$ is componentwise a weak $G$-equivalence.
\end{enumerate}
\end{theorem}

\begin{proof}
The natural transformation $\kbgg$ is the composite\label{not:kbgg}
\begin{equation}\label{kbgg_composite}
\begin{tikzpicture}[vcenter]
\def\h{2.2}
\draw[0cell]
(0,0) node (a1) {\Kfg\Bc}
(a1)++(\h,0) node (a2) {\Kfg} 
(a2)++(.9*\h,0) node (a3) {\phantom{\Kgg\ifgl}} 
(a3)++(0,-.04) node (a3') {\Kgg\ifgl}
;
\draw[1cell=.9]
(a1) edge node {\Kfg\retn} (a2)
(a2) edge node {\kggik^{-1}} node[swap] {\iso} (a3)
;
\draw[1cell=.9]
(a1) [rounded corners=2pt] |- ($(a2)+(-1,.6)$) -- node {\kbgg} ($(a2)+(1,.6)$) -| (a3')
 ;
\end{tikzpicture}
\end{equation}
of 
\begin{itemize}
\item the whiskered natural transformation $\Kfg\retn$ and
\item the natural isomorphism $\kggik^{-1} \cn \Kfg \fiso \Kgg\ifgl$ in \cref{kfgkgg_compare}.
\end{itemize}  
For \cref{kfgbkgg_compare_ii}, suppose $X$ is a proper $\FGG$-space for a finite group $G$ \pcref{def:proper_fgg}.  The retraction \cref{ret_FG}
\[\Bc X \fto{\retn_X} X \inspace \FGTopg\]
is componentwise a pointed $G$-homotopy equivalence.  By \cref{bar_proper}, the $\FGG$-space $\Bc X$ is always proper.  By \cref{thm:Kfg_inv} and the properness of $X$, the $G$-morphism of orthogonal $G$-spectra
\begin{equation}\label{kfg_retn_x}
\Kfg \Bc X \fto{\Kfg\retn_X} \Kfg X
\end{equation}
is componentwise a weak $G$-equivalence.  Thus, its composite with the $G$-homeomorphism $\kggik^{-1}_X$ \cref{kfgkgg_isoxv},
\begin{equation}\label{kbggX_weakg}
\begin{tikzpicture}[vcenter]
\def\h{2.7} \def\u{.7}
\draw[0cell]
(0,0) node (a1) {\Kfg\Bc X}
(a1)++(\h,0) node (a2) {\Kfg X} 
(a2)++(.87*\h,0) node (a3) {\phantom{\Kgg\ifgl X}}
(a3)++(0,-.04) node (a3') {\Kgg\ifgl X,}
;
\draw[1cell=.9]
(a1) edge node {\Kfg\retn_X} (a2)
(a2) edge node {\kggik^{-1}_X} node[swap] {\iso} (a3)
;
\draw[1cell=.9]
(a1) [rounded corners=2pt] |- ($(a2)+(-1,\u)$) -- node {\kbgg_X} ($(a2)+(1,\u)$) -| (a3')
 ;
\end{tikzpicture}
\end{equation}
is also componentwise a weak $G$-equivalence.  Assertion \cref{kfgbkgg_compare_iii} follows from \cref{kfgbkgg_compare_ii} and the fact that $|X|_* \in \FGTopg$ is proper \pcref{reast_proper}.
\end{proof}

\cref{kfgbkgg_compare} describes Shimakawa's homotopical functor $\Kfg\Bc$ in terms of $\Kgg$ for proper $\FGG$-spaces.   To obtain the converse that describes $\Kgg$ in terms of $\Kfg\Bc$, we first discuss proper $\GGG$-spaces.

\subsection*{Proper $\GGG$-Spaces}

Recall from \cref{ex:smai_adj} the adjunction
\begin{equation}\label{smash_adjunction}
\begin{tikzpicture}[vcenter]
\def\u{.7} \def\h{2.6} \def\b{.6ex} \def\c{-.4ex}
\draw[0cell]
(0,0) node (a1) {\GGTopg}
(a1)++(\h,0) node (a2) {\FGTopg}
;
\draw[1cell=.9]
(a1) edge[transform canvas={yshift=\b}] node {\smal} (a2)
(a2) edge[transform canvas={yshift=\c}] node {\smast} (a1)
;
\end{tikzpicture}
\end{equation}
induced by the smash functor $\sma \cn \GG \to \FG$ \pcref{def:smashFGGG}.

\begin{definition}\label{def:proper_ggg}
For a finite group $G$, a $\GGG$-space $Y \cn \GG \to \Topgst$ \cref{GGTopg_obj} is \emph{proper}\index{proper!GGG-space@$\GGG$-space}\index{GGG-space@$\GGG$-space!proper} if the $\FGG$-space $\smal Y$ is proper in the sense of \cref{def:proper_fgg}.
\end{definition}

To see that \cref{def:proper_ggg} is reasonable, recall that $\ssetgst$ is the pointed $G$-category of pointed simplicial $G$-sets and pointed simplicial morphisms with the conjugation $G$-action \pcref{def:Gssetst}.  \cref{def:GGssetgst} is the $\GG$-analogue of \cref{def:Gssetst} \cref{def:Gssetst_fgssetgst}.  

\begin{definition}\label{def:GGssetgst}
For a group $G$, the category $\GGssetgst$ has pointed $G$-functors
\[(\GG,\vstar) \to (\ssetgst,*)\]
as objects, called \index{GG-simplicial G-set@$\GG$-simplicial $G$-set}\emph{$\GG$-simplicial $G$-sets}, and $G$-natural transformations as morphisms.
\end{definition}

With $\Topgst$ replaced by $\ssetgst$, the description of $\GGTopg$ in \cref{expl:GGTopg} also applies to $\GGssetgst$.  There is a factorization of the functor \cref{clast_GGCatg}
\begin{equation}\label{clast_factorsGG}
\begin{tikzpicture}[vcenter]
\def\h{2.5}
\draw[0cell]
(0,0) node (a1) {\GGCatg}
(a1)++(\h,0) node (a2) {\GGssetgst} 
(a2)++(\h,0) node (a3) {\GGTopg}
;
\draw[1cell=.9]
(a1) edge node {\Nerst} (a2)
(a2) edge node {\Reast} (a3)
;
\draw[1cell=.9]
(a1) [rounded corners=2pt] |- ($(a2)+(-1,.6)$) -- node {\clast} ($(a2)+(1,.6)$) -| (a3)
 ;
\end{tikzpicture}
\end{equation}
where $\Nerst$ and $\Reast$ postcompose with, respectively, the nerve $\Ner$ and  the realization $\Rea$ \pcref{expl:GGCatg_GGTopg}.  \cref{reast_properGG} \cref{reast_properGG_ii} proves that the image of the componentwise realization functor $\Reast$ is always proper \pcref{def:proper_ggg}.

\begin{lemma}\label{reast_properGG}
Suppose $G$ is a group in \cref{reast_properGG_i} and a finite group in \cref{reast_properGG_ii}.
\begin{enumerate}
\item\label{reast_properGG_i} The diagram
\begin{equation}\label{smash_realize}
\begin{tikzpicture}[vcenter]
\def\v{-1.4}
\draw[0cell]
(0,0) node (a11) {\GGssetgst}
(a11)++(2.6,0) node (a12) {\FGssetgst}
(a11)++(0,\v) node (a21) {\GGTopg}
(a12)++(0,\v) node (a22) {\FGTopg}
;
\draw[1cell=.9]
(a11) edge node {\smal} (a12)
(a21) edge node {\smal} (a22)
(a11) edge node[swap] {\Reast} (a21)
(a12) edge node[inner sep=4pt] {\Reast} (a22)
;
\end{tikzpicture}
\end{equation}
commutes up to a natural isomorphism.  
\item\label{reast_properGG_ii}  
For each $\GG$-simplicial $G$-set $Y$, the $\GGG$-space $|Y|_*$ is proper.
\end{enumerate}
\end{lemma}

\begin{proof}
In the diagram in \cref{reast_properGG_i}, the top functor $\smal$ is the left adjoint of the adjunction 
\begin{equation}\label{sma_adjunction}
\begin{tikzpicture}[vcenter]
\def\u{.7} \def\h{2.8} \def\b{.6ex} \def\c{-.4ex}
\draw[0cell]
(0,0) node (a1) {\GGssetgst}
(a1)++(\h,0) node (a2) {\FGssetgst}
;
\draw[1cell=.9]
(a1) edge[transform canvas={yshift=\b}] node {\smal} (a2)
(a2) edge[transform canvas={yshift=\c}] node {\smast} (a1)
;
\end{tikzpicture}
\end{equation}
induced by the smash functor $\sma \cn \GG \to \FG$ \pcref{def:smashFGGG}.  This adjunction is established by the same proof as \cref{fun_adj} by replacing $\Topgst$ with $\ssetgst$.  The desired natural isomorphism in \cref{reast_properGG_i} is given by the pointed $G$-homeomorphisms in \cref{reast_properGG_iso} for a $\GG$-simplicial $G$-set $Y$ and a pointed finite $G$-set $\mal \in \FG$, using \cref{real_simpgspace,smaly_mal}.  The wedge $\bigvee$ is indexed by the $G$-set $\FGpunc(\sma\nbe; \mal)$ of nonzero morphisms $\sma\nbe \to \mal$ in $\FG$.
\begin{equation}\label{reast_properGG_iso}
\begin{split}
& \big(\!\smal\! |Y|_*\big)\mal \\
& = \txint^{\nbe \in \GG} \bigvee |Y\nbe| \\
&=  \txint^{\nbe \in \GG} \bigvee \big[\txint^{\ordr \in \Delta} (Y\nbe)_r \sma \simppl^r \big] \\
&\iso \txint^{\nbe \in \GG} \txint^{\ordr \in \Delta} \bigvee \big[(Y\nbe)_r \sma \simppl^r \big] \\
&\iso \txint^{\ordr \in \Delta} \txint^{\nbe \in \GG} \bigvee \big[(Y\nbe)_r \sma \simppl^r \big] \\
&\iso \txint^{\ordr \in \Delta} \big[\txint^{\nbe \in \GG} \bigvee (Y\nbe)_r \big] \sma \simppl^r \\
&\iso \txint^{\ordr \in \Delta} \big[\txint^{\nbe \in \GG} \bigvee Y\nbe \big]_r \sma \simppl^r \\
&= \big|\! \txint^{\nbe \in \GG} \bigvee Y\nbe  \big| \\
&= |\!\smal\! Y|_* \mal
\end{split}
\end{equation}
The pointed homeomorphisms in \cref{reast_properGG_iso} are $G$-equivariant by \cref{funlx_gaction} and the trivial $G$-action on the topological $r$-simplex $\simp^r$ \cref{top_simplex}.  They are natural in $\mal \in \FG$ by \cref{funlx_mor}.  They are natural in $Y$ by \cref{funl_tha_d}.

For \cref{reast_properGG_ii}, by \cref{def:proper_ggg,reast_properGG_i}, the properness of $|Y|_* \in \GGTopg$ means the properness of 
\[\smal |Y|_* \iso |\!\smal\! Y|_* \inspace \FGTopg.\]
The $\FGG$-space $|\!\smal\! Y|_*$ is proper by \cref{reast_proper}.
\end{proof}

\subsection*{Comparison of $\Kgg$ and $\Kfg\Bc$}

\cref{kggkfgb_compare} describes the prolongation functor $\Kgg$ \pcref{def:Kgg_functor} in terms of Shimakawa's homotopical functor $\Kfg\Bc$ \pcref{def:Kfg_functor,def:gtop_retn_secn}, the retraction $\retn \cn \Bc \to 1_{\FGTopg}$ \cref{retn_barFG_id}, and the natural isomorphism $\kfgsk \cn \Kfg\!\smal\! \fiso \Kgg$ in \cref{kfgkgg_compare}.

\begin{theorem}\label{kggkfgb_compare}\index{Shimakawa K-theory@Shimakawa $K$-theory!comparison with multifunctorial $K$-theory}\index{multifunctorial K-theory@multifunctorial $K$-theory!comparison with Shimakawa}
Suppose $G$ is a compact Lie group in \cref{kggkfgb_compare_i} and a finite group in \cref{kggkfgb_compare_ii,kggkfgb_compare_iii}.
\begin{enumerate}
\item\label{kggkfgb_compare_i} 
The retraction $\retn$ and the natural isomorphism $\kfgsk$ induce a natural transformation $\kgbk$ as follows.
\begin{equation}\label{kgbk_composite}
\begin{tikzpicture}[vcenter]
\def\t{25} \def\h{2.6} \def\v{-1.4} \def\c{1em}
\draw[0cell]
(0,0) node (a11) {\GGTopg}
(a11)++(\h,0) node (a12) {\FGTopg}
(a11)++(0,\v) node (a21) {\Gspec}
(a12)++(0,\v) node (a22) {\FGTopg}
;
\draw[1cell=.9]
(a11) edge node {\smal} (a12)
(a12) edge[transform canvas={xshift=\c}, bend left=\t] node {\Bc} (a22)
(a12) edge[transform canvas={xshift=-\c}, bend right=\t] node[swap] {1} (a22)
(a22) edge node[pos=.4] {\Kfg} (a21)
(a11) edge node[swap] {\Kgg} (a21)
;
\draw[2cell=1]
node[between=a12 and a22 at .6, rotate=180, 2labelmed={below,\retn}] {\Rightarrow}
node[between=a12 and a22 at .6, shift={(-.63*\h,0)}, rotate=180, 2labelmed={below,\kfgsk}] {\Rightarrow}
;
\end{tikzpicture}
\end{equation}
\item\label{kggkfgb_compare_ii} 
The component of $\kgbk$ at each proper $\GGG$-space is componentwise a weak $G$-equivalence.
\item\label{kggkfgb_compare_iii} 
For each $\GG$-simplicial $G$-set $Y$, the component of $\kgbk$ at $|Y|_* \in \GGTopg$ is componentwise a weak $G$-equivalence.
\end{enumerate}
\end{theorem}

\begin{proof}
The natural transformation $\kgbk$ is the composite\label{not:kgbk}
\begin{equation}\label{kgbk_def}
\begin{tikzpicture}[vcenter]
\def\h{2.7} \def\u{.65}
\draw[0cell]
(0,0) node (a1) {\Kfg\Bc\smal}
(a1)++(\h,0) node (a2) {\Kfg\smal} 
(a2)++(.7*\h,0) node (a3) {\phantom{\Kgg}}
(a3)++(0,.04) node (a3') {\Kgg}
;
\draw[1cell=.9]
(a1) edge node {\Kfg\retn_{\smal}} (a2)
(a2) edge node {\kfgsk} node[swap] {\iso} (a3)
;
\draw[1cell=.9]
(a1) [rounded corners=2pt] |- ($(a2)+(-.3*\h,\u)$) -- node {\kgbk} ($(a2)+(0,\u)$) -| (a3')
 ;
\end{tikzpicture}
\end{equation}
of 
\begin{itemize}
\item the whiskered natural transformation $\Kfg\retn_{\smal}$ and
\item the natural isomorphism $\kfgsk \cn \Kfg\!\smal\! \fiso \Kgg$ given by \cref{kfgkgg_compare}.
\end{itemize}  
For \cref{kggkfgb_compare_ii}, suppose $Y$ is a proper $\GGG$-space for a finite group $G$ \pcref{def:proper_ggg}, meaning $\smal Y \in \FGTopg$ is proper.  The retraction \cref{ret_FG}
\[\Bc\!\smal\! Y \fto{\retn_{\smal Y}} \!\smal\! Y \inspace \FGTopg\]
is componentwise a pointed $G$-homotopy equivalence.  By \cref{bar_proper}, the $\FGG$-space $\Bc\!\smal\! Y$ is proper.  By \cref{thm:Kfg_inv} and the properness of $\Bc\!\smal\! Y$ and $\smal Y$, the $G$-morphism of orthogonal $G$-spectra
\[\Kfg \Bc\!\smal\! Y \fto{\Kfg\retn_{\smal Y}} \Kfg \!\smal\! Y\]
is componentwise a weak $G$-equivalence.  Thus, its composite with the $G$-homeomorphism $\kfgsk_Y$ \cref{kfgsmalkgg_iso},
\begin{equation}\label{kgbkY}
\begin{tikzpicture}[vcenter]
\def\h{3.2} \def\u{.7}
\draw[0cell]
(0,0) node (a1) {\Kfg\Bc\!\smal\! Y}
(a1)++(\h,0) node (a2) {\Kfg\!\smal\! Y} 
(a2)++(.72*\h,0) node (a3) {\Kgg Y,}
;
\draw[1cell=.9]
(a1) edge node {\Kfg\retn_{\smal Y}} (a2)
(a2) edge node {\kfgsk_Y} node[swap] {\iso} (a3)
;
\draw[1cell=.9]
(a1) [rounded corners=2pt] |- ($(a2)+(-.3*\h,\u)$) -- node {\kgbk_Y} ($(a2)+(0,\u)$) -| (a3)
 ;
\end{tikzpicture}
\end{equation}
is also componentwise a weak $G$-equivalence.  Assertion \cref{kggkfgb_compare_iii} follows from \cref{kggkfgb_compare_ii} and the fact that $|Y|_* \in \GGTopg$ is proper (\cref{reast_properGG} \cref{reast_properGG_ii}).
\end{proof}

%% file: chap/h_comparison.tex
\cref{part:kgo_shi_comp} compares the categorical parts of our equivariant $K$-theory \cref{Khgo_functors} and of the homotopical Shimakawa $K$-theory \cref{khsho_khshosg}.  For a $\Tinf$-operad $\Op$ \pcref{as:OpA}, the main goals of \cref{part:kgo_shi_comp} are
\begin{enumerate}
\item\label{hthy_comp_i} 
to compute Shimakawa $H$-theory \pcref{Sgo_twofunctor}
\[\AlglaxO \fto{\Sgo} \FGCatg\]
in terms of $H$-theory \pcref{Hgo_twofunctor}
\[\AlglaxO \fto{\Hgo} \GGCatg\]
and, conversely,
\item\label{hthy_comp_ii}
to compute $\Hgo$ in terms of $\Sgo$.
\end{enumerate}
In the topological comparison in \cref{ch:shim_top}, there is a conceptual symmetry between computing Shimakawa's homotopical functor $\Kfg\Bc$ in terms of $\Kgg$ and its converse \pcref{kfgbkgg_compare,kggkfgb_compare}.  The categorical comparisons are vastly different, with \cref{hthy_comp_i} being much simpler than \cref{hthy_comp_ii}.  On the other hand, the categorical comparisons in \cref{hthy_comp_i,hthy_comp_ii} have the advantage that they hold for arbitrary groups $G$, not just finite groups.
\begin{description}
\item[Computing $\Sgo$] 
For \cref{hthy_comp_i}, Shimakawa $H$-theory $\Sgo$ factors as the composite
\[\AlglaxO \fto{\Hgo} \GGCatg \fto{\ifgst} \FGCatg\]
with $\ifgst$ the pullback functor along the length-1 inclusion functor $\ifg \cn \FG \to \GG$ \cref{ifgGG}.  At the object level, this means that for each $\Op$-pseudoalgebra $\A$, the $\FGG$-category \pcref{sys_FGcat} 
\[\FG \fto{\Sgo\A} \Catgst\]
is the restriction of the $\GGG$-category \pcref{A_ptfunctorGG} 
\[\GG \fto{\Hgo\A} \Catgst\]
to length-1 objects in $\GG$.  The strong and $J$-theory variants are also true.  This is discussed in \cref{sec:factor_hj}.
\item[Computing $\Hgo$]
The majority of \cref{part:kgo_shi_comp} deals with the categorical comparison in \cref{hthy_comp_ii}, which computes $\Hgo$ in terms of $\Sgo$.  As the first step, after \cref{sec:factor_hj}, the rest of this chapter constructs the (strong) $H$-theory comparison 2-natural transformations
\begin{equation}\label{Pist_chi}
\begin{tikzpicture}[vcenter]
\def\h{2.3} \def\v{.7} \def\s{15}
\draw[0cell]
(0,0) node (a1) {\AlglaxO}
(a1)++(0,.2*\v) node (a1') {\phantom{\AlglaxO}}
(a1)++(\h,\v) node (a2) {\FGCatg}
(a1)++(\h,-\v) node (a3) {\GGCatg}
;
\draw[1cell=.9]
(a1) edge[bend left=\s] node {\Sgo} (a2)
(a2) edge[bend left=\s, shorten <=-.5ex] node {\smast} (a3)
(a1') edge[bend right=\s] node[swap] {\Hgo} (a3)
;
\draw[2cell]
node[between=a2 and a3 at .55, shift={(-.3*\h,0)}, rotate=-120, 2labelw={below,\Pist,1pt}] {\Rightarrow}
;
\begin{scope}[shift={(2*\h,0)}]
\draw[0cell]
(0,0) node (a1) {\AlgpspsO}
(a1)++(\h,\v) node (a2) {\FGCatg}
(a1)++(\h,-\v) node (a3) {\GGCatg}
;
\draw[1cell=.9]
(a1) edge[bend left=\s] node {\Sgosg} (a2)
(a2) edge[bend left=\s, shorten <=-.5ex] node {\smast} (a3)
(a1) edge[bend right=\s] node[swap] {\Hgosg} (a3)
;
\draw[2cell]
node[between=a2 and a3 at .55, shift={(-.3*\h,0)}, rotate=-120, 2labelw={below,\Pistsg,1pt}] {\Rightarrow}
;
\end{scope}
\end{tikzpicture}
\end{equation}
and the commutative diagrams of pointed $G$-functors
\begin{equation}\label{Pist_fact_chi}
\begin{tikzpicture}[vcenter]
\def\h{3}
\draw[0cell]
(0,0) node (a11) {\Asmaangordnbe}
(a11)++(\h,0) node (a12) {\Aangordnbe}
(a11)++(\h/2,-1.1) node (a2) {\phantom{\proAnbe}}
(a2)++(0,.15) node (a2') {\proAnbe}
;
\draw[1cell=.9]
(a11) edge node {\Pist_{\A,\angordnbe}} (a12)
(a11) edge node[swap] {\zb} (a2)
(a12) edge node {\zd} (a2)
;
\begin{scope}[shift={(1.75*\h,0)}]
\draw[0cell]
(0,0) node (a11) {\Asgsmaangordnbe}
(a11)++(\h,0) node (a12) {\Asgangordnbe}
(a11)++(\h/2,-1.1) node (a2) {\phantom{\proAnbe}}
(a2)++(0,.15) node (a2') {\proAnbe}
;
\draw[1cell=.9]
(a11) edge node {\Pistsg_{\A,\angordnbe}} (a12)
(a11) edge node[swap] {\zbsg} (a2)
(a12) edge node {\zdsg} (a2)
;
\end{scope}
\end{tikzpicture}
\end{equation}
for an arbitrary group $G$, an $\Op$-pseudoalgebra $\A$, and an object $\angordnbe \in \GG$, where
\[\begin{split}
\Asmaangordnbe &= (\smast\Sgo\A)\nbe,\\
\Aangordnbe &= (\Hgo\A)\nbe,\\
\Asgsmaangordnbe &= (\smast\Sgosg\A)\nbe, \andspace\\
\Asgangordnbe &= (\Hgosg\A)\nbe.
\end{split}\] 
\end{description}

\connection
Further properties of $\Pist$ and $\Pistsg$ are discussed in \crefrange{ch:sgoprod}{ch:special}.  The strong variant $\Pistsg$ has nicer properties than $\Pist$, and $\Pistsg$ is the desired categorical comparison.  \cref{ch:sgoprod} proves that $\zbsg$ is part of a nonequivariant adjoint equivalence \pcref{thm:zbsg_eq}.  Its left adjoint inverse $\zbsgad$ is \emph{pseudo} $G$-equivariant, but not $G$-equivariant in general \pcref{thm:zbsgad_pseudo}.  The functor $\zb$ admits a pseudo $G$-equivariant left adjoint $\zbad$, but $(\zbad,\zb)$ is not an adjoint equivalence in general \pcref{expl:sys_adjunction,expl:zbad_pseudo}.  \cref{ch:hgoprod} proves analogous results for $\zdsg$ and $\zd$: $\zdsg$ admits a pseudo $G$-equivariant left adjoint inverse \pcref{thm:zdsg_eq,thm:zdsgad_pseudo}, and $\zd$ admits a pseudo $G$-equivariant left adjoint \pcref{expl:system_adjunction,expl:zdad_pseudo}.  

For an arbitrary group $G$, a $\Uinf$-operad $\Op$ \pcref{as:OpA'}, and an $\Op$-pseudoalgebra $\A$, the strong $H$-theory comparison 
\[\smast \Sgosg\A \fto{\Pistsg_\A} \Hgosg\A\]
is componentwise a nonequivariant equivalence of categories \pcref{thm:PistAequivalence}.  \cref{ch:compgen} proves the main result of \cref{part:kgo_shi_comp}: For the $\Uinf$-operad $\Oph = \Catg(\EG,\Op)$ and an $\Oph$-pseudoalgebra of the form $\Ah = \Catg(\EG,\A)$, the strong $H$-theory comparison
\[\smast \Sgohsg\Ah \fto{\Pistsg_{\Ah}} \Hgohsg\Ah\]
is componentwise a \emph{categorical weak $G$-equivalence} \pcref{thm:pistweakgeq}.  This means that, for each subgroup $H \subseteq G$ and each object $\nbe \in \GG$, the $H$-fixed subfunctor of $\Pistsg_{\Ah,\nbe}$ is a nonequivariant equivalence of categories.  \cref{ch:special} proves that the domain and the codomain of $\Pistsg_{\Ah}$ are \emph{special} $\GGG$-categories \pcref{thm:h_special,cor:h_special_shi}.  This means that each of their Segal functors is a categorical weak $G$-equivalence \pcref{def:sp_gggcat}.  The variant $\Pist$ does not have these properties \pcref{expl:Pist_not_eq,expl:pistweq_necessity}.

\organization
This chapter consists of the following sections.

\secname{sec:factor_hj}
This section proves that Shimakawa $H$-theory $\Sgo$ is equal to $\ifgst\Hgo$.  The strong variant and the $J$-theory variant are also true.

\secname{sec:pistar}
This section constructs the (strong) $H$-theory comparison 2-natural transformations $\Pist$ and $\Pistsg$.

\secname{sec:pistar_proof}
This section proves several Lemmas used in the construction in \cref{sec:pistar}.

\secname{sec:hcomp_prod}
This section defines the $\angordnbe$-twisted product $\proAnbe$ and the zigzag factorizations of the pointed $G$-functors $\Pist_{\A,\angordnbe}$ and $\Pistsg_{\A,\angordnbe}$.

\secname{sec:twprod_fixed}
This section proves that, for a pointed $G$-category $\C$, the $G$-fixed subcategory of the $\nbe$-twisted product $\proCnbe$ splits into a product of $G_t$-fixed subcategories of $\C$ \pcref{twprod_fixed}.  The product runs over the set of $G$-orbits of the finite $G$-set $\ufs{n_1 \Cdots n_q}$, and each $G_t$ is the stabilizer of a chosen element in that $G$-orbit.  In addition to its intrinsic value, this observation is used in \cref{sec:gggcat_weq} to study special $\GGG$-categories and weak $G$-equivalences in $\GGCatg$.

\section{Computing Shimakawa $H$-Theory and $J$-Theory}
\label{sec:factor_hj}

This section proves that Shimakawa $H$-theory \pcref{Sgo_twofunctor}
\[\AlglaxO \fto{\Sgo} \FGCatg\]
factors as the composite
\[\AlglaxO \fto{\Hgo} \GGCatg \fto{\ifgst} \FGCatg\]
of $H$-theory $\Hgo$ \pcref{Hgo_twofunctor} and the pullback functor $\ifgst$ along the length-1 inclusion functor $\ifg \cn \FG \to \GG$ \cref{ifgGG}.  Thus, for each $\Op$-pseudoalgebra $\A$, the $\FGG$-category $\Sgo\A$ \pcref{sys_FGcat} can also be obtained from the $\GGG$-category $\Hgo\A$ \pcref{A_ptfunctorGG} by pulling back along $\ifg$.  The strong variant and the $J$-theory variant are also true.  Throughout this section, $G$ denotes an arbitrary group.

\secoutline
\begin{itemize}
\item \cref{smashstar} proves that smash product $\sma \cn \GG \to \FG$ induces a pullback 2-functor
\[\FGCatg \fto{\smast} \GGCatg.\]
This 2-functor admits a retraction $\ifgst$ induced by the length-1 inclusion functor $\ifg \cn \FG \to \GG$.
\item \cref{clast_smast} observes that $\smast$ and $\ifgst$ are compatible with the classifying space functor $\cla$.
\item \cref{HgoSgo} proves that Shimakawa $H$-theory $\Sgo$ is the restriction of $H$-theory $\Hgo$ along the length-1 inclusion functor $\ifg \cn \FG \to \GG$, and likewise for the strong variant.
\item \cref{JgoJgos} proves that Shimakawa $J$-theory $\Jgos$ is the restriction of $J$-theory $\Jgo$ \pcref{thm:Jgo_twofunctor} along the length-1 inclusion functor $\ifg \cn \Fsk \to \Gsk$, and likewise for the strong variant.
\end{itemize}

\subsection*{$\FGG$-Categories and $\GGG$-Categories}

To compare Shimakawa $H$-theory $\Sgo$ and $H$-theory $\Hgo$, we first connect their codomains, $\FGCatg$ and $\GGCatg$.  Recall from \cref{FG_permutative,GG_permutative,sma_symmon} the strict symmetric monoidal pointed $G$-functor
\[\big(\GG,\oplus,\ang{},\xi\big) \fto{\sma} \big(\FG,\sma,\ord{1},\xi\big)\]
between the naive permutative $G$-categories $\FG$ and $\GG$.  Also recall the length-1 inclusion pointed $G$-functor \cref{ifgGG}
\[\FG \fto{\ifg} \GG.\]
It satisfies $\sma \ifg = 1_{\FG}$ by the definition of $\sma$.  The (2-)categories in \cref{smashstar} are defined in \cref{def:ggcatg,def:GGTopg,def:fgcatg,def:ggtopg}.

\begin{lemma}\label{smashstar}
Precomposing with the pointed $G$-functors $\sma$ and $\ifg$ induces 2-functors
\[\FGCatg \fto{\smast} \GGCatg \fto{\ifgst} \FGCatg\]
and functors
\[\FGTopg \fto{\smast} \GGTopg \fto{\ifgst} \FGTopg\]
such that $\ifgst\smast = 1$ in each case.
\end{lemma}

\begin{proof}
The topological case is a part of \cref{ex:smai_adj}, restated here for convenience.  In the categorical case, the objects, 1-cells, and 2-cells of $\FGCatg$ are, respectively, pointed $G$-functors $\FG \to \Catgst$, $G$-natural transformations, and $G$-modifications.  Precomposing them with the pointed $G$-functor $\sma$
\begin{equation}\label{smashst_diag}
\begin{tikzpicture}[baseline={(a.base)}]
\def\t{23}
\draw[0cell]
(0,0) node (a) {\FG}
(a)++(-1.7,0) node (g) {\GG}
(a)++(2.6,0) node (b) {\phantom{\FG}}
(b)++(.15,0) node (b') {\Catgst}
;
\draw[1cell]
(g) edge node {\sma} (a)
(a) edge[bend left=\t] node {} (b)
(a) edge[bend right=\t] node[swap] {} (b)
;
\draw[2cell=1.2]
node[between=a and b at .3, rotate=-90] {\Rightarrow}
node[between=a and b at .67, rotate=-90] {\Rightarrow}
;
\draw[2cell=1.2]
node[between=a and b at .5, shift={(0,0)}] {\Rrightarrow}
;
\end{tikzpicture}
\end{equation}
yields pointed $G$-functors $\GG \to \Catgst$, $G$-natural transformations, and $G$-modifications.  These assignments on objects, 1-cells, and 2-cells comprise a 2-functor $\smast$ because, in each of $\GGCatg$ and $\FGCatg$, the 2-categorical structure is defined componentwise in the 2-category $\Catgst$ \cref{Catgst_iicat}.  Replacing $\sma$ by $\ifg$ proves that $\ifgst$ is a 2-functor.  The equality $\ifgst\smast = 1$ follows from the equality $\sma\ifg = 1_{\FG}$.  
\end{proof}

\cref{clast_smast} records the fact that the functors $\smast$ and $\ifgst$ \pcref{smashstar} commute with the functors $\clast$ \pcref{GGCatg_GGTopg,ggcatg_ggtopg} induced by the classifying space functor $\cla \cn \Cat \to \Top$ \cref{classifying_space}.

\begin{lemma}\label{clast_smast}
In the diagram
\begin{equation}\label{Bsmai}
\begin{tikzpicture}[vcenter]
\def\v{-1.4} \def\s{.5ex}
\draw[0cell]
(0,0) node (a11) {\FGCatg}
(a11)++(2.7,0) node (a12) {\FGTopg}
(a11)++(0,\v) node (a21) {\GGCatg}
(a12)++(0,\v) node (a22) {\GGTopg}
;
\draw[1cell=.9]
(a11) edge node {\clast} (a12)
(a21) edge node {\clast} (a22)
(a11) edge[transform canvas={xshift=-\s}] node[swap] {\smast} (a21)
(a12) edge[transform canvas={xshift=-\s}] node[swap] {\smast} (a22)
(a21) edge[transform canvas={xshift=\s}] node[swap] {\ifgst} (a11)
(a22) edge[transform canvas={xshift=\s}] node[swap] {\ifgst} (a12)
;
\end{tikzpicture}
\end{equation}
of functors, the equalities
\[\begin{split}
\clast \smast &= \smast \clast \andspace\\ 
\clast \ifgst &= \ifgst \clast
\end{split}\]
hold.
\end{lemma}

\begin{proof}
The desired equalities follow from the fact that each $\clast$ is given by postcomposition with $\cla$, while $\smast$ and $\ifgst$ are given by precomposition with $\sma$ and $\ifg$.
\end{proof}

\subsection*{Factoring $\Sgo$ through $\Hgo$}

Recall
\begin{itemize}
\item the 2-categories $\AlglaxO$ and $\AlgpspsO$ \pcref{oalgps_twocat},
\item (strong) $H$-theory $\Hgo$ and $\Hgosg$ \pcref{Hgo_twofunctor}, and
\item Shimakawa (strong) $H$-theory $\Sgo$ and $\Sgosg$ \pcref{Sgo_twofunctor}. 
\end{itemize} 
\cref{HgoSgo} proves that $H$-theory $\Hgo$ restricts to Shimakawa $H$-theory $\Sgo$ along the pullback functor $\ifgst$ \pcref{smashstar}.  Thus, each $\FGG$-category in the image of Shimakawa (strong) $H$-theory can also be obtained from a $\GGG$-category in the image of (strong) $H$-theory by pulling back along $\ifg$.

\begin{proposition}\label{HgoSgo}\index{Shimakawa H-theory@Shimakawa $H$-theory!factors through H-theory@factors through $H$-theory}
In each of the diagrams
\begin{equation}\label{HgoSgo_diag}
\begin{tikzpicture}[vcenter]
\def\e{-.4} \def\h{1.5} \def\g{2.2} \def\v{1.4} \def\s{.8}
\draw[0cell=.9]
(0,0) node (a1) {\phantom{A^l_l}}
(a1)++(0,-.04) node (a1') {\AlglaxO}
(a1)++(\e,0) node (a1'') {\phantom{A^l_l}}
(a1)++(\h,\v/2) node (a2) {\FGCatg}
(a1)++(\h,-\v/2) node (a3) {\GGCatg}
(a2)++(\g,0) node (a4) {\FGTopg}
(a3)++(\g,0) node (a5) {\GGTopg}
;
\draw[1cell=.85]
(a1'') [rounded corners=2pt] |- node[pos=\s] {\Sgo} (a2)
;
\draw[1cell=.85]
(a1'') [rounded corners=2pt] |- node[pos=\s] {\Hgo} (a3)
;
\draw[1cell=.85]
(a3) edge node[swap] {\ifgst} (a2)
(a2) edge node {\clast} (a4)
(a3) edge node {\clast} (a5)
(a5) edge node {\ifgst} (a4)
;
\draw[0cell=.9]
(a1)++(2.3*\h+\g,0) node (b1) {\AlgpspsO}
(b1)++(\e,0) node (b1') {\phantom{A^l_l}}
(b1)++(\h,\v/2) node (b2) {\FGCatg}
(b1)++(\h,-\v/2) node (b3) {\GGCatg}
(b2)++(\g,0) node (b4) {\FGTopg}
(b3)++(\g,0) node (b5) {\GGTopg}
;
\draw[1cell=.85]
(b1') [rounded corners=2pt] |- node[pos=\s] {\Sgosg} (b2)
;
\draw[1cell=.85]
(b1') [rounded corners=2pt] |- node[pos=\s] {\Hgosg} (b3)
;
\draw[1cell=.85]
(b3) edge node[swap] {\ifgst} (b2)
(b2) edge node {\clast} (b4)
(b3) edge node {\clast} (b5)
(b5) edge node {\ifgst} (b4)
;
\end{tikzpicture}
\end{equation}
with $\Op$ a $\Tinf$-operad \pcref{as:OpA}, the left square of 2-functors commutes, and the right square of functors commutes.
\end{proposition}

\begin{proof}
We consider the left diagram; the strong variant is proved in the same way by restricting to strong systems and $\Op$-pseudomorphisms.  The right square commutes by \cref{clast_smast}.  The left square commutes for the following reasons.
\begin{description}
\item[Objects]
For an $\Op$-pseudoalgebra $\A$ \pcref{def:pseudoalgebra}, the pointed $G$-category $\Anbe$ of $\nbe$-systems \pcref{def:nbeta_gcat} reduces to the pointed $G$-category $\Aordnbe$ of $\nbeta$-systems \pcref{def:nsys_gcat} if $\nbe \in \GG$ is an object of length 1.  As we point out in \cref{expl:nsys} \eqref{expl:nsys_ii}, the key point is that the commutativity axiom \cref{system_commutativity} of $\nbe$-systems only applies if $\nbe$ has length $>1$.  On morphisms between objects of length 1 in $\GG$, \cref{def:AfangpsiGG} reduces to \cref{def:Apsi}.  Thus, there is an equality of $\FGG$-categories.
\[\ifgst \Hgo\A = \Sgo\A.\]
\item[1-cells]
For a lax $\Op$-morphism $f$ between $\Op$-pseudoalgebras \pcref{def:laxmorphism}, the $G$-natural transformation $\Hgo f$ \cref{hgo_f} reduces to $\Sgo f$ \cref{sgo_1cell} on length-1 objects of $\GG$.  
\item[2-cells]
For an $\Op$-transformation $\omega$ \pcref{def:algtwocells}, the $G$-modification $\Hgo\omega$ \cref{hgo_omega} reduces to $\Sgo\omega$ \cref{sgo_2cell} on length-1 objects of $\GG$.
\end{description}
This proves the equality 
\[\ifgst \Hgo = \Sgo\]
of 2-functors.  
\end{proof}

\subsection*{Factoring $\Jgos$ through $\Jgo$}

Recall the 2-categories $\GGCatii$ and $\FGCat$ \pcref{def:GGCat,def:fgcat} and the categories $\GGTopii$ and $\FGTop$ \pcref{def:ggtop_smc}.  \cref{JgoJgos} is the analogue of \cref{HgoSgo} involving $J$-theory $\Jgo$ \pcref{thm:Jgo_twofunctor}, Shimakawa $J$-theory $\Jgos$ \pcref{def:jgos}, and the length-1 inclusion functor $\ifg \cn \Fsk \to \Gsk$ \cref{ifg}.

\begin{proposition}\label{JgoJgos}\index{Shimakawa J-theory@Shimakawa $J$-theory!factors through J-theory@factors through $J$-theory}
In each of the diagrams
\begin{equation}\label{JBi}

\end{equation}
with $\Op$ a $\Tinf$-operad \pcref{as:OpA}, the left square of 2-functors commutes, and the right square of functors commutes.
\end{proposition}

\begin{proof}
The equality $\Jgos = \ifgst\Jgo$ follows from the following commutative diagram of 2-functors.
\begin{equation}\label{JHi}
%
\end{equation}
\begin{itemize}
\item The two $\igst$ are the 2-equivalences in \cref{thm:ggcat_ggcatg_iieq,thm:fgcat_fgcatg_iieq} induced by the pointed full subcategory inclusions $\Gsk \to \GG$ \cref{ig} and $\Fsk \to \FG$ \cref{ig_FG}.
\item The bottom and top regions commute by \cref{jgohgoigst,jgos_jgossg}.  
\item The left triangle commutes by \cref{HgoSgo}.  
\item The right square commutes because it is induced by precomposition with the commutative diagram
\begin{equation}\label{ii_diagram}
%
\end{equation}
involving the full subcategory inclusions defined in \cref{ifg,ifgGG,ig,ig_FG}.
\end{itemize}
The equality of functors 
\[\clast\ifgst = \ifgst\clast\] 
follows from the fact that each $\clast$ is given by postcomposition with $\cla$, while each $\ifgst$ is given by precomposition with $\ifg$.  The strong variant is proved in the same way by replacing $\Hgo$ and $\Sgo$ with their strong variants $\Hgosg$ and $\Sgosg$.
\end{proof}

\section{$H$-Theory Comparison 2-Natural Transformations}
\label{sec:pistar}

This section defines the $H$-theory comparison 2-natural transformation $\Pist$ that compares Shimakawa $H$-theory \pcref{Sgo_twofunctor} and $H$-theory \pcref{Hgo_twofunctor}:
\[\begin{split}
& \AlglaxO \fto{\Sgo} \FGCatg \andspace\\
& \AlglaxO \fto{\Hgo} \GGCatg.
\end{split}\]
The strong variant $\Pistsg$ compares Shimakawa strong $H$-theory $\Sgosg$ and strong $H$-theory $\Hgosg$:
\[\begin{split}
& \AlgpspsO \fto{\Sgosg} \FGCatg \andspace\\
& \AlgpspsO \fto{\Hgosg} \GGCatg.
\end{split}\]
Several statements used in \cref{def:Pistar} are proved in \cref{sec:pistar_proof}.  Recall the 2-functor $\smast$ \pcref{smashstar}. 

\begin{definition}[$H$-Theory Comparison]\label{def:Pistar}
For a $\Tinf$-operad $\Op$ \pcref{as:OpA}, the \emph{$H$-theory comparison}\index{H-theory comparison@$H$-theory comparison} is the 2-natural transformation
\begin{equation}\label{Pistar_twonat}
\begin{tikzpicture}[vcenter]
\def\h{2.3} \def\v{.7} \def\s{15}
\draw[0cell]
(0,0) node (a1) {\AlglaxO}
(a1)++(0,.2*\v) node (a1') {\phantom{\AlglaxO}}
(a1)++(\h,\v) node (a2) {\FGCatg}
(a1)++(\h,-\v) node (a3) {\GGCatg}
;
\draw[1cell=.9]
(a1) edge[bend left=\s] node {\Sgo} (a2)
(a2) edge[bend left=\s, shorten <=-.5ex] node {\smast} (a3)
(a1') edge[bend right=\s] node[swap] {\Hgo} (a3)
;
\draw[2cell]
node[between=a2 and a3 at .55, shift={(-.3*\h,0)}, rotate=-120, 2labelw={below,\Pist,1pt}] {\Rightarrow}
;
\end{tikzpicture}
\end{equation}
defined as follows.  A strong variant is defined in \cref{Pistarsg_twonat}.  For each $\Op$-pseudoalgebra $(\A,\gaA,\phiA)$ \pcref{def:pseudoalgebra}, the $\A$-component $G$-natural transformation
\begin{equation}\label{PistA}
\begin{tikzpicture} [vcenter]
\def\s{25}
\draw[0cell]
(0,0) node (a1) {\GG}
(a1)++(2.2,0) node (a2) {\phantom{\GG}}
(a2)++(.15,0) node (a2') {\Catgst}
;
\draw[1cell=.8]
(a1) edge[bend left=\s] node {\smast \Sgo\A} (a2)
(a1) edge[bend right=\s] node[swap] {\Hgo\A} (a2) 
;
\draw[2cell]
node[between=a1 and a2 at .35, rotate=-90, 2label={above,\Pist_\A}] {\Rightarrow}
;
\end{tikzpicture}
\end{equation}
has, for each object $\angordnbe \in \GG$, $\angordnbe$-component pointed $G$-functor
\begin{equation}\label{PistAnbe}
(\smast\Sgo\A)\angordnbe = \Asmaangordnbe
\fto{\Pist_{\A,\angordnbe}} \Aangordnbe = (\Hgo\A)\angordnbe
\end{equation}
defined as follows.  
\begin{itemize}
\item In the domain in \cref{PistAnbe}, the pointed finite $G$-set $\sma\angordnbe \in \FG$ is defined in \cref{smash_GGobjects}.  The category $\Asmaangordnbe$ of $(\sma\angordnbe)$-systems is defined in \cref{def:nsys,def:nsys_morphism,def:nsys_gcat}.
\item In the codomain in \cref{PistAnbe}, the category $\Aangordnbe$ of $\angordnbe$-systems is defined in \cref{def:nsystem,def:nsystem_morphism,def:nbeta_gcat}.
\end{itemize}
\begin{description}
\item[Base cases] If $\angordnbe = \vstar$, the initial-terminal basepoint in $\GG$, then $\sma\vstar = \ordz$.  The $\vstar$-component of $\Pist_\A$ is the identity functor on the terminal $G$-category $\bone$:
\begin{equation}\label{PistA_vstar}
\bone = \Aordz \fto{\Pist_{\A,\vstar} \,=\, 1} \sys{\A}{\vstar} = \bone.
\end{equation}
If $\angordnbe = \ang{}$, the empty tuple in $\GG$, then $\sma\ang{} = \ord{1}$.  The $\ang{}$-component
\begin{equation}\label{PistA_empty}
\Aordone \fto[\iso]{\Pist_{\A,\ang{}}} \sys{\A}{\ang{}} = \A
\end{equation}
 of $\Pist_\A$ is given by the inverse of the pointed $G$-isomorphism in \cref{ex:nsys_zero_one} \eqref{ex:nsys_one}.  In other words, $\Pist_{\A,\ang{}}$ sends
\begin{itemize}
\item a $\ord{1}$-system to its $\{1\}$-component object and 
\item a morphism of $\ord{1}$-systems to its $\{1\}$-component morphism.
\end{itemize}
\end{description}

In the rest of this definition, we assume that the object
\[\angordnbe = \sordi{n}{\be}{j}_{j\in \ufs{q}} \in \GG \setminus \{\vstar, \ang{}\}\] 
has length $q>0$.  Recall from \cref{mal-sma-nbe,smash_GGobjects} that the smash product of pointed finite $G$-sets is given by 
\[\sma\angordnbe = \ord{n_1 \Cdots n_q}^{\sma_{j \in \ufs{q}}\, \be_j} \in \FG.\]
\begin{description}
\item[Component objects]
Given a $(\sma\angordnbe)$-system $(a,\gl)$ in $\A$, the desired $\angordnbe$-system
\begin{equation}\label{agltimes}
\Pist_{\A,\angordnbe} (a,\gl) = (a^\stimes,\gl^\stimes)
\end{equation}
is defined as follows.  For each marker $\angs = \ang{s_j \subseteq \ufs{n}_j}_{j \in \ufs{q}}$, we first use the lexicographic ordering \cref{lex_bijection} to define the product subset
\begin{equation}\label{angstimes}
\angstimes = \txprod_{j\in \ufs{q}}\, s_j \subseteq 
\txprod_{j\in \ufs{q}}\, \ufs{n}_j = \ufs{n_1 \Cdots n_q}.
\end{equation}
The $\angs$-component object \cref{a_angs} of the $\angordnbe$-system $(a^\stimes,\gl^\stimes)$ is defined as the $\angstimes$-component object of $(a,\gl)$:
\begin{equation}\label{PistAnbe_aangs}
a^\stimes_{\angs} = a_{\angstimes} \in \A.
\end{equation}
\item[Gluing]
Given an object $x \in \Op(r)$ with $r \geq 0$, an index $k \in \ufs{q}$, a marker $\angs = \ang{s_j \subseteq \ufs{n}_j}_{j \in \ufs{q}}$, and a partition
\[s_k = \coprod_{i \in \ufs{r}}\, s_{k,i} \subseteq \ufs{n}_k,\]
the gluing morphism \cref{gluing-morphism} of $(a^\stimes,\gl^\stimes)$ at $(x; \angs, k, \ang{s_{k,i}}_{i \in \ufsr})$ is defined by the following commutative diagram.
\begin{equation}\label{gltimes_component}
\begin{tikzpicture}[vcenter]
\def\v{-1.2}
\draw[0cell]
(0,0) node (a1) {\gaA_r\big(x \sscs \ang{a^\stimes_{\angs \,\compk\, s_{k,i}}}_{i \in \ufs{r}} \big)}
(a1)++(5.5,0) node (a2) {a^\stimes_{\angs}}
(a1)++(0,\v) node (b1) {\gaA_r\big(x \sscs \ang{a_{(\angs \,\compk\, \ski)^\stimes}}_{i \in \ufs{r}} \big)}
(a2)++(0,\v) node (b2) {a_{\angstimes}}
;
\draw[1cell=.9]
(a1) edge[equal,shorten >=-1ex] (b1)
(a2) edge[equal] (b2)
(a1) edge node {\gl^\stimes_{x;\, \angs,\, k, \ang{s_{k,i}}_{i \in \ufs{r}}}} (a2)
(b1) edge node {\gl_{x;\, \angstimes,\, \ang{(\angs \,\compk\, \ski)^\stimes}_{i \in \ufs{r}}}} (b2)
;
\end{tikzpicture}
\end{equation}
The bottom horizontal arrow in \cref{gltimes_component} is the gluing morphism of $(a,\gl)$ at the object $x \in \Op(r)$, the subset $\angstimes$ \cref{angstimes}, and the partition
\[\angstimes = \coprod_{i \in \ufs{r}}\, \angscompkskitimes\]
of $\angstimes$ by the subsets
\begin{equation}\label{angscompkskitimes}
\angscompkskitimes = 
s_1 \times \Cdots \times s_{k-1} \times s_{k,i} \times s_{k+1} \times \Cdots \times s_q.
\end{equation}
\cref{atimes_system} proves that the pair $(a^\stimes,\gl^\stimes)$ is an $\angordnbe$-system in $\A$.
\item[Morphisms]
Given a morphism of $(\sma\angordnbe)$-systems in $\A$
\[(a,\gl^a) \fto{\theta} (b,\gl^b),\]
the morphism of $\angordnbe$-systems
\[\Pist_{\A,\angordnbe} (a,\gl^a) 
\fto{\Pist_{\A,\angordnbe} \theta = \theta^\stimes} 
\Pist_{\A,\angordnbe} (b,\gl^b)\]
has, for each marker $\angs = \ang{s_j \subseteq \ufs{n}_j}_{j \in \ufs{q}}$, $\angs$-component morphism \cref{theta_angs} defined as the $\angstimes$-component morphism of $\theta$:
\begin{equation}\label{thetatimes_angs}
a^\stimes_{\angs} = a_{\angstimes} 
\fto{\theta^\stimes_{\angs} = \theta_{\angstimes}}
b^\stimes_{\angs} = b_{\angstimes}.
\end{equation}
The unity axiom \cref{nsystem_mor_unity} and the compatibility axiom \cref{nsystem_mor_compat} for a morphism of $\angordnbe$-systems hold for $\theta^\stimes$ by the corresponding axioms, \cref{nsys_mor_unity,nsys_mor_compat}, for $\theta$, \cref{gltimes_component}, and \cref{thetatimes_angs}.
\item[Pointed functoriality]  
The assignment \cref{agltimes} is pointed, sending the base $(\sma\angordnbe)$-system $(\zero,1_\zero)$ to the base $\angordnbe$-system $(\zero,1_\zero)$, by \cref{PistAnbe_aangs,gltimes_component}.  Functoriality of $\Pist_{\A,\angordnbe}$ follows from \cref{thetatimes_angs} and the fact that identities and composition of morphisms of $(\sma\angordnbe)$-systems and $\angordnbe$-systems are defined componentwise in $\A$.
\end{description}
\cref{PistA_natural,PistA_gnat} show that, as $\angordnbe$ varies in $\GG$, $\Pist_\A$ \cref{PistA} is a $G$-natural transformation.  \cref{Pist_iinat} shows that, as $\A$ varies in $\AlglaxO$, $\Pist$ \cref{Pistar_twonat} is a 2-natural transformation.

\parhead{Strong variant}.  The \emph{strong $H$-theory comparison}\index{H-theory comparison@$H$-theory comparison!strong} is the 2-natural transformation
\begin{equation}\label{Pistarsg_twonat}
\begin{tikzpicture}[vcenter]
\def\h{2.3} \def\v{.7} \def\s{15}
\draw[0cell]
(0,0) node (a1) {\AlgpspsO}
(a1)++(\h,\v) node (a2) {\FGCatg}
(a1)++(\h,-\v) node (a3) {\GGCatg}
;
\draw[1cell=.9]
(a1) edge[bend left=\s] node {\Sgosg} (a2)
(a2) edge[bend left=\s, shorten <=-.5ex] node {\smast} (a3)
(a1) edge[bend right=\s] node[swap] {\Hgosg} (a3)
;
\draw[2cell]
node[between=a2 and a3 at .55, shift={(-.3*\h,0)}, rotate=-120, 2labelw={below,\Pistsg,1pt}] {\Rightarrow}
;
\end{tikzpicture}
\end{equation}
defined by restricting \Crefrange{PistAnbe}{thetatimes_angs}  
to strong systems.  For each $\Op$-pseudoalgebra $\A$, the $\A$-component $G$-natural transformation
\begin{equation}\label{PistsgA}
\begin{tikzpicture} [vcenter]
\def\s{25}
\draw[0cell]
(0,0) node (a1) {\GG}
(a1)++(2.2,0) node (a2) {\phantom{\GG}}
(a2)++(.15,0) node (a2') {\Catgst}
;
\draw[1cell=.85]
(a1) edge[bend left=\s] node {\smast \Sgosg\A} (a2)
(a1) edge[bend right=\s] node[swap] {\Hgosg\A} (a2) 
;
\draw[2cell]
node[between=a1 and a2 at .35, rotate=-90, 2label={above,\Pistsg_\A}] {\Rightarrow}
;
\end{tikzpicture}
\end{equation}
has, for each object $\angordnbe \in \GG$, $\angordnbe$-component pointed $G$-functor
\begin{equation}\label{PistsgAnbe}
(\smast\Sgosg\A)\angordnbe = \Asgsmaangordnbe 
\fto{\Pistsg_{\A,\angordnbe}} \Asgangordnbe = (\Hgosg\A)\angordnbe.
\end{equation}
This is well defined because, if each gluing morphism $\gl$ is invertible, then the same is true for $\gl^\stimes$ in the diagram \cref{gltimes_component}.  This finishes the definition of the 2-natural transformation $\Pistsg$.
\end{definition}

\section{Proof of 2-Naturality}
\label{sec:pistar_proof}

This section proves \cref{atimes_system,PistA_natural,PistA_gnat,Pist_iinat}, which are used in \cref{def:Pistar} for the $H$-theory comparison 2-natural transformation $\Pist$ \cref{Pistar_twonat} and the strong variant $\Pistsg$ \cref{Pistarsg_twonat}.  Each of these lemmas is stated for $\Pist$.  The strong variant using $\Pistsg$ is also true with the same proof, after restricting to strong systems.  We remind the reader that $\Op$ denotes a $\Tinf$-operad \pcref{as:OpA}.

\secoutline
\begin{itemize}
\item \cref{atimes_system} proves that $\Pist_{\A,\angordnbe}$ \cref{PistAnbe} is a pointed functor.
\item \cref{PistA_natural} proves that $\Pist_\A$ \cref{PistA} is a natural transformation.
\item \cref{PistA_gnat} proves that $\Pist_\A$ is $G$-equivariant.
\item \cref{Pist_iinat} proves that $\Pist$ \cref{Pistar_twonat} is a 2-natural transformation.
\end{itemize}

\begin{lemma}\label{atimes_system}
For each $\Op$-pseudoalgebra $\A$ and object $\angordnbe \in \GG$, the assignment in \cref{PistAnbe}
\[\Asmaangordnbe \fto{\Pist_{\A,\angordnbe}} \Aangordnbe\]
is a pointed functor.
\end{lemma}

\begin{proof}
Once we prove that the object assignment \cref{agltimes}
\[\Asmaangordnbe \ni (a,\gl) \mapsto \Pist_{\A,\angordnbe} (a,\gl) = (a^\stimes,\gl^\stimes) \in \Aangordnbe\]
is well defined, the explanation between \cref{thetatimes_angs,Pistarsg_twonat} proves the pointed functoriality of  $\Pist_{\A,\angordnbe}$.  To prove that $(a^\stimes, \gl^\stimes)$, defined in \cref{PistAnbe_aangs,gltimes_component}, is an $\angordnbe$-system in $\A$, we verify the axioms \crefrange{system_obj_unity}{system_commutativity} in \cref{def:nsystem}. 
\begin{description}
\item[Object unity, naturality, unity, equivariance, and associativity]  
The axioms \crefrange{system_obj_unity}{system_associativity} for $(a^\stimes,\gl^\stimes)$ follow from the corresponding axioms \crefrange{nsys_obj_unity}{nsys_associativity} for the $(\sma\angordnbe)$-system $(a,\gl)$.
\item[Commutativity]  
To verify the commutativity axiom \cref{system_commutativity} for $(a^\stimes,\gl^\stimes)$, suppose we are given a pair of objects
\[(x,y) \in \Op(r) \times \Op(t)\]
with $r,t \geq 0$, a marker $\angs = \ang{s_j \subseteq \ufsn_j}_{j \in \ufsq}$, two distinct indices $k,\ell \in \ufs{q}$, and partitions
\[s_k = \coprod_{i \in \ufs{r}} s_{k,i} \subseteq \ufs{n}_k \andspace 
s_\ell = \coprod_{p \in \ufs{t}} s_{\ell,p} \subseteq \ufs{n}_\ell.\]
We use the notation
\[\begin{split}
\ang{s_{(k,i)}} &= \angs \compk s_{k,i} = \left(s_1, \ldots, s_{k-1}, s_{k,i}, s_{k+1}, \ldots, s_q\right)\\
\ang{s^{(\ell,p)}} &= \angs \compell s_{\ell,p} = \big(s_1, \ldots, s_{\ell-1}, s_{\ell,p}, s_{\ell+1}, \ldots, s_q\big)\\
\ang{s_{(k,i)}^{(\ell,p)}} &= \angs \compk s_{k,i} \compell s_{\ell,p} =
\big(s_1, \ldots, s_{k,i}, \ldots, s_{\ell,p}, \ldots, s_q\big)
\end{split}\]
with the subscript $(k,i)$ indicating $s_{k,i}$ in the $k$-th entry and the superscript $(\ell,p)$ indicating $s_{\ell,p}$ in the $\ell$-th entry.  The subset $\angstimes$ \cref{angstimes} admits the following partitions.
\[\begin{split}
\angstimes &= \coprod_{i \in \ufs{r}} \ang{s_{(k,i)}}^\stimes
= \coprod_{i \in \ufs{r}} \coprod_{p \in \ufs{t}} \ang{s_{(k,i)}^{(\ell,p)}}^\stimes\\ 
&= \coprod_{p \in \ufs{t}} \ang{s^{(\ell,p)}}^\stimes
= \coprod_{p \in \ufs{t}} \coprod_{i \in \ufs{r}} \ang{s_{(k,i)}^{(\ell,p)}}^\stimes
\end{split}\]
The commutativity diagram \cref{system_commutativity} for $(a^\stimes,\gl^\stimes)$ is the boundary of the following diagram in $\A$.
\begin{equation}\label{commutativity_cross}
\begin{tikzpicture}[vcenter]
\def\s{2em} \def\h{3.5} \def\g{.8} \def\u{-1.2} \def\v{-1.5}
\draw[0cell=.7]
(0,0) node (a1) {\gaA_{tr} \big((x \intr y)\twist_{t,r} \sscs \ang{\ang{a_{\ang{s_{(k,i)}^{(\ell,p)}}^\stimes}}_{i \in \ufs{r}}}_{p \in \ufs{t}}\big)}
(a1)++(-\h,\u) node (a21) {\gaA_{tr} \big(y \intr x \sscs \ang{\ang{a_{\ang{s_{(k,i)}^{(\ell,p)}}^\stimes}}_{i \in \ufs{r}}}_{p \in \ufs{t}}\big)}
(a1)++(\h,\u) node (a22) {\gaA_{rt} \big(x \intr y \sscs \ang{\ang{a_{\ang{s_{(k,i)}^{(\ell,p)}}^\stimes}}_{p \in \ufs{t}}}_{i \in \ufs{r}}\big)}
(a21)++(0,\v) node (a31) {\gaA_t\big(y \sscs \bang{\gaA_r\big(x \sscs \ang{a_{\ang{s_{(k,i)}^{(\ell,p)}}^\stimes}}_{i \in \ufs{r}}\big)}_{p \in \ufs{t}}\big)}
(a22)++(0,\v) node (a32) {\gaA_r\big(x \sscs \bang{\gaA_t\big(y \sscs \ang{a_{\ang{s_{(k,i)}^{(\ell,p)}}^\stimes}}_{p \in \ufs{t}}\big)}_{i \in \ufs{r}}\big)}
(a31)++(\g,\v) node (a41) {\gaA_t\big(y \sscs \ang{a_{\ang{s^{(\ell,p)}}^\stimes}}_{p \in \ufs{t}}\big)}
(a32)++(-\g,\v) node (a42) {\gaA_r\big(x \sscs \ang{a_{\ang{s_{(k,i)}}^\stimes}}_{i \in \ufs{r}}\big)}
(a31)++(\h,0) node (a5) {a_{\angstimes}}
;
\draw[1cell=.7]
(a1) edge[transform canvas={xshift=-\s}, shorten >=.5ex] node[swap] {\gaA_{tr}(\pcom_{r,t} \sscs 1^{tr})} (a21)
(a21) edge node[swap] {\big(\phiA_{(t;\, r,\ldots,r)} \big)^{-1}} (a31)
(a31) edge node[swap,pos=.2] {\gaA_t\big(y \sscs \ang{\gl_{x;\, \ang{s^{(\ell,p)}}^\stimes}}_{p \in \ufs{t}}\big)} (a41)
(a41) edge node[swap] {\gl_{y;\, \angstimes}} (a5)
(a1) edge[equal, transform canvas={xshift=\s}, shorten >=.5ex] (a22)
(a22) edge node {\big(\phiA_{(r;\, t,\ldots,t)} \big)^{-1}} (a32)
(a32) edge node[pos=.2] {\gaA_r\big(x \sscs \ang{\gl_{y;\, \ang{s_{(k,i)}}^\stimes}}_{i \in \ufs{r}}\big)} (a42)
(a42) edge node {\gl_{x;\, \angstimes}} (a5)
(a1) edge node[pos=.3] {\gl_{(x \intr y)\twist_{t,r}}} (a5)
(a21) edge node {\gl_{y \intr x}} (a5)
(a22) edge node[swap] {\gl_{x \intr y}} (a5)
;
\end{tikzpicture}
\end{equation}
In the diagram above, the seven gluing morphisms of $(a,\gl)$, with their full subscripts, are given as follows.
\[\begin{gathered}
\gl_{(x \intr y)\twist_{t,r}} = 
\gl_{(x \intr y)\twist_{t,r};\, \angstimes,\, \ang{\ang{ \ang{s_{(k,i)}^{(\ell,p)}}^\stimes }_{i \in \ufs{r}}}_{p \in \ufs{t}}}\\
\begin{aligned}
\gl_{y \intr x} &= 
\gl_{y \intr x;\, \angstimes,\, \ang{\ang{ \ang{s_{(k,i)}^{(\ell,p)}}^\stimes }_{i \in \ufs{r}}}_{p \in \ufs{t}}} & \gl_{x \intr y} &= 
\gl_{x \intr y;\, \angstimes,\, \ang{\ang{ \ang{s_{(k,i)}^{(\ell,p)}}^\stimes }_{p \in \ufs{t}}}_{i \in \ufs{r}}}\\
\gl_{x;\, \ang{s^{(\ell,p)}}^\stimes} &=
\gl_{x;\, \ang{s^{(\ell,p)}}^\stimes ,\, \ang{\ang{s_{(k,i)}^{(\ell,p)}}^\stimes}_{i \in \ufs{r}}} &
\gl_{y;\, \ang{s_{(k,i)}}^\stimes} &= 
\gl_{y;\, \ang{s_{(k,i)}}^\stimes ,\, \ang{\ang{s_{(k,i)}^{(\ell,p)}}^\stimes}_{p \in \ufs{t}}}\\
\gl_{y;\, \angs^\stimes} &= 
\gl_{y;\, \angs^\stimes,\, \ang{\ang{s^{(\ell,p)}}^\stimes}_{p \in \ufs{t}}} &
\gl_{x;\, \angs^\stimes} &= 
\gl_{x;\, \angs^\stimes,\, \ang{\ang{s_{(k,i)}}^\stimes}_{i \in \ufs{r}}}
\end{aligned}
\end{gathered}\]
In the previous diagram, the upper left and upper right triangles commute by, respectively, the naturality axiom \cref{nsys_naturality} and the equivariance axiom \cref{nsys_equivariance} for the $(\sma\angordnbe)$-system $(a,\gl)$.  The lower left and lower right quadrilaterals commute by the associativity axiom \cref{nsys_associativity} for $(a,\gl)$.\qedhere
\end{description}
\end{proof}

\cref{PistA_natural,PistA_gnat} together show that $\Pist_\A$ is a $G$-natural transformation.

\begin{lemma}\label{PistA_natural}
For each $\Op$-pseudoalgebra $\A$, the assignment in \cref{PistA}
\begin{equation}\label{PistA_iicell}
\begin{tikzpicture}[vcenter]
\def\s{25}
\draw[0cell]
(0,0) node (a1) {\GG}
(a1)++(2.2,0) node (a2) {\phantom{\GG}}
(a2)++(.15,0) node (a2') {\Catgst}
;
\draw[1cell=.9]
(a1) edge[bend left=\s] node {\smast \Sgo\A} (a2)
(a1) edge[bend right=\s] node[swap] {\Hgo\A} (a2) 
;
\draw[2cell]
node[between=a1 and a2 at .35, rotate=-90, 2label={above,\Pist_\A}] {\Rightarrow}
;
\end{tikzpicture}
\end{equation}
is a natural transformation.
\end{lemma}

\begin{proof}
The components of $\Pist_\A$ are pointed functors by \cref{atimes_system}.  We verify that, for each morphism \cref{GG_morphisms}
\[\angordmal \fto{\upom} \angordnbe \inspace \GG,\]
the naturality diagram \cref{PistA_nat} of pointed functors commutes, where $\sma\upom$, $\Aupom$, $\Asmaupom$, and $\Pist_{\A,\angordnbe}$ are defined in, respectively, \cref{smash_fpsiFG,Aupom,psitil_f,PistAnbe}.
\begin{equation}\label{PistA_nat}
\begin{tikzpicture}[vcenter]
\def\v{-1.5}
\draw[0cell]
(0,0) node (a11) {\Asmaangordmal}
(a11)++(3,0) node (a12) {\Aangordmal}
(a11)++(0,\v) node (a21) {\Asmaangordnbe}
(a12)++(0,\v) node (a22) {\Aangordnbe}
;
\draw[1cell=.9]
(a11) edge node {\Pist_{\A,\angordmal}} (a12)
(a12) edge node {\Aupom} (a22)
(a11) edge[transform canvas={xshift=1ex}] node[swap] {\Asmaupom} (a21)
(a21) edge node {\Pist_{\A,\angordnbe}} (a22)
;
\end{tikzpicture}
\end{equation}

\parhead{Base cases}.  First, we consider the following two base cases.
\begin{enumerate}
\item If $\upom$ is the 0-morphism, factoring through the initial-terminal basepoint $\vstar \in \GG$, then $\sma\upom$ is also the 0-morphism, factoring through $\ordz \in \FG$.  In this case, each of $\Asmaupom$ and $\Aupom$ is the constant functor at the basepoint in its codomain.  The diagram \cref{PistA_nat} commutes because $\Pist_{\A,\angordnbe}$ is also pointed.  Thus, we assume that $\upom$ is not the 0-morphism.  In particular, neither $\angordmal$ nor $\angordnbe$ is the basepoint $\vstar \in \GG$.
\item If $\angordnbe = \ang{}$, the empty tuple in $\GG$, then $\angordmal = \ang{}$, and $\upom$ is the identity morphism, since we are assuming that $\upom$ is not the 0-morphism.  In this case, each of $\Asmaupom$ and $\Aupom$ is the identity functor, and the diagram \cref{PistA_nat} commutes.
\end{enumerate}
With the two base cases taken care of, we assume that $\upom$ has the form \cref{fangpsiGG}
\[\angordmal \fto{\upom \,=\, (f,\angpsi)} \angordnbe,\]
where $\angordmal$ and $\angordnbe$ have lengths $p \geq 0$ and $q>0$, respectively.  The case $p=0$ is covered in the following argument with $\sma\ang{} = \ord{1}$ and $\Pist_{\A,\ang{}}$ given by the isomorphism in \cref{PistA_empty}.
\begin{description}
\item[Component objects]  
To prove that the diagram \cref{PistA_nat} commutes on objects, suppose we are given a $(\sma \angordmal)$-system $(a,\gl)$ in $\A$ and a marker $\angs = \ang{s_j \subseteq \ufs{n}_j}_{j \in \ufs{q}}$.  Using \cref{smash_fpsiFG,A_fangpsi,psitil_f,apsitil_s,nsys_obj_unity,PistAnbe_aangs}, the following computation shows that the two composites in \cref{PistA_nat} yield the same $\angs$-component object.
\begin{equation}\label{PistA_nat_obj}
\begin{split}
&\big((\Aupom) \Pist_{\A,\angordmal} (a,\gl)\big)_{\angs}\\
&= \scalebox{.9}{$\begin{cases}
\zero & \text{if $\psiinv_j s_j = \emptyset$ for some $j \in \ufs{q}$,}\\
a^\stimes_{\ang{\psiinv_{f(i)} s_{f(i)}}_{i \in \ufs{p}}} 
= a_{\ang{\psiinv_{f(i)} s_{f(i)}}_{i \in \ufs{p}}^\stimes} & \text{if $\psiinv_j s_j \neq \emptyset$ for each $j \in \ufs{q}$,}
\end{cases}$}\\
&= a_{(\sma\upom)^{-1} \angstimes}\\
&= \big((\Asmaupom) (a,\gl)\big)_{\angstimes}\\
&= \big(\Pist_{\A,\angordnbe} (\Asmaupom) (a,\gl)\big)_{\angs}.
\end{split}
\end{equation}
\item[Gluing]  
Suppose we are given an object $x \in \Op(r)$ with $r \geq 0$, a marker $\angs = \ang{s_j \subseteq \ufs{n}_j}_{j \in \ufs{q}}$, an index $k \in \ufs{q}$, and a partition $s_k = \coprod_{\ell \in \ufs{r}} s_{k,\ell} \subseteq \ufs{n}_k$.  If either
\begin{itemize}
\item $\psiinv_j s_j = \emptyset$ for some $j \in \ufs{q}$ or
\item $\finv(k) = \emptyset$, 
\end{itemize}
then each composite in \cref{PistA_nat} yields an identity gluing morphism at $(x; \angs, k, \ang{s_{k,\ell}}_{\ell \in \ufs{r}})$ by \cref{glutil_component_ii,glutil_component_iii,nsys_unity_empty,nsys_unity_one}.  Thus, we assume that $\psiinv_j s_j \neq \emptyset$ for each $j \in \ufs{q}$ and $\finv(k) \neq \emptyset$.  

Using \cref{smash_fpsiFG,A_fangpsi,glpsitil,gltimes_component}, the following computation shows that the two composites in \cref{PistA_nat} yield the same gluing morphism at $(x; \angs, k, \ang{s_{k,\ell}}_{\ell \in \ufs{r}})$.
\[\begin{split}
&\big((\Aupom) \Pist_{\A,\angordmal} (a,\gl)\big)_{x;\, \angs,\, k, \ang{s_{k,\ell}}_{\ell \in \ufs{r}}}\\
&= \big(\Pist_{\A,\angordmal} (a,\gl)\big)_{x;\, \ang{\psiinv_{f(i)} s_{f(i)}}_{i \in \ufs{p}},\, \finv(k),\, \ang{\psiinv_k s_{k,\ell}}_{\ell \in \ufs{r}}}\\
&= \gl_{x;\, \ang{\psiinv_{f(i)} s_{f(i)}}_{i \in \ufs{p}}^\stimes,\, \ang{( \ang{\psiinv_{f(i)} s_{f(i)}}_{i \in \ufs{p}} \,\comp_{\finv(k)}\, (\psiinv_k s_{k,\ell}))^\stimes}_{\ell \in \ufs{r}}}\\
&= \gl_{x;\, (\sma\upom)^{-1} \angstimes,\, \ang{(\sma\upom)^{-1} (\angs \,\compk\, s_{k,\ell})^\stimes}_{\ell \in \ufs{r}}}\\
&= \big((\Asmaupom) (a,\gl)\big)_{x;\, \angstimes,\, \ang{(\angs \,\compk\, s_{k,\ell})^\stimes}_{\ell \in \ufs{r}}}\\
&= \big(\Pist_{\A,\angordnbe} (\Asmaupom) (a,\gl)\big)_{x;\, \angs,\, k, \ang{s_{k,\ell}}_{\ell \in \ufs{r}}}
\end{split}\]
This proves that the diagram \cref{PistA_nat} commutes on $(\sma\angordmal)$-systems in $\A$.
\item[Morphisms]  
To prove that the diagram \cref{PistA_nat} commutes on morphisms of $(\sma\angordmal)$-systems, we reuse the computation \cref{PistA_nat_obj}, along with  \cref{thapsitil_comp,thetatimes_angs,nsys_mor_unity} for morphisms.  
\end{description}
This proves that $\Pist_\A$ is a natural transformation.
\end{proof}

\begin{lemma}\label{PistA_gnat}
For each $\Op$-pseudoalgebra $\A$, the natural transformation in \cref{PistA_natural}
\begin{equation}\label{PistA_Geq_diag}
\begin{tikzpicture}[vcenter]
\def\s{25}
\draw[0cell]
(0,0) node (a1) {\GG}
(a1)++(2.2,0) node (a2) {\phantom{\GG}}
(a2)++(.15,0) node (a2') {\Catgst}
;
\draw[1cell=.9]
(a1) edge[bend left=\s] node {\smast \Sgo\A} (a2)
(a1) edge[bend right=\s] node[swap] {\Hgo\A} (a2) 
;
\draw[2cell]
node[between=a1 and a2 at .35, rotate=-90, 2label={above,\Pist_\A}] {\Rightarrow}
;
\end{tikzpicture}
\end{equation}
is $G$-equivariant.
\end{lemma}

\begin{proof}
By \cref{ggcatg_icell_geq}, the $G$-equivariance of $\Pist_\A$ means that the following diagram of pointed functors commutes for each object $\angordnbe \in \GG$ and $g \in G$.
\begin{equation}\label{PistA_Geq}
\begin{tikzpicture}[vcenter]
\def\v{-1.5}
\draw[0cell]
(0,0) node (a11) {\Asmaangordnbe}
(a11)++(3,0) node (a12) {\Aangordnbe}
(a11)++(0,\v) node (a21) {\Asmaangordnbe}
(a12)++(0,\v) node (a22) {\Aangordnbe}
;
\draw[1cell=.9]
(a11) edge node {\Pist_{\A,\angordnbe}} (a12)
(a12) edge node {g} (a22)
(a11) edge node[swap] {g} (a21)
(a21) edge node {\Pist_{\A,\angordnbe}} (a22)
;
\end{tikzpicture}
\end{equation}

\parhead{Bases cases}.  First, we consider the following two base cases.
\begin{enumerate}
\item If $\angordnbe = \vstar$, the basepoint in $\GG$, then $\sys{\A}{\vstar} = \bone$, the terminal category.  Thus, the diagram \cref{PistA_Geq} commutes. 
\item Suppose $\angordnbe = \ang{}$, the empty tuple in $\GG$.  Then $\sma\ang{} = \ord{1}$, and the isomorphism \cref{PistA_empty} 
\[\Asmaang = \Aordone \fto[\iso]{\Pist_{\A,\ang{}}} \Aang = \A\]
sends a $\ord{1}$-system $(a,\gl)$ to its $\{1\}$-component object $a_{\{1\}} \in \A$.  Using \cref{ga_scomp} and the trivial $G$-action on $\ord{1}$, the object equalities
\[(ga)_{\{1\}} = ga_{\ginv\{1\}} = ga_{\{1\}}\]
show that the diagram \cref{PistA_Geq} commutes on objects.  The same computation holds for morphisms of $\{1\}$-systems, using \cref{gtha_s} in place of \cref{ga_scomp}.  Thus, the diagram \cref{PistA_Geq} commutes. 
\end{enumerate}
In the rest of this proof, we assume that $\angordnbe \in \GG\setminus \{\vstar,\ang{}\}$ has length $q > 0$.
\begin{description}
\item[Component objects]  
Using \cref{al-sma-be,ga_scomponentGG,ga_scomp,PistAnbe_aangs}, the computation \cref{PistA_Geq_obj} shows that the two composites in \cref{PistA_Geq} yield the same $\angs$-component object for each $(\sma\angordnbe)$-system $(a,\gl)$ in $\A$ and each marker $\angs = \ang{s_j \subseteq \ufs{n}_j}_{j \in \ufs{q}}$.
\begin{equation}\label{PistA_Geq_obj}
\begin{split}
&\big(g \Pist_{\A,\angordnbe} (a,\gl)\big)_{\angs}\\
&= g\big(\Pist_{\A,\angordnbe} (a,\gl)\big)_{\ginv\angs}\\
&= ga_{(\ginv\angs)^\stimes}\\
&= ga_{\ginv \angstimes}\\
&= (ga)_{\angstimes}\\
&= \big(\Pist_{\A,\angordnbe}\, g (a,\gl)\big)_{\angs}
\end{split}
\end{equation}
\item[Gluing]  
Using \cref{al-sma-be,ga_gluingGG,ga_gl,gltimes_component}, the following computation shows that the two composites in \cref{PistA_Geq} yield the same gluing morphism for each object $x \in \Op(r)$ with $r \geq 0$, marker $\angs = \ang{s_j \subseteq \ufs{n}_j}_{j \in \ufs{q}}$, index $k \in \ufs{q}$, and partition $s_k = \coprod_{\ell \in \ufs{r}} s_{k,\ell} \subseteq \ufs{n}_k$.
\[\begin{split}
&\big(g \Pist_{\A,\angordnbe} (a,\gl)\big)_{x;\, \angs,\, k, \ang{s_{k,\ell}}_{\ell \in \ufs{r}}}\\
&= g\big( \Pist_{\A,\angordnbe} (a,\gl) \big)_{\ginv x;\, \ginv\angs,\, k, \ang{\ginv s_{k,\ell}}_{\ell \in \ufs{r}} }\\
&= g \gl_{\ginv x;\, (\ginv\angs)^\stimes,\, \ang{((\ginv\angs) \,\compk\, (\ginv s_{k,\ell}))^\stimes}_{\ell \in \ufs{r}}}\\
&= g \gl_{\ginv x;\, \ginv\angs^\stimes,\, \ang{\ginv(\angs \,\compk\, s_{k,\ell})^\stimes}_{\ell \in \ufs{r}}}\\
&= (g\gl)_{x;\, \angstimes,\, \ang{(\angs \,\compk\, s_{k,\ell})^\stimes}_{\ell \in \ufs{r}}}\\
&= \big(\Pist_{\A,\angordnbe} \, g (a,\gl)\big)_{x;\, \angs,\, k, \ang{s_{k,\ell}}_{\ell \in \ufs{r}}}
\end{split}\]
This proves that the diagram \cref{PistA_Geq} commutes on $(\sma\angordnbe)$-systems in $\A$.
\item[Morphisms]  
To prove that the diagram \cref{PistA_Geq} commutes on morphisms of $(\sma\angordnbe)$-systems in $\A$, we reuse the computation \cref{PistA_Geq_obj}, along with \cref{gtheta_angsGG,gtha_s,thetatimes_angs} for morphisms.  
\end{description}
This proves that $\Pist_\A$ is $G$-equivariant.
\end{proof}

\begin{lemma}\label{Pist_iinat}
The assignment in \cref{Pistar_twonat}
\begin{equation}\label{Pist_iinat_diagram}
\begin{tikzpicture}[vcenter]
\def\h{2.3} \def\v{.7} \def\s{15}
\draw[0cell]
(0,0) node (a1) {\AlglaxO}
(a1)++(0,.2*\v) node (a1') {\phantom{\AlglaxO}}
(a1)++(\h,\v) node (a2) {\FGCatg}
(a1)++(\h,-\v) node (a3) {\GGCatg}
;
\draw[1cell=.9]
(a1) edge[bend left=\s] node {\Sgo} (a2)
(a2) edge[bend left=\s, shorten <=-.5ex] node {\smast} (a3)
(a1') edge[bend right=\s] node[swap] {\Hgo} (a3)
;
\draw[2cell]
node[between=a2 and a3 at .55, shift={(-.3*\h,0)}, rotate=-120, 2labelw={below,\Pist,1pt}] {\Rightarrow}
;
\end{tikzpicture}
\end{equation}
is a 2-natural transformation.
\end{lemma}

\begin{proof}
The component 1-cells of $\Pist$ are the $G$-natural transformations $\Pist_\A$ in \cref{PistA_gnat} for $\Op$-pseudoalgebras $\A$.  We verify the 1-cell and 2-cell naturality of $\Pist$ \pcref{def:twonaturaltr}. 
\begin{description}
\item[1-cell naturality]  
The naturality of $\Pist$ with respect to 1-cells in $\AlglaxO$ means that, for each lax $\Op$-morphism between $\Op$-pseudoalgebras \pcref{def:laxmorphism}
\[(\A,\gaA,\phiA) \fto{(f,\actf)} (\B,\gaB,\phiB)\]
and each object $\angordnbe \in \GG$, the diagram \cref{Pist_1nat} of pointed functors commutes, where $(\Hgo f)_{\angordnbe}$, $(\Sgo f)_{\sma \angordnbe}$, and $\Pist_{\A,\angordnbe}$ are defined in, respectively, \cref{hgo_f_nbe,sgosf_nbe,PistAnbe}.
\begin{equation}\label{Pist_1nat}
\begin{tikzpicture}[vcenter]
\def\v{-1.5}
\draw[0cell]
(0,0) node (a11) {\Asmaangordnbe}
(a11)++(3,0) node (a12) {\Aangordnbe}
(a11)++(0,\v) node (a21) {\Bsmaangordnbe}
(a12)++(0,\v) node (a22) {\Bangordnbe}
;
\draw[1cell=.9]
(a11) edge node {\Pist_{\A,\angordnbe}} (a12)
(a12) edge node {(\Hgo f)_{\angordnbe}} (a22)
(a11) edge node[swap] {(\Sgo f)_{\sma \angordnbe}} (a21)
(a21) edge node {\Pist_{\B,\angordnbe}} (a22)
;
\end{tikzpicture}
\end{equation}
\begin{description}
\item[Base cases]  
First, we consider the following two base cases.
\begin{enumerate}
\item If $\angordnbe = \vstar$, the basepoint in $\GG$, then $\sys{\B}{\vstar} = \bone$, the terminal category.  Thus, the diagram \cref{Pist_1nat} commutes.
\item Suppose $\angordnbe = \ang{}$, the empty tuple in $\GG$.  Then $\sys{\B}{\ang{}} = \B$ and $\sma\ang{} = \ord{1}$.  Each of $\Pist_{\A,\ang{}}$ and $\Pist_{\B,\ang{}}$ \cref{PistA_empty} sends a $\ord{1}$-system to its $\{1\}$-component object.  Thus, each of the two composites in \cref{Pist_1nat} sends a $\ord{1}$-system $(a,\gl)$ in $\A$ to $fa_{\{1\}} \in \B$.  Similarly, each of those two composites sends a morphism $\tha$ of $\ord{1}$-systems in $\A$ to the morphism $f\tha_{\{1\}}$ in $\B$.
\end{enumerate}
In the rest of this proof, we assume that $\angordnbe \in \GG\setminus \{\vstar,\ang{}\}$ has length $q > 0$.
\item[Component objects]  
Using \cref{hgof_aglu_comp,sgof_nbe_a,PistAnbe_aangs}, the computation \cref{Pist_1nat_obj} shows that the two composites in \cref{Pist_1nat} yield the same $\angs$-component object in $\B$ for each $(\sma\angordnbe)$-system $(a,\gl)$ in $\A$ and marker $\angs = \ang{s_j \subseteq \ufs{n}_j}_{j \in \ufs{q}}$.
\begin{equation}\label{Pist_1nat_obj}
\begin{split}
&\big( (\Hgo f)_{\angordnbe} \Pist_{\A,\angordnbe} (a,\gl) \big)_{\angs}\\
&= f \big( \Pist_{\A,\angordnbe} (a,\gl) \big)_{\angs}\\
&= f a_{\angstimes}\\
&= \big((\Sgo f)_{\sma \angordnbe} (a,\gl) \big)_{\angstimes}\\
&= \big(\Pist_{\B,\angordnbe} (\Sgo f)_{\sma \angordnbe} (a,\gl) \big)_{\angs}
\end{split}
\end{equation}
\item[Gluing]  
Using \cref{hgof_gluing,sgof_nbe_gl,gltimes_component}, the following computation shows that the two composites in \cref{Pist_1nat} yield the same gluing morphism for each object $x \in \Op(r)$ with $r \geq 0$, marker $\angs = \ang{s_j \subseteq \ufs{n}_j}_{j \in \ufs{q}}$, index $k \in \ufs{q}$, and partition $s_k = \coprod_{\ell \in \ufs{r}} s_{k,\ell} \subseteq \ufs{n}_k$.
\[\begin{split}
&\big( (\Hgo f)_{\angordnbe} \Pist_{\A,\angordnbe} (a,\gl) \big)_{x;\, \angs,\, k, \ang{s_{k,\ell}}_{\ell \in \ufs{r}}}\\
&= \Big(f \big(\Pist_{\A,\angordnbe} (a,\gl)\big)_{x;\, \angs,\, k, \ang{s_{k,\ell}}_{\ell \in \ufs{r}}}\Big) \circ \actf_r\\
&= \big(f \gl_{x;\, \angstimes, \ang{(\angs \,\compk\, s_{k,\ell})^\stimes}_{\ell \in \ufs{r}}}\big) \circ \actf_r\\
&= \big((\Sgo f)_{\sma \angordnbe} (a,\gl) \big)_{x;\, \angstimes, \ang{(\angs \,\compk\, s_{k,\ell} )^\stimes}_{\ell \in \ufs{r}}}\\
&= \big(\Pist_{\B,\angordnbe} (\Sgo f)_{\sma \angordnbe} (a,\gl) \big)_{x;\, \angs,\, k, \ang{s_{k,\ell}}_{\ell \in \ufs{r}}}
\end{split}\]
This proves that the diagram \cref{Pist_1nat} commutes on $(\sma\angordnbe)$-systems in $\A$.
\item[Morphisms]  
To prove that the diagram \cref{Pist_1nat} commutes on morphisms of $(\sma\angordnbe)$-systems in $\A$, we reuse the computation \cref{Pist_1nat_obj}, along with \cref{hgof_tha_angs,sgof_tha_s,thetatimes_angs} for morphisms.  This proves that the diagram \cref{Pist_1nat} commutes, proving the 1-cell naturality of $\Pist$.
\end{description}
\item[2-cell naturality]  
The naturality of $\Pist$ with respect to 2-cells in $\AlglaxO$ means that, for each $\Op$-transformation $\omega$ \pcref{def:algtwocells} between lax $\Op$-morphisms between $\Op$-pseudoalgebras 
\begin{equation}\label{ome_Otr}
\begin{tikzpicture}[vcenter]
\def\s{25}
\draw[0cell]
(0,0) node (a) {\phantom{A}}
(a)++(2,0) node (b) {\phantom{A}}
(a)++(-.7,0) node (a') {(\A,\gaA,\phiA)}
(b)++(.7,0) node (b') {(\B,\gaB,\phiB)}
;
\draw[1cell=.9]
(a) edge[bend left=\s] node {(f,\actf)} (b)
(a) edge[bend right=\s] node[swap] {(h,\acth)} (b)
;
\draw[2cell]
node[between=a and b at .42, rotate=-90, 2label={above,\omega}] {\Rightarrow}
;
\end{tikzpicture}
\end{equation}
and each object $\angordnbe \in \GG$, the two whiskered natural transformations in \cref{Pist_2nat} are equal.
\begin{equation}\label{Pist_2nat}
\begin{tikzpicture}[vcenter]
\def\v{-1.5} \def\p{10} \def\s{25} \def\t{1.2em} \def\u{1em}
\draw[0cell=.85]
(0,0) node (a11) {\Asmaangordnbe}
(a11)++(3.5,0) node (a12) {\Aangordnbe}
(a11)++(0,\v) node (a21) {\Bsmaangordnbe}
(a12)++(0,\v) node (a22) {\Bangordnbe}
;
\draw[1cell=.8]
(a11) edge[bend left=\p] node {\Pist_{\A,\angordnbe}} (a12)
(a21) edge[bend right=\p] node[swap] {\Pist_{\B,\angordnbe}} (a22)
(a11) edge[bend right=\s, transform canvas={xshift=-\t}] node[swap,pos=.4] {(\Sgo f)_{\sma \angordnbe}} (a21)
(a11) edge[bend left=\s, transform canvas={xshift=\t}] node[pos=.6] {(\Sgo h)_{\sma \angordnbe}} (a21)
(a12) edge[bend right=\s, transform canvas={xshift=-\u}] node[swap,pos=.4] {(\Hgo f)_{\angordnbe}} (a22)
(a12) edge[bend left=\s, transform canvas={xshift=\u}] node[pos=.6] {(\Hgo h)_{\angordnbe}} (a22)
;
\draw[2cell]
node[between=a11 and a21 at .65, 2label={above,\Sgo\omega}] {\Rightarrow}
node[between=a12 and a22 at .65, 2label={above,\Hgo\omega}] {\Rightarrow}
;
\end{tikzpicture}
\end{equation}
In \cref{Pist_2nat}, the natural transformations 
\[\Hgo\omega = (\Hgo\omega)_{\angordnbe} \andspace 
\Sgo\omega = (\Sgo\omega)_{\sma\angordnbe}\]
are defined in, respectively, \cref{hgo_om_nbe,sgo_omega_nbe}.    There are two base cases as follows.
\begin{enumerate}
\item If $\angordnbe = \vstar$, then $\sys{\B}{\vstar} = \bone$.  Thus, the two whiskered natural transformations in \cref{Pist_2nat} are equal.
\item Suppose $\angordnbe = \ang{}$, so $\sys{\B}{\ang{}} = \B$ and $\sma\ang{} = \ord{1}$.  Each of the two whiskered natural transformations in \cref{Pist_2nat} sends a $\ord{1}$-system $(a,\gl)$ in $\A$ to the morphism 
\[fa_{\{1\}} \fto{\omega_{a_{\{1\}}}} ha_{\{1\}} \inspace \B.\]
\end{enumerate}
If $\angordnbe \in \GG \setminus \{\vstar,\ang{}\}$, then we assume that $(a,\gl)$ is a $(\sma\angordnbe)$-system in $\A$ and $\angs = \ang{s_j \subseteq \ufs{n}_j}_{j \in \ufs{q}}$ is a marker.  To show that the two whiskered natural transformations in \cref{Pist_2nat} yield the same $\angs$-component morphism in $\B$, we reuse the computation \cref{Pist_1nat_obj} with $\omega$ in place of $f$, along with \cref{omega_aangs,sgo_omega_nbe_s,thetatimes_angs}.\qedhere
\end{description}
\end{proof}

\section{$H$-Theory Comparison and Twisted Products}
\label{sec:hcomp_prod}

This section shows that, for each $\Tinf$-operad $\Op$ \pcref{as:OpA}, $\Op$-pseudoalgebra $\A$ \pcref{def:pseudoalgebra}, and object $\angordnbe \in \GG \setminus \{\vstar,\ang{}\}$, the (strong) $H$-theory comparison pointed $G$-functors $\Pist_{\A,\angordnbe}$ and $\Pistsg_{\A,\angordnbe}$ admit the following zigzag factorizations.
\begin{equation}\label{Pistfact_seci}
\begin{tikzpicture}[vcenter]
\def\h{3}
\draw[0cell]
(0,0) node (a11) {\Asmaangordnbe}
(a11)++(\h,0) node (a12) {\Aangordnbe}
(a11)++(\h/2,-1.1) node (a2) {\phantom{\proAnbe}}
(a2)++(0,.15) node (a2') {\proAnbe}
;
\draw[1cell=.9]
(a11) edge node {\Pist_{\A,\angordnbe}} (a12)
(a11) edge node[swap] {\zb} (a2)
(a12) edge node {\zd} (a2)
;
\begin{scope}[shift={(1.75*\h,0)}]
\draw[0cell]
(0,0) node (a11) {\Asgsmaangordnbe}
(a11)++(\h,0) node (a12) {\Asgangordnbe}
(a11)++(\h/2,-1.1) node (a2) {\phantom{\proAnbe}}
(a2)++(0,.15) node (a2') {\proAnbe}
;
\draw[1cell=.9]
(a11) edge node {\Pistsg_{\A,\angordnbe}} (a12)
(a11) edge node[swap] {\zbsg} (a2)
(a12) edge node {\zdsg} (a2)
;
\end{scope}
\end{tikzpicture}
\end{equation}
These factorizations are used in subsequent chapters to establish further properties of the (strong) $H$-theory comparisons $\Pist$ and $\Pistsg$.

\secoutline
\begin{itemize}
\item \cref{def:proCnbe} defines the twisted product $\proCnbe$ for each pointed $G$-category $\C$ and object $\angordnbe \in \GG \setminus \{\vstar,\ang{}\}$.
\item \cref{def:Htoprod} constructs the pointed $G$-functors $\zb$, $\zbsg$, $\zd$, and $\zdsg$.
\item \cref{Pist_to_prod} proves that each of $\zb$, $\zbsg$, $\zd$, and $\zdsg$ is actually a pointed $G$-functor and establishes the desired zigzag factorizations of $\Pist_{\A,\angordnbe}$ and $\Pistsg_{\A,\angordnbe}$.
\end{itemize}

\subsection*{Twisted Products}
Recall from \cref{PistA_vstar,PistA_empty} that $\Pist_{\A,\vstar} = 1_\bone$, and $\Pist_{\A,\ang{}}$ is a pointed $G$-isomorphism.  Thus, \cref{def:proCnbe} only deals with objects in $\GG \setminus \{\vstar,\ang{}\}$.  Recall that $\ufsn = \{1,2,\ldots,n\}$ denotes an unpointed finite set with $n$ elements \cref{ufsn}.
 
\begin{definition}[Twisted Products]\label{def:proCnbe}
For a group $G$, suppose $\C$ is a pointed $G$-category with basepoint $\zero$ \pcref{def:GCat,def:ptGcat}, and  
\[\angordnbe = \sordi{n}{\be}{j}_{j \in \ufs{q}} \in \GG \setminus \{\vstar,\ang{}\}\]
is an object of length $q > 0$ \cref{GG_objects}.  The \index{twisted product}\emph{$\angordnbe$-twisted product} is the pointed $G$-category \label{not:proCnbe}$\proCnbe$ defined as follows.
\begin{description}
\item[Underlying category] The category $\proCnbe$ is the Cartesian product
\[\proCnbe = \C^{n_1 n_2 \cdots n_q}.\]
An object or a morphism $c \in \proCnbe$ is denoted by its components as
\begin{equation}\label{proCnbe_object}
c = \ang{c_{\bdi} \in \C}_{\bdi \sins \ufs{n_1n_2\cdots n_q}}
\end{equation}
with the indices ordered lexicographically \cref{lex_bijection}:
\begin{equation}\label{bdi}
\bdi = \ang{i_j}_{j \in \ufs{q}} \in \txprod_{j \in \ufs{q}}\, \ufs{n}_j = \ufs{\txprod_{j \in \ufs{q}}\, n_j}.
\end{equation}
The object or morphism $c_{\bdi} \in \C$ is called the \emph{$\bdi$-th coordinate of $c$}.  For an element $g \in G$, the diagonal $g$-action \cref{al-sma-be} on a $q$-tuple $\bdi = \ang{i_j}_{j \in \ufs{q}}$ is denoted by
\begin{equation}\label{gbdi}
g\bdi = \ang{g i_j}_{j \in \ufs{q}} = \ang{(\be_j g)(i_j) \in \ufs{n}_j}_{j \in \ufsq}.
\end{equation} 
\item[$G$-action] For $g \in G$, the $g$-action on $c = \ang{c_{\bdi}} \in \proCnbe$ \cref{proCnbe_object} is defined as
\begin{equation}\label{proCnbe_gaction}
gc = \bang{(gc)_{\bdi} = g c_{\ginv\bdi}}_{\bdi \sins \ufs{n_1n_2\cdots n_q}} \in \proCnbe
\end{equation}
with the index given by 
\begin{equation}\label{ginvbdi}
\ginv\bdi = \ang{\ginv i_j}_{j \in \ufs{q}} \in \ufs{\txprod_{j \in \ufs{q}}\, n_j}.
\end{equation}
The object or morphism $c_{\ginv\bdi}$ is the $\ginv\bdi$-th coordinate of $c$, and $g c_{\ginv\bdi}$ is its image under the $g$-action on $\C$.
\item[Basepoint] The $G$-fixed basepoint $\zero \in \proCnbe$ has each coordinate given by the $G$-fixed basepoint $\zero \in \C$.
\end{description}
The pointed $G$-category axioms for $\proCnbe$ follow from \cref{proCnbe_gaction}, \cref{ginvbdi}, the pointed $G$-category axioms for $\C$, and the axioms for the pointed finite $G$-sets $\ordi{n}{\be}{j}$ for $j \in \ufsq$.
\end{definition}

\begin{explanation}\label{expl:proCnbe}
Unless each $\be_j$ is the trivial $G$-action on $\ord{n}_j$, the pointed $G$-category $\proCnbe$ is \emph{not} generally the product $G$-category of $\txprod_{j\in \ufs{q}}\, n_j$ copies of $\C$.  By \cref{proCnbe_gaction}, the $G$-action on $\proCnbe$ is twisted by the $G$-action on each pointed finite $G$-set $\ordi{n}{\be}{j}$.
\end{explanation}

\subsection*{Comparing $H$-Theory and Twisted Products}

Recall the pointed $G$-categories of (strong) $\ordnbe$-systems 
\[\Aordnbe \andspace \Asgordnbe\]
for $\ordnbe \in \FG$ \pcref{def:nsys,def:nsys_morphism,def:nsys_gcat} and the pointed $G$-categories of (strong) $\angordnbe$-systems 
\[\Aangordnbe \andspace \Asgangordnbe\] 
for $\angordnbe \in \GG$ \pcref{def:nsystem,def:nsystem_morphism,def:nbeta_gcat}.  \cref{def:Htoprod} compares these pointed $G$-categories of (strong) systems to twisted products \pcref{def:proCnbe}.

\begin{definition}\label{def:Htoprod}
For a $\Tinf$-operad $\Op$ \pcref{as:OpA}, an $\Op$-pseudoalgebra $\A$ \pcref{def:pseudoalgebra}, and an object $\angordnbe \in \GG \setminus \{\vstar,\ang{}\}$ of length $q > 0$ \cref{GG_objects}, we define pointed $G$-functors $\zb = \zb_{\A,\angordnbe}$ \cref{sgotoprod}, $\zd = \zd_{\A,\angordnbe}$ \cref{hgotoprod}, and their strong variants \cref{sgosgtoprod,hgosgtoprod} as follows.
\begin{description}
\item[Shimakawa $H$-theory to twisted products] 
The functor
\begin{equation}\label{sgotoprod}
\Asmaangordnbe \fto{\zb} \proAnbe
\end{equation}
sends a $(\sma\angordnbe)$-system $(a,\gl)$ \cref{nsys} and a morphism $\theta$ \cref{nsys_mor} in $\Asmaangordnbe$ to, respectively,
\begin{equation}\label{sgotoprod_def}
\begin{split}
\zb(a,\gl) &= \bang{\zb(a,\gl)_{\bdi} = a_{\{\bdi\}}}_{\bdi \sins \ufs{n_1 n_2 \cdots n_q}} \andspace\\
\zb\theta &= \bang{(\zb\theta)_{\bdi} = \theta_{\{\bdi\}}}_{\bdi \sins \ufs{n_1 n_2 \cdots n_q}}.
\end{split}
\end{equation}
Here\label{not:bdiangij} 
\[\{\bdi\} = \{\ang{i_j}_{j \in \ufs{q}}\} \subseteq \txprod_{j \in \ufs{q}}\, \ufs{n}_j = \ufs{\txprod_{j \in \ufs{q}}\, n_j}\] 
is the one-element subset consisting of only the $q$-tuple $\bdi = \ang{i_j}_{j \in \ufsq}$.
\item[Shimakawa strong $H$-theory to twisted products]   
The functor
\begin{equation}\label{sgosgtoprod}
\Asgsmaangordnbe \fto{\zbsg} \proAnbe
\end{equation}
is the restriction of $\zb$ to the full subcategory of strong $(\sma\angordnbe)$-systems in $\A$ \cref{Asgordnbe}.
\item[$H$-theory to twisted products]
The functor
\begin{equation}\label{hgotoprod}
\Aangordnbe \fto{\zd} \proAnbe
\end{equation}
sends an $\angordnbe$-system $(a,\glu)$ \cref{nsystem} and a morphism $\theta$ \cref{nsystem_mor} in $\Aangordnbe$ to, respectively,
\begin{equation}\label{hgotoprod_def}
\begin{split}
\zd(a,\glu) &= \bang{a_{\ang{\{i_j\}}_{j \in \ufsq}}}_{\bdi \sins \ufs{n_1 n_2 \cdots n_q}} \andspace\\
\zd\theta &= \bang{\theta_{\ang{\{i_j\}}_{j \in \ufsq}}}_{\bdi \sins \ufs{n_1 n_2 \cdots n_q}}.
\end{split}
\end{equation}
Here\label{not:angij}
\[\ang{\{i_j\}}_{j \in \ufsq} = \big(\{i_1\}, \{i_2\}, \ldots, \{i_q\}\big)\] 
is the marker \cref{marker} whose $j$-th entry is the one-element subset $\{i_j\} \subseteq \ufs{n}_j$. 
\item[Strong $H$-theory to twisted products]  
The functor
\begin{equation}\label{hgosgtoprod}
\Asgangordnbe \fto{\zdsg} \proAnbe
\end{equation}
is the restriction of $\zd$ to the full subcategory of strong $\angordnbe$-systems in $\A$ \cref{Asgangordnbe}.
\end{description}
Note that if $\nbe$ has length $q=1$, meaning $\nbeta \in \FG$ is a pointed finite $G$-set, then $\zb$ and $\zbsg$ coincide with, respectively, $\zd$ and $\zdsg$.
\end{definition}

\cref{Pist_to_prod} \eqref{Pist_to_prod_i} proves that each of $\zb$, $\zbsg$, $\zd$, and $\zdsg$ is a pointed $G$-functor.  For \cref{Pist_to_prod} \eqref{Pist_to_prod_ii}, recall the $H$-theory comparison pointed $G$-functor $\Pist_{\A,\angordnbe}$ \cref{PistAnbe} and its strong variant $\Pistsg_{\A,\angordnbe}$ \cref{PistsgAnbe}.

\begin{lemma}\label{Pist_to_prod}
In the setting of \cref{def:Htoprod}, the following statements hold.
\begin{enumerate}
\item\label{Pist_to_prod_i}
Each of $\zb$, $\zbsg$, $\zd$, and $\zdsg$ is a pointed $G$-functor.
\item\label{Pist_to_prod_ii}
The diagrams
\begin{equation}\label{Pist_prod_diag}
\begin{tikzpicture}[vcenter]
\def\h{3.5}
\draw[0cell]
(0,0) node (a11) {\Asmaangordnbe}
(a11)++(\h,0) node (a12) {\Aangordnbe}
(a11)++(\h/2,-1.1) node (a2) {\phantom{\proAnbe}}
(a2)++(0,.15) node (a2') {\proAnbe}
;
\draw[1cell=.9]
(a11) edge node {\Pist_{\A,\angordnbe}} (a12)
(a11) edge node[swap] {\zb} (a2)
(a12) edge node {\zd} (a2)
;
\end{tikzpicture}
\end{equation}
and 
\begin{equation}\label{Pistsg_prod_diag}
\begin{tikzpicture}[vcenter]
\def\h{3.5}
\draw[0cell]
(0,0) node (a11) {\Asgsmaangordnbe}
(a11)++(\h,0) node (a12) {\Asgangordnbe}
(a11)++(\h/2,-1.1) node (a2) {\phantom{\proAnbe}}
(a2)++(0,.15) node (a2') {\proAnbe}
;
\draw[1cell=.9]
(a11) edge node {\Pistsg_{\A,\angordnbe}} (a12)
(a11) edge node[swap] {\zbsg} (a2)
(a12) edge node {\zdsg} (a2)
;
\end{tikzpicture}
\end{equation}
of pointed $G$-functors commute.
\end{enumerate}
\end{lemma}

\begin{proof}
\begin{description}
\item[\cref{Pist_to_prod_i}]  
We first prove that $\zb$, defined in \cref{sgotoprod_def}, is a pointed $G$-functor.  
\begin{itemize}
\item Functoriality of $\zb$ follows from the fact that identity morphisms and composition are defined componentwise in the domain, $\Asmaangordnbe$, and coordinatewise in the $\angordnbe$-twisted product $\proAnbe$.
\item The functor $\zb$ sends the base $(\sma\angordnbe)$-system $(\zero,1_\zero)$ to the object in $\proAnbe$ with each coordinate given by $\zero \in \A$, which is the basepoint of $\proAnbe$.
\item Using \cref{ga_scomp,proCnbe_gaction,sgotoprod_def}, the computation \cref{Pist_prod_zb} shows that,  for each element $g \in G$, $\zb$ commutes with the $g$-action on each $(\sma\angordnbe)$-system $(a,\gl) \in \Asmaangordnbe$.
\begin{equation}\label{Pist_prod_zb}
\begin{split}
&\zb\big(g(a,\gl)\big)\\
&= \ang{(ga)_{\{\bdi\}}}_{\bdi \sins \ufs{n_1 n_2 \cdots n_q}}\\
&= \ang{ga_{\ginv \{\bdi\}}}_{\bdi \sins \ufs{n_1 n_2 \cdots n_q}}\\
&= \ang{g a_{\{\ginv\bdi\}}}_{\bdi \sins \ufs{n_1 n_2 \cdots n_q}}\\
&= g\ang{a_{\{\bdi\}}}_{\bdi \sins \ufs{n_1 n_2 \cdots n_q}}\\
&= g\big(\zb(a,\gl)\big)
\end{split}
\end{equation}
With a morphism $\theta$ in place of an object $(a,\gl)$ in $\Asmaangordnbe$ and \cref{gtha_s} in place of \cref{ga_scomp}, the computation \cref{Pist_prod_zb} also proves that $\zb$ commutes with the $g$-action on morphisms. 
\end{itemize}
This proves that $\zb$ is a pointed $G$-functor.  The strong variant $\zbsg$ \cref{sgosgtoprod} is a pointed $G$-functor because it is the precomposition of $\zb$ with the full subcategory inclusion 
\[\Asgsmaangordnbe \to \Asmaangordnbe,\]
which is a pointed $G$-functor.  

To prove that $\zd$ \cref{hgotoprod} and $\zdsg$ \cref{hgosgtoprod} are pointed $G$-functors, we reuse the proof for $\zb$ and $\zbsg$, using \cref{ga_scomponentGG,gtheta_angsGG,hgotoprod_def} instead of \cref{ga_scomp,gtha_s,sgotoprod_def}.  In the computation \cref{Pist_prod_zb}, the one-element subsets
\[\begin{split}
\{\bdi\} &= \{\ang{i_j}_{j \in \ufsq}\} \andspace \\
\ginv\{\bdi\} &= \{\ginv\bdi\} \\
&= \{\ang{\ginv i_j}_{j \in \ufsq}\} \subseteq \ufs{n_1 n_2 \Cdots n_q}
\end{split}\]
are replaced by, respectively, the $q$-tuples of one-element subsets
\[\begin{split}
\ang{\{i_j\}}_{j \in \ufsq} & \andspace\\ 
\ginv\ang{\{i_j\}}_{j \in \ufsq} &= \ang{\{\ginv i_j\}}_{j \in \ufsq}\\
&= \big(\{\ginv i_1\}, \ldots, \{\ginv i_q\} \big).
\end{split}\]
This proves assertion \cref{Pist_to_prod_i}.
\item[\cref{Pist_to_prod_ii}]  
Once we show that the diagram \cref{Pist_prod_diag} commutes, the diagram \cref{Pistsg_prod_diag} also commutes by restricting to strong $(\sma\angordnbe)$-systems.  Using \cref{PistAnbe_aangs,sgotoprod_def,hgotoprod_def}, the computation \cref{dPistb} shows that the two composites in \cref{Pist_prod_diag} yield the same $\bdi$-th coordinate object for each $(\sma\angordnbe)$-system $(a,\gl)$ in $\A$ and $q$-tuple $\bdi = \ang{i_j \in \ufs{n}_j}_{j \in \ufs{q}}$.
\begin{equation}\label{dPistb}
\begin{split}
&\big(\zd \Pist_{\A,\angordnbe} (a,\gl)\big)_{\bdi}\\
&= \big(\Pist_{\A,\angordnbe} (a,\gl)\big)_{\ang{\{i_j\}}_{j \in \ufsq}}\\
&= a_{\{i_1\} \times \{i_2\} \times \cdots \times \{i_q\}}\\
&= a_{\{(i_1,\, i_2,\, \ldots,\, i_q)\}}\\
&= \zb(a,\gl)_{\bdi}
\end{split}
\end{equation}
This proves that the diagram \cref{Pist_prod_diag} commutes on objects.  Using \cref{thetatimes_angs} instead of \cref{PistAnbe_aangs}, the computation \cref{dPistb} also proves that the diagram \cref{Pist_prod_diag} commutes on morphisms.  This proves assertion \cref{Pist_to_prod_ii}.\qedhere
\end{description}
\end{proof}

\section{Fixed Points of Twisted Products}
\label{sec:twprod_fixed}

This section identifies the $G$-fixed subcategory of a twisted product $\proCnbe$ \pcref{twprod_fixed}.  This result is used in \cref{sec:gggcat_weq} to study special $\GGG$-categories and weak $G$-equivalences in $\GGCatg$; see \cref{f_nbe,ggg_weq_sp,ggg_sp_weq}.  \cref{as:twprod_fixed} states the context for \cref{twprod_fixed} and fixes notation.

\begin{assumption}\label{as:twprod_fixed}
For a group $G$ with group unit $e$, suppose
\[\angordnbe = \sordi{n}{\be}{j}_{j \in \ufs{q}} \in \GG \setminus \{\vstar,\ang{}\}\]
is an object of length $q > 0$ \cref{GG_objects}.  Using the lexicographic ordering \cref{lex_bijection} and the diagonal $G$-action \cref{gbdi}, suppose the ordered finite $G$-set 
\[\ufs{n_1 \Cdots n_q} = \txprod_{j \in \ufsq}\, \ufsn_j\]
has $r>0$ $G$-orbits.  There is a $G$-equivariant decomposition\label{not:bdkt}
\[\ufs{n_1 \Cdots n_q} = \txcoprod_{t \in \ufsr}\, G\bdk^t\]
such that $\bdk^t \in \ufs{n_1 \Cdots n_q}$ is the least element in its $G$-orbit $G\bdk^t$ with the ordering inherited from $\ufs{n_1 \Cdots n_q}$.  The set of $G$-orbits is ordered by their least elements $\{\bdk^t\}_{t \in \ufsr}$.  For each $t \in \ufsr$, denote by
\begin{itemize}
\item $G_t$\label{not:Gt} the stabilizer of $\bdk^t$ and
\item $N_t$ the cardinality of the $G$-orbit $G\bdk^t$. 
\end{itemize} 
The inherited ordering of $G\bdk^t$ equips the unpointed finite set $\ufsN_t = \{1,2,\ldots,N_t\}$ with a $G$-action, along with an ordered $G$-bijection
\begin{equation}\label{Gkt_Nt}
G\bdk^t \iso \ufsN_t
\end{equation}
that sends $\bdk^t$ to 1 and a $G$-bijection
\begin{equation}\label{nonenqNt}
\ufs{n_1 \Cdots n_q} \iso \txcoprod_{t \in \ufsr}\, \ufsN_t.
\end{equation}
Denoting an element $d \in \ufsN_t$ also by $(t,d)$, the $G$-bijection \cref{nonenqNt} sends 
\[\bdk^t \in \ufs{n_1 \Cdots n_q} \tospace (t,1) \in \txcoprod_{t \in \ufsr}\, \ufsN_t.\]  
For each $d \in \ufsN_t$, choose an element $g_{t,d} \in G$ such that
\begin{equation}\label{gtd}
g_{t,d}(1) = d \andspace g_{t,1} = e.
\end{equation}
The element $g_{t,d}$ exists because $G$ acts transitively on the $G$-orbit $G\bdk^t$ and hence also on $\ufsN_t$.

We consider a pointed $G$-category $\C$ with basepoint $\zero$ and the $\angordnbe$-twisted product $\proCnbe$ \pcref{def:proCnbe}.  The $G$-bijection \cref{nonenqNt} and the elements $g_{t,d} \in G$ \cref{gtd} yield pointed functors
\begin{equation}\label{uphuphinv}
\begin{tikzpicture}[vcenter]
\draw[0cell]
(0,0) node (a1) {\proCnbe}
(a1)++(2.2,0) node (a2) {\phantom{\txprod_{t \in \ufsr}\, \C}}
(a2)++(0,-.07) node (a2') {\txprod_{t \in \ufsr}\, \C}
;
\draw[1cell=.9]
(a1) edge[transform canvas={yshift=.4ex}] node {\uph} (a2)
(a2) edge[transform canvas={yshift=-.5ex}] node {\uphinv} (a1)
;
\end{tikzpicture}
\end{equation}
defined by
\begin{equation}\label{uph}
\begin{split}
\uph\ang{c_{\bdi}}_{\bdi \in \ufs{n_1\cdots n_q}} &= \ang{c_{\bdk^t}}_{t \in \ufsr} \andspace\\
\uphinv\ang{c_t}_{t \in \ufsr} &= \ang{g_{t,d} c_t}_{d \in \ufsN_t,\, t \in \ufsr}
\end{split}
\end{equation}
on both objects and morphisms.  Using the $G$-bijection \cref{nonenqNt}, $\uph$ is also given by
\begin{equation}\label{uph_ctd}
\uph\ang{c_{t,d}}_{d \in \ufsN_t,\, t \in \ufsr} = \ang{c_{t,1}}_{t \in \ufsr}.
\end{equation}
Denote by
\begin{itemize}
\item $\Cnbeg$ the $G$-fixed subcategory of $\proCnbe$ and
\item $\C^{G_t}$ the $G_t$-fixed subcategory of $\C$ for each $t \in \ufsr$.\defmark
\end{itemize}
\end{assumption}

\begin{lemma}\label{twprod_fixed}\index{twisted product!fixed point}
Under \cref{as:twprod_fixed}, $\uph$ restricts to a pointed isomorphism of categories
\begin{equation}\label{uph_iso}
\Cnbeg \fto[\iso]{\uph} \txprod_{t \in \ufsr}\, \C^{G_t}
\end{equation}
with inverse given by the restriction of $\uphinv$.  Both $\uph$ and $\uphinv$ are natural in pointed $G$-functors on $\C$.
\end{lemma}

\begin{proof}
We first verify that the functors
\begin{equation}\label{uph_res}
\begin{tikzpicture}[vcenter]
\draw[0cell]
(0,0) node (a1) {\Cnbeg}
(a1)++(2.7,0) node (a2) {\phantom{\txprod_{t \in \ufsr}\, \C^{G_t}}}
(a2)++(0,0) node (a2') {\txprod_{t \in \ufsr}\, \C^{G_t}}
;
\draw[1cell=.9]
(a1) edge[transform canvas={yshift=.5ex}] node {\uph} (a2)
(a2) edge[transform canvas={yshift=-.4ex}] node {\uphinv} (a1)
;
\end{tikzpicture}
\end{equation}
are well defined.  Then we verify that they are mutual inverses.  The naturality of $\uph$ and $\uphinv$ with respect to pointed $G$-functors on $\C$ follows from the fact that the resulting functors on $\Cnbeg$ and $\txprod_{t\in \ufsr}\, \C^{G_t}$ are defined entrywise.
\begin{description}
\item[$\uph$ is well defined]  
The functor $\uph$ is defined in \cref{uph}, restricted to the $G$-fixed subcategory of the $\nbe$-twisted product $\proCnbe$.  We abbreviate $\ang{\Cdots}_{\bdi \in \ufs{n_1 \cdots n_q}}$ to $\ang{\Cdots}_{\bdi}$.  Given an object or a morphism $\ang{c_{\bdi}}_{\bdi} \in \Cnbeg$, we verify that $\ang{c_{\bdk^t}}_{t\in \ufsr}$ belongs to $\txprod_{t \in \ufsr}\, \C^{G_t}$, meaning that $c_{\bdk^t}$ is $G_t$-fixed for each $t \in \ufsr$.  Given an element $h \in G_t \subseteq G$, the fact that $\ang{c_{\bdi}}_{\bdi}$ is fixed by $h$ means that
\[\ang{c_{\bdi}}_{\bdi} = h \cdot \ang{c_{\bdi}}_{\bdi} = \ang{hc_{\hinv \bdi}}_{\bdi}.\]
The $\bdk^t$-th coordinate of $\ang{c_{\bdi}}_{\bdi}$ is given by
\[c_{\bdk^t} = hc_{\hinv \bdk^t} = hc_{\bdk^t}.\]
The second equality uses the fact that $\hinv \in G_t$, the stabilizer of $\bdk^t$.  Thus, $c_{\bdk^t}$ is fixed by $G_t$, proving that $\uph$ in \cref{uph_res} is a well-defined pointed functor.
\item[$\uphinv$ is well defined]  
The functor $\uphinv$ is defined in \cref{uph}, restricted in the $t$-th entry to the $G_t$-fixed subcategory of $\C$ for each $t \in \ufsr$.  Given an object or a morphism $\ang{c_t}_{t \in \ufsr}$ in $\txprod_{t\in \ufsr}\, \C^{G_t}$, we verify that its image under $\uphinv$ belongs to $\Cnbeg$, meaning that it is $G$-fixed.  Given an element $h \in G$, the $h$-action on $\uphinv\ang{c_t}_{t \in \ufsr}$ is given by
\begin{equation}\label{h_uphinv}
\begin{split}
h \cdot \uphinv\ang{c_t}_{t \in \ufsr} 
&= h \cdot \ang{g_{t,d} c_t}_{d \in \ufsN_t,\, t \in \ufsr} \\
&= \ang{h g_{t,\hinv d} c_t}_{d \in \ufsN_t,\, t \in \ufsr}.
\end{split}
\end{equation}
To show that the last entry in \cref{h_uphinv} is equal to $\uphinv\ang{c_t}_{t \in \ufsr}$, note that \cref{gtd} implies that, for each element $d \in \ufsN_t$,
\[h g_{t,\hinv d}(1) = h (\hinv d) = d = g_{t,d}(1).\]
These equalities imply that $\ginv_{t,d} h g_{t,\hinv d}$ fixes the element $1 \in \ufsN_t$, which corresponds to $\bdk^t \in G\bdk^t$ under the $G$-bijection \cref{Gkt_Nt}.  Thus, there exists an element $w_{t,d} \in G_t$, the stabilizer of $\bdk^t$, such that
\[h g_{t,\hinv d} = g_{t,d} w_{t,d}.\]
Continuing \cref{h_uphinv}, we have the following equalities in $\proCnbe$.
\[\begin{split}
h \cdot \uphinv\ang{c_t}_{t \in \ufsr} 
&= \ang{h g_{t,\hinv d} c_t}_{d \in \ufsN_t,\, t \in \ufsr} \\
&= \ang{g_{t,d} w_{t,d} c_t}_{d \in \ufsN_t,\, t \in \ufsr} \\
&= \ang{g_{t,d} c_t}_{d \in \ufsN_t,\, t \in \ufsr} \\
&= \uphinv\ang{c_t}_{t \in \ufsr}
\end{split}\]
The third equality uses the fact that $c_t \in \C^{G_t}$, so it is fixed by $w_{t,d}$.  Thus, $\uphinv\ang{c_t}_{t \in \ufsr}$ is fixed by $G$, proving that $\uphinv$ in \cref{uph_res} is a well-defined pointed functor.
\item[Mutual inverses]
For an object or a morphism $\ang{c_t}_{t \in \ufsr} \in \txprod_{t\in \ufsr}\, \C$, the following equalities prove that the composite $\uph \uphinv$ is the identity functor.
\[\begin{aligned}
& \uph \uphinv\ang{c_t}_{t \in \ufsr} &&\\
&= \uph \ang{g_{t,d} c_t}_{d \in \ufsN_t,\, t \in \ufsr} && \text{by \cref{uph}}\\
&= \ang{g_{t,1} c_t}_{t \in \ufsr} && \text{by \cref{uph_ctd}}\\
&= \ang{e c_t}_{t \in \ufsr} && \text{by \cref{gtd}}\\
&= \ang{c_t}_{t \in \ufsr} && 
\end{aligned}\]

For the composite $\uphinv\uph$, we use the $G$-bijection \cref{nonenqNt} and note that, for each $u \in \ufsr$ and $b \in \ufsN_u$, the $g_{u,b}$-action on an object or a morphism 
\[\ang{c_{t,d}}_{d \in \ufsN_t,\, t \in \ufsr} \in \proCnbe\] 
yields
\begin{equation}\label{gub_ctd}
g_{u,b} \cdot \ang{c_{t,d}}_{d \in \ufsN_t,\, t \in \ufsr}
= \ang{g_{u,b} c_{t, \ginv_{u,b}(d)}}_{d \in \ufsN_t,\, t \in \ufsr}.
\end{equation}
If $\ang{c_{t,d}}_{d \in \ufsN_t,\, t \in \ufsr}$ is $G$-fixed, then, by \cref{gtd,gub_ctd}, its $(u,b)$-th coordinate is given by
\begin{equation}\label{cub_gc}
c_{u,b} = g_{u,b} c_{u, \ginv_{u,b}(b)} = g_{u,b} c_{u,1}.
\end{equation}
The following equalities prove that the composite $\uphinv\uph$ is the identity functor.
\[\begin{aligned}
& \uphinv\uph\ang{c_{t,d}}_{d \in \ufsN_t,\, t \in \ufsr} &&\\
&= \uphinv\ang{c_{t,1}}_{t \in \ufsr} && \text{by \cref{uph_ctd}}\\
&= \ang{g_{t,d} c_{t,1}}_{d \in \ufsN_t,\, t \in \ufsr} && \text{by \cref{uph}}\\
&= \ang{c_{t,d}}_{d \in \ufsN_t,\, t \in \ufsr} && \text{by \cref{cub_gc}}
\end{aligned}\]
This proves that $\uph$ and $\uphinv$ are mutual inverses.\qedhere
\end{description}
\end{proof}

\begin{remark}\label{rk:mmo2.15}
\cref{twprod_fixed} is analogous to \cite[Lemma 2.15]{mmo}, which concern $G$-spaces and the indexing $G$-category $\FG$ \pcref{def:FG}.
\end{remark}

%% file: chap/sgoprod.tex
This chapter proves that there is an adjoint equivalence \pcref{thm:zbsg_eq}
\begin{equation}\label{zbsg_chi}
\begin{tikzpicture}[vcenter]
\def\s{22}
\draw[0cell]
(0,0) node (a1) {\phantom{X}}
(a1)++(1.8,0) node (a2) {\phantom{X}}
(a1)++(-.15,0) node (a1') {\proAnbe}
(a2)++(.7,-.04) node (a2') {\Asgsmaangordnbe}
node[between=a1 and a2 at .5] {\sim}
;
\draw[1cell=.9]
(a1) edge[bend left=\s] node {\zbsgad} (a2)
(a2) edge[bend left=\s] node {\zbsg} (a1)
;
\end{tikzpicture}
\end{equation}
between the twisted product $\proAnbe$ \pcref{def:proCnbe} and Shimakawa strong $H$-theory $\Asgsmaangordnbe$ \pcref{sys_FGcat} at the pointed finite $G$-set $\sma\angordnbe\in \FG$ \cref{smash_GGobjects} for each object $\angordnbe \in \GG \setminus \{\vstar,\ang{}\}$ and $\Op$-pseudoalgebra $\A$ \pcref{def:pseudoalgebra}.  \cref{thm:zbsg_eq} is used in \cref{thm:PistAequivalence} to prove that the strong $H$-theory comparison $\Pistsg_\A$ is componentwise an equivalence of categories.

While the right adjoint $\zbsg$ \cref{sgosgtoprod} is $G$-equivariant, its adjoint inverse $\zbsgad$ is only \emph{pseudo} $G$-equivariant in general \pcref{thm:zbsgad_pseudo}.  This means that $\zbsgad$ commutes with the $G$-actions of its domain and codomain up to coherent natural isomorphisms.  Throughout this chapter, $\Op$ is assumed to be a \emph{$\Uinf$-operad} \pcref{as:OpA'}, meaning that it is a 1-connected $\Gcat$-operad that is levelwise a nonempty translation category.

\organization
This chapter consists of the following sections.

\secname{sec:sgo_tprod}
This section constructs the pointed functor $\zbsgad$ from the $\angordnbe$-twisted product $\proAnbe$ to the category $\Asgsmaangordnbe$ of strong $(\sma\angordnbe)$-systems in $\A$.

\secname{sec:sgo_prod_eq}
This section proves that the functors $(\zbsgad,\zbsg)$ form an adjoint equivalence.  The unit and counit,
\[1_{\proAnbe} \fto{\unis} \zbsg\zbsgad \andspace \zbsgad\zbsg \fto[\iso]{\cous} 1_{\Asgsmaangordnbe},\]
are, respectively, the identity and a natural isomorphism.  Replacing the category $\Asgsmaangordnbe$ by the larger category $\Asmaangordnbe$ of all $(\sma\angordnbe)$-systems yields an adjunction
\begin{equation}\label{zb_chi}
\begin{tikzpicture}[vcenter]
\def\s{22}
\draw[0cell]
(0,0) node (a1) {\phantom{X}}
(a1)++(1.8,0) node (a2) {\phantom{X}}
(a1)++(-.15,0) node (a1') {\proAnbe}
(a2)++(.6,-.04) node (a2') {\Asmaangordnbe}
;
\draw[0cell=.8]
node[between=a1 and a2 at .5] {\perp}
;
\draw[1cell=.9]
(a1) edge[bend left=\s] node {\zbad} (a2)
(a2) edge[bend left=\s] node {\zb} (a1)
;
\end{tikzpicture}
\end{equation}
whose counit is not a natural isomorphism in general.

\secname{sec:pseudoequiv}
This section proves that the left adjoint $\zbsgad$ is a pseudo $G$-equivariant functor.  For a nontrivial group $G$, the functor $\zbsgad$ is \emph{not} $G$-equivariant even for the Barratt-Eccles $\Gcat$-operad $\BE$, on which $G$ acts trivially \pcref{ex:al_g_s}.  The left adjoint $\zbad$ is also pseudo $G$-equivariant.

\section{Shimakawa Strong $H$-Theory from Twisted Products}
\label{sec:sgo_tprod}

This section constructs a pointed functor 
\[\proAnbe \fto{\zbsgad} \Asgsmaangordnbe\]
from the twisted product to Shimakawa strong $H$-theory for each $\Uinf$-operad $\Op$, $\Op$-pseudoalgebra $\A$, and object $\angordnbe \in \GG \setminus \{\vstar,\ang{}\}$.  \cref{sec:sgo_prod_eq} shows that the functors $(\zbsgad,\zbsg)$ form a nonequivariant adjoint equivalence, where
\[\Asgsmaangordnbe \fto{\zbsg} \proAnbe\]
is the pointed $G$-functor from Shimakawa strong $H$-theory to the twisted product \cref{sgosgtoprod}.

\secoutline
\begin{itemize}
\item \cref{as:OpA'} fixes a $\Uinf$-operad $\Op$, which strengthens a $\Tinf$-operad \pcref{as:OpA}. 
\item \cref{expl:OpA'} discusses the reasons for assuming a $\Uinf$-operad $\Op$.  The Barratt-Eccles operad $\BE$ and the $G$-Barratt-Eccles operad $\GBE$ are $\Uinf$-operads \pcref{ex:OpA'}.
\item \cref{def:zbsgad} defines the pointed functor $\zbsgad$.
\item \cref{zbsgad_obj_welldef,zbsgad_mor_welldef} prove that $\zbsgad$ is well defined on objects and morphisms.
\end{itemize}

\subsection*{$\Uinf$-Operads}
\cref{as:OpA'} is in effect throughout this chapter.

\begin{assumption}[$\Uinf$-Operads]\label{as:OpA'}\index{operad!U-infinity@$\Uinf$}
For a group $G$, we assume that $(\Op,\ga,\opu)$ is a $\Gcat$-operad that is
\begin{itemize}
\item 1-connected \cref{i_connected}, meaning $\Op(0) = \{*\}$ and $\Op(1) = \{\opu\}$, and 
\item levelwise a nonempty translation category \pcref{def:translation_cat}.  
\end{itemize}
Such a $\Gcat$-operad is called a \emph{$\Uinf$-operad}.  Given a $\Uinf$-operad $\Op$, we choose, once and for all, an arbitrary object 
\begin{equation}\label{Opn_object}
\ep_n \in \Op(n)
\end{equation}
for each $n \geq 0$. 
\end{assumption}

\begin{explanation}[$\Uinf$]\label{expl:OpA'}
By \cref{translation_pseudocom}, a $\Uinf$-operad $\Op$ has a unique pseudo-commutative structure \pcref{def:pseudocom_operad} and, in particular, is a $\Tinf$-operad \pcref{as:OpA}.  The $\Uinf$ assumption on $\Op$ is used in the following ways.
\begin{enumerate}
\item To define the functor $\zbsgad$---from the twisted product $\proAnbe$ \pcref{def:proCnbe} to the category $\Asgsmaangordnbe$ of strong $(\sma\angordnbe)$-systems \cref{Asgordnbe}---we need to be able to glue together objects in an $\Op$-pseudoalgebra $\A$ parametrized by arbitrary subsets of $\txprod_{j \in \ufs{q}}\, \ufs{n}_j$.  The objects $\ep_n \in \Op(n)$ \cref{Opn_object} are used for this purpose; see \cref{zbsgad_comp_obj,zbsgad_gl,zbsgad_comp_mor} in \cref{def:zbsgad}.
\item In the course of proving that $\zbsgad$ is a well-defined functor, it is necessary to show that the objects $\ep_n \in \Op(n)$ are appropriately compatible with the axioms of $(\sma\angordnbe)$-systems and their morphisms.  The assumption that each $\Op(n)$ is a translation category is used for this purpose; see \cref{al_Op} and the proofs of \cref{zbsgad_obj_welldef,zbsgad_mor_welldef}.
\item The construction and proof that $\zbsgad$ is pseudo $G$-equivariant also use the $\Uinf$ assumption on $\Op$.  See \cref{al_g_s} and the proofs of \cref{pse_g_a_welldef,thm:zbsgad_pseudo}.
\end{enumerate}
The objects $\ep_n \in \Op(n)$ can be chosen arbitrarily because $\Op(n)$ is a translation category.
\end{explanation}

\begin{example}\label{ex:OpA'}
Here are some examples of $\Uinf$-operads \pcref{as:OpA'}.
\begin{enumerate}
\item\label{ex:OpA'_i} For the Barratt-Eccles $\Gcat$-operad $\BE$ \pcref{def:BE}, we can choose the identity permutation
\begin{equation}\label{BE_epn}
\ep_n = \id_n \in \ESigma_n \forspace n \geq 0.
\end{equation}
\item\label{ex:OpA'_ii} For the $G$-Barratt-Eccles operad $\GBE$ \pcref{def:GBE}, we can choose
\begin{equation}\label{GBE_epn}
\big(\EG \fto{\ep_n} \ESigma_n\big) \in \GBE(n) = \Catg(\EG,\ESigma_n)
\end{equation}
to be the constant functor at $\id_n$ for each $n \geq 0$.
\item\label{ex:OpA'_iii}
For a $\Uinf$-operad $\Op$, the $\Gcat$-operad $\Oph = \Catg(\EG,\Op)$ in \cref{catgego} \eqref{ng_i} is also a $\Uinf$-operad, since
\[\Oph(n) = \Catg(\EG,\Op(n)) \forspace n \geq 0.\]
Given the choice of an object $\ep_n \in \Op(n)$, we can choose the functor 
\[\big(\EG \fto{\ep_n} \Op(n)\big) \in \Catg(\EG,\Op(n)) = \Catg(\EG,\Op)(n)\]
to be the constant functor at $\ep_n \in \Op(n)$.  This is how the constant functor in \cref{GBE_epn} arises from the permutation in \cref{BE_epn}.
\end{enumerate}
We reiterate that $\ep_n$ can be any object in $\Op(n)$.  In the examples $\BE$, $\GBE$, and $\Oph$, we point out the simplest choices of $\ep_n$. 
\end{example}

\subsection*{Inverse of $\zbsg$}

We now define an adjoint inverse $\zbsgad$ of the pointed $G$-functor $\zbsg$ \cref{sgosgtoprod}.  

\begin{definition}\label{def:zbsgad}
Given a $\Uinf$-operad $\Op$ \pcref{as:OpA'}, an $\Op$-pseudoalgebra $(\A,\gaA,\phiA)$ \pcref{def:pseudoalgebra}, and an object \cref{GG_objects}
\[\angordnbe = \sordi{n}{\be}{j}_{j \in \ufs{q}} \in \GG \setminus \{\vstar,\ang{}\}\]
of length $q>0$, we define a pointed functor
\begin{equation}\label{zbsgad}
\proAnbe \fto{\zbsgad} \Asgsmaangordnbe
\end{equation}
as follows.  
\begin{description}
\item[Domain]
The $\angordnbe$-twisted product $\proAnbe$ is defined in \cref{def:proCnbe}.  
\item[Codomain]
The category $\Asgsmaangordnbe$ of strong $(\sma\angordnbe)$-systems is defined in \cref{def:smashFGGG,def:nsys,def:nsys_morphism,def:nsys_gcat}.  
\item[Component objects]
Given an object \cref{proCnbe_object}
\[a = \ang{a_{\bdi} \in \A}_{\bdi \sins \ufs{n_1 n_2 \cdots n_q}} \in \proAnbe\]
with $\bdi = \ang{i_j \in \ufs{n}_j}_{j \in \ufs{q}}$, the strong $(\sma\angordnbe)$-system \cref{nsys} 
\begin{equation}\label{zbsgad_a}
(\zbsgad a, \gl) \in \Asgsmaangordnbe
\end{equation} 
has, for each lexicographically ordered subset \cref{lex_bijection}
\begin{equation}\label{subset_s}
s \subseteq \txprod_{j \in \ufs{q}}\, \ufs{n}_j = \ufs{\txprod_{j \in \ufs{q}}\, n_j}
\end{equation}
of cardinality $\cards$, $s$-component object \cref{nsys_s} defined as 
\begin{equation}\label{zbsgad_comp_obj}
(\zbsgad a)_s = \gaA_{\cards}\big(\ep_{\cards} \sscs \ang{a_{\bdi}}_{\bdi \in s}\big) \in \A.
\end{equation}
The right-hand side of \cref{zbsgad_comp_obj} is given as follows.
\begin{itemize}
\item The functor
\[\Op(\cards) \times \A^{\cards} \fto{\gaA_{\cards}} \A\]
is the $\cards$-th $\Op$-action $G$-functor of $\A$ \cref{gaAn}.
\item $\ep_{\cards} \in \Op(\cards)$ is the object chosen in \cref{Opn_object}.
\item The $\cards$-tuple $\ang{a_{\bdi}}_{\bdi \in s} \in \A^{\cards}$ uses the lexicographic ordering of $s$.
\end{itemize}
\item[Gluing]
Given an object $x \in \Op(r)$ with $r \geq 0$, a subset $s \subseteq \txprod_{j \in \ufs{q}}\, \ufs{n}_j$, and a partition 
\[s = \txcoprod_{\ell \in \ufs{r}}\, s_\ell \subseteq \txprod_{j \in \ufs{q}}\, \ufs{n}_j,\]
we first define the following permutation and isomorphism.
\begin{itemize}
\item We define the permutation 
\begin{equation}\label{sigma_angsi}
\Big(s \fto[\iso]{\sig_{s,\ang{s_\ell}_{\ell \in \ufs{r}}}} \txcoprod_{\ell \in \ufs{r}}\, s_\ell\Big) \in \Si_{\cards}
\end{equation}
that reorders $s$ to $\coprod_{\ell \in \ufs{r}}\, s_\ell$, where $s$ and each $s_\ell$ are equipped with the lexicographic ordering inherited from $\txprod_{j \in \ufs{q}}\, \ufs{n}_j$. 
\item We define the unique isomorphism 
\begin{equation}\label{al_Op}
\ga\big(x \sscs \ang{\ep_{\cardsl}}_{\ell \in \ufs{r}}\big) \sig_{s,\ang{s_\ell}_{\ell \in \ufs{r}}} 
\fto[\iso]{\alp_{x;\, s, \ang{s_\ell}_{\ell \in \ufs{r}}}} \ep_{\cards} 
\end{equation}
in the translation category $\Op(\cards)$ with the indicated domain and codomain.
\end{itemize}
The gluing isomorphism \cref{gl-morphism} of the strong $(\sma\angordnbe)$-system $(\zbsgad a, \gl)$ at $(x; s, \ang{s_\ell}_{\ell \in \ufs{r}})$ is defined as the following composite isomorphism in $\A$.
\begin{equation}\label{zbsgad_gl}
\begin{tikzpicture}[vcenter]
\def\h{1} \def\u{-1.2} \def\v{-1.5}
\draw[0cell=.9]
(0,0) node (a11) {\gaA_r\big(x \sscs \ang{(\zbsgad a)_{s_\ell}}_{\ell \in \ufs{r}} \big)}
(a11)++(4.2,0) node (a12) {(\zbsgad a)_s}
(a11)++(0,\u) node (a21) {\gaA_r\Big(x \sscs \bang{\gaA_{\cardsl} \big( \ep_{\cardsl} \sscs \ang{a_{\bdi}}_{\bdi \in s_\ell}\big)}_{\ell \in \ufs{r}} \Big)}
(a12)++(0,\u) node (a22) {\gaA_{\cards}\big(\ep_{\cards} \sscs \ang{a_{\bdi}}_{\bdi \in s}\big)}
(a21)++(0,\v-.3) node (a3) {\gaA_{\cards} \Big(\ga(x \sscs \ang{\ep_{\cardsl}}_{\ell \in \ufs{r}}) \sscs \ang{\ang{a_{\bdi}}_{\bdi \in s_\ell}}_{\ell \in \ufs{r}} \Big)}
(a3)++(0,\v+.3) node (a4) {\gaA_{\cards} \Big(\ga(x \sscs \ang{\ep_{\cardsl}}_{\ell \in \ufs{r}})\, \sig_{s,\ang{s_\ell}_{\ell \in \ufs{r}}} \spsc \ang{a_{\bdi}}_{\bdi \in s} \Big)}
;
\draw[1cell=.9]
(a11) edge node {\gl_{x;\, s, \ang{s_\ell}_{\ell \in \ufs{r}}}} (a12)
(a12) edge[equal,transform canvas={xshift=0ex}] (a22)
(a11) edge[equal,transform canvas={xshift=0em}, shorten >=-.5ex] (a21)
(a21) edge node[swap] {\phiA_{(r;\, \card{s_1},\ldots,\card{s_r})}} node {\iso} (a3)
(a3) edge[equal,shorten <=-.5ex,shorten >=-.5ex] node[swap,pos=.6] {\mathbf{eq}} (a4)
(a4) [rounded corners=2pt, shorten <=-0ex] -| node[swap,pos=.65] {\gaA_{\cards}(\alp_{x;\, s, \ang{s_\ell}_{\ell \in \ufs{r}}} \sscs 1^{\cards})} node[pos=.65] {\iso} (a22)
;
\end{tikzpicture}
\end{equation}
The arrows in \cref{zbsgad_gl} are given as follows.
\begin{itemize}
\item The top left and right equalities follow from \cref{zbsgad_comp_obj}.
\item The isomorphism $\phiA_{(r;\, \card{s_1},\ldots,\card{s_r})}$ is a component of the associativity constraint of $\A$ \cref{phiA_component}.
\item The equality labeled $\mathbf{eq}$ follows from the action equivariance axiom \cref{pseudoalg_action_sym} for $\A$ and the object equality 
\[\sig_{s,\ang{s_\ell}_{\ell \in \ufs{r}}} \ang{a_{\bdi}}_{\bdi \in s} 
= \ang{\ang{a_{\bdi}}_{\bdi \in s_\ell}}_{\ell \in \ufs{r}}\]
using the permutation defined in \cref{sigma_angsi}.
\item $\alp_{x;\, s, \ang{s_\ell}_{\ell \in \ufs{r}}}$ is the unique isomorphism defined in \cref{al_Op}.
\end{itemize}
\cref{zbsgad_obj_welldef} proves that $(\zbsgad a,\gl)$ is a strong $(\sma\angordnbe)$-system in $\A$.
\item[Morphisms]
For a morphism 
\[f = \ang{f_{\bdi} \in \A}_{\bdi \sins \ufs{n_1 n_2 \cdots n_q}} \cn a \to b \inspace \proAnbe\]
and a subset $s \subseteq \txprod_{j \in \ufs{q}}\, \ufs{n}_j$, the morphism
\begin{equation}\label{zbsgad_f}
\zbsgad a \fto{\zbsgad f} \zbsgad b \inspace \Asgsmaangordnbe
\end{equation} 
has $s$-component morphism \cref{theta_s} defined as 
\begin{equation}\label{zbsgad_comp_mor}
(\zbsgad f)_s = \gaA_{\cards} \big(1_{\ep_{\cards}} \sscs \ang{f_{\bdi}}_{\bdi \in s}\big) \inspace \A.
\end{equation}
\cref{zbsgad_mor_welldef} proves that $\zbsgad f$ is a morphism of strong $(\sma\angordnbe)$-systems in $\A$.
\item[Functoriality]
The functoriality of $\zbsgad$ follows from the functoriality of $\gaA$ and the fact that identities and composition are defined coordinatewise in the $\angordnbe$-twisted product $\proAnbe$ and componentwise in $\Asgsmaangordnbe$.
\item[Basepoint preservation] 
The basepoint $\zero \in \proAnbe$ is given coordinatewise by $\zero = \gaA_0(*) \in \A$.  It is sent by $\zbsgad$ to the base $(\sma\angordnbe)$-system $(\zero,1_\zero)$ in $\A$ for the following reasons. 
\begin{itemize}
\item The $s$-component object \cref{zbsgad_comp_obj} 
\[(\zbsgad \zero)_s = \gaA_{\cards}\big(\ep_{\cards} \sscs \ang{\zero}_{\bdi \in s}\big)\] 
is equal to $\zero \in \A$ by $\cards$ applications of the basepoint axiom \cref{pseudoalg_basept_axiom} for $\A$.
\item The gluing isomorphism $\gl_{x;\, s,\ang{s_\ell}_{\ell \in \ufs{r}}}$ \cref{zbsgad_gl} is equal to the identity morphism $1_\zero$ for the following reasons.
\begin{itemize}
\item The associativity constraint $\phiA_{(r;\, \card{s_1},\ldots,\card{s_r})}$ is equal to $1_\zero$ by the composition axiom \cref{pseudoalg_comp_axiom} for $\A$ and \cref{phi_id}.
\item The isomorphism $\gaA_{\cards}(\alp ;\, 1^{\cards})$ is equal to $1_\zero$ by the basepoint axiom \cref{pseudoalg_basept_axiom} for $\A$ and the naturality of the associativity constraint $\phiA$.
\end{itemize}
\end{itemize}
\end{description}
This finishes the definition of the pointed functor $\zbsgad$ in \cref{zbsgad}.
\end{definition}

\subsection*{Well Definedness}
The rest of this section proves \cref{zbsgad_obj_welldef,zbsgad_mor_welldef}, which ensure that $\zbsgad$ is a well-defined pointed functor.  Recall that $\Op$ is a $\Uinf$-operad \pcref{as:OpA'} and $(\A,\gaA,\phiA)$ is an $\Op$-pseudoalgebra \pcref{def:pseudoalgebra}.

\begin{lemma}\label{zbsgad_obj_welldef}
The pair $(\zbsgad a,\gl)$ defined in \crefrange{zbsgad_comp_obj}{zbsgad_gl} is a strong $(\sma\angordnbe)$-system in $\A$.
\end{lemma}

\begin{proof}
Each component of $\gl$ \cref{zbsgad_gl} is an isomorphism.  We verify that the pair $(\zbsgad a,\gl)$ satisfies the axioms \crefrange{nsys_obj_unity}{nsys_associativity} in \cref{def:nsys} for a $(\sma\angordnbe)$-system, where, by definition \cref{smash_GGobjects},
\[\sma\angordnbe = \ord{n_1 \cdots n_q}^{\sma_{j \in \ufs{q}}\, \be_j} \in \FG.\]
In the rest of this proof, $s$ is a subset of $\ufs{n_1 \cdots n_q}$, and $s = \coprod_{\ell \in \ufsr} s_\ell$ is a partition.
\begin{description}
\item[Object unity]  
The axiom \cref{nsys_obj_unity} follows from \cref{pseudoalg_zero,i_connected,zbsgad_comp_obj} as follows:
\[(\zbsgad a)_\emptyset = \gaA_0(\ep_0) = \gaA_0(*) = \zero \in \A.\]
\item[Naturality]  
Given a morphism $h \cn x \to y$ in $\Op(r)$ with $r \geq 0$, the naturality diagram \cref{nsys_naturality} for $(\zbsgad a,\gl)$ is the boundary diagram in \cref{zbsgad_naturality}, where 
\[\phiA = \phiA_{(r;\, \card{s_1}, \ldots, \card{s_r})} \cq 
\sig = \sig_{s,\ang{s_\ell}_{\ell \in \ufs{r}}}, \andspace 
\alp_? = \alp_{?;\, s, \ang{s_\ell}_{\ell \in \ufs{r}}}.\]
\begin{equation}\label{zbsgad_naturality}

\end{equation}
The diagram \cref{zbsgad_naturality} commutes for the following reasons.
\begin{itemize}
\item The left quadrilateral commutes by the naturality of the associativity constraint $\phiA$ \cref{phiA}.
\item The middle quadrilateral commutes by the action equivariance axiom \cref{pseudoalg_action_sym} for $\A$.
\item The right quadrilateral commutes by the functoriality of $\gaA_{\cards}$ and the commutativity of the diagram
\begin{equation}\label{zb_natdiag}
%
\end{equation}
in the translation category $\Op(\cards)$.
\end{itemize}
\item[Unity]  
The first unity axiom \cref{nsys_unity_empty} states that, if $s=\emptyset$---which implies $s_\ell = \emptyset$ for each $\ell \in \ufsr$---then 
\[\gl_{x;\, \emptyset, \ang{\emptyset}_{\ell \in \ufsr}} = 1_\zero \inspace \A.\]
This axiom holds for the following reasons.
\begin{itemize}
\item In \cref{zbsgad_gl}, the associativity constraint $\phiA_{(r;\, 0,\ldots,0)}$ is equal to $1_\zero$ by repeated applications of the basepoint axiom \cref{pseudoalg_basept_axiom} for $\A$.
\item The isomorphism $\alp_{x;\, \emptyset, \ang{\emptyset}_{\ell \in \ufsr}}$ \cref{al_Op} in $\Op(0) = \{*\}$ is the identity morphism $1_*$.
\end{itemize}
The second unity axiom \cref{nsys_unity_one} states that, if $r=1$---which implies $x = \opu \in \Op(1)$---then 
\[\gl_{\opu;\, s,\,s} = 1 \inspace \A.\]
This axiom holds for the following reasons.
\begin{itemize}
\item In \cref{zbsgad_gl}, the associativity constraint $\phiA_{(1;\,|s|)}$ is an identity morphism by the bottom half of the unity axiom \cref{pseudoalg_unity} for $\A$.
\item The permutation $\sig_{s,s} \cn s \to s$ \cref{sigma_angsi} is the identity.  Thus, the isomorphism \cref{al_Op}
\[\ga(\opu;\, \ep_{\cards}) \sig_{s,s} = \ep_{\cards} \fto{\alp_{\opu;\, s,\,s}} \ep_{\cards}\]
is the identity morphism $1_{\ep_{\cards}}$ in the translation category $\Op(\cards)$.
\end{itemize}
\item[Equivariance]  
For each permutation $\tau \in \Si_r$, the equivariance diagram \cref{nsys_equivariance} for $(\zbsgad a,\gl)$ is the boundary diagram in \cref{zbsgad_equiv}, where 
\[\ang{\Cdots}_{\ell} = \ang{\Cdots}_{\ell \in \ufs{r}} \cq 
\ang{\Cdots\,_{\bdi \in ?}} = \ang{\Cdots\,_{\bdi}}_{\bdi \in ?}, \andspace 
\ell' = \tauinv\ell\]
for $\ell \in \ufs{r}$.
\setlength{\fboxsep}{1pt}
\setlength{\fboxrule}{0.1pt}
\begin{equation}\label{zbsgad_equiv}

\end{equation}
The arrows in the diagram \cref{zbsgad_equiv} are given as follows.
\begin{itemize}
\item Starting from the upper left corner and using \cref{zbsgad_gl}, the composites along the left and right boundaries are, respectively, $\gl_{x\tau;\, s, \ang{s_\ell}_\ell}$ and $\gl_{x;\, s, \ang{s_{\ell'}}_\ell}$.
\item Each equality labeled $\mathbf{eq}$ holds by the action equivariance axiom \cref{pseudoalg_action_sym} for $\A$.
\item Each equality labeled $\mathbf{t}$ holds by the top equivariance axiom for the $\Gcat$-operad $\Op$ \pcref{def:voperad}, where 
$\taubar \in \Si_s$ is the block permutation induced by $\tau \in \Si_r$ that permutes $r$ consecutive blocks of lengths $\card{s_1}, \ldots, \card{s_r}$.
\item The equality labeled $\mathbf{p}$ holds because, by \cref{sigma_angsi}, the following diagram of permutations commutes.
\begin{equation}\label{sigma_taubar}
\begin{tikzpicture}[vcenter]
\def\h{1.5} \def\v{-1}
\draw[0cell]
(0,0) node (a) {s}
(a)++(-\h,\v) node (b) {\txcoprod_{\ell \in \ufs{r}}\, s_\ell}
(a)++(\h,\v) node (c) {\txcoprod_{\ell \in \ufs{r}}\, s_{\ell'}}
;
\draw[1cell=.9]
(a) edge node[swap] {\sig_{s,\ang{s_\ell}_\ell}} (b)
(b) edge node {\taubar} (c)
(a) edge node {\sig_{s,\ang{s_{\ell'}}_\ell}} (c)
;
\end{tikzpicture}
\end{equation}
\end{itemize}
Regions $\fbox{a}$ through $\fbox{c}$ in the diagram \cref{zbsgad_equiv} commute for the following reasons.
\begin{enumerate}
\item Region $\fbox{a}$ commutes by the top equivariance axiom \cref{pseudoalg_topeq} for $\A$.
\item Region $\fbox{b}$ commutes because it consists entirely of equalities.
\item Region $\fbox{c}$ commutes by the functoriality of $\gaA_{\cards}$ \cref{gaAn} and the commutativity of the diagram
\begin{equation}\label{region_c}
\begin{tikzpicture}[vcenter]
\def\h{4} 
\draw[0cell]
(0,0) node (a) {\ga(x\tau \sscs \ang{\ep_{\cardsl}}_{\ell}) \sig_{s,\ang{s_\ell}_\ell}}
(a)++(\h,0) node (b) {\phantom{\ga(x \sscs \ang{\ep_{\card{s_{\ell'}}}}_{\ell}) \sig_{s,\ang{s_{\ell'}}_\ell}}}
(b)++(0,-.04) node (b') {\ga(x \sscs \ang{\ep_{\card{s_{\ell'}}}}_{\ell}) \sig_{s,\ang{s_{\ell'}}_\ell}}
(a)++(\h/2,-1) node (c) {\phantom{\ep_{\cards}}}
(c)++(0,-.05) node (c') {\ep_{\cards}}
;
\draw[1cell=.9]
(a) edge[equal] (b)
(b') [rounded corners=2pt] |- node[swap,pos=.25] {\alp_{x;\, s, \ang{s_{\ell'}}_\ell}} (c)
;
\draw[1cell=.9]
(a) [rounded corners=2pt] |- node[swap,pos=.25] {\alp_{x\tau;\, s, \ang{s_\ell}_\ell}} (c)
;
\end{tikzpicture}
\end{equation}
in the translation category $\Op(\cards)$.
\end{enumerate}
\item[Associativity]  
To prove the associativity axiom \cref{nsys_associativity} for $(\zbsgad a,\gl)$, we consider objects
\[(x \sscs \ang{x_\ell}_{\ell \in \ufs{r}}) \in \Op(r) \times \txprod_{\ell\in \ufs{r}}\, \Op(t_\ell),\]
partitions
\[s = \txcoprod_{\ell \in \ufs{r}}\, s_\ell \andspace s_\ell = \txcoprod_{p \in \ufs{t}_\ell}\, s_{\ell,p}\]
for each $\ell \in \ufsr$, and the following notation, with $\sig_{s,\ang{s_\ell}_\ell}$ the permutation in \cref{sigma_angsi} and $\alp_{x;\, s,\ang{s_\ell}_\ell}$ the unique isomorphism in \cref{al_Op}.
\[\left\{\begin{aligned}
t &= \txsum_{\ell \in \ufs{r}}\, t_\ell & \ang{\Cdots}_\ell &= \ang{\Cdots}_{\ell \in \ufs{r}} & \bx &= \ga(x \sscs \ang{x_\ell}_\ell) \in \Op(t)\\
\ang{\Cdots\,_{\bdi \in ?}} &= \ang{\Cdots\,_{\bdi}}_{\bdi \in ?} & \ang{\Cdots}_p &= \ang{\Cdots}_{p \in \ufs{t}_\ell} & y_\ell &= \ga(x_\ell \sscs \ang{\ep_{\card{s_{\ell,p}}}}_p) \in \Op(\cardsl)\\
&& \ang{\Cdots}_{p,\ell} &= \ang{\ang{\Cdots}_p}_\ell & y &= \ga(\bx \sscs \ang{\ep_{\card{s_{\ell,p}}}}_{p,\ell}) \in \Op(\cards)
\end{aligned}\right.\]
\[\left\{\scalebox{.85}{$\begin{split}
\sig &= \sig_{s, \ang{s_\ell}_\ell} \cn s \fiso \txcoprod_{\ell \in \ufs{r}}\, s_\ell\\
\sig' &= \sig_{s, \ang{s_{\ell,p}}_{p,\ell}} \cn s \fiso \txcoprod_{\ell \in \ufs{r}}\, \txcoprod_{p \in \ufs{t}_\ell}\, s_{\ell,p}\\
\sig_\ell &= \sig_{s_\ell, \ang{s_{\ell,p}}_{p}} \cn s_\ell \fiso \txcoprod_{p \in \ufs{t}_\ell}\, s_{\ell,p}\\
\sig^\stimes &= \txcoprod_{\ell \in \ufs{r}}\, \sig_\ell \cn \txcoprod_{\ell \in \ufs{r}}\, s_\ell \fiso \txcoprod_{\ell \in \ufs{r}}\, \txcoprod_{p \in \ufs{t}_\ell}\, s_{\ell,p}
\end{split}$}
\right.
\quad
\left\{\scalebox{.8}{$\begin{split}
\alp_x &= \alp_{x;\, s, \ang{s_\ell}_\ell} \cn \ga(x \sscs \ang{\ep_{\card{s_\ell}}}_\ell)\sig \fiso \ep_{\cards}\\
\alp_{\bx} &= \alp_{\bx;\, s, \ang{s_{\ell,p}}_{p,\ell}} \cn y\sig' \fiso \ep_{\cards}\\
\alp_{x_\ell} &= \alp_{x_\ell;\, s_\ell, \ang{s_{\ell,p}}_{p}} \cn y_\ell \sig_\ell \fiso \ep_{\card{s_\ell}}
\end{split}$}
\right.\]
The associativity diagram \cref{nsys_associativity} for $(\zbsgad a,\gl)$ is the boundary diagram in \cref{zbsgad_assoc}, where subscripts of the associativity constraint $\phiA$ are omitted to save space.
\setlength{\fboxsep}{1pt}
\setlength{\fboxrule}{0.1pt}
\begin{equation}\label{zbsgad_assoc}

\end{equation}
Regions $\fbox{I}$ through $\fbox{VI}$ in the diagram \cref{zbsgad_assoc} commute for the following reasons.
\begin{enumerate}
\item In region $\fbox{I}$, the equality labeled $\mathbf{a}$ holds by the associativity axiom for the $\Gcat$-operad $\Op$, which yields
\begin{equation}\label{gaxyelly}
\begin{split}
\ga\big(x\sscs \ang{y_\ell}_\ell\big) 
&= \ga\big(x\sscs \ang{\ga(x_\ell \sscs \ang{\ep_{\card{s_{\ell,p}}}}_p)}_\ell\big)\\
&= \ga\big(\ga\big(x \sscs \ang{x_\ell}_\ell\big) \sscs \ang{\ep_{\card{s_{\ell,p}}}}_{p,\ell} \big)\\
&= \ga\big(\bx \sscs \ang{\ep_{\card{s_{\ell,p}}}}_{p,\ell} \big) = y.
\end{split}
\end{equation}
Region $\fbox{I}$ commutes by the composition axiom \cref{pseudoalg_comp_axiom} for $\A$.
\item In region $\fbox{II}$, each equality labeled $\mathbf{eq}$ holds by the action equivariance axiom \cref{pseudoalg_action_sym} for $\A$.  The equality labeled $\mathbf{b}$ holds by the bottom equivariance axiom for the $\Gcat$-operad $\Op$.  The equality labeled $\mathbf{b}'$ holds by the equalities
\[\ga\big(x \sscs \ang{y_\ell \sig_\ell}_\ell \big)\sig 
= \ga\big(x \sscs \ang{y_\ell}\big) (\sig^\stimes \sig)= y\sig'.\]
\begin{itemize}
\item The first equality holds by the bottom equivariance axiom and the symmetric group action axiom for the $\Gcat$-operad $\Op$.
\item The second equality holds by \cref{gaxyelly} and the commutativity of the following diagram of permutations, which, in turn, follows from the definition of $\sig_{s,\ang{s_\ell}_\ell}$ \cref{sigma_angsi}.
\begin{equation}\label{sigmas_diag}
\begin{tikzpicture}[vcenter]
\def\h{1.5} \def\v{-1}
\draw[0cell]
(0,0) node (a) {s}
(a)++(-\h,\v) node (b) {\txcoprod_{\ell \in \ufs{r}}\, s_\ell}
(a)++(\h,\v) node (c) {\phantom{\txcoprod_{\ell \in \ufs{r}}\, s_\ell}}
(c)++(.6,0) node (c') {\txcoprod_{\ell \in \ufs{r}}\, \txcoprod_{p \in \ufs{t}_\ell}\, s_{\ell,p}}
;
\draw[1cell=.9]
(a) edge node[swap] {\sig} (b)
(b) edge node {\sig^\stimes} (c)
(a) edge node {\sig'} (c)
;
\end{tikzpicture}
\end{equation}
\end{itemize}
Region $\fbox{II}$ commutes because it consists entirely of equalities.
\item Region $\fbox{III}$ commutes by the functoriality of $\gaA_{\cards}$ and the commutativity of the diagram
\begin{equation}\label{region_iii}
\begin{tikzpicture}[vcenter]
\def\v{-1.5}
\draw[0cell]
(0,0) node (a11) {\ga(x\sscs \ang{y_\ell \sig_\ell}_\ell)\sig}
(a11)++(3,0) node (a12) {y\sig'}
(a11)++(0,\v) node (a21) {\ga(x \sscs \ang{\ep_{\cardsl}}_\ell)\sig}
(a12)++(0,\v) node (a22) {\phantom{\ep_{\cards}}}
(a22)++(0,-.05) node (a22') {\ep_{\cards}}
;
\draw[1cell=.9]
(a11) edge[equal] (a12)
(a12) edge node {\alp_{\bx}} (a22')
(a11) edge[transform canvas={xshift=1.75em}] node[swap] {\ga(x; \ang{\alp_{x_\ell}}_\ell)\sig} (a21)
(a21) edge node {\alp_x} (a22)
;
\end{tikzpicture}
\end{equation}
in the translation category $\Op(\cards)$.
\item Region $\fbox{IV}$ commutes by the action equivariance axiom \cref{pseudoalg_action_sym} for $\A$.
\item Region $\fbox{V}$ commutes by the naturality of $\phiA$ \cref{phiA}.
\item In region $\fbox{VI}$, the left equality is given by
\[\gaA_{\cardsl}\big(y_\ell \sscs \ang{\ang{a_{\bdi \in s_{\ell,p}}}}_p \big) 
= \gaA_{\cardsl}\big(y_\ell \sig_\ell \sscs \ang{a_{\bdi \in s_{\ell}}} \big),\]
which follows from the action equivariance axiom \cref{pseudoalg_action_sym} for $\A$.  Region $\fbox{VI}$ commutes by the bottom equivariance axiom \cref{pseudoalg_boteq} for $\A$.
\end{enumerate}
This proves the associativity axiom \cref{nsys_associativity} for $(\zbsgad a,\gl)$, proving that it is a strong $(\sma\angordnbe)$-system in $\A$.\qedhere
\end{description}
\end{proof}

\begin{lemma}\label{zbsgad_mor_welldef}
For each morphism $f \cn a \to b$ in $\proAnbe$, the collection 
\[\zbsgad a \fto{\zbsgad f} \zbsgad b\]
defined in \cref{zbsgad_f} is a morphism in $\Asgsmaangordnbe$.
\end{lemma}

\begin{proof}
We verify that $\zbsgad f$ satisfies the axioms \cref{nsys_mor_unity,nsys_mor_compat} for a morphism of $(\sma\angordnbe)$-systems.
\begin{description}
\item[Unity]  
The unity axiom \cref{nsys_mor_unity} for $\zbsgad f$ holds by \cref{zbsgad_comp_mor}, \cref{i_connected}, and the functoriality of $\gaA_0$ as follows:
\[(\zbsgad f)_\emptyset = \gaA_0(1_{\ep_0}) = \gaA_0(1_*) = 1_\zero \inspace \A.\]
\item[Compatibility]  
We consider the gluing morphism $\gl_{x;\, s,\ang{s_\ell}_{\ell \in \ufs{r}}}$ \cref{zbsgad_gl} for both $\zbsgad a$ and $\zbsgad b$, along with the following abbreviation.
\[\begin{aligned}
\ang{\Cdots}_\ell &= \ang{\Cdots}_{\ell \in \ufs{r}} & \sig &= \sig_{s,\ang{s_\ell}_\ell} 
& \phiA &= \phiA_{(r;\, \card{s_1}, \ldots, \card{s_r})}\\
\ang{\Cdots\,_{\bdi \in ?}} &= \ang{\Cdots\,_{\bdi}}_{\bdi \in ?} & \alp_x &= \alp_{x;\, s, \ang{s_{\ell}}_\ell} &&
\end{aligned}\]
The compatibility diagram \cref{nsys_mor_compat} for $\zbsgad f$ is the boundary diagram in \cref{zbsgad_mor_compat}.
\begin{equation}\label{zbsgad_mor_compat}
\begin{tikzpicture}[vcenter]
\def\h{2.5} \def\u{1.2} \def\v{-2.7}
\draw[0cell=.7]
(0,0) node (a1) {\gaA_r\big(x \sscs \bang{\gaA_{\cardsl}\big( \ep_{\cardsl} \sscs \ang{a_{\bdi \in s_\ell}}\big)}_{\ell} \big)}
(a1)++(\h,\u) node (a2) {\gaA_{\cards}\big(\ga(x \sscs \ang{\ep_{\cardsl}}_{\ell}) \sscs \ang{\ang{a_{\bdi \in s_\ell}}}_{\ell} \big)}
(a2)++(\h,-\u) node (a3) {\gaA_{\cards}\big(\ga(x \sscs \ang{\ep_{\cardsl}}_{\ell}) \sig \sscs \ang{a_{\bdi \in s}} \big)}
(a3)++(\h,\u) node (a4) {\gaA_{\cards}\big(\ep_{\cards} \sscs \ang{a_{\bdi \in s}}\big)}
(a1)++(0,\v) node (b1) {\gaA_r\big(x \sscs \bang{\gaA_{\cardsl}\big( \ep_{\cardsl} \sscs \ang{b_{\bdi \in s_\ell}}\big)}_{\ell} \big)}
(a2)++(0,\v) node (b2) {\gaA_{\cards}\big(\ga(x \sscs \ang{\ep_{\cardsl}}_{\ell}) \sscs \ang{\ang{b_{\bdi \in s_\ell}}}_{\ell} \big)}
(a3)++(0,\v) node (b3) {\gaA_{\cards}\big(\ga(x \sscs \ang{\ep_{\cardsl}}_{\ell}) \sig \sscs \ang{b_{\bdi \in s}} \big)}
(a4)++(0,\v) node (b4) {\gaA_{\cards}\big(\ep_{\cards} \sscs \ang{b_{\bdi \in s}}\big)}
;
\draw[1cell=.7]
(a1) edge node[pos=.2] {\phiA} (a2)
(a2) edge[equal] node[pos=.6] {\mathbf{eq}} (a3)
(a3) edge node[pos=.4] {\gaA_{\cards}(\alp_x; 1^{\cards})} (a4)
(b1) edge node[swap,pos=.8] {\phiA} (b2)
(b2) edge[equal] node[swap,pos=.4] {\mathbf{eq}} (b3)
(b3) edge node[swap,pos=.75] {\gaA_{\cards}(\alp_x; 1^{\cards})} (b4)
(a1) edge[transform canvas={xshift=1em}] node[swap] {\gaA_r(1; \ang{\gaA_{\cardsl}(1 ; \ang{f_{\bdi \in s_\ell}})}_\ell)} (b1)
(a2) edge node[swap,pos=.75] {\gaA_{\cards}(1; \ang{\ang{f_{\bdi \in s_\ell}}}_\ell)} (b2)
(a3) edge node[swap,pos=.25] {\gaA_{\cards}(1; \ang{f_{\bdi \in s}})} (b3)
(a4) edge[transform canvas={xshift=-1em}] node {\gaA_{\cards}(1 ; \ang{f_{\bdi \in s}})} (b4)
;
\end{tikzpicture}
\end{equation}
The diagram \cref{zbsgad_mor_compat} commutes for the following reasons.
\begin{itemize}
\item The left quadrilateral commutes by the naturality of $\phiA$ \cref{phiA}.
\item The middle quadrilateral commutes by the action equivariance axiom \cref{pseudoalg_action_sym} for $\A$.
\item The right quadrilateral commutes by the functoriality of $\gaA_{\cards}$.
\end{itemize}
\end{description}
This proves that $\zbsgad f$ is a morphism of strong $(\sma\angordnbe)$-systems.
\end{proof}

\section{Adjoint Equivalence}
\label{sec:sgo_prod_eq}

Throughout this section, we consider a $\Uinf$-operad $(\Op,\ga,\opu)$ \pcref{as:OpA'}, an $\Op$-pseudoalgebra $(\A,\gaA,\phiA)$ \pcref{def:pseudoalgebra}, and an object $\angordnbe = \ang{\ordn_j^{\be_j}}_{j \in \ufsq} \in \GG \setminus \{\vstar,\ang{}\}$ of length $q>0$ \cref{GG_objects}.  This section proves that there is an adjoint equivalence
\begin{equation}\label{zbsg_adjeq_seci}
\begin{tikzpicture}[vcenter]
\def\s{22}
\draw[0cell]
(0,0) node (a1) {\phantom{X}}
(a1)++(1.8,0) node (a2) {\phantom{X}}
(a1)++(-.15,0) node (a1') {\proAnbe}
(a2)++(.7,-.04) node (a2') {\Asgsmaangordnbe}
node[between=a1 and a2 at .5] {\sim}
;
\draw[1cell=.9]
(a1) edge[bend left=\s] node {\zbsgad} (a2)
(a2) edge[bend left=\s] node {\zbsg} (a1)
;
\end{tikzpicture}
\end{equation}
between the twisted product $\proAnbe$ \pcref{def:proCnbe} and Shimakawa strong $H$-theory \pcref{sys_FGcat} at the pointed finite $G$-set $\sma\angordnbe\in \FG$ \cref{smash_GGobjects}.  The left and right adjoints are the pointed functors $\zbsgad$ \cref{zbsgad} and $\zbsg$ \cref{sgosgtoprod}.  See \cref{thm:zbsg_eq}.

\secoutline
\begin{itemize}
\item \cref{def:zbsgunit} defines the unit $\unis \cn 1 \to \zbsg\zbsgad$.  The composite $\zbsg\zbsgad$ is the identity functor on the twisted product $\proAnbe$, and $\unis$ is the identity natural transformation of the identity functor.
\item \cref{def:zbsgcounit} defines the counit $\cous \cn \zbsgad\zbsg \fiso 1$.  The objects $\ep_n \in \Op(n)$ \cref{Opn_object} play a crucial role in this definition; see \cref{cous_a_s}.
\item \cref{cous_welldef} proves that each component of $\cous$ is an isomorphism of strong $(\sma\angordnbe)$-systems.
\item \cref{thm:zbsg_eq} proves that the quadruple $(\zbsgad,\zbsg,\unis,\cous)$ is an adjoint equivalence.
\item \cref{expl:sys_adjunction} clarifies that the invertibility of the counit $\cous$ requires \emph{strong} systems.  Then it discusses the variant adjunction
\begin{equation}\label{zb_zbad_seci}
\begin{tikzpicture}[vcenter]
\def\s{22}
\draw[0cell]
(0,0) node (a1) {\phantom{X}}
(a1)++(1.8,0) node (a2) {\phantom{X}}
(a1)++(-.15,0) node (a1') {\proAnbe}
(a2)++(.6,-.04) node (a2') {\Asmaangordnbe}
;
\draw[0cell=.8]
node[between=a1 and a2 at .5] {\perp}
;
\draw[1cell=.9]
(a1) edge[bend left=\s] node {\zbad} (a2)
(a2) edge[bend left=\s] node {\zb} (a1)
;
\end{tikzpicture}
\end{equation}
that involves all $(\sma\angordnbe)$-systems, not just the strong ones.
\end{itemize}

\subsection*{Unit}
First, we define the unit of the adjunction $(\zbsgad, \zbsg)$.  

\begin{definition}[Unit]\label{def:zbsgunit}
For the pointed functors $\zbsgad$ \cref{zbsgad} and $\zbsg$ \cref{sgosgtoprod}, we define $\unis$ as the identity natural transformation
\begin{equation}\label{unit_zbsg}
\begin{tikzpicture}[vcenter]
\def\h{1.8} \def\v{-.7} \def\s{15}
\draw[0cell]
(0,0) node (a1) {\proAnbe}
(a1)++(\h,\v) node (a2) {\phantom{A}}
(a2)++(.6,0) node (a2') {\Asgsmaangordnbe}
(a1)++(0,2*\v) node (a3) {\proAnbe}
(a1)++(0,\v) node (a) {}
;
\draw[1cell=.9]
(a1) edge[bend right=\s] node[swap] {1} (a3)
(a1) edge[bend left=\s] node {\zbsgad} (a2)
(a2) edge[bend left=\s, shorten <=-0ex] node {\zbsg} (a3)
;
\draw[2cell]
node[between=a and a2 at .4, shift={(0,.3*\v)}, rotate=0, 2labelw={above,\unis = 1_1,2pt}] {\Rightarrow}
;
\end{tikzpicture}
\end{equation}
of the identity functor on $\proAnbe$.  To verify that $\unis \cn 1 \to \zbsg\zbsgad$ is well defined, we consider an object or a morphism \cref{proCnbe_object}
\[a = \ang{a_{\bdi}}_{\bdi \in \ufs{n_1 n_2 \cdots n_q}} \in \proAnbe.\]
Using the action unity axiom \cref{pseudoalg_action_unity} for $\A$, \cref{i_connected,sgotoprod_def,Opn_object,zbsgad_comp_obj,zbsgad_comp_mor}, the $\bdi$-th coordinate of $\zbsg \zbsgad a$ is given by
\begin{equation}\label{unis_cod_obj}
(\zbsg \zbsgad a)_{\bdi} 
= (\zbsgad a)_{\{\bdi\}} 
= \gaA_1\big(\ep_1 \sscs a_{\bdi}) 
= \gaA_1\big(\opu \sscs a_{\bdi}) 
= a_{\bdi}.
\end{equation}
Thus, $\zbsg \zbsgad$ is equal to the identity functor on $\proAnbe$.  This finishes the definition of the identity natural transformation $\unis \cn 1 \to \zbsg\zbsgad$.
\end{definition}

\subsection*{Counit}
Next, we define the counit of the adjunction $(\zbsgad,\zbsg)$.

\begin{definition}[Counit]\label{def:zbsgcounit}
For the pointed functors $\zbsgad$ \cref{zbsgad} and $\zbsg$ \cref{sgosgtoprod}, we define $\cous \cn \zbsgad\zbsg \to 1$ as the pointed natural isomorphism
\begin{equation}\label{counit_zbsg}
\begin{tikzpicture}[vcenter]
\def\h{2} \def\v{.7} \def\s{15}
\draw[0cell]
(0,0) node (a1) {\phantom{A}}
(a1)++(-.1,0) node (a1') {\proAnbe}
(a1)++(\h,\v) node (a2) {\Asgsmaangordnbe}
(a1)++(\h,-\v) node (a3) {\Asgsmaangordnbe}
(a1)++(\h,0) node (a) {}
;
\draw[1cell=.9]
(a1) edge[bend left=\s, shorten >=0ex] node {\zbsgad} (a2)
(a2) edge[bend left=\s, transform canvas={xshift=-1.5em}, shorten >=-.3ex, shorten <=-.3ex] node {1} (a3)
(a3) edge[bend left=\s] node {\zbsg} (a1)
;
\draw[2cell]
node[between=a1 and a at .5, shift={(0,-.2*\v)}, rotate=0, 2label={above,\!\cous}] {\Rightarrow}
;
\end{tikzpicture}
\end{equation}
with, for each strong $(\sma\angordnbe)$-system $(a,\gl)$ in $\A$ \cref{nsys}, $(a,\gl)$-component isomorphism
\begin{equation}\label{cous_a}
\zbsgad\zbsg(a,\gl) \fto[\iso]{\cous_{(a,\gl)}} (a,\gl)
\end{equation}
in $\Asgsmaangordnbe$ defined as follows.  For a subset $s \subseteq \txprod_{j \in \ufs{q}}\, \ufs{n}_j$, the $s$-component isomorphism \cref{theta_s} of $\cous_{(a,\gl)}$ is defined by the commutative diagram \cref{cous_a_s}.
\begin{equation}\label{cous_a_s}
\begin{tikzpicture}[vcenter]
\def\v{-1}
\draw[0cell]
(0,0) node (a11) {\big(\zbsgad\zbsg(a,\gl)\big)_s}
(a11)++(4,0) node (a12) {(a,\gl)_s}
(a11)++(0,\v) node (a21) {\gaA_{\cards}\big(\ep_{\cards} \sscs \ang{a_{\{\bdi\}}}_{\bdi \in s}\big)}
(a12)++(0,\v) node (a22) {a_s}
;
\draw[1cell=.9]
(a11) edge node {\cous_{(a,\gl),\, s}} (a12)
(a21) edge node {\gl_{\ep_{\cards} ;\, s, \ang{\{\bdi\}}_{\bdi \in s}}} node[swap] {\iso} (a22)
(a11) edge[equal,shorten <=-.5ex,shorten >=-.5ex] (a21)
(a12) edge[equal] (a22)
;
\end{tikzpicture}
\end{equation}
\begin{itemize}
\item Using \cref{sgotoprod_def,zbsgad_comp_obj}, the left vertical equality in \cref{cous_a_s} is given by
\begin{equation}\label{cous_dom_obj}
\begin{split}
\big(\zbsgad\zbsg(a,\gl)\big)_s 
&= \gaA_{\cards}\big(\ep_{\cards} \sscs \ang{(\zbsg(a,\gl))_{\bdi}}_{\bdi \in s}\big)\\
&= \gaA_{\cards}\big(\ep_{\cards} \sscs \ang{a_{\{\bdi\}}}_{\bdi \in s}\big).
\end{split}
\end{equation}
The computation \cref{cous_dom_obj} is also valid for morphisms in $\Asgsmaangordnbe$ \pcref{def:nsys_morphism}, using \cref{zbsgad_comp_mor} instead of \cref{zbsgad_comp_obj}.
\item $\gl_{\ep_{\cards} ;\, s, \ang{\{\bdi\}}_{\bdi \in s}}$ is the gluing isomorphism \cref{gl-morphism} of $(a,\gl)$ at $(\ep_{\cards}; s, \ang{\{\bdi\}}_{\bdi \in s})$, using the partition 
\[s = \coprod_{\bdi \in s}\, \{\bdi\}\]
of $s$ into one-element subsets ordered lexicographically.
\end{itemize}
\cref{cous_welldef} proves that, as $s$ varies among the subsets of $\txprod_{j \in \ufs{q}}\, \ufs{n}_j$, $\cous_{(a,\gl)}$ is an isomorphism of strong $(\sma\angordnbe)$-systems.  

\parhead{Pointed naturality}.  As $(a,\gl)$ varies in the category $\Asgsmaangordnbe$, $\cous$ is a pointed natural transformation for the following reasons.
\begin{itemize}
\item Using \cref{sgotoprod_def,zbsgad_comp_mor}, the naturality of $\cous$ with respect to morphisms of strong $(\sma\angordnbe)$-systems follows from the unity axiom \cref{nsys_mor_unity} and the compatibility axiom \cref{nsys_mor_compat} for such morphisms.  
\item The pointedness of $\cous$ means that, for the base $(\sma\angordnbe)$-system $(\zero,1_\zero)$, the component $\cous_{(\zero,1_\zero)}$ is given by $1_\zero$ for each subset $s \subseteq \txprod_{j \in \ufs{q}}\, \ufs{n}_j$.  This holds by \cref{cous_a_s}.
\end{itemize}
This finishes the definition of the pointed natural isomorphism $\cous \cn \zbsgad\zbsg \fiso 1$.
\end{definition}

\cref{cous_welldef}, which is used in \cref{def:zbsgcounit}, proves that the components of the counit $\cous$ are well defined.

\begin{lemma}\label{cous_welldef}
For each strong $(\sma\angordnbe)$-system $(a,\gl)$ in $\A$, the collection in \cref{cous_a}
\[\zbsgad\zbsg(a,\gl) \fto{\cous_{(a,\gl)}} (a,\gl)\]
is an isomorphism of strong $(\sma\angordnbe)$-systems.
\end{lemma}

\begin{proof}
Each component of $\cous_{(a,\gl)}$ \cref{cous_a_s} is an isomorphism.  We verify that $\cous_{(a,\gl)}$ satisfies the two axioms in \cref{def:nsys_morphism} for a morphism of $(\sma\angordnbe)$-systems.
\begin{description}
\item[Unity]  
The unity axiom \cref{nsys_mor_unity} holds by \cref{i_connected} and the unity property \cref{nsys_unity_zero} of $(a,\gl)$:
\[\cous_{(a,\gl),\emptyset} 
= \gl_{\ep_0 \scsp \emptyset, \ang{}} 
= \gl_{* \scsp \emptyset, \ang{}} = 1_\zero.\]
\item[Compatibility]  
The compatibility diagram \cref{nsys_mor_compat} for $\cous_{(a,\gl)}$ is the boundary diagram in \cref{cous_compat_diagram}, using the notation $\ang{\Cdots}_\ell = \ang{\Cdots}_{\ell \in \ufs{r}}$, $\phiA = \phiA_{(r;\, \card{s_1}, \ldots, \card{s_r})}$ \cref{phiA_component}, $\sig = \sig_{s, \ang{s_\ell}_\ell}$ \cref{sigma_angsi}, $\alp_x = \alp_{x;\, s, \ang{s_\ell}_\ell}$ \cref{al_Op},
\[\begin{split}
\gl^1 &= \gl_{\ga(x;\, \ang{\ep_{\cardsl}}_\ell) \sscs \, s, \ang{\ang{\{\bdi\}}_{\bdi \in s_\ell}}_{\ell}}, \andspace\\
\gl^2 &= \gl_{\ga(x;\, \ang{\ep_{\cardsl}}_\ell) \sig \sscs \, s, \ang{\{\bdi\}}_{\bdi \in s}}.
\end{split}\]
\begin{equation}\label{cous_compat_diagram}
\begin{tikzpicture}[vcenter]
\def\u{-1} \def\v{-1.5} \def\h{1} 
\draw[0cell=.8]
(0,0) node (a11) {\gaA_r\big(x\sscs \bang{\gaA_{\cardsl}\big(\ep_{\cardsl}\sscs \ang{a_{\{\bdi\}}}_{\bdi \in s_\ell} \big)}_{\ell} \big)}
(a11)++(0,\u) node (a21) {\gaA_r\big(x\sscs \bang{(\zbsgad \zbsg a)_{s_\ell}}_{\ell}\big)}
(a21)++(\h,\v) node (a31) {\gaA_r\big(x \sscs \ang{a_{s_\ell}}_\ell\big)}
(a11)++(\h/2,-\v) node (b1) {\gaA_{\cards}\big(\ga(x\sscs \ang{\ep_{\cardsl}}_\ell) \sscs \ang{\ang{a_{\{\bdi\}}}_{\bdi \in s_\ell}}_\ell \big)}
(b1)++(5,0) node (a12) {\gaA_{\cards}\big(\ga(x\sscs \ang{\ep_{\cardsl}}_\ell) \sig \sscs \ang{a_{\{\bdi\}}}_{\bdi \in s} \big)}
(b1)++(4,0) node (a12') {\phantom{X}}
(a12)++(\h/2,\v) node (a22) {\gaA_{\cards}\big(\ep_{\cards} \sscs \ang{a_{\{\bdi\}}}_{\bdi \in s} \big)}
(a22)++(0,\u) node (a) {(\zbsgad \zbsg a)_s}
(a)++(-\h,\v) node (a32) {a_s}
;
\draw[1cell=.8]
(a11) edge[equal] (a21)
(a21) edge node[swap,pos=.4] {\gaA_r(1; \ang{\gl_{\ep_{\cardsl};\, s_\ell, \ang{\{\bdi\}}_{\bdi \in s_\ell}}}_\ell)} (a31)
(a31) edge node[pos=.45] {\gl_{x;\, s, \ang{s_\ell}_\ell}} (a32)
(a11) edge node[pos=.3] {\phiA} (b1)
(b1) edge[equal] node {\mathbf{eq}} (a12)
(a12) edge node[pos=.8] {\gaA_{\cards}(\alp_x; 1^{\cards})} (a22)
(a22) edge[equal] (a)
(a) edge node[pos=.4] {\gl_{\ep_{\cards};\, s, \ang{\{\bdi\}}_{\bdi \in s}}} (a32)
(b1) edge[bend left=10] node[pos=.4] {\gl^1} (a32)
(a12') edge[bend right=10] node[swap,pos=.3] {\gl^2} (a32)
;
\end{tikzpicture}
\end{equation}
\begin{itemize}
\item Using \cref{zbsgad_gl,cous_dom_obj}, the top composite in \cref{cous_compat_diagram} from $\gaA_r\big(x\sscs \ang{(\zbsgad \zbsg a)_{s_\ell}}_{\ell}\big)$ to $(\zbsgad \zbsg a)_s$ is the gluing isomorphism of the strong $(\sma\angordnbe)$-system $\zbsgad\zbsg(a,\gl)$ at $(x; s, \ang{s_\ell}_\ell)$.  
\item From left to right, the three regions in \cref{cous_compat_diagram} commute by, respectively, the associativity axiom \cref{nsys_associativity}, the equivariance axiom \cref{nsys_equivariance}, and the naturality axiom \cref{nsys_naturality} for the strong $(\sma\angordnbe)$-system $(a,\gl)$.
\end{itemize}
This proves the compatibility axiom for $\cous_{(a,\gl)}$.\qedhere
\end{description}
\end{proof}

\subsection*{Proof of Equivalence}

We now prove that the quadruple $(\zbsgad,\zbsg, \unis,\cous)$ is an adjoint equivalence.

\begin{theorem}\label{thm:zbsg_eq}
For a $\Uinf$-operad $(\Op,\ga,\opu)$ \pcref{as:OpA'}, an $\Op$-pseudoalgebra $(\A,\gaA,\phiA)$ \pcref{def:pseudoalgebra}, and an object $\angordnbe = \ang{\ordn_j^{\be_j}}_{j \in \ufsq} \in \GG \setminus \{\vstar,\ang{}\}$, there is an adjoint equivalence of categories
\begin{equation}\label{zbsg_zbsgad_thm}
\begin{tikzpicture}[vcenter]
\def\s{22}
\draw[0cell]
(0,0) node (a1) {\phantom{X}}
(a1)++(1.8,0) node (a2) {\phantom{X}}
(a1)++(-.17,0) node (a1') {\proAnbe}
(a2)++(.7,-.04) node (a2') {\Asgsmaangordnbe}
node[between=a1 and a2 at .5] {\sim}
;
\draw[1cell=.9]
(a1) edge[bend left=\s] node {\zbsgad} (a2)
(a2) edge[bend left=\s] node {\zbsg} (a1)
;
\end{tikzpicture}
\end{equation}
given by the following data.
\begin{itemize}
\item $\proAnbe$ is the $\angordnbe$-twisted product \pcref{def:proCnbe}.
\item $\Asgsmaangordnbe$ is the category of strong $(\sma\angordnbe)$-systems defined in \cref{smash_GGobjects,def:nsys,def:nsys_morphism,def:nsys_gcat}.
\item The left adjoint is the pointed functor $\zbsgad$ \cref{zbsgad}.
\item The right adjoint is the pointed functor $\zbsg$ \cref{sgosgtoprod}.
\item The unit is the identity natural transformation \cref{unit_zbsg}
\[1_{\proAnbe} \fto{\unis = 1_1} \zbsg\zbsgad.\]
\item The counit is the pointed natural isomorphism \cref{counit_zbsg}
\[\zbsgad \zbsg \fto[\iso]{\cous} 1_{\Asgsmaangordnbe}.\]
\end{itemize}
\end{theorem}

\begin{proof}
The unit $\unis$ is the identity natural transformation, and the counit $\cous$ is a natural isomorphism.  We verify the two triangle identities for an adjunction \pcref{def:adjunction}.
\begin{description}
\item[Left triangle identity] 
This triangle identity states that, for each object \cref{proCnbe_object}
\[a = \ang{a_{\bdi}}_{\bdi \in \ufs{n_1 n_2 \cdots n_q}} \in \proAnbe,\] 
the following composite in $\Asgsmaangordnbe$ is the identity morphism.  
\begin{equation}\label{zbsg_eq_lefttriangle}
\zbsgad a \fto{\zbsgad \unis_a} \zbsgad\zbsg\zbsgad a \fto{\cous_{\zbsgad a}} \zbsgad a
\end{equation}
We observe that each of these two morphisms is the identity.  The first morphism $\zbsgad \unis_a$ in \cref{zbsg_eq_lefttriangle} is the identity by the functoriality of $\zbsgad$, since $\unis_a = 1_a$.

For the second morphism $\cous_{\zbsgad a}$ in \cref{zbsg_eq_lefttriangle}, we consider a subset $s \subseteq \txprod_{j \in \ufs{q}}\, \ufs{n}_j$.  By \cref{sgotoprod_def,zbsgad_comp_obj,zbsgad_gl,cous_dom_obj,cous_a_s}, the $s$-component morphism $\cous_{\zbsgad a, s}$ is the following composite in $\A$.
\begin{equation}\label{zbsg_eq_l}
\begin{tikzpicture}[vcenter]
\def\v{-1.4}
\draw[0cell=.85]
(0,0) node (a1) {(\zbsgad\zbsg\zbsgad a)_s}
(a1)++(0,\v) node (a2) {\gaA_{\cards}\big(\ep_{\cards} \sscs \ang{(\zbsg\zbsgad a)_{\bdi}}_{\bdi \in s} \big)}
(a2)++(0,\v) node (a3) {\gaA_{\cards}\big(\ep_{\cards} \sscs \ang{(\zbsgad a)_{\{\bdi\}}}_{\bdi \in s} \big)}
(a3)++(0,\v) node (a4) {\gaA_{\cards}\big(\ep_{\cards} \sscs \bang{\gaA_1(\ep_1 \sscs a_{\bdi})}_{\bdi \in s}\big)}
(a1)++(5.2,0) node (b1) {(\zbsgad a)_s}
(b1)++(0,\v) node (b2) {\gaA_{\cards}\big(\ep_{\cards} \sscs \ang{a_{\bdi}}_{\bdi \in s}\big)}
(b2)++(0,\v) node (b3) {\phantom{\gaA_{\cards}}}
(b3)++(-.6,0) node (b3') {\gaA_{\cards}\big(\ga(\ep_{\cards} \sscs \ep_1^{\cards}) \sig_{s,\ang{\{\bdi\}}_{\bdi \in s}} \sscs \ang{a_{\bdi}}_{\bdi \in s} \big)}
(b3)++(0,\v) node (b4) {\gaA_{\cards}\big(\ga(\ep_{\cards} \sscs \ep_1^{\cards}) \sscs \ang{a_{\bdi}}_{\bdi \in s} \big)}
;
\draw[1cell=.85]
(a1) edge node {\cous_{\zbsgad a, s}} (b1)
(a1) edge[equal] (a2)
(a2) edge[equal] (a3)
(a3) edge[equal] (a4)
(a4) edge node {\phiA_{(\cards;\, 1^{\cards})}} (b4)
(b4) edge[equal] node[swap] {\mathbf{eq}} (b3)
(b3) edge node[swap] {\gaA_{\cards}(\alp_{\ep_{\cards};\, s, \ang{\{\bdi\}}_{\bdi \in s}} \sscs 1^{\cards})} (b2)
(b2) edge[equal] (b1)
;
\end{tikzpicture}
\end{equation}
Each of the two constituent arrows in the diagram \cref{zbsg_eq_l} is an identity morphism for the following reasons.
\begin{itemize}
\item By \cref{as:OpA'}, $\ep_1$ is the operadic unit $\opu \in \Op(1)$.  The associativity constraint $\phiA_{(\cards;\, 1^{\cards})}$ is an identity morphism by the top half of the unity axiom \cref{pseudoalg_unity} for $\A$.
\item The permutation \cref{sigma_angsi}
\[s \fto[\iso]{\sig_{s,\ang{\{\bdi\}}_{\bdi \in s}}} \txcoprod_{\bdi \in s}\, \{\bdi\}\]
is the identity.  The isomorphism \cref{al_Op}
\[\ga(\ep_{\cards} \sscs \ep_1^{\cards}) \sig_{s,\ang{\{\bdi\}}_{\bdi \in s}} 
= \ga(\ep_{\cards} \sscs \opu^{\cards}) \id_{\cards} 
= \ep_{\cards} \fto{\alp_{\ep_{\cards};\, s, \ang{\{\bdi\}}_{\bdi \in s}}} \ep_{\cards}\]
is the identity morphism on $\ep_{\cards}$ in the translation category $\Op(\cards)$.  Thus, the functoriality of $\gaA_{\cards}$ implies that the arrow $\gaA_{\cards}(\alp_{\ep_{\cards};\, s, \ang{\{\bdi\}}_{\bdi \in s}} \sscs 1^{\cards})$ is an identity morphism.
\end{itemize}
Since $s \subseteq \txprod_{j \in \ufs{q}}\, \ufs{n}_j$ is arbitrary, this proves that the arrow $\cous_{\zbsgad a}$ in \cref{zbsg_eq_lefttriangle} is an identity morphism, proving the left triangle identity.
\item[Right triangle identity]  
This triangle identity states that, for each object $(a,\gl) \in \Asgsmaangordnbe$, the following composite in $\proAnbe$ is the identity morphism.
\begin{equation}\label{zbsg_eq_righttriangle}
\zbsg(a,\gl) \fto{\unis_{\zbsg(a,\gl)}} \zbsg\zbsgad\zbsg(a,\gl) \fto{\zbsg \cous_{(a,\gl)}} \zbsg(a,\gl)
\end{equation}
The first morphism $\unis_{\zbsg(a,\gl)}$ is the identity by definition \cref{unit_zbsg}.

Using the unity axiom \cref{nsys_unity_one} for $(a,\gl)$, \cref{i_connected,sgotoprod_def,cous_a_s}, the following computation proves that the morphism $\zbsg \cous_{(a,\gl)}$ in \cref{zbsg_eq_righttriangle} has $\bdi$-th coordinate given by the identity morphism for each index $\bdi \in \txprod_{j \in \ufs{q}}\, \ufs{n}_j$.
\[\big(\zbsg \cous_{(a,\gl)}\big)_{\bdi} 
= \cous_{(a,\gl), \{\bdi\}} 
= \gl_{\ep_1;\, \{\bdi\}, \{\bdi\}} 
= \gl_{\opu;\, \{\bdi\}, \{\bdi\}} = 1_{a_{\{\bdi\}}}.\]
This proves that $\zbsg \cous_{(a,\gl)}$ is the identity morphism, proving the right triangle identity.\qedhere
\end{description}
\end{proof}

\begin{explanation}[Systems and Adjunction]\label{expl:sys_adjunction}
In \cref{cous_welldef,thm:zbsg_eq}, we need $(a,\gl)$ to be a \emph{strong} $(\sma\angordnbe)$-system because we want an adjoint equivalence, not just an adjunction.  The components of the counit $\cous$ \cref{cous_a_s} are gluing morphisms of $(a,\gl)$.  Thus, in order for $\cous$ to be a natural isomorphism, we need the gluing morphisms of $(a,\gl)$ to be isomorphisms.  The next paragraph discusses what happens if the strong requirement is dropped.

The unit and counit in \cref{def:zbsgunit,def:zbsgcounit} are still defined even if the category $\Asgsmaangordnbe$ is replaced by the category $\Asmaangordnbe$ of all $(\sma\angordnbe)$-systems \pcref{def:nsys,def:nsys_morphism}.
\begin{description}
\item[Unit] There is an identity natural transformation
\begin{equation}\label{unit_zb}

\end{equation}
of the identity functor on $\proAnbe$.  The functor $\zb$ is the pointed $G$-functor defined in \cref{sgotoprod}.  The pointed functor $\zbad$ is the composite
\begin{equation}\label{zbad}
%
\end{equation}
of the pointed functor $\zbsgad$ \cref{zbsgad} and the full subcategory inclusion $\iota$ into $\Asmaangordnbe$.  The computation \cref{unis_cod_obj} shows that this variant of the unit $\unis = 1_1$ is also well defined.
\item[Counit] There is a pointed natural transformation
\begin{equation}\label{counit_zb}
%
\end{equation}
whose components are defined as gluing morphisms of $(a,\gl)$ as in \cref{cous_a_s}.  This natural transformation is well defined by the computation \cref{cous_dom_obj} and the proof of \cref{cous_welldef}, applied to $(\sma\angordnbe)$-systems.  This variant of the counit $\cous$ is a pointed natural transformation, but not a natural isomorphism.
\end{description}
The proof of the triangle identities in \cref{thm:zbsg_eq} is still valid in this setting.  Thus, there is an adjunction
\begin{equation}\label{nbe_sys_adj}
%
\end{equation}
whose unit and counit are the natural transformations in \cref{unit_zb,counit_zb}.  \cref{expl:zbad_pseudo} discusses the pseudo $G$-equivariance of the left adjoint $\zbad$.
\end{explanation}

\section{Pseudo Equivariance of Left Adjoint}
\label{sec:pseudoequiv}

In the adjoint equivalence 
\begin{equation}\label{zbsg_adeq_seci}
%
\end{equation}
in \cref{thm:zbsg_eq}, the right adjoint $\zbsg$ is a pointed $G$-functor by \cref{Pist_to_prod} \eqref{Pist_to_prod_i}.  The main result of this section, \cref{thm:zbsgad_pseudo}, proves that the left adjoint $\zbsgad$ is \emph{pseudo $G$-equivariant}.  This means that the lack of $G$-equivariance of $\zbsgad$ is controlled by a family of natural isomorphisms, called pseudo $G$-equivariant constraints, that satisfy two coherence axioms of their own.

\secoutline
\begin{itemize}
\item \cref{def:pseudoG} defines pseudo $G$-equivariant functors.
\item \cref{def:zbsgad_pseudoG} defines the pseudo $G$-equivariant constraints for the functor $\zbsgad$.  The definition of its components \cref{pse_gas} shows precisely why $\zbsgad$ is not generally $G$-equivariant.  The corresponding components of $\zbsgad g$ and $g\zbsgad$ differ by an isomorphism that is generally not the identity.
\item \cref{ex:al_g_s} illustrates the nontriviality of the pseudo $G$-equivariant constraints for $\zbsgad$ when $\Op$ is either the Barratt-Eccles $\Gcat$-operad $\BE$ or the $G$-Barratt-Eccles operad $\GBE$. 
\item \cref{pse_g_a_welldef} proves that the components of the pseudo $G$-equivariant constraints for $\zbsgad$ are well-defined morphisms of strong $(\sma\angordnbe)$-systems.
\item \cref{thm:zbsgad_pseudo} proves that $\zbsgad$, equipped with the pseudo $G$-equivariant constraints in \cref{def:zbsgad_pseudoG}, is a pseudo $G$-equivariant functor.
\item \cref{expl:zbad_pseudo} discusses the pseudo $G$-equivariance of the pointed functor \cref{zbad}
\[\proAnbe \fto{\zbad = \iota \zbsgad} \Asmaangordnbe.\]
\end{itemize}

\subsection*{Pseudo $G$-Equivariance}
The notion of a pseudo $G$-equivariant functor is due to Merling \cite[Def.\ 3.1]{merling}.  We use the definition of a $G$-category in \cref{expl:GCat}.

\begin{definition}[Pseudo $G$-Equivariance]\label{def:pseudoG}
Suppose $\C$ and $\D$ are $G$-categories for an arbitrary group $G$.  A \emph{pseudo $G$-equivariant functor}\index{pseudo G-equivariant@pseudo $G$-equivariant!functor}\index{G-functor@$G$-functor!pseudo equivariant} from $\C$ to $\D$ is a pair
\[\C \fto{(\func,\pse)} \D\]
consisting of the following data.
\begin{description}
\item[Functor] $\func \cn \C \to \D$ is a functor.
\item[Constraints] For each element $g \in G$, $\pse^g$ is a natural isomorphism, called a \index{pseudo G-equivariant@pseudo $G$-equivariant!constraint}\emph{pseudo $G$-equivariant constraint}, as follows.
\begin{equation}\label{pseg}

\end{equation}
\end{description}
The pair $(\func,\pse)$ is required to satisfy the axioms \cref{pseudoG_unity,pseudoG_mult}.
\begin{description}
\item[Unity] For the identity element $e \in G$, there is an equality
\begin{equation}\label{pseudoG_unity}
\pse^e = 1_{\func}, 
\end{equation}
where $1_{\func}$ is the identity natural transformation of $\func$.
\item[Multiplicativity] For elements $g,h \in G$, the following equality of natural isomorphisms holds.
\begin{equation}\label{pseudoG_mult}
%
\end{equation}
\end{description}
This finishes the definition of a pseudo $G$-equivariant functor.
\end{definition}

\begin{explanation}\label{expl:pseudoG_mult}
A pseudo $G$-equivariant functor is called a \emph{pseudo equivariant functor} in \cite[Def.\ 3.1]{merling}, where $\pse^g$ is denoted by $\theta_g$.  
For each object $x \in \C$, the $x$-component of $\pse^g$ is an isomorphism
\[\func(gx) \fto[\iso]{\pse^g_x} g\func(x) \inspace \D.\]
The multiplicativity axiom \cref{pseudoG_mult} means that, for each object $x \in \C$, the diagram
\begin{equation}\label{pseudoG_mult_c}
\begin{tikzpicture}[vcenter]
\def\h{2.7}
\draw[0cell]
(0,0) node (a1) {\func(hgx)}
(a1)++(\h,0) node (a2) {h\func(gx)}
(a2)++(\h,0) node (a3) {hg\func(x)}
;
\draw[1cell=.9]
(a1) edge node {\pse^h_{gx}} (a2)
(a2) edge node {h\pse^g_x} (a3)
;
\draw[1cell=.9]
(a1) [rounded corners=2pt] |- ($(a2)+(-1,.7)$) -- node {\pse^{hg}_x} ($(a2)+(1,.7)$) -| (a3)
 ;
\end{tikzpicture}
\end{equation}
in $\D$ commutes.
\end{explanation}

\begin{example}[$G$-Functors]\label{ex:pseudoG}
Pseudo $G$-equivariant functors generalize $G$-functors \cref{Gfunctor} in the following sense.  A $G$-functor $F \cn \C \to \D$ becomes a pseudo $G$-equivariant functor $(F,1)$ with each pseudo $G$-equivariant constraint $\pse^g$ given by the identity natural transformation of $Fg = gF$.  Conversely, if $(\func,1)$ is a pseudo $G$-equivariant functor, then $\func$ is a $G$-functor.
\end{example}

\subsection*{Constraints for $\zbsgad$}
Next, we define pseudo $G$-equivariant constraints for the pointed functor \pcref{def:zbsgad}
\[\proAnbe \fto{\zbsgad} \Asgsmaangordnbe,\]
which is part of an adjoint equivalence of categories \pcref{thm:zbsg_eq}.
\cref{thm:zbsgad_pseudo} proves that $\zbsgad$, equipped with the natural isomorphisms in \cref{def:zbsgad_pseudoG}, is a pseudo $G$-equivariant functor.  Recall that $\Asgsmaangordnbe$ and $\proAnbe$ \pcref{def:proCnbe,def:nsys_gcat} are $G$-categories.

\begin{definition}\label{def:zbsgad_pseudoG}
For a $\Uinf$-operad $(\Op,\ga,\opu)$ \pcref{as:OpA'}, an $\Op$-pseudoalgebra $(\A,\gaA,\phiA)$ \pcref{def:pseudoalgebra}, an object $\angordnbe = \ang{\ordn_j^{\be_j}}_{j \in \ufsq} \in \GG \setminus \{\vstar,\ang{}\}$, and an element $g \in G$, we define a natural isomorphism
\begin{equation}\label{zbsgad_constraint}
\begin{tikzpicture}[vcenter]
\def\v{-1.5}
\draw[0cell=.9]
(0,0) node (a1) {\proAnbe}
(a1)++(2.8,0) node (b1) {\proAnbe}
(a1)++(0,\v) node (a2) {\Asgsmaangordnbe}
(b1)++(0,\v) node (b2) {\Asgsmaangordnbe}
;
\draw[1cell=.9]
(a1) edge node {g} (b1)
(a2) edge node[swap] {g} (b2)
(a1) edge node[swap] {\zbsgad} (a2)
(b1) edge node {\zbsgad} (b2)
;
\draw[2cell]
node[between=a1 and b2 at .55, rotate=-135, 2labelalt={below,\zbpg}] {\Rightarrow}
;
\end{tikzpicture}
\end{equation}
as follows.  For each subset $s \subseteq \txprod_{j \in \ufs{q}}\, \ufs{n}_j$, we first define a permutation $\si^{g,s} \in \Si_{\cards}$ and an isomorphism $\al^{g,s} \in \Op(\cards)$ as follows.
\begin{itemize}
\item We define the permutation 
\begin{equation}\label{si_g_s}
\Big(\ginv s \fto[\iso]{\si^{g,s}} s\Big) \in \Si_{\cards}
\end{equation}
that sends an element $\ginv\bdi \in \ginv s$ to $\bdi \in s$, where each of $s$ and $\ginv s$ is equipped with the lexicographic ordering inherited from $\txprod_{j \in \ufs{q}}\, \ufs{n}_j$, on which $g$ acts diagonally \cref{al-sma-be}.  In other words, $\si^{g,s}$ is the $g$-action on $\txprod_{j \in \ufs{q}}\, \ufs{n}_j$ restricted to the subset $\ginv s$.
\item We define the unique isomorphism 
\begin{equation}\label{al_g_s}
\ginv \ep_{\cards} \si^{g,s} \fto[\iso]{\al^{g,s}} \ep_{\cards}
\end{equation}
in the translation category $\Op(\cards)$ with the indicated domain and codomain, where $\ep_{\cards}$ is the chosen object in \cref{Opn_object}.  Since $\Op$ is a $\Gcat$-operad, the symmetric group action $- \cdot \si^{g,s}$ on $\Op(\cards)$ is $G$-equivariant.  Thus, the object $\ginv \ep_{\cards} \si^{g,s} \in \Op(\cards)$ is unambiguous.
\end{itemize}
\begin{description}
\item[Component isomorphisms] 
For each object \cref{proCnbe_object}
\[a = \ang{a_{\bdi}}_{\bdi \in \ufs{n_1 n_2 \cdots n_q}} \in \proAnbe,\]
the $a$-component of $\zbpg$ is the isomorphism of strong $(\sma\angordnbe)$-systems \pcref{def:nsys_morphism}
\begin{equation}\label{pse_g_a}
(\zbsgad g)a \fto[\iso]{\zbpg_a} (g \zbsgad) a
\end{equation}
with $s$-component isomorphism defined by the following commutative diagram in $\A$.
\begin{equation}\label{pse_gas}
\begin{tikzpicture}[vcenter]
\def\v{-1.2}
\draw[0cell=.85]
(0,0) node (a1) {\big((\zbsgad g)a\big)_s}
(a1)++(0,\v) node (a2) {\gaA_{\cards}\big(\ep_{\cards}\sscs \ang{(ga)_{\bdi}}_{\bdi \in s} \big)}
(a2)++(0,\v) node (a3) {\gaA_{\cards}\big(\ep_{\cards}\sscs \ang{ga_{\ginv\bdi}}_{\bdi \in s}\big)}
(a3)++(0,\v) node (a4) {g\gaA_{\cards}\big(\ginv \ep_{\cards}\sscs \ang{a_{\ginv \bdi}}_{\bdi \in s}\big)}
(a4)++(0,\v) node (a5) {g\gaA_{\cards}\big(\ginv \ep_{\cards} \si^{g,s} \sscs \ang{a_{\bdj}}_{\bdj \in \ginv s}\big)}
(a1)++(4.5,0) node (b1) {\big((g \zbsgad) a\big)_s}
(b1)++(0,\v) node (b2) {g(\zbsgad a)_{\ginv s}}
(b2)++(0,\v) node (b3) {g\gaA_{\card{\ginv s}}\big(\ep_{\card{\ginv s}}\sscs \ang{a_{\bdj}}_{\bdj \in \ginv s}\big)}
(b3)++(0,\v) node (b4) {g\gaA_{\cards}\big(\ep_{\cards}\sscs \ang{a_{\bdj}}_{\bdj \in \ginv s}\big)}
;
\draw[1cell=.8]
(a1) edge node {\zbpg_{a,s}} (b1)
(a1) edge[equal] node[swap] {\mathbf{d}} (a2)
(a2) edge[equal] node[swap] {\mathbf{g}'} (a3)
(a3) edge[equal] node[swap] {\mathbf{eq}'} (a4)
(a4) edge[equal] node[swap] {\mathbf{eq}} (a5)
(b1) edge[equal] node {\mathbf{g}} (b2)
(b2) edge[equal] node {\mathbf{d}} (b3)
(b3) edge[equal] node {\mathbf{c}} (b4)
(a5) [rounded corners=2pt] -| node[pos=.23] {g\gaA_{\cards}(\al^{g,s}; 1^{\cards})} node[swap,pos=.23] {\iso} (b4)
;
\end{tikzpicture}
\end{equation}
In the bottom arrow in \cref{pse_gas}, the isomorphism $\al^{g,s}$ is defined in \cref{al_g_s}.  The seven equalities in \cref{pse_gas} are given as follows.
\begin{itemize}
\item The two equalities labeled $\mathbf{d}$ follow from the definition \cref{zbsgad_comp_obj} of $\zbsgad$.
\item The equality labeled $\mathbf{g}$ follows from the definition \cref{ga_scomp} of the $g$-action on systems.
\item The equality labeled $\mathbf{g}'$ follows from the definition \cref{proCnbe_gaction} of the $g$-action on the $\angordnbe$-twisted product $\proAnbe$.
\item The equality labeled $\mathbf{eq}'$ follows from the $G$-equivariance of the functor $\gaA_{\cards}$ \cref{gaAn}.
\item The equality labeled $\mathbf{eq}$ holds by the action equivariance axiom \cref{pseudoalg_action_sym} for $\A$, applied to the permutation $\si^{g,s} \in \Si_{\cards}$ \cref{si_g_s}.
\item The equality labeled $\mathbf{c}$ holds because $\ginv s$ has the same cardinality as $s$.
\end{itemize}
\cref{pse_g_a_welldef} proves that, as $s$ varies among the subsets of $\txprod_{j \in \ufs{q}}\, \ufs{n}_j$, $\zbps_a$ is an isomorphism of strong $(\sma\angordnbe)$-systems.  
\item[Naturality]
The naturality of $\zbps$ with respect to morphisms in $\proAnbe$ follows from the naturality of $g \gaA_{\cards}(\al^{g,s}; -)$ for each subset $s \subseteq \txprod_{j \in \ufs{q}}\, \ufs{n}_j$, \cref{gtha_s,proCnbe_gaction,zbsgad_comp_mor,pse_gas}.  
\end{description}
This finishes the definition of the natural isomorphism $\zbps$ \cref{zbsgad_constraint}.
\end{definition}

\begin{example}[Nontriviality of Pseudo Equivariant Constraints]\label{ex:al_g_s}
The component $\zbps_{a,s}$ \cref{pse_gas} is defined as $g\gaA_{\cards}(\al^{g,s}; 1^{\cards})$, where $\al^{g,s}$ is the unique isomorphism \cref{al_g_s}
\[\ginv \ep_{\cards} \si^{g,s} \fto[\iso]{\al^{g,s}} \ep_{\cards} \inspace \Op(\cards).\]
Here are two examples that illustrate the nontriviality of $\al^{g,s}$ and $\zbpg_{a,s}$.
\begin{enumerate}
\item For the Barratt-Eccles $\Gcat$-operad $\BE$ \pcref{def:BE}, on which $G$ acts trivially, $\al^{g,s}$ is the unique isomorphism
\[\ginv \ep_{\cards} \si^{g,s} 
= \ep_{\cards} \si^{g,s} \fto[\iso]{\al^{g,s}} \ep_{\cards} \inspace \ESigma_{\cards}.\]
\item For the $G$-Barratt-Eccles operad $\GBE$ \pcref{def:GBE} and each element $h \in G$, the $h$-component of $\al^{g,s}$ is the unique isomorphism
\[(\ginv \ep_{\cards} \si^{g,s})(h) = \ep_{\cards}(gh) \cdot \si^{g,s} 
\fto[\iso]{\al^{g,s}_h} \ep_{\cards}(h) \inspace \ESigma_{\cards}.\]
\end{enumerate}
The nonidentity permutation $\si^{g,s} \cn \ginv s \fiso s$ \cref{si_g_s} changes according to $g \in G$, the subset $s \subseteq \txprod_{j \in \ufs{q}}\, \ufs{n}_j$, and the $G$-action $\be_j$ on $\ordn_j$ for $j \in \ufsq$.  Thus, for each of $\BE$ and $\GBE$, the isomorphism $\al^{g,s}$ is not generally the identity, and neither is $\zbpg_{a,s}$.
\end{example}

\cref{pse_g_a_welldef} is used in \cref{def:zbsgad_pseudoG}.

\begin{lemma}\label{pse_g_a_welldef}
For each object $a \in \proAnbe$, the collection in \cref{pse_g_a}
\[(\zbsgad g)a \fto[\iso]{\zbpg_a} (g \zbsgad) a\]
is an isomorphism of strong $(\sma\angordnbe)$-systems in $\A$.
\end{lemma}

\begin{proof}
For each subset $s \subseteq \txprod_{j \in \ufs{q}}\, \ufs{n}_j$, the $s$-component morphism of $\zbpg_a$ is the isomorphism \cref{pse_gas} 
\[\zbpg_{a,s} = g\gaA_{\cards}(\al^{g,s}; 1^{\cards}) \inspace \A.\]
We verify the axioms \cref{nsys_mor_unity,nsys_mor_compat} in \cref{def:nsys_morphism} for a morphism of $(\sma\angordnbe)$-systems.
\begin{description}
\item[Unity]  
For the empty subset $\emptyset \subseteq \txprod_{j \in \ufs{q}}\, \ufs{n}_j$, by \cref{si_g_s,al_g_s}, there are equalities
\[\si^{g,\emptyset} = \id_0 \in \Si_0 \andspace \al^{g,\emptyset} = 1_* \in \Op(0).\]
By \cref{pse_gas}, the functoriality of $\gaA_0$, and the fact that the identity morphism $1_\zero \in \A$ is $G$-fixed, the $\emptyset$-component of $\zbpg_a$ is given by
\[\zbpg_{a,\emptyset} = g\gaA_0(1_*) = g1_\zero = 1_\zero,\] 
proving the unity axiom \cref{nsys_mor_unity}.
\item[Compatibility]  
Given an object $x \in \Op(r)$ with $r \geq 0$, a subset $s \subseteq \txprod_{j \in \ufs{q}}\, \ufs{n}_j$, and a partition $s = \txcoprod_{\ell \in \ufs{r}}\, s_\ell$, the compatibility diagram \cref{nsys_mor_compat} for $\zbps_a$ is the diagram \cref{zbps_compat} in $\A$, where $\ang{\Cdots}_\ell$ means $\ang{\Cdots}_{\ell \in \ufs{r}}$.  The top and bottom horizontal arrows are the gluing isomorphisms of, respectively, $\zbsgad ga$ and $g\zbsgad a$ at $(x; s, \ang{s_\ell}_\ell)$.  
\begin{equation}\label{zbps_compat}
\begin{tikzpicture}[vcenter]
\def\v{-1.5}
\draw[0cell=.9]
(0,0) node (a11) {\gaA_r\big(x\sscs \ang{(\zbsgad ga)_{s_\ell}}_{\ell} \big)}
(a11)++(3,0) node (a12) {(\zbsgad ga)_s}
(a11)++(0,\v) node (a21) {\gaA_r\big(x\sscs \ang{(g\zbsgad a)_{s_\ell}}_{\ell} \big)}
(a12)++(0,\v) node (a22) {(g\zbsgad a)_s}
;
\draw[1cell=.9]
(a11) edge (a12)
(a21) edge (a22)
(a11) edge[transform canvas={xshift={2em}}] node[swap] {\gaA_r(1;\, \ang{\zbps_{a,s_\ell}}_{\ell})} (a21)
(a12) edge node {\zbps_{a,s}} (a22)
;
\end{tikzpicture}
\end{equation}
Using \cref{ga_gl}, \cref{sigma_angsi,al_Op,zbsgad_gl}, \cref{si_g_s}, \cref{al_g_s}, \cref{pse_gas}, and the notation
\[\begin{split}
y &= \ga\big(x \sscs \ang{\ep_{\cardsl}}_\ell\big) \in \Op(\cards), \\
z &= \ga\big(\ginv x \sscs \ang{\ep_{\cardsl}}_\ell\big) \in \Op(\cards),\text{ and}\\
\alp' &= \alp_{\ginv x;\, \ginv s, \ang{\ginv s_\ell}_\ell} \cn z\sig_{\ginv s, \ang{\ginv s_\ell}_\ell} \fiso \ep_{\cards},
\end{split}\]
the compatibility diagram \cref{zbps_compat} unravels to the boundary diagram in \cref{zbps_compat_ii}.
\setlength{\fboxsep}{1pt}
\setlength{\fboxrule}{0.1pt}
\begin{equation}\label{zbps_compat_ii}

\end{equation}
In the diagram \cref{zbps_compat_ii}, each equality labeled $\mathbf{eq}$ holds by the action equivariance axiom \cref{pseudoalg_action_sym} for $\A$, and $\mathbf{eq}'$ indicates the $G$-equivariance of $\gaA$ \cref{gaAn}.
The objects $A_0$ through $A_5$ in the interior of \cref{zbps_compat_ii} are defined as follows.
\[\begin{split}
A_0 &= g\gaA_r\Big(\ginv x; \bang{\gaA_{\cardsl}\big(\ginv \ep_{\cardsl} ; \ang{a_{\ginv\bdi}}_{\bdi \in s_\ell}\big)}_\ell\Big)\\
A_1 &= g\gaA_r\Big(\ginv x; \bang{\gaA_{\cardsl}\big(\ginv \ep_{\cardsl} \si^{g, s_\ell} ; \ang{a_{\bdj}}_{\bdj \in \ginv s_\ell}\big)}_\ell\Big)\\
A_2 &= g\gaA_{\cards}\Big(\ginv \ga(x; \ang{\ep_{\cardsl}}_\ell) ; \ang{\ang{a_{\ginv \bdi}}_{\bdi \in s_\ell}}_{\ell}\Big)\\
A_3 &= g\gaA_{\cards}\Big(\ginv \ga(x; \ang{\ep_{\cardsl}}_{\ell}) \sig_{s, \ang{s_\ell}_\ell} \si^{g,s} ; \ang{a_{\bdj}}_{\bdj \in \ginv s}\Big)\\
A_4 &= g\gaA_{\cards}\Big(\ga\big(\ginv x; \ang{\ginv \ep_{\cardsl} \si^{g,s_\ell}}_\ell\big) ; \ang{\ang{a_{\bdj}}_{\bdj \in \ginv s_\ell}}_\ell \Big)\\
A_5 &= g\gaA_{\cards}\Big(\ginv \ga(x; \ang{\ep_{\cardsl}}_\ell) (\txcoprod_{\ell \in \ufs{r}}\, \si^{g,s_\ell}) ; \ang{ \ang{a_{\bdj}}_{\bdj \in \ginv s_\ell}}_\ell \Big)
\end{split}\]
Each of the two unlabeled regions in \cref{zbps_compat_ii} commutes because it consists entirely of equalities.  

Regions $\fbox{I}$ through $\fbox{VIII}$ in \cref{zbps_compat_ii} commute for the following reasons.
\begin{enumerate}
\item In region $\fbox{I}$, the arrow $g\phiA \cn A_0 \to A_2$ is the image under the $g$-action of the associativity constraint \cref{phiA}
\begin{equation}\label{regionI}
\begin{tikzpicture}[vcenter]
\draw[0cell=.9]
(0,0) node (a) {\gaA_r\Big(\ginv x; \bang{\gaA_{\cardsl}\big(\ginv \ep_{\cardsl} ; \ang{a_{\ginv\bdi}}_{\bdi \in s_\ell}\big)}_\ell\Big)}
(a)++(0,-1.4) node (b) {\gaA_{\cards}\Big(\ga\big(\ginv x; \ang{\ginv \ep_{\cardsl}}_\ell\big) ; \ang{\ang{a_{\ginv \bdi}}_{\bdi \in s_\ell}}_{\ell}\Big).}
;
\draw[1cell=.9]
(a) edge[shorten <=-1ex, shorten >=-1ex] node {\phiA} (b)
;
\end{tikzpicture}
\end{equation}
The codomain of $g\phiA$ is equal to $A_2$ by the $G$-equivariance of the operadic composition $\ga$ of $\Op$.  Region $\fbox{I}$ commutes by the $G$-equivariance of the associativity constraint $\phiA$.
\item In region $\fbox{II}$, the equality $A_4 = A_5$ labeled $\mathbf{b}$ follows from the following equalities in $\Op(\cards)$.
\begin{equation}\label{afourfive}
\begin{split}
&\ga\big(\ginv x; \ang{\ginv \ep_{\cardsl} \si^{g,s_\ell}}_\ell\big)\\
&= \ginv \ga\big(x; \ang{\ep_{\cardsl} \si^{g,s_\ell}}_\ell \big)\\
&= \ginv \ga(x; \ang{\ep_{\cardsl}}_\ell) (\txcoprod_{\ell \in \ufs{r}}\, \si^{g,s_\ell}) 
\end{split}
\end{equation}
The first equality holds by the $G$-equivariance of the operadic composition $\ga$ of $\Op$.  The second equality holds by the bottom equivariance axiom for $\Op$.  Note that the symmetric group action on $\Op$ is $G$-equivariant, so it commutes with the $\ginv$-action.  Region $\fbox{II}$ commutes by the bottom equivariance axiom \cref{pseudoalg_boteq} for $\A$
\item In region $\fbox{III}$, the unlabeled arrow is induced by the isomorphism \cref{al_g_s}
\[\ginv \ep_{\cardsl} \si^{g,s_\ell} \fto[\iso]{\al^{g,s_\ell}} \ep_{\cardsl}\]
for $\ell \in \ufs{r}$:
\begin{equation}\label{regionIII}

\end{equation}
Region $\fbox{III}$ commutes by the $G$-equivariance of $\gaA_r$ \cref{gaAn}.
\item In region $\fbox{IV}$, recalling $z = \ga(\ginv x; \ang{\ep_{\cardsl}}_\ell)$, the unlabeled arrow from $A_4$ to the bottom object is induced by $\al^{g,s_\ell}$ for $\ell \in \ufs{r}$:
\begin{equation}\label{regionIV}
%
\end{equation}
Region $\fbox{IV}$ commutes by the naturality of $\phiA$.
\item In region $\fbox{V}$, the unlabeled arrow from $A_5$ to the bottom object is induced by the unique isomorphism
\begin{equation}\label{alpha_one}
\ginv \ga(x; \ang{\ep_{\cardsl}}_\ell) (\txcoprod_{\ell \in \ufs{r}}\, \si^{g,s_\ell}) 
\fto[\iso]{\al_1} z =\ga(\ginv x; \ang{\ep_{\cardsl}}_\ell)
\end{equation}
in the translation category $\Op(\cards)$:
\begin{equation}\label{afive_to_stuff}
%
\end{equation}
Using \cref{afourfive}, region $\fbox{V}$ is obtained by applying $g \gaA_{\cards}\big(-; \ang{\ang{a_{\bdj}}_{\bdj \in \ginv s_\ell}}_\ell\big)$ to the commutative diagram
\begin{equation}\label{regionV}
%
\end{equation}
in the translation category $\Op(\cards)$.
\item In region $\fbox{VI}$, the unlabeled vertical arrow from $A_3$ is induced by the unique isomorphism
\begin{equation}\label{alpha_two}
\ginv \ga(x; \ang{\ep_{\cardsl}}_{\ell}) \sig_{s, \ang{s_\ell}_\ell} \si^{g,s} \fto[\iso]{\al_2} 
z\sig_{\ginv s, \ang{\ginv s_\ell}_\ell}
\end{equation}
in the translation category $\Op(\cards)$:
\begin{equation}\label{athree_to_stuff}
%
\end{equation}
By \cref{sigma_angsi,si_g_s}, there is a commutative diagram of permutations
\begin{equation}\label{regionVI}
%
\end{equation}
in which each composite sends an element $\ginv \bdi \in \ginv s$ with $\bdi \in s_\ell$ to $\bdi \in s_\ell$.  By uniqueness of morphisms in $\Op(\cards)$ with specified domain and codomain, this implies the morphism equality
\begin{equation}\label{alpha_two_one_si}
\al_2 = \al_1 \sig_{\ginv s, \ang{\ginv s_\ell}_\ell} \inspace \Op(\cards),
\end{equation}
where $\al_1$ and $\al_2$ are the unique isomorphisms in \cref{alpha_one,alpha_two}.  The equality \cref{alpha_two_one_si} and the action equivariance axiom \cref{pseudoalg_action_sym} for $\A$, applied to the permutation $\sig_{\ginv s, \ang{\ginv s_\ell}_\ell}$, imply that the two isomorphisms 
\[g\gaA_{\cards}\big(\al_1; 1^{\cards}\big) \andspace g\gaA_{\cards}\big(\al_2; 1^{\cards}\big)\]
in, respectively, \cref{afive_to_stuff,athree_to_stuff} are equal.  Since the boundary of region $\fbox{VI}$ consists of only these two isomorphisms and equalities, this proves that region $\fbox{VI}$ commutes.
\item In region $\fbox{VII}$, the unlabeled slanted arrow from $A_3$ is induced by the isomorphism \cref{al_Op}
\[\ga\big(x; \ang{\ep_{\cardsl}}_\ell\big) \sig_{s,\ang{s_\ell}_\ell} 
\fto[\iso]{\alp_{x;\, s,\ang{s_\ell}_\ell}} \ep_{\cards}\]
in the translation category $\Op(\cards)$:
\begin{equation}\label{region_vii}
\begin{tikzpicture}[vcenter]
\draw[0cell=.9]
(0,0) node (a) {g\gaA_{\cards}\Big(\ginv \ga(x; \ang{\ep_{\cardsl}}_{\ell}) \sig_{s, \ang{s_\ell}_\ell} \si^{g,s} ; \ang{a_{\bdj}}_{\bdj \in \ginv s}\Big)}
(a)++(-3.8,0) node (a') {A_3}
(a)++(0,-1.4) node (b) {g\gaA_{\cards}\big(\ginv \ep_{\cards} \si^{g,s} \sscs \ang{a_{\bdj}}_{\bdj \in \ginv s}\big).}
;
\draw[1cell=.9]
(a') edge[equal] (a)
(a) edge[transform canvas={xshift={-4.5em}}, shorten <=-.5ex, shorten >=0ex] node {g\gaA_{\cards}\big(\ginv \alp_{x;\, s, \ang{s_\ell}_\ell} \si^{g,s} ; 1^{\cards} \big)} (b)
;
\end{tikzpicture}
\end{equation}
Region $\fbox{VII}$ is obtained by applying $g\gaA_{\cards}(-; \ang{a_{\bdj}}_{\bdj \in \ginv s})$ to the commutative diagram
\begin{equation}\label{regionVII}
\begin{tikzpicture}[vcenter]
\def\v{-1.5}
\draw[0cell=.9]
(0,0) node (a11) {\ginv \ga(x; \ang{\ep_{\cardsl}}_{\ell}) \sig_{s, \ang{s_\ell}_\ell} \si^{g,s}} 
(a11)++(0,\v) node (a21) {z\sig_{\ginv s, \ang{\ginv s_\ell}_\ell}}
(a11)++(5,0) node (a12) {\ginv \ep_{\cards} \si^{g,s}}
(a12)++(0,\v) node (a22) {\phantom{\ep_{\cards}}}
(a22)++(0,-.05) node (a22') {\ep_{\cards}}
;
\draw[1cell=.8]
(a11) edge node {\ginv \alp_{x;\, s, \ang{s_\ell}_\ell} \si^{g,s}} (a12)
(a12) edge node {\al^{g,s}} (a22')
(a11) edge node[swap] {\al_2} (a21)
(a21) edge node {\alp'} (a22)
;
\end{tikzpicture}
\end{equation}
in the translation category $\Op(\cards)$.
\item Both arrows in the boundary of region $\fbox{VIII}$ are induced by the unique isomorphism $\alp_{x;\, s, \ang{s_\ell}_\ell}$ in $\Op(\cards)$.  Region $\fbox{VIII}$ commutes by the $G$-equivariance of $\gaA_{\cards}$ and the action equivariance axiom \cref{pseudoalg_action_sym} for $\A$, applied to the permutation $\si^{g,s}$ \cref{si_g_s}.
\end{enumerate}
This proves that the compatibility diagram \cref{zbps_compat_ii} commutes.\qedhere
\end{description}
\end{proof}

\subsection*{Pseudo $G$-Equivariance of $\zbsgad$}

We now prove that the equivalence $\zbsgad$ \pcref{thm:zbsg_eq} is a pseudo $G$-equivariant functor \pcref{def:pseudoG}.

\begin{theorem}\label{thm:zbsgad_pseudo}
Under the same assumptions as \cref{def:zbsgad_pseudoG}, the pair
\[\proAnbe \fto{(\zbsgad,\zbp)} \Asgsmaangordnbe\]
consisting of
\begin{itemize}
\item the pointed functor $\zbsgad$ \pcref{def:zbsgad} and
\item the natural isomorphisms $\{\zbpg\}_{g \in G}$ \cref{zbsgad_constraint}
\end{itemize}
is a pseudo $G$-equivariant functor.
\end{theorem}

\begin{proof}
We verify the axioms \cref{pseudoG_unity,pseudoG_mult} for a pseudo $G$-equivariant functor.
\begin{description}
\item[Unity]  
For the unity axiom \cref{pseudoG_unity}, we consider the identity element $e \in G$ and an arbitrary subset $s \subseteq \txprod_{j \in \ufs{q}}\, \ufs{n}_j$.  The permutation \cref{si_g_s}
\[\big(e^{-1}s = s \fto[\iso]{\si^{e,s}} s\big) \in \Si_{|s|}\]
is the identity, and the unique isomorphism \cref{al_g_s}
\[e^{-1} \ep_{\cards} \si^{e,s} = \ep_{\cards} \fto[\iso]{\al^{e,s}} \ep_{\cards}\]
is the identity morphism in the translation category $\Op(\cards)$.  By \cref{pse_gas}, for each object $a \in \proAnbe$, the $s$-component isomorphism 
\[\zbpe_{a,s} = e \gaA_{\cards}\big(1_{\ep_{\cards}} ; 1^{\cards}\big)\]
is the identity morphism by the functoriality of $\gaA_{\cards}$ and the fact that the $e$-action on $\A$ is the identity.  Since $a$ and $s$ are arbitrary, this proves that $\zbpe$ is the identity natural transformation of $\zbsgad$.
\item[Multiplicativity]  
Morphisms in the category $\Asgsmaangordnbe$ \pcref{def:nsys_morphism} are determined by their components.  The multiplicativity axiom \cref{pseudoG_mult} for $\zbp$ means that, for elements $g,h \in G$, each object $a \in \proAnbe$ \cref{proCnbe_object}, and each subset $s \subseteq \txprod_{j \in \ufs{q}}\, \ufs{n}_j$, the diagram 
\begin{equation}\label{zbp_mult}
\begin{tikzpicture}[vcenter]
\def\h{3.5} \def\t{15}
\draw[0cell]
(0,0) node (a1) {\big(\zbsgad hga\big)_s}
(a1)++(\h/2,-1) node (a2) {\big(h \zbsgad ga\big)_s}
(a1)++(\h,0) node (a3) {\big(hg\zbsgad a\big)_s}
;
\draw[1cell=.85]
(a1) edge node {\zbphg_{a,s}} (a3)
(a1) [rounded corners=2pt] |- node[swap,pos=.25] {\zbph_{ga,s}} (a2)
;
\draw[1cell]
(a2) [rounded corners=2pt] -| node[swap,pos=.75] {(h\zbpg_a)_s} (a3)
;
\end{tikzpicture}
\end{equation}
in $\A$ commutes.  To verify that \cref{zbp_mult} commutes, we first unravel each of its three arrows using  \cref{si_g_s,al_g_s,pse_gas}.
\begin{enumerate}
\item The top arrow in \cref{zbp_mult} is given by
\begin{equation}\label{zbphg_as}
\zbphg_{a,s} = hg \gaA_{\cards}\big(\al^{hg,s} ; \ang{a_{\bdj}}_{\bdj \in (hg)^{-1} s} \big)
\end{equation}
in which
\[(hg)^{-1} \ep_{\cards} \si^{hg,s} \fto[\iso]{\al^{hg,s}} \ep_{\cards}\]
is the unique isomorphism in $\Op(\cards)$ with the indicated domain and codomain.  The permutation \cref{si_g_s}
\[\ginv\hinv s = (hg)^{-1} s \fto[\iso]{\si^{hg,s}} s\]
is the $(hg)$-action on $\txprod_{j \in \ufs{q}}\, \ufs{n}_j$ restricted to the subset $(hg)^{-1} s$.
\item Using \cref{proCnbe_gaction}, the $G$-equivariance of $\gaA_{\cards}$, and the action equivariance axiom \cref{pseudoalg_action_sym} for $\A$, the lower left arrow in \cref{zbp_mult} is given as follows.
\begin{equation}\label{zbph_gas}
\begin{split}
\zbph_{ga,s} 
&= h\gaA_{\cards}\big(\al^{h,s}; \ang{(ga)_{\bdj}}_{\bdj \in \hinv s}\big)\\
&= h\gaA_{\cards}\big(\al^{h,s}; \ang{ga_{\ginv \bdj}}_{\bdj \in \hinv s}\big)\\
&= hg\gaA_{\cards}\big(\ginv \al^{h,s}; \ang{a_{\ginv \bdj}}_{\bdj \in \hinv s} \big)\\
&= hg\gaA_{\cards}\big(\ginv \al^{h,s} \si^{g,\hinv s}; \ang{a_{\bdj}}_{\bdj \in \ginv\hinv s} \big)\\
\end{split}
\end{equation}
The last equality in \cref{zbph_gas} uses the permutation $\si^{g,\hinv s}$ defined by the following commutative diagram of permutations.
\begin{equation}\label{sig_hinvs}
\begin{tikzpicture}[vcenter]
\def\h{2.5} \def\t{15}
\draw[0cell]
(0,0) node (a1) {(hg)^{-1} s}
(a1)++(-1.8,0) node (a) {\ginv\hinv s}
(a1)++(\h/2,-1) node (a2) {\hinv s}
(a1)++(\h,0) node (a3) {s}
;
\draw[1cell=.85]
(a) edge[equal] (a1)
(a1) edge node {\si^{hg,s}} (a3)
(a1) [rounded corners=2pt] |- node[swap,pos=.25] {\si^{g,\hinv s}} (a2)
;
\draw[1cell=.85]
(a2) [rounded corners=2pt] -| node[swap,pos=.75] {\si^{h,s}} (a3)
;
\end{tikzpicture}
\end{equation}
In other words, $\si^{g,\hinv s}$ is the restriction of the $g$-action on $\txprod_{j \in \ufs{q}}\, \ufs{n}_j$ to $(hg)^{-1}s$.
\item Using \cref{gtha_s}, the lower right arrow in \cref{zbp_mult} is given as follows.
\begin{equation}\label{hzbpg_ahinvs}
\begin{split}
(h\zbpg_a)_s &= h\zbpg_{a,\hinv s}\\
&= hg \gaA_{\card{\hinv s}}\big(\al^{g,\hinv s}; \ang{a_{\bdj}}_{\bdj \in \ginv\hinv s}\big)\\
&= hg \gaA_{\cards}\big(\al^{g,\hinv s}; \ang{a_{\bdj}}_{\bdj \in \ginv\hinv s}\big)
\end{split}
\end{equation}
\end{enumerate}
Each of $\zbphg_{a,s}$ \cref{zbphg_as}, $\zbph_{ga,s}$ \cref{zbph_gas}, and $(h\zbpg_a)_s$ \cref{hzbpg_ahinvs} has the form 
\[hg\gaA_{\cards}\big(-; \ang{a_{\bdj}}_{\bdj \in (hg)^{-1} s}\big)\]
for some morphism $-$ in $\Op(\cards)$.  Thus, the diagram \cref{zbp_mult} commutes because the following diagram in the translation category $\Op(\cards)$ commutes.
\begin{equation}\label{alghinvs}
\begin{tikzpicture}[vcenter]
\def\h{5} \def\v{-1.5} \def\t{15}
\draw[0cell]
(0,0) node (a11) {(hg)^{-1} \ep_{\cards} \si^{hg,s}}
(a11)++(\h,0) node (a12) {\phantom{\ep_{\cards} = \ep_{\card{\hinv s}}}}
(a12)++(0,-.1) node (a12') {\ep_{\cards} = \ep_{\card{\hinv s}}}
(a11)++(0,\v) node (a21) {\ginv\hinv \ep_{\cards} \si^{h,s} \si^{g,\hinv s}}
(a12)++(0,\v) node (a22) {\ginv \ep_{\cards} \si^{g,\hinv s}}
;
\draw[1cell=.85]
(a11) edge node {\al^{hg,s}} (a12)
(a11) edge[equal] (a21)
(a21) edge node {\ginv \al^{h,s} \si^{g,\hinv s}} (a22)
(a22) edge[shorten <=-.5ex, transform canvas={xshift=0ex}] node[swap] {\al^{g,\hinv s}} (a12') 
;
\end{tikzpicture}
\end{equation}
In the previous diagram, the equality on the left uses \cref{sig_hinvs}.  This proves the multiplicativity axiom for $(\zbsgad,\zbp)$.\qedhere
\end{description}
\end{proof}

\begin{explanation}[Pseudo $G$-Equivariance of $\zbad$]\label{expl:zbad_pseudo}
Recall the adjunction \pcref{expl:sys_adjunction} 
\begin{equation}\label{zb_adjun}
\begin{tikzpicture}[vcenter]
\def\s{22}
\draw[0cell]
(0,0) node (a1) {\phantom{X}}
(a1)++(1.8,0) node (a2) {\phantom{X}}
(a1)++(-.15,0) node (a1') {\proAnbe}
(a2)++(.6,-.04) node (a2') {\Asmaangordnbe}
;
\draw[0cell=.8]
node[between=a1 and a2 at .5] {\perp}
;
\draw[1cell=.9]
(a1) edge[bend left=\s] node {\zbad} (a2)
(a2) edge[bend left=\s] node {\zb} (a1)
;
\end{tikzpicture}
\end{equation}
involving the pointed $G$-functor $\zb$ \cref{sgotoprod} and the pointed functor $\zbad = \iota\zbsgad$ \cref{zbad}.  As a consequence of \cref{thm:zbsgad_pseudo}, the pair 
\[\proAnbe \fto{(\zbad, \iota * \zbp)} \Asmaangordnbe\]
is a pseudo $G$-equivariant functor, where $\zbp = \{\zbpg\}_{g \in G}$ is the pseudo $G$-equivariant constraint for $\zbsgad$ \cref{zbsgad_constraint}. 
\begin{description}
\item[Constraints] For each $g \in G$, its pseudo $G$-equivariant constraint $\iota * \zbpg$ is given by whiskering the natural isomorphism $\zbpg$ \cref{zbsgad_constraint} with the inclusion $\iota$, which is a $G$-functor, as follows.
\begin{equation}\label{zbps_diag}
\begin{tikzpicture}[vcenter]
\def\v{-1.5} \def\u{-1.2}
\draw[0cell=.9]
(0,0) node (a1) {\proAnbe}
(a1)++(0,\v) node (a2) {\Asgsmaangordnbe}
(a2)++(0,\u) node (a3) {\Asmaangordnbe}
(a1)++(2.8,0) node (b1) {\proAnbe}
(b1)++(0,\v) node (b2) {\Asgsmaangordnbe}
(b2)++(0,\u) node (b3) {\Asmaangordnbe}
;
\draw[1cell=.9]
(a1) edge node {g} (b1)
(a2) edge node[swap] {g} (b2)
(a3) edge node {g} (b3)
(a1) edge node[swap] {\zbsgad} (a2)
(a2) edge node[swap] {\iota} (a3)
(b1) edge node {\zbsgad} (b2)
(b2) edge node {\iota} (b3)
;
\draw[2cell]
node[between=a1 and b2 at .55, rotate=-135, 2labelalt={below,\zbps}] {\Rightarrow}
;
\end{tikzpicture}
\end{equation}
\item[Unity] The unity axiom \cref{pseudoG_unity} for $\iota * \zbpe$ follows from the unity axiom for $\zbpe$ because the inclusion $\iota$ preserves identity morphisms.
\item[Multiplicativity] The multiplicativity axiom \cref{pseudoG_mult} for $\iota * \zbp$ is obtained from the multiplicativity axiom for $\zbp$ by whiskering with $\iota$.
\end{description}
In summary, the right adjoint $\zb$ is $G$-equivariant (\cref{Pist_to_prod} \eqref{Pist_to_prod_i}), and the left adjoint $\zbad$ is pseudo $G$-equivariant, with pseudo $G$-equivariant constraint given by $\iota * \zbp$.
\end{explanation}

%% file: chap/hgoprod.tex
This chapter proves that there is an adjoint equivalence \pcref{thm:zdsg_eq}
\begin{equation}\label{zdsg_chi}
\begin{tikzpicture}[vcenter]
\def\s{22}
\draw[0cell]
(0,0) node (a1) {\phantom{X}}
(a1)++(1.8,0) node (a2) {\phantom{X}}
(a1)++(-.15,0) node (a1') {\proAnbe}
(a2)++(.4,-.04) node (a2') {\Asgangordnbe}
node[between=a1 and a2 at .5] {\sim}
;
\draw[1cell=.9]
(a1) edge[bend left=\s] node {\zdsgad} (a2)
(a2) edge[bend left=\s] node {\zdsg} (a1)
;
\end{tikzpicture}
\end{equation}
between the twisted product $\proAnbe$ \pcref{def:proCnbe} and strong $H$-theory $\Asgangordnbe$ \pcref{A_ptfunctorGG} for each object $\angordnbe \in \GG \setminus \{\vstar,\ang{}\}$, $\Uinf$-operad $\Op$ \pcref{as:OpA'}, and $\Op$-pseudoalgebra $\A$ \pcref{def:pseudoalgebra}.  While the right adjoint $\zdsg$ \cref{hgosgtoprod} is $G$-equivariant, its adjoint inverse $\zdsgad$ is only \emph{pseudo} $G$-equivariant in general \pcref{thm:zdsgad_pseudo}.  Using \cref{thm:zdsg_eq}, \cref{thm:PistAequivalence} proves that the strong $H$-theory comparison pointed $G$-functor \cref{PistsgAnbe} 
\[(\smast\Sgosg\A)\angordnbe = \Asgsmaangordnbe 
\fto[\sim]{\Pistsg_{\A,\angordnbe}} 
\Asgangordnbe = (\Hgosg\A)\angordnbe\]
is an equivalence of categories.  In other words, Shimakawa strong $H$-theory and strong $H$-theory are nonequivariantly equivalent. 
In \cref{ch:compgen}, the equivalence $\Pistsg$ is upgraded to a categorical weak $G$-equivalence for a $\Catg(\EG,\Op)$-pseudoalgebra of the form $\Ah = \Catg(\EG,\A)$ \pcref{thm:pistweakgeq}.

\summary
\Cref{table.sgosg} summarizes Shimakawa strong $H$-theory $\Sgosg$ and strong $H$-theory $\Hgosg$.
\begin{figure}[H] 
\centering
\resizebox{.9\columnwidth}{!}{%
{\renewcommand{\arraystretch}{1.4}%
{\setlength{\tabcolsep}{1ex}
\begin{tabular}{c|cr|cr}
strong $H$-theories & $\AlgpspsO \fto{\Sgosg} \FGCatg$ & \eqref{def:sgosg} & $\AlgpspsO \fto{\Hgosg} \GGCatg$ & \eqref{Hgo_twofunctor} \\[.3ex] \hline
strong systems & \raisebox{8pt}{\phantom{I}}$\Asgsmaangordnbe = (\smast\Sgosg\A)\angordnbe$\raisebox{8pt}{\phantom{I}} & \eqref{def:nsys_gcat} & $\Asgangordnbe = (\Hgosg\A)\angordnbe$ & \eqref{def:nbeta_gcat}\\
right adjoints & $\Asgsmaangordnbe \fto{\zbsg} \proAnbe$ & \eqref{sgosgtoprod} & $\Asgangordnbe \fto{\zdsg} \proAnbe$ & \eqref{hgosgtoprod}\\
left adjoints & $\proAnbe \fto{\zbsgad} \Asgsmaangordnbe$ & \eqref{def:zbsgad} & $\proAnbe \fto{\zdsgad} \Asgangordnbe$ & \eqref{def:zdsgad}\\
units & $1 \fto{\unis} \zbsg\zbsgad$ & \eqref{def:zbsgunit} & $1 \fto{\unit} \zdsg\zdsgad$ & \eqref{def:zdsgunit} \\
counits & $\zbsgad\zbsg \fto{\cous} 1$ & \eqref{def:zbsgcounit} & $\zdsgad\zdsg \fto{\counit} 1$ & \eqref{def:zdsgcounit} \\
adjoint equivalences & $(\zbsgad,\zbsg,\unis,\cous)$ & \eqref{thm:zbsg_eq} & $(\zdsgad,\zdsg,\unit,\counit)$ & \eqref{thm:zdsg_eq} \\
pseudo $G$-equivariance & $(\zbsgad, \zbp)$ & \eqref{thm:zbsgad_pseudo} & $(\zdsgad, \zdp)$ & \eqref{thm:zdsgad_pseudo} \\
comparison & $\Asgsmaangordnbe \fto[\sim]{\Pistsg_{\A,\angordnbe}} \Asgangordnbe$ & \eqref{thm:PistAequivalence} & -- & \\
\end{tabular}}}}
\caption{Shimakawa strong $H$-theory and strong $H$-theory.}
\label{table.sgosg}
\end{figure}
\cref{expl:zdsgad_zbsgad,expl:zdp_zbp} provide a more detailed comparison between $(\zbsgad, \zbp)$ and $(\zdsgad, \zdp)$.  The following diagram summaries the five main functors related to the comparison $\Pistsg$.
\begin{equation}\label{Pist_main}
%
\end{equation}
The three functors in the inner commutative triangle \cref{Pistsg_prod_diag}---$\Pistsg_{\A,\angordnbe}$, $\zbsg$, and $\zdsg$---are $G$-equivariant \pcref{PistA_gnat,Pist_to_prod}.  The left adjoint $\zbsgad$ of $\zbsg$ and the left adjoint $\zdsgad$ of $\zdsg$ are pseudo $G$-equivariant \pcref{thm:zbsgad_pseudo,thm:zdsgad_pseudo}.

\organization
This chapter consists of the following sections.

\secname{sec:hgo_tprod}
This section constructs the pointed functor $\zdsgad$ from the $\angordnbe$-twisted product $\proAnbe$ to the category $\Asgangordnbe$ of strong $\angordnbe$-systems in $\A$.

\secname{sec:hgo_unit_counit}
This section constructs the unit and the counit,
\[1_{\proAnbe} \fto{\unit = 1_1} \zdsg\zdsgad \andspace 
\zdsgad\zdsg \fto[\iso]{\counit} 1_{\Asgangordnbe},\]
for the functors $(\zdsgad,\zdsg)$.

\secname{sec:hgo_prod_eq}
This section proves that the quadruple $(\zdsgad,\zdsg,\unit,\counit)$ is an adjoint equivalence and that the strong $H$-theory comparison $\Pistsg$ is objectwise a level equivalence of categories.  Replacing the category $\Asgangordnbe$ by the larger category $\Aangordnbe$ of all $\angordnbe$-systems yields an adjunction
\begin{equation}\label{zd_chi}
\begin{tikzpicture}[vcenter]
\def\s{20}
\draw[0cell]
(0,0) node (a1) {\phantom{X}}
(a1)++(1.8,0) node (a2) {\phantom{X}}
(a1)++(-.15,0) node (a1') {\proAnbe}
(a2)++(.3,-.04) node (a2') {\Aangordnbe}
;
\draw[0cell=.8]
node[between=a1 and a2 at .5] {\perp}
;
\draw[1cell=.9]
(a1) edge[bend left=\s] node {\zdad} (a2)
(a2) edge[bend left=\s] node {\zd} (a1)
;
\end{tikzpicture}
\end{equation}
whose counit is not a natural isomorphism in general.

\secname{sec:pseudoequiv_ii}
This section proves that the left adjoint $\zdsgad$ is a pseudo $G$-equivariant functor.  For a nontrivial group $G$, the functor $\zdsgad$ is \emph{not} $G$-equivariant even for the Barratt-Eccles $\Gcat$-operad $\BE$, on which $G$ acts trivially \pcref{ex:ual_g_angs}.  The left adjoint $\zdad$ is also pseudo $G$-equivariant.

\section{Strong $H$-Theory from Twisted Products}
\label{sec:hgo_tprod}

This section constructs a pointed functor 
\[\proAnbe \fto{\zdsgad} \Asgangordnbe\]
from the twisted product to strong $H$-theory for each $\Uinf$-operad $\Op$, $\Op$-pseudoalgebra $\A$, and object $\angordnbe \in \GG \setminus \{\vstar,\ang{}\}$. 
\cref{sec:hgo_unit_counit,sec:hgo_prod_eq} show that the functors $(\zdsgad,\zdsg)$ form a nonequivariant adjoint equivalence, where
\[\Asgangordnbe \fto{\zdsg} \proAnbe\]
is the pointed $G$-functor from strong $H$-theory to the twisted product \cref{hgosgtoprod}.

\secoutline
\begin{itemize}
\item \cref{def:zdsgad} defines the pointed functor $\zdsgad$.
\item \cref{expl:zdsgad_zbsgad} discusses a useful correspondence between $\zbsgad$ \cref{zbsgad} and $\zdsgad$.  This correspondence allows us to reuse some of the proofs for $\zbsgad$ to prove analogous statements for $\zdsgad$ simply by modifying the notation. 
\item \cref{zdsgad_obj_welldef,zdsgad_mor_welldef} prove that $\zdsgad$ is well defined on objects and morphisms.
\end{itemize}

\subsection*{Inverse of $\zdsg$}

We now define an adjoint inverse $\zdsgad$ of the pointed $G$-functor $\zdsg$ \cref{hgosgtoprod}.  

\begin{definition}\label{def:zdsgad}
Given a $\Uinf$-operad $\Op$ \pcref{as:OpA'}, an $\Op$-pseudoalgebra $(\A,\gaA,\phiA)$ \pcref{def:pseudoalgebra}, and an object \cref{GG_objects}
\[\angordnbe = \sordi{n}{\be}{j}_{j \in \ufs{q}} \in \GG \setminus \{\vstar,\ang{}\}\]
of length $q>0$, we define a pointed functor
\begin{equation}\label{zdsgad}
\proAnbe \fto{\zdsgad} \Asgangordnbe
\end{equation}
as follows.  
\begin{description}
\item[Domain]
The $\angordnbe$-twisted product $\proAnbe$ is defined in \cref{def:proCnbe}.  
\item[Codomain]
The category $\Asgangordnbe$ of strong $\angordnbe$-systems is defined in \cref{def:nsystem,def:nsystem_morphism,def:nbeta_gcat}.
\item[Component objects]
Given an object \cref{proCnbe_object}
\[a = \ang{a_{\bdi} \in \A}_{\bdi \sins \ufs{n_1 n_2 \cdots n_q}} \in \proAnbe,\]
the strong $\angordnbe$-system 
\begin{equation}\label{zdsgad_a}
(\zdsgad a, \glu) \in \Asgangordnbe
\end{equation} 
has, for each marker $\angs = \ang{s_j \subseteq \ufs{n}_j}_{j \in \ufs{q}}$ \cref{marker}, $\angs$-component object \cref{a_angs} defined as 
\begin{equation}\label{zdsgad_comp_obj}
(\zdsgad a)_{\angs} 
= \gaA_{\card{\angstimes}}\big( \ep_{\card{\angstimes}} \sscs \ang{a_{\bdi}}_{\bdi \in \angstimes}\big) \in \A.
\end{equation}
The right-hand side of \cref{zdsgad_comp_obj} is given as follows.
\begin{itemize}
\item $\angstimes$ is the lexicographically ordered subset \cref{angstimes}
\[\angstimes = \txprod_{j\in \ufs{q}}\, s_j \subseteq 
\txprod_{j\in \ufs{q}}\, \ufs{n}_j = \ufs{n_1 \cdots n_q}\]
with cardinality $\card{\angstimes}$.
\item The functor
\[\Op(\card{\angstimes}) \times \A^{\card{\angstimes}} \fto{\gaA_{\card{\angstimes}}} \A\]
is the $\card{\angstimes}$-th $\Op$-action $G$-functor of $\A$ \cref{gaAn}.
\item $\ep_{\card{\angstimes}} \in \Op(\card{\angstimes})$ is the object chosen in \cref{Opn_object}.
\item The $\card{\angstimes}$-tuple $\ang{a_{\bdi}}_{\bdi \in \angstimes} \in \A^{\card{\angstimes}}$ uses the lexicographic ordering of $\angstimes$.
\end{itemize}
\item[Gluing]
Given an object $x \in \Op(r)$ with $r \geq 0$, a marker $\ang{s} = \ang{s_j \subseteq \ufs{n}_j}_{j \in \ufs{q}}$, an index $k \in \ufs{q}$, and a partition
\[s_k = \txcoprod_{\ell \in \ufs{r}} \, s_{k,\ell} \subseteq \ufs{n}_k,\]
we first define the following marker, permutation, and isomorphism.
\begin{itemize}
\item We define the marker 
\begin{equation}\label{bskl}
\bskl = \angscompkskell
\end{equation}
obtained from $\ang{s}$ by replacing the $k$-th subset, $s_k$, by $s_{k,\ell} \subseteq \ufs{n}_k$ \cref{angs_compk_ski}.
\item We define the permutation
\begin{equation}\label{usi_perm}
\big(\angstimes \fto[\iso]{\usi_{\angs, k, \ang{s_{k,\ell}}_{\ell \in \ufs{r}}}} 
\txcoprod_{\ell \in \ufs{r}}\, \bskltimes \big) \in \Si_{\angstimesc}
\end{equation}
that reorders $\angstimes$ to $\coprod_{\ell \in \ufs{r}}\, \bskltimes$, where $\angstimes$ and each 
\begin{equation}\label{bskltimes}
\bskltimes = s_1 \times \cdots \times s_{k-1} \times s_{k,\ell} \times s_{k+1} \times \cdots \times s_q
\end{equation}
are equipped with the lexicographic ordering inherited from $\txprod_{j \in \ufs{q}}\, \ufs{n}_j$.  
\item We define the unique isomorphism 
\begin{equation}\label{ual_iso}
\ga\big(x \sscs \ang{\ep_{\bskltimesc}}_{\ell \in \ufs{r}}\big) \usi_{\angs, k, \ang{s_{k,\ell}}_{\ell \in \ufs{r}}} 
\fto[\iso]{\ual_{x; \angs, k, \ang{s_{k,\ell}}_{\ell \in \ufs{r}}}} \ep_{\angstimesc}
\end{equation}
in the translation category $\Op(\angstimesc)$ with the indicated domain and codomain.
\end{itemize}
The gluing isomorphism of the strong $\angordnbe$-system $(\zdsgad a, \glu)$ at $(x; \angs, k, \ang{s_{k,\ell}}_{\ell \in \ufs{r}})$ \cref{gluing-morphism} is defined as the composite isomorphism in $\A$ in \cref{zdsgad_glu}, using the following abbreviations.
\[\ang{\cdots}_{\ell} = \ang{\cdots}_{\ell \in \ufs{r}} \qquad
\usi = \usi_{\angs, k, \ang{s_{k,\ell}}_{\ell \in \ufs{r}}} \qquad
\ual = \ual_{x; \angs, k, \ang{s_{k,\ell}}_{\ell \in \ufs{r}}}\]
\begin{equation}\label{zdsgad_glu}
\begin{tikzpicture}[vcenter]
\def\h{1} \def\u{-1.2} \def\v{-1.8}
\draw[0cell=.9]
(0,0) node (a11) {\gaA_r\big(x \sscs \ang{(\zdsgad a)_{\bskl}}_{\ell} \big)}
(a11)++(4,0) node (a12) {(\zdsgad a)_{\angs}}
(a11)++(0,\u) node (a21) {\gaA_r\big(x \sscs \bang{\gaA_{\bskltimesc} ( \ep_{\bskltimesc} \sscs \ang{a_{\bdi}}_{\bdi \in \bskltimes} ) }_{\ell} \big)}
(a12)++(0,\u+.5*\v) node (a22) {\gaA_{\angstimesc} \big(\ep_{\angstimesc} \sscs \ang{a_{\bdi}}_{\bdi \in \angstimes}\big)}
(a21)++(0,\v) node (a3) {\gaA_{\angstimesc}\big(\ga\big( x \sscs \ang{\ep_{\bskltimesc}}_{\ell} \big) \sscs \ang{\ang{a_{\bdi}}_{\bdi \in \bskltimes}}_{\ell} \big)}
(a3)++(0,\u) node (a4) {\gaA_{\card{\angstimes}}\big(\ga( x \sscs \ang{\ep_{\bskltimesc}}_{\ell}) \usi \sscs \ang{a_{\bdi}}_{\bdi \in \angstimes} \big)}
;
\draw[1cell=.9]
(a11) edge node {\glu_{x; \angs, k, \ang{\skell}_{\ell}}} (a12)
(a12) edge[equal,transform canvas={xshift=0ex}] (a22)
(a11) edge[equal,transform canvas={xshift=0em}, shorten >=-.5ex] (a21)
(a21) edge node[swap] {\phiA_{(r;\, \ang{\bskltimesc}_{\ell})}} node {\iso} (a3)
(a3) edge[equal,shorten <=-.5ex,shorten >=-.5ex] node[swap,inner sep=3pt,pos=.6] {\mathbf{eq}} (a4)
(a4) [rounded corners=2pt] -| node[swap,pos=.75] {\gaA_{\angstimesc}(\ual \sscs 1^{\angstimesc})} node[pos=.75] {\iso} (a22)
;
\end{tikzpicture}
\end{equation}
The arrows in \cref{zdsgad_glu} are given as follows.
\begin{itemize}
\item The top left and right equalities follow from \cref{zdsgad_comp_obj}.
\item The isomorphism $\phiA_{(r;\, \ang{\bskltimesc}_{\ell})}$ is a component of the associativity constraint of $\A$ \cref{phiA_component}.
\item The equality labeled $\mathbf{eq}$ follows from the action equivariance axiom \cref{pseudoalg_action_sym} for $\A$, applied to the permutation $\usi$ \cref{usi_perm}.
\item In the lower right arrow, $\ual$ is the unique isomorphism defined in \cref{ual_iso}.
\end{itemize}
\cref{zdsgad_obj_welldef} proves that $(\zdsgad a,\glu)$ is a strong $\angordnbe$-system in $\A$.
\item[Morphisms]
For a morphism 
\[f = \ang{f_{\bdi} \in \A}_{\bdi \sins \ufs{n_1 n_2 \cdots n_q}} \cn a \to b \inspace \proAnbe\]
and a marker $\angs = \ang{s_j \subseteq \ufs{n}_j}_{j \in \ufs{q}}$, the morphism
\begin{equation}\label{zdsgad_f}
\zdsgad a \fto{\zdsgad f} \zdsgad b \inspace \syssg{\A}{\sord{n}{\be}}
\end{equation} 
has $\angs$-component morphism \cref{theta_angs} defined as 
\begin{equation}\label{zdsgad_comp_mor}
(\zdsgad f)_{\angs} = \gaA_{\angstimesc}\big(1_{\ep_{\angstimesc}} \sscs \ang{f_{\bdi}}_{\bdi \in \angstimes}\big) \inspace \A.
\end{equation}
\cref{zdsgad_mor_welldef} proves that $\zdsgad f$ is a morphism of strong $\angordnbe$-systems in $\A$.
\item[Functoriality]
The functoriality of $\zdsgad$ follows from the functoriality of $\gaA$ and the fact that identities and composition are defined coordinatewise in the $\angordnbe$-twisted product $\proAnbe$ and componentwise in $\Asgangordnbe$.
\item[Basepoint preservation]
The basepoint $\zero \in \proAnbe$ is given coordinatewise by $\zero = \gaA_0(*) \in \A$.  It is sent by $\zdsgad$ to the base $\angordnbe$-system $(\zero,1_\zero)$ in $\A$ for the following reasons. 
\begin{itemize}
\item The $\angs$-component object \cref{zdsgad_comp_obj} 
\[(\zdsgad \zero)_{\angs} = \gaA_{\angstimesc}\big(\ep_{\angstimesc} \sscs \ang{\zero}_{\bdi \in \angstimes}\big)\] 
is equal to $\zero \in \A$ by $\angstimesc$ applications of the basepoint axiom \cref{pseudoalg_basept_axiom} for $\A$.
\item The gluing isomorphism $\glu_{x; \angs, k, \ang{\skell}_{\ell}}$ \cref{zdsgad_glu} is equal to the identity $1_\zero$ for the following reasons.
\begin{itemize}
\item The associativity constraint $\phiA_{(r;\, \ang{\bskltimesc}_{\ell})}$ in \cref{zdsgad_glu} is equal to $1_\zero$ by the composition axiom \cref{pseudoalg_comp_axiom} for $\A$ and \cref{phi_id}.
\item The isomorphism $\gaA_{\angstimesc}(\ual; 1^{\angstimesc})$ in \cref{zdsgad_glu} is equal to $1_\zero$ by the basepoint axiom \cref{pseudoalg_basept_axiom} for $\A$ and the naturality of the associativity constraint $\phiA$.
\end{itemize}
\end{itemize}
\end{description}
This finishes the definition of the pointed functor $\zdsgad$ \cref{zdsgad}.
\end{definition}

\begin{explanation}[Comparing $\zbsgad$ and $\zdsgad$]\label{expl:zdsgad_zbsgad}
\Cref{table.zbsgad} compares the pointed functors $\zbsgad$ and $\zdsgad$ \pcref{def:zbsgad,def:zdsgad}.
\begin{figure}[H] 
\centering
{\renewcommand{\arraystretch}{1.3}%
{\setlength{\tabcolsep}{1ex}
\begin{tabular}{cr|cr}
$\proAnbe \fto{\zbsgad} \Asgsmaangordnbe$ & \cref{zbsgad} 
& $\proAnbe \fto{\zdsgad} \Asgangordnbe$ & \cref{zdsgad}\\[.3ex] \hline 
$s \subseteq \txprod_{j \in \ufs{q}}\, \ufs{n}_j$ & \cref{subset_s}
& $\angstimes \subseteq \txprod_{j \in \ufs{q}}\, \ufs{n}_j$ & \cref{angstimes}\\ 
$\gaA_{\cards}\big(\ep_{\cards} \sscs \ang{a_{\bdi}}_{\bdi \in s}\big)$ & \cref{zbsgad_comp_obj} 
& $\gaA_{\angstimesc}\big( \ep_{\angstimesc} \sscs \ang{a_{\bdi}}_{\bdi \in \angstimes}\big)$ & \cref{zdsgad_comp_obj}\\
$\sig_{s,\ang{s_\ell}_\ell} \in \Si_{\cards}$ & \cref{sigma_angsi} 
& $\usi_{\angs, k, \ang{\skell}_\ell} \in \Si_{\angstimesc}$ & \cref{usi_perm}\\
$\alp_{x;\, s, \ang{s_\ell}_\ell} \in \Op(\cards)$ & \cref{al_Op}
& $\ual_{x; \angs, k, \ang{\skell}_\ell} \in \Op(\angstimesc)$ & \cref{ual_iso}\\
$\gl_{x;\, s, \ang{s_\ell}_\ell} $ & \cref{zbsgad_gl}
& $\glu_{x; \angs, k, \ang{\skell}_\ell} $ & \cref{zdsgad_glu} \\
$\gaA_{\cards}(\alp; 1^{\cards}) \circ \mathbf{eq} \circ \phiA$ & \cref{zbsgad_gl} &
$\gaA_{\angstimesc}(\ual; 1^{\angstimesc}) \circ \mathbf{eq} \circ \phiA$ & \cref{zdsgad_glu} \\
$\gaA_{\cards}\big(1_{\ep_{\cards}} \sscs \ang{f_{\bdi}}_{\bdi \in s}\big)$ & \cref{zbsgad_comp_mor} 
& $\gaA_{\angstimesc}\big(1_{\ep_{\angstimesc}} \sscs \ang{f_{\bdi}}_{\bdi \in \angstimes}\big)$ & \cref{zdsgad_comp_mor}
\end{tabular}}}
\caption{A comparison of two pointed functors.}
\label{table.zbsgad}
\end{figure}
\cref{zdsgad_obj_welldef,zdsgad_mor_welldef} reuse some of the proofs in \cref{ch:sgoprod} by using this dictionary between \cref{def:zbsgad,def:zdsgad}.
\end{explanation}

\subsection*{Well Definedness}
The rest of this section proves \cref{zdsgad_obj_welldef,zdsgad_mor_welldef}, which ensure that $\zdsgad$ is a well-defined pointed functor.  Recall that $\Op$ is a $\Uinf$-operad \pcref{as:OpA'} and $(\A,\gaA,\phiA)$ is an $\Op$-pseudoalgebra \pcref{def:pseudoalgebra}.

\begin{lemma}\label{zdsgad_obj_welldef}
The pair $(\zdsgad a,\glu)$ defined in \crefrange{zdsgad_comp_obj}{zdsgad_glu} is a strong $\angordnbe$-system in $\A$.
\end{lemma}

\begin{proof}
Each component of the gluing morphism $\glu$ \cref{zdsgad_glu} is an isomorphism.  We verify that the pair $(\zdsgad a,\glu)$ satisfies the axioms \crefrange{system_obj_unity}{system_commutativity} in \cref{def:nsystem} for an $\angordnbe$-system, which consists of the same data as an $\angordn$-system.
\begin{description}
\item[All but commutativity]
For the first unity axiom \cref{system_unity_i}, we assume $s_j = \emptyset$ for some $j \in \ufs{q}$.  We want to show that the gluing isomorphism $\glu_{x; \angs, k, \ang{\skell}_\ell}$ \cref{zdsgad_glu} is equal to $1_\zero$ in $\A$.  In fact, each of its two constituent arrows is equal to $1_\zero$ for the following reasons.  
\begin{enumerate}
\item The assumption $s_j = \emptyset$ implies
\[\angstimes = \emptyset = \bskltimes \forspace \ell \in \ufs{r}\]
and
\[\ep_{\bskltimesc} = \ep_0 = * \in \Op(0).\]
The associativity constraint $\phiA_{(r;\,\ang{0}_\ell)}$ in \cref{zdsgad_glu} is equal to $1_\zero$ by \cref{phi_id}.
\item The unique isomorphism \cref{ual_iso}
\[\ual_{x;\angs, k, \ang{\skell}_\ell} \inspace \Op(\angstimesc) = \Op(0) = \{*\}\]
is $1_*$, so the other arrow in \cref{zdsgad_glu} is $\gaA_0(1_*)$.  This is equal to $1_\zero$ by the functoriality of $\gaA_0$.
\end{enumerate}
This proves the first unity axiom \cref{system_unity_i}.  Each of the five axioms  \cref{system_obj_unity,system_naturality,system_unity_iii,system_equivariance,system_associativity} is proved by the corresponding proof in \cref{zbsgad_obj_welldef} after changing the notation as discussed in \cref{expl:zdsgad_zbsgad}.
\item[Commutativity]
For the commutativity axiom \cref{system_commutativity} for $(\zdsgad a,\glu)$, we consider a pair of objects
\[(x,y) \in \Op(r) \times \Op(t),\]
a marker $\angs = \ang{s_j \subseteq \ufs{n}_j}_{j \in \ufs{q}}$, two distinct indices $k \neq h \in \ufs{q}$, partitions
\[s_k = \txcoprod_{\ell \in \ufs{r}}\, \skell \subseteq \ufs{n}_k \andspace 
s_h = \txcoprod_{p \in \ufs{t}}\, \shp \subseteq \ufs{n}_h,\]
and the following notation involving $\intr$ \cref{intr_jk}, $\twist_{t,r}$ \cref{eq:transpose_perm}, $\gaA$ \cref{gaAn}, $\phiA$ \cref{phiA}, $\compk$ \cref{compk}, $\angstimes$ \cref{angstimes}, $\ep$ \cref{Opn_object}, $\usi$ \cref{usi_perm}, and $\ual$ \cref{ual_iso}.
\[\left\{
\begin{aligned}
\ang{\cdots}_\ell &= \ang{\cdots}_{\ell \in \ufs{r}} & \ang{\cdots}_{\ell,p} &= \ang{\ang{\cdots}_{\ell \in \ufs{r}}}_{p \in \ufs{t}} & \ang{a_{\bdi \in ?}} &= \ang{a_{\bdi}}_{\bdi \in ?}\\
\ang{\cdots}_p &= \ang{\cdots}_{p \in \ufs{t}} & \ang{\cdots}_{p,\ell} &= \ang{\ang{\cdots}_{p \in \ufs{t}}}_{\ell \in \ufs{r}} & 
\end{aligned}
\right.\]
\[\left\{\scalebox{.85}{$
\begin{aligned}
u &= \angstimes & \bda &= \ang{a_{\bdi \in u}}\\
u_\ell &= \angscompkskelltimes & \bda_\ell &= \ang{a_{\bdi \in u_\ell}}\\
u^p &= (\angs \comph \shp)^\stimes & \bda^p &= \ang{a_{\bdi \in u^p}}\\
\uellp &= \big(\angs \compk \skell \comph \shp\big)^\stimes \phantom{M} & \bdaellp &= \ang{a_{\bdi \in \uellp}}
\end{aligned}$}
\right.\qquad
\left\{\scalebox{.85}{$
\begin{aligned}
\vphi &= \phiA & \bdaellpp &= \gaA_{|\uellp|}\big( \ep_{|\uellp|}; \bdaellp \big)\\
xy &= x \intr y \phantom{M} & (xy)' &= \ga\big( xy; \ang{\ep_{|\uellp|}}_{p,\ell} \big)\\
xy\twist &= xy\twist_{t,r} & (xy\twist)' &= \ga\big( (xy)\twist; \ang{\ep_{|\uellp|}}_{\ell,p} \big)\\
yx &= y \intr x & (yx)' &= \ga\big( yx; \ang{\ep_{|\uellp|}}_{\ell,p} \big)
\end{aligned}$}
\right.\]
\begin{equation}\label{usigma_permutations}
\left\{
\begin{aligned}
\usi^p &= \usi_{\angs \comph \shp,\, k, \ang{\skell}_\ell} \cn u^p \fiso \txcoprod_{\ell \in \ufs{r}}\, \uellp
\phantom{M} & \usi^\stimes &= \txcoprod_{p \in \ufs{t}}\, \usi^p\\
\usi_\ell &= \usi_{\angscompkskell,\, h, \ang{\shp}_p} \cn u_\ell \fiso \txcoprod_{p \in \ufs{t}}\, \uellp
& \usi_\stimes &= \txcoprod_{\ell \in \ufs{r}}\, \usi_\ell\\
\usi' &= \usi_{\angs,\, k, \ang{\skell}_\ell} \cn u \fiso \txcoprod_{\ell \in \ufs{r}}\, u_\ell\\
\usi'' &= \usi_{\angs,\, h, \ang{\shp}_p} \cn u \fiso \txcoprod_{p \in \ufs{t}}\, u^p
\end{aligned}
\right.
\end{equation}
\[\left\{\scalebox{.9}{$
\begin{aligned}
x_p' &= \ga\big(x; \ang{\ep_{|\uellp|}}_\ell \big) \phantom{M} &
\ual^p &= \ual_{x; \angs \comph \shp,\, k, \ang{\skell}_\ell} \cn x_p' \usi^p \fiso \ep_{|u^p|}\\
y_\ell' &= \ga\big(y; \ang{\ep_{|\uellp|}}_p \big) &
\ual_\ell &= \ual_{y; \angscompkskell,\, h, \ang{\shp}_p} \cn y_\ell' \usi_\ell \fiso \ep_{|u_\ell|}\\
x' &= \ga\big(x; \ang{\ep_{|u_\ell|}}_\ell \big) &
\ual' &= \ual_{x; \angs,\, k, \ang{\skell}_\ell} \cn x'\usi' \fiso \ep_{|u|}\\
y' &= \ga\big(y; \ang{\ep_{|u^p|}}_p \big) &
\ual'' &= \ual_{y; \angs,\, h, \ang{\shp}_p} \cn y'\usi'' \fiso \ep_{|u|}
\end{aligned}$}
\right.\]
In addition, we use juxtaposition to denote $\gaA$,
\[w\ang{b_i}_{i \in \ufs{m}} = \gaA_m\big(w ; \ang{b_i}_{i \in \ufs{m}} \big),\]
for any objects or morphisms $w \in \Op(m)$ and $\ang{b_i}_{i \in \ufs{m}} \in \A^m$, and similarly for $\ga$.  For example, the following objects in $\A$ appear in the diagram \cref{zdsgad_com_diagram}.
\[\begin{split}
(xy\twist)\ang{\bdaellpp}_{\ell,p} 
&= \gaA_{tr} \Big((x \intr y)\twist_{t,r}; \bang{\bang{\gaA_{|\uellp|}\big( \ep_{|\uellp|}; \ang{a_{\bdi}}_{\bdi \in \uellp} \big) }_{\ell \in \ufs{r}} }_{p \in \ufs{t}} \Big)\\
y\ang{x \ang{\bdaellpp }_\ell }_p 
&= \gaA_t\Big(y; \Bang{\gaA_r\Big(x; \bang{\gaA_{|\uellp|}\big( \ep_{|\uellp|}; \ang{a_{\bdi}}_{\bdi \in \uellp} \big) }_{\ell \in \ufs{r}} \Big) }_{p \in \ufs{t}} \Big)\\
y\ang{\ep_{|u^p|} \bda^p}_p
&= \gaA_t\Big(y; \bang{\gaA_{|u^p|}\big(\ep_{|u^p|}; \ang{a_{\bdi}}_{\bdi \in u^p} \big) }_{p \in \ufs{t}}\Big)\\
(y' \usi'')\bda &= \gaA_{|u|}\Big(\ga\big(y; \ang{\ep_{|u^p|}}_{p \in \ufs{t}} \big) \usi_{\angs,\, h, \ang{\shp}_{p \in \ufs{t}}} ; \ang{a_{\bdi}}_{\bdi \in u} \Big)
\end{split}\]
We define the block permutation  
\begin{equation}\label{twist_block}
\ol{\twist} \in \Si_{|u|}
\end{equation}
induced by the transpose permutation $\twist_{t,r} \in \Si_{tr}$ \cref{eq:transpose_perm} that permutes $tr$ blocks of lengths $\ang{\ang{|\uellp|}_{\ell \in \ufs{r}}}_{p \in \ufs{t}}$.  Moreover, a typical isomorphism in the translation category $\Op(|u|)$ is denoted by $!$. 

Using the notation above and the definition of $\glu$ \cref{zdsgad_glu}, the commutativity diagram \cref{system_commutativity} for $(\zdsgad a,\glu)$ is the boundary diagram in \cref{zdsgad_com_diagram}.
\setlength{\fboxrule}{0.2pt}
\setlength{\fboxsep}{1.5pt}
\begin{equation}\label{zdsgad_com_diagram}

\end{equation}
In \cref{zdsgad_com_diagram}, the equalities and arrows are given as follows.
\begin{itemize}
\item Each equality labeled $\mathbf{eq}$ is given by the action equivariance axiom \cref{pseudoalg_action_sym} for $\A$.
\item The equalities labeled $\mathbf{as}$, $\mathbf{top}$, and $\mathbf{b}$ are given by, respectively, the associativity axiom, the top equivariance axiom, and the bottom equivariance axiom for the $\Gcat$-operad $\Op$.
\item Each $\vphi = \phiA$ is a component of the associativity constraint for $\A$ \cref{phiA_component}.
\item $\pcom_{r,t} \cn xy\twist \fiso yx$ is a component of the $(r,t)$-pseudo-commutativity isomorphism of $\Op$ \cref{pseudocom_isos}.  It is the unique isomorphism with the indicated domain and codomain, since $\Op(tr)$ is a translation category.  
\item Each of the three arrows
\[\begin{split}
(x\ang{y_\ell'}_\ell) \ol{\twist} & \fto{!} y\ang{x_p'}_p,\\
(x\ang{y_\ell'}_\ell) \usi_\stimes \usi' & \fto{!} (y\ang{x_p'}_p) \usi^\stimes \usi'', \andspace\\
(x \ang{y_\ell' \usi_\ell}_\ell) \usi' & \fto{!} (y \ang{x_p' \usi^p}_p) \usi''
\end{split}\]
is the unique isomorphism in the translation category $\Op(|u|)$ with the indicated domain and codomain.
\item $\ual^p \in \Op(|u^p|)$, $\ual_\ell \in \Op(|u_\ell|)$, $\ual' \in \Op(|u|)$, and $\ual'' \in \Op(|u|)$ are the unique isomorphisms defined before the diagram \cref{zdsgad_com_diagram}.
\end{itemize}

The regions inside the diagram \cref{zdsgad_com_diagram} commute for the following reasons.
\begin{itemize}
\item Each of the three unlabeled regions commutes because it consists entirely of equalities.
\item Each of the three regions labeled $\fbox{Tr}$ commutes because $\Op(|u|)$ is a translation category.
\item The upper right region labeled $\fbox{T}$ commutes by the top equivariance axiom \cref{pseudoalg_topeq} for $\A$.
\item Each of the three regions labeled $\fbox{N}$ commutes by the naturality of the associativity constraint $\phiA$.
\item Each of the two regions labeled $\fbox{C}$ commutes by the composition axiom \cref{pseudoalg_comp_axiom} for $\A$.
\item Each of the two regions labeled $\fbox{B}$ commutes by the bottom equivariance axiom \cref{pseudoalg_boteq} for $\A$.
\item Each of the two regions labeled $\fbox{E}$ commutes by the action equivariance axiom \cref{pseudoalg_action_sym} for $\A$.
\end{itemize}

To show that the central region $\fbox{E'}$ commutes, we first observe that the permutations that appear in $\fbox{E'}$ yield the commutative pentagon \cref{perm_pentagon} by  \cref{usi_perm,usigma_permutations,twist_block}.
\begin{equation}\label{perm_pentagon}
\begin{tikzpicture}[vcenter]
\def\e{1em} \def\g{.25} \def\h{2.2} \def\u{-.8} \def\v{-1.4}
\draw[0cell]
(0,0) node (a0) {u}
(a0)++(-\h,\u) node (a11) {\txcoprod_{p \in \ufs{t}}\, u^p}
(a0)++(\h,\u) node (a12) {\txcoprod_{\ell \in \ufs{r}}\, u_\ell}
(a11)++(\g,\v) node (a21) {\txcoprod_{p \in \ufs{t}}\, \txcoprod_{\ell \in \ufs{r}}\, \uellp}
(a12)++(-\g,\v) node (a22) {\txcoprod_{\ell \in \ufs{r}}\, \txcoprod_{p \in \ufs{t}}\, \uellp}
;
\draw[1cell=.9]
(a0) edge node {\usi'} (a12)
(a12) edge[transform canvas={xshift=-\e}] node[pos=.3] {\usi_\stimes} (a22)
(a0) edge node[swap] {\usi''} (a11)
(a11) edge[transform canvas={xshift=\e}, shorten <=-.5ex] node[swap,pos=.2] {\usi^\stimes} (a21)
(a21) edge node {\ol{\twist}} (a22)
;
\end{tikzpicture}
\end{equation}
The commutative pentagon \cref{perm_pentagon} implies that the unique isomorphism in the bottom boundary of $\fbox{E'}$ has the form
\[(x\ang{y_\ell'}_\ell) \ol{\twist} \usi^\stimes \usi'' 
= (x\ang{y_\ell'}_\ell) \usi_\stimes \usi' \to (y\ang{x_p'}_p) \usi^\stimes \usi''.\]
Since $\Op(|u|)$ is a translation category, the above isomorphism is equal to the $(\usi^\stimes \usi'')$-action on the unique isomorphism in the top boundary of $\fbox{E'}$,
\[(x\ang{y_\ell'}_\ell) \ol{\twist} \to y\ang{x_p'}_p.\]
Thus, the region $\fbox{E'}$ commutes by the action equivariance axiom \cref{pseudoalg_action_sym} for $\A$.  This proves the commutativity axiom \cref{system_commutativity} for $(\zdsgad a,\glu)$.\qedhere
\end{description}
\end{proof}

\begin{lemma}\label{zdsgad_mor_welldef}
For each morphism $f \cn a \to b$ in $\proAnbe$, the collection 
\[\zdsgad a \fto{\zdsgad f} \zdsgad b\]
defined in \cref{zdsgad_f} is a morphism in $\Asgangordnbe$.
\end{lemma}

\begin{proof}
To check that $\zdsgad f$ satisfies the unity axiom \cref{nsystem_mor_unity} and the compatibility axiom \cref{nsystem_mor_compat}, we reuse the proof of \cref{zbsgad_mor_welldef} 
by changing the notation as discussed in \cref{expl:zdsgad_zbsgad}.
\end{proof}

\section{Unit and Counit}
\label{sec:hgo_unit_counit}

Throughout this section, we consider a $\Uinf$-operad $(\Op,\ga,\opu)$ \pcref{as:OpA'}, an $\Op$-pseudoalgebra $(\A,\gaA,\phiA)$ \pcref{def:pseudoalgebra}, and an object $\angordnbe = \ang{\ordn_j^{\be_j}}_{j \in \ufsq} \in \GG \setminus \{\vstar,\ang{}\}$ of length $q>0$ \cref{GG_objects}.  This section constructs the unit $\unit \cn 1 \to \zdsg \zdsgad$ and the counit $\counit \cn \zdsgad\zdsg \to 1$ for the functors
\begin{equation}\label{zdsg_seci}
\begin{tikzpicture}[vcenter]
\def\s{20}
\draw[0cell]
(0,0) node (a1) {\phantom{X}}
(a1)++(1.6,0) node (a2) {\phantom{X}}
(a1)++(-.15,0) node (a1') {\proAnbe}
(a2)++(.4,-.02) node (a2') {\Asgangordnbe}
;
\draw[1cell=.9]
(a1) edge[bend left=\s] node {\zdsgad} (a2)
(a2) edge[bend left=\s] node {\zdsg} (a1)
;
\end{tikzpicture}
\end{equation}
defined in \cref{hgosgtoprod,zdsgad}.  \cref{thm:zdsg_eq} proves that $(\zdsgad,\zdsg,\unit,\counit)$ is an adjoint equivalence.

\secoutline
\begin{itemize}
\item \cref{def:zdsgunit} defines the unit $\unit \cn 1 \to \zdsg \zdsgad$ as the identity natural transformation of the identity functor on $\proAnbe$. 
\item \cref{def:zdsgcounit} defines the counit $\counit \cn \zdsgad\zdsg \to 1$.
\item \cref{ex:counit_two} illustrates the counit $\counit$ for objects in $\GG$ of length 2.
\item \cref{counit_agl_welldef} proves that each component of $\counit$ is a well-defined isomorphism of strong $\angordnbe$-systems.  This proof uses the full extent of the axioms of an $\Op$-pseudoalgebra and of an $\angordnbe$-system \pcref{def:pseudoalgebra,def:nsystem,def:nbeta_gcat}.
\end{itemize}

\subsection*{Unit}
First, we define the unit for the adjunction $(\zdsgad, \zdsg)$.   This definition is analogous to \cref{def:zbsgunit}, which defines the unit for the adjoint equivalence $(\zbsgad,\zbsg)$.

\begin{definition}[Unit]\label{def:zdsgunit}
For the pointed functors $\zdsgad$ \cref{zdsgad} and $\zdsg$ \cref{hgosgtoprod}, we define $\unit$ as the identity natural transformation
\begin{equation}\label{unit_zdsg}
\begin{tikzpicture}[vcenter]
\def\h{1.8} \def\v{-.7} \def\s{17}
\draw[0cell]
(0,0) node (a1) {\proAnbe}
(a1)++(\h,\v) node (a2) {\phantom{\proAnbe}}
(a2)++(.2,0) node (a2') {\Asgangordnbe}
(a1)++(0,2*\v) node (a3) {\proAnbe}
(a1)++(0,\v) node (a) {}
;
\draw[1cell=.9]
(a1) edge[bend right=\s] node[swap] {1} (a3)
(a1) edge[bend left=\s] node {\zdsgad} (a2)
(a2) edge[bend left=\s, shorten <=-0ex] node {\zdsg} (a3)
;
\draw[2cell]
node[between=a and a2 at .4, shift={(0,.3*\v)}, rotate=0, 2labelw={above,\unit = 1_1,2pt}] {\Rightarrow}
;
\end{tikzpicture}
\end{equation}
of the identity functor on $\proAnbe$.  To verify that $\unit \cn 1 \to \zdsg\zdsgad$ is well defined, we consider an object or a morphism
\[a = \ang{a_{\bdi}}_{\bdi \in \ufs{n_1 n_2 \cdots n_q}} \in \proAnbe.\]
Using the action unity axiom \cref{pseudoalg_action_unity} for $\A$, \cref{i_connected,hgotoprod_def,Opn_object,zdsgad_comp_obj,zdsgad_comp_mor}, the $\bdi$-th coordinate of $\zdsg \zdsgad a$ is given by
\begin{equation}\label{uni_cod_obj}
(\zdsg \zdsgad a)_{\bdi} 
= (\zdsgad a)_{\ang{\{i_j\}}_{j \in \ufs{q}}} 
= \gaA_1\big(\ep_1 \sscs a_{\bdi}) 
= \gaA_1\big(\opu \sscs a_{\bdi}) 
= a_{\bdi}.
\end{equation}
The second equality in \cref{uni_cod_obj} uses the fact that
\[\big(\ang{\{i_j\}}_{j \in \ufs{q}}\big)^\stimes 
= \txprod_{j \in \ufs{q}}\, \{i_j\} 
= \{(i_1,\ldots,i_q)\} 
= \{\bdi\},\]
the one-element set containing only $\bdi$.  Thus, $\zdsg \zdsgad$ is equal to the identity functor on $\proAnbe$.  This finishes the definition of the identity natural transformation $\unit \cn 1 \to \zdsg\zdsgad$.  
\end{definition}

\subsection*{Counit}
Next, we define the counit $\counit$ for the adjunction $(\zdsgad,\zdsg)$.  Although this definition is the conceptual counterpart of \cref{def:zbsgcounit} defining the counit $\cous$ for $(\zbsgad,\zbsg)$, the current situation is more complicated.  The counit $\cous$ is componentwise a gluing isomorphism.  On the other hand, since objects in $\GG$ can have any finite length, the counit $\counit$ is componentwise an iterated composite of gluing isomorphisms, isomorphisms in $\Op$, and the associativity constraint $\phiA$.  See \cref{counit_as_gen}.

\begin{definition}[Counit]\label{def:zdsgcounit}
For the pointed functors $\zdsgad$ \cref{zdsgad} and $\zdsg$ \cref{hgosgtoprod}, we define the pointed natural isomorphism
\begin{equation}\label{counit_zdsg}
\begin{tikzpicture}[vcenter]
\def\h{2} \def\v{.7} \def\s{17}
\draw[0cell]
(0,0) node (a1) {\phantom{\proAnbe}}
(a1)++(.2,0) node (a1') {\proAnbe}
(a1)++(\h,\v) node (a2) {\Asgangordnbe}
(a1)++(\h,-\v) node (a3) {\Asgangordnbe}
(a1)++(\h,0) node (a) {}
;
\draw[1cell=.9]
(a1) edge[bend left=\s, shorten >=0ex] node {\zdsgad} (a2)
(a2) edge[bend left=\s, transform canvas={xshift=-.8em}, shorten >=-.3ex, shorten <=-.3ex] node {1} (a3)
(a3) edge[bend left=\s] node {\zdsg} (a1)
;
\draw[2cell]
node[between=a1 and a at .6, shift={(0,-.2*\v)}, rotate=0, 2label={above,\!\counit}] {\Rightarrow}
;
\end{tikzpicture}
\end{equation}
as follows.  For each $q$-tuple $\bdi = \ang{i_j}_{j \in \ufs{q}} \in \txprod_{j \in \ufs{q}}\, \ufs{n}_j$, we denote by
\begin{equation}\label{bdip}
\bdip = \bang{\{i_j\}}_{j \in \ufs{q}} = \big(\{i_1\}, \{i_2\}, \ldots, \{i_q\}\big)
\end{equation}
the marker whose $j$-th entry is the one-element subset $\{i_j\} \subseteq \ufs{n}_j$. 

\parhead{Domain of $\counit$}.  
Suppose $(a,\glu)$ is a strong $\angordnbe$-system in $\A$ \pcref{def:nsystem} and $\angs = \ang{s_j \subseteq \ufs{n}_j}_{j \in \ufs{q}}$ is a marker.  By \cref{hgotoprod_def,zdsgad_comp_obj}, the $\angs$-component object of $\zdsgad\zdsg(a,\glu)$ \cref{a_angs} is given as follows, where $\angstimes  = \txprod_{j \in \ufs{q}}\, s_j$.
\begin{equation}\label{cou_dom_obj}
\begin{split}
&\big(\zdsgad\zdsg(a,\glu)\big)_{\angs} \\
&= \gaA_{\angstimesc}\big(\ep_{\angstimesc} \sscs \ang{(\zdsg(a,\glu))_{\bdi}}_{\bdi \in \angstimes}\big)\\
&= \gaA_{\angstimesc}\big(\ep_{\angstimesc} \sscs \ang{a_{\bdip}}_{\bdi \in \angstimes}\big).
\end{split}
\end{equation}
The computation in \cref{cou_dom_obj} is also valid for morphisms in $\Asgangordnbe$ \pcref{def:nsystem_morphism}, using \cref{zdsgad_comp_mor} instead of \cref{zdsgad_comp_obj}.  The $(a,\glu)$-component of $\counit$ is an isomorphism of $\angordnbe$-systems 
\begin{equation}\label{cou_a}
\zdsgad\zdsg(a,\glu) \fto[\iso]{\counit_{(a,\glu)}} (a,\glu).
\end{equation}
Its $\angs$-component isomorphism \cref{theta_angs} is defined in \cref{cou_sj_empty,cou_a_s,counit_as_gen}.  Its naturality and pointedness are explained after that.
\begin{description}
\item[Base case]
Suppose $s_j = \emptyset$ for some $j \in \ufs{q}$.  Then 
\[\angstimes = \emptyset \andspace \ep_0 = * \in \Op(0).\] 
We define the $\angs$-component of $\counit_{(a,\glu)}$ to be the identity of $\zero \in \A$:
\begin{equation}\label{cou_sj_empty}
\big(\zdsgad\zdsg(a,\glu)\big)_{\angs} = \gaA_0(*) = \zero 
\fto{\counit_{(a,\glu), \angs} = 1_\zero} 
(a,\glu)_{\angs} = \zero.
\end{equation}
The equality in the codomain holds by the object unity axiom \cref{system_obj_unity} for $(a,\glu)$.  Assuming $s_j \neq \emptyset$ for all $j \in \ufs{q}$, the $\angs$-component of $\counit_{(a,\glu)}$ is defined inductively on $q > 0$ as follows. 
\item[Initial case]
If $q=1$, then 
\[\angs = \ang{s_j}_{j \in \ufs{q}} = (s_1) \andspace \angstimes = s_1.\]  
The $\angs$-component of $\counit_{(a,\glu)}$ is defined as the gluing isomorphism $\glu$ at $(\ep_{\card{s_1}}; (s_1), 1, \ang{\{i\}}_{i \in s_1})$ \cref{gluing-morphism}, as in the diagram \cref{cou_a_s}.
\begin{equation}\label{cou_a_s}
\begin{tikzpicture}[vcenter]
\def\v{-1}
\draw[0cell=.9]
(0,0) node (a11) {\big(\zdsgad\zdsg(a,\glu)\big)_{\angs}}
(a11)++(4.7,0) node (a12) {\phantom{a_{\angs}}}
(a12)++(0,-.05) node (a12') {a_{\angs}}
(a11)++(0,\v) node (a21) {\gaA_{\card{s_1}}\big(\ep_{\card{s_1}} \sscs \ang{a_{\{i\}}}_{i \in s_1}\big)}
(a12)++(0,\v) node (a22) {\phantom{a_{(s_1)}}}
(a22)++(0,-.05) node (a22') {a_{(s_1)}}
;
\draw[1cell=.9]
(a11) edge node {\counit_{(a,\glu), \angs}} (a12)
(a21) edge node {\glu_{\ep_{|s_1|} ;\, (s_1), 1, \ang{\{i\}}_{i \in s_1}}} node[swap] {\iso} (a22)
(a11) edge[equal,shorten <=-.5ex,shorten >=-.5ex] (a21)
(a12') edge[equal] (a22')
;
\end{tikzpicture}
\end{equation}
This gluing isomorphism of $(a,\glu)$ uses the partition of $s_1$ into one-element subsets, 
\begin{equation}\label{sone_singletons}
s_1 = \txcoprod_{i \in s_1} \{i\},
\end{equation}
with its natural ordering inherited from $\ufs{n}_1$.
\item[General case]
Inductively, for $q \geq 2$, we define the component $\counit_{(a,\glu),\angs}$ in \cref{counit_as_gen} using the following notation.
\begin{itemize}
\item We denote by
\begin{equation}\label{nbe_prime}
\nbehat = \bordi{n}{\be}{j}_{j=2}^q \in \GG
\end{equation}
the object of length $q-1$ obtained from $\sord{n}{\be}$ by removing the pointed finite $G$-set $\ordi{n}{\be}{1}$.
\item For each marker $\angs = \ang{s_j \subseteq \ufs{n}_j}_{j \in \ufs{q}}$, we denote by
\begin{equation}\label{s_prime}
\angshat = \ang{s_j \subseteq \ufs{n}_j}_{j=2}^q
\end{equation}
the marker obtained from $\angs$ by removing the first entry $s_1$.  Thus, we have
\begin{equation}\label{sone_angshat}
\begin{split}
\angs &= \big(s_1, \angshat\big),\\
\angshattimes &= \txprod_{j=2}^q\, s_j \subseteq \txprod_{j=2}^q\, \ufs{n}_j, \andspace\\
\angstimes &= s_1 \times \angshattimes \subseteq \txprod_{j \in \ufs{q}}\, \ufs{n}_j.
\end{split}
\end{equation}
\item For each $i_1 \in \ufs{n}_1$, we define the strong $\nbehat$-system 
\begin{equation}\label{aglu_ione}
\big(a_{\{i_1\}}, \glu_{\{i_1\}}\big) \in \Asgnbehat,
\end{equation}
called the \emph{$i_1$-restriction of $(a,\glu)$}, as follows.
\begin{description}
\item[Components] For each marker $\angsiiq = \ang{s_j \subseteq \ufs{n}_j}_{j=2}^q$, its $\angsiiq$-component object \cref{a_angs} is defined as the $(\{i_1\},\angsiiq)$-component object of $a$:
\begin{equation}\label{aglu_ione_comp}
(a_{\{i_1\}})_{\angsiiq} = a_{(\{i_1\},\, \angsiiq)} \in \A.
\end{equation}
\item[Gluing] For each object $x \in \Op(r)$ with $r \geq 0$, $k \in \ufs{q-1}$, and partition
\[s_{k+1} = \txcoprod_{\ell \in \ufs{r}}\, s_{k+1,\ell} \subseteq \ufs{n}_{k+1},\]
we define the gluing isomorphism \cref{gluing-morphism} as
\begin{equation}\label{aglu_ione_gluing}
(\glu_{\{i_1\}})_{x;\, \angsiiq, \,k, \ang{s_{k+1,\ell}}_{\ell \in \ufs{r}}} 
= \glu_{x;\, (\{i_1\},\, \angsiiq), \, k+1, \ang{s_{k+1,\ell}}_{\ell \in \ufs{r}}}.
\end{equation}
This gluing isomorphism is well defined because
\[(a_{\{i_1\}})_{(\angsiiq \,\compk\, s_{k+1,\ell})} 
= a_{((\{i_1\},\, \angsiiq)\, \comp_{k+1}\, (s_{k+1,\ell}))}.\]
\end{description}
Each axiom in \cref{def:nsystem} for $(a_{\{i_1\}}, \glu_{\{i_1\}})$ to be a strong $\nbehat$-system follows from the same axiom for the strong $\nbe$-system $(a,\glu)$.
\end{itemize}

For each marker $\angs = \ang{s_j \subseteq \ufs{n}_j}_{j \in \ufs{q}}$, the $\angs$-component of $\counit_{(a,\glu)}$ is defined as the composite of isomorphisms in $\A$ in \cref{counit_as_gen}.
\begin{equation}\label{counit_as_gen}
\begin{tikzpicture}[vcenter]
\def\h{3.4} \def\u{-1} \def\v{-1.4}
\draw[0cell=.8]
(0,0) node (a1) {\big(\zdsgad \zdsg(a,\glu)\big)_{\angs}}
(a1)++(0,\u) node (a2) {\gaA_{\angstimesc}\big(\ep_{\angstimesc} \sscs \ang{a_{\bdip}}_{\bdi \in \angstimes} \big)}
(a2)++(0,\v) node (a3) {\gaA_{\angstimesc}\big(\ep_{\cardsi} \intr \ep_{\angshattimesc} \sscs \ang{\ang{a_{(\{i_1\},\, \bdjp)}}_{\bdj \in \angshattimes} }_{i_1 \in s_1} \big)} 
(a3)++(0,\v) node (a4) {\gaA_{\cardsi}\big(\ep_{\cardsi} \sscs \bang{\gaA_{\angshattimesc} \big(\ep_{\angshattimesc} \sscs \ang{a_{(\{i_1\},\, \bdjp)}}_{\bdj \in \angshattimes} \big) }_{i_1 \in s_1} \big)}
(a3)++(0,\v) node (a4') {\phantom{\gaA_{|s_1|}\Big(\Big)}}
(a1)++(\h,0) node (b1) {\phantom{a_{\angs}}}
(b1)++(0,-.05) node (b1') {a_{\angs}}
(b1)++(0,\u+.5*\v) node (b2) {\gaA_{\cardsi}\big(\ep_{\cardsi} \sscs \ang{a_{(\{i_1\},\, \angshat)}}_{i_1 \in s_1} \big)}
;
\draw[1cell=.8]
(a1) edge node {\counit_{(a,\glu), \angs}} (b1)
(a1) edge[equal, shorten <=-.5ex] (a2)
(a2) edge node[swap] {\gaA_{|\angstimes|}(!\,; 1^{\angstimesc})} (a3)
(a3) edge[shorten >=-.5ex] node[swap] {(\phiA)^{-1}} (a4)
(b2) edge node[swap,pos=.5] {\glu_{\ep_{|s_1|}; \angs,\, 1, \ang{\{i_1\}}_{i_1 \in s_1}}} (b1')
(a4) [rounded corners=2pt] -| node[swap,pos=.7] {\gaA_{|s_1|}\big(1; \ang{\counit_{(a_{\{i_1\}}, \glu_{\{i_1\}}),\, \angshat}}_{i_1 \in s_1}\big)} (b2)
;
\end{tikzpicture}
\end{equation}
Going counterclockwise, the arrows in \cref{counit_as_gen} are given as follows.
\begin{itemize}
\item The upper left equality is given by the computation in \cref{cou_dom_obj}, where $\bdip = \ang{\{i_j\}}_{j \in \ufs{q}}$ is the $q$-tuple of one-element subsets \cref{bdip}.
\item In the middle left arrow, the objects
\[\ep_{\angstimesc} \in \Op(\angstimesc), \quad \ep_{\cardsi} \in \Op(\cardsi), \andspace 
\ep_{\angshattimesc} \in \Op(\angshattimesc)\]
are the chosen objects in \cref{Opn_object}.  The operation $\intr$ is the intrinsic pairing of $\Op$ \cref{intr_jk}.  Using \cref{sone_angshat}, $!$ denotes the unique isomorphism
\[\ep_{\angstimesc} \fto[\iso]{!} \ep_{\cardsi} \intr \ep_{\angshattimesc}\]
in the translation category $\Op(\angstimesc)$ with the indicated domain and codomain.  The object equality
\begin{equation}\label{anga_bdip}
\ang{a_{\bdip}}_{\bdi \in \angstimes} 
= \bang{\bang{a_{(\{i_1\},\, \bdjp)}}_{\bdj \in \angshattimes} }_{i_1 \in s_1} \in \A^{\angstimesc}
\end{equation}
follows from \cref{bdip}, \cref{sone_angshat}, and the lexicographic ordering on $\angstimes$ and $\angshattimes$.
\item The lower left arrow $(\phiA)^{-1}$ is the inverse of the associativity constraint
\[\phiA_{(\cardsi; \angshattimesc, \ldots, \angshattimesc)}\]
of $\A$ \cref{phiA}.
\item In the lower right arrow, for each $i_1 \in s_1$, the isomorphism
\[\gaA_{\angshattimesc}\big(\ep_{\angshattimesc} \sscs \ang{a_{(\{i_1\},\, \bdjp)}}_{\bdj \in \angshattimes} \big) 
\fto[\iso]{\counit_{(a_{\{i_1\}}, \glu_{\{i_1\}}),\, \angshat}} a_{(\{i_1\},\, \angshat)}\]
is the $\angshat$-component of the counit $\counit_{(a_{\{i_1\}}, \glu_{\{i_1\}})}$ for the $i_1$-restriction $(a_{\{i_1\}}, \glu_{\{i_1\}})$ of $(a,\glu)$ defined in \crefrange{aglu_ione}{aglu_ione_gluing}.  This component of the counit exists by the induction hypothesis and the fact that the $i_1$-restriction of $(a,\glu)$ is a strong $\nbehat$-system, where the object $\nbehat \in \GG$ \cref{nbe_prime} has length $q-1$.
\item The upper right arrow $\glu_{\ep_{\cardsi}; \angs,\, 1, \ang{\{i_1\}}_{i_1 \in s_1}}$ is the indicated gluing isomorphism of $(a,\glu)$, involving the object $\ep_{\cardsi} \in \Op(\cardsi)$ and the partition of $s_1$ into its one-element subsets \cref{sone_singletons}.
\end{itemize}
\cref{counit_agl_welldef} proves that $\counit_{(a,\glu)}$ is a well-defined isomorphism in $\Asgangordnbe$.
\item[Naturality]  
Morphisms in $\Asgangordnbe$ are determined by their $\angs$-component morphisms for markers $\angs = \ang{s_j \subseteq \ufsn_j}_{j \in \ufsq}$.  The naturality of $\counit$ with respect to morphisms in $\Asgangordnbe$ is proved as follows. 
\begin{description}
\item[Base case] 
For a marker $\angs$ with some $s_j = \emptyset$, the naturality of $\counit_{(a,\glu),\angs}$ with respect to $(a,\glu)$ follows from
\begin{itemize}
\item the definition \cref{cou_sj_empty} and
\item the unity axiom \cref{nsystem_mor_unity} for morphisms of $\angordnbe$-systems.
\end{itemize}
Assuming $s_j \neq \emptyset$ for all $j \in \ufs{q}$, the naturality of $\counit_{(a,\glu),\angs}$ with respect to $(a,\glu)$ is proved by an induction on $q>0$ as follows.
\item[Initial case]
If $q=1$, then the naturality of $\counit$ follows from 
\begin{itemize}
\item the definition of $\counit_{(a,\glu),\angs}$ in \cref{cou_a_s} and
\item the compatibility axiom \cref{nsystem_mor_compat} for morphisms of $\nbe$-systems.
\end{itemize}
\item[General case]
Inductively, if $q \geq 2$, then the naturality of $\counit$ follows from 
\begin{itemize}
\item the definition of $\counit_{(a,\glu),\angs}$ in \cref{counit_as_gen},
\item the functoriality of $\gaA$ \cref{gaAn},
\item the naturality of the associativity constraint $\phiA$ \cref{phiA}, 
\item the induction hypothesis, and
\item the compatibility axiom \cref{nsystem_mor_compat} for morphisms of $\nbe$-systems.
\end{itemize}
\end{description}
\item[Pointedness]
The basepoint of $\Asgangordnbe$ is the base $\angordnbe$-system $(\zero,1_\zero)$, with each component object given by the basepoint $\zero \in \A$ and each gluing morphism given by the identity $1_\zero$.  The pointedness of $\counit$ means that each component of $\counit_{(\zero,1_\zero)}$ is equal to $1_\zero$.
\begin{description}
\item[Base case]  
For a marker $\angs$ with some $s_j = \emptyset$, $\counit_{(\zero,1_\zero), \angs}$ is $1_\zero$ by \cref{cou_sj_empty}.  Assuming $s_j \neq \emptyset$ for all $j \in \ufs{q}$, the equality $\counit_{(\zero,1_\zero), \angs} = 1_\zero$
is proved by an induction on $q>0$ as follows.
\item[Initial case]
If $q=1$, then, by \cref{cou_a_s}, 
\[\counit_{(\zero,1_\zero),\angs} = \glu_{\ep_{\cardsi} ;\, (s_1), 1, \ang{\{i\}}_{i \in s_1}} = 1_\zero.\]
\item[General case]
Inductively, suppose $q \geq 2$.  Then each of the four constituent arrows of $\counit_{(\zero, 1_\zero),\angs}$ \cref{counit_as_gen} is equal to $1_\zero$ for the following reasons.
\begin{enumerate}
\item The first isomorphism $\gaA_{\angstimesc}(!\,; 1^{\angstimesc})$ is equal to $1_\zero$ by 
\begin{itemize}
\item the basepoint axiom \cref{pseudoalg_basept_axiom} for $\A$ and
\item the naturality of the associativity constraint $\phiA$.
\end{itemize}
\item The second isomorphism $(\phiA)^{-1}$ is equal to $1_\zero$ by
\begin{itemize}
\item the composition axiom \cref{pseudoalg_comp_axiom} for $\A$, with each $k_{j,i}$ there equal to 0,
\item \cref{phi_id}, and
\item the functoriality of $\gaA$.
\end{itemize}
\item The third isomorphism
\[\gaA_{\cardsi}\big(1; \bang{\counit_{(a_{\{i_1\}}, \glu_{\{i_1\}}), \angshat}}_{i_1 \in s_1}\big)\]
is equal to $1_\zero$ by 
\begin{itemize}
\item the induction hypothesis, which applies because the $i_1$-restriction of $(\zero,1_\zero)$ \cref{aglu_ione} is the base $\nbehat$-system, and 
\item the functoriality of $\gaA_{\cardsi}$.
\end{itemize}
\item The fourth isomorphism is
\[\glu_{\ep_{\cardsi}; \angs,\, 1, \ang{\{i_1\}}_{i_1 \in s_1}} = 1_\zero.\]
\end{enumerate}
\end{description}
\end{description}
This finishes the definition of the pointed natural isomorphism $\counit \cn \zdsgad\zdsg \to 1$ \cref{counit_zdsg}.
\end{definition}

\begin{example}[Counit for $q=2$]\label{ex:counit_two}
We unravel the $\angs$-component $\counit_{(a,\glu),\angs}$ \cref{counit_as_gen} in the case $q=2$.  In this case, we have
\[\begin{split}
\angordnbe &= \big(\ordi{n}{\be}{1}, \ordi{n}{\be}{2}\big) \in \GG,\\
\angs &= \big(s_1 \subseteq \ufs{n}_1, s_2 \subseteq \ufs{n}_2\big), \\
\angstimes &= s_1 \times s_2, \andspace \angshat = (s_2).
\end{split}\]
The $\angs$-component $\counit_{(a,\glu),\angs}$ is the following composite isomorphism in $\A$.
\begin{equation}\label{qtwo_diag}
\begin{tikzpicture}[vcenter]
\def\h{3.4} \def\u{-1} \def\v{-1.4}
\draw[0cell=.8]
(0,0) node (a1) {\big(\zdsgad \zdsg(a,\glu)\big)_{\angs}}
(a1)++(0,\u) node (a2) {\gaA_{\angstimesc}\big(\ep_{\angstimesc} \sscs \ang{\ang{a_{(\{i_1\}, \{i_2\})}}_{i_2 \in s_2}}_{i_1 \in s_1} \big)}
(a2)++(0,\v) node (a3) {\gaA_{\angstimesc}\big(\ep_{\cardsi} \intr \ep_{\cardsii} \sscs \ang{\ang{a_{(\{i_1\}, \{i_2\})}}_{i_2 \in s_2}}_{i_1 \in s_1} \big)} 
(a3)++(0,\v) node (a4) {\gaA_{\cardsi}\big(\ep_{\cardsi} \sscs \bang{\gaA_{\cardsii}\big(\ep_{\cardsii} \sscs \ang{a_{(\{i_1\}, \{i_2\})}}_{i_2 \in s_2} \big) }_{i_1 \in s_1} \big)}
(a1)++(\h,0) node (b1) {\phantom{a_{\angs}}}
(b1)++(0,-.05) node (b1') {a_{\angs}}
(b1)++(0,\u+\v/2) node (b2) {\gaA_{\cardsi}\big(\ep_{\cardsi} \sscs \ang{a_{(\{i_1\},\, s_2)}}_{i_1 \in s_1} \big)}
;
\draw[1cell=.8]
(a1) edge node {\counit_{(a,\glu), \angs}} (b1)
(a1) edge[equal, shorten <=-.5ex] (a2)
(a2) edge node[swap] {\gaA_{|\angstimes|}(!\,; 1^{|\angstimes|})} (a3)
(a3) edge[shorten >=-.3ex] node[swap] {(\phiA)^{-1}} (a4)
(b2) edge node[swap] {\glu_{\ep_{\cardsi};\, \angs,\, 1, \ang{\{i_1\}}_{i_1 \in s_1}}} (b1')
(a4) [rounded corners=2pt] -| node[swap,pos=.7] {\gaA_{\cardsi}\big(1; \ang{\glu_{\ep_{\cardsii};\, (\{i_1\},\, s_2),\, 2, \ang{\{i_2\}}_{i_2 \in s_2}}}_{i_1 \in s_1}\big)} (b2)
;
\end{tikzpicture}
\end{equation}
The essence of $\counit_{(a,\glu),\angs}$ is the iterated gluing given by its last two constituent arrows.  Before these gluing isomorphisms can be used, the object $\ep_{\angstimesc}$ needs to be replaced by $\ep_{\cardsi}$ and $\ep_{\cardsii}$ using the translation category structure of $\Op(\angstimesc)$ and the associativity constraint $\phiA$.
\end{example}

\subsection*{Well Definedness}
The rest of this section proves \cref{counit_agl_welldef}, which is used in \cref{def:zdsgcounit} to ensure that each component of the counit $\counit$ is well defined.

\begin{lemma}\label{counit_agl_welldef}
For each strong $\angordnbe$-system $(a,\glu)$ in $\A$, the collection 
\[\zdsgad\zdsg(a,\glu) \fto[\iso]{\counit_{(a,\glu)}} (a,\glu)\]
defined in \crefrange{cou_a}{counit_as_gen} is an isomorphism in $\Asgangordnbe$.
\end{lemma}

\begin{proof}
Each component of $\counit_{(a,\glu)}$ is an isomorphism.  The unity axiom \cref{nsystem_mor_unity} for $\counit_{(a,\glu)}$ holds by \cref{cou_sj_empty}.

\parhead{Compatibility}.  The compatibility axiom \cref{nsystem_mor_compat} for $\counit_{(a,\glu)}$ states that the diagram \cref{counit_aglu_compat} in $\A$ commutes for each object $x \in \Op(r)$ with $r \geq 0$, marker $\angs = \ang{s_j \subseteq \ufs{n}_j}_{j \in \ufs{q}}$, index $k \in \ufs{q}$, and partition $s_k = \txcoprod_{\ell \in \ufs{r}}\, s_{k,\ell}$, where $\glu'$ denotes the gluing isomorphism of the strong $\angordnbe$-system $\zdsgad\zdsg(a,\glu)$ \cref{zdsgad_glu}.
\begin{equation}\label{counit_aglu_compat}
\begin{tikzpicture}[vcenter]
\def\v{-1.5} \def\u{-.01} \def\w{-.06}
\draw[0cell=.85]
(0,0) node (a1) {\gaA_r\big(x \sscs \bang{\big(\zdsgad\zdsg(a,\glu)\big)_{\ang{s} \compk\, s_{k,\ell}}}_{\ell \in \ufs{r}} \big)}
(a1)++(0,\v) node (b1) {\gaA_r\big(x \sscs \ang{a_{\ang{s} \compk\, s_{k,\ell}}}_{\ell \in \ufs{r}} \big)}
(a1)++(5.2,0) node (a2) {\phantom{\big(\zdsgad\zdsg(a,\glu)\big)_{\ang{s}}}}
(a2)++(0,\u) node (a2') {\big(\zdsgad\zdsg(a,\glu)\big)_{\ang{s}}}
(a2)++(0,\v) node (b2) {\phantom{a_{\ang{s}}}}
(b2)++(0,\w) node (b2') {a_{\ang{s}}}
;
\draw[1cell=.8]
(a1) edge node {\glu'_{x;\, \ang{s},\, k, \ang{s_{k,\ell}}_{\ell \in \ufs{r}}}} (a2)
(b1) edge node {\glu_{x;\, \ang{s},\, k, \ang{s_{k,\ell}}_{\ell \in \ufs{r}}}} (b2)
(a1) edge[transform canvas={xshift=3em}] node[swap] {\gaA_r\big(1_x \sscs \bang{\counit_{(a,\glu), \ang{s} \compk\, s_{k,\ell}}}_{\ell \in \ufs{r}} \big)} (b1)
(a2') edge node {\counit_{(a,\glu), \ang{s}}} (b2')
;
\end{tikzpicture}
\end{equation}
\begin{description}
\item[Base case] 
If $s_j = \emptyset$ for some $j \in \ufs{q}$, then, by \cref{cou_sj_empty},
\[\counit_{(a,\glu),\angs} = 1_\zero = \counit_{(a,\glu), \angs \compk \skell}.\]
Moreover, by the first unity axiom \cref{system_unity_i} for the $\nbe$-systems $(a,\glu)$ and $\zdsgad\zdsg(a,\glu)$, both horizontal arrows in the diagram \cref{counit_aglu_compat} are equal to $1_\zero$.  In this case, the diagram \cref{counit_aglu_compat} commutes because each of the four arrows is equal to $1_\zero$.

Thus, we assume that $s_j \neq \emptyset$ for all $j \in \ufs{q}$.  The commutativity of the diagram  \cref{counit_aglu_compat} is proved by an induction on $q>0$ as follows.
\item[Initial case]  
If $q=1$, then $k=1$, $\angs = (s_1)$, and $s_1 = \txcoprod_{\ell \in \ufs{r}}\, s_{1,\ell}$.  By \cref{cou_a_s,zdsgad_glu}, we have the following gluing isomorphisms.
\[\begin{split}
\counit_{(a,\glu), \ang{s}} &= \glu_{\ep_{\cardsi};\, (s_1), 1, \ang{\{i\}}_{i \in s_1}}\\
\counit_{(a,\glu), \angscompkskell} 
&= \glu_{\ep_{|s_{1,\ell}|};\, (s_{1,\ell}), 1, \ang{\{i\}}_{i \in s_{1,\ell}}}\\
\glu'_{x; \angs,\, k,\ang{\skell}_{\ell \in \ufs{r}}} 
&= \gaA_{\cardsi}\big(\ual_{x;\, (s_1), 1, \ang{s_{1,\ell}}_{\ell \in \ufs{r}}} ; 1^{\cardsi} \big) \circ \mathbf{eq} \circ \phiA
\end{split}\]
Reusing the diagram \cref{cous_compat_diagram} by changing the notation as discussed in \cref{expl:zdsgad_zbsgad}, the compatibility diagram \cref{counit_aglu_compat} factors into three regions.  These regions commute by the associativity axiom \cref{system_associativity}, the equivariance axiom \cref{system_equivariance}, and the naturality axiom \cref{system_naturality} for the strong $\nbe$-system $(a,\glu)$.
\item[General case with $k=1$]
Inductively, suppose $q \geq 2$.  We first deal with the case $k=1$, so $s_1 = \txcoprod_{\ell \in \ufs{r}}\, s_{1,\ell}$.  We use the following notation, along with \cref{as_gamma,sigma_angsi,usi_perm,bdip,s_prime,sone_angshat,aglu_ione}.
\[\left\{\scalebox{.8}{$
\begin{aligned}
 \ang{\cdots}_{j_\ell} &= \ang{\cdots}_{j_\ell \in s_{1,\ell}}
& \ang{\cdots}_{i_1} &= \ang{\cdots}_{i_1 \in s_1} 
& \ang{\cdots}_\ell &= \ang{\cdots}_{\ell \in \ufs{r}} 
\\
v &= \angs & \ang{\cdots}_{\bdi} &= \ang{\cdots}_{\bdi \in v^\times}
& v_\ell &= \angs \comp_1 s_{1,\ell}  
\\
\ang{\cdots}_{p,t} &= \ang{\ang{\cdots}_p}_t
& \ang{\cdots}_{\bdj} &= \ang{\cdots}_{\bdj \in \angshattimes}
& \usi &= \usi_{\angs, 1, \ang{s_{1,\ell}}_{\ell}}  
\\
a_{j_\ell \bdj} &= a_{(\{j_\ell\},\, \bdjp)} & a_{j_\ell \angshat} &= a_{(\{j_\ell\},\, \angshat)}
& \sig_1 &= \sig_{s_1, \ang{s_{1,\ell}}_\ell}
\\
a_{i_1 \bdj} &= a_{(\{i_1\},\, \bdjp)} & a_{i_1 \angshat} &= a_{(\{i_1\},\, \angshat)}
& \sigbar_1 &= \sig_1\ang{\angshattimesc, \ldots, \angshattimesc}
\end{aligned}$}
\right.\]
\[\left\{\scalebox{.8}{$
\begin{aligned}
x_1 &= \ga\big(x; \ang{\ep_{|v_\ell^\stimes|}}_{\ell} \big)
& x_2 &= \ga\big(x; \ang{\ep_{|s_{1,\ell}|} }_{\ell} \big)
& x_3 &= \ga\big(x; \ang{\ep_{|s_{1,\ell}|} \intr \ep_{\angshattimesc}}_\ell \big)
\\
\counit_1 &= \counit_{(a_{\{j_\ell\}}, \glu_{\{j_\ell\}}),\, \angshat}
& \counit_2 &= \counit_{(a_{\{i_1\}}, \glu_{\{i_1\}}),\, \angshat} &&
\\
\glu_1 &= \glu_{\ep_{|s_{1,\ell}|};\, v_\ell, 1, \ang{\{j_\ell\}}_{j_\ell}} 
& \glu_2 &= \glu_{x;\, v, 1, \ang{s_{1,\ell}}_\ell}
& \glu_3 &= \glu_{\ep_{\cardsi};\, v, 1, \ang{\{i_1\}}_{i_1}}
\\
\glu_4 &= \glu_{x_2;\, v, 1, \ang{\ang{\{j_\ell\}}_{j_\ell}}_\ell}
& \glu_5 &= \glu_{x_2 \sig_1;\, v, 1, \ang{\{i_1\}}_{i_1}} && 
\end{aligned}$}
\right.\]
We abbreviate $\gaA$ \cref{gaAn}, $\ga$, and $\intr$ \cref{intr_jk} using concatenation, which yield, for example,
\[\begin{split}
x \ang{a_{v_\ell}}_\ell &= \gaA_r\big(x; \ang{a_{v_\ell}}_{\ell \in \ufs{r}} \big),\\
\ep_{|v^\stimes|}\ang{a_{\bdip}}_{\bdi} 
&= \gaA_{|v^\stimes|} \big(\ep_{|v^\stimes|} ; \ang{a_{\bdip}}_{\bdi \in v^\times} \big), \andspace\\
\ep_{|s_{1,\ell}|} \ep_{\angshattimesc} &= \ep_{|s_{1,\ell}|} \intr \ep_{\angshattimesc}.
\end{split}\]
In each level of $\Op$, which is a translation category, each unique isomorphism with given domain and codomain is denoted by $!$.  Using \cref{zdsgad_glu,cou_dom_obj,counit_as_gen}, the compatibility diagram \cref{counit_aglu_compat} is the boundary diagram in \cref{counit_compat_i}, where subscripts of $\vphi = \phiA$ are omitted.
\setlength{\fboxrule}{0.2pt}
\setlength{\fboxsep}{1.5pt}
\begin{equation}\label{counit_compat_i}

\end{equation}
The equalities in the diagram \cref{counit_compat_i} are given as follows.
\begin{itemize}
\item Each equality labeled $\mathbf{eq}$ holds by the action equivariance axiom \cref{pseudoalg_action_sym} for $\A$.  
\item The equality labeled $\mathbf{top}$ holds by the top equivariance axiom for the $\Gcat$-operad $\Op$.
\item In the region $\fbox{E'}$, the left equality labeled $\mathbf{as}$ holds by the associativity axiom for $\Op$ and the definition \cref{intr_jk} of the intrinsic pairing $\intr$:
\begin{equation}\label{xtwoeps_xthree}
\begin{split}
 x_2 \ep_{\angshattimesc}
&= \ga\big(\ga\big(x; \ang{\ep_{|s_{1,\ell}|}}_{\ell \in \ufs{r}}  \big); \ang{\ep_{\angshattimesc}}_{j_\ell \in s_{1,\ell},\, \ell \in \ufs{r}} \big)\\
&= \ga\big(x; \bang{\ga\big(\ep_{|s_{1,\ell}|} ; \ang{\ep_{\angshattimesc}}_{j_\ell \in s_{1,\ell}} \big) }_{\ell \in \ufs{r}} \big)\\
&= \ga\big(x; \ang{\ep_{|s_{1,\ell}|} \intr \ep_{\angshattimesc} }_{\ell \in \ufs{r}}  \big)\\
&= x_3.
\end{split}
\end{equation}
The right equality labeled $\mathbf{as}$ uses \cref{xtwoeps_xthree} and the permutation equality
\[\scalebox{.9}{$
\usi_{\angs, 1, \ang{s_{1,\ell}}_\ell} 
= \sig_{s_1, \ang{s_{1,\ell}}_\ell} \ang{\angshattimesc, \ldots, \angshattimesc}  
\cn \angstimes \fiso \txcoprod_{\ell \in \ufs{r}}\, (\angs \comp_1 s_{1,\ell})^\stimes$},\]
which holds by \cref{as_gamma,sigma_angsi,usi_perm}.  The object equality
\[\ang{a_{\bdip}}_{\bdi \in v^\times} 
= \ang{a_{i_1 \bdj}}_{\bdj \in \angshattimes,\, i_1 \in s_1}\]
holds by \cref{anga_bdip}.
\end{itemize}
The regions in the diagram \cref{counit_compat_i} commute for the following reasons.
\begin{itemize}
\item The region $\fbox{E'}$ commutes because it consists entirely of equalities.
\item The upper right region $\fbox{Tr}$ commutes by the fact that $\Op(|v^\stimes|)$ is a translation category and
the functoriality of $\gaA_{|v^\stimes|}$.
\item The regions $\fbox{F}$, $\fbox{N}$, $\fbox{E}$, $\fbox{T}$, and $\fbox{C}$ commute by, respectively, the functoriality of $\gaA_{\cardsi}$, the naturality of the associativity constraint $\phiA$ \cref{phiA}, the action equivariance axiom \cref{pseudoalg_action_sym}, the top equivariance axiom \cref{pseudoalg_topeq}, and the composition axiom \cref{pseudoalg_comp_axiom} for $\A$.
\item The bottom regions $\fbox{A}$, $\fbox{Eq}$, and $\fbox{Na}$ commute by, respectively, the associativity axiom \cref{system_associativity}, the equivariance axiom \cref{system_equivariance}, and the naturality axiom \cref{system_naturality} for the strong $\nbe$-system $(a,\glu)$.
\end{itemize}
This proves that the compatibility diagram \cref{counit_compat_i} commutes.
\item[General case with $k>1$]
Next, suppose $q  \geq 2$ and $k \geq 2$.  Similar to the notation in the diagram \cref{counit_compat_i}, we abbreviate $\gaA$, $\ga$, and $\intr$ \cref{intr_jk} using concatenation, and we abbreviate $\ang{\ang{\cdots}_p}_t$ to $\ang{\cdots}_{p,t}$.  We denote each isomorphism in each level of $\Op$ by $!$.  We omit subscripts of $\vphi = \phiA$.  Moreover, we use the following notation, along with \cref{as_gamma,eq:transpose_perm,Opn_object,usi_perm,bdip,s_prime,sone_angshat,aglu_ione}.
\[\left\{\scalebox{.9}{$
\begin{aligned}
\ang{\cdots}_{i_1} &= \ang{\cdots}_{i_1 \in s_1} 
& \ang{\cdots}_\ell &= \ang{\cdots}_{\ell \in \ufs{r}} 
& v_\ell &= \angs \compk s_{k,\ell}  
\\
v &= \angs 
& \ang{\cdots}_{\bdi} &= \ang{\cdots}_{\bdi \in v^\stimes} 
& \angshatl &= \angshat \comp_{k-1} \skell
\\
&
& \ang{\cdots}_{\bdj} &= \ang{\cdots}_{\bdj \in \angshattimes} 
& \ang{\cdots}_{\bdp_\ell} &= \ang{\cdots}_{\bdp_\ell \in \angshatltimes}  
\end{aligned}$}
\right.\]
\[\left\{\scalebox{.9}{$
\begin{aligned}
\usi &= \usi_{v,\, k, \ang{s_{k,\ell}}_{\ell}}  \phantom{M} 
& \usi' &= \usi_{\angshat,\, k-1, \ang{\skell}_\ell} \phantom{M} 
& \usi'_\stimes &= \txcoprod_{i_1 \in s_1} \usi'
\\
\twist &= \twist_{\cardsi,\, r} & \ol{\twist} &= \twist\ang{\angshatltimesc}_{\ell,\, i_1}
\end{aligned}$}
\right.\]
\[\left\{\scalebox{.9}{$
\begin{aligned}
\ep^v &= \ep_{|v^\stimes|} \phantom{M} 
& \ep^{v_\ell} &= \ep_{|v_\ell^\stimes|} \phantom{M} 
& \ep^{s_1} &= \ep_{\cardsi} &&
\\
\ep^\angshat &= \ep_{\angshattimesc} 
& \ep^{\angshatl} &= \ep_{\angshatltimesc} 
& \ep^{s_1 \angshat} &= \ep^{s_1} \ep^\angshat \phantom{M} 
& \ep^{s_1 \angshatl} &= \ep^{s_1} \ep^{\angshatl}
\\
\end{aligned}$}
\right.\]
\begin{equation}\label{wxy}
\left\{\scalebox{.8}{$
\begin{aligned}
x_1 &= x\ang{\ep^{v_\ell}}_{\ell} 
& x_2 &= x \ep^{s_1} 
& x_3 &= \ep^{s_1} x
\\ 
x_4 &= x\ang{\ep^{\angshatl}}_{\ell}
& x_5 &= (x\ang{\ep^{\angshatl}}_{\ell}) \usi'
& x_6 &= x\ang{\ep^{s_1 \angshatl}}_{\ell}
\\
x_7 &= (x \ep^{s_1}) \ang{\ep^{\angshatl}}_{i_1,\, \ell}
& x_8 &= (((x \ep^{s_1}) \ang{\ep^{\angshatl}}_{i_1,\, \ell}) \ol{\twist})\usi'_\stimes 
& x_9 &= (\ep^{s_1} (x\ang{\ep^{\angshatl}}_{\ell})) \usi'_\stimes
\\
x_{10} &= \ep^{s_1} ((x\ang{\ep^{\angshatl}}_{\ell}) \usi')
& x_{11} &= ((x \ep^{s_1}) \twist) \ang{\ep^{\angshatl}}_{\ell,\, i_1} 
& x_{12} &= (\ep^{s_1} x) \ang{\ep^{\angshatl}}_{\ell,\, i_1} 
\\
x_{13} &= \ep^{s_1} (x\ang{\ep^{\angshatl}}_{\ell}) &&&&
\end{aligned}$}
\right.
\end{equation}
\[\left\{\scalebox{.9}{$
\begin{aligned}
a_{i_1 \angshatl} &= a_{(\{i_1\},\, \angshatl)} \phantom{M} 
& a_{i_1 \bdp_\ell} &= a_{(\{i_1\}, \bdpp_\ell)} \phantom{M} 
& \ang{a_{i_1 \bdp_\ell}} &= \ang{a_{i_1 \bdp_\ell}}_{\bdp_\ell}
\\
a_{i_1 \angshat} &= a_{(\{i_1\},\, \angshat)} & a_{i_1 \bdj} &= a_{(\{i_1\},\, \bdjp)} 
& b_{i_1 \ell} &= \ep^{\angshatl} \ang{a_{i_1 \bdp_\ell}}
\end{aligned}$}
\right.\]
\[\left\{\scalebox{.9}{$
\begin{aligned}
\counit_1 &= \counit_{(a_{\{i_1\}}, \glu_{\{i_1\}}),\, \angshatl} \phantom{M} 
& \counit_2 &= \counit_{(a_{\{i_1\}}, \glu_{\{i_1\}}),\, \angshat}
\\
\glu_1 &= \glu_{\ep^{s_1};\, v_\ell, 1, \ang{\{i_1\}}_{i_1}} 
& \glu_2 &= \glu_{x;\, v,\, k, \ang{\skell}_\ell}
\\
\glu_3 &= \glu_{\ep^{s_1};\, v, 1, \ang{\{i_1\}}_{i_1}}
& \glu_4 &= \glu_{x;\, (\{i_1\},\, \angshat),\, k, \ang{\skell}_\ell} 
\end{aligned}$}
\right.\]
Among the objects in \cref{wxy}, there are equalities in $\Op(|v^\stimes|)$ 
\begin{equation}\label{x_sixseven}
x_6 = x_7 \andspace x_{12} = x_{13}
\end{equation}
by the associativity axiom for $\Op$, 
\begin{equation}\label{x_seveneleven}
x_{11} = x_7 \ol{\twist}
\end{equation}
by the top equivariance axiom for $\Op$, and
\begin{equation}\label{x_nineten}
x_9 = x_{10}
\end{equation}
by the bottom equivariance axiom for $\Op$.  

Using \cref{zdsgad_glu,cou_dom_obj,counit_as_gen}, the compatibility diagram \cref{counit_aglu_compat} is the boundary diagram in \cref{counit_compat_ii}.
\setlength{\fboxsep}{1.5pt}
\begin{equation}\label{counit_compat_ii}

\end{equation}
The equalities in the diagram \cref{counit_compat_ii} are given as follows.
\begin{itemize}
\item The equality labeled $\blacktriangle$ uses the object equalities
\[x_8 = x_7(\ol{\twist} \usi'_\stimes) = x_7 \usi \inspace \Op(|v^\stimes|).\]
\begin{itemize}
\item The first equality holds by the symmetric group action axiom for $\Op$.
\item The second equality holds because the following diagram of permutations in $\Si_{|v^\stimes|}$ commutes.
\begin{equation}\label{secondeq}
%
\end{equation}
\end{itemize}
The object equality
\[\ang{a_{\bdip}}_{\bdi \in v^\stimes} 
= \ang{a_{i_1 \bdj}}_{\bdj \in \angshattimes,\, i_1 \in s_1}\]
holds by \cref{anga_bdip}.
\item The equalities labeled $\mathbf{eq}$ hold by the action equivariance axiom \cref{pseudoalg_action_sym} for $\A$.
\item The equalities labeled $\mathbf{as}$, $\mathbf{t}$, and $\mathbf{b}$ hold by, respectively, the associativity axiom, the top equivariance axiom, and the bottom equivariance axiom for the $\Gcat$-operad $\Op$, as noted in \cref{x_sixseven,x_seveneleven,x_nineten}.
\end{itemize}
The regions in the diagram \cref{counit_compat_ii} commute for the following reasons.
\begin{itemize}
\item The region $\fbox{P}$ and the unlabeled region above it commute because each of them consists entirely of equalities.
\item The two regions $\fbox{Tr}$ commute by the assumption that $\Op(|v^\stimes|)$ is a translation category and the functoriality of $\gaA_{|v^\stimes|}$.
\item The regions $\fbox{F}$, $\fbox{N}$, $\fbox{E}$, $\fbox{T}$, $\fbox{B}$, and $\fbox{C}$ commute by, respectively, the functoriality of $\gaA_{|s_1|r}$, the naturality of the associativity constraint $\phiA$ \cref{phiA}, the action equivariance axiom \cref{pseudoalg_action_sym}, the top equivariance axiom \cref{pseudoalg_topeq}, the bottom equivariance axiom \cref{pseudoalg_boteq}, and the composition axiom \cref{pseudoalg_comp_axiom} for $\A$.
\item The bottom region $\fbox{Com}$ commutes by the commutativity axiom \cref{system_commutativity} for the strong $\nbe$-system $(a,\glu)$.
\item The lower right region $\fbox{IH}$ commutes by the induction hypothesis, since the object $\angshat = \ang{s_j}_{j=2}^q \in \GG$ has length $q-1$.
\end{itemize}
This proves that the compatibility diagram \cref{counit_compat_ii} commutes.\qedhere
\end{description}
\end{proof}

\section{Adjoint Equivalence}
\label{sec:hgo_prod_eq}

Throughout this section, we consider a $\Uinf$-operad $(\Op,\ga,\opu)$ \pcref{as:OpA'}, an $\Op$-pseudoalgebra $(\A,\gaA,\phiA)$ \pcref{def:pseudoalgebra}, and an object $\angordnbe = \ang{\ordn_j^{\be_j}}_{j \in \ufsq} \in \GG \setminus \{\vstar,\ang{}\}$ of length $q>0$ \cref{GG_objects}.  This section first proves that there is an adjoint equivalence \pcref{thm:zdsg_eq}
\begin{equation}\label{zdsg_adeq_seci}
\begin{tikzpicture}[vcenter]
\def\s{22}
\draw[0cell]
(0,0) node (a1) {\phantom{X}}
(a1)++(1.8,0) node (a2) {\phantom{X}}
(a1)++(-.15,0) node (a1') {\proAnbe}
(a2)++(.4,-.04) node (a2') {\Asgangordnbe}
node[between=a1 and a2 at .5] {\sim}
;
\draw[1cell=.9]
(a1) edge[bend left=\s] node {\zdsgad} (a2)
(a2) edge[bend left=\s] node {\zdsg} (a1)
;
\end{tikzpicture}
\end{equation}
between the twisted product $\proAnbe$ \pcref{def:proCnbe} and strong $H$-theory at $\angordnbe$ \pcref{A_ptfunctorGG}.  The left and right adjoints are the pointed functors $\zdsgad$ \cref{zdsgad} and $\zdsg$ \cref{hgosgtoprod}.  The unit is $\unit$ \cref{unit_zdsg}, and the counit is $\counit$ \cref{counit_zdsg}.  Combining this result with \cref{thm:zbsg_eq} proves that the strong $H$-theory comparison pointed $G$-functor \cref{PistsgAnbe}
\[\Asgsmaangordnbe \fto[\sim]{\Pistsg_{\A,\angordnbe}} \Asgangordnbe\]
is an equivalence of categories \pcref{thm:PistAequivalence}.  Thus, Shimakawa strong $H$-theory at $\sma\angordnbe$ is nonequivariantly equivalent to strong $H$-theory at $\angordnbe$.

\secoutline
\begin{itemize}
\item \cref{zdsg_left_triangle,zdsg_right_triangle} prove, respectively, the left triangle identity and the right triangle identity for the quadruple $(\zdsgad,\zdsg,\unit,\counit)$, proving \cref{thm:zdsg_eq}.
\item \cref{thm:PistAequivalence} proves that the strong $H$-theory comparison $\Pistsg_\A$ is componentwise an equivalence of categories.
\item \cref{expl:Pist_not_eq} discusses the necessity of using strong systems in \cref{thm:zdsg_eq,thm:PistAequivalence}.
\item \cref{expl:system_adjunction} discusses the variant adjunction
\begin{equation}\label{zd_adj_seci}
\begin{tikzpicture}[vcenter]
\def\s{20}
\draw[0cell]
(0,0) node (a1) {\phantom{X}}
(a1)++(1.8,0) node (a2) {\phantom{X}}
(a1)++(-.15,0) node (a1') {\proAnbe}
(a2)++(.3,-.04) node (a2') {\Aangordnbe}
;
\draw[0cell=.8]
node[between=a1 and a2 at .5] {\perp}
;
\draw[1cell=.9]
(a1) edge[bend left=\s] node {\zdad} (a2)
(a2) edge[bend left=\s] node {\zd} (a1)
;
\end{tikzpicture}
\end{equation}
that involves all $\angordnbe$-systems, not just the strong ones.
\end{itemize}

\subsection*{Triangle Identities}
For the quadruple $(\zdsgad,\zdsg,\unit,\counit)$, \cref{zdsg_left_triangle,zdsg_right_triangle} verify the two triangle identities that define an adjunction \pcref{def:adjunction}. 

\begin{lemma}\label{zdsg_left_triangle}
In the context of \cref{def:zdsgunit,def:zdsgcounit}, the unit \cref{unit_zdsg} 
\[1_{\proAnbe} \fto{\unit = 1_1} \zdsg\zdsgad\]
and the counit \cref{counit_zdsg} 
\[\zdsgad \zdsg \fto[\iso]{\counit} 1_{\Asgangordnbe}\]
satisfy the left triangle identity.
\end{lemma}

\begin{proof}
The left triangle identity states that, for each object 
\[a = \ang{a_{\bdi}}_{\bdi \in \ufs{n_1 n_2 \cdots n_q}} \in \proAnbe,\] 
the diagram \cref{zdsg_eq_lefttriangle} in $\Asgangordnbe$ commutes.  
\begin{equation}\label{zdsg_eq_lefttriangle}
\begin{tikzpicture}[vcenter]
\def\h{2.5}
\draw[0cell]
(0,0) node (a1) {\zdsgad a}
(a1)++(\h,0) node (a2) {\zdsgad\zdsg\zdsgad a}
(a2)++(\h,0) node (a3) {\zdsgad a}
;
\draw[1cell=.9]
(a1) edge node {\zdsgad \unit_a} (a2)
(a2) edge node {\counit_{\zdsgad a}} (a3)
;
\draw[1cell=.9]
(a1) [rounded corners=2pt] |- ($(a2)+(-1,.6)$) -- node {1} ($(a2)+(1,.6)$) -| (a3)
 ;
\end{tikzpicture}
\end{equation}
The morphism $\zdsgad \unit_a$ in \cref{zdsg_eq_lefttriangle} is the identity by the fact that $\unit_a = 1_a$ and the functoriality of $\zdsgad$ \cref{zdsgad}.

To show that the morphism $\counit_{\zdsgad a}$ in \cref{zdsg_eq_lefttriangle} is the identity, it suffices to show that, for each marker $\angs = \ang{s_j \subseteq \ufs{n}_j}_{j \in \ufs{q}}$, the $\angs$-component $\counit_{\zdsgad a, \angs}$ is the identity morphism in $\A$.  
\begin{description}
\item[Base case]
If $s_j = \emptyset$ for some $j \in \ufs{q}$, then, by \cref{cou_sj_empty}, 
\[\counit_{\zdsgad a, \angs} = 1_\zero.\]
Assuming $s_j \neq \emptyset$ for all $j \in \ufs{q}$, we prove that $\counit_{\zdsgad a, \angs}$ is the identity by an induction on $q>0$.
\item[Initial case]
If $q=1$, then $\angs = (s_1)$ and $\angstimes = s_1$.  The permutation \cref{usi_perm}
\[\big(s_1 \fto[\iso]{\usi_{(s_1);\, 1, \ang{\{i\}}_{i \in s_1}}} \txcoprod_{i \in s_1} \{i\}\big) \in \Si_{\cardsi}\]
is the identity.  Since $\ep_1 = \opu \in \Op(1)$, the right unity axiom for the $\Gcat$-operad $\Op$ implies
\[\ga\big(\ep_{\cardsi} ; \ang{\ep_1}_{i \in s_1}\big) \usi_{(s_1);\, 1, \ang{\{i\}}_{i \in s_1}} 
= \ga\big(\ep_{\cardsi} ; \ang{\opu}_{i \in s_1}\big) \id 
= \ep_{\cardsi}.\]
Thus, the unique isomorphism \cref{ual_iso}
\[\ga\big(\ep_{\cardsi} ; \ang{\ep_1}_{i \in s_1}\big) \usi_{(s_1);\, 1, \ang{\{i\}}_{i \in s_1}}
= \ep_{\cardsi} 
\fto[\iso]{\ual_{\ep_{\cardsi;\, (s_1),\, 1,\, \ang{\{i\}}_{i \in s_1}}}} \ep_{\cardsi}\]
in the translation category $\Op(\cardsi)$ is the identity.  Denoting by
\[\begin{split}
\usi &= \usi_{(s_1);\, 1, \ang{\{i\}}_{i \in s_1}} = \id \andspace\\
\ual &= \ual_{\ep_{\cardsi;\, (s_1), 1, \ang{\{i\}}_{i \in s_1}}} = 1_{\ep_{\cardsi}}
\end{split}\]
and using \cref{zdsgad_comp_obj,zdsgad_glu,cou_dom_obj,cou_a_s}, the $\angs$-component $\counit_{\zdsgad a, \angs}$ is the following composite in $\A$, where $\ang{\Cdots}_i = \ang{\Cdots}_{i \in s_1}$.
\begin{equation}\label{zdsg_ltri_initial}
\begin{tikzpicture}[vcenter]
\def\h{4.5} \def\u{-1}
\draw[0cell=.85]
(0,0) node (a11) {(\zdsgad\zdsg\zdsgad a)_{\angs}}
(a11)++(\h,0) node (a12) {(\zdsgad a)_{\angs}}
(a11)++(0,\u) node (a21) {\gaA_{\cardsi}\big(\ep_{\cardsi}; \ang{(\zdsgad a)_{\{i\}}}_{i} \big)}
(a12)++(0,\u) node (a22) {\gaA_{\cardsi}\big(\ep_{\cardsi}; \ang{a_i}_{i} \big)}
(a21)++(0,\u) node (a31) {\gaA_{\cardsi}\big(\ep_{\cardsi}; \bang{\gaA_1(\opu; a_i)}_{i} \big)}
(a31)++(0,-1.5) node (a41) {\gaA_{\cardsi}\big(\ga(\ep_{\cardsi}; \ang{\opu}_{i}) ; \ang{a_i}_{i} \big)}
(a41)++(\h,0) node (a42) {\gaA_{\cardsi}\big(\ga(\ep_{\cardsi}; \ang{\opu}_i) \usi ; \ang{a_i}_i \big)}
;
\draw[1cell=.8]
(a11) edge node {\counit_{\zdsgad a, \angs}} (a12)
(a11) edge[equal] (a21)
(a21) edge[equal] (a31)
(a31) edge node[swap] {\phiA} (a41)
(a41) edge[equal] node{\mathbf{eq}} (a42)
(a42) edge node[swap] {\gaA_{\cardsi}(\ual; 1^{\cardsi})} (a22)
(a22) edge[equal] (a12)
;
\end{tikzpicture}
\end{equation}
In the diagram \cref{zdsg_ltri_initial}, the associativity constraint $\phiA = \phiA_{(\cardsi;\, 1^{\cardsi})}$ is the identity by the unity axiom \cref{pseudoalg_unity} for $\A$.  The arrow $\gaA_{\cardsi}(\ual; 1^{\cardsi})$ is the identity by the functoriality of $\gaA_{\cardsi}$ and the fact that $\ual = 1_{\ep_{\cardsi}}$.  This proves that $\counit_{\zdsgad a, \angs}$ is the identity when $q=1$.
\item[General case]
Inductively, suppose $q \geq 2$.  We use the following notation, along with \cref{intr_jk,s_prime,sone_angshat}.
\[\left\{
\begin{gathered}
\begin{aligned}
v &= \angs & \ang{\Cdots}_{\bdi} &= \ang{\Cdots}_{\bdi \in v^\stimes} & \ang{\Cdots}_{i_1} &= \ang{\Cdots}_{i_1 \in s_1}\\
&& \ang{\Cdots}_{\bdj} &= \ang{\Cdots}_{\bdj \in \angshattimes} & \ep^{s_1} &= \ep_{\cardsi}\\
\ep^v &= \ep_{|v^\stimes|} & \ep^\angshat &= \ep_{\angshattimesc} & \ep^{s_1 \angshat} &= \ep^{s_1} \intr \ep^\angshat 
\end{aligned}\\
\usi = \usi_{v, 1, \ang{\{i_1\}}_{i_1}} \cn v^\stimes \fiso \txcoprod_{i_1 \in s_1} (\{i_1\} \times \angshattimes)
\end{gathered}
\right.\]
The permutation $\usi$ \cref{usi_perm} is the identity because each of $\angshattimes$ and $v^\stimes = s_1 \times \angshattimes$ has the lexicographic ordering.  For each $q$-tuple $\bdi = \ang{i_j}_{j \in \ufs{q}} \in \txprod_{j \in \ufs{q}}\, \ufs{n}_j$, recall from \cref{bdip} the $q$-tuple of one-element subsets
\[\bdip = \big(\{i_1\}, \{i_2\}, \ldots, \{i_q\}\big).\]
Since
\[\bdip^\stimes = \txprod_{j \in \ufs{q}}\, \{i_j\} = \{\bdi\} \andspace |\bdip^\stimes| = 1,\]
by \cref{pseudoalg_action_unity,i_connected,zdsgad_comp_obj}, there are object equalities
\begin{equation}\label{zdsgad_a_bdip}
(\zdsgad a)_{\bdip} = \gaA_1(\ep_1; a_\bdi) = \gaA_1(\opu; a_\bdi) = a_\bdi \inspace \A.
\end{equation}
By \cref{zdsgad_comp_obj,zdsgad_glu,cou_dom_obj,counit_as_gen,zdsgad_a_bdip}, the $\angs$-component $\counit_{\zdsgad a, \angs}$ is the following composite in $\A$, where $\ang{\Cdots}_{\bdj,\, i_1} = \ang{\ang{\Cdots}_{\bdj}}_{i_1}$.  
\begin{equation}\label{zdsg_ltri_diagram}
\begin{tikzpicture}[vcenter]
\def\u{-1} \def\v{-1.4} \def\h{4.7} \def\w{-1}
\draw[0cell=.85]
(0,0) node (a11) {(\zdsgad \zdsg \zdsgad a)_{\angs}}
(a11)++(\h,0) node (a12) {(\zdsgad a)_{\angs}}
(a11)++(0,\u) node (a21) {\gaA_{|v^\stimes|}\big(\ep^v; \ang{(\zdsgad a)_{\bdip}}_{\bdi} \big)}
(a12)++(0,\u) node (a22) {\gaA_{|v^\stimes|}\big(\ep^v; \ang{a_\bdi}_{\bdi} \big)}
(a21)++(0,\u) node (a31) {\gaA_{|v^\stimes|}\big(\ep^v; \ang{a_\bdi}_{\bdi} \big)}
(a22)++(0,\v) node (a32) {\gaA_{|v^\stimes|}\big(\ep^{s_1 \angshat} \usi ; \ang{a_\bdi}_{\bdi} \big)}
(a31)++(0,\v) node (a41) {\gaA_{|v^\stimes|}\big(\ep^{s_1 \angshat} ; \ang{a_{(i_1, \bdj)}}_{\bdj, i_1}  \big)}
(a32)++(0,\u) node (a42) {\gaA_{|v^\stimes|}\big(\ep^{s_1 \angshat} ; \ang{a_{(i_1, \bdj)}}_{\bdj, i_1}  \big)}
(a41)++(0,\v) node (a51) {\gaA_{\cardsi}\big(\ep^{s_1}; \bang{\gaA_{\angshattimesc}\big(\ep^\angshat ; \ang{a_{(i_1, \bdj)}}_{\bdj} \big) }_{i_1} \big)}
(a42)++(0,\v) node (a52) {\gaA_{\cardsi}\big(\ep^{s_1}; \bang{\gaA_{\angshattimesc}\big(\ep^\angshat ; \ang{a_{(i_1, \bdj)}}_{\bdj} \big) }_{i_1} \big)}
(a51)++(\h/2,0) node (a5) {}
;
\draw[1cell=.8]
(a11) edge node {\counit_{\zdsgad a, \angs}} (a12)
(a11) edge[equal] (a21)
(a21) edge[equal] (a31)
(a31) edge node[swap] {\gaA_{|v^\stimes|}(!\,; 1^{|v^\stimes|})} (a41)
(a41) edge node[swap] {(\phiA)^{-1}} (a51)
(a52) edge node[swap] {\phiA} (a42)
(a42) edge[equal, shorten >=-1ex] node[swap] {\mathbf{eq}} (a32)
(a32) edge node[swap] {\gaA_{|v^\stimes|}(!\,; 1^{|v^\stimes|})} (a22)
(a22) edge[equal] (a12)
(a51) [rounded corners=2pt] |- ($(a5)+(-1,\w)$)
-- node {\gaA_{\cardsi} (1 ; \ang{\counit_{(\zdsgad a)_{\{i_1\}},\, \angshat}}_{i_1})} ($(a5)+(1,\w)$) -| (a52)
;
\end{tikzpicture}
\end{equation}
In the bottom arrow in the diagram \cref{zdsg_ltri_diagram}, $(\zdsgad a)_{\{i_1\}}$ denotes the $i_1$-restriction of $\zdsgad a$ \cref{aglu_ione}.  Since the object $\angshat$ \cref{s_prime} has length $q-1$, the induction hypothesis implies that $\counit_{(\zdsgad a)_{\{i_1\}},\, \angshat}$ is the identity for each $i_1 \in s_1$.  The functoriality of $\gaA_{\cardsi}$ implies that the bottom arrow in \cref{zdsg_ltri_diagram} is the identity.  

Since $\usi = \id$ and since $\Op(|v^\stimes|)$ is a translation category, the isomorphism in the upper right of \cref{zdsg_ltri_diagram},
\[\ep^{s_1 \angshat}\usi = \ep^{s_1 \angshat} \fto{!} \ep^v,\]
is the inverse of the isomorphism
\[\ep^v \fto{!} \ep^{s_1 \angshat}\]
in the middle left of \cref{zdsg_ltri_diagram}.  Thus, the long composite in \cref{zdsg_ltri_diagram}, which is $\counit_{\zdsgad a, \angs}$, is the identity.  
\end{description}
This finishes the induction and proves the left triangle identity \cref{zdsg_eq_lefttriangle}.
\end{proof}

\begin{lemma}\label{zdsg_right_triangle}
In the context of \cref{def:zdsgunit,def:zdsgcounit}, the unit \cref{unit_zdsg} 
\[1_{\proAnbe} \fto{\unit = 1_1} \zdsg\zdsgad\]
and the counit \cref{counit_zdsg} 
\[\zdsgad \zdsg \fto[\iso]{\counit} 1_{\Asgangordnbe}\]
satisfy the right triangle identity.
\end{lemma}

\begin{proof}
 The right triangle identity states that, for each object $(a,\glu) \in \Asgangordnbe$, the diagram \cref{zdsg_eq_righttriangle} in $\proAnbe$ commutes.
\begin{equation}\label{zdsg_eq_righttriangle}
\begin{tikzpicture}[vcenter]
\def\h{3}
\draw[0cell]
(0,0) node (a1) {\zdsg(a,\glu)}
(a1)++(\h,0) node (a2) {\zdsg\zdsgad\zdsg (a,\glu)}
(a2)++(\h,0) node (a3) {\zdsg (a,\glu)}
;
\draw[1cell=.9]
(a1) edge node {\unit_{\zdsg(a,\glu)}} (a2)
(a2) edge node {\zdsg \counit_{(a,\glu)}} (a3)
;
\draw[1cell=.9]
(a1) [rounded corners=2pt] |- ($(a2)+(-1,.7)$) -- node {1} ($(a2)+(1,.7)$) -| (a3)
 ;
\end{tikzpicture}
\end{equation}
The morphism $\unit_{\zdsg(a,\glu)}$ is the identity by definition \cref{unit_zdsg}.

To show that $\zdsg \counit_{(a,\glu)}$ is the identity morphism, we show that, for each index $\bdi = \ang{i_j}_{j \in \ufs{q}} \in \txprod_{j \in \ufs{q}}\, \ufs{n}_j$, the $\bdi$-th coordinate morphism of $\zdsg \counit_{(a,\glu)}$ is the identity.  By \cref{hgotoprod_def}, this $\bdi$-th coordinate morphism is
\[\big(\zdsg \counit_{(a,\glu)}\big)_{\bdi} 
= \counit_{(a,\glu), \bdip},\]
where $\bdip = \ang{\{i_j\}}_{j \in \ufs{q}}$ is the $q$-tuple of one-element subsets \cref{bdip}.  We prove that this is the identity morphism by an induction on $q>0$.
\begin{description}
\item[Initial case]  
If $q=1$, then \cref{i_connected,cou_a_s,system_unity_iii} imply that
\[\counit_{(a,\glu), \bdip} = \glu_{\ep_1;\, \{i_1\}, 1, \{i_1\}} = \glu_{\opu;\, \{i_1\}, 1, \{i_1\}} = 1.\]
\item[General case]  
Inductively, suppose $q \geq 2$.  By \cref{counit_as_gen}, $\counit_{(a,\glu), \bdip}$ is the composite of four isomorphisms.  We observe that each of these four constituent arrows is the identity.
\begin{enumerate}
\item Since $\bdip^\stimes = \{\bdi\}$, the first constituent arrow of $\counit_{(a,\glu), \bdip}$ is
\[\gaA_{|\bdip^\stimes|}\big(!\, ; 1^{|\bdip^\stimes|}\big) = \gaA_1(!\, ; 1),\] 
where the unique isomorphism is given by
\[\ep_{|\bdip^\stimes|} = \ep_1 = \opu \fto{!} \ep_1 \intr \ep_1 = \opu \intr \opu = \opu\]
in $\Op(1) = \{\opu\}$.  Since this unique isomorphism is the identity morphism of $\opu$, the first constituent arrow of $\counit_{(a,\glu), \bdip}$ is the identity by the functoriality of $\gaA_1$.
\item The second constituent arrow of $\counit_{(a,\glu), \bdip}$ is $(\phiA_{(1;1)})^{-1}$.  This is the identity by the unity axiom \cref{pseudoalg_unity} for $\A$.
\item The third constituent arrow of $\counit_{(a,\glu), \bdip}$ is 
\[\gaA_1\big(1; \counit_{(a_{\{i_1\}}, \glu_{\{i_1\}}) ,\, \what{\bdi}} \big),\]
where
\begin{itemize}
\item $\what{\bdi} = \ang{\{i_j\}}_{j=2}^q$ and
\item $(a_{\{i_1\}}, \glu_{\{i_1\}})$ is the $i_1$-restriction of $(a,\glu)$ \cref{aglu_ione}.
\end{itemize}
Since $\what{\bdi}$ has length $q-1$, the induction hypothesis implies that $\counit_{(a_{\{i_1\}}, \glu_{\{i_1\}}) ,\, \what{\bdi}}$ is the identity.  Thus, the functoriality of $\gaA_1$ implies that the third constituent arrow of $\counit_{(a,\glu), \bdip}$ is the identity.
\item The fourth constituent arrow of $\counit_{(a,\glu), \bdip}$ is 
\[\glu_{\ep_1;\, \bdip, 1, \{i_1\}} = \glu_{\opu;\, \bdip, 1, \{i_1\}}.\]
This is the identity by the second unity axiom \cref{system_unity_iii} for $(a,\glu)$.
\end{enumerate}
\end{description}
This finishes the induction and proves the right triangle identity \cref{zdsg_eq_righttriangle}.
\end{proof}

\begin{theorem}\label{thm:zdsg_eq}
For a $\Uinf$-operad $(\Op,\ga,\opu)$ \pcref{as:OpA'}, an $\Op$-pseudoalgebra $(\A,\gaA,\phiA)$ \pcref{def:pseudoalgebra}, and an object $\angordnbe = \ang{\ordn_j^{\be_j}}_{j \in \ufsq} \in \GG \setminus \{\vstar,\ang{}\}$, there is an adjoint equivalence of categories
\begin{equation}\label{zdsg_thm}
\begin{tikzpicture}[vcenter]
\def\s{22}
\draw[0cell]
(0,0) node (a1) {\phantom{X}}
(a1)++(1.8,0) node (a2) {\phantom{X}}
(a1)++(-.15,0) node (a1') {\proAnbe}
(a2)++(.4,-.04) node (a2') {\Asgangordnbe}
node[between=a1 and a2 at .5] {\sim}
;
\draw[1cell=.9]
(a1) edge[bend left=\s] node {\zdsgad} (a2)
(a2) edge[bend left=\s] node {\zdsg} (a1)
;
\end{tikzpicture}
\end{equation}
given by the following data.
\begin{itemize}
\item $\proAnbe$ is the $\angordnbe$-twisted product \pcref{def:proCnbe}.
\item $\Asgangordnbe$ is the category of strong $\angordnbe$-systems \pcref{def:nsystem,def:nsystem_morphism,def:nbeta_gcat}.
\item The left adjoint is the pointed functor $\zdsgad$ \cref{zdsgad}.
\item The right adjoint is the pointed functor $\zdsg$ \cref{hgosgtoprod}.
\item The unit is the identity natural transformation \cref{unit_zdsg}
\[1_{\proAnbe} \fto{\unit = 1_1} \zdsg\zdsgad.\]
\item The counit is the natural isomorphism \cref{counit_zdsg}
\[\zdsgad \zdsg \fto[\iso]{\counit} 1_{\Asgangordnbe}.\]
\end{itemize}
\end{theorem}

\begin{proof}
The unit $\unit$ is the identity natural transformation, and the counit $\counit$ is a natural isomorphism.   The left and right triangle identities for an adjunction are verified in \cref{zdsg_left_triangle,zdsg_right_triangle}.
\end{proof}

\subsection*{Strong $H$-Theory Comparison is an Equivalence}

Next, we observe that the strong $H$-theory comparison pointed $G$-functor \cref{PistsgAnbe}
\[\Asgsmaangordnbe \fto{\Pistsg_{\A,\angordnbe}} \Asgangordnbe\]
is an equivalence of categories.  Thus, Shimakawa strong $H$-theory is nonequivariantly equivalent to strong $H$-theory.

\begin{theorem}\label{thm:PistAequivalence}\index{H-theory comparison@$H$-theory comparison!equivalence}
For a $\Uinf$-operad $(\Op,\ga,\opu)$ \pcref{as:OpA'} and an $\Op$-pseudoalgebra $(\A,\gaA,\phiA)$ \pcref{def:pseudoalgebra}, the strong $H$-theory comparison $G$-natural transformation \cref{PistsgA}
\begin{equation}\label{PistsgA_thm}
\begin{tikzpicture} [vcenter]
\def\s{25}
\draw[0cell]
(0,0) node (a1) {\GG}
(a1)++(2.5,0) node (a2) {\phantom{\GG}}
(a2)++(.15,0) node (a2') {\Catgst}
;
\draw[1cell=.9]
(a1) edge[bend left=\s] node {\smast \Sgosg\A} (a2)
(a1) edge[bend right=\s] node[swap] {\Hgosg\A} (a2) 
;
\draw[2cell]
node[between=a1 and a2 at .4, rotate=-90, 2label={above,\Pistsg_\A}] {\Rightarrow}
;
\end{tikzpicture}
\end{equation}
is componentwise an equivalence of categories.
\end{theorem}

\begin{proof}
By \cref{PistA_vstar,PistA_empty}, $\Pistsg_{\A,\vstar} = 1_\boldone$, and $\Pistsg_{\A,\ang{}}$ is a pointed $G$-isomorphism.  For each object $\angordnbe \in \GG \setminus \{\vstar,\ang{}\}$, there is a commutative diagram \cref{Pistsg_prod_diag}
\begin{equation}\label{PistsgA_proof}
\begin{tikzpicture}[vcenter]
\def\h{3.5}
\draw[0cell]
(0,0) node (a11) {\Asgsmaangordnbe}
(a11)++(\h,0) node (a12) {\Asgangordnbe}
(a11)++(\h/2,-1.1) node (a2) {\phantom{\proAnbe}}
(a2)++(0,.15) node (a2') {\proAnbe}
;
\draw[1cell=.9]
(a11) edge node {\Pistsg_{\A,\angordnbe}} (a12)
(a11) edge node[swap] {\zbsg} (a2)
(a12) edge node {\zdsg} (a2)
;
\end{tikzpicture}
\end{equation}
of pointed $G$-functors.  The functors $\zbsg$ and $\zdsg$ are equivalences of categories by, respectively, \cref{thm:zbsg_eq,thm:zdsg_eq}.  Thus, $\Pistsg_{\A,\angordnbe}$ is also an equivalence of categories.
\end{proof}

\begin{explanation}[Necessity of Strong Systems]\label{expl:Pist_not_eq}
Both \cref{thm:zdsg_eq,thm:PistAequivalence} involve \emph{strong} $\angordnbe$-systems, with invertible gluing morphisms.  The invertibility of the gluing morphisms are used in \cref{cou_a_s,counit_as_gen} to make sure that the counit $\counit$ is componentwise an isomorphism of $\angordnbe$-systems.  \cref{expl:system_adjunction} discusses an adjunction analogous to the one in \cref{thm:zdsg_eq} that uses the larger category $\Aangordnbe$ of all $\angordnbe$-systems.  That variant adjunction is not generally an adjoint equivalence.  

On the other hand, there is no analogue of \cref{thm:PistAequivalence} for the $G$-natural transformation \cref{PistA}
\begin{equation}\label{PistA_expl}
\begin{tikzpicture}[vcenter]
\def\s{25}
\draw[0cell]
(0,0) node (a1) {\GG}
(a1)++(2.5,0) node (a2) {\phantom{\GG}}
(a2)++(.15,0) node (a2') {\Catgst}
;
\draw[1cell=.9]
(a1) edge[bend left=\s] node {\smast \Sgo\A} (a2)
(a1) edge[bend right=\s] node[swap] {\Hgo\A} (a2) 
;
\draw[2cell]
node[between=a1 and a2 at .4, rotate=-90, 2label={above,\Pist_\A}] {\Rightarrow}
;
\end{tikzpicture}
\end{equation}
that compares Shimakawa $H$-theory and $H$-theory.  For an object $\angordnbe \in \GG \setminus \{\vstar,\ang{}\}$, there is still an equality \cref{Pist_prod_diag}
\[\zb = \zd \circ \Pist_{\A,\angordnbe}.\]
However, $\zb$ and $\zd$ are not equivalences in general \pcref{expl:sys_adjunction,expl:system_adjunction}, so neither is $\Pist_{\A,\nbe}$.
\end{explanation}

\begin{explanation}[Systems and Adjunction]\label{expl:system_adjunction}
Analogous to the discussion in \cref{expl:sys_adjunction}, the unit $\unit$ and the counit $\counit$ \pcref{def:zdsgunit,def:zdsgcounit} are still defined even if the category $\Asgangordnbe$ is replaced by the larger category $\Aangordnbe$ of all $\angordnbe$-systems \pcref{def:nsystem,def:nsystem_morphism,def:nbeta_gcat}.  
\begin{description}
\item[Unit] There is an identity natural transformation
\begin{equation}\label{unit_zd}

\end{equation}
of the identity functor on $\proAnbe$.  The functor $\zd$ is the pointed $G$-functor defined in \cref{hgotoprod}.  The pointed functor $\zdad$ denotes the composite
\begin{equation}\label{zdad}
%
\end{equation}
of the pointed functor $\zdsgad$ \cref{zdsgad} and the full subcategory inclusion $\iota$.  The computation  \cref{uni_cod_obj} shows that this variant of the unit $\unit = 1_1$ is also well defined.
\item[Counit] There is a pointed natural transformation
\begin{equation}\label{counit_zd}
%
\end{equation}
whose components are defined as in \cref{cou_sj_empty,cou_a_s,counit_as_gen}.  This is well defined by the computation \cref{cou_dom_obj} and the proof of \cref{counit_agl_welldef}, applied to $\angordnbe$-systems.  By \cref{cou_a_s,counit_as_gen}, if $s_j \neq \emptyset$ for all $j \in \ufs{q}$, then the $\angs$-component $\counit_{(a,\glu),\angs}$ is a gluing morphism or involves at least two gluing morphisms of $(a,\glu) \in \Aangordnbe$.  Thus, $\counit_{(a,\glu),\angs}$ is not generally an isomorphism.  This variant of the counit $\counit$ is a pointed natural transformation, but not a natural isomorphism.
\end{description}
The proofs of the triangle identities in \cref{zdsg_left_triangle,zdsg_right_triangle} are still valid in this setting.  Thus, there is an adjunction, but not generally an adjoint equivalence,
\begin{equation}\label{nbe_system_adj}
\begin{tikzpicture}[vcenter]
\def\s{20}
\draw[0cell]
(0,0) node (a1) {\phantom{X}}
(a1)++(1.8,0) node (a2) {\phantom{X}}
(a1)++(-.15,0) node (a1') {\proAnbe}
(a2)++(.3,-.04) node (a2') {\Aangordnbe}
;
\draw[0cell=.8]
node[between=a1 and a2 at .5] {\perp}
;
\draw[1cell=.9]
(a1) edge[bend left=\s] node {\zdad} (a2)
(a2) edge[bend left=\s] node {\zd} (a1)
;
\end{tikzpicture}
\end{equation}
whose unit and counit are the natural transformations in \cref{unit_zd,counit_zd}.  \cref{expl:zdad_pseudo} discusses the pseudo $G$-equivariance of the left adjoint $\zdad$.
\end{explanation}

\section{Pseudo $G$-Equivariance of Left Adjoint}
\label{sec:pseudoequiv_ii}

In the adjoint equivalence 
\begin{equation}\label{zdsgadj_seci}
\begin{tikzpicture}[vcenter]
\def\s{22}
\draw[0cell]
(0,0) node (a1) {\phantom{X}}
(a1)++(1.8,0) node (a2) {\phantom{X}}
(a1)++(-.15,0) node (a1') {\proAnbe}
(a2)++(.4,-.04) node (a2') {\Asgangordnbe}
node[between=a1 and a2 at .5] {\sim}
;
\draw[1cell=.9]
(a1) edge[bend left=\s] node {\zdsgad} (a2)
(a2) edge[bend left=\s] node {\zdsg} (a1)
;
\end{tikzpicture}
\end{equation}
in \cref{thm:zdsg_eq}, the right adjoint $\zdsg$ is a pointed $G$-functor by \cref{Pist_to_prod} \eqref{Pist_to_prod_i}.  The main result of this section, \cref{thm:zdsgad_pseudo}, proves that the left adjoint $\zdsgad$ is pseudo $G$-equivariant \pcref{def:pseudoG}.  Thus, the lack of $G$-equivariance of $\zdsgad$ is controlled by its pseudo $G$-equivariant constraints, which satisfy two coherence axioms of their own.  

\secoutline
\begin{itemize}
\item \cref{def:zdsgad_pseudoG} defines the pseudo $G$-equivariant constraints for the functor $\zdsgad$.  The definition of its components \cref{zdp_gas} shows why $\zdsgad$ is not generally $G$-equivariant.  The corresponding components of $\zdsgad g$ and $g\zdsgad$ differ by an isomorphism that is generally not the identity.
\item \cref{ex:ual_g_angs} illustrates the nontriviality of the pseudo $G$-equivariant constraints for $\zdsgad$ when $\Op$ is either the Barratt-Eccles $\Gcat$-operad $\BE$ or the $G$-Barratt-Eccles operad $\GBE$. 
\item \cref{expl:zdp_zbp} discusses the formal similarity between the pseudo $G$-equivariant constraints for $\zbsgad$ and $\zdsgad$.
\item \cref{zdpg_a_welldef} proves that the components of the pseudo $G$-equivariant constraints for $\zdsgad$ are well-defined morphisms of strong $\angordnbe$-systems.
\item \cref{thm:zdsgad_pseudo} proves that $\zdsgad$, equipped with the pseudo $G$-equivariant constraints in \cref{def:zdsgad_pseudoG}, is a pseudo $G$-equivariant functor.
\item \cref{expl:zdad_pseudo} discusses the pseudo $G$-equivariance of the pointed functor \cref{zdad}
\[\proAnbe \fto{\zdad = \iota \zdsgad} \Aangordnbe.\]
\end{itemize}

\subsection*{Constraints for $\zdsgad$}
First, we define the pseudo $G$-equivariant constraints for the pointed functor \pcref{def:zdsgad}
\[\proAnbe \fto{\zdsgad} \Asgangordnbe,\]
which is part of an adjoint equivalence of categories \pcref{thm:zdsg_eq}. 
Recall that $\Asgangordnbe$ and $\proAnbe$ \pcref{def:proCnbe,def:nbeta_gcat} are $G$-categories.
\cref{def:zdsgad_pseudoG} is analogous to \cref{def:zbsgad_pseudoG}, which defines the pseudo $G$-equivariant constraints for $\zbsgad$.

\begin{definition}\label{def:zdsgad_pseudoG}
Given a $\Uinf$-operad $(\Op,\ga,\opu)$ \pcref{as:OpA'}, an $\Op$-pseudoalgebra $(\A,\gaA,\phiA)$ \pcref{def:pseudoalgebra}, an object $\angordnbe = \ang{\ordn_j^{\be_j}}_{j \in \ufsq} \in \GG \setminus \{\vstar,\ang{}\}$ of length $q>0$ \cref{GG_objects}, the pointed functor \cref{zdsgad}
\[\proAnbe \fto{\zdsgad} \Asgangordnbe,\]
and an element $g \in G$, we define a natural isomorphism
\begin{equation}\label{zdsgad_constraint}
\begin{tikzpicture}[vcenter]
\def\v{-1.5}
\draw[0cell=.9]
(0,0) node (a1) {\proAnbe}
(a1)++(2.8,0) node (b1) {\proAnbe}
(a1)++(0,\v) node (a2) {\Asgangordnbe}
(b1)++(0,\v) node (b2) {\Asgangordnbe}
;
\draw[1cell=.9]
(a1) edge node {g} (b1)
(a2) edge node[swap] {g} (b2)
(a1) edge node[swap] {\zdsgad} (a2)
(b1) edge node {\zdsgad} (b2)
;
\draw[2cell]
node[between=a1 and b2 at .55, rotate=-135, 2labelalt={below,\zdpg}] {\Rightarrow}
;
\end{tikzpicture}
\end{equation}
as follows.  For each marker $\angs = \ang{s_j \subseteq \ufs{n}_j}_{j \in \ufs{q}}$, we recall the marker \cref{ginvs}
\[\ginv\angs = \ang{\ginv s_j \subseteq \ufs{n}_j}_{j \in \ufs{q}}.\]
We first define a permutation $\usi^{g,\angs} \in \Si_{\angstimesc}$, with $\angstimes = \txprod_{j \in \ufs{q}}\, s_j$, and an isomorphism $\ual^{g,\angs} \in \Op(\angstimesc)$ as follows.  The components of $\zdpg$ are defined in \cref{zdp_gas} after these preliminary constructions.
\begin{itemize}
\item We define the permutation
\begin{equation}\label{usi_g_angs}
\big(\ginvangstimes \fto[\iso]{\usi^{g,\angs}} \angstimes\big) \in \Si_{|\angstimes|}
\end{equation}
given by
\[\begin{split}
\ginv\bdi = \ang{\ginv i_j}_{j \in \ufs{q}} & \in \ginvangstimes = \txprod_{j \in \ufs{q}}\, \ginv s_j = \ginv\angstimes\\
\mapsto \bdi = \ang{i_j}_{j \in \ufs{q}} & \in \angstimes.
\end{split}\] 
Each of $\angstimes$ and $\ginvangstimes$ is equipped with the lexicographic ordering inherited from $\txprod_{j \in \ufs{q}}\, \ufs{n}_j$, on which $g$ acts diagonally \cref{al-sma-be}.  In other words, $\usi^{g,\angs}$ is the $g$-action on $\txprod_{j \in \ufs{q}}\, \ufs{n}_j$ restricted to the subset $\ginvangstimes$.
\item Using \cref{Opn_object,usi_g_angs}, we define the unique isomorphism
\begin{equation}\label{ual_g_angs}
\ginv \ep_{\angstimesc} \usi^{g,\angs} \fto[\iso]{\ual^{g,\angs}} 
\ep_{\ginvangstimesc} = \ep_{\angstimesc}
\end{equation}
in the translation category $\Op(\angstimesc)$ with the indicated domain and codomain.   Since $\Op$ is a $\Gcat$-operad, the symmetric group action $- \cdot \usi^{g,\angs}$ on $\Op(\angstimesc)$ is $G$-equivariant.  Thus, the object $\ginv \ep_{\angstimesc} \usi^{g,\angs}$ is unambiguous.
\end{itemize}
\begin{description}
\item[Component isomorphisms] 
For each object \cref{proCnbe_object}
\[a = \ang{a_{\bdi}}_{\bdi \in \ufs{n_1 n_2 \cdots n_q}} \in \proAnbe,\]
the $a$-component of $\zdps$ is the isomorphism of strong $\angordnbe$-systems \pcref{def:nsystem_morphism}
\begin{equation}\label{zdp_g_a}
(\zdsgad g)a \fto[\iso]{\zdpg_a} (g \zdsgad) a
\end{equation}
with $\angs$-component isomorphism defined by the following commutative diagram in $\A$.
\begin{equation}\label{zdp_gas}
\begin{tikzpicture}[vcenter]
\def\v{-1.2}
\draw[0cell=.8]
(0,0) node (a1) {\big((\zdsgad g)a\big)_{\angs}}
(a1)++(0,\v) node (a2) {\gaA_{\angstimesc}\big(\ep_{\angstimesc}\sscs \ang{(ga)_{\bdi}}_{\bdi \in \angstimes} \big)}
(a2)++(0,\v) node (a3) {\gaA_{\angstimesc}\big(\ep_{\angstimesc}\sscs \ang{ga_{\ginv\bdi}}_{\bdi \in \angstimes}\big)}
(a3)++(0,\v) node (a4) {g\gaA_{\angstimesc}\big(\ginv \ep_{\angstimesc}\sscs \ang{a_{\ginv \bdi}}_{\bdi \in \angstimes}\big)}
(a4)++(0,\v) node (a5) {g\gaA_{\angstimesc}\big(\ginv \ep_{\angstimesc} \usi^{g,\angs} \sscs \ang{a_{\bdj}}_{\bdj \in \ginvangstimes}\big)}
(a1)++(5.5,0) node (b1) {\big((g \zdsgad) a\big)_{\angs}}
(b1)++(0,\v) node (b2) {g(\zdsgad a)_{\ginv \angs}}
(b2)++(0,\v) node (b3) {\phantom{\gaA_{\angstimes}}}
(b3)++(-1,0) node (b3') {g\gaA_{\ginvangstimesc}\big(\ep_{\ginvangstimesc}\sscs \ang{a_{\bdj}}_{\bdj \in \ginvangstimes}\big)}
(b3)++(0,\v) node (b4) {\phantom{\gaA_{\angstimes}}}
(b4)++(-.5,0) node (b4') {g\gaA_{\angstimesc}\big(\ep_{\angstimesc}\sscs \ang{a_{\bdj}}_{\bdj \in \ginvangstimes}\big)}
;
\draw[1cell=.8]
(a1) edge node {\zdpg_{a,\angs}} (b1)
(a1) edge[equal] node[swap] {\mathbf{d}} (a2)
(a2) edge[equal] node[swap] {\mathbf{g}'} (a3)
(a3) edge[equal] node[swap] {\mathbf{eq}'} (a4)
(a4) edge[equal] node[swap] {\mathbf{eq}} (a5)
(b1) edge[equal] node {\mathbf{g}} (b2)
(b2) edge[equal] node {\mathbf{d}} (b3)
(b3) edge[equal] node {\mathbf{c}} (b4)
;
\draw[1cell=.75]
(a5) [rounded corners=2pt] -| node[pos=.25] {g\gaA_{\angstimesc}(\ual^{g,\angs}; 1^{\angstimesc})} node[swap,pos=.25] {\iso} (b4)
;
\end{tikzpicture}
\end{equation}
In the bottom arrow in \cref{zdp_gas}, the isomorphism $\ual^{g,\angs}$ is defined in \cref{ual_g_angs}.  The seven equalities in \cref{zdp_gas} are given as follows.
\begin{itemize}
\item The two equalities labeled $\mathbf{d}$ follow from the definition \cref{zdsgad_comp_obj} of $\zdsgad$.
\item The equality labeled $\mathbf{g}$ follows from the definition \cref{ga_scomponentGG} of the $g$-action on $\angordnbe$-systems.
\item The equality labeled $\mathbf{g}'$ follows from the definition \cref{proCnbe_gaction} of the $g$-action on $\proAnbe$.
\item The equality labeled $\mathbf{eq}'$ follows from the $G$-equivariance of the functor $\gaA_{\angstimesc}$ \cref{gaAn}.
\item The equality labeled $\mathbf{eq}$ holds by the action equivariance axiom \cref{pseudoalg_action_sym} for $\A$, applied to the permutation $\usi^{g,\angs} \in \Si_{\angstimesc}$ \cref{usi_g_angs}.
\item The equality labeled $\mathbf{c}$ holds because the set
\[\ginvangstimes = \txprod_{j \in \ufs{q}}\, \ginv s_j\]
has the same cardinality as the set $\angstimes = \txprod_{j \in \ufs{q}}\, s_j$.
\end{itemize}
\cref{zdpg_a_welldef} proves that, as the marker $\angs$ varies, $\zdpg_a$ is an isomorphism of strong $\angordnbe$-systems.  
\item[Naturality] 
The naturality of $\zdpg$ with respect to morphisms in $\proAnbe$ follows from the naturality of $g \gaA_{\angstimesc}(\ual^{g,\angs}; -)$ for each marker $\angs$, \cref{gtheta_angsGG,proCnbe_gaction,zdsgad_comp_mor,zdp_gas}.  
\end{description}
This finishes the definition of the natural isomorphism $\zdpg$.
\end{definition}

\begin{example}[Nontriviality of Pseudo Equivariant Constraints]\label{ex:ual_g_angs}
By \cref{zdp_gas}, the component $\zdps_{a,\angs}$ is defined as $g\gaA_{\angstimesc}(\ual^{g,\angs}; 1^{\angstimesc})$, where $\ual^{g,\angs}$ is the unique isomorphism \cref{ual_g_angs}
\[\ginv \ep_{\angstimesc} \usi^{g,\angs} \fto[\iso]{\ual^{g,\angs}} \ep_{\angstimesc}\]
in the translation category $\Op(\angstimesc)$.  Here are two examples that illustrate the nontriviality of $\ual^{g,\angs}$.
\begin{enumerate}
\item For the Barratt-Eccles $\Gcat$-operad $\BE$ \pcref{def:BE}, on which $G$ acts trivially, $\ual^{g,\angs}$ is the unique isomorphism
\[\ginv \ep_{\angstimesc} \usi^{g,\angs}
= \ep_{\angstimesc} \usi^{g,\angs} 
\fto[\iso]{\ual^{g,\angs}} \ep_{\angstimesc} \inspace \ESigma_{\angstimesc}.\]
\item For the $G$-Barratt-Eccles operad $\GBE$ \pcref{def:GBE} and each element $h \in G$, the $h$-component of $\ual^{g,\angs}$ is the unique isomorphism
\[(\ginv \ep_{\angstimesc} \usi^{g,\angs})(h)
= \ep_{\angstimesc}(gh) \cdot \usi^{g,\angs} 
\fto[\iso]{\ual^{g,\angs}_h} \ep_{\angstimesc}(h) \inspace \ESigma_{\angstimesc}.\]
\end{enumerate}
The nonidentity permutation \cref{usi_g_angs}
\[\ginvangstimes \fto{\usi^{g,\angs}} \angstimes \inspace \Si_{\angstimesc}\]
changes according to $g \in G$, the marker $\angs$, and the $G$-action $\be_j$ on $\ordn_j$ for $j \in \ufsq$.    Thus, for each of $\BE$ and $\GBE$, the isomorphism $\ual^{g,\angs}$ is not generally the identity.
\end{example}

\begin{explanation}[Comparing $\zbp$ and $\zdp$]\label{expl:zdp_zbp}
\cref{expl:zdsgad_zbsgad} compares the pointed functors $\zbsgad$ \cref{zbsgad} and $\zdsgad$ \cref{zdsgad}.  Adding to that discussion, \cref{table.pGec} compares their pseudo $G$-equivariant constraints $\zbpg$ and $\zdpg$ \pcref{def:zbsgad_pseudoG,def:zdsgad_pseudoG}.
\begin{figure}[H] 
\centering
\resizebox{.9\width}{!}{%
{\renewcommand{\arraystretch}{1.3}%
{\setlength{\tabcolsep}{1ex}
\begin{tabular}{cr|cr}
$\zbsgad g \fto{\zbpg} g\zbsgad$ & \cref{zbsgad_constraint} & $\zdsgad g \fto{\zdpg} g\zdsgad$ & \cref{zdsgad_constraint} \\[.2ex] \hline
$\big(\ginv s \fto{\si^{g,s}} s\big) \in \Si_{\cards}$ & \cref{si_g_s} 
& $\big(\ginvangstimes \fto{\usi^{g,\angs}} \angstimes\big) \in \Si_{\angstimesc}$ & \cref{usi_g_angs}
\\
$\big(\ginv \ep_{\cards} \si^{g,s} \fto{\al^{g,s}} \ep_{\cards}\big) \in \Op(\cards)$ & \cref{al_g_s} 
& $\big(\ginv \ep_{\angstimesc} \usi^{g,\angs} \fto{\ual^{g,\angs}} \ep_{\angstimesc}\big) \in \Op(\angstimesc)$ & \cref{ual_g_angs}
\\
$((\zbsgad g)a)_s \fto{\zbpg_{a,s}} ((g\zbsgad)a)_s$ & \cref{pse_gas} 
& $((\zdsgad g)a)_{\angs} \fto{\zdpg_{a,\angs}} ((g\zdsgad)a)_{\angs}$ & \cref{zdp_gas}
\\[.4ex]
$\zbpg_{a,s} = g\gaA_{\cards} (\al^{g,s}; 1^{\cards})$ & \cref{pse_gas} 
& $\zdpg_{a,\angs} = g\gaA_{\angstimesc}(\ual^{g,\angs}; 1^{\angstimesc})$ & \cref{zdp_gas}
\end{tabular}}}}
\caption{Comparison of pseudo $G$-equivariant constraints.}
\label{table.pGec}
\end{figure}
\cref{zdpg_a_welldef,thm:zdsgad_pseudo} reuse some of the proofs in \cref{sec:pseudoequiv} by using this dictionary between \cref{def:zbsgad_pseudoG,def:zdsgad_pseudoG}. 
\end{explanation}

\cref{zdpg_a_welldef} is used in \cref{def:zdsgad_pseudoG}.

\begin{lemma}\label{zdpg_a_welldef}
For each object $a \in \proAnbe$, the collection \cref{zdp_g_a}
\[(\zdsgad g)a \fto[\iso]{\zdps_a} (g \zdsgad) a\]
is an isomorphism of strong $\angordnbe$-systems in $\A$.
\end{lemma}

\begin{proof}
For a marker $\angs = \ang{s_j \subseteq \ufs{n}_j}_{j \in \ufs{q}}$, the $\angs$-component \cref{zdp_gas}
\[\zdpg_{a,\angs} = g\gaA_{\angstimesc} \big(\ual^{g,\angs}; 1^{\angstimesc} \big)\]
is an isomorphism in $\A$, where $\ual^{g,\angs} \in \Op(\angstimesc)$ is the unique isomorphism defined in \cref{ual_g_angs}.  We verify the unity axiom \cref{nsystem_mor_unity} and the compatibility axiom \cref{nsystem_mor_compat} for $\zdps_a$.
\begin{description}
\item[Unity]
If $s_j = \emptyset$ for some $j \in \ufs{q}$, then 
\[\begin{split}
\angstimes &= \emptyset = \ginvangstimes,\\
\usi^{g,\angs} &= \id_0 \in \Si_0,\andspace\\
\ual^{g,\angs} &= 1_* \cn \ep_0 = * \to \ep_0 = * \in \Op(0).
\end{split}\]
The functoriality of $\gaA_0$, the functoriality of the $g$-action on $\A$, and the fact that $\gaA_0(*) = \zero \in \A$ is $G$-fixed \cref{pseudoalg_zero} imply
\[\zdpg_{a,\angs} = g\gaA_0( 1_* ) = g 1_\zero = 1_\zero \inspace \A.\]
This proves the unity axiom \cref{nsystem_mor_unity} for $\zdps_a$.
\item[Compatibility]
The compatibility diagram \cref{nsystem_mor_compat} for $\zdps_a$ is the diagram \cref{zdps_compat} for an object $x \in \Op(r)$ with $r \geq 0$, an index $k \in \ufs{q}$, and a partition 
\[s_k = \txcoprod_{\ell \in \ufs{r}}\, \skell \subseteq \ufs{n}_k,\]  
where $\ang{\Cdots}_\ell$ means $\ang{\Cdots}_{\ell \in \ufs{r}}$.  The top and bottom horizontal arrows are the gluing isomorphisms of, respectively, $\zdsgad ga$ and $g\zdsgad a$ at $(x; \angs, k,  \ang{\skell}_\ell)$.
\begin{equation}\label{zdps_compat}
\begin{tikzpicture}[vcenter]
\def\v{-1.5}
\draw[0cell=.9]
(0,0) node (a11) {\gaA_r\big(x\sscs \bang{(\zdsgad ga)_{\angs \compk\, s_{k,\ell}}}_{\ell} \big)}
(a11)++(3.5,0) node (a12) {(\zdsgad ga)_{\angs}}
(a11)++(0,\v) node (a21) {\gaA_r\big(x\sscs \bang{(g\zdsgad a)_{\angs \compk\, s_{k,\ell}}}_{\ell} \big)}
(a12)++(0,\v) node (a22) {(g\zdsgad a)_{\angs}}
;
\draw[1cell=.9]
(a11) edge (a12)
(a21) edge (a22)
(a11) edge[transform canvas={xshift={2.5em}}, shorten >=-.3ex] node[swap] {\gaA_r(1;\, \ang{\zdps_{a,\angs \compk s_{k,\ell}}}_{\ell})} (a21)
(a12) edge node {\zdps_{a,\angs}} (a22)
;
\end{tikzpicture}
\end{equation}
The diagram \cref{zdps_compat} is the $\zdpg_a$-analogue of the compatibility diagram \cref{zbps_compat} for $\zbpg_a$, which factors into the commutative diagram \cref{zbps_compat_ii}.  Starting from the commutative diagram \cref{zbps_compat_ii} and changing the notation as discussed in \cref{expl:zdsgad_zbsgad,expl:zdp_zbp}, the diagram \cref{zdps_compat} factors into a number of commutative regions.  In addition to the dictionaries in \cref{expl:zdsgad_zbsgad,expl:zdp_zbp}, we replace 
\begin{itemize}
\item $s_\ell$ by $\angscompkskelltimes$ in $|s_\ell|$, $\bdi \in s_\ell$, and $\bdj \in \ginv s_\ell$;
\item the permutation \cref{si_g_s} 
\[\ginv s_\ell \fto{\si^{g,s_\ell}} s_\ell\]
by the permutation
\[\ginv \angscompkskelltimes \fto{\usi^{g,\angscompkskell}} \angscompkskelltimes\]
in \cref{usi_g_angs};
\item the isomorphism \cref{al_g_s} 
\[\ginv \ep_{|s_\ell|} \si^{g,s_\ell} \fto{\al^{g,s_\ell}} \ep_{|s_\ell|}\]
by the isomorphism \cref{ual_g_angs}
\[\ginv \ep_{|\angscompkskelltimes|} \usi^{g,\angscompkskell} \fto{\ual^{g,\angscompkskell}} \ep_{|\angscompkskelltimes|}\]
in $\Op(|\angscompkskelltimes|)$;
\item the permutation \cref{sigma_angsi} 
\[\ginv s \fto{\sig_{\ginv s, \ang{\ginv s_\ell}_\ell}} \txcoprod_{\ell \in \ufs{r}}\, \ginv s_\ell\] 
by the permutation
\[\ginv\angstimes \fto{\usi_{\ginv\angs,\, k, \ang{\ginv \skell}_\ell}} 
\txcoprod_{\ell \in \ufs{r}}\, \ginv \angscompkskelltimes\]
in \cref{usi_perm}; and
\item the isomorphism \cref{al_Op} 
\[\ga(\ginv x; \ang{\ep_{|\ginv s_\ell|}}_\ell) \sig_{\ginv s, \ang{\ginv s_\ell}_\ell} 
\fto{\alp_{\ginv x;\, \ginv s, \ang{\ginv s_\ell}_\ell}} \ep_{|\ginv s|}\] 
by the isomorphism 
\[\scalebox{.9}{$
\ga\big(\ginv x; \ang{\ep_{|\ginv \angscompkskelltimes|}}_\ell \big) \usi_{\ginv \angs,\, k, \ang{\ginv \skell}_\ell}
\fto{\ual_{\ginv x;\, \ginv \angs,\, k, \ang{\ginv \skell}_\ell}}
\ep_{|\ginv \angstimes|}
$}\]
in \cref{ual_iso}.
\end{itemize} 
This proves the compatibility axiom \cref{nsystem_mor_compat} for $\zdps_a$.\qedhere
\end{description}
\end{proof}

\subsection*{Pseudo $G$-Equivariance of $\zdsgad$}

We now prove that the equivalence $\zdsgad$ \pcref{thm:zdsg_eq} is a pseudo $G$-equivariant functor \pcref{def:pseudoG}.

\begin{theorem}\label{thm:zdsgad_pseudo}
Under the same assumptions as \cref{def:zdsgad_pseudoG}, the pair
\[\proAnbe \fto{(\zdsgad,\zdp)} \Asgangordnbe\]
consisting of
\begin{itemize}
\item the pointed functor $\zdsgad$ \pcref{def:zdsgad} and
\item the natural isomorphisms $\{\zdps\}_{g \in G}$ \cref{zdsgad_constraint}
\end{itemize}
is a pseudo $G$-equivariant functor.
\end{theorem}

\begin{proof}
We verify the axioms \cref{pseudoG_unity,pseudoG_mult} of a pseudo $G$-equivariant functor for $(\zdsgad,\zdp)$.
\begin{description}
\item[Unity]  
For the identity element $e \in G$ and a marker $\angs = \ang{s_j \subseteq \ufs{n}_j}_{j \in \ufs{q}}$, the permutation \cref{usi_g_angs}
\[\big(e^{-1} \angstimes = \angstimes \fto[\iso]{\usi^{e,\angs}} \angstimes\big) \in \Si_{\angstimesc}\]
is the identity.  Thus, the unique isomorphism \cref{ual_g_angs}
\[e^{-1} \ep_{\angstimesc} \usi^{e,\angs} = \ep_{\angstimesc} 
\fto[\iso]{\ual^{e,\angs}} \ep_{\angstimesc}\]
in the translation category $\Op(\angstimesc)$ is the identity morphism.  For each object $a \in \proAnbe$, the $\angs$-component isomorphism \cref{zdp_gas}
\[\zdpe_{a,\angs} = e \gaA_{\angstimesc}\big(1_{\ep_{\angstimesc}} ; 1^{\angstimesc}\big)\]
is the identity morphism by the functoriality of $\gaA_{\angstimesc}$ and the fact that the $e$-action on $\A$ is the identity.  This proves that $\zdpe$ is the identity natural transformation of $\zdsgad$, proving the unity axiom \cref{pseudoG_unity}.
\item[Multiplicativity]   
Morphisms in the category $\Asgangordnbe$ \pcref{def:nsystem_morphism} are determined by their components.  The multiplicativity axiom \cref{pseudoG_mult} for $\zdp$ means that, for elements $g,h \in G$, an object $a \in \proAnbe$ \cref{proCnbe_object}, and a marker $\angs = \ang{s_j \subseteq \ufs{n}_j}_{j \in \ufs{q}}$, the diagram \cref{zdp_mult} in $\A$ commutes.
\begin{equation}\label{zdp_mult}
\begin{tikzpicture}[vcenter]
\def\h{3.5} \def\t{15}
\draw[0cell]
(0,0) node (a1) {\big(\zdsgad hga\big)_{\angs}}
(a1)++(\h/2,-1) node (a2) {\big(h \zdsgad ga\big)_{\angs}}
(a1)++(\h,0) node (a3) {\big(hg\zdsgad a\big)_{\angs}}
;
\draw[1cell=.9]
(a1) edge node {\zdphg_{a,\angs}} (a3)
(a1) [rounded corners=2pt,shorten <=-.5ex] |- node[swap,pos=.25] {\zdph_{ga,\angs}} (a2)
;
\draw[1cell=.9]
(a2) [rounded corners=2pt,shorten >=-.5ex] -| node[swap,pos=.7] {(h\zdpg_a)_{\angs}} (a3)
;
\end{tikzpicture}
\end{equation}
The diagram \cref{zdp_mult} is the $\zdp$-analogue of the multiplicativity diagram \cref{zbp_mult} for $\zbp$.  By  changing the notation as discussed in \cref{expl:zdsgad_zbsgad,expl:zdp_zbp}, the proof of the commutativity of \cref{zbp_mult} also proves that the diagram \cref{zdp_mult} commutes.  

More precisely, similar to \cref{zbphg_as,zbph_gas,hzbpg_ahinvs}, we use \cref{usi_g_angs,ual_g_angs,zdp_gas} to unravel the three morphisms in \cref{zdp_mult} as follows.
\[\begin{split}
\zdphg_{a,\angs} 
&= hg \gaA_{\angstimesc}\big(\ual^{hg,\angs} ; \ang{a_{\bdj}}_{\bdj \in (hg)^{-1} \angstimes} \big)\\
\zdph_{ga,\angs} 
&= hg\gaA_{\angstimesc}\big(\ginv \ual^{h,\angs} \usi^{g,\hinv \angs}; \ang{a_{\bdj}}_{\bdj \in \ginv\hinv \angstimes} \big)\\
(h\zdpg_a)_{\angs}
&= hg \gaA_{\angstimesc}\big(\ual^{g,\hinv \angs}; \ang{a_{\bdj}}_{\bdj \in \ginv\hinv \angstimes}\big)
\end{split}\]
Each of these three morphisms has the form
\[hg\gaA_{\angstimesc}\big(-; \ang{a_{\bdj}}_{\bdj \in (hg)^{-1} \angstimes}\big)\]
for some morphism $-$ in $\Op(\angstimesc)$.  Thus, the diagram \cref{zdp_mult} commutes because the following diagram in the translation category $\Op(\angstimesc)$ commutes.
\begin{equation}\label{zdp_mult_diag}

\end{equation}
The left vertical equality uses the following commutative diagram of permutations, which is analogous to \cref{sig_hinvs}.
\begin{equation}\label{zdp_mult_perm}
%
\end{equation}
The permutations $\usi^{hg,\angs}$, $\usi^{g,\hinv \angs}$, and $\usi^{h,\angs}$ are restrictions of, respectively, the $hg$-action, the $g$-action, and the $h$-action on $\txprod_{j \in \ufs{q}}\, \ufs{n}_j$.  This proves the multiplicativity axiom \cref{pseudoG_mult} for $(\zdsgad,\zdp)$.\qedhere
\end{description}
\end{proof}

\begin{explanation}[Pseudo $G$-Equivariance of $\zdad$]\label{expl:zdad_pseudo}
Recall the adjunction \pcref{expl:system_adjunction} 
\begin{equation}\label{zd_zdad_adj}
%
\end{equation}
involving the pointed $G$-functor $\zd$ \cref{hgotoprod} and the pointed functor $\zdad = \iota\zdsgad$ \cref{zdad}.  As a consequence of \cref{thm:zdsgad_pseudo}, the pair 
\[\proAnbe \fto{(\zdad, \iota * \zdp)} \Aangordnbe\]
is a pseudo $G$-equivariant functor, where $\zdp = \{\zdpg\}_{g \in G}$ is the pseudo $G$-equivariant constraint for $\zdsgad$ \cref{zdsgad_constraint}. 
\begin{description}
\item[Constraints] For each $g \in G$, its pseudo $G$-equivariant constraint $\iota * \zdpg$ is given by whiskering the natural isomorphism $\zdpg$ \cref{zdsgad_constraint} with the inclusion $\iota$, which is a $G$-functor, as follows.
\begin{equation}\label{zdpg_whisker}
\begin{tikzpicture}[vcenter]
\def\v{-1.5} \def\u{-1.2}
\draw[0cell=.9]
(0,0) node (a1) {\proAnbe}
(a1)++(0,\v) node (a2) {\Asgangordnbe}
(a2)++(0,\u) node (a3) {\Aangordnbe}
(a1)++(2.8,0) node (b1) {\proAnbe}
(b1)++(0,\v) node (b2) {\Asgangordnbe}
(b2)++(0,\u) node (b3) {\Aangordnbe}
;
\draw[1cell=.9]
(a1) edge node {g} (b1)
(a2) edge node[swap] {g} (b2)
(a3) edge node {g} (b3)
(a1) edge node[swap] {\zdsgad} (a2)
(a2) edge node[swap] {\iota} (a3)
(b1) edge node {\zdsgad} (b2)
(b2) edge node {\iota} (b3)
;
\draw[2cell]
node[between=a1 and b2 at .55, rotate=-135, 2labelalt={below,\zdpg}] {\Rightarrow}
;
\end{tikzpicture}
\end{equation}
\item[Unity] The unity axiom \cref{pseudoG_unity} for $\iota * \zdpe$ follows from the unity axiom for $\zdpe$ because the inclusion $\iota$ preserves identity morphisms.
\item[Multiplicativity] The multiplicativity axiom \cref{pseudoG_mult} for $\iota * \zdp$ is obtained from the multiplicativity axiom for $\zdp$ by whiskering with $\iota$.
\end{description}
In summary, the right adjoint $\zd$ is $G$-equivariant (\cref{Pist_to_prod} \eqref{Pist_to_prod_i}), and the left adjoint $\zdad$ is pseudo $G$-equivariant, with pseudo $G$-equivariant constraint given by $\iota * \zdp$.
\end{explanation}

%% file: chap/compgen.tex
\cref{thm:PistAequivalence} proves that the strong $H$-theory comparison pointed $G$-functor \cref{PistsgAnbe} 
\[(\smast\Sgosg\A)\angordnbe = \Asgsmaangordnbe 
\fto[\sim]{\Pistsg_{\A,\angordnbe}} 
\Asgangordnbe = (\Hgosg\A)\angordnbe\]
is a nonequivariant equivalence of categories for each object $\angordnbe \in \GG$ \cref{GG_objects}, $\Uinf$-operad $\Op$ \pcref{as:OpA'}, and $\Op$-pseudoalgebra $\A$ \pcref{def:pseudoalgebra}.  The main result of this chapter, \cref{thm:pistweakgeq}, upgrades \cref{thm:PistAequivalence} for the $\Uinf$-operad $\Oph =  \Catg(\EG,\Op)$ and an $\Oph$-pseudoalgebra of the form $\Ah = \Catg(\EG,\A)$ \pcref{catgego}.  \cref{thm:pistweakgeq} proves that the strong $H$-theory comparison pointed $G$-functor 
\[(\smast\Sgohsg\Ah)\angordnbe = \Ahsgsmaangordnbe 
\fto[\eqg]{\Pistsg_{\Ah,\angordnbe}} 
\Ahsgangordnbe = (\Hgohsg\Ah)\angordnbe\]
is a \emph{categorical weak $G$-equivalence} \pcref{def:cat_weakg}.  This means that, for each subgroup $H \subseteq G$, the $H$-fixed subfunctor $\big(\Pistsg_{\Ah,\nbe}\big)^H$ is an equivalence of categories.  Thus, for each $\Oph$-pseudoalgebra of the form $\Ah = \Catg(\EG,\A)$, Shimakawa strong $H$-theory and strong $H$-theory are componentwise categorically weakly $G$-equivalent via $\Pistsg$.

\subsection*{Application to $G$-Spaces}
As a consequence of \cref{thm:pistweakgeq}, applying the classifying space functor $\cla \cn \Cat \to \Top$ yields a weak $G$-equivalence between $G$-spaces
\[\cla\Ahsgsmaangordnbe
\fto[\eqg]{\cla\Pistsg_{\Ah,\angordnbe}} \cla\Ahsgangordnbe.\]
Thus, the strong $H$-theory comparison pointed $G$-functor $\Pistsg_{\Ah,\nbe}$ is a \emph{topological weak $G$-equivalence} in the sense of Merling \cite{merling} \pcref{def:weakG}.  See \cref{ex:catweakg_weakg,ex:pistsg_weakg}.

\strategy
\cref{thm:pistweakgeq} is proved by a 2-out-of-3 argument applied to the following commutative diagram of $G$-functors, where $\Pistsg_{\Ah,\nbe}$ is abbreviated to $\Pistsg$.
\begin{equation}\label{Pist_strategy}
\begin{tikzpicture}[vcenter]
\def\v{-1.4}
\draw[0cell=.85]
(0,0) node (a11) {\sgAhsmanbe}
(a11)++(5,0) node (a12) {\sgAhnbe}
(a11)++(0,\v) node (a21) {\Catg(\EG,\sgAhsmanbe)}
(a12)++(0,\v) node (a22) {\Catg(\EG,\sgAhnbe)}
;
\draw[1cell=.8]
(a11) edge node {\Pistsg} (a12)
(a12) edge node {\inc^{\sgAhnbe}} node[swap] {\gsim} (a22)
(a11) edge node {\gsim} node[swap] {\inc^{\sgAhsmanbe}} (a21)
(a21) edge node {\Catg(\EG,\Pistsg)} node[swap] {\eqg} (a22)
;
\end{tikzpicture}
\end{equation}
For each small $G$-category $\C$, the inclusion $G$-functor $\inc \cn \C \to \Catg(\EG,\C)$ is a nonequivariant equivalence \pcref{inc_eq}.  \cref{thm:SgoAh,thm:HgoAh} prove that, for each of the $G$-categories $\sgAhsmanbe$ and $\sgAhnbe$, the inclusion $G$-functor $\inc$ is the left adjoint of an \emph{adjoint $G$-equivalence} \pcref{def:adjointGeq}.  This means that each of the two vertical $G$-functors $\inc$ in \cref{Pist_strategy} admits a $G$-equivariant inverse and invertible $G$-equivariant unit and counit that satisfy the triangle identities for a $G$-adjunction.  Thus, each vertical $G$-functor $\inc$ is a categorical weak $G$-equivalence \pcref{def:cat_weakg}.  Moreover, since the $G$-functor $\Pistsg$ is a nonequivariant equivalence of categories \pcref{thm:PistAequivalence}, the $G$-functor $\Catg(\EG,\Pistsg)$ is a categorical weak $G$-equivalence \pcref{merling_2.16}.  Since the left, right, and bottom $G$-functors in \cref{Pist_strategy} are categorical weak $G$-equivalences, so is the $G$-functor $\Pistsg$.

\summary 
\Cref{table.equivalen} summarizes the $G$-functors mentioned in the preceding discussion, from the weakest notion (nonequivariant equivalence) to the strongest one (adjoint $G$-equivalence).
\begin{figure}[H] 
\centering
\resizebox{.9\width}{!}{%
{\renewcommand{\arraystretch}{1.3}%
{\setlength{\tabcolsep}{1.5ex}
\begin{tabular}{c|cr}
\multirow{2}{*}{\makecell{nonequivariant \\ equivalences}} &
$\Pistsg$ for $\Op$-pseudoalgebras $\A$ & \ref{thm:PistAequivalence} \\
& $\inc$ for $G$-categories & \ref{inc_eq} \\ \hline
\multirow{2}{*}{\makecell{categorical weak \\ $G$-equivalences \eqref{def:cat_weakg}}} &
$\Catg(\EG,\fun)$ for equivalences $\fun$ & \ref{merling_2.16} \\
& $\Pistsg$ for $\Oph$-pseudoalgebras $\Ah$ & \ref{thm:pistweakgeq} \\ \hline
\multirow{2}{*}{\makecell{adjoint \\ $G$-equivalences \eqref{def:adjointGeq}}} &
\raisebox{6pt}{\phantom{M}}$\inc$ for $\Ahmal$, $\sgAhmal$\raisebox{6pt}{\phantom{M}} & \ref{thm:SgoAh} \\
& $\inc$ for $\Ahnbe$, $\sgAhnbe$ & \ref{thm:HgoAh}
\end{tabular}}}}
\caption{Three notions of equivalences for equivariant functors.}
\label{table.equivalen}
\end{figure}
\noindent
The following diagram summarizes the relationships between various kinds of equivalences.  The letters and arrows are explained after the diagram.
\begin{equation}\label{equiv_summary}
\begin{tikzpicture}[vcenter]
\def\v{1.3} \def\h{4}
\draw[0cell]
(0,0) node (a) {A}
(a)++(0,\v) node (c) {C}
(c)++(0,\v) node (n) {N}
(c)++(\h,0) node (s) {S}
(c)++(\h/2,-\v) node (m) {T}
;
\draw[1cell]
(a) edge (c)
(c) edge (n)
(n) edge[bend right=40] node[swap] {\Catg(\EG,-)} (c)
(c) edge node {\cla} (s)
(c) edge (m)
(m) edge node[swap] {\cla} (s)
;
\end{tikzpicture}
\end{equation}
\begin{itemize}
\item $N$: $G$-functors that are nonequivariant equivalences of categories.
\item $C$: categorical weak $G$-equivalences \pcref{def:cat_weakg}.
\item $A$: adjoint $G$-equivalences \pcref{def:adjointGeq}.
\item $S$: weak $G$-equivalences between $G$-spaces \pcref{def:weakG_top}.
\item $T$: topological weak $G$-equivalences in the sense of Merling \cite{merling} \pcref{def:weakG}.
\item Each of the three unlabeled arrows indicates that something in its source is also in the target.  For example, the arrow $A \to C$ means that a $G$-functor that is the left or right adjoint of an adjoint $G$-equivalence is also a categorical weak $G$-equivalence \pcref{expl:cat_weakg,expl:weakG}.
\item $\cla$ is the classifying space functor \pcref{ex:catweakg_weakg}.
\item The arrow labeled $\Catg(\EG,-)$ is given by \cref{merling_2.16}.
\end{itemize}

\organization
This chapter consists of the following sections.

\secname{sec:inc_gfunctor}
This section discusses the inclusion $G$-functor $\inc \cn \C \to \Catg(\EG,\C)$ for a small $G$-category $\C$.  After discussing the fact that it is an equivalence of categories, we observe that its adjoint inverse
\[\Catg(\EG,\C) \fto{\pn} \C\]
is a pseudo $G$-equivariant functor \pcref{ppn_pseudoG}.  This observation provides a useful contrast with some $G$-functors in \cref{sec:sgo_geq,sec:hgo_geq} from $\Catg(\EG,A)$ to $A$, where $A$ is either Shimakawa (strong) $H$-theory or (strong) $H$-theory for an $\Oph$-pseudoalgebra of the form $\Ah = \Catg(\EG,\A)$.

\secname{sec:sgo_geq}
To show that the inclusion $G$-functor $\inc$ for each of the $G$-categories $\Ahmal$ and $\sgAhmal$ is an adjoint $G$-equivalence, this section constructs a $G$-functor
\[\Catg(\EG, \Ahmal) \fto{\cni} \Ahmal\]
and its strong variant
\[\Catg(\EG, \sgAhmal) \fto{\cnisg} \sgAhmal\]
that serve as $G$-equivariant inverses of $\inc$ for, respectively, $\Ahmal$ and $\sgAhmal$.  These constructions are based on the diagonal.  See \cref{cni_fsg,cnif_gl_g,cni_theta_sg}, where $g \in G$ is used twice.  They are different from the pseudo $G$-equivariant adjoint inverse $\pn$.

\secname{sec:sgo_geq_unit}
This section constructs the invertible $G$-equivariant unit and counit for each pair of $G$-functors $(\inc,\cni)$ and $(\inc,\cnisg)$.

\secname{sec:sgosg_gequiv}
The first half of this section discusses \emph{adjoint $G$-equivalences}.  As an example, \cref{EHEG} shows that a subgroup inclusion $H \subseteq G$ induces the left adjoint of an adjoint $H$-equivalence between the translation categories $\EH$ and $\EG$.  Using the unit and counit in \cref{sec:sgo_geq_unit}, the rest of this section proves that the $G$-functors
\begin{equation}\label{inc_cni_chi}
\begin{tikzpicture}[vcenter]
\def\h{1.8} \def\t{25}
\draw[0cell]
(0,0) node (a1) {\phantom{A}}
(a1)++(\h,0) node (a2) {\phantom{A}}
(a1)++(-.2,0) node (a1') {\Ahmal}
(a2)++(1,0) node (a2') {\Catg(\EG, \Ahmal)}
(a1)++(\h/2,0) node () {\gsim}
;
\draw[1cell=.9]
(a1) edge[bend left=\t] node {\inc} (a2)
(a2) edge[bend left=\t] node {\cni} (a1)
;
\end{tikzpicture}
\end{equation}
form an adjoint $G$-equivalence, and likewise for the strong variant \pcref{thm:SgoAh}.

\secname{sec:hgo_geq}
This section is the $H$-theory analogue of \cref{sec:sgo_geq}.  It constructs a $G$-functor
\[\Catg(\EG, \Ahnbe) \fto{\ci} \Ahnbe\]
and its strong variant
\[\Catg(\EG, \sgAhnbe) \fto{\cisg} \sgAhnbe\]
that serve as $G$-equivariant inverses of $\inc$ for, respectively, $\Ahnbe$ and $\sgAhnbe$.  \cref{expl:cnici_compare} compares the $G$-functors $\cni$ and $\ci$.

\secname{sec:hgo_gequiv}
This section is the $H$-theory analogue of \cref{sec:sgo_geq_unit,sec:sgosg_gequiv}.  It constructs the invertible $G$-equivariant unit and counit for each pair of $G$-functors $(\inc,\ci)$ and $(\inc,\cisg)$.  Then it proves that the $G$-functors
\begin{equation}\label{inc_ci_seci}
\begin{tikzpicture}[vcenter]
\def\h{1.8} \def\t{25}
\draw[0cell]
(0,0) node (a1) {\phantom{A}}
(a1)++(\h,0) node (a2) {\phantom{A}}
(a1)++(\h/2,0) node () {\gsim}
(a1)++(-.3,0) node (a1') {\Ahnbe}
(a2)++(1.1,0) node (a2') {\Catg(\EG, \Ahnbe)}
;
\draw[1cell=.9]
(a1) edge[bend left=\t] node {\inc} (a2)
(a2) edge[bend left=\t] node {\ci} (a1)
;
\end{tikzpicture}
\end{equation}
form an adjoint $G$-equivalence, and likewise for the strong variant \pcref{thm:HgoAh}.

\secname{sec:cat_weakg}
The first half of this section introduces \emph{categorical weak $G$-equivalences} and discusses how they are related to nonequivariant equivalences and adjoint $G$-equivalences, with illustrative examples from earlier sections.  \cref{merling_2.16} provides a self-contained treatment of a result of Merling \cite{merling}.  This result states that, for each $G$-functor $\fun$ between small $G$-categories that is also a nonequivariant equivalence of categories, the $G$-functor $\Catg(\EG,\fun)$ is a categorical weak $G$-equivalence.  The proof of \cref{thm:pistweakgeq} uses this result for $\Pistsg_{\Ah,\nbe}$.

\secname{sec:levelgeq}
This section proves the main result of this chapter, \cref{thm:pistweakgeq}.  It states that the comparison pointed $G$-functor $\Pistsg_{\Ah}$ is componentwise a categorical weak $G$-equivalence.  \cref{ex:pistsg_weakg,ex:pistweq} are applications to weak $G$-equivalences and genuine symmetric monoidal $G$-categories.

\section{Inclusion $G$-Functors and Pseudo $G$-Equivariant Inverses}
\label{sec:inc_gfunctor}

This section discusses the inclusion $G$-functor 
\[\C \fto{\inc} \Catg(\EG,\C)\] 
for a group $G$ and a small $G$-category $\C$.  The inclusion $G$-functor $\inc$ embeds $\C$ into the $G$-thickening $\Catg(\EG,\C)$, which is the $G$-category of functors $\EG \to \C$ and natural transformations, with $G$ acting by conjugation.  It is a nonequivariant equivalence of categories.  Its adjoint inverse 
\[\Catg(\EG,\C) \fto{\pn} \C,\]
given by evaluating at the group unit of $G$, is pseudo $G$-equivariant \pcref{def:pseudoG}.  The discussion about $\pn$ in this section provides an instructive contrast with some key constructions in later sections; see \cref{def:cni}.  In subsequent sections, it is shown that, for (Shimakawa) $H$-theory of a $\Catg(\EG,\Op)$-pseudoalgebra of the form $\Ah = \Catg(\EG,\A)$ with $\A$ an $\Op$-pseudoalgebra, the inclusion $G$-functor $\inc$ is part of an \emph{adjoint $G$-equivalence}, and similarly for the strong variant. 

\secoutline
\begin{itemize}
\item \cref{def:gcat_inc} defines the inclusion $G$-functor $\inc \cn \C \to \Catg(\EG,\C)$.
\item \cref{inc_eq} proves that $\inc$ is a nonequivariant equivalence of categories.
\item \cref{expl:inc_eq} discusses the adjoint equivalence $\inc \dashv \pn$ in detail.
\item \cref{def:pn_pseudo} defines the pseudo $G$-equivariant constraints $\ppn$ for $\pn$.
\item \cref{ppn_pseudoG} proves that $(\pn,\ppn)$ is a pseudo $G$-equivariant functor.
\end{itemize}

\subsection*{Inclusion $G$-Functors}

Recall that, for a group $G$ and small $G$-categories $\C$ and $\D$, $\Catg(\C,\D)$ is the small $G$-category with functors $\C \to \D$ as objects, natural transformations as morphisms, and the conjugation $G$-action \pcref{def:Catg}.  Also recall the translation category $\EG$ of $G$ with the regular $G$-action \pcref{def:translation_cat}. 

\begin{definition}[Inclusion $G$-Functors]\label{def:gcat_inc}
Suppose $\C$ is a small $G$-category for an arbitrary group $G$.  The \emph{inclusion $G$-functor}\index{inclusion G-functor@inclusion $G$-functor}\index{G-functor@$G$-functor!inclusion}
\begin{equation}\label{incC}
\C \fto{\inc} \Catg(\EG,\C)
\end{equation}
is defined by sending
\begin{itemize}
\item an object $c \in \C$ to the constant functor $\inc c \cn \EG \to \C$ at $c$ and
\item a morphism $d \cn c \to c'$ in $\C$ to the constant natural transformation $\inc d \cn \inc c \to \inc c'$ with each component given by $d$.
\end{itemize}
The inclusion $G$-functor $\inc$ is also denoted by $\incC$.  The $G$-category $\Catg(\EG,\C)$ is called the \emph{$G$-thickening}\index{G-thickening@$G$-thickening} of $\C$.
\end{definition}

\begin{explanation}[Equivariance of $\inc$]\label{expl:inc_equivariant}
For an element $g \in G$ and an object or a morphism $c \in \C$, the conjugation $g$-action on $\inc c$ 
\begin{equation}\label{ic_conjugation}
\begin{tikzpicture}[vcenter]
\def\h{2}
\draw[0cell]
(0,0) node (a1) {\EG}
(a1)++(\h,0) node (a2) {\EG}
(a2)++(\h,0) node (a3) {\C}
(a3)++(\h,0) node (a4) {\C}
;
\draw[1cell=.9]
(a1) edge node {\ginv} (a2)
(a2) edge node {\inc c} (a3)
(a3) edge node {g} (a4)
;
\draw[1cell=.9]
(a1) [rounded corners=2pt] |- ($(a2)+(0,.6)$) -- node {g(\inc c)\ginv} ($(a3)+(0,.6)$) -| (a4)
 ;
\end{tikzpicture}
\end{equation}
is equal to $\inc(gc)$---that is, constant at $gc \in \C$---because $\inc c$ is constant at $c$.  
\end{explanation}

\subsection*{Nonequivariant Equivalences}
\cref{inc_eq} proves that each inclusion $G$-functor is a nonequivariant equivalence.

\begin{lemma}\label{inc_eq}
For each small $G$-category $\C$, the inclusion $G$-functor \cref{incC}
\[\C \fto{\inc} \Catg(\EG,\C)\]
is an equivalence of categories. 
\end{lemma}

\begin{proof}
The unique functor $\EG \to \bone$ to the terminal category is an equivalence of categories, since it is fully faithful on morphisms and essentially surjective on objects.  The inclusion $G$-functor $\inc$ is obtained from the equivalence $\EG \to \bone$ by applying the contravariant functor $\Catg(-,\C)$ and using the canonical isomorphism $\C \iso \Catg(\bone,\C)$.  
\end{proof}

\begin{explanation}[Equivalence]\label{expl:inc_eq}
To better understand the equivalence $\inc \cn \C \to \Catg(\EG,\C)$ in \cref{inc_eq}, we discuss in detail an adjoint inverse of $\inc$, the unit and counit for the adjunction, and the triangle identities \pcref{def:adjunction}.
\begin{description}
\item[Inverse] 
A right adjoint inverse of $\inc$ is given by the functor
\begin{equation}\label{proC}
\Catg(\EG,\C) \fto{\pn} \C
\end{equation}
that sends
\begin{itemize}
\item a functor $f \cn \EG \to \C$ to the object 
\begin{equation}\label{pn_h}
\pn(f) = f(e) \in \C
\end{equation}
with $e \in G$ denoting the group unit and
\item a natural transformation $\theta \cn f \to f'$ to its $e$-component morphism 
\begin{equation}\label{pn_theta}
\pn(f) = f(e) \fto{\pn(\theta) = \theta_e} \pn(f') = f'(e).
\end{equation}
\end{itemize}
In other words, $\pn$ is given by evaluating at the group unit $e \in G$.
\item[Unit] 
The composite $\pn\inc$ is equal to the identity functor on $\C$.  The unit 
\begin{equation}\label{iunit}
1_\C \fto{\iuni = 1_{1_\C}} \pn\inc = 1_\C
\end{equation}
is the identity natural transformation on $1_\C$.  
\item[Counit]
The counit is the natural isomorphism
\begin{equation}\label{icounit}
\inc\pn \fto[\iso]{\icou} 1_{\Catg(\EG,\C)}
\end{equation}
whose component at a functor $f \cn \EG \to \C$ is a natural isomorphism 
\[\inc\pn f \fto[\iso]{\icou_f} f \inspace \Catg(\EG,\C).\] 
By \cref{incC,pn_h}, the functor $\inc\pn f \cn \EG \to \C$ is constant at $f(e)$.  The component of $\icou_f$ at an object $g \in \EG$ is given by the isomorphism
\begin{equation}\label{icou_hg}
(\inc\pn f)(g) = f(e) \fto[\iso]{\icou_{f,g} = f[g,e]} f(g) \inspace \C.
\end{equation}
This is the image under the functor $f$ of the unique isomorphism 
\[ e \fto[\iso]{[g,e]} g \inspace \EG(e,g).\] 
\begin{itemize}
\item The naturality of $\icou_{f,g}$ with respect to isomorphisms $[g',g] \cn g \to g'$ in $\EG$ follows from the morphism equality
\[[g',e] = [g',g] \circ [g,e] \cn e \fiso g' \inspace \EG\]
and the functoriality of $f$.
\item The naturality of $\icou_f$ with respect to $f$ follows from the fact that morphisms in $\Catg(\EG,\C)$ are natural transformations.
\end{itemize}
\item[Triangle Identities]  
The left triangle identity for the quadruple $(\inc,\pn,\iuni,\icou)$ states that, for each object $c \in \C$, the following composite is the identity natural transformation of $\inc c \in \Catg(\EG,\C)$.
\[\inc c \fto{\inc \iuni_c} \inc\pn\inc c \fto{\icou_{\inc c}} \inc c\]
This composite is equal to $1_{\inc c}$ because, by \cref{incC,iunit,icou_hg}, 
\[\begin{split}
\iuni_c &= 1_c \andspace\\
\icou_{\inc c} &= (\inc c)[-,e] = 1_c.
\end{split}\]
The right triangle identity states that, for each functor $f \cn \EG \to \C$, the following composite is the identity morphism of $\pn f = f(e) \in \C$.
\[\pn f \fto{\iuni_{\pn f}} \pn\inc\pn f \fto{\pn \icou_f} \pn f\]
This composite is equal to $1_{f(e)}$ because, by \cref{pn_h,pn_theta,iunit,icou_hg},
\[\begin{split}
\iuni_{\pn f} &= 1_{\pn f} = 1_{f(e)} \andspace\\
\pn \icou_f &= \icou_{f,e} = f[e,e] = f 1_e = 1_{f(e)}. 
\end{split}\]
The last two equalities use, respectively, the morphism equality 
\begin{equation}\label{eeonee}
e \fto{[e,e] = 1_e} e \inspace \EG
\end{equation}
and the functoriality of $f$.  
\end{description}
In summary, for each small $G$-category $\C$, there is an adjoint equivalence
\begin{equation}\label{inc_pn_adeq}
\begin{tikzpicture}[vcenter]
\def\h{1.8} \def\b{20}
\draw[0cell]
(0,0) node (a1) {\C}
(a1)++(\h,0) node (a2) {\phantom{\C}}
(a1)++(\h/2,0) node () {\sim}
(a2)++(.8,-.03) node (a2') {\Catg(\EG,\C)}
;
\draw[1cell=1]
(a1) edge[bend left=\b] node {\inc} (a2)
(a2) edge[bend left=\b] node {\pn} (a1)
;
\end{tikzpicture}
\end{equation}
with left adjoint $\inc$ \cref{incC}, right adjoint $\pn$ \cref{proC}, unit $\iuni$ \cref{iunit}, and counit $\icou$ \cref{icounit}.  We emphasize that the quadruple $(\inc,\pn,\iuni,\icou)$ is an adjoint equivalence in the nonequivariant sense.  The left adjoint $\inc$ is $G$-equivariant, but its adjoint inverse $\pn$ is not $G$-equivariant in general.
\end{explanation}

\subsection*{Pseudo $G$-Equivariance of $\pn$}
The rest of this section explains that the functor $\pn$ \cref{proC}, the adjoint inverse of the inclusion $G$-functor $\inc$, is \emph{pseudo} $G$-equivariant \pcref{def:pseudoG}.  Its pseudo $G$-equivariant constraints are given in \cref{def:pn_pseudo}.  The composites $\pn g$  \cref{ppng_h_dom} and $g \pn$ \cref{ppng_h_cod} differ by a nonidentity isomorphism in general.  

\begin{definition}\label{def:pn_pseudo}
For a group $G$, an element $g \in G$, a small $G$-category $\C$, and the functor $\pn$ \cref{proC}, we define a natural isomorphism
\begin{equation}\label{ppng}
\begin{tikzpicture}[vcenter]
\def\v{-1.3}
\draw[0cell]
(0,0) node (a1) {\Catg(\EG,\C)}
(a1)++(3,0) node (b1) {\Catg(\EG,\C)}
(a1)++(0,\v) node (a2) {\C}
(b1)++(0,\v) node (b2) {\C}
;
\draw[1cell=.9]
(a1) edge node {g} (b1)
(a2) edge node[swap] {g} (b2)
(a1) edge node[swap] {\pn} (a2)
(b1) edge node {\pn} (b2)
;
\draw[2cell]
node[between=a1 and b2 at .55, rotate=-135, 2labelw={below,\ppng,-1pt}] {\Rightarrow}
;
\end{tikzpicture}
\end{equation}
as follows.
\begin{description}
\item[Domain]
For a functor $f \cn \EG \to \C$, by \cref{pn_h} and the conjugation $G$-action on $\Catg(\EG,\C)$, the domain of $\ppng_f$ is the object
\begin{equation}\label{ppng_h_dom}
\pn(gf\ginv) = (gf\ginv)(e) = gf(\ginv) \in \C.
\end{equation}
This is the $g$-action on the object $f(\ginv) \in \C$.
\item[Codomain]
The codomain of $\ppng_f$ is the object
\begin{equation}\label{ppng_h_cod}
g\pn(f) = gf(e) \in \C.
\end{equation}
This is the $g$-action on the object $f(e) \in \C$.
\item[Components]
The $f$-component of $\ppng$ is defined as the isomorphism
\begin{equation}\label{ppng_h}
\pn(gf\ginv) = gf(\ginv) \fto[\iso]{\ppng_f = gf[e,\ginv]} gf(e) = g\pn(f)
\end{equation}
in $\C$.  This is the $g$-action on the isomorphism
\[f(\ginv) \fto[\iso]{f[e,\ginv]} f(e) \inspace \C,\]
which, in turn, is the $f$-image of the unique isomorphism 
\[\ginv \fto[\iso]{[e,\ginv]} e \inspace \EG.\]
\end{description}
This finishes the definition of $\ppng$.  We denote the collection $\{\ppng\}_{g \in G}$ by $\ppn$.
\cref{ppn_pseudoG} proves that $(\pn,\ppn)$ is a pseudo $G$-equivariant functor \pcref{def:pseudoG}. 
\end{definition}

\begin{explanation}[Nontrivial Constraints]\label{expl:ppng}
The pseudo $G$-equivariant constraints $\ppng$ \pcref{def:pn_pseudo} are not identities in general.  Indeed, for a nontrivial group $G$ and an element $g \in G$, the isomorphism $[e,\ginv] \cn \ginv \fiso e$ in $\EG$ is not an identity morphism.  Thus, for a general small $G$-category $\C$ and a functor $f \cn \EG \to \C$, the isomorphism $\ppng_f$ \cref{ppng_h}, which is defined as $gf[e,\ginv]$, is not an identity morphism.  
\end{explanation}

\begin{proposition}\label{ppn_pseudoG}
For a small $G$-category $\C$, the pair \pcref{def:pn_pseudo}
\[\Catg(\EG,\C) \fto{(\pn,\ppn)} \C\]
is a pseudo $G$-equivariant functor.
\end{proposition}

\begin{proof}
We verify that $\ppng$ is a natural isomorphism for each $g \in G$ and that $(\pn,\ppn)$ satisfies the two axioms in \cref{def:pseudoG}.
\begin{description}
\item[Naturality]
Each component of $\ppng$ \cref{ppng_h} is an isomorphism.  By \crefrange{ppng_h_dom}{ppng_h}, the naturality of $\ppng$ means the commutativity of the diagram \cref{ppng_nat_diag} in $\C$ for each natural transformation $\theta \cn f \to f'$ between functors $f,f' \cn \EG \to \C$.
\begin{equation}\label{ppng_nat_diag}
\begin{tikzpicture}[vcenter]
\def\v{-1.4}
\draw[0cell]
(0,0) node (a11) {gf(\ginv)}
(a11)++(3.5,0) node (a12) {gf(e)}
(a11)++(0,\v) node (a21) {gf'(\ginv)}
(a12)++(0,\v) node (a22) {gf'(e)}
;
\draw[1cell=.9]
(a11) edge node {gf[e,\ginv]} (a12)
(a21) edge node {gf'[e,\ginv]} (a22)
(a11) edge node[swap] {g\theta_{\ginv}} (a21)
(a12) edge node {g\theta_e} (a22)
;
\end{tikzpicture}
\end{equation}
The diagram \cref{ppng_nat_diag} commutes by the naturality of $\theta$ and the functoriality of the $g$-action on $\C$.
\item[Unity]
By \cref{pn_h,eeonee,ppng_h}, for the unit element $e \in G$ and a functor $f \cn \EG \to \C$, the $f$-component of $\ppne$ is given by
\[\ppne_f = ef[e,e] = f1_e = 1_{f(e)} = 1_{\pn(f)}.\]
The second equality uses the fact that the $e$-action on $\C$ is equal to $1_\C$.  The third equality uses the functoriality of $f$.  Thus, $\ppne$ is equal to the identity natural transformation $1_\pn$, proving the unity axiom \cref{pseudoG_unity} for $(\pn,\ppn)$.
\item[Multiplicativity]
The multiplicativity axiom \cref{pseudoG_mult} for $(\pn,\ppn)$ states that, for any two elements $g,k \in G$, the following pasting is equal to $\ppnkg$.
\begin{equation}\label{ppn_mult}
\begin{tikzpicture}[vcenter]
\def\h{3} \def\v{-1.3} 
\draw[0cell=.9]
(0,0) node (a11) {\Catg(\EG,\C)}
(a11)++(\h,0) node (a12) {\Catg(\EG,\C)}
(a12)++(\h,0) node (a13) {\Catg(\EG,\C)}
(a11)++(0,\v) node (a21) {\C}
(a12)++(0,\v) node (a22) {\C}
(a13)++(0,\v) node (a23) {\C}
;
\draw[1cell=.9]
(a11) edge node {g} (a12)
(a12) edge node {k} (a13)
(a21) edge node[swap] {g} (a22)
(a22) edge node[swap] {k} (a23)
(a11) edge node[swap] {\pn} (a21)
(a12) edge node {\pn} (a22)
(a13) edge node {\pn} (a23)
;
\draw[2cell=.9]
node[between=a11 and a22 at .55, rotate=-135, 2labelw={below,\ppng,0pt}] {\Rightarrow}
node[between=a12 and a23 at .55, rotate=-135, 2labelw={below,\ppnk,-1pt}] {\Rightarrow}
;
\end{tikzpicture}
\end{equation}
For a functor $f \cn \EG \to \C$, by \cref{ppng_h}, the pasting \cref{ppn_mult} yields the following composite isomorphism in $\C$.
\begin{equation}\label{ppn_mult_left}
kgf(\ginv\kinv) \fto{k(gf\ginv [e,\kinv])} kgf(\ginv) \fto{kgf[e,\ginv]} kgf(e)
\end{equation}
In the translation category $\EG$, there are morphism equalities as follows.
\begin{equation}\label{ginv_ekinv}
\begin{split}
\ginv[e,\kinv] &= [\ginv, \ginv\kinv] \cn \ginv\kinv \to \ginv e = \ginv\\
[e,\kginv] &= [e,\ginv] \circ [\ginv,\ginv\kinv] \cn \kginv \to e
\end{split}
\end{equation}
The equalities \cref{ginv_ekinv} and the functoriality of $f$ and of the $kg$-action on $\C$ imply that the composite in \cref{ppn_mult_left} is equal to 
\[\begin{split}
& kg\big(f[e,\ginv] \circ f[\ginv, \ginv\kinv]\big)\\
&= kgf\big( [e,\ginv] \circ [\ginv,\ginv\kinv] \big)\\
&= kgf[e,\kginv]\\
&= \ppnkg_f.
\end{split}\]
This proves the multiplicativity axiom \cref{pseudoG_mult} for $(\pn,\ppn)$.\qedhere
\end{description}
\end{proof}

\section{$G$-Thickening to Shimakawa $H$-Theory}
\label{sec:sgo_geq}

This section constructs a $G$-functor 
\[\Catg\big(\EG, \Catg(\EG,\A)\mal \big) \fto{\cni} \Catg(\EG,\A)\mal\]
to Shimakawa $H$-theory of the $\Catg(\EG,\Op)$-pseudoalgebra $\Catg(\EG,\A)$ at a pointed finite $G$-set $\mal$, along with its strong variant $\cnisg$.  \cref{thm:SgoAh} proves that $\cni$ is a $G$-equivariant inverse of the inclusion $G$-functor $\inc$ \cref{incC} for Shimakawa $H$-theory of the $\Catg(\EG,\Op)$-pseudoalgebra $\Catg(\EG,\A)$.  The strong variant is also true.  

\secoutline
\begin{itemize}
\item \cref{as:Ahat} states the assumptions for this section and \cref{sec:sgo_geq_unit,sec:sgosg_gequiv}.
\item \cref{def:cni} constructs the $G$-functors $\cni$ and $\cnisg$.
\item \cref{cnifs_welldef,glcnif_welldef,cnif_welldef,cnitha_welldef,cni_equivariant} prove that $\cni$ and $\cnisg$ are well defined.
\end{itemize}

\cref{as:Ahat} is in effect in \cref{sec:sgo_geq,sec:sgo_geq_unit,sec:sgosg_gequiv}.

\begin{assumption}\label{as:Ahat}
We consider a 1-connected $\Gcat$-operad $(\Op,\ga,\opu)$ \cref{i_connected} for a group $G$, an $\Op$-pseudoalgebra $(\A,\gaA,\phiA)$ \pcref{def:pseudoalgebra}, the 1-connected $\Gcat$-operad $\Oph = \Catg(\EG,\Op)$, and the $\Oph$-pseudoalgebra \pcref{catgego}
\begin{equation}\label{as:Ah}
\big(\Ah = \Catg(\EG,\A), \gaAh, \phiAh\big),
\end{equation}
where $\EG$ is the translation category of $G$ \pcref{def:translation_cat}.
\end{assumption}

\begin{notation}\label{not:fgsh}
For each pointed finite $G$-set $\mal \in \FG$ \pcref{def:ptGset,def:FG}, recall the small pointed $G$-categories  \pcref{def:nsys,def:nsys_morphism,def:nsys_gcat}  
\begin{equation}\label{sgohAhmal}
\begin{split}
(\Sgoh\Ah)(\mal) &= \Ahmal \andspace\\
(\Sgohsg\Ah)(\mal) &= \sgAhmal
\end{split}
\end{equation}
of (strong) $\mal$-systems in $\Ah$, where $\Sgoh$ and $\Sgohsg$ denote, respectively, Shimakawa $H$-theory for $\Oph$ and its strong variant \pcref{def:sgo}.  For a functor 
\begin{equation}\label{fEGAhmal}
\EG \fto{f} \Ahmal,
\end{equation}
elements $g,h \in G$, and a subset $s \subseteq \ufs{m} = \{1,2,\ldots,m\}$, we denote by
\begin{equation}\label{f_sub_g}
(f_g, \gl^{f_g}) = f(g) \in \Ahmal
\end{equation}
the image of $g$ under $f$; by
\begin{equation}\label{fgs}
f_{g,s} = (f_g)_s \in \Ah = \Catg(\EG,\A)
\end{equation}
the $s$-component object of $f_g$ \cref{nsys_s}; and by
\begin{equation}\label{fgsh}
f_{g,s,h} = (f_{g,s})(h) \in \A
\end{equation}
the image of $h$ under the functor $f_{g,s} \cn \EG \to \A$.  We use similar notation for morphisms and for the full subcategory $\sgAhmal \subseteq \Ahmal$ of strong $\mal$-systems in $\Ah$.
\end{notation}

The pseudo $G$-equivariant adjoint inverse $\pn$ \cref{proC} of the inclusion $G$-functor $\inc$ is defined by evaluating at the group unit $e \in G$.  In contrast, the $G$-functors $\cni$ and $\cnisg$ are given by the diagonal.  See \cref{cni_fsg,cnif_gl_g,cni_theta_sg}, where an element $g \in G$ is used twice in each case.

\begin{definition}\label{def:cni}
Under \cref{as:Ahat}, for each pointed finite $G$-set $\mal \in \FG$, the $G$-functor
\begin{equation}\label{cni_functor}
\Catg(\EG, \Ahmal) \fto{\cni} \Ahmal
\end{equation}
is defined as follows.  A strong variant is defined in \cref{cnisg_functor}.
\begin{description}
\item[Objects] For a functor $f \cn \EG \to \Ahmal$, the $\mal$-system in $\Ah$
\begin{equation}\label{cni_f}
(\cni f, \gl^{\cni f}) \in \Ahmal
\end{equation}
has component objects defined in \cref{cni_fs} and gluing morphisms defined in \cref{cnif_gl}.
\begin{description}
\item[Component objects]
For each subset $s \subseteq \ufs{m}$, the $s$-component object \cref{nsys_s}
\begin{equation}\label{cni_fs}
(\cni f)_s \in \Ah = \Catg(\EG,\A)
\end{equation}
is the functor $\EG \to \A$ with object assignment \cref{cni_fsg} and morphism assignment \cref{cni_fs_mor}.
\begin{description}
\item[Component objects on objects]
Using \cref{fgsh}, the functor $(\cni f)_s$ sends an object $g \in \EG$ to the object
\begin{equation}\label{cni_fsg}
(\cni f)_s(g) = f_{g,s,g} \in \A.
\end{equation}
\item[Component objects on morphisms]
For an isomorphism $[h,g] \cn g \fiso h$ in $\EG$, the isomorphism
\[(\cni f)_s(g) = f_{g,s,g} \fto[\iso]{(\cni f)_s [h,g]} (\cni f)_s(h) = f_{h,s,h} \inspace \A\]
is defined by the commutative diagram \cref{cni_fs_mor}.
\begin{equation}\label{cni_fs_mor}
\begin{tikzpicture}[vcenter]
\def\h{3.5} \def\v{1} \def\t{20}
\draw[0cell]
(0,0) node (a11) {f_{g,s,g}}
(a11)++(\h/2,\v) node (a12) {f_{h,s,g}}
(a11)++(\h/2,-\v) node (a21) {f_{g,s,h}}
(a11)++(\h,0) node (a22) {f_{h,s,h}}
;
\draw[1cell=.9]
(a11) edge[bend left=\t] node[pos=.3] {f_{[h,g],s,g}} (a12)
(a12) edge[bend left=\t] node[pos=.7] {f_{h,s,[h,g]}} (a22)
(a11) edge[bend right=\t] node[swap,pos=.3] {f_{g,s,[h,g]}} (a21)
(a21) edge[bend right=\t] node[swap,pos=.7] {f_{[h,g],s,h}} (a22)
(a11) edge node {(\cni f)_s [h,g]} (a22)
;
\end{tikzpicture}
\end{equation}
Using \cref{fgs}, the boundary diagram in \cref{cni_fs_mor} commutes by the naturality of the isomorphism
\[f_{g,s} \fto[\iso]{f_{[h,g],s}} f_{h,s} \inspace \Ah = \Catg(\EG,\A)\]
with respect to the isomorphism $[h,g]$ in $\EG$.
\end{description}
\cref{cnifs_welldef} proves that $(\cni f)_s \cn \EG \to \A$ is a functor.
\item[Gluing]
Given an object 
\[x \in \Oph(r) = \Catg(\EG,\Op(r))\] 
with $r \geq 0$, a subset $s \subseteq \ufsm$, and a partition 
\[s = \txcoprod_{\ell \in \ufs{r}}\, s_\ell \subseteq \ufs{m},\]
the gluing morphism \cref{gl-morphism} of $\cni f$ at $(x; s, \ang{s_\ell}_{\ell})$, where $\ang{\Cdots}_\ell = \ang{\Cdots}_{\ell \in \ufs{r}}$, is a morphism 
\begin{equation}\label{cnif_gl}
\gaAh_r\big(x; \ang{(\cni f)_{s_\ell}}_{\ell} \big) 
\fto{\gl^{\cni f}_{x;\, s, \ang{s_\ell}_{\ell}}} (\cni f)_s
\end{equation}
in $\Ah = \Catg(\EG,\A)$, meaning a natural transformation.  For each object $g \in \EG$, the $g$-component morphism of $\gl^{\cni f}_{x;\, s, \ang{s_\ell}_{\ell}}$ is defined by the diagram \cref{cnif_gl_g} in $\A$.
\begin{equation}\label{cnif_gl_g}
\begin{tikzpicture}[vcenter]
\def\h{4.5} \def\v{-1}
\draw[0cell=.9]
(0,0) node (a11) {\gaAh_r\big(x; \ang{(\cni f)_{s_\ell}}_{\ell} \big) (g)}
(a11)++(\h,0) node (a12) {(\cni f)_s (g)} 
(a12)++(0,3*\v) node (a22) {f_{g,s,g}}
(a11)++(0,\v) node (a21) {\gaA_r\big(x(g) ; \ang{(\cni f)_{s_\ell} (g)}_\ell \big)}
(a21)++(0,\v) node (a31) {\gaA_r\big(x(g) ; \ang{f_{g, s_\ell, g}}_\ell \big)}
(a31)++(0,\v) node (a32) {\gaAh_r\big(x; \ang{f_{g,s_\ell}}_\ell \big) (g)}
;
\draw[1cell=.9]
(a11) edge node {\gl^{\cni f}_{x;\, s, \ang{s_\ell}_{\ell}, g}} (a12)
(a32) edge node {\gl^{f_g}_{x;\, s, \ang{s_\ell}_\ell, g}} (a22)
(a12) edge[equal] node {\mathbf{d}^\cni} (a22)
(a11) edge[equal] node[swap] {\mathbf{d}} (a21)
(a21) edge[equal] node[swap] {\mathbf{d}^\cni} (a31)
(a31) edge[equal] node[swap] {\mathbf{d}} (a32)
;
\end{tikzpicture}
\end{equation}
The diagram \cref{cnif_gl_g} is defined as follows.
\begin{itemize}
\item Each of the two equalities labeled $\mathbf{d}$ holds by \cref{as:Ah}.
\item Each of the two equalities labeled $\mathbf{d}^\cni$ holds by \cref{cni_fsg}.
\item Using \cref{f_sub_g,fgs,fgsh}, the gluing morphism of the $\mal$-system $(f_g, \gl^{f_g}) \in \Ahmal$ at $(x; s, \ang{s_\ell}_\ell)$ is a morphism 
\[\gaAh_r\big(x; \ang{f_{g,s_\ell}}_\ell \big) \fto{\gl^{f_g}_{x;\, s, \ang{s_\ell}_\ell}} f_{g,s} \inspace \Ah,\]
meaning a natural transformation.  Its $g$-component morphism in $\A$ is the bottom arrow in \cref{cnif_gl_g}.
\end{itemize}
\cref{glcnif_welldef} proves that $\gl^{\cni f}_{x;\, s, \ang{s_\ell}_{\ell}}$ is a natural transformation.  \cref{cnif_welldef} proves that the pair $(\cni f, \gl^{\cni f})$ is an $\mal$-system in $\Ah$.  This finishes the definition of $\cni$ on objects.  
\end{description}
\item[Morphisms]
Suppose $\theta \cn f \to f'$ is a morphism in $\Catg(\EG,\Ahmal)$, meaning a natural transformation as follows.
\begin{equation}\label{theta_mor_EG}
\begin{tikzpicture}[vcenter]
\def\t{27}
\draw[0cell]
(0,0) node (a1) {\phantom{A}}
(a1)++(1.8,0) node (a2) {\phantom{A}}
(a1)++(-.1,0) node (a1') {\EG}
(a2)++(.2,0) node (a2') {\Ahmal}
;
\draw[1cell=.9]
(a1) edge[bend left=\t] node {f} (a2)
(a1) edge[bend right=\t] node[swap] {f'} (a2)
;
\draw[2cell]
node[between=a1 and a2 at .42, rotate=-90, 2label={above,\theta}] {\Rightarrow}
;
\end{tikzpicture}
\end{equation}
The morphism
\begin{equation}\label{cni_theta}
\cni f \fto{\cni\theta} \cni f' \inspace \Ahmal
\end{equation}
has, for each subset $s \subseteq \ufs{m}$, an $s$-component morphism \cref{theta_s} in $\Ah$, meaning a natural transformation as follows.
\begin{equation}\label{cni_theta_s}
\begin{tikzpicture}[vcenter]
\def\t{27}
\draw[0cell]
(0,0) node (a1) {\phantom{A}}
(a1)++(2.3,0) node (a2) {\A}
(a1)++(-.1,0) node (a1') {\EG}
;
\draw[1cell=.9]
(a1) edge[bend left=\t] node {(\cni f)_s} (a2)
(a1) edge[bend right=\t] node[swap] {(\cni f')_s} (a2)
;
\draw[2cell]
node[between=a1 and a2 at .33, rotate=-90, 2label={above,(\cni\theta)_s}] {\Rightarrow}
;
\end{tikzpicture}
\end{equation}
Using \cref{fgsh,cni_fsg}, for each object $g \in \EG$, the $g$-component of $(\cni\theta)_s$ is defined as the morphism 
\begin{equation}\label{cni_theta_sg}
(\cni f)_s(g) = f_{g,s,g} \fto{(\cni\theta)_{s,g} = \theta_{g,s,g}} (\cni f')_s(g) = f'_{g,s,g}
\end{equation}
in $\A$.  \cref{cnitha_welldef} proves that $(\cni\theta)_s$ is a natural transformation and that $\cni\theta$ is a morphism in $\Ahmal$.

\item[Functoriality]
The assignments
\[f \mapsto (\cni f, \gl^{\cni f}) \andspace \theta \mapsto \cni\theta\]
in \cref{cni_f,cni_theta} define a functor by \cref{cni_theta_sg} because identity morphisms and composition in $\Catg(\EG,-)$ and $\Ahmal$ \pcref{def:Catg,def:nsys_morphism} are defined componentwise.
\end{description}
This finishes the definition of the functor $\cni$ \cref{cni_functor}.  \cref{cni_equivariant} proves that $\cni$ is a $G$-functor.

\parhead{Strong variant}.  For the full subcategory $\sgAhmal \subseteq \Ahmal$ of strong $\mal$-systems in $\Ah$, the $G$-functor
\begin{equation}\label{cnisg_functor}
\Catg(\EG, \sgAhmal) \fto{\cnisg} \sgAhmal
\end{equation}
is defined by 
\begin{itemize}
\item \crefrange{cni_f}{cnif_gl_g} on objects and
\item \cref{cni_theta,cni_theta_s,cni_theta_sg} on morphisms.
\end{itemize}
This is well defined because, for each functor $f \cn \EG \to \sgAhmal$ and each object $g \in \EG$, $f_g \in \sgAhmal$ is now a strong $\mal$-system in $\Ah$.  Each component of its gluing morphism $\gl^{f_g}$ is an isomorphism in $\Ah = \Catg(\EG,\A)$.  Thus, the bottom arrow $\gl^{f_g}_{x;\, s, \ang{s_\ell}_\ell, g}$ in \cref{cnif_gl_g} is an isomorphism in $\A$.  This shows that $(\cnisg f, \gl^{\cnisg f})$ is a strong $\mal$-system in $\Ah$. 
\end{definition}

\subsection*{Proofs}
The rest of this section proves \cref{cnifs_welldef,glcnif_welldef,cnif_welldef,cnitha_welldef,cni_equivariant}, which are used in \cref{def:cni} to ensure that $\cni$ is a well-defined $G$-functor.

\begin{lemma}\label{cnifs_welldef}
In the context of \cref{def:cni}, the data
\[\EG \fto{(\cni f)_s} \A\]
in \cref{cni_fs,cni_fsg,cni_fs_mor} define a functor.
\end{lemma}

\begin{proof}
We verify that $(\cni f)_s$ preserves identity morphisms and composition.
\begin{description}
\item[Preservation of identities]  
Suppose $g \in \EG$ is an object with identity morphism $1_g$.  By \cref{cni_fs_mor}, the isomorphism $(\cni f)_s 1_g$ is defined as the composite
\[f_{g,s,g} \fto{f_{1_g,s,g}} f_{g,s,g} \fto{f_{g,s,1_g}} f_{g,s,g} \inspace \A.\]
The functoriality of $f \in \Catg(\EG,\Ahmal)$ implies that 
\[f_{1_g} = 1_{f_g} \inspace \Ahmal,\]
so $f_{1_g,s,g}$ is the identity morphism.  Similarly, the functoriality of $f_{g,s} \in \Ah = \Catg(\EG,\A)$ implies that $f_{g,s,1_g}$ is the identity morphism.
\item[Preservation of composition]  
For objects $g,h,k \in \EG$, applying $(\cni f)_s$ to the commutative diagram
\begin{equation}\label{hkg_diag}
\begin{tikzpicture}[vcenter]
\def\h{2}
\draw[0cell]
(0,0) node (a1) {g}
(a1)++(\h,0) node (a2) {h}
(a2)++(\h,0) node (a3) {k}
;
\draw[1cell=.9]
(a1) edge node {[h,g]} (a2)
(a2) edge node {[k,h]} (a3)
;
\draw[1cell=.9]
(a1) [rounded corners=2pt] |- ($(a2)+(-1,.6)$) -- node {[k,g]} ($(a2)+(1,.6)$) -| (a3)
 ;
\end{tikzpicture}
\end{equation}
in $\EG$ yields the following boundary diagram in $\A$.
\begin{equation}\label{fgskg}
\begin{tikzpicture}[vcenter]
\def\h{3} \def\v{-1.2} \def\t{40}
\draw[0cell]
(0,0) node (a11) {f_{g,s,g}}
(a11)++(\h,0) node (a12) {f_{g,s,k}}
(a12)++(\h,0) node (a13) {f_{k,s,k}}
(a11)++(\h/2,\v) node (a21) {f_{g,s,h}}
(a21)++(\h,0) node (a22) {f_{h,s,k}}
(a21)++(\h/2,\v) node (a3) {f_{h,s,h}}
;
\draw[1cell=.8]
(a11) edge node {f_{g,s,[k,g]}} (a12)
(a12) edge node {f_{[k,g],s,k}} (a13)
(a11) edge node[pos=.6, inner sep=0pt] {f_{g,s,[h,g]}} (a21)
(a21) edge node[swap, pos=.15, inner sep=0pt] {f_{g,s,[k,h]}} (a12)
(a12) edge node[swap, pos=.2, inner sep=0pt] {f_{[h,g],s,k}} (a22)
(a22) edge node[pos=.4, inner sep=0pt] {f_{[k,h],s,k}} (a13)
(a21) edge node[pos=.8, inner sep=0pt] {f_{[h,g],s,h}} (a3)
(a3) edge node[pos=.8, inner sep=0pt] {f_{h,s,[k,h]}} (a22)
;
\draw[1cell=.9]
(a11) [rounded corners=2pt] |- ($(a12)+(-1,.6)$) -- node {(\cni f)_s[k,g]} ($(a12)+(1,.6)$) -| (a13)
;
\draw[1cell=.9]
(a11) [rounded corners=2pt] |- node[pos=.7] {(\cni f)_s[h,g]} (a3)
;
\draw[1cell=.9]
(a3) [rounded corners=2pt] -| node[pos=.3] {(\cni f)_s[k,h]} (a13)
;
\end{tikzpicture}
\end{equation}
This diagram commutes for the following reasons.
\begin{itemize}
\item The three boundary regions commute by the definition \cref{cni_fs_mor} of $(\cni f)_s$ at a morphism.
\item The two upper triangles commute by the functoriality of $f_{g,s} \in \Ah$ and $f \in \Catg(\EG,\Ahmal)$.
\item The central quadrilateral commutes by the naturality of 
\[f_{g,s} \fto{f_{[h,g],s}} f_{h,s} \inspace \Ah\]
with respect to the isomorphism $[k,h] \cn h \fiso k$ in $\EG$.
\end{itemize}
This proves that $(\cni f)_s$ preserves composition.\qedhere
\end{description}
\end{proof}

\begin{lemma}\label{glcnif_welldef}
In the context of \cref{def:cni}, the data
\[\gaAh_r\big(x; \ang{(\cni f)_{s_\ell}}_{\ell} \big) 
\fto{\gl^{\cni f}_{x;\, s, \ang{s_\ell}_{\ell}}} (\cni f)_s\]
in \cref{cnif_gl,cnif_gl_g} define a morphism in $\Ah = \Catg(\EG,\A)$.
\end{lemma}

\begin{proof}
We verify that the morphism \cref{cnif_gl_g}
\[\gaAh_r\big(x; \ang{(\cni f)_{s_\ell}}_{\ell} \big) (g)
\fto{\gl^{\cni f}_{x;\, s, \ang{s_\ell}_{\ell}, g}} (\cni f)_s(g)\]
is natural in $g \in \EG$.  By \cref{cni_fs_mor}, \cref{cnif_gl_g}, and the functoriality of $\gaA_r$, the naturality diagram of $\gl^{\cni f}_{x;\, s, \ang{s_\ell}_{\ell}}$ with respect to an isomorphism $[h,g] \cn g \fiso h$ in $\EG$ is the following boundary diagram in $\A$.
\begin{equation}\label{glcnif_diag}
\begin{tikzpicture}[vcenter]
\def\v{-1.4} \def\t{2em}
\draw[0cell=.9]
(0,0) node (a11) {\gaA_r\big(x(g) ; \ang{f_{g,s_\ell,g}}_\ell \big)}
(a11)++(4,0) node (a12) {f_{g,s,g}}
(a11)++(0,\v) node (a21) {\gaA_r\big(x(h) ; \ang{f_{g,s_\ell,h}}_\ell \big)}
(a12)++(0,\v) node (a22) {f_{g,s,h}}
(a21)++(0,\v) node (a31) {\gaA_r\big(x(h) ; \ang{f_{h,s_\ell,h}}_\ell \big)}
(a22)++(0,\v) node (a32) {f_{h,s,h}}
;
\draw[1cell=.8]
(a11) edge node {\gl^{f_g}_{x;\, s, \ang{s_\ell}_\ell, g}} (a12)
(a21) edge node {\gl^{f_g}_{x;\, s, \ang{s_\ell}_\ell, h}} (a22)
(a31) edge node {\gl^{f_h}_{x;\, s, \ang{s_\ell}_\ell, h}} (a32)
(a12) edge node {f_{g,s,[h,g]}} (a22)
(a22) edge node {f_{[h,g],s,h}} (a32)
(a11) edge[transform canvas={xshift=\t}] node[swap] {\gaA_r(x[h,g] ; \ang{f_{g,s_\ell,[h,g]}}_\ell)} (a21)
(a21) edge[transform canvas={xshift=\t}] node[swap] {\gaA_r(1 ; \ang{f_{[h,g], s_\ell, h}}_\ell)} (a31)
;
\end{tikzpicture}
\end{equation}
The top rectangle commutes by the naturality of the gluing morphism
\[\gaAh_r\big(x; \ang{f_{g,s_\ell}}_\ell\big) \fto{\gl^{f_g}_{x;\, s, \ang{s_\ell}_\ell}} f_{g,s} \inspace \Ah\]
with respect to the isomorphism $[h,g]$.  The bottom rectangle commutes because it is the $h$-component of the compatibility diagram \cref{nsys_mor_compat} for the morphism
\[f_g \fto{f_{[h,g]}} f_h \inspace \Ahmal.\]
This proves that $\gl^{\cni f}_{x;\, s, \ang{s_\ell}_{\ell}}$ is a morphism in $\Ah$.
\end{proof}

\begin{lemma}\label{cnif_welldef}
In the context of \cref{def:cni}, the pair \cref{cni_f} 
\[(\cni f, \gl^{\cni f})\]
is an $\mal$-system in $\Ah = \Catg(\EG,\A)$.
\end{lemma}

\begin{proof}
We verify that the pair $(\cni f, \gl^{\cni f})$ satisfies the axioms \crefrange{nsys_obj_unity}{nsys_associativity} for an $\mal$-system in $\Ah$.
\begin{description}
\item[Object unity]  
By \cref{cni_fsg}, for the empty subset $\emptyset \subseteq \ufs{m}$, the functor
\[\EG \fto{(\cni f)_\emptyset} \A\]
sends an object $g \in \EG$ to the object 
\[(\cni f)_{\emptyset}(g) = f_{g,\emptyset,g} = \zero \in \A.\]
The last equality follows from
\begin{itemize}
\item the object unity axiom \cref{nsys_obj_unity} for the $\mal$-system $(f_g, \gl^{f_g}) \in \Ahmal$ and 
\item the fact that the basepoint of $\Ah = \Catg(\EG,\A)$ is the constant functor at the basepoint $\zero = \gaA_0(*) \in \A$.
\end{itemize}
Thus, $(\cni f)_\emptyset$ is the constant functor at $\zero \in \A$, which is the basepoint of $\Ah$.  This proves the object unity axiom \cref{nsys_obj_unity} for $\cni f$.
\item[Other axioms]  
Each of the remaining axioms,  \crefrange{nsys_naturality}{nsys_associativity}, for $(\cni f, \gl^{\cni f})$ asserts the commutativity of a diagram in $\Ah = \Catg(\EG,\A)$, which is verified objectwise for $g \in \EG$.  By \cref{cnif_gl_g}, each of these axioms, when evaluated at an object $g \in \EG$, is the $g$-component of the same axiom for the $\mal$-system $(f_g, \gl^{f_g}) \in \Ahmal$.\qedhere
\end{description}
\end{proof}

\begin{lemma}\label{cnitha_welldef}
In the context of \cref{def:cni}, the following statements hold.
\begin{enumerate}
\item\label{cnitha_i} The data \cref{cni_theta_s} 
\begin{equation}\label{cnitha_wd}
\begin{tikzpicture}[vcenter]
\def\t{27}
\draw[0cell]
(0,0) node (a1) {\phantom{A}}
(a1)++(2.3,0) node (a2) {\A}
(a1)++(-.1,0) node (a1') {\EG}
;
\draw[1cell=.9]
(a1) edge[bend left=\t] node {(\cni f)_s} (a2)
(a1) edge[bend right=\t] node[swap] {(\cni f')_s} (a2)
;
\draw[2cell]
node[between=a1 and a2 at .33, rotate=-90, 2label={above,(\cni\theta)_s}] {\Rightarrow}
;
\end{tikzpicture}
\end{equation}
define a natural transformation.
\item\label{cnitha_ii} The data \cref{cni_theta} 
\[\cni f \fto{\cni\theta} \cni f'\]
define a morphism in $\Ahmal$.
\end{enumerate}
\end{lemma}

\begin{proof}
\begin{description}
\item[\cref{cnitha_i}]  
By \cref{cni_fs_mor,cni_theta_sg}, the naturality diagram of $(\cni\theta)_s$ with respect to an isomorphism $[h,g] \cn g \fiso h$ in $\EG$ is the following boundary diagram in $\A$.
\begin{equation}\label{thaffprime}
\begin{tikzpicture}[vcenter]
\def\v{-1.4}
\draw[0cell=1]
(0,0) node (a11) {f_{g,s,g}}
(a11)++(3,0) node (a12) {f'_{g,s,g}}
(a11)++(0,\v) node (a21) {f_{g,s,h}}
(a12)++(0,\v) node (a22) {f'_{g,s,h}}
(a21)++(0,\v) node (a31) {f_{h,s,h}}
(a22)++(0,\v) node (a32) {f'_{h,s,h}}
;
\draw[1cell=.9]
(a11) edge node {\theta_{g,s,g}} (a12)
(a21) edge node {\theta_{g,s,h}} (a22)
(a31) edge node {\theta_{h,s,h}} (a32)
(a11) edge node[swap] {f_{g,s,[h,g]}} (a21)
(a21) edge node[swap] {f_{[h,g],s,h}} (a31)
(a12) edge node {f'_{g,s,[h,g]}} (a22)
(a22) edge node {f'_{[h,g],s,h}} (a32)
;
\end{tikzpicture}
\end{equation}
The top rectangle is the naturality diagram of the morphism
\[f_{g,s} \fto{\theta_{g,s}} f'_{g,s} \inspace \Ah = \Catg(\EG,\A)\]
with respect to $[h,g]$.  The bottom rectangle is the $(s,h)$-component of the naturality diagram of the morphism
\[f \fto{\theta} f' \inspace \Catg(\EG,\Ahmal)\]
with respect to $[h,g]$.  This proves assertion \eqref{cnitha_i}.
\item[\cref{cnitha_ii}]  
We verify that $\cni\theta \cn \cni f \to \cni f'$ satisfies the unity axiom \cref{nsys_mor_unity} and the compatibility axiom \cref{nsys_mor_compat} for a morphism of $\mal$-systems in $\Ah$.
\begin{description}
\item[Unity]  
For the empty subset $\emptyset \subseteq \ufs{m}$ and each object $g \in \EG$, by \cref{cni_theta_sg}, there are morphism equalities
\[(\cni\theta)_{\emptyset,g} = \theta_{g,\emptyset,g} = 1_\zero \inspace \A.\]
The second equality follows from
\begin{itemize}
\item the unity axiom \cref{nsys_mor_unity} for the morphism
\[f_g \fto{\theta_g} f'_g \inspace \Ahmal\]
and
\item the fact that the basepoint of $\Ah = \Catg(\EG,\A)$ is the constant functor at the basepoint $\zero = \gaA_0(*) \in \A$.
\end{itemize}
This proves the unity axiom \cref{nsys_mor_unity} for $\cni\theta$.
\item[Compatibility]  
For an object $x \in \Oph(r)$ with $r \geq 0$, a subset $s \subseteq \ufsm$, and a partition $s = \txcoprod_{\ell \in \ufs{r}}\, s_\ell \subseteq \ufs{m}$, the compatibility diagram \cref{nsys_mor_compat} for $\cni\theta$ is the following diagram in $\Ah$, where $\ang{\Cdots}_\ell = \ang{\Cdots}_{\ell \in \ufs{r}}$.
\begin{equation}\label{gaAh_diag}
\begin{tikzpicture}[vcenter]
\def\v{-1.4} \def\t{2em}
\draw[0cell=.9]
(0,0) node (a11) {\gaAh_r\big(x; \ang{(\cni f)_{s_\ell}}_\ell \big)}
(a11)++(4,0) node (a12) {(\cni f)_s}
(a11)++(0,\v) node (a21) {\gaAh_r\big(x; \ang{(\cni f')_{s_\ell}}_\ell \big)}
(a12)++(0,\v) node (a22) {(\cni f')_s}
;
\draw[1cell=.85]
(a11) edge node {\gl^{\cni f}_{x;\, s, \ang{s_\ell}_\ell}} (a12)
(a21) edge node {\gl^{\cni f'}_{x;\, s, \ang{s_\ell}_\ell}} (a22)
(a11) edge[transform canvas={xshift=\t}] node[swap] {\gaAh_r(1; \ang{(\cni\theta)_{s_\ell}}_\ell)} (a21)
(a12) edge node {(\cni\theta)_s} (a22)
;
\end{tikzpicture}
\end{equation}
By \cref{cnif_gl_g,cni_theta_sg}, evaluating the previous diagram at an object $g \in \EG$ yields the following diagram in $\A$.
\begin{equation}\label{gaA_diag}
\begin{tikzpicture}[vcenter]
\def\v{-1.4} \def\t{2em}
\draw[0cell=.9]
(0,0) node (a11) {\gaA_r\big(x(g); \ang{f_{g,s_\ell,g}}_\ell \big)}
(a11)++(4,0) node (a12) {f_{g,s,g}}
(a11)++(0,\v) node (a21) {\gaA_r\big(x(g); \ang{f'_{g,s_\ell,g}}_\ell \big)}
(a12)++(0,\v) node (a22) {f'_{g,s,g}}
;
\draw[1cell=.85]
(a11) edge node {\gl^{f_g}_{x;\, s, \ang{s_\ell}_\ell, g}} (a12)
(a21) edge node {\gl^{f'_g}_{x;\, s, \ang{s_\ell}_\ell, g}} (a22)
(a11) edge[transform canvas={xshift=\t}] node[swap] {\gaA_r(1; \ang{\theta_{g,s_\ell,g}}_\ell)} (a21)
(a12) edge node {\theta_{g,s,g}} (a22)
;
\end{tikzpicture}
\end{equation}
This diagram is the compatibility diagram \cref{nsys_mor_compat} for the morphism $\theta_g \cn f_g \to f'_g$ in $\Ahmal$, evaluated at the object $g \in \EG$.  This proves the compatibility axiom for $\cni\theta$ and finishes the proof of assertion \eqref{cnitha_ii}.\qedhere
\end{description}
\end{description}
\end{proof}

\begin{lemma}\label{cni_equivariant}
In the context of \cref{def:cni}, the functor \cref{cni_functor}
\[\Catg(\EG, \Ahmal) \fto{\cni} \Ahmal\]
is $G$-equivariant.
\end{lemma}

\begin{proof}
For each $g \in G$, we verify the commutativity of the diagram
\begin{equation}\label{cni_eq_diag}
\begin{tikzpicture}[vcenter]
\def\v{-1.4} 
\draw[0cell=1]
(0,0) node (a11) {\Catg(\EG,\Ahmal)}
(a11)++(3,0) node (a12) {\Ahmal}
(a11)++(0,\v) node (a21) {\Catg(\EG,\Ahmal)}
(a12)++(0,\v) node (a22) {\Ahmal}
;
\draw[1cell=.9]
(a11) edge node {\cni} (a12)
(a21) edge node {\cni} (a22)
(a11) edge node[swap] {g} (a21)
(a12) edge node {g} (a22)
;
\end{tikzpicture}
\end{equation}
on objects and morphisms of $\Catg(\EG,\Ahmal)$.  
\begin{description}
\item[Objects]  
For a functor $f \cn \EG \to \Ahmal$, the object equality
\begin{equation}\label{g_cnif}
g \cdot \cni f = \cni(g \cdot f) \inspace \Ahmal
\end{equation}
is proved by the following statements \eqref{g_cnif_i} and \eqref{g_cnif_ii}.
\begin{enumerate}
\item\label{g_cnif_i} For each subset $s \subseteq \ufs{m}$, there is an equality of $s$-component objects
\begin{equation}\label{g_cnif_s}
(g \cdot \cni f)_s = \big(\cni(g \cdot f)\big)_s \inspace \Ah = \Catg(\EG,\A).
\end{equation}
This means that the functors $(g \cdot \cni f)_s$ and $\big(\cni(g \cdot f)\big)_s$ are equal on both the objects and morphisms of $\EG$.
\item\label{g_cnif_ii} The corresponding gluing morphisms of $g \cdot \cni f$ and $\cni(g \cdot f)$ are equal.
\end{enumerate}
\begin{description}
\item[\cref{g_cnif_i}: objects of $\EG$] 
The desired equality \cref{g_cnif_s} at an object $h \in \EG$ is proved by the following object equalities in $\A$, which are explained further below.
\setlength{\fboxrule}{.1pt}
\setlength{\fboxsep}{1.3pt}
\begin{equation}\label{cni_eq_objects}
\begin{split}
(g \cdot \cni f)_s (h) 
&\seq{1} \big(g \cdot (\cni f)_{\ginv s} \big) (h) 
\seq{2} g (\cni f)_{\ginv s} (\ginv h) \\
&\seq{3} g f_{\ginv h, \ginv s, \ginv h} 
\seq{4} (g \cdot f_{\ginv h, \ginv s})_h \\
&\seq{5} (g \cdot f_{\ginv h})_{s,h}
\seq{6} (g \cdot f)_{h,s,h} \\
&\seq{7} \big(\cni(g \cdot f)\big)_s (h)
\end{split}
\end{equation}
The equalities in \cref{cni_eq_objects} hold for the following reasons.
\begin{itemize}
\item $\fbox{1}$ and $\fbox{5}$ hold by the definition \cref{ga_scomp} of the $G$-action on $\mal$-systems.
\item $\fbox{2}$, $\fbox{4}$, and $\fbox{6}$ hold by the definition \cref{conjugation-gaction} of the conjugation $G$-action on $\Catg(\EG,-)$.
\item $\fbox{3}$ and $\fbox{7}$ hold by the definition \cref{cni_fsg} of $(\cni f)_s(g)$.
\end{itemize}
This proves that $(g \cdot \cni f)_s$ and $(\cni (g \cdot f))_s$ are equal at each object of $\EG$.
\item[\cref{g_cnif_i}: morphisms of $\EG$] 
The desired equality \cref{g_cnif_s} at an isomorphism $[k,h] \cn h \fiso k$ in $\EG$ is proved by the following morphism equalities in $\A$, which are analogous to those in \cref{cni_eq_objects}.
\setlength{\fboxrule}{.1pt}
\setlength{\fboxsep}{1.3pt}
\begin{equation}\label{cni_eq_mor}
\begin{split}
& (g \cdot \cni f)_s  [k,h] \\ 
&\seq{1} \big(g \cdot (\cni f)_{\ginv s} \big) [k,h] \\
&\seq{2} g (\cni f)_{\ginv s} [\ginv k, \ginv h] \\
&\seq{3} \big(g f_{\ginv[k,h], \ginv s, \ginv k} \big) \circ \big( g f_{\ginv h, \ginv s, \ginv[k,h]} \big) \\
&\seq{4} (g \cdot f_{\ginv[k,h], \ginv s} )_k \circ (g \cdot f_{\ginv h, \ginv s} )_{[k,h]} \\
&\seq{5} (g \cdot f_{\ginv [k,h]})_{s,k} \circ (g \cdot f_{\ginv h})_{s, [k,h]} \\
&\seq{6} (g \cdot f)_{[k,h], s, k} \circ (g \cdot f)_{h,s,[k,h]} \\
&\seq{7} \big(\cni(g \cdot f)\big)_s [k,h]
\end{split}
\end{equation}
The equalities in \cref{cni_eq_mor} hold for the following reasons.
\begin{itemize}
\item $\fbox{1}$ and $\fbox{5}$ hold by the definitions \cref{ga_scomp,gtha_s} of the $G$-action on $\mal$-systems and their morphisms.
\item $\fbox{2}$, $\fbox{4}$, and $\fbox{6}$ hold by the definition \cref{conjugation-gaction} of the conjugation $G$-action on $\Catg(\EG,-)$.  $\fbox{2}$ also uses the equality
\begin{equation}\label{ginvkh}
\ginv[k,h] = [\ginv k, \ginv h] \cn \ginv h \fiso \ginv k
\end{equation}
in the translation category $\EG$ \pcref{def:translation_cat}.
\item $\fbox{3}$ and $\fbox{7}$ hold by the definition \cref{cni_fs_mor} of $(\cni f)_s [k,h]$.  $\fbox{3}$ also uses the functoriality of the $g$-action on $\A$ and \cref{ginvkh}.
\end{itemize}
This proves that $(g \cdot \cni f)_s$ and $(\cni (g \cdot f))_s$ are equal at each morphism of $\EG$, proving the equality of $s$-component objects in \cref{g_cnif_s}.
\item[\cref{g_cnif_ii}: gluing]
It suffices to verify that the gluing morphisms of $g \cdot \cni f$ and $\cni(g \cdot f)$ in $\Ah = \Catg(\EG,\A)$ are equal at each object of $\EG$.  For each object $x \in \Oph(r)$ with $r \geq 0$, subset $s \subseteq \ufsm$, partition $s = \txcoprod_{\ell \in \ufs{r}}\, s_\ell \subseteq \ufs{m}$, and object $h \in \EG$, there are morphism equalities in $\A$ as follows, where $\ang{\Cdots}_\ell = \ang{\Cdots}_{\ell \in \ufs{r}}$.
\setlength{\fboxrule}{.1pt}
\setlength{\fboxsep}{1.3pt}
\begin{equation}\label{cni_eq_gl}
\begin{split}
\gl^{(g \cdot \cni f)}_{x;\, s, \ang{s_\ell}_\ell, h} 
&\seq{1} \big(g \cdot \gl^{\cni f}_{\ginv \cdot x;\, \ginv s, \ang{\ginv s_\ell}_\ell} \big)_h \\
&\seq{2} g \gl^{\cni f}_{\ginv \cdot x;\, \ginv s, \ang{\ginv s_\ell}_\ell, \ginv h} \\
&\seq{3} g \gl^{f_{\ginv h}}_{\ginv \cdot x;\, \ginv s, \ang{\ginv s_\ell}_\ell, \ginv h} \\
&\seq{4} \big(g \cdot \gl^{f_{\ginv h}}_{\ginv \cdot x;\, \ginv s, \ang{\ginv s_\ell}_\ell} \big)_h \\
&\seq{5} \gl^{g \cdot f_{\ginv h}}_{x;\, s, \ang{s_\ell}_\ell, h} \\
&\seq{6} \gl^{(g \cdot f)_h}_{x;\, s, \ang{s_\ell}_\ell, h} 
\seq{7} \gl^{\cni(g \cdot f)}_{x;\, s, \ang{s_\ell}_\ell, h} 
\end{split}
\end{equation}
The equalities in \cref{cni_eq_gl} hold for the following reasons.
\begin{itemize}
\item $\fbox{1}$ and $\fbox{5}$ hold by the definition \cref{ga_gl} of the gluing morphisms of the $G$-action on $\mal$-systems.
\item $\fbox{2}$, $\fbox{4}$, and $\fbox{6}$  hold by the definition \cref{conjugation-gaction} of the conjugation $G$-action on $\Catg(\EG,-)$.  
\item $\fbox{3}$ and $\fbox{7}$ hold by the definition \cref{cnif_gl_g} of $\gl^{\cni f}$.
\end{itemize}
This proves that $g \cdot \cni f$ and $\cni(g \cdot f)$ have the same gluing morphisms, proving the object equality in \cref{g_cnif}.
\end{description}
\item[Morphisms]  
Two morphisms in $\Ahmal$ are equal if they are equal in $\A$ after evaluating at each subset $s \subseteq \ufs{m}$ and each object $h \in \EG$.  The commutativity of the diagram \cref{cni_eq_diag} on morphisms is proved by the computation \cref{cni_eq_objects} with the following adjustment.
\setlength{\fboxrule}{.1pt}
\setlength{\fboxsep}{1.3pt}
\begin{itemize}
\item $f$ is replaced by a morphism $\theta \cn f \to f'$ in $\Catg(\EG,\Ahmal)$.
\item For $\fbox{1}$ and $\fbox{5}$, \cref{gtha_s} is used for the $G$-action on morphisms of $\mal$-systems.
\item For $\fbox{3}$ and $\fbox{7}$, \cref{cni_theta_sg} is used to interpret $(\cni \theta)_{s,h}$ as $\theta_{h,s,h}$.
\end{itemize}
\end{description}
This proves that the diagram \cref{cni_eq_diag} commutes.
\end{proof}

\section{$G$-Equivariant Unit and Counit}
\label{sec:sgo_geq_unit}

This section constructs the $G$-equivariant unit $\ucni$ and counit $\ccni$ for the $G$-functors
\begin{equation}\label{inc_cni_seci}
\begin{tikzpicture}[vcenter]
\def\h{1.8} \def\t{20}
\draw[0cell]
(0,0) node (a1) {\phantom{A}}
(a1)++(\h,0) node (a2) {\phantom{A}}
(a1)++(-.2,0) node (a1') {\Ahmal}
(a2)++(1,0) node (a2') {\Catg(\EG, \Ahmal)}
;
\draw[1cell=.9]
(a1) edge[bend left=\t] node {\inc^{\Ahmal}} (a2)
(a2) edge[bend left=\t] node {\cni} (a1)
;
\end{tikzpicture}
\end{equation}
in \cref{def:gcat_inc,def:cni}.  \cref{thm:SgoAh} proves that the quadruple $(\inc^{\Ahmal}, \cni, \ucni, \ccni)$ is an adjoint $G$-equivalence.  The strong variant, involving $\sgAhmal$ and $\cnisg$, is also true.  \cref{as:Ahat} is in effect throughout this section, so $\A$ is an $\Op$-pseudoalgebra for a 1-connected $\Gcat$-operad $\Op$, $\Oph = \Catg(\EG,\Op)$, and $\Ah = \Catg(\EG,\A)$.

\secoutline
\begin{itemize}
\item \cref{inccni_unit_dom} proves that $\cni\inc^{\Ahmal}$ and $\cnisg\inc^{\sgAhmal}$ are identity functors of, respectively, $\Ahmal$ and $\sgAhmal$.
\item \cref{def:inccni_unit} defines the unit 
\[1_{\Ahmal} \fto{\ucni} \cni\inc^{\Ahmal}\]
as the identity natural transformation, and similarly for the strong variant.
\item \cref{def:inccni_counit} defines the counit
\[\inc^{\Ahmal}\cni \fto{\ccni} 1_{\Catg(\EG,\,\Ahmal)}\]
and its strong variant.
\item \cref{expl:ccni_domain} describes the domain $\inc^{\Ahmal}\cni$ of $\ccni$ explicitly.
\item \cref{ccni_welldef} proves that $\ccni$ and its strong variant are $G$-natural isomorphisms.
\end{itemize}

\subsection*{Unit}
\cref{inccni_unit_dom} is needed to define the $G$-equivariant unit for the $G$-functors $(\inc^{\Ahmal}, \cni)$ in \cref{def:gcat_inc,def:cni}.

\begin{lemma}\label{inccni_unit_dom}
The following two diagrams of $G$-functors commute.
\begin{equation}\label{inc_cnisg_diag}
\begin{tikzpicture}[vcenter]
\def\h{2.5} \def\t{20}
\draw[0cell]
(0,0) node (a1) {\Ahmal}
(a1)++(\h/2,-1.2) node (a2) {\Catg(\EG,\Ahmal)}
(a1)++(\h,0) node (a3) {\Ahmal}
;
\draw[1cell=.9]
(a1) edge[bend left=0] node {1} (a3)
(a1) edge[bend right=\t, shorten >=-.5ex] node[swap,pos=.2] {\inc^{\Ahmal}} (a2)
(a2) edge[bend right=\t, shorten <=-.5ex] node[swap,pos=.8] {\cni} (a3)
;
\begin{scope}[shift={(5,0)}]
\draw[0cell]
(0,0) node (a1) {\sgAhmal}
(a1)++(\h/2,-1.2) node (a2) {\Catg(\EG,\sgAhmal)}
(a1)++(\h,0) node (a3) {\sgAhmal}
;
\draw[1cell=.9]
(a1) edge[bend left=0] node {1} (a3)
(a1) edge[bend right=\t, shorten >=-.5ex] node[swap,pos=.2] {\inc^{\sgAhmal}} (a2)
(a2) edge[bend right=\t, shorten <=-.5ex] node[swap,pos=.8] {\cnisg} (a3)
;
\end{scope}
\end{tikzpicture}
\end{equation}
\end{lemma}

\begin{proof}
We prove the equality $\cni\inc^{\Ahmal} = 1$.  The proof for the equality $\cnisg\inc^{\sgAhmal} = 1$ is the same after restricting to strong $\mal$-systems in $\Ah = \Catg(\EG,\A)$.  To show that $\cni \inc^{\Ahmal}$ is the identity on objects, suppose $(a,\gl) \in \Ahmal$ is an object and $s \subseteq \ufs{m}$ is a subset.
\begin{description}
\item[Components and objects] 
For each object $g \in \EG$, there are object equalities in $\A$
\begin{equation}\label{ucni_cod_obj}
\big(\cni \inc^{\Ahmal} (a,\gl) \big)_s(g)
= \big(\inc^{\Ahmal} (a,\gl)\big)_{g,s,g}
= a_s(g).
\end{equation}
In \cref{ucni_cod_obj}, the first equality holds by \cref{cni_fsg}.  The second equality holds because, by \cref{def:gcat_inc}, $\inc^{\Ahmal} (a,\gl)$ is the constant functor at $(a,\gl)$, so $(\inc^{\Ahmal} (a,\gl))_g$ is $(a,\gl)$.  
\item[Components and morphisms] 
For an isomorphism $[h,g] \cn g \fiso h$ in $\EG$, \cref{cni_fs_mor} implies the following morphism equalities in $\A$.
\begin{equation}\label{ucni_cod_mor}
\begin{split}
&\big(\cni \inc^{\Ahmal} (a,\gl) \big)_s [h,g]\\
&= \big(\inc^{\Ahmal} (a,\gl) \big)_{h,s,[h,g]} \circ \big(\inc^{\Ahmal} (a,\gl) \big)_{[h,g], s, g} \\
&= a_s [h,g] \circ 1_{(a,\gl), s, g} \\
& = a_s [h,g]
\end{split}
\end{equation}
By \cref{ucni_cod_obj,ucni_cod_mor}, the $\mal$-system $\cni \inc^{\Ahmal} (a,\gl)$ has the same $s$-component object as $(a,\gl)$ for each subset $s \subseteq \ufs{m}$.
\item[Gluing] 
For each object $x \in \Oph(r)$ with $r \geq 0$, subset $s \subseteq \ufsm$, partition $s = \txcoprod_{\ell \in \ufs{r}}\, s_\ell \subseteq \ufs{m}$, and object $g \in \EG$, \cref{cnif_gl_g} implies the morphism equalities in $\A$
\begin{equation}\label{ucni_cod_gl}
\gl^{\cni\inc^{\Ahmal} (a,\gl)}_{x;\, s, \ang{s_\ell}_\ell, g}
= \gl^{(\inc^{\Ahmal} (a,\gl))_g}_{x;\, s, \ang{s_\ell}_\ell, g}
=  \gl_{x;\, s, \ang{s_\ell}_\ell, g}.
\end{equation}
By \cref{ucni_cod_gl}, $\cni \inc^{\Ahmal} (a,\gl)$ has the same gluing morphisms as $(a,\gl)$.  Thus, the composite $\cni \inc^{\Ahmal}$ is the identity on objects of $\Ahmal$.
\item[Morphisms]
For a morphism $\theta$ in $\Ahmal$, $\inc^{\Ahmal}\theta$ is the constant natural transformation at $\theta$, so $(\inc^{\Ahmal}\theta)_g = \theta$ for each object $g \in \EG$.  Replacing $(a,\gl)$ by $\theta$, the computation \cref{ucni_cod_obj} and \cref{cni_theta_sg} imply that the composite $\cni \inc^{\Ahmal}$ is the identity on morphisms of $\Ahmal$.
\end{description}
This proves that $\cni \inc^{\Ahmal}$ is the identity functor.
\end{proof}

\begin{definition}[Unit]\label{def:inccni_unit}
In the context of \cref{def:gcat_inc,def:cni}, the $G$-natural transformation
\begin{equation}\label{ucni}

\end{equation}
is defined as the identity natural transformation of the identity $G$-functor on $\Ahmal$, using \cref{inccni_unit_dom} for the equality $\cni \inc^{\Ahmal} = 1_{\Ahmal}$.

\parhead{Strong variant}.
The $G$-natural transformation
\begin{equation}\label{ucnisg}
%
\end{equation}
is defined as the identity natural transformation of the identity $G$-functor on $\sgAhmal$, using \cref{inccni_unit_dom} for the equality $\cnisg \inc^{\sgAhmal} = 1_{\sgAhmal}$.  
\end{definition}

\subsection*{Counit}
Next, we define the $G$-equivariant counit for the $G$-functors $(\inc^{\Ahmal}, \cni)$.

\begin{definition}[Counit]\label{def:inccni_counit}
In the context of \cref{def:gcat_inc,def:cni}, the $G$-natural isomorphism
\begin{equation}\label{ccni}
%
\end{equation}
is defined as follows.  A strong variant is defined in \cref{ccnisg}.  For a functor $f \cn \EG \to \Ahmal$, by \cref{def:gcat_inc}, the functor
\[\EG \fto{\inc^{\Ahmal} \cni f} \Ahmal\]
is constant at $\cni f \in \Ahmal$ \cref{cni_f}.  By \cref{cni_fsg}, for a subset $s \subseteq \ufs{m}$ and objects $g,h \in \EG$, its $(g,s,h)$-component object is given by
\begin{equation}\label{ccni_dom_obj}
(\inc^{\Ahmal} \cni f)_{g,s,h} = (\cni f)_s(h) = f_{h,s,h} \in \A.
\end{equation} 
The $f$-component of $\ccni$ is the natural isomorphism
\begin{equation}\label{ccni_f}
\begin{tikzpicture}[vcenter]
\def\t{35}
\draw[0cell]
(0,0) node (a1) {\phantom{A}}
(a1)++(2,0) node (a2) {\phantom{A}}
(a1)++(-.1,0) node (a1') {\EG}
(a2)++(.2,0) node (a2') {\Ahmal}
;
\draw[1cell=.9]
(a1) edge[bend left=\t] node {\inc^{\Ahmal}\cni f} (a2)
(a1) edge[bend right=\t] node[swap] {f} (a2)
;
\draw[2cell]
node[between=a1 and a2 at .4, rotate=-90, 2label={above,\ccni_f}] {\Rightarrow}
;
\end{tikzpicture}
\end{equation}
whose $(g,s,h)$-component isomorphism is defined as
\begin{equation}\label{ccni_f_gsh}
(\inc^{\Ahmal} \cni f)_{g,s,h} = f_{h,s,h}
\fto[\iso]{\ccni_{f,g,s,h} = f_{[g,h],s,h}} f_{g,s,h}
\end{equation}
in $\A$, with $[g,h] \cn h \fiso g$ the unique isomorphism in the translation category $\EG$.  \cref{ccni_welldef} proves that $\ccni$ is a $G$-natural isomorphism.

\parhead{Strong variant}.
The $G$-natural isomorphism
\begin{equation}\label{ccnisg}
\begin{tikzpicture}[vcenter]
\def\h{4.2} \def\t{25} \def\b{.8}
\draw[0cell]
(0,0) node (a1) {\Catg(\EG,\sgAhmal)}
(a1)++(\h/2,1.2) node (a2) {\sgAhmal}
(a1)++(\h,0) node (a3) {\Catg(\EG,\sgAhmal)}
(a1)++(\b,0) node (a1') {\phantom{\Ah}}
(a3)++(-\b,0) node (a3') {\phantom{\Ah}}
;
\draw[1cell=.9]
(a1) edge node[swap] {1} (a3)
(a1') edge[bend left=\t] node {\cnisg} (a2)
(a2) edge[bend left=\t] node {\inc^{\sgAhmal}} (a3')
;
\draw[2cell]
node[between=a1 and a3 at .45, shift={(0,.5)}, rotate=-90, 2label={above,\ccnisg}] {\Rightarrow}
;
\end{tikzpicture}
\end{equation}
is defined by applying \cref{ccni_f_gsh} to functors $f \cn \EG \to \sgAhmal$.  This is well defined because $\sgAhmal$ is a full subcategory of $\Ahmal$.
\end{definition}

\begin{explanation}[Domain of $\ccni$]\label{expl:ccni_domain}
We unravel the domain of $\ccni$ \cref{ccni}, which is the composite 
\[\Catg(\EG,\Ahmal) \fto{\cni} \Ahmal \fto{\inc^{\Ahmal}} \Catg(\EG,\Ahmal)\]
of the $G$-functors in \cref{def:gcat_inc,def:cni}.  
\begin{description}
\item[Objects]
For a functor $f \cn \EG \to \Ahmal$, the functor
\[\EG \fto{\inc^{\Ahmal} \cni f} \Ahmal\]
is constant at $\cni f \in \Ahmal$, which is defined in \crefrange{cni_f}{cnif_gl_g}. 
\begin{description}
\item[Component objects]  For a subset $s \subseteq \ufs{m}$ and objects $g,h \in \EG$, the $(g,s,h)$-component object of $\inc^{\Ahmal} \cni f$ is $f_{h,s,h}$ \cref{ccni_dom_obj}.
\item[Component morphisms] For an isomorphism $[k,h] \cn h \fiso k$ in $\EG$, its $(g,s,[k,h])$-component isomorphism in $\A$ is given by
\begin{equation}\label{ccni_dom_mor}
\begin{split}
& (\inc^{\Ahmal} \cni f)_{g,s,[k,h]} \\
&= (\cni f)_s [k,h] \\
&= f_{k,s,[k,h]} \circ f_{[k,h],s,h} \cn f_{h,s,h} \to f_{k,s,h} \to f_{k,s,k}.
\end{split}
\end{equation}
\item[Gluing] Given an object $x \in \Oph(r)$ with $r \geq 0$, a subset $s \subseteq \ufsm$, a partition $s = \txcoprod_{\ell \in \ufs{r}}\, s_\ell \subseteq \ufs{m}$, and objects $g,h \in \EG$, the gluing morphism of the $\mal$-system
\[(\inc^{\Ahmal} \cni f)_g  = \cni f\in \Ahmal\] 
at $(x; s, \ang{s_\ell}_\ell)$ has the following $h$-component morphism in $\A$, where $\ang{\Cdots}_\ell = \ang{\Cdots}_{\ell \in \ufsr}$. 
\begin{equation}\label{ccni_dom_gl}

\end{equation}
to the constant natural transformation
\begin{equation}\label{EGiptha}
%
\end{equation}
at $\cni\theta \cn \cni f \to \cni f'$, which is defined in \cref{cni_theta,cni_theta_s,cni_theta_sg}.  This means that, for each subset $s \subseteq \ufs{m}$ and objects $g,h \in \EG$, its $(g,s,h)$-component morphism in $\A$ is given as follows.
\begin{equation}\label{ccni_dom_nt}
%
\end{equation}
\end{description}
This finishes the unraveling of the domain of $\ccni \cn \inc^{\Ahmal} \cni \to 1$.

\parhead{Strong variant}. 
The domain of the strong variant $\ccnisg$ \cref{ccnisg}, which is the composite $G$-functor
\[\Catg(\EG,\sgAhmal) \fto{\cnisg} \sgAhmal \fto{\inc^{\sgAhmal}} \Catg(\EG,\sgAhmal),\]
admits the same description as \crefrange{ccni_dom_mor}{ccni_dom_nt} after changing the notation from $(\Ahmal, \cni)$ to $(\sgAhmal,\cnisg)$.
The gluing morphism $\gl^{f_h}_{x;\, s, \ang{s_\ell}_\ell, h}$ in \cref{ccni_dom_gl} is now an isomorphism because, for $f \in \Catg(\EG,\sgAhmal)$ and $h \in \EG$, $f_h$ is a strong $\mal$-system in $\Ah$.
\end{explanation}

The rest of this section proves \cref{ccni_welldef}, which is used in \cref{def:inccni_counit}.

\begin{lemma}\label{ccni_welldef}
In the context of \cref{def:inccni_counit}, the data 
\begin{equation}\label{ccni_statement}
\begin{tikzpicture}[vcenter]
\def\h{3} \def\t{28}
\draw[0cell]
(0,0) node (a1) {\phantom{A}}
(a1)++(2,0) node (a2) {\phantom{A}}
(a1)++(-1,0) node (a1') {\Catg(\EG,\Ahmal)}
(a2)++(1,0) node (a2') {\Catg(\EG,\Ahmal)}
;
\draw[1cell=.9]
(a1) edge[bend left=\t] node {\inc^{\Ahmal} \cni} (a2)
(a1) edge[bend right=\t] node[swap] {1} (a2)
;
\draw[2cell]
node[between=a1 and a2 at .42, rotate=-90, 2label={above,\ccni}] {\Rightarrow}
;
\end{tikzpicture}
\end{equation}
define a $G$-natural isomorphism.  Moreover, the strong variant $\ccnisg$ \cref{ccnisg} is also a $G$-natural isomorphism.
\end{lemma}

\begin{proof}
We only need to consider $\ccni$.  The proof for the strong variant $\ccnisg$ is obtained by changing the notation from $(\Ahmal, \cni)$ to $(\sgAhmal,\cnisg)$.  The $(f,g,s,h)$-component of $\ccni$ \cref{ccni_f_gsh} is the isomorphism
\[(\inc^{\Ahmal} \cni f)_{g,s,h} = f_{h,s,h}
\fto[\iso]{\ccni_{f,g,s,h} = f_{[g,h],s,h}} f_{g,s,h} \inspace \A.\]
To prove that $\ccni$ is a $G$-natural isomorphism, we need to prove the following statements \crefrange{ccni_def_i}{ccni_def_v}.
\begin{enumerate}
\item\label{ccni_def_i} As $h \in \EG$ varies, 
\begin{equation}\label{ccnifgs_iicell}
\begin{tikzpicture}[vcenter]
\def\t{30}
\draw[0cell]
(0,0) node (a1) {\phantom{A}}
(a1)++(2.5,0) node (a2) {\A}
(a1)++(-.1,0) node (a1') {\EG}
;
\draw[1cell=.9]
(a1) edge[bend left=\t] node {(\inc^{\Ahmal}\cni f)_{g,s}} (a2)
(a1) edge[bend right=\t] node[swap] {f_{g,s}} (a2)
;
\draw[2cell]
node[between=a1 and a2 at .35, rotate=-90, 2label={above,\ccni_{f,g,s}}] {\Rightarrow}
;
\end{tikzpicture}
\end{equation}
is a natural transformation.
\item\label{ccni_def_ii} As $s \subseteq \ufs{m}$ varies,
\[(\inc^{\Ahmal} \cni f)_g \fto{\ccni_{f,g}} f_g\]
is a morphism of $\mal$-systems in $\Ah = \Catg(\EG,\A)$ \pcref{def:nsys_morphism}.
\item\label{ccni_def_iii} As $g \in \EG$ varies, 
\begin{equation}\label{ccnif_iicell}
\begin{tikzpicture}[vcenter]
\def\t{33}
\draw[0cell]
(0,0) node (a1) {\phantom{A}}
(a1)++(2,0) node (a2) {\phantom{A}}
(a1)++(-.1,0) node (a1') {\EG}
(a2)++(.2,0) node (a2') {\Ahmal}
;
\draw[1cell=.9]
(a1) edge[bend left=\t] node {\inc^{\Ahmal}\cni f} (a2)
(a1) edge[bend right=\t] node[swap] {f} (a2)
;
\draw[2cell]
node[between=a1 and a2 at .4, rotate=-90, 2label={above,\ccni_f}] {\Rightarrow}
;
\end{tikzpicture}
\end{equation}
is a natural transformation.
\item\label{ccni_def_iv} As $f \in \Catg(\EG,\Ahmal)$ varies, 
\[\inc^{\Ahmal}\cni \fto{\ccni} 1_{\Catg(\EG,\,\Ahmal)}\]
is a natural transformation.
\item\label{ccni_def_v} The natural transformation $\ccni$ is $G$-equivariant.
\end{enumerate} 
In the rest of this proof, we abbreviate $\inc^{\Ahmal}$ to $\inc$.
\begin{description}
\item[\cref{ccni_def_i}] 
Given $f \in \Catg(\EG,\Ahmal)$, $g \in \EG$, and $s \subseteq \ufs{m}$, the naturality of $\ccni_{f,g,s}$ means that, for each isomorphism $[k,h] \cn h \fiso k$ in $\EG$, the following boundary diagram in $\A$ commutes.
\begin{equation}\label{ccni_i_diag}
\begin{tikzpicture}[vcenter]
\def\g{2.5} \def\v{-1.4} \def\u{-2.8}
\draw[0cell]
(0,0) node (a11) {(\inc\cni f)_{g,s,h}}
(a11)++(\g,0) node (a12) {f_{h,s,h}}
(a12)++(3,0) node (a13) {f_{g,s,h}}
(a12)++(0,\v) node (a22) {f_{h,s,k}}
(a13)++(0,\v) node (a23) {f_{g,s,k}} 
(a11)++(0,\u) node (a31) {(\inc\cni f)_{g,s,k}}
(a12)++(0,\u) node (a32) {f_{k,s,k}}
(a13)++(0,\u) node (a33) {f_{g,s,k}}
;
\draw[1cell=.9]
(a11) edge[equal] (a12)
(a12) edge node {f_{[g,h],s,h}} (a13)
(a13) edge node {f_{g,s,[k,h]}} (a23)
(a23) edge[equal] (a33)
(a11) edge node[swap] {(\inc\cni f)_{g,s,[k,h]}} (a31)
(a31) edge[equal] (a32)
(a32) edge node {f_{[g,k],s,k}} (a33)
(a12) edge node[swap] {f_{h,s,[k,h]}} (a22)
(a22) edge node[swap] {f_{[k,h],s,k}} (a32)
(a22) edge node {f_{[g,h],s,k}} (a23)
;
\end{tikzpicture}
\end{equation}
The three regions in the previous diagram commutes for the following reasons.
\begin{itemize}
\item The left rectangle commutes by the equality (\cref{def:gcat_inc})
\[(\inc\cni f)_{g,s,[k,h]} = (\cni f)_{s,[k,h]}\]
and \cref{cni_fs_mor}.
\item The upper right rectangle commutes by the naturality of 
\[f_{h,s} \fto{f_{[g,h],s}} f_{g,s} \inspace \Ah = \Catg(\EG,\A)\]
with respect to the isomorphism $[k,h]$.
\item The lower right rectangle commutes by the commutativity of the diagram
\begin{equation}\label{hkg}
\begin{tikzpicture}[vcenter]
\def\h{2}
\draw[0cell]
(0,0) node (a1) {h}
(a1)++(\h,0) node (a2) {k}
(a2)++(\h,0) node (a3) {g}
;
\draw[1cell=.9]
(a1) edge node {[k,h]} (a2)
(a2) edge node {[g,k]} (a3)
;
\draw[1cell=.9]
(a1) [rounded corners=2pt] |- ($(a2)+(-1,.6)$) -- node {[g,h]} ($(a2)+(1,.6)$) -| (a3)
 ;
\end{tikzpicture}
\end{equation}
in the translation category $\EG$ and the functoriality of $f$.
\end{itemize}
\item[\cref{ccni_def_ii}] 
To show that 
\[(\inc\cni f)_g \fto{\ccni_{f,g}} f_g\]
is a morphism of $\mal$-systems in $\Ah$, we verify the unity axiom \cref{nsys_mor_unity} and the compatibility axiom \cref{nsys_mor_compat}.
\begin{description}
\item[Unity]
For the unity axiom \cref{nsys_mor_unity}, the $(\emptyset,h)$-component of $\ccni_{f,g}$
\[(\inc\cni f)_{g,\emptyset,h} = f_{h,\emptyset,h}
\fto{\ccni_{f,g,\emptyset,h} = f_{[g,h],\emptyset,h}} f_{g,\emptyset,h}\]
is equal to $1_\zero$ in $\A$ by
\begin{itemize}
\item the unity axiom \cref{nsys_mor_unity} for the morphism
\[f_h \fto{f_{[g,h]}} f_g \inspace \Ahmal\]
and
\item the fact that the basepoint of $\Ah = \Catg(\EG,\A)$ is the constant functor at the basepoint $\zero \in \A$.
\end{itemize}
This proves the unity axiom \cref{nsys_mor_unity} for $\ccni_{f,g}$.
\item[Compatibility]
A diagram in $\Ah$ commutes if and only if it commutes in $\A$ after evaluating at each object $h \in \EG$.  For an object $x \in \Oph(r)$ with $r \geq 0$, a subset $s \subseteq \ufsm$, a partition $s = \txcoprod_{\ell \in \ufs{r}}\, s_\ell \subseteq \ufs{m}$, and an object $h \in \EG$, by \cref{ccni_dom_gl}, the compatibility diagram \cref{nsys_mor_compat} for $\ccni_{f,g}$ at $(x; s, \ang{s_\ell}_\ell, h)$ is the following diagram in $\A$.
\begin{equation}\label{gaAgl_compat}

\end{equation}
This diagram commutes because it is the compatibility diagram \cref{nsys_mor_compat} for the morphism $f_{[g,h]} \cn f_h \to f_g$ in $\Ahmal$ at $(x; s, \ang{s_\ell}_\ell, h)$.
\end{description}
\item[\cref{ccni_def_iii}] 
The naturality diagram of $\ccni_f \cn \inc\cni f \to f$ at an isomorphism $[g,k] \cn k \fiso g$ in $\EG$ is the diagram
\begin{equation}\label{ccnifk_diag}
%
\end{equation}
in $\Ahmal$.  This diagram commutes if and only if its $(s,h)$-component in $\A$ commutes for each subset $s \subseteq \ufs{m}$ and object $h \in \EG$.  Recalling that $\inc\cni f$ is constant at $\cni f$, the $(s,h)$-component of the previous diagram is the diagram
\begin{equation}\label{ipfgsh_diag}
%
\end{equation}
in $\A$.  This diagram commutes by \cref{hkg} and the functoriality of $f$.
\item[\cref{ccni_def_iv}]
The naturality diagram of $\ccni \cn \inc\cni \to 1$ at a natural transformation
\begin{equation}\label{EGthaffprime}
%
\end{equation}
is the diagram
\begin{equation}\label{ccninat_diagram}
%
\end{equation}
in $\Catg(\EG,\Ahmal)$.  This diagram commutes if and only if its $(g,s,h)$-component in $\A$, which is the following diagram, commutes for each subset $s \subseteq \ufs{m}$ and objects $g,h \in \EG$.
\begin{equation}\label{inccni_tha_diag}
%
\end{equation}
This diagram commutes by the naturality of $\theta$ with respect to the isomorphism $[g,h] \cn h \fiso g$ in $\EG$.
\item[\cref{ccni_def_v}]  
Since $G$ acts on $\Catg(\EG,-)$ by conjugation \cref{conjugation-gaction}, the $G$-equivariance of the natural transformation $\ccni \cn \inc\cni \to 1$ means that, for each functor $f \cn \EG \to \Ahmal$ and object $g \in \EG$, there is a morphism equality
\[\ccni_{g \cdot f \cdot \ginv} = g \cdot \ccni_f \cdot \ginv \inspace \Catg(\EG,\Ahmal).\]
For each subset $s \subseteq \ufs{m}$ and objects $h,k \in \EG$, the previous equality is verified at the $(h,s,k)$-component by the following morphism equalities in $\A$.
\setlength{\fboxrule}{0.4pt}
\setlength{\fboxsep}{1.3pt}
\[\begin{split}
\ccni_{g \cdot f \cdot \ginv,h,s,k} 
&\seq{1} (g \cdot f \cdot \ginv)_{[h,k],s,k} 
\seq{2} (g \cdot f_{[\ginv h, \ginv k]})_{s,k} \\
&\seq{3} (g \cdot f_{[\ginv h, \ginv k], \ginv s})_k 
\seq{4} g f_{[\ginv h, \ginv k], \ginv s, \ginv k} \\
&\seq{5} g \ccni_{f, \ginv h, \ginv s, \ginv k} 
\seq{6} (g \cdot \ccni_{f, \ginv h, \ginv s})_k \\
&\seq{7} (g \cdot \ccni_{f, \ginv h})_{s,k}
\seq{8} (g \cdot \ccni_f \cdot \ginv)_{h,s,k}
\end{split}\]
These morphism equalities in $\A$ hold for the following reasons.
\begin{itemize}
\item $\fbox{1}$ and $\fbox{5}$ hold by the definition \cref{ccni_f_gsh} of $\ccni$.
\item $\fbox{2}$ and $\fbox{8}$ hold by the fact that $G$ acts on $\EG$ by left multiplication.  $\fbox{2}$ also uses \cref{ginvkh}.
\item $\fbox{3}$ and $\fbox{7}$ hold by the definition \cref{gtha_s} of the $g$-action on a morphism of $\mal$-systems in $\Ah$.
\item $\fbox{4}$ and $\fbox{6}$ hold by the conjugation $g$-action on $\Ah = \Catg(\EG,\A)$.
\end{itemize}
\end{description}
This finishes the proof that $\ccni$ is a $G$-natural isomorphism.
\end{proof}

\section{$G$-Equivalence to $G$-Thickening for Shimakawa $H$-Theory}
\label{sec:sgosg_gequiv}

Under \cref{as:Ahat}, this section proves that the $G$-functors
\begin{equation}\label{ip_aGeq_seci}
\begin{tikzpicture}[vcenter]
\def\h{1.8} \def\t{25}
\draw[0cell]
(0,0) node (a1) {\phantom{A}}
(a1)++(\h,0) node (a2) {\phantom{A}}
(a1)++(-.2,0) node (a1') {\Ahmal}
(a2)++(1,0) node (a2') {\Catg(\EG, \Ahmal)}
(a1)++(\h/2,0) node () {\gsim}
;
\draw[1cell=.9]
(a1) edge[bend left=\t] node {\inc^{\Ahmal}} (a2)
(a2) edge[bend left=\t] node {\cni} (a1)
;
\end{tikzpicture}
\end{equation}
in \cref{def:gcat_inc,def:cni}, along with the $G$-equivariant unit and counit in \cref{def:inccni_unit,def:inccni_counit}, form an adjoint $G$-equivalence.  The strong version, involving $(\inc^{\sgAhmal}, \cnisg, \ucnisg, \ccnisg)$, is also true.  See \cref{thm:SgoAh}.  Thus, for Shimakawa (strong) $H$-theory of the $\Oph$-pseudoalgebra $\Ah$, the inclusion $G$-functor $\inc$ is a $G$-equivalence into the $G$-thickening $\Catg(\EG, \Ahmal)$.  This result improves \cref{inc_eq,ppn_pseudoG}, which assert that the inclusion $G$-functor $\inc$ for a general small $G$-category is a nonequivariant equivalence with a pseudo $G$-equivariant adjoint inverse $\pn$.

\secoutline
\begin{itemize}
\item \cref{def:adjointGeq} defines adjoint $G$-equivalences, with further elaboration given in \cref{expl:adjointGeq}.
\item As an example, \cref{EHEG} shows that, for each subgroup $H \subseteq G$, the inclusion functor from $\EH$ to $\EG$ is part of an adjoint $H$-equivalence.
\item \cref{thm:SgoAh} is the first main result of this chapter, proving that $(\inc^{\Ahmal}, \cni, \ucni, \ccni)$ and its strong variant are adjoint $G$-equivalences.
\item \cref{ex:SgoAh} is an application of \cref{thm:SgoAh} to the Barratt-Eccles $\Gcat$-operad $\BE$ and its pseudoalgebras.  The conclusion is that, for a genuine symmetric monoidal $G$-category of the form $\Ah = \Catg(\EG,\A)$ with $\A$ a $\BE$-pseudoalgebra, the $G$-category of (strong) $\mal$-systems is adjoint $G$-equivalent to its $G$-thickening.
\end{itemize}

\subsection*{Adjoint $G$-Equivalences}

Recall the 2-category $\Gcat$ \pcref{def:GCat} of small $G$-categories, $G$-functors, and $G$-natural transformations.

\begin{definition}\label{def:adjointGeq}
For a group $G$, an \emph{adjoint $G$-equivalence}\index{adjoint G-equivalence@adjoint $G$-equivalence} means an adjoint equivalence in the 2-category $\Gcat$ \pcref{def:internal_adj}.
\end{definition}

\begin{explanation}[Unraveling]\label{expl:adjointGeq}
For $G$-categories $\C$ and $\D$, an adjoint $G$-equivalence from $\C$ to $\D$ is given by an adjoint equivalence of categories 
\begin{equation}\label{aGeq_diag}
\begin{tikzpicture}[vcenter]
\def\h{1.8} \def\t{25}
\draw[0cell]
(0,0) node (a1) {\C}
(a1)++(\h,0) node (a2) {\D}
(a1)++(\h/2,0) node () {\gsim}
;
\draw[1cell=.9]
(a1) edge[bend left=\t] node {\ladt} (a2)
(a2) edge[bend left=\t] node {\radt} (a1)
;
\end{tikzpicture}
\end{equation}
such that 
\begin{itemize}
\item the left adjoint $\ladt$, 
\item the right adjoint $\radt$,
\item the unit $\uadt \cn 1_\C \fiso \radt\ladt$, and
\item the counit $\vadt \cn \ladt\radt \fiso 1_\D$
\end{itemize}
 are $G$-equivariant \pcref{expl:GCat}.  An adjoint $G$-equivalence is denoted by $\gsim$.
\end{explanation}

For a subgroup $H \subseteq G$, restricting the $G$-action to $H$ turns each $G$-category into an $H$-category with the same underlying category.  Thus, the translation category $\EG$ \pcref{def:translation_cat} can be regarded as an $H$-category.  \cref{EHEG} contains the essential content of \cite[Obs.\ 2.10]{merling} and is used in the proof of \cref{catgegc_h}.

\begin{proposition}\label{EHEG}
For each subgroup $H$ of a group $G$, there is an adjoint $H$-equivalence
\begin{equation}\label{ionpon_statement}
\begin{tikzpicture}[vcenter]
\def\h{1.8} \def\t{25} \def\d{.15}
\draw[0cell]
(0,0) node (a1) {\phantom{E}}
(a1)++(\h,0) node (a2) {\phantom{E}}
(a1)++(\h/2,0) node () {\hsim}
(a1)++(-\d,0) node (a1') {\EH}
(a2)++(\d,0) node (a2') {\EG}
;
\draw[1cell=1]
(a1) edge[bend left=\t] node {\ion} (a2)
(a2) edge[bend left=\t] node {\pon} (a1)
;
\end{tikzpicture}
\end{equation}
with the left adjoint $\ion$ induced by the inclusion $H \subseteq G$.
\end{proposition}

\begin{proof}
Since the $H$-action on $\EG$ is the restriction of the $G$-action, the functor $\ion$ is $H$-equivariant.  We construct the right adjoint $\pon$, the unit $\uon \cn 1_{\EH} \fiso \pon\ion$, and the counit $\von \cn \ion\pon \fiso 1_{\EG}$.  Then we prove the triangle identities for $(\ion,\pon,\uon,\von)$, as stated in \cref{def:adjunction}.
\begin{description}
\item[Right adjoint]  
We choose a representative $g_t \in G$ for each right $H$-coset 
\[Hg_t = \big\{hg_t \cn h \in H\big\} \subseteq G,\]
so $G$ is the disjoint union of $\{Hg_t\}_{t \in T}$ for some indexing set $T$.  For the right $H$-coset $H$, we choose the group unit $e \in H$ as the representative.  The right adjoint $\pon \cn \EG \to \EH$ is defined on objects by
\[\pon(hg_t) = h \in \EH \forspace (h,t) \in H \times T,\]
and on morphisms by
\[\pon[kg_r,hg_t] = [k,h] \cn h \fiso k \inspace \EH\]
for $h,k \in H$ and $t,r \in T$.  Since $\EH$ and $\EG$ are translation categories, this defines an $H$-functor $\pon \cn \EG \to \EH$, which is, furthermore, the identity on $\EH$.
\item[Unit]  
The composite $\pon\ion$ is equal to the identity functor $1_{\EH}$.  We define the $H$-equivariant unit $\uon$ as the identity natural transformation
\[1_{\EH} \fto{\uon = 1_1} \pon \ion.\]
\item[Counit]  
For each object $hg_t \in Hg_t$ with $(h,t) \in H \times T$, the counit
\[\ion\pon \fto{\von} 1_{\EG}\]
has $(hg_t)$-component
\[\ion\pon(hg_t) = \ion(h) = h \fto[\iso]{\von_{hg_t}} hg_t\]
given by the unique isomorphism in $\EG$ from $h$ to $hg_t$.  This defines an $H$-natural isomorphism $\von$ because $\EG$ is a translation category.
\item[Triangle identities]  
The unit $\uon$ is the identity on $1_{\EH}$.  The left triangle identity holds because, for each object $h \in \EH$, the $(\ion h)$-component of the counit
\[h \fto{\von_{\ion h} = \von_h} h\]
is the identity morphism on $h$.  The right triangle identity holds because, for each object $hg_t \in Hg_t$ with $(h,t) \in H \times T$, there are morphism equalities
\[\pon \von_{hg_t} = \pon[hg_t,h] = [h,h] = 1_h\]
in $\EH$.\qedhere
\end{description}
\end{proof}

\subsection*{$G$-Thickening of Shimakawa $H$-Theory}

\cref{thm:SgoAh} is the main result of this section.  It improves \cref{inc_eq,ppn_pseudoG} when the $G$-category in question comes from Shimakawa (strong) $H$-theory of the $\Oph$-pseudoalgebra $\Ah = \Catg(\EG,\A)$ for an $\Op$-pseudoalgebra $\A$.  Recall \cref{not:fgsh}.

\begin{theorem}\label{thm:SgoAh}
Under \cref{as:Ahat}, statements \eqref{thm:SgoAh_i} and \eqref{thm:SgoAh_ii} hold for each pointed finite $G$-set $\mal \in \FG$.
\begin{enumerate}
\item\label{thm:SgoAh_i}
There is an adjoint $G$-equivalence
\begin{equation}\label{thm:SgoAh_i_diag}
\begin{tikzpicture}[vcenter]
\def\h{1.8} \def\t{25}
\draw[0cell]
(0,0) node (a1) {\phantom{A}}
(a1)++(\h,0) node (a2) {\phantom{A}}
(a1)++(\h/2,0) node () {\gsim}
(a1)++(-.2,0) node (a1') {\Ahmal}
(a2)++(1,0) node (a2') {\Catg(\EG, \Ahmal)}
;
\draw[1cell=.9]
(a1) edge[bend left=\t] node {\inc^{\Ahmal}} (a2)
(a2) edge[bend left=\t] node {\cni} (a1)
;
\end{tikzpicture}
\end{equation}
given by the following data.
\begin{itemize}
\item $\Ahmal$ is the $G$-category of $\mal$-systems in $\Ah = \Catg(\EG,\A)$ \pcref{def:nsys,def:nsys_morphism,def:nsys_gcat}.
\item $\Catg(\EG,\Ahmal)$ is the $G$-category defined in \cref{def:Catg,def:translation_cat}.
\item The left adjoint is the inclusion $G$-functor $\inc^{\Ahmal}$ (\cref{def:gcat_inc}).
\item The right adjoint is the $G$-functor $\cni$ (\cref{def:cni}).
\item The unit is the identity natural transformation (\cref{def:inccni_unit})
\[1_{\Ahmal} \fto{\ucni = 1_1} \cni\inc^{\Ahmal}.\]
\item The counit is the $G$-natural isomorphism (\cref{def:inccni_counit})
\[\inc^{\Ahmal}\cni \fto[\iso]{\ccni} 1_{\Catg(\EG,\, \Ahmal)}.\]
\end{itemize}
\item\label{thm:SgoAh_ii}
There is an adjoint $G$-equivalence
\begin{equation}\label{thm:SgoAh_ii_diag}
\begin{tikzpicture}[vcenter]
\def\h{1.8} \def\t{25}
\draw[0cell]
(0,0) node (a1) {\phantom{A}}
(a1)++(\h,0) node (a2) {\phantom{A}}
(a1)++(\h/2,0) node () {\gsim}
(a1)++(-.3,0) node (a1') {\sgAhmal}
(a2)++(1.1,0) node (a2') {\Catg(\EG, \sgAhmal)}
;
\draw[1cell=.9]
(a1) edge[bend left=\t] node {\inc^{\sgAhmal}} (a2)
(a2) edge[bend left=\t] node {\cnisg} (a1)
;
\end{tikzpicture}
\end{equation}
involving the strong variants $\sgAhmal$ \cref{Asgordnbe}, $\cnisg$ \cref{cnisg_functor}, $\ucnisg$ \cref{ucnisg}, and $\ccnisg$ \cref{ccnisg}.
\end{enumerate}
\end{theorem}

\begin{proof}
We prove statement \eqref{thm:SgoAh_i}.  The proof for statement \eqref{thm:SgoAh_ii} is the same after a change of notation for the strong variants.  In this proof, $\inc^{\Ahmal}$ is abbreviated to $\inc$.  Since $\inc$, $\cni$, $\ucni$, and $\ccni$ are $G$-equivariant, by \cref{expl:adjointGeq}, it remains to verify that $(\inc,\cni,\ucni,\ccni)$ satisfies the two triangle identities for an adjunction \pcref{def:adjunction}.  
\begin{description}
\item[Left triangle identity]
Since the unit $\ucni$ is the identity natural transformation, the left triangle identity means that, for each $\mal$-system $(a,\gl) \in \Ahmal$, the $\inc(a,\gl)$-component of the counit
\[\inc\cni\inc(a,\gl) \fto{\ccni_{\inc(a,\gl)}} \inc(a,\gl) \inspace \Catg(\EG,\Ahmal)\]
is the identity morphism.  For each subset $s \subseteq \ufs{m}$ and objects $g,h \in \EG$, the $(g,s,h)$-component morphism of $\ccni_{\inc(a,\gl)}$ in $\A$ is given by
\begin{equation}\label{SgoAh_left_tri}
\ccni_{\inc(a,\gl),g,s,h} = \big(\inc(a,\gl)\big)_{[g,h],s,h} = 1_{(a,\gl),s,h} = 1_{a_s(h)}
\end{equation}
\begin{itemize}
\item The first equality holds by the definition \cref{ccni_f_gsh} of $\ccni$.
\item The second equality holds because $\inc(a,\gl)$ is constant at $(a,\gl)$ by \cref{def:gcat_inc}. 
\item The third equality holds because identity morphisms in $\Ahmal$ and $\Catg(-,-)$ are defined componentwise.
\end{itemize}
This proves the left triangle identity for $(\inc,\cni,\ucni,\ccni)$.
\item[Right triangle identity]
Since the unit $\ucni$ is the identity natural transformation, the right triangle identity means that, for each functor $f \cn \EG \to \Ahmal$, the morphism
\[\cni\inc\cni f \fto{\cni\ccni_{f}} \cni f \inspace \Ahmal\]
is the identity.  For each subset $s \subseteq \ufs{m}$ and object $g \in \EG$, the $(s,g)$-component morphism of $\cni\ccni_{f}$ in $\A$ is given by
\begin{equation}\label{SgoAh_right_tri}
(\cni\ccni_f)_{s,g} = \ccni_{f,g,s,g} = f_{[g,g],s,g} = f_{1_g,s,g} = 1_{f_{g,s,g}}
\end{equation}
\begin{itemize}
\item The first equality holds by the definition \cref{cni_theta_sg} of $\cni$ on morphisms.
\item The second equality holds by the definition \cref{ccni_f_gsh} of $\ccni$.
\item The third equality holds by the morphism equality 
\[ [g,g] = 1_g \cn g \fiso g\]
in the translation category $\EG$.
\item The last equality holds by the functoriality of $f$ and the fact that identity morphisms in $\Ahmal$ and $\Catg(-,-)$ are defined componentwise.
\end{itemize}
This proves the right triangle identity for $(\inc,\cni,\ucni,\ccni)$.\qedhere
\end{description}
\end{proof}

\begin{example}[Genuine Symmetric Monoidal $G$-Categories]\label{ex:SgoAh}
\cref{thm:SgoAh} applies to any $\BE$-pseudoalgebra $\A$, where $\BE$ is the Barratt-Eccles $\Gcat$-operad \pcref{def:BE}.  We recall that $\BE$-pseudoalgebras correspond to naive symmetric monoidal $G$-categories under any one of the three 2-equivalences in \cref{thm:BEpseudoalg}.  Applying $\Catg(\EG,-)$ to $\BE$ and the $\BE$-pseudoalgebra $\A$ yields, respectively, 
\begin{itemize}
\item the $G$-Barratt-Eccles operad $\GBE = \Catg(\EG,\BE)$ \pcref{def:GBE} and
\item the genuine symmetric monoidal $G$-category $\Ah = \Catg(\EG,\A)$ \pcref{naive_genuine_smgcat}.
\end{itemize}
In this context, \cref{thm:SgoAh} states that the quadruple $(\inc^{\Ahmal}, \cni, \ucni, \ccni)$ is an adjoint $G$-equivalence between
\begin{itemize}
\item the $G$-category $\Ahmal$ of $\mal$-systems in the genuine symmetric monoidal $G$-category $\Ah = \Catg(\EG,\A)$ and
\item its $G$-thickening $\Catg(\EG,\Ahmal)$.
\end{itemize}
Moreover, the strong variant involving the quadruple $(\inc^{\sgAhmal}, \cnisg, \ucnisg, \ccnisg)$ is also true.
\end{example}

\section{$G$-Thickening to $H$-Theory}
\label{sec:hgo_geq}

This section constructs a $G$-functor 
\[\Catg(\EG, \Ahnbe) \fto{\ci} \Ahnbe\]
to the $H$-theory of the $\Catg(\EG,\Op)$-pseudoalgebra $\Ah = \Catg(\EG,\A)$ at an object $\nbe \in \GG \setminus \{\vstar, \ang{}\}$ \cref{GG_objects}, along with its strong variant $\cisg$.  \cref{thm:HgoAh} proves that $\ci$ is a $G$-equivariant inverse of the inclusion $G$-functor $\inc$ \cref{incC} for the $H$-theory of the $\Catg(\EG,\Op)$-pseudoalgebra $\Ah$.  The strong variant is also true.  

\secoutline
\begin{itemize}
\item \cref{as:OpAh} states the assumptions for this section and \cref{sec:hgo_gequiv}.
\item \cref{def:ci} constructs the $G$-functor $\ci$ and its strong variant $\cisg$ involving $\sgAhnbe$. 
\item \cref{expl:cnici_compare} compares $\ci$ and $\cni$ \pcref{def:cni}, which involves Shimakawa $H$-theory.
\item \cref{ci_welldef} proves that $\ci$ is a well-defined $G$-functor.
\end{itemize}

\cref{as:OpAh} is in effect in \cref{sec:hgo_geq,sec:hgo_gequiv}.

\begin{assumption}\label{as:OpAh}
We consider a $\Tinf$-operad $(\Op,\ga,\opu,\pcom)$ \pcref{as:OpA} for a group $G$, an $\Op$-pseudoalgebra $(\A,\gaA,\phiA)$ \pcref{def:pseudoalgebra}, the $\Tinf$-operad $\Oph = \Catg(\EG,\Op)$, and the $\Oph$-pseudoalgebra \pcref{catgego}
\begin{equation}\label{as:Ah'}
\big(\Ah = \Catg(\EG,\A), \gaAh, \phiAh\big),
\end{equation}
where $\EG$ is the translation category of $G$ \pcref{def:translation_cat}.
\end{assumption}

\begin{notation}\label{not:Ahnbe}
For each object \cref{GG_objects}
\[\nbe = \nbej \in \GG \setminus \{\vstar, \ang{}\}\]
of length $q>0$, recall the small pointed $G$-categories \pcref{def:nsystem,def:nsystem_morphism,def:nbeta_gcat} 
\begin{equation}\label{JgohAhnbe}
\begin{split}
(\Hgoh\Ah)\nbe &= \Ahnbe \andspace\\
(\Hgohsg\Ah)\nbe &= \sgAhnbe
\end{split}
\end{equation}
of (strong) $\nbe$-systems in $\Ah$, where $\Hgoh$ and $\Hgohsg$ denote, respectively, the $H$-theory for $\Oph$ and its strong variant \pcref{Hgo_twofunctor}.  For a functor
\begin{equation}\label{fEGAhnbe}
\EG \fto{f} \Ahnbe,
\end{equation}
elements $g,h \in G$, a marker $\angs = \ang{s_j \subseteq \ufs{n}_j}_{j \in \ufs{q}}$ \cref{marker}, we denote by
\begin{equation}\label{f_g}
(f_g, \glu^{f_g}) = f(g) \in \Ahnbe
\end{equation}
the image of $g$ under $f$; by
\begin{equation}\label{fgangs}
f_{g,\angs} = (f_g)_{\angs} \in \Ah = \Catg(\EG,\A)
\end{equation}
the $\angs$-component object \cref{a_angs} of $f_g$; and by
\begin{equation}\label{fgangsh}
f_{g,\angs,h} = (f_{g,\angs})(h) \in \A
\end{equation}
the image of $h$ under the functor $f_{g,\angs} \cn \EG \to \A$.  We use similar notation for morphisms and the full subcategory $\sgAhnbe \subseteq \Ahnbe$ of strong $\nbe$-systems in $\Ah$.
\end{notation}

\begin{definition}\label{def:ci}
Under \cref{as:OpAh}, th $G$-functor
\begin{equation}\label{ci_functor}
\Catg(\EG, \Ahnbe) \fto{\ci} \Ahnbe
\end{equation}
is defined as follows.  A strong variant is defined in \cref{cisg_functor}.
\begin{description}
\item[Objects] For a functor $f \cn \EG \to \Ahnbe$, the $\nbe$-system in $\Ah$
\begin{equation}\label{ci_f}
(\ci f, \glu^{\ci f}) \in \Ahnbe
\end{equation}
has component objects defined in \cref{ci_fs} and gluing morphisms defined in \cref{cif_glu}.
\begin{description}
\item[Component objects]
For each marker $\angs = \ang{s_j \subseteq \ufs{n}_j}_{j \in \ufs{q}}$, the $\angs$-component object
\begin{equation}\label{ci_fs}
(\ci f)_{\angs} \in \Ah = \Catg(\EG,\A)
\end{equation}
is the functor $\EG \to \A$ with object assignment \cref{ci_fsg} and morphism assignment \cref{ci_fs_mor}.
\begin{description}
\item[Component objects on objects]
$(\ci f)_{\angs}$ sends an object $g \in \EG$ to the object
\begin{equation}\label{ci_fsg}
(\ci f)_{\angs}(g) = f_{g,\angs,g} \in \A.
\end{equation}
\item[Component objects on morphisms]
For an isomorphism $[h,g] \cn g \fiso h$ in $\EG$, the isomorphism
\[(\ci f)_{\angs}(g) = f_{g,\angs,g} \fto[\iso]{(\ci f)_{\angs} [h,g]} (\ci f)_{\angs}(h) = f_{h,\angs,h} \inspace \A\]
is defined by the commutative diagram \cref{ci_fs_mor}.
\begin{equation}\label{ci_fs_mor}
\begin{tikzpicture}[vcenter]
\def\h{3.5} \def\v{1} \def\t{20}
\draw[0cell]
(0,0) node (a11) {\phantom{f_{g,\angs,g}}}
(a11)++(-.15,0) node (a11') {f_{g,\angs,g}}
(a11)++(\h/2,\v) node (a12) {f_{h,\angs,g}}
(a11)++(\h/2,-\v) node (a21) {f_{g,\angs,h}}
(a11)++(\h,0) node (a22) {\phantom{f_{h,\angs,h}}}
(a22)++(.15,0) node (a22') {f_{h,\angs,h}}
;
\draw[1cell=.9]
(a11) edge[bend left=\t] node[pos=.3] {f_{[h,g],\angs,g}} (a12)
(a12) edge[bend left=\t] node[pos=.7] {f_{h,\angs,[h,g]}} (a22)
(a11) edge[bend right=\t] node[swap,pos=.3] {f_{g,\angs,[h,g]}} (a21)
(a21) edge[bend right=\t] node[swap,pos=.7] {f_{[h,g],\angs,h}} (a22)
(a11') edge[transform canvas={yshift=-.2ex}] node {(\ci f)_{\angs} [h,g]} (a22')
;
\end{tikzpicture}
\end{equation}
The boundary diagram in \cref{ci_fs_mor} commutes by the naturality of the isomorphism
\[f_{g,\angs} \fto[\iso]{f_{[h,g],\angs}} f_{h,\angs} \inspace \Ah = \Catg(\EG,\A)\]
with respect to the isomorphism $[h,g]$.
\end{description}
\item[Gluing]
For an object 
\[x \in \Oph(r) = \Catg(\EG,\Op(r))\] 
with $r \geq 0$, a marker $\angs = \ang{s_j \subseteq \ufs{n}_j}_{j \in \ufs{q}}$, an index $p \in \ufs{q}$, and a partition 
\[s_p = \txcoprod_{\ell \in \ufs{r}}\, \spell \subseteq \ufs{n}_p,\]  
the gluing morphism \cref{gluing-morphism} of $\ci f$ at $(x; \angs, p, \ang{\spell}_{\ell})$, where $\ang{\Cdots}_\ell = \ang{\Cdots}_{\ell \in \ufsr}$, is a morphism 
\begin{equation}\label{cif_glu}
\gaAh_r\big(x; \ang{(\ci f)_{\angscomppspell}}_{\ell} \big) 
\fto{\glu^{\ci f}_{x;\, \angs, p, \ang{\spell}_{\ell}}} (\ci f)_{\angs}
\end{equation}
in $\Ah = \Catg(\EG,\A)$, meaning a natural transformation.  For each object $g \in \EG$, the $g$-component morphism of $\glu^{\ci f}_{x;\, \angs, p, \ang{\spell}_{\ell}}$ is defined as the following composite in $\A$.
\begin{equation}\label{cif_glu_g}
\begin{tikzpicture}[vcenter]
\def\h{5.3} \def\v{-1.2}
\draw[0cell=.9]
(0,0) node (a11) {\gaAh_r\big(x; \ang{(\ci f)_{\angscomppspell}}_{\ell} \big) (g)}
(a11)++(\h,0) node (a12) {(\ci f)_{\angs} (g)} 
(a12)++(0,3*\v) node (a22) {f_{g,\angs,g}}
(a11)++(0,\v) node (a21) {\gaA_r\big(x(g) ; \ang{(\ci f)_{\angscomppspell} (g)}_\ell \big)}
(a21)++(0,\v) node (a31) {\gaA_r\big(x(g) ; \ang{f_{g, \angscomppspell, g}}_\ell \big)}
(a31)++(0,\v) node (a32) {\gaAh_r\big(x; \ang{f_{g,\angscomppspell}}_\ell \big) (g)}
;
\draw[1cell=.9]
(a11) edge node {\glu^{\ci f}_{x;\, \angs, p, \ang{\spell}_{\ell}, g}} (a12)
(a32) edge node {\glu^{f_g}_{x;\, \angs, p, \ang{\spell}_\ell, g}} (a22)
(a12) edge[equal] node {\mathbf{d}^\ci} (a22)
(a11) edge[equal] node[swap] {\mathbf{d}} (a21)
(a21) edge[equal] node[swap] {\mathbf{d}^\ci} (a31)
(a31) edge[equal] node[swap] {\mathbf{d}} (a32)
;
\end{tikzpicture}
\end{equation}
The diagram \cref{cif_glu_g} is defined as follows.
\begin{itemize}
\item Each of the two equalities labeled $\mathbf{d}$ holds by \cref{as:Ah'}, since $\Ah = \Catg(\EG,\A)$.
\item Each of the two equalities labeled $\mathbf{d}^\ci$ holds by the definition \cref{ci_fsg} of $(\ci f)_{\angs}(g)$.
\item Using \cref{f_g,fgangs,fgangsh}, the gluing morphism of $(f_g, \glu^{f_g}) \in \Ahnbe$ at $(x; \angs, p, \ang{\spell}_\ell)$ is a morphism 
\[\gaAh_r\big(x; \ang{f_{g,\angscomppspell}}_\ell \big) \fto{\glu^{f_g}_{x;\, \angs, p, \ang{\spell}_\ell}} f_{g,\angs} \inspace \Ah,\]
meaning a natural transformation.  Its $g$-component morphism in $\A$ is the bottom arrow in \cref{cif_glu_g}.
\end{itemize}
This finishes the definition of $\ci$ on objects.  
\end{description}
\item[Morphisms]
Suppose $\theta \cn f \to f'$ is a morphism in $\Catg(\EG,\Ahnbe)$, meaning a natural transformation as follows.
\begin{equation}\label{thaffprime_iicell}
\begin{tikzpicture}[vcenter]
\def\t{28}
\draw[0cell]
(0,0) node (a1) {\phantom{A}}
(a1)++(1.8,0) node (a2) {\phantom{A}}
(a1)++(-.1,0) node (a1') {\EG}
(a2)++(.3,0) node (a2') {\Ahnbe}
;
\draw[1cell=.9]
(a1) edge[bend left=\t] node {f} (a2)
(a1) edge[bend right=\t] node[swap] {f'} (a2)
;
\draw[2cell]
node[between=a1 and a2 at .45, rotate=-90, 2label={above,\theta}] {\Rightarrow}
;
\end{tikzpicture}
\end{equation}
The morphism
\begin{equation}\label{ci_theta}
\ci f \fto{\ci\theta} \ci f' \inspace \Ahnbe
\end{equation}
has, for each marker $\angs = \ang{s_j \subseteq \ufs{n}_j}_{j \in \ufs{q}}$, an $\angs$-component morphism \cref{theta_angs} in $\Ah$, meaning a natural transformation as follows.
\begin{equation}\label{ci_theta_angs}
\begin{tikzpicture}[vcenter]
\def\t{28}
\draw[0cell]
(0,0) node (a1) {\phantom{A}}
(a1)++(2.5,0) node (a2) {\A}
(a1)++(-.1,0) node (a1') {\EG}
;
\draw[1cell=.9]
(a1) edge[bend left=\t] node {(\ci f)_{\angs}} (a2)
(a1) edge[bend right=\t] node[swap] {(\ci f')_{\angs}} (a2)
;
\draw[2cell]
node[between=a1 and a2 at .31, rotate=-90, 2label={above,(\ci\theta)_{\angs}}] {\Rightarrow}
;
\end{tikzpicture}
\end{equation}
Using \cref{fgangsh,ci_fsg}, for each object $g \in \EG$, the $g$-component of $(\ci\theta)_{\angs}$ is defined as the morphism
\begin{equation}\label{ci_theta_angsg}
(\ci f)_{\angs}(g) = f_{g,\angs,g} \fto{(\ci\theta)_{\angs,g} = \theta_{g,\angs,g}} 
(\ci f')_{\angs}(g) = f'_{g,\angs,g}
\end{equation}
in $\A$.  This finishes the definition of $\ci$ on morphisms.

\item[Functoriality]
The assignments
\[f \mapsto (\ci f, \glu^{\ci f}) \andspace \theta \mapsto \ci\theta\]
in \cref{ci_f,ci_theta} define a functor by \cref{ci_theta_angsg} because identity morphisms and composition in $\Catg(\EG,-)$ and $\Ahnbe$ \pcref{def:Catg,def:nsystem_morphism} are defined componentwise.
\end{description}
This finishes the definition of the functor $\ci$ \cref{ci_functor}.  \cref{ci_welldef} proves that $\ci$ is a $G$-functor.

\parhead{Strong variant}.  For the full subcategory $\sgAhnbe \subseteq \Ahnbe$ of strong $\nbe$-systems in $\Ah$, the $G$-functor
\begin{equation}\label{cisg_functor}
\Catg(\EG, \sgAhnbe) \fto{\cisg} \sgAhnbe
\end{equation}
is defined by
\begin{itemize}
\item \cref{ci_f,ci_fs,ci_fsg,ci_fs_mor,cif_glu,cif_glu_g} on objects and
\item \cref{ci_theta,ci_theta_angs,ci_theta_angsg} on morphisms.
\end{itemize}
This is well defined because, for each functor $f \cn \EG \to \sgAhnbe$ and each object $g \in \EG$, $f_g \in \sgAhnbe$ is now a strong $\nbe$-system in $\Ah$.  Each component of its gluing morphism $\glu^{f_g}$ is an isomorphism in $\Ah = \Catg(\EG,\A)$.  Thus, the bottom arrow $\glu^{f_g}_{x;\, \angs, p, \ang{\spell}_\ell, g}$ in \cref{cif_glu_g} is an isomorphism in $\A$.  This shows that $(\cisg f, \glu^{\cisg f})$ is a strong $\nbe$-system in $\Ah$. 
\end{definition}

\begin{explanation}[Comparing $\cni$ and $\ci$]\label{expl:cnici_compare}
Some proofs for $\cni$ in \crefrange{sec:sgo_geq}{sec:sgosg_gequiv} can be reused for $\ci$ by changing the notation according to the following two tables.  \Cref{table.cnici} compares \cref{not:fgsh,not:Ahnbe}.
\begin{figure}[H] 
\centering
{\renewcommand{\arraystretch}{1.3}%
{\setlength{\tabcolsep}{1em}
\begin{tabular}{cr|cr}
\multicolumn{2}{c|}{Shimakawa $H$-theory} & \multicolumn{2}{c}{ $H$-theory} \\ \hline
$\mal \in \FG$ & \eqref{def:FG} & $\nbe \in \GG \setminus \{\vstar,\ang{}\}$ & \eqref{def:GG} \\
$\Ahmal$, $\sgAhmal$ & \eqref{def:sgo} & $\Ahnbe$, $\sgAhnbe$ & \eqref{Hgo_twofunctor} \\
$s \subseteq \ufs{m}$ & \cref{ufsn} & $\angs = \ang{s_j \subseteq \ufs{n}_j}_{j \in \ufs{q}}$ & \cref{marker} \\[.4ex] \hline
\raisebox{6pt}{\phantom{M}}$f \cn \EG \to \Ahmal$\raisebox{6pt}{\phantom{M}} 
& \cref{fEGAhmal} & $f \cn \EG \to \Ahnbe$ & \cref{fEGAhnbe} \\
$(f_g,\gl^{f_g}) \in \Ahmal$ & \cref{f_sub_g} & $(f_g,\glu^{f_g}) \in \Ahnbe$ & \cref{f_g} \\
$f_{g,s} \in \Ah$ & \cref{fgs} & $f_{g,\angs} \in \Ah$ & \cref{fgangs} \\
$f_{g,s,h} \in \A$ & \cref{fgsh} & $f_{g,\angs,h} \in \A$ & \cref{fgangsh} \\
\end{tabular}}}
\caption{A comparison of \cref{not:fgsh,not:Ahnbe}.}
\label{table.cnici}
\end{figure}
\noindent
\Cref{table.pq} compares $\cni$ and $\ci$ \pcref{def:cni,def:ci}.   
\begin{figure}[H] 
\centering
\resizebox{.9\columnwidth}{!}{%
{\renewcommand{\arraystretch}{1.3}%
{\setlength{\tabcolsep}{1ex}
\begin{tabular}{cr|cr}
$\Catg(\EG, \Ahmal) \fto{\cni} \Ahmal$ & \cref{cni_functor} & $\Catg(\EG, \Ahnbe) \fto{\ci} \Ahnbe$ & \cref{ci_functor} \\[.3ex] \hline
\raisebox{6pt}{\phantom{M}}$(\cni f, \gl^{\cni f}) \in \Ahmal$\raisebox{6pt}{\phantom{M}} 
& \cref{cni_f} & $(\ci f, \glu^{\ci f}) \in \Ahnbe$ & \cref{ci_f} \\
$(\cni f)_s(g) = f_{g,s,g}$ & \cref{cni_fsg} & $(\ci f)_{\angs}(g) = f_{g,\angs,g}$ & \cref{ci_fsg} \\
$(\cni f)_s[h,g] = f_{h,s,[h,g]} \circ f_{[h,g],s,g}$ & \cref{cni_fs_mor} & 
$(\ci f)_{\angs}[h,g] = f_{h,\angs,[h,g]} \circ f_{[h,g],\angs,g}$ & \cref{ci_fs_mor} \\
$\phantom{(\cni f)_s[h,g]} = f_{[h,g],s,h} \circ f_{g,s,[h,g]}$ & \cref{cni_fs_mor} & 
$\phantom{(\ci f)_{\angs}[h,g]} = f_{[h,g],\angs,h} \circ f_{g,\angs,[h,g]}$ & \cref{ci_fs_mor} \\[.5ex] \hline
\raisebox{5pt}{\phantom{M}}$s = \txcoprod_{\ell \in \ufs{r}}\, s_\ell \subseteq \ufs{m}$\raisebox{5pt}{\phantom{M}} 
& \cref{partition} 
& $s_p = \txcoprod_{\ell \in \ufs{r}}\, \spell \subseteq \ufs{n}_p$ & \cref{partition} \\
$\gaAh_r\big(x; \ang{(\cni f)_{s_\ell}}_{\ell} \big) 
\fto{\gl^{\cni f}_{x;\, s, \ang{s_\ell}_{\ell}}} (\cni f)_s$ & \cref{cnif_gl} & 
$\gaAh_r\big(x; \ang{(\ci f)_{\angscomppspell}}_{\ell} \big) 
\fto{\glu^{\ci f}_{x;\, \angs, p, \ang{\spell}_{\ell}}} (\ci f)_{\angs}$ & \cref{cif_glu} \\[.3ex] 
$\gl^{\cni f}_{x;\, s, \ang{s_\ell}_{\ell}, g} = \gl^{f_g}_{x;\, s, \ang{s_\ell}_{\ell}, g}$ & \cref{cnif_gl_g} &
$\glu^{\ci f}_{x;\, \angs, p, \ang{\spell}_{\ell}, g} = \glu^{f_g}_{x;\, \angs, p, \ang{\spell}_\ell, g}$ & \cref{cif_glu_g} \\
\end{tabular}}}}
\caption{A comparison of \cref{def:cni,def:ci}.}
\label{table.pq}
\end{figure}
\noindent
\cref{expl:ucni_uci_compare} further compares the units and counits for the pairs of $G$-functors $(\inc^{\Ahmal},\cni)$ and $(\inc^{\Ahnbe},\ci)$.
\end{explanation}

\begin{lemma}\label{ci_welldef}
In \cref{def:ci}, the data
\[\Catg(\EG, \Ahnbe) \fto{\ci} \Ahnbe\]
define a $G$-functor.
\end{lemma}

\begin{proof}
It suffices to prove the following statements \crefrange{ci_welldef_i}{ci_welldef_iii}, where changing the notation means using the correspondence between $\cni$ and $\ci$ discussed in \cref{expl:cnici_compare}.
\begin{enumerate}
\item\label{ci_welldef_i} The pair \cref{ci_f}
\[(\ci f, \glu^{\ci f})\]
is an $\nbe$-system in $\Ah = \Catg(\EG,\A)$.  This statement is proved by reusing the proofs of \cref{cnifs_welldef,glcnif_welldef,cnif_welldef} and changing the notation.
\item\label{ci_welldef_ii} The data \cref{ci_theta}
\[\ci f \fto{\ci\theta} \ci f'\] 
define a morphism in $\Ahnbe$.  This statement is proved by reusing the proof of \cref{cnitha_welldef} and changing the notation.  Given \eqref{ci_welldef_i} and \eqref{ci_welldef_ii}, it is already explained in \cref{def:ci} that $\ci$ \cref{ci_functor} is a functor.
\item\label{ci_welldef_iii} The functor
\[\Catg(\EG, \Ahnbe) \fto{\ci} \Ahnbe\]
is $G$-equivariant.  This statement is proved by reusing the proof of \cref{cni_equivariant} and changing the notation.  Instead of the $G$-action on $\mal$-systems and their morphisms \pcref{def:nsys_gcat}, we use the $G$-action on $\nbe$-systems and their morphisms \pcref{def:nbeta_gcat}.
\end{enumerate}
This proves that $\ci$ is a $G$-functor.
\end{proof}

\section{$G$-Equivalence to $G$-Thickening for $H$-Theory}
\label{sec:hgo_gequiv}

Under \cref{as:OpAh}, this section proves that the $G$-functors
\begin{equation}\label{incci_seci}
\begin{tikzpicture}[vcenter]
\def\h{1.8} \def\t{25}
\draw[0cell]
(0,0) node (a1) {\phantom{A}}
(a1)++(\h,0) node (a2) {\phantom{A}}
(a1)++(-.3,0) node (a1') {\Ahnbe}
(a2)++(1.13,0) node (a2') {\Catg(\EG, \Ahnbe)}
(a1)++(\h/2,0) node () {\gsim}
;
\draw[1cell=.9]
(a1) edge[bend left=\t] node {\inc^{\Ahnbe}} (a2)
(a2) edge[bend left=\t] node {\ci} (a1)
;
\end{tikzpicture}
\end{equation}
in \cref{def:gcat_inc,def:ci} are part of an adjoint $G$-equivalence.  The strong version, involving the $G$-functors $(\inc^{\sgAhnbe}, \cisg)$, is also true.  See \cref{thm:HgoAh}.  Thus, for (strong) $H$-theory of the $\Oph$-pseudoalgebra $\Ah = \Catg(\EG,\A)$, the inclusion $G$-functor $\inc$ is a $G$-equivalence into its $G$-thickening $\Catg(\EG,\Ahnbe)$.  This result is analogous to \cref{thm:SgoAh}, which proves the same statement for Shimakawa $H$-theory.

\secoutline
\begin{itemize}
\item \cref{def:incci_unit} defines the unit 
\[1_{\Ahnbe} \fto{\uci = 1_1} \ci\inc^{\Ahnbe}\]
and its strong variant $\ucisg = 1_1$.
\item \cref{def:incci_counit} defines the counit 
\[\inc^{\Ahnbe}\ci \fto[\iso]{\cci} 1_{\Catg(\EG,\,\Ahnbe)}\]
and its strong variant $\ccisg$.
\item \cref{expl:cci_domain} describes the domain of $\cci$ explicitly.
\item \cref{expl:ucni_uci_compare} compares $(\ucni,\ccni)$ \pcref{def:inccni_unit,def:inccni_counit} with $(\uci,\cci)$.
\item \cref{cci_welldef} proves that $\cci$ and its strong variant are $G$-natural isomorphisms.
\item \cref{thm:HgoAh} proves that $(\inc^{\Ahnbe}, \ci, \uci, \cci)$ and its strong variant are adjoint $G$-equivalences.
\item \cref{ex:HgoAh} is an application of \cref{thm:HgoAh} to the Barratt-Eccles $\Gcat$-operad $\BE$ and its pseudoalgebras.  The conclusion is that, for a genuine symmetric monoidal $G$-category of the form $\Catg(\EG,\A)$ with $\A$ a $\BE$-pseudoalgebra, the $G$-category of (strong) $\nbe$-systems is adjoint $G$-equivalent to its $G$-thickening.
\end{itemize}

\subsection*{Unit and Counit}

\begin{definition}[Unit]\label{def:incci_unit}
In the context of \cref{def:gcat_inc,def:ci}, we define the identity $G$-natural transformations
\begin{equation}\label{uci}

\end{equation}
of, respectively, $1_{\Ahnbe}$ and $1_{\sgAhnbe}$.  Using the correspondence between $\cni$ and $\ci$ \pcref{expl:cnici_compare}, the proof of \cref{inccni_unit_dom} also proves the equalities 
\[\ci \inc^{\Ahnbe} = 1_{\Ahnbe} \andspace \cisg \inc^{\sgAhnbe} = 1_{\sgAhnbe}.\]
Thus, $\uci$ and $\ucisg$ are well defined.
\end{definition}

\begin{definition}[Counit]\label{def:incci_counit}
In the context of \cref{def:gcat_inc,def:ci}, the $G$-natural isomorphism
\begin{equation}\label{cci}
%
\end{equation}
is defined as follows.  A strong variant is defined in \cref{ccisg}.  For a functor $f \cn \EG \to \Ahnbe$, by \cref{def:gcat_inc}, the functor
\[\EG \fto{\inc^{\Ahnbe} \ci f} \Ahnbe\]
is constant at $\ci f \in \Ahnbe$ \cref{ci_f}.  By \cref{ci_fsg}, for a marker $\angs = \ang{s_j \subseteq \ufs{n}_j}_{j \in \ufs{q}}$ and objects $g,h \in \EG$, its $(g,\angs,h)$-component object is given by
\begin{equation}\label{cci_dom_obj}
(\inc^{\Ahnbe} \ci f)_{g,\angs,h} = (\ci f)_{\angs}(h) = f_{h,\angs,h} \in \A.
\end{equation} 
The $f$-component of $\cci$ is the natural isomorphism
\begin{equation}\label{cci_f}
\begin{tikzpicture}[vcenter]
\def\t{32}
\draw[0cell]
(0,0) node (a1) {\phantom{A}}
(a1)++(2,0) node (a2) {\phantom{A}}
(a1)++(-.1,0) node (a1') {\EG}
(a2)++(.3,0) node (a2') {\Ahnbe}
;
\draw[1cell=.9]
(a1) edge[bend left=\t] node {\inc^{\Ahnbe}\ci f} (a2)
(a1) edge[bend right=\t] node[swap] {f} (a2)
;
\draw[2cell]
node[between=a1 and a2 at .4, rotate=-90, 2label={above,\cci_f}] {\Rightarrow}
;
\end{tikzpicture}
\end{equation}
whose $(g,\angs,h)$-component isomorphism is defined as
\begin{equation}\label{cci_f_gsh}
(\inc^{\Ahnbe} \ci f)_{g,\angs,h} = f_{h,\angs,h}
\fto[\iso]{\cci_{f,g,\angs,h} = f_{[g,h],\angs,h}} f_{g,\angs,h}
\end{equation}
in $\A$, with $[g,h] \cn h \fiso g$ the unique isomorphism in the translation category $\EG$.  \cref{cci_welldef} proves that $\cci$ is a well-defined $G$-natural isomorphism.

\parhead{Strong variant}.
For the category $\sgAhnbe$ of strong $\nbe$-systems in $\Ah$ \cref{JgohAhnbe} and the $G$-functor $\cisg$ \cref{cisg_functor}, the $G$-natural isomorphism
\begin{equation}\label{ccisg}
\begin{tikzpicture}[vcenter]
\def\h{4.2} \def\t{25} \def\b{.8}
\draw[0cell]
(0,0) node (a1) {\Catg(\EG,\sgAhnbe)}
(a1)++(\h/2,1.2) node (a2) {\phantom{\sgAhnbe}}
(a2)++(0,-.15) node (a2') {\sgAhnbe}
(a1)++(\h,0) node (a3) {\Catg(\EG,\sgAhnbe)}
(a1)++(\b,0) node (a1') {\phantom{\Ah}}
(a3)++(-\b,0) node (a3') {\phantom{\Ah}}
;
\draw[1cell=.9]
(a1) edge node[swap] {1} (a3)
(a1') edge[bend left=\t] node {\cisg} (a2)
(a2) edge[bend left=\t] node {\inc^{\sgAhnbe}} (a3')
;
\draw[2cell]
node[between=a1 and a3 at .45, shift={(0,.4)}, rotate=-90, 2label={above,\ccisg}] {\Rightarrow}
;
\end{tikzpicture}
\end{equation}
is defined by \cref{cci_f_gsh}, applied to functors $f \cn \EG \to \sgAhnbe$.  This is well defined because $\sgAhnbe$ is a full subcategory of $\Ahnbe$.
\end{definition}

\begin{explanation}[Domain of $\cci$]\label{expl:cci_domain}
We unravel the domain of $\cci$ \cref{cci}, which is the composite 
\[\Catg(\EG,\Ahnbe) \fto{\ci} \Ahnbe \fto{\inc^{\Ahnbe}} \Catg(\EG,\Ahnbe)\]
of the $G$-functors in \cref{def:gcat_inc,def:ci}.  
\begin{description}
\item[Objects]
For a functor $f \cn \EG \to \Ahnbe$, the functor
\[\EG \fto{\inc^{\Ahnbe} \ci f} \Ahnbe\]
is constant at $\ci f \in \Ahnbe$, which is defined in \cref{ci_f,ci_fs,ci_fsg,ci_fs_mor,cif_glu,cif_glu_g}. 
\begin{description}
\item[Component objects]  For a marker $\angs = \ang{s_j \subseteq \ufs{n}_j}_{j \in \ufs{q}}$ and objects $g,h \in \EG$, the $(g,\angs,h)$-component object of $\inc^{\Ahnbe} \ci f$ is $f_{h,\angs,h}$ \cref{cci_dom_obj}.
\item[Component morphisms] For an isomorphism $[k,h] \cn h \fiso k$ in $\EG$, its $(g,\angs,[k,h])$-component isomorphism in $\A$ is given by
\begin{equation}\label{cci_dom_mor}
\begin{split}
& (\inc^{\Ahnbe} \ci f)_{g,\angs,[k,h]} \\
&= (\ci f)_{\angs} [k,h] \\
&= f_{k,\angs,[k,h]} \circ f_{[k,h],\angs,h} \cn f_{h,\angs,h} \to f_{k,\angs,h} \to f_{k,\angs,k}.
\end{split}
\end{equation}
\item[Gluing] Given an object $x \in \Oph(r)$ with $r \geq 0$, a marker $\angs = \ang{s_j \subseteq \ufsn_j}_{j \in \ufsq}$, an index $p \in \ufs{q}$, a partition $s_p = \txcoprod_{\ell \in \ufs{r}}\, \spell \subseteq \ufs{n}_p$, and objects $g,h \in \EG$, the gluing morphism of the $\nbe$-system
\[(\inc^{\Ahnbe} \ci f)_g  = \ci f \in \Ahnbe\] 
at $(x; \angs, p, \ang{\spell}_\ell)$, where $\ang{\Cdots}_\ell = \ang{\Cdots}_{\ell \in \ufsr}$, has the following $h$-component morphism in $\A$. 
\begin{equation}\label{cci_dom_gl}

\end{equation}
to the constant natural transformation
\begin{equation}\label{inccitha}
%
\end{equation}
at $\ci\theta \cn \ci f \to \ci f'$, which is defined in \cref{ci_theta,ci_theta_angs,ci_theta_angsg}.  For each marker $\angs = \ang{s_j \subseteq \ufs{n}_j}_{j \in \ufs{q}}$ and objects $g,h \in \EG$, its $(g,\angs,h)$-component morphism in $\A$ is given as follows.
\begin{equation}\label{cci_dom_nt}
%
\end{equation}
\end{description}
This finishes the unraveling of the domain of the counit $\cci$.

\parhead{Strong variant}. 
Changing the notation from $(\Ahnbe, \ci)$ to $(\sgAhnbe,\cisg)$, the domain of the strong variant $\ccisg$ \cref{ccisg}, which is the composite $G$-functor
\[\Catg(\EG,\sgAhnbe) \fto{\cisg} \sgAhnbe \fto{\inc^{\sgAhnbe}} \Catg(\EG,\sgAhnbe),\]
admits the same description as \crefrange{cci_dom_mor}{cci_dom_nt}.
The gluing morphism $\glu^{f_h}_{x;\, \angs, p, \ang{\spell}_\ell, h}$ in \cref{cci_dom_gl} is an isomorphism because, for $f \in \Catg(\EG,\sgAhnbe)$ and $h \in \EG$, $f_h$ is a strong $\nbe$-system in $\Ah$.
\end{explanation}

\begin{explanation}[Comparing $(\ucni,\ccni)$ and $(\uci,\cci)$]\label{expl:ucni_uci_compare}
\Cref{table.units} compares the units \pcref{def:inccni_unit,def:incci_unit}. 
\begin{figure}[H] 
\centering
{\renewcommand{\arraystretch}{1.3}%
{\setlength{\tabcolsep}{1ex}
\begin{tabular}{cr|cr}
\multicolumn{2}{c|}{Shimakawa $H$-theory} & \multicolumn{2}{c}{$H$-theory} \\ \hline
$1_{\Ahmal} \fto{\ucni = 1_1} \cni\inc^{\Ahmal}$ & \cref{ucni} & 
$1_{\Ahnbe} \fto{\uci = 1_1} \ci\inc^{\Ahnbe}$ & \cref{uci} \\
$1_{\sgAhmal} \fto{\ucnisg = 1_1} \cnisg\inc^{\sgAhmal}$ & \cref{ucnisg} & 
$1_{\sgAhnbe} \fto{\ucisg = 1_1} \cisg\inc^{\sgAhnbe}$ & \cref{ucisg}
\end{tabular}}}
\caption{A comparison of \cref{def:inccni_unit,def:incci_unit}.}
\label{table.units}
\end{figure}
\noindent
\Cref{table.counits} compares the counits \pcref{def:inccni_counit,def:incci_counit}.
\begin{figure}[H] 
\centering
\resizebox{.9\width}{!}{%
{\renewcommand{\arraystretch}{1.3}%
{\setlength{\tabcolsep}{1ex}
\begin{tabular}{cr|cr}
$\inc^{\Ahmal}\cni \fto{\ccni} 1_{\Catg(\EG,\,\Ahmal)}$ & \cref{ccni} & 
$\inc^{\Ahnbe}\ci \fto{\cci} 1_{\Catg(\EG,\,\Ahnbe)}$ & \cref{cci} \\ [1ex]\hline
\raisebox{8pt}{\phantom{M}}$(\inc^{\Ahmal}\cni f)_{g,s,h} = f_{h,s,h}$\raisebox{8pt}{\phantom{M}} 
& \cref{ccni_dom_obj} & 
$(\inc^{\Ahnbe}\ci f)_{g,\angs,h} = f_{h,\angs,h}$ & \cref{cci_dom_obj} \\
$\inc^{\Ahmal}\cni f \fto{\ccni_f} f$ & \cref{ccni_f} & $\inc^{\Ahnbe}\ci f \fto{\cci_f} f$ & \cref{cci_f} \\
$\ccni_{f,g,s,h} = f_{[g,h],s,h}$ & \cref{ccni_f_gsh} & 
$\cci_{f,g,\angs,h} = f_{[g,h],\angs,h}$ & \cref{cci_f_gsh} \\ [1ex]\hline
\raisebox{8pt}{\phantom{M}}$\inc^{\sgAhmal}\cnisg \fto{\ccnisg} 1_{\Catg(\EG,\,\sgAhmal)}$\raisebox{8pt}{\phantom{M}} 
& \cref{ccnisg} & 
$\inc^{\sgAhnbe}\cisg \fto{\ccisg} 1_{\Catg(\EG,\,\sgAhnbe)}$ & \cref{ccisg}
\end{tabular}}}}
\caption{A comparison of \cref{def:inccni_counit,def:incci_counit}.}
\label{table.counits}
\end{figure}
\noindent
\Cref{table.counit_domains} compares the domains of the counits \pcref{expl:ccni_domain,expl:cci_domain}.
\begin{figure}[H] 
\centering
\resizebox{.9\width}{!}{%
{\renewcommand{\arraystretch}{1.8}%
{\setlength{\tabcolsep}{1ex}
\begin{tabular}{cr|cr}
$(\inc^{\Ahmal}\cni f)_{g,s,[k,h]} = f_{k,s,[k,h]} \circ f_{[k,h],s,h}$ & \cref{ccni_dom_mor} & 
$(\inc^{\Ahnbe}\ci f)_{g,\angs,[k,h]} = f_{k,\angs,[k,h]} \circ f_{[k,h],\angs,h}$ & \cref{cci_dom_mor} \\ [.5ex]
$\gl^{(\inc^{\Ahmal} \cni f)_g}_{x;\, s, \ang{s_\ell}_\ell, h} = \gl^{f_h}_{x;\, s, \ang{s_\ell}_\ell, h}$ & \cref{ccni_dom_gl} & 
$\glu^{(\inc^{\Ahnbe} \ci f)_g}_{x;\, \angs, p, \ang{\spell}_\ell, h} = \glu^{f_h}_{x;\, \angs, p, \ang{\spell}_\ell, h}$ & \cref{cci_dom_gl} \\
$(\inc^{\Ahmal}\cni \theta)_{g,s,h} = \theta_{h,s,h}$ & \cref{ccni_dom_nt} & 
$(\inc^{\Ahnbe}\ci \theta)_{g,\angs,h} = \theta_{h,\angs,h}$ & \cref{cci_dom_nt}
\end{tabular}}}}
\caption{A comparison of the domains of the counits.}
\label{table.counit_domains}
\end{figure}
\noindent
Some proofs for $(\ucni,\ccni)$ in \cref{sec:sgo_geq_unit,sec:sgosg_gequiv} can be reused for $(\uci,\cci)$ by changing the notation according to the preceding tables and those in \cref{expl:cnici_compare}.
\end{explanation}

\begin{lemma}\label{cci_welldef}
In \cref{def:incci_counit}, the data 
\begin{equation}\label{cci_statement}
%
\end{equation}
define a $G$-natural isomorphism.  Moreover, the strong variant $\ccisg$ \cref{ccisg} is also a $G$-natural isomorphism.
\end{lemma}

\begin{proof}
This \namecref{cci_welldef} is proved by reusing the proof of \cref{ccni_welldef}, changing the notation as described in \cref{expl:ucni_uci_compare}, and using the following remarks.
\begin{itemize}
\item When we reuse the proof of statement \eqref{ccni_def_ii} in \cref{ccni_welldef}, instead of \cref{def:nsys_morphism} for a morphism of $\mal$-systems, we use \cref{def:nsystem_morphism} for a morphism of $\nbe$-systems. 
\item When we reuse the proof of statement \eqref{ccni_def_v} in \cref{ccni_welldef}, instead of the $g$-action on a morphism of $\mal$-systems \cref{gtha_s}, we use the $g$-action on a morphism of $\nbe$-systems \cref{gtheta_angs}.
\end{itemize}
This finishes the proof.
\end{proof}

\subsection*{Adjoint $G$-Equivalences}

\cref{thm:HgoAh} is the $H$-theory analogue of \cref{thm:SgoAh}.  Recall from \cref{expl:adjointGeq} that an \emph{adjoint $G$-equivalence} is an adjoint equivalence of categories in which the categories, functors, unit, and counit are all $G$-equivariant.

\begin{theorem}\label{thm:HgoAh}
Under \cref{as:OpAh}, statements \eqref{thm:HgoAh_i} and \eqref{thm:HgoAh_ii} hold for each object $\nbe \in \GG \setminus \{\vstar, \ang{}\}$ \cref{GG_objects}.
\begin{enumerate}
\item\label{thm:HgoAh_i}
There is an adjoint $G$-equivalence
\begin{equation}\label{thm:HgoAh_i_diag}
\begin{tikzpicture}[vcenter]
\def\h{1.8} \def\t{25}
\draw[0cell]
(0,0) node (a1) {\phantom{A}}
(a1)++(\h,0) node (a2) {\phantom{A}}
(a1)++(\h/2,0) node () {\gsim}
(a1)++(-.3,0) node (a1') {\Ahnbe}
(a2)++(1.1,0) node (a2') {\Catg(\EG, \Ahnbe)}
;
\draw[1cell=.9]
(a1) edge[bend left=\t] node {\inc^{\Ahnbe}} (a2)
(a2) edge[bend left=\t] node {\ci} (a1)
;
\end{tikzpicture}
\end{equation}
given by the following data.
\begin{itemize}
\item $\Ahnbe$ is the $G$-category of $\nbe$-systems in $\Ah = \Catg(\EG,\A)$ \pcref{def:nsystem,def:nsystem_morphism,def:nbeta_gcat}.
\item $\Catg(\EG,\Ahnbe)$ is the $G$-category defined in \cref{def:Catg,def:translation_cat}.
\item The left adjoint is the inclusion $G$-functor $\inc^{\Ahnbe}$ \pcref{def:gcat_inc}.
\item The right adjoint is the $G$-functor $\ci$ \pcref{def:ci}.
\item The unit is the identity natural transformation \pcref{def:incci_unit}
\[1_{\Ahnbe} \fto{\uci = 1_1} \ci\inc^{\Ahnbe}.\]
\item The counit is the $G$-natural isomorphism \pcref{def:incci_counit}
\[\inc^{\Ahnbe}\ci \fto[\iso]{\cci} 1_{\Catg(\EG,\, \Ahnbe)}.\]
\end{itemize}
\item\label{thm:HgoAh_ii}
There is an adjoint $G$-equivalence
\begin{equation}\label{thm:HgoAh_ii_diag}
\begin{tikzpicture}[vcenter]
\def\h{1.8} \def\t{25}
\draw[0cell]
(0,0) node (a1) {\phantom{A}}
(a1)++(\h,0) node (a2) {\phantom{A}}
(a1)++(\h/2,0) node () {\gsim}
(a1)++(-.4,0) node (a1') {\sgAhnbe}
(a2)++(1.2,0) node (a2') {\Catg(\EG, \sgAhnbe)}
;
\draw[1cell=.9]
(a1) edge[bend left=\t] node {\inc^{\sgAhnbe}} (a2)
(a2) edge[bend left=\t] node {\cisg} (a1)
;
\end{tikzpicture}
\end{equation}
involving the strong variants $\sgAhnbe$ \cref{Asgangordnbe}, $\cisg$ \cref{cisg_functor}, $\ucisg$ \cref{ucisg}, and $\ccisg$ \cref{ccisg}.
\end{enumerate}
\end{theorem}

\begin{proof}
It remains to verify the two triangle identities for an adjunction \pcref{def:adjunction}.  These triangle identities are proved by reusing the proof of \cref{thm:SgoAh} and changing the notation as described in \cref{expl:cnici_compare,expl:ucni_uci_compare}.
\end{proof}

\begin{example}[Genuine Symmetric Monoidal $G$-Categories]\label{ex:HgoAh}
The Barratt-Eccles $\Gcat$-operad $\BE$ \pcref{def:BE} is a $\Tinf$-operad \pcref{as:OpA}.  In this case, \cref{thm:HgoAh} states that, for each $\BE$-pseudoalgebra $\A$, the quadruple $(\inc^{\Ahnbe}, \ci, \uci, \cci)$ is an adjoint $G$-equivalence between
\begin{itemize}
\item the $G$-category $\Ahnbe$ of $\nbe$-systems in the $\GBE$-pseudoalgebra $\Ah = \Catg(\EG,\A)$ (which, by \cref{def:GBE_pseudoalg}, is a genuine symmetric monoidal $G$-category) and
\item its $G$-thickening $\Catg(\EG,\Ahnbe)$.
\end{itemize}
Moreover, the strong variant involving $(\inc^{\sgAhnbe}, \cisg, \ucisg, \ccisg)$ is also true.
\end{example}

\section{Categorical Weak $G$-Equivalences}
\label{sec:cat_weakg}

This section introduces the concept of a categorical weak $G$-equivalence between $G$-categories, which is used in \cref{thm:pistweakgeq}.  It is a categorical analogue of a weak $G$-equivalence between $G$-spaces.   A categorical weak $G$-equivalence is stronger than a nonequivariant equivalence and weaker than an adjoint $G$-equivalence.  Categorical weak $G$-equivalences are automatically topological weak $G$-equivalences.  Throughout this section, $G$ denotes an arbitrary group.

\secoutline
\begin{itemize}
\item \cref{def:cat_weakg} defines categorical weak $G$-equivalences between $G$-categories.
\item \cref{expl:cat_weakg} discusses the relationships between nonequivariant equivalences, categorical weak $G$-equivalences, and adjoint $G$-equivalences, along with some examples.
\item \cref{ex:catweakg_weakg} shows that, for each categorical weak $G$-equivalence $\fun$, $\cla \fun$ is a weak $G$-equivalence between $G$-spaces, where $\cla$ is the classifying space functor.
\item \cref{def:weakG} recalls the concept of a topological weak $G$-equivalence between $G$-categories from \cite{merling}.
\item \cref{expl:weakG} observes that each categorical weak $G$-equivalence is also a topological weak $G$-equivalence.
\item \cref{merling_2.16} recalls a result of Merling \cite{merling} that asserts that, for a $G$-functor $\fun$ that is also a nonequivariant equivalence of categories, the induced $G$-functor $\Catg(\EG,\fun)$ is a categorical weak $G$-equivalence. 
\end{itemize}

\cref{def:cat_weakg} defines the categorical analogue of a weak $G$-equivalence between $G$-spaces \pcref{def:weakG_top}.

\begin{definition}\label{def:cat_weakg}
A \emph{categorical weak $G$-equivalence}\index{categorical weak G-equivalence@categorical weak $G$-equivalence}\index{weak G-equivalence@weak $G$-equivalence!categorical} is a $G$-functor $\fun \cn \C \to \D$ between $G$-categories such that, for each subgroup $H \subseteq G$, the $H$-fixed subfunctor
\[\C^H \fto[\sim]{\fun^H} \D^H\]
is an equivalence of categories.  Here, $\C$ and $\D$ are regarded as $H$-categories by restricting their $G$-actions, and $\fun$ is regarded as an $H$-functor.  A categorical weak $G$-equivalence is also denoted by \label{not:cateqg}$\eqg$.
\end{definition}

\begin{explanation}[Three Notions of Equivalences for $G$-Functors]\label{expl:cat_weakg}
A functor is an equivalence of categories if and only if it is essentially surjective on objects and fully faithful on morphisms \pcref{thm:cat_equiv}.  Thus, a $G$-functor $\fun \cn \C \to \D$ between $G$-categories is a categorical weak $G$-equivalence if and only if, for each subgroup $H \subseteq G$, the $H$-fixed subfunctor $\fun^H \cn \C^H \to \D^H$ is essentially surjective on objects and fully faithful on morphisms.  Since $\fun^H$ is a restriction of $\fun$, faithfulness ($=$ injectivity on morphism sets) only needs to be checked for $\fun^{\{e\}} = \fun$, where $\{e\} \subseteq G$ denotes the subgroup  consisting of only the unit $e \in G$.

There are three progressively stronger notions of equivalences for a $G$-functor:
\begin{enumerate}
\item\label{expl:cat_weakg_i} a nonequivariant equivalence of categories;
\item\label{expl:cat_weakg_ii} a categorical weak $G$-equivalence \pcref{def:cat_weakg}; and
\item\label{expl:cat_weakg_iii} an adjoint $G$-equivalence \pcref{def:adjointGeq}.
\end{enumerate}
If $G$ is the trivial group, then these three notions are equivalent to each other.  For a general group $G$, \eqref{expl:cat_weakg_iii} implies \eqref{expl:cat_weakg_ii}, which, in turn, implies \eqref{expl:cat_weakg_i}.  We provide more details and examples next.
\begin{description}
\item[\cref{expl:cat_weakg_i,expl:cat_weakg_ii}]  
A categorical weak $G$-equivalence $\fun \cn \C \to \D$ \pcref{def:cat_weakg} is, in particular, an equivalence of categories, since $\fun^{\{e\}} = \fun$.  On the other hand, a $G$-functor $\fun'$ that is an equivalence of categories is not generally a categorical weak $G$-equivalence, since its $H$-fixed subfunctor $(\fun')^H$ may not be essentially surjective on objects or full on morphisms.  Here are two examples.
\begin{itemize}
\item For a small $G$-category $\C$, the inclusion $G$-functor \cref{incC}
\[\C \fto{\inc} \Catg(\EG,\C)\]
is a nonequivariant equivalence \pcref{inc_eq}, but it is not generally a categorical weak $G$-equivalence.
\item The strong $H$-theory comparison \cref{PistsgAnbe}
\[\Asgsmaangordnbe \fto{\Pistsg_{\A,\nbe}} \Asgangordnbe\]
is a $G$-functor \cref{PistA_Geq} and a nonequivariant equivalence \pcref{thm:PistAequivalence}.  However, we do not expect it to be a categorical weak $G$-equivalence in general.
\end{itemize}
\item[\cref{expl:cat_weakg_ii,expl:cat_weakg_iii}]  
An adjoint $G$-equivalence  $(\ladt,\radt,\uadt,\vadt)$ \pcref{expl:adjointGeq} is also an adjoint $H$-equivalence for each subgroup $H \subseteq G$.  Passing to $H$-fixed subcategories yields an adjoint equivalence.  Thus, $\ladt$ and $\radt$ are categorical weak $G$-equivalences.  Here are some examples.
\begin{itemize}
\item For each subgroup $H \subseteq G$, in the adjoint $H$-equivalence
\begin{equation}\label{ionpon_aHeq}

\end{equation}
in \cref{EHEG}, both $\ion$ and $\pon$ are categorical weak $H$-equivalences.
\item In the adjoint $G$-equivalences
\begin{equation}\label{inccni_adjG}
%
\end{equation}
in \cref{thm:SgoAh}, the $G$-functors $\inc^{\Ahmal}$, $\cni$, $\inc^{\sgAhmal}$, and $\cnisg$ are categorical weak $G$-equivalences.
\item In the adjoint $G$-equivalences
\begin{equation}\label{ci_diag}
%
\end{equation}
in \cref{thm:HgoAh}, the $G$-functors $\inc^{\Ahnbe}$, $\ci$, $\inc^{\sgAhnbe}$, and $\cisg$ are categorical weak $G$-equivalences.
\end{itemize}
On the other hand, a general categorical weak $G$-equivalence cannot be expected to admit $G$-equivariant inverse, unit, and counit.\defmark
\end{description}
\end{explanation}

\begin{example}\label{ex:catweakg_weakg}
Suppose $\fun \cn \C \to \D$ is a categorical weak $G$-equivalence \pcref{def:cat_weakg} with $\C$ and $\D$ small.  Thus, for each subgroup $H \subseteq G$, the $H$-fixed subfunctor $\fun^H \cn \C^H \fsim \D^H$ is an equivalence of categories.  The classifying space functor \cref{classifying_space}
\[\Cat \fto{\cla} \Top\]
commutes with $G$-action and taking $H$-fixed points and sends adjunctions to homotopy equivalences.  Thus, the morphism
\[(\cla\C)^H \iso \cla(\C^H) \fto[\sim]{(\cla \fun)^H \iso \cla(\fun^H)} (\cla\D)^H \iso \cla(\D^H)\]
is a weak homotopy equivalence, and the $G$-morphism 
\[\cla\C \fto[\eqg]{\cla \fun} \cla\D\] 
is a weak $G$-equivalence between $G$-spaces \pcref{def:weakG_top}.
\end{example}

\subsection*{Topological Weak $G$-Equivalences}

The following notion is due to Merling \cite[Def.\ 2.1]{merling}.

\begin{definition}\label{def:weakG}
A $G$-functor $\fun \cn \C \to \D$ between small $G$-categories is called a \emph{topological weak $G$-equivalence}\index{topological weak G-equivalence@topological weak $G$-equivalence}\index{weak G-equivalence@weak $G$-equivalence!topological} if 
\[\cla\C \fto{\cla \fun} \cla\D\]
is a weak $G$-equivalence between $G$-spaces \pcref{def:weakG_top}.  
\end{definition}

\begin{explanation}\label{expl:weakG}
A $G$-functor $\fun \cn \C \to \D$ is a topological weak $G$-equivalence if, for each subgroup $H \subseteq G$, the $H$-fixed morphism
\[(\cla\C)^H \fto[\sim]{(\cla \fun)^H} (\cla\D)^H\]
is a weak homotopy equivalence of spaces.  By \cref{ex:catweakg_weakg}, each categorical weak $G$-equivalence \pcref{def:cat_weakg} is also a topological weak $G$-equivalence \pcref{def:weakG}.
\end{explanation}

A categorical weak $G$-equivalence is, in particular, a $G$-functor that is also a nonequivariant equivalence, but the converse is not generally true \pcref{expl:cat_weakg}.  The following partial converse is due to Merling \cite[2.10 and 2.16]{merling}.  It is used in \cref{thm:pistweakgeq,thm:h_special}.  We provide a proof here for the reader's convenience.  Recall that $\Catg(\C,\D)$ is the $G$-category of functors $\C \to \D$ and natural transformations with the conjugation $G$-action \pcref{def:Catg}.

\begin{proposition}[Merling]\label{merling_2.16}
Suppose $\fun \cn \C \to \D$ is a $G$-functor between small $G$-categories that is also a nonequivariant equivalence of categories.  Then the $G$-functor given by postcomposing with $\fun$,
\[\Catg(\EG,\C) \fto{\Catg(\EG,\fun)} \Catg(\EG,\D),\]
is a categorical weak $G$-equivalence, hence also a topological weak $G$-equivalence.
\end{proposition}

\begin{proof}
The $G$-functor $\Catg(\EG,\fun)$ is a categorical weak $G$-equivalence if, for each subgroup $H \subseteq G$, the $H$-fixed subfunctor
\begin{equation}\label{catgegfh}
\Catg(\EG,\C)^H \fto{\Catg(\EG,\fun)^H} \Catg(\EG,\D)^H
\end{equation}
is an equivalence of categories.  The rest of this proof consists of three main steps.
\begin{description}
\item[Reduction]
The $H$-fixed subcategory $\Catg(\EG,\C)^H$ can be rewritten as follows.
\begin{equation}\label{catgegc_h}
\begin{split}
\Catg(\EG,\C)^H &= \Cath(\EG,\C)^H \\
&= \Hcat(\EG,\C) \\
&\simeq \Hcat(\EH,\C)
\end{split}
\end{equation}
\begin{itemize}
\item The first equality follows from the fact that the $H$-fixed subcategory of $\Catg(\EG,\C)$ involves only the $H$-action.
\item The second equality holds by \cref{catg_gfixed}. 
\item The last equivalence of categories is induced by the adjoint $H$-equivalence $\ion \cn \EH \to \EG$ in \cref{EHEG} by applying $\Hcat(-,\C)$.
\end{itemize}
Using \cref{catgegc_h}, to show that the functor $\Catg(\EG,\fun)^H$ in \cref{catgegfh} is an equivalence of categories, it suffices to show that, for each $H$-functor $\hun \cn \C \to \D$ that is also an equivalence, the induced functor
\[\Hcat(\EH,\C) \fto{\Hcat(\EH,\hun)} \Hcat(\EH,\D)\]
is an equivalence of categories.  Since $\hun$ is faithful on morphisms, so is $\Hcat(\EH,\hun)$.  It remains to show that $\Hcat(\EH,\hun)$ is (i) essentially surjective on objects and (ii) full on morphisms.
\item[Essential surjectivity]  
Given an $H$-functor $\urho \cn \EH \to \D$, we construct an essential $\Hcat(\EH,\hun)$-preimage as follows.  
\begin{description}
\item[Objects]
Denoting the unit element by $e \in H$, the object assignment of $\urho$ is determined by  
\begin{equation}\label{urho_eh}
\begin{split}
\urho e &= d \in \D \andspace\\
\urho\ell &= \ell(\urho e) = \ell d \forspace \ell \in \EH.
\end{split}
\end{equation}
By the essential surjectivity of $\hun$, there exist an object $c \in \C$ and an isomorphism 
\begin{equation}\label{durhoe_hc}
d = \urho e \fto[\iso]{f} \hun c \inspace \D.
\end{equation}  
We define the $H$-functor $\urho' \cn \EH \to \C$ with the object assignment
\begin{equation}\label{urhoph_hc}
\urho' \ell = \ell c \forspace \ell \in \EH.
\end{equation}
This object assignment is $H$-equivariant by definition.
\item[Morphisms]
For an isomorphism $[\ell,e] \cn e \fiso \ell$ in $\EH$, by the fully faithfulness of $\hun$ and the $H$-equivariance of $\hun$ and $\urho$, there is a unique isomorphism 
\[c = \urho' e \fto[\iso]{c_\ell} \ell c = \urho' \ell \inspace \C\]
such that the following diagram in $\D$ commutes.
\begin{equation}\label{Fch}
\begin{tikzpicture}[vcenter]
\def\v{-1.3}
\draw[0cell]
(0,0) node (a11) {\urho e}
(a11)++(2.7,0) node (a12) {\urho\ell}
(a11)++(0,\v) node (a21) {\hun c}
(a12)++(0,\v) node (a22) {\hun(\ell c)}
;
\draw[1cell=.9]
(a11) edge node {\urho[\ell,e]} node[swap] {\iso} (a12)
(a12) edge node {\ell f} node[swap] {\iso} (a22)
(a11) edge node[swap] {f} node {\iso} (a21)
(a21) edge node {\hun c_\ell} (a22)
;
\end{tikzpicture}
\end{equation}
The $H$-equivariant morphism assignment of $\urho' \cn \EH \to \C$ is defined as
\begin{equation}\label{urhop_mor}
\begin{split}
\urho'[\ell,e] &= c_\ell \cn \urho' e \fiso \urho' \ell \andspace\\
\urho'[k,\ell] &= \ell \urho'[\ellinv k,e] = \ell c_{\ellinv k} \forspace k,\ell \in \EH.
\end{split}
\end{equation}
\item[$H$-functoriality]
The morphism assignment of $\urho' \cn \EH \to \C$ preserves identity morphisms because 
\[ [e,e] = 1_e \andspace c_e = 1_c.\]
Since the morphism assignment of $\urho'$ is $H$-equivariant, it suffices to show that $\urho'$ preserves composition for morphisms
\[e \to \ell \to \ell k\]
for $k,\ell \in \EH$.  The desired morphism equality
\[\urho'[\ell k,e] = \big(\ell \urho'[k,e]\big) \circ \urho'[\ell,e] \inspace \C\]
follows from the same condition for $\urho$ using \cref{Fch,urhop_mor}, the functoriality and faithfulness of $\hun$, and the $H$-equivariance of $\hun$ and $\urho$.  Thus, $\urho' \cn \EH \to \C$ is an $H$-functor.
\item[Essential preimage]
Using \crefrange{urho_eh}{urhoph_hc}, the $H$-natural isomorphism $\utha \cn \urho \fiso \hun\urho'$ is defined by the components
\[\urho\ell \fto[\iso]{\utha_\ell = \ell f} \hun(\ell c) = \hun\urho'\ell \forspace \ell \in \EH.\]
The naturality diagram of $\utha$ for an isomorphism $[\ell,e] \cn e \fiso \ell$ in $\EH$ is the commutative diagram \cref{Fch}.  By \cref{urhop_mor} and $H$-equivariance, this implies the naturality of $\utha$ for a general isomorphism $[k,\ell] \cn \ell \fiso k$ in $\EH$.  
\end{description}
This proves that $\Hcat(\EH,\hun)$ is essentially surjective on objects.
\item[Fullness]  
To show that the functor $\Hcat(\EH,\hun)$ is full on morphisms, suppose $\utha \cn \hun\urho_1 \to \hun\urho_2$ is an $H$-natural transformation for $H$-functors $\urho_1, \urho_2 \cn \EH \to \C$.   For each object $\ell \in \EH$, the fully faithfulness of $\hun$ implies that there exists a unique morphism
\[\urho_1\ell \fto{\utha'_\ell} \urho_2\ell \inspace \C\]
such that 
\begin{equation}\label{hthaph}
\hun\urho_1\ell \fto{\hun\utha'_\ell = \utha_\ell} \hun\urho_2\ell \inspace \D.
\end{equation}
We define an $H$-natural transformation $\utha' \cn \urho_1 \to \urho_2$ with $\ell$-component given by $\utha'_\ell$ for $\ell \in \EH$.  
\begin{itemize}
\item The $H$-equivariance of $\utha'$ follows from \cref{hthaph}, the faithfulness of $\hun$, and the $H$-equivariance of $\hun$ and $\utha$.
\item The naturality of $\utha'$ follows from \cref{hthaph}, the functoriality and faithfulness of $\hun$, and the naturality of $\utha$.
\end{itemize}
The equality
\[1_\hun * \utha' = \utha\] 
proves that $\Hcat(\EH,\hun)$ is full on morphisms.\qedhere
\end{description}
\end{proof}

\section{Comparison Weak $G$-Equivalence}
\label{sec:levelgeq}

For a $\Uinf$-operad \pcref{as:OpA'} and an $\Op$-pseudoalgebra $\A$, the strong $H$-theory comparison \cref{PistsgAnbe} 
\[(\smast\Sgosg\A)\nbe = \Asgsmaangordnbe
\fto{\Pistsg_{\A,\nbe}} 
\Asgangordnbe = (\Hgosg\A)\nbe\]
is a $G$-functor \cref{PistA_Geq} and a nonequivariant equivalence of categories \pcref{thm:PistAequivalence}.  The main result of this section, \cref{thm:pistweakgeq}, proves that, for the $\Catg(\EG,\Op)$-pseudoalgebra $\Ah = \Catg(\EG,\A)$, the strong $H$-theory comparison $\Pistsg_{\Ah}$ is componentwise a categorical weak $G$-equivalence \pcref{def:cat_weakg} and, therefore, also a topological weak $G$-equivalence \pcref{def:weakG}.  Thus, Shimakawa strong $H$-theory and strong $H$-theory for $\Ah$ are componentwise categorically weakly $G$-equivalent.

\secoutline
\begin{itemize}
\item \cref{thm:pistweakgeq} proves that each component of $\Pistsg_{\Ah}$ is a categorical weak $G$-equivalence.
\item \cref{ex:pistsg_weakg} applies \cref{thm:pistweakgeq} to obtain weak $G$-equivalences between $G$-spaces.
\item \cref{ex:pistweq} is an application to the Barratt-Eccles $\Gcat$-operad $\BE$ and a $\BE$-pseudoalgebra $\A$.  In this case, \cref{thm:pistweakgeq} states that the strong $H$-theory comparison for the genuine symmetric monoidal $G$-category $\Catg(\EG,\A)$ is componentwise a categorical weak $G$-equivalence.
\item \cref{expl:pistweq_necessity} discusses why \cref{thm:pistweakgeq} uses the \emph{strong} $H$-theory comparison.
\end{itemize}

For \cref{thm:pistweakgeq}, recall that a $\Uinf$-operad $\Op$ \pcref{as:OpA'} is a 1-connected $\Gcat$-operad that is levelwise a nonempty translation category.  It is, in particular, a $\Tinf$-operad \pcref{expl:OpA'}.  The $\Gcat$-operad $\Oph = \Catg(\EG,\Op)$ is also a $\Uinf$-operad (\cref{ex:OpA'} \eqref{ex:OpA'_iii}).  Applying $\Catg(\EG,-)$ to an $\Op$-pseudoalgebra $\A$ \pcref{def:pseudoalgebra} yields an $\Oph$-pseudoalgebra $\Ah = \Catg(\EG,\A)$ \pcref{catgego}.

\begin{theorem}\label{thm:pistweakgeq}\index{H-theory comparison@$H$-theory comparison!weak G-equivalence@weak $G$-equivalence}\index{weak G-equivalence@weak $G$-equivalence!H-theory comparison@$H$-theory comparison}\index{Shimakawa K-theory@Shimakawa $K$-theory!comparison with multifunctorial $K$-theory}\index{multifunctorial K-theory@multifunctorial $K$-theory!comparison with Shimakawa}
For a $\Uinf$-operad $\Op$, an $\Op$-pseudoalgebra $\A$, and the $\Oph$-pseudoalgebra $\Ah$, the strong $H$-theory comparison $G$-natural transformation \cref{PistsgA}
\begin{equation}\label{PistsgAh_thm}
\begin{tikzpicture} [vcenter]
\def\s{25}
\draw[0cell]
(0,0) node (a1) {\GG}
(a1)++(2.5,0) node (a2) {\phantom{\GG}}
(a2)++(.15,0) node (a2') {\Catgst}
;
\draw[1cell=.9]
(a1) edge[bend left=\s] node {\smast \Sgohsg\Ah} (a2)
(a1) edge[bend right=\s] node[swap] {\Hgohsg\Ah} (a2) 
;
\draw[2cell]
node[between=a1 and a2 at .4, rotate=-90, 2label={above,\Pistsg_{\Ah}}] {\Rightarrow}
;
\end{tikzpicture}
\end{equation}
is componentwise a categorical weak $G$-equivalence, hence also a topological weak $G$-equivalence.
\end{theorem}

\begin{proof}
For an object $\nbe \in \GG$ \cref{GG_objects}, we abbreviate the $\nbe$-component pointed $G$-functor \cref{PistsgAnbe} 
\[(\smast \Sgohsg \Ah)\nbe = \Ahsgsmaangordnbe 
\fto{\Pistsg_{\Ah,\nbe}} \Ahsgangordnbe = (\Hgohsg\Ah)\nbe\]
to $\Pistsg$.  To verify that $\Pistsg$ is a categorical weak $G$-equivalence \pcref{def:cat_weakg}, we consider the following three cases.
\begin{itemize}
\item If $\nbe$ is the basepoint $\vstar$, then $\Pistsg = 1_{\boldone}$ \cref{PistA_vstar}. 
\item If $\nbe$ is the empty tuple $\ang{}$, then $\Pistsg$ is a pointed $G$-isomorphism \cref{PistA_empty}, hence also a categorical weak $G$-equivalence. 
\end{itemize} 
In the rest of this proof, we assume that $\nbe \in \GG \setminus \{\vstar,\ang{}\}$.
\begin{description}
\item[Commutative diagram]
We use the diagram \cref{pistsgahnbe_inc} of $G$-functors and a 2-out-of-3 argument to show that $\Pistsg$ is a categorical weak $G$-equivalence.
\begin{equation}\label{pistsgahnbe_inc}
\begin{tikzpicture}[vcenter]
\def\v{-1.4}
\draw[0cell=.85]
(0,0) node (a11) {\sgAhsmanbe}
(a11)++(5,0) node (a12) {\sgAhnbe}
(a11)++(0,\v) node (a21) {\Catg(\EG,\sgAhsmanbe)}
(a12)++(0,\v) node (a22) {\Catg(\EG,\sgAhnbe)}
;
\draw[1cell=.8]
(a11) edge node {\Pistsg} (a12)
(a12) edge node {\inc^{\sgAhnbe}} node[swap] {\gsim} (a22)
(a11) edge node {\gsim} node[swap] {\inc^{\sgAhsmanbe}} (a21)
(a21) edge node {\Catg(\EG,\Pistsg)} node[swap] {\eqg} (a22)
;
\end{tikzpicture}
\end{equation}
The four $G$-functors in the diagram \cref{pistsgahnbe_inc} are given as follows.
\begin{description}
\item[Top] $\Pistsg$ is the strong $H$-theory comparison \cref{PistsgAnbe} for the $\Oph$-pseudoalgebra $\Ah = \Catg(\EG,\A)$.  It is a $G$-functor by \cref{PistA_Geq} applied to $\Ah$.
\item[Left] $\inc^{\sgAhsmanbe}$ is the inclusion $G$-functor \cref{incC} for the small $G$-category $\sgAhsmanbe$ \cref{Asgordnbe} of strong $(\sma\nbe)$-systems in $\Ah$, where $\sma\nbe \in \FG$ is the pointed finite $G$-set defined in \cref{smash_GGobjects}.
\item[Right] $\inc^{\sgAhnbe}$ is the inclusion $G$-functor \cref{incC} for the small $G$-category $\sgAhnbe$ \cref{Asgangordnbe} of strong $\nbe$-systems in $\Ah$.
\item[Bottom] $\Catg(\EG,\Pistsg)$ is obtained from the top $G$-functor by applying $\Catg(\EG,-)$.
\item[Commutativity] 
For each object or morphism $b \in \sgAhsmanbe$, under either composite in \cref{pistsgahnbe_inc}, the image of $b$ in $\Catg(\EG,\sgAhnbe)$ is constant at $\Pistsg b$.  
\end{description}
Thus, the diagram \cref{pistsgahnbe_inc} of $G$-functors commutes.
\item[Weak $G$-equivalences]
The left, right, and bottom $G$-functors in the diagram \cref{pistsgahnbe_inc} are categorical weak $G$-equivalences for the following reasons.
\begin{description}
\item[Left] By \cref{thm:SgoAh} \eqref{thm:SgoAh_ii} applied to the object $\sma\nbe \in \FG$, the $G$-functor $\inc^{\sgAhsmanbe}$ is the left adjoint of an adjoint $G$-equivalence, hence also a categorical weak $G$-equivalence \pcref{expl:cat_weakg}.
\item[Right] By \cref{thm:HgoAh} \eqref{thm:HgoAh_ii}, the $G$-functor $\inc^{\sgAhnbe}$ is the left adjoint of an adjoint $G$-equivalence, hence also a categorical weak $G$-equivalence.
\item[Bottom] By \cref{thm:PistAequivalence} applied to the $\Uinf$-operad $\Oph = \Catg(\EG,\Op)$ and the $\Oph$-pseudoalgebra $\Ah$, $\Pistsg$ is a $G$-functor \cref{PistA_Geq} that is also a nonequivariant equivalence of categories.  \cref{merling_2.16} implies that $\Catg(\EG,\Pistsg)$ is a categorical weak $G$-equivalence.
\end{description}
Since equivalences of categories satisfy the 2-out-of-3 property, so do categorical weak $G$-equivalences.  Thus, the commutative diagram \cref{pistsgahnbe_inc} implies that $\Pistsg$ is a categorical weak $G$-equivalence.\qedhere
\end{description}
\end{proof}

\begin{example}[Weak $G$-Equivalences]\label{ex:pistsg_weakg}
\cref{thm:pistweakgeq} shows that, for each object $\nbe \in \GG$ \cref{GG_objects}, the strong $H$-theory comparison pointed $G$-functor $\Pistsg_{\Ah,\nbe}$ \cref{PistsgAnbe} is a categorical weak $G$-equivalence \pcref{def:cat_weakg} from Shimakawa strong $H$-theory \cref{Asgordnbe} at $\sma\nbe \in \FG$ \cref{smash_GGobjects} to strong $H$-theory at $\nbe$ \cref{Asgangordnbe} for the $\Oph$-pseudoalgebra $\Ah = \Catg(\EG,\A)$.  By  \cref{ex:catweakg_weakg}, applying the classifying space functor $\cla$ to $\Pistsg_{\Ah,\nbe}$ yields a weak $G$-equivalence 
\[\cla\Ahsgsmaangordnbe
\fto[\eqg]{\cla\Pistsg_{\Ah,\nbe}} \cla\Ahsgangordnbe\]
between $G$-spaces.
\end{example}

\begin{example}[Genuine Symmetric Monoidal $G$-Categories]\label{ex:pistweq}
The Barratt-Eccles $\Gcat$-operad $\BE$ \pcref{def:BE} is a $\Uinf$-operad \pcref{ex:OpA'}.  Each $\BE$-pseudoalgebra $\A$ yields a $\GBE$-pseudoalgebra $\Ah = \Catg(\EG,\A)$ \pcref{naive_genuine_smgcat}, where $\GBE = \Catg(\EG,\BE)$ is the $G$-Barratt-Eccles operad \pcref{def:GBE}.  In this case, \cref{thm:pistweakgeq} asserts that, for each object $\nbe \in \GG$, the strong $H$-theory comparison \cref{PistsgAnbe} 
\[(\smast \Sgosggbe \Ah)\nbe = \Ahsgsmaangordnbe 
\fto[\eqg]{\Pistsg_{\Ah,\nbe}} \Ahsgangordnbe = (\Hgosggbe\Ah)\nbe\]
is a categorical weak $G$-equivalence, hence also a topological weak $G$-equivalence.  The domain $\Ahsgsmaangordnbe$ is Shimakawa strong $H$-theory \cref{Asgordnbe} for the genuine symmetric monoidal $G$-category $\Ah$ \pcref{def:GBE_pseudoalg} at $\sma\nbe \in \FG$.  The codomain $\Ahsgangordnbe$ is the strong $H$-theory of $\Ah$ at $\nbe$ \cref{Asgangordnbe}.
\end{example}

\begin{explanation}[Necessity of Strong Variant]\label{expl:pistweq_necessity}
\cref{thm:pistweakgeq} is about the \emph{strong} $H$-theory comparison $\Pistsg$ for the $\Oph$-pseudoalgebra $\Ah = \Catg(\EG,\A)$.  To see why we use the strong variant $\Pistsg$ instead of the $H$-theory comparison
\[\smast\Sgoh \fto{\Pist} \Hgoh\]
defined in \cref{Pistar_twonat}, consider the diagram \cref{pistsgahnbe_inc} in the proof of \cref{thm:pistweakgeq}.
\begin{itemize}
\item Using the category $\Ahsmanbe$ of all $(\sma\nbe)$-systems in $\Ah$ \cref{A_nbeta}, the inclusion $G$-functor
\[\Ahsmanbe \fto[\gsim]{\inc^{\Ahsmanbe}} \Catg(\EG,\Ahsmanbe)\]
is the left adjoint of an adjoint $G$-equivalence by \cref{thm:SgoAh} \eqref{thm:SgoAh_i} applied to the object $\sma\nbe \in \FG$.  Thus, $\inc^{\Ahsmanbe}$ is a categorical weak $G$-equivalence.
\item Using the category $\Ahnbe$ of all $\nbe$-systems in $\Ah$ \pcref{def:nbeta_gcat}, the inclusion $G$-functor
\[\Ahnbe \fto[\gsim]{\inc^{\Ahnbe}} \Catg(\EG,\Ahnbe)\]
is the left adjoint of an adjoint $G$-equivalence by \cref{thm:HgoAh} \eqref{thm:HgoAh_i}.  Thus, $\inc^{\Ahnbe}$ is a categorical weak $G$-equivalence.
\end{itemize}
However, the $H$-theory comparison pointed $G$-functor \cref{PistAnbe}
\[(\smast\Sgoh\Ah)\nbe = \sys{\Ah}{(\sma \nbe)} 
\fto{\Pist_{\Ah,\nbe}} \sys{\Ah}{\nbe} = (\Hgoh\Ah)\nbe\]
is \emph{not} generally an equivalence of categories \pcref{expl:Pist_not_eq}.  Thus, \cref{merling_2.16} does \emph{not} imply that $\Catg(\EG,\Pist_{\Ah,\nbe})$ is a categorical weak $G$-equivalence.  This is where the proof of \cref{thm:pistweakgeq} would fail if $\Pist$ is used instead of the strong variant $\Pistsg$.
\end{explanation}

%% file: chap/special.tex
This chapter develops the closely related concepts of special objects and weak $G$-equivalences for both $\GGG$-categories and $\Gskg$-categories.  Some of the results in this chapter hold for arbitrary groups, while some only hold for finite groups.  The notions defined in this chapter extend simpler ones for $G$-spaces and the indexing categories $\Fsk$ and $\FG$, from the works of Segal \cite{segal_eht}, Shimakawa \cite{shimakawa91}, and May, Merling, and Osorno \cite{mmo}.  See \cref{expl:sp_gggcat,expl:gph_subgrp,expl:sp_ggcat,rk:shi_cor,expl:ggcat_weq,expl:ggcat_weq}.

\applications

The main applications of special objects and weak $G$-equivalences are strong $H$-theory, strong $J$-theory, and their comparison for $\Oph$-pseudoalgebras of the form $\Ah = \Catg(\EG,\A)$ with $\Oph = \Catg(\EG,\Op)$.  
\begin{description}
\item[$H$-theory] The strong $H$-theory comparison $G$-natural transformation in \cref{thm:pistweakgeq} is a weak $G$-equivalence
\begin{equation}\label{PistAh_chpi}
\begin{tikzpicture}[vcenter]
\def\s{25}
\draw[0cell]
(0,0) node (a1) {\GG}
(a1)++(2.5,0) node (a2) {\phantom{\GG}}
(a2)++(.15,0) node (a2') {\Catgst}
;
\draw[1cell=.9]
(a1) edge[bend left=\s] node {\smast \Sgohsg\Ah} (a2)
(a1) edge[bend right=\s] node[swap] {\Hgohsg\Ah} (a2) 
;
\draw[2cell]
node[between=a1 and a2 at .4, rotate=-90, 2label={above,\Pistsg_{\Ah}}] {\Rightarrow}
;
\end{tikzpicture}
\end{equation}
in $\GGCatg$ from Shimakawa strong $H$-theory to strong $H$-theory \pcref{ex:ggcat_weq}.  Moreover, (Shimakawa) strong $H$-theory for $\Ah$ is a special $\GGG$-category \pcref{thm:h_special,cor:h_special_shi}.  A key part of the proof of \cref{thm:h_special} is also used in \cref{thm:gmmo_shi} \cref{thm:gmmo_shi_iii} to prove the equivalence between Shimakawa and GMMO $K$-theories.  See \cref{rk:zdsg_gmmo}.
\item[$J$-theory] For a \emph{finite} group $G$, there is a weak $G$-equivalence
\begin{equation}\label{igstpist_jthy}
\begin{tikzpicture} [vcenter]
\def\s{27}
\draw[0cell]
(0,0) node (a1) {\Gsk}
(a1)++(2.5,0) node (a2) {\phantom{\Gsk}}
(a2)++(.27,0) node (a2') {\Gcatst}
;
\draw[1cell=.9]
(a1) edge[bend left=\s] node {\smast \Jgohssg\Ah} (a2)
(a1) edge[bend right=\s] node[swap] {\Jgohsg\Ah} (a2) 
;
\draw[2cell=.9]
node[between=a1 and a2 at .33, rotate=-90, 2label={above,\igst\Pistsg_{\Ah}}] {\Rightarrow}
;
\end{tikzpicture}
\end{equation}
in $\GGCatii$ from Shimakawa strong $J$-theory to strong $J$-theory \pcref{thm:pistsgah_weq}.  Moreover, (Shimakawa) strong $J$-theory for $\Ah$ is a special $\Gskg$-category \pcref{thm:j_special,cor:j_special_shi}.
\end{description}  
Thus, the strong $H$-theory and strong $J$-theory comparisons are weak $G$-equivalences between special objects.

\subtleties

It is useful to keep in mind the following points.
\begin{description}
\item[Empty tuple]
The definitions of special $\GGG$-categories and special $\Gskg$-categories do not involve the empty tuple $\ang{}$ because the Segal functor at $\ang{}$ is simply the identity \pcref{expl:sp_gggcat_mar,expl:sp_ggcat}.  On the other hand, the definitions of weak $G$-equivalences in $\GGCatg$ and $\GGCatii$ involve the empty tuple because the $\ang{}$-component of a 1-cell is not generally a categorical weak $G$-equivalence \pcref{expl:ggcatg_weq,expl:ggcat_weq}.
\item[Graph subgroups] For $\GGG$-categories, special objects and weak $G$-equivalences \pcref{def:sp_gggcat,def:ggg_weq} are defined in terms of categorical weak $G$-equivalences \pcref{def:cat_weakg}.  On the other hand, for $\Gskg$-categories, special objects and weak $G$-equivalences \pcref{def:sp_ggcat,def:gg_weq} are defined in terms of categorical weak $K$-equivalences for graph subgroups $K \subseteq \GSin$ \pcref{def:graph_subgrp}.  This asymmetry is reconciled by the fact that, for a finite group $G$, these notions correspond to each other under the adjoint 2-equivalence \pcref{thm:ggcat_ggcatg_iieq}
\begin{equation}\label{Lgigst_chpi}
\begin{tikzpicture}[vcenter]
\draw[0cell]
(0,0) node (a1) {\GGCatii}
(a1)++(2.5,0) node (a2) {\GGCatg}
;
\draw[1cell=.9]
(a1) edge[transform canvas={yshift=.5ex}] node {\Lg} (a2)
(a2) edge[transform canvas={yshift=-.4ex}] node {\igst} (a1)
;
\end{tikzpicture}
\end{equation}
relating $\Gskg$-categories and $\GGG$-categories.  See \cref{ggcat_sp_LX,gggcat_sp_iX,Lg_weq,igst_weq}. 
\item[Finiteness of $G$] The assertions that $\Lg$ and $\igst$ preserve and reflect special objects and weak $G$-equivalences require the group $G$ to be finite.  \cref{rk:shi_cor} discusses why the finiteness of $G$ is necessary in these results.
\item[$H$ vs. $J$] There is an asymmetry between $H$-theory and $J$-theory.  Results about strong $H$-theory---\cref{ex:ggcat_weq,thm:h_special,cor:h_special_shi}---hold for arbitrary groups $G$.  On the other hand, results about strong $J$-theory---\cref{thm:pistsgah_weq,thm:j_special,cor:j_special_shi}---hold for \emph{finite} groups $G$ because they involve \cref{igst_weq,gggcat_sp_iX}.
\end{description}

\organization
This chapter consists of the following sections.

\secname{sec:sp_gggcat}
This section defines Segal functors and special $\GGG$-categories.

\secname{sec:sp_ggcat}
This section discusses graph subgroups and defines Segal functors and special $\Gskg$-categories.

\secname{sec:char_special}
This section characterizes special $\Gskg$-categories in terms of subgroups $H \subseteq G$ and the left 2-adjoint $\Lh$ in the adjoint 2-equivalence \pcref{thm:ggcat_ggcatg_iieq} 
\begin{equation}\label{Lhihst_chpi}
\begin{tikzpicture}[vcenter]
\draw[0cell]
(0,0) node (a1) {\GHCatii}
(a1)++(2.5,0) node (a2) {\GHCath}
;
\draw[1cell=.9]
(a1) edge[transform canvas={yshift=.5ex}] node {\Lh} (a2)
(a2) edge[transform canvas={yshift=-.4ex}] node {\ihst} (a1)
;
\end{tikzpicture}
\end{equation}
applied to $H$.

\secname{sec:specials}
This section proves that, for a finite group $G$, the 2-equivalences $\Lg$ and $\igst$ preserve and reflect special objects.

\secname{sec:gggcat_weq}
This section defines weak $G$-equivalences between $\GGG$-categories.  The strong $H$-theory comparison $\Pistsg_{\Ah}$ in \cref{thm:pistweakgeq} is a weak $G$-equivalence in $\GGCatg$.  Weak $G$-equivalences preserve and reflect special objects.  Between special $\GGG$-categories, weak $G$-equivalences are detected at objects $(\ordone,\ldots,\ordone) \in \GG$ consisting of copies of $\ordone = \{0,1\}$.

\secname{sec:ggcat_weq}
This section defines weak $G$-equivalences between $\Gskg$-categories and proves results analogous to those in \cref{sec:gggcat_weq}.

\secname{sec:comparing_weq}
This section proves that, for a finite group $G$, the 2-equivalences $\Lg$ and $\igst$ preserve and reflect weak $G$-equivalences.

\secname{sec:weq_app}
This section applies \cref{igst_weq} to show that the strong $J$-theory comparison \cref{igstpist_jthy} for an $\Oph$-pseudoalgebra $\Ah$ is a weak $G$-equivalence in $\GGCatii$.

\secname{sec:special_theory}
This section proves that (Shimakawa) strong $H$-theory and strong $J$-theory for an $\Oph$-pseudoalgebra $\Ah$ are special.

\section{Special $\GGG$-Categories}
\label{sec:sp_gggcat}

This section discusses Segal functors and special $\GGG$-categories for an arbitrary group $G$.

\secoutline
\begin{itemize}
\item \cref{def:angordone} defines $\bdi$-characteristic morphisms $\segbdi$ for $\bdi \in \ufs{n_1 \Cdots n_q}$.  These morphisms are used to define $\nbe$-Segal functors.
\item \cref{def:sp_gggcat} defines $\nbe$-Segal functors and special $\GGG$-categories.
\item \cref{segnbe_welldef} verifies that the $\nbe$-Segal functor is $G$-equivariant.  
\end{itemize}

\begin{definition}[Characteristic Morphisms]\label{def:angordone}
Suppose $q \geq 0$. 
\begin{itemize}
\item The object \cref{Gsk_objects}
\begin{equation}\label{angordone}
\angordone = \ang{\ordone}_{j \in \ufsq} \in \Gsk
\end{equation}
consists of $q$ copies of the pointed finite set $\ordone = \{0,1\}$.  Given any group $G$, equipping $\ordone$ with the trivial $G$-action, $\angordone$ is also regarded as an object in $\GG$ \cref{GG_objects}.  
\item Suppose $\nbe = \sordi{n}{\be}{j}_{j \in \ufsq} \in \GG \setminus \{\vstar,\ang{}\}$ is an object of length $q>0$, and 
\[\bdi = \ang{i_j \in \ufsn_j}_{j \in \ufs{q}} \in \ufs{n_1 \Cdots n_q}\]
is a $q$-tuple of positive integers \cref{bdi}.  The \emph{$\bdi$-characteristic morphism}\index{characteristic morphism} \cref{GG_morphisms}
\begin{equation}\label{i-char}
\nbe \fto{\segbdi = (1_{\ufsq}, \ang{\seg_{i_j}}_{j \in \ufsq})} \angordone = \ang{\ordone}_{j \in \ufsq} \inspace \GG
\end{equation}
consists of
\begin{itemize}
\item the identity reindexing function $1_{\ufsq}$ and
\item for each $j \in \ufsq$, the \emph{$i_j$-characteristic function}\index{characteristic function}
\begin{equation}\label{ij-char}
\ordn_j \fto{\seg_{i_j}} \ordone \definedby 
\seg_{i_j}(k) = 
\begin{cases} 1 & \text{if $k=i_j$},\\ 
0 & \text{if $k \neq i_j$.}
\end{cases}
\end{equation}
\end{itemize}
\item For an object $\angordn \in \Gsk \setminus \{\vstar,\ang{}\}$, the $\bdi$-characteristic morphism \cref{Gsk_morphisms}
\begin{equation}\label{i_char}
\angordn \fto{\segbdi = (1_{\ufsq}, \ang{\seg_{i_j}}_{j \in \ufsq})} \angordone \inspace \Gsk
\end{equation}
is defined in the same way as \cref{i-char}.\defmark
\end{itemize}
\end{definition}

Recall that a $\GGG$-category is a pointed $G$-functor $\GG \to \Catgst$ \cref{ggcatg_obj}; see \cref{expl:GGcategory}.

\begin{definition}[Segal Functors and Specialness]\label{def:sp_gggcat}
Suppose 
\[(\GG,\vstar) \fto{X} (\Catgst,\bone)\]
is a pointed $G$-functor for a group $G$, and $\nbe \in \GG \setminus \{\vstar,\ang{}\}$ is an object of length $q>0$.  
\begin{itemize}
\item For the object $\angordone \in \GG$ of length $q$ \cref{angordone}, 
\begin{equation}\label{Xonenbe}
\Xonenbe
\end{equation} 
denotes the $\nbe$-twisted product \pcref{def:proCnbe} of the pointed $G$-category $X\angordone$. 
\item The \emph{$\nbe$-Segal functor}\index{Segal functor} for $X$ is the pointed $G$-functor
\begin{equation}\label{nbe_segal}
X\nbe \fto{\segnbe} \Xonenbe
\end{equation}
whose $\bdi$-th coordinate is the pointed functor
\begin{equation}\label{Xsegbdi}
X\nbe \fto{X\segbdi} X\angordone
\end{equation}
with $\segbdi$ the $\bdi$-characteristic morphism \cref{i-char} for each index $\bdi \in \ufs{n_1 \Cdots n_q}$.  The $G$-equivariance of $\segnbe$ is proved in \cref{segnbe_welldef}.  We write $\segnbe$ as $\segnbex$ if we want to emphasize $X$.
\end{itemize}
A $\GGG$-category $X$ is \emph{special}\index{special!GGG-category@$\GGG$-category} if, for each object $\nbe \in \GG \setminus \{\vstar,\ang{}\}$, the $\nbe$-Segal functor $\segnbe$ \cref{nbe_segal} is a categorical weak $G$-equivalence \pcref{def:cat_weakg}.
\end{definition}

\begin{lemma}\label{segnbe_welldef}
In \cref{def:sp_gggcat}, the $\nbe$-Segal functor \cref{nbe_segal}
\[X\nbe \fto{\segnbe} \Xonenbe\]
is a pointed $G$-functor.
\end{lemma}

\begin{proof}
For each $\bdi \in \ufs{n_1 \Cdots n_q}$, $X\segbdi$ \cref{Xsegbdi} is a pointed functor, so $\segnbe$ is a pointed functor.

The $G$-equivariance of the pointed functor $\segnbe$ is proved in \cref{segnbe_geq} after some preliminary observations.  The following equalities show that, for each $g \in G$, the $\ginv$-action on the $\bdi$-characteristic morphism $\segbdi$ \cref{i-char} yields the $(\ginv\bdi)$-characteristic morphism in $\GG$.
\begin{equation}\label{ginv_segbdi}
\begin{aligned}
& \ginv \cdot \segbdi &&\\
&= \ginv \cdot (1_{\ufsq}, \ang{\seg_{i_j}}_{j \in \ufsq}) && \text{by \cref{i-char}}\\
&= \big(1_{\ufsq}, \ang{\ginv \seg_{i_j} g}_{j \in \ufsq} \big) && \text{by \cref{GG_Gaction}}\\
&= \big(1_{\ufsq}, \ang{\seg_{i_j} g}_{j \in \ufsq} \big) && \text{by \cref{angordone}}\\
&= \big(1_{\ufsq}, \ang{\seg_{\ginv i_j}}_{j \in \ufsq} \big) && \text{by \cref{ij-char}}\\
&= \seg_{\ginv\bdi} && \text{by \cref{ginvbdi}}
\end{aligned}
\end{equation}
By \cref{GGGcat_Gequiv,ginv_segbdi}, there are equalities of functors
\[X\seg_{\ginv\bdi} = X(\ginv \cdot \segbdi) = \ginv (X\segbdi) g,\]
which can be rewritten as
\begin{equation}\label{gXsegginvbdi}
g(X\seg_{\ginv\bdi}) = (X\segbdi) g.
\end{equation}
For an object or a morphism $x \in X\nbe$, the following equalities in the $\nbe$-twisted product $\Xonenbe$ \cref{Xonenbe} prove that $\segnbe$ is $G$-equivariant.
\begin{equation}\label{segnbe_geq}
\begin{aligned}
& \segnbe(gx) && \\
&= \ang{(X\segbdi)(gx)}_{\bdi \in \ufs{n_1\cdots n_q}} && \text{by \cref{Xsegbdi}}\\
&= \ang{g(X\seg_{\ginv\bdi}) (x)}_{\bdi \in \ufs{n_1\cdots n_q}} && \text{by \cref{gXsegginvbdi}}\\
&= g\ang{(X\segbdi)x}_{\bdi \in \ufs{n_1\cdots n_q}} && \text{by \cref{proCnbe_gaction}}\\
&= g(\segnbe x) && \text{by \cref{Xsegbdi}}
\end{aligned}
\end{equation}
This finishes the proof.
\end{proof}

\begin{explanation}[Specialness]\label{expl:sp_gggcat}
A $\GGG$-category $X \cn \GG \to \Catgst$ is special if and only if, for each object $\nbe \in \GG \setminus \{\vstar,\ang{}\}$ and each subgroup $H \subseteq G$, the $H$-fixed subfunctor
\begin{equation}\label{segnbeH}
(X\nbe)^H \fto{\segnbe^H} (\Xonenbe)^H
\end{equation}
is an equivalence of categories.  Special $\GGG$-categories are analogous to Shimakawa's \emph{special $\GaG$-spaces} \cite[p.\ 225]{shimakawa91} and Segal's \emph{$G\mh\Ga$-spaces} \cite[p.\ 19]{segal_eht}, which concern $G$-spaces and the indexing $G$-category $\FG$ \pcref{def:FG}.  See also \cite[Def.\ 2.34]{mmo}.
\end{explanation}

\begin{explanation}[Marginal Cases]\label{expl:sp_gggcat_mar}
The definition of a special $\GGG$-category is about the $\nbe$-Segal functors $\segnbe$ for objects $\nbe \in \GG \setminus \{\vstar,\ang{}\}$.  The basepoint $\vstar$ and the empty tuple $\ang{}$ are excluded because they lead to empty conditions.
\begin{itemize}
\item For the basepoint $\vstar \in \GG$, the only reasonable definition of $X\angordone^{\vstar}$ is the terminal $G$-category $\bone$.  The $\vstar$-Segal functor 
\[X\vstar = \bone \fto{\seg_{\vstar}} X\angordone^{\vstar} = \bone\]
is the identity functor on $\bone$.
\item For the empty tuple $\ang{} \in \GG$, the only reasonable definition of $X\ang{}^{\ang{}}$ is the pointed $G$-category $X\ang{}$, and the $\ang{}$-Segal functor
\[X\ang{} \fto{\seg_{\ang{}}} X\ang{}^{\ang{}}\]
should be defined as the identity functor on $X\ang{}$.
\end{itemize}
Thus, for the objects $\vstar$ and $\ang{}$, the Segal functors are identities, which are automatically categorical weak $G$-equivalences.
\end{explanation}

\section{Special $\Gskg$-Categories}
\label{sec:sp_ggcat}

This section discusses Segal functors and special $\Gskg$-categories for an arbitrary group $G$.

\secoutline
\begin{itemize}
\item \cref{def:graph_subgrp} defines graph subgroups of $\GSin$, where $\Siangordn$ is the product $\txprod_{j=1}^q \Si_{n_j}$.  Graph subgroups are used to define special $\Gskg$-categories.
\item \cref{gph_subgrp} characterizes graph subgroups in terms of subgroups of $G$ equipped with a homomorphism to $\Siangordn$.
\item \cref{def:sp_ggcat} defines $\angordn$-Segal functors and special $\Gskg$-categories.
\item \cref{segn_welldef} verifies that the $\angordn$-Segal functor is $\GSinp$-equivariant.
\end{itemize}

\subsection*{Graph Subgroups}

\begin{definition}[Graph Subgroups]\label{def:graph_subgrp}
Suppose $\angordn = \ang{\ordn_j}_{j \in \ufsq} \in \Gsk \setminus \{\vstar\}$ is an object of length $q \geq 0$ \cref{Gsk_objects}. 
\begin{itemize}
\item The product group of the $n_j$-th symmetric group $\Si_{n_j}$ for $j \in \ufsq$ is denoted by 
\begin{equation}\label{Siangordn}
\Siangordn = \txprod_{j=1}^q \Si_{n_j}.
\end{equation}
For the empty tuple $\ang{}$, $\Siang$ is defined to be the trivial group $\{e\}$.
\item For a group $G$, a \emph{graph subgroup}\index{graph subgroup} of $\GSin$ is a subgroup $K$ such that
\[K \cap \Siangordn = \{e\}\]
with $\Siangordn$ denoting the subgroup $\{e\} \ttimes \Siangordn \subseteq \GSin$ and $e$ denoting the group unit.  
\item The family of graph subgroups of $\GSin$ is denoted by \label{not:FGSin}$\FGSin$.\defmark
\end{itemize}
\end{definition}

\begin{example}\label{ex:FGSione}
Suppose $\angordone = \ang{\ordone}_{j \in \ufsq} \in \Gsk$ consists of $q>0$ copies of $\ordone = \{0,1\}$.  Then 
\[\Siangordone = \txprod_{j=1}^q \Si_1\]
is the trivial group.  Under the first-factor projection, graph subgroups of 
\[G \ttimes \Siangordone \iso G \ttimes \{e\}\]
are precisely subgroups of $G$.
\end{example}

The name graph subgroup comes from the following basic observation.

\begin{lemma}\label{gph_subgrp}
In \cref{def:graph_subgrp}, a subgroup $K \subseteq \GSin$ is a graph subgroup if and only if there exist a unique subgroup $H \subseteq G$ and a unique homomorphism $\be \cn H \to \Siangordn$ such that
\[K = \{(h, \be(h)) \tmid h \in H \}.\]
\end{lemma}

\begin{proof}
A subgroup of $\GSin$ of the form $\{(h, \be(h)) \tmid h \in H \}$ intersects $\Siangordn$ at $\{e\}$ because $\be(e) = e$.  Conversely, if $K \subseteq \GSin$ satisfies $K \cap \Siangordn = \{e\}$, then the composite homomorphism
\[K \subseteq \GSin \fto{\pr_1} G\]
is injective, where $\pr_1$ denotes the first-factor projection.  Thus, it is an isomorphism onto its image, say, $H \subseteq G$.  Under the projection $\pr_1$, each $h \in H$ has a unique preimage $(h,\be(h)) \in K$.  The assignment $h \mapsto \be(h)$ defines a homomorphism $\be \cn H \to \Siangordn$.
\end{proof}

\begin{explanation}\label{expl:gph_subgrp}
In the context of \cref{gph_subgrp}, via the $j$-th projection $\pr_j \cn \Siangordn \to \Si_{n_j}$, a homomorphism $\be \cn H \to \Siangordn$ is equivalent to a family of homomorphisms 
\begin{equation}\label{betaj}
H \fto{\be_j} \Si_{n_j} \forspace j \in \ufsq.
\end{equation}  
Such a family is equivalent to an object \cref{GG_objects}
\begin{equation}\label{nbe_GH}
\nbe = \sordi{n}{\be}{j}_{j \in \ufsq} \in \GH \setminus \{\vstar\},
\end{equation}
meaning a $q$-tuple of pointed finite $H$-sets.  \cref{def:graph_subgrp} is a generalization to the case $q \geq 0$ of a definition due to Shimakawa \cite[p.\ 226]{shimakawa91}, which also appears in \cite[Lemma 4.5]{gmm17} and \cite[Def.\ 1.2]{mmo}.  The name graph subgroup follows the usage in \cite[Def.\ 1.1.25]{schwede_ght}.
\end{explanation}

\subsection*{Special Objects}

Recall that a $\Gskg$-category is a pointed functor $\Gsk \to \Gcatst$ \cref{ggcat_obj}; see \cref{expl:GGCat_iicat}. 

\begin{definition}[Segal Functors and Specialness]\label{def:sp_ggcat}
Suppose 
\[(\Gsk,\vstar) \fto{X} (\Gcatst,\bone)\]
is a pointed functor, and $\angordn = \ang{\ordn_j}_{j \in \ufsq} \in \Gsk \setminus \{\vstar,\ang{}\}$ is an object of length $q>0$. 
\begin{itemize}
\item For each $q$-tuple of permutations \cref{Siangordn}
\begin{equation}\label{Siangord_obj}
\bsi = \ang{\si_j \in \Si_{n_j}}_{j \in \ufsq} \in \Siangordn,
\end{equation}
there is an isomorphism \cref{Gsk_morphisms}
\[\angordn \fto[\iso]{(1_{\ufsq}, \bsi)} \angordn \inspace \Gsk,\]
where each $\si_j \cn \ordn_j \fiso \ordn_j$ is the pointed bijection determined by the permutation $\si_j$.  Thus,
\[X\angordn \fto[\iso]{X(1_{\ufsq}, \bsi)} X\angordn\]
is a pointed $G$-isomorphism \cref{f_upom}.  The pointed $G$-category $X\angordn$ is also regarded as a pointed $\GSinp$-category with $(g,\bsi)$-action functor given by the composite isomorphism
\begin{equation}\label{Xn_gsin_act}
X\angordn \fto[\iso]{g} X\angordn \fto[\iso]{X(1_{\ufsq}, \bsi)} X\angordn
\end{equation}
for $(g,\bsi) \in \GSin$.  The composite in \cref{Xn_gsin_act} is also equal to $g \circ X(1_{\ufsq}, \bsi)$ because $X(1_{\ufsq}, \bsi)$ is $G$-equivariant.
\item For the object $\angordone \in \Gsk$ of length $q$ \cref{angordone}, the product $G$-category 
\[X\angordone^{n_1\cdots n_q} = \prod_{\bdi \sins \ufs{n_1 \cdots n_q}} X\angordone\]
becomes a pointed $\GSinp$-category with $(g,\bsi)$-action given by
\begin{equation}\label{xonen_gsin_act}
(g,\bsi) \cdot \ang{x_{\bdi}}_{\bdi \in \ufs{n_1\cdots n_q}} 
= \ang{gx_{\bsiinv\bdi}}_{\bdi \in \ufs{n_1\cdots n_q}} 
\end{equation}
for objects or morphisms $\ang{x_{\bdi}}_{\bdi \in \ufs{n_1\cdots n_q}} \in X\angordone^{n_1 \cdots n_q}$.  The index on the right-hand side of \cref{xonen_gsin_act} is given by
\begin{equation}\label{bsiinv_i}
\bsiinv\bdi = \ang{\sigmainv_j (i_j)}_{j \in \ufsq} \in \ufs{n_1\Cdots n_q}
\end{equation}
for each $q$-tuple of positive integers $\bdi = \ang{i_j \in \ufsn_j}_{j \in \ufsq}$ \cref{bdi}.
\item Using \cref{Xn_gsin_act,xonen_gsin_act}, the \emph{$\angordn$-Segal functor}\index{Segal functor} for $X$ is the pointed $\GSinp$-functor
\begin{equation}\label{nsegal}
X\angordn \fto{\segn} X\angordone^{n_1\cdots n_q}
\end{equation}
whose $\bdi$-th coordinate is the pointed $G$-functor
\begin{equation}\label{X_segbdi}
X\angordn \fto{X\segbdi} X\angordone
\end{equation}
with $\segbdi \cn \angordn \to \angordone$ the $\bdi$-characteristic morphism \cref{i_char} for each index $\bdi \in \ufs{n_1 \Cdots n_q}$.  The $\GSinp$-equivariance of $\segn$ is proved in \cref{segn_welldef}.  We write $\segn$ as $\segnx$ if we want to emphasize $X$.
\end{itemize}
A $\Gskg$-category $X$ is \emph{special}\index{special!G-G-category@$\Gskg$-category} if, for each object $\angordn \in \Gsk \setminus \{\vstar,\ang{}\}$ and each graph subgroup $K \in \FGSin$ \pcref{def:graph_subgrp}, the $\angordn$-Segal functor $\segn$ \cref{nsegal} is a categorical weak $K$-equivalence \pcref{def:cat_weakg}.
\end{definition}

\begin{explanation}[Specialness]\label{expl:sp_ggcat}\
\begin{enumerate}
\item\label{expl:sp_ggcat_i} The family $\FGSin$ of graph subgroups of $\GSin$ \pcref{def:graph_subgrp} is closed under passage to subgroups.  Thus, a $\Gskg$-category $X$ is special if and only if, for each object $\angordn \in \Gsk \setminus \{\vstar,\ang{}\}$ and each graph subgroup $K \in \FGSin$, the $K$-fixed subfunctor 
\begin{equation}\label{segnK}
X\angordn^K \fto{\segn^K} (X\angordone^{n_1\cdots n_q})^K
\end{equation}
is an equivalence of categories.  The basepoint $\vstar$ and the empty tuple $\ang{}$ are excluded in the definition of a special $\Gskg$-category because they do not impose any condition at all.  See \cref{expl:sp_gggcat_mar}.
\item\label{expl:sp_ggcat_ii} 
Special $\Gskg$-categories are analogous to Shimakawa's \emph{special $\Ga\mh G$-spaces} \cite[p.\ 226]{shimakawa91}, which concern $G$-spaces and the indexing category $\Fsk$ \pcref{def:Fsk}.  Special $\Ga\mh G$-spaces are called \emph{$\mathbbm{F}_{\!\bullet}$-special} $\Fskg$-spaces in \cite[Def.\ 2.14 (i)]{mmo}.\defmark
\end{enumerate}
\end{explanation}

\begin{lemma}\label{segn_welldef}
The $\angordn$-Segal functor \cref{nsegal}
\[X\angordn \fto{\segn} X\angordone^{n_1\cdots n_q}\]
is $\GSinp$-equivariant.
\end{lemma}

\begin{proof}
For each $q$-tuple of permutations $\bsi \in \Siangordn$ and each index $\bdi \in \ufs{n_1\Cdots n_q}$, by \cref{Gsk_composite,ij-char,bsiinv_i}, the diagram 
\begin{equation}\label{seg_bsiinv_i}
\begin{tikzpicture}[vcenter]
\def\v{-1.3}
\draw[0cell]
(0,0) node (a11) {\angordn}
(a11)++(2.3,0) node (a12) {\angordn}
(a11)++(0,\v) node (a21) {\angordone}
(a12)++(0,\v) node (a22) {\angordone}
;
\draw[1cell=.9]
(a11) edge node {(1_{\ufsq}, \bsi)} (a12)
(a12) edge node {\segbdi} (a22)
(a11) edge node[swap] {\seg_{\bsiinv\bdi}} (a21)
(a21) edge[equal] (a22)
;
\end{tikzpicture}
\end{equation}
in $\Gsk$ commutes.  For an element $g \in G$ and an object or a morphism $x \in X\angordn$, the following equalities in $X\angordone^{n_1\cdots n_q}$ prove that $\segn$ commutes with the $(g,\bsi)$-action functors, where $\ang{\Cdots}_{\bdi} = \ang{\Cdots}_{\bdi \in \ufs{n_1\cdots n_q}}$.
\[\begin{aligned}
& (g,\bsi) \cdot (\segn x) &&\\
&= (g,\bsi) \cdot \ang{(X\segbdi)x}_{\bdi} && \text{by \cref{X_segbdi}}\\
&= \ang{g ((X\seg_{\bsiinv\bdi}) x)}_{\bdi} && \text{by \cref{xonen_gsin_act}}\\
&= \ang{(X\seg_{\bsiinv\bdi}) (gx)}_{\bdi} && \text{by \cref{f_upom}}\\
&= \ang{(X\segbdi)(X(1_{\ufsq}, \bsi)) (gx)}_{\bdi} && \text{by \cref{seg_bsiinv_i}}\\
&= \segn\big((X(1_{\ufsq}, \bsi)) (gx)\big) && \text{by \cref{X_segbdi}}\\
&= \segn\big((g,\bsi) \cdot x\big) && \text{by \cref{Xn_gsin_act}}
\end{aligned}\]
This proves that $\segn$ is $\GSinp$-equivariant.
\end{proof}

\section{Characterizations of Special $\Gskg$-Categories}
\label{sec:char_special}

This section provides alternative characterizations of special $\Gskg$-categories \pcref{def:sp_ggcat} for an arbitrary group $G$.

\secoutline
\begin{itemize}
\item \cref{spec_ggcat} describes special $\Gskg$-categories using subgroups of $G$ instead of graph subgroups.
\item \cref{spec_ggcat_ii} shows that, for a special $\Gskg$-category $X$, its $H$-restriction $\Xh \cn \Gsk \to \Hcatst$ is a special $\Gskh$-category for each subgroup $H$ of $G$.
\item \cref{ggcat_sp_Xh} proves that a $\Gskg$-category $X$ is special if and only if the $\GHH$-category $\Lh\Xh \cn \GH \to \Cathst$ is special for each subgroup $H$ of $G$, where 
\[\GHCatii \fto{\Lh} \GHCath\]
is the 2-equivalence in \cref{thm:ggcat_ggcatg_iieq} applied to the group $H$.  For a finite group $G$, an improvement of \cref{ggcat_sp_Xh} is given in \cref{ggcat_sp_LX}.
\end{itemize}

\begin{lemma}\label{spec_ggcat}
A $\Gskg$-category 
\[\Gsk \fto{X} \Gcatst\]
is special if and only if, for each object $\angordn \in \Gsk \setminus \{\vstar,\ang{}\}$, subgroup $H \subseteq G$, and homomorphism $\be \cn H \to \Siangordn$, the $\angordn$-Segal functor
\begin{equation}\label{xnbeta_segn}
\Xnbeta \fto{\segn} \Xonenbe
\end{equation}
is a categorical weak $H$-equivalence.
\end{lemma}

\begin{proof}
First, we clarify the domain and codomain $H$-categories of the $\angordn$-Segal functor $\segn$ in \cref{xnbeta_segn}.  
\begin{description}
\item[Domain]
The pointed $G$-category $X\angordn$ is regarded as a pointed $H$-category by restricting the $G$-action to the subgroup $H$.  The pointed $H$-category $\Xnbeta$ \cref{Xnbeta} is the pointed category $X\angordn$ with, for each $h \in H$, $h$-action functor given by the following composite \cref{Xnbeta_gaction}.
\begin{equation}\label{Xnbeta_haction}
\begin{tikzpicture}[vcenter]
\def\u{.65}
\draw[0cell]
(0,0) node (a1) {X\angordn}
(a1)++(1.8,0) node (a2) {X\angordn}
(a2)++(2,0) node (a3) {X\angordn}
;
\draw[1cell=.9]
(a1) edge node {h} (a2)
(a2) edge node {X\betah} (a3)
(a1) [rounded corners=2pt] |- ($(a2)+(-1,\u)$) -- node {h \cdot -} ($(a2)+(1,\u)$) -| (a3)
;
\end{tikzpicture}
\end{equation}
In \cref{Xnbeta_haction}, $\betah$ is the isomorphism \cref{betag}
\begin{equation}\label{betah}
\angordn \fto[\iso]{\betah = (1_{\ufsq}, \ang{\be_j h}_{j \in \ufsq})} \angordn \inspace \Gsk,
\end{equation}
and $\be_j$ is the composite of $\be$ with the $j$-th projection \cref{betaj}:
\[H \fto{\be} \Siangordn \fto{\pr_j} \Si_{n_j}.\]
\item[Codomain]
For the object $\angordone \in \Gsk$ of length $q$ \cref{angordone}, the pointed $G$-category $X\angordone$ becomes a pointed $H$-category by restricting the $G$-action to the subgroup $H$.  The pointed $H$-category $\Xonenbe$ is the $\nbe$-twisted product \pcref{def:proCnbe} for the object $\nbe \in \GH \setminus \{\vstar,\ang{}\}$ \cref{nbe_GH}.
\item[Equivariance]
Since the underlying pointed categories of $\Xnbeta$ and $\Xonenbe$ are, respectively, $X\angordn$ and $X\angordone^{n_1\cdots n_q}$, the $\angordn$-Segal functor $\segn$ is a well-defined pointed functor.  Its coordinates are the pointed $G$-functors $X\segbdi \cn X\angordn \to X\angordone$ \cref{X_segbdi}.  The pair $(H,\be)$ specifies a unique graph subgroup \pcref{gph_subgrp}
\[K = \{(h, \be(h)) \tmid h \in H \}.\]
By \cref{Xn_gsin_act,xonen_gsin_act,Xnbeta_haction,proCnbe_gaction}, the $K$-equivariance of the $\angordn$-Segal functor \pcref{segn_welldef}
\begin{equation}\label{X_Kequiv}
X\angordn \fto{\segn} X\angordone^{n_1\cdots n_q}
\end{equation}
is precisely the $H$-equivariance of the $\angordn$-Segal functor \cref{xnbeta_segn}
\begin{equation}\label{X_Hequiv}
\Xnbeta \fto{\segn} \Xonenbe.
\end{equation}
\item[Assertion]
For the \emph{only if} assertion, recall that $\segn$ \cref{xnbeta_segn} is a categorical weak $H$-equivalence if, for each subgroup $L \subseteq H$, the $L$-fixed subfunctor 
\[\Xnbeta^L \fto{\segn^L} (\Xonenbe)^L\]
is an equivalence of categories \pcref{def:cat_weakg}.  Since $\segn^L$ is defined by the $L$-action, it remains the same if $\be \cn H \to \Siangordn$ is replaced by its restriction $\be|_L \cn L \to \Siangordn$ to $L$.  The pair $(L,\be|_L)$ specifies a graph subgroup $K' \subseteq \GSin$ \pcref{gph_subgrp}.  Thus, $\segn^L = \segn^{K'}$ is an equivalence by the assumption that $X$ is special.

For the \emph{if} assertion, recall that a $\Gskg$-category $X$ is special if, for each graph subgroup $K \subseteq \GSin$, the $K$-fixed subfunctor \cref{segnK}
\[X\angordn^K \fto{\segn^K} (X\angordone^{n_1\cdots n_q})^K\] 
is an equivalence. 
Using \cref{X_Kequiv,X_Hequiv}, $\segn^K = \segn^H$ is an equivalence by assumption.\qedhere
\end{description}
\end{proof}

Using \cref{def:ghcat}, \cref{spec_ggcat_ii} shows that specialness is preserved under passage to subgroups.

\begin{definition}[Restrictions]\label{def:ghcat}
For a $\Gskg$-category $X \cn \Gsk \to \Gcatst$ \cref{ggcat_obj} and a subgroup $H \subseteq G$, the $\Gskh$-category
\begin{equation}\label{Xh}
(\Gsk,\vstar) \fto{\Xh} (\Hcatst,\bone),
\end{equation}
called the \emph{$H$-restriction}\index{restriction} of $X$, is obtained from $X$ by
\begin{itemize}
\item restricting the $G$-action on the small pointed $G$-category $X\angordm$ to the subgroup $H$ for each object $\angordm \in \Gsk$ and
\item regarding the pointed $G$-functor $X\upom$ as a pointed $H$-functor for each morphism $\upom$ in $\Gsk$.
\end{itemize}
Here, $\Hcatst$ is the 2-category of small pointed $H$-categories, pointed $H$-functors, and pointed $H$-natural transformations \pcref{def:gcatst}.  Note that the domain and codomain of the $\angordn$-Segal functor $\segn \cn \Xnbeta \to \Xonenbe$ \cref{xnbeta_segn} are $\Xhnbeta$ and $\Xhonenbe$.
\end{definition}

\begin{lemma}\label{spec_ggcat_ii}
If 
\[\Gsk \fto{X} \Gcatst\] 
is a special $\Gskg$-category, then 
\[\Gsk \fto{\Xh} \Hcatst\] 
is a special $\Gskh$-category for each subgroup $H \subseteq G$.
\end{lemma}

\begin{proof}
This follows from \cref{spec_ggcat} because a subgroup of $H$ is also a subgroup of $G$.
\end{proof}

\begin{notation}\label{not:Lh_ihst}
For a subgroup $H \subseteq G$, we denote by
\begin{equation}\label{Lhihst_not_diag}
\begin{tikzpicture}[vcenter]
\draw[0cell]
(0,0) node (a1) {\GHCatii}
(a1)++(2.5,0) node (a2) {\GHCath}
;
\draw[1cell=.9]
(a1) edge[transform canvas={yshift=.5ex}] node {\Lh} (a2)
(a2) edge[transform canvas={yshift=-.4ex}] node {\ihst} (a1)
;
\end{tikzpicture}
\end{equation}
the adjoint 2-equivalence in \cref{thm:ggcat_ggcatg_iieq} applied to the group $H$.  The right 2-adjoint $\ihst$ is induced by the full subcategory inclusion $\ih \cn \Gsk \to \GH$ \cref{ig}.  The left 2-adjoint $\Lh$ is componentwise a coend \cref{Lg_f_nbe}.
\end{notation}

Applying the left 2-adjoint $\Lh$ to the $H$-restriction $\Xh \cn \Gsk \to \Hcatst$ \cref{Xh} yields a $\GHH$-category 
\[\GH \fto{\Lh\Xh} \Cathst.\]  
\cref{ggcat_sp_Xh} compares specialness in \cref{def:sp_gggcat,def:sp_ggcat} by showing that the relevant Segal functors correspond to each other under the units of the 2-adjunctions.

\begin{lemma}\label{ggcat_sp_Xh}
A $\Gskg$-category 
\[\Gsk \fto{X} \Gcatst\] 
is special if and only if, for each subgroup $H \subseteq G$, the $\GHH$-category 
\[\GH \fto{\Lh\Xh} \Cathst\] 
is special.
\end{lemma}

\begin{proof}
We consider a triple $(\angordn, H, \be)$ consisting of
\begin{itemize}
\item an object $\angordn \in \Gsk \setminus \{\vstar,\ang{}\}$,
\item a subgroup $H \subseteq G$, and
\item a homomorphism $\be \cn H \to \Siangordn$.
\end{itemize}
By \cref{expl:gph_subgrp}, such a triple is equivalent to an object $\nbe \in \GH \setminus \{\vstar,\ang{}\}$.  We first show that, for each such triple $(\angordn, H, \be)$, there is a commutative diagram of pointed $H$-functors
\begin{equation}\label{ggcat_sp_Xh_dia}
\begin{tikzpicture}[vcenter]
\def\v{-1.4} \def\h{6.5}
\draw[0cell=.9]
(0,0) node (a11) {\Xhnbeta}
(a11)++(\h,0) node (a13) {\Xhonenbe}
(a11)++(0,\v) node (a21) {(\Lh \Xh)\angordn_{\be}}
(a13)++(0,\v) node (a23) {\LXhonenbe}
(a21)++(.52*\h,0) node (a22) {\phantom{(\Lh \Xh)\nbe}}
(a22)++(0,.03) node (a22') {(\Lh \Xh)\nbe}
;
\draw[1cell=.85]
(a11) edge node {\segn} (a13)
(a13) edge node {\ug_{\Xh,\angordone}^{\nbe}} node[swap] {\iso} (a23)
(a11) edge node {\iso} node[swap] {\ug_{\Xh,\angordn}} (a21)
(a21) edge node {(\Lh \Xh)\jin} node[swap] {\iso} (a22)
(a22) edge node {\segnbe} (a23)
;
\end{tikzpicture}
\end{equation}
constructed as follows.
\begin{description}
\item[Top] $\segn$ is the $\angordn$-Segal functor \cref{xnbeta_segn} for the $\Gskh$-category $\Xh$.  The proof of \cref{spec_ggcat}---see \cref{X_Kequiv,X_Hequiv}---shows that $\segn$ is a pointed $H$-functor.
\item[Bottom left] The isomorphism $\jin \cn \angordn \to \nbe$ in $\GH$ is given by $(1_{\ufsq}, \ang{1})$ \cref{jin_nbe}.  By \cref{Xjininv}, the functor $(\Lh \Xh)\jin$ is inverse to an $H$-isomorphism, so it is also an $H$-isomorphism.
\item[Bottom right] The $\nbe$-Segal functor $\segnbe$ \cref{nbe_segal} for the $\GHH$-category $\Lh\Xh$ is $H$-equivariant by \cref{segnbe_welldef}.
\item[Left] $\ug_{\Xh,\angordn}$ is the unit isomorphism \cref{ugf_angordm} for the adjoint 2-equivalence $(\Lh,\ihst)$ \pcref{not:Lh_ihst}.  For $h \in H$ and $\be_j = \pr_j \circ \be \cn H \to \Si_{n_j}$, recall the isomorphism \cref{betah}
\[\angordn \fto[\iso]{\betah = (1_{\ufsq}, \ang{\be_j h}_{j \in \ufsq})} \angordn \inspace \Gsk.\]
For $x \in \Xhnbeta$, the following equalities in $(\Lh \Xh)\angordn_{\be}$ prove that $\ug_{\Xh,\angordn}$ is $H$-equivariant.
\begin{equation}\label{uXh_Hequiv}
\begin{aligned}
& h \cdot (\ug_{\Xh,\angordn})(x) &&\\
&= h \cdot (1_{\angordn}; x) && \text{by \cref{ugx_angordmx}}\\
&= [(\Lh \Xh) \betah] [h(1_{\angordn}; x)] && \text{by \cref{Xnbeta_gaction}}\\
&= [(\Lh \Xh) \betah] [(1_{\angordn}; hx)] && \text{by \cref{LXnbe_gaction}}\\
&= (\betah ; hx) && \text{by \cref{Lg_f_nbe}}\\
&= \big(1_{\angordn} ; (\Xh\betah)(hx) \big) && \text{by \cref{Lg_f_nbe}}\\
&= (\ug_{\Xh,\angordn}) [(\Xh\betah)(hx)] && \text{by \cref{ugx_angordmx}}\\
&= (\ug_{\Xh,\angordn})(h \cdot x) && \text{by \cref{Xnbeta_haction}}
\end{aligned}
\end{equation}
Thus, the left vertical functor $\ug_{\Xh,\angordn}$ in \cref{ggcat_sp_Xh_dia} is an $H$-isomorphism.
\item[Right] Each entry of $\ug_{\Xh,\angordone}^{\nbe}$ is the unit isomorphism \cref{ugf_angordm}
\[\Xh\angordone \fto[\iso]{\ug_{\Xh,\angordone}} (\Lh \Xh)\angordone\]
for the adjoint 2-equivalence $(\Lh,\ihst)$.  For $\ang{x_\bdi}_{\bdi} \in \Xhonenbe$, the following equalities in $\LXhonenbe$ prove that $\ug_{\Xh,\angordone}^{\nbe}$ is $H$-equivariant, where $\ang{\Cdots}_\bdi = \ang{\Cdots}_{\bdi \in \ufs{n_1 \cdots n_q}}$.
\[\begin{aligned}
& h \cdot \big(\ug_{\Xh,\angordone}^{\nbe}\big) \ang{x_\bdi}_{\bdi} &&\\
&= h \cdot \ang{(1_{\angordone} ; x_{\bdi})}_{\bdi} && \text{by \cref{ugx_angordmx}}\\
&= \ang{h(1_{\angordone} ; x_{\hinv\bdi})}_{\bdi} && \text{by \cref{proCnbe_gaction}}\\
&= \ang{(1_{\angordone} ; hx_{\hinv\bdi})}_{\bdi} && \text{by \cref{LXnbe_gaction}}\\
&= \big(\ug_{\Xh,\angordone}^{\nbe}\big) \ang{hx_{\hinv\bdi}}_{\bdi} && \text{by \cref{ugx_angordmx}}\\
&= \big(\ug_{\Xh,\angordone}^{\nbe}\big) (h \cdot \ang{x_\bdi}_{\bdi}) && \text{by \cref{proCnbe_gaction}}
\end{aligned}\]
Thus, the right vertical functor $\ug_{\Xh,\angordone}^{\nbe}$ in \cref{ggcat_sp_Xh_dia} is an $H$-isomorphism.
\item[Commutativity] By \cref{jin_nbe,i-char,ij-char,i_char}, the diagram
\begin{equation}\label{segbdi_jin}
\begin{tikzpicture}[vcenter]
\def\h{1.8} \def\u{.6}
\draw[0cell]
(0,0) node (a1) {\angordn}
(a1)++(\h,0) node (a2) {\phantom{\nbe}}
(a2)++(0,.03) node (a2') {\nbe}
(a2)++(\h,0) node (a3) {\angordone}
;
\draw[1cell=.9]
(a1) edge node {\jin} (a2)
(a2) edge node {\segbdi} (a3)
(a1) [rounded corners=2pt] |- ($(a2)+(-1,\u)$) -- node {\segbdi} ($(a2)+(1,\u)$) -| (a3)
;
\end{tikzpicture}
\end{equation}
in $\GH$ commutes.  For $x \in \Xhnbeta$, the following equalities in $\LXhonenbe$ prove that the diagram \cref{ggcat_sp_Xh_dia} commutes, where $\ang{\Cdots}_\bdi = \ang{\Cdots}_{\bdi \in \ufs{n_1 \cdots n_q}}$.
\[\begin{aligned}
& \segnbe [(\Lh \Xh)\jin] (\ug_{\Xh,\angordn}) (x) &&\\
&= \segnbe [(\Lh \Xh)\jin] (1_{\angordn} ; x) && \text{by \cref{ugx_angordmx}}\\
&= \segnbe (\jin ; x) && \text{by \cref{Lg_f_nbe}}\\
&= \ang{((\Lh \Xh)\segbdi) (\jin ; x)}_{\bdi} && \text{by \cref{Xsegbdi}}\\
&= \ang{(\segbdi\jin ; x)}_{\bdi} && \text{by \cref{Lg_f_nbe}}\\
&= \ang{(\segbdi ; x)}_{\bdi} && \text{by \cref{segbdi_jin}}\\
&= \ang{(1_{\angordone} ; (\Xh\segbdi) x)}_{\bdi} && \text{by \cref{Lg_f_nbe}}\\
&= \big(\ug_{\Xh,\angordone}^{\nbe}\big) \ang{(\Xh\segbdi) x}_{\bdi} && \text{by \cref{ugx_angordmx}}\\
&= \big(\ug_{\Xh,\angordone}^{\nbe}\big) (\segn x) && \text{by \cref{X_segbdi}}
\end{aligned}\]
\end{description}
This proves that \cref{ggcat_sp_Xh_dia} is a commutative diagram of pointed $H$-functors.

By \cref{spec_ggcat}, the $\Gskg$-category $X$ is special if and only if the $\angordn$-Segal functor 
\[\Xhnbeta \fto{\segn} \Xhonenbe\]
is a categorical weak $H$-equivalence for each triple $(\angordn, H, \be)$.  By \cref{def:sp_gggcat}, for each subgroup $H \subseteq G$, the $\GHH$-category $\Lh\Xh$ is special if and only if  the $\nbe$-Segal functor \cref{nbe_segal}
\[(\Lh \Xh)\nbe \fto{\segnbe} \LXhonenbe\]
is a categorical weak $H$-equivalence for each object $\nbe \in \GH \setminus \{\vstar,\ang{}\}$.  The commutative diagram \cref{ggcat_sp_Xh_dia} of $H$-functors consists of $\segn$, $\segnbe$, and three $H$-isomorphisms.  Thus, $\segn$ is a categorical weak $H$-equivalence if and only if $\segnbe$ is so.
\end{proof}

\section{Equivalence of Specialness}
\label{sec:specials}

For a \emph{finite} group $G$, this section proves that special $\GGG$-categories correspond to special $\Gskg$-categories \pcref{def:sp_gggcat,def:sp_ggcat} under the adjoint 2-equivalence \pcref{thm:ggcat_ggcatg_iieq}
\begin{equation}\label{Lgigst_secti}
\begin{tikzpicture}[vcenter]
\draw[0cell]
(0,0) node (a1) {\GGCatii}
(a1)++(2.5,0) node (a2) {\GGCatg}
;
\draw[1cell=.9]
(a1) edge[transform canvas={yshift=.5ex}] node {\Lg} (a2)
(a2) edge[transform canvas={yshift=-.4ex}] node {\igst} (a1)
;
\end{tikzpicture}
\end{equation}
between the 2-categories in \cref{def:GGCat,def:ggcatg}.  

\secoutline
\begin{itemize}
\item \cref{weq_retract} proves that categorical weak $G$-equivalences are closed under retracts for any group $G$.  This observation is used in the proof of \cref{ggcat_sp_LX}.
\item \cref{ggcat_sp_LX} proves that $\Lg$ preserves and reflects specialness for a finite group $G$.  At the cost of restricting to finite groups, \cref{ggcat_sp_LX} improves \cref{ggcat_sp_Xh} by eliminating the condition about the $\GHH$-category $\Lh\Xh$ for proper subgroups $H \subsetneq G$.
\item \cref{gggcat_sp_iX} proves that $\igst$ preserves and reflects specialness for a finite group $G$.
\end{itemize}

Recall that a $G$-functor is a categorical weak $G$-equivalence if, for each subgroup of $G$, the fixed-point subfunctor is an equivalence of categories \pcref{def:cat_weakg}.  

\begin{lemma}\label{weq_retract}
For an arbitrary group $G$, suppose
\begin{equation}\label{weq_retract_dia}

\end{equation}
is a commutative diagram of $G$-functors such that $\fun \cn \C \to \D$ is a categorical weak $G$-equivalence.  Then $\fuk \cn \A \to \B$ is a categorical weak $G$-equivalence.
\end{lemma}

\begin{proof}
For a subgroup $H \subseteq G$, we verify that the $H$-fixed subfunctor $\fuk^H \cn \A^H \to \B^H$ is an equivalence of categories.  Restricting the commutative diagram \cref{weq_retract_dia} to the $H$-fixed subcategories and subfunctors yields the following commutative diagram of solid-arrow functors.
\begin{equation}\label{weq_retract_H}
%
\end{equation}
By the assumption on $\fun$, the $H$-fixed subfunctor $\fun^H$ is an equivalence of categories.  Thus, it admits an inverse functor $\Jrp$, depicted as the dashed arrow in \cref{weq_retract_H}, and natural isomorphisms
\begin{equation}\label{JrpfH}
\Jrp \circ \fun^H \iso 1_{\C^H} \andspace \fun^H \circ \Jrp \iso 1_{\D^H}.
\end{equation}
The commutative diagram \cref{weq_retract_H} and the natural isomorphisms in \cref{JrpfH} imply that the composite functor
\begin{equation}\label{inverse_kH}
%
\end{equation}
is an inverse of $\fuk^H$.  This proves that $\fuk^H$ is an equivalence of categories and that $\fuk$ is a categorical weak $G$-equivalence.
\end{proof}

\begin{theorem}\label{ggcat_sp_LX}
For a finite group $G$, a $\Gskg$-category 
\[\Gsk \fto{X} \Gcatst\] 
is special if and only if the $\GGG$-category 
\[\GG \fto{\Lg X} \Catgst\] 
is special.
\end{theorem}

\begin{proof}
The \emph{only if} assertion follows from the \emph{only if} assertion in \cref{ggcat_sp_Xh}, applied to the group $G$ itself.  This assertion does not require $G$ to be finite.

For the \emph{if} assertion, assuming that the $\GGG$-category $\Lg X$ is special, we use \cref{ggcat_sp_Xh} and verify that, for each proper subgroup $H \subsetneq G$, the $\GHH$-category 
\[\GH \fto{\Lh\Xh} \Cathst\]
is special.   By \cref{def:sp_gggcat}, $\Lh\Xh$ is special if and only if the $\nbe$-Segal functor \cref{nbe_segal}
\begin{equation}\label{segnbe_LXh}
(\Lh \Xh)\nbe \fto{\segnbe} \LXhonenbe
\end{equation}
is a categorical weak $H$-equivalence for each object 
\begin{equation}\label{ggcat_sp_nbe}
\nbe = \sordi{n}{\be}{j}_{j\in \ufsq} \in \GH \setminus \{\vstar,\ang{}\}
\end{equation}
with length $q>0$.  The rest of this proof shows that the $\nbe$-Segal functor $\segnbe$ is the retract of a categorical weak $H$-equivalence.  \cref{weq_retract} then implies that $\segnbe$ is a categorical weak $H$-equivalence.

\newcounter{ggcatspstep}
\def\ggcatspstep{Step \stepcounter{ggcatspstep}\arabic{ggcatspstep}}
\begin{description}
\item[\ggcatspstep. Enlarging $\ordn_j^{\be_j}$]  
Using the finiteness of $G$, this step enlarges the pointed finite $H$-set $\ordn_j^{\be_j}$ to a pointed finite $G$-set $\ordm_j^{\al_j}$, along with some auxiliary objects and morphisms.  Choosing an ordering of $G$, we denote by\label{not:Gpl}
\[\Gpl = G \sqcup {*}\] 
the pointed finite $G$-set \cref{ordn_be} obtained from $G$ by adjoining a disjoint $G$-fixed basepoint $*$, with $G$ acting on itself by left multiplication.  For each $j \in \ufsq$, starting from the pointed finite $H$-set $\ordn_j^{\be_j}$, we define the pointed finite $G$-set
\begin{equation}\label{ordmj_alphaj}
\ordm_j^{\al_j} = \Gpl \sma_H \ordn_j^{\be_j}
\end{equation}
as the quotient of the smash product $\Gpl \sma \ordn_j^{\be_j}$ by the identification
\[(gh,a) \sim (g,ha)\]
for $g \in G$, $h \in H$, and $a \in \ordn$.  With $e \in G$ denoting the group unit, the $H$-invariant subset $\{e\} \sma_H \ordn_j^{\be_j}$ of $\ordm_j^{\al_j}$ is canonically isomorphic to $\ordn_j^{\be_j}$.  Using this canonical isomorphism and restricting the $G$-action on $\ordm_j^{\al_j}$ to the subgroup $H$, there is a wedge decomposition of pointed finite $H$-sets
\begin{equation}\label{mnr_decomp}
\ordm_j^{\al_j} = \ordm_{j,1}^{\al_j} \wed \ordm_{j,2}^{\al_j} \withspace 
\ordm_{j,1}^{\al_j} = \ordn_j^{\be_j}.
\end{equation}
For each $\epz_j \in \{1,2\}$, there are pointed $H$-equivariant inclusion and projection
\begin{equation}\label{wprinc_onetwo}
\ordm_{j,\, \epz_j}^{\al_j} \fto{\winc^{\epz_j}_j} \ordm_j^{\al_j} \fto{\wpr^{\epz_j}_j} \ordm_{j,\, \epz_j}^{\al_j}
\end{equation}
whose composite $\wpr^{\epz_j}_j \winc^{\epz_j}_j$ is the identity function.
\item[\ggcatspstep. Enlarging $\nbe$]  
This step enlarges the object $\nbe \in \GH$ \cref{ggcat_sp_nbe} to an object $\malp \in \GG$, along with some auxiliary objects and morphisms.  Using the abbreviation $\ang{\Cdots}_j = \ang{\Cdots}_{j \in \ufsq}$, we define the object of length $q>0$
\begin{equation}\label{malp_GG}
\malp = \sordi{m}{\al}{j}_{j} \in \GG
\end{equation}
whose $j$-th entry is the pointed finite $G$-set $\ordm_j^{\al_j}$ \cref{ordmj_alphaj}.  We also regard $\malp$ as an object in $\GH$ with each $\ordm_j^{\al_j}$ regarded as a pointed finite $H$-set.  For each $q$-tuple of indices
\begin{equation}\label{epz_tuple}
\epz = \ang{\epz_j}_{j} \witheachspace \epz_j \in \{1,2\},
\end{equation}
we define the object of length $q>0$
\begin{equation}\label{malpepz}
\malpepz = \ang{\ordm_{j,\,\epz_j}^{\al_j}}_{j} \in \GH
\end{equation}
whose $j$-th entry is the pointed finite $H$-set $\ordm_{j,\,\epz_j}^{\al_j}$ \cref{mnr_decomp}.  Using the $H$-equivariant inclusion and projection in \cref{wprinc_onetwo}, there are $H$-equivariant morphisms
\begin{equation}\label{ipepz}
\malpepz \fto[= (1_{\ufsq}, \ang{\winc^{\epz_j}_j}_j)]{\winc^\epz} \malp 
\fto[= (1_{\ufsq}, \ang{\wpr^{\epz_j}_j}_j)]{\wpr^\epz} \malpepz \inspace \GH
\end{equation}
whose composite $\wpr^\epz \winc^\epz$ is the identity morphism.

For example, the $q$-tuple $\eps = (1,1,\ldots,1)$ of 1's yields the object
\begin{equation}\label{malp_ones}
\ang{\ordm^\al_{\eps}} = \ang{\ordm^{\al_j}_{j,1}}_j = \sordi{n}{\be}{j}_j = \nbe \in \GH
\end{equation}
and the $H$-equivariant morphisms
\begin{equation}\label{ip_ones}
\nbe \fto[= (1_{\ufsq}, \ang{\winc^1_j}_j)]{\winc^{\eps}} \malp 
\fto[= (1_{\ufsq}, \ang{\wpr^1_j}_j)]{\wpr^{\eps}} \nbe \inspace \GH
\end{equation}
whose composite is the identity morphism.
\item[\ggcatspstep. Embedding $\segnbe$]  
This step constructs the following commutative diagram of pointed $H$-functors.
\begin{equation}\label{ggcat_sp_LX_dia}
\begin{tikzpicture}[vcenter]
\def\v{-1.4} \def\h{3} \def\c{.6} \def\d{-.6}
\draw[0cell=.9]
(0,0) node (a11) {(\Lh\Xh)\nbe}
(a11)++(\h,0) node (a12) {(\Lg X)\malp}
(a12)++(\h,0) node (a13) {(\Lh\Xh)\nbe}
(a11)++(0,\v) node (a21) {\LXhonenbe}
(a12)++(0,\v) node (a22) {\LXonemalp}
(a13)++(0,\v) node (a23) {\LXhonenbe}
;
\draw[1cell=.8]
(a11) edge node {\winc^\eps_*} (a12)
(a12) edge node {\wpr^\eps_*} (a13)
(a21) edge node {\rinc} (a22)
(a22) edge node {\rpr} (a23)
(a11) edge node[swap] {\segnbe} (a21)
(a12) edge node [swap] {\segmalp} (a22)
(a13) edge node {\segnbe} (a23)
(a11) [rounded corners=2pt] |- ($(a12)+(-1,\c)$) -- node {1} ($(a12)+(1,\c)$) -| (a13)
;
\draw[1cell=.8]
(a21) [rounded corners=2pt] |- ($(a22)+(-1,\d)$) -- node {1} node[swap] {\phantom{x}} ($(a22)+(1,\d)$) -| (a23)
;
\end{tikzpicture}
\end{equation}
\begin{description}
\item[Top categories] 
By \cref{Lg_f_nbe,malp_GG}, the top middle entry in \cref{ggcat_sp_LX_dia} is the pointed $G$-category
\begin{equation}\label{LXmalp}
(\Lg X)\malp = \ecint^{\angorda \in \Gsk} \bigvee_{\GGpunc(\angorda; \malp)} X\angorda
\end{equation}
with $\GGpunc(-;-)$ denoting the set of nonzero morphisms in $\GG$ \cref{GG_morphisms}.  The top left entry, which is equal to the top right entry, is the pointed $H$-category
\begin{equation}\label{LhXhnbe}
(\Lh\Xh)\nbe = \ecint^{\angordb \in \Gsk} \bigvee_{\GHpunc(\angordb; \nbe)} \Xh\angordb
\end{equation}
with $\Xh$ denoting the $H$-restriction of $X$ \cref{Xh} and $\Lh$ denoting the adjoint 2-inverse of $\ihst$ \pcref{not:Lh_ihst}. 
\item[Top functors]
Denote by $\eps = (1,1,\ldots,1)$ the $q$-tuple of 1's and by $\winc^\eps \cn \nbe \to \malp$ the $H$-equivariant morphism in \cref{ip_ones}.  Using \cref{LXmalp,LhXhnbe}, the pointed $H$-functor $\winc^\eps_*$ in \cref{ggcat_sp_LX_dia} is defined on representatives \cref{LX_reps} by the assignment
\begin{equation}\label{wincepz_star}
\begin{split}
& \big(\angordb \fto{\upom} \nbe ; x\big) \in \GH(\angordb; \nbe) \ttimes \Xh\angordb\\
& \mapsto \big(\angordb \fto{\upom} \nbe \fto{\winc^\eps} \malp ; x\big) \in \GG(\angordb; \malp) \ttimes X\angordb.
\end{split}
\end{equation}
This assignment preserves the relations defining the coends \cref{LX_relations}, identity morphisms, and composition.  The group $H$ acts trivially on entries of $\angordb \in \Gsk$ and diagonally \cref{LXnbe_gaction} on representatives of $(\Lh\Xh)\nbe$ and $(\Lg X)\malp$.  Thus, the $H$-equivariance of the functor $\winc^\eps_*$ follows from the $H$-equivariance of the morphism $\winc^\eps$. 

Similarly, with $\wpr^\eps \cn \malp \to \nbe$ denoting the $H$-equivariant morphism in \cref{ip_ones}, the pointed $H$-functor $\wpr^\eps_*$ in \cref{ggcat_sp_LX_dia} is defined on representatives \cref{LX_reps} by the assignment
\begin{equation}\label{wprepz_star}
\begin{split}
& \big(\angorda \fto{\upom} \malp ; x\big) \in \GG(\angorda; \malp) \ttimes X\angorda\\
& \mapsto \big(\angorda \fto{\upom} \malp \fto{\wpr^\eps} \nbe ; x\big) \in \GH(\angorda; \nbe) \ttimes \Xh\angorda.
\end{split}
\end{equation}
The composite $H$-functor $\wpr^\eps_* \winc^\eps_*$ is the identity because the composite morphism $\wpr^\eps \winc^\eps$ is the identity.
\item[Bottom categories]
By the trivial $H$-action on the entries of $\angorda, \angordb, \angordone \in \Gsk$, \cref{GG_Gaction,angordone,Xh,LXmalp,LhXhnbe}, there are equalities of pointed $H$-categories
\begin{equation}\label{LXone_LXhone}
\begin{split}
(\Lg X)\angordone 
&= \txint^{\angorda \in \Gsk} \bigvee_{\GGpunc(\angorda; \angordone)} X\angorda \\
&= \txint^{\angordb \in \Gsk} \bigvee_{\GHpunc(\angordb; \angordone)} \Xh\angordb \\
&= (\Lh\Xh)\angordone.
\end{split}
\end{equation}
By \cref{Xonenbe}, the following two statements hold.
\begin{itemize}
\item The bottom middle entry $\LXonemalp$ in \cref{ggcat_sp_LX_dia} is the $\malp$-twisted product of the pointed $G$-category $(\Lg X)\angordone$ \pcref{def:proCnbe}.
\item The bottom left entry $\LXhonenbe$ in \cref{ggcat_sp_LX_dia}, which is equal to the bottom right entry, is the $\nbe$-twisted product of the pointed $H$-category $(\Lh\Xh)\angordone$.
\end{itemize}

To facilitate the rest of the proof, we construct a pointed $H$-isomorphism
\begin{equation}\label{uptau_Hiso}
\LXonemalp \fto[\iso]{\uptau} \txprod_{\epz} \LXhonemalpepz
\end{equation}
with the product indexed by the lexicographically ordered set of $q$-tuples $\epz = \ang{\epz_j \in \{1,2\}}_j$ \cref{epz_tuple}.  By \cref{ufsn,mnr_decomp}, there are an equality and an isomorphism of unpointed finite $H$-sets
\begin{equation}\label{ufsmalp}
\begin{split}
\prod_{j \in \ufsq} \ufsm_j 
&= \prod_{j \in \ufsq} \big( \ufsm_{j,1} \sqcup \ufsm_{j,2}\big)\\
&\iso \coprod_{\epz = \ang{\epz_j \in \{1,2\}}_{j \in \ufsq}} \Big(\prod_{j \in \ufsq} \ufsm_{j,\,\epz_j} \Big).
\end{split}
\end{equation}
Using \cref{malp_GG,malpepz,LXone_LXhone,ufsmalp}, the pointed $H$-isomorphism $\uptau$ \cref{uptau_Hiso} sends an object or a morphism
\[\ang{y_{\bdi}}_{\bdi \sins \ufs{m_1 m_2 \cdots m_q}} \in \LXonemalp\]
to
\begin{equation}\label{tau_y}
\bang{\ang{y_{\bdi}}_{\bdi \sins \ufs{m_{1,\, \epz_1} m_{2,\, \epz_2} \cdots m_{q,\, \epz_q}}}}_{\epz} 
\in \txprod_{\epz} \LXhonemalpepz.
\end{equation}
In other words, $\uptau$ permutes the coordinates according to the $H$-bijection \cref{ufsmalp}.  The $H$-equivariance of $\uptau$ follows from \cref{proCnbe_gaction,ufsmalp}.
\item[Bottom functors]
Using \crefrange{LXone_LXhone}{tau_y}, the pointed $H$-functor
\begin{equation}\label{rinc_Hfunctor}
\LXhonenbe \fto{\rinc} \LXonemalp
\end{equation}
in \cref{ggcat_sp_LX_dia} is given by
\begin{itemize}
\item the identity for coordinates 
\[\bdi \in \ufs{m_{1,1} m_{2,1} \Cdots m_{q,1}} = \ufs{n_1 n_2 \Cdots n_q}\]
and
\item the basepoint of $(\Lg X)\angordone$ for all other coordinates
\[\bdi \in \ufs{m_{1,\, \epz_1} m_{2,\, \epz_2} \Cdots m_{q,\, \epz_q}} \withspace \epz \neq (1,1,\ldots,1).\]
\end{itemize}
The $H$-equivariance of $\rinc$ follows from \cref{proCnbe_gaction}, \cref{ufsmalp}, and the fact that the basepoint of any pointed $G$-category, including $(\Lg X)\angordone$, is $G$-fixed \pcref{def:ptGcat}, hence also $H$-fixed.

Similarly, the pointed $H$-functor
\begin{equation}\label{rpr_Hfunctor}
\LXonemalp \fto{\rpr} \LXhonenbe
\end{equation}
in \cref{ggcat_sp_LX_dia} projects to the coordinates in $\ufs{m_{1,1} \Cdots m_{q,1}} = \ufs{n_1 \Cdots n_q}$.  By construction, the composite $H$-functor $\rpr \rinc$ is the identity functor.
\item[Vertical functors]
In \cref{ggcat_sp_LX_dia}, $\segnbe$ is the $\nbe$-Segal functor of $\Lh\Xh$ \cref{segnbe_LXh}.  The middle vertical functor $\segmalp$ in \cref{ggcat_sp_LX_dia} is the $\malp$-Segal functor \cref{nbe_segal} of the $\GGG$-category $\Lg X$.  It is a pointed $G$-functor, hence also a pointed $H$-functor.
\item[Right commutativity]
To prove that the right rectangle in the diagram \cref{ggcat_sp_LX_dia} commutes, we consider the $q$-tuple $\eps = (1,1,\ldots,1)$ of 1's and a $q$-tuple of positive integers \cref{bdi}
\[\bdi = \ang{i_j}_j \in \ufs{m_{1,1} m_{2,1} \Cdots m_{q,1}} = \ufs{n_1 n_2 \Cdots n_q}.\]
The diagram in $\GH$
\begin{equation}\label{wprepz_segbdi}
\begin{tikzpicture}[vcenter]
\def\h{1.8} \def\u{.6}
\draw[0cell]
(0,0) node (a1) {\malp}
(a1)++(\h,0) node (a2) {\phantom{\nbe}}
(a2)++(0,.03) node (a2') {\nbe}
(a2)++(\h,0) node (a3) {\angordone}
;
\draw[1cell=.9]
(a1) edge node {\wpr^\eps} (a2)
(a2) edge node {\segbdi} (a3)
(a1) [rounded corners=2pt] |- ($(a2)+(-1,\u)$) -- node {\segbdi} ($(a2)+(1,\u)$) -| (a3)
;
\end{tikzpicture}
\end{equation}
commutes because, for each $j \in \ufsq$, the diagram of pointed finite sets
\begin{equation}\label{wprj_diag}
\begin{tikzpicture}[vcenter]
\def\v{-1.3}
\draw[0cell]
(0,0) node (a11) {\ordm_{j,1}^{\al_j} \wed \ordm_{j,2}^{\al_j}}
(a11)++(2.5,0) node (a12) {\ordn_j^{\be_j}}
(a11)++(0,\v) node (a21) {\ordm_j^{\al_j}}
(a12)++(0,\v) node (a22) {\ordone}
;
\draw[1cell=.9]
(a11) edge node {\wpr_j^1} (a12)
(a12) edge[shorten <=-.5ex] node {\seg_{i_j}} (a22)
(a11) edge[equal, shorten <=-.5ex] (a21)
(a21) edge node {\seg_{i_j}} (a22)
;
\end{tikzpicture}
\end{equation}
commutes by \cref{ij-char,mnr_decomp,wprinc_onetwo,ipepz}.  For each representative pair
\[\big(\angorda \fto{\upom} \malp; x\big) \in \GG(\angorda; \malp) \ttimes X\angorda\]
in $(\Lg X)\malp$, the following equalities in $\LXhonenbe$ prove that the right rectangle in \cref{ggcat_sp_LX_dia} commutes.
\[\begin{aligned}
& \segnbe \wpr^\eps_* \,(\upom; x) &&\\
&= \segnbe (\wpr^\eps \upom ; x) && \text{by \cref{wprepz_star}}\\
&= \ang{(\segbdi \wpr^\eps \upom ; x)}_{\bdi \sins \ufs{m_{1,1} \cdots m_{q,1}}}
&& \text{by \cref{nbe_segal}}\\
&= \ang{(\segbdi \upom ; x)}_{\bdi \sins \ufs{m_{1,1} \cdots m_{q,1}}}
&& \text{by \cref{wprepz_segbdi}}\\
&= \rpr \ang{(\segbdi \upom ; x)}_{\bdi \sins \ufs{m_1 \cdots m_q}} && \text{by \cref{rpr_Hfunctor}}\\
&= \rpr \segmalp (\upom; x) && \text{by \cref{nbe_segal}}
\end{aligned}\]
\item[Left commutativity]
To prove that the left rectangle in the diagram \cref{ggcat_sp_LX_dia} commutes, we consider the $q$-tuple $\eps = (1,1,\ldots,1)$ of 1's, a $q$-tuple of positive integers \cref{bdi}
\[\bdi = \ang{i_j}_j \in \ufs{m_1 m_2 \Cdots m_q},\]
and the following diagram in $\GH$.
\begin{equation}\label{wincepz_segbdi}
\begin{tikzpicture}[vcenter]
\def\h{1.8} \def\u{.6} \def\t{-.6}
\draw[0cell]
(0,0) node (a1) {\phantom{\nbe}}
(a1)++(0,.03) node (a1') {\nbe}
(a1)++(\h,0) node (a2) {\malp}
(a2)++(\h,0) node (a3) {\angordone}
;
\draw[1cell=.9]
(a1) edge node {\winc^\eps} (a2)
(a2) edge node {\segbdi} (a3)
(a1) [rounded corners=2pt] |- ($(a2)+(-1,\u)$) -- node {\segbdi} ($(a2)+(1,\u)$) -| (a3)
;
\draw[1cell=.9]
(a1) [rounded corners=2pt] |- ($(a2)+(-1,\t)$) -- node {0} node[swap] {\phantom{x}} ($(a2)+(1,\t)$) -| (a3)
;
\end{tikzpicture}
\end{equation}
By \cref{ij-char,mnr_decomp,wprinc_onetwo,ipepz}, the following two statements hold regarding the diagram \cref{wincepz_segbdi}.
\begin{description}
\item[Top] If
\[\bdi \in \ufs{m_{1,1} m_{2,1} \Cdots m_{q,1}} = \ufs{n_1 n_2 \Cdots n_q},\]
then the top half of the diagram \cref{wincepz_segbdi} commutes.
\item[Bottom] If
\[\bdi \not\in \ufs{m_{1,1} m_{2,1} \Cdots m_{q,1}} = \ufs{n_1 n_2 \Cdots n_q},\]
then the bottom half of the diagram \cref{wincepz_segbdi} commutes, where 0 denotes the 0-morphism in $\GH$ \cref{GG_morphisms}.
\end{description}
For each representative pair
\[\big(\angordb \fto{\upom} \nbe; x\big) \in \GH(\angordb; \nbe) \ttimes \Xh\angordb\]
in $(\Lh\Xh)\nbe$, the following equalities in $\LXonemalp$ prove that the left rectangle in \cref{ggcat_sp_LX_dia} commutes.
\[\begin{aligned}
& \segmalp \winc^\eps_* (\upom; x) &&\\
&= \segmalp (\winc^\eps \upom ; x) && \text{by \cref{wincepz_star}}\\
&= \ang{(\segbdi \winc^\eps \upom ; x)}_{\bdi \in \ufs{m_1 \cdots m_q}} && \text{by \cref{nbe_segal}}\\
&= \rinc \ang{(\segbdi \upom ; x)}_{\bdi \in \ufs{m_{1,1} \cdots m_{q,1}}} && \text{by \cref{rinc_Hfunctor,wincepz_segbdi}}\\
&= \rinc \segnbe (\upom; x) && \text{by \cref{nbe_segal}}
\end{aligned}\]
\end{description}
This finishes the construction of the commutative diagram \cref{ggcat_sp_LX_dia} of pointed $H$-functors.
\item[\ggcatspstep]  
By the assumption that $\Lg X$ is a special $\GGG$-category, its $\malp$-Segal functor $\segmalp$ is a categorical weak $G$-equivalence, hence also a categorical weak $H$-equivalence.  \cref{weq_retract}, applied to the group $H$ and the diagram \cref{ggcat_sp_LX_dia}, implies that the $\nbe$-Segal functor $\segnbe$ for $\Lh\Xh$ is also a categorical weak $H$-equivalence.  
\end{description}
This proves that the $\GHH$-category $\Lh\Xh$ is special and completes the proof of the \emph{if} assertion of the \namecref{ggcat_sp_LX}.
\end{proof}

\cref{ggcat_sp_LX} shows that the left 2-adjoint $\Lg$ \cref{Lg} preserves and reflects specialness.  The next result shows that its 2-adjoint inverse $\igst$ \pcref{thm:ggcat_ggcatg_iieq} also preserves and reflects specialness.

\begin{theorem}\label{gggcat_sp_iX}
For a finite group $G$, a $\GGG$-category
\[\GG \fto{X} \Catgst\]
is special if and only if the $\Gskg$-category
\[\Gsk \fto{\igst X} \Gcatst\]
is special.
\end{theorem}

\begin{proof}
The counit $\vg$ of the adjoint 2-equivalence $(\Lg,\igst)$ has component pointed $G$-isomorphisms \cref{vgf_nbe}
\[(\Lg\igst X)\nbe \fto[\iso]{\vg_{X,\nbe}} X\nbe \forspace \nbe \in \GG,\]
where $G$-equivariance is proved in \cref{vg_gnatural}.  For each object $\nbe \in \GG \setminus \{\vstar,\ang{}\}$, the naturality of $\vg_X$ with respect to the $\bdi$-characteristic morphisms $\segbdi$ \cref{i-char} implies that the diagram of pointed $G$-functors
\begin{equation}\label{vgxnbe_diag}
\begin{tikzpicture}[vcenter]
\def\v{-1.4}
\draw[0cell]
(0,0) node (a11) {(\Lg\igst X)\nbe}
(a11)++(3.75,0) node (a12) {X\nbe}
(a11)++(0,\v) node (a21) {\LiXonenbe}
(a12)++(0,\v) node (a22) {\Xonenbe}
;
\draw[1cell=.85]
(a11) edge node {\vg_{X,\nbe}} node[swap] {\iso} (a12)
(a12) edge[transform canvas={xshift=-.8em}] node {\segnbex} (a22)
(a11) edge[transform canvas={xshift=1em}] node[swap] {\segnbelix} (a21)
(a21) edge node {\pro{(\vg_{X,\angordone})}{\nbe}} node[swap] {\iso} (a22)
;
\end{tikzpicture}
\end{equation}
commutes.  Thus, the $\nbe$-Segal functor $\segnbex$ for $X$ is a categorical weak $G$-equivalence if and only if the $\nbe$-Segal functor $\segnbelix$ for $\Lg\igst X$ is a categorical weak $G$-equivalence.  By \cref{def:sp_gggcat}, this means that $X$ is special if and only if $\Lg\igst X$ is special.  By \cref{ggcat_sp_LX}, $\Lg\igst X$ is special if and only if $\igst X$ is special.
\end{proof}

\begin{explanation}[Finiteness of $G$]\label{rk:shi_cor}
\cref{ggcat_sp_LX,gggcat_sp_iX} are analogous to \cite[p.\ 227, Corollary]{shimakawa91}, which concerns $G$-spaces and the indexing categories $\Fsk$ and $\FG$ \pcref{def:Fsk,def:FG}.  We emphasize that \cref{gggcat_sp_iX} requires the group $G$ to be finite because its proof uses \cref{ggcat_sp_LX}, which requires $G$ to be finite.  Specifically, the finiteness of $G$ is used in \cref{ordmj_alphaj} to construct the pointed finite $G$-set $\ordm_j^{\al_j}$ from the pointed finite $H$-set $\ordn_j^{\be_j}$.  If $G$ is not finite, then there is no guarantee that $\ordm_j^{\al_j}$ is finite and that $\malp$ \cref{malp_GG} is an object in $\GG$.
\end{explanation}

\section{Weak $G$-Equivalences of $\GGG$-Categories}
\label{sec:gggcat_weq}

This section introduces weak $G$-equivalences between $\GGG$-categories for an arbitrary group $G$.  In this work, the most important weak $G$-equivalences are the strong $H$-theory comparison $G$-natural transformations in \cref{thm:pistweakgeq}.

\secoutline
\begin{itemize}
\item \cref{def:ggg_weq} defines weak $G$-equivalences between $\GGG$-categories.
\item \cref{ex:ggcat_weq} observes that the strong $H$-theory comparison $\Pistsg_{\Ah}$ in \cref{thm:pistweakgeq} is a weak $G$-equivalence.
\item \cref{f_nbe} proves that twisted products preserve categorical weak $G$-equivalences.
\item \cref{ggg_weq_sp} proves that the domain of a weak $G$-equivalence is special if and only if the codomain is special.
\item \cref{ggg_sp_weq} proves that, between special $\GGG$-categories, weak $G$-equivalences are detected at objects of the form $(\ordone,\ldots,\ordone) \in \GG$.
\end{itemize}

Recall from \cref{def:ggcatg} that a 1-cell in $\GGCatg$ is a $G$-natural transformation.  Its components are pointed $G$-functors \cref{ggcatg_icell_geq}.

\begin{definition}[Weak $G$-Equivalences]\label{def:ggg_weq}
A 1-cell
\begin{equation}\label{weakGeq_iicell}
\begin{tikzpicture}[vcenter]
\def\t{28}
\draw[0cell]
(0,0) node (a1) {\phantom{\Gsk}}
(a1)++(1.8,0) node (a2) {\phantom{\Gsk}}
(a1)++(-.08,0) node (a1') {\GG}
(a2)++(.2,0) node (a2') {\Catgst}
;
\draw[1cell=.9]
(a1) edge[bend left=\t] node {X} (a2)
(a1) edge[bend right=\t] node[swap] {X'} (a2)
;
\draw[2cell]
node[between=a1 and a2 at .45, rotate=-90, 2label={above,\tha}] {\Rightarrow}
;
\end{tikzpicture}
\end{equation}
in $\GGCatg$ is called a \emph{weak $G$-equivalence}\index{weak G-equivalence@weak $G$-equivalence!GGG-category@$\GGG$-category}\index{GGG-category@$\GGG$-category!weak G-equivalence@weak $G$-equivalence} if, for each object $\malp \in \GG \setminus \{\vstar\}$, its $\malp$-component pointed $G$-functor
\[X\malp \fto{\tha_{\malp}} X'\malp\]
is a categorical weak $G$-equivalence \pcref{def:cat_weakg}.
\end{definition}

\begin{explanation}\label{expl:ggcatg_weq}
Consider \cref{def:ggg_weq}.
\begin{enumerate}
\item\label{expl:ggcatg_weq_i} 
A 1-cell $\tha \cn X \to X'$ in $\GGCatg$ is a weak $G$-equivalence if and only if, for each object $\malp \in \GG \setminus \{\vstar\}$ and subgroup $H \subseteq G$, the $H$-fixed subfunctor
\[X\malp^H \fto{\tha_{\malp}^H} X'\malp^H\]
is an equivalence of categories.  Weak $G$-equivalences are analogous to \emph{level $G$-equivalences} between $\FGG$-spaces \cite[Def.\ 2.34]{mmo}, which concern $G$-spaces and the indexing $G$-category $\FG$ \pcref{def:FG}.
\item\label{expl:ggcatg_weq_ii} 
There is no requirement regarding the basepoint $\vstar \in \GG$ because the $\vstar$-component of each 1-cell $\tha$ in $\GGCatg$ is the identity functor on the terminal $G$-category $\bone$.  On the other hand, the empty tuple $\ang{} \in \GG$ is included in \cref{def:ggg_weq} because, for a general 1-cell $\tha$, there is no guarantee that the $\ang{}$-component
\[X\ang{} \fto{\tha_{\ang{}}} X'\ang{}\]
is a categorical weak $G$-equivalence.  This is a subtle conceptual difference with special $\GGG$-categories \pcref{expl:sp_gggcat_mar}.\defmark
\end{enumerate}
\end{explanation}

\begin{example}\label{ex:ggcat_weq}
By \cref{thm:pistweakgeq}, for a $\Uinf$-operad $\Op$ \pcref{as:OpA'}, an $\Op$-pseudoalgebra $\A$ \pcref{def:pseudoalgebra}, and the $\Oph$-pseudoalgebra $\Ah = \Catg(\EG,\A)$ with $\Oph = \Catg(\EG,\Op)$ \pcref{catgego}, the strong $H$-theory comparison $G$-natural transformation \cref{PistsgA}
\begin{equation}\label{PistsgAh_ex_diag}
\begin{tikzpicture} [vcenter]
\def\s{25}
\draw[0cell]
(0,0) node (a1) {\GG}
(a1)++(2.5,0) node (a2) {\phantom{\GG}}
(a2)++(.15,0) node (a2') {\Catgst}
;
\draw[1cell=.9]
(a1) edge[bend left=\s] node {\smast \Sgohsg\Ah} (a2)
(a1) edge[bend right=\s] node[swap] {\Hgohsg\Ah} (a2) 
;
\draw[2cell]
node[between=a1 and a2 at .4, rotate=-90, 2label={above,\Pistsg_{\Ah}}] {\Rightarrow}
;
\end{tikzpicture}
\end{equation}
is a weak $G$-equivalence in $\GGCatg$.
\end{example}

The rest of this section discusses properties of weak $G$-equivalences and how they relate to special $\GGG$-categories.  \cref{f_nbe} is a key observation needed for the rest of this section.  It states that categorical weak $G$-equivalences \pcref{def:cat_weakg} are closed under twisted products \pcref{def:proCnbe}.

\begin{lemma}\label{f_nbe}
Suppose $\fun \cn \C \to \D$ is a pointed categorical weak $G$-equivalence between pointed $G$-categories, and $\nbe \in \GG \setminus \{\vstar,\ang{}\}$ is an object of length $q>0$.  Then the induced pointed $G$-functor between the $\nbe$-twisted products
\[\Cnbe \fto{\funnbe} \Dnbe\]
is a categorical weak $G$-equivalence, where each entry of $\funnbe$ is given by $\fun$.
\end{lemma}

\begin{proof}
Given a subgroup $H \subseteq G$, we verify that the $H$-fixed subfunctor 
\[(\Cnbe)^H \fto{(\funnbe)^H} (\Dnbe)^H\]
is an equivalence of categories.  The object $\nbe$ is uniquely determined by
\begin{itemize}
\item an object $\angordn \in \Gsk \setminus \{\vstar,\ang{}\}$ and 
\item a homomorphism $\be \cn G \to \Siangordn$ \pcref{expl:gph_subgrp}.
\end{itemize}  
Denoting by $\be' \cn H \to \Siangordn$ the restriction of $\be$ to $H$, we consider the object $\nbep \in \GH \setminus \{\vstar,\ang{}\}$.  There is an equality of categories
\[(\Cnbe)^H = (\Cnbep)^H\]
because, by \cref{proCnbe_gaction}, the $H$-action on $\Cnbe$ involves only the $H$-actions on $\C$ and on each entry $\ordi{n}{\be}{j}$ in $\nbe$.  In the $\nbep$-twisted product $\Cnbep$, $\C$ is regarded as a pointed $H$-category.  The same discussion also yields the equality of $H$-fixed subfunctors
\[(\funnbe)^H = (\funnbep)^H,\]
where $\fun$ in $\funnbep$ is regarded as a pointed $H$-functor.  \cref{twprod_fixed}, applied to
\begin{itemize}
\item the group $H$,
\item the object $\nbep \in \GH \setminus \{\vstar,\ang{}\}$, and
\item the pointed $H$-functor $\fun \cn \C \to \D$,
\end{itemize}
states that the $H$-fixed subfunctor $(\funnbep)^H$ splits into a product $\txprod_{t \in \ufsr}\, \fun^{H_t}$ for some $r>0$ and some subgroups $H_t \subseteq H$.  Since each $H_t$ is also a subgroup of $G$, the $H_t$-fixed subfunctor $\fun^{H_t}$ is an equivalence of categories by the assumption that the $G$-functor $\fun$ is a categorical weak $G$-equivalence.  Thus, the product $\txprod_{t \in \ufsr}\, \fun^{H_t}$ and also $(\funnbep)^H = (\funnbe)^H$ are equivalences of categories.
\end{proof}

\subsection*{Weak $G$-Equivalences and Specialness}

\cref{ggg_weq_sp} shows that weak $G$-equivalences \pcref{def:ggg_weq} preserve and reflect special $\GGG$-categories \pcref{def:sp_gggcat}.

\begin{lemma}\label{ggg_weq_sp}
Suppose $\tha \cn X \to X'$ is a weak $G$-equivalence in $\GGCatg$.  Then the following two statements hold.
\begin{enumerate}
\item\label{ggg_weq_sp_i} For each object $\nbe \in \GG \setminus \{\vstar,\ang{}\}$, the $\nbe$-Segal functor for $X$ \cref{nbe_segal} 
\[X\nbe \fto{\segnbex} \Xonenbe\]
is a categorical weak $G$-equivalence if and only if the $\nbe$-Segal functor for $X'$ 
\[\Xpnbe \fto{\segnbexp} \Xponenbe\]
is a categorical weak $G$-equivalence.
\item\label{ggg_weq_sp_ii} $X$ is special if and only if $X'$ is special.
\end{enumerate}
\end{lemma}

\begin{proof}
By \cref{def:sp_gggcat}, \eqref{ggg_weq_sp_ii} follows from \eqref{ggg_weq_sp_i}.  To prove \eqref{ggg_weq_sp_i}, observe that the naturality of $\tha$ for the $\bdi$-characteristic morphisms $\segbdi \cn \nbe \to \angordone$ \cref{i-char} implies that the diagram of pointed $G$-functors
\begin{equation}\label{segnbe_thanbe}
\begin{tikzpicture}[vcenter]
\def\v{-1.4}
\draw[0cell]
(0,0) node (a11) {X\nbe}
(a11)++(3,0) node (a12) {\Xonenbe}
(a11)++(0,\v) node (a21) {\Xpnbe}
(a12)++(0,\v) node (a22) {\Xponenbe}
;
\draw[1cell=.85]
(a11) edge node {\segnbex} (a12)
(a12) edge[transform canvas={xshift=-1ex}] node {\tha_{\angordone}^{\nbe}} (a22)
(a11) edge node[swap] {\tha_{\nbe}} (a21)
(a21) edge node {\segnbexp} (a22)
;
\end{tikzpicture}
\end{equation}
commutes, where $\tha_{\angordone}^{\nbe}$ is given entrywise by $\tha_{\angordone}$.  The components $\tha_{\nbe}$ and $\tha_{\angordone}$ of $\tha$ are categorical weak $G$-equivalences by assumption.  By \cref{f_nbe}, $\tha_{\angordone}^{\nbe}$ is a categorical weak $G$-equivalence.  Thus, $\segnbex$ is a categorical weak $G$-equivalence if and only if $\segnbexp$ is so.
\end{proof}

\cref{ggg_sp_weq} shows that, between special $\GGG$-categories \pcref{def:sp_gggcat}, weak $G$-equivalences \pcref{def:ggg_weq} are detected at the objects $\ang{\ordone} = (\ordone,\ldots,\ordone) \in \GG$ consisting of $q \geq 0$ copies of $\ordone$.

\begin{lemma}\label{ggg_sp_weq}
Suppose $\tha \cn X \to X'$ is a 1-cell in $\GGCatg$ such that $X$ and $X'$ are special $\GGG$-categories.  Then the following two statements are equivalent.
\begin{enumerate}
\item\label{ggg_sp_weq_i} $\tha$ is a weak $G$-equivalence.
\item\label{ggg_sp_weq_ii} $\tha_{\angordone}$ is a categorical weak $G$-equivalence for each object $\ang{\ordone} = (\ordone,\ldots,\ordone) \in \GG$ of length $q \geq 0$.
\end{enumerate}
\end{lemma}

\begin{proof}
The implication $(1) \Rightarrow (2)$ holds by \cref{def:ggg_weq}, regardless of whether $X$ and $X'$ are special or not.

To prove $(2) \Rightarrow (1)$, note that $\angordone_{j \in \ufs{0}}$ is the empty tuple $\ang{}$.  To show that $\tha$ is a weak $G$-equivalence, it suffices to prove that $\tha_{\nbe}$ is a categorical weak $G$-equivalence for each object $\nbe \in \GG \setminus \{\vstar,\ang{}\}$ of length $q>0$.  Reusing the commutative diagram \cref{segnbe_thanbe}, the $\nbe$-Segal functors $\segnbex$ and $\segnbexp$ are categorical weak $G$-equivalences by the assumption that $X$ and $X'$ are special.  By the assumption on $\tha_{\angordone}$ and \cref{f_nbe}, $\tha_{\angordone}^{\nbe}$ is a categorical weak $G$-equivalence.  It follows from the commutative diagram \cref{segnbe_thanbe} that $\tha_{\nbe}$ is also a categorical weak $G$-equivalence.
\end{proof}

\section{Weak $G$-Equivalences of $\Gskg$-Categories}
\label{sec:ggcat_weq}

This section introduces weak $G$-equivalences between $\Gskg$-categories for an arbitrary group $G$.

\secoutline
\begin{itemize}
\item \cref{def:gg_weq} defines weak $G$-equivalences between $\Gskg$-categories.
\item \cref{gg_weq_char} describes weak $G$-equivalences in terms of subgroups of $G$
\item \cref{gg_weq_sp} proves that the domain of a weak $G$-equivalence is special if and only if the codomain is special.
\item \cref{gg_sp_weq} proves that, between special $\Gskg$-categories, weak $G$-equivalences are detected at objects of the form $(\ordone,\ldots,\ordone) \in \Gsk$.
\end{itemize}

For a $\Gskg$-category $X$ \cref{ggcat_obj} and an object $\angordn \in \Gsk \setminus \{\vstar\}$ \cref{Gsk_objects}, recall that $X\angordn$ is a pointed $\GSinp$-category \cref{Xn_gsin_act}.  For the empty tuple $\ang{}$, $\Siang$ is the trivial group $\{e\}$, and $X\ang{}$ is the original pointed $G$-category.  Also recall that a 1-cell in $\GGCatii$ is a natural transformation \pcref{def:GGCat}.  

\begin{definition}[Weak $G$-Equivalences]\label{def:gg_weq}
A 1-cell
\begin{equation}\label{GGCatii_weq_iicell}
\begin{tikzpicture}[vcenter]
\def\t{28}
\draw[0cell]
(0,0) node (a1) {\Gsk}
(a1)++(1.8,0) node (a2) {\phantom{\Gskel}}
(a2)++(.3,0) node (a2') {\Gcatst}
;
\draw[1cell=.9]
(a1) edge[bend left=\t] node {X} (a2)
(a1) edge[bend right=\t] node[swap] {X'} (a2)
;
\draw[2cell]
node[between=a1 and a2 at .45, rotate=-90, 2label={above,\theta}] {\Rightarrow}
;
\end{tikzpicture}
\end{equation}
in $\GGCatii$ is called a \emph{weak $G$-equivalence}\index{weak G-equivalence@weak $G$-equivalence!G-G-category@$\Gskg$-category}\index{G-G-category@$\Gskg$-category!weak G-equivalence@weak $G$-equivalence} if, for each object $\angordn \in \Gsk \setminus \{\vstar\}$ and each graph subgroup $K \in \FGSin$ \pcref{def:graph_subgrp}, the $\angordn$-component pointed $\GSinp$-functor
\[X\angordn \fto{\tha_{\angordn}} X'\angordn\]
is a categorical weak $K$-equivalence \pcref{def:cat_weakg}.
\end{definition}

\begin{explanation}\label{expl:ggcat_weq}
Consider \cref{def:gg_weq}.
\begin{enumerate}
\item\label{expl:ggcat_weq_i} 
For a 1-cell $\tha$ in $\GGCatii$ and an object $\angordn \in \Gsk \setminus \{\vstar\}$, the $\angordn$-component $\tha_{\angordn}$ is a pointed $\GSinp$-functor by $G$-equivariance \cref{ggcat_mor_component}, naturality \cref{ggcat_mor_naturality}, and \cref{Xn_gsin_act}.
\item\label{expl:ggcat_weq_ii} 
Graph subgroups are closed under passage to subgroups.  Thus, a 1-cell $\tha \cn X \to X'$ in $\GGCatii$ is a weak $G$-equivalence if and only if, for each object $\angordn \in \Gsk \setminus \{\vstar\}$ and each graph subgroup $K \in \FGSin$, the $K$-fixed subfunctor
\[X\angordn^K \fto{\tha_{\angordn}^K} X'\angordn^K\]
is an equivalence of categories.   Weak $G$-equivalences in \cref{def:gg_weq} are analogous to \emph{$\mathbbm{F}_{\!\bullet}$-level equivalences} between $\Fskg$-spaces \cite[Def.\ 2.14 (iii)]{mmo}, which concern $G$-spaces and the indexing category $\Fsk$ \pcref{def:Fsk}.
\item\label{expl:ggcat_weq_iii} 
There is no requirement regarding the basepoint $\vstar \in \GG$ because the $\vstar$-component of each 1-cell in $\GGCatii$ is the identity functor on the terminal $G$-category $\bone$.  On the other hand, the empty tuple $\ang{} \in \Gsk$ is included in \cref{def:gg_weq} because, for a general 1-cell $\tha$, there is no guarantee that the $\ang{}$-component
\[X\ang{} \fto{\tha_{\ang{}}} X'\ang{}\]
is a categorical weak $G$-equivalence.  This is a subtle conceptual difference with special $\Gskg$-categories \pcref{expl:sp_ggcat}.\defmark
\end{enumerate}
\end{explanation}

By \cref{def:gg_weq}, weak $G$-equivalences in $\GGCatii$ are defined in terms of graph subgroups of $\GSin$.  \cref{gg_weq_char} characterizes weak $G$-equivalences in $\GGCatii$ in terms of subgroups of $G$.  For a $\Gskg$-category $X$, a subgroup $H \subseteq G$, and a homomorphism $\be \cn H \to \Siangordn$, recall from \cref{Xnbeta_haction} the pointed $H$-category $\Xnbeta$ obtained from $X\angordn$ by twisting the $H$-action by $\beta$.

\begin{lemma}\label{gg_weq_char}
A 1-cell $\tha \cn X \to X'$ in $\GGCatii$ is a weak $G$-equivalence if and only if, for each object $\angordn \in \Gsk \setminus \{\vstar\}$, subgroup $H \subseteq G$, and homomorphism $\be \cn H \to \Siangordn$, the $\angordn$-component
\begin{equation}\label{than_xxp}
\Xnbeta \fto{\tha_{\angordn}} \Xpnbeta
\end{equation}
is a categorical weak $H$-equivalence.
\end{lemma}

\begin{proof}
The $\angordn$-component $\tha_{\angordn}$ in \cref{than_xxp} is $H$-equivariant by \cref{Xnbeta_haction}, the $G$-equivariance of 
\begin{equation}\label{than_xnxpn}
X\angordn \fto{\tha_{\angordn}} X'\angordn
\end{equation}
with respect to the original $G$-actions on $X\angordn$ and $X'\angordn$ \cref{ggcat_mor_component}, and naturality \cref{ggcat_mor_naturality}.  By \cref{gph_subgrp}, a graph subgroup $K \in \FGSin$ is uniquely determined by a subgroup $H \subseteq G$ and a homomorphism $\be \cn H \to \Siangordn$ such that 
\[K = \{(h, \be(h)) \tmid h \in H \}.\]
By \cref{Xn_gsin_act,Xnbeta_haction}, the $K$-equivariance of $\tha_{\angordn}$ \cref{than_xnxpn} is precisely the $H$-equivariance of $\tha_{\angordn}$ \cref{than_xxp}.  Moreover, subgroups of $K$ correspond to subgroups of $H$ together with the corresponding restriction of $\be$.  Thus, $\tha_{\angordn}$ \cref{than_xnxpn} is a categorical weak $K$-equivariance if and only if $\tha_{\angordn}$ \cref{than_xxp} is a categorical weak $H$-equivalence.
\end{proof}

\subsection*{Weak $G$-Equivalences and Specialness}

\cref{gg_weq_sp} shows that weak $G$-equivalences \pcref{def:gg_weq} preserve and reflect special $\Gskg$-categories \pcref{def:sp_ggcat}.

\begin{lemma}\label{gg_weq_sp}
Suppose $\tha \cn X \to X'$ is a weak $G$-equivalence in $\GGCatii$.  Then the following two statements hold.
\begin{enumerate}
\item\label{gg_weq_sp_i} For each object $\angordn \in \Gsk \setminus \{\vstar,\ang{}\}$ and each graph subgroup $K \in \FGSin$, the $\angordn$-Segal functor for $X$ \cref{nsegal} 
\[X\angordn \fto{\segnx} \Xonenoneq\]
is a categorical weak $K$-equivalence if and only if the $\angordn$-Segal functor for $X'$ 
\[X'\angordn \fto{\segnxp} \Xponenoneq\]
is a categorical weak $K$-equivalence.
\item\label{gg_weq_sp_ii} $X$ is special if and only if $X'$ is special.
\end{enumerate}
\end{lemma}

\begin{proof}
By \cref{def:sp_ggcat}, \eqref{gg_weq_sp_ii} follows from \eqref{gg_weq_sp_i}.  To prove \eqref{gg_weq_sp_i}, observe that the naturality of $\tha$ for the $\bdi$-characteristic morphisms $\segbdi \cn \angordn \to \angordone$ \cref{i_char} implies that the diagram of pointed $\GSinp$-functors
\begin{equation}\label{segn_thanbe}
\begin{tikzpicture}[vcenter]
\def\v{-1.4}
\draw[0cell]
(0,0) node (a11) {X\angordn}
(a11)++(3,0) node (a12) {\Xonenoneq}
(a11)++(0,\v) node (a21) {X'\angordn}
(a12)++(0,\v) node (a22) {\Xponenoneq}
;
\draw[1cell=.85]
(a11) edge node {\segnx} (a12)
(a12) edge[transform canvas={xshift=-1ex}] node {\tha_{\angordone}^{n_1 \cdots n_q}} (a22)
(a11) edge node[swap] {\tha_{\angordn}} (a21)
(a21) edge node {\segnxp} (a22)
;
\end{tikzpicture}
\end{equation}
commutes, where $\tha_{\angordone}^{n_1 \cdots n_q}$ is given entrywise by $\tha_{\angordone}$.  For each graph subgroup $K \in \FGSin$, the component $\tha_{\angordn}$ is a categorical weak $K$-equivalence by the assumption on $\tha$.  It remains to prove that $\tha_{\angordone}^{n_1 \cdots n_q}$ is a categorical weak $K$-equivalence.

By \cref{gph_subgrp}, the graph subgroup $K$ has the form 
\[K = \{(h, \be(h)) \tmid h \in H \}\]
for a unique subgroup $H \subseteq G$ and a unique homomorphism $\be \cn H \to \Siangordn$.  Subgroups of $K$ correspond to subgroups of $H$ together with the corresponding restriction of $\be$. By \cref{proCnbe_gaction,xonen_gsin_act}, $\tha_{\angordone}^{n_1 \cdots n_q}$ is a categorical weak $K$-equivalence if and only if
\[\Xonenbe \fto{\tha_{\angordone}^{\nbe}} \Xponenbe\]
is a categorical weak $H$-equivalence, in which $\tha_{\angordone} \cn X\angordone \to X'\angordone$ is regarded as an $H$-functor.  By the assumption on $\tha$, $\tha_{\angordone}$ is a categorical weak $K'$-equivalence for each graph subgroup $K' \subseteq G \ttimes \Siangordone$.  By \cref{ex:FGSione}, graph subgroups of $G \ttimes \Siangordone$ are precisely subgroups of $G$.  Thus, $\tha_{\angordone}$ is a categorical weak $H$-equivalence.  By \cref{f_nbe}, applied to $\tha_{\angordone}$ and the group $H$, $\tha_{\angordone}^{\nbe}$ is a categorical weak $H$-equivalence.
\end{proof}

\cref{gg_sp_weq} shows that, between special $\Gskg$-categories \pcref{def:sp_ggcat}, weak $G$-equivalences \pcref{def:gg_weq} are detected at the objects $\ang{\ordone} = (\ordone,\ldots,\ordone) \in \Gsk$ consisting of $q \geq 0$ copies of $\ordone$.  Recall that graph subgroups of $G \ttimes \Siangordone$ are precisely subgroups of $G$ \pcref{ex:FGSione}.

\begin{lemma}\label{gg_sp_weq}
Suppose $\tha \cn X \to X'$ is a 1-cell in $\GGCatii$ such that $X$ and $X'$ are special $\Gskg$-categories.  Then the following two statements are equivalent.
\begin{enumerate}
\item\label{gg_sp_weq_i} $\tha$ is a weak $G$-equivalence.
\item\label{gg_sp_weq_ii} $\tha_{\angordone}$ is a categorical weak $G$-equivalence for each object $\ang{\ordone} = (\ordone,\ldots,\ordone) \in \Gsk$ of length $q \geq 0$.
\end{enumerate}
\end{lemma}

\begin{proof}
The implication $(1) \Rightarrow (2)$ holds by \cref{ex:FGSione,def:gg_weq}, regardless of whether $X$ and $X'$ are special or not.

To prove $(2) \Rightarrow (1)$, note that $\angordone_{j \in \ufs{0}}$ is the empty tuple $\ang{}$.  To show that $\tha$ is a weak $G$-equivalence, it suffices to prove that $\tha_{\angordn}$ is a categorical weak $K$-equivalence for each object $\angordn \in \Gsk \setminus \{\vstar,\ang{}\}$ of length $q>0$ and each graph subgroup $K \in \FGSin$.  Reusing the commutative diagram \cref{segn_thanbe}, the $\angordn$-Segal functors $\segnx$ and $\segnxp$ are categorical weak $K$-equivalences by the assumption that $X$ and $X'$ are special.  By assumption, $\tha_{\angordone}$ is a categorical weak $H$-equivalence for each subgroup $H \subseteq G$.  The last paragraph of the proof of \cref{gg_weq_sp}---which uses \cref{proCnbe_gaction}, \cref{xonen_gsin_act}, and \cref{gph_subgrp,f_nbe}---proves that $\tha_{\angordone}^{n_1 \cdots n_q}$ is a categorical weak $K$-equivalence.  It follows from the commutative diagram \cref{segn_thanbe} that $\tha_{\angordn}$ is also a categorical weak $K$-equivalence.
\end{proof}

\section{Comparing Weak $G$-Equivalences}
\label{sec:comparing_weq}

For a \emph{finite} group $G$, this section proves that weak $G$-equivalences in $\GGCatg$ and $\GGCatii$ \pcref{def:ggg_weq,def:gg_weq} correspond to each other under the adjoint 2-equivalence \pcref{thm:ggcat_ggcatg_iieq}
\begin{equation}\label{Li_secti}
\begin{tikzpicture}[vcenter]
\draw[0cell]
(0,0) node (a1) {\GGCatii}
(a1)++(2.5,0) node (a2) {\GGCatg}
;
\draw[1cell=.9]
(a1) edge[transform canvas={yshift=.5ex}] node {\Lg} (a2)
(a2) edge[transform canvas={yshift=-.4ex}] node {\igst} (a1)
;
\end{tikzpicture}
\end{equation}
between the 2-categories in \cref{def:GGCat,def:ggcatg}.

\secoutline
\begin{itemize}
\item For an arbitrary group $G$, \cref{gg_ggg_weq} describes weak $G$-equivalences in $\GGCatii$ in terms of weak $H$-equivalences in $\GHCath$ for subgroups $H \subseteq G$.  This preliminary result is used in the proof of \cref{Lg_weq}.
\item \cref{Lg_weq} proves that $\Lg$ preserves and reflects weak $G$-equivalences.
\item \cref{igst_weq} proves that $\igst$ preserves and reflects weak $G$-equivalences.
\end{itemize}

For a subgroup $H \subseteq G$, recall the $H$-restriction $\Xh \in \GHCatii$ of a $\Gskg$-category $X$ \pcref{def:ghcat}.

\begin{definition}[Restrictions]\label{def:tha_h}
For a 1-cell $\tha \cn X \to X'$ in $\GGCatii$ \pcref{def:GGCat} and a subgroup $H$ of a group $G$, the 1-cell
\begin{equation}\label{thah}
\begin{tikzpicture}[vcenter]
\def\t{30}
\draw[0cell]
(0,0) node (a1) {\Gsk}
(a1)++(1.8,0) node (a2) {\phantom{\Gskel}}
(a2)++(.3,0) node (a2') {\Hcatst}
;
\draw[1cell=.9]
(a1) edge[bend left=\t] node {\Xh} (a2)
(a1) edge[bend right=\t] node[swap] {\Xh'} (a2)
;
\draw[2cell]
node[between=a1 and a2 at .4, rotate=-90, 2label={above,\thah}] {\Rightarrow}
;
\end{tikzpicture}
\end{equation}
in $\GHCatii$, called the \emph{$H$-restriction}\index{restriction} of $\tha$, is defined by regarding each pointed $G$-functor $\tha_{\angordn}$ for $\angordn \in \Gsk$ as a pointed $H$-functor.
\end{definition}

Using the adjoint 2-equivalence \pcref{not:Lh_ihst}
\begin{equation}\label{Lhih_diag}
\begin{tikzpicture}[vcenter]
\draw[0cell]
(0,0) node (a1) {\GHCatii}
(a1)++(2.5,0) node (a2) {\GHCath}
;
\draw[1cell=.9]
(a1) edge[transform canvas={yshift=.5ex}] node {\Lh} (a2)
(a2) edge[transform canvas={yshift=-.4ex}] node {\ihst} (a1)
;
\end{tikzpicture}
\end{equation}
in \cref{thm:ggcat_ggcatg_iieq} applied to the group $H$, \cref{gg_ggg_weq} compares weak $G$-equivalences in \cref{def:ggg_weq,def:gg_weq}.

\begin{lemma}\label{gg_ggg_weq}
For an arbitrary group $G$, a 1-cell 
\[X \fto{\tha} X' \inspace \GGCatii\]
is a weak $G$-equivalence if and only if, for each subgroup $H \subseteq G$, the 1-cell 
\[\Lh\Xh \fto{\Lh\thah} \Lh\Xh' \inspace \GHCath\] 
is a weak $H$-equivalence.
\end{lemma}

\begin{proof}
We consider a triple $(\angordn, H, \be)$ consisting of
\begin{itemize}
\item an object $\angordn \in \Gsk \setminus \{\vstar\}$,
\item a subgroup $H \subseteq G$, and
\item a homomorphism $\be \cn H \to \Siangordn$.
\end{itemize}
By \cref{expl:gph_subgrp}, such a triple is equivalent to an object $\nbe \in \GH \setminus \{\vstar\}$.  For each such triple $(\angordn, H, \be)$, there is a commutative diagram of pointed $H$-functors
\begin{equation}\label{gg_ggg_weq_dia}
\begin{tikzpicture}[vcenter]
\def\h{3.3} \def\v{-1.5}
\draw[0cell=.9]
(0,0) node (a11) {\Xhnbeta}
(a11)++(\h,0) node (a12) {(\Lh \Xh)\angordn_{\be}}
(a12)++(1.15*\h,0) node (a13) {\phantom{(\Lh \Xh)\nbe}} 
(a13)++(0,.03) node (a13') {(\Lh \Xh)\nbe}
(a11)++(0,\v) node (a21) {\Xphnbeta}
(a12)++(0,\v) node (a22) {(\Lh \Xh')\angordn_{\be}}
(a13)++(0,\v) node (a23) {\phantom{(\Lh \Xh')\nbe}}
(a23)++(0,.02) node (a23') {(\Lh \Xh')\nbe}
;
\draw[1cell=.85]
(a11) edge node {\ug_{\Xh,\angordn}} node[swap] {\iso} (a12)
(a12) edge node {(\Lh \Xh)\jin} node[swap] {\iso} (a13)
(a21) edge node {\ug_{\Xh',\angordn}} node[swap] {\iso} (a22)
(a22) edge node {(\Lh \Xh')\jin} node[swap] {\iso} (a23)
(a11) edge node[swap] {\thahn} (a21)
(a12) edge node[swap] {\Lthahn} (a22)
(a13') edge[transform canvas={xshift=1em}] node[swap] {\Lthahnbe} (a23')
;
\end{tikzpicture}
\end{equation}
constructed as follows.
\begin{description}
\item[Top left]
The pointed $H$-category $\Xhnbeta$ is denoted by $\Xnbeta$ in \cref{Xnbeta_haction}.  Its $H$-action is the original one on $X\angordn$ twisted by $\be$.  The pointed $H$-functor $\ug_{\Xh,\angordn}$ is the unit isomorphism \cref{ugf_angordm} for the adjoint 2-equivalence $(\Lh,\ihst)$ \pcref{not:Lh_ihst}.  It is an $H$-isomorphism by \cref{uXh_Hequiv}.
\item[Top right] 
The isomorphism $\jin \cn \angordn \to \nbe$ in $\GH$ is given by $(1_{\ufsq}, \ang{1})$ \cref{jin_nbe}.  The functor $(\Lh \Xh)\jin$ is inverse to an $H$-isomorphism \pcref{Xjininv}, so it is an $H$-isomorphism.
\item[Bottom] 
The $H$-isomorphisms $\ug_{\Xh',\angordn}$ and $(\Lh \Xh')\jin$ are defined in the same way as the top horizontal functors, using $\Xh'$ instead of $\Xh$.
\item[Vertical]
The pointed $H$-functor $\thahn$ is the $\angordn$-component of $\thah$; it is denoted by $\tha_{\angordn}$ in \cref{than_xxp}.  The pointed $H$-functors $\Lthahn$ and $\Lthahnbe$ are the components of the 1-cell $\Lh\thah$ in $\GHCath$ at the objects $\angordn$ and $\nbe$.
\item[Commutativity]
The left rectangle commutes by the naturality of the unit $\ug$ with respect to the 1-cell $\thah$ in $\GHCatii$.  The right rectangle commutes by the naturality of $\Lh\thah$ with respect to the morphism $\jin$ in $\GH$.
\end{description}
By \cref{gg_weq_char}, $\tha$ is a weak $G$-equivalence if and only if $\thahn$ is a categorical weak $H$-equivalence for each triple $(\angordn, H, \be)$.  By the commutative diagram \cref{gg_ggg_weq_dia}, this happens if and only if $\Lthahnbe$ is a categorical weak $H$-equivalence for each subgroup $H \subseteq G$ and each object $\nbe \in \GH \setminus \{\vstar\}$.  This, in turn, means that $\Lh\thah$ is a weak $H$-equivalence \pcref{def:ggg_weq} for each subgroup $H \subseteq G$.
\end{proof}

\cref{Lg_weq} shows that the 2-equivalence \pcref{thm:ggcat_ggcatg_iieq} 
\[\GGCatii \fto{\Lg} \GGCatg\]
preserves and reflects weak $G$-equivalences in the sense of \cref{def:ggg_weq,def:gg_weq}.  It improves \cref{gg_ggg_weq} when the group $G$ is finite.

\begin{theorem}\label{Lg_weq}
For a finite group $G$, a 1-cell 
\[X \fto{\tha} X' \inspace \GGCatii\]
is a weak $G$-equivalence if and only if the 1-cell 
\[\Lg X \fto{\Lg\tha} \Lg X' \inspace \GGCatg\] 
is a weak $G$-equivalence.
\end{theorem}

\begin{proof}
The \emph{only if} assertion follows from the \emph{only if} assertion of \cref{gg_ggg_weq}, applied to the group $G$ itself.  This assertion does not require $G$ to be finite.

For the \emph{if} assertion, we assume that the 1-cell $\Lg\tha$ is a weak $G$-equivalence \pcref{def:ggg_weq}.  To prove that $\tha$ is a weak $G$-equivalence, we use \cref{gg_ggg_weq} and verify that, for each proper subgroup $H \subsetneq G$, the 1-cell  
\[\Lh\Xh \fto{\Lh\thah} \Lh\Xh' \inspace \GHCath\] 
is a weak $H$-equivalence.  In other words, we verify that the pointed $H$-functor
\begin{equation}\label{Lthah_nbe}
(\Lh \Xh)\nbe \fto{\Lthahnbe} (\Lh \Xh')\nbe
\end{equation}
is a categorical weak $H$-equivalence for each object $\nbe \in \GH \setminus \{\vstar\}$ with length $q \geq 0$.  There are two cases, depending on whether $\nbe$ is the empty sequence or not.
\begin{description}
\item[Empty sequence]
Suppose $\nbe \in \GH$ is the empty sequence $\ang{}$.  The only nonzero morphism in $\GG$ with codomain $\ang{}$ is the identity morphism of $\ang{}$ \cref{GG_empty_mor}.  By the definition of $\Lg$ \cref{Lg_f_nbe} applied to $G$ and $H$, there are an equality of pointed $G$-functors
\[(\Lg X)\ang{} = X\ang{} \fto{(\Lg\tha)_{\ang{}} = \tha_{\ang{}}} (\Lg X')\ang{} = X'\ang{}\]
and an equality of pointed $H$-functors
\[(\Lh \Xh)\ang{} = \Xh\ang{} \fto{(\Lh\thah)_{\ang{}} = \thahang} (\Lh \Xh')\ang{} = \Xh'\ang{}.\]
By the assumption that the 1-cell $\Lg\tha$ is a weak $G$-equivalence, the $\ang{}$-component $(\Lg\tha)_{\ang{}} = \tha_{\ang{}}$ is a categorical weak $G$-equivalence, hence also a categorical weak $H$-equivalence.  Since $\thahang$ is $\tha_{\ang{}}$ regarded as an $H$-functor \pcref{def:tha_h}, it follows that $(\Lh\thah)_{\ang{}} = \thahang$ is a categorical weak $H$-equivalence.
\item[Nonempty sequences]
For an object 
\[\nbe = \sordi{n}{\be}{j}_{j \in \ufsq} \in \GH \setminus \{\vstar,\ang{}\}\]
of length $q>0$, we prove that $\Lthahnbe$ \cref{Lthah_nbe} is the retract of a categorical weak $H$-equivalence.  Using the finiteness of $G$ and the constructions in \crefrange{ordmj_alphaj}{ip_ones}, there is a commutative diagram of pointed $H$-functors
\begin{equation}\label{Ltha_diag}
\begin{tikzpicture}[vcenter]
\def\v{-1.4} \def\h{3} \def\c{.6} \def\d{-.6}
\draw[0cell=.9]
(0,0) node (a11) {(\Lh\Xh)\nbe}
(a11)++(\h,0) node (a12) {(\Lg X)\malp}
(a12)++(\h,0) node (a13) {(\Lh\Xh)\nbe}
(a11)++(0,\v) node (a21) {(\Lh\Xh')\nbe}
(a12)++(0,\v) node (a22) {(\Lg X')\malp}
(a13)++(0,\v) node (a23) {(\Lh\Xh')\nbe}
;
\draw[1cell=.85]
(a11) edge node {\winc^\eps_*} (a12)
(a12) edge node {\wpr^\eps_*} (a13)
(a21) edge node {\winc^\eps_*} (a22)
(a22) edge node {\wpr^\eps_*} (a23)
(a11) edge node[swap] {\Lthahnbe} (a21)
(a12) edge node [swap] {\Lthamalp} (a22)
(a13) edge node {\Lthahnbe} (a23)
(a11) [rounded corners=2pt] |- ($(a12)+(-1,\c)$) -- node {1} ($(a12)+(1,\c)$) -| (a13)
;
\draw[1cell=.8]
(a21) [rounded corners=2pt] |- ($(a22)+(-1,\d)$) -- node {1} node[swap] {\phantom{x}} ($(a22)+(1,\d)$) -| (a23)
;
\end{tikzpicture}
\end{equation}
constructed as follows.
\begin{description}
\item[Top and bottom] 
The top row in \cref{Ltha_diag} is the same as the top row in \cref{ggcat_sp_LX_dia}, constructed in \crefrange{LXmalp}{wprepz_star}.  The bottom row in \cref{Ltha_diag} is constructed in the same way as the top row, using $\Xh'$ and $X'$ instead of $\Xh$ and $X$.
\item[Vertical]
The left and right vertical functors are given by the pointed $H$-functor $\Lthahnbe$ \cref{Lthah_nbe}.  The middle vertical functor $\Lthamalp$ is the $\malp$-component pointed $G$-functor of the 1-cell $\Lg\tha$ in $\GGCatg$.  
\item[Commutativity]
In the left rectangle in \cref{Ltha_diag}, by \cref{Lg_f_nbe,wincepz_star}, each of the two composites sends a representative pair \cref{LX_reps}
\[\big(\angordb \fto{\upom} \nbe; x\big) \in \GH(\angordb; \nbe) \ttimes \Xh\angordb\]
in $(\Lh\Xh)\nbe$ to the representative pair
\[\big(\angordb \fto{\upom} \nbe \fto{\winc^\eps} \malp; \tha_{\angordb} x\big) 
\in \GG(\angordb; \malp) \ttimes X'\angordb\]
in $(\Lg X')\malp$.  Similarly, by \cref{wprepz_star}, the right rectangle in \cref{Ltha_diag} commutes because each of the two composites sends a representative pair 
\[\big(\angorda \fto{\upom} \malp; x\big) \in \GG(\angorda; \malp) \ttimes X\angorda\]
in $(\Lg X)\malp$ to the representative pair
\[\big(\angorda \fto{\upom} \malp \fto{\wpr^\eps} \nbe ; \tha_{\angorda} x\big) 
\in \GH(\angorda; \nbe) \ttimes \Xh'\angorda\]
in $(\Lh\Xh')\nbe$.
\end{description}
This finishes the construction of the commutative diagram \cref{Ltha_diag} of pointed $H$-functors.

By the assumption that $\Lg\tha$ is a weak $G$-equivalence, its $\malp$-component $\Lthamalp$ is a categorical weak $G$-equivalence, hence also a categorical weak $H$-equivalence.  \cref{weq_retract}, applied to the group $H$ and the diagram \cref{Ltha_diag}, implies that $\Lthahnbe$ \cref{Lthah_nbe} is a categorical weak $H$-equivalence.  This proves that the 1-cell $\Lh\thah$ in $\GHCath$ is a weak $H$-equivalence, thereby proving the \emph{if} assertion of the \namecref{Lg_weq}.\qedhere
\end{description}
\end{proof}

\cref{Lg_weq} shows that the 2-equivalence $\Lg$ preserves and reflects weak $G$-equivalences.  \cref{igst_weq} shows that the 2-adjoint inverse \pcref{thm:ggcat_ggcatg_iieq}
\[\GGCatg \fto{\igst} \GGCatii\]
of $\Lg$ also preserves and reflects weak $G$-equivalences in the sense of \cref{def:ggg_weq,def:gg_weq}.

\begin{theorem}\label{igst_weq}
For a finite group $G$, a 1-cell 
\[X \fto{\tha} X' \inspace \GGCatg\]
is a weak $G$-equivalence if and only if the 1-cell 
\[\igst X \fto{\igst\tha} \igst X' \inspace \GGCatii\] 
is a weak $G$-equivalence.
\end{theorem}

\begin{proof}
For each object $\nbe \in \GG \setminus \{\vstar\}$, the naturality of the counit $\vg \cn \Lg\igst \to 1$ \cref{vg} with respect to the 1-cell $\tha$ implies that the diagram of pointed $G$-functors
\begin{equation}\label{vgnat_diag}
\begin{tikzpicture}[vcenter]
\def\v{-1.4}
\draw[0cell]
(0,0) node (a11) {(\Lg\igst X)\nbe}
(a11)++(3.5,0) node (a12) {X\nbe}
(a11)++(0,\v) node (a21) {(\Lg\igst X')\nbe}
(a12)++(0,\v) node (a22) {X'\nbe}
;
\draw[1cell=.9]
(a11) edge node {\vg_{X,\nbe}} node[swap] {\iso} (a12)
(a12) edge[transform canvas={xshift=-.8em}] node {\tha_{\nbe}} (a22)
(a11) edge[transform canvas={xshift=1em}] node[swap] {(\Lg\igst\tha)_{\nbe}} (a21)
(a21) edge node {\vg_{X',\nbe}} node[swap] {\iso} (a22)
;
\end{tikzpicture}
\end{equation}
commutes.  Thus, $\tha_{\nbe}$ is a categorical weak $G$-equivalence if and only if $(\Lg\igst\tha)_{\nbe}$ is a categorical weak $G$-equivalence.  By \cref{def:ggg_weq}, this means that $\tha$ is a weak $G$-equivalence if and only if $\Lg\igst \tha$ is a weak $G$-equivalence.  By the finiteness of $G$ and \cref{Lg_weq}, $\Lg\igst \tha$ is a weak $G$-equivalence if and only if $\igst \tha$ is a weak $G$-equivalence.
\end{proof}

\section{Strong $J$-Theory Comparison Weak $G$-Equivalence}
\label{sec:weq_app}

This section shows that, for an $\Oph$-pseudoalgebra $\Ah = \Catg(\EG,\A)$ with $\Oph = \Catg(\EG,\Op)$ \pcref{catgego} and $\Op$ a $\Uinf$-operad \pcref{as:OpA'}, Shimakawa strong $J$-theory \cref{jgos_jgossg} and strong $J$-theory \pcref{thm:Jgo_twofunctor} are connected by a weak $G$-equivalence in $\GGCatii$ \pcref{def:gg_weq}.

\begin{theorem}\label{thm:pistsgah_weq}
With the same assumptions as \cref{thm:pistweakgeq} and $G$ a finite group, there is a weak $G$-equivalence 
\begin{equation}\label{igst_pistsg_j}

\end{equation}
is a weak $G$-equivalence in $\GGCatg$ \pcref{def:ggg_weq}.  Applying \cref{igst_weq} to $\Pistsg_{\Ah}$ yields the weak $G$-equivalence
\begin{equation}\label{igst_pistsg}
%
\end{equation}
in $\GGCatii$.  To identify $\igst\! \smast \Sgohsg$ and $\igst\Hgohsg$ with, respectively, $\smast \Jgohssg$ and $\Jgohsg$, we first consider the following diagram of 2-functors for a $\Tinf$-operad $\Op$ \pcref{as:OpA} and an arbitrary group $G$.
\begin{equation}\label{j_comparison}
%
\end{equation}
\begin{itemize}
\item The 2-categories in \cref{j_comparison} are defined in \cref{oalgps_twocat,def:GGCat,def:ggcatg,def:fgcat,def:fgcatg}.
\item Each $\igst$ or $\smast$ is given by precomposing with $\ig$ or $\sma$.  Each $\igst$ is a 2-equivalence by \cref{thm:ggcat_ggcatg_iieq,thm:fgcat_fgcatg_iieq}.  The right square in \cref{j_comparison} commutes because the diagram of functors
\begin{equation}\label{FG_isma}
%
\end{equation}
commutes by \cref{def:smashFskGsk,def:smashFGGG,def:gskel_gg,def:fsk_fg}. 
\item The 2-natural transformation $\Pistsg$ in \cref{j_comparison} is the strong $H$-theory comparison for $\Op$ \cref{Pistarsg_twonat}, from Shimakawa strong $H$-theory $\Sgosg$ \cref{def:sgosg} to strong $H$-theory $\Hgosg$ \pcref{Hgo_twofunctor}.  
\end{itemize}
\begin{description}
\item[Domain]
Along the top of \cref{j_comparison}, the equality
\[\AlgpspsO \fto{\Jgossg = \igst \Sgosg} \FGCat\]
is the definition of Shimakawa strong $J$-theory for $\Op$ \cref{jgos_jgossg}.  Applying this to $\Oph = \Catg(\EG,\Op)$ for a $\Uinf$-operad $\Op$ and using the commutativity of the right square in \cref{j_comparison}, the domain of the weak $G$-equivalence $\igst\Pistsg_{\Ah}$ \cref{igst_pistsg} is the $\Gskg$-category given by
\begin{equation}\label{igst_pistsg_dom}
\begin{split}
\igst\! \smast \Sgohsg\Ah 
&= \smast \igst \Sgohsg\Ah \\
&= \smast \Jgohssg \Ah.
\end{split}
\end{equation}
\item[Codomain] 
Along the bottom of \cref{j_comparison}, by \cref{jgohgoigst}, there is a factorization of the strong $J$-theory 2-functor for $\Op$ \pcref{thm:Jgo_twofunctor} as
\[\AlgpspsO \fto{\Jgosg = \igst \Hgosg} \GGCatii.\]
Applying this factorization to $\Oph$, the codomain of the weak $G$-equivalence $\igst\Pistsg_{\Ah}$ \cref{igst_pistsg} is the strong $J$-theory
\begin{equation}\label{igst_pistsg_cod}
\Gsk \fto{\Jgohsg\Ah = \igst\Hgohsg\Ah} \Gcatst
\end{equation}
for the $\Oph$-pseudoalgebra $\Ah$ \pcref{A_ptfunctor}.  
\end{description}
Combining \cref{igst_pistsg,igst_pistsg_cod,igst_pistsg_dom} yields the desired weak $G$-equivalence \cref{igst_pistsg_j}.
\end{proof}

\begin{explanation}[Strong $J$-Theory Comparison]\label{expl:pistsgah_weq}
By \cref{expl:ggcat_weq} \eqref{expl:ggcat_weq_ii}, the fact that the strong $J$-theory comparison \cref{igst_pistsg_j} 
\[\smast \Jgohssg\Ah \fto{\igst\Pistsg_{\Ah}} \Jgohsg\Ah\]
is a weak $G$-equivalence means that, for each object $\angordn \in \Gsk \setminus \{\vstar\}$ \cref{Gsk_objects} and each graph subgroup $K \in \FGSin$ \pcref{def:graph_subgrp}, the $K$-fixed subfunctor
\begin{equation}\label{pistsgahn_k}
\begin{tikzpicture}[vcenter]
\def\v{-1}
\draw[0cell=.9]
(0,0) node (a11) {((\smast \Jgohssg\Ah)\angordn)^K}
(a11)++(4.3,0) node (a12) {((\Jgohsg\Ah)\angordn)^K}
(a11)++(0,\v) node (a21) {(\Ahsgsman)^K}
(a12)++(0,\v) node (a22) {(\Ahsgn)^K}
;
\draw[1cell=.85]
(a11) edge[equal] (a21)
(a12) edge[equal] (a22)
(a11) edge node {(\igst\Pistsg_{\Ah})_{\!\angordn}^K} (a12)
(a21) edge node {(\Pistsg_{\Ah,\angordn})^K} node[swap] {\sim} (a22)
;
\end{tikzpicture}
\end{equation}
is an equivalence of categories. 
\begin{itemize}
\item In the domain in \cref{pistsgahn_k}, $\Ahsgsman$ is the pointed $G$-category of strong $(\sma\angordn)$-systems in $\Ah$ \cref{Asgordnbe}.  Its $G$-action is described in \crefrange{g_Aordn}{gtha_scomp}.  This $G$-action extends to a $\GSinp$-action using \cref{psitilsg_f,Xn_gsin_act}.
\item In the codomain in \cref{pistsgahn_k}, $\Ahsgn$ is the pointed $G$-category of strong $\angordn$-systems in $\Ah$ \cref{sgAordnbe}.  This $G$-action extends to a $\GSinp$-action using \cref{AF_sg,Xn_gsin_act}.
\item The $\angordn$-component $\Pistsg_{\Ah,\angordn}$ is defined in \crefrange{PistA_empty}{thetatimes_angs}.
\end{itemize}
More detailed description of the $\GSinp$-action on $\Ahsgsman$ and $\Ahsgn$ is given next.
\begin{description}
\item[Domain $\GSinp$-action]
Using \cref{psitil_f,g_Aordn,gtha_scomp,Xn_gsin_act}, the $\GSinp$-action on the category $\Ahsgsman$, the domain of $\Pistsg_{\Ah,\angordn}$, is given as follows.  Suppose
\begin{itemize}
\item $(a,\gl)$ is a strong $(\sma\angordn)$-system in $\Ah$ \cref{nsys} with $\angordn = \ang{\ordn_j}_{j \in \ufsq} \in \Gsk$ and
\item $(g,\bsi) \in \GSin$ is an element with $\bsi = \ang{\si_j \in \Si_{n_j}}_{j \in \ufsq}$ \cref{Siangord_obj}.
\end{itemize}  
For a subset 
\[s \subseteq \ufs{n_1 \Cdots n_q} = \txprod_{j \in \ufsq}\, \ufsn_j,\] 
the $s$-component object \cref{nsys_s} of the strong $(\sma\angordn)$-system $(g,\bsi) \cdot (a,\gl)$ in $\Ah$ is
\begin{equation}\label{gbsi_as}
\big((g,\bsi) \cdot (a,\gl) \big)_s = g a_{\bsi^{-1} s} \in \Ah
\end{equation}
with the subscript of $a$ given by the subset
\begin{equation}\label{bsiinvs}
\bsi^{-1} s = (\si_1 \ttimes \Cdots \ttimes \si_q)^{-1} s \subseteq \ufs{n_1 \Cdots n_q}.
\end{equation}
For an object $x \in \Oph(r)$ with $r \geq 0$ and a partition $s = \txcoprod_{i \in \ufsr}\, s_i$, the gluing isomorphism \cref{gl-morphism} of the strong $(\sma\angordn)$-system $(g,\bsi) \cdot (a,\gl)$ at $(x; s, \ang{s_i}_{i \in \ufsr})$ is the following isomorphism in $\Ah$.
\begin{equation}\label{gbsi_gl}
\begin{tikzpicture}[vcenter]
\draw[0cell=.9]
(0,0) node (a1) {\gaAh_r\big(x; \ang{g a_{\bsi^{-1} s_i}}_{i \in \ufsr} \big)}
(a1)++(0,-1) node (a2) {g \gaAh_r\big(\ginv x ; \ang{a_{\bsi^{-1} s_i}}_{i \in \ufsr} \big)}
(a1)++(2.5,-.05) node (a3) {g a_{\bsi^{-1} s}}
(a3)++(-.3,0) node (a3') {\phantom{g a_{\bsi^{-1} s}}}
;
\draw[1cell=.85]
(a1) edge[equal] (a2)
(a2) [rounded corners=2pt, shorten >=-.5ex] -| node[swap,pos=.75] {g \gl_{\ginv x ; \bsi^{-1}s, \ang{\bsi^{-1} s_i}_{i \in \ufsr}}} (a3')
;
\end{tikzpicture}
\end{equation}
For a morphism $\tha$ in $\Ahsgsman$, the $s$-component morphism \cref{theta_s} of $(g,\bsi) \cdot \tha$ is
\begin{equation}\label{gbsi_thas}
\big((g,\bsi) \cdot \tha) \big)_s = g \tha_{\bsi^{-1} s} \inspace \Ah.
\end{equation}
\item[Codomain $\GSinp$-action]
Using \cref{nsystem_gaction,gtheta,psitil_functor,Xn_gsin_act}, the $\GSinp$-action on the category $\Ahsgn$, the codomain of $\Pistsg_{\Ah,\angordn}$, is given as follows.  Suppose $(a,\glu)$ is a strong $\angordn$-system in $\Ah$ \cref{nsystem}.  For a marker $\angs = \ang{s_j \subseteq \ufsn_j}_{j \in \ufsq}$, the $\angs$-component object \cref{a_angs} of the strong $\angordn$-system $(g,\bsi) \cdot (a,\glu)$ in $\Ah$ is
\begin{equation}\label{gbsi_aangs}
\big((g,\bsi) \cdot (a,\glu) \big)_{\angs} = g a_{\ang{\siinv s}} \in \Ah
\end{equation}
with the subscript of $a$ given by the marker
\begin{equation}\label{bsiinvangs}
\ang{\siinv s} = \ang{\siinv_j s_j \subseteq \ufsn_j}_{j \in \ufsq}.
\end{equation}
For an object $x \in \Oph(r)$ with $r \geq 0$, an index $k \in \ufsq$, and a partition $s_k = \txcoprod_{i \in \ufsr}\, s_{k,i}$, the gluing isomorphism \cref{gluing-morphism} of the strong $\angordn$-system $(g,\bsi) \cdot (a,\glu)$ at $(x; \angs, k, \ang{s_{k,i}}_{i \in \ufsr})$ is the following isomorphism in $\Ah$.
\begin{equation}\label{gbsi_glu}
\begin{tikzpicture}[vcenter]
\draw[0cell=.9]
(0,0) node (a1) {\gaAh_r\big(x; \ang{g a_{\ang{\siinv s} \compk (\siinv_k s_{k,i})}}_{i \in \ufsr} \big)}
(a1)++(0,-1) node (a2) {g \gaAh_r\big(\ginv x ; \ang{a_{\ang{\siinv s} \compk (\siinv_k s_{k,i})}}_{i \in \ufsr} \big)}
(a1)++(3,-.05) node (a3) {g a_{\ang{\siinv s}}}
(a3)++(-.3,0) node (a3') {\phantom{g a_{\ang{\siinv s}}}}
;
\draw[1cell=.85]
(a1) edge[equal, shorten >=-1ex] (a2)
(a2) [rounded corners=2pt, shorten >=-.5ex] -| node[swap,pos=.75] {g \glu_{\ginv x ; \ang{\siinv s},\, k, \ang{\siinv_k s_{k,i}}_{i \in \ufsr}}} (a3')
;
\end{tikzpicture}
\end{equation}
For a morphism $\tha$ in $\Ahsgn$, the $\angs$-component morphism \cref{theta_angs} of $(g,\bsi) \cdot \tha$ is
\begin{equation}\label{gbsi_thaangs}
\big((g,\bsi) \cdot \tha) \big)_{\angs} = g \tha_{\ang{\siinv s}}
\end{equation}
in $\Ah$.\defmark
\end{description}
\end{explanation}

\section{Specialness of Strong $H$-Theory and $J$-Theory}
\label{sec:special_theory}

This section proves that strong $H$-theory and strong $J$-theory for $\Oph$-pseudoalgebras of the form $\Ah = \Catg(\EG,\A)$, where $\Oph = \Catg(\EG,\Op)$, are special in the sense of \cref{def:sp_gggcat,def:sp_ggcat}.

\secoutline
\begin{itemize}
\item \cref{i_bbr} proves that, for the twisted product of a $G$-category of the form $\Catg(\EG,\C)$, the inclusion $G$-functor is the left adjoint of an adjoint $G$-equivalence.  \cref{i_bbr} is used in the proof of \cref{thm:h_special}; see \cref{zdsg_catweakg}.
\item \cref{thm:h_special} proves that, for an $\Oph$-pseudoalgebra of the form $\Ah$, the strong $H$-theory $\Hgohsg\Ah$ is a special $\GGG$-category.
\item \cref{cor:h_special_shi} proves that the $\GGG$-category $\smast\Sgohsg \Ah$ is special, where $\Sgohsg$ is Shimakawa strong $H$-theory.
\item \cref{expl:segnbe_zd} discusses why it is necessary to use \emph{strong} $H$-theory in \cref{thm:h_special,cor:h_special_shi}. 
\item \cref{thm:j_special} proves that, for a finite group $G$, the strong $J$-theory $\Jgohsg\Ah$ of an $\Oph$-pseudoalgebra $\Ah$ is a special $\Gskg$-category.
\item \cref{cor:j_special_shi} proves that, for a finite group $G$, the $\Gskg$-category $\smast\Jgohssg \Ah$ is special, where $\Jgohssg$ is Shimakawa strong $J$-theory.
\item \cref{rk:hj_special} discusses why $G$ is assumed to be finite in \cref{thm:j_special,cor:j_special_shi}.
\end{itemize}

\subsection*{Inclusion $G$-Functors for Twisted Products of $G$-Thickening}

For \cref{i_bbr}, recall
\begin{itemize}
\item the inclusion $G$-functor \cref{incC}
\[\C \fto{\inc} \Chat = \Catgegc\]
from a small $G$-category $\C$ to its $G$-thickening $\Catgegc$ \pcref{def:Catg,def:translation_cat},
\item the $\nbe$-twisted product $\proCnbe$ \pcref{def:proCnbe}, and
\item an adjoint $G$-equivalence \pcref{def:adjointGeq}.
\end{itemize}

\begin{lemma}\label{i_bbr}
Suppose $\C$ is a small $G$-category for a group $G$, and $\nbe \in \GG \setminus \{\vstar,\ang{}\}$ is an object.  Then there is an adjoint $G$-equivalence\index{adjoint G-equivalence@adjoint $G$-equivalence}
\begin{equation}\label{incbbr_statement}
\begin{tikzpicture}[vcenter]
\def\h{1.8} \def\t{25}
\draw[0cell]
(0,0) node (a1) {\phantom{C}}
(a1)++(-.2,0) node (a1') {\Chnbe}
(a1)++(\h,0) node (a2) {\phantom{C}}
(a2)++(1.05,0) node (a2') {\Catgegchnbe}
(a1)++(\h/2,0) node () {\gsim}
;
\draw[1cell=.9]
(a1) edge[bend left=\t] node {\inc} (a2)
(a2) edge[bend left=\t] node {\bbr} (a1)
;
\end{tikzpicture}
\end{equation}
with $\inc$ denoting the inclusion $G$-functor for the $\nbe$-twisted product of $\Chat = \Catgegc$.
\end{lemma}

\begin{proof}
This proof is obtained from the proof of \cref{thm:SgoAh} via a change of notation as follows.  Suppose $\nbe = \sordi{n}{\be}{j}_{j \in \ufsq}$ has length $q>0$.  We reuse the notation in \crefrange{f_sub_g}{fgsh} by replacing $\Ahmal$ and $s \subseteq \ufsm$ with, respectively, $\Chnbe$ and $\bdi \in \ufs{n_1 \Cdots n_q}$ \cref{bdi}.  We construct the $G$-functor $\bbr$, the $G$-equivariant unit $\iunit \cn 1 \fto{=} \bbr\inc$, and the $G$-equivariant counit $\icounit \cn \inc\bbr \fiso 1$.
\begin{description}
\item[Inverse $G$-functor]
The $G$-functor $\bbr$ is defined in the same way as the $G$-functor \cref{cni_functor}
\[\Catg(\EG, \Ahmal) \fto{\cni} \Ahmal\]
by changing the notation from $(\Ahmal,s)$ to $(\Chnbe,\bdi)$ and ignoring gluing morphisms.
\begin{description}
\item[Objects] 
More precisely, for a functor $f \cn \EG \to \Chnbe$, the image $\bbr f \in \Chnbe$ has, for each $q$-tuple $\bdi \in \ufs{n_1 \Cdots n_q}$, $\bdi$-th coordinate functor
\[\EG \fto{(\bbr f)_{\bdi}} \C\]
that sends an object $g \in \EG$ to the object
\[(\bbr f)_{\bdi}(g) = f_{g,\bdi,g} \in \C.\]
The functor $(\bbr f)_{\bdi}$ sends an isomorphism $[h,g] \cn g \fiso h$ in $\EG$ to the isomorphism in $\C$ defined by the following commutative diagram.
\begin{equation}\label{rfihg_diag}
\begin{tikzpicture}[vcenter]
\def\h{3.5} \def\v{1} \def\t{20}
\draw[0cell]
(0,0) node (a11) {f_{g,\bdi,g}}
(a11)++(\h/2,\v) node (a12) {f_{h,\bdi,g}}
(a11)++(\h/2,-\v) node (a21) {f_{g,\bdi,h}}
(a11)++(\h,0) node (a22) {f_{h,\bdi,h}}
;
\draw[1cell=.9]
(a11) edge[bend left=\t] node[pos=.3] {f_{[h,g],\bdi,g}} (a12)
(a12) edge[bend left=\t] node[pos=.7] {f_{h,\bdi,[h,g]}} (a22)
(a11) edge[bend right=\t] node[swap,pos=.3] {f_{g,\bdi,[h,g]}} (a21)
(a21) edge[bend right=\t] node[swap,pos=.7] {f_{[h,g],\bdi,h}} (a22)
(a11) edge node {(\bbr f)_\bdi [h,g]} (a22)
;
\end{tikzpicture}
\end{equation}
With $(\bbr,\bdi)$ instead of $(\cni,s)$, the proof of \cref{cnifs_welldef} proves that the object $\bbr f \in \Chnbe$ is well defined.  In other words, each $(\bbr f)_{\bdi}$ is a functor.
\item[Morphisms] 
The functor $\bbr$ sends a morphism $\tha \cn f \to f'$ in $\Catgegchnbe$ to the morphism 
\[\bbr f \fto{\bbr\tha} \bbr f' \inspace \Chnbe\]
with, for each $\bdi \in \ufs{n_1 \Cdots n_q}$, $\bdi$-th coordinate natural transformation as follows.
\begin{equation}\label{rtha_diag}
\begin{tikzpicture}[vcenter]
\def\t{27}
\draw[0cell]
(0,0) node (a1) {\phantom{A}}
(a1)++(2.3,0) node (a2) {\C}
(a1)++(-.1,0) node (a1') {\EG}
;
\draw[1cell=.9]
(a1) edge[bend left=\t] node {(\bbr f)_{\bdi}} (a2)
(a1) edge[bend right=\t] node[swap] {(\bbr f')_{\bdi}} (a2)
;
\draw[2cell]
node[between=a1 and a2 at .33, rotate=-90, 2label={above,(\bbr\tha)_{\bdi}}] {\Rightarrow}
;
\end{tikzpicture}
\end{equation}
For each object $g \in \EG$, the $g$-component of $(\bbr\tha)_{\bdi}$ is defined as the morphism 
\[(\bbr f)_{\bdi}(g) = f_{g,\bdi,g} \fto{(\bbr\tha)_{\bdi,g} = \tha_{g,\bdi,g}} (\bbr f')_{\bdi}(g) = f'_{g,\bdi,g}\]
in $\C$.  With $(\bbr,\bdi)$ instead of $(\cni,s)$, the proof of \cref{cnitha_welldef} \eqref{cnitha_i} proves that the morphism $\bbr \tha$ is well defined.  In other words, each $(\bbr\tha)_{\bdi}$ is a natural transformation.
\item[$G$-functoriality]
The functoriality of $\bbr$ follows from the fact that identity morphisms and composition are defined componentwise in $\Catg(\EG,-)$ \pcref{def:Catg} and coordinatewise in $\Chnbe$ \pcref{def:proCnbe}.  Using the $G$-action on twisted products \cref{proCnbe_gaction} and $(\bbr,\bdi)$ instead of $(\cni,s)$, the computation in \cref{cni_eq_objects,cni_eq_mor} proves that the functor $\bbr$ is $G$-equivariant.
\end{description}
This finishes the construction of the $G$-functor $\bbr$.
\item[Unit]
With $(\bbr,\bdi)$ instead of $(\cni,s)$, the computation in \cref{ucni_cod_obj,ucni_cod_mor} proves that the composite $G$-functor
\[\Chnbe \fto{\inc} \Catgegchnbe \fto{\bbr} \Chnbe\]
is equal to the identity.  Analogous to \cref{ucni}, we define the unit $\iunit$ as the identity natural transformation
\begin{equation}\label{irunit_diag}
\begin{tikzpicture}[vcenter]
\def\h{3} \def\t{20}
\draw[0cell]
(0,0) node (a1) {\Chnbe}
(a1)++(-.3,0) node (a1') {\phantom{\Chnbe}}
(a1)++(\h/2,-1) node (a2) {\Catgegchnbe}
(a1)++(\h,0) node (a3) {\Chnbe}
(a3)++(.3,0) node (a3') {\phantom{\Chnbe}}
;
\draw[1cell=.9]
(a1) edge[bend left=0] node {1} (a3)
(a1') edge[bend right=\t, shorten >=-.5ex] node[swap,pos=.4] {\inc} (a2)
(a2) edge[bend right=\t, shorten <=-.5ex] node[swap,pos=.6] {\bbr} (a3')
;
\draw[2cell]
node[between=a1 and a3 at .32, shift={(0,-.45)}, rotate=-90, 2label={above,\iunit = 1_1}] {\Rightarrow}
;
\end{tikzpicture}
\end{equation}
of the identity $G$-functor on $\Chnbe$.
\item[Counit]
For a functor $f \cn \EG \to \Chnbe$, objects $g,h \in \EG$, and $\bdi \in \ufs{n_1 \Cdots n_q}$, the $(g,\bdi,h)$-component object of $\inc\bbr f$ is given by
\[(\inc\bbr f)_{g,\bdi,h} = (\bbr f)_{\bdi,h} = f_{h,\bdi,h}.\]
Analogous to \cref{ccni}, we define the counit $\icounit$ as the natural isomorphism
\begin{equation}\label{ricounit_diag}
\begin{tikzpicture}[vcenter]
\def\h{4} \def\t{25} \def\b{.7}
\draw[0cell]
(0,0) node (a1) {\Catgegchnbe}
(a1)++(\h/2,1) node (a2) {\Chnbe}
(a1)++(\h,0) node (a3) {\Catgegchnbe}
(a1)++(\b,0) node (a1') {\phantom{\Ah}}
(a3)++(-\b,0) node (a3') {\phantom{\Ah}}
;
\draw[1cell=.9]
(a1) edge node[swap] {1} (a3)
(a1') edge[bend left=\t, shorten <=.3ex] node {\bbr} (a2)
(a2) edge[bend left=\t] node {\inc} (a3')
;
\draw[2cell]
node[between=a1 and a3 at .45, shift={(0,.4)}, rotate=-90, 2label={above,\icounit}] {\Rightarrow}
;
\end{tikzpicture}
\end{equation}
with $(f,g,\bdi,h)$-component isomorphism given by
\[(\inc \bbr f)_{g,\bdi,h} = f_{h,\bdi,h}
\fto[\iso]{\icounit_{f,g,\bdi,h} = f_{[g,h],\bdi,h}} f_{g,\bdi,h}\]
in $\C$.  Using the $G$-action on twisted products \cref{proCnbe_gaction} and $(\bbr,\bdi)$ instead of $(\cni,s)$, the proof of \cref{ccni_welldef} proves that $\icounit \cn \inc\bbr \to 1$ is a $G$-natural isomorphism.
\item[Triangle identities]
The left and right triangle identities for an adjunction are proved by the computation in \cref{SgoAh_left_tri,SgoAh_right_tri} with $(\bbr,\bdi)$ instead of $(\cni,s)$.  
\end{description}
By \cref{expl:adjointGeq}, the quadruple $(\inc,\bbr,\iunit,\icounit)$ is an adjoint $G$-equivalence.
\end{proof}

\subsection*{Specialness of Strong $H$-Theory}

\begin{assumption}\label{as:uinf_ah}
We consider a $\Uinf$-operad $\Op$ \pcref{as:OpA'} for an arbitrary group $G$, an $\Op$-pseudoalgebra $(\A,\gaA,\phiA)$ \pcref{def:pseudoalgebra}, the $\Uinf$-operad $\Oph = \Catg(\EG,\Op)$, and the $\Oph$-pseudoalgebra \pcref{catgego}
\[\big(\Ah = \Catg(\EG,\A), \gaAh, \phiAh\big),\]
where $\EG$ is the translation category of $G$ \pcref{def:translation_cat}.
\end{assumption}

For \cref{thm:h_special}, recall strong $H$-theory \pcref{Hgo_twofunctor}
\[\AlgpspsO \fto{\Hgosg} \GGCatg.\]
With the $\Gcat$-operad $\Oph = \Catg(\EG,\Op)$ in place of $\Op$, this 2-functor sends each $\Oph$-pseudoalgebra, such as $\Ah = \Catg(\EG,\A)$ for an $\Op$-pseudoalgebra $\A$, to a $\GGG$-category \pcref{A_ptfunctorGG}.  Recall that a $\GGG$-category is \emph{special} if each $\nbe$-Segal functor $\segnbe$ is a categorical weak $G$-equivalence \pcref{def:sp_gggcat}.

\begin{theorem}\label{thm:h_special}\index{H-theory@$H$-theory!special}\index{special!GGG-category@$\GGG$-category!strong $H$-theory}
Under \cref{as:uinf_ah}, the strong $H$-theory of the $\Oph$-pseudoalgebra $\Ah$
\[\GG \fto{\Hgohsg \Ah = \Ahsgdash} \Catgst\]
is a special $\GGG$-category.
\end{theorem}

\begin{proof}
We verify that, for each object $\nbe \in \GG \setminus \{\vstar,\ang{}\}$ of length $q>0$ and the strong $H$-theory of $\Ah = \Catg(\EG,\A)$, the $\nbe$-Segal functor \cref{nbe_segal} 
\begin{equation}\label{segal_Ahsg}
\Ahsgnbe \fto{\segnbe} \Ahsgonenbe
\end{equation}
is a categorical weak $G$-equivalence \pcref{def:cat_weakg}.  We use a 2-out-of-3 argument similar to the proof of \cref{thm:pistweakgeq}.
\begin{description}
\item[Segal functor] 
For the strong $H$-theory of the $\Op$-pseudoalgebra $\A$, we first identify the $\nbe$-Segal functor
\[\Asgnbe \fto{\segnbe} \Asgonenbe\]
with the pointed $G$-functor \cref{hgosgtoprod} 
\[\Asgnbe \fto{\zdsg} \proAnbe\]
as follows.  With the object $\angordone = \ang{\ordone}_{j \in \ufsq} \in \GG$ \cref{angordone}, there is a pointed $G$-isomorphism \cref{Aordtuone}
\begin{equation}\label{upkap}
\Asgone \fto[\iso]{\upkap} \A
\end{equation}
that sends
\begin{itemize}
\item a strong $\angordone$-system to its $\ang{\{1\}}_{j \in \ufsq}$-component object and
\item a morphism of $\angordone$-systems to its $\ang{\{1\}}_{j \in \ufsq}$-component morphism.  
\end{itemize}
Consider the $\bdi$-characteristic morphism $\segbdi \cn \nbe \to \angordone$ \cref{i-char} for a $q$-tuple $\bdi = \ang{i_j}_{j \in \ufsq} \in \ufs{n_1 \Cdots n_q}$.  By the axioms \cref{system_obj_unity,nsystem_mor_unity}, the pointed $G$-functor \cref{psitil_functor}
\begin{equation}\label{Ahsg_segbdi}
\Asgnbe \fto{\Asgsegbdi} \Asgone
\end{equation}
sends a strong $\nbe$-system $(a,\glu)$ in $\A$ to the strong $\angordone$-system with 
\begin{itemize}
\item $\ang{\{1\}}_{j \in \ufsq}$-component object given by $a_{\ang{\{i_j\}}_{j \in \ufsq}}$,
\item all other component objects given by the basepoint $\zero \in \A$, and
\item all gluing morphisms given by identities,
\end{itemize}
and likewise for morphisms in $\Asgnbe$.  By \cref{hgotoprod_def,Xsegbdi,upkap,Ahsg_segbdi}, there is a commutative diagram of pointed $G$-functors
\begin{equation}\label{segnbe_zdsg}
\begin{tikzpicture}[vcenter]
\def\h{2.8} \def\c{.7}
\draw[0cell=1]
(0,0) node (a11) {\Asgnbe}
(a11)++(\h,0) node (a12) {\Asgonenbe}
(a12)++(\h,0) node (a13) {\phantom{\proAnbe}}
(a13)++(0,.055) node (a13') {\proAnbe}
;
\draw[1cell=.85]
(a11) edge node {\segnbe} (a12)
(a12) edge node {\pro{\upkap}{\nbe}} node [swap] {\iso} (a13)
(a11) [rounded corners=2pt] |- ($(a12)+(-1,\c)$) -- node {\zdsg} ($(a12)+(1,\c)$) -| (a13')
;
\end{tikzpicture}
\end{equation}
with the $G$-isomorphism $\pro{\upkap}{\nbe}$ given coordinatewise by $\upkap$.  Since the commutative diagram \cref{segnbe_zdsg} is defined for an arbitrary $\Op$-pseudoalgebra $\A$, it is also defined for the $\Oph$-pseudoalgebra $\Ah$.
\item[2-out-of-3] 
Using \cref{segnbe_zdsg} for $(\Oph,\Ah)$, to show that the $\nbe$-Segal functor $\segnbe$ \cref{segal_Ahsg} is a categorical weak $G$-equivalence, it suffices to show that $\zdsg$ is a categorical weak $G$-equivalence.  To prove that, we consider the commutative diagram of pointed $G$-functors 
\begin{equation}\label{zdsg_catweakg}
\begin{tikzpicture}[vcenter]
\def\v{-1.4}
\draw[0cell=1]
(0,0) node (a11) {\Ahsgnbe}
(a11)++(5,0) node (a12) {\phantom{\proAhnbe}}
(a12)++(0,.04) node (a12') {\proAhnbe}
(a11)++(0,\v) node (a21) {\Catg(\EG,\Ahsgnbe)}
(a12)++(0,\v) node (a22) {\Catg(\EG,\proAhnbe)}
;
\draw[1cell=.85]
(a11) edge node {\zdsg} (a12)
(a12') edge node {\inc} node[swap] {\gsim} (a22)
(a11) edge node {\gsim} node[swap] {\inc} (a21)
(a21) edge node {\Catg(\EG,\zdsg)} node[swap] {\eqg} (a22)
;
\end{tikzpicture}
\end{equation}
constructed as follows.
\begin{description}
\item[Top] $\zdsg$ is the pointed $G$-functor in \cref{segnbe_zdsg} for $(\Oph,\Ah)$.
\item[Left] The functor $\inc$ on the left is the inclusion $G$-functor \cref{incC} for the small $G$-category $\Ahsgnbe$ of strong $\nbe$-systems in $\Ah$.  By \cref{thm:HgoAh} \eqref{thm:HgoAh_ii}, it is part of an adjoint $G$-equivalence, hence also a categorical weak $G$-equivalence.
\item[Right] The functor $\inc$ on the right is the inclusion $G$-functor for the $\nbe$-twisted product $\proAhnbe$ of $\Ah$.  By \cref{i_bbr}, it is part of an adjoint $G$-equivalence, hence also a categorical weak $G$-equivalence.
\item[Bottom] By \cref{Pist_to_prod} \eqref{Pist_to_prod_i} and \cref{thm:zdsg_eq} applied to $(\Oph,\Ah)$, $\zdsg$ is a pointed $G$-functor and a nonequivariant equivalence of categories.  By \cref{merling_2.16}, $\Catg(\EG,\zdsg)$ is a categorical weak $G$-equivalence.
\item[Commutativity] 
For each object or morphism $b \in \Ahsgnbe$, under either composite in \cref{zdsg_catweakg}, the image of $b$ in $\Catg(\EG,\proAhnbe)$ is constant at $\zdsg b$.  
\end{description}
\end{description}
Thus, the commutative diagram \cref{zdsg_catweakg} implies that $\zdsg$ is a categorical weak $G$-equivalence.
\end{proof}

\begin{remark}\label{rk:zdsg_gmmo}
The categorical weak $G$-equivalence $\zdsg$ in \cref{zdsg_catweakg} is also used in the equivalence between Shimakawa and GMMO $K$-theories.  See \cref{thm:gmmo_shi} \cref{thm:gmmo_shi_iii}.  That proof uses $\zdsg$ when $\nbe$ has length one, meaning that it is a pointed finite $G$-set.  In that context, the pointed $G$-functor $\zdsg$ is denoted by $\cgs_{\Ah,\mal}$.
\end{remark}

Recall Shimakawa strong $H$-theory $\Sgosg$ \pcref{Sgo_twofunctor} and the 2-functor $\smast$ \pcref{smashstar}:
\[\AlgpspsO \fto{\Sgosg} \FGCatg \fto{\smast} \GGCatg.\]
With the $\Gcat$-operad $\Oph = \Catg(\EG,\Op)$ in place of $\Op$, the preceding composite 2-functor sends each $\Oph$-pseudoalgebra to a $\GGG$-category.

\begin{corollary}\label{cor:h_special_shi}\index{Shimakawa H-theory@Shimakawa $H$-theory!special}\index{special!GGG-category@$\GGG$-category!Shimakawa strong $H$-theory}
Under \cref{as:uinf_ah}, the $\GGG$-category
\[\GG \fto{\smast\Sgohsg \Ah} \Catgst\]
is special.
\end{corollary}

\begin{proof}
By \cref{thm:pistweakgeq}, there is a weak $G$-equivalence
\[\smast\Sgohsg\Ah \fto{\Pistsg_{\Ah}} \Hgohsg\Ah\]
in $\GGCatg$ \pcref{def:ggg_weq}.  By \cref{thm:h_special}, the $\GGG$-category $\Hgohsg\Ah$ is special.  Thus, $\smast\Sgohsg\Ah$ is special by \cref{ggg_weq_sp} \eqref{ggg_weq_sp_ii}.
\end{proof}

\begin{explanation}[Necessity of Strong $H$-Theory]\label{expl:segnbe_zd}
We explain why \cref{thm:h_special,cor:h_special_shi} use \emph{strong} $H$-theory instead of $H$-theory \pcref{A_ptfunctorGG}.  The discussion from \crefrange{upkap}{segnbe_zdsg} also yields a commutative diagram of pointed $G$-functors
\begin{equation}\label{zd_factor}
\begin{tikzpicture}[vcenter]
\def\h{2.6} \def\c{.7}
\draw[0cell=1]
(0,0) node (a11) {\Anbe}
(a11)++(\h,0) node (a12) {\phantom{\Aonenbe}}
(a12)++(0,.01) node (a12') {\Aonenbe}
(a12)++(\h,0) node (a13) {\phantom{\proAnbe}}
(a13)++(0,.055) node (a13') {\proAnbe}
;
\draw[1cell=.85]
(a11) edge node {\segnbe} (a12)
(a12) edge node {\pro{\upkap}{\nbe}} node [swap] {\iso} (a13)
(a11) [rounded corners=2pt] |- ($(a12)+(-1,\c)$) -- node {\zd} ($(a12)+(1,\c)$) -| (a13')
;
\end{tikzpicture}
\end{equation}
that identifies
\begin{itemize}
\item the $\nbe$-Segal functor $\segnbe$ for $\Hgo\A = \Adash$ and 
\item the pointed $G$-functor $\zd$ \cref{hgotoprod}.
\end{itemize}  
Analogous to \cref{zdsg_catweakg}, there is a commutative diagram of pointed $G$-functors
\begin{equation}\label{zd_twothree}
\begin{tikzpicture}[vcenter]
\def\v{-1.4}
\draw[0cell=1]
(0,0) node (a11) {\Ahnbe}
(a11)++(5,0) node (a12) {\phantom{\proAhnbe}}
(a12)++(0,.04) node (a12') {\proAhnbe}
(a11)++(0,\v) node (a21) {\Catg(\EG,\Ahnbe)}
(a12)++(0,\v) node (a22) {\Catg(\EG,\proAhnbe)}
;
\draw[1cell=.85]
(a11) edge node {\zd} (a12)
(a12') edge node {\inc} node[swap] {\gsim} (a22)
(a11) edge node {\gsim} node[swap] {\inc} (a21)
(a21) edge node {\Catg(\EG,\zd)} (a22)
;
\end{tikzpicture}
\end{equation}
in which the left inclusion $G$-functor $\inc$ is part of an adjoint $G$-equivalence by \cref{thm:HgoAh} \eqref{thm:HgoAh_i}.  However, while $\zd$ admits a left adjoint \cref{nbe_system_adj}, it is not generally an equivalence of categories \pcref{expl:Pist_not_eq,expl:system_adjunction}.  Thus, we \emph{cannot} use \cref{merling_2.16} to infer that $\Catg(\EG,\zd)$, and hence also $\zd$ and the $\nbe$-Segal functor $\segnbe$ for $\Adash$, are categorical weak $G$-equivalences.  This is where the proof of \cref{thm:h_special} would break down had we used $H$-theory instead of strong $H$-theory.
\end{explanation}

\subsection*{Specialness of Strong $J$-Theory}

For \cref{thm:j_special}, recall strong $J$-theory \pcref{thm:Jgo_twofunctor}
\[\AlgpspsO \fto{\Jgosg} \GGCatii.\]
With the $\Gcat$-operad $\Oph = \Catg(\EG,\Op)$ in place of $\Op$, this 2-functor sends each $\Oph$-pseudoalgebra to a $\Gskg$-category \pcref{A_ptfunctor}.  Recall that a $\Gskg$-category is \emph{special} if the $\angordn$-Segal functor $\segn$ is a categorical weak $K$-equivalence for each object $\angordn \in \Gsk \setminus \{\vstar,\ang{}\}$ and each graph subgroup $K \in \FGSin$ \pcref{def:sp_ggcat}.

\begin{theorem}\label{thm:j_special}\index{J-theory@$J$-theory!special}\index{special!G-G-category@$\Gskg$-category!strong $J$-theory}
Under \cref{as:uinf_ah} with $G$ a finite group, the strong $J$-theory of the $\Oph$-pseudoalgebra $\Ah$
\[\Gsk \fto{\Jgohsg \Ah = \Ahsgdash} \Gcatst\]
is a special $\Gskg$-category.
\end{theorem}

\begin{proof}
By \cref{jgohgoigst}, the strong $J$-theory of $\Ah$ factors through strong $H$-theory as
\begin{equation}\label{jgohsgah_fact}
\Jgohsg \Ah = \igst\Hgohsg\Ah.
\end{equation}
By \cref{gggcat_sp_iX} and the finiteness of $G$, the $\Gskg$-category $\Jgohsg \Ah$ is special if and only if the $\GGG$-category $\Hgohsg\Ah$ is special, which is true by \cref{thm:h_special}.
\end{proof}

Recall Shimakawa strong $J$-theory \cref{jgos_jgossg} and the 2-functor $\smast$ \cref{j_comparison}:
\[\AlgpspsO \fto{\Jgossg = \igst\Sgosg} \FGCat \fto{\smast} \GGCatii.\]
With the $\Gcat$-operad $\Oph = \Catg(\EG,\Op)$ in place of $\Op$, the preceding composite 2-functor sends each $\Oph$-pseudoalgebra to a $\Gskg$-category.

\begin{corollary}\label{cor:j_special_shi}\index{Shimakawa J-theory@Shimakawa $J$-theory!special}\index{special!G-G-category@$\Gskg$-category!Shimakawa strong $J$-theory}
Under \cref{as:uinf_ah} with $G$ a finite group, the $\Gskg$-category
\[\Gsk \fto{\smast\Jgohssg \Ah} \Gcatst\]
is special.
\end{corollary}

\begin{proof}
By \cref{thm:pistsgah_weq}, there is a weak $G$-equivalence
\begin{equation}\label{igst_smast_sgohsg_ah}
\igst\!\smast\Sgohsg\Ah \fto{\igst\Pistsg_{\Ah}} \igst\Hgohsg\Ah
\end{equation}
in $\GGCatii$ \pcref{def:gg_weq}.  By \cref{thm:j_special,jgohsgah_fact}, the $\Gskg$-category 
\[\Jgohsg \Ah = \igst\Hgohsg\Ah\]
is special.  By \cref{gg_weq_sp} \eqref{gg_weq_sp_ii}, \cref{igst_pistsg_dom}, and \cref{igst_smast_sgohsg_ah}, the $\Gskg$-category 
\[\igst\!\smast\Sgohsg\Ah = \smast\Jgohssg \Ah\]
is special.
\end{proof}

\begin{explanation}[Finiteness of $G$]\label{rk:hj_special}
\cref{thm:h_special,cor:h_special_shi} about (Shimakawa) strong $H$-theory hold for an arbitrary group $G$.  In contrast, \cref{thm:j_special,cor:j_special_shi} about (Shimakawa) strong $J$-theory hold for a \emph{finite} group $G$.  The reason for this difference is that the proof of \cref{thm:j_special} uses \cref{gggcat_sp_iX}, which requires $G$ to be finite.  See \cref{rk:shi_cor}.
\end{explanation}

%% file: chap/kgmmo.tex
This chapter compares
\begin{itemize}
\item the homotopical Shimakawa strong $K$-theory $\Khshosg$ \pcref{def:khsho} and
\item the equivariant $K$-theory functor $\Kgmmo$ of Guillou, May, Merling, and Osorno \cite{gmmo23}, which we call GMMO $K$-theory.
\end{itemize} 
For a finite group $G$ and a chaotic $\Einfg$-operad $\Op$ in $\Gcat$, each of these two equivariant $K$-theory machines sends $\Op$-algebras to orthogonal $G$-spectra.  \cref{thm:gmmo_shi} establishes a natural zigzag between $\Khshosg$ and $\Kgmmo$.  Moreover, for a chaotic $\Einfg$-operad of the form $\Oph = \Catg(\EG,\Op)$ and an $\Oph$-algebra $\Ah = \Catg(\EG,\A)$ with $\A$ an $\Op$-algebra, the zigzag comparison between $\Khshohsg\Ah$ and $\Kgmmoh\Ah$ consists of componentwise weak $G$-equivalences of orthogonal $G$-spectra.  In particular, for a naive permutative $G$-category $\A$ and the genuine permutative $G$-category $\Ah$, Shimakawa and GMMO $K$-theories produce weakly $G$-equivalent orthogonal $G$-spectra $\Khshgbesg\Ah$ and $\Kgmmogbe\Ah$, where $\GBE$ is the $G$-Barratt-Eccles operad.  See \cref{ex:gmmo_shi_gbe}.  Throughout this chapter, $G$ denotes a finite group.

\summary
The following diagram summarizes $\Khshosg$ and $\Kgmmo$ along the top and left-bottom-right boundaries.
\begin{equation}\label{Kgmmo_diag}
\begin{tikzpicture}[vcenter]
\def\v{1.4} \def\h{2.8} \def\u{-1.2} \def\w{1}
\draw[0cell]
(0,0) node (a1) {\AlgstpsO}
(a1)++(\h,0) node (b) {\AlgpspsO}
(a1)++(0,-\v) node (a2) {\DGCatg}
(a2)++(0,.2) node (a2') {\phantom{\DGCatg}}
(a2)++(0,-\v) node (a3) {\FGCatgps}
(a3)++(\h,0) node (a4) {\FGCatg}
(a4)++(\h,0) node (a5) {\FGTopg}
(a5)++(0,\v) node (a6) {\FGTopg}
(a6)++(0,\v) node (a7) {\Gspec}
;
\draw[1cell=.9]
(a1) edge[right hook->] node {\algi} (b)
(b) edge node {\Khshosg} (a7)
(b) edge node {\Sgosg} (a4)
(a1) edge node[swap] {\Rg} (a2)
(a2) edge node[swap] {\gzest} (a3)
(a3) edge[transform canvas={yshift=-.3ex}] node[swap] {\str} (a4)
(a4) edge[transform canvas={yshift=.6ex}, left hook->] node[swap] {\Incj} (a3)
(a4) edge node {\clast} (a5)
(a5) edge node {\Bc} (a6)
(a6) edge node {\Kfg} (a7)
(a4) edge[bend right=10,shorten >=.2ex] node[swap] {\gxist} (a2')
(a1) [rounded corners=2pt] -| ($(a2)+(\u,0)$) -- node[pos=0,swap] {\Kgmmo} ($(a3)+(\u,0)$) |- ($(a5)+(0,-.6)$) -| node[pos=0,swap] {\phantom{x}} ($(a5)+(\w,0)$) |- (a7)
;
\draw[2cell=.9]
node[between=a1 and a2 at .55, shift={(.4*\h,0)}, rotate=180, 2label={below,\cgs}] {\Rightarrow}
;
\end{tikzpicture}
\end{equation}
Each of $\Khshosg$ and $\Kgmmo$ ends with the composite functor 
\[\FGCatg \fto{\Kfg\Bc\clast} \Gspec.\]  
The main comparison between $\Khshosg$ and $\Kgmmo$ is the 2-natural transformation $\cgs$ \pcref{sec:cgs} that compares their first steps, namely, $\Sgosg$ and $\Rg$.  The comparison $\cgs$ occurs at the categorical level before passing to topological data.  For each $\Oph$-algebra of the form $\Ah = \Catg(\EG,\A)$, each component of $\cgs_{\Ah}$ is a categorical weak $G$-equivalence \pcref{def:cat_weakg}; see \cref{cgs_ah_wgeq}.  Along with the counit of the 2-adjunction $(\str,\Incj)$, these data yield the desired weak $G$-equivalence between $\Khshohsg\Ah$ and $\Kgmmoh\Ah$; see \cref{gmmo_shi_zig}.

\subtlety
Enriched pseudonaturality plays a key role in the comparison natural transformation $\cgs$.  At the object level, each of $\Rg$ and $\gxist$ yields pointed $\Gcatst$-functors $\DG \to \Gcatst$.  However, for a 1-connected $\Gcat$-operad $\Op$ and an $\Op$-algebra $\A$, the comparison \cref{cgsa} 
\[\gxist\Sgosg\algi\A \fto{\cgs_\A} \Rg\A\] 
is a $\Pig$-strict $\Gcatst$-pseudotransformation that is generally \emph{not} $\Gcatst$-natural.  Its lack of strict $\Gcatst$-naturality is controlled by the natural isomorphisms in \cref{cgsax}.  Thus, a crucial aspect of the comparison between $\Khshosg$ and $\Kgmmo$ is that GMMO $K$-theory goes through the 2-category $\DGCatg$, whose 1-cells are $\Pig$-strict $\Gcatst$-pseudotransformations instead of $\Gcatst$-natural transformations.  As we discuss in the introduction of \cref{ch:kgl_gmmo}, there is a similar subtlety in the equivalence between GMMO and global $K$-theories at a finite group $G$.

\organization
\crefrange{sec:dgcatg}{sec:kgmmo} review GMMO $K$-theory.  \cref{sec:cgs,sec:gmmo_shi_weq} compare $\Khshosg$ and $\Kgmmo$.

\secname{sec:dgcatg}
This section reviews the $\Gcatst$-category $\DG$ associated to a reduced $\Gcat$-operad $\Op$ and the 2-category $\DGCatg$.  The objects of $\DGCatg$ are pointed $\Gcatst$-functors 
\[\DG \to \Gcatst\] 
called $\DGG$-categories.  The subscript in the notation $\DGCatg$ refers to the 1-cells, which are $\Gcatst$-pseudotransformations that are strictly natural on a pointed $G$-subcategory $\Pig$.  The 2-cells of $\DGCatg$ are $\Gcatst$-modifications.

\secname{sec:rg}
This section reviews the first step $\Rg$ of $\Kgmmo$.  For a reduced $\Gcat$-operad $\Op$, $\AlgstpsO$ is the 2-category of $\Op$-algebras, $\Op$-pseudomorphisms, and $\Op$-transformations.  The 2-functor $\Rg$ sends them to $\DGG$-categories, $\Pig$-strict $\Gcatst$-pseudotransformations, and $\Gcatst$-modifications.  There is an inclusion 2-functor 
\[\AlgstpsO \fto{\algi} \AlgpspsO\] 
into the domain $\AlgpspsO$ \pcref{oalgps_twocat} of $\Khshosg$.

\secname{sec:fgcatgps}
This section reviews the 2-category $\FGCatgps$.  Its objects, called pseudo $\FGG$-categories, are $\Gcatst$-pseudofunctors 
\[\FG \to \Gcatst.\]  
Its 1-cells and 2-cells are $\Gcatst$-pseudotransformations and $\Gcatst$-modifications.  The 2-category $\FGCatg$ of $\FGG$-categories \pcref{def:fgcatg} is a non-full sub-2-category of $\FGCatgps$.

\secname{sec:gzest}
This section reviews the second step $\gzest$ of $\Kgmmo$.  The 2-functor $\gzest$ is given by pulling back along a strictly unital $\Gcatst$-pseudofunctor 
\[\FG \fto{\gze} \DG.\]
Here, $\DG$ is the $\Gcatst$-category from \cref{sec:dgcatg} associated to a chaotic $\Einfg$-operad $\Op$ \pcref{def:chaotic_einf}. 

\secname{sec:str}
This section reviews the third step $\str$ of $\Kgmmo$.  It is a left 2-adjoint of the inclusion 2-functor 
\[\FGCatg \fto{\Incj} \FGCatgps.\]

\secname{sec:kgmmo}
This section summarizes the functor $\Kgmmo$.  The functors $\Rg$, $\gzest$, $\str$, $\clast$, and $\Kfg\Bc$ preserve weak $G$-equivalences, and so does $\Kgmmo$ \pcref{kgmmo_weq}.

\secname{sec:cgs}
This section constructs a 2-natural transformation
\begin{equation}\label{cgs_chpi}
\begin{tikzpicture}[vcenter]
\def\t{30} \def\h{2.7} \def\v{1.4}
\draw[0cell]
(0,0) node (a1) {\AlgstpsO}
(a1)++(0,\v) node (a2) {\AlgpspsO}
(a2)++(\h,0) node (a3) {\FGCatg}
(a3)++(.2,0) node (a3') {\phantom{\FGCatg}}
(a3')++(0,-\v) node (a4) {\DGCatg}
;
\draw[1cell=.85]
(a1) edge[right hook->] node {\algi} (a2)
(a2) edge node {\Sgosg} (a3)
(a3') edge[transform canvas={xshift=-1ex}] node {\gxist} (a4)
(a1) edge node[swap] {\Rg} (a4)
;
\draw[2cell=1]
node[between=a2 and a3 at .42, shift={(0,-.5*\v)}, rotate=-90, 2label={above,\cgs}] {\Rightarrow}
;
\end{tikzpicture}
\end{equation}
that compares $\Sgosg$ and $\Rg$, which are the first steps of, respectively, $\Khshosg$ and $\Kgmmo$.

\secname{sec:gmmo_shi_weq}
This section proves that there is a natural zigzag connecting $\Khshosg$ and $\Kgmmo$ for each chaotic $\Einfg$-operad $\Op$.  Moreover, for a chaotic $\Einfg$-operad of the form $\Oph = \Catg(\EG,\Op)$ and an $\Oph$-algebra $\Ah = \Catg(\EG,\A)$ with $\A$ an $\Op$-algebra, \cref{thm:gmmo_shi} proves that there is a zigzag 
\[\Khshohsg\algi\Ah \fot[]{\eqg} \bdot \fto{\eqg} \Kgmmoh\Ah\]
of weak $G$-equivalences between orthogonal $G$-spectra.  \cref{ex:gmmo_shi_gbe} discusses this zigzag when $\Ah$ is a genuine permutative $G$-category.

\section{The 2-Category of $\DGG$-Categories}
\label{sec:dgcatg}

This section recalls the 2-category $\DGCatg$ from \cite[Sections 4--6]{gmmo23}.  For brief reviews of enriched categories and  2-categories, see \cref{sec:enrichedcat,sec:twocategories}. 

\secoutline
\begin{itemize}
\item \cref{def:dgo} defines the $\Gcatst$-category $\DG$ associated to each reduced $\Gcat$-operad $\Op$.
\item \cref{def:dgcatg} defines the objects of $\DGCatg$, called $\DGG$-categories, which are further elaborated in \cref{expl:dgcatg}.
\item \cref{def:dgcatg_onecell,def:Pig,def:DGcoop,def:pig_strict} define the 1-cells of $\DGCatg$, called $\Pig$-strict $\Gcatst$-pseudotransformations.
\item \cref{def:dgcatg_iicell} defines the 2-cells of $\DGCatg$, called $\Gcatst$-modifications.
\item \cref{def:dgcatg_iicat} defines the 2-category $\DGCatg$, which is further unpacked in \cref{expl:dgcatg_iicat}.
\item \cref{expl:gxist} unpacks the pullback functor $\gxist$ induced by the $\Gcatst$-functor $\gxi \cn \DG \to \FG$ \cref{dgcoop}.
\end{itemize}

\subsection*{The $\Gcatst$-Category $\DG$}

Recall the following.
\begin{itemize}
\item For $n \geq 0$, $\ordn$ is the pointed finite set $\{0,1,\ldots,n\}$ with basepoint $0$, and $\ufsn = \{1,2,\ldots,n\}$ is the unpointed finite set with $n$ elements \cref{ufsn}.
\item $\FG$ is the indexing $G$-category of pointed finite $G$-sets of the form $\mal$ for $m \geq 0$ \cref{ordn_be} and pointed morphisms with the conjugation $G$-action \pcref{def:FG,def:ctopgst}.  Its initial-terminal basepoint is the object $\ordz$.
\item $(\Gcatst,\sma,\bonep)$ is the complete and cocomplete symmetric monoidal closed category with small pointed $G$-categories as objects, pointed $G$-functors as morphisms, and $\Catgst$ as the internal hom \pcref{def:gcatst,expl:Gcatst}.
\item A $\Gcat$-operad $\Op$ is \emph{reduced} if $\Op(0)$ is the terminal $G$-category $\bone$.   We also denote $\Op(n)$ by $\Op_n$.\label{not:Opn}
\end{itemize}

\begin{definition}\label{def:dgo}
For a reduced $\Gcat$-operad $(\Op,\ga,\opu)$, the $\Gcatst$-category $\DG$ is defined as follows.  Note that $\Op$ is suppressed from the notation $\DG$ to simplify the typography.
\begin{description}
\item[Objects]
$\DG$ has the same objects as $\FG$, namely, pointed finite $G$-sets of the form $\mal$ for $m \geq 0$.
\item[Hom pointed $G$-categories]
For pointed finite $G$-sets $\mal$ and $\nbeta$, the pointed $G$-category $\DG(\mal,\nbeta)$ has underlying pointed category given by the coproduct of product categories
\begin{equation}\label{dg_mn}
\DG(\mal,\nbeta) = \coprod_{\psi \in \FG(\mal,\,\nbeta)} \,\prod_{j \in \ufsn}\, \Op_{|\psiinv j|}.
\end{equation}
An object or a morphism in $\DG(\mal,\nbeta)$ is denoted by
\begin{equation}\label{dgmn_obj}
\objx = \big(\psi; \ang{x_j}_{j \in \ufsn}\big).
\end{equation}
\begin{itemize}
\item $\psi \in \FG(\mal,\nbeta)$ is a pointed morphism $\mal \to \nbeta$.
\item Depending on whether $\objx$ is an object or a morphism, each $x_j \in \Op_{|\psiinv j|}$ is an object or a morphism.
\end{itemize}
\begin{description}
\item[Basepoint] 
The basepoint of $\DG(\mal,\nbeta)$ is the object 
\begin{equation}\label{dgmn_basept}
\objzero = (0; \ang{*}_{j \in \ufsn})
\end{equation} 
consisting of
\begin{itemize}
\item the 0-morphism $0 \cn \mal \to \nbeta$ and
\item $n$ copies of the unique object $* \in \Op(0) = \bone$.
\end{itemize}
\item[Identity 2-cells]
In \cref{dgmn_obj}, if each $x_j \in \Op_{|\psiinv j|}$ is an object, then the identity morphism of the object $\objx = (\psi; \ang{x_j}_{j \in \ufsn})$ is 
\begin{equation}\label{dg_idtwocell}
\objone_{\objx} = (\psi; \ang{1_{x_j}}_{j \in \ufsn}).
\end{equation}
\item[$G$-action]
For an element $g \in G$, the $g$-action on $\DG(\mal,\nbeta)$ is given by
\begin{equation}\label{dgmn_gact}
g \big(\psi; \ang{x_j}_{j \in \ufsn}\big) 
= \big(g\psi\ginv; \ang{g x_{\ginv j} \tau_{\ginv}}_{j \in \ufsn} \big).
\end{equation}
\begin{itemize}
\item $g\psi\ginv \cn \mal \to \nbeta$ is the conjugation $g$-action on $\psi$ \cref{gpsi}.
\item $gx_{\ginv j}$ is the $g$-action on $x_{\ginv j} \in \Op_{|\psiinv\ginv j|}$.
\item $gx_{\ginv j} \tau_{\ginv}$ is the right symmetric group action on $gx_{\ginv j} \in \Op_{|\psiinv\ginv j|}$ for the permutation
\begin{equation}\label{tauginv}
\ufsm \supset (g\psi\ginv)^{-1}j \fto[\iso]{\tau_{\ginv}} (g\psi)^{-1}j \subset \ufsm
\end{equation}
given by the restriction of the $\ginv$-action on $\mal$, using the equalities 
\[|\psiinv\ginv j| = |(g\psi)^{-1}j| = |(g\psi\ginv)^{-1}j|.\]
\end{itemize}
\end{description}
\item[Identity 1-cells]
The identity 1-cell of a pointed finite $G$-set $\mal \in \DG$ is the object 
\begin{equation}\label{dg_idonecell}
\objone_{\mal} = \big(1_{\mal}; \ang{\opu}_{i \in \ufsm}\big) \in \DG(\mal,\mal),
\end{equation}
where $\opu \in \Op(1)$ is the operadic unit.
\item[Composition]
For pointed finite $G$-sets $\ldea$, $\mal$, and $\nbeta$, the composition pointed $G$-functor
\[\DG(\mal,\nbeta) \sma \DG(\ldea,\mal) \fto{\comp} \DG(\ldea,\nbeta)\]
is defined by
\begin{equation}\label{dg_comp}
\big(\psi; \ang{x_j}_{j \in \ufsn}\big) \comp \big(\phi; \ang{v_i}_{i \in \ufsm} \big) 
= \big(\psi\phi; \ang{\ga(x_j; \ang{v_i}_{i \in \psiinv j}) \tau^j_{\psi,\phi}}_{j \in \ufsn} \big).
\end{equation}
For each $j \in \ufsn$, the shuffle
\begin{equation}\label{tauj_psiphi}
(\psi\phi)^{-1}(j) \fto[\iso]{\tau^j_{\psi,\phi}} \coprod_{i \in \psiinv j} \phiinv i
\end{equation}
is determined by the inherited orderings of 
\[(\psi\phi)^{-1}(j), \phiinv i \subseteq \ufsl \andspace \psiinv j \subseteq \ufsm.\]
\end{description}
This finishes the definition of $\DG$.
\end{definition}

\begin{explanation}\label{expl:dgo}
Consider \cref{def:dgo}.
\begin{enumerate}
\item The $\Catst$-category axioms of $\DG$ follow from the reduced $\Cat$-operad axioms of $\Op$. A detailed proof is given in \cite[Prop.\ 5.3]{johnson-yau-permmult}. 
\item The $G$-equivariance of the identity 1-cells \cref{dg_idonecell} follows from the fact that the operadic unit $\opu \in \Op(1)$ is $G$-fixed.
\item The $G$-equivariance of the composition \cref{dg_comp} follows from
\begin{itemize}
\item the $G$-equivariance of the operadic composition $\ga$ and the right symmetric group action on $\Op$, 
\item the symmetric group action equivariance axioms of $\Op$, and
\item the uniqueness of permutations.
\end{itemize}
\item The object $\ordz = \{0\} \in \DG$ is a \emph{0-object} in the sense that $\DG(\ordz,\nbeta)$ and $\DG(\mal,\ordz)$ are terminal $G$-categories.  Indeed, each of $\FG(\ordz,\nbeta)$ and $\FG(\mal,\ordz)$ consists of only the 0-morphism.  If $m = 0$, then we use the fact that $\Op$ is reduced, meaning $\Op(0) = \bone$.  If $\ufsn = \ufs{0} = \emptyset$, then an empty product is, by definition, the terminal $G$-category $\bone$.\defmark
\end{enumerate}
\end{explanation}

\subsection*{$\DGG$-Categories}

For the rest of this section, $\Op$ is a reduced $\Gcat$-operad, and $\DG$ is the $\Gcatst$-category in \cref{def:dgo}.  

\begin{definition}\label{def:dgcatg}
A \emph{$\DGG$-category}\index{DGG-category@$\DGG$-category} is a pointed $\Gcatst$-functor
\begin{equation}\label{dgcatg_obj}
(\DG,\ordz) \fto{X} (\Gcatst,\bone),
\end{equation}
where pointedness means $X\ordz = \bone$.
\end{definition}

\begin{explanation}[Unpacking]\label{expl:dgcatg}
A $\DGG$-category $X \cn \DG \to \Gcatst$ \cref{dgcatg_obj} consists of the following data.
\begin{description}
\item[Pointed $G$-categories] 
$X$ sends each pointed finite $G$-set $\mal$ to a small pointed $G$-category $X\mal$ such that $X\ordz = \bone$.
\item[Pointed $G$-functors] 
For each pair $(\mal,\nbeta)$ of pointed finite $G$-sets, $X$ is equipped with an $(\mal,\nbeta)$-component pointed $G$-functor
\begin{equation}\label{dgcatg_xmn}
\DG(\mal,\nbeta) \fto{X_{\mal,\,\nbeta}} \Catgst(X\mal, X\nbeta)
\end{equation}
that preserves identity 1-cells \cref{dg_idonecell} and composition of $\DG$ \cref{dg_comp}.  We often abbreviate $X_{\mal,\,\nbeta}$ to $X$.   We further unpack the pointed $G$-functor $X_{\mal,\, \nbeta}$.
\begin{description}
\item[Objects] 
$X_{\mal,\, \nbeta}$ sends each object $\objx = (\psi; \ang{x_j}_{j \in \ufsn})$ in $\DG(\mal,\nbeta)$ \cref{dgmn_obj} to a pointed functor 
\begin{equation}\label{xmn_psix}
X\mal \fto{X\objx} X\nbeta.
\end{equation}
Note that the pointed functor $X\objx$ is \emph{not} required to be $G$-equivariant. 
\item[Pointedness]
There is an equality of functors
\begin{equation}\label{xzero_basept}
X\mal \fto{X\objzero = *} X\nbeta
\end{equation}
for the basepoints $\objzero \in \DG(\mal,\nbeta)$ \cref{dgmn_basept} and $* \in X\nbeta$.
\item[Morphisms]
$X_{\mal,\, \nbeta}$ sends each morphism 
\begin{equation}\label{dgmn_morphism}
\objx = \big(\psi; \ang{x_j}_{j \in \ufsn}\big)
\fto{\objf = (\psi; \ang{f_j}_{j \in \ufsn})}
\objy = \big(\psi; \ang{y_j}_{j \in \ufsn}\big)
\end{equation}
in $\DG(\mal,\nbeta)$ \cref{dg_mn} to a pointed natural transformation 
\begin{equation}\label{xobjf}
\begin{tikzpicture}[vcenter]
\def\t{28}
\draw[0cell]
(0,0) node (a1) {\phantom{A}}
(a1)++(1.8,0) node (a2) {\phantom{A}}
(a1)++(-.2,0) node (a1') {X\mal}
(a2)++(.2,.02) node (a2') {X\nbeta}
;
\draw[1cell=.9]
(a1) edge[bend left=\t] node {X\objx} (a2)
(a1) edge[bend right=\t] node[swap] {X\objy} (a2)
;
\draw[2cell=.9]
node[between=a1 and a2 at .38, rotate=-90, 2label={above,X\objf}] {\Rightarrow}
;
\end{tikzpicture}
\end{equation}
such that $X$ preserves identity 2-cells \cref{dg_idtwocell} and composition of morphisms in $\DG(\mal,\nbeta)$ \cref{dg_mn}.
\item[$G$-Equivariance] 
By \cref{conjugation-gaction,dgmn_gact}, the $G$-equivariance of $X_{\mal,\, \nbeta}$ means that, for $g \in G$ and $\objx = (\psi; \ang{x_j}_{j \in \ufsn})$ in $\DG(\mal,\nbeta)$, the following diagram commutes.
\begin{equation}\label{xmn_geq}
\begin{tikzpicture}[vcenter]
\def\v{-1.3}
\draw[0cell]
(0,0) node (a11) {X\mal}
(a11)++(2.3,0) node (a12) {X\nbeta}
(a11)++(0,\v) node (a21) {X\mal}
(a12)++(0,\v) node (a22) {X\nbeta}
;
\draw[1cell=.9]
(a11) edge node {X(g\objx)} (a12)
(a11) edge node[swap] {\ginv} (a21)
(a21) edge node {X\objx} (a22)
(a22) edge node[swap] {g} (a12)
;
\end{tikzpicture}
\end{equation}
Depending on whether $\objx$ is an object or a morphism, \cref{xmn_geq} is a diagram of pointed functors or pointed natural transformations.
\item[Preservation of identity 1-cells]
There is an equality of functors
\begin{equation}\label{xmm_one}
X\mal \fto{X\objone_{\mal} = 1_{X\mal}} X\mal
\end{equation}
for each identity 1-cell $\objone_{\mal} \in \DG(\mal,\mal)$ \cref{dg_idonecell}.
\item[Preservation of composition]
For $\objx = (\psi; \ang{x_j}_{j \in \ufsn})$ in $\DG(\mal,\nbeta)$ and $\objv = (\phi; \ang{v_i}_{i \in \ufsm})$ in $\DG(\ldea,\mal)$ as displayed in \cref{dg_comp}, the diagram
\begin{equation}\label{xlmn}
\begin{tikzpicture}[vcenter]
\def\h{2} \def\u{.7}
\draw[0cell]
(0,0) node (a1) {X\ldea}
(a1)++(\h,0) node (a2) {\phantom{X\mal}}
(a2)++(0,-.03) node (a2') {X\mal}
(a2)++(\h,0) node (a3) {X\nbeta}
;
\draw[1cell=.9]
(a1) edge node {X\objv} (a2)
(a2) edge node {X\objx} (a3)
(a1) [rounded corners=2pt] |- ($(a2)+(-1,\u)$) -- node {X(\objx\comp\objv)} ($(a2)+(1,\u)$) -| (a3)
;
\end{tikzpicture}
\end{equation}
commutes.  Depending on whether $\objx$ and $\objv$ are objects or morphisms, \cref{xlmn} is a diagram of pointed functors or pointed natural transformations.
\end{description}
\end{description}
The adjoint of $X_{\mal,\,\nbeta}$ is a pointed $G$-functor
\begin{equation}\label{dgcatg_xmn_adj}
\DG(\mal,\nbeta) \sma X\mal \fto{X_{\mal,\,\nbeta}} X\nbeta
\end{equation}
that we denote by the same notation.
\end{explanation}

\subsection*{$\Gcatst$-Pseudotransformations}

\begin{definition}\label{def:dgcatg_onecell}
Suppose $X,X' \cn \DG \to \Gcatst$ are pointed $\Gcatst$-functors \cref{dgcatg_obj}.  A \emph{$\Gcatst$-pseudotransformation}\index{Gcatst-pseudotransformation@$\Gcatst$-pseudotransformation} $\tha \cn X \to X'$ consists of the following data.
\begin{description}
\item[$G$-functors] For each object $\mal \in \DG$, $\tha$ is equipped with an $\mal$-component pointed $G$-functor
\begin{equation}\label{dgcatg_one_comp}
X\mal \fto{\tha_{\mal}} X'\mal.
\end{equation}
\item[Natural isomorphisms] For each pair $(\mal,\nbeta)$ of objects in $\DG$ and each object $\objx \in \DG(\mal,\nbeta)$ \cref{dgmn_obj}, $\tha$ is equipped with an $\objx$-component pointed natural isomorphism $\tha_{\objx}$ as follows, where $X\objx$ and $X'\objx$ are the pointed functors in \cref{xmn_psix}.
\begin{equation}\label{dgcatg_one_compii}

\end{equation}
\end{description}
These data are required to satisfy the axioms \crefrange{thaobjx_nat}{tha_xv}.
\begin{description}
\item[Naturality] $\tha_{\objx}$ is natural in $\objx$.  This means that, for each morphism $\objf \cn \objx \to \objy$ in $\DG(\mal,\nbeta)$ \cref{dgmn_morphism}, the following two pasting diagrams of natural transformations are equal.
\begin{equation}\label{thaobjx_nat}
%
\end{equation}
\item[$G$-equivariance]
The assignment $\objx \mapsto \tha_{\objx}$ is $G$-equivariant, meaning
\begin{equation}\label{thaobjx_gequiv}
\tha_{g\objx} = g \tha_{\objx} \ginv
\end{equation}
for $g \in G$ and $\objx \in \DG(\mal,\nbeta)$. 
\begin{itemize}
\item The object $g\objx$ is the $G$-action on $\objx$ \cref{dgmn_gact}. 
\item The right-hand side is the conjugation $G$-action on $\tha_{\objx}$ \cref{conjugation-gaction}, using the $\ginv$-action on $X\mal$ and the $g$-action on $X'\nbeta$. 
\end{itemize} 
The axiom \cref{thaobjx_gequiv} is well defined by \cref{xmn_geq} and the $G$-equivariance of $\tha_{\mal}$ and $\tha_{\nbeta}$ \cref{dgcatg_one_comp}.
\item[Basepoint]
The assignment $\objx \mapsto \tha_{\objx}$ preserves basepoints, meaning 
\begin{equation}\label{tha_objzero}
\tha_{\objzero} = 1_*
\end{equation}
for the basepoints $\objzero \in \DG(\mal,\nbeta)$ \cref{dgmn_basept} and $* \in X'\nbeta$.  This axiom is well defined by \cref{xzero_basept} and the pointedness of $\tha_{\nbeta}$.
\item[Unity]
The assignment $\objx \mapsto \tha_{\objx}$ preserves identity 1-cells, meaning
\begin{equation}\label{tha_onex}
\tha_{\objone_{\mal}} = 1_{\tha_{\mal}}
\end{equation}
for the identity 1-cell $\objone_{\mal} \in \DG(\mal,\mal)$ \cref{dg_idonecell}.  This axiom is well defined by \cref{xmm_one}.
\item[Compositionality]
Using the notation in \cref{xlmn}, the following two pasting diagrams of natural transformations are equal, where $\objx\objv$ means $\objx \comp \objv$.
\begin{equation}\label{tha_xv}
\begin{tikzpicture}[vcenter]
\def\h{2} \def\v{-1.3} \def\s{.9} \def\q{.8} \def\p{.85} \def\r{-120}
\def\boundary{
\draw[0cell=\s]
(0,0) node (a11) {X\ldea}
(a11)++(\h,0) node (a12) {X'\ldea}
(a11)++(0,2*\v) node (a31) {X\nbeta}
(a12)++(0,2*\v) node (a32) {X'\nbeta}
;
\draw[1cell=\q]
(a11) edge node {\tha_{\ldea}} (a12)
(a31) edge node[swap] {\tha_{\nbeta}} (a32)
;}
\boundary
\draw[0cell=\s]
(a11)++(0,\v) node (a21) {X\mal}
(a12)++(0,\v) node (a22) {X'\mal}
;
\draw[1cell=\q]
(a21) edge node {\tha_{\mal}} (a22)
(a11) edge node[swap] {X\objv} (a21)
(a21) edge node[swap] {X\objx} (a31)
(a12) edge node {X'\objv} (a22)
(a22) edge node {X'\objx} (a32)
;
\draw[2cell=\p]
node[between=a12 and a21 at .5, shift={(0,.1)}, rotate=\r, 2labelalt={below,\tha_{\objv}}] {\Rightarrow}
node[between=a22 and a31 at .5, shift={(0,-.15)}, rotate=\r, 2labelalt={below,\tha_{\objx}}] {\Rightarrow}
;
\begin{scope}[shift={(2*\h,0)}]
\boundary
\draw[1cell=\q]
(a11) edge node[swap] {X(\objx\objv)} (a31)
(a12) edge node {X'(\objx\objv)} (a32)
;
\draw[2cell=\p]
node[between=a12 and a31 at .5, shift={(0,0)}, rotate=\r, 2labelalt={below,\tha_{\objx\objv}}] {\Rightarrow}
;
\end{scope}
\end{tikzpicture}
\end{equation}
This axiom is well defined by \cref{xlmn}.
\end{description}
A $\Gcatst$-pseudotransformation $\tha$ is called a \emph{weak $G$-equivalence}\index{Gcatst-pseudotransformation@$\Gcatst$-pseudotransformation!weak G-equivalence@weak $G$-equivalence}\index{weak G-equivalence@weak $G$-equivalence!Gcatst-pseudotransformation@$\Gcatst$-pseudotransformation} if each $\tha_{\mal}$ \cref{dgcatg_one_comp} is a categorical weak $G$-equivalence \pcref{def:cat_weakg}.  Note that a \emph{$\Gcatst$-natural transformation} $\tha \cn X \to X'$ yields a $\Gcatst$-pseudotransformation with each $\tha_{\objx}$ \cref{dgcatg_one_compii} given by the identity natural transformation and vice versa.
\end{definition}

\subsection*{$\Pig$-Strictness}

The 1-cells in the 2-category $\DGCatg$ are $\Gcatst$-pseudotransformations that are strict on a $G$-subcategory $\Pig \subseteq \FG$ defined next. 

\begin{definition}\label{def:Pig}
Denote by $\Pig$ the pointed $G$-subcategory of $\FG$ consisting of all the objects and pointed morphisms $\psi \cn \mal \to \nbeta$ such that
\[|\psiinv j| \in \{0,1\} \forspace j \in \ufsn.\]
Each hom set in each of $\Pig$ and $\FG$ is a pointed finite $G$-set with the conjugation $G$-action \cref{gpsi} and basepoint given by the 0-morphism.  Each such pointed $G$-set is regarded as a small pointed $G$-category with only identity morphisms.  The pointed $G$-categories $\Pig$ and $\FG$ are regarded as $\Gcatst$-categories with only identity 2-cells.
\end{definition}

\begin{definition}\label{def:DGcoop}
For the $\Gcatst$-categories $\FG$, $\DG$, and $\Pig$ in \cref{def:FG,def:dgo,def:Pig}, the $\Gcatst$-functors
\begin{equation}\label{dgcoop}
\Pig \fto{\gio} \DG \fto{\gxi} \FG
\end{equation}
are defined as the identity on objects.  On hom pointed $G$-categories, they are defined as follows.
\begin{itemize}
\item $\gxi$ sends $(\psi; \ang{x_j}_{j \in \ufsn}) \in \DG(\mal,\nbeta)$ to $\psi \in \FG(\mal,\nbeta)$.
\item $\gio$ sends a pointed morphism $\psi \cn \mal \to \nbeta$ in $\Pig$ to the object
\begin{equation}\label{iota_psi}
\gio(\psi) = \big(\psi; \ang{x_j}_{j \in \ufsn}\big) \in \DG(\mal,\nbeta)
\end{equation}
such that 
\[x_j = \begin{cases} 
* \in \Op(0) & \text{if $|\psiinv j| = 0$ and}\\
\opu \in \Op(1) & \text{if $|\psiinv j| = 1$.}
\end{cases}\]
\end{itemize}
Note that the composite $\gxi\gio \cn \Pig \to \FG$ is the inclusion.
\end{definition}

\begin{example}\label{ex:DGFG}
Suppose $\Op$ is the terminal $\Gcat$-operad with $\Op(n) = \bone$ for each $n \geq 0$.  Then $\DG$ is isomorphic to $\FG$ via $\gxi$.
\end{example}

\begin{definition}\label{def:pig_strict}\index{Gcatst-pseudotransformation@$\Gcatst$-pseudotransformation!Pig-strict@$\Pig$-strict}\index{Pig-strict@$\Pig$-strict}
A $\Gcatst$-pseudotransformation $\tha \cn X \to X'$ between $\DGG$-categories \pcref{def:dgcatg,def:dgcatg_onecell} is \emph{$\Pig$-strict} if $\tha_{\objx}$ \cref{dgcatg_one_compii} is the identity natural transformation for each $\objx \in \DG(\mal,\nbeta)$ in the image of $\gio \cn \Pig \to \DG$ \cref{dgcoop}.
\end{definition}

\subsection*{2-Category Structure}

\begin{definition}\label{def:dgcatg_iicell}
Suppose $\tha, \ups \cn X \to X'$ are $\Gcatst$-pseudotransformations between $\DGG$-categories \pcref{def:dgcatg,def:dgcatg_onecell}.  A \emph{$\Gcatst$-modification}\index{Gcatst-modification@$\Gcatst$-modification}\index{modification!Gcatst@$\Gcatst$} $\Theta \cn \tha \to \ups$ consists of, for each object $\mal \in \DG$, an $\mal$-component pointed $G$-natural transformation
\begin{equation}\label{dgcatg_iicell_comp}
\begin{tikzpicture}[vcenter]
\def\t{30}
\draw[0cell]
(0,0) node (a1) {\phantom{A}}
(a1)++(2,0) node (a2) {\phantom{A}}
(a1)++(-.2,0) node (a1') {X\mal}
(a2)++(.25,0) node (a2') {X'\mal}
;
\draw[1cell=.9]
(a1) edge[bend left=\t] node {\tha_{\mal}} (a2)
(a1) edge[bend right=\t] node[swap] {\ups_{\mal}} (a2)
;
\draw[2cell=.9]
node[between=a1 and a2 at .38, rotate=-90, 2label={above,\Theta_{\mal}}] {\Rightarrow}
;
\end{tikzpicture}
\end{equation}
such that, for each object $\objx \in \DG(\mal,\nbeta)$ \cref{dgmn_obj}, the following two pasting diagrams of natural transformations are equal.
\begin{equation}\label{dgcatg_iicell_modax}
\begin{tikzpicture}[vcenter]
\def\h{2} \def\v{-1.4} \def\u{1.4} \def\s{.8} \def\c{.4} \def\d{.6} \def\t{23} 
\def\boundary{
\draw[0cell=.9]
(0,0) node (a11) {X\mal}
(a11)++(\h,0) node (a12) {X'\mal}
(a11)++(0,\v) node (a21) {X\nbeta}
(a12)++(0,\v) node (a22) {X'\nbeta}
;
\draw[1cell=\s]
(a11) edge node[swap] {X\objx} (a21)
(a12) edge node {X'\objx} (a22)
(a11) edge[bend left=\t] node[pos=\c] {\tha_{\mal}} (a12)
(a21) edge[bend right=\t] node[swap,pos=\d] {\ups_{\nbeta}} (a22)
;}
\boundary
\draw[1cell=\s]
(a11) edge[bend right=\t] node[swap,pos=\d] {\ups_{\mal}} (a12)
;
\draw[2cell=\s]
node[between=a11 and a12 at .4, shift={(0,0)}, rotate=-90, 2label={above,\Theta_{\mal}}] {\Rightarrow}
node[between=a11 and a12 at .4, shift={(0,-\u)}, rotate=-90, 2label={above,\ups_{\objx}}] {\Rightarrow}
;
\begin{scope}[shift={(2*\h,0)}]
\boundary
\draw[1cell=\s]
(a21) edge[bend left=\t] node[pos=\c] {\tha_{\nbeta}} (a22)
;
\draw[2cell=\s]
node[between=a21 and a22 at .4, shift={(0,0)}, rotate=-90, 2label={above,\Theta_{\nbeta}}] {\Rightarrow}
node[between=a21 and a22 at .4, shift={(0,\u)}, rotate=-90, 2label={above,\tha_{\objx}}] {\Rightarrow}
;
\end{scope}
\end{tikzpicture}
\end{equation}
Note that $\Theta_{\mal}$ is not required to be invertible.
\end{definition}

\begin{definition}\label{def:dgcatg_iicat}
Suppose $\Op$ is a reduced $\Gcat$-operad, and $\DG$ is the $\Gcatst$-category in \cref{def:dgo}. 
\begin{enumerate}
\item\label{def:dgcatg_iicat_i} Define the 2-category $\DGCatg$ with
\begin{itemize}
\item $\DGG$-categories \pcref{def:dgcatg} as objects,
\item $\Pig$-strict $\Gcatst$-pseudotransformations \pcref{def:dgcatg_onecell,def:pig_strict} as 1-cells, and
\item $\Gcatst$-modifications \pcref{def:dgcatg_iicell} as 2-cells.
\end{itemize}
Other 2-categorical structures---identity and horizontal composition of 1-cells and 2-cells and vertical composition of 2-cells---are defined componentwise; see \cref{expl:dgcatg_iicat}.  The underlying 1-category of $\DGCatg$ is denoted by the same notation.
\item\label{def:dgcatg_iicat_ii} Define the sub-2-category 
\[\begin{tikzpicture}[vcenter]
\draw[0cell]
(0,0) node (a1) {\dgcatg}
(a1)++(2.5,0) node (a2) {\DGCatg}
;
\draw[1cell=.9]
(a1) edge[right hook->] node {\dgi} (a2)
;
\end{tikzpicture}\]
with
\begin{itemize}
\item the same objects and 2-cells as $\DGCatg$ and 
\item $\Gcatst$-natural transformations as 1-cells.
\end{itemize}
\item\label{def:dgcatg_iicat_iii} Denote by
\begin{equation}\label{gxist}
\FGCatg \fto{\gxist} \dgcatg
\end{equation} 
the pullback 2-functor given by precomposing with the $\Gcatst$-functor $\gxi \cn \DG \to \FG$ \cref{dgcoop}, where $\FGCatg$ is the 2-category of $\FGG$-categories \pcref{def:fgcatg}.  See \cref{expl:gxist}.  The composite 
\[\FGCatg \fto{\gxist} \dgcatg \fto{\dgi} \DGCatg\]
is also denoted by $\gxist$.  \defmark
\end{enumerate}
\end{definition}

\begin{explanation}[Unpacking $\DGCatg$]\label{expl:dgcatg_iicat}
Compositions and identities in the 2-category $\DGCatg$ are given explicitly as follows.
\begin{description}
\item[Identity 1-cells]
For a pointed $\Gcatst$-functor $X \cn \DG \to \Catgst$ \cref{dgcatg_obj}, the identity 1-cell $1_X$ consists of
\begin{itemize}
\item the identity functor $1_{X\mal}$ for each object $\mal \in \DG$ and
\item the identity natural transformation $1_{X\objx}$ for each object $\objx \in \DG(\mal,\nbeta)$.
\end{itemize}
\item[Horizontal composition of 1-cells]
For pointed $\Gcatst$-functors $X,X',X'' \cn \DG \to \Catgst$ \cref{dgcatg_obj}, the horizontal composition of $\Gcatst$-pseudotransformations \pcref{def:dgcatg_onecell}
\[X \fto{\tha} X' \fto{\tha'} X'',\]
denoted by $\tha'\tha \cn X \to X''$, has $\mal$-component pointed $G$-functor \cref{dgcatg_one_comp} given by the composite
\[X\mal \fto{\tha_{\mal}} X'\mal \fto{\tha'_{\mal}} X''\mal\]
for each object $\mal \in \DG$.  For each object $\objx \in \DG(\mal,\nbeta)$, the $\objx$-component pointed natural isomorphism $(\tha'\tha)_{\objx}$ \cref{dgcatg_one_compii} is given by the following pasting diagram.
\begin{equation}\label{dgcatg_icell_paste}
\begin{tikzpicture}[vcenter]
\def\h{2.3}\def\v{-1.4} \def\u{-.1} \def\r{240}
\draw[0cell]
(0,0) node (a11) {X\mal}
(a11)++(\h,0) node (a12) {X'\mal}
(a12)++(\h,0) node (a13) {X''\mal}
(a11)++(0,\v) node (a21) {X\nbeta}
(a12)++(0,\v) node (a22) {X'\nbeta}
(a13)++(0,\v) node (a23) {X''\nbeta}
;
\draw[1cell=.9]
(a11) edge node {\tha_{\mal}} (a12)
(a12) edge node {\tha'_{\mal}} (a13)
(a21) edge node[swap] {\tha_{\nbeta}} (a22)
(a22) edge node[swap] {\tha'_{\nbeta}} (a23)
(a11) edge node[swap] {X\objx} (a21)
(a12) edge node[swap] {X'\objx} (a22)
(a13) edge node {X''\objx} (a23)
;
\draw[2cell]
node[between=a12 and a21 at .5, shift={(0,\u)}, rotate=\r, 2labelmed={below,\tha_{\objx}}] {\Rightarrow}
node[between=a13 and a22 at .5, shift={(0,\u)}, rotate=\r, 2labelmed={below,\tha'_{\objx}}] {\Rightarrow}
;
\end{tikzpicture}
\end{equation}
If $\tha$ and $\tha'$ are both $\Pig$-strict \pcref{def:pig_strict}, then so is the horizontal composite $\tha'\tha$.  
\item[Identity 2-cells]
For a $\Pig$-strict $\Gcatst$-pseudotransformation $\tha \cn X \to X'$, the identity 2-cell $1_{\tha}$ consists of the identity natural transformation $1_{\tha_{\mal}}$ for each object $\mal \in \DG$.
\item[Compositions of 2-cells]  
Vertical and horizontal compositions of $\Gcatst$-modifications are defined, for each object $\mal \in \DG$, as the corresponding compositions of the $\mal$-component $G$-natural transformations \cref{dgcatg_iicell_comp}.
\end{description}
The subscript in the notation $\DGCatg$ is a reminder that the 1-cells---$\Pig$-strict $\Gcatst$-pseudotransformations---are not strictly natural.
\end{explanation}

\begin{explanation}[Unpacking $\gxist$]\label{expl:gxist}
Using \cref{expl:dgcatg,expl:fgcatg_iicat,dgcoop}, the 2-functor \cref{gxist}
\[\FGCatg \fto{\gxist} \dgcatg\] 
is given as follows.
\begin{description}
\item[Objects]
For an $\FGG$-category $X \cn \FG \to \Catgst$ \cref{fgcatg_obj}, the $\DGG$-category 
\[(\DG,\ordz) \fto{\gxist X} (\Gcatst,\bone)\]
is given by
\begin{itemize}
\item the pointed $G$-category 
\[(\gxist X)\mal = X\mal\]
for each pointed finite $G$-set $\mal$;
\item the pointed functor
\[X\mal \fto{(\gxist X)\objx \,=\, X\psi} X\nbeta\]
for each object $\objx = (\psi; \ang{x_j}_{j \in \ufsn})$ in $\DG(\mal,\nbeta)$ \cref{dgmn_obj}; and
\item the identity natural transformation
\[(\gxist X)\objf = 1_{X\psi}\]
for each morphism $\objf = (\psi; \ang{f_j}_{j \in \ufsn})$ in $\DG(\mal,\nbeta)$.
\end{itemize}
\item[1-cells] 
For a 1-cell $\tha \cn X \to X'$ in $\FGCatg$ \cref{fgcatg_icell}, the $\Gcatst$-natural transformation 
\[\gxist X \fto{\gxist\tha} \gxist X'\] 
is given by the pointed $G$-functor
\[X\mal \fto{(\gxist \tha)_{\mal} \,=\, \tha_{\mal}} X'\mal\]
for each pointed finite $G$-set $\mal$.
\item[2-cells]
For a 2-cell $\Theta \cn \tha \to \ups$ in $\FGCatg$ \cref{fgcatg_iicell}, the $\Gcatst$-modification 
\[\gxist\tha \fto{\gxist\Theta} \gxist\ups\] 
is given by the pointed $G$-natural transformation
\[\tha_{\mal} \fto{(\gxist\Theta)_{\mal} \,=\, \Theta_{\mal}} \ups_{\mal}\]
for each pointed finite $G$-set $\mal$.\defmark
\end{description}
\end{explanation}

\section{The 2-Functor $\Rg$}
\label{sec:rg}

This section recalls the 2-functor
\[\AlgstpsO \fto{\Rg} \DGCatg\]
for a reduced $\Gcat$-operad $\Op$, where $\DGCatg$ is the 2-category in \cref{def:dgcatg_iicat}.  The 2-functor $\Rg$ is the first step of GMMO $K$-theory \cref{Kgmmo_diag}, sending $\Op$-algebras to $\DGG$-categories at the object level.  The definitions in this section correspond to those in \cite[Sections 5 and 6]{gmmo23} without multifunctoriality.

\secoutline
\begin{itemize}
\item \cref{def:algstpsO} defines the 2-category $\AlgstpsO$.
\item \cref{def:rg_obj,def:rg_onecell,def:rg_twocell} define the object, 1-cell, and 2-cell assignments of $\Rg$, which are further elaborated in \cref{expl:rg_obj,expl:rg_onecell}.
\end{itemize}

\subsection*{The 2-Category $\AlgstpsO$}

Recall the following.
\begin{itemize}
\item $\AlgpspsO$ is the 2-category \pcref{oalgps_twocat} with $\Op$-pseudoalgebras as objects, $\Op$-pseudomorphisms as 1-cells, and $\Op$-transformations as 2-cells \pcref{def:pseudoalgebra,def:laxmorphism,def:algtwocells}. 
\item An \emph{$\Op$-algebra} is an $\Op$-pseudoalgebra $(\A,\gaA,\phiA)$ such that each component of the associativity constraint $\phiA$ \cref{phiA} is the identity natural transformation \pcref{def:pseudoalgebra,def:operadalg}.
\item An \emph{$\Op$-algebra morphism} is a lax $\Op$-morphism $(f,\actf)$ between $\Op$-algebras such that each action constraint $\actf_n$ is the identity \pcref{def:laxmorphism}.
\end{itemize}

\begin{definition}\label{def:algstpsO}
For a reduced $\Gcat$-operad $\Op$, define the sub-2-categories 
\[\begin{tikzpicture}
\draw[0cell]
(0,0) node (a1) {\AlgstpsO}
(a1)++(2.5,0) node (a2) {\AlgpspsO}
(a1)++(-2.3,0) node (a0) {\AlgO}
;
\draw[1cell=.9]
(a0) edge[right hook->] node {\algi} (a1)
(a1) edge[right hook->] node {\algi} (a2)
;
\end{tikzpicture}\]
of $\AlgpspsO$ as follows. 
\begin{itemize}
\item $\AlgO$ has $\Op$-algebras as objects and $\Op$-algebra morphisms as 1-cells.
\item $\AlgstpsO$ has $\Op$-algebras as objects and $\Op$-pseudomorphisms as 1-cells.
\end{itemize}
In each of $\AlgO$ and $\AlgstpsO$, the 2-cells are $\Op$-transformations.
\end{definition}

\begin{explanation}\label{expl:algstpsO}
For an $\Op$-pseudomorphism \pcref{def:laxmorphism}
\[(\A,\gaA) \fto{(f,\actf)} (\B,\gaB)\]
between $\Op$-algebras, each action constraint $\actf_n$ \cref{laxmorphism_constraint} is a $G$-natural isomorphism.  In the associativity axiom \cref{laxmorphism_associativity}, the arrows $\phiB_{(n;\, m_\crdot)}$ and $f\phiA_{(n;\, m_\crdot)}$ are identity morphisms.
\end{explanation}

\subsection*{$\Rg$ on Objects}

Recall that, for a reduced $\Gcat$-operad $(\Op,\ga,\opu)$, an object in $\DGCatg$ is a $\DGG$-category, meaning a pointed $\Gcatst$-functor $\DG \to \Gcatst$ \pcref{def:dgo,def:dgcatg}.

\begin{definition}\label{def:rg_obj}
For a reduced $\Gcat$-operad $(\Op,\ga,\opu)$ and an $\Op$-algebra $(\A,\gaA)$, the pointed $\Gcatst$-functor
\begin{equation}\label{rga}
\DG \fto{\Rg\A} \Gcatst
\end{equation}
is defined as follows.
\begin{description}
\item[Objects]
$\Rg\A$ sends each pointed finite $G$-set $\mal \in \DG$ to the $\mal$-twisted product \pcref{def:proCnbe}
\begin{equation}\label{rga_mal}
(\Rg\A)\mal = \proAmal
\end{equation}
with $\mal$ regarded as a length-1 object in $\GG$ and $\pro{\A}{\ordz} = \bone$.
\item[1-cells]
For each pair $(\mal,\nbeta)$ of pointed finite $G$-sets, the pointed $G$-functor
\begin{equation}\label{rga_mn}
\DG(\mal,\nbeta) \fto{\Rg\A} \Catgst(\proAmal,\proAnbeta)
\end{equation}
sends an object $\objx = (\psi; \ang{x_j}_{j \in \ufsn})$ in $\DG(\mal,\nbeta)$ \cref{dgmn_obj} to the pointed functor
\begin{equation}\label{rga_objx}
\proAmal \fto{(\Rg\A)\objx} \proAnbeta
\end{equation}
defined by the assignment
\begin{equation}\label{rgaxa}
\big((\Rg\A)\objx\big) \ang{a_i}_{i \in \ufsm} = 
\bang{\gaA\big(x_j; \ang{a_i}_{i \in \psiinv j} \big)}_{j \in \ufsn} \in \proAnbeta
\end{equation}
for each $m$-tuple $\ang{a_i}_{i \in \ufsm} \in \proAmal$ of all objects or all morphisms.  In \cref{rgaxa}, each subset $\psiinv j \subseteq \ufsm$ has the inherited ordering. 
\item[2-cells]
The functor $\Rg\A$ \cref{rga_mn} sends each morphism $\objf \cn \objx \to \objy$ in $\DG(\mal,\nbeta)$ \cref{dgmn_morphism} to the pointed natural transformation
\begin{equation}\label{rgaf}
\begin{tikzpicture}[vcenter]
\def\t{28}
\draw[0cell]
(0,0) node (a1) {\phantom{A}}
(a1)++(2.5,0) node (a2) {\phantom{A}}
(a1)++(-.1,0) node (a1') {\proAmal}
(a2)++(.1,.03) node (a2') {\proAnbeta}
;
\draw[1cell=.8]
(a1) edge[bend left=\t] node {(\Rg\A)\objx} (a2)
(a1) edge[bend right=\t] node[swap] {(\Rg\A)\objy} (a2)
;
\draw[2cell=.9]
node[between=a1 and a2 at .3, rotate=-90, 2label={above,(\Rg\A)\objf}] {\Rightarrow}
;
\end{tikzpicture}
\end{equation}
whose component at an object $\ang{a_i}_{i \in \ufsm} \in \proAmal$ is the morphism
\begin{equation}\label{rgafa}
\scalebox{.9}{$\bang{\gaA\big(x_j; \ang{a_i}_{i \in \psiinv j} \big)}_{j \in \ufsn}
\fto{\ang{\gaA(f_j; \ang{a_i}_{i \in \psiinv j})}_{j \in \ufsn}}
\bang{\gaA\big(y_j; \ang{a_i}_{i \in \psiinv j} \big)}_{j \in \ufsn}$}
\end{equation}
in $\proAnbeta$.
\end{description}
This finishes the definition of the pointed $\Gcatst$-functor $\Rg\A$.
\end{definition}

\begin{explanation}[Enriched Functoriality]\label{expl:rg_obj}
Using \cref{expl:dgcatg}, we discuss the pointed $\Gcatst$-functoriality of $\Rg\A$ in \cref{def:rg_obj}.
\begin{description}
\item[Well-defined 1-cells and 2-cells] 
$(\Rg\A)\objx$ \cref{rga_objx} and $(\Rg\A)\objf$ \cref{rgaf} are pointed by
\begin{itemize}
\item the associativity axiom \pcref{phi_id} and
\item the fact that the basepoint of $\proAmal$ is the $m$-tuple $\ang{\zero}_{i \in \ufsm}$ of copies of the basepoint $\zero = \gaA_0(*) \in \A$.
\end{itemize}  
The functoriality of $(\Rg\A)\objx$, the naturality of $(\Rg\A)\objf$, and the functoriality of the assignment $\objf \mapsto (\Rg\A)\objf$ follow from the functoriality of the $\Op$-action $\gaA$ \cref{gaAn}.
\item[Identity 1-cells] $(\Rg\A)\objx$ sends each identity 1-cell $\objone_{\mal}$ \cref{dg_idonecell} to the identity functor $1_{\proAmal}$ by the action unity axiom \cref{pseudoalg_action_unity} for $\A$.
\item[Composition] $(\Rg\A)\objx$ preserves composition in the sense of \cref{xlmn} by the associativity axiom \cref{phiA} and the action equivariance axiom \cref{pseudoalg_action_sym} for $\A$.
\item[$G$-equivariance]
By \cref{xmn_geq}, the $G$-equivariance of $\Rg\A$ \cref{rga_mn} means the commutativity of the diagram
\begin{equation}\label{rga_geq_diag}
\begin{tikzpicture}[vcenter]
\def\v{-1.3}
\draw[0cell]
(0,0) node (a11) {\proAmal}
(a11)++(2.5,0) node (a12) {\phantom{\proAnbeta}}
(a12)++(0,.02) node (a12') {\proAnbeta}
(a11)++(0,\v) node (a21) {\proAmal}
(a12)++(0,\v) node (a22) {\phantom{\proAnbeta}}
(a22)++(0,.02) node (a22') {\proAnbeta}
;
\draw[1cell=.9]
(a11) edge node {(\Rg\A)(g\objx)} (a12)
(a11) edge node[swap] {\ginv} (a21)
(a21) edge node {(\Rg\A)\objx} (a22)
(a22') edge node[swap] {g} (a12')
;
\end{tikzpicture}
\end{equation}
for each $g \in G$ and $\objx = (\psi; \ang{x_j}_{j \in \ufsn})$ in $\DG(\mal,\nbeta)$.  For each $\obja = \ang{a_i}_{i \in \ufsm} \in \proAmal$, the following equalities prove that the diagram \cref{rga_geq_diag} commutes.
\[\begin{aligned}
& g\big((\Rg\A)\objx \big) (\ginv\obja) &&\\
&= g\big((\Rg\A)\objx \big) \ang{\ginv a_{gi}}_{i \in \ufsm} && \text{by \cref{proCnbe_gaction}}\\
&= g\bang{\gaA\big(x_j; \ang{\ginv a_{gi}}_{i \in \psiinv j} \big) }_{j \in \ufsn} && \text{by \cref{rgaxa}}\\
&= \bang{g \gaA\big(x_{\ginv j}; \ang{\ginv a_{gi}}_{i \in \psiinv\ginv j} \big) }_{j \in \ufsn} && \text{by \cref{proCnbe_gaction}}\\
&= \bang{\gaA\big(gx_{\ginv j}; \ang{a_{gi}}_{i \in (g\psi)^{-1} j} \big) }_{j \in \ufsn} && \text{by \cref{gaAn}}\\
&= \bang{\gaA\big(gx_{\ginv j} ; \tau_{\ginv} \ang{a_{i}}_{i \in (g\psi\ginv)^{-1} j} \big) }_{j \in \ufsn} && \text{by \cref{tauginv}}\\
&= \bang{\gaA\big(gx_{\ginv j} \tau_{\ginv} ; \ang{a_{i}}_{i \in (g\psi\ginv)^{-1} j} \big) }_{j \in \ufsn} && \text{by \cref{pseudoalg_action_sym}}\\
&= (\Rg\A)(g\objx)\obja && \text{by \cref{dgmn_gact,rgaxa}}
\end{aligned}\]
The preceding computation makes sense if
\begin{itemize}
\item $\objx$ is an object and $\obja$ is either an object or a morphism, or
\item $\objx$ is a morphism and $\obja$ is an object.\defmark
\end{itemize}
\end{description}
\end{explanation}

\subsection*{$\Rg$ on 1-Cells}

Recall that a 1-cell in $\DGCatg$ is a $\Pig$-strict $\Gcatst$-pseudotransformation \pcref{def:dgcatg_onecell,def:pig_strict}.

\begin{definition}\label{def:rg_onecell}
For an $\Op$-pseudomorphism \pcref{def:laxmorphism} between $\Op$-algebras
\[(\A,\gaA) \fto{(f,\actf)} (\B,\gaB),\]
the $\Pig$-strict $\Gcatst$-pseudotransformation 
\begin{equation}\label{rgf}
(\Rg\A) \fto{(\Rg f)} (\Rg\B)
\end{equation}
has, for each object $\mal \in \DG$, $\mal$-component pointed $G$-functor defined as
\begin{equation}\label{rgf_mal}
(\Rg\A)\mal = \proAmal \fto{(\Rg f)_{\mal} = f^m} (\Rg\B)\mal = \proBmal.
\end{equation}
For each object $\objx = (\psi; \ang{x_j}_{j \in \ufsn})$ in $\DG(\mal,\nbeta)$ \cref{dgmn_obj}, the $\objx$-component pointed natural isomorphism
\begin{equation}\label{rgf_objx}
\begin{tikzpicture}[vcenter]
\def\v{-1.4}
\draw[0cell]
(0,0) node (a11) {\proAmal}
(a11)++(2.3,0) node (a12) {\proBmal}
(a11)++(0,\v) node (a21) {\proAnbeta}
(a12)++(0,\v) node (a22) {\proBnbeta}
;
\draw[1cell=.9]
(a11) edge node {f^m} (a12)
(a12) edge node {(\Rg\B)\objx} (a22)
(a11) edge node[swap] {(\Rg\A)\objx} (a21)
(a21) edge node[swap] {f^n} (a22)
;
\draw[2cell]
node[between=a12 and a21 at .5, shift={(.4,-.1)}, rotate=240, 2labelalt={below,(\Rg f)_{\objx}}] {\Rightarrow}
;
\end{tikzpicture}
\end{equation}
sends an object $\obja = \ang{a_i}_{i \in \ufsm} \in \proAmal$ to the isomorphism
\begin{equation}\label{rgfxa}
\begin{tikzpicture}[vcenter]
\def\v{-1}
\draw[0cell=.9]
(0,0) node (a11) {\big((\Rg\B)\objx\big) f^m \obja}
(a11)++(5.2,0) node (a12) {f^n\big((\Rg\A)\objx \big) \obja}
(a11)++(0,\v) node (a21) {\bang{\gaB\big(x_j; \ang{fa_i}_{i \in \psiinv j}\big)}_{j \in \ufsn}}
(a12)++(0,\v) node (a22) {\bang{f \gaA\big(x_j ; \ang{a_i}_{i \in \psiinv j}\big)}_{j \in \ufsn}}
;
\draw[1cell=.85]
(a11) edge[equal] (a21)
(a12) edge[equal] (a22)
(a11) edge node {(\Rg f)_{\objx, \obja}} (a12)
(a21) edge node {\ang{\actf_{|\psiinv j|}}_{j \in \ufsn}} node[swap] {\iso} (a22)
;
\end{tikzpicture}
\end{equation}
in $\proBnbeta$ given by the action constraint $\actf$ \cref{laxmorphism_constraint}.  Note that if $f$ is an $\Op$-algebra morphism, then $\Rg f$ is a $\Gcatst$-natural transformation, with each $(\Rg f)_{\objx}$ given by the identity.
\end{definition}

\begin{explanation}\label{expl:rg_onecell}
We discuss the $\Pig$-strict $\Gcatst$-pseudotransformation axioms of $\Rg f$ in \cref{def:rg_onecell}.
\begin{description}
\item[Pointedness] 
$(\Rg f)_{\objx}$ is pointed---meaning $(\Rg f)_{\objx,\obja}$ is the identity morphism of the basepoint in $\proBnbeta$ when $\obja \in \proAmal$ is the basepoint---by the basepoint axiom \cref{laxmorphism_basepoint} and the associativity axiom \cref{laxmorphism_associativity} for $f$. 
\item[Naturality] 
The naturality of $(\Rg f)_{\objx, \obja}$ in $\obja \in \proAmal$ and in $\objx \in \DG(\mal,\nbeta)$ \cref{thaobjx_nat} follows from the naturality of $\actf$ \cref{laxmorphism_constraint}.
\item[Unity] 
For the identity 1-cell $\objone_{\mal}$ \cref{dg_idonecell}, $(\Rg f)_{\objone_{\mal}}$ is the identity natural transformation by the unity axiom \cref{laxmorphism_unity} for $f$.
\item[Compositionality] 
$\Rg f$ preserves composition \cref{tha_xv} by \cref{dg_comp}, the equivariance axiom \cref{laxmorphism_equiv}, and the associativity axiom \cref{laxmorphism_associativity} for $f$.
\item[$G$-equivariance] 
The axiom \cref{thaobjx_gequiv} for $(\Rg f)_{\objx}$ means the morphism equality
\begin{equation}\label{rgfx_geq}
(\Rg f)_{g\objx, \obja} = g(\Rg f)_{\objx, \ginv\obja} \inspace \proBnbeta
\end{equation}
for each $g \in G$, object $\objx = (\psi; \ang{x_j}_{j \in \ufsn})$ in $\DG(\mal,\nbeta)$ \cref{dgmn_obj}, and object $\obja = \ang{a_i}_{i \in \ufsm} \in \proAmal$.  The equality \cref{rgfx_geq} is proved by the commutative diagram \cref{rgfx_geq_diag}, where $\ang{\Cdots}_j = \ang{\Cdots}_{j \in \ufsn}$.
\begin{equation}\label{rgfx_geq_diag}
\begin{tikzpicture}[vcenter]
\def\v{-1}
\draw[0cell=.8]
(0,0) node (a11) {g\ang{\gaB(x_j; \ang{f\ginv a_{gi}}_{i \in \psiinv j})}_j}
(a11)++(6,0) node (a12) {g\ang{f\gaA(x_j; \ang{\ginv a_{gi}}_{i \in \psiinv j})}_j}
(a11)++(0,\v) node (a21) {\ang{\gaB(gx_{\ginv j}; \ang{fa_{gi}}_{i \in (g\psi)^{-1} j})}_j}
(a12)++(0,\v) node (a22) {\ang{f\gaA(gx_{\ginv j}; \ang{a_{gi}}_{i \in (g\psi)^{-1} j})}_j}
(a21)++(0,\v) node (a31) {\ang{\gaB(gx_{\ginv j}; \tau_{\ginv} \ang{fa_{i}}_{i \in (g\psi\ginv)^{-1} j})}_j}
(a22)++(0,\v) node (a32) {\ang{f\gaA(gx_{\ginv j}; \tau_{\ginv} \ang{a_{i}}_{i \in (g\psi\ginv)^{-1} j})}_j}
(a31)++(0,\v) node (a41) {\ang{\gaB(gx_{\ginv j} \tau_{\ginv} ; \ang{fa_{i}}_{i \in (g\psi\ginv)^{-1} j})}_j}
(a32)++(0,\v) node (a42) {\ang{f\gaA(gx_{\ginv j} \tau_{\ginv} ; \ang{a_{i}}_{i \in (g\psi\ginv)^{-1} j})}_j}
;
\draw[1cell=.7]
(a11) edge[equal] (a21)
(a21) edge[equal] (a31)
(a31) edge[equal] (a41)
(a12) edge[equal] (a22)
(a22) edge[equal] (a32)
(a32) edge[equal] (a42)
(a11) edge node {g\ang{\actf_{|\psiinv j|}}_j} (a12)
(a21) edge node {\ang{\actf_{|(g\psi)^{-1} j|}}_j} (a22)
(a31) edge node {\ang{\actf_{|(g\psi)^{-1} j|}}_j} (a32)
(a41) edge node {\ang{\actf_{|(g\psi\ginv)^{-1} j|}}_j} (a42)
;
\end{tikzpicture}
\end{equation}
Details of the diagram \cref{rgfx_geq_diag} are given next.
\begin{itemize}
\item By \cref{proCnbe_gaction,rgfxa}, the top horizontal arrow in \cref{rgfx_geq_diag} is $g(\Rg f)_{\objx, \ginv\obja}$, the right-hand side of the desired equality \cref{rgfx_geq}.  By \cref{dgmn_gact,rgfxa}, the bottom horizontal arrow in \cref{rgfx_geq_diag} is $(\Rg f)_{g\objx, \obja}$, the left-hand side of \cref{rgfx_geq}.
\item The top rectangle in \cref{rgfx_geq_diag} commutes by \cref{proCnbe_gaction} and the $G$-equivariance of $\gaA$, $\gaB$, $f$, and $\actf$.
\item The middle rectangle in \cref{rgfx_geq_diag} commutes by the definition of $\tau_{\ginv}$ \cref{tauginv}.
\item The bottom rectangle in \cref{rgfx_geq_diag} commutes by the equivariance axiom \cref{laxmorphism_equiv} for $f$.
\end{itemize}
\item[$\Pig$-strictness]
Suppose $\psi \cn \mal \to \nbeta$ is a pointed morphism in $\Pig$, and $\iota(\psi) \in \DG(\mal,\nbeta)$ is its image under $\iota \cn \Pig \to \DG$ \pcref{def:Pig,def:DGcoop}.  Since $|\psiinv j|$ is either 0 or 1 for each $j \in \ufsn$, 
\[(\Rg f)_{\iota(\psi), \obja} = \bang{\actf_{|\psiinv j|}}_{j \in \ufsn}\] 
is the identity morphism by the basepoint axiom \cref{laxmorphism_basepoint} and the unity axiom \cref{laxmorphism_unity} for $f$.\defmark
\end{description}
\end{explanation}

\subsection*{$\Rg$ on 2-Cells}

Recall that a 2-cell in $\DGCatg$ is a $\Gcatst$-modification \pcref{def:dgcatg_iicell}.

\begin{definition}\label{def:rg_twocell}
Suppose 
\begin{equation}\label{omega_fh_iicell}

\end{equation}
is an $\Op$-transformation \pcref{def:algtwocells} between $\Op$-pseudomorphisms $(f,\actf)$ and $(h,\acth)$ between $\Op$-algebras $(\A,\gaA)$ and $(\B,\gaB)$.  The $\Gcatst$-modification
\begin{equation}\label{rgom}
\Rg f \fto{\Rg\omega} \Rg h
\end{equation}
has $\mal$-component pointed $G$-natural transformation defined as
\begin{equation}\label{rgom_mal}
%
\end{equation}
for each object $\mal \in \DG$.  The pointed $G$-naturality of 
\[(\Rg\omega)_{\mal} = \omega^m\] 
follows from that of $\omega$ \pcref{expl:Otr_pointed}.  The modification axiom \cref{dgcatg_iicell_modax} holds by \cref{Otransformation_ax,rgfxa}.
\end{definition}

\subsection*{The 2-Functors}
The object, 1-cell, and 2-cell assignments in \cref{def:rg_obj,def:rg_onecell,def:rg_twocell} define the 2-functors
\begin{equation}\label{rg_twofunctor}
%
\end{equation}
between the 2-categories in \cref{def:algstpsO,def:dgcatg_iicat}.  For example, $\Rg$ preserves horizontal composition of 1-cells by \cref{Omorphism_paste,dgcatg_icell_paste,rgfxa}.  The composite 
\[\AlgO \fto{\algi\Rg = \Rg\algi} \DGCatg\]
is also denoted by $\Rg$.

\section{The 2-Category of Pseudo $\FGG$-Categories}
\label{sec:fgcatgps}

This section recalls the 2-category $\FGCatgps$ from \cite[Section 4]{gmmo23}.  The objects, 1-cells, and 2-cells are defined in \cref{def:psfggcat,def:psfggcat_icell,def:psfggcat_iicell}, which are further elaborated in \cref{expl:psfggcat,expl:psfggcat_iicat}.

\subsection*{Pseudo $\FGG$-Categories}
Recall the indexing $G$-category $\FG$ \pcref{def:FG} with pointed finite $G$-sets as objects and pointed morphisms with the conjugation $G$-action.  As discussed in \cref{def:Pig}, $\FG$ is also regarded as a $\Gcatst$-category with only identity 2-cells.  \cref{def:psfggcat} defines the pseudo analogue of an $\FGG$-category \cref{fgcatg_obj}.

\begin{definition}\label{def:psfggcat}
A \emph{pseudo $\FGG$-category}\index{pseudo FGG-category@pseudo $\FGG$-category}\index{FGG-category@$\FGG$-category!pseudo} is a pointed $\Gcatst$-pseudofunctor 
\[(\FG, \ordz) \fto{X} (\Gcatst,\bone)\]
that is strictly unital.
\end{definition}

\begin{explanation}[Unpacking]\label{expl:psfggcat}
A pseudo $\FGG$-category $X$ consists of the following data.
\begin{description}
\item[Pointed $G$-categories] 
$X$ sends each pointed finite $G$-set $\mal$ to a small pointed $G$-category $X\mal$ such that $X\ordz = \bone$, a terminal $G$-category.
\item[Pointed $G$-functors]
For each pair $(\mal,\nbeta)$ of pointed finite $G$-sets, $X$ is equipped with a pointed $G$-functor
\begin{equation}\label{psfggcat_xmn}
\FG(\mal,\nbeta) \fto{X_{\mal,\,\nbeta}} \Catgst(X\mal, X\nbeta).
\end{equation}
We often abbreviate $X_{\mal,\,\nbeta}$ to $X$.  Since $\FG$ has only identity 2-cells, the functor $X_{\mal,\,\nbeta}$ is determined by a pointed functor
\begin{equation}\label{psfggcat_xpsi}
X\mal \fto{X\psi} X\nbeta
\end{equation}
for each pointed morphism $\psi \cn \mal \to \nbeta$ such that the following two conditions are satisfied.
\begin{description}
\item[Pointedness] There is an equality of functors
\begin{equation}\label{psfggcat_basept}
X\mal \fto{X0 = *} X\nbeta
\end{equation}
for the 0-morphism $0 \cn \mal \to \nbeta$ and the basepoint $* \in X\nbeta$.
\item[$G$-equivariance] For each $g \in G$, there is an equality of functors
\begin{equation}\label{psfggcat_geq}
X\mal \fto{X(g\psi\ginv) = g(X\psi)\ginv} X\nbeta.
\end{equation}
\end{description}
\item[Pointed $G$-natural isomorphisms]
For each triple $(\ldea, \mal, \nbeta)$ of pointed finite $G$-sets, $X$ is equipped with a pointed $G$-natural isomorphism $X^2$ as follows.
\begin{equation}\label{psfggcat_xtwo}
\begin{tikzpicture}[vcenter]
\def\v{-1.4}
\draw[0cell=.85]
(0,0) node (a11) {\FG(\mal,\nbeta) \sma \FG(\ldea,\mal)}
(a11)++(5,0) node (a12) {\Catgst(X\mal,X\nbeta) \sma \Catgst(X\ldea,X\mal)}
(a11)++(0,\v) node (a21) {\FG(\ldea,\nbeta)}
(a12)++(0,\v) node (a22) {\Catgst(X\ldea,X\nbeta)}
;
\draw[1cell=.8]
(a11) edge node {X \sma X} (a12)
(a12) edge node {\comp} (a22)
(a11) edge node[swap] {\comp} (a21)
(a21) edge node[swap] {X} (a22)
;
\draw[2cell=.9]
node[between=a21 and a12 at .5, shift={(-.2,-.1)}, rotate=-120, 2labelalt={below,X^2}, 2label={above,\scalebox{.8}{$\iso$}}] {\Rightarrow}
;
\end{tikzpicture}
\end{equation}
Since $\FG$ has only identity 2-cells, $X^2$ is determined by a pointed natural isomorphism
\begin{equation}\label{psfggcat_xtwo_comp}
\begin{tikzpicture}[vcenter]
\def\t{26}
\draw[0cell]
(0,0) node (a1) {X\ldea}
(a1)++(2.3,0) node (a2) {X\nbeta}
;
\draw[1cell=.8]
(a1) edge[bend left=\t] node {(X\psi)(X\phi)} (a2)
(a1) edge[bend right=\t] node[swap] {X(\psi\phi)} (a2)
;
\draw[2cell=.9]
node[between=a1 and a2 at .36, rotate=-90, 2label={above, X^2_{\psi,\phi}}] {\Rightarrow}
;
\end{tikzpicture}
\end{equation}
for each pair $(\psi \cn \mal \to \nbeta, \phi \cn \ldea \to \mal)$ of pointed morphisms such that the following two conditions are satisfied.
\begin{description}
\item[Pointedness] 
There are equalities of natural isomorphisms
\begin{equation}\label{xtwo_basept}
X^2_{0,\phi} = 1_* = X^2_{\psi,0}
\end{equation}
for the 0-morphisms $0 \cn \mal \to \nbeta$ and $0 \cn \ldea \to \mal$ and the basepoint $* \in X\nbeta$.
\item[$G$-equivariance] For each $g \in G$, there is an equality of natural isomorphisms
\begin{equation}\label{xtwo_gequiv}
X^2_{g\psi\ginv, g\phi\ginv} = g(X^2_{\psi,\phi})\ginv.
\end{equation}
\end{description}
\end{description}
The preceding data are required to satisfy the axioms \crefrange{psfggcat_unity}{psfggcat_as}.
\begin{description}
\item[Unity]
For each pointed finite $G$-set $\mal$, there is an equality of functors
\begin{equation}\label{psfggcat_unity}
X\mal \fto{X1_{\mal} = 1_{X\mal}} X\mal.
\end{equation}
Moreover, in the context of \cref{psfggcat_xtwo_comp}, there are equalities of natural isomorphisms
\begin{equation}\label{psfggcat_xtwo_u}
X^2_{\psi,1_{\mal}} = 1_{X\psi} \andspace X^2_{1_{\mal},\phi} = 1_{X\phi}.
\end{equation}
\item[Associativity]
For each triple of pointed morphisms 
\[\ldea \fto{\phi} \mal \fto{\psi} \nbeta \fto{\rho} \rka\]
between pointed finite $G$-sets, the diagram of natural isomorphisms
\begin{equation}\label{psfggcat_as}
\begin{tikzpicture}[vcenter]
\def\v{-1.4}
\draw[0cell=.9]
(0,0) node (a11) {(X\rho)(X\psi)(X\phi)}
(a11)++(3.75,0) node (a12) {(X\rho)X(\psi\phi)}
(a11)++(0,\v) node (a21) {\big(X(\rho\psi)\big)(X\phi)}
(a12)++(0,\v) node (a22) {X(\rho\psi\phi)}
;
\draw[1cell=.8]
(a11) edge node {1_{X\rho} * X^2_{\psi,\phi}} (a12)
(a12) edge[transform canvas={xshift=-1ex}] node {X^2_{\rho,\psi\phi}} (a22)
(a11) edge[transform canvas={xshift=1em}] node[swap] {X^2_{\rho,\psi} * 1_{X\phi}} (a21)
(a21) edge node {X^2_{\rho\psi,\phi}} (a22)
;
\end{tikzpicture}
\end{equation}
commutes.
\end{description}
Note that an $\FGG$-category \cref{fgcatg_obj} yields a pseudo $\FGG$-category with $X^2 = 1$ and vice versa.
\end{explanation}

\subsection*{$\Gcatst$-Pseudotransformations}

\cref{def:psfggcat_icell} defines the pseudo analogue of a 1-cell \cref{fgcatg_icell} in the 2-category $\FGCatg$.

\begin{definition}\label{def:psfggcat_icell}
Suppose $X, Y \cn \FG \to \Gcatst$ are pseudo $\FGG$-categories \pcref{def:psfggcat}.  A \emph{$\Gcatst$-pseudotransformation}\index{Gcatst-pseudotransformation@$\Gcatst$-pseudotransformation} $\tha \cn X \to Y$ consists of the following data.
\begin{description}
\item[$G$-functors] For each object $\mal \in \FG$, $\tha$ is equipped with an $\mal$-component pointed $G$-functor
\begin{equation}\label{psfggcat_one_comp}
X\mal \fto{\tha_{\mal}} Y\mal.
\end{equation}
\item[Natural isomorphisms] For each pair $(\mal,\nbeta)$ of pointed finite $G$-sets and each pointed morphism $\psi \cn \mal \to \nbeta$, $\tha$ is equipped with a $\psi$-component pointed natural isomorphism $\tha_{\psi}$ as follows, where $X\psi$ and $Y\psi$ are the pointed functors in \cref{psfggcat_xpsi}.
\begin{equation}\label{psfggcat_one_compii}
\begin{tikzpicture}[vcenter]
\def\v{-1.4}
\draw[0cell]
(0,0) node (a11) {X\mal}
(a11)++(2.2,0) node (a12) {Y\mal}
(a11)++(0,\v) node (a21) {X\nbeta}
(a12)++(0,\v) node (a22) {Y\nbeta}
;
\draw[1cell=.9]
(a11) edge node {\tha_{\mal}} (a12)
(a12) edge node {Y\psi} (a22)
(a11) edge node[swap] {X\psi} (a21)
(a21) edge node[swap] {\tha_{\nbeta}} (a22)
;
\draw[2cell]
node[between=a12 and a21 at .5, shift={(0,-.05)}, rotate=240, 2label={above,\scalebox{.8}{$\iso$}}, 2labelw={below,\tha_{\psi},0pt}] {\Rightarrow}
;
\end{tikzpicture}
\end{equation}
\end{description}
The preceding data are required to satisfy the axioms \crefrange{thapsi_geq}{tha_psiphi}.
\begin{description}
\item[$G$-equivariance]
The assignment $\psi \mapsto \tha_{\psi}$ is $G$-equivariant, meaning
\begin{equation}\label{thapsi_geq}
\tha_{g\psi\ginv} = g \tha_{\psi} \ginv
\end{equation}
for $g \in G$ and $\psi \in \FG(\mal,\nbeta)$. 
\begin{itemize}
\item $g\psi\ginv$ is the conjugation $G$-action on $\psi$ \cref{gpsi}. 
\item The right-hand side is the conjugation $G$-action on $\tha_{\psi}$ \cref{conjugation-gaction}, using the $\ginv$-action on $X\mal$ and the $g$-action on $Y\nbeta$. 
\end{itemize} 
The axiom \cref{thapsi_geq} is well defined by \cref{psfggcat_geq} and the $G$-equivariance of $\tha_{\mal}$ and $\tha_{\nbeta}$ \cref{psfggcat_one_comp}.
\item[Basepoint]
The assignment $\psi \mapsto \tha_{\psi}$ preserves basepoints, meaning
\begin{equation}\label{thazero_istar}
\tha_0 = 1_*
\end{equation}
for the 0-morphism $0 \cn \mal \to \nbeta$ and the basepoint $* \in Y\nbeta$.  This axiom is well defined by \cref{psfggcat_basept}.
\item[Unity]
The assignment $\psi \mapsto \tha_{\psi}$ preserves identity 1-cells, meaning
\begin{equation}\label{tha_onemal}
\tha_{1_{\mal}} = 1_{\tha_{\mal}}
\end{equation}
for the identity morphism $1_{\mal} \in \FG(\mal,\mal)$.  This axiom is well defined by \cref{psfggcat_unity}.
\item[Compositionality]
Using the notation in \cref{psfggcat_xtwo_comp}, the following two pasting diagrams of natural isomorphisms are equal.
\begin{equation}\label{tha_psiphi}
\begin{tikzpicture}[vcenter]
\def\h{1.8} \def\d{1.4} \def\v{-1.3} \def\s{.8} \def\q{.7} \def\p{.8} \def\r{200}
\def\boundary{
\draw[0cell=\s]
(0,0) node (a11) {X\ldea}
(a11)++(\h,0) node (a12) {Y\ldea}
(a11)++(0,2*\v) node (a31) {X\nbeta}
(a12)++(0,2*\v) node (a32) {Y\nbeta}
;
\draw[1cell=\q]
(a11) edge node {\tha_{\ldea}} (a12)
(a31) edge node[swap] {\tha_{\nbeta}} (a32)
;}
\boundary
\draw[0cell=\s]
(a11)++(0,\v) node (a21) {X\mal}
(a12)++(0,\v) node (a22) {Y\mal}
;
\draw[1cell=\q]
(a21) edge node {\tha_{\mal}} (a22)
(a11) edge node[swap] {X\phi} (a21)
(a21) edge node[swap] {X\psi} (a31)
(a12) edge node {Y\phi} (a22)
(a22) edge node {Y\psi} (a32)
(a11) [rounded corners=2pt] -| ($(a21)+(-\d,1)$) -- node[swap] {X(\psi\phi)} ($(a21)+(-\d,-1)$) |- (a31)
;
\draw[2cell=\p]
node[between=a12 and a21 at .5, shift={(0,0)}, rotate=\r, 2labelmed={below,\tha_{\phi}}] {\Rightarrow}
node[between=a22 and a31 at .5, shift={(0,-.15)}, rotate=\r, 2labelmed={below,\tha_{\psi}}] {\Rightarrow}
node[between=a11 and a31 at .5, shift={(-.7*\d,0)}, rotate=\r, 2label={below,\phantom{a}}] {\Rightarrow}
node[between=a11 and a31 at .48, shift={(-.7*\d,0)}, rotate=180, 2labelmed={below,X^2_{\psi,\phi}}] {\phantom{\Rightarrow}}
;
\begin{scope}[shift={(2*\h,0)}]
\boundary
\draw[0cell=\s]
(a12)++(\d,\v) node (a22) {Y\mal}
;
\draw[1cell=\q]
(a11) edge node[swap] {X(\psi\phi)} (a31)
(a12) edge node[swap,pos=.15] {Y(\psi\phi)} (a32)
(a12) [rounded corners=2pt] -| node[pos=.8] {Y\phi} (a22)
;
\draw[1cell=\q]
(a22) [rounded corners=2pt] |- node[pos=.2] {Y\psi} (a32)
;
\draw[2cell=\p]
node[between=a12 and a31 at .5, shift={(0,0)}, rotate=\r, 2labelmed={below,\tha_{\psi\phi}}] {\Rightarrow}
node[between=a12 and a32 at .5, shift={(.5*\d,0)}, rotate=\r, 2label={below,\phantom{a}}] {\Rightarrow}
node[between=a12 and a32 at .48, shift={(.5*\d,0)}, rotate=180, 2labelmed={below,Y^2_{\psi,\phi}}] {\phantom{\Rightarrow}}
;
\end{scope}
\end{tikzpicture}
\end{equation}
\end{description}
Moreover, a $\Gcatst$-pseudotransformation $\tha$ is called a \emph{weak $G$-equivalence}\index{Gcatst-pseudotransformation@$\Gcatst$-pseudotransformation!weak G-equivalence@weak $G$-equivalence}\index{weak G-equivalence@weak $G$-equivalence!Gcatst-pseudotransformation@$\Gcatst$-pseudotransformation} if each $\tha_{\mal}$ \cref{psfggcat_one_comp} is a categorical weak $G$-equivalence \pcref{def:cat_weakg}.  This means that the $H$-fixed subfunctor 
\[(X\mal)^H \fto{\tha_{\mal}^H} (Y\mal)^H\]
is an equivalence of categories for each subgroup $H \subseteq G$.  Note that a 1-cell $\tha$ \cref{fgcatg_icell} in the 2-category $\FGCatg$ yields a $\Gcatst$-pseudotransformation between $\FGG$-categories with each $\tha_\psi = 1$ and vice versa.  Thus, weak $G$-equivalences also make sense for 1-cells in $\FGCatg$.
\end{definition}

\subsection*{2-Category Structure}

\begin{definition}\label{def:psfggcat_iicell}
Suppose $\tha, \ups \cn X \to Y$ are $\Gcatst$-pseudotransformations between pseudo $\FGG$-categories $X$ and $Y$ \pcref{def:psfggcat,def:psfggcat_icell}.  A \emph{$\Gcatst$-modification}\index{Gcatst-modification@$\Gcatst$-modification} $\Theta \cn \tha \to \ups$ consists of, for each object $\mal \in \FG$, an $\mal$-component pointed $G$-natural transformation
\begin{equation}\label{psfggcat_iicell_comp}
\begin{tikzpicture}[vcenter]
\def\t{30}
\draw[0cell]
(0,0) node (a1) {\phantom{A}}
(a1)++(2,0) node (a2) {\phantom{A}}
(a1)++(-.2,0) node (a1') {X\mal}
(a2)++(.2,0) node (a2') {Y\mal}
;
\draw[1cell=.9]
(a1) edge[bend left=\t] node {\tha_{\mal}} (a2)
(a1) edge[bend right=\t] node[swap] {\ups_{\mal}} (a2)
;
\draw[2cell=.9]
node[between=a1 and a2 at .38, rotate=-90, 2label={above,\Theta_{\mal}}] {\Rightarrow}
;
\end{tikzpicture}
\end{equation}
such that, for each pointed morphism $\psi \cn \mal \to \nbeta$, the following two pasting diagrams of natural transformations are equal.
\begin{equation}\label{psfggcat_iicell_modax}
\begin{tikzpicture}[vcenter]
\def\h{2} \def\v{-1.4} \def\u{1.4} \def\s{.8} \def\c{.4} \def\d{.6} \def\t{23} 
\def\boundary{
\draw[0cell=.9]
(0,0) node (a11) {X\mal}
(a11)++(\h,0) node (a12) {Y\mal}
(a11)++(0,\v) node (a21) {X\nbeta}
(a12)++(0,\v) node (a22) {Y\nbeta}
;
\draw[1cell=\s]
(a11) edge node[swap] {X\psi} (a21)
(a12) edge node {Y\psi} (a22)
(a11) edge[bend left=\t] node[pos=\c] {\tha_{\mal}} (a12)
(a21) edge[bend right=\t] node[swap,pos=\d] {\ups_{\nbeta}} (a22)
;}
\boundary
\draw[1cell=\s]
(a11) edge[bend right=\t] node[swap,pos=\d] {\ups_{\mal}} (a12)
;
\draw[2cell=\s]
node[between=a11 and a12 at .4, shift={(0,0)}, rotate=-90, 2label={above,\Theta_{\mal}}] {\Rightarrow}
node[between=a11 and a12 at .4, shift={(0,-\u)}, rotate=-90, 2label={above,\ups_{\psi}}] {\Rightarrow}
;
\begin{scope}[shift={(2*\h,0)}]
\boundary
\draw[1cell=\s]
(a21) edge[bend left=\t] node[pos=\c] {\tha_{\nbeta}} (a22)
;
\draw[2cell=\s]
node[between=a21 and a22 at .4, shift={(0,0)}, rotate=-90, 2label={above,\Theta_{\nbeta}}] {\Rightarrow}
node[between=a21 and a22 at .4, shift={(0,\u)}, rotate=-90, 2label={above,\tha_{\psi}}] {\Rightarrow}
;
\end{scope}
\end{tikzpicture}
\end{equation}
Note that $\Theta_{\mal}$ is not required to be invertible.
\end{definition}

\begin{definition}\label{def:psfggcat_iicat}
Define the 2-category $\FGCatgps$ with
\begin{itemize}
\item pseudo $\FGG$-categories \pcref{def:psfggcat} as objects,
\item $\Gcatst$-pseudotransformations \pcref{def:psfggcat_icell} as 1-cells, and
\item $\Gcatst$-modifications \pcref{def:psfggcat_iicell} as 2-cells.
\end{itemize}
Other 2-categorical structures---identities and horizontal composition of 1-cells and 2-cells and vertical composition of 2-cells---are defined componentwise.
\end{definition}

\begin{explanation}[Unpacking]\label{expl:psfggcat_iicat}
The description of the 2-category $\DGCatg$ in \cref{expl:dgcatg_iicat} also applies to the 2-category $\FGCatgps$ after replacing
\begin{itemize}
\item $\DG$ with $\FG$ and
\item the objects, 1-cells, and 2-cells in $\DGCatg$ with those in $\FGCatgps$.  
\end{itemize}
There is an inclusion 2-functor
\begin{equation}\label{fgcatg_psfgcatg}
\begin{tikzpicture}
\draw[0cell]
(0,0) node (a1) {\FGCatg}
(a1)++(2.5,0) node (a2) {\FGCatgps}
;
\draw[1cell=.9]
(a1) edge[right hook->] node {\Incj} (a2)
;
\end{tikzpicture}
\end{equation}
between the 2-categories in \cref{def:fgcatg,def:psfggcat_iicat}.  The superscript and subscript in the notation $\FGCatgps$ are reminders that the objects and 1-cells are pseudo structures.
\end{explanation}

\section{The 2-Functor $\gzest$}
\label{sec:gzest}

This section recalls the 2-functor
\[\DGCatg \fto{\gzest} \FGCatgps\]
between the 2-categories in \cref{def:dgcatg_iicat,def:psfggcat_iicat}.  The 2-functor $\gzest$ is the second step of GMMO $K$-theory \cref{Kgmmo_diag}, sending $\DGG$-categories to pseudo $\FGG$-categories at the object level.   The definitions in this section correspond to those in \cite[Section 7]{gmmo23} without multifunctoriality.

\secoutline
\begin{itemize}
\item \cref{def:chaotic_einf} defines chaotic $\Einfg$-operads.
\item \cref{ex:gbe_chaotic} shows that the $G$-Barratt-Eccles operad $\GBE$ is a chaotic $\Einfg$-operad.
\item \cref{def:ze_psfunctor} defines the strictly unital $\Gcatst$-pseudofunctor $\gze \cn \FG \to \DG$.  It is a section of $\gxi \cn \DG \to \FG$ \pcref{expl:gze}.
\item \cref{def:gzest} defines the pullback 2-functor $\gzest$.
\item \cref{expl:gzest,expl:gzest_icell,expl:gzest_iicell} further unpack the object, 1-cell, and 2-cell assignments of $\gzest$.
\end{itemize}

\subsection*{Chaotic $\Einfg$-Operads}

Recall the following.
\begin{itemize}
\item A $\Gcat$-operad $\Op$ is \emph{reduced} if $\Op(0)$ is a terminal $G$-category. 
\item $\cla$ denotes the classifying space functor \cref{classifying_space}.
\item For $n \geq 0$, a \emph{graph subgroup} of $\gsin$ is a subgroup $K$ such that $K \cap \Si_n = \{e\}$ \pcref{def:graph_subgrp}.  Each subgroup $H \subseteq G$ is also regarded as the graph subgroup $H \ttimes \{\id_n\}$ of $\gsin$.
\end{itemize}

\begin{definition}\label{def:chaotic_einf}
A \emph{chaotic $\Einfg$-operad}\index{chaotic Einfg-operad@chaotic $\Einfg$-operad}\index{operad!chaotic $\Einfg$} is a reduced $\Gcat$-operad $\Op$ that satisfies the following two conditions.
\begin{enumerate}
\item\label{chaotic_einf_i} $\Op$ is levelwise a translation category \pcref{def:translation_cat}.
\item\label{chaotic_einf_ii} For each $n \geq 0$ and subgroup $K \subseteq \gsin$, the $K$-fixed point subspace $\cla\Op(n)^K$ is contractible if $K$ is a graph subgroup and is empty otherwise.\defmark
\end{enumerate}
\end{definition}

\begin{explanation}\label{expl:gsin_action}
In \cref{def:chaotic_einf} \cref{chaotic_einf_ii}, the case $n=0$ is automatically true because, by the reduced assumption on $\Op$, 
\[\Op(0)^K = \bone^K = \bone\]
is a terminal $G$-category for each subgroup $K \subseteq G\ttimes\Si_0 \iso G$.  The $G$-category $\Op(n)$ with its right $\Si_n$-action becomes a left $(\gsin)$-category with the action defined as
\begin{equation}\label{six_left}
(g,\si) x = gx\sigmainv
\end{equation}
for $g \in G$, $\si \in \Si_n$, and $x \in \Op(n)$.  \cref{def:chaotic_einf} \cref{chaotic_einf_ii} implies that, for a subgroup $K \subseteq \gsin$, the $K$-fixed subcategory $\Op(n)^K$ is nonempty if $K$ is a graph subgroup and is empty otherwise.  Thus, 
\begin{equation}\label{OpnH}
\Op(n)^H \neq \emptyset
\end{equation}
for each $n \geq 0$ and each subgroup $H \subseteq G$.  
\end{explanation}

\begin{example}[$G$-Barratt-Eccles]\label{ex:gbe_chaotic}
The $G$-Barratt-Eccles operad $\GBE$ \pcref{def:GBE}, whose algebras are genuine permutative $G$-categories, is a chaotic $\Einfg$-operad.  We provide a proof here by specializing \cite[Lemma 3.7]{gmm17}.  It is reduced and levelwise a translation category because
\[\GBE(n) \iso \tn[G,\Si_n]\]
for $n \geq 0$.

To prove \cref{def:chaotic_einf} \cref{chaotic_einf_ii} for $\GBE$, note that for each $g \in G$ and each function $f \cn G \to \Si_n$, the $g$-action on $f$ is defined as
\[gf = f(\ginv(-)).\]
By \cref{six_left}, for each element $(g,\si) \in \gsin$, the $(g,\si)$-action on $f$ yields the function
\begin{equation}\label{gsif}
(g,\si)f = f(\ginv(-))\sigmainv.
\end{equation}
\begin{description}
\item[Nongraph subgroups]
If $K \subseteq \gsin$ is not a graph subgroup, then there is an element $(e,\si) \in K$ with $e \in G$ the group unit and $\si \not= \id_n \in \Si_n$.  Thus, we have that
\[(e,\si)f = f(-)\sigmainv \neq f,\]
and the $K$-fixed subcategory $\GBE(n)^K$ is empty.  
\item[Graph subgroups]
Next, suppose $K \subseteq \gsin$ is a graph subgroup.  By \cref{gph_subgrp} for a length-1 object in $\Gsk$, $K$ has the form
\[K = \{(h, \be h) \tmid h \in H \}\]
for some subgroup $H \subseteq G$ and homomorphism $\be \cn H \to \Si_n$.  To show that $\cla\GBE(n)^K$ is contractible, first note that 
\[\GBE(n)^K = \Catg(\EG,\ESigma_n)^K\] 
is a translation category because $\ESigma_n$ is so.  The classifying space $\cla\GBE(n)^K$ is contractible if the translation category $\GBE(n)^K$ is nonempty.  Thus, it suffices to construct a $K$-fixed function $f \cn G \to \Si_n$.  Regarding $G$ as an $H$-set by restricting the regular $G$-action, the $H$-set $G$ decomposes into its $H$-orbits as
\[G = \txcoprod_{i \in \ufsr}\, G_i.\]
Choosing an element $g_i \in G_i$ in each $H$-orbit of $G$, we define a function $f \cn G \to \Si_n$ by setting
\begin{equation}\label{fhgi}
f(hg_i) = (\be h)^{-1}
\end{equation}
for $h \in H$ and $i \in \ufsr$.  The function $f$ is well defined because, if $hg_i = h'g_i \in G_i$, then $h = h'$ and $\be h = \be h'$.  For each $(h,\be h) \in K$ and $h'g_i \in G_i$ for any $h, h' \in H$ and $i \in \ufsr$, the following equalities in $\Si_n$ prove that $f \in \GBE(n)$ is $K$-fixed.
\[\begin{aligned}
&\big((h,\be h)f\big)(h'g_i) &&\\
&= f(\hinv h'g_i) (\be h)^{-1} && \text{by \cref{gsif}}\\
&= \be(\hinv h')^{-1} (\be h)^{-1} && \text{by \cref{fhgi}}\\
&= (\be h')^{-1} && \text{by multiplicativity of $\be$}\\
&= f(h'g_i) && \text{by \cref{fhgi}}
\end{aligned}\]
Thus, the $K$-fixed subcategory $\GBE(n)^K$ is nonempty.\defmark
\end{description}
\end{example}

\subsection*{The $\Gcatst$-Pseudofunctor $\gze$}

\begin{definition}\label{def:ze_psfunctor}
Consider the $\Gcatst$-categories $\FG$ \pcref{def:FG,def:Pig} and $\DG$ associated to a chaotic $\Einfg$-operad $\Op$ \pcref{def:dgo,def:chaotic_einf}.  Define the strictly unital $\Gcatst$-pseudofunctor
\begin{equation}\label{gze}
\FG \fto{\gze} \DG
\end{equation}
as follows.
\begin{description}
\item[Objects] $\gze$ is the identity function on objects.
\item[Pointed $G$-functors] For each pair $(\mal,\nbeta)$ of pointed finite $G$-sets, the pointed $G$-functor
\begin{equation}\label{ze_mn}
\FG(\mal,\nbeta) \fto{\gze} \DG(\mal,\nbeta)
\end{equation}
is defined as the pointed $G$-functor
\begin{equation}\label{gze_pig}
\Pig(\mal,\nbeta) \fto{\gio} \DG(\mal,\nbeta)
\end{equation}
on the pointed $G$-subcategory $\Pig(\mal,\nbeta)$ \pcref{def:Pig,def:DGcoop}.  In particular, $\gze$ sends
\begin{itemize}
\item the 0-morphism $0 \cn \mal \to \nbeta$ to the basepoint $\objzero \in \DG(\mal,\nbeta)$ \cref{dgmn_basept} and
\item the identity morphism $1_{\mal} \in \Pig(\mal,\mal)$ to the identity 1-cell $\objone_{\mal} \in \DG(\mal,\mal)$ \cref{dg_idonecell}.
\end{itemize}

To define the rest of $\gze$, we choose a point in each $G$-orbit of the finite $G$-set 
\[\FG(\mal,\nbeta) \setminus \Pig(\mal,\nbeta).\]
For each such choice $\phi \cn \mal \to \nbeta$, suppose $H_\phi \subseteq G$ is the stabilizer of $\phi$.  By \cref{OpnH}, the $H_\phi$-fixed subcategory
\[\Big(\prod_{j \in \ufsn}\, \Op_{|\phiinv j|} \Big)^{H_\phi} 
\iso \prod_{j \in \ufsn}\, (\Op_{|\phiinv j|})^{H_\phi} \]
is nonempty.  We choose an $H_\phi$-fixed object \cref{dgmn_obj}
\[\gze(\phi) \in \DG(\mal,\nbeta)^{H_\phi} 
= \Big(\coprod_{\psi \in \FG(\mal,\,\nbeta)} \,\prod_{j \in \ufsn}\, \Op_{|\psiinv j|}\Big)^{H_\phi}\]
in the $\phi$-component, meaning
\begin{equation}\label{gxi_gze_phi}
\gxi\big(\gze(\phi)\big) = \phi \in \FG(\mal,\nbeta)
\end{equation}
for the $\Gcatst$-functor $\gxi \cn \DG \to \FG$ \cref{dgcoop}.  On the rest of the $G$-orbit of $\phi$, we use the $G$-actions on $\FG$ \cref{gpsi} and $\DG$ \cref{dgmn_gact} to define
\begin{equation}\label{gze_gphi}
\gze(g \cdot \phi) = g \cdot \gze(\phi).
\end{equation}
This finishes the construction of the pointed $G$-functor $\gze$ \cref{ze_mn}.
\item[Pointed $G$-natural isomorphisms]
For each triple $(\ldea, \mal, \nbeta)$ of pointed finite $G$-sets, $\gze$ is equipped with a pointed $G$-natural isomorphism $\gze^2$ as follows.
\begin{equation}\label{gze_two}
\begin{tikzpicture}[vcenter]
\def\v{-1.4}
\draw[0cell=.85]
(0,0) node (a11) {\FG(\mal,\nbeta) \sma \FG(\ldea,\mal)}
(a11)++(4.5,0) node (a12) {\DG(\mal,\nbeta) \sma \DG(\ldea,\mal)}
(a11)++(0,\v) node (a21) {\FG(\ldea,\nbeta)}
(a12)++(0,\v) node (a22) {\DG(\ldea,\nbeta)}
;
\draw[1cell=.9]
(a11) edge node {\gze \sma \gze} (a12)
(a12) edge node {\comp} (a22)
(a11) edge node[swap] {\comp} (a21)
(a21) edge node[swap] {\gze} (a22)
;
\draw[2cell=.9]
node[between=a21 and a12 at .5, shift={(0,-.1)}, rotate=-120, 2labelalt={below,\gze^2}, 2label={above,\scalebox{.8}{$\iso$}}] {\Rightarrow}
;
\end{tikzpicture}
\end{equation}
Recall the assumption that $\Op$ is levelwise a translation category (\cref{def:chaotic_einf} \cref{chaotic_einf_i}).  For each pair 
\[(\psi \cn \mal \to \nbeta, \phi \cn \ldea \to \mal)\] 
of pointed morphisms in $\FG$, the $(\psi,\phi)$-component of $\gze^2$ is defined as the unique isomorphism
\begin{equation}\label{gzetwo_psiphi}
\begin{split}
& \big((\gze\psi)(\gze\phi) \fto[\iso]{\gze^2_{\psi,\phi}} \gze(\psi\phi)\big) \\
& \in \DG(\ldea,\nbeta) = \coprod_{\rho \in \FG(\ldea,\nbeta)} \prod_{j \in \ufsn}\, \Op_{|\rho^{-1}j|}
\end{split}
\end{equation}
in the $(\psi\phi)$-component of $\DG(\ldea,\nbeta)$ \cref{dg_mn} with the indicated domain and codomain.  All the axioms for $\gze^2$---pointed $G$-naturality, unity of $\gze^2_{-,1_{\mal}}$ and $\gze^2_{1_{\mal},-}$, and associativity---follow from the assumption that each level of $\Op$ is a translation category.
\end{description}
This finishes the definition of the strictly unital $\Gcatst$-pseudofunctor $\gze$.
\end{definition}

\subsection*{Defining $\gzest$}

Recall the 2-categories $\DGCatg$ and $\FGCatgps$ \pcref{def:dgcatg_iicat,def:psfggcat_iicat}.

\begin{definition}\label{def:gzest}
For a chaotic $\Einfg$-operad $\Op$ \pcref{def:chaotic_einf}, the 2-functor
\begin{equation}\label{gzest}
\DGCatg \fto{\gzest} \FGCatgps
\end{equation}
is defined by precomposing with the strictly unital $\Gcatst$-pseudofunctor $\gze \cn \FG \to \DG$ \cref{gze}.
\end{definition}

\begin{explanation}\label{expl:gze}
By \cref{gze_pig,gxi_gze_phi,gze_gphi} and the $G$-equivariance of $\gxi$, the composite \pcref{def:DGcoop,def:ze_psfunctor}
\[\FG \fto{\gze} \DG \fto{\gxi} \FG\]
is the identity $\Gcatst$-functor on $\FG$.  Thus, there is a commutative diagram 
\begin{equation}\label{gxist_gzest}
\begin{tikzpicture}[vcenter]
\def\h{2.5} \def\u{.7}
\draw[0cell]
(0,0) node (a1) {\FGCatg}
(a1)++(\h,0) node (a2) {\DGCatg}
(a2)++(1.1*\h,0) node (a3) {\FGCatgps}
;
\draw[1cell=.9]
(a1) edge node {\gxist} (a2)
(a2) edge node {\gzest} (a3)
(a1) [rounded corners=2pt] |- ($(a2)+(-1,\u)$) -- node {\Incj} ($(a2)+(1,\u)$) -| (a3)
;
\end{tikzpicture}
\end{equation}
consisting of the pullback 2-functors $\gxist$ \cref{gxist} and $\gzest$ \cref{gzest} and the inclusion 2-functor $\Incj$ \cref{fgcatg_psfgcatg}.
\end{explanation}

The rest of this section unpacks the object, 1-cell, and 2-cell assignments of $\gzest$.

\begin{explanation}[$\gzest$ on Objects]\label{expl:gzest}
The 2-functor $\gzest$ \cref{gzest} sends a $\DGG$-category \cref{dgcatg_obj} 
\[\DG \fto{X} \Gcatst\]
to the pseudo $\FGG$-category \pcref{def:psfggcat}
\[\FG \fto{\gzest X} \Gcatst\]
given as follows.
\begin{description}
\item[Pointed $G$-categories]
On objects, $\gzest X$ is given by the small pointed $G$-categories
\begin{equation}\label{gzest_xmal}
(\gzest X)\mal = X\mal \forspace \mal \in \FG.
\end{equation} 
\item[Pointed $G$-functors]
For each pair $(\mal,\nbeta)$ of pointed finite $G$-sets, the $(\mal,\nbeta)$-component pointed $G$-functor of $\gzest X$ \cref{psfggcat_xmn} is the composite
\begin{equation}\label{gzestx_mn}
\begin{tikzpicture}[vcenter]
\def\u{.6}
\draw[0cell]
(0,0) node (a1) {\FG(\mal,\nbeta)}
(a1)++(2.7,0) node (a2) {\DG(\mal,\nbeta)}
(a2)++(3.2,0) node (a3) {\Catgst(X\mal,X\nbeta)}
;
\draw[1cell=.9]
(a1) edge node {\gze} (a2)
(a2) edge node {X} (a3)
(a1) [rounded corners=2pt, shorten >=-.5ex] |- ($(a2)+(-1,\u)$) -- node {\gzest X} ($(a2)+(1,\u)$) -| (a3)
;
\end{tikzpicture}
\end{equation}
of the pointed $G$-functors in \cref{ze_mn,dgcatg_xmn}.  The unity axiom \cref{psfggcat_unity} for $\gzest X$ holds by \cref{xmm_one,gze_pig}.
\item[Pointed $G$-natural isomorphisms]
For each triple $(\ldea, \mal, \nbeta)$ of pointed finite $G$-sets, the pointed $G$-natural isomorphism $(\gzest X)^2$ \cref{psfggcat_xtwo} is the whiskering of $\gze^2$ \cref{gze_two} with $X$, as displayed in the following diagram.
\begin{equation}\label{gzestx_two}
\begin{tikzpicture}[vcenter]
\def\h{3} \def\v{1} \def\u{-1.4}
\draw[0cell=.85]
(0,0) node (a11) {\FG(\mal,\nbeta) \sma \FG(\ldea,\mal)}
(a11)++(\h,\v) node (a12) {\DG(\mal,\nbeta) \sma \DG(\ldea,\mal)}
(a12)++(\h,-\v) node (a13) {\Catgst(X\mal,X\nbeta) \sma \Catgst(X\ldea,X\mal)}
(a11)++(0,\u) node (a21) {\FG(\ldea,\nbeta)}
(a21)++(\h,0) node (a22) {\DG(\ldea,\nbeta)}
(a22)++(\h,0) node (a23) {\Catgst(X\ldea,X\nbeta)}
;
\draw[1cell=.8]
(a11) edge node[swap] {\comp} (a21)
(a12) edge node[pos=.6] {\comp} (a22)
(a13) edge node {\comp} (a23)
(a21) edge node {\gze} (a22)
(a22) edge node {X} (a23)
(a11) [rounded corners=2pt] |- node[pos=.3] {\gze \sma \gze} (a12)
;
\draw[1cell=.8]
(a12) [rounded corners=2pt] -| node[pos=.7] {X \sma X} (a13)
;
\draw[2cell=.9]
node[between=a12 and a22 at .7, shift={(-.3*\h,0)}, rotate=-135, 2labelalt={below,\gze^2}] {\Rightarrow}
;
\end{tikzpicture}
\end{equation}
The right region commutes by \cref{xlmn}.  The unity axiom \cref{psfggcat_xtwo_u} for $(\gzest X)^2$ holds by
\begin{itemize}
\item the unity of $\gze^2_{-,1_{\mal}}$ and $\gze^2_{1_{\mal},-}$ and
\item the fact that $X$ preserves identity 2-cells \cref{xobjf}.
\end{itemize}  
The associativity axiom \cref{psfggcat_as} for $(\gzest X)^2$ holds by the associativity of $\gze^2$ and the fact that $X$ preserves composition \cref{xlmn}.\defmark
\end{description}
\end{explanation}

\begin{explanation}[$\gzest$ on 1-Cells]\label{expl:gzest_icell}
The 2-functor $\gzest$ \cref{gzest} sends a 1-cell $\tha \cn X \to X'$ in $\DGCatg$ \pcref{def:pig_strict} to the 1-cell \pcref{def:psfggcat_icell}
\[(\gzest X) \fto{\gzest\tha} (\gzest X') \inspace \FGCatgps\]
given as follows.
\begin{description}
\item[$G$-functors] For each pointed finite $G$-set $\mal$, the $\mal$-component pointed $G$-functor of $\gzest\tha$ \cref{psfggcat_one_comp} is given by the $\mal$-component of $\tha$ \cref{dgcatg_one_comp}:
\begin{equation}\label{gzest_icell_comp}
(\gzest X)\mal = X\mal \fto{(\gzest\tha)_{\mal} = \tha_{\mal}} (\gzest X')\mal = X'\mal.
\end{equation}
\item[Natural isomorphisms] For each pointed morphism $\psi \cn \mal \to \nbeta$ between pointed finite $G$-sets, the $\psi$-component pointed natural isomorphism of $\gzest\tha$ \cref{psfggcat_one_compii} is given by 
\[(\gzest\tha)_\psi = \tha_{\gze\psi}.\]
This is the component of $\tha$ \cref{dgcatg_one_compii} at the object $\gze\psi \in \DG(\mal,\nbeta)$ in the image of $\gze$ \cref{ze_mn}, as displayed in the following diagram.
\begin{equation}\label{gzest_icell_compii}
\begin{tikzpicture}[vcenter]
\def\v{-1.4}
\draw[0cell]
(0,0) node (a11) {X\mal}
(a11)++(2.2,0) node (a12) {X'\mal}
(a11)++(0,\v) node (a21) {X\nbeta}
(a12)++(0,\v) node (a22) {X'\nbeta}
;
\draw[1cell=.9]
(a11) edge node {\tha_{\mal}} (a12)
(a12) edge node {X'(\gze\psi)} (a22)
(a11) edge node[swap] {X(\gze\psi)} (a21)
(a21) edge node[swap] {\tha_{\nbeta}} (a22)
;
\draw[2cell]
node[between=a12 and a21 at .5, shift={(.2,-.1)}, rotate=240, 2labelalt={below,\tha_{\gze\psi}}] {\Rightarrow}
;
\end{tikzpicture}
\end{equation}
\end{description}
\begin{itemize}
\item The $G$-equivariance axiom \cref{thapsi_geq} for $\gzest\tha$ follows from the $G$-equivariance of $\tha$ \cref{thaobjx_gequiv} and $\gze$ \cref{ze_mn}.
\item The basepoint axiom \cref{thazero_istar} for $\gzest\tha$ holds by
\begin{itemize}
\item the basepoint axiom \cref{tha_objzero} for $\tha$ and
\item the fact that $\gze$ sends the 0-morphism $0 \cn \mal \to \nbeta$ to the basepoint $\objzero \in \DG(\mal,\nbeta)$ \cref{gze_pig}.
\end{itemize}
\item The unity axiom \cref{tha_onemal} for $\gzest\tha$ holds by
\begin{itemize}
\item the unity axiom \cref{tha_onex} for $\tha$ and
\item the fact that $\gze$ sends the identity morphism $1_{\mal}$ in $\FG$ to the identity 1-cell $\objone_{\mal}$ \cref{dg_idonecell} in $\DG$ \cref{gze_pig}.
\end{itemize}
\item The compositionality axiom \cref{tha_psiphi} for $\gzest\tha$ follows from 
\begin{itemize}
\item the definition \cref{gzestx_two} of $(\gzest X)^2$ and $(\gzest X')^2$;
\item  the compositionality axiom \cref{tha_xv} for $\tha$ applied to the objects 
\[\objv = \gze\phi \in \DG(\ldea,\mal) \andspace \objx = \gze\psi \in \DG(\mal,\nbeta);\]
and
\item the naturality of $\tha$ \cref{thaobjx_nat} for the isomorphism 
\[(\gze\psi)(\gze\phi) \fto[\iso]{\gze^2_{\psi,\phi}} \gze(\psi\phi)\]
in $\DG(\ldea,\nbeta)$ \cref{gzetwo_psiphi}.
\end{itemize}
\end{itemize}
The $\Pig$-strictness of $\tha$ and \cref{gze_pig} imply that
\begin{equation}\label{tha_giopsi}
\tha_{\gze\psi} = \tha_{\gio\psi} = 1 \forspace \psi \in \Pig(\mal,\nbeta).
\end{equation}
Moreover, by \cref{gzest_icell_comp}, if $\tha$ is a weak $G$-equivalence, then so is $\gzest\tha$.
\end{explanation}

\begin{explanation}[$\gzest$ on 2-Cells]\label{expl:gzest_iicell}
The 2-functor $\gzest$ \cref{gzest} sends a 2-cell $\Theta \cn \tha \to \ups$ between 1-cells $\tha, \ups \cn X \to X'$ in $\DGCatg$ \pcref{def:dgcatg_iicell} to the 2-cell \pcref{def:psfggcat_iicell}
\[\gzest\tha \fto{\gzest\Theta} \gzest\ups \inspace \FGCatgps\]
given by
\begin{equation}\label{gzest_iicell_comp}
(\gzest\Theta)_{\mal} = \Theta_{\mal} \forspace \mal \in \FG.
\end{equation}
This is well defined by \cref{gzest_icell_comp}.  The modification axiom \cref{psfggcat_iicell_modax} for $\gzest\Theta$ follows from the modification axiom \cref{dgcatg_iicell_modax} for $\Theta$ and \cref{gzest_icell_compii}.
\end{explanation}

\section{The Strictification 2-Functor}
\label{sec:str}

This section recalls the 2-functor
\[\FGCatgps \fto{\str} \FGCatg\]
between the 2-categories in \cref{def:fgcatg,def:psfggcat_iicat}.  The 2-functor $\str$ is the third step of GMMO $K$-theory \cref{Kgmmo_diag}, sending pseudo $\FGG$-categories to $\FGG$-categories at the object level.  The definitions in this section correspond to those in \cite[Section 8.1]{gmmo23} restricted to (pseudo) $\FGG$-categories.

\secoutline
\begin{itemize}
\item \cref{def:str_iifunctor} defines the 2-functor $\str$.
\item \cref{expl:str_obj,expl:str_icell,expl:str_iicell} describe the object, 1-cell, and 2-cell assignments of $\str$.
\item \cref{expl:str_unit} describes the unit and counit of the 2-adjunction $(\str,\Incj)$.
\end{itemize}


\begin{definition}\label{def:str_iifunctor}
Define the \emph{strictification 2-functor}\index{strictification 2-functor}\index{2-functor!strictification}
\begin{equation}\label{str_iifunctor}
\FGCatgps \fto{\str} \FGCatg
\end{equation}
as a left 2-adjoint of the inclusion 2-functor $\Incj \cn \FGCatg \hookrightarrow \FGCatgps$ \cref{fgcatg_psfgcatg}.
\end{definition}

The rest of this section unpacks the 2-functor $\str$ and the 2-adjunction $(\str,\Incj)$.

\begin{explanation}[$\str$ on Objects]\label{expl:str_obj}
The strictification 2-functor $\str$ \cref{str_iifunctor} sends a pseudo $\FGG$-category \pcref{def:psfggcat}
\[\FG \fto{X} \Gcatst\]
to the $\FGG$-category \cref{fgcatg_obj}
\begin{equation}\label{str_x}
\FG \fto{\str X} \Catgst
\end{equation}
given as follows.
\begin{description}
\item[Objects]
$\str X$ sends each pointed finite $G$-set $\mal$ to the small pointed $G$-category $(\str X)\mal$ with object pointed $G$-set
\begin{equation}\label{obstrx_mal}
\Ob\!\big((\str X)\mal\big) = \bigvee_{\ldea \in \FG} \FG(\ldea,\mal) \sma \Ob(X\ldea).
\end{equation}
The wedge in \cref{obstrx_mal} is indexed by the set of objects in $\FG$.  An object in $(\str X)\mal$ has the form
\begin{equation}\label{strx_mal_obj}
\objx = (\phi; x)
\end{equation}
with $\phi \in \FG(\ldea,\mal)$ and $x \in \Ob(X\ldea)$.  The object $\objx$ is the basepoint of $(\str X)\mal$ if either
\begin{itemize}
\item $\phi$ is the 0-morphism or
\item $x$ is the basepoint in $X\ldea$.
\end{itemize}  
\begin{description}
\item[$G$-action] The group $G$ acts diagonally, meaning 
\begin{equation}\label{strx_objx_gact}
g\objx = (g\phi\ginv; gx) \forspace g \in G.
\end{equation}
\item[Morphisms]
For an object 
\[\objy = (\rho; y) \in (\str X)\mal\]
with $\rho \in \FG(\rka, \mal)$ and $y \in \Ob(X\rka)$, the morphism $G$-set from $\objx$ to $\objy$ is given by
\begin{equation}\label{morstrx}
\big((\str X)\mal\big)(\objx, \objy) = (X\mal)\big((X\phi)x, (X\rho)y \big).
\end{equation}
Identities, composition, and the $G$-action on morphisms of $(\str X)\mal$ are induced by those of the small pointed $G$-category $X\mal$.  The $G$-action on morphisms is well defined by \cref{strx_objx_gact} and the $G$-equivariance of $X$ \cref{psfggcat_geq}.  More precisely, for $g \in G$ and a morphism $\objd \cn \objx \to \objy$ in $(\str X)\mal$ \cref{morstrx} given by a morphism
\[(X\phi)x \fto{\objd} (X\rho)y \inspace X\mal,\]
the $g$-action yields the following morphism in $((\str X)\mal)(g\objx,g\objy)$.
\begin{equation}\label{objd_gact}
\begin{tikzpicture}[vcenter]
\def\u{-.9}
\draw[0cell=.9]
(0,0) node (a11) {g(X\phi)x}
(a11)++(3,0) node (a12) {g(X\rho)y}
(a11)++(0,\u) node (a21) {(g(X\phi)\ginv)(gx)}
(a12)++(0,\u) node (a22) {(g(X\rho)\ginv)(gy)}
(a21)++(0,\u) node (a31) {X(g\phi\ginv)(gx)}
(a22)++(0,\u) node (a32) {X(g\rho\ginv)(gy)}
;
\draw[1cell=.9]
(a11) edge node {g\objd} (a12)
(a11) edge[equal] (a21)
(a21) edge[equal] (a31)
(a12) edge[equal] (a22)
(a22) edge[equal] (a32)
;
\end{tikzpicture}
\end{equation}
\item[Pointedness]
Since $\FG(-,\ordz) = *$ and $X\ordz = \bone$, $(\str X)\ordz$ is a terminal $G$-category, as required for an $\FGG$-category.
\end{description}
\item[Morphisms]
Suppose $\psi \cn \mal \to \nbeta$ is a pointed morphism between pointed finite $G$-sets.  The pointed functor
\begin{equation}\label{strx_psi}
(\str X)\mal \fto{(\str X)\psi} (\str X)\nbeta
\end{equation}
sends an object $\objx = (\phi;x)$ in $(\str X)\mal$ \cref{strx_mal_obj} to the object
\begin{equation}\label{strx_psi_objx}
\big((\str X)\psi\big) (\phi; x) = (\psi\phi; x).
\end{equation}
The pointed functor $(\str X)\psi$ sends a morphism $\objd \cn \objx \to \objy$ in $(\str X)\mal$ \cref{morstrx} to the morphism
\[((\str X)\psi)\objx = (\psi\phi; x) \fto{((\str X)\psi)\objd} ((\str X)\psi)\objy = (\psi\rho; y)\]
in $(\str X)\nbeta$ given by the following composite morphism in $X\nbeta$.
\begin{equation}\label{strx_psi_mor}
\begin{tikzpicture}[vcenter]
\def\v{-1.3}
\draw[0cell=.9]
(0,0) node (a11) {X(\psi\phi)x}
(a11)++(3.2,0) node (a12) {X(\psi\rho)y}
(a11)++(0,\v) node (a21) {(X\psi)(X\phi)x}
(a12)++(0,\v) node (a22) {(X\psi)(X\rho)y}
;
\draw[1cell=.85]
(a11) edge node {((\str X)\psi)\objd} (a12)
(a11) edge node[swap] {(X^2_{\psi,\phi,x})^{-1}} node {\iso} (a21)
(a21) edge node {(X\psi)\objd} (a22)
(a22) edge node[swap] {X^2_{\psi,\rho,y}} node {\iso} (a12)
;
\end{tikzpicture}
\end{equation}
\begin{itemize}
\item $(X\psi)\objd$ is the image of the morphism $\objd \cn (X\phi)x \to (X\rho)y$ under the pointed functor $X\psi \cn X\mal \to X\nbeta$ \cref{psfggcat_xpsi}.
\item $(X^2_{\psi,\phi,x})^{-1}$ is the inverse of the $x$-component of the pointed natural isomorphism  \cref{psfggcat_xtwo_comp}
\[(X\psi)(X\phi) \fto{X^2_{\psi,\phi}} X(\psi\phi).\]
\item $X^2_{\psi,\rho,y}$ is the $y$-component of the pointed natural isomorphism 
\[(X\psi)(X\rho) \fto{X^2_{\psi,\rho}} X(\psi\rho).\] 
\end{itemize}
The functoriality of $(\str X)\psi$ follows from \cref{strx_psi_mor} and the functoriality of $X\psi$.  It preserves the basepoints by \cref{strx_psi_objx} and the fact that $\psi\phi$ is the 0-morphism whenever $\phi$ is so.
\item[Functoriality]
To see that the assignment $\psi \mapsto (\str X)\psi$ preserves identities and composition, we consider objects and morphisms separately.
\begin{description}
\item[Objects]
$\str X$ preserves identities and composition on objects of $(\str X)\mal$ \cref{obstrx_mal} by \cref{strx_psi_objx}, since $\FG$ is a category.  
\item[Morphisms]
Consider morphisms $\objd \cn \objx \to \objy$ in $(\str X)\mal$ \cref{morstrx}.
\begin{description}
\item[Identities] $\str X$ preserves identities, meaning 
\[\big((\str X)1_{\mal}\big) \objd = \objd,\]
by \cref{strx_psi_mor} and the unity axioms \cref{psfggcat_unity,psfggcat_xtwo_u} of $X$. 
\item[Composition] $\str X$ preserves composition if, given pointed morphisms $\psi \cn \mal \to \nbeta$ and $\varphi \cn \nbeta \to \tpi$, there is an equality of morphisms
\[(\str X)(\varphi\psi)\objd = \big((\str X)\varphi\big) \big((\str X)\psi\big)\objd \inspace (\str X)\tpi.\]
The preceding equality is proved by the following commutative diagram in $X\tpi$, where $X^{-2}$ means $(X^2)^{-1}$.
\begin{equation}\label{StX_diag}
\begin{tikzpicture}[vcenter]
\def\h{5.2} \def\b{1} \def\v{-1.4} \def\s{.7} \def\t{.3}
\draw[0cell=.8]
(0,0) node (a11) {(X\varphi)(X\psi)(X\phi)x}
(a11)++(\h,0) node (a12) {(X\varphi)(X\psi)(X\rho)y}
(a11)++(-\b,\v) node (a21) {(X\varphi) X(\psi\phi)x}
(a11)++(\b,\v) node (a22) {X(\varphi\psi) (X\phi)x}
(a12)++(-\b,\v) node (a23) {X(\varphi\psi) (X\rho)y}
(a12)++(\b,\v) node (a24) {(X\varphi) X(\psi\rho)y}
(a11)++(0,2*\v) node (a31) {X(\varphi\psi\phi)x}
(a12)++(0,2*\v) node (a32) {X(\varphi\psi\rho)y}
;
\draw[1cell=.75]
(a31) edge node[pos=\s] {X^{-2}_{\varphi,\psi\phi,x}} (a21)
(a21) edge node[pos=\t] {(X\varphi) X^{-2}_{\psi,\phi,x}} (a11)
(a11) edge node {(X\varphi)(X\psi)\objd} (a12)
(a12) edge node[pos=\s] {(X\varphi)X^2_{\psi,\rho,y}} (a24)
(a24) edge node[pos=\t] {X^2_{\varphi,\psi\rho,y}} (a32)
(a31) edge node {(\str X)(\varphi\psi)\objd} (a32)
(a31) edge node[swap,pos=\s] {X^{-2}_{\varphi\psi,\phi,x}} (a22)
(a22) edge node {X(\varphi\psi)\objd} (a23)
(a23) edge node[swap,pos=\t] {X^2_{\varphi\psi,\rho,y}} (a32)
(a22) edge node[swap,pos=\t] {X^{-2}_{\varphi,\psi,(X\phi)x}} (a11)
(a12) edge node[swap,pos=\s] {X^2_{\varphi,\psi,(X\rho)y}} (a23)
;
\end{tikzpicture}
\end{equation}
\begin{itemize}
\item The outer boundary composite along the left, top, and right is $\big((\str X)\varphi\big) \big((\str X)\psi\big)\objd$ by \cref{strx_psi_mor} for $\varphi$ and $\psi$, together with the functoriality of $X\varphi$.
\item The left and right quadrilaterals commute by the associativity axiom \cref{psfggcat_as} for $X$.
\item The top middle trapezoid commutes by the naturality of $X^2_{\varphi,\psi}$.
\item The bottom middle trapezoid commutes by \cref{strx_psi_mor} for the composite $\varphi\psi$.
\end{itemize}
\end{description}
\end{description}
\item[$G$-equivariance]
The $G$-equivariance \cref{FGGcat_Gequiv} of $\str X$ means the equality of functors
\begin{equation}\label{strx_geq}
(\str X)\mal \fto[ = g((\str X)\psi)\ginv]{(\str X)(g\psi\ginv)} (\str X)\nbeta
\end{equation}
for each $g \in G$ and pointed morphism $\psi \cn \mal \to \nbeta$ in $\FG$.  
\begin{description}
\item[Objects]
For an object $(\phi; x) \in (\str X)\mal$ \cref{obstrx_mal}, the desired equality \cref{strx_geq} is proved by the following object equalities in $(\str X)\nbeta$.
\[\begin{aligned}
&\big(g((\str X)\psi)\ginv\big) (\phi; x) &&\\
&= \big(g((\str X)\psi)\big) (\ginv\phi g; \ginv x) && \text{by \cref{strx_objx_gact}}\\
&= g\big(\psi\ginv\phi g; \ginv x\big) && \text{by \cref{strx_psi_objx}}\\
&= \big(g\psi\ginv\phi g\ginv; g\ginv x  \big) && \text{by \cref{strx_objx_gact}}\\
&= (g\psi\ginv\phi; x) && \text{by $g\ginv = 1$} \\
&= \big((\str X)(g\psi\ginv)\big) (\phi; x) && \text{by \cref{strx_psi_objx}}
\end{aligned}\]
\item[Morphisms]
For a morphism $\objd \cn \objx \to \objy$ in $(\str X)\mal$ \cref{morstrx}, the $\ginv$-action yields a morphism
\begin{equation}\label{objd_ginv}
\ginv\objx = (\ginv\phi g; \ginv x) \fto{\ginv\objd} \ginv\objy = (\ginv\rho g; \ginv y)
\end{equation}
by \cref{objd_gact}.  The desired equality \cref{strx_geq} is proved by the following morphism equalities in $X\nbeta$, where $X^{-2}$ means $(X^2)^{-1}$.
{\small
\[\begin{aligned}
& g[((\str X)\psi) (\ginv\objd)] &&\\
&= g\big[\big(X^2_{\psi,\ginv\rho g,\ginv y}\big) \big((X\psi)(\ginv\objd)\big) \big(X^{-2}_{\psi,\ginv\phi g, \ginv x}\big) \big] && \text{by \cref{strx_psi_mor,objd_ginv}}\\
&= \big(gX^2_{\psi,\ginv\rho g,\ginv y}\big) \big(g(X\psi)(\ginv\objd)\big) \big(gX^{-2}_{\psi,\ginv\phi g, \ginv x}\big) && \text{by functoriality of $g \cdot -$}\\
&= \big(X^2_{g\psi\ginv,\rho,y} \big) \big(X(g\psi\ginv)\objd\big) \big(X^{-2}_{g\psi\ginv,\phi,x}\big) && \text{by \cref{psfggcat_geq,xtwo_gequiv}}\\
&= (\str X)(g\psi\ginv)\objd && \text{by \cref{strx_psi_mor}}
\end{aligned}\]
\par}
\end{description}
\end{description}
This finishes the description of the $\FGG$-category $\str X$ \cref{str_x}.
\end{explanation}

\begin{explanation}[$\str$ on 1-Cells]\label{expl:str_icell}
Suppose $\tha \cn X \to Y$ is a $\Gcatst$-pseudotransformation between pseudo $\FGG$-categories $X,Y \cn \FG \to \Gcatst$ \pcref{def:psfggcat_icell}.  The strictification 2-functor $\str$ \cref{str_iifunctor} sends $\tha$ to the $G$-natural transformation \cref{fgcatg_icell}
\begin{equation}\label{str_tha}
\begin{tikzpicture}[vcenter]
\def\t{30}
\draw[0cell]
(0,0) node (a1) {\phantom{A}}
(a1)++(2,0) node (a2) {\phantom{A}}
(a1)++(-.05,0) node (a1') {\FG}
(a2)++(.2,0) node (a2') {\Catgst}
;
\draw[1cell=.8]
(a1) edge[bend left=\t] node {\str X} (a2)
(a1) edge[bend right=\t] node[swap] {\str Y} (a2)
;
\draw[2cell=1]
node[between=a1 and a2 at .37, rotate=-90, 2label={above,\str\tha}] {\Rightarrow}
;
\end{tikzpicture}
\end{equation}
with, for each pointed finite $G$-set $\mal$, $\mal$-component pointed $G$-functor
\begin{equation}\label{str_tha_mal}
(\str X)\mal \fto{(\str\tha)_{\mal}} (\str Y)\mal
\end{equation}
given as follows.
\begin{description}
\item[Objects]
$(\str\tha)_{\mal}$ sends an object $\objx = (\phi; x)$ in $(\str X)\mal$ \cref{strx_mal_obj} to the object
\begin{equation}\label{strtha_mal_obj}
(\str\tha)_{\mal} (\phi; x) = \big(\phi; \tha_{\ldea}x\big) \inspace (\str Y)\mal, 
\end{equation}
where $\tha_{\ldea} \cn X\ldea \to Y\ldea$ is the $\ldea$-component pointed $G$-functor of $\tha$ \cref{psfggcat_one_comp}.  The pointed $G$-equivariance of $(\str\tha)_{\mal}$ on objects follows from \cref{strx_objx_gact} and the pointed $G$-equivariance of $\tha_{\ldea}$.
\item[Morphisms]
$(\str\tha)_{\mal}$ sends a morphism $\objd \cn \objx \to \objy$ in $(\str X)\mal$ \cref{morstrx} to the morphism
\[(\str\tha)_{\mal}\objx = \big(\phi; \tha_{\ldea}x\big) \fto{(\str\tha)_{\mal}\objd} 
(\str\tha)_{\mal}\objy = \big(\rho; \tha_{\rka}y\big)\]
in $(\str Y)\mal$ given by the following composite morphism in $Y\mal$.
\begin{equation}\label{strx_tha_d}
\begin{tikzpicture}[vcenter]
\def\v{-1.3}
\draw[0cell=.9]
(0,0) node (a11) {(Y\phi)(\tha_{\ldea} x)}
(a11)++(3.2,0) node (a12) {(Y\rho)(\tha_{\rka} y)}
(a11)++(0,\v) node (a21) {\tha_{\mal}(X\phi)x}
(a12)++(0,\v) node (a22) {\tha_{\mal}(X\rho)y}
;
\draw[1cell=.85]
(a11) edge node {(\str\tha)_{\mal}\objd} (a12)
(a11) edge node[swap] {\tha_{\phi,x}} node {\iso} (a21)
(a21) edge node {\tha_{\mal}\objd} (a22)
(a22) edge node[swap] {\tha^{-1}_{\rho,y}} node {\iso} (a12)
;
\end{tikzpicture}
\end{equation}
\begin{itemize}
\item $\tha_{\mal}\objd$ is the image of the morphism $\objd \cn (X\phi)x \to (X\rho)y$ under the pointed $G$-functor $\tha_{\mal} \cn X\mal \to Y\mal$ \cref{psfggcat_one_comp}.
\item $\tha_{\phi,x}$ is the $x$-component of the pointed natural isomorphism \cref{psfggcat_one_compii}
\[(Y\phi)\tha_{\ldea} \fto{\tha_{\phi}} \tha_{\mal}(X\phi).\]
\item $\tha^{-1}_{\rho,y}$ is the inverse of the $y$-component of the pointed natural isomorphism 
\[(Y\rho)\tha_{\rka} \fto{\tha_{\rho}} \tha_{\mal}(X\rho).\]
\end{itemize}
The functoriality of the assignment $\objd \mapsto (\str\tha)_{\mal}\objd$ follows from the functoriality of $\tha_{\mal}$.  The following equalities for $g \in G$ prove that it is $G$-equivariant.
{\small
\[\begin{aligned}
& (\str\tha)_{\mal}(g\objd) &&\\
&= \big(\thainv_{g\rho\ginv, gy} \big) \big(\tha_{\mal}(g\objd)\big) \big(\tha_{g\phi\ginv,gx}\big) && \text{by \cref{objd_gact,strx_tha_d}}\\
&= \big(g\thainv_{\rho,y}\big) \big(g\tha_{\mal}\objd\big) \big(g\tha_{\phi,x}\big) && \text{by \cref{thapsi_geq} and $G$-equivariance of $\tha_{\mal}$}\\
&= g\big((\str \tha)_{\mal}\objd \big) && \text{by \cref{strx_tha_d} and functoriality of $g \cdot -$}
\end{aligned}\]
\par}
\item[Naturality]
The naturality \cref{fgcatg_icell_nat} of $\str\tha$ means the commutativity of the diagram 
\begin{equation}\label{strtha_nat}
\begin{tikzpicture}[vcenter]
\def\v{-1.3} \def\s{1em}
\draw[0cell=.9]
(0,0) node (a11) {(\str X)\mal}
(a11)++(2.8,0) node (a12) {(\str Y)\mal}
(a11)++(0,\v) node (a21) {(\str X)\nbeta}
(a12)++(0,\v) node (a22) {(\str Y)\nbeta}
;
\draw[1cell=.85]
(a11) edge node {(\str\tha)_{\mal}} (a12)
(a12) edge[transform canvas={xshift=-\s}] node {(\str Y)\psi} (a22)
(a11) edge[transform canvas={xshift=\s}] node[swap] {(\str X)\psi} (a21)
(a21) edge node {(\str\tha)_{\nbeta}} (a22)
;
\end{tikzpicture}
\end{equation}
for each pointed morphism $\psi \cn \mal \to \nbeta$ between pointed finite $G$-sets. 
\begin{description}
\item[Objects] By \cref{strx_psi_objx,strtha_mal_obj}, each composite in the diagram \cref{strtha_nat} sends an object $(\phi;x)$ in $(\str X)\mal$ \cref{strx_mal_obj} to the object $(\psi\phi; \tha_{\ldea} x)$ in $(\str Y)\nbeta$.
\item[Morphisms]
The following commutative diagram in $(\str Y)\nbeta$ proves that the diagram \cref{strtha_nat} commutes on each morphism $\objd \cn \objx \to \objy$ in $(\str X)\mal$ \cref{morstrx}, where $X^{-2} = (X^2)^{-1}$ and $Y^{-2} = (Y^2)^{-1}$.
\begin{equation}\label{StYnbeta_diag}
\begin{tikzpicture}[vcenter]
\def\v{1.3} \def\h{2.5} \def\t{22}
\draw[0cell=.8]
(0,0) node (a11) {Y(\psi\phi)(\tha_{\ldea}x)}
(a11)++(0,\v) node (a12) {(Y\psi)(Y\phi)(\tha_{\ldea}x)}
(a12)++(\h,0) node (a13) {(Y\psi)\tha_{\mal}(X\phi)x}
(a13)++(\h,0) node (a14) {(Y\psi)\tha_{\mal}(X\rho)y}
(a14)++(\h,0) node (a15) {(Y\psi)(Y\rho)(\tha_{\rka}y)}
(a15)++(0,-\v) node (a16) {Y(\psi\rho) (\tha_{\rka}y)}
(a11)++(0,-\v) node (a21) {\tha_{\nbeta} X(\psi\phi)x}
(a21)++(\h,0) node (a22) {\tha_{\nbeta} (X\psi)(X\phi)x}
(a22)++(\h,0) node (a23) {\tha_{\nbeta} (X\psi)(X\rho)y}
(a23)++(\h,0) node (a24) {\tha_{\nbeta} X(\psi\rho)y}
;
\draw[1cell=.75]
(a11) edge node {Y^{-2}_{\psi,\phi,\tha_{\ldea}x}} (a12)
(a12) edge[bend left=\t] node {(Y\psi)\tha_{\phi,x}} (a13)
(a13) edge[bend left=\t] node {(Y\psi)\tha_{\mal}\objd} (a14)
(a14) edge[bend left=\t] node {(Y\psi)\thainv_{\rho,y}} (a15)
(a15) edge node {Y^2_{\psi,\rho,\tha_{\rka}y}} (a16)
(a11) edge node[swap] {\tha_{\psi\phi,x}} (a21)
(a21) edge[bend right=\t] node[swap] {\tha_{\nbeta} X^{-2}_{\psi,\phi,x}} (a22)
(a22) edge[bend right=\t] node[swap] {\tha_{\nbeta} (X\psi)\objd} (a23)
(a23) edge[bend right=\t] node[swap] {\tha_{\nbeta} X^2_{\psi,\rho,y}} (a24)
(a24) edge node[swap] {\thainv_{\psi\rho,y}} (a16)
(a13) edge node {\tha_{\psi,(X\phi)x}} (a22)
(a14) edge node {\tha_{\psi,(X\rho)y}} (a23)
;
\end{tikzpicture}
\end{equation}
\begin{itemize}
\item By \cref{strx_psi_mor}, \cref{strx_tha_d}, and the functoriality of $Y\psi$, the top boundary composite is the morphism $\big((\str Y)\psi\big) (\str\tha)_{\mal}\objd$.  Similarly, by the functoriality of $\tha_{\nbeta}$, the bottom boundary composite is the morphism $(\str\tha)_{\nbeta} \big((\str X)\psi\big)\objd$.
\item The left and right regions commute by the compositionality axiom \cref{tha_psiphi} for $\tha$ at, respectively, the objects $x \in X\ldea$ and $y \in X\rka$.  The middle region commutes by the naturality of $\tha_{\psi}$ \cref{psfggcat_one_compii} at the morphism $\objd \cn (X\phi)x \to (X\rho)y$.
\end{itemize}
\end{description}
This finishes the description of the $G$-natural transformation $\str\tha$ \cref{str_tha}.
\item[Preservation of weak $G$-equivalences]
The strictification 2-functor $\str$ preserves weak $G$-equivalences of 1-cells \pcref{def:psfggcat_icell}.  In more detail, suppose $\tha \cn X \to Y$ is a weak $G$-equivalence, meaning that each $\tha_{\mal} \cn X\mal \to Y\mal$ is a categorical weak $G$-equivalence \pcref{def:cat_weakg}.  To show that $\str\tha$ \cref{str_tha} is a weak $G$-equivalence, we need to show that, for each pointed finite $G$-set $\mal$ and subgroup $H \subseteq G$, the $H$-fixed subfunctor
\begin{equation}\label{strtham_H}
\big((\str X)\mal\big)^H \fto{(\str\tha)_{\mal}^H} \big((\str Y)\mal\big)^H
\end{equation}
of $(\str\tha)_{\mal}$ \cref{str_tha_mal} is an equivalence of categories.
\begin{description}
\item[Essential surjectivity]
The functor $(\str\tha)_{\mal}^H$ is essentially surjective on objects by \cref{strx_objx_gact}, \cref{strtha_mal_obj}, the $G$-equivariance of $\tha_{\ldea}$, and the essential surjectivity of the $H$-fixed subfunctor 
\[(X\ldea)^H \fto{\tha_{\ldea}^H} (Y\ldea)^H.\]
\item[Fully faithfulness]
For objects $\objx = (\phi;x)$ and $\objy = (\rho;y)$ in $((\str X)\mal)^H$ \cref{morstrx}, by \cref{strx_tha_d}, the morphism assignment of $(\str\tha)_{\mal}^H$ is the following composite bijection.
\begin{equation}\label{fullfaith_diag}
\begin{tikzpicture}[vcenter]
\def\h{5.2} \def\u{-1} \def\s{.8}
\draw[0cell=.85]
(0,0) node (a11) {\big((\str X)\mal\big)^H(\objx, \objy)}
(a11)++(\h,0) node (a12) {\big((\str Y)\mal\big)^H\big((\str\tha)_{\mal}\objx, (\str\tha)_{\mal}\objy \big)}
(a11)++(0,\u) node (a21) {(X\mal)^H\big((X\phi)x, (X\rho)y\big)}
(a12)++(0,\u) node (a22) {(Y\mal)^H\big((Y\phi)(\tha_{\ldea}x), (Y\rho)(\tha_{\rka}y) \big)}
(a21)++(\h/2,\u) node (a3) {(Y\mal)^H\big(\tha_{\mal}(X\phi)x, \tha_{\mal}(X\rho)y \big)}
;
\draw[1cell=\s]
(a11) edge node {(\str\tha)_{\mal}^H} (a12)
(a11) edge[equal] (a21)
(a12) edge[equal] (a22)
(a21) [rounded corners=2pt] |- node[swap,pos=.3] {\tha_{\mal}^H} node[pos=.3] {\iso} (a3)
;
\draw[1cell=\s]
(a3) [rounded corners=2pt] -| node[swap,pos=.7] {\thainv_{\rho,y} (-) \tha_{\phi,x}} node [pos=.7] {\iso} (a22)
;
\end{tikzpicture}
\end{equation}
The morphism assignment of the equivalence of categories 
\[(X\mal)^H \fto{\tha_{\mal}^H} (Y\mal)^H\] 
is a bijection.  The function $\thainv_{\rho,y} (-) \tha_{\phi,x}$ is a bijection because $\thainv_{\rho,y}$ and $\tha_{\phi,x}$ are isomorphisms in $Y\mal$.
\end{description}
Thus, $(\str\tha)_{\mal}^H$ is an equivalence of categories for each object $\mal \in \FG$ and subgroup $H \subseteq G$, proving that $\str\tha$ is a weak $G$-equivalence.\defmark
\end{description}
\end{explanation}

\begin{explanation}[$\str$ on 2-Cells]\label{expl:str_iicell}
Suppose $\Theta \cn \tha \to \ups$ is a $\Gcatst$-modification \pcref{def:psfggcat_iicell} between $\Gcatst$-pseudotransformations $\tha, \ups \cn X \to Y$ between pseudo $\FGG$-categories $X$ and $Y$.  The strictification 2-functor $\str$ \cref{str_iifunctor} sends $\Theta$ to the $G$-modification \cref{fgcatg_iicell}
\begin{equation}\label{str_iicell}
\begin{tikzpicture}[vcenter]
\def\t{25} \def\u{.6} \def\h{3.5}
\draw[0cell]
(0,0) node (a) {\FG}
(a)++(\h,0) node (b) {\Catgst}
;
\draw[1cell=.8]
(a) [rounded corners=2pt] |- ($(a)+(1,\u)$) -- node[pos=.5] {\str X} ($(b)+(-1,\u)$) -| (b)
;
\draw[1cell=.8]
(a) [rounded corners=2pt] |- ($(a)+(1,-\u)$) -- node[pos=.5] {\str Y} ($(b)+(-1,-\u)$) -| (b)
;
\draw[2cell=.9]
node[between=a and b at .34, rotate=-90, 2label={below,\str\tha}] {\Rightarrow}
node[between=a and b at .65, rotate=-90, 2label={above,\str\ups}] {\Rightarrow}
;
\draw[2cell=.9]
node[between=a and b at .52, rotate=0, shift={(0,-.05)}, 2labelmed={above,\!\str\Theta}] {\Rrightarrow}
;
\end{tikzpicture}
\end{equation}
with an $\mal$-component pointed $G$-natural transformation \cref{fgcatg_iicell_comp}
\begin{equation}\label{str_iicell_comp}
\begin{tikzpicture}[vcenter]
\def\t{25}
\draw[0cell]
(0,0) node (a1) {\phantom{X}}
(a1)++(2.5,0) node (a2) {\phantom{X}}
(a1)++(-.52,0) node (a1') {(\str X)\mal}
(a2)++(.52,0) node (a2') {(\str Y)\mal}
;
\draw[1cell=.8]
(a1) edge[bend left=\t] node {(\str\tha)_{\mal}} (a2)
(a1) edge[bend right=\t] node[swap] {(\str\ups)_{\mal}} (a2)
;
\draw[2cell=.8]
node[between=a1 and a2 at .3, rotate=-90, 2label={above, (\str\Theta)_{\mal}}] {\Rightarrow}
;
\end{tikzpicture}
\end{equation}
for each pointed finite $G$-set $\mal$.  Using \cref{morstrx,strtha_mal_obj}, for each object $\objx = (\phi; x)$ in $(\str X)\mal$ \cref{strx_mal_obj}, the $\objx$-component morphism
\[(\str\tha)_{\mal}\objx = \big(\phi; \tha_{\ldea}x\big) \fto{(\str\Theta)_{\mal,\objx}} 
(\str\ups)_{\mal}\objx = \big(\phi; \ups_{\ldea}x \big)\]
of $(\str\Theta)_{\mal}$ is given by the morphism
\begin{equation}\label{str_iicell_mx}
(Y\phi)(\tha_{\ldea}x) \fto{(Y\phi) \Theta_{\ldea,x}} (Y\phi)(\ups_{\ldea}x) \inspace Y\mal.
\end{equation}
In the preceding morphism, $\Theta_{\ldea,x}$ is the $x$-component of the $\ldea$-component pointed $G$-natural transformation $\Theta_{\ldea} \cn \tha_{\ldea} \to \ups_{\ldea}$ \cref{psfggcat_iicell_comp}.
\begin{description}
\item[Equivariance] The pointed $G$-equivariance of $(\str\Theta)_{\mal}$ follows from the pointed $G$-equivariance of $\Theta_{\ldea}$, the pointed functoriality of $Y\phi$, \cref{psfggcat_basept}, and \cref{psfggcat_geq}.
\item[Naturality]
The following commutative diagram in $Y\mal$ proves that $(\str\Theta)_{\mal}$ is natural in morphisms $\objd \cn \objx \to \objy$ in $(\str X)\mal$ \cref{morstrx}.
\begin{equation}\label{StTheta_natdiag}

\end{equation}
The top and bottom regions are the definitions \cref{strx_tha_d} of $(\str\tha)_{\mal}\objd$ and $(\str\ups)_{\mal}\objd$.  The left and right regions commute by the modification axiom \cref{psfggcat_iicell_modax} for $\Theta$.  The middle region commutes by the naturality of $\Theta_{\mal}$ \cref{psfggcat_iicell_comp}.
\item[Modification axiom]
The axiom \cref{fgcatg_iicell_modax} for $\str\Theta$ states that the two whiskered natural transformations in the diagram
\begin{equation}\label{str_iicell_modax}
%
\end{equation}
are equal for each pointed morphism $\psi \cn \mal \to \nbeta$ in $\FG$.  To see that \cref{str_iicell_modax} commutes, for each object $\objx = (\phi; x)$ in $(\str X)\mal$ \cref{strx_mal_obj}, we consider the following commutative diagram in $Y\nbeta$, where $Y^{-2} = (Y^2)^{-1}$.
\begin{equation}\label{str_modax_x}
\def\u{.2ex}
%
\end{equation}
By \cref{strx_psi_objx,strx_psi_mor,str_iicell_mx}, the bottom horizontal arrow and the left-top-right composite in \cref{str_modax_x} are, respectively, the morphisms 
\[\big(\str\Theta\big)_{\nbeta,\, ((\str X)\psi) \objx} \andspace 
\big((\str Y)\psi\big) (\str\Theta)_{\mal,\objx}.\]
By the naturality of $Y^2_{\psi,\phi}$ \cref{psfggcat_xtwo_comp}, the diagram \cref{str_modax_x} commutes, proving the modification axiom \cref{str_iicell_modax} for $\str\Theta$.
\end{description}
This finishes the description of the $G$-modification $\str\Theta$ \cref{str_iicell}.
\end{explanation}

\begin{explanation}[Unit and Counit]\label{expl:str_unit}
The unit 
\begin{equation}\label{stru}
\begin{tikzpicture}[vcenter]
\def\t{27}
\draw[0cell]
(0,0) node (a1) {\phantom{X}}
(a1)++(1.8,0) node (a2) {\phantom{X}}
(a1)++(-.7,0) node (a1') {\FGCatgps}
(a2)++(.67,0) node (a2') {\FGCatgps}
;
\draw[1cell=.8]
(a1) edge[bend left=\t] node {1} (a2)
(a1) edge[bend right=\t] node[swap] {\Incj\str} (a2)
;
\draw[2cell=.9]
node[between=a1 and a2 at .43, rotate=-90, 2label={above, \stru}] {\Rightarrow}
;
\end{tikzpicture}
\end{equation}
of the 2-adjunction \cref{str_iifunctor}
\begin{equation}\label{StIn_expl_diag}
\begin{tikzpicture}[vcenter]
\draw[0cell]
(0,0) node (a1) {\FGCatgps}
(a1)++(2.7,0) node (a2) {\FGCatg}
;
\draw[1cell=.9]
(a1) edge[transform canvas={yshift=.5ex}] node {\str} (a2)
(a2) edge[transform canvas={yshift=-.4ex}] node {\Incj} (a1)
;
\end{tikzpicture}
\end{equation}
sends a pseudo $\FGG$-category $X$ \pcref{def:psfggcat} to the $\Gcatst$-pseudotransformation \pcref{def:psfggcat_icell}
\begin{equation}\label{strux}
X \fto{\stru_X} \Incj\str X
\end{equation}
defined by the data \cref{strux_mal,strux_psi}.
\begin{description}
\item[$G$-functors]
For each pointed finite $G$-set $\mal$, the $\mal$-component pointed $G$-functor
\begin{equation}\label{strux_mal}
X\mal \fto{\stru_{X,\mal}} (\Incj\str X)\mal = (\str X)\mal
\end{equation}
sends an object $x \in X\mal$ to the object \cref{strx_mal_obj}
\[\stru_{X,\mal}(x) = (1_{\mal}; x) \in (\str X) \mal.\]
Using \cref{psfggcat_unity,morstrx}, $\stru_{X,\mal}$ sends a morphism $d \in (X\mal)(x,y)$ to the morphism
\[\begin{split}
d & \in (X\mal)\big((X1_{\mal})x, (X1_{\mal})y \big)\\
&= \big((\str X)\mal\big)\big(\stru_{X,\mal}(x), \stru_{X,\mal}(y)\big).
\end{split}\]
\item[Natural isomorphisms]
For each pointed morphism $\psi \cn \mal \to \nbeta$, the $\psi$-component pointed natural isomorphism
\begin{equation}\label{strux_psi}
\begin{tikzpicture}[vcenter]
\def\v{-1.4}
\draw[0cell]
(0,0) node (a11) {X\mal}
(a11)++(2.5,0) node (a12) {(\str X)\mal}
(a11)++(0,\v) node (a21) {X\nbeta}
(a12)++(0,\v) node (a22) {(\str X)\nbeta}
;
\draw[1cell=.9]
(a11) edge node {\stru_{X,\mal}} (a12)
(a12) edge[transform canvas={xshift=-1em}] node {(\str X)\psi} (a22)
(a11) edge node[swap] {X\psi} (a21)
(a21) edge node[swap] {\stru_{X,\nbeta}} (a22)
;
\draw[2cell=.9]
node[between=a12 and a21 at .5, shift={(0,-.05)}, rotate=240, 2label={above,\scalebox{.8}{$\iso$}}, 2labelw={below,\stru_{X,\psi},0pt}] {\Rightarrow}
;
\end{tikzpicture}
\end{equation}
sends an object $x \in X\mal$ to the isomorphism
\[\begin{tikzpicture}
\draw[0cell]
(0,0) node (a1) {\big((\str X)\psi\big) \big(\stru_{X,\mal} x\big) = (\psi 1_{\mal}; x)}
(a1)++(0,-1.3) node (a2) {(\stru_{X,\nbeta}) (X\psi)x = \big(1_{\nbeta}; (X\psi)x\big) }
;
\draw[1cell=.9]
(a1) edge[shorten >=-.4ex] node {\stru_{X,\psi,x}} (a2)
;
\end{tikzpicture}\]
in $(\str X)\nbeta$ given by the identity morphism 
\[X(\psi 1_{\mal})x =  (X\psi)x \fto{1_{(X\psi)x}} (X\psi)x = (X1_{\nbeta})(X\psi)x\]
in $X\nbeta$.
\item[Inverse $G$-functors]
The pointed $G$-functor $\stru_{X,\mal}$ \cref{strux_mal} is a pointed $G$-equivalence.  Its inverse pointed $G$-equivalence
\begin{equation}\label{strvx_mal}
(\str X)\mal \fto{\strv_{X,\mal}} X\mal
\end{equation}
sends an object $\objx = (\phi;x) \in (\str X)\mal$ \cref{strx_mal_obj} to the object
\begin{equation}\label{strvxm_obj}
\strv_{X,\mal}(\phi;x) = (X\phi)x \in X\mal.
\end{equation}
On morphisms, $\strv_{X,\mal}$ is defined by the equality \cref{morstrx}. 
\begin{itemize}
\item The composite $\strv_{X,\mal}\stru_{X,\mal}$ is the identity on $X\mal$. 
\item There is a $G$-natural isomorphism 
\begin{equation}\label{strwXmal_diag}
\begin{tikzpicture}[vcenter]
\def\t{30}
\draw[0cell]
(0,0) node (a1) {\phantom{X}}
(a1)++(2,0) node (a2) {\phantom{X}}
(a1)++(-.52,0) node (a1') {(\str X)\mal}
(a2)++(.52,0) node (a2') {(\str X)\mal}
;
\draw[1cell=.9]
(a1) edge[bend left=\t] node {1} (a2)
(a1) edge[bend right=\t] node[swap] {\stru_{\mal} \strv_{\mal}} (a2)
;
\draw[2cell=.85]
node[between=a1 and a2 at .33, rotate=-90, 2label={above, \strw_{X,\mal}}] {\Rightarrow}
;
\end{tikzpicture}
\end{equation}
that sends an object $\objx = (\phi;x) \in (\str X)\mal$ to the isomorphism
\[(\phi;x) \fto{\strw_{X,\mal,\objx}} \stru_{X,\mal}\strv_{X,\mal}\objx = \big(1_{\mal}; (X\phi)x\big)\]
in $(\str X)\mal$ given by the identity morphism
\[(X\phi)x \fto{1_{(X\phi)x}} (X\phi)x = (X1_{\mal})(X\phi)x\]
in $X\mal$.
\end{itemize}
\item[Counit]
The counit 
\begin{equation}\label{strv}
\begin{tikzpicture}[vcenter]
\def\t{27}
\draw[0cell]
(0,0) node (a1) {\phantom{X}}
(a1)++(1.8,0) node (a2) {\phantom{X}}
(a1)++(-.4,0) node (a1') {\FGCatg}
(a2)++(.4,0) node (a2') {\FGCatg}
;
\draw[1cell=.8]
(a1) edge[bend left=\t] node {\str\Incj} (a2)
(a1) edge[bend right=\t] node[swap] {1} (a2)
;
\draw[2cell=.9]
node[between=a1 and a2 at .43, rotate=-90, 2label={above, \strv}] {\Rightarrow}
;
\end{tikzpicture}
\end{equation}
of the 2-adjunction $(\str,\Incj)$ sends an $\FGG$-category $X$ \cref{fgcatg_obj} to the $G$-natural transformation
\begin{equation}\label{strvx}
\begin{tikzpicture}[vcenter]
\def\t{27}
\draw[0cell]
(0,0) node (a1) {\phantom{X}}
(a1)++(2,0) node (a2) {\phantom{X}}
(a1)++(-.08,0) node (a1') {\FG}
(a2)++(.22,0) node (a2') {\Catgst}
;
\draw[1cell=.8]
(a1) edge[bend left=\t] node {\str\Incj X} (a2)
(a1) edge[bend right=\t] node[swap] {X} (a2)
;
\draw[2cell=.9]
node[between=a1 and a2 at .4, rotate=-90, 2label={above, \strv_X}] {\Rightarrow}
;
\end{tikzpicture}
\end{equation}
whose $\mal$-component pointed $G$-equivalence is $\strv_{X,\mal}$ \cref{strvx_mal} for each object $\mal \in \FG$.  To apply the construction $\strv_{X,\mal}$ to an $\FGG$-category $X$, we regard $X$ as a pseudo $\FGG$-category via the inclusion $\Incj$ with $X^2 = 1$ \cref{psfggcat_xtwo} and abbreviate $\Incj X$ to $X$.
\end{description}
This finishes the description of the 2-adjunction $(\str,\Incj,\stru,\strv)$.
\end{explanation}

\section{GMMO $K$-Theory}
\label{sec:kgmmo}

This section reviews GMMO $K$-theory $\Kgmmo$ \cite[(1.2)]{gmmo23} and observes that it preserves weak $G$-equivalences \pcref{kgmmo_weq}.

\begin{definition}\label{def:kgmmo}
For a finite group $G$ and a chaotic $\Einfg$-operad $\Op$ \pcref{def:chaotic_einf}, \emph{GMMO $K$-theory}\index{GMMO K-theory@GMMO $K$-theory}\index{K-theory@$K$-theory!GMMO} $\Kgmmo$ is the composite functor
\begin{equation}\label{kgmmo_diag}
\begin{tikzpicture}[vcenter]
\def\v{1.4} \def\h{2.8}
\draw[0cell]
(0,0) node (a1) {\AlgstpsO}
(a1)++(0,-\v) node (a2) {\DGCatg}
(a2)++(0,-\v) node (a3) {\FGCatgps}
(a3)++(\h,0) node (a4) {\FGCatg}
(a4)++(\h,0) node (a5) {\FGTopg}
(a5)++(0,\v) node (a6) {\FGTopg}
(a6)++(0,\v) node (a7) {\Gspec}
;
\draw[1cell=.9]
(a1) edge node {\Kgmmo} (a7)
(a1) edge node[swap] {\Rg} (a2)
(a2) edge node[swap] {\gzest} (a3)
(a3) edge node {\str} (a4)
(a4) edge node {\clast} (a5)
(a5) edge node[swap] {\Bc} (a6)
(a6) edge node[swap] {\Kfg} (a7)
;
\end{tikzpicture}
\end{equation}
defined as follows.
\begin{enumerate}
\item The 2-categories $\DGCatg$ and $\AlgstpsO$ are defined in \cref{def:dgcatg_iicat,def:algstpsO}.  The 2-functor $\Rg$ is defined in \cref{def:rg_obj,def:rg_onecell,def:rg_twocell} for objects, 1-cells, and 2-cells.
\item The 2-category $\FGCatgps$ is defined in \cref{def:psfggcat,def:psfggcat_icell,def:psfggcat_iicell}.  The 2-functor $\gzest$ is defined in \cref{def:gzest}.
\item The strictification 2-functor $\str$ is defined in \cref{def:str_iifunctor}.
\item The functor $\clast$ is induced by the classifying space functor $\cla$ \pcref{ggcatg_ggtopg}.
\item The bar functor $\Bc$ on $\FGTopg$ is defined in \cref{bar_functor_FG}.
\item The prolongation functor $\Kfg$ is defined in \cref{def:Kfg_functor}.\defmark
\end{enumerate}
\end{definition}

\begin{explanation}[Unpacking]\label{expl:kgmmo}
The last three constituent functors of $\Kgmmo$---$\clast$, $\Bc$, and $\Kfg$---are the same ones as in the homotopical Shimakawa (strong) $K$-theory \pcref{def:khsho}.  By \cref{Kfgxv,real_FG}, for each $\Op$-algebra $\A$, the orthogonal $G$-spectrum $\Kgmmo\A$ sends each object $V \in \IU$ \pcref{def:IU_spaces} to the pointed $G$-space
\[(\Kgmmo\A)_V = 
\int^{\mal \sins \FG} (S^V)^{\mal} \sma \Bc\big(\FG(-,\mal),\FG,\clast\str\gzest\Rg\A\big).\]
The sphere action on $\Kgmmo\A$ is defined in \cref{Kfgx_action_vw}.  The $G$-action on $(\Kgmmo\A)_V$ is defined in \cref{Kfgxv_rep_gact,Br_FGX,strx_objx_gact,morstrx,proCnbe_gaction}.  By \cref{rga_mal,gzest_xmal}, the pseudo $\FGG$-category $\gzest\Rg\A \in \FGCatgps$ sends each pointed finite $G$-set $\nbeta$ to the $\nbeta$-twisted product \pcref{def:proCnbe}
\[(\gzest\Rg\A)\nbeta = (\Rg\A)\nbeta = \proAnbeta.\]
The strictified $\FGG$-category $\str\gzest\Rg\A$ is discussed in \cref{expl:str_obj} with $X = \gzest\Rg\A$.  The $\FGG$-space $\clast\str\gzest\Rg\A$ sends $\nbeta$ to the pointed $G$-space
\[(\clast\str\gzest\Rg\A)\nbeta = \cla (\str\gzest\Rg\A)\nbeta,\]
which is the classifying space \cref{classifying_space} of the pointed $G$-category $(\str\gzest\Rg\A)\nbeta$.
\end{explanation}

\begin{definition}[Weak $G$-Equivalences]\label{def:kgmmo_weakg}
Each category in the diagram \cref{kgmmo_diag} has a natural notion of \emph{weak $G$-equivalence}\index{weak G-equivalence@weak $G$-equivalence} for 1-cells or morphisms.
\begin{itemize}
\item A 1-cell in $\AlgstpsO$ is an $\Op$-pseudomorphism \pcref{def:laxmorphism} between $\Op$-algebras.  It is called a \emph{weak $G$-equivalence} if its underlying pointed $G$-functor \cref{fAB} is a categorical weak $G$-equivalence \pcref{def:cat_weakg}.
\item A 1-cell in $\DGCatg$ is a $\Pig$-strict $\Gcatst$-pseudotransformation \pcref{def:dgcatg_onecell,def:pig_strict}.  It is called a \emph{weak $G$-equivalence} if its component pointed $G$-functors \cref{dgcatg_one_comp} are categorical weak $G$-equivalences.  
\item A 1-cell in $\FGCatgps$ is a $\Gcatst$-pseudotransformation \pcref{def:psfggcat_icell}.  It is called a \emph{weak $G$-equivalence} if its component pointed $G$-functors \cref{psfggcat_one_comp} are categorical weak $G$-equivalences.  The same definition also applies to $\FGCatg$.
\item A morphism in $\FGTopg$ \pcref{def:ggtopg} is a $G$-natural transformation.  It is called a \emph{weak $G$-equivalence} if its components \cref{ggtopg_icell_comp} are weak $G$-equivalences of $G$-spaces \pcref{def:weakG_top}.
\item A morphism in $\Gspec$ is a $G$-morphism of orthogonal $G$-spectra (\cref{def:gsp_morphism} \cref{def:gsp_morphism_ii}).  It is called a \emph{weak $G$-equivalence} if its component pointed $G$-morphisms \cref{eqiu_mor_component} are weak $G$-equivalences of $G$-spaces.\defmark
\end{itemize}
\end{definition}

The prolongation functor $\Kfg$ \pcref{def:Kfg_functor} does \emph{not} generally preserve weak $G$-equivalences, but everything else in $\Kgmmo$ does.

\begin{lemma}\label{kgmmo_weq}
In the diagram \cref{kgmmo_diag}, each of the five functors 
\[\Rg \csp \gzest \csp \str \csp \clast \csp \text{and}\,\, \Kfg\Bc\]
preserves weak $G$-equivalences.  Thus, their composite $\Kgmmo$ also preserves weak $G$-equivalences.
\end{lemma}

\begin{proof}
The functors $\gzest$, $\str$, and $\clast$ preserve weak $G$-equivalences by, respectively, \cref{gzest_icell_comp}, \cref{strtham_H}, and \cref{ex:catweakg_weakg}.
\begin{description}
\item[$\Rg$]
Suppose $f \cn \A \to \B$ is a weak $G$-equivalence in $\AlgstpsO$.  To show that $\Rg f$ is a weak $G$-equivalence in $\DGCatg$, recall from \cref{rgf_mal} that, for each object $\mal \in \DG$, the $\mal$-component $(\Rg f)_{\mal}$ is $f^m \cn \proAmal \to \proBmal$.   We need to show that $f^m$ is a categorical weak $G$-equivalence.  By \cref{twprod_fixed}, for each subgroup $H \subseteq G$, there is a commutative diagram 
\begin{equation}\label{RgweakG_diag}
\begin{tikzpicture}[vcenter]
\def\v{-1.4}
\draw[0cell]
(0,0) node (a11) {(\proAmal)^H}
(a11)++(3.3,0) node (a12) {(\proBmal)^H}
(a11)++(0,\v) node (a21) {\txprod_{t\in \ufsr}\, \A^{H_t}}
(a12)++(0,\v) node (a22) {\txprod_{t\in \ufsr}\, \B^{H_t}}
;
\draw[1cell=.9]
(a11) edge node {(f^m)^H} (a12)
(a12) edge node {\uph} node[swap] {\iso} (a22)
(a11) edge node[swap] {\uph} node {\iso} (a21)
(a21) edge node {\txprod_{t \in \ufsr}\, f^{H_t}} (a22)
;
\end{tikzpicture}
\end{equation}
of functors, decomposing the $H$-fixed subfunctor $(f^m)^H$, up to isomorphisms, into a finite product $\prod_{t \in \ufsr} f^{H_t}$ for some subgroups $H_t \subseteq H$.  Since $f$ is a categorical weak $G$-equivalence, each $H_t$-fixed subfunctor 
\[\A^{H_t} \fto{f^{H_t}} \B^{H_t}\]
is an equivalence of categories.  Thus, $\prod_{t \in \ufsr} f^{H_t}$ and $(f^m)^H$ are also equivalences of categories, proving that $(\Rg f)_{\mal} = f^m$ is a categorical weak $G$-equivalence.
\item[$\Kfg\Bc$]
For each morphism $\tha \cn X \to Y$ in $\FGTopg$, the naturality of the retraction $\retn \cn \Bc \to 1_{\FGTopg}$ \cref{retn_barFG_id} yields a commutative diagram
\begin{equation}\label{Barnat_diag}
\begin{tikzpicture}[vcenter]
\def\v{-1.3}
\draw[0cell]
(0,0) node (a11) {\Bc X}
(a11)++(2,0) node (a12) {\Bc Y}
(a11)++(0,\v) node (a21) {X}
(a12)++(0,\v) node (a22) {Y}
;
\draw[1cell=.9]
(a11) edge node {\Bc\tha} (a12)
(a12) edge node {\retn_Y} (a22)
(a11) edge node[swap] {\retn_X} (a21)
(a21) edge node {\tha} (a22)
;
\end{tikzpicture}
\end{equation}
in $\FGTopg$.  Each of $\retn_X$ and $\retn_Y$ \cref{ret_FG} is a weak $G$-equivalence in $\FGTopg$.   By \cref{bar_proper}, for each $\FGG$-space $Z$, the bar construction $\Bc Z$ \cref{bar_FG} is proper in the sense of \cref{def:proper_fgg}.  Thus, if $\tha$ is a weak $G$-equivalence, then $\Bc\tha$ is a weak $G$-equivalence between proper $\FGG$-spaces.  \cref{thm:Kfg_inv} implies that $\Kfg\Bc\tha$ is a weak $G$-equivalence in $\Gspec$.
\end{description}
This finishes the proof.
\end{proof}

\begin{remark}[Multifunctoriality]\label{rk:gmmo_mult}
This chapter and this work in general do not use any of the claims of multifunctoriality in \cite{gmmo23}.  However, we note that some assertions of multifunctoriality in \cite{gmmo23} require corrections.  It is claimed in \cite[Theorem A]{gmmo23} that the functor $\Kgmmo$ \cref{kgmmo_diag} extends to a nonsymmetric multifunctor between symmetric multicategories.  However, the intermediate constructions 
\[\gmmofgps \andspace \gmmocalg\]
are \emph{not} symmetric multicategories, as claimed in \cite[7.14 and 8.12]{gmmo23}, but only nonsymmetric multicategories.  In each case, the symmetric group action on $k$-ary morphisms for $k>1$, as constructed in the bottom diagram in \cite[p.\ 46]{gmmo23}, is not well defined.  Specifically, using the notation in \cite[(5.18)]{gmmo23}, for a $k$-ary morphism 
\[F \cn (\mathcal{X}_1,\ldots, \mathcal{X}_k) \to \mathcal{Y}\] 
and permutations $\si,\si' \in \Si_k$, the $k$-ary morphisms $(F\si)\si'$ and $F(\si\si')$ are not generally equal because $\mathcal{Y}$ is not strictly functorial even on permutations.
\end{remark}

%% file: chap/gmmo_shi.tex
\section{Comparison 2-Natural Transformation}
\label{sec:cgs}

This section constructs a comparison 2-natural transformation
\begin{equation}\label{cgs_secti}
\begin{tikzpicture}[vcenter]
\def\t{30} \def\h{2.7} \def\v{1.4}
\draw[0cell]
(0,0) node (a1) {\AlgstpsO}
(a1)++(0,\v) node (a2) {\AlgpspsO}
(a2)++(\h,0) node (a3) {\FGCatg}
(a3)++(.2,0) node (a3') {\phantom{\FGCatg}}
(a3')++(0,-\v) node (a4) {\DGCatg}
;
\draw[1cell=.85]
(a1) edge[right hook->] node {\algi} (a2)
(a2) edge node {\Sgosg} (a3)
(a3') edge[transform canvas={xshift=-1ex}] node {\gxist} (a4)
(a1) edge node[swap] {\Rg} (a4)
;
\draw[2cell=1]
node[between=a2 and a3 at .42, shift={(0,-.5*\v)}, rotate=-90, 2label={above,\cgs}] {\Rightarrow}
;
\end{tikzpicture}
\end{equation}
for a finite group $G$ and a 1-connected $\Gcat$-operad $(\Op,\ga,\opu)$ \pcref{as:Op_iconn}.  The comparison 2-natural transformation $\cgs$ is used in \cref{sec:gmmo_shi_weq} to compare the homotopical Shimakawa strong $K$-theory and GMMO $K$-theory.

\secoutline
\begin{itemize}
\item \cref{def:cgs} constructs $\cgs$.
\item \cref{cgs_welldef} proves that components of $\cgs$ are 1-cells in $\DGCatg$.  
\item \cref{cgs_iinat} proves that $\cgs$ is 2-natural.
\end{itemize}

For \cref{def:cgs}, recall the following.
\begin{itemize}
\item The first step of $\Kgmmo$ is the 2-functor \pcref{sec:rg} 
\[\AlgstpsO \fto{\Rg} \DGCatg\]
between the 2-categories in \cref{def:dgcatg_iicat,def:algstpsO}.  There is a full sub-2-category inclusion
\[\AlgstpsO \fto{\algi} \AlgpspsO\]
with codomain $\AlgpspsO$ \pcref{oalgps_twocat}.
\item Shimakawa strong $H$-theory for $\Op$ is the 2-functor \pcref{def:sgo}
\[\AlgpspsO \fto{\Sgosg} \FGCatg\]
with codomain $\FGCatg$ \pcref{def:fgcatg}.  Its object assignment $\A \mapsto \Sgosg\A = \Asgdash$ is given by \cref{sys_FGcat}.
\item Precomposing with the $\Gcatst$-functor $\gxi \cn \DG \to \FG$ \cref{dgcoop} yields the pullback 2-functor \cref{gxist}
\[\FGCatg \fto{\gxist} \dgcatg \subseteq \DGCatg.\]
\end{itemize}

\begin{definition}\label{def:cgs}
For a finite group $G$, a 1-connected $\Gcat$-operad $(\Op,\ga,\opu)$, and an $\Op$-algebra $(\A,\gaA)$ \pcref{def:pseudoalgebra}, the $\Gcatst$-pseudotransformation \pcref{def:dgcatg_onecell} between $\DGG$-categories 
\begin{equation}\label{cgsa}
\gxist\Sgosg\algi\A \fto{\cgs_\A} \Rg\A
\end{equation}
is defined as follows.
\begin{description}
\item[$G$-functors]
Using \cref{rga_mal} and the fact that $\gxi \cn \DG \to \FG$ is the identity on objects, the $\mal$-component pointed $G$-functor \cref{dgcatg_one_comp} of $\cgs_\A$ is defined by the following commutative diagram for each object $\mal \in \DG$ \pcref{def:dgo}.
\begin{equation}\label{cgsam}
\begin{tikzpicture}[vcenter]
\def\v{-1}
\draw[0cell]
(0,0) node (a11) {(\gxist\Sgosg\algi\A)\mal}
(a11)++(3.3,0) node (a12) {(\Rg\A)\mal}
(a11)++(0,\v) node (a21) {\Asgmal}
(a12)++(0,\v) node (a22) {\phantom{\proAmal}}
(a22)++(0,.05) node (a22') {\proAmal}
;
\draw[1cell=.9]
(a11) edge node {\cgs_{\A,\mal}} (a12)
(a21) edge node {\zbsg} (a22)
(a11) edge[equal] (a21)
(a12) edge[equal] (a22')
;
\end{tikzpicture}
\end{equation}
\begin{itemize}
\item $\Asgmal$ is the small pointed $G$-category of strong $\mal$-systems \cref{Asgordnbe}.  It is the terminal $G$-category $\bone$ if $\mal = \ordz$ (\cref{ex:nsys_zero_one} \cref{ex:nsys_zero}).
\item $\proAmal$ is the $\mal$-twisted product \pcref{def:proCnbe}, regarding $\mal$ as a length-1 object in $\GG$.  By definition, $\proAzero = \bone$, so $\cgs_{\A,\ordz}$ is necessarily the identity on $\bone$.
\item $\cgs_{\A,\mal} = \zbsg$ is the pointed $G$-functor defined in \cref{sgosgtoprod}.  It sends a strong $\mal$-system $(a,\gl)$ \cref{nsys} and a morphism $\tha$ \cref{nsys_mor} in $\Asgmal$ to, respectively, the $m$-tuples
\begin{equation}\label{cgsam_def}
\begin{split}
\cgs_{\A,\mal}(a,\gl) &= \ang{a_{\{i\}}}_{i \in \ufsm} \andspace\\
\cgs_{\A,\mal}\tha &= \ang{\tha_{\{i\}}}_{i \in \ufsm} \quad \text{in $\proAmal$.}
\end{split}
\end{equation}
Note that $\cgs_{\A,\mal}$ is the same as the pointed $G$-functor $\zdsg$ \cref{hgosgtoprod} because $\mal$ is a length-1 object in $\GG$.  
\end{itemize}
\cref{Pist_to_prod} \cref{Pist_to_prod_i} proves that $\cgs_{\A,\mal}$ is a pointed $G$-functor.
\item[Natural isomorphisms]
For each object $\objx = (\psi; \ang{x_j}_{j \in \ufsn})$ in $\DG(\mal,\nbeta)$ \cref{dgmn_obj}, $\gxi\objx$ is the pointed morphism $\psi \cn \mal \to \nbeta$.  The $\objx$-component pointed natural isomorphism \cref{dgcatg_one_compii} of $\cgs_\A$ is defined as follows, where $\Asgpsi$ and $(\Rg\A)\objx$ are the pointed functors in \cref{psitilsg_f,rga_objx}.
\begin{equation}\label{cgsax}
\begin{tikzpicture}[vcenter]
\def\v{-1.4}
\draw[0cell]
(0,0) node (a11) {\Asgmal}
(a11)++(2.7,0) node (a12) {\phantom{\proAmal}}
(a12)++(0,.05) node (a12') {\proAmal}
(a11)++(0,\v) node (a21) {\Asgnbeta}
(a12)++(0,\v) node (a22) {\phantom{\proAnbeta}}
(a22)++(0,.05) node (a22') {\proAnbeta}
;
\draw[1cell=.9]
(a11) edge node {\cgs_{\A,\mal}} (a12)
(a12') edge node {(\Rg\A)\objx} (a22')
(a11) edge node[swap] {\Asgpsi} (a21)
(a21) edge node[swap] {\cgs_{\A,\nbeta}} (a22)
;
\draw[2cell]
node[between=a12 and a21 at .5, shift={(.2,-.05)}, rotate=-120, 2label={above,\scalebox{.8}{$\iso$}}, 2labelalt={below,\cgs_{\A,\objx}}] {\Rightarrow}
;
\end{tikzpicture}
\end{equation}
Using \cref{apsitil_s,rgaxa,cgsam_def}, the component of $\cgs_{\A,\objx}$ at a strong $\mal$-system $(a,\gl) \in \Asgmal$ is the following isomorphism in $\proAnbeta$ defined by the gluing isomorphism $\gl$ \cref{gl-morphism}.
\begin{equation}\label{cgsax_comp}
\begin{tikzpicture}[vcenter]
\def\v{-1.4}
\draw[0cell=.9]
(0,0) node (a11) {\big((\Rg\A)\objx\big) (\cgs_{\A,\mal}) (a,\gl)}
(a11)++(4.2,0) node (a12) {\bang{\gaA\big(x_j; \ang{a_{\{i\}}}_{i \in \psiinv j} \big)}_{j \in \ufsn}}
(a11)++(0,\v) node (a21) {(\cgs_{\A,\nbeta}) (\Asgpsi) (a,\gl)}
(a12)++(0,\v) node (a22) {\ang{a_{\psiinv j}}_{j \in \ufsn}}
;
\draw[1cell=.85]
(a11) edge[equal] (a12)
(a21) edge[equal] (a22)
(a11) edge node[swap] {\cgs_{\A,\objx,(a,\gl)}} (a21)
(a12) edge node[swap] {\ang{\gl_{x_j;\, \psiinv j,\, \ang{\{i\}}_{i \in \psiinv j}}}_{j \in \ufsn}} node {\iso} (a22)
;
\end{tikzpicture}
\end{equation}
\begin{itemize}
\item The pointedness of $\cgs_{\A,\objx}$ follows from the fact that each gluing morphism of the base $\mal$-system $(\zero, 1_{\zero})$ is the identity morphism of the basepoint $\zero \in \A$ \cref{pseudoalg_zero}. 
\item The naturality of $\cgs_{\A,\objx}$ follows from the compatibility axiom \cref{nsys_mor_compat} for morphisms in $\Asgmal$.
\end{itemize}
\end{description}
This finishes the definition of $\cgs_\A$.
\end{definition}

Recall that a 1-cell in $\DGCatg$ is a $\Pig$-strict $\Gcatst$-pseudotransformation \pcref{def:dgcatg_onecell,def:pig_strict}.

\begin{lemma}\label{cgs_welldef}
$\cgs_\A$ \cref{cgsa} is a 1-cell in $\DGCatg$.
\end{lemma}

\begin{proof}
In this proof, $(a,\gl) \in \Asgmal$ denotes an arbitrary strong $\mal$-system in $\A$.
We check each 1-cell axiom for $\cgs_\A$.
\begin{description}
\item[$\Pig$-strictness]
If $\objx \in \DG(\mal,\nbeta)$ is in the image of $\gio \cn \Pig \to \DG$ \cref{iota_psi}, then each gluing morphism $\gl$ in \cref{cgsax_comp} is the identity by the unity axioms \cref{nsys_unity_empty,nsys_unity_one} for $(a,\gl)$. 
\item[Naturality]
By \cref{cgsax_comp}, \cref{rgafa}, and the fact that $\FG$ has only identity 2-cells \pcref{def:Pig}, the naturality axiom \cref{thaobjx_nat} for $\cgs_\A$ means that the following diagram in $\A$ commutes for each $j \in \ufsn$ and each morphism $\objf \cn \objx \to \objy$ in $\DG(\mal,\nbeta)$ \cref{dgmn_morphism}.
\begin{equation}\label{cgsnat_diag}
\begin{tikzpicture}[vcenter]
\def\v{-1.4}
\draw[0cell=.9]
(0,0) node (a11) {\gaA\big(x_j; \ang{a_{\{i\}}}_{i \in \psiinv j} \big)}
(a11)++(4.5,0) node (a12) {a_{\psiinv j}}
(a11)++(0,\v) node (a21) {\gaA\big(y_j; \ang{a_{\{i\}}}_{i \in \psiinv j} \big)}
(a12)++(0,\v) node (a22) {a_{\psiinv j}}
;
\draw[1cell=.85]
(a11) edge node {\gl_{x_j;\, \psiinv j,\, \ang{\{i\}}_{i \in \psiinv j}}} (a12)
(a12) edge[equal] (a22)
(a11) edge[transform canvas={xshift=-2.5em}] node {\gaA(f_j; \ang{a_{\{i\}}}_{i \in \psiinv j})} (a21)
(a21) edge node {\gl_{y_j;\, \psiinv j,\, \ang{\{i\}}_{i \in \psiinv j}}} (a22)
;
\end{tikzpicture}
\end{equation}
The preceding diagram commutes by the naturality axiom \cref{nsys_naturality} for $(a,\gl)$.
\item[$G$-equivariance]
The axiom \cref{thaobjx_gequiv} for $\cgs_\A$ means that, for each $g \in G$ and each object $\objx = (\psi; \ang{x_j}_{j \in \ufsn})$ in $\DG(\mal,\nbeta)$ \cref{dgmn_obj}, there is an equality of natural transformations
\begin{equation}\label{cgsa_gequiv}
\cgs_{\A,g\objx} = g \cgs_{\A,\objx} \ginv.
\end{equation}
The following commutative diagram in $\proAnbeta$ proves the desired equality \cref{cgsa_gequiv}, where $\ang{\Cdots}_j = \ang{\Cdots}_{j \in \ufsn}$.
\begin{equation}\label{cgsgeq_diag}
\begin{tikzpicture}[vcenter]
\def\v{-1.2}
\draw[0cell=.7]
(0,0) node (a11) {\ang{\gaA(g x_{\ginv j} \tau_{\ginv} ; \ang{a_{\{i\}}}_{i \in (g\psi\ginv)^{-1}j} )}_j}
(a11)++(7,0) node (a12) {\ang{a_{(g\psi\ginv)^{-1} j}}_j}
(a11)++(0,\v) node (a21) {\ang{\gaA(g x_{\ginv j} ; \ang{a_{\{gi\}}}_{i \in \psiinv \ginv j})}_j}
(a12)++(0,\v) node (a22) {\ang{a_{g\psiinv\ginv j} }_j}
(a21)++(0,\v) node (a31) {\ang{g\gaA(x_{\ginv j} ; \ang{(\ginv a)_{\{i\}}}_{i \in \psiinv\ginv j})}_j}
(a22)++(0,\v) node (a32) {\ang{g(\ginv a)_{\psiinv\ginv j}}_j}
;
\draw[1cell=.7]
(a11) edge[equal] (a21)
(a21) edge[equal] (a31)
(a12) edge[equal] (a22)
(a22) edge[equal] (a32)
(a11) edge node {\ang{\gl_{g x_{\ginv j} \tau_{\ginv} ;\, (g\psi\ginv)^{-1}j ,\, \ang{\{i\}}_{i \in (g\psi\ginv)^{-1}j}}}_j} (a12)
(a21) edge node {\ang{\gl_{g x_{\ginv j} ;\, g\psiinv\ginv j ,\, \ang{\{gi\}}_{i \in \psiinv\ginv j}}}_j} (a22)
(a31) edge node {\ang{g(\ginv\gl)_{x_{\ginv j} ;\, \psiinv\ginv j ,\, \ang{\{i\}}_{i \in \psiinv\ginv j}}}_j} (a32)
;
\end{tikzpicture}
\end{equation}
\begin{itemize}
\item By \cref{dgmn_gact,cgsax_comp}, the top horizontal arrow is $\cgs_{\A,g\objx,(a,\gl)}$.  By \cref{nsys_gaction,proCnbe_gaction,cgsax_comp}, the bottom horizontal arrow is $g\cgs_{\A,\objx,\ginv(a,\gl)}$. 
\item The top rectangle commutes by the equivariance axiom \cref{nsys_equivariance} for $(a,\gl)$ applied to the permutation $\tau_{\ginv}$ \cref{tauginv}.
\item The bottom rectangle commutes by the equality $g\ginv = 1$, the $G$-equivariance of $\gaA$ \cref{gaAn}, and the definitions of $\ginv a$ \cref{ga_scomp} and $\ginv\gl$ \cref{ga_gl}.
\end{itemize}
\item[Basepoint]
The basepoint axiom \cref{tha_objzero} for $\cgs_\A$ follows from the unity axiom \cref{nsys_unity_empty} for $(a,\gl)$ and the definition of the basepoint $\objzero \in \DG(\mal,\nbeta)$ \cref{dgmn_basept}.
\item[Unity]
The unity axiom \cref{tha_onex} for $\cgs_\A$ follows from the unity axiom \cref{nsys_unity_one} for $(a,\gl)$ and the definition of the identity 1-cell $\objone_{\mal} \in \DG(\mal,\mal)$ \cref{dg_idonecell}.
\item[Compositionality]
The axiom \cref{tha_xv} for $\cgs_\A$ means that the two pasting diagrams of natural transformations
\begin{equation}\label{cgsa_compax}
\begin{tikzpicture}[vcenter]
\def\h{2.2} \def\v{-1.3} \def\s{.8} \def\q{.8} \def\p{.8} \def\r{240} \def\t{.3}
\def\boundary{
\draw[0cell=\s]
(0,0) node (a11) {\Asgldea}
(a11)++(\h,0) node (a12) {\phantom{\proAldea}}
(a12)++(0,.05) node (a12') {\proAldea}
(a11)++(0,2*\v) node (a31) {\Asgnbeta}
(a12)++(0,2*\v) node (a32) {\phantom{\proAnbeta}}
(a32)++(0,.05) node (a32') {\proAnbeta}
;
\draw[1cell=\q]
(a11) edge node {\cgs_{\A,\ldea}} (a12)
(a31) edge node[swap] {\cgs_{\A,\nbeta}} (a32)
;}
\boundary
\draw[0cell=\s]
(a11)++(0,\v) node (a21) {\Asgmal}
(a12)++(0,\v) node (a22) {\phantom{\proAmal}}
(a22)++(0,.05) node (a22') {\proAmal}
;
\draw[1cell=\q]
(a21) edge node {\cgs_{\A,\mal}} (a22)
(a11) edge node[swap] {\Asgphi} (a21)
(a21) edge node[swap] {\Asgpsi} (a31)
(a12') edge node {(\Rg\A)\objv} (a22')
(a22') edge node {(\Rg\A)\objx} (a32')
;
\draw[2cell=\p]
node[between=a12 and a21 at .5, shift={(\t,0)}, rotate=\r, 2labelmed={below,\cgs_{\A,\objv}}] {\Rightarrow}
node[between=a22 and a31 at .5, shift={(\t,-.15)}, rotate=\r, 2labelmed={below,\cgs_{\A,\objx}}] {\Rightarrow}
;
\begin{scope}[shift={(2.2*\h,0)}]
\boundary
\draw[1cell=\q]
(a11) edge node[swap] {\Asgpsiphi} (a31)
(a12') edge node {(\Rg\A)(\objx\objv)} (a32')
;
\draw[2cell=\p]
node[between=a12 and a31 at .5, shift={(\t,0)}, rotate=\r, 2labelmed={below,\cgs_{\A,\objx\objv}}] {\Rightarrow}
;
\end{scope}
\end{tikzpicture}
\end{equation}
are equal for objects 
\[\begin{split}
\objx &= (\psi; \ang{x_j}_{j \in \ufsn}) \in \DG(\mal,\nbeta) \andspace\\
\objv &= (\phi; \ang{v_i}_{i \in \ufsm}) \in \DG(\ldea, \mal) 
\end{split}\]
with composite \cref{dg_comp}
\[\objx\objv = \big(\psi\phi; \ang{\bx_j \tau^j_{\psi,\phi} }_{j \in \ufsn}\big) \inspace \DG(\ldea,\nbeta),\]
where
\[\bx_j = \ga\big(x_j; \ang{v_i}_{i \in \psiinv j}\big) \in \Op_{\sum_{i \in \psiinv j} |\phiinv i|}.\]
To prove this axiom, we consider the following diagram in $\A$ for $j \in \ufsn$ and a strong $\ldea$-system $(a,\gl) \in \Asgldea$.
\begin{equation}\label{cgscompos_diag}
\begin{tikzpicture}[vcenter]
\def\h{6.7} \def\s{.7} \def\v{-1.4} \def\u{-1}
\draw[0cell=\s]
(0,0) node (a11) {\gaA(x_j; \ang{\gaA(v_i; \ang{a_{\{k\}}}_{k \in \phiinv i})}_{i \in \psiinv j})}
(a11)++(\h,0) node (a12) {\gaA(x_j; \ang{a_{\phiinv i}}_{i \in \psiinv j})}
(a11)++(0,\v) node (a21) {\gaA(\bx_j ; \ang{a_{\{k\}}}_{k \in \phiinv i,\, i \in \psiinv j})}
(a12)++(0,\v) node (a22) {a_{\phiinv\psiinv j}}
(a21)++(0,\u) node (a31) {\gaA(\bx_j \tau^j_{\psi,\phi} ; \ang{a_{\{i\}}}_{i \in \psiphiinv j})}
(a22)++(0,\u) node (a32) {a_{\psiphiinv j}}
;
\draw[1cell=\s]
(a11) edge node {\gaA(x_j; \ang{\gl_{v_i;\, \phiinv i ,\, \ang{\{k\}}_{k \in \phiinv i}}}_{i \in \psiinv j})} (a12)
(a21) edge node {\gl_{\bx_j ;\, \psiphiinv j ,\, \ang{\{k\}}_{k \in \phiinv i,\, i \in \psiinv j}}} (a22)
(a31) edge node {\gl_{\bx_j \tau^j_{\psi,\phi} ;\, \psiphiinv j ,\, \ang{\{i\}}_{i \in \psiphiinv j}}} (a32)
(a12) edge node[swap,pos=.4] {\gl_{x_j ;\, \phiinv\psiinv j ,\, \ang{\phiinv i}_{i \in \psiinv j}}} (a22)
(a11) edge[equal] (a21)
(a21) edge[equal] (a31)
(a22) edge[equal] (a32)
;
\end{tikzpicture}
\end{equation}
\begin{itemize}
\item By \cref{glpsitil,rgaxa,cgsax_comp}, the top-right boundary composite is the $j$-th component morphism of the left pasting diagram in \cref{cgsa_compax} applied to $(a,\gl)$.  By \cref{dg_comp,cgsax_comp}, the bottom horizontal arrow is the $j$-th component morphism of $\cgs_{\A,\objx\objv,(a,\gl)}$.
\item The top rectangle commutes by the associativity axiom \cref{nsys_associativity} for $(a,\gl)$ and the fact that the associativity constraint $\phiA$ \cref{phiA} is the identity, since $\A$ is an $\Op$-algebra. 
\item The bottom rectangle commutes by the equivariance axiom \cref{nsys_equivariance} for $(a,\gl)$ applied to the permutation $\tau^j_{\psi,\phi}$ \cref{tauj_psiphi}.
\end{itemize}
Thus, the two pasting diagrams in \cref{cgsa_compax} are equal.
\end{description}
This proves that $\cgs_\A$ is a $\Pig$-strict $\Gcatst$-pseudotransformation.
\end{proof}

\begin{lemma}\label{cgs_iinat}
There is a 2-natural transformation\label{not:cgs}
\begin{equation}\label{cgs_nat_lemma}
\begin{tikzpicture}[vcenter]
\def\t{30} \def\h{2.7} \def\v{1.4}
\draw[0cell]
(0,0) node (a1) {\AlgstpsO}
(a1)++(0,\v) node (a2) {\AlgpspsO}
(a2)++(\h,0) node (a3) {\FGCatg}
(a3)++(.2,0) node (a3') {\phantom{\FGCatg}}
(a3')++(0,-\v) node (a4) {\DGCatg}
;
\draw[1cell=.85]
(a1) edge[right hook->] node {\algi} (a2)
(a2) edge node {\Sgosg} (a3)
(a3') edge[transform canvas={xshift=-1ex}] node {\gxist} (a4)
(a1) edge node[swap] {\Rg} (a4)
;
\draw[2cell=1]
node[between=a2 and a3 at .42, shift={(0,-.5*\v)}, rotate=-90, 2label={above,\cgs}] {\Rightarrow}
;
\end{tikzpicture}
\end{equation}
given by the assignment $\A \mapsto \cgs_\A$ in \cref{def:cgs}.
\end{lemma}

\begin{proof}
\cref{cgs_welldef} proves that each $\cgs_\A$ is a 1-cell in $\DGCatg$.  We check that $\cgs$ is natural with respect to 1-cells and 2-cells.
\begin{description}
\item[1-naturality]
The 1-cell naturality of $\cgs$ means that, for each $\Op$-pseudomorphism between $\Op$-algebras \pcref{def:laxmorphism}
\[(\A,\gaA) \fto{(f,\actf)} (\B,\gaB),\]
the following diagram of 1-cells in $\DGCatg$ \pcref{def:dgcatg_onecell,def:pig_strict} commutes.
\begin{equation}\label{cgs_onenat}
\begin{tikzpicture}[vcenter]
\def\v{-1.4}
\draw[0cell]
(0,0) node (a11) {\gxist\Sgosg\algi\A}
(a11)++(2.5,0) node (a12) {\Rg\A}
(a11)++(0,\v) node (a21) {\gxist\Sgosg\algi\B}
(a12)++(0,\v) node (a22) {\Rg\B}
;
\draw[1cell=.9]
(a11) edge node {\cgs_\A} (a12)
(a21) edge node {\cgs_\B} (a22)
(a11) edge[transform canvas={xshift=1em}] node[swap] {\gxist\Sgosg\algi f} (a21)
(a12) edge[transform canvas={xshift=-1ex}] node {\Rg f} (a22)
;
\end{tikzpicture}
\end{equation}
We verify that the two composite 1-cells in \cref{cgs_onenat} have the same component pointed $G$-functors \cref{dgcatg_one_comp} and natural isomorphisms \cref{dgcatg_one_compii}. 
\begin{description}
\item[$G$-functors]
For each object $\mal \in \DG$, we need to show that the diagram
\begin{equation}\label{cgsgfunctor_diag}
\begin{tikzpicture}[vcenter]
\def\v{-1.3}
\draw[0cell]
(0,0) node (a11) {\Asgmal}
(a11)++(2.5,0) node (a12) {\phantom{\proAmal}}
(a12)++(0,.05) node (a12') {\proAmal}
(a11)++(0,\v) node (a21) {\Bsgmal}
(a12)++(0,\v) node (a22) {\phantom{\proBmal}}
(a22)++(0,.05) node (a22') {\proBmal}
;
\draw[1cell=.85]
(a11) edge node {\cgs_{\A,\mal}} (a12)
(a21) edge node {\cgs_{\B,\mal}} (a22)
(a11) edge[transform canvas={xshift=1ex}] node[swap] {(\Sgosg f)_{\mal}} (a21)
(a12') edge[transform canvas={xshift=-.8ex}] node {(\Rg f)_{\mal}} (a22')
;
\end{tikzpicture}
\end{equation}
of pointed $G$-functors commutes.  By \cref{sgof_nbe_a,rgf_mal,cgsam_def}, each of the two composites in the preceding diagram has object assignment 
\begin{equation}\label{cgs_inat_obj}
\Asgmal \ni (a,\gl) \mapsto \ang{fa_{\{i\}}}_{i \in \ufsm} \in \proBmal
\end{equation}
and likewise for morphisms using \cref{sgof_tha_s}.
\item[Natural isomorphisms]
Recall that $\gxist$ \cref{gxist} actually lands in $\dgcatg$.  By \cref{dgcatg_icell_paste}, we need to show that, for each object $\objx = (\psi; \ang{x_j}_{j \in \ufsn})$ in $\DG(\mal,\nbeta)$ \cref{dgmn_obj}, the following two pasting diagrams of natural isomorphisms are equal.
\begin{equation}\label{cgsnatiso_diag}

\end{equation}
By \cref{sgof_nbe_gl,cgsax_comp}, the left pasting diagram applied to a strong $\mal$-system $(a,\gl) \in \Asgmal$ is the following composite isomorphism in $\proBnbeta$, where $\ang{\Cdots}_j = \ang{\Cdots}_{j \in \ufsn}$.
\begin{equation}\label{cgsnatiso_left}
%
\end{equation}
By  \cref{rgfxa,cgsam_def,cgsax_comp}, the right pasting diagram applied to $(a,\gl)$ is also the preceding composite isomorphism.  This proves that the diagram \cref{cgs_onenat} commutes.
\end{description}
\item[2-naturality]
The 2-cell naturality of $\cgs$ is verified in the same way as \cref{cgs_inat_obj} by replacing $f$ with an $\Op$-transformation between $\Op$-pseudomorphisms \pcref{def:algtwocells} and using \cref{sgo_omega_nbe_s,rgom_mal,cgsam_def}.
\end{description}
This proves that $\cgs \cn \gxist\Sgosg\algi \to \Rg$ is a 2-natural transformation.
\end{proof}

\section{Comparison Weak $G$-Equivalences}
\label{sec:gmmo_shi_weq}

This section compares three equivariant $K$-theory machines:
\begin{itemize}
\item the prolongation $\Kgg$ \pcref{def:Kgg_functor},
\item the homotopical Shimakawa strong $K$-theory $\Khshosg$ \cref{khsho_khshosg}, and
\item GMMO $K$-theory $\Kgmmo$ \cref{kgmmo_diag}.
\end{itemize}  
The main result, \cref{thm:gmmo_shi}, proves that $\Khshosg$ and $\Kgmmo$ are connected by a natural zigzag that is a weak $G$-equivalence in important cases.  For example, these machines yield weakly $G$-equivalent orthogonal $G$-spectra for each genuine permutative $G$-category of the form $\Ah = \Catg(\EG,\A)$ for some naive permutative $G$-category $\A$.  See \cref{ex:gmmo_shi_gbe}.

\recollection
For \cref{thm:gmmo_shi}, recall
\begin{itemize}
\item chaotic $\Einfg$-operads \pcref{def:chaotic_einf};
\item the equality of 2-functors \cref{gxist_gzest} 
\[\FGCatg \fto{\gzest\gxist = \Incj} \FGCatgps;\]
\item the counit $\strv \cn \str\Incj \to 1_{\FGCatg}$ \cref{strv} of the 2-adjunction $(\str,\Incj)$;
\item the construction $\Catg(\EG,-)$ \pcref{catgego};
\item the functors 
\[\FGCatg \fto{\clast} \FGTopg \fto{\Bc} \FGTopg \fto{\Kfg} \Gspec\]
in \cref{ggcatg_ggtopg}, \cref{bar_functor_FG}, and \cref{def:Kfg_functor}; 
\item the 2-natural transformation $\cgs \cn \gxist\Sgosg\algi \to \Rg$ \pcref{cgs_iinat};
\item the natural transformation $\kbgg \cn \Kfg\Bc \to \Kgg\ifgl$ \pcref{kfgbkgg_compare}; and 
\item that weak $G$-equivalences between orthogonal $G$-spectra are defined componentwise \pcref{def:kgmmo_weakg}. 
\end{itemize} 
The diagram \cref{gmmo_shi_diag} summarizes the homotopical Shimakawa strong $K$-theory $\Khshosg$ \cref{khsho_khshosg} and GMMO $K$-theory $\Kgmmo$ \cref{kgmmo_diag} along its top and left-bottom-right boundaries.

\begin{equation}\label{gmmo_shi_diag}

\end{equation}

\begin{theorem}\label{thm:gmmo_shi}\index{GMMO K-theory@GMMO $K$-theory!comparison with Shimakawa}\index{Shimakawa K-theory@Shimakawa $K$-theory!comparison with GMMO}
For \cref{thm:gmmo_shi_i,thm:gmmo_shi_ii}, suppose $\Op$ is a chaotic $\Einfg$-operad for a finite group $G$.
\begin{enumerate}
\item\label{thm:gmmo_shi_i}
There are natural transformations
\begin{equation}\label{gmmo_shi_zigzag}
%
\end{equation}
between functors $\AlgstpsO \to \Gspec$.
\item\label{thm:gmmo_shi_ii}
Each component of the left arrow $\Kfg\Bc\clast \strv_{\Sgosg\algi}$ in \cref{gmmo_shi_zigzag} is a weak $G$-equivalence of orthogonal $G$-spectra.
\enumeratext{For \cref{thm:gmmo_shi_iii,thm:gmmo_shi_iv}, suppose $\Op$ is a $\Uinf$-operad \pcref{as:OpA'} such that $\Oph = \Catg(\EG,\Op)$ is a chaotic $\Einfg$-operad, and $\Ah = \Catg(\EG,\A)$ is an $\Oph$-algebra for some $\Op$-algebra $\A$.}
\item\label{thm:gmmo_shi_iii}
Applying \cref{gmmo_shi_zigzag} to $(\Oph,\Ah)$ yields weak $G$-equivalences 
\begin{equation}\label{gmmo_shi_zig}
\begin{tikzpicture}[vcenter]
\def\u{-1.4} \def\v{-.9}
\draw[0cell=.9]
(0,0) node (a11) {\Kfg\Bc\clast\str\Incj\Sgohsg\algi\Ah}
(a11)++(0,\u) node (a21) {\Kfg\Bc\clast\Sgohsg\algi\Ah}
(a21)++(0,\v) node (a31) {\Khshohsg\algi\Ah}
(a11)++(3.5,0) node (a12) {\Kfg\Bc\clast\str\gzest\gxist\Sgohsg\algi\Ah}
(a12)++(0,\u) node (a22) {\Kfg\Bc\clast\str\gzest\Rg\Ah}
(a22)++(0,\v) node (a32) {\phantom{\Kgmmoh\Ah}}
(a32)++(0,-.01) node (a32') {\Kgmmoh\Ah}
;
\draw[1cell=.85]
(a11) edge node[swap] {\Kfg\Bc\clast \strv_{\Sgohsg\algi\Ah}} (a21)
(a12) edge node[swap] {\Kfg\Bc\clast\str\gzest\cgs_{\Ah}} (a22)
(a11) edge[equal] (a12)
(a21) edge[equal] (a31)
(a22) edge[equal] (a32')
;
\end{tikzpicture}
\end{equation}
between orthogonal $G$-spectra.
\item\label{thm:gmmo_shi_iv}
Combining \cref{kbggX_weakg,gmmo_shi_zig} yields weak $G$-equivalences 
\begin{equation}\label{gmmo_kgg_zigzag}
\begin{tikzpicture}[vcenter]
\def\u{-1.4} \def\v{-.9}
\draw[0cell=.9]
(0,0) node (a11) {\Kfg\Bc\clast\str\Incj\Sgohsg\algi\Ah}
(a11)++(0,\u) node (a21) {\Kfg\Bc\clast\Sgohsg\algi\Ah}
(a21)++(0,\v) node (a31) {\Khshohsg\algi\Ah}
(a31)++(-3.2,0) node (a30) {\Kgg\ifgl\clast\Sgohsg\algi\Ah}
(a11)++(3.3,0) node (a12) {\Kfg\Bc\clast\str\gzest\gxist\Sgohsg\algi\Ah}
(a12)++(0,\u) node (a22) {\Kfg\Bc\clast\str\gzest\Rg\Ah}
(a22)++(0,\v) node (a32) {\phantom{\Kgmmoh\Ah}}
(a32)++(0,-.01) node (a32') {\Kgmmoh\Ah}
;
\draw[1cell=.85]
(a31) edge node[swap] {\kbgg_{\clast\Sgohsg\algi\Ah}} (a30)
(a11) edge node[swap] {\Kfg\Bc\clast \strv_{\Sgohsg\algi\Ah}} (a21)
(a12) edge node[swap] {\Kfg\Bc\clast\str\gzest\cgs_{\Ah}} (a22)
(a11) edge[equal] (a12)
(a21) edge[equal] (a31)
(a22) edge[equal] (a32')
;
\end{tikzpicture}
\end{equation}
between orthogonal $G$-spectra.
\end{enumerate}
\end{theorem}

\begin{proof}
\begin{description}
\item[\cref{thm:gmmo_shi_i}]
The two natural transformations in the diagram \cref{gmmo_shi_zigzag} exist by the 2-naturality of
\begin{itemize}
\item the counit $\strv \cn \str\Incj \to 1$ \cref{strv} and
\item $\cgs \cn \gxist\Sgosg\algi \to \Rg$ \pcref{cgs_iinat}.  
\end{itemize}
The top horizontal equality in \cref{gmmo_shi_zigzag} follows from the equality $\gzest\gxist = \Incj$ \cref{gxist_gzest}.  
\item[\cref{thm:gmmo_shi_ii}]
This assertion follows from the following two facts.
\begin{itemize}
\item The composite functor $\Kfg\Bc\clast$ preserves weak $G$-equivalences \pcref{kgmmo_weq}.
\item Each component of the counit $\strv$ \cref{strvx} is a weak $G$-equivalence in $\FGCatg$.  In fact, for each $\FGG$-category $X$, each $\strv_{X,\mal}$ \cref{strvx_mal} is a pointed $G$-equivalence between pointed $G$-categories that is inverse to $\stru_{X,\mal}$ \cref{strux_mal}.  Thus, each $\strv_{X,\mal}$ is a categorical weak $G$-equivalence \pcref{def:cat_weakg}.
\end{itemize}
\item[\cref{thm:gmmo_shi_iii}]
Recall that there is an equality of pointed $G$-functors 
\[\Asgmal \fto{\cgs_{\A,\mal} \,=\,\zdsg} \proAmal\]
between $\cgs_{\A,\mal}$ \cref{cgsam_def} and $\zdsg$ \cref{hgosgtoprod}.  Regarding a pointed finite $G$-set $\mal$ as a length-1 object in $\GG$, the diagram \cref{zdsg_catweakg} yields a categorical weak $G$-equivalence
\[\Ahsgmal \fto{\cgs_{\Ah,\mal} \,=\, \zdsg} \proAhmal.\]
Thus, the 1-cell
\begin{equation}\label{cgs_ah_wgeq}
\gxist\Sgohsg\algi\Ah \fto{\cgs_{\Ah}} \Rg\Ah
\end{equation}
in $\DGCatg$ is a weak $G$-equivalence.  Since the functor $\Kfg\Bc\clast\str\gzest$ preserves weak $G$-equivalences \pcref{kgmmo_weq}, the $G$-morphism
\begin{equation}\label{gmmo_shi_zig_i}
\Kfg\Bc\clast\str\gzest\gxist\Sgohsg\algi\Ah \fto{\Kfg\Bc\clast\str\gzest\cgs_{\Ah}} 
\Kfg\Bc\clast\str\gzest\Rg\Ah = \Kgmmoh\Ah
\end{equation}
is a weak $G$-equivalence in $\Gspec$.  Moreover, by assertion \cref{thm:gmmo_shi_ii} for $(\Oph,\Ah)$, the $G$-morphism
\begin{equation}\label{gmmo_shi_zig_ii}
\Kfg\Bc\clast\str\Incj\Sgohsg\algi\Ah \fto{\Kfg\Bc\clast \strv_{\Sgohsg\algi\Ah}} 
\Kfg\Bc\clast\Sgohsg\algi\Ah = \Khshohsg\algi\Ah
\end{equation}
is a weak $G$-equivalence in $\Gspec$.  The desired zigzag \cref{gmmo_shi_zig} is given by combining \cref{gmmo_shi_zig_i,gmmo_shi_zig_ii}. 
\item[\cref{thm:gmmo_shi_iv}]
The zigzag \cref{gmmo_kgg_zigzag} is given by combining \cref{gmmo_shi_zig} and the weak $G$-equivalence (\cref{kfgbkgg_compare} \cref{kfgbkgg_compare_iii})
\[\Kfg\Bc X \fto{\kbgg_X} \Kgg\ifgl X\]
for the $\FGG$-space 
\[X = \clast\Sgohsg\algi\Ah\]
in the image of $\clast = \Reast \comp \Nerst$ \cref{clast_factors}.\qedhere
\end{description}
\end{proof}

\begin{example}[Permutative $G$-Categories]\label{ex:gmmo_shi_gbe}
The Barratt-Eccles $\Gcat$-operad $\BE$ \pcref{def:BE} is a $\Uinf$-operad, which means that it is 1-connected and levelwise a nonempty translation category \pcref{as:OpA'}.  Its algebras in $\Gcat$ are \emph{naive} permutative $G$-categories.  Applying $\Catg(\EG,-)$ to $\BE$ yields the $G$-Barratt-Eccles operad $\GBE$ \pcref{def:GBE}, which is a chaotic $\Einfg$-operad \pcref{ex:gbe_chaotic}.  Its algebras are \emph{genuine} permutative $G$-categories.  Each naive permutative $G$-category $\A$ yields a genuine permutative $G$-category $\Ah = \Catg(\EG,\A)$ \pcref{naive_genuine_pGcat}.  By \cref{gmmo_kgg_zigzag}, for each genuine permutative $G$-category $\Ah$, there are weak $G$-equivalences of orthogonal $G$-spectra
\begin{equation}\label{KweakG_ex}
\begin{tikzpicture}[vcenter]
\def\u{-1.4} \def\v{-.9}
\draw[0cell=.9]
(0,0) node (a11) {\Kfg\Bc\clast\str\Incj\Sgosggbe\algi\Ah}
(a11)++(0,\u) node (a21) {\Kfg\Bc\clast\Sgosggbe\algi\Ah}
(a21)++(0,\v) node (a31) {\Khshgbesg\algi\Ah}
(a31)++(-3.2,0) node (a30) {\Kgg\ifgl\clast\Sgosggbe\algi\Ah}
(a11)++(3.3,0) node (a12) {\Kfg\Bc\clast\str\gzest\gxist\Sgosggbe\algi\Ah}
(a12)++(0,\u) node (a22) {\Kfg\Bc\clast\str\gzest\Rg\Ah}
(a22)++(0,\v) node (a32) {\phantom{\Kgmmogbe\Ah}}
(a32)++(0,0) node (a32') {\Kgmmogbe\Ah}
;
\draw[1cell=.85]
(a31) edge node[swap] {\kbgg_{\clast\Sgosggbe\algi\Ah}} (a30)
(a11) edge node[swap] {\Kfg\Bc\clast \strv_{\Sgosggbe\algi\Ah}} (a21)
(a12) edge node[swap] {\Kfg\Bc\clast\str\gzest\cgs_{\Ah}} (a22)
(a11) edge[equal] (a12)
(a21) edge[equal] (a31)
(a22) edge[equal] (a32')
;
\end{tikzpicture}
\end{equation}
connecting the prolongation $\Kgg$ \pcref{def:Kgg_functor}, the homotopical Shimakawa strong $K$-theory $\Khshgbesg$ \cref{khsho_khshosg}, and GMMO $K$-theory $\Kgmmogbe$ \cref{kgmmo_diag}.
\end{example}

%% file: chap/kgl.tex
This chapter reviews Schwede global algebraic $K$-theory $\Ksc$ \cite{schwede_global}, which is the following composite functor. 
\begin{equation}\label{ksc_intro}
\begin{tikzpicture}[vcenter]
\def\h{2.5} \def\v{1.4}
\draw[0cell]
(0,0) node (a1) {\Parcat}
(a1)++(0,-\v) node (a2) {\FMCat}
(a2)++(\h,0) node (a3) {\FICat}
(a3)++(\h,0) node (a4) {\phantom{\FITop}}
(a4)++(0,-.035) node (a4') {\FITop}
(a4)++(0,\v) node (a5) {\Sptop}
;
\draw[1cell=.9]
(a1) edge node {\Ksc} (a5)
(a1) edge node[swap] {\Jsc} (a2) 
(a2) edge node {\uprst} (a3)
(a3) edge node {\clast} (a4)
(a4') edge node[swap] {\Kfi} (a5)
;
\end{tikzpicture}
\end{equation}
Global $K$-theory keeps track of $G$-actions for all finite groups $G$ simultaneously.  The input data for global $K$-theory are parsummable categories.  They are commutative monoids in the symmetric monoidal category of tame $\schm$-categories, and $\schm$ is the translation category of the monoid of injections 
\[\{0,1,2,\ldots\} = \ome \to \ome.\]  
Global $K$-theory produces symmetric spectra based on topological spaces.  Each such symmetric spectrum yields an underlying $G$-spectrum with the trivial $G$-action for each finite group $G$.  \cref{ch:kgl_gmmo} proves that the underlying $G$-spectra produced by global $K$-theory are $G$-stably equivalent to those produced by GMMO $K$-theory.  In this chapter, $G$ denotes a finite group.

\organization
This chapter consists of the following sections.

\secname{sec:mcat}
This section discusses universal $G$-sets, the injection category $\schm$, and $\schm$-categories.  Each $\schm$-category can be reparametrized to take into account countably infinite sets other than $\ome = \{0,1,2,\ldots\}$.  Reparametrization with respect to universal $G$-sets leads to the notion of a global equivalence between $\schm$-categories.

\secname{sec:jgl}
This section defines the box product for $\schm$-categories, parsummable categories, and the functor 
\[\Parcat \fto{\Jsc} \FMCat\]
that sends parsummable categories to $\FM$-categories.

\secname{sec:kfm}
This section discusses the passage from $\FM$-categories to symmetric spectra using the functors
\[\FMCat \fto{\uprst} \FICat \fto{\clast} \FITop \fto{\Kfi} \Sptop.\]
In the intermediate categories, $\bdI$ is the category of finite sets and injections.

\secname{sec:ksc}
This section discusses Schwede global $K$-theory 
\[\Parcat \fto{\Ksc = \Kfi \clast \uprst \Jsc} \Sptop .\]  
For a parsummable category $\C$ and a finite group $G$, the underlying $G$-spectrum of the global $K$-theory $\Ksc\C$ of $\C$ is computed by more accessible $G$-symmetric spectra $\kcg$ \pcref{def:sch4.15} and $\Kscgb\C$ \pcref{def:kscgb}.  The $G$-symmetric spectrum $\Kscgb\C$ is used in \cref{ch:kgl_gmmo} to compare global $K$-theory with GMMO $K$-theory.

\section{$\schm$-Categories}
\label{sec:mcat}

This section reviews universal $G$-sets, $\schm$-categories, their reparametrizations, and global equivalences between $\schm$-categories.  Most of this section is adapted from \cite[Section 2]{schwede_global}.

\secoutline
\begin{itemize}
\item \cref{def:univ_gset} defines universal $G$-sets.  \cref{schwede2.17} proves some useful properties of universal $G$-sets.
\item \cref{def:mcat} defines the injection category $\schm$, $\schm$-categories, and supports for objects in an $\schm$-category. \cref{sch2.13} recalls some basic properties of $\schm$-categories and supports.
\item \cref{def:repar} defines reparametrizations of $\schm$-categories.  \cref{def:sch2.21} defines the $G$-fixed $\schm$-category of an $\schm$-category.  It is used in \cref{def:sch2.26} to define global equivalences between $\schm$-categories.
\end{itemize}

\subsection*{Universal $G$-Sets}

\begin{definition}\label{def:univ_gset}
For a finite group $G$, a \emph{universal $G$-set}\index{universal G-set@universal $G$-set}\index{G-set@$G$-set!universal} is a countable $G$-set such that each subgroup of $G$ is the stabilizer group of infinitely many elements.
\end{definition}

\begin{example}\label{ex:univ_gset}
For a finite group $G$ and a subgroup $H \subseteq G$, the \emph{orbit $G$-set}\index{orbit G-set@orbit $G$-set}\index{G-set@$G$-set!orbit} 
\[G/H = \{gH \tmid g \in G\}\] 
has $G$ acting by left multiplication, meaning $g' \cdot gH = g'gH$.  Equipping \label{not:ome}$\ome = \{0,1,2,\ldots\}$ with the trivial $G$-action and letting $H \subseteq G$ run through all the subgroups of $G$, the countable $G$-set
\begin{equation}\label{ugsetw}
\ugsetw = \coprod_{H \subseteq G} \, \ome \times G/H
\end{equation}
is a universal $G$-set.  Indeed, for each $n \in \ome$, the element $(n,eH) \in \ugsetw$ has stabilizer group $H$, where $e \in G$ is the group unit.
\end{example}

\cref{schwede2.17} is stated in \cite[2.17]{schwede_global}.  We include a proof here for future reference.

\begin{lemma}\label{schwede2.17}
The following statements hold for a finite group $G$.
\begin{enumerate}
\item\label{sch2.17_i}
A countable $G$-set $\ugsetu$ is a universal $G$-set if and only if each finite $G$-set admits a $G$-injection into $\ugsetu$.
\item\label{sch2.17_ii}
Any two universal $G$-sets are $G$-isomorphic.
\item\label{sch2.17_iii}
Suppose $\ugsetu$ is a universal $G$-set, and $\ugsetv$ is a countable $G$-set that contains $\ugsetu$ as a $G$-subset.  Then $\ugsetv$ is a universal $G$-set.
\item\label{sch2.17_iv}
Suppose $\ugsetu$ is a universal $G$-set.  For each subgroup $H \subseteq G$, the underlying $H$-set of $\ugsetu$ is a universal $H$-set.
\end{enumerate}
\end{lemma}

\begin{proof}
Recall that $\ufsn = \{1,2,\ldots,n\}$ denotes an unpointed finite set with $n$ elements \cref{ufsn}.
\begin{description}
\item[\cref{sch2.17_i}]
To prove the \emph{if} assertion, we consider a subgroup $H \subseteq G$, an integer $n \geq 1$, and the finite $G$-set $\ufsn \times G/H$ with $G$ acting trivially on $\ufsn$.  By assumption, there is a $G$-injection 
\[\ufsn \times G/H \to \ugsetu.\]  
For each $i \in \ufsn$, $H$ is the stabilizer group of $(i,eH)$.  Since $n \geq 1$ is arbitrary, $H$ is the stabilizer group of infinitely many elements in $\ugsetu$, proving that $\ugsetu$ is a universal $G$-set.

To prove the \emph{only if} assertion, suppose $\ugsetu$ is a universal $G$-set, and $A$ is a finite $G$-set.  We need to construct a $G$-injection $f \cn A \to \ugsetu$.  Suppose $a \in A$ is an element with stabilizer group $H_a \subseteq G$.  We choose distinct left $H_a$-coset representatives 
\[\{e = g^a_1,g^a_2,\ldots,g^a_r\}\] 
such that
\[G/H_a = \{g^a_i H_a\}_{i \in \ufsr}.\]
Then the $G$-orbit of $a$ is given by the $r$ distinct elements
\[Ga = \{g^a_i a\}_{i \in \ufsr}.\]
Choosing an element $u_a \in \ugsetu$ with stabilizer group $H_a$, we define $f$ on $Ga$ by
\begin{equation}\label{fgia}
f(g^a_i a) = g^a_i u_a \forspace i \in \ufsr.
\end{equation}
The $G$-orbit of $u_a$ is given by the $r$ distinct elements 
\[Gu_a = \{g^a_i u_a\}_{i \in \ufsr}.\]  
The assignment \cref{fgia} defines a $G$-bijection 
\[Ga \fiso Gu_a.\]  
More explicitly, for $g \in G$ and $i \in \ufsr$, suppose 
\[gg^a_i a = g^a_j a.\]  
Then 
\[(g^a_j)^{-1} gg^a_i \in H_a,\] 
the stabilizer group of $u_a$, and
\[f(gg^a_i a)= f(g^a_j a)= g^a_j u_a = gg^a_i u_a = gf(g^a_i a).\]

To define $f$ on the rest of $A$, we consider an element $b \in A \setminus Ga$ with stabilizer group $H_b \subseteq G$.  Reusing the previous paragraph, we choose distinct left $H_b$-coset representatives $\{g^b_k\}_{k \in \ufst}$ and an element $u_b \in \ugsetu \setminus Gu_a$ with stabilizer group $H_b$.  Such an element $u_b$ exists by the assumption that $\ugsetu$ is a universal $G$-set.  The assignment
\[f(g^b_k b) = g^b_k u_b \forspace k \in \ufst\]
defines a $G$-bijection 
\[Gb \fiso Gu_b.\]  
Moreover, the $G$-orbits 
\[Gu_a \andspace Gu_b = \{g^b_k u_b\}_{k \in \ufst}\] 
are disjoint in $\ugsetu$.  Continuing inductively over the finite set of $G$-orbits in $A$, we obtain a $G$-injection $f \cn A \to \ugsetu$.  This finishes the proof of the \emph{only if} assertion.
\item[\cref{sch2.17_ii}]
It suffices to show that each universal $G$-set $\ugsetu$ is $G$-isomorphic to the universal $G$-set $\ugsetw$ \cref{ugsetw}.  For each subgroup $H \subseteq G$, we choose distinct left $H$-coset representatives $\{g_i^H\}_{i \in \ufsr_H}$ such that 
\[G/H = \{g_i^H H\}_{i \in \ufsr_H}.\]
We choose an element $u_j \in \ugsetu$ in each $G$-orbit in $\ugsetu$, decomposing $\ugsetu$ into a countable disjoint union of $G$-orbits
\[\ugsetu = \coprod_{j \geq 0}\, Gu_j,\]
such that each subgroup $H \subseteq G$ is the stabilizer group of infinitely many $u_j$'s.  This is possible by the finiteness of $G$ and the universality of the $G$-set $\ugsetu$.  We define a $G$-isomorphism $f \cn \ugsetu \fiso \ugsetw$ by induction on $j \geq 0$.  

If $u_0$ has stabilizer group $H \subseteq G$, then, analogous to \cref{fgia}, we define the $G$-injection 
\[\{g_i^H u_0\}_{i \in \ufsr_H} = Gu_0 \fto[\iso]{f} (0, G/H) \subseteq \ugsetw\]
by
\[f(g_i^H u_0) = (0, g_i^H H) \forspace i \in \ufsr_H.\]  
Inductively, suppose $u_j$ has stabilizer group $K \subseteq G$.  We define the $G$-injection
\[\{g_i^K u_j\}_{i \in \ufsr_K} = Gu_j \fto[\iso]{f} (\ell, G/K) \subseteq \ugsetw\]
by 
\[f(g_i^K u_j) = (\ell, g_i^K K) \forspace i \in \ufsr_K,\] 
where $\ell$ is the smallest integer that has not been used in the previous $j-1$ steps.  For example, if the stabilizer group $K$ of $u_1$ is not $H$, the stabilizer group of $u_0$, then we define the $G$-injection
\[Gu_1 \fto[\iso]{f} (0, G/K) \subseteq \ugsetw\]
that sends $g_i^K u_1$ to $(0, g_i^K K)$.  If $u_1$ also has stabilizer group $H$, then we define the $G$-injection
\[Gu_1 \fto[\iso]{f} (1, G/H) \subseteq \ugsetw\]
that sends $g_i^H u_1$ to $(1, g_i^H H)$.  This defines a $G$-injection $f \cn \ugsetu \to \ugsetw$ that is also surjective by the choices of the $u_j$'s.  Thus, $f$ is a $G$-isomorphism.
\item[\cref{sch2.17_iii}]
By the \emph{only if} part of \cref{sch2.17_i}, each finite $G$-set admits a $G$-injection into $\ugsetu$, hence also into $\ugsetv$.  Thus, $\ugsetv$ is a universal $G$-set by the \emph{if} part of \cref{sch2.17_i}.
\item[\cref{sch2.17_iv}]
Each subgroup $K \subseteq H$ is also a subgroup of $G$.  Thus, $K$ is the stabilizer group of infinitely many elements in $\ugsetu$.\qedhere
\end{description}
\end{proof}

\begin{definition}\label{def:omea}
For the countable set $\ome = \{0,1,2,\ldots\}$, a finite group $G$, and a finite $G$-set $A$, define the countable $G$-set $\omea$ whose elements are functions $A \to \ome$.  For $g \in G$, the $g$-action on a function $f \cn A \to \ome$ is defined as
\begin{equation}\label{g_omea}
(gf)(a) = f(\ginv a) \forspace a \in A.
\end{equation}
Note that $gf$ is the conjugation $g$-action on $f$ \cref{ginv_h_g}, regarding $\ome$ as a trivial $G$-set.
\end{definition}

\begin{example}\label{ex:omea}
Suppose $G$ is a finite group, and $A$ is a finite $G$-set with a free $G$-orbit.  Then $\omea$ is a universal $G$-set.  In particular, equipping $G$ with the regular $G$-action, $\omeg$ is a universal $G$-set.  These examples are verified in \cite[2.19]{schwede_global} using \cref{schwede2.17} \cref{sch2.17_iii}.
\end{example}

\subsection*{Modules over the Injection Category}

Monoids and modules are reviewed in \cref{def:monoid,def:modules}.

\begin{definition}\label{def:mcat}
Denote by $\ome = \{0,1,2,\ldots\}$ the countable set of nonnegative integers.
\begin{enumerate}
\item\label{def:mcat_i}
The \emph{injection monoid}\index{injection monoid} $\Injm$ is the monoid in $(\Set,\times,*)$ of injections $\ome \to \ome$ under composition.  
\item\label{def:mcat_ii}
The \emph{injection category}\index{injection category} $\schm$ is the translation category $\tn\Injm$ \pcref{def:translation_cat}.  Thus, the objects of $\schm$ are injections $\ome \to \ome$, and each hom set in $\schm$ is a one-element set.  For objects $u,v \in \schm$, the unique isomorphism from $u$ to $v$ is denoted by \label{not:mor_vu}$[v,u] \cn u \to v$.  The monoid structure on $\Injm$ extends to a strict monoidal category structure on $\schm$.  The latter is also regarded as a monoid $(\schm,\comp,1_\ome)$ in the Cartesian monoidal category $(\Cat,\times, \bone)$ of small categories and functors. 
\item\label{def:mcat_iii}
An \emph{$\schm$-category}\index{M-category@$\schm$-category} is a left $\schm$-module in $\Cat$.  An \emph{$\schm$-functor}\index{M-functor@$\schm$-functor} is a morphism of left $\schm$-modules.  The category of $\schm$-categories and $\schm$-functors is denoted by $\Mcat$.  For an object or an isomorphism $u \in \schm$, the $u$-action functor or natural isomorphism on an $\schm$-category is denoted by \label{not:ustar}$u_\mstar$.  An object $x$ in an $\schm$-category $\C$ is \emph{$\schm$-fixed}\index{M-fixed@$\schm$-fixed} if $u_{\mstar}(x) = x$ for each injection $u \cn \ome \to \ome$.  Denote by $\C^{\suppempty}$\label{not:Csuppempty} the full subcategory of $\C$ consisting of $\schm$-fixed objects.
\item\label{def:mcat_iv}
An object $x$ in an $\schm$-category is \emph{supported} on a subset $A \subseteq \ome$ if $u_\mstar(x) = x$ for each injection $u \cn \ome \to \ome$ whose restriction to $A$ is the identity.  An object $x$ is \emph{finitely supported}\index{finitely supported} if it is supported on some finite subset of $\ome$.    The \emph{support}\index{support} of $x$, denoted by \label{not:suppx}$\supp(x)$, is defined as the intersection of all the finite subsets of $\ome$ on which $x$ is supported.  By convention, $\supp(x) = \ome$ if $x$ is not finitely supported.
\item\label{def:mcat_v}
An $\schm$-category is \emph{tame}\index{tame M-category@tame $\schm$-category}\index{M-category@$\schm$-category!tame} if each of its objects is finitely supported.  The category of tame $\schm$-categories and $\schm$-functors is denoted by \label{not:Mcatt}$\Mcatt$.\defmark
\end{enumerate}
\end{definition}

\begin{explanation}[Unpacking]\label{expl:mcat}
An $\schm$-category consists of
\begin{itemize}
\item a small category $\C$ and
\item an $\schm$-action functor $\mmu \cn \schm \ttimes \C \to \C$
\end{itemize} 
such that the following associativity and unity diagrams commute, where iterated products are interpreted using \cref{conv:left_normal}.
\begin{equation}\label{mcat_axioms}
\begin{tikzpicture}[vcenter]
\def\h{3} \def\d{2} \def\v{-1.4}
\draw[0cell]
(0,0) node (a11) {\schm \ttimes \schm \ttimes \C}
(a11)++(\h,0) node (a12) {\schm \ttimes \C}
(a11)++(0,\v) node (a21) {\schm \ttimes \C}
(a12)++(0,\v) node (a22) {\C}
(a12)++(2.5,0) node (b1) {\bone \ttimes \C}
(b1)++(0,.1) node (b1') {\phantom{\bone \ttimes \C}}
(b1)++(0,\v) node (b2) {\schm \ttimes \C}
(b2)++(\d,0) node (b3) {\C}
(b3)++(.15,0) node (b3') {\phantom{\C}}
;
\draw[1cell=.9]
(a11) edge node {1_{\schm} \ttimes \mmu} (a12)
(a12) edge node {\mmu} (a22)
(a11) edge[transform canvas={xshift=.9em}] node[swap] {\comp \ttimes 1_{\C}} (a21)
(a21) edge node {\mmu} (a22)
(b1) edge[transform canvas={xshift=1ex}] node[swap] {1_\ome \ttimes 1_{\C}} (b2)
(b2) edge node[pos=.45] {\mmu} (b3)
(b1') edge[bend left=15] node {\iso} (b3')
;
\end{tikzpicture}
\end{equation}
For each object $u \in \schm$, the $u$-action
\begin{equation}\label{u_action}
u_{\mstar} = \mmu(u,-) \cn \C \to \C
\end{equation}
is a functor.  An object $x \in \C$ is supported on the empty subset $\emptyset \subseteq \ome$ if and only if $u_{\mstar}(x) = x$ for each injection $u \cn \ome \to \ome$, meaning that $x$ is an $\schm$-fixed object.  For each isomorphism $[v,u] \cn u \to v$ in $\schm$, the natural isomorphism 
\begin{equation}\label{vu_star}
[v,u]_{\mstar} = \mmu([v,u],-) \cn u_{\mstar} \to v_{\mstar}
\end{equation}
between functors $\C \to \C$ sends each object $x \in \C$ to an isomorphism
\begin{equation}\label{vustar_x}
u_{\mstar}(x) \fto{[v,u]_{\mstar}^x} v_{\mstar}(x) \inspace \C.
\end{equation}
An $\schm$-functor $\fun \cn (\C,\mmu) \to (\D,\mmu')$ is a functor $\fun \cn \C \to \D$ such that the diagram
\begin{equation}\label{mfunctor_axiom}
\begin{tikzpicture}[vcenter]
\def\h{2.2} \def\v{-1.4}
\draw[0cell]
(0,0) node (a11) {\schm \ttimes \C}
(a11)++(\h,0) node (a12) {\C}
(a11)++(0,\v) node (a21) {\schm \ttimes \D}
(a12)++(0,\v) node (a22) {\D}
;
\draw[1cell=.9]
(a11) edge node {\mmu} (a12)
(a12) edge node {\fun} (a22)
(a11) edge[transform canvas={xshift=1em}] node[swap] {1_{\schm} \ttimes \fun} (a21)
(a21) edge node {\mmu'} (a22)
;
\end{tikzpicture}
\end{equation}
commutes.  
\end{explanation}

\begin{example}[Finite Sets]\label{ex:sch2.14}
There is an $\schm$-category $\Fome$ \cite[2.14]{schwede_global} with
\begin{itemize}
\item finite subsets of $\ome = \{0,1,2,\ldots\}$ as objects and
\item bijections as morphisms.
\end{itemize}  
For an injection $u \cn \ome \to \ome$ and a finite subset $A \subseteq \ome$, the $u$-action on $A$ is given by the image of $A$ under $u$:
\[u_{\mstar}(A) = u(A) \subseteq \ome.\]
For a bijection $p \cn A \fiso B$ between finite subsets of $\ome$, the $u$-action on $p$ is the bijection
\begin{equation}\label{ustarp}
\begin{tikzpicture}[vcenter]
\def\h{2} \def\u{.7}
\draw[0cell]
(0,0) node (a1) {u_{\mstar}(A)}
(a1)++(\h,0) node (a2) {\phantom{A}}
(a2)++(0,.03) node (a2') {A}
(a2)++(.8*\h,0) node (a3) {\phantom{B}}
(a3)++(0,.03) node (a3') {B}
(a3)++(\h,0) node (a4) {u_{\mstar}(B).}
;
\draw[1cell=.9]
(a1) edge node {(u|_A)^{-1}} (a2)
(a2) edge node {p} (a3)
(a3) edge node {u|_B} (a4)
(a1) [rounded corners=2pt] |- ($(a2)+(0,\u)$) -- node {u_{\mstar}(p)} ($(a3)+(0,\u)$) -| (a4)
;
\end{tikzpicture}
\end{equation}
For an isomorphism $[v,u] \cn u \to v$ in $\schm$, the $A$-component of the natural isomorphism $[v,u]_{\mstar}$ is the bijection
\begin{equation}\label{vustarA}
\begin{tikzpicture}[vcenter]
\def\h{2.2} \def\u{.7}
\draw[0cell]
(0,0) node (a1) {u_{\mstar}(A)}
(a1)++(\h,0) node (a2) {\phantom{A}}
(a2)++(0,.03) node (a2') {A}
(a2)++(\h,0) node (a3) {v_{\mstar}(A).}
;
\draw[1cell=.9]
(a1) edge node {(u|_A)^{-1}} (a2)
(a2) edge node {v|_A} (a3)
(a1) [rounded corners=2pt] |- ($(a2)+(-1,\u)$) -- node {[v,u]_{\mstar}^A} ($(a2)+(1,\u)$) -| (a3)
;
\end{tikzpicture}
\end{equation}
The support of $A$ is $A$ itself, so $\Fome$ is a tame $\schm$-category.
\end{example}

\begin{example}[Limits and Colimits]\label{ex:sch2.10}
The category $\Mcat$ has all small limits and colimits, and they are created by the forgetful functor $\Mcat \to \Cat$.  In particular, for $\schm$-categories $\C$ and $\D$, the Cartesian product $\C \ttimes \D$ is an $\schm$-category with the diagonal $\schm$-action.  It is the product in the category $\Mcat$.  
\end{example}

\cref{sch2.13} is \cite[2.13]{schwede_global}.

\begin{lemma}\label{sch2.13}
Suppose $x$ is an object in an $\schm$-category $\C$, and $t,u,v \cn \ome \to \ome$ are injections.
\begin{enumerate}
\item\label{sch2.13_i}
The object $x$ is supported on $\supp(x)$.
\item\label{sch2.13_ii}
If $u$ and $v$ agree on $\supp(x)$, then the following equalities hold in $\C$.
\[\begin{split}
u_{\mstar}(x) &= v_{\mstar}(x)\\
[t,u]_{\mstar}(x) &= [t,v]_{\mstar}(x) \cn u_{\mstar}(x) = v_{\mstar}(x) \to t_{\mstar}(x)\\
[u,t]_{\mstar}(x) &= [v,t]_{\mstar}(x) \cn t_{\mstar}(x) \to u_{\mstar}(x) = v_{\mstar}(x)
\end{split}\]
\item\label{sch2.13_iii}
If $x$ is supported on a subset $A \subseteq \ome$, then $u_{\mstar}(x)$ is supported on $u(A) \subseteq \ome$.  If $x$ is finitely supported, then so is $u_{\mstar}(x)$, and
\[\supp\big(u_{\mstar}(x)\big) = u\big(\supp(x)\big).\]
\item\label{sch2.13_iv}
Suppose $p \cn x \to y$ is a morphism in $\C$, and $u$ and $v$ agree on $\supp(x) \cup \supp(y)$.  Then there is a morphism equality
\[u_{\mstar}(x) = v_{\mstar}(x) \fto{u_{\mstar}(p) = v_{\mstar}(p)} u_{\mstar}(y) = v_{\mstar}(y).\]
\item\label{sch2.13_v}
For each $\schm$-functor $\fun \cn \C \to \D$, there is a set inclusion
\[\supp(\fun x) \subseteq \supp(x).\]
\end{enumerate}
\end{lemma}

\subsection*{Reparametrization}

Next, we recall reparametrization of $\schm$-categories \cite[2.20]{schwede_global}.

\begin{definition}[Reparametrization]\label{def:repar}
Recall the countable set $\ome = \{0,1,2,\ldots\}$.
\begin{enumerate}
\item\label{def:repar_i} 
For each countably infinite set $S$, we choose a bijection 
\begin{equation}\label{rep_S}
\ome \fto[\iso]{\rep_S} S
\end{equation}
such that $\rep_\ome = 1_\ome$.
\item\label{def:repar_ii} 
Denote by $\schj$ the 2-category with
\begin{itemize}
\item countably infinite sets as objects,
\item injections as 1-cells, 
\item a unique 2-cell for each pair of parallel 1-cells, and
\item horizontal composition of 1-cells given by composition of injections.
\end{itemize}
Note that for each pair $(U,V)$ of countably infinite sets, the hom-category $\schj(U,V)$ is a translation category with injections $U \to V$ as objects \pcref{def:translation_cat}.
\item\label{def:repar_iiii} 
For each $\schm$-category $\C$, the \emph{reparametrization}\index{reparametrization} of $\C$ is the 2-functor
\begin{equation}\label{Cbrac}
\schj \fto{\Cbrac} \Cat
\end{equation}
defined as follows.
\begin{description}
\item[Objects] $\Cbrac$ sends each countably infinite set $U$ to the category $\CbracU = \C$.
\item[1-cells] $\Cbrac$ sends an injection $p \cn U \to V$ between countably infinite sets to the $\schm$-action functor on $\C$ 
\begin{equation}\label{Cbrac_p}
\begin{tikzpicture}[vcenter]
\def\v{-1}
\draw[0cell]
(0,0) node (a11) {\CbracU}
(a11)++(2.5,0) node (a12) {\CbracV}
(a11)++(0,\v) node (a21) {\C}
(a12)++(0,\v) node (a22) {\C}
;
\draw[1cell=.9]
(a11) edge[equal] (a21)
(a12) edge[equal] (a22)
(a11) edge node {\Cbracp} (a12)
(a21) edge node {p^\ome _{\mstar}} (a22)
;
\end{tikzpicture}
\end{equation}
for the injection
\begin{equation}\label{repU_p_repV}
\begin{tikzpicture}[vcenter]
\def\h{1.6} \def\u{.7}
\draw[0cell]
(0,0) node (a1) {\ome}
(a1)++(\h,0) node (a2) {U}
(a2)++(.9*\h,0) node (a3) {V}
(a3)++(\h,0) node (a4) {\ome.}
;
\draw[1cell=.9]
(a1) edge node {\rep_U} node[swap] {\iso} (a2)
(a2) edge node {p} (a3)
(a3) edge node {\rep_V^{-1}} node[swap] {\iso} (a4)
(a1) [rounded corners=2pt] |- ($(a2)+(0,\u)$) -- node {p^\ome} ($(a3)+(0,\u)$) -| (a4)
;
\end{tikzpicture}
\end{equation}
\item[2-cells] For injections $p,q \cn U \to V$, $\Cbrac$ sends the unique 2-cell $[q,p] \cn p \to q$ to the $\schm$-action natural isomorphism on $\C$ 
\begin{equation}\label{Cbrac_iicell}
\Cbracp = p^\ome_{\mstar} \fto{[q^\ome,p^\ome]_{\mstar}} q^\ome_{\mstar} = \Cbracq
\end{equation}
for the unique morphism $[q^\ome,p^\ome] \cn p^\ome \to q^\ome$ in $\schm$ \cref{vu_star}.  The 2-functor axioms for $\Cbrac$ follow from the $\schm$-category axioms for $\C$ \cref{mcat_axioms}.
\end{description}
An object $x \in \CbracU$ is \emph{supported} on a subset $A \subseteq U$ if $p^\ome_{\mstar}(x) = x$ for each injection $p \cn U \to U$ that is the identity on $A$.  An object $x \in \CbracU$ is \emph{finitely supported} if it is supported on some finite subset of $U$.
\item\label{def:repar_iv}
For each $\schm$-functor $\fun \cn \C \to \D$, the \emph{reparametrization} of $\fun$ is the 2-natural transformation
\begin{equation}\label{fbrac}
\begin{tikzpicture}[vcenter]
\def\t{29}
\draw[0cell]
(0,0) node (a1) {\schj}
(a1)++(2,0) node (a2) {\phantom{\schj}}
(a2)++(.15,0) node (a2') {\Cat}
;
\draw[1cell=.8]
(a1) edge[bend left=\t] node {\Cbrac} (a2)
(a1) edge[bend right=\t] node[swap] {\Dbrac} (a2)
;
\draw[2cell=.9]
node[between=a1 and a2 at .35, rotate=-90, 2label={above,\fbrac}] {\Rightarrow}
;
\end{tikzpicture}
\end{equation}
with each component functor $\fU \cn \CbracU \to \DbracU$ given by $\fun$.  The 2-naturality of $\fbrac$ follows from the $\schm$-functor axiom for $\fun$ \cref{mfunctor_axiom}.\defmark
\end{enumerate}
\end{definition}

\begin{explanation}[Preservation of Finitely Supportedness]\label{expl:cu_support}
Consider an $\schm$-category $\C$, an object $x \in \C$ supported on a subset $A \subseteq \ome$, and a countably infinite set $U$.  If an injection $p \cn U \to U$ is the identity on $\rep_U(A) \subseteq U$, then the injection $p^\ome \cn \ome \to \ome$ \cref{repU_p_repV} is the identity on $A$, implying $p^\ome_{\mstar}(x) = x$.  Thus, the object $x \in \CbracU$ is supported on the subset $\rep_U(A) \subseteq U$.  In particular, if $x \in \C$ is finitely supported, then $x \in \CbracU$ is also finitely supported.
\end{explanation}

Recall the universal $G$-set $\omeg$ of functions $G \to \ome$ \pcref{ex:omea}, with $g \in G$ sending a function $f \cn G \to \ome$ to $gf = f(\ginv \cdot -)$.  \cref{def:sch2.21} uses the reparametrization $\Cbrac$ \cref{Cbrac} for the countably infinite set $\omeg$.

\begin{definition}[$G$-Fixed $\schm$-Categories]\label{def:sch2.21}
Suppose $G$ is a finite group, and $\C$ is an $\schm$-category.
\begin{enumerate}
\item\label{def:sch2.21_i}
The \emph{$\omeg$-reparametrization}\index{reparametrization!omeg@$\omeg$} of $\C$ is the $\gschm$-category $\Comeg$ defined as follows.
\begin{description}
\item[Category] The underlying category of $\Comeg$ is $\C$.
\item[$G$-action] For each $g \in G$, the $g$-action functor on $\Comeg$ is the $\schm$-action functor on $\C$ 
\begin{equation}\label{Cbrac_g}

\end{equation}
for the bijection
\begin{equation}\label{rep_g_rep}
%
\end{equation}
\item[$\schm$-action] 
For the injection monoid $\Injm$ of injections $\ome \to \ome$ (\cref{def:mcat} \cref{def:mcat_i}), there is an $\Injm$-action on $\omeg$ given by postcomposition: 
\begin{equation}\label{u_omeg}
(uf)(h) = u(f(h))
\end{equation}
for each injection $u \cn \ome \to \ome$, $f \in \omeg$, and $h \in G$.  The $u$-action functor on $\Comeg$ is the $\schm$-action functor on $\C$ 
\begin{equation}\label{Cbrac_u}
%
\end{equation}
for the injection
\begin{equation}\label{rep_u_rep}
%
\end{equation}
For another injection $v \cn \ome \to \ome$ and the isomorphism $[v,u] \cn u \to v$ in $\schm$, the $\schm$-action natural isomorphism on $\Comeg$ is the $\schm$-action natural isomorphism on $\C$ 
\[\Cbracu = u^\ome_{\mstar} \fto{[v^\ome,u^\ome]_{\mstar}} v^\ome_{\mstar} = \Cbracv\]
for the unique morphism $[v^\ome,u^\ome] \cn u^\ome \to v^\ome$ in $\schm$.  See \cref{Cbrac_iicell}.  
\end{description}
The $G$-action commutes with the $\Injm$-action \cref{u_omeg} on $\omeg$, with
\begin{equation}\label{guf}
(g,u)f = uf(\ginv \cdot -).
\end{equation}
Thus, the $G$-action \cref{Cbrac_g} commutes with the $\schm$-action \cref{Cbrac_u} on $\Comeg$, giving it the structure of a $\gschm$-category.
\item\label{def:sch2.21_ii}
Define the \emph{$G$-fixed $\schm$-category}\index{M-category@$\schm$-category!G-fixed@$G$-fixed}
\begin{equation}\label{fixgc}
\Fixg\C = \Comeg^G
\end{equation}
as the $G$-fixed $\schm$-subcategory of the $\gschm$-category $\Comeg$.
\item\label{def:sch2.21_iii}
Suppose $\fun \cn \C \to \D$ is an $\schm$-functor.
\begin{itemize}
\item The \emph{$\omeg$-reparametrization} of $\fun$ is the $\gschm$-functor
\begin{equation}\label{fomeg}
\Comeg \fto{\fomeg} \Domeg
\end{equation}
given by the functor $\fun$.
\item The \emph{$G$-fixed $\schm$-functor}\index{M-functor@$\schm$-functor!G-fixed@$G$-fixed}
\begin{equation}\label{fixgf}
\Fixg\C \fto{\Fixg\fun} \Fixg\D
\end{equation}
is the $G$-fixed $\schm$-subfunctor $\fomeg^G$ of the $\gschm$-functor $\fomeg$.\defmark
\end{itemize}
\end{enumerate}
\end{definition}

\begin{explanation}[Preservation of Tameness]\label{expl:comeg_tame}
Suppose $\C$ is a tame $\schm$-category. 
Then the $\schm$-categories $\Comeg$ and $\Fixg\C = \Comeg^G$ are tame \cite[2.22]{schwede_global}.  Indeed, under the $\schm$-action on $\C$, any given object $x \in \Comeg$ is supported on some finite subset $A \subseteq \ome$.  Suppose $u \cn \ome \to \ome$ is an injection that is the identity on the finite subset
\[A' = \bigcup_{a \in A}\, \Img(\rep_{\omeg}(a)) \subseteq \ome\]
consisting of the images of $\rep_{\omeg}(a) \cn G \to \ome$ for $a \in A$.  Then the injection $u^\ome \cn \ome \to \ome$ \cref{rep_u_rep} is the identity on $A$ by \cref{u_omeg}, so $u^\ome_{\mstar}(x) = x$.  Thus, under the $\schm$-action on $\Comeg$ \cref{Cbrac_u}, $x$ is supported on the finite subset $A' \subseteq \ome$.
\end{explanation}

\subsection*{Global Equivalences}

A functor $\fun$ between small categories is called a \emph{weak equivalence} if its image $\cla\fun$ under the classifying space functor $\cla \cn \Cat \to \Top$ \cref{classifying_space} is a weak homotopy equivalence of spaces.  \cref{def:sch2.26} uses \cref{schwede2.17} \cref{sch2.17_ii}---that any two universal $G$-sets are $G$-isomorphic---and the universal $G$-set $\omeg$ \pcref{ex:omea}.

\begin{definition}\label{def:sch2.26}
An $\schm$-functor $\fun \cn \C \to \D$ is called a \emph{global equivalence}\index{global equivalence} if the following three equivalent conditions hold for each finite group $G$.
\begin{enumerate}
\item\label{gleq_i}
The functor $\Fixg\fun \cn \Fixg\C \to \Fixg\D$ \cref{fixgf} is a weak equivalence.
\item\label{gleq_ii}
For some universal $G$-set $\ugsetu$, the $G$-fixed subfunctor
\[\Cugsetu^G \fto{\fugsetu^G} \Dugsetu^G\]
of $\fugsetu$ \cref{fbrac} is a weak equivalence.
\item\label{gleq_iii}
For each universal $G$-set $\ugsetu$, the functor $\fugsetu^G$ is a weak equivalence.
\end{enumerate}
Moreover, $\fun$ is called a \emph{global categorical equivalence}\index{global categorical equivalence} if the equivalent conditions \crefrange{gleq_i}{gleq_iii} hold with weak equivalence replaced by equivalence of categories.
\end{definition}

\begin{example}\label{ex:gleq}
Each global categorical equivalence is also a global equivalence.  By \cite[2.31]{schwede_global}, if an $\schm$-functor $\fun \cn \C \to \D$ is a global (categorical) equivalence, then so is the $G$-fixed $\schm$-functor $\Fixg\fun \cn \Fixg\C \to \Fixg\D$ \cref{fixgf} for each finite group $G$.
\end{example}

\section{Parsummable Categories to $\FM$-Categories}
\label{sec:jgl}

This section reviews parsummable categories and their associated $\FM$-categories.  Parsummable categories are the input data for global $K$-theory.  The passage to $\FM$-categories is the first step of global $K$-theory.

\secoutline
\begin{itemize}
\item \cref{def:bxtimes} defines the box product $\Bxtimes$ for $\schm$-categories.
\item \cref{def:parcat} defines parsummable categories as commutative monoids in the symmetric monoidal category $\Mcatt$ of tame $\schm$-categories equipped with the box product.
\item \cref{def:fmcat} defines $\FM$-categories as pointed functors $\Fsk \to \Mcat$.
\item \cref{def:sch4.3} defines the first step of global $K$-theory: the functor
\[\Parcat \fto{\Jsc} \FMCat\]
called \emph{global $J$-theory}.
\end{itemize}

\subsection*{Box Product}

Parsummable categories are defined using the box product in \cref{def:bxtimes}.

\begin{definition}\label{def:bxtimes}
Suppose $\C$ and $\D$ are $\schm$-categories \pcref{def:mcat}.
\begin{enumerate}
\item\label{def:bxtimes_i}
An object $(x,y) \in \C \ttimes \D$ is \emph{disjointly supported}\index{disjointly supported} if there exist disjoint subsets 
\[A, B \subseteq \ome = \{0,1,2,\ldots\}\] 
such that $x$ is supported on $A$ and $y$ is supported on $B$.  In this case, we also say that $x$ and $y$ are disjointly supported.
\item\label{def:bxtimes_ii}
The \emph{box product}\index{box product} $\C\bxtimes\D$ is the full $\schm$-subcategory of the product $\schm$-category $\C\ttimes\D$ consisting of disjointly supported objects.  The box product is well defined by \cref{sch2.13} \cref{sch2.13_iii}.  Each object $(x,y) \in \C \ttimes \D$ is supported on $\supp(x) \cup \supp(y)$ by \cref{sch2.13} \cref{sch2.13_i}.  Thus, if $\C$ and $\D$ are tame $\schm$-categories, then so is $\C\bxtimes\D$ \cite[2.34]{schwede_global}.
\item\label{def:bxtimes_iii}
Given $\schm$-functors $\fun \cn \C \to \D$ and $\fun' \cn \C' \to \D'$, the $\schm$-functor
\[\C\bxtimes\C' \fto{\fun\bxtimes\fun'} \D\bxtimes\D'\]
is given by the restriction of $\fun\ttimes\fun'$ to $\C\bxtimes\C'$.  This is well defined by \cref{sch2.13} \cref{sch2.13_v}.
\item\label{def:bxtimes_iv}
Equip the category $\Mcatt$ of tame $\schm$-categories and $\schm$-functors \pcref{def:mcat} with the symmetric monoidal structure $(\Bxtimes, \bone, \bxi)$.  The unit object $\bone$, unit isomorphisms, associativity isomorphism, and braiding $\bxi$ are inherited from the corresponding structures for the symmetric monoidal category $(\Cat,\times,\bone,\bxi)$.\defmark
\end{enumerate}
\end{definition}

\begin{example}\label{ex:sch2.33}
If either $x \in \C$ or $y\in \D$ is $\schm$-fixed, meaning the support is the empty set, then $(x,y) \in \C\ttimes\D$ is disjointly supported.  By \cite[2.33]{schwede_global}, for any $\schm$-categories $\C$ and $\D$, the inclusion $\schm$-functor\label{not:bxinc} 
\[\C\bxtimes\D \fto{\iota} \C\ttimes\D\]
is a global categorical equivalence, hence also a global equivalence \pcref{def:sch2.26}.
\end{example}

\subsection*{Parsummable Categories}

\begin{definition}\label{def:parcat}
A \emph{parsummable category}\index{parsummable category} is a commutative monoid in the symmetric monoidal category $(\Mcatt,\Bxtimes,\bone,\bxi)$ (\cref{def:bxtimes} \cref{def:bxtimes_iv}).  A \emph{parsummable functor}\index{parsummable functor} between parsummable categories is a morphism of commutative monoids in $\Mcatt$.  The category of parsummable categories and parsummable functors is denoted by $\Parcat$.
\end{definition}

\begin{explanation}[Unpacking]\label{expl:parcat}
A parsummable category $(\C,\psum,\pzero)$ consists of
\begin{itemize}
\item a tame $\schm$-category $\C$, 
\item an $\schm$-functor $\psum \cn \C\bxtimes\C \to \C$ called the \index{parsummable category!sum}\emph{sum}, and
\item an $\schm$-fixed object $\pzero \in \C$ called the \index{parsummable category!unit}\emph{unit}
\end{itemize}
such that the following associativity, unity, and symmetry diagrams commute, where iterated box products are interpreted using \cref{conv:left_normal}.
\begin{equation}\label{parcat_axioms}
\begin{tikzpicture}[vcenter]
\def\h{2.7} \def\d{2.2} \def\v{-1.4}
\draw[0cell]
(0,0) node (a11) {\C\bxtimes\C\bxtimes\C}
(a11)++(\h,0) node (a12) {\C\bxtimes\C}
(a11)++(0,\v) node (a21) {\C\bxtimes\C}
(a12)++(0,\v) node (a22) {\C}
;
\draw[1cell=.9]
(a11) edge node {1_\C \bxtimes \psum} (a12)
(a12) edge node {\psum} (a22)
(a11) edge node[swap] {\psum \bxtimes 1_\C} (a21)
(a21) edge node {\psum} (a22)
;
\begin{scope}[shift={(\h+1.5,0)}]
\draw[0cell]
(0,0) node (b11) {\bone\bxtimes\C}
(b11)++(\d,0) node (b12) {\C\bxtimes\C}
(b12)++(\d,0) node (b13) {\C\bxtimes\C}
(b12)++(0,\v) node (b2) {\C}
;
\draw[1cell=.9]
(b11) edge node {\pzero \bxtimes 1_{\C}} (b12)
(b11) edge node[swap,pos=.4] {\iso} (b2)
(b12) edge node {\psum} (b2)
(b13) edge node[swap] {\bxi} (b12)
(b13) edge node[pos=.4] {\psum} (b2)
;
\end{scope}
\end{tikzpicture}
\end{equation}
Applying \cref{sch2.13} \cref{sch2.13_v} to the $\schm$-functor $\psum$, for each object $(x,y) \in \C\bxtimes\C$, there is an inclusion\index{subadditivity}
\begin{equation}\label{supp_xsumy}
\supp(x \psum y) \subseteq \supp(x,y) = \supp(x) \sqcup \supp(y)
\end{equation}
called the \emph{subadditivity} of support.  For parsummable categories $\C$ and $\D$, a parsummable functor $\fun \cn \C \to \D$ is an $\schm$-functor such that the compatibility diagrams 
\begin{equation}\label{parfun_axioms}
\begin{tikzpicture}[vcenter]
\def\h{2} \def\d{1.5} \def\v{-1.4} \def\u{-.7}
\draw[0cell]
(0,0) node (a11) {\C\bxtimes\C}
(a11)++(\h,0) node (a12) {\C}
(a11)++(0,\v) node (a21) {\D\bxtimes\D}
(a12)++(0,\v) node (a22) {\D}
;
\draw[1cell=.9]
(a11) edge node {\psum} (a12)
(a12) edge node {\fun} (a22)
(a11) edge node[swap] {\fun\bxtimes\fun} (a21)
(a21) edge node {\psum} (a22)
;
\begin{scope}[shift={(\h+1.5,\u)}]
\draw[0cell]
(0,0) node (b1) {\bone}
(b1)++(\d,-\u) node (b2) {\C}
(b1)++(\d,\u) node (b3) {\D}
;
\draw[1cell=.9]
(b1) edge node {\pzero} (b2)
(b2) edge node {\fun} (b3)
(b1) edge node[swap] {\pzero} (b3)
;
\end{scope}
\end{tikzpicture}
\end{equation}
commute.
\end{explanation}

\begin{example}[Limits and Coproducts]\label{ex:sch4.10}
The category $\Parcat$ of parsummable categories has all small limits, and they are created by the forgetful functor $\Parcat \to \Mcatt$ \cite[4.11]{schwede_global}.  For parsummable categories $\C$ and $\D$, the box product $\C\bxtimes\D$ is also a parsummable category with unit $(\pzero,\pzero)$ and the sum given by the composite $\schm$-functor
\[(\C\bxtimes\D)\bxtimes(\C\bxtimes\D) \fto{1\bxtimes\bxi\bxtimes 1}
 (\C\bxtimes\C)\bxtimes(\D\bxtimes\D) \fto{\psum\bxtimes\psum} \C\bxtimes\D.\]
The parsummable functors
\[\C \fto{i_1 \,=\, (-,\pzero)} \C\bxtimes\D \fot{i_2 \,=\, (\pzero,-)} \D\]
make $\C\bxtimes\D$ into a coproduct of $\C$ and $\D$ in $\Parcat$.  See \cite[4.10]{schwede_global}.
\end{example}

\begin{example}[Finite Sets]\label{ex:sch4.5}
Recall the tame $\schm$-category $\Fome$ of finite subsets of $\ome$ and bijections \pcref{ex:sch2.14}.  An object in $\Fome \bxtimes \Fome$ is a pair of disjoint finite subsets of $\ome$.  There is a parsummable structure on $\Fome$ with unit given by the empty set $\emptyset$.  The $\schm$-functor
\[\Fome \bxtimes \Fome \fto{\psum} \Fome\]
sends a pair $(A,B) \in \Fome \bxtimes \Fome$ to their disjoint union $A \sqcup B \subseteq \ome$.  For bijections 
\[A \fto[\iso]{p} B \andspace A' \fto[\iso]{p'} B'\] 
in $\Fome$ such that $A \cap A' = \emptyset = B \cap B'$, the sum 
\[A \sqcup A' \fto{p \psum p'} B \sqcup B'\]
is the disjoint union of $p$ and $p'$.  See \cite[4.5]{schwede_global}.
\end{example}

\subsection*{$\FM$-Categories}

Recall the pointed category $(\Fsk,\ordz)$ of pointed finite sets $\ordn = \{0,1,\ldots,n\}$ with basepoint 0 for $n \geq 0$ and pointed morphisms \pcref{def:Fsk}.  The category $\Mcat$ of $\schm$-categories and $\schm$-functors (\cref{def:mcat} \cref{def:mcat_iii}) becomes a pointed category with the terminal $\schm$-category $\bone$ as the basepoint.

\begin{definition}\label{def:fmcat}
An \emph{$\FM$-category}\index{FM-category@$\FM$-category} is a pointed functor
\[(\Fsk,\ordz) \fto{X} (\Mcat,\bone).\]
A \emph{morphism} of $\FM$-categories is a natural transformation.  The category of $\FM$-categories and morphisms is denoted by $\FMCat$.
\end{definition}

\begin{explanation}[Unpacking]\label{expl:fmcat}
An $\FM$-category $X$ consists of
\begin{itemize}
\item an $\schm$-category $X\ordn$ for each $n \geq 0$ and
\item an $\schm$-functor
\[X\ordm \fto{X\psi} X\ordn\]
for each pointed morphism $\psi \cn \ordm \to \ordn$ 
\end{itemize}
such that
\begin{equation}\label{fmcat_axioms}
\begin{split}
X\ordz &= \bone,\\
X1_{\ordm} &= 1_{X\ordm}, \andspace\\ 
X(\phi\psi) &= (X\phi)(X\psi)
\end{split}
\end{equation}
for composable pointed morphisms $\psi \cn \ordm \to \ordn$ and $\phi \cn \ordn \to \ordr$.  For the unique pointed morphism $\ordz \to \ordn$, the $\schm$-functor 
\[\bone = X\ordz \to X\ordn\]
equips each $X\ordn$ with an $\schm$-fixed basepoint object 0.  By the functoriality of $X$, each $\schm$-functor $X\psi$ preserves the $\schm$-fixed basepoint 0.

A morphism $\tha \cn X \to Y$ of $\FM$-categories assigns to each object $\ordn \in \Fsk$ an $\schm$-functor
\[X\ordn \fto{\tha_{\ordn}} Y\ordn\]
such that, for each pointed morphism $\psi \cn \ordm \to \ordn$, the naturality diagram
\begin{equation}\label{fmcat_mor_nat}
\begin{tikzpicture}[vcenter]
\def\v{-1.4}
\draw[0cell]
(0,0) node (a11) {X\ordm}
(a11)++(2.,0) node (a12) {Y\ordm}
(a11)++(0,\v) node (a21) {X\ordn}
(a12)++(0,\v) node (a22) {Y\ordn}
;
\draw[1cell=.9]
(a11) edge node {\tha_{\ordm}} (a12)
(a12) edge node {Y\psi} (a22)
(a11) edge node[swap] {X\psi} (a21)
(a21) edge node {\tha_{\ordn}} (a22)
;
\end{tikzpicture}
\end{equation}
in $\Mcat$ commutes.  Each morphism $\tha$ is necessarily pointed, meaning $\tha_{\ordz} = 1_{\bone}$.  Moreover, the naturality diagram \cref{fmcat_mor_nat} for the unique pointed morphism $\ordz \to \ordn$ implies that each $\schm$-functor $\tha_{\ordn}$ preserves the $\schm$-fixed basepoint 0.
\end{explanation}

\subsection*{Schwede Global $J$-Theory}

Next, we recall the $\FM$-categories associated to parsummable categories, denoted by $\ga(-)$ in \cite[4.3]{schwede_global}.  Since this construction is analogous to $J$-theory \pcref{thm:Jgo_twofunctor} and Shimakawa $J$-theory \cref{jgos_jgossg}, we use the notation $\Jsc$.

\begin{definition}\label{def:sch4.3}
The \emph{global $J$-theory}\index{global J-theory@global $J$-theory}\index{J-theory@$J$-theory!global} functor
\begin{equation}\label{Jsc}
\Parcat \fto{\Jsc} \FMCat
\end{equation}
is defined as follows.
\begin{description}
\item[Objects] 
For a parsummable category $(\C,\psum,\pzero)$, the pointed functor
\begin{equation}\label{JscC}
(\Fsk,\ordz) \fto{\Jsc\C} (\Mcat,\bone)
\end{equation}
sends each object $\ordn \in \Fsk$ to the full $\schm$-subcategory
\begin{equation}\label{JscCn}
(\Jsc\C)\ordn = \C^{\Bxtimes n} \subseteq \C^n
\end{equation}
of disjointly supported $n$-tuples $\ang{x_j}_{j \in \ufsn} \in \C^n$.  The $\schm$-category $(\Jsc\C)\ordz$ is defined as the terminal $\schm$-category $\bone$.  For a pointed morphism $\psi \cn \ordm \to \ordn$, the $\schm$-functor
\[(\Jsc\C)\ordm = \C^{\Bxtimes m} \fto{(\Jsc\C)\psi} (\Jsc\C)\ordn = \C^{\Bxtimes n}\]
sends an $m$-tuple $\ang{x_i}_{i \in \ufsm} \in \C^{\Bxtimes m}$ of all objects or all morphisms to the $n$-tuple
\begin{equation}\label{JscCpsi}
\big((\Jsc\C)\psi\big) \ang{x_i}_{i \in \ufsm} 
= \bang{\,\txsum_{\,i \in \psiinv j} x_i}_{j \in \ufsn} \in \C^{\Bxtimes n}.
\end{equation}
When $\psiinv j = \emptyset$, an empty sum means the unit $\pzero \in \C$ or its identity morphism.  The $n$-tuple in \cref{JscCpsi} is disjointly supported by subadditivity \cref{supp_xsumy}.  The functoriality of $\Jsc\C$ follows from the parsummable category axioms \cref{parcat_axioms} for $\C$.
\item[Morphisms] 
For a parsummable functor $\fun \cn \C \to \D$ between parsummable categories, the natural transformation
\begin{equation}\label{Jscf}
\begin{tikzpicture}[vcenter]
\def\t{30}
\draw[0cell]
(0,0) node (a1) {\Fsk}
(a1)++(2,0) node (a2) {\phantom{\Fsk}}
(a2)++(.3,0) node (a2') {\Mcat}
;
\draw[1cell=.85]
(a1) edge[bend left=\t] node {\Jsc\C} (a2)
(a1) edge[bend right=\t] node[swap] {\Jsc\D} (a2)
;
\draw[2cell]
node[between=a1 and a2 at .35, rotate=-90, 2label={above,\Jsc\fun}] {\Rightarrow}
;
\end{tikzpicture}
\end{equation}
sends each object $\ordn \in \Fsk$ to the $\schm$-functor
\[(\Jsc\C)\ordn = \C^{\Bxtimes n} \fto{(\Jsc\fun)_{\ordn} \,=\, \fun^{\Bxtimes n}} 
(\Jsc\D)\ordn = \D^{\Bxtimes n}.\]
The preceding $\schm$-functor is well defined by \cref{sch2.13} \cref{sch2.13_v}.  The naturality of $\Jsc\fun$ follows from the parsummable functor axioms \cref{parfun_axioms} for $\fun$.  The functoriality of $\Jsc$ follows from the fact that $\fun^{\Bxtimes n}$ is given entrywise by $\fun$.
\end{description}
This finishes the definition of $\Jsc$.
\end{definition}

\begin{example}[Commutative Monoids]\label{ex:sch4.4}
Given a commutative monoid $(\A,\psum,\pzero)$, regard it as a discrete $\schm$-category with the trivial $\schm$-action.  Thus, each object $a \in \A$ is $\schm$-fixed and supported on the empty set.  The global $J$-theory of $\A$ has the discrete trivial $\schm$-category 
\[(\Jsc\A)\ordn = \A^{\Bxtimes n} = \A^n\]
for each $n \geq 0$.  This example is from \cite[4.4 and 4.29]{schwede_global}.
\end{example}

\section{$\FM$-Categories to Symmetric Spectra}
\label{sec:kfm}

This section reviews the functors
\[\FMCat \fto{\uprst} \FICat \fto{\clast} \FITop \fto{\Kfi} \Sptop\]
from the category $\FMCat$ of $\FM$-categories to the category $\Sptop$ of symmetric spectra.

\secoutline
\begin{itemize}
\item \cref{def:sch8.2} defines the categories $\ICat$ and $\ITop$ of $\bdI$-categories and $\bdI$-spaces, where $\bdI$ is the category of finite sets and injections.
\item \cref{def:upr} defines the functor $\upr \cn \Mcat \to \ICat$ from $\schm$-categories to $\bdI$-categories.  \cref{def:ficat} defines the postcomposition functors $\uprst$ and $\clast$.
\item \cref{def:sptop} defines symmetric spectra, (underlying) $G$-symmetric spectra, $G$-stable equivalences, and global equivalences of symmetric spectra.  \cref{expl:sptop} briefly discusses global and $G$-equivariant stable model structures.
\item \cref{def:kfi} defines the functor $\Kfi$.
\end{itemize}

\subsection*{$\FM$-Categories to $\FI$-Spaces}

Recall that $\Cat$ is the category of small categories and functors \pcref{def:GCat} and that $\Top$ is the category of compactly generated weak Hausdorff spaces and continuous morphisms \pcref{def:Gtop}.

\begin{definition}[$\bdI$-Categories and $\bdI$-Spaces]\label{def:sch8.2}\
\begin{enumerate}
\item\label{def:sch8.2_i} Denote by $\bdI$ the category of finite sets and injections.  The category $\bdI$ has a small skeleton consisting of the finite sets $\ufsn = \{1,2,\ldots,n\}$ for $n \geq 0$ and injections among them.  
\item\label{def:sch8.2_ii} An \emph{$\bdI$-category}\index{I-category@$\bdI$-category} is a functor $X \cn \bdI \to \Cat$.  A \emph{morphism} of $\bdI$-categories is a natural transformation.  The category of $\bdI$-categories and morphisms is denoted by $\ICat$.  Denote by $\bone \in \ICat$ the constant functor at the terminal category $\bone$.
\item\label{def:sch8.2_iii} An \emph{$\bdI$-space}\index{I-space@$\bdI$-space} is a functor $X \cn \bdI \to \Top$.  A \emph{morphism} of $\bdI$-spaces is a natural transformation.  The category of $\bdI$-spaces and morphisms is denoted by $\ITop$.  Denote by $* \in \ITop$ the constant functor at the one-point space $*$.
\end{enumerate}
Both $\ICat$ and $\ITop$ are well-defined categories with small hom sets because $\bdI$ has a small skeleton.
\end{definition}

\begin{definition}[$\schm$-Categories to $\bdI$-Categories]\label{def:upr}
Define the functor
\begin{equation}\label{upr_functor}
\Mcat \fto{\upr} \ICat
\end{equation}
as follows.
\begin{description}
\item[Objects]
For an $\schm$-category $\C$ (\cref{def:mcat} \cref{def:mcat_iii}), the $\bdI$-category $\upr\C \cn \bdI \to \Cat$ sends a finite set $A$ to the category
\begin{equation}\label{uprca}
(\upr\C)A = \begin{cases}
\Comea & \text{if $A \neq \emptyset$ and}\\
\C^{\suppempty} & \text{if $A=\emptyset$.}
\end{cases}
\end{equation}
\begin{itemize}
\item For $A \neq \emptyset$, $\Comea$ is the reparametrization of $\C$ \cref{Cbrac} at the countably infinite set $\omea$ of functions $A \to \ome = \{0,1,2,\ldots\}$ \pcref{def:omea}.  Recall that the underlying category of $\Comea$ is $\C$.
\item For $A = \emptyset$, $\C^{\suppempty}$ is the full subcategory of $\C$ consisting of $\schm$-fixed objects, meaning objects with empty support.  
\end{itemize}
For an injection $i \cn A \to B$ between finite sets, the functor
\begin{equation}\label{uprci}
(\upr\C)A \fto{(\upr\C)i} (\upr\C)B
\end{equation}
is defined as follows.
\begin{itemize}
\item If $A = \emptyset = B$, then $(\upr\C)i$ is the identity functor on $\C^{\suppempty}$.
\item If $A = \emptyset \neq B$, then $(\upr\C)i$ is the full subcategory inclusion 
\begin{equation}\label{uprci_emptyb}
\C^{\suppempty} \fto{\iota} \Comeb = \C.
\end{equation} 
\item If $A \neq \emptyset$, then we consider the injection
\begin{equation}\label{extzero}
\omea \fto{i_!} \omeb
\end{equation}
defined by
\[(i_! f)(b) = \begin{cases}
f(a) & \text{if $b = i(a)$ and}\\
0 & \text{if $b \not\in i(A)$}
\end{cases}\]
for $f \in \omea$, $a \in A$, and $b \in B$.  The functor $(\upr\C)i$ is defined as the $\schm$-action functor 
\begin{equation}\label{uprci_nonempty}
\begin{tikzpicture}[vcenter]
\def\v{-1} \def\h{1.2}
\draw[0cell]
(0,0) node (a11) {(\upr\C)A}
(a11)++(4.3,0) node (a12) {(\upr\C)B}
(a11)++(0,\v) node (a21) {\Comea}
(a12)++(0,\v) node (a22) {\Comeb}
(a21)++(\h,0) node (b1) {\C}
(a22)++(-\h,0) node (b2) {\C}
;
\draw[1cell=.9]
(a11) edge node {(\upr\C)i} (a12)
(a11) edge[equal] (a21)
(a12) edge[equal] (a22)
(a21) edge[equal] (b1)
(b2) edge[equal] (a22)
(b1) edge node {(i_!^\ome)_{\mstar}} (b2)
;
\end{tikzpicture}
\end{equation}
on $\C$ for the injection
\begin{equation}\label{iome_ext}
\begin{tikzpicture}[vcenter]
\def\h{1.8} \def\g{1.6} \def\u{.7} \def\d{.05}
\draw[0cell]
(0,0) node (a1) {\ome}
(a1)++(\h,0) node (a2) {\phantom{\omeg}}
(a2)++(0,\d) node (a2') {\omea}
(a2)++(\g,0) node (a3) {\phantom{\omeg}}
(a3)++(0,\d) node (a3') {\omeb}
(a3)++(\h,0) node (a4) {\ome.}
;
\draw[1cell=.9]
(a1) edge node {\rep_{\omea}} node[swap] {\iso} (a2)
(a2) edge node {i_!} (a3)
(a3) edge node {\rep_{\omeb}^{-1}} node[swap] {\iso} (a4)
(a1) [rounded corners=2pt] |- ($(a2)+(0,\u)$) -- node {i_!^\ome} ($(a3)+(0,\u)$) -| (a4)
;
\end{tikzpicture}
\end{equation}
The functor $(i_!^\ome)_{\mstar}$ is an instance of \cref{Cbrac_p}. 
\end{itemize}
\item[Morphisms]
For an $\schm$-functor $\fun \cn \C \to \D$ between $\schm$-categories, the natural transformation
\begin{equation}\label{uprf}
\begin{tikzpicture}[vcenter]
\def\t{29}
\draw[0cell]
(0,0) node (a1) {\bdI}
(a1)++(2,0) node (a2) {\phantom{\bdI}}
(a2)++(.15,0) node (a2') {\Cat}
;
\draw[1cell=.85]
(a1) edge[bend left=\t] node {\upr\C} (a2)
(a1) edge[bend right=\t] node[swap] {\upr\D} (a2)
;
\draw[2cell=1]
node[between=a1 and a2 at .4, rotate=-90, 2label={above,\upr\fun}] {\Rightarrow}
;
\end{tikzpicture}
\end{equation}
has, for each finite set $A$, $A$-component given by
\begin{itemize}
\item the functor $\fun \cn \Comea \to \Domea$ if $A \neq \emptyset$ and
\item the restriction $\fun|_{\C^{\suppempty}} \cn \C^{\suppempty} \to \D^{\suppempty}$ if $A = \emptyset$.
\end{itemize}
The naturality of $\upr\fun$ follows from the $\schm$-functor axiom \cref{mfunctor_axiom} for $\fun$.
\end{description}
This finishes the definition of the functor $\upr$.
\end{definition}

\begin{definition}[$\FI$-Categories and $\FI$-Spaces]\label{def:ficat}\
\begin{enumerate}
\item\label{def:ficat_i} An \emph{$\FI$-category}\index{FI-category@$\FI$-category} is a pointed functor
\[(\Fsk,\ordz) \fto{X} (\ICat,\bone).\]
A \emph{morphism} of $\FI$-categories is a natural transformation.  The category of $\FI$-categories and morphisms is denoted by $\FICat$.
\item\label{def:ficat_ii} An \emph{$\FI$-space}\index{FI-space@$\FI$-space} is a pointed functor
\[(\Fsk,\ordz) \fto{X} (\ITop,*).\]
A \emph{morphism} of $\FI$-spaces is a natural transformation.  The category of $\FI$-spaces and morphisms is denoted by $\FITop$.  For an $\FI$-category or $\FI$-space $X$, $n \geq 0$, and a finite set $A$, the category or space $(X\ordn)(A)$ is also denoted by $X(\ordn,A)$.  For a morphism $\tha$ of $\FI$-categories or $\FI$-spaces, $(\tha_{\ordn})_A$ is also denoted by $\tha_{\ordn,A}$.
\item\label{def:ficat_iii} 
The functor
\begin{equation}\label{fmcat_ficat}
\FMCat \fto{\uprst} \FICat
\end{equation}
is defined by postcomposition with the functor $\upr \cn \Mcat \to \ICat$ \cref{upr_functor}.  Thus, for an $\FM$-category $X \cn \Fsk \to \Mcat$ \pcref{def:fmcat}, 
\[(\uprst X)(\ordn,A) = \big(\upr(X\ordn)\big)(A).\]
The functor $\uprst$ is well defined because $X\ordz$ is the terminal $\schm$-category $\bone$ \cref{fmcat_axioms} and $\upr\bone$ is the terminal $\bdI$-category $\bone$.
\item\label{def:ficat_iv} 
The functor
\begin{equation}\label{ficat_fitop}
\FICat \fto{\clast} \FITop
\end{equation} 
is defined by componentwise postcomposition with the classifying space functor $\cla \cn \Cat \to \Top$ \cref{classifying_space}.  Thus, for an $\FI$-category $X$, 
\[(\clast X)(\ordn,A) = \cla X(\ordn,A),\]
the classifying space of $X(\ordn,A)$.\defmark
\end{enumerate}
\end{definition}

\subsection*{Symmetric Spectra}

Recall the category $\Topst$ of pointed spaces and pointed morphisms \pcref{def:gtopst}.

\begin{definition}\label{def:asphere}
For each finite set $A$, denote by $\bR[A]$ the real vector space of functions $A \to \bR$.  The \emph{$A$-sphere}\index{sphere} $S^A$\label{not:asphere} is the one-point compactification $\bR[A] \sqcup \{\infty\}$ with basepoint $\infty$.  The $\emptyset$-sphere \label{not:emptysp}$S^\emptyset = \{*,\infty\}$ is the two-point space, where $* \cn \emptyset \to \bR$ is the unique function.  Thus, for each pointed space $X$, there is a canonical pointed homeomorphism 
\begin{equation}\label{zero_sphere}
X \sma S^\emptyset \fiso X
\end{equation}
that sends $(x,*)$ to $x \in X$.
\end{definition}

Global equivariant algebraic $K$-theory uses symmetric spectra based on topological spaces and injections between finite sets \cite[1.1]{schwede_global}.  

\begin{definition}\label{def:sptop}\
\begin{enumerate}
\item\label{def:sptop_i} 
A \emph{symmetric spectrum}\index{symmetric spectrum} $X$ consists of 
\begin{itemize}
\item a pointed space $X_A$ for each finite set $A$ and
\item a pointed morphism\label{not:istarab}
\[X_A \sma S^{B \setminus i(A)} \fto{i_*} X_B\]
for each injection $i \cn A \to B$ between finite sets
\end{itemize}
such that the following two axioms hold.
\begin{description}
\item[Unity]
For each finite set $A$, the pointed morphisms
\begin{equation}\label{sptop_unity}
\begin{tikzpicture}[vcenter]
\draw[0cell]
(0,0) node (a1) {X_A \sma S^\emptyset}
(a1)++(2.5,0) node (a2) {\phantom{X_A}}
(a2)++(0,-.05) node (a2') {X_A}
;
\draw[1cell=.9]
(a1) edge[transform canvas={yshift=.4ex}] node {(1_A)_*} (a2)
(a1) edge[transform canvas={yshift=-.5ex}] node[swap] {\iso} (a2) 
;
\end{tikzpicture}
\end{equation}
are equal, where $\iso$ is the homeomorphism in \cref{zero_sphere}.
\item[Associativity]
For composable injections $i \cn A \to B$ and $j \cn B \to C$ between finite sets, the following diagram commutes.
\begin{equation}\label{sptop_assoc}
\begin{tikzpicture}[vcenter]
\def\v{-1.4}
\draw[0cell]
(0,0) node (a11) {X_A \sma S^{B \setminus i(A)} \sma S^{C \setminus j(B)}}
(a11)++(4.5,0) node (a12) {X_A \sma S^{C \setminus ji(A)}}
(a11)++(0,\v) node (a21) {X_B \sma S^{C \setminus j(B)}}
(a12)++(0,\v) node (a22) {\phantom{X_C}}
(a22)++(0,-.05) node (a22') {X_C}
;
\draw[1cell=.9]
(a11) edge node {1_{X_A} \sma \iso} (a12)
(a12) edge node {(ji)_*} (a22')
(a11) edge[transform canvas={xshift=1em}] node[swap] {i_* \sma 1_{S^{C \setminus j(B)}}} (a21)
(a21) edge node {j_*} (a22)
;
\end{tikzpicture}
\end{equation}
\end{description}
\item\label{def:sptop_ii}
A \emph{morphism} $f \cn X \to Y$ between symmetric spectra consists of a pointed morphism
\[X_A \fto{f_A} Y_A\]
for each finite set $A$ such that the following diagram commutes for each injection $i \cn A \to B$ between finite sets.
\begin{equation}\label{sptop_mor}
\begin{tikzpicture}[vcenter]
\def\v{-1.4}
\draw[0cell]
(0,0) node (a11) {X_A \sma S^{B \setminus i(A)}}
(a11)++(2.5,0) node (a12) {\phantom{X_B}}
(a12)++(0,-.05) node (a12') {X_B}
(a11)++(0,\v) node (a21) {Y_A \sma S^{B \setminus i(A)}}
(a12)++(0,\v) node (a22) {\phantom{Y_B}}
(a22)++(0,-.05) node (a22') {Y_B}
;
\draw[1cell=.9]
(a11) edge node {i_*} (a12)
(a12') edge node {f_B} (a22')
(a11) edge[transform canvas={xshift=1em}] node[swap] {f_A \sma 1_{S^{B \setminus i(A)}}} (a21)
(a21) edge node {i_*} (a22)
;
\end{tikzpicture}
\end{equation}
Identity morphisms and composition are defined levelwise at each finite set $A$.  Denote by $\Sptop$ the category of symmetric spectra and morphisms.  The category $\Sptop$ is well defined with small hom sets because the category $\bdI$ of finite sets and injections has a small skeleton \pcref{def:sch8.2}.
\item\label{def:sptop_iii}
For a finite group $G$, the category $\Sptopg$ of \emph{$G$-symmetric spectra}\index{G-symmetric spectrum@$G$-symmetric spectrum}\index{symmetric spectrum!G@$G$-}\index{G-spectrum@$G$-spectrum!symmetric} \cite[2.2]{hausmann17} has functors $G \to \Sptop$ as objects and natural transformations as morphisms.  In other words, a $G$-symmetric spectrum is a symmetric spectrum $X$ such that the following two statements hold.
\begin{enumerate}
\item\label{def:sptop_iiia} For each finite set $A$, the pointed space $X_A$ is equipped with the structure of a pointed $G$-space.
\item\label{def:sptop_iiib} For each injection $i \cn A \to B$ between finite sets, the structure morphism $i_*$ is a pointed $G$-morphism with $G$ acting trivially on the sphere $S^{B \setminus i(A)}$.
\end{enumerate}
A morphism of $G$-symmetric spectra $f \cn X \to Y$ is a morphism of underlying symmetric spectra such that $f_A \cn X_A \to Y_A$ is a pointed $G$-morphism for each finite set $A$.  
\item\label{def:sptop_iv}
For a finite group $G$ and a symmetric spectrum $X$, the \emph{underlying $G$-spectrum}\index{underlying G-spectrum@underlying $G$-spectrum} $X_G$ is the $G$-symmetric spectrum $X$ with $G$ acting trivially.  Fixing an arbitrary universal $G$-set $\ugsetu$ \pcref{def:univ_gset}, a \emph{$G$-stable equivalence}\index{G-stable equivalence@$G$-stable equivalence} between $G$-symmetric spectra is a $G^{\ugsetu}$-stable equivalence in the sense of \cite[2.35]{hausmann17}.  A morphism of symmetric spectra $f \cn X \to Y$ is a \emph{global equivalence}\index{global equivalence!symmetric spectra} \cite[2.9]{hausmann19} if it induces a $G$-stable equivalence  $f_G \cn X_G \to Y_G$ between underlying $G$-spectra for each finite group $G$.
\item\label{def:sptop_v}
Consider a finite group $G$, a subgroup $H \subseteq G$, and an integer $n$.  The \emph{$n$-th $H$-equivariant homotopy group}\index{equivariant homotopy group} of a $G$-symmetric spectrum $X$ is defined as the colimit\label{not:pinhu}
\[\pinhu(X) = \colimover{A \subset \ugsetu} \,[S^{\ufsn \sqcup A}, X_A]^H\]
with
\begin{itemize}
\item $A$ running through finite $G$-subsets of the given universal $G$-set $\ugsetu$,
\item $X_A$ denoting the pointed $G$-space in \cref{def:sptop_iiia},
\item $\ufsn$ denoting the unpointed finite set $\{1,2,\ldots,n\}$, and
\item $[S^{\ufsn \sqcup A}, X_A]^H$ denoting the set of $H$-homotopy classes of $H$-morphisms $S^{\ufsn \sqcup A} \to X_A$.
\end{itemize}  
See \cite[3.1]{hausmann17} for the precise meaning of $\pinhu(X)$ when $n$ is negative.  A morphism $f \cn X \to Y$ of $G$-symmetric spectra is a \emph{$\pistu$-isomorphism} if, for each integer $n$ and subgroup $H \subseteq G$, the induced morphism
\[\pinhu(X) \fto{\pinhu(f)} \pinhu(Y)\]
is an isomorphism \cite[3.3]{hausmann17}.\defmark
\end{enumerate}
\end{definition}

By \cite[3.36]{hausmann17}, each $\pistu$-isomorphism between $G$-symmetric spectra is a $G$-stable equivalence.

\begin{explanation}[Stable Model Structures]\label{expl:sptop}
The category $\Sptop$ of symmetric spectra in \cref{def:sptop} admits a stable model structure that is Quillen equivalent to the original one defined by Hovey, Shipley, and Smith \cite[3.4.4]{hss}, which is based on pointed simplicial sets and indexed on nonnegative integers.  See \cite[2.18]{hausmann19} for an explanation of this Quillen equivalence and \cite[Ch.\ 7]{cerberusIII} for an elementary discussion of symmetric spectra.  By \cite[2.17]{hausmann19}, the category of symmetric spectra admits a \emph{global stable model structure}\index{global stable model structure} with global equivalences (\cref{def:sptop} \cref{def:sptop_iv}) as weak equivalences.

By the work of Hausmann \cite[4.8]{hausmann17}, for each finite group $G$ and universal $G$-set $\ugsetu$, the category $\Sptopg$ of $G$-symmetric spectra (\cref{def:sptop} \cref{def:sptop_iii}) admits a stable model structure with $G$-stable equivalences as weak equivalences.  
\begin{itemize}
\item By \cite[7.4]{hausmann17}, $\Sptopg$ is Quillen equivalent to a slight variant of the Mandell-May stable model category of orthogonal $G$-spectra with $\pi_*$-isomorphisms as weak equivalences \cite[III.4.2]{mandell_may}.  The right Quillen equivalence is the forgetful functor 
\begin{equation}\label{gsp_forget}
\Gspec \fto{\Gspu} \Sptopg
\end{equation}
from orthogonal $G$-spectra to $G$-symmetric spectra.
\item By \cite[7.7]{hausmann17}, $\Sptopg$ is Quillen equivalent to Mandell's stable model category of equivariant symmetric spectra \cite[4.1]{mandell04}.  The right Quillen equivalence is the forgetful functor from $G$-symmetric spectra to Mandell's equivariant symmetric spectra.\defmark
\end{itemize}
\end{explanation}

\subsection*{$\FI$-Spaces to Symmetric Spectra}

\cref{def:kfi} defines the passage from $\FI$-spaces to symmetric spectra \pcref{def:ficat,def:sptop}.  Recall that for an $\FI$-space $X$ and a morphism $\tha$ of $\FI$-spaces, we also denote the pointed space $(X\ordn)(A)$ by $X(\ordn,A)$ and the pointed morphism $(\tha_{\ordn})_A$ by $\tha_{\ordn,A}$.

\begin{definition}\label{def:kfi}
The functor
\begin{equation}\label{kfi_functor}
\FITop \fto{\Kfi} \Sptop
\end{equation}
is defined as follows.
\begin{description}
\item[Objects]
For an $\FI$-space $X \cn \Fsk \to \ITop$, the symmetric spectrum $\Kfi X$ sends a finite set $A$ to the pointed space
\begin{equation}\label{kfixa}
(\Kfi X)_A = \begin{cases}
\int^{\ordn \in \Fsk} (S^A)^n \sma X(\ordn,A) & \text{if $A \neq \emptyset$ and}\\
X(\ordone, \emptyset) & \text{if $A = \emptyset$.}
\end{cases}
\end{equation}
\begin{description}
\item[Coend]
For $A \neq \emptyset$, $(S^A)^n$ is the $n$-fold Cartesian product of the $A$-sphere $S^A$ \pcref{def:asphere} with basepoint $\ang{\infty}_{j \in \ufsn}$, and $(S^A)^0 = *$ if $n=0$.  The coend in \cref{kfixa} is obtained from the wedge
\[\bigvee_{\ordn \in \Fsk} (S^A)^n \sma X(\ordn,A)\]
by identifying the points
\[\begin{cases}
\big(y; X(\psi,A)(x) \big) \in (S^A)^n \sma X(\ordn,A) & \text{and}\\
\big(\ang{y_{\psi(i)}}_{i \in \ufsm} \,; x \big) \in (S^A)^m \sma X(\ordm,A) &
\end{cases}\]
for each pointed morphism $\psi \cn \ordm \to \ordn$ in $\Fsk$, $x \in X(\ordm,A)$, and $y = \ang{y_j}_{j \in \ufsn} \in (S^A)^n$.  If $\psi(i) = 0$, then $y_{\psi(i)} \in S^A$ means the basepoint $\infty$.
\item[Structure morphisms]
For an injection $i \cn A \to B$ between finite sets, the pointed morphism 
\[(\Kfi X)_A \sma S^{B \setminus i(A)} \fto{i_*} (\Kfi X)_B\]
is defined as follows.
\begin{itemize}
\item If $A = \emptyset = B$, then 
\begin{equation}\label{kfix_emptysets}
X(\ordone,\emptyset) \sma S^{\emptyset} \fto{i_*} X(\ordone,\emptyset)
\end{equation}
is the canonical homeomorphism \cref{zero_sphere}.
\item If $A = \emptyset \neq B$, then $i_*$ is defined as the following composite.
\begin{equation}\label{kfix_istar_empty}
\begin{tikzpicture}[vcenter]
\def\v{-1.4} \def\h{3}
\draw[0cell=.9]
(0,0) node (a11) {X(\ordone,\emptyset) \sma S^B}
(a11)++(\h,0) node (a12) {(\Kfi X)_B}
(a11)++(0,\v) node (a21) {X(\ordone,B) \sma S^B}
(a12)++(0,\v) node (a22) {S^B \sma X(\ordone,B)}
;
\draw[1cell=.85]
(a11) edge node {i_*} (a12)
(a11) edge[transform canvas={xshift=1.7em}] node[swap] {X(\ordone,i) \sma 1_{S^B}} (a21)
(a21) edge node {\xi} node[swap] {\iso} (a22)
(a22) edge node[swap] {\eta_{\ordone}} (a12)
;
\end{tikzpicture}
\end{equation}
The morphism $\eta_{\ordone}$ is part of the coend $(\Kfi X)_B$ at the object $\ordone \in \Fsk$, and $\xi$ is the braiding for the smash product $\sma$.
\item If $A \neq \emptyset$, then $i_*$ is defined as the following composite.
\begin{equation}\label{kfix_istar}
\begin{tikzpicture}[vcenter]
\def\u{-1} \def\v{-1.4} \def\h{5}
\draw[0cell=.9]
(0,0) node (a11) {(\Kfi X)_A \sma S^{B \setminus i(A)}}
(a11)++(\h,0) node (a12) {(\Kfi X)_B}
(a11)++(0,\u) node (a21) {\big(\txint^{\ordn \in \Fsk} (S^A)^n \sma X(\ordn,A) \big) \sma S^{B \setminus i(A)}}
(a12)++(0,\u) node (a22) {\txint^{\ordn \in \Fsk} (S^B)^n \sma X(\ordn,B) }
(a21)++(0,\v) node (a31) {\txint^{\ordn \in \Fsk} (S^A)^n \sma S^{B \setminus i(A)} \sma X(\ordn,A)}
;
\draw[1cell=.85]
(a11) edge node {i_*} (a12)
(a11) edge[equal] (a21)
(a12) edge[equal] (a22)
(a21) edge node[swap] {\iso} (a31)
(a31) [rounded corners=2pt] -| node [pos=.25]{\txint^{\ordn} \assm_n \sma X(\ordn,i)} (a22)
;
\end{tikzpicture}
\end{equation}
The pointed homeomorphism $\iso$ is given by the commutation of $-\sma S^{B \setminus i(A)}$ with coends, along with the associativity and braiding for $\sma$.  Using the pointed homeomorphisms 
\[S^A \sma S^{B \setminus i(A)} \fto[\iso]{(- \comp i^{-1}) \sma 1} S^{i(A)} \sma S^{B \setminus i(A)} \iso S^B,\] 
the pointed morphism\label{not:assmn}
\[(S^A)^n \sma S^{B \setminus i(A)} \fto{\assm_n} (S^B)^n\]
is defined by
\begin{equation}\label{kfi_assm}
\assm_n\big(\ang{x_j}_{j \in \ufsn} \, ; y\big) 
= \ang{(x_j \comp i^{-1}, y)}_{j \in \ufsn} \in (S^B)^n
\end{equation}
for $\ang{x_j}_{j \in \ufsn}  \in (S^A)^n$ and $y \in S^{B \setminus i(A)}$.
\end{itemize}
\item[Axioms]
The unity axiom \cref{sptop_unity} for $\Kfi X$ follows from \cref{kfix_emptysets}, \cref{kfix_istar}, and the fact that $\assm_n$ is the canonical homeomorphism 
\begin{equation}\label{assm_onea}
(S^A)^n \sma S^\emptyset \iso (S^A)^n \ifspace i = 1_A.
\end{equation} 
The associativity axiom \cref{sptop_assoc} for $\Kfi X$ follows from \cref{kfix_istar_empty}, \cref{kfix_istar}, the universal property of coends, the functoriality of each $X\ordn \in \ITop$, and the following commutative diagram for composable injections $i \cn A \to B$ and $j \cn B \to C$ between finite sets.
\begin{equation}\label{assm_assoc}
\begin{tikzpicture}[vcenter]
\def\v{-1.4}
\draw[0cell]
(0,0) node (a11) {(S^A)^n \sma S^{B \setminus i(A)} \sma S^{C \setminus j(B)}}
(a11)++(4.7,0) node (a12) {(S^A)^n \sma S^{C \setminus ji(A)}}
(a11)++(0,\v) node (a21) {(S^B)^n \sma S^{C \setminus j(B)}}
(a12)++(0,\v) node (a22) {(S^C)^n}
;
\draw[1cell=.9]
(a11) edge node {1\, \sma \iso} (a12)
(a12) edge node {\assm_n} (a22)
(a11) edge node[swap] {\assm_n \sma 1} (a21)
(a21) edge node {\assm_n} (a22)
;
\end{tikzpicture}
\end{equation}
\end{description}
This finishes the definition of the symmetric spectrum $\Kfi X$.
\item[Morphisms]
For a morphism $\tha \cn X \to Y$ between $\FI$-spaces, the morphism of symmetric spectra
\[\Kfi X \fto{\Kfi\tha} \Kfi Y\]
consists of the pointed morphisms
\begin{equation}\label{Kfi_tha}
(\Kfi \tha)_A = \begin{cases}
\txint^{\ordn \in \Fsk} 1_{(S^A)^n} \sma \tha_{\ordn,A} & \text{if $A \neq \emptyset$ and}\\
\tha_{\ordone, \emptyset} & \text{if $A = \emptyset$.}
\end{cases}
\end{equation}
For $A \neq \emptyset$, the pointed morphism $(\Kfi \tha)_A$ is well defined by the naturality of $\tha$ in $\ordn \in \Fsk$.  The symmetric spectrum morphism axiom \cref{sptop_mor} holds for $\Kfi\tha$ by the universal property of coends and the naturality of each morphism $\tha_{\ordn} \cn X\ordn \to Y\ordn$ in $\ITop$.
\end{description}
The functoriality of $\Kfi$ follows from \cref{Kfi_tha} and the universal property of coends.
\end{definition}

\section{Global $K$-Theory and $G$-Stable Homotopy Type}
\label{sec:ksc}

This section reviews Schwede global $K$-theory \cite[4.14]{schwede_global}, which sends parsummable categories to symmetric spectra, and its underlying $G$-stable homotopy type.

\secoutline
\begin{itemize}
\item \cref{def:sch3.3} defines global $K$-theory $\Ksc$ as the composite functor $\Kfi\clast\uprst\Jsc$.  \cref{expl:ksc} unpacks the functor $\Ksc$.
\item \cref{sch4.16} records the fact that global $K$-theory sends global equivalences of parsummable categories to global equivalences of symmetric spectra.
\item \cref{def:sch4.15} defines a $G$-symmetric spectrum $\kcg$ that has the same pointed spaces as global $K$-theory $\Ksc\C$, simpler structure morphisms, and a nontrivial $G$-action.  
\cref{sch4.15} records the fact that $\kcg$ computes the underlying $G$-spectrum of the global $K$-theory of $\C$.
\item Using the functors in \cref{def:jscg,def:kfgsi,def:kscgb}, \cref{Isc_geq} provides a simpler description of the $G$-symmetric spectrum $\kcg$.
\end{itemize}

\subsection*{Parsummable Categories to Symmetric Spectra}

\begin{definition}\label{def:sch3.3}
Define \emph{global $K$-theory}\index{global K-theory@global $K$-theory}\index{K-theory@$K$-theory!global} $\Ksc$ as the composite 
\begin{equation}\label{ksc_def}
\begin{tikzpicture}[vcenter]
\def\h{2.5} \def\v{1.4}
\draw[0cell]
(0,0) node (a1) {\Parcat}
(a1)++(0,-\v) node (a2) {\FMCat}
(a2)++(\h,0) node (a3) {\FICat}
(a3)++(\h,0) node (a4) {\phantom{\FITop}}
(a4)++(0,-.035) node (a4') {\FITop}
(a4)++(0,\v) node (a5) {\Sptop}
;
\draw[1cell=.9]
(a1) edge node {\Ksc} (a5)
(a1) edge node[swap] {\Jsc} (a2) 
(a2) edge node {\uprst} (a3)
(a3) edge node {\clast} (a4)
(a4') edge node[swap] {\Kfi} (a5)
;
\end{tikzpicture}
\end{equation}
of the functors $\Jsc$ \cref{Jsc}, $\uprst$ \cref{fmcat_ficat}, $\clast$ \cref{ficat_fitop}, and $\Kfi$ \cref{kfi_functor}.
\end{definition}

\begin{explanation}[Unpacking]\label{expl:ksc}
For a parsummable category $(\C,\psum,\pzero)$ \pcref{def:parcat}, the symmetric spectrum $\Ksc\C$ is given explicitly as follows.  
\begin{description}
\item[Pointed spaces]
By \cref{JscCn,uprca,kfixa}, $\Ksc\C$ sends a finite set $A$ to the pointed space
\begin{equation}\label{KscCa}
(\Ksc\C)_A = \begin{cases}
\int^{\ordn \in \Fsk} (S^A)^n \sma \cla\Cboxnomea & \text{if $A \neq \emptyset$ and}\\
\cla\C^{\suppempty} & \text{if $A = \emptyset$.}
\end{cases}
\end{equation}
\begin{itemize}
\item $\C^{\suppempty} \subseteq \C$ is the full subcategory of $\schm$-fixed objects, and $\cla \cn \Cat \to \Top$ is the classifying space functor.
\item For $A \neq \emptyset$, the underlying category of $\Cboxnomea$ is the full subcategory $\Cboxn \subseteq \C^n$ of disjointly supported $n$-tuples.  In defining the coend in \cref{KscCa}, for a pointed morphism $\psi \cn \ordm \to \ordn$ in $\Fsk$, the functor 
\[\Cboxmomea \fto{\psi_*} \Cboxnomea\]
sends an $m$-tuple $\ang{x_i}_{i \in \ufsm}$ to the $n$-tuple $\ang{\,\sum_{i \in \psiinv j} x_i}_{j \in \ufsn}$ \cref{JscCpsi}.
\end{itemize}
\item[Structure morphisms]
For an injection $i \cn A \to B$ between finite sets, the structure pointed morphism 
\[(\Ksc\C)_A \sma S^{B \setminus i(A)} \fto{i_*} (\Ksc\C)_B\]
is given as follows.
\begin{description}
\item[Empty domain and codomain]
If $A = \emptyset = B$, then, by \cref{JscCn,uprca,kfix_emptysets},
\begin{equation}\label{ksc_empty_empty}
\cla\C^{\suppempty} \sma S^\emptyset \fto{i_*} \cla\C^{\suppempty}
\end{equation}
is the canonical pointed homeomorphism \cref{zero_sphere}.
\item[Empty domain]
If $A = \emptyset \neq B$, then, by \cref{kfix_istar_empty,uprci_emptyb}, $i_*$ is the following composite.
\begin{equation}\label{ksc_empty_dom}
\begin{tikzpicture}[vcenter]
\def\v{-1.4} \def\h{6.5}
\draw[0cell]
(0,0) node (a11) {(\Ksc\C)_\emptyset \sma S^B}
(a11)++(\h,0) node (a12) {(\Ksc\C)_B}
(a11)++(0,\v) node (a21) {\cla\C^{\suppempty} \sma S^B}
(a21)++(3.5,0) node (a22) {\cla\Comeb \sma S^B}
(a12)++(0,\v) node (a23) {S^B \sma \cla\Comeb}
;
\draw[1cell=.9]
(a11) edge node {i_*} (a12)
(a11) edge[equal] (a21)
(a21) edge node {\cla\iota \sma 1} (a22)
(a22) edge node {\xi} node[swap] {\iso} (a23)
(a23) edge node[swap] {\eta_{\ordone}} (a12)
;
\end{tikzpicture}
\end{equation}
The functor $\iota \cn \C^{\suppempty} \to \Comeb = \C$ is the full subcategory inclusion, and $\xi$ is the braiding for $\sma$.  The pointed morphism $\eta_{\ordone}$ is part of the coend $(\Ksc\C)_B$ at the object $\ordone \in \Fsk$.
\item[Nonempty domain]
If $A \neq \emptyset$, then, by \cref{uprci_nonempty,kfix_istar}, $i_*$ is the following composite.
\begin{equation}\label{ksc_nonempty_dom}
\begin{tikzpicture}[vcenter]
\def\u{-1} \def\v{-1.4} \def\h{4.6} \def\d{.8}
\draw[0cell=.9]
(0,0) node (a11) {(\Ksc\C)_A \sma S^{B \setminus i(A)}}
(a11)++(\h+\d,0) node (a12) {(\Ksc\C)_B}
(a11)++(0,\u) node (a21) {\big(\txint^{\ordn \in \Fsk} (S^A)^n \sma \cla\Cboxnomea\big) \sma S^{B \setminus i(A)}}
(a21)++(\h,0) node (a22) {\txint^{\ordn \in \Fsk} (S^B)^n \sma \cla\Cboxnomeb}
(a22)++(\d,0) node (a22') {\phantom{\txint^{\ordn \in \Fsk}}}
(a21)++(0,\v) node (a31) {\txint^{\ordn \in \Fsk} (S^A)^n \sma S^{B \setminus i(A)} \sma \cla\Cboxnomea}
;
\draw[1cell=.85]
(a11) edge node {i_*} (a12)
(a11) edge[equal] (a21)
(a12) edge[equal] (a22')
(a21) edge node[swap] {\iso} (a31)
(a31) [rounded corners=2pt] -| node [pos=.25]{\txint^{\ordn} \assm_n \sma \cla(i_!^\ome)_\mstar^{\Bxtimes n}} (a22')
;
\end{tikzpicture}
\end{equation}
The functor
\[\Cboxnomea = \Cboxn \fto{(i_!^\ome)_\mstar^{\Bxtimes n}} \Cboxn = \Cboxnomeb\]
is given entrywise by the $\schm$-action functor $(i_!^\ome)_\mstar \cn \C \to \C$ for the injection $i_!^\ome \cn \ome \to \ome$ \cref{iome_ext}.  The pointed morphism 
\[(S^A)^n \sma S^{B \setminus i(A)} \fto{\assm_n} (S^B)^n\]
is defined in \cref{kfi_assm}.  This finishes the description of the symmetric spectrum $\Ksc\C$.
\end{description}
\item[$\Ksc$ on morphisms]
Suppose $\fun \cn \C \to \D$ is a parsummable functor.  By \cref{Jscf,uprf,Kfi_tha}, the morphism of symmetric spectra
\[\Ksc\C \fto{\Ksc\fun} \Ksc\D\]
has, for each finite set $A$, $A$-component pointed morphism
\[(\Ksc\fun)_A = \begin{cases}
\int^{\ordn \in \Fsk} 1_{(S^A)^n} \sma \cla\fboxn & \text{if $A \neq \emptyset$ and}\\
\cla\fun|_{\C^{\suppempty}} & \text{if $A = \emptyset$.}
\end{cases}\]
\begin{itemize}
\item For $A \neq \emptyset$, the functor 
\[\Cboxnomea \fto{\fboxn} \Dboxnomea\] 
is given entrywise by $\fun$. 
\item For $A = \emptyset$, the functor 
\[\C^{\suppempty} \fto{\fun|_{\C^{\suppempty}}} \D^{\suppempty}\] 
is the restriction of $\fun$ to the full subcategory $\C^{\suppempty} \subseteq \C$ of $\schm$-fixed objects.  
\end{itemize}
\end{description}
This finishes the description of the global $K$-theory functor $\Ksc$.
\end{explanation}

\begin{remark}[Notation]\label{rk:sch_notation}
The subscripts in $\Jsc$ and $\Ksc$ refer to Schwede.  In \cite[4.3]{schwede_global}, global $J$-theory $\Jsc$ is denoted by $\ga(-)$.  In \cite[3.3]{schwede_global}, the composite functor 
\[\FMCat \fto{\uprst} \FICat \fto{\clast} \FITop \fto{\Kfi} \Sptop\]
is presented as a single construction and denoted by $\Kfm$.  The classifying space functor is denoted there by $\Rea$, with the nerve suppressed from the notation.  In \cite[4.14]{schwede_global}, global $K$-theory $\Ksc$ is denoted by $\Schk$.
\end{remark}

The following result from \cite[4.16]{schwede_global} shows that global $K$-theory preserves global equivalences, as defined in \cref{def:sch2.26,def:sptop}.

\begin{theorem}[Schwede]\label{sch4.16}
Suppose $\fun \cn \C \to \D$ is a parsummable functor whose underlying $\schm$-functor is a global equivalence.  Then the morphism
\[\Ksc\C \fto{\Ksc\fun} \Ksc\D\]
is a global equivalence of symmetric spectra.
\end{theorem}

\subsection*{Underlying $G$-Spectrum of Global $K$-Theory}

For each finite group $G$, the underlying $G$-spectrum of the global $K$-theory of a parsummable category \pcref{def:sch3.3} can be computed more directly as follows.  For \cref{def:sch4.15}, recall the following.
\begin{itemize}
\item A \emph{$G$-symmetric spectrum} is a symmetric spectrum equipped with a $G$-action by automorphisms (\cref{def:sptop} \cref{def:sptop_iii}).  The \emph{underlying $G$-spectrum} of a symmetric spectrum is equipped with the trivial $G$-action (\cref{def:sptop} \cref{def:sptop_iv}). 
\item The universal $G$-set $\omeg$ consists of functions $G \to \ome$, with $G$-action $gf = f(\ginv \cdot -)$ \pcref{def:omea}.
\end{itemize}

\begin{definition}\label{def:sch4.15}
For a parsummable category $(\C,\psum,\pzero)$ and a finite group $G$, the $G$-symmetric spectrum $\kcg$ is defined as follows.
\begin{description}
\item[Pointed spaces]
$\kcg$ sends a finite set $A$ to the pointed space
\begin{equation}\label{kcga}
(\kcg)_A = \begin{cases}
\int^{\ordn \in \Fsk} (S^A)^n \sma \cla\Cboxnomeg & \text{if $A \neq \emptyset$ and}\\
\cla\C^{\suppempty} & \text{if $A = \emptyset$.}
\end{cases}
\end{equation}
\begin{itemize}
\item For $A=\emptyset$, $\cla\C^{\suppempty}$ is the classifying space of the full subcategory $\C^{\suppempty} \subseteq \C$ of $\schm$-fixed objects.
\item For $A \neq \emptyset$, $\Cboxnomeg$ is the $\omeg$-reparametrization of the $\schm$-category $\Cboxn$ \pcref{def:sch2.21,def:bxtimes}.  The underlying category of $\Cboxnomeg$ is the full subcategory $\Cboxn \subseteq \C^n$ of disjointly supported $n$-tuples.  In defining the coend in \cref{kcga}, for a pointed morphism $\psi \cn \ordm \to \ordn$ in $\Fsk$, the functor 
\begin{equation}\label{cmg_cng}
\Cboxmomeg \fto{\psi_*} \Cboxnomeg
\end{equation}
sends an $m$-tuple $\ang{x_i}_{i \in \ufsm}$ to the $n$-tuple $\ang{\,\sum_{i \in \psiinv j} x_i}_{j \in \ufsn}$ \cref{JscCpsi}.  
\end{itemize}
\item[Structure morphisms]
For an injection $i \cn A \to B$ between finite sets, the structure pointed morphism 
\[(\kcg)_A \sma S^{B \setminus i(A)} \fto{i_*} (\kcg)_B\]
is defined as follows.
\begin{description}
\item[Empty domain and codomain]
If $A = \emptyset = B$, then 
\begin{equation}\label{kcg_empty_empty}
\cla\C^{\suppempty} \sma S^\emptyset \fto{i_*} \cla\C^{\suppempty}
\end{equation}
is the canonical pointed homeomorphism \cref{zero_sphere}. 
\item[Empty domain]
If $A = \emptyset \neq B$, then $i_*$ is the following composite.
\begin{equation}\label{kcg_empty_dom}
\begin{tikzpicture}[vcenter]
\def\v{-1.4} \def\h{6.6}
\draw[0cell]
(0,0) node (a11) {(\kcg)_\emptyset \sma S^B}
(a11)++(\h,0) node (a12) {(\kcg)_B}
(a11)++(0,\v) node (a21) {\cla\C^{\suppempty} \sma S^B}
(a21)++(3.5,0) node (a22) {\cla\Comeg \sma S^B}
(a12)++(0,\v) node (a23) {S^B \sma \cla\Comeg}
;
\draw[1cell=.9]
(a11) edge node {i_*} (a12)
(a11) edge[equal] (a21)
(a21) edge node {\cla\iota \sma 1} (a22)
(a22) edge node {\xi} node[swap] {\iso} (a23)
(a23) edge node[swap] {\eta_{\ordone}} (a12)
;
\end{tikzpicture}
\end{equation}
The functor $\iota \cn \C^{\suppempty} \to \Comeg = \C$ is the full subcategory inclusion, and $\xi$ is the braiding for $\sma$.  The pointed morphism $\eta_{\ordone}$ is part of the coend $(\kcg)_B$ at the object $\ordone \in \Fsk$. 
\item[Nonempty domain]
If $A \neq \emptyset$, then $i_*$ is the following composite.
\begin{equation}\label{kcg_nonempty_dom}
\begin{tikzpicture}[vcenter]
\def\u{-1} \def\v{-1.4} \def\h{4.6} \def\d{.8}
\draw[0cell=.9]
(0,0) node (a11) {(\kcg)_A \sma S^{B \setminus i(A)}}
(a11)++(\h+\d,0) node (a12) {(\kcg)_B}
(a11)++(0,\u) node (a21) {\big(\txint^{\ordn \in \Fsk} (S^A)^n \sma \cla\Cboxnomeg \big) \sma S^{B \setminus i(A)}}
(a21)++(\h,0) node (a22) {\txint^{\ordn \in \Fsk} (S^B)^n \sma \cla\Cboxnomeg}
(a22)++(\d,0) node (a22') {\phantom{\txint^{\ordn \in \Fsk}}}
(a21)++(0,\v) node (a31) {\txint^{\ordn \in \Fsk} (S^A)^n \sma S^{B \setminus i(A)} \sma \cla\Cboxnomeg}
;
\draw[1cell=.85]
(a11) edge node {i_*} (a12)
(a11) edge[equal] (a21)
(a12) edge[equal] (a22')
(a21) edge node[swap] {\iso} (a31)
(a31) [rounded corners=2pt] -| node [pos=.25]{\txint^{\ordn} \assm_n \sma 1} (a22')
;
\end{tikzpicture}
\end{equation}
The pointed morphism
\[(S^A)^n \sma S^{B \setminus i(A)} \fto{\assm_n} (S^B)^n\]
is defined in \cref{kfi_assm}, and 1 denotes the identity morphism of the pointed space $\cla\Cboxnomeg$.
\end{description}
\item[Axioms]
The unity axiom \cref{sptop_unity} for $\kcg$ follows from \cref{kcg_empty_empty,kcg_nonempty_dom,assm_onea}.  The associativity axiom \cref{sptop_assoc} for $\kcg$ follows from the universal property of coends, \cref{assm_assoc,kcg_empty_dom,kcg_nonempty_dom}.  Thus, $\kcg$ is a symmetric spectrum.
\item[$G$-action]
For each finite set $A$, the pointed space $(\kcg)_A$ \cref{kcga} is equipped with the following $G$-action.
\begin{description}
\item[Empty set]
For $A = \emptyset$, $G$ acts trivially on $(\kcg)_\emptyset = \cla\C^{\suppempty}$.
\item[Nonempty sets]
For $A \neq \emptyset$ and $n \geq 0$, the $\schm$-category $\Cboxn$ yields a pointed $G$-category structure on $\Cboxnomeg$, as defined in \cref{Cbrac_g}.  In more detail, for each $g \in G$, the $g$-action functor on $\Cboxnomeg$ is given entrywise by the $\schm$-action functor $g^\ome_\mstar \cn \C \to \C$ for the bijection $g^\ome \cn \ome \to \ome$ \cref{rep_g_rep}:
\begin{equation}\label{kcg_gaction}
\Cboxnomeg = \Cboxn \fto[\iso]{\ang{g^\ome_\mstar}_{j \in \ufsn}} \Cboxn = \Cboxnomeg.
\end{equation}
For each pointed morphism $\psi \cn \ordm \to \ordn$ in $\Fsk$, the functor \cref{cmg_cng}
\[\Cboxmomeg \fto{\psi_*} \Cboxnomeg\]
is $G$-equivariant because $g^\ome_\mstar$ is an $\schm$-action functor and the sum is an $\schm$-functor.  Applying the classifying space functor $\cla$ gives $\cla\Cboxnomeg$ the structure of a pointed $G$-space.  Letting $G$ act trivially on $(S^A)^n$, the universal property of coends implies that
\[(\kcg)_A  = \txint^{\ordn \in \Fsk} (S^A)^n \sma \cla\Cboxnomeg\]
is a pointed $G$-space.
\item[Compatibility with structure morphisms]
The $G$-action on $\kcg$ commutes with the structure morphisms of $\kcg$ in the sense of \cref{sptop_mor}.  
\begin{itemize}
\item For an injection $i \cn \emptyset \to B$ with $B$ a finite set, the $G$-action commutes with $i_*$, defined in \cref{kcg_empty_empty,kcg_empty_dom}, because $\C^{\suppempty}$ consists of $\schm$-fixed objects in $\C$ and the functor $g^\ome_\mstar$ is part of the $\schm$-action on $\C$.
\item For an injection $i \cn A \to B$ between finite sets with $A \neq \emptyset$, the $G$-action commutes with $i_*$ \cref{kcg_nonempty_dom} because $G$ acts trivially on spheres.
\end{itemize}
\end{description}
\end{description}
This finishes the construction of the $G$-symmetric spectrum $\kcg$.  

Moreover, for a parsummable functor $\fun \cn \C \to \D$, the morphism of $G$-symmetric spectra
\[\kcg \fto{\Kscg\fun} \kdg\]
is defined by the component pointed $G$-morphisms
\[(\Kscg\fun)_A = \begin{cases}
\txint^{\ordn \in \Fsk} 1_{(S^A)^n} \sma \cla \fboxn & \text{if $A\neq\emptyset$ and}\\
\fun|_{\cla\C^{\suppempty}} \cn \cla\C^{\suppempty} \to \cla\D^{\suppempty} & \text{if $A = \emptyset$.}
\end{cases}\]
This defines a functor $\Kscg \cn \Parcat \to \Sptopg$.
\end{definition}

\begin{explanation}[$\Ksc$ versus $\Kscg$]\label{expl:kcg_ksc}
For a parsummable category $(\C,\psum,\pzero)$ and a finite group $G$, we compare the definitions of $\Ksc\C$ and $\kcg$.
\begin{description}
\item[Similarities]
For each finite set $A$, by \cref{KscCa,kcga}, there is an equality of pointed spaces 
\[(\Ksc\C)_A = (\kcg)_A\]
because the underlying categories of $\Cboxnomea$ and $\Cboxnomeg$ are both given by the full subcategory $\Cboxn \subseteq \C^n$ of disjointly supported $n$-tuples.  Moreover, for an injection $i \cn \emptyset \to B$ with $B$ a finite set, the structure morphism $i_*$ for $\Ksc\C$, defined in \cref{ksc_empty_empty,ksc_empty_dom}, is equal to $i_*$ for $\kcg$, defined in \cref{kcg_empty_empty,kcg_empty_dom}.  
\item[Differences]
For an injection $i \cn A \to B$ between finite sets with $A \neq \emptyset$, the structure morphism $i_*$ for $\Ksc\C$ \cref{ksc_nonempty_dom} involves the pointed morphism 
\[\cla\Cboxnomea \fto{\cla(i_!^\ome)_\mstar^{\Bxtimes n}} \cla\Cboxnomeb.\] 
On the other hand, $i_*$ for $\kcg$ \cref{kcg_nonempty_dom} involves the identity morphism on $\cla\Cboxnomeg$.  Furthermore, $\kcg$ has a nontrivial $G$-action \cref{kcg_gaction}, while the underlying $G$-spectrum of $\Ksc\C$ has the trivial $G$-action.\defmark
\end{description}
\end{explanation}

Recall $G$-stable equivalences and $\pistu$-isomorphisms between $G$-symmetric spectra \pcref{def:sptop}.  By \cite[3.36]{hausmann17}, each $\pistu$-isomorphism is also a $G$-stable equivalence.  The following result from \cite[3.14 and 4.15 (ii)]{schwede_global} shows that the underlying $G$-spectrum $(\Ksc\C)_G$ of global $K$-theory $\Ksc\C$ \pcref{def:sch3.3} is naturally $\pistu$-isomorphic to the $G$-symmetric spectrum $\kcg$ \pcref{def:sch4.15}.

\begin{theorem}[Schwede]\label{sch4.15}\index{global K-theory@global $K$-theory!underlying G-spectrum@underlying $G$-spectrum}\index{underlying G-spectrum@underlying $G$-spectrum!global K-theory@global $K$-theory}
For each parsummable category $\C$ and finite group $G$, there are natural $\pistu$-isomorphisms 
\[(\Ksc\C)_G \fto{\ascgc} \gacomegs \fot{\bscgc} \kcg.\]
Thus, the $G$-symmetric spectra $(\Ksc\C)_G$ and $\kcg$ are naturally $G$-stably equivalent.
\end{theorem}

\begin{explanation}\label{expl:sch3.13}
In \cite[3.13 and 4.15]{schwede_global}, $\kcg$ is denoted by $|\ga(\C)[\omeg]|(\mathbbm{S})$.  Analogous to \cref{KscCa}, the $G$-symmetric spectrum $\gacomegs$ sends a finite set $A$ to the pointed space
\begin{equation}\label{gacomeg_a}
\gacomegs_A = \begin{cases}
\int^{\ordn \in \Fsk} (S^A)^n \sma \cla\Cboxnomega & \text{if $A \neq \emptyset$ and}\\
\cla\C^{\suppempty} & \text{if $A = \emptyset$.}
\end{cases}
\end{equation}
The group $G$ acts trivially on $\cla\C^{\suppempty}$.  For $A \neq \emptyset$ and $g \in G$, the $g$-action on $\gacomegs_A$ is given by the $\schm$-action functor on $\Cboxn$ for the bijection
\[\ome \fto{\rep_{\omeg \sqcup \omea}} \omeg \sqcup \omea \fto{g \sqcup 1_{\omea}} 
\omeg \sqcup \omea \fto{\rep_{\omeg \sqcup \omea}^{-1}} \ome\]
involving
\begin{itemize}
\item the $g$-action $g \cn \omeg \to \omeg$ \cref{g_omea} and
\item the bijection $\rep_{\omeg \sqcup \omea}$ \cref{rep_S}. 
\end{itemize} 
The structure morphisms of $\gacomegs$ are similar to those for $\Ksc\C$, as explained in \crefrange{ksc_empty_empty}{ksc_nonempty_dom}, with the following modification.  In \cref{ksc_nonempty_dom}, instead of $(i_!^\ome)_\mstar^{\Bxtimes n}$, here we use the $\schm$-action functor on $\Cboxn$ for the injection
\[\ome \fto{\rep_{\omeg \sqcup \omea}} \omeg \sqcup \omea \fto{1_{\omeg} \sqcup i_!} 
\omeg \sqcup \omeb \fto{\rep_{\omeg \sqcup \omeb}^{-1}} \ome.\]
By \cref{KscCa,kcga,gacomeg_a}, all three $G$-symmetric spectra in \cref{sch4.15} have the same pointed space at each finite set $A$, including the empty set, but different $G$-actions for nonempty finite sets $A$. 
\begin{itemize}
\item At the empty set, $\ascgc$ and $\bscgc$ are given by the identity morphism on the pointed space $\cla\C^{\suppempty}$. 
\item At each nonempty finite set $A$, $\ascgc$ and $\bscgc$ are induced by the $\schm$-action functors on $\Cboxn$ for, respectively, the following inclusions.
\begin{equation}\label{ionetwoome}
\begin{tikzpicture}[vcenter]
\def\h{1.8} \def\g{2} \def\f{2.75} \def\u{.7} \def\d{.05} \def\s{.85}
\draw[0cell]
(0,0) node (a1) {\ome}
(a1)++(\h,0) node (a2) {\phantom{\omeg}}
(a2)++(0,\d) node (a2') {\omea}
(a2)++(\g,0) node (a3) {\phantom{\omeg\sqcup \omea}}
(a3)++(0,\d) node (a3') {\omeg \sqcup \omea}
(a3)++(\f,0) node (a4) {\ome}
;
\draw[1cell=\s]
(a1) edge node {\rep_{\omea}} (a2)
(a2) edge node {i_2} (a3)
(a3) edge node {\rep^{-1}_{\omeg \sqcup \omea}} (a4)
(a1) [rounded corners=2pt] |- ($(a2)+(0,\u)$) -- node[pos=.7] {i_2^\ome} ($(a3)+(0,\u)$) -| (a4)
;
\begin{scope}[shift={(0,-1.6)}]
\draw[0cell]
(0,0) node (a1) {\ome}
(a1)++(\h,0) node (a2) {\phantom{\omeg}}
(a2)++(0,\d) node (a2') {\omeg}
(a2)++(\g,0) node (a3) {\phantom{\omeg\sqcup \omea}}
(a3)++(0,\d) node (a3') {\omeg \sqcup \omea}
(a3)++(\f,0) node (a4) {\ome}
;
\draw[1cell=\s]
(a1) edge node {\rep_{\omeg}} (a2)
(a2) edge node {i_1} (a3)
(a3) edge node {\rep^{-1}_{\omeg \sqcup \omea}} (a4)
(a1) [rounded corners=2pt] |- ($(a2)+(0,\u)$) -- node[pos=.7] {i_1^\ome} ($(a3)+(0,\u)$) -| (a4)
;
\end{scope}
\end{tikzpicture}
\end{equation}
\end{itemize}
Here, $i_1$ and $i_2$ denote the first and the second summand inclusions.
\end{explanation}

\subsection*{Another Description of $\kcg$}

There is a simpler description of the $G$-symmetric spectrum $\kcg$ \pcref{def:sch4.15} using the following constructions.  Recall
\begin{itemize}
\item the category $\Gcatst$ of small pointed $G$-categories and pointed $G$-functors \pcref{def:GCat} and
\item the category $\FGCat$ of $\Fskg$-categories \pcref{def:fgcat}. 
\end{itemize} 
\cref{def:jscg} is the $G$-equivariant analogue of global $J$-theory $\Jsc$ \cref{Jsc}.

\begin{definition}\label{def:jscg}
For a finite group $G$, \emph{Schwede $J$-theory}\index{Schwede J-theory@Schwede $J$-theory}\index{J-theory@$J$-theory!Schwede} for $G$ is the functor
\begin{equation}\label{Jscg}
\Parcat \fto{\Jscg} \FGCat
\end{equation}
defined as follows.
\begin{description}
\item[Objects] 
For a parsummable category $(\C,\psum,\pzero)$, the pointed functor
\begin{equation}\label{JscgC}
(\Fsk,\ordz) \fto{\Jscg\C} (\Gcatst,\bone)
\end{equation}
sends each object $\ordn \in \Fsk$ to the pointed $G$-category \cref{Cbrac_g}
\begin{equation}\label{JscgCn}
(\Jscg\C)\ordn = \Cboxnomeg.
\end{equation}
Its underlying pointed category is the full $\schm$-subcategory $\Cboxn \subseteq \C^n$ of disjointly supported $n$-tuples.  Its basepoint is the $n$-tuple $\ang{\pzero}_{j \in \ufsn}$.  For each $g \in G$, the $g$-action functor on $\Cboxnomeg$ is given entrywise by the $\schm$-action functor $g^\ome_\mstar \cn \C \to \C$ for the bijection $g^\ome \cn \ome \to \ome$ \cref{rep_g_rep}:
\begin{equation}\label{jscg_gaction}
\Cboxnomeg = \Cboxn \fto[\iso]{\ang{g^\ome_\mstar}_{j \in \ufsn}} \Cboxn = \Cboxnomeg.
\end{equation}
This functor is well defined by \cref{sch2.13} \cref{sch2.13_iii}.

For each pointed morphism $\psi \cn \ordm \to \ordn$ in $\Fsk$, the pointed $G$-functor 
\begin{equation}\label{jscg_psist}
(\Jscg\C)\ordm = \Cboxmomeg \fto{(\Jscg\C)\psi \,=\, \psi_*} (\Jscg\C)\ordn = \Cboxnomeg 
\end{equation}
sends an $m$-tuple $\ang{x_i}_{i \in \ufsm} \in \Cboxmomeg$ to the $n$-tuple
\[\psi_* \ang{x_i}_{i \in \ufsm} = \bang{\,\txsum_{i \in \psiinv j} x_i}_{j \in \ufsn} \in \Cboxnomeg.\]  
The preceding $n$-tuple is disjointly supported by subadditivity \cref{supp_xsumy}.  When $\psiinv j = \emptyset$, an empty sum means the unit $\pzero \in \C$ or its identity morphism.  The $G$-equivariance of $\psi_*$ follows from the $\schm$-equivariance of the sum on $\C$ and the fact that $g^\ome_\mstar$ is an $\schm$-action functor on $\C$.  The functoriality of $\Jscg\C$ follows from the parsummable category axioms \cref{parcat_axioms} for $\C$.
\item[Morphisms] 
For a parsummable functor $\fun \cn \C \to \D$ between parsummable categories, the natural transformation
\begin{equation}\label{Jscgf}
\begin{tikzpicture}[vcenter]
\def\t{30}
\draw[0cell]
(0,0) node (a1) {\Fsk}
(a1)++(2,0) node (a2) {\phantom{\Fsk}}
(a2)++(.3,-.03) node (a2') {\Gcatst}
;
\draw[1cell=.85]
(a1) edge[bend left=\t] node {\Jscg\C} (a2)
(a1) edge[bend right=\t] node[swap] {\Jscg\D} (a2)
;
\draw[2cell]
node[between=a1 and a2 at .35, rotate=-90, 2label={above,\Jscg\fun}] {\Rightarrow}
;
\end{tikzpicture}
\end{equation}
sends each object $\ordn \in \Fsk$ to the pointed $G$-functor
\[(\Jscg\C)\ordn = \Cboxnomeg \fto{(\Jscg\fun)_{\ordn} \,=\, \fun^{\Bxtimes n}} 
(\Jscg\D)\ordn = \Dboxnomeg.\]
The preceding $G$-functor is well defined by \cref{sch2.13} \cref{sch2.13_v}.  The naturality of $\Jscg\fun$ follows from the parsummable functor axioms \cref{parfun_axioms} for $\fun$.  The functoriality of $\Jscg$ follows from the fact that $\fun^{\Bxtimes n}$ is given entrywise by $\fun$.
\end{description}
This finishes the definition of the functor $\Jscg$.
\end{definition}

Recall the category $\Gtopst$ of pointed $G$-spaces and pointed $G$-morphisms \pcref{def:gtopst}.  \cref{def:fgtop} is the topological analogue of $\FGCat$.

\begin{definition}\label{def:fgtop}
For a group $G$, denote by $\FGTop$ the category with
\begin{itemize}
\item pointed functors $(\Fsk,\ord{0}) \to (\Gtopst,*)$, called \emph{$\Fskg$-spaces}, as objects and
\item natural transformations between such functors as morphisms.
\end{itemize}
Identity morphisms and composition are defined componentwise.
\end{definition}

\cref{expl:fgcat_iicat}, with $\Gcatst$ replaced by $\Gtopst$, also applies to $\FGTop$.  \cref{def:kfgsi} defines the functor from $\Fskg$-spaces to $G$-symmetric spectra (\cref{def:sptop} \cref{def:sptop_iii}).

\begin{definition}\label{def:kfgsi}
For a finite group $G$, define the functor
\begin{equation}\label{kfgsi_functor}
\FGTop \fto{\Kfgsi} \Sptopg
\end{equation}
as follows.  
\begin{description}
\item[Objects]
For an $\Fskg$-space $X \cn \Fsk \to \Gtopst$, the $G$-symmetric spectrum $\Kfgsi X$ is defined as follows.
\begin{description}
\item[Pointed $G$-spaces]
$\Kfgsi X$ sends each finite set $A$ to the coend
\begin{equation}\label{Kfgsixa}
(\Kfgsi X)_A = \int^{\ordn\in \Fsk} (S^A)^n \sma X\ordn
\end{equation}
taken in the category $\Gtopst$.  The group $G$ acts trivially on the $A$-sphere $S^A$ \pcref{def:asphere}. 
\item[Structure $G$-morphisms]
For each injection $i \cn A \to B$ between finite sets, the structure morphism $i_*$ for $\Kfgsi X$ is the following composite of pointed $G$-morphisms.
\begin{equation}\label{Kfgsix_istar}

\end{equation}
The pointed morphism
\[(S^A)^n \sma S^{B \setminus i(A)} \fto{\assm_n} (S^B)^n\]
is defined in \cref{kfi_assm}.  The unity axiom \cref{sptop_unity} and the associativity axiom \cref{sptop_assoc} for a $G$-symmetric spectrum follow from the corresponding properties for $\assm_n$ in \cref{assm_onea,assm_assoc}.  This finishes the definition of the $G$-symmetric spectrum $\Kfgsi X$.
\end{description}
\item[Morphisms]
For a natural transformation 
\begin{equation}\label{thaXY}
%
\end{equation}
between $\Fskg$-spaces $X$ and $Y$, the morphism between $G$-symmetric spectra
\begin{equation}\label{Kfgsi_tha}
\Kfgsi X \fto{\Kfgsi\tha} \Kfgsi Y
\end{equation}
sends each finite set $A$ to the following pointed $G$-morphism.
\begin{equation}\label{Kfgsi_tha_a}
%
\end{equation}
The $G$-equivariance of $\Kfgsi\tha$ follows from the $G$-equivariance of each component $\tha_{\ordn}$.  The compatibility \cref{sptop_mor} of $\Kfgsi\tha$ with the structure $G$-morphisms of $\Kfgsi X$ and $\Kfgsi Y$ follows from \cref{Kfgsix_istar}, \cref{Kfgsi_tha}, the universal property of coends, and the functoriality of $\sma$.  The functoriality of $\Kfgsi$ follows from the fact that identity morphisms and composition in $\FGTop$ are defined componentwise.
\end{description}
This finishes the definition of the functor $\Kfgsi$.
\end{definition}

\begin{definition}\label{def:kscgb}
For a finite group $G$, \emph{Schwede $K$-theory}\index{Schwede K-theory@Schwede $K$-theory}\index{K-theory@$K$-theory!Schwede} $\Kscgb$ is defined as the composite functor
\begin{equation}\label{kscgb_functor}
\begin{tikzpicture}[vcenter]
\def\h{2.3} \def\g{1.6} \def\u{.7} \def\d{.05}
\draw[0cell]
(0,0) node (a1) {\phantom{\Parcat}}
(a1)++(0,.02) node (a1') {\Parcat}
(a1)++(\h,0) node (a2) {\FGCat}
(a2)++(\h,0) node (a3) {\FGTop}
(a3)++(\h,0) node (a4) {\Sptopg}
;
\draw[1cell=.9]
(a1) edge node {\Jscg} (a2)
(a2) edge node {\clast} (a3)
(a3) edge node {\Kfgsi} (a4)
(a1') [rounded corners=2pt] |- ($(a2)+(0,\u)$) -- node {\Kscgb} ($(a3)+(0,\u)$) -| (a4)
;
\end{tikzpicture}
\end{equation}
of
\begin{itemize}
\item Schwede $J$-theory $\Jscg$ \cref{Jscg},
\item the functor $\clast$ given by postcomposing with the classifying space functor $\cla$ \cref{cla_gcat_gtop}, and
\item the functor $\Kfgsi$ \cref{kfgsi_functor}.\defmark
\end{itemize}
\end{definition}

\begin{explanation}[$\Kscg$ versus $\Kscgb$]\label{expl:kscgb}
By \cref{JscgCn,Kfgsixa}, for a parsummable category $(\C,\psum,\pzero)$ and a finite group $G$, the $G$-symmetric spectrum $\Kscgb\C$ sends each finite set $A$ to the pointed $G$-space
\begin{equation}\label{kscgbca}
(\Kscgb\C)_A = \int^{\ordn \in \Fsk} (S^A)^n \sma \cla\Cboxnomeg.
\end{equation}
The $G$-action on $\Cboxnomeg$ is given in \cref{jscg_gaction}.  The structure pointed $G$-morphisms are defined in \cref{Kfgsix_istar}.  

Comparing \cref{def:sch4.15,def:kscgb}, the only differences between the $G$-symmetric spectra $\kcg$ and $\Kscgb\C$ are their values at the empty set and the structure morphisms for inclusions $\emptyset \to B$.  More explicitly, there is a pointed $G$-morphism
\begin{equation}\label{Isc_empty}
\begin{tikzpicture}[vcenter]
\def\v{-1.4} \def\h{5.7}
\draw[0cell]
(0,0) node (a11) {(\kcg)_\emptyset}
(a11)++(\h,0) node (a12) {(\Kscgb\C)_\emptyset}
(a11)++(0,\v) node (a21) {\cla\C^{\suppempty}}
(a21)++(.4*\h,0) node (a22) {S^\emptyset \sma \cla\C^{\suppempty}}
(a12)++(0,\v) node (a23) {\phantom{S^\emptyset \sma \cla\Comeg}} 
(a23)++(0,-.03) node (a23') {S^\emptyset \sma \cla\Comeg}
;
\draw[1cell=.9]
(a11) edge node {\Isc^\C_\emptyset} (a12)
(a11) edge[equal] (a21)
(a21) edge node {\iso} (a22)
(a22) edge node {1 \sma \cla\iota} (a23)
(a23) edge node[swap] {\eta_{\ordone}} (a12)
;
\end{tikzpicture}
\end{equation}
defined by
\begin{itemize}
\item the canonical pointed homeomorphism $\iso$ \cref{zero_sphere},
\item the full subcategory inclusion $\iota \cn \C^{\suppempty} \to \Comeg = \C$, and
\item the pointed $G$-morphism $\eta_{\ordone}$ for the coend $(\Kscgb\C)_\emptyset$ at the object $\ordone \in \Fsk$.
\end{itemize}
There is a morphism of $G$-symmetric spectra
\begin{equation}\label{IscC}
\kcg \fto{\Isc^\C} \Kscgb\C
\end{equation}
with $A$-component pointed $G$-morphism $\Isc^\C_A$ given by
\begin{itemize}
\item $\Isc^\C_\emptyset$ \cref{Isc_empty} if $A=\emptyset$ and
\item the identity morphism for each nonempty finite set $A$.
\end{itemize}
There is  a natural transformation 
\begin{equation}\label{Isc_nt}
\begin{tikzpicture}[vcenter]
\def\t{28}
\draw[0cell]
(0,0) node (a1) {\phantom{A}}
(a1)++(1.8,0) node (a2) {\phantom{A}}
(a1)++(-.4,0) node (a1') {\Parcat}
(a2)++(.27,-.03) node (a2') {\Sptopg}
;
\draw[1cell=.9]
(a1) edge[bend left=\t] node {\Kscg} (a2)
(a1) edge[bend right=\t] node[swap] {\Kscgb} (a2)
;
\draw[2cell]
node[between=a1 and a2 at .42, rotate=-90, 2label={above,\Isc}] {\Rightarrow}
;
\end{tikzpicture}
\end{equation}
with component morphisms $\Isc^\C$ \cref{IscC}.
\end{explanation}

\begin{lemma}\label{Isc_geq}
For each finite group $G$, the natural transformation 
\[\Kscg \fto{\Isc} \Kscgb\]
in \cref{Isc_nt} is componentwise a $\pistu$-isomorphism between $G$-symmetric spectra.  Thus, $\Isc$ is componentwise a $G$-stable equivalence.
\end{lemma}

\begin{proof}
For each parsummable category $\C$, the $A$-component pointed $G$-morphism $\Isc^\C_A$ is the identity for each nonempty finite set $A$.  Thus, the morphism $\Isc^\C$ of $G$-symmetric spectra is a $\pistu$-isomorphism (\cref{def:sptop} \cref{def:sptop_v}), hence also a $G$-stable equivalence by \cite[3.36]{hausmann17}.
\end{proof}

\begin{explanation}\label{expl:Isc_geq}
Combining \cref{sch4.15,Isc_geq}, for each finite group $G$, there are natural $\pistu$-isomorphisms, hence natural $G$-stable equivalences, between $G$-symmetric spectra
\begin{equation}\label{IscC_Gequiv}
(\Ksc\C)_G \fto{\ascgc} \gacomegs \fot{\bscgc} \kcg \fto{\Isc^\C} \Kscgb\C.
\end{equation}
To study the underlying $G$-spectrum of global $K$-theory $\Ksc\C$, it suffices to consider the $G$-symmetric spectrum $\Kscgb\C$ \pcref{def:kscgb}. 
\end{explanation}

%% file: chap/kgl_gmmo.tex
This chapter proves that, for each finite group $G$ and parsummable category $\C$, there is a natural zigzag of $\pistu$-isomorphisms 
\begin{equation}\label{gl_gmmo_intro}
(\Ksc\C)_G \to \cdots \leftarrow \Kiopgb\Ig\C
\end{equation}
connecting
\begin{itemize}
\item the underlying $G$-spectrum $(\Ksc\C)_G$ of the global $K$-theory $\Ksc\C$ of $\C$ \pcref{def:sch3.3} and
\item the GMMO $K$-theory $\Kiopgb\Ig\C$ of the $\Einfg$-category $\Ig\C = \Comeg$ associated to $\C$ \cref{gacomeg}. 
\end{itemize} 
See \cref{gl_gmmo_geq}.  In particular, \cref{gl_gmmo_intro} is a natural zigzag of $G$-stable equivalences between $G$-symmetric spectra.  The $\Einfg$-category $\Ig\C$ and GMMO $K$-theory $\Kiopgb$ are defined using the chaotic $\Einfg$-operad $\Iopg$ \cref{iopgn}, called the categorical injection $G$-operad.  \cref{gl_gmmo_geq} is independent of \cref{thm:gmmo_shi}, which compares Shimakawa $K$-theory and GMMO $K$-theory.

\summary
\cref{gl_gmmo_geq} is proved as follows.    By \cref{sch4.15,Isc_geq}, there is a natural zigzag of $\pistu$-isomorphisms connecting global $K$-theory at $G$ and Schwede $K$-theory:
\[(\Ksc\C)_G \fto{\ascgc} \gacomegs \fot{\bscgc} \kcg \fto{\Isc^\C} \Kscgb\C.\]
With this natural zigzag of $\pistu$-isomorphisms in mind, this chapter focuses on comparing Schwede $K$-theory $\Kscgb\C$ and GMMO $K$-theory $\Kiopgb\Ig\C$, as summarized in the following diagram.
\begin{equation}\label{kiopg_ksc_intro}

\end{equation}
The boundary of \cref{kiopg_ksc_intro} consists of the following functors.
\begin{itemize}
\item The top composite $\Kscgb$ \cref{kscgb_functor} is Schwede $K$-theory for a finite group $G$. 
\item The upper-left vertical functor $\Ig$ \cref{Ig_functor} sends parsummable categories and parsummable functors to $\Iopg$-algebras and $\Iopg$-algebra morphisms for the chaotic $\Einfg$-operad $\Iopg$.
\item The left-bottom-right composite $\Kiopgb$ \cref{kiopgb_def} is a version of GMMO $K$-theory that produces $G$-symmetric spectra from $\Iopg$-algebras. 
\end{itemize} 
The interior of \cref{kiopg_ksc_intro} is given as follows.
\begin{itemize}
\item The natural transformation $\retn \cn \Bc \to 1$ is the retraction for the bar functor \cref{retn_barFG_id}.
\item Each instance of $\Lg$ is the left adjoint of an adjoint equivalence of categories \pcref{thm:fgcat_fgcatg_iieq,fgtop_fgtopg_eq}.
\item The two regions decorated by $\iso$ commute up to natural isomorphisms \pcref{kfgs_kfgsi_eq,Lg_clast_com}.   
\item The functor $\gxist$ is the pullback along $\gxi \cn \DG \to \FG$ \cref{gxist}, where $\DG$ is defined using $\Iopg$ \pcref{expl:dgiopg}. 
\item The inclusion functor $\Incj$ \cref{fgcatg_psfgcatg} is equal to the composite $\gzest\gxist$ \cref{gxist_gzest}.  Each component of the counit \cref{strv} 
\[\str\Incj \fto{\strv} 1_{\FGCatg}\]
of the 2-adjunction $(\str,\Incj)$ is componentwise a pointed $G$-equivalence between pointed $G$-categories.
\item The main comparison between Schwede and GMMO $K$-theories is the natural transformation \cref{cglg_nat} 
\[\gxist\Lg\Jscg \fto{\cglg} \Rg\Ig\]
that compares their first steps $\Jscg$ and $\Rg$.  \cref{cglgc_equiv} proves that each component of $\cglg$ is componentwise the left adjoint of an adjoint pointed $G$-equivalence between pointed $G$-categories.
\end{itemize}
\cref{gl_gmmo_geq} explains how these data produce the desired natural zigzag of $\pistu$-isomorphisms \cref{gl_gmmo_intro}. 

\subtleties
There are two subtle points about the comparison natural transformation $\cglg$.
\begin{description}
\item[Pseudonaturality] In \cref{kiopg_ksc_intro}, the functors $\gxist$ and $\Rg$ land in the subcategory (\cref{def:dgcatg_iicat} \cref{def:dgcatg_iicat_ii})
\[\dgcatg \subseteq \DGCatg\]
with $\Gcatst$-natural transformations as morphisms.  However, the comparison natural transformation $\cglg$ is only defined for $\DGCatg$ and \emph{not} $\dgcatg$.  The components of $\cglg$ are $\Pig$-strict $\Gcatst$-pseudotransformations \pcref{cglgc_welldef}, but generally \emph{not} $\Gcatst$-natural transformations.  Their lack of strict $\Gcatst$-naturality is controlled by the natural isomorphisms in \cref{cglgcx}.  Thus, a crucial aspect of the comparison $\cglg$ is that GMMO $K$-theory involves the category $\DGCatg$, whose morphisms are $\Pig$-strict $\Gcatst$-pseudotransformations instead of $\Gcatst$-natural transformations.  As we discuss in the introduction of \cref{ch:kgmmo}, the comparison between Shimakawa and GMMO $K$-theories involves a similar subtlety.
\item[Direction of comparison] By \cref{cglgc_equiv}, for each pointed finite $G$-set $\nbeta$ and each parsummable category $\C$, the comparison pointed $G$-functor 
\[(\gxist\Lg\Jscg\C)\nbeta \fto{\cglgcnbeta} (\Rg\Ig\C)\nbeta\]
admits a $G$-adjoint inverse $\cgglcnbeta$ \cref{cgglcnbeta_functor}.  As $\nbeta$ varies, the pointed $G$-functors $\cgglcnbeta$ assemble into a $\Gcatst$-pseudotransformation $\cgglc$.  However, $\cgglc$ is \emph{not} a morphism in the category $\DGCatg$ because it is not $\Pig$-strict \pcref{rk:cggl_notpigst}.  Thus, the direction of the comparison $\cglg$ cannot be reversed.
\end{description}

\organization
This chapter consists of the following sections.  Throughout this chapter, $G$ is a finite group.

\secname{sec:iopg}
This section defines the chaotic $\Einfg$-operad $\Iopg$ and the functor
\[\Parcat \fto{\Ig} \Algiopg\]
that sends parsummable categories to $\Iopg$-algebras.

\secname{sec:cglg}
This section defines the natural transformation $\cglg$ that compares the first step $\Jscg$ of Schwede $K$-theory and the first step $\Rg$ of GMMO $K$-theory for $\Iopg$.

\secname{sec:cglg_geq}
This section proves that each component of $\cglg$ is componentwise a pointed $G$-equivalence between pointed $G$-categories.

\secname{sec:gsp_fggspaces}
This section defines the prolongation functor 
\[\FGTopg \fto{\Kfgs} \Sptopg\] 
used in the comparison of GMMO and global $K$-theories.  The second half of this section records several properties of $\Kfgs$ needed in \cref{sec:kgl_kgmmo_eq}.  The functor $\Kfgs$ is homotopical when restricted to proper $\FGG$-spaces \pcref{thm:Kfgs_ht}.

\secname{sec:kgl_kgmmo_eq}
Using $\cglg$, this section constructs the natural zigzag of $\pistu$-isomorphisms \cref{gl_gmmo_intro} that connects the $G$-symmetric spectra $(\Ksc\C)_G$ and $\Kiopgb\Ig\C$.   \cref{expl:gl_gmmo_geq} further elaborates this zigzag at the point-set level.

\section{$\Einfg$-Categories from Parsummable Categories}
\label{sec:iopg}

To relate the domain categories of GMMO and global $K$-theories, this section constructs a chaotic $\Einfg$-operad $\Iopg$, called the categorical injection $G$-operad, and a functor
\[\Parcat \fto{\Ig} \Algiopg\]
from the category of parsummable categories and parsummable functors \pcref{def:parcat} to the category of $\Iopg$-algebras and morphisms.  At the object level, $\Ig$ sends each parsummable category $\C$ to an $\Iopg$-algebra $\Comeg$.  Since $\Iopg$ is a chaotic $\Einfg$-operad, each $\Iopg$-algebra, including $\Comeg$, is an $\Einfg$-category.  This fact is mentioned in \cite[4.20]{schwede_global}; this section provides a detailed construction.  Throughout this section, $G$ is a finite group.

\secoutline
\begin{itemize}
\item \cref{def:inop} defines the categorical injection operad $\Iop$.  \cref{parcat_iopalg} proves that parsummable categories are $\Iop$-algebras with underlying tame $\schm$-categories.  \cref{def:inop,parcat_iopalg} are not actually needed, but they motivate the constructions in the rest of this section.
\item \cref{def:global_einf} defines the categorical injection $G$-operad $\Iopg$.
\item \cref{iopg_einf} proves that $\Iopg$ is a chaotic $\Einfg$-operad.
\item \cref{def:comeg_iopg} defines the $\Iopg$-algebra structure on the $\omeg$-reparametrization $\Comeg$.
\item \cref{comeg_iopgalg} verifies the $\Iopg$-algebra axioms for $\Comeg$.
\item \cref{fomeg_iopgmor} shows that each parsummable functor $\fun$ yields an $\Iopg$-algebra morphism $\fomeg$. 
\item \cref{def:Ig} defines the functor $\Ig$.
\end{itemize}

\subsection*{Parsummable Categories as Operadic Algebras}
\label{sec:parcat_operad}

Recall the unpointed finite set $\ufsn = \{1,2,\ldots,n\}$ with $\ufs{0} = \emptyset$ and the countable set $\ome = \{0,1,2,\ldots\}$ of nonnegative integers.  

\begin{definition}\label{def:inop}
Define the operads $\Inop$ and $\Iop$ as follows.
\begin{description}
\item[Injection operad]
The \emph{injection operad}\index{injection operad}\index{operad!injection} $\Inop$ is the operad whose $n$-th set $\Inop(n)$ consists of injections $\iphi \cn \ufsn \itimes \ome \to \ome$.  
\begin{description}
\item[Unit] The operadic unit is the identity function 
\begin{equation}\label{iop_unit}
\ome \fto{1_{\ome}} \ome.
\end{equation}
\item[Symmetric group action]
The symmetric group $\Si_n$ permutes the first coordinate $\ufsn$, meaning
\begin{equation}\label{iop_symmetry}
(\iphi^\si)(j,-) = \iphi(\si j, -)
\end{equation}
for $\iphi \in \Inop(n)$, $\si \in \Si_n$, and $j \in \ufsn$.  
\item[Operadic composition]
For $n \geq 1$, $k_1,\ldots,k_n \geq 0$, and $k = \sum_{j\in \ufsn} k_j$, the composition
\[\Inop(n) \times \prod_{j \in \ufsn}\, \Inop(k_j) \fto{\ga} \Inop(k)\]
sends an $(n+1)$-tuple of injections $(\iphi; \ang{\iphi_j}_{j \in \ufsn})$ to the composite injection
\begin{equation}\label{iop_comp}
\ufsk \itimes \ome = \coprod_{j \in \ufsn}\, \ufsk_j \itimes \ome \fto{\sqcup_j\, \iphi_j} 
\coprod_{j \in \ufsn}\,\ome = \ufsn \itimes \ome \fto{\iphi} \ome.
\end{equation}
\end{description}
\item[Categorical injection operad]
Applying the translation category construction $\tn$ \pcref{def:translation_cat} to the injection operad $\Inop$ yields the \emph{categorical injection operad}\index{categorical injection operad}\index{operad!categorical injection} $\Iop$.  Its $n$-th category is the translation category
\begin{equation}\label{iopn}
\Iop(n) = \tn\Inop(n)
\end{equation}
with
\begin{itemize}
\item injections $\ufsn \itimes \ome \to \ome$ as objects and
\item each hom set given by a one-element set. 
\end{itemize} 
The operad structure on objects is given by \cref{iop_unit,iop_symmetry,iop_comp}, and it is uniquely determined on morphisms.
\end{description}
The unary category $\Iop(1)$ is the injection category $\schm$ (\cref{def:mcat} \cref{def:mcat_ii}).  The \emph{underlying $\schm$-category}\index{M-category@$\schm$-category!underlying} of an $\Iop$-algebra is obtained by restricting the $\Iop$-algebra structure to $\Iop(1) = \schm$.
\end{definition}

Recall the category $\Parcat$ of parsummable categories and parsummable functors \pcref{def:parcat}.  \cref{parcat_iopalg} is the categorical analogue of \cite[Appendix A]{sagave_schwede}, which concerns the injection operad $\Inop$ instead of the categorical injection operad $\Iop$.  It is analogous to the fact that small naive permutative $G$-categories are precisely algebras over the Barratt-Eccles operad in $\Gcat$ \pcref{def:BE}.  \cref{parcat_iopalg} is not used anywhere in this work.  It is included here to help motivate \cref{def:global_einf,def:comeg_iopg}.

\begin{lemma}\label{parcat_iopalg}
There is an isomorphism between $\Parcat$ and the category of $\Iop$-algebras with underlying tame $\schm$-categories.
\end{lemma}

\begin{proof}
The passage from parsummable categories to $\Iop$-algebras with underlying tame $\schm$-categories is given as follows.
\begin{description}
\item[Parsummable categories as $\Iop$-algebras] 
Suppose $(\C,\psum,\pzero)$ is a \index{parsummable category!as $\Iop$-algebra}parsummable category.  For each injection $\iphi \cn \ufsn \itimes \ome \to \ome$, the injections
\begin{equation}\label{iphi_j}
\iphi^j = \iphi(j,-) \cn \ome \to \ome \forspace j \in \ufsn
\end{equation}
have disjoint images.  Using the notation in \cref{u_action}, the $\iphi$-action $\iphi_* \cn \C^n \to \C$ is defined as the composite functor
\begin{equation}\label{iphist_caction}
\begin{tikzpicture}[vcenter]
\def\u{.7}
\draw[0cell]
(0,0) node (a1) {\C^n}
(a1)++(2.3,0) node (a2) {\phantom{\Cboxn}}
(a2)++(0,.03) node (a2') {\Cboxn}
(a2)++(1.5,0) node (a3) {\C}
;
\draw[1cell=.85]
(a1) edge node {\ang{\iphist^j}_{j \in \ufsn}} (a2)
(a2) edge node {\psum} (a3)
(a1) [rounded corners=2pt] |- ($(a2)+(-1,\u)$) -- node {\iphi_*} ($(a2)+(.5,\u)$) -| (a3) 
;
\end{tikzpicture}
\end{equation}
that sends an object or a morphism $\objx = \ang{x_j}_{j \in \ufsn} \in \C^n$ to the sum
\begin{equation}\label{iphi_action}
\iphi_* \objx = \sum_{j \in \ufsn}\, \iphist^j x_j.
\end{equation}
\begin{itemize}
\item By \cref{sch2.13} \cref{sch2.13_iii}, for each object $\ang{x_j}_{j \in \ufsn} \in \C^n$, the objects $\iphist^j x_j \in \C$ for $j \in \ufsn$ have disjoint supports.  Thus, the object $\ang{\iphist^j x_j}_{j \in \ufsn} \in \C^n$ lies in the full subcategory $\Cboxn$ of disjointly supported objects \pcref{def:bxtimes}.  The iterated sum functor $\psum$ is well defined by the parsummable category axioms \cref{parcat_axioms} and the subadditivity of supports \cref{supp_xsumy}. 
\item For the case $n=0$, $\iphi \cn \emptyset \itimes \ome \to \ome$ is the unique injection, and $\iphi_* \cn \bone \to \C$ is given by the unit $\pzero \in \C$.
\end{itemize}
Using the notation in \cref{vu_star,vustar_x}, given an injection $\ipsi \cn \ufsn \itimes \ome \to \ome$ and the isomorphism $[\ipsi, \iphi] \cn \iphi \to \ipsi$ in $\Iop(n)$, the natural isomorphism
\begin{equation}\label{ipsiiphi_st}
\begin{tikzpicture}[vcenter]
\def\t{28}
\draw[0cell]
(0,0) node (a1) {\phantom{\C}}
(a1)++(-.07,0) node (a1') {\C^n}
(a1)++(2.5,0) node (a2) {\C}
;
\draw[1cell=.9]
(a1) edge[bend left=\t] node {\iphi_*} (a2)
(a1) edge[bend right=\t] node[swap] {\ipsi_*} (a2)
;
\draw[2cell=.9]
node[between=a1 and a2 at .33, rotate=-90, 2label={above,[\ipsi, \iphi]_*}] {\Rightarrow}
;
\end{tikzpicture}
\end{equation}
sends an object $\objx = \ang{x_j}_{j \in \ufsn} \in \C^n$ to the isomorphism
\begin{equation}\label{ipsiiphi_x}
[\ipsi, \iphi]_*^\objx = \sum_{j \in \ufsn}\, [\ipsi^j, \iphi^j]_{\mstar} ^{x_j} \cn
\sum_{j \in \ufsn}\, \iphist^j x_j \fiso \sum_{j \in \ufsn}\, \ipsist^j x_j
\end{equation}
in $\C$.  For the case $n=0$, $[\ipsi,\iphi]_*$ is given by the identity morphism of the unit $\pzero \in \C$.  The $\Iop$-algebra axioms hold by the $\schm$-category axioms \cref{mcat_axioms} and the parsummable category axioms \cref{parcat_axioms} for $\C$.  The restriction of the $\Iop$-action to $\Iop(1) = \schm$ is the given tame $\schm$-category structure on $\C$.  Each parsummable functor becomes an $\Iop$-algebra morphism by the $\schm$-functor axiom \cref{mfunctor_axiom} and the parsummable functor axioms \cref{parfun_axioms}.
\item[$\Iop$-algebras as parsummable categories] 
For the other direction, suppose $\C$ is an $\Iop$-algebra whose underlying $\schm$-category is tame (\cref{def:mcat} \cref{def:mcat_v}).  The tame $\schm$-category $\C$ becomes a parsummable category $(\C,\psum,\pzero)$ as follows.
\begin{description}
\item[Unit] The unit $\pzero \in \C$ is the image of the $\iphi$-action $\iphi_* \cn \bone \to \C$ for the unique injection $\iphi \cn \emptyset \itimes \ome \to \ome$ in $\Iop(0) = \bone$.
\item[Sum] Using the $\Iop$-algebra structure on $\C$, the functor $\psum \cn \C \bxtimes \C \to \C$ sends each object $(x_1,x_2) \in \C \bxtimes \C$ to the object
\begin{equation}\label{xone_sum_xtwo}
x_1 \psum x_2 = \iphi_*(x_1,x_2) \in \C
\end{equation}
for any choice of an injection $\iphi \cn \ufstwo \itimes \ome \to \ome$ such that
\[\iphi^j = \iphi(j,-) \cn \ome \to \ome\]
is the identity on $\supp(x_j)$ for each $j \in \{1,2\}$.  At least one such $\iphi$ exists because $x_1$ and $x_2$ have disjoint supports.  Moreover, the sum is independent of the choice of $\iphi$ by \cite[A.6 (i)]{sagave_schwede}.  On morphisms of $\C\bxtimes\C$, the sum is also defined by the formula \cref{xone_sum_xtwo}.  
\end{description}
Each morphism between $\Iop$-algebras with underlying tame $\schm$-categories preserves the unit $\pzero$ and the sum $\psum$ because they are parts of the $\Iop$-algebra structure. 
\end{description}
Minor extensions of \cite[A.11 and A.13]{sagave_schwede} show that $(\C,\psum,\pzero)$ is a parsummable category and that the preceding constructions define mutually inverse isomorphisms between $\Parcat$ and the category of $\Iop$-algebras with underlying tame $\schm$-categories.
\end{proof}

\subsection*{Injection $G$-Operads}

Recall that, for a finite group $G$, the universal $G$-set $\omeg$ \pcref{ex:omea} consists of functions $f \cn G \to \ome$ with the $G$-action $gf = f(\ginv \cdot -)$.  The group $G$ acts trivially on each unpointed finite set $\ufsn = \{1,2,\ldots,n\}$ and the countable set $\ome = \{0,1,2,\ldots\}$.  Recall the category $\Gcat$ of small $G$-categories and $G$-functors \pcref{def:GCat}.  \cref{def:global_einf} defines the $G$-equivariant analogues of the injection operad $\Inop$ and the categorical injection operad $\Iop$ \pcref{def:inop}.

\begin{definition}\label{def:global_einf}
For a finite group $G$, define the operads $\Inopg$ and $\Iopg$ as follows.
\begin{description}
\item[Injection $G$-operad]
The \emph{injection $G$-operad}\index{injection G-operad@injection $G$-operad}\index{operad!injection G-@injection $G$-} $\Inopg$ is the operad in the category of $G$-sets and $G$-morphisms whose $n$-th $G$-set $\Inopg(n)$ consists of injections $\iphi \cn \ufsn \itimes \omeg \to \omeg$.  The group $G$ acts by conjugation, meaning
\begin{equation}\label{inopg_gact}
(g\iphi\ginv)(j,f) = g\big(\iphi(j,\ginv f) \big)
\end{equation}
for $\iphi \in \Inopg(n)$, $g \in G$, $j \in \ufsn$, and $f \in \omeg$.
\begin{description}
\item[Unit]
The operadic unit is the identity function 
\begin{equation}\label{inopg_unit}
\omeg \fto{1_{\omeg}} \omeg.
\end{equation} 
\item[Symmetric group action] The symmetric group $\Si_n$ permutes the first coordinate $\ufsn$, meaning
\begin{equation}\label{inopg_symmetry}
(\iphi^\si)(j,-) = \iphi(\si j, -)
\end{equation}
for $\iphi \in \Inopg(n)$, $\si \in \Si_n$, and $j \in \ufsn$.  
\item[Operadic composition] 
For $n \geq 1$, $k_1,\ldots,k_n \geq 0$, and $k = \sum_{j\in \ufsn} k_j$, the composition $G$-morphism
\[\Inopg(n) \times \prod_{j \in \ufsn}\, \Inopg(k_j) \fto{\ga} \Inopg(k)\]
sends an $(n+1)$-tuple of injections $(\iphi; \ang{\iphi_j}_{j \in \ufsn})$ to the composite injection
\begin{equation}\label{inopg_comp}
\ufsk \itimes \omeg = \coprod_{j \in \ufsn}\, \ufsk_j \itimes \omeg \fto{\sqcup_j\, \iphi_j} 
\coprod_{j \in \ufsn}\,\omeg = \ufsn \itimes \omeg \fto{\iphi} \omeg.
\end{equation}
\end{description}
\item[Categorical injection $G$-operad]
Applying the translation category construction $\tn$ \pcref{def:translation_cat} to $\Inopg$ yields the \emph{categorical injection $G$-operad}\index{categorical injection G-operad@categorical injection $G$-operad}\index{operad!categorical injection G-@categorical injection $G$-} $\Iopg$.  It is the $\Gcat$-operad whose $n$-th $G$-category is the translation category
\begin{equation}\label{iopgn}
\Iopg(n) = \tn\Inopg(n)
\end{equation}
with
\begin{itemize}
\item injections $\ufsn \itimes \omeg \to \omeg$ as objects and
\item each hom set given by a one-element set. 
\end{itemize} 
The group $G$ acts by conjugation on objects \cref{inopg_gact} and on isomorphisms $[\ipsi,\iphi] \cn \iphi \to \ipsi$ in $\Iopg(n)$, meaning
\begin{equation}\label{inopg_gmor}
g[\ipsi,\iphi]\ginv = [g\ipsi\ginv, g\iphi\ginv] \cn g\iphi\ginv \to g\ipsi\ginv
\end{equation}
for $g \in G$.  The operad structure on objects is given by \cref{inopg_unit,inopg_symmetry,inopg_comp}.  On morphisms, the symmetric group action is given by
\begin{equation}\label{inopg_sym_mor}
[\ipsi,\iphi]^\si = [\ipsi^\si,\iphi^\si] \cn \iphi^\si \to \ipsi^\si.
\end{equation}
The composition is given by
\begin{equation}\label{inopg_comp_mor}
\ga\big([\ipsi,\iphi] ; \ang{[\ipsi_j,\iphi_j]}_{j \in \ufsn} \big) 
= \big[\ga(\ipsi; \ang{\ipsi_j}_{j \in \ufsn}), \ga(\iphi; \ang{\iphi_j}_{j \in \ufsn})\big]
\end{equation}
for injections $\iphi,\ipsi \in \Iopg(n)$ and $\iphi_j, \ipsi_j \in \Iopg(k_j)$.
\end{description}
Denote by $\Algiopg$ the category of $\Iopg$-algebras and $\Iopg$-algebra morphisms \pcref{def:operadalg}.
\end{definition}

\begin{explanation}[Unpacking]\label{expl:iopg}
The operadic composition $\ga$ in \cref{inopg_comp} is given by
\[\big(\ga(\iphi; \ang{\iphi_j}_{j \in \ufsn})\big) \big(k_1 + \Cdots + k_{i-1} + t , f \big)
= \iphi\big(i, \iphi_i(t,f)\big)\]
for $i \in \ufsn$, $t \in \ufsk_i$, and $f \in \omeg$.
\end{explanation}

Recall chaotic $\Einfg$-operads \pcref{def:chaotic_einf}.

\begin{lemma}\label{iopg_einf}
For each finite group $G$, the categorical injection $G$-operad $\Iopg$ is a chaotic $\Einfg$-operad.
\end{lemma}

\begin{proof}
The $\Gcat$-operad $\gmvg$ in \cite[7.4]{gm17} is defined in the same way as $\Iopg$ \pcref{def:global_einf} using the universal $G$-set $\ugsetw$ \cref{ugsetw} instead of $\omeg$.  By \cite[7.8]{gm17}, $\gmvg$ is a chaotic $\Einfg$-operad.  Since the universal $G$-sets $\omeg$ and $\ugsetw$ are $G$-isomorphic (\cref{schwede2.17} \cref{sch2.17_ii}), the $\Gcat$-operads $\Iopg$ and $\gmvg$ are isomorphic.  Thus, $\Iopg$ is also a chaotic $\Einfg$-operad.
\end{proof}

\subsection*{Reparametrizations as $\Einfg$-Categories}

\cref{def:comeg_iopg} extends the $\Iop$-algebra structure on a parsummable category \pcref{parcat_iopalg} to an $\Iopg$-algebra structure on the $G$-category $\Comeg$ (\cref{def:sch2.21} \cref{def:sch2.21_i}).  Recall that $\Comeg$ has underlying category $\C$.  By \cref{rep_S,Cbrac_g}, for each $g \in G$, the $g$-action on $\Comeg$ is the $\schm$-action functor $g^\ome_{\mstar} \cn \C \to \C$ for the bijection
\begin{equation}\label{gome_mstar}
\begin{tikzpicture}[vcenter]
\def\h{1.8} \def\g{1.5} \def\u{.7} \def\d{.05}
\draw[0cell]
(0,0) node (a1) {\ome}
(a1)++(\h,0) node (a2) {\phantom{\omeg}}
(a2)++(0,\d) node (a2') {\omeg}
(a2)++(\h,0) node (a3) {\phantom{\omeg}}
(a3)++(0,\d) node (a3') {\omeg}
(a3)++(\h,0) node (a4) {\ome.}
;
\draw[1cell=.9]
(a1) edge node {\rep_{\omeg}} (a2)
(a2) edge node {g} (a3)
(a3) edge node {\rep_{\omeg}^{-1}} (a4)
(a1) [rounded corners=2pt] |- ($(a2)+(0,\u)$) -- node {\gome} ($(a3)+(0,\u)$) -| (a4)
;
\end{tikzpicture}
\end{equation}

\begin{definition}\label{def:comeg_iopg}\index{reparametrization!as $\Iopg$-algebra}
For a finite group $G$, the categorical injection $G$-operad $\Iopg$ \cref{iopgn}, and a parsummable category $(\C,\psum,\pzero)$ \pcref{def:parcat}, define the structure of an $\Iopg$-algebra on the $G$-category $\Comeg$ with $\Iopg$-action $G$-functor
\begin{equation}\label{gacomeg}
\Iopg(n) \times \Comeg^n \fto{\gacomeg} \Comeg
\end{equation}
given as follows for $n \geq 0$.  
\begin{description}
\item[$\Iopg$-action on objects] 
For an injection $\iphi \cn \ufsn \itimes \omeg \to \omeg$, the injections\label{not:iphij}
\[\iphi^j = \iphi(j,-) \cn \omeg \to \omeg \forspace j \in \ufsn\]
have disjoint images.  The $n$ injections
\begin{equation}\label{iphijome}
\begin{tikzpicture}[vcenter]
\def\h{1.8} \def\g{1.5} \def\u{.7} \def\d{.05}
\draw[0cell]
(0,0) node (a1) {\ome}
(a1)++(\h,0) node (a2) {\phantom{\omeg}}
(a2)++(0,\d) node (a2') {\omeg}
(a2)++(\h,0) node (a3) {\phantom{\omeg}}
(a3)++(0,\d) node (a3') {\omeg}
(a3)++(\h,0) node (a4) {\ome}
;
\draw[1cell=.9]
(a1) edge node {\rep_{\omeg}} (a2)
(a2) edge node {\iphi^j} (a3)
(a3) edge node {\rep_{\omeg}^{-1}} (a4)
(a1) [rounded corners=2pt] |- ($(a2)+(0,\u)$) -- node {\iphijome} ($(a3)+(0,\u)$) -| (a4)
;
\end{tikzpicture}
\end{equation}
also have disjoint images, and 
\begin{equation}\label{iphijomest}
\C \fto{\iphijomest} \C 
\end{equation}
denotes the $\iphijome$-action functor on $\C$ \cref{u_action}.  For an object or a morphism $\objx = \ang{x_j}_{j \in \ufsn} \in \Comeg^n$, define the object or morphism
\begin{equation}\label{iphi_comeg}
\gacomeg(\iphi; \objx) = \iphi_* \objx = \sum_{j \in \ufsn}\, \iphijomest x_j \inspace \C.
\end{equation}
This iterated sum is given by
\begin{itemize}
\item the sum on $\C$ if $n \geq 1$ and
\item the unit $\pzero \in \C$ if $n=0$.
\end{itemize}
It is well defined by \cref{sch2.13} \cref{sch2.13_iii}, \cref{parcat_axioms,supp_xsumy,iphijome}.  The functoriality of 
\[\Comeg^n \fto{\gacomeg(\iphi;-)} \Comeg\]
follows from the functoriality of the sum $\psum$ and each $\schm$-action functor $\iphijomest$.
\item[$\Iopg$-action on morphisms] 
For an isomorphism $[\ipsi,\iphi] \cn \iphi \to \ipsi$ in $\Iopg(n)$, the natural isomorphism
\begin{equation}\label{phipsi_comeg}
\begin{tikzpicture}[vcenter]
\def\h{4} \def\u{.6} \def\s{.85}
\draw[0cell]
(0,0) node (a1) {\Comeg^n}
(a1)++(\h,0) node (a2) {\Comeg}
(a1)++(.5*\h,0) node (a0) {\phantom{a}}
;
\draw[1cell=\s]
(a1) [rounded corners=2pt] |- ($(a0)+(-1,\u)$) -- node {\iphi_* = \gacomeg(\iphi;-)} ($(a0)+(1,\u)$) -| (a2)
;
\draw[1cell=\s]
(a1) [rounded corners=2pt] |- ($(a0)+(-1,-\u)$) -- node[swap] {\ipsi_* = \gacomeg(\ipsi;-)} ($(a0)+(1,-\u)$) -| (a2)
;
\draw[2cell=1]
node[between=a1 and a0 at .5, rotate=-90, 2label={above,\gacomeg([\ipsi,\iphi];-)}] {\Rightarrow}
;
\end{tikzpicture}
\end{equation}
sends an object $\objx \in \Comeg^n$ to the isomorphism
\begin{equation}\label{phipsi_comeg_x}
\gacomeg\big([\ipsi,\iphi]; \objx\big) 
= \sum_{j \in \ufsn}\, [\ipsijome, \iphijome]_{\mstar}^{x_j} \cn
\sum_{j \in \ufsn} \iphijomest x_j \fiso \sum_{j \in \ufsn} \ipsijomest x_j
\end{equation}
in $\C$.  For each $j \in \ufsn$,\label{not:ipsiphijome} 
\[\iphijomest \fto{[\ipsi^{j,\ome}, \iphi^{j,\ome}]_{\mstar}} \ipsijomest\]
is the $\schm$-action natural isomorphism on $\C$ \cref{vu_star} for the isomorphism 
\[\iphijome \fto{[\ipsijome, \iphijome]} \ipsijome \inspace \schm.\]
For $n=0$, the isomorphism in \cref{phipsi_comeg_x} is interpreted as the identity morphism $1_\pzero$.
\end{description}
This finishes the definition of the functor $\gacomeg$ \cref{gacomeg}.
\end{definition}

\begin{lemma}\label{comeg_iopgalg}
In the context of \cref{def:comeg_iopg}, $(\Comeg,\gacomeg)$ is an $\Iopg$-algebra in $\Gcat$.
\end{lemma}

\begin{proof}
We verify the $\Iopg$-algebra axioms \pcref{def:operadalg} for $\Comeg$.
\begin{description}
\item[$G$-equivariance] 
The $\Iopg$-action functor $\gacomeg$ \cref{gacomeg} is $G$-equivariant on objects if and only if
\begin{equation}\label{gacomeg_geq}
\gacomeg(g\iphi\ginv; g\objx) = g\gacomeg(\iphi; \objx) \inspace \Comeg
\end{equation}
for $g \in G$, $\iphi \in \Iopg(n)$, and $\objx = \ang{x_j}_{j\in \ufsn} \in \Comeg^n$.  Using \cref{mcat_axioms,inopg_gact,gome_mstar,iphijome,iphi_comeg}, the desired equality \cref{gacomeg_geq} is proved as follows.
\[\begin{split}
& \gacomeg(g\iphi\ginv; g\objx) \\
&= \txsum_{j \in \ufsn}\, (g\iphi\ginv)^{j,\ome}_{\mstar} (g^\ome_{\mstar} x_j) \\
&= \txsum_{j \in \ufsn} \big((\rep_{\omeg}^{-1} g \iphi^j \ginv \rep_{\omeg}) (\rep_{\omeg}^{-1} g \rep_{\omeg})\big)_{\mstar} (x_j) \\
&= \txsum_{j \in \ufsn} \big(\rep_{\omeg}^{-1} g \iphi^j \rep_{\omeg} \big)_{\mstar} (x_j) \\
&= \txsum_{j \in \ufsn} \big((\rep_{\omeg}^{-1} g \rep_{\omeg}) (\rep_{\omeg}^{-1} \iphi^j \rep_{\omeg}) \big)_{\mstar} (x_j) \\
&= \txsum_{j \in \ufsn}\, g^{\ome}_{\mstar} \iphi^{j,\ome}_{\mstar} x_j \\
&= g^{\ome}_{\star} \big(\txsum_{j \in \ufsn}\, \iphi^{j,\ome}_{\mstar} x_j \big) \\
&= g\gacomeg(\iphi; \objx)
\end{split}\]
The second-to-last equality uses the $\schm$-equivariance of the sum on $\C$.  The $G$-equivariance of $\gacomeg$ on morphisms is proved by the preceding equalities using \cref{inopg_gmor,phipsi_comeg_x}.
\item[Unity]
For the identity function $1_{\omeg} \cn \omeg \to \omeg$, the $G$-functor
\[\Comeg \fto{\gacomeg(1_{\omeg}; -)} \Comeg\]
is the identity because $(1_{\ome})_{\mstar} = 1_{\C}$ by the unity axiom of the $\schm$-category $\C$ \cref{mcat_axioms}.
\item[Symmetry]
Using \cref{inopg_symmetry}, \cref{iphi_comeg}, and the symmetry axiom \cref{parcat_axioms} for the sum, the following equalities for permutations $\si \in \Si_n$ prove the symmetry axiom for $\gacomeg$ on objects.
\[\begin{split}
\gacomeg(\iphi^\si; \objx) 
&= \txsum_{j \in \ufsn}\, (\iphi^\si)^{j,\ome}_{\mstar} x_j 
= \txsum_{j \in \ufsn}\, \iphi^{\si j,\ome}_{\mstar} x_j \\
&= \txsum_{j \in \ufsn}\, \iphi^{j,\ome}_{\mstar} x_{\sigmainv j} 
= \gacomeg(\iphi; \si \objx)
\end{split}\]
The symmetry axiom for $\gacomeg$ on morphisms is proved by the preceding equalities using \cref{phipsi_comeg_x,inopg_sym_mor}.
\item[Associativity]
Using \cref{inopg_comp}, \cref{iphi_comeg}, the associativity of the $\schm$-action on $\C$ \cref{mcat_axioms}, the associativity of the sum \cref{parcat_axioms}, and the $\schm$-equivariance of the sum, the following equalities prove the associativity axiom for $\gacomeg$ on objects, where $\iphi \in \Iopg(n)$ with $n \geq 1$, $\iphi_j \in \Iopg(k_j)$ for $j \in \ufsn$, and $x_{j,t} \in \Comeg$ for $j \in \ufsn$ and $t \in \ufsk_j$.
\[\begin{split}
& \gacomeg\big( \ga(\iphi; \ang{\iphi_j}_{j \in \ufsn}) ; \ang{x_{j,t}}_{t\in \ufsk_j \csp j \in \ufsn}\big) \\
&= \txsum_{j \in \ufsn} \txsum_{t \in \ufsk_j}\, (\iphi^j \iphi_j^t)^\ome_{\mstar} (x_{j,t}) \\
&= \txsum_{j \in \ufsn} \txsum_{t \in \ufsk_j} \big(\rep_{\omeg}^{-1} \iphi^j \iphi_j^t \rep_{\omeg} \big)_{\mstar} (x_{j,t})\\
&= \txsum_{j \in \ufsn} \txsum_{t \in \ufsk_j} \, (\rep_{\omeg}^{-1} \iphi^j \rep_{\omeg})_{\mstar} (\rep_{\omeg}^{-1} \iphi_j^t \rep_{\omeg})_{\mstar} (x_{j,t}) \\
&= \txsum_{j \in \ufsn} \txsum_{t \in \ufsk_j} \, \iphi^{j,\ome}_{\mstar} (\iphi_j^t)^\ome_{\mstar} (x_{j,t})\\
&= \txsum_{j \in \ufsn}\, \iphi^{j,\ome}_{\mstar} \big(\sum_{t \in \ufsk_j} \, (\iphi_j^t)^\ome_{\mstar} (x_{j,t})\big) \\
&= \gacomeg\big(\iphi; \bang{\gacomeg(\iphi_j ; \ang{x_{j,t}}_{t\in \ufsk_j})}_{j \in \ufsn} \big)
\end{split}\]
\end{description}
The associativity axiom for $\gacomeg$ on morphisms is proved by the preceding equalities using \cref{phipsi_comeg_x,inopg_comp_mor}.  This finishes the proof that $\Comeg$ is an $\Iopg$-algebra.
\end{proof}

\cref{fomeg_iopgmor} extends \cref{comeg_iopgalg} to parsummable functors \cref{parfun_axioms}.

\begin{lemma}\label{fomeg_iopgmor}
For a finite group $G$, each parsummable functor $\fun \cn \C \to \D$ yields an $\Iopg$-algebra morphism
\begin{equation}\label{fomeg_iopalg}
(\Comeg,\gacomeg) \fto{\fomeg} (\Domeg,\gadomeg).
\end{equation}
\end{lemma}

\begin{proof}
The $G$-functor $\fomeg = \fun$ \cref{fomeg} is an $\Iopg$-algebra morphism \pcref{def:laxmorphism} if the diagram
\begin{equation}\label{fomeg_ga}
\def\v{-1.4}
\begin{tikzpicture}[vcenter]
\draw[0cell]
(0,0) node (a11) {\Iopg(n) \times \Comeg^n}
(a11)++(3,0) node (a12) {\Comeg}
(a11)++(0,\v) node (a21) {\Iopg(n) \times \Domeg^n}
(a12)++(0,\v) node (a22) {\Domeg}
;
\draw[1cell=.85]
(a11) edge node {\gacomeg} (a12)
(a21) edge node {\gadomeg} (a22)
(a11) edge[transform canvas={xshift=1em}] node[swap] {1 \times \fun^n} (a21)
(a12) edge node {\fun} (a22)
;
\end{tikzpicture}
\end{equation}
commutes for each $n \geq 0$.  The diagram \cref{fomeg_ga} commutes by \cref{iphi_comeg}, \cref{phipsi_comeg_x}, the parsummable axioms for $\fun$ \cref{parfun_axioms}, and the $\schm$-equivariance of $\fun$ \cref{mfunctor_axiom}.
\end{proof}

Recall that $\Algiopg$ is the category of $\Iopg$-algebras and $\Iopg$-algebra morphisms \pcref{def:global_einf}.

\begin{definition}\label{def:Ig}
For a finite group $G$, define the functor
\begin{equation}\label{Ig_functor}
\Parcat \fto{\Ig} \Algiopg
\end{equation}
with object assignment \label{not:IgC}\pcref{comeg_iopgalg} 
\[\C \mapsto \Ig\C = (\Comeg,\gacomeg)\] 
and morphism assignment  \pcref{fomeg_iopgmor}
\[\fun \mapsto \Ig\fun = \fomeg.\]
The functoriality of $\Ig$ follows from the fact that the underlying $G$-functor of $\fomeg$ is $\fun$ \cref{fomeg}.
\end{definition}

\section{Comparison Natural Transformation}
\label{sec:cglg}

To compare GMMO and global $K$-theories for a finite group $G$, this section constructs the comparison natural transformation
\begin{equation}\label{cglg_diag_secti}
\begin{tikzpicture}[vcenter]
\def\h{2.5} \def\v{1.4}
\draw[0cell]
(0,0) node (a1) {\Parcat}
(a1)++(.95*\h,0) node (a2) {\Algiopg}
(a2)++(1.1*\h,0) node (a3) {\DGCatg}
(a1)++(0,\v) node (b1) {\FGCat}
(a3)++(0,\v) node (b2) {\FGCatg}
;
\draw[1cell=.85]
(a1) edge node {\Ig} (a2)
(a2) edge node {\Rg} (a3)
(a1) edge node {\Jscg} (b1)
(b1) edge node {\Lg} (b2)
(b2) edge node {\gxist} (a3)
;
\draw[2cell]
node[between=b1 and b2 at .45, shift={(0,-.5*\v)}, rotate=-90, 2label={above,\cglg}] {\Rightarrow}
;
\end{tikzpicture}
\end{equation}
involving
\begin{itemize}
\item the categories $\FGCat$ \pcref{def:fgcat}, $\FGCatg$ \pcref{def:fgcatg}, $\DGCatg$ \pcref{def:dgcatg_iicat}, $\Parcat$ \pcref{def:parcat}, and $\Algiopg$ \pcref{def:global_einf}; and
\item the functors $\Jscg$ \cref{Jscg}, $\Lg$ \pcref{thm:fgcat_fgcatg_iieq}, $\gxist$ \cref{gxist}, $\Ig$ \cref{Ig_functor}, and $\Rg$ \cref{rg_twofunctor}.
\end{itemize}
The functor $\Rg$ is the first step of GMMO $K$-theory $\Kiopg$ \cref{kgmmo_diag} restricted to the subcategory $\Algiopg$ for the chaotic $\Einfg$-operad $\Iopg$ \cref{iopgn}.  Schwede $J$-theory $\Jscg$ is the first step of Schwede $K$-theory $\Kscgb$ \cref{kscgb_functor}.  The categories $\Algiopg$ and $\DGCatg$, as well as the functors $\gxist$ and $\Rg$, are defined using $\Iopg$.  Throughout this section, $G$ is a finite group.

\secoutline
\begin{itemize}
\item \cref{expl:dgiopg} describes the $\Gcatst$-category $\DG$ and the category $\DGCatg$ defined using $\Iopg$.
\item \cref{expl:LgJscg} discusses the composite functor $\Lg\Jscg$.
\item \cref{expl:cglg_dom} discusses the composite functor $\gxist\Lg\Jscg$, which is the domain of $\cglg$.
\item \cref{expl:cglg_cod} discusses the composite functor $\Rg\Ig$, which is the codomain of $\cglg$.
\item \cref{def:cglg} constructs the natural transformation $\cglg$.
\item \cref{cglgc_welldef} proves that each component of $\cglg$ is a morphism in $\DGCatg$.
\item \cref{cglg_natural} proves that $\cglg$ is a natural transformation.
\end{itemize}

\subsection*{Domain and Codomain of $\cglg$}
To facilitate the construction of the comparison natural transformation $\cglg$, \cref{expl:dgiopg,expl:cglg_dom,expl:cglg_cod} describe $\DG$, $\DGCatg$, the domain functor, and the codomain functor of $\cglg$.  Recall the unpointed finite set $\ufsn = \{1,2,\ldots,n\}$ with $\ufs{0} = \emptyset$.

\begin{explanation}[$\DG$ and $\DGCatg$]\label{expl:dgiopg}
Using the chaotic $\Einfg$-operad $\Iopg$ \cref{iopgn}, the $\Gcatst$-category $\DG$ \pcref{def:dgo} has pointed finite $G$-sets of the form $\mal$ \cref{ordn_be} as objects.  For each pair $(\mal,\nbeta)$ of pointed finite $G$-sets, it is equipped with the hom pointed $G$-category
\begin{equation}\label{dgmn_iopg}
\DG(\mal,\nbeta) = \coprod_{\psi \in \FG(\mal,\,\nbeta)} \,\prod_{j \in \ufsn}\, \Iopg(|\psiinv j|).
\end{equation}
\begin{description}
\item[Objects]
An object 
\begin{equation}\label{dgmn_iopg_obj}
\objx = \big(\psi; \ang{\iphi_j}_{j \in \ufsn}\big) \inspace \DG(\mal,\nbeta)
\end{equation}
consists of
\begin{itemize}
\item a pointed morphism $\psi \cn \mal \to \nbeta$ and
\item an object $\iphi_j \in \Iopg(|\psiinv j|)$ for each $j \in \ufsn$, meaning an injection 
\[\{1,2,\ldots,|\psiinv j|\} \itimes \omeg = \ufs{|\psiinv j|} \itimes \omeg \fto{\iphi_j} \omeg.\]
\end{itemize}
Recall that $\omeg$ is the universal $G$-set of functions $f \cn G \to \ome = \{0,1,2,\ldots\}$ with the $G$-action $gf = f(\ginv \cdot -)$ \pcref{def:omea}.  The basepoint of $\DG(\mal,\nbeta)$ is the object
\begin{equation}\label{dgiopg_bp}
\objzero = (0; \ang{*}_{j \in \ufsn})
\end{equation}
consisting of
\begin{itemize}
\item the 0-morphism $0 \cn \mal \to \nbeta$ and
\item $n$ copies of the unique object $* \in \Iopg(0) = \bone$.
\end{itemize}  
The identity 1-cell of an object $\mal \in \DG$ is the object 
\begin{equation}\label{dgiopg_idonecell}
\objone_{\mal} = \big(1_{\mal}; \ang{1_{\omeg}}_{i \in \ufsm}\big) \in \DG(\mal,\mal),
\end{equation}
where the identity function $(1_{\omeg} \cn \omeg \to \omeg) \in \Iopg(1)$ is the operadic unit.
\item[Morphisms]
A morphism in $\DG(\mal,\nbeta)$ has the form
\begin{equation}\label{dgmn_iopg_mor}
\objx = \big(\psi; \ang{\iphi_j}_{j \in \ufsn}\big) 
\fto{\objp \,=\, (\psi; \ang{[\ivphi_j,\iphi_j]}_{j \in \ufsn})}
\objy = \big(\psi; \ang{\ivphi_j}_{j \in \ufsn}\big)
\end{equation}
with each $[\ivphi_j,\iphi_j] \cn \iphi_j \to \ivphi_j$ an isomorphism in $\Iopg(|\psiinv j|)$.  The identity of $\objx$ is the morphism
\begin{equation}\label{dgiopg_idtwo}
\objone_{\objx} = \big(\psi; \ang{1_{\iphi_j}}_{j \in \ufsn}\big).
\end{equation}
\item[$G$-action]
For $g \in G$, the $g$-action on an object $\objx$ \cref{dgmn_iopg_obj} is given by
\begin{equation}\label{dgmn_iopg_gact}
g \big(\psi; \ang{\iphi_j}_{j \in \ufsn}\big) 
= \big(g\psi\ginv; \bang{g \iphi_{\ginv j}^{\tau_{\ginv}}\ginv}_{j \in \ufsn} \big).
\end{equation}
\begin{itemize}
\item $g\psi\ginv \cn \mal \to \nbeta$ is the conjugation $g$-action on $\psi$ \cref{gpsi}.
\item The object 
\[\iphi_{\ginv j}^{\tau_{\ginv}} \in \Iopg(|(g\psi\ginv)^{-1}j|)\] 
is the symmetric group action \cref{inopg_symmetry} on the object 
\[\iphi_{\ginv j} \in \Iopg(|(g\psi)^{-1} j|)\] 
for the permutation $\tau_{\ginv}$ \cref{tauginv} induced by the $\ginv$-action on $\mal$:
\[\ufsm \supset (g\psi\ginv)^{-1}j \fto[\iso]{\tau_{\ginv}} (g\psi)^{-1}j \subset \ufsm.\]
The object 
\[g \iphi_{\ginv j}^{\tau_{\ginv}}\ginv \in \Iopg(|(g\psi\ginv)^{-1}j|)\] 
is the conjugation $g$-action \cref{inopg_gact} on $\iphi_{\ginv j}^{\tau_{\ginv}}$.  Using the order-preserving bijections
\begin{equation}\label{ord_bij}
\begin{split}
\ufs{|s|} = \{1,2,\ldots,|s|\} \fto[\iso]{\obij_s} s  \subseteq \ufsm,
\end{split}
\end{equation}
$g \iphi_{\ginv j}^{\tau_{\ginv}}\ginv$ is given explicitly by
\begin{equation}\label{dgiopg_g}
\begin{split}
& \big(g \iphi_{\ginv j}^{\tau_{\ginv}}\ginv\big)(r,f) \\
&= g\big(\iphi_{\ginv j} \big(\obijinv_{(g\psi)^{-1}j} \ginv \obij_{(g\psi\ginv)^{-1} j} (r) , \ginv f\big)\big)
\end{split}
\end{equation}
for $r \in \ufs{|(g\psi\ginv)^{-1} j|}$ and  $f \in \omeg$.  The right-hand side of \cref{dgiopg_g} involves the following permutation.
\begin{equation}\label{obij_ginv_obijinv}
\begin{tikzpicture}[vcenter]
\def\v{-1.4}
\draw[0cell]
(0,0) node (a1) {\ufs{|(g\psi\ginv)^{-1} j|}}
(a1)++(0,\v) node (a2) {(g\psi\ginv)^{-1}j}
(a2)++(3,0) node (a3) {(g\psi)^{-1}j}
(a3)++(0,-\v) node (a4) {\ufs{|(g\psi)^{-1}j|}}
;
\draw[1cell=.85]
(a1) edge node[swap] {\obij_{(g\psi\ginv)^{-1}j}} (a2)
(a2) edge node {\ginv} (a3)
(a3) edge node {\obijinv_{(g\psi)^{-1}j}} (a4)
;
\end{tikzpicture}
\end{equation}
\end{itemize}
Using the conjugation $G$-action \cref{inopg_gmor} and the symmetric group action \cref{inopg_sym_mor} on morphisms in $\Iopg$, the $g$-action on a morphism $\objp$ \cref{dgmn_iopg_mor} is given by
\begin{equation}\label{dgmn_iopg_gmor}
\begin{split}
& g \big(\psi; \ang{[\ivphi_j,\iphi_j]}_{j \in \ufsn}\big) \\
&= \big(g\psi\ginv; \bang{[g \ivphi_{\ginv j}^{\tau_{\ginv}}\ginv,g \iphi_{\ginv j}^{\tau_{\ginv}}\ginv]}_{j \in \ufsn} \big).
\end{split}
\end{equation}
\item[Composition]
The composition pointed $G$-functor\label{not:DGcomp}
\[\DG(\mal,\nbeta) \sma \DG(\kdea,\mal) \fto{\comp} \DG(\kdea,\nbeta)\]
is given on objects by
\begin{equation}\label{dg_iopg_comp}
\big(\psi; \ang{\iphi_j}_{j \in \ufsn}\big) \comp \big(\phi; \ang{\ivphi_i}_{i \in \ufsm} \big) 
= \big(\psi\phi; \ang{\ga(\iphi_j; \ang{\ivphi_i}_{i \in \psiinv j})^{\tau^j_{\psi,\phi}}}_{j \in \ufsn} \big).
\end{equation}
It uses the operadic composition of objects in $\Iopg$ \cref{inopg_comp} and the shuffle \cref{tauj_psiphi}
\[(\psi\phi)^{-1}(j) \fto[\iso]{\tau^j_{\psi,\phi}} \coprod_{i \in \psiinv j} \phiinv i\]
determined by the inherited orderings
\[(\psi\phi)^{-1}(j) , \phiinv i \subseteq \ufsk \andspace \psiinv j \subseteq \ufsm.\]
The composition on morphisms is given by
\begin{equation}\label{dg_iopg_mor}
\begin{split}
& \big(\psi; \ang{[\iphi_j',\iphi_j]}_{j \in \ufsn}\big) \comp 
\big(\phi; \ang{[\ivphi_i',\ivphi_i]}_{i \in \ufsm} \big) \\
&= \big(\psi\phi; \bang{[\ga(\iphi_j'; \ang{\ivphi_i'}_{i \in \psiinv j})^{\tau^j_{\psi,\phi}} , \ga(\iphi_j; \ang{\ivphi_i}_{i \in \psiinv j})^{\tau^j_{\psi,\phi}}]}_{j \in \ufsn} \big)
\end{split}
\end{equation}
using the operadic composition of morphisms in $\Iopg$ \cref{inopg_comp_mor}.
\end{description}
The objects of the category $\DGCatg$ \pcref{def:dgcatg_iicat} are $\DGG$-categories \cref{dgcatg_obj}, meaning pointed $\Gcatst$-functors
\[(\DG,\ordz) \fto{X} (\Gcatst,\bone).\]
Morphisms in $\DGCatg$ are $\Pig$-strict $\Gcatst$-pseudotransformations \pcref{def:dgcatg_onecell,def:DGcoop,def:pig_strict}.
\end{explanation}

\begin{explanation}[The Functor $\Lg\Jscg$]\label{expl:LgJscg}
Up to a natural isomorphism, the composite 
\[\Parcat \fto{\Jscg} \FGCat \fto{\Lg} \FGCatg\]
of Schwede $J$-theory $\Jscg$ \cref{Jscg} and the equivalence $\Lg$ \pcref{thm:fgcat_fgcatg_iieq} sends a parsummable category $(\C,\psum,\pzero)$ \pcref{expl:parcat} to the pointed $G$-functor
\begin{equation}\label{LgJscgC}
(\FG,\ordz) \fto{\Lg\Jscg\C} (\Catgst,\bone)
\end{equation}
given as follows.   Recall that a parsummable category $(\C,\psum,\pzero)$ is strictly associative, unital, and symmetric \cref{parcat_axioms}.
\begin{description}
\item[Pointed $G$-categories]
For a pointed finite $G$-set $\nbeta$ \cref{ordn_be}, by \cref{LX_nbeta,Xjin_giso,JscgCn}, the underlying pointed category 
\begin{equation}\label{LgJscg_nbeta}
(\Lg\Jscg\C)\nbeta = \Cboxn
\end{equation}
is the full subcategory of $\C^n$ consisting of disjointly supported $n$-tuples.  Its basepoint is the $n$-tuple $\ang{\pzero}_{j \in \ufsn}$.  For each $g \in G$, by \cref{Xordnbe_g,jscg_gaction,jscg_psist}, the $g$-action on $(\Lg\Jscg\C)\nbeta$ sends an $n$-tuple $\ang{x_j}_{j \in \ufsn} \in \Cboxn$ to the $n$-tuple
\begin{equation}\label{LgJscg_gact}
g\ang{x_j}_{j \in \ufsn} = \ang{\gomest x_{\ginv j}}_{j \in \ufsn} \in \Cboxn.
\end{equation}
Here, $\gomest \cn \C \to \C$ is the $\schm$-action functor for the bijection $g^\ome \cn \ome \to \ome$ \cref{rep_g_rep}, and $\ginv j$ means $(\be g)^{-1}j$.
\item[Pointed functors]
For a pointed morphism $\psi \cn \mal \to \nbeta$ between pointed finite $G$-sets, by \cref{Xjin_giso,jscg_psist}, the pointed functor
\begin{equation}\label{LgJscgC_psi}
(\Lg\Jscg\C)\mal = \Cboxm \fto{(\Lg\Jscg\C)\psi \,=\, \psi_*} (\Lg\Jscg\C)\nbeta = \Cboxn
\end{equation}
sends an $m$-tuple $\ang{x_i}_{i \in \ufsm} \in \Cboxm$ to the $n$-tuple
\begin{equation}\label{LgJscgC_psix}
\psi_*\ang{x_i}_{i \in \ufsm} = \bang{\,\txsum_{i \in \psiinv j} x_i}_{j \in \ufsn} \in \Cboxn.
\end{equation}
\end{description}
For a parsummable functor $\fun \cn \C \to \D$, by \cref{Xjin_giso,Jscgf}, the $G$-natural transformation
\begin{equation}\label{LgJscgf}
\begin{tikzpicture}[vcenter]
\def\t{27}
\draw[0cell]
(0,0) node (a1) {\phantom{\Fsk}}
(a1)++(-.05,0) node (a1') {\FG}
(a1)++(2.2,0) node (a2) {\phantom{\Fsk}}
(a2)++(.2,-.02) node (a2') {\Catgst}
;
\draw[1cell=.85]
(a1) edge[bend left=\t] node {\Lg\Jscg\C} (a2)
(a1) edge[bend right=\t] node[swap] {\Lg\Jscg\D} (a2)
;
\draw[2cell=.9]
node[between=a1 and a2 at .35, rotate=-90, 2label={above,\Lg\Jscg\fun}] {\Rightarrow}
;
\end{tikzpicture}
\end{equation}
sends each pointed finite $G$-set $\nbeta$ to the pointed $G$-functor
\begin{equation}\label{LgJscgf_nbeta}
(\Lg\Jscg\C)\nbeta = \Cboxn \fto{(\Lg\Jscg\fun)_{\nbeta} \,=\, \fboxn} (\Lg\Jscg\D)\nbeta = \Dboxn
\end{equation}
given entrywise by $\fun$.
\end{explanation}

\begin{explanation}[The Functor $\gxist\Lg\Jscg$]\label{expl:cglg_dom}
 Using \cref{expl:dgcatg,expl:gxist,expl:parcat,expl:LgJscg,expl:dgiopg}, the composite functor
\[\Parcat \fto{\Jscg} \FGCat \fto{\Lg} \FGCatg \fto{\gxist} \dgcatg \subseteq \DGCatg\]
is given as follows.
\begin{description}
\item[Objects]
For a parsummable category $(\C,\psum,\pzero)$, the $\DGG$-category \cref{dgcatg_obj}
\[(\DG,\ordz) \fto{\gxist\Lg\Jscg\C} (\Gcatst,\bone)\]
sends an object $\nbeta \in \DG$ to the pointed $G$-category
\begin{equation}\label{cglgdom_nbeta}
(\gxist\Lg\Jscg\C)\nbeta = (\Lg\Jscg\C)\nbeta = \Cboxn
\end{equation}
of disjointly supported $n$-tuples in $\C^n$ \pcref{def:bxtimes}.  Its basepoint is $\ang{\pzero}_{j \in \ufsn}$.  The $G$-action is given by
\begin{equation}\label{cglgdom_gact}
g\ang{a_j}_{j \in \ufsn} = \ang{\gomest a_{\ginv j}}_{j \in \ufsn} \in \Cboxn
\end{equation}
for $g \in G$ and $\ang{a_j}_{j \in \ufsn} \in \Cboxn$, where $\gomest \cn \C \to \C$ is the $\schm$-action functor for the bijection $\gome \cn \ome \to \ome$ \cref{rep_g_rep}.  For an object $\objx = (\psi; \ang{\iphi_j}_{j \in \ufsn})$ in $\DG(\mal,\nbeta)$ \cref{dgmn_iopg_obj}, the pointed functor
\begin{equation}\label{cglgdom_psi}
(\gxist\Lg\Jscg\C)\mal = \Cboxm \fto{(\gxist\Lg\Jscg\C)\objx \,=\, \psi_*} 
(\gxist\Lg\Jscg\C)\nbeta = \Cboxn
\end{equation}
sends an $m$-tuple $\ang{a_i}_{i \in \ufsm} \in \Cboxm$ to the $n$-tuple
\begin{equation}\label{cglgdom_psix}
\psi_*\ang{a_i}_{i \in \ufsm} = \bang{\,\txsum_{i \in \psiinv j} a_i}_{j \in \ufsn} \in \Cboxn.
\end{equation}
Each morphism $(\psi; \ang{[\ivphi_j,\iphi_j]}_{j \in \ufsn})$ in $\DG(\mal,\nbeta)$ \cref{dgmn_iopg_mor} is sent by $\gxist\Lg\Jscg\C$ to the identity natural transformation:
\begin{equation}\label{cglgdom_twocell}
(\gxist\Lg\Jscg\C)(\psi; \ang{[\ivphi_j,\iphi_j]}_{j \in \ufsn}) = 1_{\psi_*}.
\end{equation}
\item[Morphisms]
For a parsummable functor $\fun \cn \C \to \D$, the $\Gcatst$-natural transformation
\begin{equation}\label{cglgdom_f}
\begin{tikzpicture}[vcenter]
\def\t{27}
\draw[0cell]
(0,0) node (a1) {\phantom{\Fsk}}
(a1)++(-.05,0) node (a1') {\DG}
(a1)++(2.5,0) node (a2) {\phantom{\Fsk}}
(a2)++(.32,-.02) node (a2') {\Gcatst}
;
\draw[1cell=.8]
(a1) edge[bend left=\t] node {\gxist\Lg\Jscg\C} (a2)
(a1) edge[bend right=\t] node[swap] {\gxist\Lg\Jscg\D} (a2)
;
\draw[2cell=.85]
node[between=a1 and a2 at .32, rotate=-90, 2label={above,\gxist\Lg\Jscg\fun}] {\Rightarrow}
;
\end{tikzpicture}
\end{equation}
sends each object $\nbeta \in \DG$ to the pointed $G$-functor
\begin{equation}\label{cglgdom_fn}
(\gxist\Lg\Jscg\C)\nbeta = \Cboxn \fto{(\gxist\Lg\Jscg\fun)_{\nbeta} \,=\, \fboxn} 
(\gxist\Lg\Jscg\D)\nbeta = \Dboxn
\end{equation}
given entrywise by $\fun$.\defmark
\end{description}
\end{explanation}

\begin{explanation}[The Functor $\Rg\Ig$]\label{expl:cglg_cod}
Using \cref{expl:dgcatg,def:rg_obj,def:rg_onecell,expl:parcat,def:Ig,expl:dgiopg}, the composite functor
\[\Parcat \fto{\Ig} \Algiopg \fto{\Rg} \dgcatg \subseteq \DGCatg\]
is given as follows.
\begin{description}
\item[Objects]
For a parsummable category $(\C,\psum,\pzero)$, the $\DGG$-category \cref{dgcatg_obj}
\[(\DG,\ordz) \fto{\Rg\Ig\C} (\Gcatst,\bone)\]
sends an object $\nbeta \in \DG$ to the $\nbeta$-twisted product \pcref{def:proCnbe}
\begin{equation}\label{ricn}
(\Rg\Ig\C)\nbeta = \Comegnbeta
\end{equation}
of the pointed $G$-category $\Comeg$ (\cref{def:sch2.21} \cref{def:sch2.21_i}).  
\begin{description}
\item[Pointed category]
The underlying pointed category of $\Comegnbeta$ is $\C^n$ with basepoint $\ang{\pzero}_{j \in \ufsn}$.  
\item[$G$-action]
The $G$-action is given by
\begin{equation}\label{ricn_gact}
g\ang{a_j}_{j \in \ufsn} = \ang{\gomest a_{\ginv j}}_{j \in \ufsn} \in \Comegnbeta
\end{equation}
for $g \in G$ and $\ang{a_j}_{j \in \ufsn} \in \Comegnbeta$, where $\gomest \cn \C \to \C$ is the $\schm$-action functor for the bijection $\gome \cn \ome \to \ome$ \cref{rep_g_rep}.
\item[Pointed functors]
For each object $\objx = (\psi; \ang{\iphi_j}_{j \in \ufsn})$ in $\DG(\mal,\nbeta)$ \cref{dgmn_iopg_obj}, the pointed functor
\begin{equation}\label{ricx}
(\Rg\Ig\C)\mal = \Comegmal \fto{(\Rg\Ig\C)\objx \,=\, \objx_*}
(\Rg\Ig\C)\nbeta = \Comegnbeta
\end{equation}
sends $\obja = \ang{a_i}_{i \in \ufsm} \in \Comegmal$ to
\begin{equation}\label{ricxa}
\begin{split}
\objx_*\obja
&= \bang{\gacomeg\big(\iphi_j; \ang{a_i}_{i \in \psiinv j} \big)}_{j \in \ufsn} \\
&= \bang{\,\txsum_{\ell \in \ufs{|\psiinv j|}} \, \iphijlomest a_{\obij_{\psiinv j} \ell}}_{j \in \ufsn} 
\end{split}
\end{equation}
in $\Comegnbeta$.  It uses the $\Iopg$-action on $\Comeg$ \cref{iphi_comeg} and the order-preserving bijection \cref{ord_bij}
\begin{equation}\label{psiinvj_ord}
\{1,2,\ldots,|\psiinv j|\} = \ufs{|\psiinv j|} \fto[\iso]{\obij_{\psiinv j}} \psiinv j \subseteq \ufsm.
\end{equation}
To explain \cref{ricxa} in more detail, recall the injections
\[\{1,2,\ldots,|\psiinv j|\} \itimes \omeg = \ufs{|\psiinv j|} \itimes \omeg \fto{\iphi_j} \omeg\]
for $j \in \ufsn$ and
\[\iphijl = \iphi_j(\ell,-) \cn \omeg \to \omeg\]
for $\ell \in \ufs{|\psiinv j|}$.  The injections
\begin{equation}\label{iphijlome}
\begin{tikzpicture}[vcenter]
\def\h{1.8} \def\g{1.5} \def\u{.7} \def\d{.05}
\draw[0cell]
(0,0) node (a1) {\ome}
(a1)++(\h,0) node (a2) {\phantom{\omeg}}
(a2)++(0,\d) node (a2') {\omeg}
(a2)++(\h,0) node (a3) {\phantom{\omeg}}
(a3)++(0,\d) node (a3') {\omeg}
(a3)++(\h,0) node (a4) {\ome}
;
\draw[1cell=.9]
(a1) edge node {\rep_{\omeg}} (a2)
(a2) edge node {\iphijl} (a3)
(a3) edge node {\rep_{\omeg}^{-1}} (a4)
(a1) [rounded corners=2pt] |- ($(a2)+(0,\u)$) -- node {\iphijlome} ($(a3)+(0,\u)$) -| (a4)
;
\end{tikzpicture}
\end{equation}
for $\ell \in \ufs{|\psiinv j|}$ have disjoint images, and 
\begin{equation}\label{iphijlomest}
\C \fto{\iphijlomest} \C 
\end{equation}
denotes the $\iphijlome$-action functor on the $\schm$-category $\C$ \cref{u_action}.  The sum in \cref{ricxa} is well defined by \cref{sch2.13} \cref{sch2.13_iii}, \cref{supp_xsumy}, and \cref{iphijlome}.
\item[Pointed natural transformations]
The $\DGG$-category $\Rg\Ig\C$ sends each morphism
\[\objx = \big(\psi; \ang{\iphi_j}_{j \in \ufsn}\big) 
\fto{\objp \,=\, (\psi; \ang{[\ivphi_j,\iphi_j]}_{j \in \ufsn})}
\objy = \big(\psi; \ang{\ivphi_j}_{j \in \ufsn}\big)\]
in $\DG(\mal,\nbeta)$ \cref{dgmn_iopg_mor} to the pointed natural isomorphism
\begin{equation}\label{ricp}

\end{equation}
whose component at an object $\obja = \ang{a_i}_{i \in \ufsm} \in \Comegmal$ is the following isomorphism in $\Comegnbeta$, with $\obij$ denoting $\obij_{\psiinv j}$ \cref{psiinvj_ord}.
\begin{equation}\label{ricpa}
%
\end{equation}
For each $\ell \in \ufs{|\psiinv j|}$,
\[\iphijlomest \fto{[\ivphijlome,\iphijlome]_{\mstar}} \ivphijlomest\]
is the $\schm$-action natural isomorphism on $\C$ \cref{vu_star} for the isomorphism
\[\iphijlome \fto{[\ivphijlome, \iphijlome]} \ivphijlome \inspace \schm.\]
\end{description}
\item[Morphisms]
For a parsummable functor $\fun \cn \C \to \D$, the $\Gcatst$-natural transformation
\begin{equation}\label{cglgcod_f}
%
\end{equation}
sends each object $\nbeta \in \DG$ to the pointed $G$-functor
\begin{equation}\label{cglgcod_fn}
(\Rg\Ig\C)\nbeta = \Comegnbeta \fto{(\Rg\Ig\fun)_{\nbeta} \,=\, \fun^n} 
(\Rg\Ig\D)\nbeta = \Domegnbeta
\end{equation}
given entrywise by $\fun$.\defmark
\end{description}
\end{explanation}

\subsection*{The Natural Transformation $\cglg$}
Using \cref{expl:dgiopg,expl:cglg_dom,expl:cglg_cod}, \cref{def:cglg} constructs the components of $\cglg$ at parsummable categories \pcref{def:parcat}.  These components of $\cglg$ compare the first step $\Jscg$ of Schwede $K$-theory and the first step $\Rg$ of GMMO $K$-theory.

\begin{definition}\label{def:cglg}
For a finite group $G$ and a parsummable category $(\C,\psum,\pzero)$, the $\Gcatst$-pseudotransformation \pcref{def:dgcatg_onecell}
\begin{equation}\label{cglgc}
\begin{tikzpicture}[vcenter]
\def\t{25}
\draw[0cell]
(0,0) node (a1) {\phantom{\Fsk}}
(a1)++(-.05,0) node (a1') {\DG}
(a1)++(2.1,0) node (a2) {\phantom{\Fsk}}
(a2)++(.32,-.02) node (a2') {\Gcatst}
;
\draw[1cell=.85]
(a1) edge[bend left=\t] node {\gxist\Lg\Jscg\C} (a2)
(a1) edge[bend right=\t] node[swap] {\Rg\Ig\C} (a2)
;
\draw[2cell=.9]
node[between=a1 and a2 at .4, rotate=-90, 2label={above,\cglgc}] {\Rightarrow}
;
\end{tikzpicture}
\end{equation}
is defined as follows.
\begin{description}
\item[$G$-functors]
Using \cref{cglgdom_nbeta,ricn}, for each object $\nbeta \in \DG$, the $\nbeta$-component pointed $G$-functor
\begin{equation}\label{cglgcn}
\begin{tikzpicture}[vcenter]
\def\v{-1}
\draw[0cell]
(0,0) node (a11) {(\gxist\Lg\Jscg\C)\nbeta}
(a11)++(3.3,0) node (a12) {(\Rg\Ig\C)\nbeta}
(a11)++(0,\v) node (a21) {\Cboxn}
(a12)++(0,\v) node (a22) {\Comegnbeta}
;
\draw[1cell=.85]
(a11) edge[equal] (a21)
(a12) edge[equal,shorten >=-.3ex] (a22)
(a11) edge node {\cglgcnbeta} (a12)
(a21) edge[right hook->] (a22)
;
\end{tikzpicture}
\end{equation}
is defined as the full subcategory inclusion.  Recall that the underlying pointed category of $\Comegnbeta$ is $\C^n$.  The basepoint in both $\Cboxn$ and $\Comegnbeta$ is $\ang{\pzero}_{j \in \ufsn}$.  The pointed functor $\cglgcnbeta$ is $G$-equivariant by \cref{cglgdom_gact,ricn_gact}.
\item[Natural isomorphisms]
Using \cref{cglgdom_psix,ricxa,psiinvj_ord,cglgcn}, for each object $\objx = (\psi; \ang{\iphi_j}_{j\in \ufsn})$ in $\DG(\mal,\nbeta)$ \cref{dgmn_iopg_obj}, the $\objx$-component pointed natural isomorphism
\begin{equation}\label{cglgcx}
\begin{tikzpicture}[vcenter]
\def\v{-1.4}
\draw[0cell]
(0,0) node (a11) {\Cboxm}
(a11)++(2.3,0) node (a12) {\Comegmal}
(a11)++(0,\v) node (a21) {\Cboxn}
(a12)++(0,\v) node (a22) {\Comegnbeta}
;
\draw[1cell=.85]
(a11) edge node {\cglgcmal} (a12)
(a12) edge[shorten >=-.3ex] node {\objx_*} (a22)
(a11) edge node[swap] {\psi_*} (a21)
(a21) edge node[swap] {\cglgcnbeta} (a22)
;
\draw[2cell]
node[between=a12 and a21 at .5, shift={(0,-.1)}, rotate=-115, 2labelw={below,\cglgcx,0pt}] {\Rightarrow}
;
\end{tikzpicture}
\end{equation}
sends an object $\obja = \ang{a_i}_{i \in \ufsm} \in \Cboxm$ to the following isomorphism in $\Comegnbeta$, where $\obij = \obij_{\psiinv j}$ \cref{psiinvj_ord}.
\begin{equation}\label{cglgcxa}
\begin{tikzpicture}[vcenter]
\def\v{-1.5} \def\d{-1.1} \def\e{-.25} \def\h{1.5} \def\g{1.5}
\draw[0cell]
(0,0) node (a11) {\objx_*\cglgcmal\obja}
(a11)++(\g,0) node (a11') {\objx_*\obja}
(a11')++(2.7,0) node (a12) {\bang{\,\txsum_{\ell \in \ufs{|\psiinv j|}} \, \iphijlomest a_{\obij\ell}}_{j \in \ufsn}}
(a11)++(0,\v) node (a21) {\cglgcnbeta\psi_*\obja}
(a11')++(0,\v) node (a21') {\psi_*\obja}
(a12)++(\e,\v) node (a22) {\bang{\,\txsum_{\ell \in \ufs{|\psiinv j|}} \, a_{\obij\ell}}_{j \in \ufsn}}
(a11')++(\h,0) node (a12') {\phantom{\sum}}
(a21')++(\h,0) node (a22') {\phantom{\sum}}
;
\draw[1cell=.8]
(a11) edge[equal] (a11')
(a11') edge[equal] (a12)
(a21) edge[equal] (a21')
(a21') edge[equal] (a22)
(a11) edge node[swap] {\cglgcxa} (a21)
(a12') edge node {\bang{\,\txsum_{\ell \in \ufs{|\psiinv j|}} \, [1_\ome,\iphijlome]_{\mstar}^{a_{\obij\ell}}}_{j \in \ufsn}} (a22')
;
\end{tikzpicture}
\end{equation}
For each $j \in \ufsn$ and $\ell \in \ufs{|\psiinv j|}$,
\[\iphijlomest \fto{[1_\ome,\iphijlome]_{\mstar}} (1_\ome)_\mstar = 1_\C\]
is the $\schm$-action natural isomorphism on $\C$ \cref{vu_star} for the isomorphism
\[\iphijlome \fto{[1_\ome,\iphijlome]} 1_\ome \inspace \schm.\]
The naturality of $\cglgcx$ in $\obja \in \Cboxm$ follows from the naturality of each $[1_\ome,\iphijlome]_{\mstar}$ and the functoriality of the sum on $\C$.  This natural isomorphism is pointed because the unit $\pzero \in \C$ is $\schm$-fixed and a strict two-sided unit for the sum.
\end{description}
This finishes the definition of $\cglgc$.  \cref{cglgc_welldef,cglg_natural} prove that the assignment $\C \mapsto \cglgc$ defines a natural transformation $\cglg$.
\end{definition}

Recall that a morphism in $\DGCatg$ is a $\Pig$-strict $\Gcatst$-pseudotransformation \pcref{def:dgcatg_onecell,def:pig_strict,def:dgcatg_iicat}.

\begin{lemma}\label{cglgc_welldef}
For each finite group $G$ and parsummable category $(\C,\psum,\pzero)$, $\cglgc$ \cref{cglgc} is a $\Pig$-strict $\Gcatst$-pseudotransformation.
\end{lemma}

\begin{proof}
We verify the axioms in \cref{def:dgcatg_onecell,def:pig_strict} for $\cglgc$, which is defined in \crefrange{cglgcn}{cglgcxa}.
\begin{description}
\item[Basepoint]
The basepoint axiom \cref{tha_objzero} states the equality
\[\cglgc_{\objzero,\obja} = 1_{\ang{\pzero}_{j \in \ufsn}} \inspace \Comegnbeta\]
for the basepoint \cref{dgiopg_bp}
\[\objzero = (0; \ang{*}_{j \in \ufsn}) \inspace \DG(\mal,\nbeta)\]
and each object $\obja \in \Cboxm$.  For the morphism $\psi = 0 \cn \mal \to \nbeta$, each subset $\psiinv j \subseteq \ufsm$ is empty.  By definition, an empty sum is the unit $\pzero$ or its identity morphism.
\item[Unity]
Recall that $\cglgcmal$ \cref{cglgcn} is the inclusion $\Cboxm \to \C^m$ for each object $\mal \in \DG$.  The unity axiom \cref{tha_onex} states the equality
\[\cglgc_{\objone_{\mal},\obja} = 1_{\obja}\]
for the identity 1-cell \cref{dgiopg_idonecell}
\[\objone_{\mal} = (1_{\mal}; \ang{1_{\omeg}}_{i \in \ufsm}) \inspace \DG(\mal,\mal)\]
and each object $\obja \in \Cboxm$.  The unity axiom holds because $(1_\ome)_\mstar = 1_\C$ and
\[1_\ome \fto{[1_\ome,1_\ome]} 1_\ome\]
is the identity morphism $1_{(1_\ome)}$ in the injection category $\schm$ \pcref{def:mcat}.
\item[Naturality]
By \cref{cglgdom_twocell,ricpa,cglgcxa}, the naturality axiom \cref{thaobjx_nat} for $\cglgc$ states that, for each morphism $\objp \cn \objx \to \objy$ in $\DG(\mal,\nbeta)$ \cref{dgmn_iopg_mor} and each object $\obja = \ang{a_i}_{i \in \ufsm} \in \Cboxm$, the following diagram in $\Comegnbeta$ commutes, where $\sum_\ell = \sum_{\ell \in \ufs{|\psiinv j|}}$, $\ang{\Cdots}_j = \ang{\Cdots}_{j \in \ufsn}$, and $\obij = \obij_{\psiinv j}$ \cref{psiinvj_ord}.
\begin{equation}\label{cglgcxa_diag}
\begin{tikzpicture}[vcenter]
\def\h{3.3} \def\v{-1} \def\c{1} \def\d{-1.1} \def\s{.8}
\draw[0cell]
(0,0) node (a1) {\objx_*\obja = \ang{\,\txsum_{\ell} \iphijlomest a_{\obij\ell}}_{j}}
(a1)++(\h,\v) node (a2) {\objy_*\obja = \ang{\,\txsum_{\ell} \ivphijlomest a_{\obij\ell}}_{j}}
(a1)++(0,2*\v) node (a3) {\phantom{\psi_*\obja = \ang{\,\txsum_{\ell} a_{\obij\ell}}_{j}}}
(a3)++(.25,0) node (a3'') {\psi_*\obja = \ang{\,\txsum_{\ell} a_{\obij\ell}}_{j}}
(a1)++(\c,0) node (a1') {\phantom{\sum}}
(a2)++(\d,0) node (a2') {\phantom{\sum}}
(a3)++(\c,0) node (a3') {\phantom{\sum}}
;
\draw[1cell=\s]
(a1') edge node[swap] {\cglgcxa \,=\, \ang{\,\txsum_\ell\, [1_\ome,\iphijlome]_\mstar^{a_{\obij\ell}}}_{j}} (a3')
(a1) [rounded corners=2pt] -| node[pos=.7] {\objp_*^\obja \,=\, \ang{\,\txsum_\ell\, [\ivphijlome,\iphijlome]_\mstar^{a_{\obij\ell}}}_j} (a2')
;
\draw[1cell=\s]
(a2') [rounded corners=2pt] |- node[pos=.3] {\cglgcya \,=\, \ang{\,\txsum_\ell\, [1_\ome,\ivphijlome]_\mstar^{a_{\obij\ell}}}_j} (a3'')
;
\end{tikzpicture}
\end{equation}
The previous diagram commutes by the functoriality of the sum on $\C$ and the morphism equality
\[[1_\ome,\iphijlome] = [1_\ome,\ivphijlome] \comp [\ivphijlome,\iphijlome] \cn \iphijlome \to 1_\ome\]
in the injection category $\schm$.
\item[$G$-equivariance]
For $g \in G$, an object $\objx = (\psi; \ang{\iphi_j}_{j \in \ufsn})$ in $\DG(\mal,\nbeta)$, and an object $\obja = \ang{a_i}_{i \in \ufsm}$ in $\Cboxm$, by \cref{dgmn_iopg_gact,cglgdom_gact}, there are objects
\begin{equation}\label{gx_ginva}
\begin{split}
g\objx &= \big(g\psi\ginv; \bang{g\iphi_{\ginv j}^{\tau_{\ginv}}\ginv}_{j \in \ufsn}\big) \in \DG(\mal,\nbeta) \andspace\\
\ginv\obja &= \bang{\ginvomest a_{gi}}_{i \in \ufsm} \in \Cboxm.
\end{split}
\end{equation}
Using \cref{cglgcxa}, \cref{gx_ginva}, and the abbreviations
\[\ang{\Cdots}_j = \ang{\Cdots}_{j \in \ufsn} \andspace 
\ell_{\ginv} = \obij^{-1}_{(g\psi)^{-1} j} \ginv \obij_{(g\psi\ginv)^{-1}j} \ell,\]
the following morphism equalities in $\Comegnbeta$ prove the $G$-equivariance axiom \cref{thaobjx_gequiv} for $\cglgc$.
\def\hsk{\hspace{-4em}}
\[\scalebox{.85}{$
\begin{aligned}
&\cglgc_{g\objx,\obja} &&\\
&= \bang{\, \txsum_{\ell \in \ufs{|(g\psi\ginv)^{-1} j|}} \, [1_\ome, (g\iphi_{\ginv j}^{\tau_{\ginv}} \ginv)^{\ell,\ome} ]_\mstar^{a_{\obij_{(g\psi\ginv)^{-1} j} \ell}}}_j &&\\
&= \bang{\, \txsum_{\ell \in \ufs{|(g\psi\ginv)^{-1} j|}} \, [1_\ome, (g\iphi_{\ginv j}^{\ell_{\ginv}} \ginv)^\ome]_\mstar^{a_{\obij_{(g\psi\ginv)^{-1} j} \ell}}}_j 
&&\hsk \text{by \cref{dgiopg_g}}\\
&= \bang{\, \txsum_{\ell \in \ufs{|(g\psi\ginv)^{-1} j|}} \, [\gome\ginvome, (\gome)(\iphi_{\ginv j}^{\ell_{\ginv}, \ome}) \ginvome]_\mstar^{a_{\obij_{(g\psi\ginv)^{-1} j} \ell}}}_j &&\\
&&&\hspace{-12em} \text{by \cref{gome_mstar}, \cref{iphijlome}, and $\repomeg\repomeginv = 1_{\omeg}$} &&\\
&= \bang{\, \txsum_{\ell \in \ufs{|(g\psi)^{-1} j|}} \, [\gome\ginvome, \gome\iphi_{\ginv j}^{\ell,\ome} \ginvome]_\mstar^{a_{g\obij_{(g\psi)^{-1} j} \ell}}}_j 
&&\hsk \text{by \cref{parcat_axioms,obij_ginv_obijinv}}\\
&= \bang{\, \txsum_{\ell \in \ufs{|(g\psi)^{-1} j|}} \, \gomest[1_\ome, \iphi_{\ginv j}^{\ell,\ome}]_\mstar^{\ginvomest a_{g\obij_{(g\psi)^{-1} j} \ell}}}_j 
&&\hsk \text{by \cref{mcat_axioms}}\\
&= \bang{\, \gomest \txsum_{\ell \in \ufs{|(g\psi)^{-1} j|}} \, [1_\ome, \iphi_{\ginv j}^{\ell,\ome}]_\mstar^{(\ginv\obja)_{\obij_{(g\psi)^{-1} j} \ell}}}_j 
&&\hsk \text{by $\schm$-equivariance of $\psum$}\\
&= g\bang{\, \txsum_{\ell \in \ufs{|\psiinv j|}} \, [1_\ome, \iphijlome]_\mstar^{(\ginv\obja)_{\obij_{\psiinv j} \ell}}}_j 
&&\hsk \text{by \cref{ricn_gact}}\\
&= g\cglgc_{\objx,\ginv\obja} 
\end{aligned}$}\]
\item[Compositionality]
The axiom \cref{tha_xv} for $\cglgc$ states that, for objects \cref{dgmn_iopg_obj}
\[\begin{split}
\objx &= \big(\psi; \ang{\iphi_j}_{j \in \ufsn} \big) \in \DG(\mal,\nbeta) \andspace\\
\objy &= \big(\phi; \ang{\ivphi_i}_{i \in \ufsm} \big) \in \DG(\kdea,\mal)
\end{split}\]
with composite \cref{dg_iopg_comp}
\[\objx\objy = \big(\psi\phi; \ang{\ga(\iphi_j; \ang{\ivphi_i}_{i \in \psiinv j})^{\tau^j_{\psi,\phi}}}_{j \in \ufsn} \big) \in \DG(\kdea,\nbeta),\]
the following two pasting diagrams of natural transformations are equal.
\begin{equation}\label{cglgc_compax}
\begin{tikzpicture}[vcenter]
\def\h{2} \def\v{-1.3} \def\s{.9} \def\q{.75} \def\p{.85} \def\r{-120}
\def\boundary{
\draw[0cell=\s]
(0,0) node (a11) {\Cboxk}
(a11)++(\h,0) node (a12) {\Comegkdea}
(a11)++(0,2*\v) node (a31) {\Cboxn}
(a12)++(0,2*\v) node (a32) {\Comegnbeta}
;
\draw[1cell=\q]
(a11) edge node {\cglgc_{\kdea}} (a12)
(a31) edge node[swap] {\cglgc_{\nbeta}} (a32)
;}
\boundary
\draw[0cell=\s]
(a11)++(0,\v) node (a21) {\Cboxm}
(a12)++(0,\v) node (a22) {\Comegmal}
;
\draw[1cell=\q]
(a21) edge node {\cglgc_{\mal}} (a22)
(a11) edge node[swap] {\phi_*} (a21)
(a21) edge node[swap] {\psi_*} (a31)
(a12) edge node {\objy_*} (a22)
(a22) edge node {\objx_*} (a32)
;
\draw[2cell=\p]
node[between=a12 and a21 at .5, shift={(0,.1)}, rotate=\r, 2labelw={below,\cglgcy,-1pt}] {\Rightarrow}
node[between=a22 and a31 at .5, shift={(0,-.15)}, rotate=\r, 2labelw={below,\cglgcx,0pt}] {\Rightarrow}
;
\begin{scope}[shift={(2*\h,0)}]
\boundary
\draw[1cell=\q]
(a11) edge node[swap] {(\psi\phi)_*} (a31)
(a12) edge node {(\objx\objy)_*} (a32)
;
\draw[2cell=\p]
node[between=a12 and a31 at .5, shift={(.1*\h,0)}, rotate=\r, 2labelalt={below,\cglgcxy}] {\Rightarrow}
;
\end{scope}
\end{tikzpicture}
\end{equation}
The equality of the pasting diagrams in \cref{cglgc_compax} means that, for each object $\obja = \ang{a_d}_{d\in \ufsk}$ in $\Cboxk$, there is a morphism equality
\begin{equation}\label{cglgc_compeq}
\big(\cglgc_{\objx,\phi_*\obja}\big) \big(\objx_*\cglgcya\big) = \cglgcxya \inspace \Comegnbeta.
\end{equation}
To prove \cref{cglgc_compeq}, recall that the first entries of $\objy$ and $\objx$ are pointed functions
\[\kdea \fto{\phi} \mal \fto{\psi} \nbeta.\]
Using the notation in \cref{ord_bij}, for $j \in \ufsn$ and $i \in \ufsm$, there are order-preserving bijections as follows.
\[\begin{split}
\ufs{|\psiinv j|} = \{1,2,\ldots,|\psiinv j|\} & \fto[\iso]{\obij_{\psiinv j}} \psiinv j \subseteq \ufsm\\
\ufs{|\phiinv i|} =\{1,2,\ldots,|\phiinv i|\} & \fto[\iso]{\obij_{\phiinv i}} \phiinv i \subseteq \ufsk\\
\ufs{|\psiphiinv j|} = \{1,2,\ldots,|\psiphiinv j|\} & \fto[\iso]{\obij_{\psiphiinv j}} \psiphiinv j \subseteq \ufsk\\
\end{split}\]
In the rest of this proof, $\obij_s(-)$ is also denoted by $(-)'$.  For example, each $t \in \ufs{|\phiinv i|}$ yields the element
\begin{equation}\label{t_prime}
t' = \obij_{\phiinv i} t \in \phiinv i \subseteq \ufsk.
\end{equation}
Similarly, each $\ell \in \ufs{|\psiinv j|}$ yields the element
\begin{equation}\label{ell_prime}
\ell' = \obij_{\psiinv j}\ell \in \psiinv j \subseteq \ufsm,
\end{equation}
the order-preserving bijection
\begin{equation}\label{obij_phiinv}
\ufs{|\phiinv \ell'|} =\{1,2,\ldots,|\phiinv \ell'|\} 
\fto[\iso]{\obij_{\phiinv \ell'}} \phiinv \ell' \subseteq (\psi\phi)^{-1} j \subseteq \ufsk,
\end{equation}
and the elements
\begin{equation}\label{r_prime}
r' = \obij_{\phiinv \ell'} r \in \phiinv \ell' \forspace r \in \ufs{|\phiinv \ell'|}.
\end{equation}
Using \crefrange{t_prime}{r_prime}, we compute the left-hand side of \cref{cglgc_compeq}.

For $j \in \ufsn$, the $j$-th entry of the leftmost morphism in \cref{cglgc_compeq} is the following morphism in $\C$.
\begin{equation}\label{cglgc_comp_i}
\scalebox{.9}{$
\begin{aligned}
& \big(\cglgc_{\objx,\phi_*\obja}\big)_j &&\\
&= \txsum_{\ell \in \ufs{|\psiinv j|}} \, [1_\ome, \iphijlome]_\mstar^{(\phi_*\obja)_{\ell'}} && \text{by \cref{cglgcxa,ell_prime}}\\
&= \txsum_{\ell \in \ufs{|\psiinv j|}} \, [1_\ome, \iphijlome]_\mstar^{\sum_{t\in \phiinv \ell'} a_t} && \text{by \cref{cglgdom_psix}}\\
&= \txsum_{\ell \in \ufs{|\psiinv j|}} \txsum_{t\in \phiinv \ell'} \, [1_\ome, \iphijlome]_\mstar^{a_t} && \text{by $\schm$-equivariance of $\psum$}\\
&= \txsum_{\ell \in \ufs{|\psiinv j|}} \txsum_{r\in \ufs{|\phiinv \ell'|}} \, [1_\ome, \iphijlome]_\mstar^{a_{r'}} && \text{by \cref{parcat_axioms,obij_phiinv,r_prime}}\\
\end{aligned}$}
\end{equation}
The $j$-th entry of the second morphism in \cref{cglgc_compeq} is the following morphism in $\C$.
\begin{equation}\label{cglgc_comp_ii}
\scalebox{.9}{$
\begin{aligned}
& \big(\objx_*\cglgcya\big)_j &&\\
&= \big(\objx_* \bang{\, \txsum_{t \in \ufs{|\phiinv i|}} \, [1_\ome, \ivphi_i^{t,\ome}]_\mstar^{a_{t'}}}_{i \in \ufsm} \big)_j && \text{by \cref{cglgcxa,t_prime}}\\
&= \txsum_{\ell \in \ufs{|\psiinv j|}} \, \iphijlomest \big(\txsum_{r \in \ufs{|\phiinv \ell'|}} \, [1_\ome,\ivphi_{\ell'}^{r,\ome}]_\mstar^{a_{r'}} \big) && \text{by \cref{ricxa,ell_prime,r_prime}}\\
&= \txsum_{\ell \in \ufs{|\psiinv j|}} \txsum_{r \in \ufs{|\phiinv \ell'|}} \, \iphijlomest [1_\ome,\ivphi_{\ell'}^{r,\ome}]_\mstar^{a_{r'}} && \text{by $\schm$-equivariance of $\psum$}\\
&= \txsum_{\ell \in \ufs{|\psiinv j|}} \txsum_{r \in \ufs{|\phiinv \ell'|}} \, [\iphijlome, \iphijlome \ivphi_{\ell'}^{r,\ome}]_\mstar^{a_{r'}} && \text{by \cref{mcat_axioms}}\\
&= \txsum_{\ell \in \ufs{|\psiinv j|}} \txsum_{r \in \ufs{|\phiinv \ell'|}} \, [\iphijlome, (\iphijl \ivphi_{\ell'}^r)^{\ome}]_\mstar^{a_{r'}} && \text{by $\repomeg\repomeginv = 1_{\omeg}$}
\end{aligned}$}
\end{equation}
Using \cref{cglgc_comp_i,cglgc_comp_ii}, the $j$-th entry of the left-hand side in \cref{cglgc_compeq} is the following morphism in $\C$.
\[\scalebox{.85}{$
\begin{aligned}
& \big(\cglgc_{\objx,\phi_*\obja}\big)_j \big(\objx_*\cglgcya\big)_j &&\\
&= \txsum_{\ell \in \ufs{|\psiinv j|}} \txsum_{r\in \ufs{|\phiinv \ell'|}} \, [1_\ome, \iphijlome]_\mstar^{a_{r'}} [\iphijlome, (\iphijl \ivphi_{\ell'}^r)^{\ome}]_\mstar^{a_{r'}} && \text{by functoriality of $\psum$}\\
&= \txsum_{\ell \in \ufs{|\psiinv j|}} \txsum_{r\in \ufs{|\phiinv \ell'|}} \, \big([1_\ome, \iphijlome] [\iphijlome, (\iphijl \ivphi_{\ell'}^r)^{\ome}]\big)_\mstar^{a_{r'}} && \text{by \cref{mcat_axioms}}\\
&= \txsum_{\ell \in \ufs{|\psiinv j|}} \txsum_{r\in \ufs{|\phiinv \ell'|}} \, [1_\ome, (\iphijl \ivphi_{\ell'}^r)^{\ome}]_\mstar^{a_{r'}} && \\
&= \txsum_{q \in \ufs{|(\psi\phi)^{-1} j|}} \, \big[1_\ome, \big(\ga(\iphi_j; \ang{\ivphi_i}_{i \in \psiinv j})^{\tau_{\psi,\phi}^j} \big)^{q,\ome} \big]_\mstar^{a_{\obij_{(\psi\phi)^{-1}j} q}} && \text{by \cref{parcat_axioms,tauj_psiphi,inopg_comp}} \\
&= (\cglgcxya)_j && \text{by \cref{dg_iopg_comp,cglgcxa}}
\end{aligned}$}\]
Since $j \in \ufsn$ is arbitrary, this proves the morphism equality \cref{cglgc_compeq}.  Thus, $\cglgc$ \cref{cglgc} is a $\Gcatst$-pseudotransformation.
\item[$\Pig$-strictness]
To prove that $\cglgc$ is $\Pig$-strict, suppose $\objx = (\psi; \ang{\iphi_j}_{j \in \ufsn}) \in \DG(\mal,\nbeta)$ is in the image of $\gio \cn \Pig \to \DG$ \cref{dgcoop}.  This means that $|\psiinv j| \in \{0,1\}$ for each $j \in \ufsn$ and
\[\iphi_j = \begin{cases}
* \in \Iopg(0) & \text{if $|\psiinv j| = 0$ and}\\
1_{\omeg} \in \Iopg(1) & \text{if $|\psiinv j| = 1$.} 
\end{cases}\]
To prove that $\cglgcx$ is the identity natural transformation, consider the component morphism \cref{cglgcxa}
\[\cglgcxa = \bang{\, \txsum_{\ell \in \ufs{|\psiinv j|}} \, [1_\ome, \iphijlome]_\mstar^{a_{\obij_{\psiinv j}\ell}}}_{j \in \ufsn} \inspace \Comegnbeta\]
for an object $\obja = \ang{a_i}_{i \in \ufsm} \in \Cboxm$. 
\begin{itemize}
\item If $|\psiinv j| = 0$, then the $j$-th entry in $\cglgcxa$ is the empty sum, which means the identity morphism $1_\pzero$.
\item If $|\psiinv j| = 1$ with $\iphi_j = 1_{\omeg}$, then there are an object equality
\[(1_{\omeg})^{1,\ome} = 1_\ome \cn \ome \to \ome\]
by \cref{iphijlome} and a morphism equality 
\[[1_\ome,1_\ome] = 1_{(1_\ome)} \inspace \schm.\]
Thus, the $j$-th entry in $\cglgcxa$ is the identity morphism.
\end{itemize}
This proves that $\cglgcx$ is the identity natural transformation.
\end{description}
Thus, $\cglgc$ is a $\Pig$-strict $\Gcatst$-pseudotransformation.
\end{proof}

The categories and functors in \cref{cglg_natural} are recalled at the beginning of this section.

\begin{lemma}\label{cglg_natural}
For a finite group $G$, there is a natural transformation
\begin{equation}\label{cglg_nat}

\end{equation}
whose component at a parsummable category $\C$ is $\cglgc$ \cref{cglgc}.
\end{lemma}

\begin{proof}
\cref{cglgc_welldef} proves that each component $\cglgc$ is a morphism in $\DGCatg$ \pcref{def:dgcatg_iicat}.  Naturality of $\cglg$ means that, for each parsummable functor $\fun \cn \C \to \D$, the diagram
\begin{equation}\label{cglg_nat_diag}
%
\end{equation}
in $\DGCatg$ commutes.  We verify that the two composites in \cref{cglg_nat_diag} have the same component pointed $G$-functors and component natural isomorphisms.
\begin{description}
\item[$G$-functors]
By \cref{cglgdom_fn,cglgcod_fn,cglgcn}, the diagram \cref{cglg_nat_diag} evaluated at an object $\nbeta \in \DG$ is the following diagram of pointed $G$-functors.
\begin{equation}\label{cglgcdf}
%
\end{equation}
This diagram commutes because $\cglgcnbeta$ and $\cglgdnbeta$ are the inclusion functors $\Cboxn \to \C^n$ and $\Dboxn \to \D^n$. 
\item[Natural isomorphisms]
By \cref{cglgdom_fn,cglgcod_fn}, the following morphism equalities in $\Domegnbeta$ prove that the two composites in \cref{cglg_nat_diag} have the same component morphisms at an object $\objx = (\psi; \ang{\iphi_j}_{j \in \ufsn})$ in $\DG(\mal,\nbeta)$ and an object $\obja = \ang{a_i}_{i \in \ufsm}$ in $\Cboxm$.
\[\begin{aligned}
\fun^n \cglgcxa
&= \bang{\fun\txsum_{\ell \in \ufs{|\psiinv j|}} \, [1_\ome,\iphijlome]_{\mstar}^{a_{\obij\ell}}}_{j \in \ufsn} && \text{by \cref{cglgcxa}}\\
&= \bang{\txsum_{\ell \in \ufs{|\psiinv j|}} \, \fun[1_\ome,\iphijlome]_{\mstar}^{a_{\obij\ell}}}_{j \in \ufsn} && \text{by \cref{parfun_axioms}}\\
&= \bang{\txsum_{\ell \in \ufs{|\psiinv j|}} \, [1_\ome,\iphijlome]_{\mstar}^{\fun a_{\obij\ell}}}_{j \in \ufsn} && \text{by \cref{mfunctor_axiom}}\\
&= \cglgd_{\objx, \fboxm \obja} && \text{by \cref{cglgcxa}}
\end{aligned}\]
\end{description}
This proves that $\cglg$ is a natural transformation.
\end{proof}

\begin{explanation}[Pseudonaturality]\label{expl:cglg_pseudo}
In the comparison natural transformation $\cglg$ \cref{cglg_nat}, the functors $\gxist$ and $\Rg$ land in the subcategory (\cref{def:dgcatg_iicat} \cref{def:dgcatg_iicat_ii})
\[\begin{tikzpicture}[vcenter]
\draw[0cell]
(0,0) node (a1) {\dgcatg}
(a1)++(2.5,0) node (a2) {\DGCatg}
;
\draw[1cell=.9]
(a1) edge[right hook->] node {\dgi} (a2)
;
\end{tikzpicture}\]
whose morphisms are $\Gcatst$-natural transformations. 
However, $\cglg$ is only defined for the category $\DGCatg$ and \emph{not} $\dgcatg$.  The reason is that the components of $\cglg$ are $\Pig$-strict $\Gcatst$-pseudotransformations but not $\Gcatst$-natural transformations.  These components of $\cglg$ are only $\Gcatst$-natural up to the natural isomorphisms $\cglgcx$ in \cref{cglgcx}.
\end{explanation}

\section{Comparison $G$-Equivalences}
\label{sec:cglg_geq}

This section proves that, for each finite group $G$, pointed finite $G$-set $\nbeta$, and parsummable category $(\C,\psum,\pzero)$, the comparison pointed $G$-functor \cref{cglgcn}
\[\Cboxn = (\gxist\Lg\Jscg\C)\nbeta \fto{\cglgcnbeta} (\Rg\Ig\C)\nbeta = \Comegnbeta\]
is a pointed $G$-equivalence.  This means that there exist a pointed $G$-functor
\[\Comegnbeta \fto{\cgglcnbeta} \Cboxn\]
and pointed $G$-natural isomorphisms
\begin{equation}\label{cglguv_secti}
\begin{tikzpicture}[vcenter]
\def\h{1.8} \def\t{28} \def\c{.15} \def\d{.43} \def\s{.8}
\draw[0cell]
(0,0) node (a1) {\phantom{A}}
(a1)++(\h,0) node (a2) {\phantom{A}}
(a1)++(-\c,0) node (a1') {\Cboxn}
(a2)++(\c,0) node (a2') {\Cboxn}
;
\draw[1cell=\s]
(a1) edge[bend left=\t] node {1} (a2)
(a1) edge[bend right=\t] node[swap] {\cgglcnbeta \cglgcnbeta} (a2)
;
\draw[2cell]
node[between=a1 and a2 at .42, rotate=-90, 2label={above,\cglgu}] {\Rightarrow}
;
\begin{scope}[shift={(4,0)}]
\draw[0cell]
(0,0) node (a1) {\phantom{A}}
(a1)++(\h,0) node (a2) {\phantom{A}}
(a1)++(-\d,0) node (a1') {\Comegnbeta}
(a2)++(\d,0) node (a2') {\Comegnbeta}
;
\draw[1cell=\s]
(a1) edge[bend left=\t] node {\cglgcnbeta\cgglcnbeta} (a2)
(a1) edge[bend right=\t] node[swap] {1} (a2)
;
\draw[2cell]
node[between=a1 and a2 at .42, rotate=-90, 2label={above,\cglgv}] {\Rightarrow}
;
\end{scope}
\end{tikzpicture}
\end{equation}
that satisfy the triangle identities for a $G$-equivariant adjunction.  Since $\cglgc_{\ordz}$ is the identity functor of the terminal $G$-category $\C^0 = \bone$, we may assume that $\nbeta \neq \ordz$.

\secoutline
\begin{itemize}
\item \cref{def:cgglc,cgglc_welldef} construct the pointed $G$-functor $\cgglcnbeta$.
\item \cref{def:cglgu,cglgu_welldef} construct the unit pointed $G$-natural isomorphism $\cglgu$.
\item \cref{def:cglgv,cglgv_welldef} construct the counit pointed $G$-natural isomorphism $\cglgv$.
\item \cref{cglgc_equiv} proves that the quadruple $(\cglgcnbeta,\cgglcnbeta,\cglgu,\cglgv)$ is an adjoint pointed $G$-equivalence.
\item \cref{rk:cggl_notpigst} discusses the fact that the pointed $G$-functors $\cgglcnbeta$ can be extended to a $\Gcatst$-pseudotransformation that is \emph{not} $\Pig$-strict. 
\end{itemize}

\subsection*{The $G$-Functor $\cgglcnbeta$}

Recall 
\begin{itemize}
\item the unpointed finite set $\ufsn = \{1,2,\ldots,n\}$ and
\item the universal $G$-set $\omeg$ consisting of functions $f \cn G \to \ome = \{0,1,2,\ldots\}$ with the $G$-action $gf = f(\ginv \cdot -)$ \pcref{def:omea}.  
\end{itemize}
The pointed $G$-functor $\cgglcnbeta$ in \cref{def:cgglc} uses the injections $\uppsjome$ in \cref{def:upps}.

\begin{definition}[$G$-Bijection $\upps$]\label{def:upps}
Given a finite group $G$ and a pointed finite $G$-set $\nbeta \neq \ordz$, there is a $G$-equivariant injection
\[\omeg \fto{i_0 = (0,-)} \nbetaomeg\]
into the product $G$-set $\nbetaomeg$.  By \cref{schwede2.17} \cref{sch2.17_iii}, $\nbetaomeg$ is a universal $G$-set.  By \cref{schwede2.17} \cref{sch2.17_ii}, any two universal $G$-sets are $G$-isomorphic.  Thus, we can choose a $G$-bijection
\begin{equation}\label{upps_iso}
\nbetaomeg \fto[\iso]{\upps} \omeg.
\end{equation}
We use the notation $\upps_{\nbeta}$ when we need to emphasize the dependency of $\upps$ on $\nbeta$.  Restricting $\upps$ to $\{j\} \ttimes \omeg$ for $j \in \ufsn$ yields $n$ injections
\begin{equation}\label{upps_j}
\omeg \fto{\upps^j \,=\, \upps(j,-)} \omeg
\end{equation}
with disjoint images.  Using the notation in \cref{repU_p_repV}, the $n$ injections
\begin{equation}\label{uppsjome}
\begin{tikzpicture}[vcenter]
\def\h{1.8} \def\u{.7} \def\d{.05}
\draw[0cell]
(0,0) node (a1) {\ome}
(a1)++(\h,0) node (a2) {\phantom{\omeg}}
(a2)++(0,\d) node (a2') {\omeg}
(a2)++(\h,0) node (a3) {\phantom{\omeg}}
(a3)++(0,\d) node (a3') {\omeg}
(a3)++(\h,0) node (a4) {\ome}
;
\draw[1cell=.9]
(a1) edge node {\rep_{\omeg}} node[swap] {\iso} (a2)
(a2) edge node {\upps^j} (a3)
(a3) edge node {\rep_{\omeg}^{-1}} node[swap] {\iso} (a4)
(a1) [rounded corners=2pt] |- ($(a2)+(0,\u)$) -- node {\uppsjome} ($(a3)+(0,\u)$) -| (a4)
;
\end{tikzpicture}
\end{equation}
for $j \in \ufsn$ also have disjoint images.  For an $\schm$-category $\C$, using the notation in \cref{u_action}, the $\schm$-action functor associated to $\uppsjome$ is denoted by
\begin{equation}\label{uppsjomest}
\C \fto{\uppsjomest} \C.
\end{equation}
Since $\upps$ is $G$-equivariant, the injections $\upps^j$ \cref{upps_j} are \emph{twisted $G$-equivariant} in the sense that
\begin{equation}\label{uppsj_g}
g\upps^j(f) = g\upps(j,f) = \upps(gj,gf) = \upps^{gj}(gf)
\end{equation}
for each $j \in \ufsn$ and $f \in \omeg$.
\end{definition}

\begin{definition}[$G$-Functor $\cgglcnbeta$]\label{def:cgglc}
For a finite group $G$, a pointed finite $G$-set $\nbeta \neq \ordz$, and a parsummable category $(\C,\psum,\pzero)$ \pcref{def:parcat}, the pointed $G$-functor
\begin{equation}\label{cgglcnbeta_functor}
\Comegnbeta = (\Rg\Ig\C)\nbeta \fto{\cgglcnbeta} (\gxist\Lg\Jscg\C)\nbeta = \Cboxn
\end{equation}
between the pointed $G$-categories in \cref{cglgdom_nbeta,ricn} is defined by
\begin{equation}\label{cgglcnbeta_def}
\cgglcnbeta\obja = \ang{\uppsjomest a_j}_{j \in \ufsn} \in \Cboxn
\end{equation}
for each object or morphism $\obja = \ang{a_j}_{j \in \ufsn} \in \Comegnbeta$.  If $\obja$ is an object, then the objects $\uppsjomest a_j \in \C$ for $j \in \ufsn$ have disjoint supports by \cref{sch2.13} \cref{sch2.13_iii} and the fact that the injections $\uppsjome$ \cref{uppsjome} have disjoint images.
\end{definition}

\begin{lemma}\label{cgglc_welldef}
In \cref{cgglcnbeta_functor}, $\cgglcnbeta \cn \Comegnbeta \to \Cboxn$ is a pointed $G$-functor.
\end{lemma}

\begin{proof}
The basepoint in each of $\Comegnbeta$ and $\Cboxn$ is $\ang{\pzero}_{j \in \ufsn}$.  Since the unit $\pzero \in \C$ is $\schm$-fixed, $\cgglcnbeta$ preserves the basepoint.  The functoriality of $\cgglcnbeta$ follows from the functoriality of each $\schm$-action functor $\uppsjomest$.  The $G$-equivariance of $\cgglcnbeta$ is proved by the following equalities in $\Cboxn$ for $g \in G$ and an object or a morphism $\obja = \ang{a_j}_{j \in \ufsn} \in \Comegnbeta$.
\begin{equation}\label{cgglc_gequi}
\begin{aligned}
\cgglcnbeta\big(g\obja\big) 
&= \ang{\uppsjomest (g\obja)_j}_{j \in \ufsn} && \text{by \cref{cgglcnbeta_def}}\\
&= \ang{\uppsjomest \gomest a_{\ginv j} }_{j \in \ufsn} && \text{by \cref{ricn_gact}}\\
&= \ang{(\uppsjome g^\ome)_\mstar a_{\ginv j}}_{j \in \ufsn} && \text{by \cref{mcat_axioms}}\\
&= \ang{(\upps^j g)^\ome_\mstar a_{\ginv j}}_{j \in \ufsn} && \text{by $\rep_{\omeg}\rep_{\omeg}^{-1} = 1$} \\
&= \ang{(g\upps^{\ginv j})^\ome_\mstar a_{\ginv j}}_{j \in \ufsn} && \text{by \cref{uppsj_g}} \\
&= \ang{(g^\ome\uppsginvjome)_\mstar a_{\ginv j}}_{j \in \ufsn} && \text{by $\rep_{\omeg}\rep_{\omeg}^{-1} = 1$} \\
&= \ang{\gomest\uppsginvjome_\mstar a_{\ginv j}}_{j \in \ufsn} && \text{by \cref{mcat_axioms}}\\
&= g\ang{\uppsjomest a_j}_{j \in \ufsn}  && \text{by \cref{cglgdom_gact}}\\
&= g(\cgglcnbeta\obja) && \text{by \cref{cgglcnbeta_def}}
\end{aligned}
\end{equation}
This proves that $\cgglcnbeta$ is a pointed $G$-functor.
\end{proof}

\subsection*{Unit}

\cref{def:cglgu} constructs the $G$-equivariant unit for the pointed $G$-functors
\begin{equation}\label{cglgcggl_diag}
\begin{tikzpicture}[vcenter]
\def\h{1.8} \def\t{28} \def\c{.15} \def\d{.5} \def\s{.8}
\draw[0cell]
(0,0) node (a1) {\phantom{A}}
(a1)++(\h,0) node (a2) {\phantom{A}}
(a1)++(-\c,0) node (a1') {\Cboxn}
(a2)++(\d,-.03) node (a2') {\Comegnbeta}
;
\draw[1cell=\s]
(a1) edge[transform canvas={yshift=.4ex}] node {\cglgcnbeta} (a2)
(a2) edge[transform canvas={yshift=-.5ex}] node {\cgglcnbeta} (a1)
;
\end{tikzpicture}
\end{equation}
in \cref{cglgcn,cgglcnbeta_functor}.

\begin{definition}\label{def:cglgu}
Under the same hypotheses as \cref{def:cgglc}, we define the pointed $G$-natural isomorphism
\begin{equation}\label{cglgu_nt}
\begin{tikzpicture}[vcenter]
\def\h{1.8} \def\t{28} \def\c{.15} \def\d{.43} \def\s{.8}
\draw[0cell]
(0,0) node (a1) {\phantom{A}}
(a1)++(\h,0) node (a2) {\phantom{A}}
(a1)++(-\c,0) node (a1') {\Cboxn}
(a2)++(\c,0) node (a2') {\Cboxn}
;
\draw[1cell=\s]
(a1) edge[bend left=\t] node {1} (a2)
(a1) edge[bend right=\t] node[swap] {\cgglcnbeta \cglgcnbeta} (a2)
;
\draw[2cell]
node[between=a1 and a2 at .42, rotate=-90, 2label={above,\cglgu}] {\Rightarrow}
;
\end{tikzpicture}
\end{equation}
whose component at an object $\objx = \ang{x_j}_{j \in \ufsn} \in \Cboxn$ is given by the commutative diagram
\begin{equation}\label{cglgu_x}
\begin{tikzpicture}[vcenter]
\def\v{-1}
\draw[0cell]
(0,0) node (a11) {\objx}
(a11)++(4.5,0) node (a12) {\phantom{\cgglcnbeta \cglgcnbeta\objx}}
(a12)++(0,.02) node (a12') {\cgglcnbeta \cglgcnbeta\objx}
(a11)++(0,\v) node (a21) {\ang{x_j}_{j \in \ufsn}}
(a12)++(0,\v) node (a22) {\ang{\uppsjomest x_j}_{j \in \ufsn}}
;
\draw[1cell=.9]
(a11) edge[equal] (a21)
(a12') edge[equal,shorten >=-.2ex] (a22)
(a11) edge node {\cglgu_\objx} (a12)
(a21) edge node {\bang{[\uppsjome,1_\ome]_\mstar^{x_j}}_{j \in \ufsn}} (a22)
;
\end{tikzpicture}
\end{equation}
in $\Cboxn$ \cref{cglgdom_nbeta}.  The right equality in \cref{cglgu_x} follows from \cref{cglgcn,cgglcnbeta_def}.  For each $j \in \ufsn$, 
\[1 = (1_\ome)_\mstar \fto{[\uppsjome,1_\ome]_\mstar} \uppsjomest\]
is the $\schm$-action natural isomorphism on $\C$ \cref{vu_star} associated to the isomorphism 
\[1_\ome \fto{[\uppsjome,1_\ome]} \uppsjome\]
in $\schm$.  
\end{definition}

\begin{lemma}\label{cglgu_welldef}
In \cref{cglgu_nt},
\[1_{\Cboxn} \fto{\cglgu} \cgglcnbeta \cglgcnbeta\]
is a pointed $G$-natural isomorphism.
\end{lemma}

\begin{proof}
For each object $\objx \in \Cboxn$, $\cglgu_\objx$ \cref{cglgu_x} is an isomorphism because  $[\uppsjome,1_\ome]_\mstar^{x_j}$ is an isomorphism for each $j \in \ufsn$.  The naturality of $\cglgu$ in $\objx \in \Cboxn$ follows from the naturality of each $[\uppsjome,1_\ome]_\mstar$.  The natural isomorphism $\cglgu$ is pointed, meaning 
\[\cglgu_{\ang{\pzero}_{j \in \ufsn}} = 1_{\ang{\pzero}_{j \in \ufsn}},\]
because the unit $\pzero \in \C$ is $\schm$-fixed.   Similar to the computation in \cref{cgglc_gequi}, for $g \in G$, the following morphism equalities in $\Cboxn$ prove that $\cglgu$ is $G$-equivariant.
\begin{equation}\label{cglgu_gequiv}
\begin{aligned}
g\cglgu_{\objx} 
&= g\ang{[\uppsjome,1_\ome]_\mstar^{x_j}}_{j \in \ufsn} && \text{by \cref{cglgu_x}}\\
&= \ang{\gomest[\uppsginvjome,1_\ome]_\mstar^{x_{\ginv j}}}_{j \in \ufsn} && \text{by \cref{cglgdom_gact}}\\
&= \ang{[\gome\uppsginvjome,\gome]_\mstar^{x_{\ginv j}}}_{j \in \ufsn} && \text{by \cref{mcat_axioms}}\\
&= \ang{[\uppsjome\gome,\gome]_\mstar^{x_{\ginv j}}}_{j \in \ufsn} && \text{by \cref{uppsj_g} and $\rep_{\omeg}\rep_{\omeg}^{-1} = 1$}\\
&= \ang{[\uppsjome,1_\ome]_\mstar^{\gomest x_{\ginv j}}}_{j \in \ufsn} && \text{by \cref{mcat_axioms}}\\
&= \ang{[\uppsjome,1_\ome]_\mstar^{(g\objx)_j}}_{j \in \ufsn} && \text{by \cref{cglgdom_gact}}\\
&= \cglgu_{g\objx} && \text{by \cref{cglgu_x}}
\end{aligned}
\end{equation}
This proves that $\cglgu$ is a pointed $G$-natural isomorphism.
\end{proof}

\subsection*{Counit}

\cref{def:cglgv} constructs the $G$-equivariant counit for the pointed $G$-functors $(\cglgcnbeta,\cgglcnbeta)$ in \cref{cglgcn,cgglcnbeta_functor}.  Recall the pointed $G$-category $\Comegnbeta$ \cref{ricn}.

\begin{definition}\label{def:cglgv}
Under the same hypotheses as \cref{def:cgglc}, we define the pointed $G$-natural isomorphism
\begin{equation}\label{cglgv_nt}
\begin{tikzpicture}[vcenter]
\def\h{1.8} \def\t{28} \def\c{.15} \def\d{.43} \def\s{.8}
\draw[0cell]
(0,0) node (a1) {\phantom{A}}
(a1)++(\h,0) node (a2) {\phantom{A}}
(a1)++(-\d,0) node (a1') {\Comegnbeta}
(a2)++(\d,0) node (a2') {\Comegnbeta}
;
\draw[1cell=\s]
(a1) edge[bend left=\t] node {\cglgcnbeta\cgglcnbeta} (a2)
(a1) edge[bend right=\t] node[swap] {1} (a2)
;
\draw[2cell]
node[between=a1 and a2 at .42, rotate=-90, 2label={above,\cglgv}] {\Rightarrow}
;
\end{tikzpicture}
\end{equation}
whose component at an object $\obja = \ang{a_j}_{j \in \ufsn} \in \Comegnbeta$ is given by the commutative diagram
\begin{equation}\label{cglgv_a}
\begin{tikzpicture}[vcenter]
\def\v{-1}
\draw[0cell]
(0,0) node (a11) {\phantom{\cglgcnbeta\cgglcnbeta\obja}}
(a11)++(0,0) node (a11') {\cglgcnbeta\cgglcnbeta\obja}
(a11)++(4.3,0) node (a12) {\obja}
(a11)++(0,\v) node (a21) {\ang{\uppsjomest a_j}_{j \in \ufsn}}
(a12)++(0,\v) node (a22) {\phantom{\ang{a_j}_{j \in \ufsn}}}
(a22)++(0,-.05) node (a22') {\ang{a_j}_{j \in \ufsn}}
;
\draw[1cell=.9]
(a11') edge[equal,shorten >=-.2ex] (a21)
(a12) edge[equal] (a22')
(a11) edge node {\cglgv_\obja} (a12)
(a21) edge node {\bang{[1_\ome,\uppsjome]_\mstar^{a_j}}_{j \in \ufsn}} (a22)
;
\end{tikzpicture}
\end{equation}
in $\Comegnbeta$.  The left equality in \cref{cglgv_a} follows from \cref{cglgcn,cgglcnbeta_def}.  For each $j \in \ufsn$, 
\[\uppsjomest \fto{[1_\ome,\uppsjome]_\mstar} (1_\ome)_\mstar = 1\]
is the $\schm$-action natural isomorphism on $\C$ \cref{vu_star} associated to the isomorphism 
\[\uppsjome \fto{[1_\ome,\uppsjome]} 1_\ome\]
in $\schm$.  
\end{definition}

\begin{lemma}\label{cglgv_welldef}
In \cref{cglgv_nt},
\[\cglgcnbeta\cgglcnbeta \fto{\cglgv} 1_{\Comegnbeta} \]
is a pointed $G$-natural isomorphism.
\end{lemma}

\begin{proof}
For each object $\obja \in \Comegnbeta$, $\cglgv_\obja$ \cref{cglgv_a} is an isomorphism because  $[1_\ome,\uppsjome]_\mstar^{a_j}$ is an isomorphism for each $j \in \ufsn$.  The naturality of $\cglgv$ in $\obja \in \Comegnbeta$ follows from the naturality of each $[1_\ome,\uppsjome]_\mstar$.  The natural isomorphism $\cglgv$ is pointed, meaning 
\[\cglgv_{\ang{\pzero}_{j \in \ufsn}} = 1_{\ang{\pzero}_{j \in \ufsn}},\]
because the unit $\pzero \in \C$ is $\schm$-fixed.   Similar to the computation in \cref{cglgu_gequiv}, for $g \in G$, the following morphism equalities in $\Comegnbeta$ prove that $\cglgv$ is $G$-equivariant.
\[\begin{aligned}
g\cglgv_{\obja} 
&= g\ang{[1_\ome,\uppsjome]_\mstar^{a_j}}_{j \in \ufsn} && \text{by \cref{cglgv_a}}\\
&= \ang{\gomest[1_\ome,\uppsginvjome]_\mstar^{a_{\ginv j}}}_{j \in \ufsn} && \text{by \cref{ricn_gact}}\\
&= \ang{[\gome,\gome\uppsginvjome]_\mstar^{a_{\ginv j}}}_{j \in \ufsn} && \text{by \cref{mcat_axioms}}\\
&= \ang{[\gome,\uppsjome\gome]_\mstar^{a_{\ginv j}}}_{j \in \ufsn} && \text{by \cref{uppsj_g} and $\rep_{\omeg}\rep_{\omeg}^{-1} = 1$}\\
&= \ang{[1_\ome,\uppsjome]_\mstar^{\gomest a_{\ginv j}}}_{j \in \ufsn} && \text{by \cref{mcat_axioms}}\\
&= \ang{[1_\ome,\uppsjome]_\mstar^{(g\obja)_j}}_{j \in \ufsn} && \text{by \cref{ricn_gact}}\\
&= \cglgv_{g\obja} && \text{by \cref{cglgv_a}}
\end{aligned}\]
This proves that $\cglgv$ is a pointed $G$-natural isomorphism.
\end{proof}

\subsection*{Adjoint $G$-Equivalence}

\begin{theorem}\label{cglgc_equiv}
For a finite group $G$, a pointed finite $G$-set $\nbeta \neq \ordz$, and a parsummable category $(\C,\psum,\pzero)$, the pointed $G$-functors
\begin{equation}\label{cglgc_thm}
\begin{tikzpicture}[vcenter]
\def\h{1.8} \def\t{28} \def\c{.15} \def\d{.5} \def\s{.8}
\draw[0cell]
(0,0) node (a1) {\phantom{A}}
(a1)++(\h,0) node (a2) {\phantom{A}}
(a1)++(-\c,.01) node (a1') {\Cboxn}
(a2)++(\d,-.03) node (a2') {\Comegnbeta}
(a1')++(-1.9,0) node (a1'') {\phantom{(\gxist\Lg\Jscg\C)\nbeta}}
(a1'')++(0,-.05) node (a0) {(\gxist\Lg\Jscg\C)\nbeta}
(a2')++(2.1,0) node (a2'') {(\Rg\Ig\C)\nbeta}
;
\draw[1cell=\s]
(a1') edge[equal,transform canvas={yshift=-.25ex}] (a1'')
(a2') edge[equal] (a2'')
(a1) edge[transform canvas={yshift=.4ex}] node {\cglgcnbeta} (a2)
(a2) edge[transform canvas={yshift=-.5ex}] node {\cgglcnbeta} (a1)
;
\end{tikzpicture}
\end{equation}
defined in \cref{cglgcn,cgglcnbeta_functor} form an adjoint pointed $G$-equivalence\index{adjoint G-equivalence@adjoint $G$-equivalence} with
\begin{itemize}
\item unit $\cglgu \cn 1_{\Cboxn} \to \cgglcnbeta\cglgcnbeta$ \cref{cglgu_nt} and
\item counit $\cglgv \cn \cglgcnbeta\cgglcnbeta \to 1_{\Comegnbeta}$ \cref{cglgv_nt}.
\end{itemize}
\end{theorem}

\begin{proof}
The pointed $G$-functors $\cglgcnbeta$ and $\cgglcnbeta$ are well defined by \cref{cglgcn,cgglc_welldef}.  The pointed $G$-natural isomorphisms $\cglgu$ and $\cglgv$ are well defined by \cref{cglgu_welldef,cglgv_welldef}.  It remains to prove the triangle identities for an adjunction \pcref{def:adjunction}.
\begin{description}
\item[Left triangle identity]
This identity states that, for each object $\objx = \ang{x_j}_{j \in \ufsn} \in \Cboxn$, the composite
\begin{equation}\label{cglg_ltri}
\begin{tikzpicture}[baseline={(a1.base)}]
\def\h{2.6} \def\u{.7} 
\draw[0cell]
(0,0) node (a1) {\cglgcnbeta\objx}
(a1)++(\h,0) node (a2) {\cglgcnbeta \cgglcnbeta \cglgcnbeta\objx}
(a2)++(\h,0) node (a3) {\cglgcnbeta\objx}
;
\draw[1cell=.9]
(a1) edge node {\cglgcnbeta\cglgu_\objx} (a2)
(a2) edge node {\cglgv_{\cglgcnbeta\objx}} (a3)
;
\end{tikzpicture}
\end{equation}
in $\Comegnbeta$ is equal to the identity morphism.  By \cref{cglgcn,cgglcnbeta_def,cglgu_x,cglgv_a}, the composite in \cref{cglg_ltri} is given as follows, where $\ang{\Cdots}_j = \ang{\Cdots}_{j \in \ufsn}$.
\[\ang{x_j}_j \fto{\ang{[\uppsjome,1_\ome]_\mstar^{x_j}}_j} \ang{\uppsjomest x_j}_j
\fto{\ang{[1_\ome,\uppsjome]_\mstar^{x_j}}_j} \ang{x_j}_{j} \]
The preceding composite is equal to the identity morphism by the functoriality of the $\schm$-action on $\C$ and the morphism equality
\[[1_\ome,\uppsjome] \comp [\uppsjome,1_\ome] = 1_{(1_\ome)} \cn 1_\ome \to 1_\ome\]
in the injection category $\schm$ \pcref{def:mcat}.
\item[Right triangle identity]
This identity states that, for each object $\obja = \ang{a_j}_{j \in \ufsn} \in \Comegnbeta$, the composite
\begin{equation}\label{cglg_rtri}
\begin{tikzpicture}[baseline={(a1.base)}]
\def\h{2.6} \def\u{.7} 
\draw[0cell]
(0,0) node (a1) {\cgglcnbeta\obja}
(a1)++(\h,0) node (a2) {\cgglcnbeta \cglgcnbeta \cgglcnbeta\obja}
(a2)++(\h,0) node (a3) {\cgglcnbeta\obja}
;
\draw[1cell=.9]
(a1) edge node {\cglgu_{\cgglcnbeta\obja}} (a2)
(a2) edge node {\cgglcnbeta\cglgv_{\obja}} (a3)
;
\end{tikzpicture}
\end{equation}
in $\Cboxn$ is equal to the identity morphism.   The composite in \cref{cglg_rtri} is given as follows.
\[\ang{\uppsjomest a_j}_j \fto{\ang{[\uppsjome, 1_\ome]_\mstar^{\uppsjomest a_j}}_j}
\ang{\uppsjomest\uppsjomest a_j}_j \fto{\ang{\uppsjomest [1_\ome,\uppsjome]_\mstar^{a_j}}_j}
\ang{\uppsjomest a_j}_j\]
The preceding composite is equal to the identity morphism by the following morphism equalities in $\C$ for each $j \in \ufsn$.
\[\begin{aligned}
& \uppsjomest [1_\ome,\uppsjome]_\mstar^{a_j} \comp [\uppsjome, 1_\ome]_\mstar^{\uppsjomest a_j} &&\\
&= [\uppsjome,\uppsjome\uppsjome]_\mstar^{a_j} \comp [\uppsjome\uppsjome, \uppsjome]_\mstar^{a_j} && \text{by \cref{mcat_axioms}}\\
&= \big([\uppsjome,\uppsjome\uppsjome] \comp [\uppsjome\uppsjome, \uppsjome]\big)_\mstar^{a_j} && \text{by functoriality of $\schm$-action}\\
&= [\uppsjome,\uppsjome]_\mstar^{a_j} &&\\
&= (1_{\uppsjome})_\mstar^{a_j} = 1_{\uppsjomest a_j} &&\\
\end{aligned}\]
This proves the right triangle identity \cref{cglg_rtri}.
\end{description}
Thus, $(\cglgcnbeta,\cgglcnbeta,\cglgu,\cglgv)$ is an adjoint pointed $G$-equivalence.
\end{proof}

\cref{rk:cggl_notpigst} discusses further properties of $\cgglcnbeta$.  It is not needed in the rest of this work, so it can be safely skipped.

\begin{explanation}[Pseudonaturality of $\cgglc$]\label{rk:cggl_notpigst}
The pointed $G$-functors $\cgglcnbeta$ \cref{cgglcnbeta_functor} can be extended to a $\Gcatst$-pseudotransformation, but it is \emph{not} $\Pig$-strict.  Thus, $\cgglc$ is not a morphism in $\DGCatg$ \pcref{expl:dgiopg}.

More precisely, suppose $\objx = (\psi; \ang{\iphi_j}_{j\in \ufsn})$ is an object in $\DG(\mal,\nbeta)$ \cref{dgmn_iopg_obj}.  Using the $\schm$-functoriality of $\psum$,  \cref{mcat_axioms,vu_star,cglgdom_psix,ricxa,psiinvj_ord,ell_prime,cgglcnbeta_def}, we define the $\objx$-component pointed natural isomorphism
\begin{equation}\label{cgglcx}
\begin{tikzpicture}[vcenter]
\def\v{-1.4}
\draw[0cell]
(0,0) node (a11) {\Comegmal}
(a11)++(2.3,0) node (a12) {\Cboxm}
(a11)++(0,\v) node (a21) {\Comegnbeta}
(a12)++(0,\v) node (a22) {\Cboxn}
;
\draw[1cell=.85]
(a11) edge node {\cgglcmal} (a12)
(a12) edge node {\psi_*} (a22)
(a11) edge[shorten >=-.3ex] node[swap] {\objx_*} (a21)
(a21) edge node[swap] {\cgglcnbeta} (a22)
;
\draw[2cell]
node[between=a12 and a21 at .5, shift={(.25,-.1)}, rotate=-115, 2labelw={below,\cgglcx,0pt}] {\Rightarrow}
;
\end{tikzpicture}
\end{equation}
whose component at an object $\obja = \ang{a_i}_{i \in \ufsm}$ in $\Comegmal$ is the following isomorphism in $\Cboxn$.
\begin{equation}\label{cgglcxa}
\begin{tikzpicture}[vcenter]
\def\h{-1.3}
\draw[0cell=.9]
(0,0) node (a1) {\psi_*\cgglcmal\obja = \bang{\,\txsum_{\ell \in \ufs{|\psiinv j|}} \, \upps^{\ell',\ome}_{\mal,\mstar} a_{\ell'}}_{j \in \ufsn}}
(a1)++(0,-1.5) node (a2) {\cgglcnbeta\objx_*\obja = \bang{\,\txsum_{\ell \in \ufs{|\psiinv j|}} \, \upps^{j,\ome}_{\nbeta,\mstar} \iphijlomest a_{\ell'}}_{j \in \ufsn}}
(a1)++(\h,0) node (a1') {\phantom{\cgglcmal}}
(a2)++(\h,0) node (a2') {\phantom{\cgglcmal}}
;
\draw[1cell=.8]
(a1') edge[shorten >=-.2ex] node {\ang{\,\txsum_{\ell \in \ufs{|\psiinv j|}} \, [\uppsjome_{\nbeta} \iphijlome , \upps^{\ell',\ome}_{\mal}]_\mstar^{a_{\ell'}}}_{j \in \ufsn}} node[swap] {\cgglc_{\objx,\obja} \,=} (a2')
;
\end{tikzpicture}
\end{equation}
In \cref{cgglcxa}, the $G$-bijections
\[\malomeg \fto[\iso]{\upps_{\mal}} \omeg \andspace \nbetaomeg \fto[\iso]{\upps_{\nbeta}} \omeg\]
are defined in \cref{upps_iso}.  An argument similar to the proof of \cref{cglgc_welldef} proves that the data $\cgglc$---consisting of the pointed $G$-functors $\cgglcnbeta$ \cref{cgglcnbeta_functor} and the pointed natural isomorphisms $\cgglcx$ \cref{cgglcx}---define a $\Gcatst$-pseudotransformation \pcref{def:dgcatg_onecell}.  

However, $\cgglc$ is not $\Pig$-strict \pcref{def:pig_strict}.  The reason is that, even if $|\psiinv j|=1$ and $\iphi_j = 1_{\omeg} \in \Iopg(1)$, the injections \cref{uppsjome}
\[\ome \fto{\upps^{\psiinv j,\ome}_{\mal}} \ome \andspace 
\ome \fto{\uppsjome_{\nbeta}} \ome\]
are not generally equal.  Thus, the domain and codomain of $\cgglc_{\objx,\obja}$---namely, the objects $\psi_*\cgglcmal\obja$ and $\cgglcnbeta\objx_*\obja$ in \cref{cgglcxa}---are not generally equal even if $\objx$ is in the image of $\gio \cn \Pig \to \DG$ \cref{iota_psi}.
\end{explanation}

\section{$G$-Symmetric Spectra from $\FGG$-Spaces}
\label{sec:gsp_fggspaces}

This section constructs the prolongation functor 
\[\FGTopg \fto{\Kfgs} \Sptopg\] 
from the category of $\FGG$-spaces \pcref{def:ggtopg} to the category of $G$-symmetric spectra (\cref{def:sptop} \cref{def:sptop_iii}).  The functor $\Kfgs$ is equal to the composite
\begin{equation}\label{KfgGspu}
\FGTopg \fto{\Kfg} \Gspec \fto{\Gspu} \Sptopg
\end{equation}
of the functors in \cref{def:Kfg_functor,gsp_forget}.  The functor $\Kfgs$ is used in \cref{sec:kgl_kgmmo_eq} to compare GMMO $K$-theory and global $K$-theory.

\secoutline
\begin{itemize}
\item \cref{def:kfgs_obj} defines the object assignment of $\Kfgs$.
\item \cref{def:kfgs_mor} defines the morphism assignment of $\Kfgs$.
\item \cref{def:kfgs_functor} defines the functor $\Kfgs$.
\item \cref{fgtop_fgtopg_eq,kfgs_kfgsi_eq} prove that $\Kfgs$ is naturally isomorphic to the prolongation functor $\Kfgsi$ with domain $\FGTop$.
\item \cref{Xjin_top} provides an alternative characterization of objects in $\FGTopg$.
\item \cref{Lg_clast_com} proves that the classifying space functor $\cla$ \cref{clast_ggcatg} is compatible with the equivalences 
\[\FGCat \fto{\Lg} \FGCatg \andspace \FGTop \fto{\Lg} \FGTopg\]
in \cref{thm:fgcat_fgcatg_iieq,fgtop_fgtopg_eq}.
\item \cref{thm:Kfgs_ht} records the fact that $\Kfgs$ preserves componentwise weak $G$-equivalences between proper $\FGG$-spaces \pcref{def:proper_fgg}.
\end{itemize}

\subsection*{The Functor $\Kfgs$}

Recall the pointed $G$-category $\Topgst$ \cref{topgst_gtopst_enr} of pointed $G$-spaces and pointed morphisms with the conjugation $G$-action \cref{ginv_h_g}.

\begin{definition}[$\Kfgs$ on Objects]\label{def:kfgs_obj}
Given a finite group $G$ and an $\FGG$-space 
\[(\FG,\ord{0}) \fto{X} (\Topgst,*),\]
the $G$-symmetric spectrum $\Kfgs X$ is defined as follows.
\begin{description}
\item[Pointed $G$-spaces]
$\Kfgs X$ sends each finite set $A$ to the pointed $G$-space
\begin{equation}\label{Kfgsxa}
(\Kfgs X)_A = \int^{\mal \sins \FG} (S^A)^{\mal} \sma X\mal.
\end{equation}
\begin{description}
\item[Coend]
The coend in \cref{Kfgsxa} is taken in the category $\Topst$ of pointed spaces and pointed morphisms.  Each pointed finite $G$-set $\mal \in \FG$ is regarded as a discrete pointed $G$-space.  The $A$-sphere $S^A$ \pcref{def:asphere} is equipped with the trivial $G$-action.  The pointed $G$-space 
\[(S^A)^{\mal} = \Topgst(\mal, S^A)\]
consists of pointed morphisms $\mal \to S^A$ \cref{Gtopst_smc}, with $G$ acting by conjugation \cref{ginv_h_g}.
\item[$G$-action]
The group $G$ acts diagonally on representatives.  This means that, for an element $g \in G$ and a representative pair
\begin{equation}\label{Kfgsxa_rep}
\big(\mal \fto{\upom} S^A ; x \in X\mal \big) \in (S^A)^{\mal} \ttimes X\mal
\end{equation}
in $(\Kfgs X)_A$, the diagonal $g$-action is given by
\begin{equation}\label{Kfgsxa_rep_gact}
g \cdot (\upom; x) = (\upom\ginv ; gx).
\end{equation}
\end{description}
\item[Structure $G$-morphisms]
For each injection $i \cn A \to B$ between finite sets, the structure morphism $i_*$ for $\Kfgs X$ is the following composite of pointed $G$-morphisms.
\begin{equation}\label{Kfgsx_istar}
\begin{tikzpicture}[vcenter]
\def\h{4.7} \def\u{-1} \def\v{-1.4}
\draw[0cell=.9]
(0,0) node (a11) {(\Kfgs X)_A \sma S^{B \setminus i(A)}}
(a11)++(\h,0) node (a12) {(\Kfgs X)_B}
(a11)++(0,\u) node (a21) {\big( \txint^{\mal \in \FG} (S^A)^{\mal} \sma X\mal \big) \sma S^{B \setminus i(A)}}
(a12)++(0,\u) node (a22) {\txint^{\mal \in \FG} (S^B)^{\mal} \sma X\mal}
(a21)++(0,\v) node (a3) {\txint^{\mal \in \FG} \big( (S^A)^{\mal} \sma S^{B \setminus i(A)} \big) \sma X\mal}
;
\draw[1cell=.8]
(a11) edge node {i_*} (a12)
(a11) edge[equal] (a21)
(a12) edge[equal] (a22)
(a21) edge node[swap] {\iso} (a3)
(a3) [rounded corners=2pt] -| node[pos=.25] {\txint^{\mal} \asm_{\mal} \sma 1_{X\mal}} (a22)
;
\end{tikzpicture}
\end{equation}
\begin{itemize}
\item The pointed $G$-homeomorphism $\iso$ uses the commutation of $- \sma S^{B \setminus i(A)}$ with coends, along with the associativity and braiding for $\sma$.
\item For $\mal \in \FG$, the pointed $G$-morphism 
\begin{equation}\label{assem_asph}
(S^A)^{\mal} \sma S^{B \setminus i(A)} \fto{\asm_{\mal}} (S^B)^{\mal}
\end{equation}
is defined by the assignment
\[\begin{split}
& \big(\mal \fto{\upom} S^A ; y \big) \in (S^A)^{\mal} \sma S^{B \setminus i(A)}\\
&\mapsto \big(\mal \fto{\upom} S^A \fto{(- \comp i^{-1}, y)} S^{i(A)} \sma S^{B \setminus i(A)} \iso S^B \big) \in (S^B)^{\mal}.
\end{split}\]
\end{itemize}
\end{description}
The unity axiom \cref{sptop_unity} and the associativity axiom \cref{sptop_assoc} for a $G$-symmetric spectrum follow from the corresponding properties for $\asm_{\mal}$, as stated in \cref{assm_onea,assm_assoc}.  This finishes the definition of the $G$-symmetric spectrum $\Kfgs X$.
\end{definition}

\begin{definition}[$\Kfgs$ on Morphisms]\label{def:kfgs_mor}
For a finite group $G$ and a $G$-natural transformation 
\begin{equation}\label{tha_XY}
\begin{tikzpicture}[vcenter]
\def\t{28}
\draw[0cell]
(0,0) node (a1) {\phantom{\Gsk}}
(a1)++(1.8,0) node (a2) {\phantom{\Gsk}}
(a1)++(-.08,0) node (a1') {\FG}
(a2)++(.23,0) node (a2') {\Topgst}
;
\draw[1cell=.9]
(a1) edge[bend left=\t] node {X} (a2)
(a1) edge[bend right=\t] node[swap] {Y} (a2)
;
\draw[2cell]
node[between=a1 and a2 at .42, rotate=-90, 2label={above,\tha}] {\Rightarrow}
;
\end{tikzpicture}
\end{equation}
between $\FGG$-spaces $X$ and $Y$, the morphism between $G$-symmetric spectra
\begin{equation}\label{Kfgs_tha}
\Kfgs X \fto{\Kfgs\tha} \Kfgs Y
\end{equation}
sends each finite set $A$ to the following pointed $G$-morphism.
\begin{equation}\label{Kfgs_tha_a}
\begin{tikzpicture}[vcenter]
\def\v{-1.4}
\draw[0cell=.9]
(0,0) node (a11) {(\Kfgs X)_A}
(a11)++(3,0) node (a12) {\txint^{\mal \sins \FG} (S^A)^{\mal} \sma X\mal}
(a11)++(0,\v) node (a21) {(\Kfgs Y)_A}
(a12)++(0,\v) node (a22) {\txint^{\mal \sins \FG} (S^A)^{\mal} \sma Y\mal}
;
\draw[1cell=.8]
(a11) edge[equal] (a12)
(a21) edge[equal] (a22)
(a11) edge[transform canvas={xshift=1em}] node[swap] {(\Kfgs\tha)_A} (a21)
(a12) edge[transform canvas={xshift=-2em}, shorten <=-.5ex] node {\txint^{\mal} 1 \sma \tha_{\mal}} (a22)
;
\end{tikzpicture}
\end{equation}
The $G$-equivariance of $\Kfgs\tha$ follows from \cref{Kfgsxa_rep_gact} and the $G$-equivariance of each component $\tha_{\mal}$ \cref{ggtopg_icell_geq}.  The compatibility \cref{sptop_mor} of $\Kfgs\tha$ with the structure $G$-morphisms of $\Kfgs X$ and $\Kfgs Y$ follows from \cref{Kfgsx_istar}, \cref{Kfgs_tha_a}, the universal property of coends, and the functoriality of $\sma$.
\end{definition}

\begin{definition}\label{def:kfgs_functor}
For a finite group $G$, the functor
\[\FGTopg \fto{\Kfgs} \Sptopg\]
is defined by
\begin{itemize}
\item the object assignment $X \mapsto \Kfgs X$ \pcref{def:kfgs_obj} and
\item the morphism assignment $\tha \mapsto \Kfgs\tha$ \pcref{def:kfgs_mor}.
\end{itemize}
The functoriality of $\Kfgs$ follows from \cref{Kfgs_tha_a} and the fact that identity morphisms and composition are defined componentwise in $\FGTopg$ and levelwise in $\Sptopg$.
\end{definition}

\subsection*{Properties of $\Kfgs$}
The rest of this section records several facts about $\Kfgs$ and its domain category.  First, the functor $\Kfgsi \cn \FGTop \to \Sptopg$ \cref{kfgsi_functor} factors through $\Kfgs$ and the left adjoint $\Lg$ in \cref{fgtop_fgtopg_eq}.  Recall the full subcategory inclusion $\ig \cn \Fsk \to \FG$ \cref{ig_FG}.  \cref{fgtop_fgtopg_eq} is \cite[Theorem 1]{shimakawa91}.

\begin{lemma}\label{fgtop_fgtopg_eq}
For each group $G$, there is an adjoint equivalence
\begin{equation}\label{LgigstTop}
\begin{tikzpicture}[vcenter]
\draw[0cell]
(0,0) node (a1) {\FGTop}
(a1)++(2.5,0) node (a2) {\FGTopg}
;
\draw[1cell=.9]
(a1) edge[transform canvas={yshift=.5ex}] node {\Lg} (a2)
(a2) edge[transform canvas={yshift=-.4ex}] node {\igst} (a1)
;
\end{tikzpicture}
\end{equation}
between the categories in \cref{def:fgtop,def:ggtopg}, where the right adjoint $\igst$ is given by precomposition with $\ig$.
\end{lemma}

\begin{proof}
The proof of \cref{thm:fgcat_fgcatg_iieq} is applicable by replacing the pair $(\Gcatst, \Catgst)$ with the pair $(\Gtopst, \Topgst)$.
\end{proof}

\begin{lemma}\label{kfgs_kfgsi_eq}
For a finite group $G$, there is a natural isomorphism
\begin{equation}\label{kiso_iicell}
\begin{tikzpicture}[vcenter]
\def\h{2} \def\v{-1.3} \def\t{15}
\draw[0cell]
(0,0) node (a1) {\FGTopg}
(a1)++(0,\v) node (a2) {\FGTop}
(a1)++(\h,\v/2) node (a3) {\phantom{\Gspec}}
(a3)++(-.08*\h,0) node (a3') {\Sptopg}
;
\draw[1cell=.9]
(a1) edge[bend left=\t] node {\Kfgs} (a3)
(a2) edge[transform canvas={xshift=-1ex}] node {\Lg} (a1)
(a2) edge[bend right=\t] node[swap] {\Kfgsi} (a3)
;
\draw[2cell=.9]
node[between=a1 and a2 at .6, shift={(.35*\h,0)}, rotate=-50, 2label={above,\kiso}] {\Rightarrow}
;
\end{tikzpicture}
\end{equation}
between the functors $\Kfgs\Lg$ and $\Kfgsi$.
\end{lemma}

\begin{proof}
The proof of \cref{KgKgg} is applicable by replacing the pair $(\Kg,\Kgg)$ with the pair $(\Kfgsi,\Kfgs)$ and restricting to length-1 objects in $\GG$ and $\Gsk$.
\end{proof}

Objects in $\FGTopg$ can be described as follows.

\begin{definition}\label{def:gnbeta_top}
Given a pointed $G$-functor $X \cn \FG \to \Topgst$ and a pointed finite $G$-set $\nbeta \in \FG$, define the pointed $G$-space $(X\ordn)_\be$ as follows.
\begin{itemize}
\item The underlying pointed space of $(X\ordn)_\be$ is the underlying pointed space of $X\ordn$, where $\ordn \in \FG$ is equipped with the trivial $G$-action.
\item For each $g \in G$, the $g$-action functor on $(X\ordn)_\be$ is the composite
\begin{equation}\label{Xnbeta_g}
\begin{tikzpicture}[vcenter]
\def\u{.65}
\draw[0cell]
(0,0) node (a1) {X\ordn}
(a1)++(1.8,0) node (a2) {X\ordn}
(a2)++(2,0) node (a3) {X\ordn}
;
\draw[1cell=.9]
(a1) edge node {g} (a2)
(a2) edge node {X(\be g)} (a3)
(a1) [rounded corners=2pt] |- ($(a2)+(-1,\u)$) -- node {g \cdot -} ($(a2)+(1,\u)$) -| (a3)
;
\end{tikzpicture}
\end{equation}
of the $g$-action on $X\ordn$ and the image under $X$ of the pointed bijection $\be g \cn \ordn \fiso \ordn$.
\end{itemize}
Moreover, denote by $\jin \cn \ordn \fiso \nbeta$\label{not:jin_nbeta} the isomorphism in $\FG$ given by the identity morphism on $\ordn$.
\end{definition}

\cref{Xjin_top} is \cite[Prop.\ 2]{shimakawa91} and the topological analogue of \cref{Xjin} with the same proof. 

\begin{lemma}\label{Xjin_top}
In the context of \cref{def:gnbeta_top}, the pointed homeomorphism
\[(X\ordn)_\be \fto[\iso]{X\jin} X\nbeta\]
is $G$-equivariant.
\end{lemma}

\cref{Lg_clast_com} shows that the equivalences $\Lg$ in \cref{thm:fgcat_fgcatg_iieq,fgtop_fgtopg_eq} are compatible with the classifying space functor $\cla$ \cref{clast_ggcatg}.

\begin{lemma}\label{Lg_clast_com}
For a finite group $G$, the diagram of functors
\begin{equation}\label{clastL}
\begin{tikzpicture}[vcenter]
\def\v{-1.4}
\draw[0cell]
(0,0) node (a11) {\FGCatg}
(a11)++(2.5,0) node (a12) {\FGTopg}
(a11)++(0,\v) node (a21) {\FGCat}
(a12)++(0,\v) node (a22) {\FGTop}
;
\draw[1cell=.9]
(a21) edge node {\Lg} (a11)
(a11) edge node {\clast} (a12)
(a21) edge node {\clast} (a22)
(a22) edge node[swap] {\Lg} (a12)
;
\end{tikzpicture}
\end{equation}
commutes up to a natural isomorphism, in which each $\clast$ is given by postcomposition with the classifying space functor $\cla$.
\end{lemma}

\begin{proof}
By \cref{LX_nbeta,LXnbeta_gaction,Xjin,Xjin_top}, up to natural isomorphisms, each of the two composites $\Lg\clast$ and $\clast\Lg$ sends an $\Fskg$-category $X \cn \Fsk \to \Gcatst$ to the $\FGG$-space whose value at a pointed finite $G$-set $\nbeta \in \FG$ is the underlying pointed space of $\cla X\ordn$ with the composite $g$-action
\[\cla X\ordn \fto{\cla g} \cla X\ordn \fto{\cla X(\be g)} \cla X\ordn\]
for $g \in G$.  This $\FGG$-space sends a pointed morphism $\psi \cn \mal \to \nbeta$ between pointed finite $G$-sets to the pointed $G$-morphism 
\[\cla X\ordm \fto{\cla X\psi} \cla X\ordn\]
in which $\psi \cn \ordm \to \ordn$ is regarded as a pointed morphism in $\Fsk$.
\end{proof}

\cref{thm:Kfgs_ht} states that $\Kfgs$ preserves componentwise weak $G$-equivalences between proper $\FGG$-spaces \pcref{def:proper_fgg}.  Its proof is the same as that of \cref{thm:Kfg_inv} with $G$-symmetric spectra in place of orthogonal $G$-spectra.

\begin{theorem}\label{thm:Kfgs_ht}
Suppose $\tha \cn X \to Y$ is a morphism of $\FGG$-spaces such that the following two conditions hold.
\begin{itemize}
\item Each of $X$ and $Y$ is a proper $\FGG$-space.
\item $\tha_{\mal} \cn X\mal \to Y\mal$ is a weak $G$-equivalence for each pointed finite $G$-set $\mal \in \FG$.
\end{itemize}
Then for each finite set $A$, the $A$-component pointed $G$-morphism \cref{Kfgs_tha_a}
\[(\Kfgs X)_A \fto{(\Kfgs\tha)_A} (\Kfgs Y)_A\]
is a weak $G$-equivalence.
\end{theorem}

\section{Comparison $G$-Stable Equivalences}
\label{sec:kgl_kgmmo_eq}

This section proves that there are natural $G$-stable equivalences connecting
\begin{itemize} 
\item the underlying $G$-spectrum $(\Ksc\C)_G$ of the global $K$-theory $\Ksc\C$ of a parsummable category $\C$ \pcref{def:parcat,def:sch3.3} and
\item the GMMO $K$-theory $\Kiopgb\Ig\C$ of the $\Einfg$-category $\Ig\C = \Comeg$ \pcref{def:global_einf,def:comeg_iopg,def:Ig,def:kiopgb}.
\end{itemize}
See \cref{gl_gmmo_geq}.  The key ingredient of the proof of \cref{gl_gmmo_geq} is the comparison natural transformation $\cglg$ \pcref{def:cglg}; see the diagram \cref{kiopg_ksc_diag}.  Throughout this section, $G$ denotes a finite group.

\secoutline
\begin{itemize}
\item \cref{def:kiopgb} defines a version of GMMO $K$-theory, denoted $\Kiopgb$, that sends $\Einfg$-categories to $G$-symmetric spectra.  \cref{expl:kiopgb} unpacks $\Kiopgb$.
\item \cref{gl_gmmo_geq} proves that the $G$-symmetric spectra $(\Ksc\C)_G$  and $\Kiopgb\Ig\C$ are naturally $G$-stably equivalent for each parsummable category $\C$.
\item \cref{expl:gl_gmmo_geq} unpacks the natural $G$-stable equivalences connecting $(\Ksc\C)_G$ and $\Kiopgb\Ig\C$ at the point-set level.  
\end{itemize}

\subsection*{GMMO $K$-Theory}
Recall that the forgetful functor $\Gspu \cn \Gspec \to \Sptopg$ \cref{gsp_forget} from the category of orthogonal $G$-spectra to the category of $G$-symmetric spectra \pcref{def:gsp_morphism,def:sptop} is the right adjoint of a Quillen equivalence.

\begin{definition}\label{def:kiopgb}\index{GMMO K-theory@GMMO $K$-theory!for $\Iopg$}
For a finite group $G$, define the functor
\begin{equation}\label{kiopgb_def}
\begin{tikzpicture}[vcenter]
\def\h{2.3} \def\u{.7}
\draw[0cell]
(0,0) node (a1) {\Algiopg}
(a1)++(\h,0) node (a2) {\Algpsiopg}
(a2)++(1.1*\h,0) node (a3) {\Gspec}
(a3)++(.8*\h,0) node (a4) {\Sptopg}
;
\draw[1cell=.85]
(a1) edge[right hook->] node {\algi} (a2)
(a2) edge node {\Kiopg} (a3)
(a3) edge node {\Gspu} (a4)
(a1) [rounded corners=2pt] |- ($(a2)+(0,\u)$) -- node {\Kiopgb} ($(a3)+(0,\u)$) -| (a4)
;
\end{tikzpicture}
\end{equation}
as the composite of
\begin{itemize}
\item the subcategory inclusion 
\[\Algiopg \fto{\algi} \Algpsiopg\]
for the chaotic $\Einfg$-operad $\Iopg$ \pcref{def:algstpsO,def:global_einf};
\item GMMO $K$-theory 
\[\Algpsiopg \fto{\Kiopg} \Gspec\]
for $\Iopg$ \cref{kgmmo_diag}; and
\item the forgetful functor $\Gspu$ \cref{gsp_forget}.
\end{itemize}
We also refer to $\Kiopgb$ as \emph{GMMO $K$-theory} for $\Iopg$. 
\end{definition}

\begin{explanation}[Unpacking]\label{expl:kiopgb}
Using \cref{kgmmo_diag}, the functor $\Kiopgb$ \cref{kiopgb_def} is the following boundary composite.
\begin{equation}\label{kiopgb_diag}
\begin{tikzpicture}[vcenter]
\def\v{1.4} \def\h{2.8} \def\u{1}
\draw[0cell]
(0,0) node (a0) {\Algiopg}
(a0)++(0,-\v) node (a1) {\Algpsiopg}
(a1)++(0,-\v) node (a2) {\DGCatg}
(a2)++(0,-\v) node (a3) {\FGCatgps}
(a3)++(\h,0) node (a4) {\FGCatg}
(a4)++(\h,0) node (a5) {\FGTopg}
(a5)++(0,\v) node (a6) {\FGTopg}
(a6)++(0,\v) node (a7) {\Gspec}
(a7)++(0,\v) node (a8) {\Sptopg}
;
\draw[1cell=.9]
(a0) edge node {\Kiopgb} (a8)
(a1) edge node {\Kiopg} (a7)
(a0) edge node[swap] {\algi} (a1)
(a1) edge node[swap] {\Rg} (a2)
(a2) edge node[swap] {\gzest} (a3)
(a3) edge node {\str} (a4)
(a4) edge node {\clast} (a5)
(a5) edge node[swap] {\Bc} (a6)
(a6) edge node[swap] {\Kfg} (a7)
(a7) edge node[swap] {\Gspu} (a8)
(a6) [rounded corners=2pt] -| ($(a7)+(\u,-1)$) -- node[swap] {\Kfgs} ($(a7)+(\u,1)$) |- (a8) 
;
\end{tikzpicture}
\end{equation}
\begin{itemize}
\item $\Kiopgb$ goes through the categories $\Algiopg$, $\Algpsiopg$ \pcref{def:algstpsO}, $\DGCatg$ \pcref{expl:dgiopg}, $\FGCatgps$ \pcref{def:psfggcat_iicat}, $\FGCatg$ \pcref{def:fgcatg}, $\FGTopg$ \pcref{def:ggtopg}, $\Gspec$ \pcref{def:gsp_morphism}, and $\Sptopg$ \pcref{def:sptop}.
\item $\Kiopgb$ consists of the functors $\algi$, $\Rg$ \cref{rg_twofunctor}, $\gzest$ \cref{gzest}, $\str$ \cref{str_iifunctor}, $\clast$ \cref{clast_ggcatg}, $\Bc$ \cref{bar_functor_FG}, and $\Kfgs = \Gspu\Kfg$ \cref{KfgGspu}.  In the rest of this section, $\Rg\algi$ is abbreviated to $\Rg$.
\end{itemize}
By \cref{Kfgsxa,real_FG}, for each $\Iopg$-algebra $\C$ \cref{iopgn}, the $G$-symmetric spectrum $\Kiopgb\C$ sends each finite set $A$ to the pointed $G$-space
\begin{equation}\label{kiopgbca}
(\Kiopgb\C)_A = 
\int^{\mal \sins \FG} (S^A)^{\mal} \sma \Bc\big(\FG(-,\mal),\FG,\clast\str\gzest\Rg\C\big).
\end{equation}
\begin{itemize}
\item The structure morphisms of $\Kiopgb\C$ involve only the factor $(S^A)^{\mal}$ \cref{Kfgsx_istar}. 
\item The $G$-action on $(\Kiopgb\C)_A$ is defined in \cref{Kfgsxa_rep_gact,Br_FGX,strx_objx_gact,morstrx,proCnbe_gaction}. 
\end{itemize} 
By \cref{rga_mal,gzest_xmal}, the pseudo $\FGG$-category 
\[\gzest\Rg\C \in \FGCatgps\] 
sends each pointed finite $G$-set $\nbeta$ to the $\nbeta$-twisted product \pcref{def:proCnbe}
\[(\gzest\Rg\C)\nbeta = (\Rg\C)\nbeta = \proCnbeta.\]
The strictified $\FGG$-category $\str\gzest\Rg\C$ is discussed in \cref{expl:str_obj} with $X = \gzest\Rg\C$.  The $\FGG$-space $\clast\str\gzest\Rg\C$ sends $\nbeta$ to the pointed $G$-space
\[(\clast\str\gzest\Rg\C)\nbeta = \cla (\str\gzest\Rg\C)\nbeta,\]
which is the classifying space \cref{classifying_space} of the pointed $G$-category $(\str\gzest\Rg\C)\nbeta$.
\end{explanation}

\subsection*{Comparing Global and GMMO $K$-Theories}
Recall that the underlying $G$-spectrum $X_G$ of a symmetric spectrum $X$ is equipped with the trivial $G$-action.  Also recall $G$-stable equivalences and $\pistu$-isomorphisms between $G$-symmetric spectra \pcref{def:sptop}.  By \cite[3.36]{hausmann17}, each $\pistu$-isomorphism is also a $G$-stable equivalence.  \cref{gl_gmmo_geq} is the main result of this chapter.  It establishes a natural zigzag of $G$-stable equivalences between global $K$-theory at a finite group $G$ and GMMO $K$-theory for the chaotic $\Einfg$-operad $\Iopg$ \cref{iopgn}.

\begin{theorem}\label{gl_gmmo_geq}\index{GMMO K-theory@GMMO $K$-theory!comparison with global}\index{global K-theory@global $K$-theory!comparison with GMMO}
For a finite group $G$ and a parsummable category $\C$, there are $\pistu$-isomorphisms
\begin{equation}\label{gl_gmmo_geqdiag}

\end{equation}
that are natural in $\C$.  Thus, the $G$-symmetric spectra $(\Ksc\C)_G$ and $\Kiopgb\Ig\C$ are naturally $G$-stably equivalent.
\end{theorem}

\begin{proof}
We explain the diagram \cref{gl_gmmo_geqdiag}, starting at the upper-left corner and proceeding counterclockwise.  The left column of \cref{gl_gmmo_geqdiag} consists of the natural $\pistu$-isomorphisms $\ascgc$, $\bscgc$, and $\Isc^\C$ in \cref{sch4.15,Isc_geq}.  The rest of \cref{gl_gmmo_geqdiag} uses the diagram 
\begin{equation}\label{kiopg_ksc_diag}
%
\end{equation}
involving
\begin{itemize}
\item the categories in \cref{def:fgcat,def:fgcatg,def:parcat,def:ggtopg,def:fgtop,def:sptop,def:global_einf,expl:dgiopg,def:psfggcat_iicat};
\item Schwede $K$-theory \cref{kscgb_functor}
\[\Kscgb = \Kfgsi\clast\Jscg \cn \Parcat \to \Sptopg;\]
\item the functor \cref{Ig_functor}
\[\Parcat \fto{\Ig} \Algiopg;\] 
and
\item GMMO $K$-theory for $\Iopg$ \cref{kiopgb_diag}
\[\Algiopg \fto{\Kiopgb} \Sptopg.\]
\end{itemize}
In the diagram \cref{gl_gmmo_geqdiag}, the bottom horizontal natural isomorphisms 
\[\Kfgsi \fot{\iso} \Kfgs\Lg \andspace \Lg\clast \fot{\iso} \clast\Lg\]
are from \cref{kfgs_kfgsi_eq,Lg_clast_com}.  They correspond to the two regions in \cref{kiopg_ksc_diag} decorated by $\iso$.  

The natural $\pistu$-isomorphisms in the right column in \cref{gl_gmmo_geqdiag} are given as follows.
\begin{enumerate}
\item
In the lower-right arrow, each component of the counit \cref{strv} 
\[\str\Incj \fto{\strv} 1_{\FGCatg}\]
is componentwise a pointed $G$-equivalence between pointed $G$-categories.  Applying the functor $\clast$ \cref{clast_ggcatg}, the morphism
\[\clast\str\Incj\Lg\Jscg\C \fto{\clast\strv_{\Lgsub\Jscg\C}} \clast\Lg\Jscg\C \inspace \FGTopg\]
is componentwise a pointed $G$-homotopy equivalence between pointed $G$-spaces, hence also a weak $G$-equivalence \pcref{def:weakG_top}.  By \cref{reast_proper}, each $\FGG$-space in the image of $\clast$ is proper \pcref{def:proper_fgg}.  By \cref{thm:Kfgs_ht}, the morphism of $G$-symmetric spectra
\[\Kfgs\clast\str\Incj\Lg\Jscg\C \fto{\Kfgs\clast\strv_{\Lgsub\Jscg\C}} \Kfgs\clast\Lg\Jscg\C\]
is levelwise a weak $G$-equivalence between pointed $G$-spaces.  Thus, $\Kfgs\clast\strv_{\Lgsub\Jscg\C}$ is a $\pistu$-isomorphism.
\item
The equality in the right column in \cref{gl_gmmo_geqdiag} comes from the equality \pcref{expl:gze}
\[\FGCatg \fto{\Incj \,=\, \gzest\gxist} \FGCatgps.\]
This equality corresponds to the lower-left triangle in \cref{kiopg_ksc_diag}.
\item
In the second-to-top arrow in the right column in \cref{gl_gmmo_geqdiag}, the $\Pig$-strict $\Gcatst$-pseudotransformation \cref{cglgc}
\[\gxist\Lg\Jscg\C \fto{\cglgc} \Rg\Ig\C \inspace \DGCatg\]
is natural in $\C$ \pcref{cglg_natural} and componentwise the left adjoint of an adjoint pointed $G$-equivalence \pcref{cglgc_equiv}.  It corresponds to the upper-left region in \cref{kiopg_ksc_diag}.  Since $\gzest$ \cref{gzest_icell_comp} and $\str$ \cref{strtham_H} preserve componentwise categorical weak $G$-equivalences \pcref{def:cat_weakg}, the morphism
\[\str\gzest\gxist\Lg\Jscg\C \fto{\str\gzest\cglgc} \str\gzest\Rg\Ig\C \inspace \FGCatg\]
is componentwise a categorical weak $G$-equivalence.  Applying the functor $\clast$ \cref{clast_ggcatg}, the morphism
\[\clast\str\gzest\gxist\Lg\Jscg\C \fto{\clast\str\gzest\cglgc} \clast\str\gzest\Rg\Ig\C \inspace \FGTopg\]
is componentwise a weak $G$-equivalence between pointed $G$-spaces.  By \cref{reast_proper,thm:Kfgs_ht}, the morphism of $G$-symmetric spectra
\[\Kfgs\clast\str\gzest\gxist\Lg\Jscg\C \fto{\Kfgs\clast\str\gzest\cglgc} \Kfgs\clast\str\gzest\Rg\Ig\C\]
is levelwise a weak $G$-equivalence between pointed $G$-spaces.  Thus, $\Kfgs\clast\str\gzest\cglgc$ is a $\pistu$-isomorphism.
\item
In the top arrow in the right column in \cref{gl_gmmo_geqdiag}, the retraction \cref{ret_FG}
\[\Bc\clast\str\gzest\Rg\Ig\C \fto{\retn_{\clast\str\gzest\Rg\Ig\C}} \clast\str\gzest\Rg\Ig\C 
\inspace \FGTopg\]
is componentwise a pointed $G$-homotopy equivalence between pointed $G$-spaces.  This retraction corresponds to the lower-right region in \cref{kiopg_ksc_diag}.  By \cref{bar_proper,reast_proper}, the domain and codomain of $\retn_{\clast\str\gzest\Rg\Ig\C}$ are proper $\FGG$-spaces.  By \cref{thm:Kfgs_ht}, the morphism of $G$-symmetric spectra
\[\Kfgs\Bc\clast\str\gzest\Rg\Ig\C 
\fto{\Kfgs\retn_{\clast\str\gzest\Rg\Ig\C}} \Kfgs\clast\str\gzest\Rg\Ig\C\]
is levelwise a weak $G$-equivalence between pointed $G$-spaces.  Thus, $\Kfgs\retn_{\clast\str\gzest\Rg\Ig\C}$ is a $\pistu$-isomorphism.
\end{enumerate}
In summary, the diagram \cref{gl_gmmo_geqdiag} consists entirely of natural $\pistu$-isomorphisms and natural isomorphisms, hence also natural $G$-stable equivalences, between $G$-symmetric spectra.
\end{proof}

\begin{explanation}[Explicit Description]\label{expl:gl_gmmo_geq}
We unpack the equivalence \cref{gl_gmmo_geqdiag} between global and GMMO $K$-theories in two steps.  
\begin{description}
\item[Left-bottom part of \cref{gl_gmmo_geqdiag}]
Evaluated at a nonempty finite set $A$, the $\pistu$-isomorphisms in the left-bottom part of the diagram \cref{gl_gmmo_geqdiag} are displayed in the left column in the diagram \cref{kgl_gmmo_geqa}, using \cref{KscCa,gacomeg_a,kcga,kscgbca,Kfgsxa,LgJscg_nbeta}.  In \cref{kgl_gmmo_geqa}, each pointed $G$-morphism in the left column is induced by the corresponding pointed morphism in the right column.
\begin{equation}\label{kgl_gmmo_geqa}
\begin{tikzpicture}[vcenter]
\def\v{-1.4} \def\s{.85} \def\t{.8}
\draw[0cell=\s]
(0,0) node (a1) {(\Ksc\C)_{G,A} = \txint^{\ordn \in \Fsk} (S^A)^n \sma \cla\Cboxnomea}
(a1)++(0,\v) node (a2) {\gacomegs_A = \txint^{\ordn \in \Fsk} (S^A)^n \sma \cla\Cboxnomega}
(a2)++(0,\v) node (a3) {(\kcg)_A = \txint^{\ordn \in \Fsk} (S^A)^n \sma \cla\Cboxnomeg}
(a3)++(0,\v) node (a4) {(\Kscgb\C)_A = \txint^{\ordn \in \Fsk} (S^A)^n \sma \cla\Cboxnomeg}
(a4)++(0,\v) node (a5) {(\Kfgs\clast\Lg\Jscg\C)_A = \txint^{\nbeta \in \FG} (S^A)^{\nbeta} \sma \cla\Cboxn}
;
\draw[0cell=\s]
(a1)++(5,0) node (c1) {\cla\Cboxnomea}
(c1)++(0,\v) node (c2) {\cla\Cboxnomega}
(c2)++(0,\v) node (c3) {\cla\Cboxnomeg}
(c3)++(0,\v) node (c4) {\cla\Cboxnomeg}
(c4)++(0,\v) node (c5) {\cla\Cboxn}
;
\draw[1cell=\t]
(a1) edge node[swap] {\ascgc} (a2)
(a3) edge node {\bscgc} (a2) 
(a3) edge node[swap] {\Isc^\C} (a4)
(a5) edge node {\iso} (a4)
;
\draw[1cell=\t]
(c1) edge node[swap] {\cla(i_2^\ome)_\mstar^{\Bxtimes n}} (c2)
(c3) edge node {\cla(i_1^\ome)_\mstar^{\Bxtimes n}} (c2)
(c3) edge node[swap] {1} (c4)
(c5) edge node {1} (c4)
;
\end{tikzpicture}
\end{equation}
From top to bottom, the pointed morphisms in the right column in \cref{kgl_gmmo_geqa} are given as follows.
\begin{itemize}
\item The pointed functors
\[\Cboxnomea \fto{(i_2^\ome)_\mstar^{\Bxtimes n}} \Cboxnomega 
\fot{(i_1^\ome)_\mstar^{\Bxtimes n}} \Cboxnomeg\]
are the $\schm$-action functors on $\Cboxn$ for the inclusions 
\[i_2^\ome, i_1^\ome \cn \ome \to \ome\] 
in \cref{ionetwoome}. 
\item The third arrow from the top denoted by 1 uses the definition \cref{IscC} of $\Isc^\C$. 
\item The fourth arrow from the top denoted by 1 uses \cref{kiso_xv} applied to $(\Kfgsi,\Kfgs)$, \cref{Xjin_top,Lg_clast_com}, and the fact that the underlying pointed spaces of
\[(S^A)^n \sma \cla\Cboxnomeg \andspace (S^A)^{\nbeta} \sma \cla\Cboxn\]
are equal. 
\end{itemize}

Suppose $A$ is empty.  Then the morphisms $\ascgc$ and $\bscgc$ in \cref{kgl_gmmo_geqa} are the identity morphism $1_{\cla\C^{\suppempty}}$.  The morphism $\Isc^\C$ for $A=\emptyset$ is induced by the full subcategory inclusion $\iota \cn \C^{\suppempty} \to \Comeg$ at the object $\ordone \in \Fsk$ \cref{Isc_empty}.  The bottom left isomorphism in \cref{kgl_gmmo_geqa} admits the same description regardless of whether $A$ is empty or nonempty.
\item[Right column of \cref{gl_gmmo_geqdiag}]
Next, the right column of the diagram \cref{gl_gmmo_geqdiag} evaluated at a finite set $A$ is displayed in the diagram \cref{kgl_gmmo_geq}, in which each arrow is a pointed $G$-morphism between pointed $G$-spaces.
\begin{equation}\label{kgl_gmmo_geq}
\begin{tikzpicture}[vcenter]
\def\u{1.5} \def\w{1} \def\s{.9} \def\t{.8}
\draw[0cell=\s]
(0,0) node (a7) {(\Kfgs\clast\Lg\Jscg\C)_A = \txint^{\nbeta \in \FG} (S^A)^{\nbeta} \sma \cla\Cboxn}
(a7)++(0,\u) node (a8) {(\Kfgs\clast\str\Incj\Lg\Jscg\C)_A = \txint^{\nbeta \in \FG} (S^A)^{\nbeta} \sma \cla(\str\Incj\Lg\Jscg\C)\nbeta}
(a8)++(0,\w) node (a9) {(\Kfgs\clast\str\gzest\gxist\Lg\Jscg\C)_A = \txint^{\nbeta \in \FG} (S^A)^{\nbeta} \sma \cla(\str\gzest\gxist\Lg\Jscg\C)\nbeta}
(a9)++(0,\u) node (a10) {(\Kfgs\clast\str\gzest\Rg\Ig\C)_A = \txint^{\nbeta \in \FG} (S^A)^{\nbeta} \sma \cla(\str\gzest\Rg\Ig\C)\nbeta}
(a10)++(0,\u) node  (a12) {(\Kiopgb\Ig\C)_A = \txint^{\nbeta \in \FG} (S^A)^{\nbeta} \sma \Bc\big(\FG(-,\nbeta), \FG, \clast\str\gzest\Rg\Ig\C\big)}
;
\draw[1cell=\t]
(a8) edge node[swap,pos=.4] {\Kfgs\clast\strv_{\Lgsub\Jscg\C}} (a7)
(a9) edge[equal,shorten <=-.7ex] (a8)
(a9) edge node {\Kfgs\clast\str\gzest\cglgc} (a10)
(a12) edge node[swap,pos=.45] {\Kfgs\retn_{\clast\str\gzest\Rg\Ig\C}} (a10)
;
\end{tikzpicture}
\end{equation}
From top to bottom, the arrows in \cref{kgl_gmmo_geq} are given as follows.  
\begin{itemize}
\item The top arrow is induced by the retraction $\retn \cn \Bc \to 1$ defined in \cref{retn,retn_cr}. 
\item The second-to-top arrow is induced by the comparison \cref{cglgc}
\[\gxist\Lg\Jscg\C \fto{\cglgc} \Rg\Ig\C\]
whose $\nbeta$-component is the full subcategory inclusion \cref{cglgcn}
\[\Cboxn \hookrightarrow \Comegnbeta.\]
\item The equal sign comes from the equality $\Incj = \gzest\gxist$ \pcref{expl:gze}.
\item The bottom arrow is induced by the counit $\strv \cn \str\Incj \to 1$ \cref{strv}.  On objects, the counit $\strv$ is defined in \cref{strvxm_obj}.  On morphisms, the counit $\strv$ is defined in \cref{morstrx} using 
\[X = \Lg\Jscg\C \andspace (\Lg\Jscg\C)\mal = \Cboxm.\]
The last equality is from \cref{LgJscg_nbeta}.
\end{itemize}
\end{description}
This finishes our explicit description of the equivalence \cref{gl_gmmo_geqdiag} between global $K$-theory and GMMO $K$-theory.
\end{explanation}

\input

%% file: chap/mon_cat.tex
This appendix reviews basic concepts of category theory \pcref{sec:category_thy}, monoidal categories \pcref{sec:sym_mon_cat}, enriched categories \pcref{sec:enrichedcat}, 2-categories \pcref{sec:twocategories}, and enriched operads \pcref{sec:operads}.

\section{Categories}
\label{sec:category_thy}

This section reviews categories, functors, natural transformations, adjunctions, (adjoint) equivalences, (co)limits, and coends.  For detailed discussion of categories, see \cite{maclane,riehl}.

\subsection*{Categories, Functors, and Natural Transformations}

\begin{definition}\label{def:category}
A \emph{category}\index{category} $(\C,\mcomp,1)$ consists of the following data and axioms.
\begin{description}
\item[Objects] $\C$ is equipped with a class of \emph{objects}\index{object} \label{not:ObC}$\Ob\C$, also denoted by $\C$.
\item[Morphisms] For each pair of objects $a,b \in \C$, $\C$ is equipped with a \emph{hom set}\index{hom set} $\C(a,b)$\label{not:Cab} of \emph{morphisms}\index{morphism} from $a$ to $b$, each denoted by $a \to b$.
\item[Composition] For objects $a,b,c \in \C$, $\C$ is equipped with a \emph{composition}\index{composition} function\label{not:mcompabc}
\[\C(b,c) \times \C(a,b) \fto{\mcomp_{a,b,c}} \C(a,c).\] 
For morphisms $h \cn b \to c$ and $f \cn a \to b$, the morphism $\mcomp_{a,b,c}(h,f)$ is denoted by $h \comp f$ or $hf$.
\item[Identities] For each object $a \in \C$, $\C$ is equipped with an \emph{identity} morphism \label{not:onea}$1_a \cn a \to a$.
\item[Axioms] These data satisfy the morphism equalities
\[\begin{split}
1_b f &= f = f 1_a \andspace\\
(f'' f')f &= f''(f' f)
\end{split}\]
whenever they are defined.
\end{description}
Suppressing $\mcomp$ and 1 from the notation, such a category is usually denoted by $\C$.  A morphism $f \cn a \to b$ is an \emph{isomorphism}\index{isomorphism} if there exists a morphism $f' \cn b \to a$ that satisfies the morphism equalities
\[f'f = 1_a \andspace ff' = 1_b.\]
The symbol $\iso$ is sometimes used to denote an isomorphism.  A category is \emph{small}\index{small category}\index{category!small} if $\Ob\C$ is a set.  It is \emph{discrete}\index{discrete category}\index{category!discrete} if it has no nonidentity morphisms.  
\end{definition}

\begin{definition}\label{def:opposites}
Given a category $\C$, the \emph{opposite category}\index{opposite category}\index{category!opposite} $(\Cop,\mcompop,\oneop)$ has the same objects as $\C$, hom set
\[\Cop(a,b) = \C(b,a) \forspace a,b \in \C,\]
identity morphisms $\oneop_a = 1_a$, and composition $\mcompop(h,f) = \mcomp(f,h)$.
\end{definition}

\begin{definition}\label{def:producat_cat}
Given categories $\C$ and $\C'$, the \emph{Cartesian product}\index{Cartesian product} $\C\times\C'$ is the category defined by 
\begin{itemize}
\item the objects $(\Ob\C) \times (\Ob\C')$;
\item the hom sets 
\[(\C\times\C')((a,a'), (b,b')) = \C(a,b) \times \C'(a',b')\] 
for objects $a,b \in \C$ and $a',b' \in \C'$;
\item the identity morphisms 
\[1_{(a,a')} = (1_a,1_{a'})\] 
for objects $a \in \C$ and $a' \in \C'$; and
\item the composition 
\[(h,h')(f,f') = (hf,h'f')\] 
for morphisms $f \in \C(a,b)$, $h \in \C(b,c)$, $f' \in \C'(a',b')$, and $h' \in \C'(b',c')$.\defmark
\end{itemize}
\end{definition}

\begin{definition}\label{def:functor}
Given categories $\C$ and $\D$, a \emph{functor}\index{functor} $F \cn \C \to \D$ consists of
\begin{itemize}
\item an object assignment $F \cn \Ob\C \to \Ob\D$ and
\item a morphism assignment $F_{a,b} \cn \C(a,b) \to \D(Fa,Fb)$ for each pair of objects $a,b \in \C$
\end{itemize}
such that the morphism equalities
\[\begin{split}
F_{a,a} 1_a &= 1_{Fa} \andspace\\
F_{a,c}(hf) &= (F_{b,c} h)(F_{a,b} f)
\end{split}\]
hold for objects $a,b,c \in \C$ and morphisms $f \cn a \to b$ and $h \cn b \to c$.  If each function $F_{a,b}$ is surjective, then $F$ is said to be \emph{full}.  If each function $F_{a,b}$ is injective, then $F$ is said to be \emph{faithful}.  Each function $F_{a,b}$ is usually abbreviated to $F$.  Identity functors and composition of functors are defined separately on object assignments and morphism assignments.  The identity functor of a category $\C$ is denoted by $1_\C$.  A functor $F \cn \C \to \D$ is an \emph{isomorphism} if there exists a functor $H \cn \D \to \C$ such that the functor equalities
\[HF = 1_{\C} \andspace FH = 1_{\D}\]
hold.
\end{definition}

\begin{definition}\label{def:nat_tr}
Given functors $F,H \cn \C\to\D$ between categories, a \index{natural transformation}\emph{natural transformation} $\phi \cn F\to H$ consists of an $a$-component\index{component} morphism $\phi_a \cn Fa \to Ha$ in $\D$ for each object $a \in \C$ such that the naturality diagram
\begin{equation}\label{nt_naturality}
\begin{tikzpicture}[vcenter]
\def\v{-1.3}
\draw[0cell]
(0,0) node (a11) {Fa}
(a11)++(2,0) node (a12) {Ha}
(a11)++(0,\v) node (a21) {Fb}
(a12)++(0,\v) node (a22) {Hb}
;
\draw[1cell=.85]
(a11) edge node {\phi_a} (a12)
(a12) edge node {Hf} (a22)
(a11) edge node[swap] {Ff} (a21)
(a21) edge node {\phi_b} (a22)
;
\end{tikzpicture}
\end{equation}
commutes for each morphism $f \cn a \to b$ in $\C$.  The \emph{identity natural transformation} $1_F \cn F \to F$ is defined by the component identity morphisms $1_a$ for $a \in \C$.  A \emph{natural isomorphism}\index{natural isomorphism} is a natural transformation $\phi$ such that each component $\phi_a$ is an isomorphism.   A natural transformation is sometimes denoted by 
\begin{equation}\label{twocellnotation}
\begin{tikzpicture}[xscale=1.8,yscale=1.7,baseline={(x1.base)}]
\draw[0cell=.9]
(0,0) node (x1) {\C}
(x1)++(1,0) node (x2) {\D}
;
\draw[1cell=.9]  
(x1) edge[bend left] node {F} (x2)
(x1) edge[bend right] node[swap] {H} (x2)
;
\draw[2cell]
node[between=x1 and x2 at .43, rotate=-90, 2label={above,\phi}] {\Rightarrow}
;
\end{tikzpicture}
\end{equation}
called the \index{2-cell!notation}2-cell notation.
\end{definition}

\begin{definition}\label{def:nt_comp}\
\begin{description}
\item[Vertical composition] 
Given functors $F,F',F'' \cn \C \to \D$ and natural transformations $\phi \cn F \to F'$ and $\phi' \cn F' \to F''$, the \emph{vertical composition}\index{vertical composition} $\phi'\phi \cn F \to F''$\label{not:vercomp} is the natural transformation with component morphism defined as the composite
\[Fa \fto{\phi_a} F'a \fto{\phi'_a} F''a\]
for each object $a \in \C$. 
\item[Horizontal composition] 
Given functors $H,H' \cn \D \to \B$ and a natural transformation $\psi \cn H \to H'$, the \emph{horizontal composition}\index{horizontal composition} $\psi \ast \phi \cn HF \to H'F'$\label{not:horcomp} is the natural transformation with component morphism defined as either composite in the commutative diagram
\begin{equation}\label{hor_comp}
\begin{tikzpicture}[vcenter]
\def\v{-1.3}
\draw[0cell]
(0,0) node (a11) {HFa}
(a11)++(2.4,0) node (a12) {\phantom{HF'a}}
(a12)++(0,.03) node (a12') {HF'a}
(a11)++(0,\v) node (a21) {H'Fa}
(a12)++(0,\v) node (a22) {H'F'a}
;
\draw[1cell=.85]
(a11) edge node {H\phi_a} (a12)
(a12') edge node {\psi_{F'a}} (a22)
(a11) edge node[swap] {\psi_{Fa}} (a21)
(a21) edge node {H'\phi_a} (a22)
;
\end{tikzpicture}
\end{equation}
for each object $a \in \C$.
\item[Whiskering]
The horizontal compositions
\[\begin{split}
H\phi &= 1_H \ast \phi \cn HF \to HF' \andspace\\
\psi F&= \psi \ast 1_F \cn HF \to H'F
\end{split}\]
are called \emph{whiskerings}.\defmark
\end{description}
\end{definition}

\subsection*{Equivalences and Adjunctions}

\begin{definition}\label{def:cat_equiv}
A functor $F \cn \C \to \D$ is an \emph{equivalence}\index{equivalence} if there exist a functor $H \cn \D \to \C$ and natural isomorphisms
\[1_{\C} \fiso HF \andspace FH \fiso 1_{\D}.\]
In this case, $H$ is called an \emph{inverse} of $F$.
\end{definition}

The following characterization of an equivalence is proved in \cite[IV.4 Theorem 1]{maclane}.  

\begin{theorem}\label{thm:cat_equiv}
A functor $F \cn \C \to \D$ is an equivalence if and only if it satisfies the following two conditions.
\begin{enumerate}
\item $F$ is essentially surjective\index{essentially surjective} on objects.  This means that, for each object $b \in \D$, there exist an object $a \in \C$ and an isomorphism $Fa \fiso b$.
\item $F$ is fully faithful\index{fully faithful} on morphisms.  This means that, for each pair of objects $a,b \in \C$, the morphism assignment $F \cn \C(a,b) \to \D(Fa,Fb)$ is a bijection.
\end{enumerate}
\end{theorem}

\begin{definition}\label{def:adjunction}
Given categories $\C$ and $\D$, an \emph{adjunction}\index{adjunction} $(L,R,u,v)$ from $\C$ to $\D$ consists of
\begin{itemize}
\item a functor $L \cn \C \to \D$ called the \index{left adjoint}\emph{left adjoint},
\item a functor $R \cn \D \to \C$ called the \index{right adjoint}\emph{right adjoint},
\item a natural transformation $u \cn 1_{\C} \to RL$ called the \index{unit}\emph{unit}, and
\item a natural transformation $v \cn LR \to 1_{\D}$ called the \index{counit}\emph{counit}
\end{itemize}
such that the equalities of natural transformations\index{triangle identities}
\[\begin{split}
(vL)(Lu) &= 1_L \andspace\\ 
(Rv)(uR) &= 1_R
\end{split}\]
hold.  These equalities are called the \emph{left triangle identity}\index{left triangle identity} and the \index{right triangle identity}\emph{right triangle identity}.  An adjunction $(L,R,u,v)$ is called an \emph{adjoint equivalence}\index{adjoint equivalence}\index{equivalence!adjoint} if $u$ and $v$ are natural isomorphisms.
\end{definition}

\subsection*{Colimits and Coends}

\begin{definition}\label{def:colimits}
Suppose $F \cn \C \to \D$ is a functor with $\C$ a small category.  A \emph{colimit}\index{colimit}, also called a \emph{small colimit}, of $F$ is a pair $(\colmt,p)$ consisting of an object $\colmt \in \D$ and a morphism $p_a \cn Fa \to \colmt$ for each object $a \in \C$ such that the following two conditions hold.
\begin{enumerate}
\item\label{colimits_i}
For each morphism $f \cn a \to b$ in $\C$, the following diagram commutes.
\begin{equation}\label{colim_i}
\begin{tikzpicture}[vcenter]
\def\h{2}
\draw[0cell]
(0,0) node (a1) {Fa}
(a1)++(\h,0) node (a2) {Fb}
(a1)++(.5*\h,-1) node (a3) {\colmt}
;
\draw[1cell=.85]
(a1) edge node {Ff} (a2)
(a2) edge node {p_b} (a3)
(a1) edge node[swap] {p_a} (a3)
;
\end{tikzpicture}
\end{equation}
\item\label{colimits_ii} 
If $(\colmt',p')$ is another such pair satisfying \cref{colimits_i}, then there exists a unique morphism $h \cn \colmt \to \colmt'$ such that the following diagram commutes for each $a \in \C$.
\begin{equation}\label{colim_ii}
\begin{tikzpicture}[vcenter]
\def\h{2}
\draw[0cell]
(0,0) node (a1) {\colmt}
(a1)++(\h,0) node (a2) {\phantom{\colmt'}}
(a2)++(0,.05) node (a2') {\colmt'}
(a1)++(.5*\h,1) node (a3) {Fa}
;
\draw[1cell=.85]
(a3) edge node[swap] {p_a} (a1)
(a1) edge node {h} (a2)
(a3) edge[shorten >=-.6ex] node[pos=.6] {p'_a} (a2)
;
\end{tikzpicture}
\end{equation}
\end{enumerate}
A \emph{limit}\index{limit}, also called a \emph{small limit}, of $F$ is defined by turning the morphisms $p_a$ and $p'_a$ for $a \in \C$ and $h$ backward in the definition of a colimit of $F$.  A \emph{(co)product}\index{product}\index{coproduct} in $\D$ is a (co)limit for a functor $F \cn \C \to \D$ with $\C$ a small discrete category.  An \emph{initial object}\index{initial object} in $\D$ is a colimit for the functor $F \cn \emptyset \to \D$ with $\emptyset$ the empty category.  A \emph{terminal object}\index{terminal object} in $\D$ is a limit for the functor $F \cn \emptyset \to \D$.  A category $\D$ is \emph{(co)complete}\index{complete}\index{cocomplete} if all (co)limits exist.
\end{definition}

\begin{definition}\label{def:coends}
Suppose $F \cn \Cop \times \C \to \D$ is a functor with $\C$ a small category.  A \emph{coend}\index{coend} of $F$ is a pair $(\colmt,p)$ consisting of an object $\colmt \in \D$ and a morphism $p_a \cn F(a,a) \to \colmt$ for each object $a \in \C$ such that the following two conditions hold.
\begin{enumerate}
\item\label{coend_i}
For each morphism $f \cn a \to b$ in $\C$, the following diagram commutes.
\begin{equation}\label{coe_i}
\begin{tikzpicture}[vcenter]
\def\h{2.8} \def\v{-1.4}
\draw[0cell]
(0,0) node (a11) {F(b,a)}
(a11)++(\h,0) node (a12) {F(b,b)}
(a11)++(0,\v) node (a21) {F(a,a)}
(a12)++(0,\v) node (a22) {\colmt}
;
\draw[1cell=.85]
(a11) edge node {F(1_b,f)} (a12)
(a12) edge node {p_b} (a22)
(a11) edge node[swap] {F(f,1_a)} (a21)
(a21) edge node {p_a} (a22)
;
\end{tikzpicture}
\end{equation}
\item\label{coend_ii} 
If $(\colmt',p')$ is another such pair satisfying \cref{coend_i}, then there exists a unique morphism $h \cn \colmt \to \colmt'$ such that the following diagram commutes for each $a \in \C$.
\begin{equation}\label{coe_ii}
\begin{tikzpicture}[vcenter]
\def\h{2}
\draw[0cell]
(0,0) node (a1) {\colmt}
(a1)++(\h,0) node (a2) {\phantom{\colmt'}}
(a2)++(0,.05) node (a2') {\colmt'}
(a1)++(.5*\h,1) node (a3) {F(a,a)}
;
\draw[1cell=.85]
(a3) edge node[swap] {p_a} (a1)
(a1) edge node {h} (a2)
(a3) edge[shorten >=-.6ex] node[pos=.6] {p'_a} (a2)
;
\end{tikzpicture}
\end{equation}
\end{enumerate}
Condition \cref{coend_ii} implies that, if a coend of $F$ exists, then it is unique up to a unique isomorphism.  In this case, a coend of $F$ is denoted by $\int^{a \in \C} F(a,a)$.\label{not:coend}
\end{definition}

\section{Monoidal Categories}
\label{sec:sym_mon_cat}

This section reviews monoidal categories, their braided and symmetric variants, symmetric monoidal closed categories, monoidal functors, monoidal natural transformations, monoids, and modules.  For detailed discussion of monoidal categories, see \cite{joyal-street,maclane,cerberusI,cerberusII}. 

\begin{definition}\label{def:monoidalcategory}\index{monoidal category}\index{category!monoidal}
A \emph{monoidal category} $(\C,\otimes,\tu,\alpha,\lambda,\rho)$ consists of a category $\C$, \label{notation:monoidal-product}a functor $\otimes \cn \C \times \C \to \C$ called the \index{monoidal category!monoidal product}\emph{monoidal product}, \label{not:monoidalunit}an object $\tu \in \C$ called the \index{monoidal category!monoidal unit}\emph{monoidal unit}, and \label{not:associativityiso}\label{not:unitisos}natural isomorphisms
\begin{equation}\label{moncat_as}
\begin{tikzcd}[column sep=huge,row sep=tiny]
(a \otimes b) \otimes c \ar{r}{\alpha_{a,b,c}}[swap]{\iso} & a \otimes (b \otimes c)
\end{tikzcd}
\end{equation}
and
\begin{equation}\label{moncat_units}
\begin{tikzcd}[column sep=large]
\tu \otimes a \ar{r}{\lambda_a}[swap]{\iso} & a & a \otimes \tu \ar{l}{\iso}[swap]{\rho_a}
\end{tikzcd}
\end{equation}
for objects $a,b,c \in \C$ called the \index{associativity isomorphism}\index{monoidal category!associativity isomorphism}\emph{associativity isomorphism}, the \index{left unit isomorphism}\index{monoidal category!left unit isomorphism}\emph{left unit isomorphism}, and the \index{right unit isomorphism}\index{monoidal category!right unit isomorphism}\emph{right unit isomorphism}.
These data are required to make the following middle unity\index{monoidal category!middle unity axiom} and pentagon\index{monoidal category!pentagon axiom}\index{pentagon axiom} diagrams commute for objects $a,b,c,d \in \C$, where $ab$ means $a \otimes b$.
\begin{equation}\label{monoidalcataxioms}
\begin{tikzpicture}[xscale=1,yscale=1,vcenter]
\tikzset{0cell/.append style={nodes={scale=.85}}}
\tikzset{1cell/.append style={nodes={scale=.85}}}
\def\h{1.6} \def\g{1.3} \def\v{.9} \def\u{2}
\draw[0cell]
(0,0) node (a) {(a \tu) b}
(a)++(0,-2) node (b) {a (\tu b)}
(a)++(1.5,-1) node (c) {ab}
;
\draw[1cell]  
(a) edge node[pos=.4] {\rho_a 1_b} (c)
(a) edge node[swap] {\alpha_{a,\tu,b}} (b)
(b) edge node[swap,pos=.4] {1_a \lambda_b} (c)
;
\begin{scope}[shift={(5.5,0)}]
\draw[0cell]
(0,0) node (x0) {(ab)(cd)}
(x0)++(-\h,-\v) node (x11) {((ab)c)d}
(x0)++(\h,-\v) node (x12) {a(b(cd))}
(x0)++(-\g,-\u) node (x21) {(a(bc))d}
(x0)++(\g,-\u) node (x22) {a((bc)d)}
;
\draw[1cell]
(x11) edge node[pos=.2] {\al_{ab,c,d}} (x0)
(x0) edge node[pos=.8] {\al_{a,b,cd}} (x12)
(x11) edge node[swap,pos=.2] {\al_{a,b,c} 1_d} (x21)
(x21) edge node {\al_{a,bc,d}} (x22)
(x22) edge node[swap,pos=.8] {1_a \al_{b,c,d}} (x12)
;
\end{scope}
\end{tikzpicture}
\end{equation}
A monoidal category is \emph{strictly unital}\index{monoidal category!strictly unital}\index{strictly unital monoidal category} if $\lambda$ and $\rho$ are identities.  It is \emph{strict}\index{monoidal category!strict}\index{strict monoidal category} if $\alpha$, $\lambda$, and $\rho$ are identities.
\end{definition}

\begin{explanation}[Unity]
The morphism equality\index{monoidal category!unity properties} 
$\lambda_{\tu} = \rho_{\tu}$ 
holds in each monoidal category.  Moreover, the unity diagrams
\begin{equation}\label{moncat-other-unit-axioms}
\begin{tikzcd}[column sep=normal]
(\tu \otimes a) \otimes b \dar[swap]{\lambda_a \otimes 1_b} \rar{\alpha_{\tu,a,b}}
& \tu \otimes (a \otimes b) \dar{\lambda_{a \otimes b}}\\ 
a \otimes b \rar[-,semithick, double distance=1.8pt,double=\pagecol] & a \otimes b
\end{tikzcd}\qquad
\begin{tikzcd}[column sep=normal]
(a \otimes b) \otimes \tu \dar[swap]{\rho_{a \otimes b}} \rar{\alpha_{a,b,\tu}}
& a \otimes (b \otimes \tu) \dar{1_a \otimes \rho_b}\\ 
a \otimes b \rar[-,semithick, double distance=1.8pt,double=\pagecol] & a \otimes b
\end{tikzcd}
\end{equation}
commute.  See \cite[Section 2.2]{johnson-yau} for the proofs.
\end{explanation}

\begin{convention}[Bracketing]\label{conv:left_normal}\index{bracketing convention}
Unless otherwise specified, iterated monoidal products are \emph{left normalized}, meaning that the left half of each pair of parentheses is at the left end of the expression.  For example, $a \otimes b \otimes c$ means $(a \otimes b) \otimes c$.  An empty monoidal product\index{monoidal product!empty} means the monoidal unit $\tu$.  Necessary coherence isomorphisms \cite[VII.2 and XI.1]{maclane} are tacitly inserted.  
\end{convention}

\begin{definition}\label{def:braidedmoncat}\index{braided monoidal category}\index{monoidal category!braided}\index{category!braided monoidal}
A \emph{braided monoidal category} $(\C,\xi)$ consists of a monoidal category $\C$ and a natural isomorphism\label{notation:symmetry-iso}
\[a \otimes b \fto[\iso]{\xi_{a,b}} b \otimes a \forspace a,b \in \C\]
called the \index{braiding}\emph{braiding} such that the hexagon diagrams\index{hexagon diagram} 
\begin{equation}\label{hexagon-braided}
\def\sb{\scalebox{.7}}
\begin{tikzpicture}[scale=.75,commutative diagrams/every diagram]
			\node (P0) at (0:2cm) {\sb{$b \otimes (c \otimes a)$}};
			\node (P1) at (60:2cm) {\makebox[3ex][l]{\sb{$b \otimes (a \otimes c)$}}};
			\node (P2) at (120:2cm) {\makebox[3ex][r]{\sb{$(b \otimes a) \otimes c$}}};
			\node (P3) at (180:2cm) {\sb{$(a \otimes b) \otimes c$}};
			\node (P4) at (240:2cm) {\makebox[3ex][r]{\sb{$a \otimes (b \otimes c)$}}};
			\node (P5) at (300:2cm) {\makebox[3ex][l]{\sb{$(b \otimes c) \otimes a$}}};
			\path[commutative diagrams/.cd, every arrow, every label]
			(P3) edge node[pos=.25] {\sb{$\xi_{a,b}\otimes 1_c$}} (P2)
			(P2) edge node {\sb{$\alpha$}} (P1)
			(P1) edge node[pos=.75] {\sb{$1_b \otimes \xi_{a,c}$}} (P0)
			(P3) edge node[swap,pos=.4] {\sb{$\alpha$}} (P4)
			(P4) edge node {\sb{$\xi_{a,b \otimes c}$}} (P5)
			(P5) edge node[swap,pos=.6] {\sb{$\alpha$}} (P0);
\end{tikzpicture}\qquad
\begin{tikzpicture}[scale=.75, commutative diagrams/every diagram]
			\node (P0) at (0:2cm) {\sb{$(c \otimes a) \otimes b$}};
			\node (P1) at (60:2cm) {\makebox[3ex][l]{\sb{$(a \otimes c) \otimes b$}}};
			\node (P2) at (120:2cm) {\makebox[3ex][r]{\sb{$a \otimes (c \otimes b)$}}};
			\node (P3) at (180:2cm) {\sb{$a \otimes (b \otimes c)$}};
			\node (P4) at (240:2cm) {\makebox[3ex][r]{\sb{$(a \otimes b) \otimes c$}}};
			\node (P5) at (300:2cm) {\makebox[3ex][l]{\sb{$c \otimes (a \otimes b)$}}};
			\path[commutative diagrams/.cd, every arrow, every label]
			(P3) edge node[pos=.25] {\sb{$1_a \otimes \xi_{b,c}$}} (P2)
			(P2) edge node {\sb{$\alpha^\inv$}} (P1)
			(P1) edge node[pos=.75] {\sb{$\xi_{a,c}\otimes 1_b$}} (P0)
			(P3) edge node[swap,pos=.4] {\sb{$\alpha^\inv$}} (P4)
			(P4) edge node {\sb{$\xi_{a \otimes b, c}$}} (P5)
			(P5) edge node[swap,pos=.6] {\sb{$\alpha^\inv$}} (P0);
\end{tikzpicture}
\end{equation}
commute for objects $a,b,c \in \C$.
\end{definition}

\begin{explanation}[Unity]
In each braided monoidal category, the unity diagram\index{braided monoidal category!unity properties} 
\begin{equation}\label{braidedunity}
\begin{tikzcd}[column sep=large]
a \otimes \tu \dar[swap]{\rho_a} \rar{\xi_{a,\tu}} & \tu \otimes a \dar{\lambda_a} \rar{\xi_{\tu,a}} & a \otimes \tu \dar{\rho_a}\\ 
a \rar[-,semithick, double distance=1.8pt,double=\pagecol] 
& a \rar[-,semithick, double distance=1.8pt,double=\pagecol] & a
\end{tikzcd}
\end{equation}
commutes for each object $a$.
See \cite[1.3.21]{cerberusII} for the proof. 
\end{explanation}

\begin{definition}\label{def:symmoncat}\index{symmetric monoidal category}\index{monoidal category!symmetric}\index{category!symmetric monoidal}
A \emph{symmetric monoidal category} $(\C,\xi)$ consists of a monoidal category $\C$ and a natural isomorphism
\[a \otimes b \fto[\iso]{\xi_{a,b}} b \otimes a \forspace a,b \in \C\]
called the \emph{braiding} such that the following symmetry\index{symmetry axiom} and \index{hexagon diagram}hexagon diagrams commute for objects $a,b,c \in \C$. 
\begin{equation}\label{symmoncatsymhexagon}
\begin{tikzpicture}[xscale=3,yscale=1.2,vcenter]
\def\h{.1}
\draw[0cell=.8] 
(0,0) node (a) {a \otimes b}
(a)++(.7,0) node (c) {a \otimes b}
(a)++(.35,-1) node (b) {b \otimes a}
;
\draw[1cell=.8] 
(a) edge node {1_{a \otimes b}} (c)
(a) edge node [swap,pos=.3] {\xi_{a,b}} (b)
(b) edge node [swap,pos=.6] {\xi_{b,a}} (c)
;
\begin{scope}[shift={(1.5,.5)}]
\draw[0cell=.8] 
(0,0) node (x11) {(b \otimes a) \otimes c}
(x11)++(.9,0) node (x12) {b \otimes (a \otimes c)}
(x11)++(-\h,-1) node (x21) {(a \otimes b) \otimes c}
(x12)++(\h,-1) node (x22) {b \otimes (c \otimes a)}
(x11)++(0,-2) node (x31) {a \otimes (b \otimes c)}
(x12)++(0,-2) node (x32) {(b \otimes c) \otimes a}
;
\draw[1cell=.8]
(x21) edge node[pos=.25] {\xi_{a,b} \otimes 1_c} (x11)
(x11) edge node {\alpha} (x12)
(x12) edge node[pos=.75] {1_b \otimes \xi_{a,c}} (x22)
(x21) edge node[swap,pos=.25] {\alpha} (x31)
(x31) edge node {\xi_{a,b \otimes c}} (x32)
(x32) edge node[swap,pos=.75] {\alpha} (x22)
;
\end{scope}
\end{tikzpicture}
\end{equation}
A \emph{permutative category}\index{permutative category}\index{category!permutative} is a strict symmetric monoidal category, meaning a symmetric monoidal category where $\alpha$, $\lambda$, and $\rho$ are identities.
\end{definition}

\begin{explanation}\label{rk:smcat}
Assuming the symmetry axiom in \cref{symmoncatsymhexagon}, the two hexagon diagrams in \cref{hexagon-braided} are equivalent.  Thus, a symmetric monoidal category is precisely a braided monoidal category that satisfies the symmetry axiom.  The unity diagrams \cref{braidedunity} commute in each symmetric monoidal category.
\end{explanation}

\begin{definition}\label{def:closedcat}\index{closed category}\index{category!closed}\index{internal hom}
A symmetric monoidal category $(\C,\otimes,\tu)$ is \emph{closed} if, for each object $a \in \C$, the functor $- \otimes a \cn \C \to \C$ admits a right adjoint \label{notation:internal-hom}$[a,-]$ called an \emph{internal hom}.  In this case, for each object $b \in \C$, the counit
\[[a,b] \otimes a \to b\]
is called the \emph{evaluation}.  A symmetric monoidal closed category where $\otimes$ is the Cartesian product is called a \emph{Cartesian closed category}.
\end{definition}

\subsection*{Monoidal Functors and Monoidal Natural Transformations}

\begin{definition}\label{def:monoidalfunctor}\index{monoidal functor}\index{functor!monoidal}
Given monoidal categories $\C$ and $\D$, a \emph{monoidal functor} $(F, F^2, F^0) \cn \C \to \D$
consists of a functor $F \cn \C \to \D$, a morphism 
\[\tu \fto{F^0} F\tu \inspace \D\]
called the \index{unit constraint}\index{constraint!unit}\emph{unit constraint}, and a natural transformation
\[\begin{tikzcd}[column sep=large]
Fa \otimes Fb \ar{r}{F^2_{a,b}} & F(a \otimes b)
\end{tikzcd}
\forspace a,b \in \C\]
called the \index{monoidal constraint}\index{constraint!monoidal}\emph{monoidal constraint} such that the diagrams 
\begin{equation}\label{monoidalfunctorunity}
\begin{tikzcd}
\tu \otimes Fa \dar[swap]{F^0 \otimes 1_{Fa}} \rar{\lambda_{Fa}} & Fa \\ 
F\tu \otimes Fa \rar{F^2_{\tu,a}} & F(\tu \otimes a) \uar[swap]{F\lambda_a}
\end{tikzcd}
\qquad
\begin{tikzcd}
Fa \otimes \tu \dar[swap]{1_{Fa} \otimes F^0} \rar{\rho_{Fa}} & Fa \\ 
Fa \otimes F\tu \rar{F^2_{a,\tu}} & F(a \otimes \tu) \uar[swap]{F\rho_a}
\end{tikzcd}
\end{equation}
and
\begin{equation}\label{monoidalfunctorassoc}
\begin{tikzcd}[column sep=large]
\bigl(Fa \otimes Fb\bigr) \otimes Fc \rar{\alpha} \dar[swap]{F^2_{a,b} \otimes 1_{Fc}} 
& Fa \otimes \bigl(Fb \otimes Fc\bigr) \dar{1_{Fa} \otimes F^2_{b,c}}\\
F(a \otimes b) \otimes Fc \dar[swap]{F^2_{a \otimes b,c}} & Fa \otimes F(b \otimes c) \dar[d]{F^2_{a,b \otimes c}}\\
F\bigl((a \otimes b) \otimes c\bigr) \rar{F\alpha} & F\bigl(a \otimes (b \otimes c)\bigr)
\end{tikzcd}
\end{equation}
commute for objects $a,b,c \in \C$.  A monoidal functor between braided monoidal categories is said to be \index{braided monoidal functor}\index{monoidal functor!braided}\emph{braided} if the diagram 
\begin{equation}\label{monoidalfunctorbraiding}
\begin{tikzcd}[column sep=large]
Fa \otimes Fb \dar[swap]{F^2_{a,b}} \rar{\xi_{Fa,Fb}} & Fb \otimes Fa \dar{F^2_{b,a}} \\ 
F(a \otimes b) \rar{F\xi_{a,b}} & F(b \otimes a)
\end{tikzcd}
\end{equation}
commutes for $a,b \in \C$.  A braided monoidal functor between symmetric monoidal categories is called a \index{symmetric monoidal functor}\index{monoidal functor!symmetric}\emph{symmetric monoidal functor}.  A monoidal functor is \index{monoidal functor!strong}\index{strong monoidal functor}\emph{strong} if $F^0$ and $F^2$ are isomorphisms; \index{monoidal functor!strictly unital}\index{strictly unital monoidal functor}\emph{strictly unital} if $F^0$ is the identity morphism; and \index{monoidal functor!strict}\index{strict monoidal functor}\emph{strict} if $F^0$ and $F^2$ are identities.
\end{definition}

\begin{definition}\label{def:monoidalnattr}
Given monoidal functors 
\[(F, F^2, F^0) \andspace (H,H^2,H^0) \cn \C \to \D\]
between monoidal categories $\C$ and $\D$, a \index{monoidal natural transformation}\index{natural transformation!monoidal}\emph{monoidal natural transformation} $\phi \cn F \to H$ is a natural transformation such that the diagrams 
\begin{equation}\label{monnattr}
\begin{tikzpicture}[xscale=2,yscale=1.3,vcenter]
\draw[0cell=.9]
(0,0) node (x11) {\tu}
(x11)++(1,0) node (x12) {F\tu}
(a11)++(0,-1) node (x21) {\tu}
(x12)++(0,-1) node (x22) {H\tu}
(x12)++(1.2,0) node (y11) {Fa \otimes Fb}
(y11)++(1.5,0) node (y12) {Ha \otimes Hb}
(y11)++(0,-1) node (y21) {F(a \otimes b)}
(y12)++(0,-1) node (y22) {H(a \otimes b)}
;
\draw[1cell=.9]  
(x11) edge node {F^0} (x12)
(x12) edge node {\phi_{\tu}} (x22)
(x11) edge[equal] (x21)
(x21) edge node {H^0} (x22)
(y11) edge node {\phi_a \otimes \phi_b} (y12)
(y12) edge node {H^2_{a,b}} (y22)
(y11) edge node[swap] {F^2_{a,b}} (y21)
(y21) edge node {\phi_{a \otimes b}} (y22)
;
\end{tikzpicture}
\end{equation}
commute for $a,b \in \C$.
\end{definition}

\subsection*{Monoids and Modules}

\begin{definition}[Monoids]\label{def:monoid}
Given a monoidal category $(\C, \otimes, \tu, \al, \la, \rho)$, a \emph{monoid}\index{monoid} $(a,\mu,\eta)$\label{notation:monoid} in $\C$ consists of an object $a \in \C$, a morphism $\mu \cn a \otimes a \to a$ called the \index{multiplication}\emph{multiplication}, and a morphism $\eta \cn \tu \to a$ called the \index{unit}\emph{unit} such that the diagrams
\begin{equation}\label{monoid_axioms}
\begin{tikzcd}[column sep=normal]
(a \otimes a) \otimes a \arrow{dd}[swap]{\mu \otimes 1_a} \rar{\alpha} 
& a \otimes (a \otimes a) \dar{1_a \otimes \mu}\\ 
& a \otimes a \dar{\mu}\\  
a  \otimes a \arrow{r}{\mu} & a
\end{tikzcd}\qquad
\begin{tikzcd}[column sep=normal]
\tu \otimes a \ar{d}[swap]{\eta \otimes 1_a} \ar{r}{\lambda_a} 
& a \ar[-,semithick, double distance=1.8pt,double=\pagecol]{d}\\ 
a \otimes a \ar{r}{\mu} 
& a \ar[-,semithick, double distance=1.8pt,double=\pagecol]{d}\\
a \otimes \tu \ar{u}{1_a \otimes \eta} \ar{r}{\rho_a} & a
\end{tikzcd}
\end{equation}
commute.  Given a symmetric monoidal category $(\C,\xi)$, a \emph{commutative monoid}\dindex{commutative}{monoid} in $\C$ is a monoid $(a,\mu,\eta)$ such that the diagram
\begin{equation}\label{com_monoid}
\begin{tikzpicture}[vcenter]
\def\h{2.2}
\draw[0cell]
(0,0) node (a1) {a \otimes a}
(a1)++(\h,0) node (a2) {a \otimes a}
(a1)++(.5*\h,-1) node (a3) {a}
;
\draw[1cell=.9]
(a1) edge node {\xi} (a2)
(a2) edge node[pos=.4] {\mu} (a3)
(a1) edge node[swap,pos=.4] {\mu} (a3)
;
\end{tikzpicture}
\end{equation}
commutes.  A \emph{morphism} of (commutative) monoids $f \cn (a,\mu^a,\eta^a) \to (b,\mu^b,\eta^b)$ is a morphism $f \cn a \to b$ in $\C$ such that the diagrams
\begin{equation}\label{mor_monoids}
\begin{tikzcd}[column sep=large]
a \otimes a \dar[swap]{\mu^a} \rar{f\otimes f} & b \otimes b \dar{\mu^b}\\ 
a \rar{f} & b
\end{tikzcd} \qquad
\begin{tikzcd}[column sep=large]
\tu \dar[-,semithick, double distance=1.8pt,double=\pagecol] \rar{\eta^a} 
& a \dar{f}\\ \tu \rar{\eta^b} & b
\end{tikzcd}
\end{equation}
commute.
\end{definition}

\begin{definition}[Modules]\label{def:modules}
Suppose $(a,\mu,\eta)$ is a monoid in a monoidal category $\C$.
\begin{enumerate}
\item\label{def:modules_i} A \emph{right $a$-module}\index{module} $(x,\umu)$ consists of an object $x \in \C$ and a morphism $\umu \cn x \otimes a \to x$ called the \emph{right $a$-action} such that the associativity and unity diagrams
\begin{equation}\label{module_axioms}
\begin{tikzcd}[column sep=normal]
(x \otimes a) \otimes a \arrow{dd}[swap]{\umu \otimes 1_a} \rar{\alpha} 
& x \otimes (a \otimes a) \dar{1_x \otimes \mu}\\ 
& x \otimes a \dar{\umu}\\  
x  \otimes a \arrow{r}{\umu} & x
\end{tikzcd}\qquad
\begin{tikzcd}[column sep=normal]
x \otimes a \ar{r}{\umu} & x \ar[-,semithick, double distance=1.8pt,double=\pagecol]{d}\\
x \otimes \tu \ar{u}{1_x \otimes \eta} \ar{r}{\rho_x} & x
\end{tikzcd}
\end{equation}
commute.
\item\label{def:modules_iv} A \emph{morphism} of right $a$-modules $f \cn (x,\umu^x) \to (y,\umu^y)$ is a morphism $f \cn x \to y$ in $\C$ such that the diagram
\begin{equation}\label{module_morphism}
\begin{tikzpicture}[vcenter]
\def\v{-1.4}
\draw[0cell]
(0,0) node (a11) {x \otimes a}
(a11)++(2,0) node (a12) {x}
(a11)++(0,\v) node (a21) {y \otimes a}
(a12)++(0,\v) node (a22) {y}
;
\draw[1cell=.9]
(a11) edge node {\umu^x} (a12)
(a12) edge node {f} (a22)
(a11) edge node[swap] {f \otimes 1_a} (a21)
(a21) edge node {\umu^y} (a22)
;
\end{tikzpicture}
\end{equation}
commutes.
\item\label{def:modules_ii} A \emph{left $a$-module} is defined analogously to a right $a$-module using a \emph{left $a$-action} $\umu \cn a \otimes x \to x$ that makes the corresponding associativity and unity diagrams commute.
\item\label{def:modules_iii} If $(a,\mu,\eta)$ is a commutative monoid in a symmetric monoidal category $(\C,\xi)$, then each right $a$-module $(x,\umu)$ becomes a left $a$-module with left $a$-action given by the composite
\[a \otimes x \fto[\iso]{\xi} x \otimes a \fto{\umu} x.\]
\end{enumerate}
This finishes the definition.
\end{definition}

\section{Enriched Categories}
\label{sec:enrichedcat}

This section reviews enriched categories, enriched functors, and enriched natural transformations.  See \cite[Section 1.3]{johnson-yau}, \cite[Chapter 1]{cerberusIII}, and \cite{kelly-enriched} for further discussion. Suppose $(\V,\otimes,\tu,\alpha,\lambda,\rho)$ is a monoidal category \pcref{def:monoidalcategory}.

\begin{definition}\label{def:enriched-category}
A \emph{$\V$-category} $(\C,\mcomp,\cone)$\index{category!enriched}\index{enriched!category} consists of a class \index{object!enriched category}$\Ob\C$ of \emph{objects}, a \index{hom object}\index{object!hom}\emph{hom $\V$-object} $\C(a,b) \in \V$ for each pair of objects $a,b \in \C$, a \index{composition!enriched category}\emph{composition} morphism\label{not:enrcomposition}
\[\begin{tikzcd}[column sep=large] 
\C(b,c) \otimes \C(a,b) \rar{\mcomp_{a,b,c}} & \C(a,c)
\end{tikzcd} \inspace \V\]
for objects $a,b,c \in \C$, and an \emph{identity} morphism\label{not:enridentity}
\[\begin{tikzcd}[column sep=normal] 
\tu \rar{\cone_a} & \C(a,a)
\end{tikzcd} \inspace \V\]
for each object $a \in \C$ such that the diagrams
\begin{equation}\label{enriched-cat-associativity}
\begin{tikzcd}[cells={nodes={scale=.9}}]
  \big(\C(c,d) \otimes \C(b,c)\big) \otimes \C(a,b) \arrow{dd}[swap]{\mcomp_{b,c,d} \otimes 1}
  \rar{\alpha}
  & \C(c,d) \otimes \big(\C(b,c) \otimes \C(a,b)\big)  \dar{1 \otimes \mcomp_{a,b,c}} \\
& \C(c,d) \otimes \C(a,c) \dar{\mcomp_{a,c,d}}\\
\C(b,d) \otimes \C(a,b) \rar{\mcomp_{a,b,d}} & \C(a,d)  
\end{tikzcd}
\end{equation}
and
\begin{equation}\label{enriched-cat-unity}
\begin{tikzcd}[cells={nodes={scale=.9}}]
\tu \otimes \C(a,b) \dar[swap]{\cone_b \otimes 1} \rar{\lambda} 
& \C(a,b) \dar[-,semithick, double distance=1.8pt,double=\pagecol] 
& \C(a,b) \otimes \tu \lar[swap]{\rho} \dar{1 \otimes \cone_a}\\
\C(b,b) \otimes \C(a,b) \rar{\mcomp_{a,b,b}} & \C(a,b) & \C(a,b) \otimes \C(a,a) \lar[swap]{\mcomp_{a,a,b}}
\end{tikzcd}
\end{equation}
commute for objects $a,b,c,d \in \C$
\end{definition}

\begin{definition}\label{def:enriched-functor}
Given $\V$-categories $\C$ and $\D$, a \index{functor!enriched}\index{enriched!functor}\emph{$\V$-functor} $F \cn \C \to \D$ consists of
\begin{itemize}
\item an object assignment $F \cn \Ob\C \to \Ob\D$ and
\item a component morphism $F_{a,b} \cn \C(a,b) \to \D(Fa,Fb)$ for each pair of objects $a,b \in \C$
\end{itemize}
such that the diagrams\index{composition!enriched functor}\index{identities!enriched functor} 
\begin{equation}\label{eq:enriched-composition}
\begin{tikzcd}[column sep=large,cells={nodes={scale=.9}}]
\C(b,c) \otimes \C(a,b) \rar{\mcomp_{a,b,c}} \dar[swap]{F_{b,c} \otimes F_{a,b}} & \C(a,c) \dar{F_{a,c}}\\
\D(Fb,Fc) \otimes \D(Fa,Fb) \rar{\mcomp_{Fa,Fb,Fc}} & \D(Fa,Fc)\end{tikzcd}
\qquad
\begin{tikzcd}[cells={nodes={scale=.9}}]
\tu \dar[-,semithick, double distance=1.8pt,double=\pagecol] \rar{\cone_a} & \C(a,a) \dar{F_{a,a}}\\
\tu \rar{\cone_{Fa}} & \D(Fa,Fa)
\end{tikzcd}
\end{equation}
commute for objects $a,b,c \in \C$.  Identity and composition of $\V$-functors are defined separately on object assignments and component morphisms.
\end{definition}

\begin{definition}\label{def:enriched-natural-transformation}
Given $\V$-functors $F,H \cn \C\to\D$ between $\V$-categories, a \index{enriched!natural transformation}\index{natural transformation!enriched}\emph{$\V$-natural transformation} $\phi \cn F\to H$ consists of an $a$-component\index{component} $\V$-morphism $\phi_a \cn \tu \to \D(Fa,Ha)$ for each object $a \in \C$ such that the naturality diagram
\begin{equation}\label{enr_naturality}
\begin{tikzpicture}[x=40mm,y=12mm,vcenter]
\def\s{.85}
\draw[0cell=\s]
(0,0) node (a) {\C(a,b)}
(.25,1) node (b1) {\tu \otimes \C(a,b)}
(b1)++(1,0) node (c1) {\D(Fb,Hb) \otimes \D(Fa,Fb)}
(.25,-1) node (b2) {\C(a,b) \otimes \tu}
(b2)++(1,0) node (c2) {\D(Ha,Hb) \otimes \D(Fa,Ha)}
(c1)++(.25,-1) node (d) {\D(Fa,Hb)}
;
\draw[1cell=.8]
(a) edge node[pos=.3] {\la^\inv} node['] {\iso}(b1)
(a) edge node[swap,pos=.2] {\rho^\inv} node {\iso} (b2)
(b1) edge node {\phi_b \otimes F_{a,b}} (c1)
(b2) edge node {H_{a,b} \otimes \phi_a} (c2)
(c1) edge node[pos=.6] {\mcomp_{Fa,Fb,Hb}} (d)
(c2) edge node[swap,pos=.7] {\mcomp_{Fa,Ha,Hb}} (d)
;
\end{tikzpicture}
\end{equation}
commutes for $a,b \in \C$.  Identity, vertical composition, and horizontal composition of $\V$-natural transformations are defined componentwise.
\end{definition}

\section{2-Categories}
\label{sec:twocategories}

This section reviews 2-categories, 2-functors, 2-natural transformations, 2-adjunctions, (adjoint) 2-equivalences, adjunctions and equivalences in a 2-category, and modifications.  See \cite{johnson-yau} for further discussion.  

\begin{definition}\label{def:twocategory}\index{2-category}\index{category!2-}
A \emph{2-category} $\A$ consists of the following data and axioms.
\begin{description}
\item[Objects] $\A$ is equipped with a class of \emph{objects}\label{not:Azero} $\A_0$, which is also denoted by $\A$.
\item[1-cells] $\A$ is equipped with a class $\A_1(a,b)$\label{not:Aoneab} of \index{1-cell}\emph{1-cells} from $a$ to $b$ for each pair of objects $a,b \in\A_0$.  Such a 1-cell is denoted by $a \to b$.
\item[2-cells] For 1-cells $f,f' \in \A_1(a,b)$, $\A$ is equipped with a set $\A_2(f,f')$\label{not:Atwo} of \index{2-cell}\emph{2-cells} from $f$ to $f'$, each denoted by
\begin{equation}\label{iicell_notation}
\begin{tikzpicture}[baseline={(x1.base)}]
\def\t{30}
\draw[0cell=1]
(0,0) node (a) {a}
(a)++(1.6,0) node (b) {b}
;
\draw[1cell=.8]  
(a) edge[bend left=\t] node {f} (b)
(a) edge[bend right=\t] node[swap] {f'} (b)
;
\draw[2cell]
node[between=a and b at .48, rotate=-90, 2label={above,}] {\Rightarrow}
;
\end{tikzpicture}
\end{equation}
or $f \to f'$.
\item[Identities]  $\A$ is equipped with
\begin{itemize}
\item an \emph{identity 1-cell} $1_a \in \A_1(a,a)$\label{not:idonecell} for each object $a$ and
\item an \emph{identity 2-cell} $1_f \in \A_2(f,f)$\label{not:idtwocell} for each 1-cell $f \in \A_1(a,b)$.
\end{itemize}
\item[Compositions]
$\A$ is equipped with the following compositions.
\begin{itemize}
\item For objects $a,b \in \A_0$ and 1-cells $f,f',f'' \in \A_1(a,b)$, $\A$ is equipped with the \index{vertical composition!2-category}\emph{vertical composition} of 2-cells\label{not:vcompiicell}
\[\begin{tikzcd}
\A_2(f',f'') \times \A_2(f,f') \rar{v} & \A_2(f,f''),
\end{tikzcd}\] 
which is denoted by $v(\alpha',\alpha) = \alpha'\alpha$.
\item For objects $a,b,c \in \A_0$, $\A$ is equipped with the \index{horizontal composition!2-category}\emph{horizontal composition} of 1-cells\label{not:hcompicell}
\[\begin{tikzcd}
\A_1(b,c) \times \A_1(a,b) \rar{h_1} & \A_1(a,c),
\end{tikzcd}\]
which is denoted by $h_1(g,f) = gf$.
\item For objects $a,b,c \in \A_0$ and 1-cells $f,f' \in \A_1(a,b)$ and $g,g' \in \A_1(b,c)$, $\A$ is equipped with the \emph{horizontal composition} of 2-cells\label{not:hcompiicell}
\[\begin{tikzcd}
\A_2(g,g') \times \A_2(f,f') \rar{h_2} & \A_2(gf,g'f'),
\end{tikzcd}\]
which is denoted by $h_2(\beta,\alpha) = \beta * \alpha$.
\end{itemize}
\item[Axioms] These data satisfy axioms \crefrange{twocat-i}{twocat-iv}.
\begin{enumerate}
\item\label{twocat-i} Vertical composition is associative and unital for identity 2-cells.
\item\label{twocat-ii} Horizontal composition preserves identity 2-cells and vertical composition.
\item\label{twocat-iii} Horizontal composition of 1-cells is associative and unital for identity 1-cells.
\item\label{twocat-iv} Horizontal composition of 2-cells is associative and unital for identity 2-cells of identity 1-cells.
\end{enumerate}
\end{description}
Moreover, we define the following.
\begin{itemize}
\item The \index{underlying 1-category}\index{1-category!underlying}\emph{underlying 1-category} of a 2-category $\A$ is the category with the same class of objects, morphisms given by 1-cells, composition given by horizontal composition of 1-cells, and identity morphisms given by identity 1-cells.
\item For objects $a$ and $b$ in a 2-category $\A$, the \index{hom category!2-category}\emph{hom category} \label{not:homcat}$\A(a,b)$ is the category with
\begin{itemize}
\item objects given by 1-cells $a \to b$,
\item morphisms $f \to f'$ given by 2-cells,
\item composition given by vertical composition of 2-cells, and
\item identity morphisms given by identity 2-cells.
\end{itemize}
 \item A 2-category is \index{locally small}\emph{locally small} if each hom category is a small category.\defmark
\end{itemize}
\end{definition}

\begin{example}[Small Categories]\label{ex:catastwocategory}\index{2-category!of small categories}
There is a locally small 2-category $\Cat$ with small categories as objects, functors as 1-cells, and natural transformations as 2-cells.  Horizontal composition of 1-cells is given by composition of functors.  Horizontal and vertical compositions of 2-cells are the corresponding compositions of natural transformations.
\end{example}

\begin{example}\label{ex:lsiicat}
The notation $\Cat$ also denotes the Cartesian closed category $(\Cat,\times,[,])$ of small categories and functors.  A $\Cat$-category in the sense of \cref{def:enriched-category} is a locally small 2-category and vice versa.
\end{example}

\subsection*{2-Functors and 2-Natural Transformations}

\begin{definition}\label{def:twofunctor}
A \index{2-functor}\index{functor!2-}\emph{2-functor} $F \cn \A \to \B$ between 2-categories consists of
\begin{itemize}
\item an \emph{object assignment} $F_0 \cn \A_0 \to \B_0$,
\item a \emph{1-cell assignment} $F_1 \cn \A_1(a,a') \to \B_1(F_0 a,F_0 a')$ for each pair of objects $a,a' \in \A$, and
\item a \emph{2-cell assignment} $F_2 \cn \A_2(f,f') \to \B_2(F_1 f,F_1 f')$ for each pair of objects $a,a'$ and 1-cells $f,f' \in \A_1(a,a')$
\end{itemize}
such that axioms \cref{twofunctor-i,twofunctor-ii,twofunctor-iii} are satisfied.
\begin{enumerate}
\item\label{twofunctor-i} $F_0$ and $F_1$ define a functor between the underlying 1-categories of $\A$ and $\B$.
\item\label{twofunctor-ii} $F$ defines a functor $F \cn \A(a,a') \to \B(F_0 a,F_0 a')$ between hom categories for each pair of objects $a,a' \in \A$.
\item\label{twofunctor-iii} $F_2$ preserves horizontal composition of 2-cells.
\end{enumerate}
We usually denote each of $F_0$, $F_1$, and $F_2$ by $F$.  Given a 2-functor $H \cn \B \to \C$ between 2-categories, the \emph{composite 2-functor} $HF \cn \A \to \C$ is defined by separately composing the object, 1-cell, and 2-cell assignments.
\end{definition}

\begin{definition}\label{def:twonaturaltr}\index{2-natural!transformation}\index{natural transformation!2-}
Given 2-functors $F,H \cn \A\to\B$ between 2-categories, a \emph{2-natural transformation} $\phi \cn F \to H$ consists of a \emph{component 1-cell} $\phi_a \cn Fa \to Ha$ in $\B$ for each object $a \in \A$ such that the following two naturality axioms are satisfied.
\begin{description}
\item[1-cell naturality] For each 1-cell $f \cn a \to b$ in $\A$, the two composite 1-cells 
\begin{equation}\label{onecellnaturality}
\begin{tikzcd}[column sep=large]
Fa \ar{d}[swap]{Ff} \ar{r}{\phi_a} & Ha \ar{d}{Hf}\\
Fb \ar{r}{\phi_b} & Hb
\end{tikzcd}
\end{equation}
in $\B(Fa,Hb)$ are equal.
\item[2-cell naturality] 
For each 2-cell $\psi \cn f \to f'$ in $\A(a,b)$, the two whiskered 2-cells
\begin{equation}\label{twocellnaturality}
\begin{tikzpicture}[xscale=3,yscale=1.4,vcenter]
\def\a{30} \def\s{.85}
\draw[0cell=\s]
(0,0) node (x11) {Fa}
(x11)++(1,0) node (x12) {Ha}
(x11)++(0,-1) node (x21) {Fb}
(x12)++(0,-1) node (x22) {Hb}
;
\draw[1cell=\s]  
(x11) edge node {\phi_a} (x12)
(x21) edge node {\phi_b} (x22)
(x11) edge[bend right] node[swap] {Ff} (x21)
(x11) edge[bend left] node {Ff'} (x21)
(x12) edge[bend right] node[swap] {Hf} (x22)
(x12) edge[bend left] node {Hf'} (x22)
;
\draw[2cell=.9]
node[between=x11 and x21 at .65, 2label={above,F\psi}] {\Rightarrow}
node[between=x12 and x22 at .65, 2label={above,H\psi}] {\Rightarrow}
;
\end{tikzpicture}
\end{equation}
are equal, meaning the 2-cell equality
\[H\psi * 1_{\phi_a} = 1_{\phi_b} * F\psi \inspace \B(Fa,Hb).\]
\end{description}
We use the 2-cell notation \cref{twocellnotation} for 2-natural transformations.  A \index{2-natural!isomorphism}\index{natural isomorphism!2-}\emph{2-natural isomorphism} is a 2-natural transformation for which each component 1-cell is an isomorphism in the underlying 1-category.
\end{definition}

\begin{definition}\label{def:twonatcomposition}
Given 2-natural transformations $\phi \cn F \to H$ and $\varphi \cn H \to K$ for 2-functors $F,H,K \cn \A \to \B$ between 2-categories $\A$ and $\B$, the \emph{horizontal composite}\index{horizontal composition!2-natural transformation} $\varphi \phi \cn F \to K$ is the 2-natural transformation defined by the horizontal composite component 1-cell 
\[\begin{tikzcd}[column sep=normal]
Fa \ar{r}{\phi_a} & Ha \ar{r}{\varphi_a} & Ka
\end{tikzcd}\]
in $\B(Fa,Ka)$ for each object $a \in \A$.
\end{definition}

\subsection*{2-Equivalences and 2-Adjunctions}

\begin{definition}\label{def:twoequivalence}
A 2-functor $F \cn \A \to \B$ is a \emph{2-equivalence}\index{2-equivalence}\index{equivalence!2-} if there exist a 2-functor $H \cn \B \to \A$ and 2-natural isomorphisms
\[1_{\A} \fiso HF \andspace FH \fiso 1_{\B}.\]
In this case, $H$ is called a \emph{2-inverse} of $F$.
\end{definition}

A proof of the following characterization of a 2-equivalence is given in \cite[7.5.8]{johnson-yau} and also in \cite[3.70]{yau-multgro}.

\begin{theorem}\label{thm:twoequivalences}
A 2-functor is a 2-equivalence if and only if it is essentially surjective on objects, fully faithful on 1-cells, and fully faithful on 2-cells.
\end{theorem}

\begin{definition}\label{def:twoadjunction}
Given 2-categories $\A$ and $\B$, a \emph{2-adjunction}\index{2-adjunction}\index{adjunction!2-} $(L,R,u,v)$ from $\A$ to $\B$ consists of
\begin{itemize}
\item a 2-functor $L \cn \A \to \B$ called the \index{left 2-adjoint}\index{left adjoint!2-}\emph{left 2-adjoint},
\item a 2-functor $R \cn \B \to \A$ called the \index{right 2-adjoint}\index{right adjoint!2-}\emph{right 2-adjoint},
\item a 2-natural transformation $u \cn 1_{\A} \to RL$ called the \index{unit}\emph{unit}, and
\item a 2-natural transformation $v \cn LR \to 1_{\B}$ called the \index{counit}\emph{counit}
\end{itemize}
such that the equalities of 2-natural transformations\index{triangle identities}
\[\begin{split}
(vL)(Lu) &= 1_L \andspace\\ 
(Rv)(uR) &= 1_R
\end{split}\]
hold.  These equalities are called the \emph{left triangle identity}\index{left triangle identity} and the \index{right triangle identity}\emph{right triangle identity}.  A 2-adjunction $(L,R,u,v)$ is called an \emph{adjoint 2-equivalence}\index{adjoint 2-equivalence}\index{2-equivalence!adjoint} if $u$ and $v$ are 2-natural isomorphisms.  In the triangle identities, the 2-natural transformations 
\[\begin{split}
Lu \cn & L \to LRL,\\
vL \cn & LRL \to L,\\
uR \cn & R \to RLR, \andspace\\
Rv \cn & RLR \to R
\end{split}\] 
are defined by the component 1-cells 
\[\begin{split}
Lu_a \cn & La \to LRLa,\\
v_{La} \cn & LRLa \to La,\\
u_{Rb} \cn & Rb \to RLRb, \andspace\\
Rv_b \cn & RLRb \to Rb
\end{split}\] 
for objects $a \in \A$ and $b \in \B$. 
\end{definition}

\begin{definition}\label{def:internal_adj}
Suppose $a$ and $b$ are objects in a 2-category $\A$.  An \emph{adjunction} $(\ell,r,u,v)$ from $a$ to $b$ consists of
\begin{itemize}
\item a 1-cell $\ell \cn a \to b$ called the \emph{left adjoint},
\item a 1-cell $r \cn b \to a$ called the \emph{right adjoint},
\item a 2-cell $u \cn 1_a \to r\ell$ called the \emph{unit}, and
\item a 2-cell $v \cn \ell r \to 1_b$ called the \emph{counit}
\end{itemize}
such that the 2-cell equalities\index{triangle identities}
\[\begin{split}
(v * 1_\ell)(1_\ell \ast u) &= 1_\ell \andspace\\ 
(1_r * v)(u * 1_r) &= 1_r
\end{split}\]
hold.  These equalities are called the \emph{left triangle identity}\index{left triangle identity} and the \index{right triangle identity}\emph{right triangle identity}.  An adjunction $(\ell, r,u,v)$ is called an \emph{adjoint equivalence}\index{adjoint equivalence}\index{equivalence!adjoint} if the 2-cells $u$ and $v$ are isomorphisms in the hom categories $\A(a,a)$ and $\A(b,b)$.
\end{definition}

\subsection*{Modifications}

\begin{definition}\label{def:modification}
Given 2-natural transformations $\phi, \varphi \cn F \to H$ for 2-functors $F,H \cn \A \to \B$ between 2-categories $\A$ and $\B$, a \index{modification}\emph{modification} $\Phi \cn \phi \to \varphi$ consists of a \emph{component 2-cell}\label{not:componentiicell}
\[\begin{tikzcd}[column sep=normal]
\phi_a \ar{r}{\Phi_a} & \varphi_a \inspace \B(Fa,Ha)
\end{tikzcd}\]
for each object $a \in \A$ such that the two whiskered 2-cells 
\begin{equation}\label{modificationaxiom}
\begin{tikzpicture}[xscale=2.5,yscale=1.5,vcenter]
\def\a{35} \def\s{.85}
\draw[0cell=\s]
(0,0) node (x11) {Fa}
(x11)++(1,0) node (x12) {Ha}
(x11)++(0,-1) node (x21) {Fb}
(x12)++(0,-1) node (x22) {Hb}
;
\draw[1cell=\s]  
(x11) edge[bend left=\a] node[pos=.4] {\phi_a} (x12)
(x11) edge[bend right=\a] node[swap,pos=.6] {\varphi_a} (x12)
(x21) edge[bend left=\a] node[pos=.4] {\phi_b} (x22)
(x21) edge[bend right=\a] node[swap,pos=.6] {\varphi_b} (x22)
(x11) edge node[swap] {Ff} (x21)
(x12) edge node {Hf} (x22)
;
\draw[2cell=.9]
node[between=x11 and x12 at .44, rotate=-90, 2label={above,\Phi_a}] {\Rightarrow}
node[between=x21 and x22 at .44, rotate=-90, 2label={above,\Phi_b}] {\Rightarrow}
;
\end{tikzpicture}
\end{equation}
in $\B(Fa,Hb)$ are equal for each 1-cell $f \cn a \to b$ in $\A$, meaning the 2-cell equality
\[1_{Hf} * \Phi_a = \Phi_b * 1_{Ff}.\]
Such a modification is also denoted as follows.
\begin{equation}\label{iiicell_notation}
\begin{tikzpicture}[vcenter]
\def\t{27}
\draw[0cell]
(0,0) node (a) {\A}
(a)++(3,0) node (b) {\B}
;
\draw[1cell=.8]
(a) edge[bend left=\t] node {F} (b)
(a) edge[bend right=\t] node[swap] {H} (b)
;
\draw[2cell=.9]
node[between=a and b at .32, rotate=-90, 2labelmed={below,\phi}] {\Rightarrow}
node[between=a and b at .65, rotate=-90, 2label={above,\varphi}] {\Rightarrow}
;
\draw[2cell]
node[between=a and b at .5, rotate=0, shift={(0,-.17)}, 2labelmed={above,\!\Phi}] {\Rrightarrow}
;
\end{tikzpicture}
\end{equation}
The axiom \cref{modificationaxiom} is called the \index{modification!axiom}\emph{modification axiom}.
\end{definition}

\begin{definition}\label{def:modcomposition}
Suppose $F,H,K \cn \A \to \B$ are 2-functors between 2-categories $\A$ and $\B$.
\begin{description}
\item[Vertical composition] Given 2-natural transformations $\phi, \varphi, \psi \cn F \to H$ and modifications $\Phi \cn \phi \to \varphi$ and $\Psi \cn \varphi \to \psi$, the \emph{vertical composite modification}\index{vertical composition!modification} $\Psi \Phi \cn \phi \to \psi$\label{not:modvcomp} is defined by the vertical composite 2-cell
\[\begin{tikzcd}[column sep=normal]
\phi_a \ar{r}{\Phi_a} & \varphi_a \ar{r}{\Psi_a} & \psi_a
\end{tikzcd}\]
in $\B(Fa,Ha)$ for each object $a \in \A$.
\item[Horizontal composition] Given a modification $\Phi' \cn \phi' \to \varphi'$ for 2-natural transformations $\phi', \varphi' \cn H \to K$, the \emph{horizontal composite modification}\index{horizontal composition!modification} $\Phi' * \Phi \cn \phi' \phi \to \varphi' \varphi$\label{not:modhcomp} is defined by the horizontal composite 2-cell
\[\begin{tikzcd}[column sep=large]
\phi'_a \phi_a \ar{r}{\Phi'_a \,*\, \Phi_a} & \varphi'_a \varphi_a
\end{tikzcd}\]
in $\B(Fa,Ka)$ for each object $a \in \A$.\defmark
\end{description}
\end{definition}

\section{Enriched Operads}
\label{sec:operads}

This section reviews enriched operads and their algebras.  For more discussion of enriched operads, see \cite{may97a,may97b,yau-operad}.  Suppose $(\V,\otimes,\tu,\al,\la,\rho,\xi)$ is a symmetric monoidal category \pcref{def:symmoncat}.  The symmetric group on $n$ letters is denoted by $\Si_n$.

\begin{definition}\label{def:voperad}\index{enriched!operad}\index{operad!enriched}
A \emph{$\V$-operad} $(\Op,\ga,\opu)$ consists of the following data and axioms.
\begin{description}
\item[$\Si_n$-objects]
$\Op$ is equipped with a right $\Si_n$-object $\Op(n) \in \V$ for each $n \geq 0$.
\item[Unit]
$\Op$ is equipped with a $\V$-morphism 
\[\tu \fto{\opu} \Op(1)\] 
called the \index{operadic unit}\emph{operadic unit}.
\item[Composition]
For integers $n \geq 1$ and $k_1, \ldots, k_n \geq 0$ with $k = \sum_{i=1}^n k_i$, $\Op$ is equipped with a $\V$-morphism
\[\Op(n) \otimes \txotimes_{i=1}^n \Op(k_i) \fto{\ga} \Op(k)\]
called the \index{operadic composition}\emph{operadic composition}.
\end{description}
These data satisfy the following unity, associativity, and equivariance axioms.
\begin{description}
\item[Unity]
The following unity diagrams in $\V$ commute for $k \geq 0$ and $n \geq 1$.
\begin{equation}\label{operad_unity}

\end{equation}
\item[Associativity]
The following associativity diagram in $\V$ commutes for integers $n \geq 1$; $k_i \geq 0$ for $1 \leq i \leq n$ with $k = \sum_{i=1}^n k_i$ and at least one $k_i \geq 1$; and $\ell_{i,r} \geq 0$ for $1 \leq i \leq n$ and $1 \leq r \leq k_i$ with $\ell_i = \sum_{r=1}^{k_i} \ell_{i,r}$ and $\ell = \sum_{i=1}^n \ell_i$.
\begin{equation}\label{operad_associativity}
%
\end{equation}
\item[Equivariance]
Suppose $n \geq 1$ and $k_1,\ldots,k_n \geq 0$ with $k = \sum_{i=1}^n k_i$.
\begin{description}
\item[Top equivariance]
The following diagram in $\V$ commutes for each permutation $\si \in \Si_n$, where $\sigmabar \in \Si_k$ permutes $n$ consecutive blocks of lengths $k_{\si(1)}, \ldots, k_{\si(n)}$ as $\si$ permutes $\{1,\ldots,n\}$.
\begin{equation}\label{operad_topeq}
%
\end{equation}
\item[Bottom equivariance]
The following diagram in $\V$ commutes for permutations $\tau_i \in \Si_{k_i}$ for $1 \leq i \leq n$, where $\tau^\times = \bigoplus_{i=1}^n \tau_i \in \Si_k$ is the block sum.
\begin{equation}\label{operad_boteq}
%
\end{equation}
\end{description}
\end{description}
An \emph{operad}\index{operad} is a $\Set$-operad for the Cartesian closed category $\Set$ of sets and functions.
\end{definition}

\begin{definition}\label{def:operadalg}
Given a $\V$-operad $(\Op,\ga,\opu)$, an \emph{$\Op$-algebra}\index{algebra}\index{operad!algebra} $(\A,\gaA)$ consists of an object $\A \in \V$ and a $\V$-morphism
\[\Op(n) \otimes \A^{\otimes n} \fto{\gaA} \A\]
for each $n \geq 0$ such that the following unity, associativity, and equivariance axioms hold.
\begin{description}
\item[Unity]
The following diagram in $\V$ commutes.
\begin{equation}\label{operadalg_unity}

\end{equation}
\item[Associativity]
For $n \geq 1$ and $k_1,\ldots,k_n \geq 0$ with $k = \sum_{i=1}^n k_i$, the following diagram in $\V$ commutes.
\begin{equation}\label{operadalg_assoc}
%
\end{equation}
\item[Equivariance]
For each permutation $\si \in \Si_n$, the following diagram in $\V$ commutes.
\begin{equation}\label{operadalg_equiv}
%
\end{equation}
\end{description}
An \emph{$\Op$-algebra morphism}\index{algebra morphism}\index{morphism!algebra} $f \cn (\A,\gaA) \to (\B,\gaB)$ between $\Op$-algebras is a morphism $f \cn \A \to \B$ in $\V$ such that the diagram 
\begin{equation}\label{operadalgmor}
%
\end{equation}
commutes for each $n \geq 0$.  The category of $\Op$-algebras and $\Op$-algebra morphisms is denoted by $\AlgO$.\label{not:AlgO}
\end{definition}